\documentclass[a4paper,12pt]{book}
\usepackage{amsthm,amsfonts,amssymb,euscript}
\usepackage{latexsym, multicol, fancybox}
\usepackage{graphicx}
\usepackage{color}
\usepackage{amsmath, amsthm, amssymb, bm}
\usepackage{epstopdf}
\usepackage{caption}
\usepackage{psfrag}

\usepackage{mathrsfs}
 \usepackage{xcolor}
 \usepackage[citebordercolor={green}]{hyperref}
 \usepackage{tikz}
 \usepackage{bbm}
 \usepackage{dsfont}
 \usepackage{bbold}

\numberwithin{equation}{section}

\newtheorem{theorem}{Theorem}[section]
\newtheorem{lemma}[theorem]{Lemma}
\newtheorem{proposition}[theorem]{Proposition}
\newtheorem{corollary}[theorem]{Corollary}
\newtheorem{definition}[theorem]{Definition}
\newtheorem{remark}[theorem]{Remark}

\newcommand{\bea}{\begin{eqnarray}}
\newcommand{\eea}{\end{eqnarray}}
\def\beaa{\begin{eqnarray*}}
\def\eeaa{\end{eqnarray*}}
\def\ba{\begin{array}}
\def\ea{\end{array}}
\def\be#1{\begin{equation} \label{#1}}
\def \eeq{\end{equation}}

\newcommand{\nn}{\nonumber}

\def\a{{\alpha}}
\def\al{{\alpha}}
\def\b{{\beta}}
\def\be{{\beta}}
\def\ga{\gamma}
\def\Ga{\Gamma}
\def\de{\delta}
\def\De{\Delta}
\def\ep{\epsilon}
\def\ka{\kappa}

\def\la{\lambda}
\def\La{\Lambda}
\def\si{\sigma}
\def\Si{\Sigma}
\def\om{\omega}

\def\vphi{\varphi}
\def\varo{{\varrho}}
\def\th{\theta}

\def\ze{\zeta}
\def\ka{\kappa}
\def\nab{\nabla}

\def\Up{\Upsilon}

\def\pr{{\partial}}
\def\les{\lesssim}
\def\c{\cdot}

\def\AA{{math\cal A}}
\def\BB{{\mathcal B}}

\def\MM{{\mathcal M}}
\def\NN{{\mathcal N}}

\def\LL{{\mathcal L}}
\def\II{{\mathcal I}}

\def\FF{{\mathcal F}}
\def\EE{{\mathcal E}}
\def\HH{{\mathcal H}}

\def\TT{{\mathcal T}}

\def\VV{{\mathcal V}}

\def\OO{{\mathcal O}}
\def\SS{{\mathcal S}}
\def\UU{{\mathcal U}}
\def\JJ{{\mathcal J}}
\def\KK{{\mathcal K}}
\def\Lie{{\mathcal L}}

\def\DD{{\mathcal D}}
\def\PP{{\mathcal P}}
\def\RR{{\mathcal R}}
\def\QQ{{\mathcal Q}}
\def\AA{{\mathcal A}}

\def\HH{{\mathcal H}}

\def\Lie{{\mathcal L}}

\def\LLb{ \, \LL     \mkern-9mu /}
\def\lap{{\triangle}}

\def\B{{\bf B}}
\def\C{{\bf C}}
\def\D{{\bf D}}

\def\G{{\bf G}}

\def\M{{\bf M}}

\def\O{{\bf O}}

\def\R{{\bf R}}

\def\S{{\bf S}}
\def\K{{\bf K}}
\def\T{{\bf T}}
\def\Z{{\bf Z}}
\def\X{\bf X}
\def\g{{\bf g}}

\def\e{{\bf e}}
\def\k{{\bold k}}
\def\q{\bold q}
\def\w{{\bf w}}

\def\SSS{{\Bbb S}}
\def\RRR{{\Bbb R}}

\def\CCC{{\Bbb C}}
\def\f12{{\frac 1 2}}

\def\dual{{\,\,^*}}
\def\div{{\mbox div\,}}
\def\curl{{\mbox curl\,}}

\def\Lb{{\,\underline{l}}}
\def\Hb{\,\underline{H}}
\def\Lb{{\,\underline{L}}}

\def\Xh{\,^{(h)}X}
\def\Yh{\,^{(h)}Y}
\def\trch{{\mbox tr}\, \chi}
\def\chih{{\widehat \chi}}
\def\chib{{\underline \chi}}
\def\chibh{{\underline{\chih}}}

\def\etab{{\underline \eta}}
\def\omb{{\underline{\om}}}
\def\bb{{\underline{\b}}}
\def\aa{\protect\underline{\a}}
\def\xib{{\underline \xi}}

\def\Xib{\underline{\Xi}}

\def\Ab{\protect\underline{A}}
\def\Bb{\protect\underline{B}}
\def\Xb{\protect\underline{X}}

\def\Xh{\widehat{X}}
\def\Xbh{\widehat{\Xb}}

\def\ub{\underline{u}}

\def\varoc{\check{\varrho}}

\newcommand{\nabb}{{\bf \nab} \mkern-13mu /\,}

\def\tr{\mbox{tr}}
\def\atr{\,^{(a)}\mbox{tr}}

\def\trchb{{\tr \,\chib}}

\def\Div{\mbox{Div}}

\def\chia{{\,^{(a)}\chi}}

\def\chiba{{\,^{(a)}\chib}}

\def\Us{{\,^{(s)}U}}
\def\Ua{{\,^{(a)}U}}
\def\Uh{{\widehat U}}

\def\Vh{{\widehat V}}

\def\atrch{\atr\chi}
\def\atrchb{\atr\chib}

\def\hot{\widehat{\otimes}}
\def\rhod{\,\dual\hspace{-2pt}\rho}

\def\GaX{\, ^{(X)}\Ga}
\def\piX{\, ^{(X)}\pi}
\def\piZ{\, ^{(\Z)}\pi}

\def\fb{\protect\underline{f}}
\def\err{{\mbox{Err}}}
\def\ov{\overline}

\newcommand{\hch}{\widehat{\chi}}

\def\f12{\frac 1 2}
\def\lab{\label}
\def\nabc{\,^{(c)}\nab}

\def\bsplit{\begin{split}}

\newcommand{\Mext}{{\,{}^{(ext)}\mathcal{M}}}

\newcommand{\Mint}{{\protect \,{}^{(int)}\mathcal{M}}}

\def\Mint{{\protect \, ^{(int)}\MM}}

\def\Mtop{{ \, ^{(top)} \MM}}

\def\qf{\frak{q}}
\def\qfb{\protect\underline{\qf}}

\def\Rk{\mathfrak{R}}

\def\Ik{\mathfrak{I}}
\def\Jk{\mathfrak{J}}

\def\sk{\mathfrak{s}}
\def\dkb{ \, \mathfrak{d}     \mkern-9mu /}
\def\dk{\mathfrak{d}}

\DeclareFontFamily{U}{mathx}{\hyphenchar\font45}
\DeclareFontShape{U}{mathx}{m}{n}{
      <5> <6> <7> <8> <9> <10>
      <10.95> <12> <14.4> <17.28> <20.74> <24.88>
      mathx10
      }{}
\DeclareSymbolFont{mathx}{U}{mathx}{m}{n}
\DeclareFontSubstitution{U}{mathx}{m}{n}
\DeclareMathAccent{\widecheck}{0}{mathx}{"71}

\def\Zc{\widecheck{Z}}
\def\Hc{\widecheck{H}}
\def\Hbc{\widecheck{\Hb}}
\def\trXc{\widecheck{\tr X}}
\def\trXbc{\widecheck{\tr\Xb}}
\def\Pc{\widecheck{P}}
\def\ombc{\widecheck{\omb}}
\def\Gac{\widecheck{\Ga}}
\def\Rc{\widecheck R}

\def\rhoc{\widecheck{\rho}}

\def\omc{\widecheck \omega}
\def\ombc{\underline{\widecheck{\omega}}}

\def\trchc{\widecheck{\tr\chi}}
\def\trchbc{\widecheck{\tr\chib}}

\def\rhodc{\widecheck{\rhod}}

\def\DDc{\,^{(c)} \DD}

\def\qfc{\widecheck{\qf}}

\newcommand{\deh}{\delta_{\mathcal{H}}}
\newcommand{\dec}{\delta_{dec}}
\newcommand{\dee}{\delta_{extra}}
\newcommand{\dt}{\delta_B}

\def\That{{{ \widehat T}}}

\def\piT{{^{(\T)}\pi}}

\newcommand{\Lieb}{\Lie \mkern-10mu /\,}

\def\DDov{\ov{\DD}}

\newcommand{\Gabb}{{\Ga \mkern-11mu /\,}}

\def\GabbX{{^{(X)} \Gabb}}

\def\DDs{ \, \DD \hspace{-2.4pt}\dual    \mkern-20mu /}
\def\DDd{ \, \DD \hspace{-2.4pt}    \mkern-8mu /}

\def\Bdot{\dot{B}}
\def\Edot{\dot{E}}
\def\Fdot{\dot{F}}

\def\Rdot{\dot{\R}}

\def\divc{\,^{(c)}\div}
\def\curlc{\,^{(c)}\curl}
\def\DDc{\,^{(c)} \DD}

\def\DDb{\ov{\DD}}
\def\DDbc{\ov{\DDc}}

\def\Lied{\dot{\Lie}}

\def\DDov{\ov{\DD}}
\def\DDcov{\ov{\DDc}}

\def\Db{\dot{\D}}
\def\Ddot{\dot{\D}}
\def\squared{\dot{\square}}

 \def\Ft{\widetilde{F}}

\def\Eext {{\,{}^{(ext)}E}}

\def\Mor{{\mbox{Mor}}}
 \def   \Mordot{{\dot{Mor}}\hspace {-2pt}}

\def\Morr{{\mbox{Morr}}}

\def\DDdS{\,^S\hspace{-3pt}\DDd}

\def\nabS{\,^S\nab}

\def\eS{\,^{S} \hspace{-1.5pt} e}

\def\Ddot{\dot{\D}}

\def\Ddot{\dot{\D}}

\def\Bext{{\,^{(ext)} B}}
\def\Next{{\,^{(ext)} \NN}}

\def\Xho{\,^{(h)}X}
\def\Kh{\,^{(h)}K}

\def\Rhat{{\widehat{R}}}

\def\nabc{\,^{(c)}\nab}

\def\lapc{\,^{(c)}\lap}

\DeclareFontFamily{U}{mathx}{\hyphenchar\font45}
\DeclareFontShape{U}{mathx}{m}{n}{
      <5> <6> <7> <8> <9> <10>
      <10.95> <12> <14.4> <17.28> <20.74> <24.88>
      mathx10
      }{}
\DeclareSymbolFont{mathx}{U}{mathx}{m}{n}
\DeclareFontSubstitution{U}{mathx}{m}{n}
\DeclareMathAccent{\widecheck}{0}{mathx}{"71}

\def\trXc{\widecheck{\tr X}}
\def\trXbc{\widecheck{\tr\Xb}}
\def\Hc{\widecheck{H}}
\def\Zc{\widecheck{Z}}
\def\Pc{\widecheck{P}}
\def\ombc{\widecheck{\omb}}

\def\psiwc{\widecheck{\psi}}
\def\ec{\widecheck{e}}

\def\aund{{\underline{a}}}
\def\bund{{\underline{b}}}
\def\cund{{\underline{c}}}

\def\SS{\mathcal{S}}

\def\Morr{\mbox{Morr}}
\def\RRt{\tilde{\RR}}

\def\RRtp{\RRt'}
\def\RRtpp{\RRt''}
 \def\SSz{{\SS,z}}

\def\aund{{\underline{a}}}
\def\bund{{\underline{b}}}
\def\cund{{\underline{c}}}

\def\Sa{S_{\underline{a}}}
\def\Sb{S_{\underline{b}}}
\def\Sc{S_{\underline{c}}}

\def\SSa{\SS_{\underline{a}}}

\def\psia{\psi_{\aund}}
\def\psib{\psi_{\bund}}
\def\psic{\psi_{\cund}}

   \def\undpsi{\underline{\psi}} 

\def\FFab{\FF^{\aund\bund}}

\def\AAa{\AA^{\aund}}
\def\VVa{\VV^{\aund}}
 \def\Qr{\mbox{Qr}}

\def\RRa {\RR^{\underline{a}}}

      \def\gadot{\dot{\ga}}
      \def\lz{{\bf \ell_z}}
      
      \def\ntrap{trap\mkern-18 mu\big/\,}
          \def\Mtrap{\,\MM_{trap}}
\def\Mntrap{{\MM_{\ntrap}}}
 
  \def\AAt{\widetilde{\AA}}

\def\ges{\gtrsim}

 \def\DDdc{\,^{(c)} \DDd}
    \def\DDsc{\,^{(c)} \DDs}

    \def\DDc{\,^{(c)} \DD}

\def\Bdot{\dot{B}}

\def\Pdot{\dot{P}}
\def\Bbdot{\dot{\Bb}}

      \def\ntrap{trap\mkern-18 mu\big/\,}
\def\Mntrap{{\MM_{\ntrap}}}
 
  \def\AAt{\widetilde{\AA}}

\def\ges{\gtrsim}

\def\Sext{{\,^{(ext)}\Si}}

\def\Sk{\frak G}

\def\Rkext{\,^{(ext)} \Rk}
\def\Skext{\,^{(ext)} \Sk}
\def\Skint{\,^{(int)} \Sk}
\def\Rkint{\,^{(int)} \Rk}

\def\Pdot{\dot{P}}
\def\Bdot{\dot{B}}
\def\Bbdot{\dot{\Bb}}

\def\Bt{\widetilde{B}}
\def\Pt{\widetilde{P}}
\def\Et{\widetilde{E}}
 \def\NNt{\widetilde{\NN}}
\def \Bbt{\widetilde{\Bb}}
\def\At{\widetilde{A}}
\def \Abt{\widetilde{\Ab}}

\def\RRa {\RR^{\underline{a}}}

 \def\DDdc{\,^{(c)} \DDd}
    \def\DDsc{\,^{(c)} \DDs}

     \def\kl{{k_L}}
    
    \def\DDc{\,^{(c)} \DD}

\def\NNmor{{\, ^{(mor)}\NN}}
\def\NNred{{\, ^{(red)}\NN}}
\def\NNen{{\, ^{(en)}\NN}}
\def\MMred{{\,^{(red)}\MM}}

\def\BEF{B\hspace{-2.5pt}E \hspace{-2.5pt} F}
\def\EF{E \hspace{-2.5pt} F}
\def\BEFdot{\dot{\BEF}\hspace{-2.5pt}}

 \def\Ptp{\Pt_{+}}
   \def\Ptm{\Pt_{-}}

   \def\Psit{\widetilde{\Psi}}
\def\Psid{\dot{\Psi}}
\def\Psib{\und{\Psi}}
\def\Fd{\dot{F}}
\def\Psidd{\ddot{\Psi}}
\def\Fdd{\ddot{F}}
\def\Pdot{\dot{P}}
\def\Bdot{\dot{B}}
\def\Bbdot{\dot{\Bb}}
\def\Abdot{\dot{\Ab}}
\def\Adot{\dot{A}}

\def\Int{\mbox{\,I\hspace{-4.5pt} n \hspace{-8.pt} t}}

\def\und{\underline}

\def\nabcheck{\widecheck{\nab}}
\def\psicheck{\widecheck{\psi}}

\def\Xone{\, ^{(1)}X}
\def\Xtwo{\, ^{(2)}X}

\def\pione{\, ^{(1)}\pi}
\def\pionet{\widetilde{ \pione}}
\def\pitwo{\, ^{(2)}\pi}

\def\ec{\widecheck{e}}

\def\Nt{\widetilde{N}}

\begin{document}

\title{Wave equations estimates and the nonlinear stability of  slowly rotating Kerr black holes}
\author{Elena  Giorgi, Sergiu Klainerman, J\'er\'emie Szeftel}

\maketitle

{\bf Abstract.} \textit{This is the last  part    of our   proof of the  nonlinear stability of the Kerr family for small angular momentum, i.e $|a|/m\ll 1$, in which we  deal with    the nonlinear wave  type estimates needed to complete the project.  More precisely  we provide  complete  proofs for  Theorems M1 and M2 as well the curvature  estimates of Theorem M8,   which were  stated without proof  in sections 3.7.1 and 9.4.7 of \cite{KS:Kerr}. Our procedure is based on a new   general interest formalism  (detailed in Part I of this work), which extends  the one used in the stability of Minkowski space.   Together with  \cite{KS:Kerr} and the GCM papers    \cite{KS-GCM1},    \cite{KS-GCM2}, \cite{Shen},    this  work   completes  proof of the Main   Theorem stated in Section 3.4 of \cite{KS:Kerr}.}

\tableofcontents


\chapter{Introduction}


This is the last  part    of our   proof of the  nonlinear stability of the Kerr family for small angular momentum, i.e $|a|/m\ll 1$, in which we  deal with    the nonlinear wave  type estimates needed to complete the project.  More precisely  we provide  complete  proofs for  Theorems M1 and M2 as well the curvature  estimates of Theorem M8,   which were  stated without proof  in sections 3.7.1 and 9.4.7 of \cite{KS:Kerr}. Our procedure is based on a new   general interest formalism  (detailed in Part I of this work), which extends  the one used in the stability of Minkowski space.  Together with  \cite{KS:Kerr} and the GCM papers    \cite{KS-GCM1},    \cite{KS-GCM2}, \cite{Shen},    this  work   completes  proof of the Main   Theorem stated in Section 3.4 of \cite{KS:Kerr}.


\section{Black hole stability problem}
  

We  give below a quick introduction to the black  hole stability problem.

  
\subsection{Einstein Vacuum equations}


 We   restrict our attention to   the Einstein vacuum equations (EVE), i.e.  spacetimes  $(\MM, \g)$  with  vanishing Ricci curvature, i.e. 
\bea
\lab{EVE-intro}
\R_{\a\b}=0.
\eea
Note that a solution of \eqref{EVE-intro}  is  in fact a class of  equivalence of solutions with respect to   diffeomorphisms  $\Phi $ of $\MM$, i.e.  $\g$ and $\Phi^*\g$ are indistinguishable  as solutions of EVE. This is precisely what is meant by the general covariance of the   Einstein equations.

\eqref{EVE-intro} corresponds to an evolution problem and an associated initial
data set $(\Si_0, g_{(0)}, k_{(0)})$  for  EVE  consists of $3$ dimensional manifold $\Si$
 together with  a  complete Riemannian metric $g_{(0)}$ and  a symmetric 2-tensor $k_{(0)}$ which verify  compatibility conditions known as the constraint equations.  A Cauchy development of an initial data set  is a globally hyperbolic  space-time   $(\MM,g)$, verifying EVE
and  an  embedding $i:\Si\to \MM$ such that 
$i_*(g_{(0)}), i_*(k_{(0)})$  are the first and second fundamental forms
of $i(\Si_{(0)})$ in $\MM$. A well known  foundational result  in GR
 associates a unique   maximal, global hyperbolic,  future development 
    to all sufficiently regular  initial data sets, see \cite{Br} \cite{Br-Ge}.   
We  further  restrict  the discussion  to asymptotically flat initial data sets,
i.e.  we assume that outside a sufficiently large compact set $K$,
 $\Si_{(0)}\setminus K$ is diffeomorphic  to the complement of the unit
 ball  in $\RRR^3$ and admits a system of coordinates in which $ g_{(0)}$ is asymptotically euclidean and $k_{(0)}$ vanishes at appropriate order.

   
\subsection{Kerr family} 
  

    (EVE) admits a remarkable   two parameter family of explicit  solutions, the  Kerr spacetimes  $\KK(a,m)$, $0\le |a|\le m$,  which are  stationary and axisymmetric. In the usual  Boyer-Lindquist coordinates  they take the form 
   \begin{equation}\label{zq1}
\g=-\frac{q^2\Delta}{\Sigma^2}(dt)^2+\frac{\Sigma^2(\sin\theta)^2}{|q|^2}\left(d\phi-\frac{2amr}{\Sigma^2}dt\right)^2 +\frac{|q|^2}{\Delta}(dr)^2+|q|^2(d\theta)^2,
\end{equation}
where
\begin{equation}\label{zq2}
\begin{cases}
&\Delta=r^2+a^2-2mr, \qquad  q=r+i a \cos \th,\\
&\Sigma^2=(r^2+a^2)|q|^2+2mra^2(\sin\theta)^2=(r^2+a^2)^2-a^2(\sin\theta)^2\Delta.
\end{cases}
\end{equation}
 Among  them 
  one distinguishes the Schwarzschild family of spherically symmetric   solutions,  of mass $m>0$,
\bea
\g=-\left(1-\frac{2m}{r}\right) dt^2+\left(1-\frac{2m}{r}\right)^{-1} dr^2 +r^2 d\si_{\SSS^2}. 
\label{eq:Schw}
\eea 
Though the metric seems singular at $r=2m$   ($r=r_+$, the largest  root   of $\De(r)=0$,  in the case of Kerr)
it turns out that one can glue together two regions $r>2m$ and two regions $r<2m$ of the Schwarzschild metric to obtain a metric which is smooth along the null hypersurface   $\EE=\{r=2m\}$
 called the Schwarzschild  event  horizon.   The portion of  $r<2m$  to the future
 of the hypersurface  $t=0$ is a black hole whose future boundary $r=0$ is 
 singular.    The region $r>2m$, free of  singularities,  is called the domain of outer communication.  
  The more general family of Kerr solutions, which   are both stationary and axially symmetric,  possesses (in addition to   well defined
  event horizons, black holes and domains of outer communication)   Cauchy horizons  ($r=r_-$, the smallest root  of $\De(r)=0$) inside the black hole region   across which predictability fails\footnote{Infinitely many smooth extensions are possible beyond the  boundary.}.  Once more, one
  can easily check, from the precise nature of the Kerr metric, that the region
  outside  the event horizon,  i.e. outside the  Kerr black hole,  is free of singularities\footnote{Consistent with  the   \textit{weak cosmic censorship conjecture (WCC)} of Penrose.}.  
  
 \begin{figure}[h!]
\centering
   \includegraphics[width=4.5in]{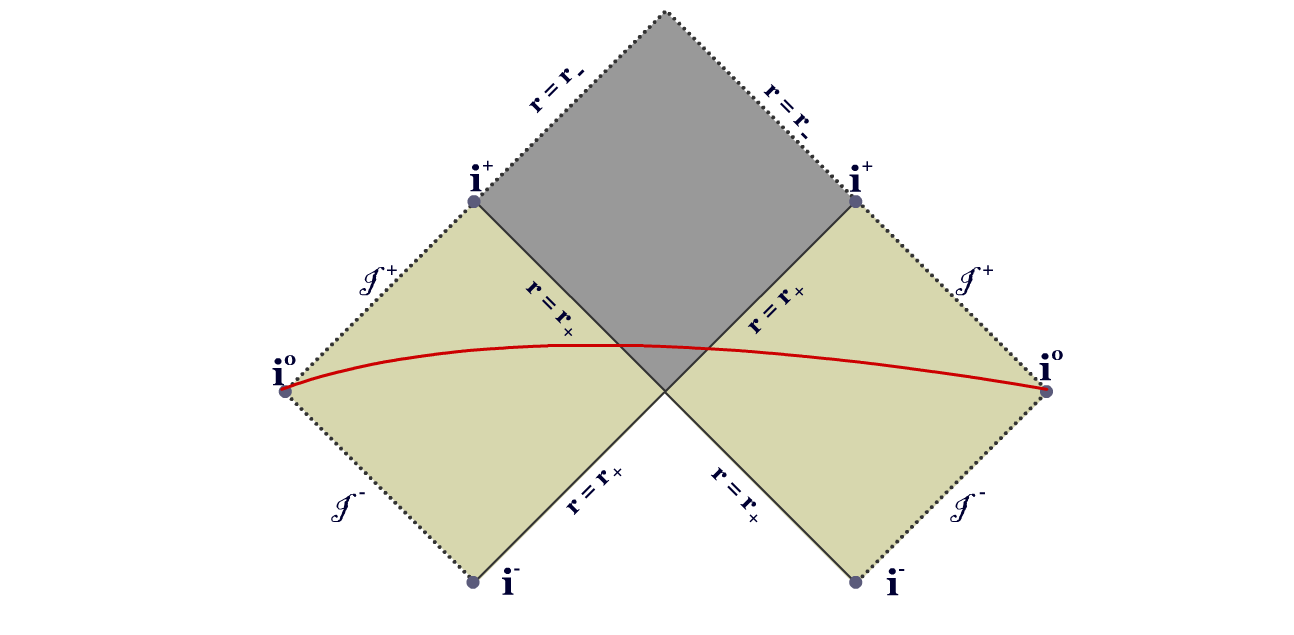}  
\caption{The Penrose diagram of Kerr for $0<|a|<m$.}
\end{figure}

   Finally we  note that    the Kerr spacetimes 
  $\KK(a,m)$ possess two Killing vectorfields: the stationary  vectorfield    $\T=\pr_t$,  which is time-like   in the 
  asymptotic region, away  from the horizon, and  the  axial symmetric Killing field $\Z=\pr_\varphi$.

  
\subsection{Stability of Kerr conjecture}  

  
   The issue of the stability\footnote{This is not only a deep mathematical 
    question but one with serious astrophysical  implications.
     Indeed, if  the Kerr family would be  unstable, 
      black holes  would be nothing more than mathematical artifacts.} of the Kerr    family 
      has been at the center of attention  of GR physics and  mathematical relativity for   more than half a century, ever since their  discovery   by Kerr in \cite{Kerr}.  Roughly the problem here is to show that all spacetime developments of  initial data sets,    sufficiently close to  that of a Kerr spacetime, behave  asymptotically  like a    Kerr solution.  Here is a more precise formulation of the conjecture.  

   {\bf Kerr stability conjecture.}\,\,{\it  Vacuum initial data sets
   sufficiently close to Kerr initial data have a maximal development with complete
   future null infinity\footnote{Thus observers  which are  far away  
    from the black hole  will never  experience  its effects.} and with a  domain of outer communication which
   approaches  (globally)  a nearby Kerr solution.}
   
    Until very recently  the only   space-time for which  full nonlinear  stability had been established was the Minkowski space. The proof is based  on  some     important  PDE  advances    of late last century:
\begin{enumerate}
\item[(i)]  Robust   vectorfield  approach  to   derive quantitative   decay  based on generalized energy estimates and commutation with (approximate)  Killing and conformal killing  vectorfields. 

\item[(ii)]   The  \textit{null condition} identifying the deep  mechanism for     nonlinear stability, 
i.e.  the specific structure of the nonlinear terms   which   enables stability  despite the low
  decay of the  perturbations. 
  
\item[(iii)]  Elaborate bootstrap argument according to which  one makes  educated assumptions 
about the behavior of  solutions  to  nonlinear wave equations  and then proceeds, by a long sequence of a-priori estimates,  to       show that   they are in fact 
 satisfied. This  amounts to a \textit{conceptual linearization}, i.e. a method  by which the equations become, essentially, linear\footnote{Note that in the context of  EVE, and other quasilinear hyperbolic systems, this    differs substantially from the usual notion of linearization around a fixed background.}  without   actually linearizing them.
 \end{enumerate}
 
  There are three, related,  major obstacles in passing from the stability of Minkowski  to that of  Kerr:
   \begin{enumerate}
 \item The first   can be understood in the general framework of nonlinear evolution equations.
 Given a nonlinear evolution equation $\NN[\phi]=0$ and a stationary solution $\phi_0$,   we have two notions of stability, \textit{orbital stability}, according to which  small perturbations    of $\phi_0$    lead to solutions $\phi$ which remain close for all time, and  \textit{asymptotic stability (AS)} according to which the perturbed solutions converge as $t\to \infty$ to $\phi_0$. In the case where $\phi_0$ is non trivial, there is a third notion of stability, which we call \textit{asymptotic orbital stability (AOS)}, to describe the fact that  the perturbed solutions may converge to a different stationary solution\footnote{This happens if $\phi_0$ belongs to a multi parameter smooth family of stationary solutions, or by applying a gauge transform to $\phi_0$ which keeps the equation invariant. In the case of Kerr, both cases are present as we shall see below.}. For quasilinear equations\footnote{Orbital stability can be  established directly (i.e. without   establishing the stronger version) only   for  Hamiltonian equations with weak  nonlinearities.}, such as  EVE,   a proof of stability means necessarily  AS  or AOS   stability.      Both    require    a    detailed  understanding of the decay properties    of the   linearized equation.  
   One is thus led to  study  the linearized equation $\LL[\phi_0]\psi =0$,  with $\LL[\phi_0]$ the Fr\'echet  derivative   $\NN'[\phi_0]$, which  is,  essentially, 
  a linear hyperbolic  system  with variable coefficients  and,  typically,  presents   instabilities. 
  In the exceptional  situation,  when  stability can ultimately  be established,   one  can  tie  all   the  instability   modes to  either  the  gauge invariance of the equation  or  the presence of a  continuum of other distinct\footnote{I.e. solutions not tied to $\phi_0$ by a gauge transformation.} stationary solutions\footnote{In the  case  of the stability of   Kerr,  there exists a $2$ parameter family of solutions $\KK(a,m)$.} near $\phi_0$. These instabilities at the linearized level are responsible for the fact   that a small perturbation of  the   fixed stationary solution $\phi_0$  may not converge to  $\phi_0$ but to  another nearby  stationary solution,  this  is  the case of  AOS. The methodology of tracking this asymptotic final state, in general different from $\phi_0$, is usually referred to as modulation.  In the case   of the Einstein equations, this problem is compounded by the presence of infinitely many instabilities related to full group of diffeomorphism, i.e. to the general covariance of the Einstein equations\footnote{Note that in  the stability of Minkowski, even though the linearized system does not contain instabilities, one must still take general covariance into account in the far $r$ region of the perturbed space-time due the presence of a non trivial mass. On the other hand, in perturbations of Kerr, the general covariance affects the entire construction of the spacetime.}. 
         
   \item A  fundamental insight   in the stability of the Minkowski space  was that the Bianchi identities  decouple at first order from the null structure equations which allows one to control curvature first, as a Maxwell type system (see \cite{Ch-Kl0}), and then proceed with the rest of the solution. This  cannot work for  perturbations of Kerr due to the fact that some of the null components\footnote{With respect to  the  \textit{principal null directions of Kerr,} i.e   a distinguished  null  pair which   diagonalizes  the  full curvature tensor, the middle component $P=\rho+ i\rhod$ is nontrivial.} of the  curvature tensor  are non-trivial in Kerr.
 
 \item    Even if  one succeeds in tackling the above mentioned issues, there are still major obstacles  in understanding the  decay properties of the solution.  Indeed,  when one  considers  the simplest,    relevant,
  linear equation  on a fixed Kerr background,  i.e. the scalar wave equation $\square_\g \psi=0$,  one 
  encounters serious difficulties  to  prove decay.   Below is a  very short description of these:
       \begin{itemize}
        \item \textit{The problem of   trapped null geodesics.}   This  concerns the existence of null geodesics\footnote{In the Schwarzschild case, these geodesics
are associated with the celebrated photon sphere $r=3m$.} neither crossing the event horizon nor escaping to 
null infinity, along which solutions   can concentrate for arbitrary long
times.  This   leads  to degenerate  energy-Morawetz estimates  which require a very delicate analysis.
       
\item \textit{The   trapping   properties of the     horizon.}       
    The horizon  itself  is ruled  by null  geodesics, which do not     communicate with  null infinity and can thus concentrate energy.         This problem was solved   by understanding      the  so called     red-shift effect associated to  the event horizon,  which     counteracts this type of  trapping.

\item \textit{The problem of superradiance.} This  is  the failure of the stationary Killing 
field $\T=\partial_t$
to be everywhere timelike in the domain of outer communication\footnote{$\T$ is timelike only outside of the so-called ergoregion.}, and thus, of the
associated  conserved energy to be positive. Note that this problem is absent in Schwarzschild and, in general,   for axially symmetric solutions of EVE.
\item   \textit{Superposition problem.}       This is  the  problem  of combining the estimates  in the near  region, close to the horizon, (including the  ergoregion  and trapping) with estimates in the asymptotic region, where the spacetime looks Minkowskian.
\end{itemize}

   \item   The full linearized system, whatever  its  formulation,  presents   many additional difficulties  due to the  huge gauge covariance of the  equations\footnote{Note that rates of decay are heavily dependent on a proper choice of gauge, thus affecting the issue of convergence.}. In particular, the full linearized system is  not  conservative and  we  thus  lack, unlike  in the case of  the  scalar wave equation
 $\square_\g \psi=0$, the most basic ingredient  in  controlling     the   solutions of the equation, i.e. energy estimates. 
\end{enumerate}

  
\section{Linear stability}


  
\subsection{Linearized gravity system}


 Historically,    two versions of linearization for  EVE  have been   considered:
   
   \begin{itemize}
   \item[(a)] At the level of the metric  itself, i.e.  with   $\G_{\a\b}:=\R_{\a\b}-\frac 1 2 \R \g_{\a\b}$, 
   \bea
\label{Lin. Grav}
\G'(\g_0) \, \de \g =0.
\eea

   \item[(b)]  At the level of curvature  via  the Newman-Penrose (NP) formalism,  based  on null frames. 
   \end{itemize}
   
In our work, we rely on a geometric variant of the second approach, see the comparison between the NP formalism, the Geroch-Held-Penrose formalism and our approach in section \ref{section:NPformalism}. In what  follows, we  describe the main known  results  concerning  solutions to linearized equations in a Kerr background.


\subsection{Formal mode analysis}


The first important  results  concerning both items  (a) and (b)  above were  obtained by mathematical  physicists  based on the classical  method 
of  separation of variables  and   formal mode analysis.
  In the particular case when $\g_0$ is the Schwarzschild metric, the  LGE equations \eqref{Lin. Grav}   
can be formally decomposed into modes,  by using Fourier transform in time and spherical harmonics. 
A similar  decomposition, using oblate  spheroidal harmonics, can be done in Kerr.
 The formal study of fixed modes from the point of view of \textit{metric perturbations}  as in \eqref{Lin. Grav}   was initiated by  Regge-Wheeler \cite{Re-W}, in perturbations of  Schwarzschild,  who
 discovered      the master Regge-Wheeler equation for  axial or odd-parity perturbations.
 This study was completed by Vishveshwara \cite{Vishev}       and Zerilli  \cite{Ze}.  A gauge-invariant formulation of \textit{metric perturbations}  was then given by Moncrief  \cite{Moncr}.  An alternative approach via the Newman-Penrose  (NP) formalism      was  first  undertaken  by Bardeen-Press \cite{Bar-Press}. This latter type of analysis was later extended to the Kerr family by Teukolsky  \cite{Teuk} who made the  important  discovery that  the extreme  curvature components, relative   to  a principal null frame, satisfy  decoupled, separable,  wave equations.   These extreme curvature components    also turn out to be gauge invariant  in the sense discussed above.     The  full  extent   of what could be done  by  mode analysis, in both approaches,    can be found in    Chandrasekhar's   book  \cite{Chand1}.
  Chandrasekhar        also introduced (see \cite{Chand2}) a  transformation theory relating  the two approaches.          More precisely, he exhibited a transformation which connects  the Teukolsky  equations  to the       Regge-Wheeler  one.  This  transformation    was further  elucidated and extended  by R. Wald \cite{Wald}. The full  mode stability, i.e. lack of exponentially growing modes,  for the Teukolsky equation in Kerr is  due  to  Whiting \cite{Whit}.


\subsection{Classical vectorfield method}  


  Mode stability     is far from establishing  even boundedness of solutions. To achieve that and, in addition, to   derive realistic decay estimates 
one needs  an entirely different approach based on  the vectorfield  method\footnote{Robust method based on the symmetries of Minkowski space   to derive decay  for nonlinear wave equations, see  \cite{Kl-vectorfield} and \cite{Kl-null}.},  used in the proof of  the nonlinear stability of Minkowski \cite{Ch-Kl}.

The  vectorfield method  was first   developed in connection   with the wave  equation  in Minkowski space.    
As well known,  solutions of the wave equation  $ \square \psi=0$ in the Minkowski space  $\RRR^{n+1}$  
  both conserve energy  and decay uniformly in time   like $t^{-\frac{n-1}{2}}$. While conservation of energy  can be established  by a simple integration by parts,  and is thus robust to perturbations of the Minkowski metric,  decay   was first derived   either using  the   Kirchhoff formula  or by Fourier methods, which are manifestly not robust.  An integrated version of local energy  decay, based on an  inspired  integration by parts argument,   was first derived by  C. Morawetz \cite{Mor1}, \cite{Mor2}.  The first    derivation of decay 
   based  on   the  commutations properties of  $\square$ with     Killing and  conformal Killing  vectorfields    of Minkowski space
    together with    energy conservation 
   appear in \cite{Kl-vectorfield} and \cite{Kl-vect}.   That method also  provides  precise  information about  the    decay properties  of  derivatives of solutions  with respect to   the standard null frame 
    of Minkowski space,   an important motivating factor in the discovery of the  null condition \cite{Kl-ICM}, \cite{Chr} and \cite{Kl-null}.    The  methodology  initiated with these papers,    to which we  refer      as the classical vectorfield method,
      has had numerous applications to nonlinear wave equations and has played an important  role in  the proof of the nonlinear stability of Minkowski space \cite{Ch-Kl}. The vectorfield method  has also been applied    to   later versions  of the stability  of Minkowski   result in  \cite{KlNi},   \cite{Lind-Rodn},  \cite{Bi}, \cite{HVas2}, and extensions of it   to  Einstein equation coupled with  various matter  fields in \cite{BiZi},  \cite{FJS}, \cite{BFJT}, \cite{Wa},   \cite{Lf-Ma},   \cite{Lind-T},   \cite{I-Pau}.


\subsection{Scalar wave equation  in  Kerr   and new vectorfield method} 
 
 
  The new vectorfield method is an extension of the classical  vectorfield method  which  compensates  for   the lack of  enough Killing and conformal Killing vectorfields   on a Schwarzschild  or  Kerr  background  by   introducing  new vectorfields  whose deformation tensors have    coercive properties    in different regions of spacetime, not necessarily causal.    The new   method   has emerged  in the last fifteen years   in connection to the study of boundedness and decay  for the    scalar wave equation in $\KK(a,m)$, $\square_{\g_{a,m}} \phi=0$. The starting and most  demanding   part of the new method, originating in \cite{B-S1}, is  the derivation of  a global,  simultaneous,  \textit{Energy-Morawetz}   estimate  which degenerates   in the trapping region.  Once an Energy-Morawetz estimate is established  one  can  commute  with the   Killing    vectorfields   of the background, and  the so called red shift  vectorfield  introduced in  \cite{DaRo1},   to derive uniform bounds  for  solutions.  The most efficient  way  to also  get decay, and  solve the\textit{ superposition  problem}, originating in \cite{Da-Ro3},  is  based on   the presence of   family of \textit{$r^p$-weighted},  quasi-conformal  vectorfields  defined in the  far $r$ region of spacetime\footnote{These replace the scaling  and inverted time translation  vectorfields used in \cite{Kl-vectorfield}  or their  corresponding  deformations used in \cite{Ch-Kl}.  A recent improvement of the method  allowing one to derive higher order decay  can be found in \cite{AnArGa}.}.

   The first Energy-Morawetz  type  results  for scalar wave  equations   in Schwarzschild   are due to  Soffer-Blue  \cite{B-S1}, \cite{B-S2} and  Blue-Sterbenz \cite{B-St},    based on a modified version of the classical  Morawetz  integral energy decay  estimate.     Further developments appear in the works of Dafermos-Rodnianski\footnote{We note in particular the red shift   vectorfield, introduced  in  \cite{DaRo1} to   deal with the degeneracy 
       of the  Morawetz-energy estimates along the horizon  and the $r^p$ weighted  estimates introduced in  \cite{Da-Ro3}  as an effective  method to derive    decay estimates   in the  asymptotic  region  of black holes, as mentioned above.}, see  \cite{DaRo1},   \cite{Da-Ro3},  and   Marzuola-Metcalfe-Tataru-Tohaneanu,    \cite{MaMeTaTo}. 
    The vectorfield method    can also be extended 	 to derive decay   for axially symmetric solutions   in Kerr, see    \cite{I-Kl1} and \cite{St}, but it is  known to fail   for  general  solutions  in Kerr, see \cite{Al}.

      In the absence of axial symmetry   the derivation of    an    Energy-Morawetz estimate   in    $\KK(a,m)$ for  $|a/m|\ll 1 $ requires  a more refined  analysis involving     both the vectorfield method and    Fourier or mode  decompositions,   see Tataru-Tohaneanu \cite{Ta-Toh}  for  the first full quantitative  decay  result (see also Dafermos-Rodnianski   \cite{DaRo2}  for boundedness of solutions).  The   derivation of    such an estimate  in the full sub-extremal case $|a|<m$  is even more subtle  and   was  achieved  by Dafermos-Rodnianski-Shlapentokh-Rothman \cite{mDiRySR2014}.    A  purely physical   space proof   the  Energy-Morawetz estimate   for small  $|a/m|$,  which extends the classical vectorfield method to  include    second order  operators (in this case  the Carter operator \cite{Carter})   was pioneered   by  Andersson-Blue in   \cite{A-B}.     Their approach  has the usual advantages of  the classical  vectorfield method, i.e it is robust with respect to perturbations, which is the very  reason  we pursue it  in this paper.

       
\subsection{Linear  stability  of Schwarzschild}  
  

    A  first  quantitative  proof   of the  linear stability of Schwarzschild  spacetime  was established  by  Dafermos-Holzegel-Rodnianski  in \cite{D-H-R}. 
   Notable in their analysis  is  
  the treatment of  the Teukolsky equation (TE)  in  a fixed Schwarzschild background.   While  TE    is separable, and  amenable to mode analysis, it  is not Lagrangian and  thus cannot  be treated   by      energy type estimates.  It is for this reason that    \cite{D-H-R}  relies  on a physical space  version of   the  Chandrasekhar transformation, which takes solutions  of  TE  to  solutions of  Regge-Wheeler (RW), an equation  which is manifestly Lagrangian, and in addition has a positive potential. Once  decay  estimates  for   the RW  equation    have   been established\footnote{Based on the technology developed earlier for the scalar wave equation in Schwarzschild.},    the authors  recover   the  expected decay for  solutions to  the original   TE.     The remaining work in \cite{D-H-R}   is to       derive decay  for the other curvature  components and most of the  linearized  Ricci coefficients associated to the double null foliation.  This   last step  requires carefully chosen gauge conditions, which the authors   make  within the framework of  a double null   foliation,  initialized by  the  given  Cauchy  data\footnote{\lab{footnote1}Note that there is also a scalar condition for the linearized lapse along the event horizon (part of what the authors call future normalized gauge), itself initialized from initial  data, see (212) and (214) in \cite{D-H-R}. This gauge fixing  from initial data leads to sub-optimal decay estimates  for some metric coefficients (see (250)--(252)  and  (254) in \cite{D-H-R}) and potentially for $\omb$,  and is  thus inapplicable to the nonlinear case. This deficiency  was  fixed in  \cite{Giorgi-T}, by relying on  a linearized version of the GCM construction  in \cite{KS}.}  of  the fixed Schwarzschild background.   
        
  
\subsection{Linear stability of  Kerr for small angular momentum}


The first  breakthrough       result on the linear stability of Kerr,      for $|a|/m\ll 1$,  is due  to  Ma \cite{Ma}, see also  \cite{D-H-R-Kerr}. Both results  are  based on a generalization of the Chandrasekhar transformation to Kerr which   takes the Teukolsky equations   to a generalized version of the  Regge-Wheeler (gRW)  equation.  Also, both methods depend on a combination of mode decomposition and vectorfield techniques  similar to those developed   for  the scalar wave equation in slowly rotating  Kerr.  These  results   were  recently partially\footnote{The analysis  of  \cite{Y-R} is still  limited  to  modes. The authors   have  however announced a full proof    of  the result.}  extended to the full subextremal range in    \cite{Y-R}.
 
The first   stability results for the full linearized   Einstein vacuum equations   near   $Kerr(a, m)$,   for $|a|/m\ll 1$,  appeared  in    \cite{ABBMa2019}  and  \cite{HHV}.  The first paper,  based on the   GHP formalism\footnote{An adapted  spinorial version  of the    NP formalism.}, see   \cite{GHP},    builds  on the results of  \cite{Ma}  while the second paper is based on an adapted  version of  the   metric formalism and builds on the seminal work of the authors on Kerr-de-sitter \cite{HVas}. Though the ultimate  relevance  of these  papers to  nonlinear stability  remains open,   they are  both  remarkable results  in so far as they deal with difficulties that looked insurmountable even ten years ago.

  
\section{Nonlinear stability}



\subsection{Nonlinear stability of Schwarzschild}  

 
      The first  nonlinear stability result of the Schwarzschild space  was established  in  \cite{KS}.    In its simplest version, the  result  states the following.

\begin{theorem}[Klainerman-Szeftel \cite{KS}] 
\lab{MainThm-firstversion-KS}
The future globally hyperbolic   development  of  an   \emph{axially symmetric, polarized},  asymptotically  flat   initial data set, sufficiently close  (in a specified  topology)  to a Schwarzschild  initial data set   of  mass $m_0>0$,  has a complete    future null infinity  $\II^+$ and converges 
in  its causal past  $\JJ^{-1}(\II^{+})$  to another  nearby  Schwarzschild solution of mass $m_{\infty}$ close to $m_0$. 
 \end{theorem}

The restriction to  axial polarized perturbations  is the  simplest  assumption   which insures  that the final state is itself Schwarzschild  and thus avoids the additional complications  of  the Kerr stability problem which we discuss below.  We refer  the reader to  the introduction in \cite{KS} for a full discussion of the result.

Recently Dafermos-Holzegel-Rodnianski-Taylor   \cite{DHRT} have   extended    the result of \cite{KS}      by properly  preparing a co-dimension 3 subset of the initial data such that the final state is  still Schwarzschild. Like in   \cite{KS},  the starting point of  \cite{DHRT} is to anchor  the entire construction   on a  far away\footnote{I.e. $r\gg u$, similar to the dominant in $r$ condition (3.3.4) of \cite{KS}.} GCM type  sphere, in the sense of \cite{KS-GCM1} \cite{KS-GCM2}, with no  direct reference to the initial data. It also  uses  the same definition of the  angular momentum as in (7.19) of \cite{KS-GCM2}. Finally, the spacetime in \cite{DHRT} is separated in an exterior region $\Mext$ and an interior region $\Mint$,   with the ingoing  foliation of $\Mint$ initialized on a  timelike hypersurface,   as in \cite{KS}. We note, however, that \cite{DHRT} does not use the geodesic foliation of \cite{KS}, but instead both  $\Mint$ and $\Mext$ are foliated  by double null    foliations, and thus, the process of estimating the gauge dependent variables is somewhat different.


\subsection{Nonlinear stability of  slowly rotating Kerr black holes}
        

        In \cite{KS:Kerr},  we  have stated the following theorem on the resolution of the 
        Kerr stability conjecture for  small angular momentum whose complete proof relies on \cite{KS:Kerr}, \cite{KS-GCM1}, \cite{KS-GCM2}, \cite{Shen} and the present paper.
        \begin{theorem}[Kerr stability for $|a|/m\ll 1$] 
\lab{MainThm-firstversion}
The future globally hyperbolic   development  of  a general,   asymptotically  flat,    initial data set, sufficiently close  (in a suitable  topology)  to a   $Kerr(a_0, m_0) $   initial data set,  
  for sufficiently  small $|a_0|/m_0$,  has a complete    future null infinity  $\II^+$ and converges 
in  its causal past  $\JJ^{-1}(\II^{+})$  to another  nearby Kerr spacetime $Kerr(a_\infty, m_\infty )$ with parameters    $(a_{\infty}, m_\infty)$ close to the initial ones $(a_0, m_0)$.
 \end{theorem}
 
 \begin{figure}[h!]
\centering
\includegraphics[scale=0.6]{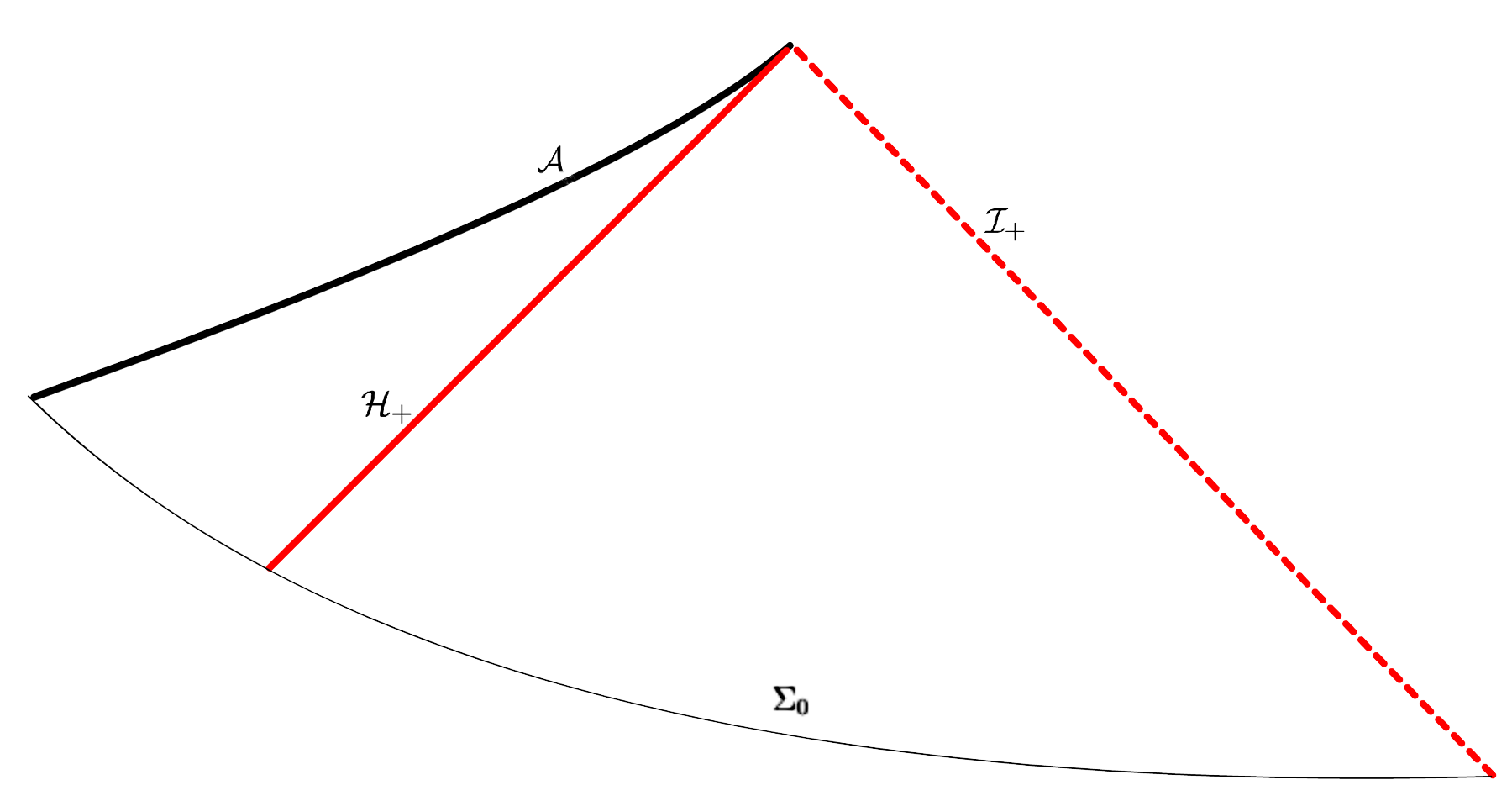}
\caption{The Penrose diagram of the  final   space-time in  the Main Theorem \ref{MainThm-firstversion}.}
\lab{fig0-introd}
\end{figure}

\begin{remark}
The analog of our Main Theorem \ref{MainThm-firstversion} in the context of the Einstein vacuum equation  with  a strictly positive cosmological constant, i.e. the nonlinear   stability   of the stationary part  of  Kerr-de Sitter with small angular momentum,  was proved in the seminal paper \cite{HVas} of Hintz and Vasy, see our discussion in section 1.1.4 in \cite{KS:Kerr}.
\end{remark}

The  proof of the Main Theorem \ref{MainThm-firstversion} rests on     the following   major    ingredients:
\begin{enumerate}
\item   A  formalism  to derive   tensorial versions of the  Teukolsky  and   generalized Regge-Wheeler type (gRW)   equations  in the   full nonlinear setting.       The formalism,   first introduced in \cite{GKS1},  is  self-contained and vastly expanded in  Part I  of this paper.

\item An analytic mechanism to        derive  combined Morawetz-energy estimates  for solutions of   these  gRW  equations,
 based on an extension of the  Andersson-Blue method, introduced in  \cite{A-B}, to  spin-2 wave equations in  suitable perturbations of Kerr.  
  This  is developed in  Part II of this paper.

\item  A   dynamical mechanism for    finding the right gauge conditions, based on GCM (generally  covariant modulated) spheres and hypersurfaces,    in which convergence to the final state  takes place. GCM spheres are codimension 2 compact surfaces, unrelated   to the initial conditions,   on which  specific geometric quantities take Schwarzschildian values  (made possible by taking into account the full  general covariance of the Einstein vacuum equations), see the   discussion in the introductions to  \cite{KS-GCM1}, \cite{KS-GCM2}. It  is hard to overstate the importance of  admissible GCM spheres,     they are  literally   the keystone of our  entire approach      to the  proof  of the nonlinear  stability  of Kerr. 
 The related  concepts of        GCM   admissible  spheres and hypersurfaces (these are codimension-1 spacelike  hypersurfaces    foliated by GCM spheres,  where additional conditions are verified)  have appeared first in \cite{KS} in the context of polarized  symmetry.   The  construction of GCM spheres,  without any  symmetries, in realistic perturbations of Kerr, is treated in \cite{KS-GCM1}, \cite{KS-GCM2}\footnote{See also chapter 16 of \cite{DHRT} in the particular case of perturbations of Schwarzschild, where they appear instead under the name  ``teleological".}, and the case of spacelike GCM hypersurfaces is treated in \cite{Shen}.

\item A dynamical mechanism  to  identify the   values of   $(a_\infty, m_\infty)$ and the  axis of rotation of the final Kerr, see  sections 3.2.4 and 8.5.2 in \cite{KS:Kerr}.  This is based on the fact that our GCM approach allows us to define the   mass $m$, as well as the  angular momentum $a$ in terms of intrinsically defined (using effective uniformization, see \cite{KS-GCM2})  $\ell=1$  modes of $\curl \b$.   This was introduced\footnote{Note that our definition of angular momentum, see (7.19) in \cite{KS-GCM2}, is  unnatural from a  physical point of view (though  very effective for our proof). A more realistic definition was introduced in \cite{Rizzi}, and another  general  definition  can be found in \cite{Chen}, see also  \cite{Sz} for a comprehensive discussion of the subject.}  in  \cite{KS-GCM2} and used in \cite{KS:Kerr}\footnote{As well as in \cite{DHRT} in the particular case of perturbations of Schwarzschild.}.

\item  As mentioned above, the main novelty of  the GCM approach   is that it    relies on  gauge conditions initialized at a far away co-dimension $2$ sphere  $S_*$ with no direct reference to the initial conditions.    
This  gauge  choice needs however  to be connected, somehow,    to the initial  conditions. This is achieved in  both \cite{KS} and  \cite{KS:Kerr} by transporting\footnote{That is, we transport the $\ell=1$ modes of some quantities from $S_*$ to $S_1$, see section 8.3.1 in \cite{KS:Kerr}.} the   sphere   $S_*$  to  a  sphere $S_1$ in the   the initial  layer        and  compare  it, using the rigidity properties 
of the GCM conditions, to a sphere of  the initial data layer. This  leads to a new foliation  of  the initial  layer  which  differs substantially from  the  original one, due to a shift  of the center of mass frame of the  final black hole, known  in the physics literature as a gravitational wave recoil. We refer the reader to  section 8.3  in \cite{KS:Kerr} for the details.

\item  A precisely formulated continuity argument,    based on a grand bootstrap   scheme,   which assigns to all geometric quantities  involved in the process   specific  decay rates,  which can   then  be   dynamically recovered  from the initial  conditions   by a long series  of estimates,   and thus ensure  convergence to a   final Kerr  state, see sections 3.5 and 3.7 in \cite{KS:Kerr}. 
\item  The continuity  argument is based  on the crucial concept of   finite  \textit{GCM admissible} spacetimes $\MM=\Mext\cup \Mint\cup \Mtop$, see Figure \ref{fig1-introd}, whose
 defining characteristic is  its spacelike \textit{GCM  boundary} $\Si_*$.   Note that the boundaries  $\Mext\cap\Mtop $ and $\Mint\cap \Mtop $ are timelike\footnote{Asymptotically null   as we pass to the limit.} and that $\Mtop$ is  needed   to have the entire space $\MM$  causal.   The regions $\Mext$ and $\Mint $ are separated by the timelike hypersurface $\TT$ and the spacelike boundary $\AA$ is beyond the future horizon $\HH_+$ of the limiting space.       Finally the region $\LL_0$,  is the initial data  layer   in which $\MM$ is  prescribed  as a solution of the Einstein vacuum equations. We   direct the reader  to section 3.2 of \cite{KS:Kerr} for more details on the construction.
\end{enumerate}

\begin{figure}[h!]
\centering
\includegraphics[scale=0.9]{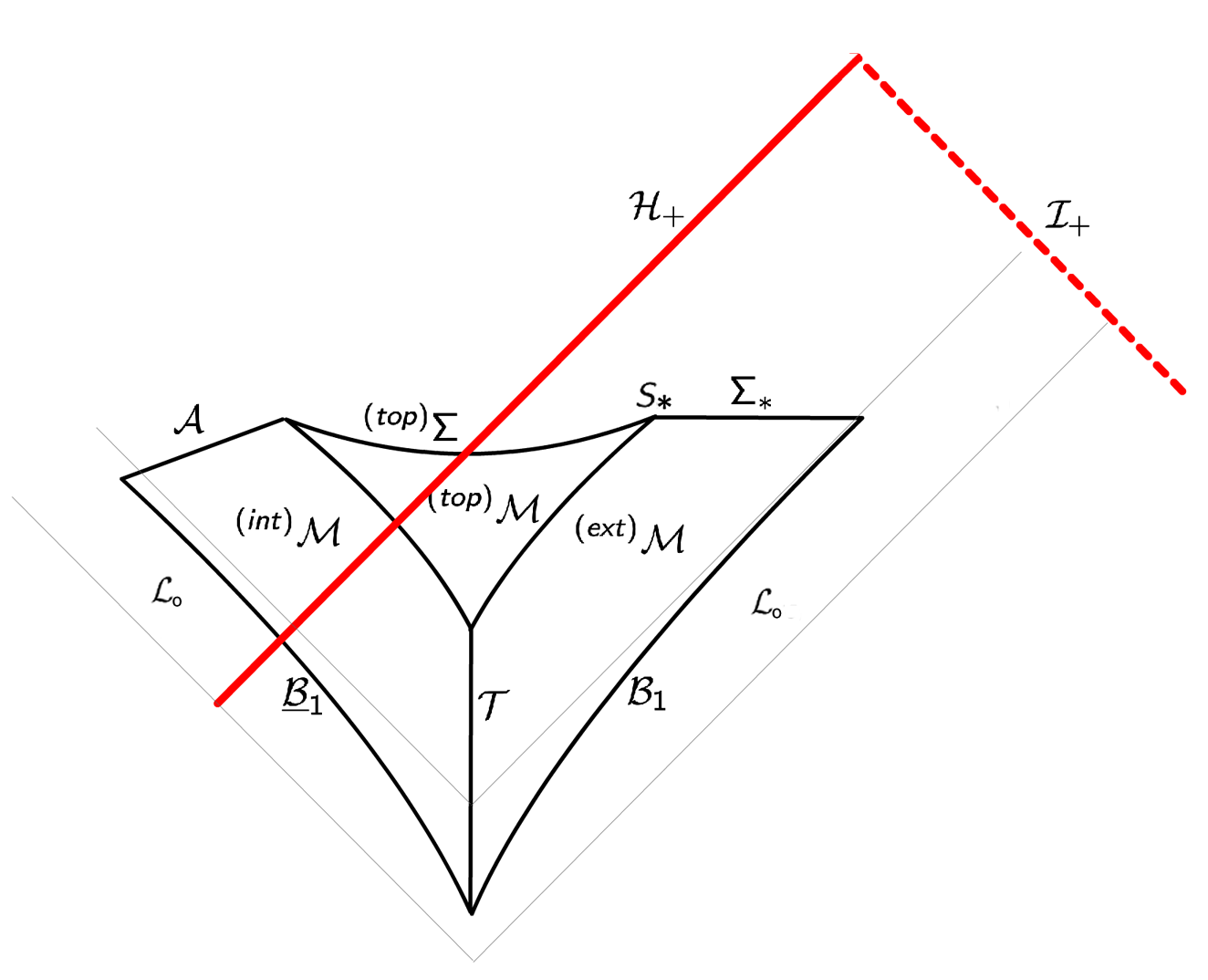}
\caption{The GCM admissible space-time $\mathcal{M}$}
\lab{fig1-introd}
\end{figure} 

Here is a short summary of where these issues are dealt with in our work:

  \begin{itemize}
 \item Papers \cite{KS-GCM1},  \cite{KS-GCM2}  provide a   framework for  
 dealing  with  the issue  (4), by  constructing general covariant modulated (GCM)  spheres, generalizing those   used in the nonlinear stability of  Schwarzschild  in the polarized case in \cite{KS}, in the asymptotic region of a general perturbation of Kerr.  Paper     \cite{KS-GCM2}   also  contains   a definition of  the   angular momentum  for GCM spheres. Relying on these GCM spheres, spacelike GCM hypersurfaces are constructed in   \cite{Shen}, generalizing the construction of GCM hypersurfaces in polarized symmetry of \cite{KS}  to the  case of general perturbations of Kerr. These results are applied in \cite{KS:Kerr}  to the construction of the crucial spacelike GCM boundary $\Si_*$, see section 8.4 and 8.5 in that paper.

   \item  Part I  of this paper     deals with  issue (1) by   developing a geometric formalism  of non-integrable  horizontal structures,  well adapted to perturbations of Kerr, and  use    it  to derive  the   generalized Regge-Wheeler (gRW)  equations   in the context of  general   perturbations of Kerr\footnote{In the linear case,  recall that  complex scalar versions  of such equations were   derived  in  \cite{Ma}, see also \cite{D-H-R-Kerr},  based on an extension of the  physical space Chandrasekhar type transformation   introduced  in \cite{Chand2}.}.    
      
    \item Paper \cite{KS:Kerr}     contains  a precise version of  Theorem  \ref{MainThm-firstversion},  definition of  the main objects  and  a roadmap for the entire proof. The proof is split into 9 distinct steps,   Theorems M0-M8,    and full proofs are given  to  all but  three  of them. Theorems M1,  M2 and half of  M8 were  carefully stated   and their proofs were assumed as a blackbox and  delayed to the   present paper.
 
 \item  Part II of the present paper    deals  with issue (2). It provides complete proofs of Theorems M1 and M2 mentioned above,   by deriving estimates for  gRW using  an     extension of the classical  vectorfield  method\footnote{The results  on decay in  \cite{Ma} and  \cite{D-H-R-Kerr}, on the other hand, depend heavily on  mode decompositions for the linearized  gRW equations  in Kerr, an approach whose
    generalization to the full nonlinear setting  seems to present substantial difficulties.},   based on commutation with   a nonlinear version of the Carter operator.  In the context of the standard scalar wave equation in Kerr,  such an approach was developed by Andersson and Blue  in their important  paper  \cite{A-B}.       
    
    \item The nonlinear terms  present  in the  full version of  the  gRW  equation derived in \cite{GKS1}, as well as those generated 
     by commutation   with vectorfields and  second order Carter operator, are treated  in a similar spirit\footnote{We note however that Theorem M2 is  treated different  here, see  discussion   below.}   as the treatment of the nonlinear terms in  \cite{KS},  by   showing that they verify  a favorable null type structure. 

\item  Part III of this paper gives a  full derivation of the  estimates  for the top  derivatives 
 of the  curvature components, i.e.   a proof of the second part of Theorem M8 stated without proof in   Theorem 9.4.15 in  \cite{KS:Kerr}.
\end{itemize}

In the  remaining part of this introduction, we  give a short presentation of  the main  ingredients of this paper and refer to the introduction in \cite{KS:Kerr} for a presentation of the other steps in the proof of Theorem  \ref{MainThm-firstversion}, i.e Theorems M0, M3--M7, and the first part of M8.

 
\section{Geometric set-up}
\lab{section:Geometricsetup-intro}
  

 
 \subsection{Spacetime $\MM$}
 \lab{sec:spacetimeMM:intro}
  

The geometric setting of this paper consists of an Einstein vacuum Lorentzian manifold 
$(\MM, \g)$  with boundaries  equipped with the following:
\begin{enumerate}
\item A regular horizontal structure defined by a null pair $(e_3, e_4)$, and the  space $\HH$ orthogonal to it. Note that the horizontal structure considered here is not integrable\footnote{In other words, the space $\HH$ forms a non integrable distribution. The formalism  was  originally   mentioned  in  \cite{IK} and   developed  in \cite{GKS1}.}.  The  formalism of non-integrable  horizontal structures, on which of our entire work is based, is developed  in full in Chapter 2. Chapter 3 contains a description of the standard   non-integrable  horizontal structure of Kerr.

\item Two constants $(a, m)$ with $|a|<m$, two scalar functions $(r, \th)$  and a time function $\tau$ on $\MM$. In addition, $\MM$ possesses a horizontal complex 1-form\footnote{By this, we mean $\Jk=j+i\dual j$ where $j$ is a real horizontal 1-form. In Kerr this quantity is specifically introduced  in  Definition \ref{def:JkandJ}.} $\Jk$,   used  to linearize  all horizontal  $1$-forms in perturbations of Kerr.

\item Boundaries given by $\pr\MM = \AA\cup\Si(\tau_*)\cup\Si_*\cup\Si(1)$
where
\begin{itemize}
\item $\AA$ is the spacelike hypersurface given by 
\beaa
\AA &:=& \MM\cap\{r=r_+(1-\deh)\}, \qquad r_+:=m+\sqrt{m^2-a^2},
\eeaa
where $\deh>0$ a sufficiently small constant.

\item $\Si(1)$ and $\Si(\tau_*)$ denote the spacelike level hypersurfaces $\tau=1$ and $\tau=\tau_*$, with  $\tau_*>1$ and $1\leq\tau\leq\tau_*$ on $\MM$.

\item $\Si_*$ is a uniformly spacelike hypersurface connecting $\Si(1)$ to $\Si(\tau_*)$.  
\end{itemize}

\item Two spacetime regions $\Mint$ and $\Mext$ such that 
\beaa
\MM=\Mint\cup\Mext, \qquad \Mext=\MM_{r\geq r_0}, \qquad \Mint=\MM_{r\leq r_0},
\eeaa
where $r_0\gg m$ is a sufficiently large constant.
\end{enumerate}

 \begin{remark}
Note that the spacetime $\MM$ considered above does not require any specific gauge conditions. 
Indeed, in this paper, we only provide gauge independent curvature estimates. The control of Ricci coefficients is provided in \cite{KS:Kerr} where specific gauge choices are made, see section  2.3 and  2.8  for the definitions of PG and PT  structures   in \cite{KS:Kerr}.  We also note that the scalar functions $r, \th$ and $\tau$ are  not   aligned with the frame, i.e. unlike in the  stability of Minkowski space, in \cite{Ch-Kl},  and all  other subsequent works\footnote{We note however  that  in the treatment of  the Regge Wheeler equation in       Chapter 10 of \cite{KS} the  foliations used are  also not aligned with the frame.},  our frames are  in no way  adapted  to foliations.
 \end{remark}
 
 The function $\tau$  is used to define the regions of integrations
 $\MM(\tau_1,\tau_2)$  where  $\tau_1\le \tau\le \tau_2$. We also define the following  significant regions of  $\MM$, see Definition \ref{def:causalregions}.  
     \begin{definition}\lab{def:causalregions-intro}
    We define the following regions of $\MM$:
    \begin{enumerate}
      \item   We define the trapping region of $\MM$ to be the set 
    \beaa
     \MM_{trap}:=\MM\cap\left\{     \frac{|\TT|}{r^3}  \le \de_{trap} \right\}, \qquad \delta_{trap} = \frac{1}{10},
    \eeaa
    where  $\TT= \TT=r^3-3mr^2+ a^2r+ma^2$.    This is the region that contains all trapped null  geodesics, for  
     sufficiently small $a/m$.
  
  \item We denote $\Mntrap$ the complement to the trapping region $\MM_{trap}$. 
 
 \item We denote $\MM_{red}:=\MM\cap\big\{r\leq r_+(1+2\de_{red})\big\}$,  for a sufficiently small  constant $\de_{red}>0$,  the  region  where the red shift effect  of the horizon is  manifest.
          \end{enumerate}
     \end{definition}

 
\subsection{Ricci and curvature coefficients}
  

 
\subsubsection{Definition of the Ricci and curvature coefficients}
  

We can  define, with respect to the horizontal structure associated to $(e_3, e_4)$, connection  and curvature coefficients   similar to   those  in the integrable case, as  in  \cite{Ch-Kl}, 
  \beaa
 \begin{split}
\chib_{ab}&=\g(\D_ae_3, e_b),\qquad \,\,\chi_{ab}=\g(\D_ae_4, e_b),\qquad\,\,\,\,
\xib_a=\frac 1 2 \g(\D_3 e_3 , e_a),\qquad \xi_a=\frac 1 2 \g(\D_4 e_4, e_a),\\
\omb&=\frac 1 4 \g(\D_3e_3 , e_4),\qquad\,\, \om=\frac 1 4 \g(\D_4 e_4, e_3),\qquad 
\etab_a=\frac 1 2 \g(\D_4 e_3, e_a),\quad\,\,\,\,\,  \eta_a=\frac 1 2 \g(\D_3 e_4, e_a),\\
 \ze_a&=\frac 1 2 \g(\D_{a}e_4,  e_3),
 \end{split}
\eeaa
\beaa
\a_{ab}=\R_{a4b4},\quad \b_a=\frac 12 \R_{a434}, \quad   \bb_a=\frac 1 2 \R_{a334},  \quad \aa_{ab}=\R_{a3b3},\quad \rho=\frac 1 4 \R_{3434} , \quad \rhod=\frac 1 4\dual \R_{3434},
\eeaa
and derive  the corresponding   null structure and  null Bianchi equations, see  Propositions 
\ref{prop-nullstr} and  \ref{prop:bianchi}. 
The non-symmetric $2$ tensors  $\chi, \chib$ are decomposed as follows.
\beaa
\chi_{ab}&=&\chih_{ab} +\frac 1 2 \de_{ab} \trch+\frac 1 2 \in_{ab}\atrch,\qquad 
\chib_{ab}=\chibh_{ab} +\frac 1 2 \de_{ab} \trchb+\frac 1 2 \in_{ab}\atrchb,
\eeaa
where   the scalars $\trch$, $\trchb$ and $\atrch$, $\atrchb$ are given by
\beaa
\trch:=\de^{ab}\chi_{ab}, \qquad \trchb:=\de^{ab}\chib_{ab}, \qquad \atrch:=\in^{ab}\chi_{ab}, \qquad \atrchb:=\in^{ab}\chib_{ab}.
\eeaa

\begin{remark}\lab{rem:nonintegrabilityandatrchatrchb}
The non integrability of $(e_3, e_4)$  corresponds to the non vanishing $\atrch$ and $\atrchb$. A  well known example of  a non integrable null frame,   is the principal null frame of Kerr for which $\atrch$ and $\atrchb$ are indeed non trivial,
 see  section \ref{section:values-Kerr}.
\end{remark}


\subsection{Frame transformations}


To start with, given an arbitrary perturbation of Kerr,  there is no reason to prefer an horizontal structure  to any other one.  It is thus   essential  that we  consider  all possible frame transformations from  one   horizontal structure  $(e_4, e_3, \HH)$ to another one  $(e_4', e_3',  \HH')$  together with the transformation  formulas $\Ga\to \Ga'$, $R\to R'$  they generate for the Ricci and curvature coefficients. The most general  transformation formulas  between two null frames is given in Lemma  2.2.1  in \cite{KS:Kerr}. It depends on two  horizontal  $1$-forms  $f, \fb$ and a  real scalar function $\la$ and is given by
 \bea
 \lab{General-frametransformation-intro}
 \bsplit
  e_4'&=\la\left(e_4 + f^b  e_b +\frac 1 4 |f|^2  e_3\right),\\
  e_a'&= \left(\de_a^b +\frac{1}{2}\fb_af^b\right) e_b +\frac 1 2  \fb_a  e_4 +\left(\frac 1 2 f_a +\frac{1}{8}|f|^2\fb_a\right)   e_3,\\
 e_3'&=\la^{-1}\left( \left(1+\frac{1}{2}f\c\fb  +\frac{1}{16} |f|^2  |\fb|^2\right) e_3 + \left(\fb^b+\frac 1 4 |\fb|^2f^b\right) e_b  + \frac 1 4 |\fb|^2 e_4 \right).
 \end{split}
 \eea
The  very important   transformation formulas  $\Ga\to \Ga'$, $R\to R'$ are given in Proposition 2.2.3 of \cite{KS:Kerr}.


\subsection{Basic  equations and complexification}  


The  null structure and null Bianchi equations   verified by the  Ricci  and curvature coefficients    are derived in   sections 2.2.    These  equations    simplify considerably, see section 2.4,  by introducing  complex notations:
\beaa
&& A:=\a+i\dual\a, \quad B:=\b+i\dual\b, \quad P:=\rho+i\dual\rho,\quad \Bb:=\bb+i\dual\bb, \quad \Ab:=\aa+i\dual\aa,
\\
&& X:=\chi+i\dual\chi, \quad \Xb:=\chib+i\dual\chib, \quad H:=\eta+i\dual \eta, \quad \Hb:=\etab+i\dual \etab, \quad Z:=\ze+i\dual\ze, \\
&&\Xi:= \xi + i \dual \xi, \quad \Xib:= \xib+ i \dual \xib,
\eeaa    
where $\dual$ denotes the Hodge  dual as defined in Definition \ref{definition-hodge-duals}. 
In particular, note that  $\tr X = \trch-i\atrch, \,   \tr\Xb = \trchb -i\atrchb$, while $\Xh$ and $\Xbh$ denote the symmetric traceless part of $X$ and $\Xb$ respectively.  Further  useful simplifications of the equations  can be obtained   with the help of conformally invariant derivative operators introduced in  section 2.2.9.


\subsection{Kerr and $O(\ep)$-perturbations  of Kerr}


The preferred (principal) null pair of Kerr is given in Boyer Lindquist coordinates by\footnote{There is an indeterminacy in the principal null pair as one may replace  $(e_3, e_4)$ with $(\la^{-1}e_3, \la e_4)$ for any $\la>0$. The formulas provided here correspond to  a  choice of $\la>0$ ensuring $\D_3e_3=0$ and thus $\omb=0$, $\xib=0$, which is regular across the horizon towards the future. This is called the  canonical incoming null frame of Kerr.}
\beaa
\bsplit
e_4&=\frac{r^2+a^2}{|q|^2} \pr_t +\frac{\De}{|q|^2} \pr_r +\frac{a}{|q|^2} \pr_\phi, \qquad
e_3  =\frac{r^2+a^2}{\De} \pr_t -\pr_r +\frac{a}{\De} \pr_\phi.
\end{split}
\eeaa 
The horizontal structure associated to  this null pair is  spanned  by the vectorfields
\bea
\lab{eq:canonicalHorizBasisKerr-intro}
e_1=\frac{1}{|q|}\pr_\th,\qquad e_2=\frac{a\sin\th}{|q|}\pr_t+\frac{1}{|q|\sin\th}\pr_\phi.
\eea

 The  complexified  Ricci and curvature coefficients   take a particularly simple form  in Kerr,  relative to the above  principal null pair, see  Section \ref{SECTION:KERR}, 
 \beaa
   \Xh=\Xbh=\Xi=\Xib=\omb=0,\qquad A=B=\Bb=\Ab=0,
 \eeaa
and
\beaa
\tr X=\frac{2\De\ov{q}}{|q|^4}, \qquad \tr\Xb=-\frac{2}{\ov{q}}, \qquad P=-\frac{2m}{q^3},\\
\Hb = -\frac{a\ov{q}}{|q|^2}\Jk, \qquad H=\frac{aq}{|q|^2}\Jk,\qquad Z=\frac{aq}{|q|^2}\Jk, 
\eeaa
where $q= r+i a \cos \th$ and $\De= r^2+a^2-2 m r$, relative to  the Boyer-Lindquist (BL) coordinates  $(r, \th)$, and where the horizontal complex 1-form $\Jk$ is given, relative to the  horizontal  basis 
 $(e_1, e_2)$, by the formula 
\bea
\lab{canonical:Jk:intro}
\Jk_1=\frac{i\sin\th}{|q|}, \qquad \Jk_2=\frac{\sin\th}{|q|}.
\eea
Relative to   the ingoing  frame $(e_3, e_4) $ and complex $1$ form  $\Jk$ the Killing vectorfields  $\T$ and $\Z$ take the form 
\bea
\lab{eq:TZinKerr-intro}
\bsplit
\T &= \frac{1}{2}\left(e_4+\frac{\Delta}{|q|^2}e_3 -2a\Re(\Jk)^be_b\right),\\
\Z &= \frac 1 2 \left(2(r^2+a^2)\Re(\Jk)^be_b -a(\sin\th)^2 e_4 -\frac{a(\sin\th)^2\De}{ |q|^2} e_3\right).
\end{split}
\eea
The same formula  is used to define  approximate Killing vectorfields  in perturbations of Kerr, see Definition \ref{Definition:vfsTZ}. Note that $\T$ becomes  spacelike   in the ergoregion $|q|^2 <2 mr$.

A spacetime $\MM$, endowed with an horizontal structure $(e_3, e_4, \HH)$  is said to be  an   $O(\ep)$  perturbation of Kerr   if  all quantities  which vanish  in Kerr are  $O(\ep)$, and  if all other quantities  stay bounded in an  $O(\ep)$ neighborhood of their corresponding\footnote{To make this precise, we also need a definition of     functions $(r, \th)$ and of a complex 1-form $\Jk$, see section \ref{sec:deflinearizedquantities:intro}.}  Kerr  values.    The definition is, of course,  ambiguous in the sense that   any   other  horizontal structure $(e'_3, e'_4, \HH')$ connected to  $(e_3, e_4, \HH)$ by the  frame transformation \eqref{General-frametransformation-intro} with $f, \fb=O(\ep)$ and $\la=1+O(\ep)$   is also an $O(\ep)$-perturbation of Kerr.  Nevertheless  the definition is useful   in that  it brings to light  the remarkable  fact that the extreme curvature components  are in fact $O(\ep^2) $ invariant. This can be easily  seen from the transformation formulas 
\beaa
\bsplit
\la^{-2} \a'&=\a +  \big(  f\hot \b  -\dual f \hot \dual  \b )+ \left( f\hot f-\frac 1 2  \dual f \hot   \dual f \right) \rho
+  \frac 3  2 \big(  f \hot  \dual  f\big) \rhod +O(\ep^3),
\\
\la^2\aa'&=\aa + \big(  \fb \hot \bb  -\dual \fb \hot \dual  \bb )+ \big( \fb \hot \fb-\frac 1 2  \dual \fb \hot   \dual \fb \big) \rho
+  \frac 3  2 \big(  \fb \hot  \dual  \fb\big) \rhod +O(\ep^3),
\end{split}
\eeaa
see Proposition 2.2.3 of \cite{KS:Kerr}. 

\begin{remark}
It is this fact that allows us to treat    $\a, \aa$   differently  from all   other  
quantities,   by choosing frames  best  adapted   to  their analysis. 
\end{remark}


\subsection{Linearization of the Ricci and curvature coefficients}
\lab{sec:deflinearizedquantities:intro}



\subsubsection{Definition of linearized quantities}


Since the quantities $\Xh$, $\Xbh$, $\Xi$, $\Xib$, $\omb$, $A$, $B$, $\Bb$, $\Ab$ all vanish in Kerr, it suffices to linearize the remaining quantities, i.e. $\tr X$, $\tr\Xb$, $\om$, $H$, $\Hb$, $Z$ and $P$. 
The linearization is given  by subtracting the Kerr values as follows, see section \ref{sec:definitionoflinearizedquantities:chap4}: 
\beaa
\bsplit
\trXc &:= \tr X-\frac{2\ov{q}\De}{|q|^4}, \qquad     \trXbc := \tr\Xb+\frac{2}{\ov{q}},\qquad \Pc:= P+\frac{2m}{q^3},\qquad \omc  := \om  + \frac{1}{2}\pr_r\left(\frac{\De}{|q|^2} \right),\\
\Hc &:= H-\frac{aq}{|q|^2}\Jk, \qquad \Hbc:=\Hb+\frac{a\ov{q}}{|q|^2}\Jk,\qquad  \Zc := Z-\frac{aq}{|q|^2}\Jk.
\end{split}
\eeaa


\subsubsection{Notation $(\Ga_g, \Ga_b)$ for Ricci coefficients}


We group the linearized Ricci coefficients in two subsets reflecting their expected decay properties, see section \ref{sec:definitionofGabandGagfirsttime}:
\beaa
\Ga_g&:=& \Big\{\trXc, \quad \Xh, \quad \trXbc, \quad \Hbc, \quad \Zc, \quad \omc, \quad \Xi\Big\},\\
\Ga_b&:=& \Big\{\Xbh, \quad \Hc, \quad \omb, \quad \Xib\Big\}. 
\eeaa

\begin{remark}
In fact, $(\Ga_g, \Ga_b)$ also include the linearization of the derivatives of the scalar functions $(r, \cos\th)$, and of the complex horizontal 1-form $\Jk$, see section \ref{sec:definitionofGabandGagfirsttime}.
\end{remark}

The justification for the above decompositions has to do with the expected  decay properties of the linearized  components in perturbations of Kerr, with respect to $\tau$ and $r$. More precisely, see  \eqref{eq:assumptionsonMextforpartII-1}, 
\bea\lab{eq:expectedbehaviorGabGag:chap2-intro}
\bsplit
\big|\dk^{\leq s}\Ga_g|&\les \ep \min\Big\{ r^{-2 }\tau^{-1/2-\dec},  \, r^{-1}\tau^{-1-\dec} \Big\}, \\
\big|\dk^{\leq s}\Ga_b\big| &\les \ep  r^{-1 }\tau^{-1-\dec},
\end{split}
\eea
for a small constant $\dec>0$, where $\dk=\{\nab_3, r\nab_4, r\nab \}$ denotes weighted derivatives, and  $\ep>0$ is   a sufficiently small bootstrap constant.
 We note also  that the curvature  components  $\Ab, r\Bb$  behave  in the same way as  $\Ga_b$,  while $ r( \Pc, B, A) $ behave like $ \Ga_g$. Moreover   $A, B$  get the optimal decay in powers of  $r$, i.e. 
 \beaa
 |A|, |B| \les \ep  r^{-7/2 -\dec}.
\eeaa


\section{Main theorems}


We refer to  section 3.4  of \cite{KS:Kerr} for   a precise  statement     of our  Main Theorem  concerning the stability of  Kerr and to section 3.7 of \cite{KS:Kerr} the main steps  in the proof.  Here  we  concentrate on  a simplified set of assumptions  needed  for the proof of Theorems M1, M2 and  the curvature estimates for Theorem M8.


\subsection{Smallness constants}


The following constants  are  involved in  the statement of Theorems M0-M8,  see section 3.4. in  \cite{KS:Kerr}:
\begin{itemize}
\item The constants $m_0>0$ and $|a_0|\ll m_0$ are the mass and the angular momentum of the Kerr solution relative to which our initial perturbation is measured. 

\item The integer $k_{large}$ which corresponds to the maximum number of derivatives of the solution.

\item The size of the initial data  perturbation  is measured by $\ep_0>0$. 

\item The size of the bootstrap assumption norms are measured by $\ep>0$.

\item $r_0>0$ is tied to  $\Mint\cap\Mext=\{r=r_0\}$. 

\item The constant $\deh$  tied to the  definition of  $\AA=\{r=r_+(1-\deh)\}$.

\item $\dec$ is tied to decay estimates  in $\tau$   for the linearized quantities of section \ref{sec:deflinearizedquantities:intro}.
\end{itemize}

These  constants are chosen such that
\bea\lab{eq:constraintsonthemainsmallconstantsepanddelta}
\bsplit
& 0<\deh,\,\, \dec  \,\, \ll \min\{m_0 -|a_0|, 1\},\\ 
& r_0\gg \max\{m_0,1\},\qquad k_{large}\gg \frac{1}{\dec}.
\end{split}
\eea
Then, $\ep$ and $\ep_0$ are chosen such that
\bea\lab{eq:constraintsonthemainsmallconstantsepanddelta:bis}
0<\ep_0, \ep\ll \min\left\{\dec,   \frac{1}{r_0}, \frac{1}{k_{large}}, m_0-|a_0|, 1\right\},
\eea
\bea\lab{eq:constraintofep0wrta0}
\ep_0, \ep\ll |a_0|\quad \textrm{ in the case }a_0\neq 0,
\eea
and
\bea\lab{eq:constraintbetweenepep*ep0}
\ep=\ep_0^{\frac{2}{3}}.
\eea

Also, we introduce the integer $k_{small}$ which corresponds to the  number of derivatives  for which the solution satisfies decay estimates. It is related to $k_{large}$ by
\bea\lab{eq:choiceksmallmaintheorem}
k_{small}=\left \lfloor\frac 1 2 k_{large}\right \rfloor +1.
\eea

From now on, in the rest of the paper, $\lesssim$ means bounded by a constant depending only on geometric universal constants (such as Sobolev embeddings, elliptic estimates,...) as well as the constants 
$$m_0,\, a_0, \, \deh,\, \dec, \, r_0, \, k_{large},$$
\textit{but not on} $\ep$ and $\ep_0$.


\subsection{Initial data assumptions}
\lab{sec:assumptioninitialdatanorm:intro}


The initial data norm denoted by $\Ik_k$, measures the size of the perturbation from Kerr at $\tau=1$,   for  the top $k$ derivatives of the  curvature tensor\footnote{The definition  used here differs slightly from the one in  Definition  9.4.9  in  \cite{KS:Kerr}, but easily follows from it by a local existence argument.}.
  
\begin{definition}\lab{def:initialdatanorm-intro}
We define the following initial data  norms on $\Si_1$
\bea
\bsplit
\Ik_{k} := &\sup_{S\subset\Si_1 }r^{\frac{5}{2}  +\de_B} \Big(  \big\| \dk^k\, (A, B) \big\|_{L^2(S)}  + \big\| \dk^k\, B\big\|_{L^2(S)}\Big)\\
&+ \sup_{S\subset\Si_1 }\Big(  r^2 \big\| \dk^k\,  \Pc  \big\|_{L^2(S)}+  r  \big\| \dk^k\, \Bb\big\|_{L^2(S)} + \big\| \dk^k\, \Ab\big\|_{L^2(S)}\Big).
\end{split}
\eea
\end{definition}

In this paper, we make the following assumption on the control of the initial data norm\footnote{The original assumption on initial data in \cite{KS:Kerr} is stated for $k_{large}+10$ derivatives, see (3.4.7) in that paper, in a given frame of an initial data layer $\LL(a_0, m_0)$. The control in the frames  used in this paper  are obtained in Theorem M0 of section 3.7.1 in \cite{KS:Kerr}, and in Theorem 9.4.12  in  \cite{KS:Kerr} for $k_{large}+7$ derivatives.}
\bea\lab{eq:controlofinitialdataforThM8-intro}
\Ik_{k_{large}+7}\leq \ep_0.
\eea
The bound \eqref{eq:controlofinitialdataforThM8-intro}   will be used   both in Part II and  Part III as   assumptions on the initial data.


\subsection{Quantitative assumptions on the spacetime $\MM$}
\lab{sec:bootstrapassumptions:intro}


The quantitative assumptions  made in this article depend on a large positive integer  $k_L$,  representing the maximal number of derivatives  for the linearized  Ricci and curvature coefficients  $(\Gac, \Rc)$ which are  required in the proof. There are in fact two types of assumptions: 
\begin{enumerate}
\item For the proof of Theorem M1 and M2 of \cite{KS:Kerr}, we rely on the following pointwise quantitative assumptions on $\Ga_b$ and $\Ga_g$, for $k\leq \kl$,  
\bea\lab{eq:assumptionsonMextforpartII-1-intro}
\bsplit
\Big(r^2\tau^{\frac{1}{2}+\dec}+r\tau^{1+\dec}\Big)|\dk^{\leq k}\Ga_g| & \leq\ep, \\ 
r\tau^{1+\dec}|\dk^{\leq k}\Ga_b| & \leq\ep.
\end{split}
\eea

\item For the proof of the curvature estimates of Theorem M8 of \cite{KS:Kerr}, we introduce weighted energy-Morawetz type norms for curvature and Ricci coefficients, denoted respectively by $\Rk_k$ and $\Sk_k$, see section \ref{subsection:MainNormsM8} for the precise definition. We then rely on 
the following quantitative assumptions on $\Rk_k$ and $\Sk_k$  
\bea\lab{eq:assumptionsonMextforpartIII-intro}
\Rk_k+\Sk_k & \leq\ep,\qquad 0\leq k\leq \kl,
\eea
as well as the following pointwise quantitative assumptions on $\Ga_b$ and $\Ga_g$
\bea\lab{eq:assumptionsonMextforpartIII-bis-intro}
r^2|\dk^k\Ga_g|+r|\dk^k\Ga_b| & \leq \frac{\ep}{\tau_{trap}^{1+\dec}}, \qquad 0\leq k\leq \frac{\kl}{2},
\eea
where the scalar function $\tau_{trap}$ defined by
\beaa
\tau_{trap} := \left\{\ba{lll}
1+\tau & \textrm{on} & \MM_{trap},\\
1& \textrm{on} & \Mntrap.
\ea\right.
\eeaa
\end{enumerate}

The integer $\kl$ is chosen as follows:
\begin{itemize}
\item For the proof of Theorem M1 and M2 of \cite{KS:Kerr} (restated in Theorem \ref{theoremM1:intro} and \ref{theoremM2:intro} below), we choose $\kl=k_{small}+120$. Then, \eqref{eq:assumptionsonMextforpartII-1-intro} follows by interpolation from the bootstrap assumptions (3.5.1) (3.5.2) in \cite{KS:Kerr} together with the construction of the global frame in section 3.6.3 of \cite{KS:Kerr}, where (3.5.1) in \cite{KS:Kerr} are bootstrap assumptions on boundedness for $k\leq k_{large}$ derivatives, and  (3.5.2) in \cite{KS:Kerr} are bootstrap assumptions on decay for $k\leq k_{small}$ derivatives.  

\item For the proof of the curvature estimates of Theorem M8 (see Theorem \ref{THEOREMM8:INTRO} below), we choose $\kl=k_{large}+7$. Then, \eqref{eq:assumptionsonMextforpartIII-intro} follows from the bootstrap assumptions  (9.4.20) of \cite{KS:Kerr} together with the construction of the global frame in section 9.4 of \cite{KS:Kerr}. Also, \eqref{eq:assumptionsonMextforpartIII-bis-intro} is a non sharp consequence of the bootstrap assumptions (9.4.22) in \cite{KS:Kerr} together with the construction of the global frame in section 9.4 of \cite{KS:Kerr}. 
\end{itemize}


\subsection{Statement of the main theorems}


Recall that the  nonlinear stability of the Kerr family for small angular momentum, i.e $|a|/m\ll 1$, is stated in the Main   Theorem in section 3.4 of \cite{KS:Kerr}. The proof is divided in a sequence of nine intermediary steps, called Theorem M0--M8, see section 3.7 in \cite{KS:Kerr}. The goal of the present paper is to provide the proof of Theorems M1 and M2 as well the curvature  estimates of Theorem M8,   which were  stated without proof  in Theorem 9.4.15 of \cite{KS:Kerr} and all involve curvature estimates of hyperbolic type.


\subsubsection{Theorems M1 and M2}


In what follows, we restate\footnote{A more precise statement is given in Theorems \ref{theoremM1:Chap11} and \ref{thm:restatementofTheoremM2}.} Theorem M1 and M2, see section 3.7.1 in \cite{KS:Kerr}.
\begin{theorem}[Theorem M1 in \cite{KS:Kerr}]
\lab{theoremM1:intro}
Assume that the spacetime $\MM$  as defined in section \ref{sec:spacetimeMM:intro}  
verifies the quantitative  assumptions \eqref{eq:assumptionsonMextforpartII-1-intro}, and the assumption \eqref{eq:controlofinitialdataforThM8-intro}  on initial data. Then, if $\ep_0>0$ is sufficiently small, there exists $\dee>\dec$ such that we have the following   estimates in  $\MM$, for all  $ k\le \kl -20$, 
\beaa
\sup_{\MM}\left(\frac{r^2(2r+\tau)^{1+\dee}}{\log(1+\tau)}+r^3(2r+\tau)^{\frac{1}{2}+\dee}\right)\Big(|\dk^kA|+r|\dk^{k-1}\nab_3A|\Big) \les \ep_0.
\eeaa
Also, the quantity $\qf$ introduced below, see section \ref{sec:TeukandRWinourapproach:intro}, satisfies, for all  $ k\le \kl -20$, 
\beaa
  \int_{\Si_*(\ge \tau)}|\nab_3 \dk^{k-1}\qf|^2 &\les& \ep^2_0 \tau^{-2-2\de_{extra}}.
\eeaa
\end{theorem}

\begin{theorem}[Theorem M2 in \cite{KS:Kerr}]
\lab{theoremM2:intro}
In addition to the assumptions of Theorem \ref{theoremM1:intro}, we make the following assumption\footnote{This is the dominant condition of $r$ on $\Si_*$, see (3.4.5) in \cite{KS:Kerr}.} on $\Si_*$
\bea
\lab{eq:behaviorofronS-star}
\min_{\Si_*}r\geq  \de_*  \ep_0^{-1}\tau_*^{1+\dec}
\eea
for some small universal constant $\de_*>0$. Then, we  have the following decay estimates for $\Ab$ along $\Si_*$
\beaa
\max_{0\leq k\leq   \kl-40}\int_{\Sigma_*}\tau^{2+2\dec}|\dk^k\Ab|^2 \les \ep_0^2.
\eeaa
\end{theorem}

Both results are proved  in Part II of this paper.

\medskip


\subsubsection{Curvature estimates in Theorem M8} 


Theorem M8 in \cite{KS:Kerr} is proved through an iteration procedure described   in section 9.4.7 of \cite{KS:Kerr}. The control of the Ricci coefficients have been derived in Chapter 9 of \cite{KS:Kerr}. In the present paper, we derive the remaining estimates for the proof of Theorem M8, i.e the estimates for curvature stated in Theorem 9.4.15 of \cite{KS:Kerr}. To this end, we introduce weighed $L^2$ type norms $\Rk_k$ and $\Sk_k$ respectively for curvature and Ricci coefficients\footnote{As well as derivatives of $(r, \cos\th)$ and $\Jk$.}, and decompose $\Rk_k$ and $\Sk_k$ in their restrictions $\Rkint$, $\Skint$ to $\Mint$ and $\Rkext$, $\Skext$ to $\Mext$, see section \ref{subsection:MainNormsM8} for the precise definition of these norms. In view of the results in Chapter 9 of \cite{KS:Kerr}, the proof of Theorem 8 reduces to the following result on the control of the curvature norm $\Rk_k$.
\begin{theorem}[Theorem 9.4.15 of \cite{KS:Kerr}]
\lab{THEOREMM8:INTRO}
Assume that the spacetime $\MM$  as defined in section \ref{sec:spacetimeMM:intro}  
verifies the quantitative assumptions \eqref{eq:assumptionsonMextforpartIII-intro} \eqref{eq:assumptionsonMextforpartIII-bis-intro} for $\kl=k_{large}+7$, and the assumption \eqref{eq:controlofinitialdataforThM8-intro}  on initial data. Let  $k_{small}-1\leq J\leq k_{large}+6$. Then,  we  have the following boundedness estimates for all components of curvature
\bea
\bsplit
    \nn\Rkint_{J+1}^2\les& r_0^{18}\Big( \ep_J(\Sk_{J+1}+\Rk_{J+1}) +\ep_J^2+\ep_0^2\Big)+|a|r_0^3\Sk^2_{J+1}\\
& +r_0^{\frac{27}{4}}\Sk_{J+1}^{\frac{3}{2}}\Big(\ep_0+\sqrt{\ep_J}\sqrt{\Sk_{J+1}+\Rk_{J+1}}\Big)^{\frac{1}{2}},\\
\Rkext_{J+1}^2 \les& r_0^{3+\de_B}\Rkint^2_{J+1}+ r_0^{-\de_B} \Skext^2_{J+1} +\ep_0^2,
\end{split}
\eea
where the constant in $\les$ is independent of $r_0$ and $\ep_J$  is such that $\Sk_J+\Rk_J\leq \ep_J$.
\end{theorem}

Part III of this paper is  entirely dedicated to the proof of  Theorem \ref{THEOREMM8:INTRO}.

 
\section{Derivation  and estimates for the gRW  equations}



\subsection{Teukolsky  and gRW  equations in  our approach}
\lab{sec:TeukandRWinourapproach:intro}


In section \ref{sec:spacetimeMM:intro}  we derive, using the formalism
 developed  in  the previous sections\footnote{This follows from the complex form of the null Bianchi identities, see Proposition \ref{prop:bianchi:complex}.},  the nonlinear version of  the Teukolsky  equations for $A$ and $\Ab$ of the form
 \bea
 \lab{eq:Teukolsky-forA-intro2}
 \LL[A]=\err[\LL[A]],\qquad \und{ \LL}[\Ab]=\err[[\und{\LL}[\Ab]],
 \eea
where $\LL,\und{ \LL} $ are  second order tensorial wave operators on our spacetime $\MM$,  
  and   where  $\err[\LL[A]]$, $\err[\LLb[\Ab]]$  are   nonlinear errors   depending on all linearized Ricci  and  curvature coefficients.
  
 Just as in linear theory, to be able to control  $A, \Ab$ we need to perform   transformations $\qf=\qf[A]$,   $\qfb=\qfb[\Ab]$,   which    take solutions $A, \Ab$ of the Teukolsky equation 
\eqref{eq:Teukolsky-forA-intro2} into solutions of  nonlinear, tensorial, versions  of Regge-Wheeler   equations,  which we   call  gRW equations.

In the  setting of polarized perturbations of Schwarzschild \cite{KS},         the  derivation of the RW  equation for\footnote{Note that  \cite{KS}    did   not rely on $\qfb$.} $\qf$   was performed  using null   frames, which had the feature   to be  both   adapted to 
  an integrable   foliation and   diagonalize the curvature tensor up to error terms.     One could thus  rely on  the geometric  formalism developed 
 in the context of the proof of the nonlinear  stability of Minkowski space \cite{Ch-Kl}. In the present paper, we rely on an extension  of the formalism of \cite{Ch-Kl} which allows for non integrable null frames and is presented in Chapter \ref{CHAPTER-NON-INTEGRABLE-STRUCTURES}. Our results on the derivation of gRW in perturbations of Kerr are obtained  in  Chapter \ref{CHAPTER-DERIVATION-MAIN-EQS} and can be summarized as follows.
         \begin{theorem}
          There exist complex  $2$ tensors  $\qf, \qfb\in \sk_2(\CCC) $ derived  from $A, \Ab$   as follows,
          \bea
          \lab{def:qf,qfb:intro}
          \bsplit
          \qf&=& q \ov{q}^{3} \left( \nabc_3\nabc_3 A + C_1  \nabc_3A + C_2   A\right),\\
          \qfb&=& \ov{q} q ^{3} \left( \nabc_4\nabc_4 A + \und{C}_1  \nabc_3A + \und{C}_2   A\right),
          \end{split}
           \eea
where $q=r+ i a \cos\th$  $\nabc_3$, $\nabc_4$ are conformal derivatives, see section \ref{section:conf.inv.Derivatives}, and  
\bea\lab{eq:C1-C2-und{C_1},und{C_2}}
\begin{split}
C_1&=2\trchb - 2\frac {\atrchb^2}{ \trchb}  -4 i \atrchb, \\
C_2  &= \frac 1 2 \trchb^2- 4\atrchb^2+\frac 3 2 \frac{\atrchb^4}{\trchb^2} +  i \left(-2\trchb\atrchb +4\frac{\atrchb^3}{\trchb}\right),\\
\und{C}_1&=2\trch - 2\frac {\atrch^2}{ \trch}  -4 i \atrch, \\
\und{C}_2  &= \frac 1 2 \trch^2- 4\atrch^2+\frac 3 2 \frac{\atrch^4}{\trch^2} +  i \left(-2\trch\atrch +4\frac{\atrch^3}{\trch}\right),
\end{split}
\eea
 which verify gRW equations of the form\footnote{Here $\squared_2$ is the covariant wave operator for horizontal $2$-tensors, see   section \ref{section:Horizwaveoperators}.} 
 \bea
 \lab{wave-equation-qf-intro}
 \bsplit
 \squared_2 \qf   -i \frac{4 a\cos\th}{|q|^2} \nab_\T \qf   - V  \qf &=  L_{\qf}[A] + \err[\squared_2 \qf],\\
  \squared_2 \qfb   +i \frac{4 a\cos\th}{|q|^2} \nab_\T \qfb   - \und{V}  \qfb &= L_{\qfb}[\Ab] + \err[\squared_2 \qfb],
  \end{split}
 \eea
 with $\T$  defined as  in \eqref{eq:TZinKerr-intro}.   The  potentials $V, \und V$ 
  are  real and positive  and the   terms  $L_{\qf}[A], L_{\qfb}[\Ab]$  are  linear in $A$, resp $\Ab$ and have  have   important  specific properties  described  in  detail in subsections \ref{subsection:DerivationgRW-eq} and \ref{subsection:DerivationgRW-eqbar}.
 Finally  the error terms  $\err[\squared_2 \qf], \, \err[\squared_2 \qfb]$
  depending on all linearized Ricci and  curvature coefficients    are acceptable  error terms, i.e.  they  verify  important   structural properties, reminiscent to the null condition. 
    \end{theorem}

 \begin{remark} 
 \lab{remark:coupledGRW-transport} 
 Due to the presence of the linear terms in $A$, resp. $\Ab$,  on the right hand side of \eqref{wave-equation-qf-intro},   one  has to  view  the  wave equations  in 
  \eqref{wave-equation-qf-intro}  as  coupled with the defining equations for $\qf$, $\qfb$  given by \eqref{def:qf,qfb:intro},  that is coupled\footnote{This is different from the case of Schwarzschild, see \cite{KS},   where these equations  decouple.}  with   second order  transport type  equations in $A$, resp. $\Ab$.  
 \end{remark}
 
 \begin{remark} 
 Note that, in the case of Kerr, the corresponding    gRW type  equations  in    \cite{Ma}  are 
  complex scalars  $\psi^{[\pm]} $ verifying the equations\footnote{ With $\square_{a,m}$  the Kerr  D'Alembertian, 
  $c, V$    are  real function of $r, \th$ and $L_{\pm}(\alpha^{[\pm2]}) $ lower order terms. } 
  \bea
\lab{Intro:Regge-Wheeler-DHR}
\square_{a,m} \psi^{[\pm]} + i a  \ c(r, \th) \partial_t \psi^{[\pm]} +V(r, \th) \psi^{[\pm]}  = a   L_{\pm}(\alpha^{[\pm2]}). 
\eea
  These scalars  are connected to our tensorial quantities $\qf, \qfb$   via the relations
$\psi^{[+]}=\qf(e_1, e_1)$,   $\psi^{[-]}=\qfb(e_1, e_1)$. The equations        \eqref{Intro:Regge-Wheeler-DHR}    can be  obtained by projecting our tensorial  equations   
\eqref{wave-equation-qf-intro}.
 Note  however  that the projection modifies the equations  by    the appearance of Christoffel symbols\footnote{Singular on the axis, i.e. at $\th=0, \pi$.} of the horizontal  frame\footnote{See Section \ref{section:projection-gRW} for  a discussion of the projection  and the relation with   equation \eqref{Intro:Regge-Wheeler-DHR}.}.  
\end{remark}


\subsection{RW model equations}    


   The most demanding part  in the analysis of the  gRW equations  \eqref{wave-equation-qf-intro} is to derive  global Energy-Morawetz
    type estimates for $(\qf, A)$ and  respectively $(\qfb, \Ab)$.  To do this,  it  helps to analyze first the reduced equations   in which the right hand side of both equations are treated as sources.    Taking also 
    $\psi=\Re(\qf)$, $\und{\psi}= \Re(\qfb)$ we are led to the real RW model equations
  \bea\lab{eq:Gen.RW-intro}
  \bsplit
\square_2 \psi -V\psi&=- \frac{4 a\cos\th}{|q|^2}\dual \nab_T  \psi+N, \qquad V= \frac{4\De}{ (r^2+a^2) |q|^2},\\
\squared_2 \underline{\psi} -V\underline{\psi}&= \frac{4 a\cos\th}{|q|^2}\dual \nab_T  \underline{\psi}+\und{N}, \qquad \,  V= \frac{4\De}{ (r^2+a^2) |q|^2}.
\end{split}
\eea
   A significant part in the proof of  Theorems \ref{theoremM1:intro}-\ref{theoremM2:intro}   is   to derive the following   result for  $\psi, \undpsi$.

 \begin{theorem}
 \lab{THEOREM:GENRW1-P-INTRO}
       The following estimates hold true  for solutions $\psi, \undpsi \in\sk_2$  of  the wave equations   \eqref{eq:Gen.RW-intro} on spacetime region $\MM(\tau_1, \tau_2)$, for all $\de\le p\le  2-\de$ and  $2\leq s\le \kl$,
       \bea
       \lab{eqtheorem:GenRW1-p-intro}
        \BEF_p^s[\psi](\tau_1, \tau_2)  \les
       E_p^s[\psi](\tau_1)+\NN_p^s[\psi, N](\tau_1, \tau_2),
       \eea
        \bea
       \lab{eqtheorem:GenRW1-p- -psib'-intro}
       \BEF_p^s[\undpsi](\tau_1, \tau_2)   \les
       E_p^s[\undpsi](\tau_1)+\NN_p^s[\undpsi, \und{N}](\tau_1, \tau_2),
       \eea
       where
         \bea
         \lab{norms:BEF-intro}
   \BEF_p^s[\psi](\tau_1, \tau_2):=  \sup_{\tau\in[\tau_1, \tau_2]} 
        E_p^s[\psi](\tau) +B_p^s[\psi](\tau_1, \tau_2) +F_p^s[\psi](\tau_1, \tau_2).
   \eea
       \end{theorem}
       
  The  energy flux norms  $E_p^s[\psi],  F_p^s[\psi]$,  bulk norms  $B_p^s[\psi]$ and  source norms  $\NN_p^s$, 
   with $p$ referring to $r^p$ weights and  $s$ to the number of derivatives, are defined  in  section \ref{subsection:basicnormsforpsi}. 
   For the sake of this introduction it suffices to  take a closer look 
    at the  crucial  bulk terms  $B_p^s$,  which degenerate 
     at the trapped set $\MM_{trap} $, see Definition \ref{def:causalregions-intro}.
    \begin{definition}
    For $0<p<2$ we define, with $\dk = (r\nab_4,  r\nab, \nab_3)$,  
     the bulk norms  $B^s_{p}[\psi](\tau_1, \tau_2) :=\sum_{k\le s} B_p[\dk^k\psi]$
\beaa
B_{p}[\psi](\tau_1, \tau_2)&:=& \Mor[\psi](\tau_1, \tau_2) + \int_{\MM_{r\geq 4m}(\tau_1,\tau_2)}  r^{-1-\de} |\nab_3 \psi|+r^{p-3}\Big(|\dk \psi|^2 +|\psi|^2  \Big),\\
 \Mor[\psi](\tau_1, \tau_2)&:=&\int_{\MM(\tau_1, \tau_2) } 
      r^{-2}  | \nab_\Rhat  \psi|^2 +r^{-3}|\psi|^2+ \int_{\Mntrap(\tau_1, \tau_2)} \left(  r^{-2}|\nab_3\psi|^2 + r^{-1}  |\nab  \psi|^2\right).
\eeaa
    \end{definition}
    
    \begin{figure}[hh!]
\centering
\includegraphics[scale=0.6]{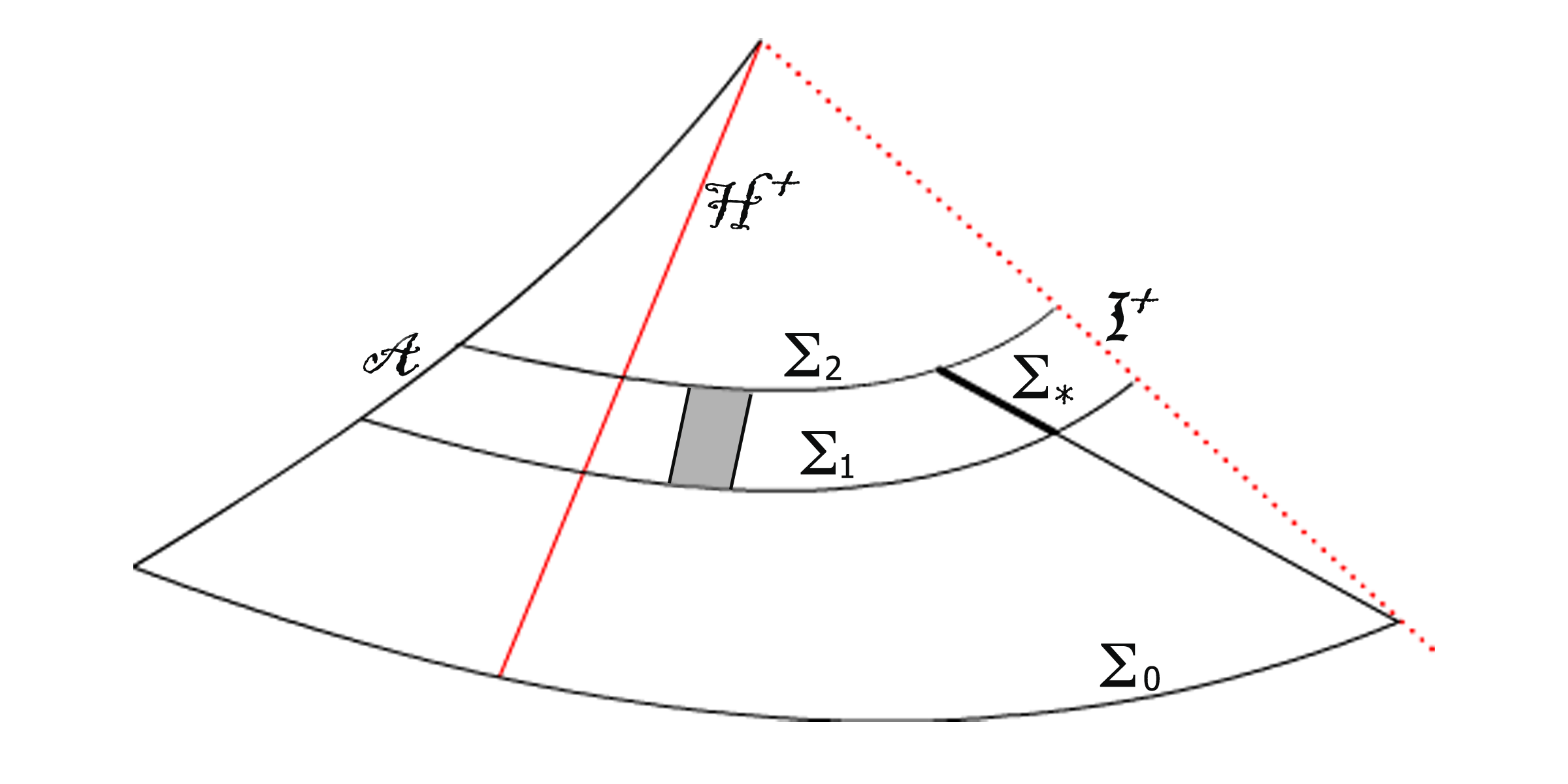}
\caption{The  spacetime region $\MM(\tau_1,\tau_2)=\MM\cap\{\tau_1\le \tau\le \tau_2\} $ between the spacelike hypersurfaces  $\Si_1=\Si(\tau_1)$ and $\Si_2=\Si(\tau_2)$, with the grey region denoting the trapped set.  }
\lab{fig00-introd}
\end{figure}

   The important thing in this definition is that  $B_p[\psi]$ controls 
   the spacetime integrals of  $|\nab_\Rhat \psi|^2$ and $|\psi|^2 $ everywhere  and  all other  derivatives  away from the trapped set.
   
   In addition, we also  derive  estimates for the quantity  $\psiwc := r^2(e_4\psi+\frac{r}{|q|^2}\psi)$  for  which   one can prove  stronger $r^p$ estimates\footnote{These results are the analog in perturbations of Kerr,  to Theorem 5.17 and  Theorem 5.18 of  \cite{KS} for perturbations of Schwarzschild. They are based  on   improved  $r^p$ weighted hierarchy first  introduced in \cite{AnArGa}.}, see Theorem \ref{THEOREM:GENRW2-Q}.

        
\subsection{Main steps in the proof of  Theorems M1 and M2}
\lab{section:estim(qfb,Ab)-intro}

    
   The proof of Theorem   \ref{THEOREM:GENRW1-P-INTRO} is by far the most demanding  part in the proof  of  Theorems \ref{theoremM1:intro}-\ref{theoremM2:intro}.  
     We provide  an introduction to  some of the main 
   ideas  of the proof in  section \ref{section:Mainideas-GenRW-p-intro}.

   To extend Theorem   \ref{THEOREM:GENRW1-P-INTRO}  to  the full  gRW equations \eqref{wave-equation-qf-intro}, 
    we have to  control the terms   $ \NN_p^s[\qf, N](\tau_1, \tau_2)$ with 
    $N= L_{\qf}[A] + \err[\squared_2 \qf]$ and  $ \NN_p^s[\qfb, \und{N}](\tau_1, \tau_2)$  with $\und{N}=   L_{\qfb}[\Ab ] + \err[\squared_2 \qfb]$  for the second equation.  
         This is  done in steps   by first eliminating the linear error terms  on the right hand side of our two gRW equations  and then eliminating the remaining nonlinear quadratic terms.    The procedure for doing this  differs substantially for the two equations.

         
\subsubsection{Estimates for $(\qf, A)$}   


The procedure of eliminating  the  error term  $N= L_{\qf}[A] + \err[\squared_2 \qf]$  requires
         the use  a global frame of $\MM$  for which we have\footnote{We remark that if  $\Hc\in \Ga_b$, we cannot even derive a Morawetz estimate.}
         \bea
         \Hc\in \Ga_g.
         \eea 
         
     {\bf Step 1.}   To  eliminate the      contribution of the    linear term $L_\qf[A]$
      one has to first  derive estimates for  $A$  using  the second  order transport   equations    which  defines  $\qf$.   In the particular   case of Kerr, we   can write, see  Proposition \ref{PROP:FACTORIZATION-QF},
\beaa
\nabc_3\nabc_3\left(\frac{\ov{q}^4}{r^2} A \right) &=&  \frac{\ov{q} }{q } r^{-2} \qf.
\eeaa
The above  factorization\footnote{In perturbations of Kerr, to avoid  the presence  of  unacceptable nonlinear   error terms,   we  use instead the following  modified  factorization, see Lemma \ref{lemma:factorizationofqfusefulfortransportequations},
\beaa
\nabc_3\left(\nabc_3\left(\frac{(\ov{\tr\Xb})^2}{(\Re(\tr\Xb))^2(\tr\Xb)^2}A\right) -\frac{r^2}{2}F A\right) 
&=& O(r^{-2})\qf  +r\dk^{\leq 1} (\Ga_b\c A),
\eeaa
 for an appropriately defined  scalar   function $F$.} is used   to  derive appropriate $\BEF_p^s$   estimates  for $A$  in terms of   $\qf$.   These can then be combined   with   the estimates  derived in Theorem \ref{THEOREM:GENRW1-P-INTRO}   to obtain  estimates for the norms 
 $\BEF_p^s[A, \qf]$  depending only on the  the nonlinear  error  terms, see Theorem \ref{Thm:Nondegenerate-Morawetz}. We note  that the smallness of $a$ and   the  specific   structure\footnote{This is needed  to  control  the corresponding  energy estimates.} of the top  terms  $L_\qf[A]$  is also essential  in  controlling\footnote{This  requires  several integration by parts.}  the terms in $ \NN_p^s[\qf, L_\qf[A] ](\tau_1, \tau_2)$.

 {\bf Step 2.}     To eliminate the  error terms $\err[\squared_2 \qf]$ in the $\qf$ equation we  proceed  as in  \cite{KS}.  As mentioned above, it is essential  that    the analysis is done in a  global frame of $\MM$ for which  $ \Hc\in \Ga_g$. 
 One can  show that   $\err[\squared_2 \qf]$  can be written in the form, see equation \eqref{eq:MaiThmParq-err}, 
    \beaa
   \bsplit
 \err[\squared_2 \qf]&= r^2 \frak{d}^{\leq 3} (\Ga_g \c (A, B))+ \nab_3 (r^2 \frak{d}^{\leq 2}( \Ga_b \c (A, B)))+\frak{d}^{\leq 1} (\Ga_g \c \qf) + r^{3}\dk^{\leq 2} \big( \Ga_b \c \Ga_g \c \Ga_g\big).
\end{split}
 \eeaa  
 
       {\bf Step 3.} Given  the above form of $\err[\squared_2 \qf]$  we are able to derive estimates   for the norms $\BEF^s_p[\qf, A]$ for $\de\leq p\leq 2-\de$. Additional improved $r^p$ estimates are then obtained for the quantity $\qfc=r^2(\nab_4\qf+\frac{r}{|q|^2}\qf)$.

    
\subsubsection{Estimates for   $(\qfb, \Ab)$}  
   
   
   We rely on a different global frame  of $\MM$ for which  we have, in $\Mext$, 
    \bea
    \lab{eq:Framecondfor-qfb}
    \Xi=0, \qquad \Hbc=0.
    \eea
    
    \begin{remark}
    \lab{remark:framecond-qfb-intro}
      We note that the temporal frame of $\Mext$ (see  Definition 9.1.1 in \cite{KS:Kerr}), in which  the  top  derivative estimates  in $\Mext$ for  the linearized Ricci coefficients  were derived,  verifies  these properties.  The conditions \eqref{eq:Framecondfor-qfb} can in fact be relaxed, i.e. $\Xi\in r^{-2} \Ga_g$, $\Hbc\in r^{-1} \Ga_g$  suffice.
    \end{remark}
    
           The structure of the  error terms   $ \err[\squared_2 \qfb]$, in that frame,  turns out to be   more subtle  than  that  of   $ \err[\squared_2 \qf]$,   as   it depends  in an essential   way on the fact that the quantities  $\Ab_4, \Bb_4, P_4$,  introduced\footnote{In the case  of integrable  $S$-foliations  these  notations  were introduced in  \cite{Ch-Kl}, Chapter 7.} in  section \ref{section:RenormBianchiIds}, have  improved decay properties in $r$.
Similar improvements appear   in the null  structure equations, see Proposition \ref{prop-nullstr:complex-conf}. Thus, for example,  the quantities  $\nabc_4\tr\Xb +\frac{1}{2}\tr X\tr\Xb $ and   $\nabc_4\widehat{\Xb} +\frac{1}{2}\tr X\, \widehat{\Xb}$ have better decay properties than respectively  $\nabc_4\tr\Xb $ and  $\nabc_4\widehat{\Xb} $.
 These improvements\footnote{Similar    improvements,  in the particular case of perturbations of Schwarzschild,    are used  to     treat    the   corresponding  error  terms   in the estimates  for the  quantity  corresponding to $\qfb$ in \cite{DHRT}.} carry over  various quadratic and cubic error terms  appearing    in $  \err[\squared_2 \qfb]$.  Using these facts we can show that
    \bea
    \lab{eq:err[squared_2qfb]-intro}
    \err[\squared_2 \qfb]&=& r^2\dk^{\le 3 }\big((A, B)\c \Ab\big)+\dk^{\le 3 }(\Ga_b\c \Ga_g).
    \eea

{\bf Step 1.}  The simplest way to factorize $\qfb$ in Kerr  is given by, see   Proposition \ref{Prop:factorization-qfb}, 
  \beaa
r \nabc_4\left( r^2   \left(   \nabc_4  ( r \frac{q^4}{r^4}  \Ab)\right)\right)= \frac{q }{\ov{q} } \qfb.
\eeaa
The above  factorization\footnote{In perturbations of Kerr, to avoid  the presence  of  unacceptable nonlinear   error terms,   we  use instead the following  modified  factorization, see
  Lemma \ref{lemma:factorizationofqfbusingpowersofqandr},
  \beaa
\qfb &=& \ov{q}q^3\left(\nabc_4+2\tr X -\frac{|\tr X|^2}{2\trch}\right)\left(\nabc_4+2\tr X -\frac{3|\tr X|^2}{2\trch}\right)\Ab +r^2\dk^{\leq 1}(\Ga_g\c\Ga_b).
\eeaa} is used   to  derive appropriate $\BEF_p^s$   estimates  for $\Ab$  in terms of   $\qfb$.   The  contribution to $ \NN_p^s[\qfb, N ](\tau_1, \tau_2)$ due to  $L_{\qfb}[\Ab]$  can then be absorbed  exactly as 
   in the case of $\qf$.
      
        {\bf Step 2.} Given  the form \eqref{eq:err[squared_2qfb]-intro} for $\err[\squared_2 \qfb]$  we can   only derive estimates   for the norms $\BEF^s_p[\qfb, \Ab]$ for $\de\leq p\leq 1-\de$.

    
\subsubsection{Proof of Theorem  \ref{theoremM1:intro}}  


Once  full  $r^p$ weighted estimates  for $(\qf, A)$ are  derived,   the  proof  of the 
estimates of Theorem M1   follow  steps similar  to those used in  Chapter 5  of \cite{KS}, see section \ref{sec:finallyproofofThmM1} for the details.

    
 \subsubsection{Proof of Theorem  \ref{theoremM2:intro}} 


Given that we  only derive estimates   for the norms $\BEF^s_p[\qfb, \Ab]$ for $\de\leq p\leq 1-\de$, we obtain at first insufficient decay estimates in $\tau$ for $(\qfb, \Ab)$. However, we can show that higher $\Lieb_\T$ derivatives of $(\qfb, \Ab)$ decay faster in powers of $\tau$. Using this observation, we obtain suitable decay in $\tau$ for the flux of $(\Lieb_\T^2\qfb, \Lieb_\T^2\Ab)$ along $\Si_*$. We then rely on this result, and on some version of Teukolsky-Starobinsky providing an identity between $\qf$ and $\Ab$ to recover the desired decay estimate for $\Ab$ on $\Si_*$ stated in Theorem M2, see section \ref{sec:finallyproofofThmM2} for the details.

   
\section{Main ideas in the proof of Theorem \ref{THEOREM:GENRW1-P-INTRO}}
\lab{section:Mainideas-GenRW-p-intro}


In the polarized  situation of  \cite{KS} where the wave equation\footnote{We note that in  \cite{KS} the quantity  $\qfb$ was not  actually used.} for $\qf$ decouples, linearly,  from the transport equation for $A$,  the proof of    the analogue of  Theorem   \ref{THEOREM:GENRW1-P-INTRO}
was done  as follows:
\begin{itemize}
\item An adapted \textit{physical space}  approach  based on Morawetz-Energy, red shift  and $r^p$-weighted  estimates which extends   the treatment of the   linearized  RW equation in Schwarzschild in \cite{D-H-R} to the  full nonlinear setting. 

\item Exploit the specific null   structure of the nonlinear error terms. 

\item Once  we  control $\qf$, we then estimate   $A$  using the corresponding second order transport equation.
\end{itemize}

To pass  to the  case of  perturbations of   Kerr,  one encounters  the following  additional problems:  
\begin{enumerate}
\item \textit{Complicated nature of the trapping region.} This is the  region of the domain of outer communication $r> r_+$ which contains trapped null geodesics.  For small $|a|/m$  one  can show\footnote{See a discussion of trapped null  geodesics in Kerr   in  section \ref{section:morawetz-geodesics}. Note that the trapped  set reduces to  the hypersurface  $r=3m$ in the case of Schwarzschild.}  that   all   such geodesics  are included  in   the  set $ \MM_{trap}$, see Definition \ref{def:causalregions-intro}.

\item \textit{Presence of  a non-trivial  ergoregion.} This is the region   of $r> r_+$ where $\T$ is spacelike.
\end{enumerate}

As mentioned  earlier, to treat these difficulties  in linear theory,    \cite{Ma} and \cite{D-H-R-Kerr}   rely on methods, first  used in the context  of the scalar wave equation in Kerr,    based  on    mode decompositions and construction of  vectorfields adapted  to different modes.  It is however not clear how to extend this method, without loss of  derivatives,  to general perturbations of Kerr. 
 
 In our work\footnote{We need in fact   to adapt the  the method of \cite{A-B} to    the far more difficult 
  case  when the metric is a perturbation of the Kerr metric, $\psi$ is a $2$-horizontal tensor rather than a    scalar and  the  gRW equation  contains additional linear terms.},   we rely instead on a  physical space  method introduced by Blue and Andersson in \cite{A-B}  in the context of   the scalar wave equation  $\square_{a,m} \psi=0$   in  $Kerr(a, m) $ for small $|a|/ m$.  In the next section, we provide a very short review of this method in the simplest case of a scalar wave equation in Kerr.

 
 \subsection{Andersson-Blue method} 
 \lab{sect.AndBluemethod}
 
 
   The crucial new  idea  in \cite{A-B} is to supplement the existing  Killing vectorfields  of $Kerr(a, m)$, i.e. $\T=\pr_t$ and  $\Z=\pr_\phi$  in Boyer Lindquist (BL)  coordinates,   with   a  second order  operator  $\KK =\D_\a(\K^{\a\b} \D_\b) $
   which commutes  with the scalar wave operator $\square_{a,m}$.
 Here   $\K_{\a\b} $ is a Killing tensor, i.e. symmetric  and verifying 
      $\D_{(\ga} \K_{\a\b)}=0$. 
   Note that  if  $X, Y$  are   Killing vectorfields,   their symmetric tensor product  $\frac 1 2 (X\otimes Y+ Y\otimes X ) $ is   automatically   a Killing tensor but the remarkable thing  about Kerr, discovered by Carter,  is that it has an additional Killing tensor which  cannot be reduced  in this manner to Killing vectorfields.  Relative to  the  basis  \eqref{eq:canonicalHorizBasisKerr-intro},  the Carter  tensor  takes  the form
    \bea
    \lab{intro:defineO}
  \K=-a^2 \cos^2\th \g+O, \qquad O=|q|^2 \big( e_1\otimes e_1+ e_2\otimes e_2\big).
  \eea
The associated  second order operator  $\OO= \D_\a(O^{\a\b} \D_\b)$   verifies  itself a  commutation  property with $\square_{a,m}$
\beaa
\,[\OO, |q|^2 \square_{a,m}]= 0.
\eeaa
 In their work,  Andersson and Blue  introduce the set of  second  order  operators  \beaa
 \SS_\aund = \D_\a(S_\aund^{\a\b} \D_\b), \qquad \aund=1,2,3,4,
 \eeaa 
 associated to the Killing   tensors
\bea
 \Big\{  S_1=\T\otimes \T,\quad    S_2= \frac 1 2a \big( \T\otimes \Z+\Z\otimes \T),\quad  S_3 = a^2 \Z\otimes \Z, \quad O\Big\}, 
 \eea
 which commute with   $ |q|^2 \square$. Thus if  $\psi$ is a solution to $\square_{a,m} \psi=0$,  so are $\psi_\aund:=\SS_\aund \psi$ for $\aund=1,2,3,4$.  
 
It is important to note that   $O$  appears naturally  in the expression of the  the inverse  Kerr metric in BL coordinates,  see Lemma \ref{lemma:inversemetricexpressioninKerr},
 \bea
 \lab{eq:inversemetric-in Kerr-intro}
 |q|^2 \g^{\a\b}&=&\De \pr_r^\a \pr_r^\b +\frac{1}{\De} \RR^{\a\b} 
 \eea
with  the $2$ tensor $\RR^{\a\b} $   a linear combination of the tensors $S_\aund$  of the form
 \beaa
          \RR^{\a\b} =\RR^\aund  S_{\aund} ^{\a\b},
          \eeaa
          with coefficients $ \RR^1=-(r^2+a^2)^2$, $\RR^2 = -2a(r^2+a^2)$, $\RR^3 =-a^2$, $\RR^4=\De$.
         
 The main idea in \cite{A-B} is  to  derive  a  Morawetz-Energy    integrated spacetime  estimates  for   an appropriate  combination
 the set of solutions  $\{\psi, \psi_\aund =\SS_\aund\psi|\,   \aund=1, 2, 3,4\}$.   This is done,  roughly,  as follows.
  
 {\bf Step 1.}   Integrate  by parts an  expression of the form
 \beaa
\int_{\MM(\tau_1, \tau_2)}  \sum_{\aund, \bund=1} ^4 \left(X^{\aund\bund}+\frac 1 2 w^{\aund\bund}\right)\psi_{\aund} \c \square_{a, m} \psi _\bund
\eeaa
 for  Morawetz type vectorfields  $X^{\aund\bund}=- z h f^{\aund\bund}  \pr_r $ and  scalars $w ^{\aund\bund}= |q|^2 \Div \big( |q|^{-2}  X^{\aund\bund} \big)
  +h \pr_r z  f^{\aund\bund}$ with a choice  for the function  $f^{\aund\bund}$ of the form
 \beaa
        f^{\underline{a}\underline{b}}&=&   \tilde{\RR}'^{(\underline{a}} \LL^{\underline{b})}= \frac 1 2 \big(\RRtp^{\aund}\LL^{\bund}+\RRtp^{\bund}\LL^{\aund} \big)
\eeaa
where $\RRtp^\aund= \pr_r\Big( \frac{z}{\De} \RRa\Big)$, and $\LL^{\aund}$ are suitable constants.   With appropriate  choices\footnote{The  scalar  functions  $z$ and $h$  are chosen by 
 $z= z_0-\de_0z_0^2$,  $z_0=\frac{\De}{(r^2+a^2)^2}$, $h= \frac{(r^2+a^2)^3}{r} $.} for $z, h$, after suitable integration by parts,  summing and absorbing terms  small in $a/m$  one  derives estimates of the form
\bea
\lab{schematic:Morawetz}
\int_{\MM(\tau_1, \tau_2)}  P  &\les& \left| \int_{\pr\MM(\tau_1,\tau_2)}  B\right|,
\eea
where  $P, B $ are quadratic quantities  involving  combinations of  $\psi_\aund$ and their  first derivatives  with $P$  everywhere positive.  One can show  that outside the trapped  set,
$\int_{\MM(\tau_1, \tau_2)}  P$ is coercive,  that is, roughly,
using the notation $|\psi|^2_{\SS}:=\sum_{\aund=1}^4\big|\psia|^2$,
we can show:
\begin{itemize}
\item The integral  $\int_{\MM(\tau_1, \tau_2)}  P$ controls, with appropriate $r$ weights, 
 the  spacetime integrals of $|\nab_r \psi|^2_{\SS}$ and $ |\psi|_{\SS}^2$ in the entire domain of outer communication   $r> r_+$.  
\item For $r> r_+$ and away from the trapping set  $\int_{\MM(\tau_1, \tau_2)}  P$ also controls  the  corresponding spacetime integrals 
of  all other derivatives, with appropriate $r$ weights,  $ |\nab_\T\psi|^2_{\SS},\,  |\nab_Z\psi|^2_{\SS}, \,  |\nab\psi|^2_{\SS}$.
\end{itemize}

 {\bf Step 2.}   The boundary term  $\int_{\pr\MM(\tau_1,\tau_2)}  B$
     can be eliminated  in  favor  of a positive   energy-flux  type integral\footnote{Here  $E[\psi](\tau)    $
 denote integrals on $\Si(\tau)$  while  $ F[\psi](\tau_1, \tau_2) $  are flux integrals  on the  spacelike boundaries  $\AA\cup \Si_* $.}   $ \sup_{\tau\in [\tau_1,\tau_2]} E[\psi](\tau) +F[\psi](\tau_1, \tau_2) $, with   the help of  an energy identity induced by   the Killing vectorfield $\T$.
This is derived by integrating by parts  in the integral
\beaa
\int_{\MM(\tau_1, \tau_2)} \T \psi \c \square_{a, m} \psi. 
\eeaa
 To avoid  the fact that $\T$ becomes  spacelike in the ergoregion one  can      replace it  by  a smooth, causal,  vectorfield  $\That_\de$     equal to $\T$ in a  small neighborhood  of size $\de$ of the trapping region  and equal to  the future  causal  Hawking  vectorfield $\That =\pr_t+\frac{a}{r^2+a^2} \pr_\phi $, see \eqref{define:That},  everywhere else.  Combining  such  an   energy  inequality with the Morawetz estimate   above we derive
estimates of the form, roughly, 
\beaa
\int_{\MM(\tau_1, \tau_2)}  P  + \sup_{\tau\in [\tau_1,\tau_2]} E[\psi](\tau) +F[\psi](\tau_1, \tau_2)  \les &E[\psi](\tau_1)
\eeaa
where the integrands in $E[\psi]$ are positive definite.

{\bf Step 3.}  Both $E$ and $P$ degenerate  at the horizon $r=r_+$.  This problem can   be easily   fixed by  using  the  standard red shift vectorfield technique, see  \cite{DaRo1},  applied to the region $\MM_{red}$.

  {\bf Step 4.}  To control   all  second   derivatives of $\psi$ away from the trapping set one needs to make  use  of   some type of    coercivity  for  the operator  $\OO$.  In that sense it is important to remark that  in  the integrable  context of   Schwarzschild, or Minkowski  space    $\OO= r^2 \lap$ with $\lap$ the   standard Laplacian on the sphere and thus $\OO\psi $ controls   $\nab^2 \psi$ by standard elliptic estimates on spheres. Such   estimates are  not possible in the non-integrable  case when there are   no  compact  $2$-surfaces  on which $\OO$ acts.  We need in fact  both $\OO$ and $\T$ to control   $\nab^2 \psi$, see  sections  \ref{Section:SpacetimeellipticI} and \ref{section:nonintegrableHodge}. 
  
{\bf Step  5.}  The method outlined above only provides  estimates
 consistent  with  the value $s=2$  in   the estimates\footnote{In  reality, the  method also generates lower order terms  in derivatives of $\psi$   so that the  estimates  cannot be closed  without  control  of  the lower derivatives.}
 of Theorem \ref{THEOREM:GENRW1-P-INTRO}.    To control the lower derivatives  one has to rely on  the  conditional   Morawetz estimates  of  Proposition \ref{proposition:Morawetz1-step1}, i.e. estimates which can only be closed   when combined with the $s=2$  estimates discussed above. These  conditional  estimates\footnote{Note that  these conditional estimates  become unconditional in the axially symmetric case.}  are  derived    using the traditional vectorfield method,    based on   a vectorfield  $X$ of the form $ \FF(r)\pr_r$.

 {\bf Step 6.}    Higher  derivative estimates, i.e. $s>2$,  can  also  be derived by making use of the  commutation properties of the  symmetry operators $\T$, $\Z$ and $\OO$.

  Once a non-degenerate combined Energy-Morawetz estimate has been derived, one can  also  derive  $r^p$-weighted   type estimates
on the large $r$ region. 
This step does not differ much from the case of Schwarzschild since the corrections in $a/r$ are small.  Then, combining the Energy-Morawetz estimate  with these $r^p$ weighted estimates, one   can derive the result of  Theorem \ref{THEOREM:GENRW1-P-INTRO}  for the particular case of the equation $\square_{m,a}\psi=0$.
   We refer the reader to section \ref{section:outline-proofs}  for 
 a more extended outline   of the  main ideas in the proof.  We also note that 
  the   case of the equation  $\square_{m,a}\psi+V\psi$ with $V= \frac{4\De}{ (r^2+a^2) |q|^2}$, which is more relevant to our  situation,    can be treated  in the same manner\footnote{The  positive potential $V$ does actually help in the estimates.}.


\subsection{Proof of Theorem  \ref{THEOREM:GENRW1-P-INTRO}}

 
Here are some of the necessary adjustments  to the Andersson-Blue method  to the case of our model  gRW  equations  \eqref{eq:Gen.RW-intro}:
 \begin{itemize}
 \item     Define approximate Killing   vectorfields  $\T, \Z$, using only the frame and the functions $r,\th$, according to  the formula 
\eqref{eq:TZinKerr-intro}. Define also  the vectorfields
 \beaa
 \That& =&\frac 1 2 \left( \frac{|q|^2}{r^2+a^2} e_4+\frac{\De}{r^2+a^2}  e_3\right),\qquad 
\Rhat=\frac 1 2 \left( \frac{|q|^2}{r^2+a^2} e_4-\frac{\De}{r^2+a^2}  e_3\right).
 \eeaa

 \item
 Define an  appropriate   version    of the $O$ tensor,  expressed   only in terms of the horizontal  structure  and  function $q= r+i a \cos \th$, i.e.
 \beaa
O^{\a\b}&:=& |q|^2\ga^{ab}e_a^\a e_b^\b,
\eeaa
where $\ga$ is the metric induced by $\g$ on the horizontal structure.

 \item Show that the  error  terms generated by  commuting $\T, \Z$ and   the operator   $\OO= \Ddot_\a(O^{\a\b} \Ddot_\b)$   with $|q|^2 \squared_2$ are acceptable error terms. It is important to note  
  that, due to the tensorial character  of the  gRW  equations,     even  in the case of   Kerr  the commutators generate linear terms  proportional to $a$, see  section \ref{section:CartertensorKerr}.    These   can  ultimately be absorbed  for sufficiently small $a/m$.
 
\item   Consider    the approximate Killing tensors $S_\aund$, $ \aund=1, 2, 3,4,$ 
\beaa
 \left\{  S_1=\T\otimes \T,\quad    S_2= \frac 1 2 a\big( \T\otimes \Z+\Z\otimes \T),\quad  S_3 =a^2 \Z\otimes \Z, \quad  S_4=O\right\}
 \eeaa
 and the corresponding second order operators $ \SS_\aund:=\Ddot_a \big(S_\aund^{\a\b} \Ddot_\b)$ such that the  commutators with $|q|^2\squared_2$ produce only acceptable  error terms. 
\end{itemize}

The proof of Theorem  \ref{THEOREM:GENRW1-P-INTRO} is the most technical part of the paper.  It is carried out first in Kerr in Chapters \ref{chapter-proof-part1} and \ref{chapter-proof-mor-2}, and is then extended to perturbations of Kerr in Chapter \ref{chapter-perturbations-Kerr}.

     
\section{Main ideas in the proof of Theorem \ref{THEOREMM8:INTRO}}


        To  derive  top derivative   estimates for all  curvature components we  need to rely on the following   facts:
        \begin{itemize}
       \item  We work in the global frame  of $\MM$ for which  the conditions     \eqref{eq:Framecondfor-qfb}  are verified in $\Mext$,   see also  Remark  \ref{remark:framecond-qfb-intro}.
               
        \item Control of the initial  data, see \ref{eq:controlofinitialdataforThM8-intro},\, $ \Ik_{k_L}\les \ep_0$.
        
        \item  Improved\footnote{By improved estimates, we mean  in terms of $\ep_0$ rather than $\ep$.}   $r^p$-estimates for $\qfb$  in the case  $p=\de$ for all derivatives  $s\le k_L-2$. This is done as in Part II, noticing that one can go to the maximum number of derivatives in the case $p=\de$. These  estimates rely heavily on the  global   frame  of $\MM$  in which  the conditions \eqref{eq:Framecondfor-qfb}   hold true. 
         
        \item The fact that the desired estimates for $\Rk_{J+1}$ are conditional on the constant $\ep_J$ where\footnote{This is based on an iteration assumption allowing to recover the control of $J+1$ derivatives from the one of $J$ derivatives. To initiate the iteration procedure, we start at $J=\frac{\kl}{2}$ for which we have full control of $\Sk_J$ and $\Rk_J$ thanks to the decay estimates established in Theorem M1 and M2 for $A$ and $\Ab$, and in Theorems M3 to M7 of \cite{KS:Kerr} for all other linearized curvature components and all linearized Ricci and metric coefficients, see Lemma 9.4.13 in \cite{KS:Kerr} for the corresponding statement.}
           \bea\lab{eq:iterationassumptiondiscussionThM8:intro}
\Sk_J+\Rk_J  & \leq&\ep_J,
\eea   
which allows to deal with lower order terms. 
 
 \item     We take into account the  quantitative assumptions \eqref{eq:assumptionsonMextforpartIII-intro} on $\Rk_k$ and $\Sk_k$  
        \bea\lab{eq:mainbootassforchapte9-intro}
\Sk_k+\Rk_k &\leq& \ep, \quad k\le k_L,
\eea
which allows to deal with nonlinear terms.
  \end{itemize}
        
        The proof of Theorem \ref{THEOREMM8:INTRO}  can be divided in two main parts:
        \begin{enumerate}
        \item Global Energy-Morawetz estimates in $ \MM$.  This step relies heavily  on the  following  main ingredients:
        \begin{enumerate} 
        \item  Derive first      $\BEF_\de^J$   estimates for $\Pc$ by relying on a suitable linearization of the wave equation of $P$.  The RHS of this scalar wave equation contains   quadratic terms,  which are dealt with  our quantitative assumptions \eqref{eq:mainbootassforchapte9-intro},  and linear terms  containing  fewer derivatives
    which can be dealt by the iterations assumption \eqref{eq:iterationassumptiondiscussionThM8:intro}.   

        \item Use  these estimates for $\Pc$ together  with the Bianchi identities 
        to derive $\BEF_\de^J$ estimates for all  other curvature components 
         $A, B, \Ab,\Bb$. The RHS of the Bianchi identities contain also  quadratic terms,  which are dealt with  our quantitative assumptions \eqref{eq:mainbootassforchapte9-intro},  and linear terms  containing  fewer derivatives
    which can be dealt by the iterations assumption \eqref{eq:iterationassumptiondiscussionThM8:intro}. The most demanding part of the proof is to deal with the estimates in the trapping\footnote{This requires in particular to differentiate between the good direction $\Rhat$ and the other directions, and make systematic use of the Bianchi identities to recover the remaining derivatives.}. Given that $\Pc$ is already under control, we rely on the triangular structure of the Bianchi identities which allow us to control first $(B, \Bb)$ from $\Pc$, and then $A$ from $B$ and $\Ab$ from $\Bb$. 
        \end{enumerate}

        \item  Weighted  estimates   in the  region  of $\Mext$, where $r\ge r_0$ for sufficiently large $r_0$, for all top derivatives of the curvature components. This is proved   using $r^p$ weighted estimates for Bianchi pairs in a similar fashion as  the corresponding  estimates in  section 8.7 of \cite{KS}. 
        \end{enumerate}

  
 \section{Organization} 
 
 
The paper  is organized  in  3  Parts  and  16 Chapters as follows:

\begin{itemize}
\item Part I:
\begin{itemize}
\item In Chapter 2, we give a full account of our geometric framework based on non-integrable horizontal structures and derive the null structure and null Bianchi equations. We also rephrase these equations using complex notations.

\item In Chapter 3, we provide the main formulas in Kerr.

\item In Chapter 4, we discuss the linearized quantities, introduce our notations $(\Ga_b, \Ga_g)$, as well as other important quantities in perturbation of Kerr. We also provide numerous useful commutators.

\item In Chapter 5, we derive, in perturbations of Kerr, the Teukolsky equations, the generalized Regge-Wheeler (gRW) equations, and the Teukolsky-Starobinsky identities.
\end{itemize}

\item Part II:
\begin{itemize}
\item In Chapter 6, we introduce a simplified RW model and provide precise statements for the corresponding Energy-Morawetz-$r^p$-weighted estimates.

\item In Chapters 7 and 8, we adapt the Blue Andersson method to derive basic  Energy-Morawetz estimates for the RW model equations \eqref{eq:Gen.RW-intro} in the particular case of Kerr.

\item In Chapter 9, we extend the results of Chapters 7 and 8 to perturbations of Kerr and also deal the case of higher order derivatives.

\item In Chapter 10, we conclude the study of the RW model \eqref{eq:Gen.RW-intro} by deriving $r^p$-weighted estimates for it, hence proving the statements of Chapter 6.

\item In Chapter 11, we extend the results of Chapter 6 on Energy-Morawetz-$r^p$-weighted estimates for the RW model \eqref{eq:Gen.RW-intro} to the full generalized Regge Wheeler equations for $(\qf, A)$. We then rely on this result to prove Theorem M1, restated here as  Theorem \ref{theoremM1:intro}.

\item In Chapter 12, we extend the results of Chapter 6 on Energy-Morawetz-$r^p$-weighted estimates for the RW model \eqref{eq:Gen.RW-intro} to the full generalized Regge Wheeler equations for $(\qfb, \Ab)$. We then rely on this result to prove Theorem M2, restated here as  Theorem \ref{theoremM2:intro}.
\end{itemize}

\item Part III:
\begin{itemize}
\item In Chapter 13, we provide the geometric set up and the precise statement of Theorem \ref{THEOREMM8:INTRO} on the control of high derivatives curvature estimates. 

\item In Chapter 14, we derive energy-Morawetz estimates for $\Pc$.

\item In Chapter 15, we derive energy-Morawetz estimates for $(B, \Bb)$ and then $(A, \Ab)$.

\item  In Chapter 16, we derive $r^p$ weighted estimates for the Bianchi pairs to recover the correct $r$ decay using the control of Energy-Morawetz for $(A, B, \Pc, \Bb, \Ab)$, hence concluding the proof of Theorem \ref{THEOREMM8:INTRO}.  
\end{itemize}
\end{itemize}


\section{Acknowledgments}


The first author is supported by the NSF grant DMS 2128386.
The second author is supported   by  the  NSF grant  DMS 180841 as well as by  the Simons grant  10011738. He would like to thank the Laboratoire Jacques-Louis Lions  of Sorbonne Universit\'e  and   IHES   for their  hospitality during  his many visits.  The third author is supported by ERC grant  ERC-2016 CoG 725589 EPGR.


\part{Formalism and derivation of the main equations}
\label{part-I}



\chapter{Non-integrable structures}\label{CHAPTER-NON-INTEGRABLE-STRUCTURES}



\section{A general formalism for non-integrable structures}\lab{section-general-formalism}


We present here a general formalism extending the one used for perturbations of Minkowski space \cite{Ch-Kl} to perturbations of Kerr spacetimes.


\subsection{Null pairs and horizontal structures}
\lab{sec:nullparisandhorizontalstruct}


Let $(\MM, \g)$ be a Lorentzian spacetime. 
Consider an arbitrary null pair $e_3=\Lb$, $e_4=L$, i.e.
\beaa
\g(e_3, e_3)=\g(e_4, e_4)=0, \qquad  \g(e_3, e_4)=-2.
\eeaa

\begin{definition}
A vectorfield $X$ is $(L,\Lb)$-horizontal, or
simply horizontal,  if
\beaa
\g(L,X)=\g(\Lb, X)=0. 
\eeaa
We denote by $\O(\MM)$ the set of horizontal vectorfields on $\MM$.  Given a fixed 
 orientation  on $\MM$,  with corresponding  volume form  $\in$,  we define  the induced 
 volume form on   $\O(\MM)$ by,
 \bea
 \in(X, Y):=\frac 1 2\in(X, Y, \Lb, L).  \label{vol.form}
 \eea
\end{definition}

Given a null pair $(L, \Lb)$, the horizontal vectorfields $\O(\MM)$ define a horizontal distribution, i.e. a sub-bundle of the tangent bundle $\T(\MM)$ of the manifold. 
In the standard terminology used in differential topology, a subbundle $E \subset \T(\MM)$ of the tangent bundle is said to be integrable if for any vectorfields $X$ and $Y$ taking values in $E$, the Lie bracket $[X, Y]$ takes values in $E$ as well. According to the Frobenius  theorem  a subbundle $E$ is integrable (or involutive) if and only if the subbundle $E$ arises from a regular foliation of $\MM$, i.e. if locally the subbundle $E$ can be realized as the tangent space of a submanifold of $\MM$. 
An  useful  example  of integrable  structures, which have played an important role 
 in  \cite{Ch-Kl},   is  that provided  by $S$-foliations, i.e. regular foliations  whose  leaves are  topological spheres orthogonal, at every point,  to  the  null pair $(L, \Lb)$.  In this work   we consider  general,
 not necessarily integrable, horizontal structures.

  Given an arbitrary vectorfield $X$ we denote by $\Xho$
its  horizontal projection,
\beaa
\Xho&=&X+ \frac 1 2 \g(X,\Lb)L+ \frac 1 2   \g(X,L)\Lb.
\eeaa

\begin{definition}
A  $k$-covariant tensor-field $U$ is said to be horizontal,  and denoted $U\in \O_k(\MM)$,
if  for any vectorfields $X_1,\ldots X_k$ we have,
\beaa
U(X_1,\ldots X_k)=U(\Xho_1,\ldots \Xho_k).
\eeaa
\end{definition}

Define the projection operator
\beaa
\Pi^{\mu\nu}&=&\g^{\mu\nu}+ \frac 1 2 (\Lb^\mu L^\nu+L^\mu\Lb^\nu).
\eeaa
Clearly $ \Pi_{\a}^{\mu}\Pi_{\mu}^\b=\Pi_{\a}^\b$.  An arbitrary tensor 
$U_{\a_1\ldots\a_m}$ is horizontal  iff
\beaa
\Pi_{\a_1}^{\b_1}\ldots \Pi_{\a_m}^{\b_m}\, U_{\b_1\ldots\b_m}=U_{\a_1\ldots\a_m}.
\eeaa

 \begin{definition}\lab{def:proxyfirstandsecondfundameentalform} 
 For any horizontal $X,Y$ we define\footnote{In the particular case where the horizontal structure is integrable, $\ga$ is the induced metric, and $\chi$ and $\chib$ are the  null second fundamental forms.} 
  \bea
 \ga(X,Y)&=&\g(X,Y)
 \eea
 and   
\begin{equation}\label{fo1}
\begin{cases}
&\chi(X,Y)=\g(\D_XL,Y),\\ 
&\chib(X,Y)=\g(\D_X\Lb,Y).
\end{cases}
\end{equation}
where $\D$ denotes the covariant derivative of $\g$.
\end{definition}

Observe that    $\chi$
 and $\chib$  are  symmetric if and 
 only if   the horizontal structure is 
 integrable. Indeed this follows easily from
 the formulas,
 \beaa
 \chi(X,Y)-\chi(Y,X)&=&\g(\D_X L, Y)-\g(\D_YL,X)=-\g(L, [X,Y]),\\
 \chib(X,Y)-\chib(Y,X)&=&\g(\D_X \Lb, Y)-\g(\D_Y\Lb,X)=-\g(\Lb, [X,Y]).
\eeaa
 We can view $\ga$,  $\chi$ and $\chib$ as horizontal 2-covariant tensor-fields
 by extending their definition to arbitrary vectorfields  $X, Y$  according to,
 \beaa
  \ga(X,Y)&=&\ga(\Xho,\Yh)
  \eeaa
  and
  \beaa
  \chi(X,Y)&=&\chi(\Xho, \Yh),\qquad \chib(X,Y)=\chib(\Xho, \Yh).
 \eeaa

  Given a  general 2-covariant horizontal  tensor $U$
  we decompose it in its symmetric and antisymmetric part as follows,
  \beaa
  \Us(X,Y)&=&\frac 1 2 \big(U(X, Y)+U(Y, X)\big),\\
   \Ua(X,Y)&=&\frac 1 2 \big(U(X, Y)-U(Y, X)\big).
  \eeaa
  Given a  horizontal structure  defined by $e_3=\Lb$,  $e_4=L$ 
  we  associate a null frame by choosing   orthonormal  horizontal  vectorfields
  $e_1, e_2$  such that $\ga(e_a, e_b)=\de_{ab}$.  By convention, we   say that 
  $(e_1, e_2)$ is positively oriented  on $\O(\MM)$  if,
  \bea
  \in(e_1, e_2)=\frac 1 2 \in(e_1, e_2,  e_3, e_4)=1.
  \eea
  
  \begin{remark}
  We note that the particular choice of an orthonormal basis  is immaterial.  All the quantities we work with are tensorial with respect to the horizontal structure.
  \end{remark}

  Given a covariant horizontal 2-tensor $U$   and an arbitrary
 orthonormal  horizontal frame $(e_a)_{a=1,2}$ we have,
\beaa
\Us_{ab}=\frac 12 (U_{ab}+U_{ba}),\qquad \Ua_{ab}=\frac 12 (U_{ab}-U_{ba}).
\eeaa
\begin{definition}
The trace  of a  horizontal  2-tensor $U$ is defined by
\bea
\tr (U):=\de^{ab}U_{ab}=\de^{ab}\Us_{ab}.
\eea
We define the anti-trace of $U$  to be
\bea
\atr (U):=\in^{ab}U_{ab}=\in^{ab}\Ua_{ab}.
\eea
Observe that  the  first    trace  is   independent  of  the particular choice  of  the  frame $e_1, e_2$.
On the other hand, for fixed $e_3, e_4$,   $\atr$  depends on the orientation of $e_1, e_2$.
Also, by interchanging  $e_3, e_4$, $\atr$  changes sign.
\end{definition}

A general  horizontal  2-tensor $U$  can be decomposed according to
\bea
U_{ab}&=&\Us_{ab}+\Ua_{ab}=\widehat{U}_{ab}+\frac 1 2 \de_{ab}\, \tr( U)+\frac 12 \in_{ab}\atr (U), \label{decompU}
\eea
where $\widehat{U}$ denotes the symmetric traceless part of $U$.

\begin{definition}
We  introduce the  notation
\bea
\trch:=\tr(\chi), \quad \atr\chi:=\atr(\chi),\quad \trchb:=\tr(\chib), \quad \atr\chib:=\atr(\chib).
\eea
The quantities $\chih, \trch$ and $\chibh, \trchb$  are called, respectively,  the shear and expansion 
of the horizontal distribution $\O(\MM)$. The scalars $ \atr\chi$ and $\atr\chib$ measure the
integrability defects of the distribution. 

Accordingly, we  decompose $\chi, \chib$ as follows
\beaa
\chi_{ab}&=&\chih_{ab} +\frac 1 2 \de_{ab} \trch+\frac 1 2 \in_{ab}\atrch,\\
\chib_{ab}&=&\chibh_{ab} +\frac 1 2 \de_{ab} \trchb+\frac 1 2 \in_{ab}\atrchb.
\eeaa
\end{definition}

In what follows we fix  a null pair $(e_3, e_4)$ and an orientation  on $\O(\MM)$.
\begin{definition}\label{definition-hodge-duals}
We define the left and right duals of a horizontal 1-form $\xi$ and a $2$-covariant  tensor-field $U$,
\beaa
&&\quad\,\,\dual \xi_{a}=\in_{ab}\xi_b,\qquad\quad\,\,\,\, \xi^*\, _{a}=\xi_b\in_{ba},\\
&&(\dual U)_{ab}=\in_{ac} U_{cb},\qquad (U^*)_{ab}=
U_{ac}\in_{cb}.
\eeaa
\end{definition}

\begin{lemma}
Given a horizontal 1-form $\xi$, we have
\beaa
\dual(\dual \xi)=-\xi, \qquad \dual\xi=-\xi^*.
\eeaa
\end{lemma}

\begin{lemma} 
Given a covariant horizontal 2-tensor $U$, we have
\begin{enumerate}
\item $\dual (\dual U)=-U$.

\item If $U$ is symmetric, then $\dual U_{ab}=- U^*_{ba}$. 

\item If  $U=\hat U$ is symmetric traceless,
then $\dual \hat U=-\hat U^*$ is also symmetric 
traceless. 
\item  In general, 
\beaa
\tr(\dual U)&=&\tr(U^*)=-\atr(U),\\
\atr(\dual U)&=&\atr(U^*)=\tr(U),\\
\widehat{\dual U}&=&\dual\hat{U}.
\eeaa
\end{enumerate}
\end{lemma}

Given a general horizontal 2-tensor $U$ we have, according to \eqref{decompU}, 
\beaa
U^*_{ab}&=&\hat{U}^*_{ab}+\frac 1 2 \in_{ab}\, \tr( U)- \frac 12 \de_{ab}\atr (U), \\
\dual U_{ab}&=&\dual \hat{U}_{ab}+\frac 1 2 \in_{ab}\, \tr( U)- \frac 12 \de_{ab}\atr (U).
\eeaa
Hence,
\beaa
U^*_{ab}&=&-\dual U_{ab}+\in_{ab}\, \tr( U)-  \de_{ab}\atr (U).
\eeaa

We note the following lemma.
\begin{lemma}
\label{le:duals}
Given two  1-forms $\xi, \eta$ we have,
\beaa
\dual\xi \c  \eta=-\xi\dual\c \eta=\xi\c\eta^*.
\eeaa
Given  a  1-form  $\xi$ and    2-tensor $U$ we have,
\beaa
 \xi_a  U_{ab} &=& \xi^a\hat{U}_{ab}+\frac 1 2 \xi_b \tr(U)-\frac 1 2 \xi_b \atr(U),\\
\dual \xi^a U^*_{ab} &=& \xi^a\hat{U}_{ab}+\frac 1 2 \xi_b \tr(U)-\frac 1 2 \xi_b \atr(U).
 \eeaa
Thus,
\beaa
\dual \xi^a U^*_{ab}+ \xi^a  U_{ab}&=& \xi_b( \tr U)-\dual \xi_b (\atr U ).
\eeaa
Also,
\beaa
\dual \xi^a U^*_{ab}- \xi^a  U_{ab}&=&-2\xi^a \hat{U}_{ab}.
\eeaa
\end{lemma}

\begin{proof} 
We have
\beaa
\dual \xi^a U^*_{ab}&=&\dual \xi^a\left(\hat{U}^*_{ab}+\frac 1 2 \in_{ab}\tr(U)-\frac 1 2 \de_{ab} \atr(U)\right)\\
&=&-\xi^a\hat{U}_{ab}+\frac 1 2 \xi_b\tr(U)-\frac 1 2 \dual \xi_b\atr(U),\\
 \xi_a  U_{ab}&=&\xi^a\left(\hat{U}_{ab}+\frac 1 2 \de_{ab} \tr(U)+\frac 1 2 \in_{ab} \atr(U)\right)\\
&=& \xi^a\hat{U}_{ab}+\frac 1 2 \xi_b \tr(U)-\frac 1 2 \xi_b \atr(U),
 \eeaa
 which implies the formulas involving $U$ and $\xi$.
\end{proof}

\begin{definition}\label{definition-SS-real}
We denote\footnote{Using the convention of raising and lowering indices  we make no distinction here between
  covariant and contravariant  tensors.}  by $\O_k(\MM)$ the set of all horizontal  tensor-fields of rank $k$ on $\MM$.
We denote by $\sk_0=\sk_0(\MM)$ the set of pairs of real scalar functions on $\MM$,  $\sk_1=\sk_1(\MM)$ the  set of real horizontal $1$-forms  on $\MM$   and for,  $ k \ge 2 $,      $\sk_k(\MM)$  the set of fully symmetric traceless  horizontal real tensors of rank $k$.  In particular $\sk_2=\sk_2(\MM)$ denotes 
  the set of symmetric traceless   horizontal real $2$-tensors on $\MM$.
\end{definition}

\begin{definition} 
Given real  $\xi, \eta\in\sk_1 $  we denote
\beaa
\xi\c \eta&:=&\de^{ab} \xi_a\eta_b,\qquad
\xi\wedge\eta:=\in^{ab} \xi_a\eta_b=\xi\c\dual \eta,\qquad
(\xi\hot \eta)_{ab}:= \xi_a \eta_b +\xi_b \eta_a-\de_{ab} \xi\c \eta.
\eeaa
 Given   $\xi\in \sk_1 $,  $U\in \sk_2$ we denote
\beaa
(\xi\c U)_a&:=&\de^{bc} \xi_b U_{ac}.
\eeaa
Given     $U, V \in \sk_2$ we denote
\beaa
(U\wedge V)_{ab} &:=& \in^{ab}U_{ac}V_{cb}.
\eeaa
\end{definition}

The following two lemmas   are   immediate.
\begin{lemma}
 Given   $\xi\in \sk_1 $,  $U\in \sk_2$, we have
\beaa
U_{ab}\xi^b&=& \left(\hat{U}\c\xi+\frac{1}{2}\tr(U)\xi+\frac{1}{2}\atr(U)\dual\xi\right)_a.
\eeaa 
\end{lemma}

\begin{lemma}
\label{le:sym-product}
Given $\hat U, \hat V \in \sk_2$
we have, with respect to an arbitrary orthonormal basis, 
\beaa
\hat U_{ac}\hat V_{cb}+\hat V_{ac}\hat U_{cb}=\de_{ab} \hat U\c\hat V
\eeaa
where
 \beaa
 \hat U\c\hat V=\de^{ac}\de^{bd} \hat U_{ab} \hat V_{cd}.
 \eeaa
In particular
\beaa
\hat V_{ac }\hat V_{cb}=\frac 1 2 \de_{ab}|\hat V|^2
\eeaa
with $|\hat V|^2=\hat V\c\hat V$.
\end{lemma}

\begin{remark}\lab{rmk:tracelesspartofproducttracelessis0}
The previous lemma implies in particular $\widehat{\hat{U}_{ac}\hat{V}_{cb}}=0$.
\end{remark}

We generalize  the lemma as follows.
\begin{lemma}
\lab{le:traces}
Given $U,V$  arbitrary 2-covariant horizontal tensor-fields, we  have
\beaa
\de^{ab}U_{ac}V_{cb}&=&\hat U\c\hat V+\frac 1 2 \big(\tr(U)\tr(V)-\atr(U)\atr(V)\big),\\
\in^{ab}U_{ac}V_{cb}&=&\hat U\wedge \hat V+\frac 12 \big(\atr(U)\tr (V)+\tr(U)\atr(V)\big),\\
\widehat{U_{ac} V_{cb}}&=&\frac 1 2 \big(\Uh_{ab}\tr(V)+\Vh_{ab}\tr(U)\big)
+\f12\big(-\dual\Uh_{ab}\atr(V)+\dual\Vh_{ab}\atr(U)\big),
\eeaa
where
\beaa
\hat U\cdot\hat V&=&\de^{ac}\de^{bd}\hat U_{ab}\hat V_{cd},\\
\hat U\wedge \hat V&=&\hat U\cdot \dual \hat V=\in^{ab}\hat U_{ac}\hat V_{cb}.
\eeaa
\end{lemma}

\begin{proof} 
In view of the decomposition \eqref{decompU}, we have
\beaa
U_{ac} V_{cb}&=&\Uh_{ac} \Vh_{cb}+\frac 1 2 \big(\tr(V) \Uh_{ab}+\tr(U) \Vh_{ab}\big)
+\frac 1 2\big(  \atr(U)\dual \Vh_{ab}+\atr(V) \Uh^*_{ab}\big)\\
&+&\frac 1 4 \big( \tr(U) \tr(V)-\atr(U)\atr(V)\big)\de_{ab}+\frac 1 4 \big( \tr(U) \atr(V)+\atr(U)\tr(V)\big)\in_{ab}
\eeaa
and the proof easily follows, using also the fact that $\widehat{\hat{U}_{ac}\hat{V}_{cb}}=0$ according to Remark \ref{rmk:tracelesspartofproducttracelessis0}.
\end{proof}

The following is an  immediate consequence of Lemma \ref{le:traces}.
\begin{corollary}
\label{special.products}
In the particular case when $U=V$, we have
\beaa
\de^{ab}U_{ac}U_{cb}
&=&|\hat U|^2+\frac 1 2 \big((\tr(U))^2-(\atr(U))^2\big),\\
\in^{ab}U_{ac}U_{cb}&=&\tr(U)\atr(U),\\
\widehat{U_{ac}U_{cb}} &=& \tr(U)\Uh_{ab}.
\eeaa
\end{corollary}

As another corollary to Lemma \ref{le:traces} we have the following.
\begin{lemma}
 \lab{le:nonsym-product}
Let   $u$ be an arbitrary $2$-horizontal tensor and $v\in \sk_2$. Then
\beaa
u_{ac} v_{cb}+ u_{bc} v_{ca}&=&\de_{ab} \hat{u}\c v   +(\tr u )v_{ab}  +\frac 1 2\Big[\big( u_{ac}-u_{ca}\big) v_{cb}+ \big( u_{bc}-u_{cb}\big)v_{ca}\Big].
\eeaa
\end{lemma}

\begin{proof}
We give below a  direct proof based on Lemma \ref{le:sym-product}
 according to which,  given  $u, v\in\sk_2$,   we have
\beaa
u_{ac} v_{cb}+ u_{bc} v_{ca}=\de_{ab} u\c v. 
\eeaa
If $u$ is only symmetric and $v\in \sk_2$   we  can write,
\beaa
u_{ac} v_{cb}+ u_{bc} v_{ca}&=&\left( \hat{u}_{ac}+\frac 1 2 \de_{ac} tr (u)\right) v_{cb}+  \left( \hat{u}_{bc}+\frac 1 2 \de_{bc} tr (u)\right) v_{ca}\\
&=&\de_{ab}\hat{u} \c v   +(\tr u )v_{ab}.
\eeaa
If $u$ is an arbitrary $2$-tensor and $v\in \sk_2$,
\beaa
u_{ac} v_{cb}+ u_{bc} v_{ca}&=&\frac 1 2 \Big( u_{ac}+u_{ca} +\big( u_{ac}-u_{ca}\big)\Big) v_{cb}+\frac 1 2  \Big( u_{bc}+u_{cb} +\big( u_{bc}-u_{cb}\big)\Big) v_{ca}\\
&=&\frac  12 \de_{ab} \Big( u_{ac}+u_{ca}\Big) v_{ac} +\frac 1 2\Big( \big( u_{ac}-u_{ca}\big) v_{cb}+ \big( u_{bc}-u_{cb}\big) v_{ca}\Big)\\
&=&\de_{ab} \hat{u}\c v   +(\tr u )v_{ab}  +\frac 1 2\Big(\big( u_{ac}-u_{ca}\big) v_{cb}+ \big( u_{bc}-u_{cb}\big) v_{ca}\Big).
\eeaa
This concludes the proof of the lemma.
\end{proof}

\begin{lemma}\label{lemma:usefulidentitiesforcomplexification}
The following formulas  hold true:
\begin{itemize}
\item Given $\xi, \eta \in \sk_1$, we have
\beaa
\bsplit
\dual\xi \c  \eta &= -\xi\c\dual\eta,\quad 
\dual\xi \c  \dual\eta = \xi\c\eta,\quad
\dual\xi\wedge\eta = -\xi\wedge\dual\eta,\quad 
\dual\xi\wedge\dual\eta = \xi\wedge\eta,\\
\dual\xi \hot  \eta &= \xi\hot\dual\eta,\quad 
\dual(\xi\hot\eta) = \dual\xi \hot  \eta,\quad 
\dual\xi \hot  \dual\eta = -\xi\hot\eta.
\end{split}
\eeaa

\item Given  $\xi \in \sk_1, U \in \sk_2$, we have
\beaa
\dual(\xi\c U) &=& \xi\c\dual U,\qquad 
\dual\xi\c U = -\xi\c\dual U,\qquad 
\dual\xi\c\dual U = \xi\c U.
\eeaa

\item Given $U,V \in \sk_2$, we  have,
\beaa
\dual U\c V = -U\c\dual V,\,\,\,\,
\dual U\c\dual V = U\c V,\,\,\,\, 
\dual U\wedge V = -U\wedge\dual V,\,\,\,\,  
\dual U\wedge\dual V = U\wedge V.
\eeaa
\end{itemize}
\end{lemma}

\begin{proof}
The statements follow from the above results except the ones involving $\hot$. To check those  we write, for an arbitrary basis $e_1, e_2$, 
      \beaa
      \dual( \xi\hot \eta)_{11}&=&( \xi\hot \eta)_{21} =   \xi_2 \eta_1+\xi_1 \eta_2, \\
      ( \xi\hot \dual \eta)_{11}&=& \xi_1( \dual \eta)_1- \xi_2 (\dual \eta)_2 =\xi_1 \eta_2 + \xi_2 \eta_1, \\
      (\dual \xi \hot \eta)_{11} &=&(\dual \xi)_1 \eta_1-(\dual \xi)_2 \eta_2= \xi_2 \eta_1+ \xi_1 \eta_2,
      \eeaa
      and
      \beaa
        \dual(\xi\hot\eta)_{12}&=&\xi_2\eta_2-\xi_1 \eta_1, \\
        (\xi\hot \dual\eta)_{12}&=&\xi_1 \dual \eta _2+\xi_2\dual  \eta_1 = -\xi_1 \eta_1 + \xi_2 \eta_2,\\
        (\dual\xi \hot\eta)_{12} &=&\dual \xi_1 \eta_2+ \dual\xi_2 \eta_1  =-\xi_1 \eta_1 + \xi_2 \eta_2.
      \eeaa
      Hence,
      \beaa
      \dual(\xi\hot\eta)=\dual \xi \hot \eta= \xi\hot \dual\eta
      \eeaa  
      as stated.
\end{proof}

\begin{lemma}
\lab{dot-hot}
Given  $\xi , \eta \in \sk_1$, $  u\in \sk_2$  we have 
\beaa
\xi \hot ( \eta  \c u) + \eta  \hot ( \xi \c u)= 2  (\xi\c \eta  ) u. 
\eeaa
\end{lemma}

\begin{proof}
Straightforward verification.
\end{proof}

  
  \subsection{Horizontal covariant derivative}
  
  
 Given $X, Y\in\O(\MM)$, the covariant derivative $\D_XY$ fails in general  to be horizontal.
 We thus define the horizontal covariant operator $\nab$ as follows
 \bea
 \nab_X Y :=^{(h)}(\D_XY)=\D_XY- \frac 1 2 \chib(X,Y)L -  \frac 1 2 \chi(X,Y) \Lb.
 \eea
 
 \begin{proposition}
 For all  $X,Y\in \O(\MM)$,
  \beaa
  \nab_X Y-\nab_Y X
  &=&[X,Y]-\chiba(X,Y)L -\chia(X,Y)\Lb\\
  &=&[X, Y]-\frac 1 2 \left(\atrchb\,  L+\atrch \, \Lb\right)\in(X, Y).
 \eeaa
 In particular,
 \bea\,^{(h)}[X, Y]&=&\frac 1 2\left(\atrchb\,  L+\atrch \, \Lb\right)\in(X, Y).
 \eea
 
 For  all  $X,Y, Z\in \O(\MM)$,
 \beaa
 Z\ga(X,Y)=\ga(\nab_Z X, Y)+\ga(X, \nab_ZY).
 \eeaa
\end{proposition}

\begin{remark} 
In the integrable case, $\nab$ coincides with the Levi-Civita connection
 of the metric induced on the integral surfaces of   $\O(\MM)$. 
\end{remark} 
 
 Given a general covariant, horizontal tensor-field  $U$
 we define its horizontal covariant derivative according to
 the formula,
 \beaa
 \nab_Z U(X_1,\ldots X_k)=Z (U(X_1,\ldots X_k))&-&U(\nab_ZX_1,\ldots X_k)-\nn
 \\
\ldots   &-& U(X_1,\ldots \nab_ZX_k).
 \eeaa
 
 Given $X$ horizontal, $\D_L X$ and $\D_\Lb X$ are in general not horizontal.
 We define $\nab_L X$ and $\nab_\Lb X$  to be the horizontal projections
 of the former.  More precisely,
 \beaa
 \nab_L X&:=&^{(h)}(\D_L X)=\D_L X- \g(X, \D_L\Lb) L- \g(X, \D_L L) \Lb,\\
 \nab_\Lb X&:=&^{(h)}(\D_\Lb X)=\D_\Lb X- \g(X, \D_\Lb\Lb) L - \g(X, \D_\Lb L) \Lb. 
 \eeaa
 
  We can extend the operators $\nab_L$ and $\nab_\Lb$ to
 arbitrary  $k$-covariant, horizontal tensor-fields  $U$  as follows,
 \beaa
 \nab_LU(X_1,\ldots, X_k)=L(U(X_1,\ldots, X_k))&-&U(\nab_L X_1,\ldots, X_k)-\\
 \ldots&-&U( X_1,\ldots, \nab_LX_k),\\
  \nab_\Lb U(X_1,\ldots, X_k)=\Lb(U(X_1,\ldots, X_k))&-&U(\nab_\Lb X_1,\ldots, X_k)-\\
  \ldots&-&U( X_1,\ldots, \nab_\Lb X_k).
 \eeaa 
 The following  proposition follows easily from the definition.
 
 \begin{proposition}
 The operators $\nab$, $\nab_L$ and $\nab_\Lb$ take horizontal tensor-fields into
 horizontal tensor-fields. We have,
 \bea
 \nab \ga=\nab_L\ga=\nab_\Lb\ga=0.
 \eea
 \end{proposition}

 We now extend the definition of horizontal covariant derivative to any $X\in \T(\MM)$ in the tangent space of $\MM$  and $Y\in\O(\MM)$.
\begin{definition}
Given  $X\in \T(\MM)$ and $Y \in \O(\MM)$ we define,
\beaa
\Ddot_X Y&:=& ^{(h)}( \D_X Y).
\eeaa
Given an orthonormal  frame  $e_1, e_2\in \O(\MM)$ we write 
\beaa
\Ddot_\mu  e_a&=&\sum_{b=1,2}(\La_\mu)_{ba}\, e_b,  \qquad  (\La_\mu)_{\a\b}:=\g(\D_\mu e_\b, e_\a).
\eeaa
\end{definition} 

 \begin{definition}
 \label{definition:S-covariantderivative}
 Given a general, covariant,  horizontal tensor-field  $U$
 we define its horizontal covariant derivative according to
 the formula
 \beaa
 \Ddot_X U(Y_1,\ldots Y_k)=X (U(Y_1,\ldots Y_k))&-&U(\Ddot_X Y_1,\ldots Y_k)-\ldots- U(Y_1,\ldots \Ddot_XY_k),
 \eeaa
 where $X\in \T(\MM)$ and $Y_1,\ldots Y_k\in \O(\MM)$.
 \end{definition}
 
  \begin{proposition}
 For  all  $X\in\T(\MM)$   and $Y_1, Y_2 \in \O(\MM)$,
 \beaa
 X \gamma (Y_1,Y_2)= \gamma (\Ddot_X Y_1, Y_2)+ \gamma(Y_1, \Ddot_X Y_2).
 \eeaa
  \end{proposition}
 
  \begin{proof}
  Indeed,
  \beaa
 X \gamma (Y_1,Y_2)&=&X \g (Y_1,Y_2)=\g (\D_X Y_1,Y_2)+\g ( Y_1,\D_XY_2)= \g (\Ddot_X Y_1,Y_2)+\g ( Y_1,\Ddot_XY_2)\\
&=&  \gamma (\Ddot_X Y_1,Y_2)+\gamma ( Y_1,\Ddot_XY_2)
  \eeaa
  as desired.
   \end{proof}

We consider tensors  $\T_k (\MM)\otimes   \O_l (\MM)$, i.e. tensors  of the form  $U_{\nu_1\ldots \nu_k,  a_1\ldots a_l}$
for which we define,
\beaa
\Ddot_\mu U_{\nu_1\ldots \nu_k,  a_1\ldots a_l}&=& e_\mu U_{\nu_1\ldots \nu_k,  a_1\ldots a_l}-U_{\D_\mu\nu_1\ldots \nu_k,  a_1\ldots a_l}-\ldots- U_{\nu_1\ldots  \D_\mu\nu_k,  a_1\ldots a_l}\\
&-& U_{\nu_1\ldots \nu_k,   \Ddot_\mu a_1\ldots a_l}-  U_{\nu_1\ldots \nu_k,   a_1 \ldots \Ddot_\mu a_l}.
\eeaa

We are now ready   to prove the following.

\begin{proposition}
\lab{Proposition:commutehorizderivatives}
For a tensor $\Psi\in \O_1 (\MM)$, we   have  the curvature formula\footnote{With an immediate generalization to tensors $\Psi\in \O_l (\MM)$.}
 \bea
( \Ddot _\mu\Ddot_\nu -\Ddot_\nu\Ddot _\mu)\Psi_a=\Rdot_{a b  \mu\nu}\Psi^b
 \eea
 where, with  connection coefficients  $(\La_\a)_{\b\ga}= \g(\D_\a e_\ga, e_\b)$,
 \bea
 \lab{eq:DefineRdot}
 \bsplit
 \Rdot_{ab   \mu\nu}&:= \R_{ab    \mu\nu}+ \frac 1 2  \B_{ab   \mu\nu},\\
  \B_{ab   \mu\nu} &:=  (\La_\mu)_{3a} (\La_\nu)_{b4}+  (\La_\mu)_{4a} (\La_\nu)_{b3}- (\La_\nu)_{3a} (\La_\mu)_{b4}-  (\La_\nu)_{4a} (\La_\mu)_{b3}.
  \end{split}
 \eea
   More  generally, for a mixed tensor $\Psi\in \T_1(\MM)\otimes   \O_1(\MM)$, we have
  \beaa
( \Ddot _\mu\Ddot_\nu -\Ddot_\nu\Ddot _\mu)\Psi_{\la a}=    \R_{\la }\, ^   \si \,_{   \mu\nu}\Psi_{\si a}+    
            \Rdot_{a}\, ^   b\,_{   \mu\nu}\Psi_{\la b}
 \eeaa
 with an immediate generalization to tensors $\Psi\in\T_k (\MM)\otimes   \O_l (\MM)$. 
 \end{proposition}

\begin{proof}
We have
 \beaa
 \D_\mu e_a=\Ddot_\mu e_a  -\frac 1 2 \g(\D_\mu e_a, e_3) e_4 -
 \frac 1 2 \g(\D_\mu e_a, e_4) e_3= \Ddot_\mu e_a-\frac 1 2 (\La_\mu)_{3a} e_4 -\frac 1 2 (\La_\mu)_{4a} e_3. 
 \eeaa
 We deduce
\beaa
\D_\mu\Psi_a
&=&\Ddot_\mu\Psi_a -\frac 1 2  (\La_\mu)_{3a} \Psi_4 -\frac 1 2  (\La_\mu)_{4a} \Psi_3= \nab_\mu\Psi_a.
\eeaa
Hence
\beaa
\D_\mu \D_\nu \Psi_a &=& e_\mu( \D_\nu \Psi_a )- \D_{\D_\mu e_\nu}\Psi_a -\D_\nu\Psi_{\D_\mu e_a} \\
 &=&\nab_\mu\nab_\nu \Psi_a - \Ddot_{\D_\mu e_\nu}\Psi_a
 +\frac 1 2   (\La_\mu)_{3a}\D_\nu \Psi_4  +\frac 1 2   (\La_\mu)_{4a}\D_\nu \Psi_3.
\eeaa
On the other hand
\beaa
\D_\nu \Psi_4 &=& -(\La_\nu)_{b4}\Psi_b, \qquad 
\D_\nu \Psi_3 = -(\La_\nu)_{b 3}\Psi_b,
\eeaa
hence
\beaa
\D_\mu \D_\nu \Psi_a &=&\nab_\mu\nab_\nu \Psi_a  - \Ddot_{\D_\mu e_\nu}\Psi_a-\frac 1 2 (\La_\mu)_{3a} (\La_\nu)_{b4}\Psi_b- \frac 1 2 (\La_\mu)_{4a} (\La_\nu)_{b3}\Psi_b.
\eeaa
By symmetry
\beaa
\D_\nu \D_\mu \Psi_a&=& \nab_\nu \nab_\mu \Psi_a  - \Ddot_{\D_\nu e_\mu}\Psi_a -\frac 1 2 (\La_\nu)_{3a} (\La_\mu)_{b4}\Psi_b- \frac 1 2 (\La_\nu)_{4a} (\La_\mu)_{b3}\Psi_b.
\eeaa
Subtracting and using  the Ricci formula and the  Lemma above we  deduce
\beaa
\R_{ab\mu\nu}\Psi^b &=&[\nab_\mu, \nab_\nu] \Psi_a 
-\nab_{[e_\mu, e_\nu] }\Psi_a \\
 &&-\frac 1 2 (\La_\mu)_{3a} (\La_\nu)_{b4}\Psi_b- \frac 1 2 (\La_\mu)_{4a} (\La_\nu)_{b3}\Psi_b+\frac 1 2 (\La_\nu)_{3a} (\La_\mu)_{b4}\Psi_b+ \frac 1 2 (\La_\nu)_{4a} (\La_\mu)_{b3}\Psi_b\\
 &=&\Ddot_\mu \Ddot_\nu\Psi_a- \Ddot_\nu \Ddot_\mu\Psi_a   -\frac 1 2 \B_{ab\mu\nu} \Psi^b,
\eeaa
 from which the desired formula follows.
\end{proof}

\begin{remark}
Note that  the tensor  \B$_{ab\mu\nu}$ is anti-symmetric in both $\mu\nu $ and $ab$.
\end{remark}

\begin{corollary}\label{Corollory:Commutatornab_anab_b}
Let $X, Y$ be arbitrary vectorfields  on $\MM$ and $U\in \O_1(\MM)$ an  horizontal  tensor. We have\footnote{Here $(\Rdot(X, Y)U)_a:= X^\mu Y^\nu  \Rdot_{ab\mu\nu} U^b$.}
\beaa
\big(\nab_X \nab_Y -\nab_Y\nab_X) U=\nab_{[X, Y]} U+    \Rdot(X, Y) U
\eeaa
with an immediate generalization to $U\in \O_l (\MM)$.
\end{corollary}

 \begin{proof}
 We have
 \beaa
 \nab_Y \nab_X U_a&=&(Y^\la \Db_\la)( X^\mu \Db_\mu )U_a =   Y^\la X^\mu \Db_\la \Db_\mu U_a+ (Y^\la \Db_\la)( X^\mu ) \Db_\mu U_a,\\
  \nab_X \nab_Y U_a&=&  X^\mu  Y^\la  \Db_\mu  \Db_\la U_a+ (X^\mu  \Db_\mu )( Y^\la  ) \Db_\la U_a.
 \eeaa
 Hence,
 \beaa
 \big(\nab_X \nab_Y -\nab_Y\nab_X) U_a&=&Y^\la X^\mu\big(  \Db_\la \Db_\mu -  \Db_\mu\Db_\la\big)    U_a+ \big(   \Db_ X( Y^\mu ) -\Db_Y (X^\mu)\big)  \Db_\mu U_a\\
 &=&    X^\mu Y^\nu  \Rdot_{ab\mu\nu} U^b+\Db_{[X, Y]}  U_a,
 \eeaa
 as stated.
 \end{proof}


\subsection{Horizontal Hodge operators}
\lab{section:HorizontalHodgeoperators}
\label{section:nonintegrableHodge}


In this section we recall the Hodge operators on  $2$-spheres as defined in \cite{Ch-Kl} and extend their properties to the case of non-integrable horizontal structures.

We first define the following operators on horizontal tensors.

\begin{definition} 
For a given horizontal   $1$-form $\xi$, we  define the frame independent   operators
\beaa
\div\xi=\de^{ab}\nab_b\xi_a,\qquad 
\curl\xi=\in^{ab}\nab_a\xi_b,\qquad 
(\nab\hot \xi)_{ba}=\nab_b\xi_a+\nab_a  \xi_b-\de_{ab}( \div \xi).
\eeaa
\end{definition}

We collect here some Leibniz rules regarding the horizontal Hodge operators.
\begin{lemma}\label{Leibniz-rule-real} 
We have for $\xi, \eta \in \sk_1$, $u\in \sk_2$,
\beaa
 (\div \eta )\xi -(\curl \eta)  \dual \xi &=& \xi \c \nab \eta +\xi \c \dual \nab \dual \eta, \\
 \xi \hot (\div u) &=& \xi \c \nab u + \xi \c \dual \nab \dual u,\\
\xi \c ( \nab\hot  \eta)&=& \xi \c \nab f -\xi \c \dual \nab \dual \eta.
\eeaa
\end{lemma}

\begin{proof} 
Define  the  tensors on the left $Z_{a}= (\div \eta)\xi_a -(\curl \eta)  \dual \xi_a$,  $Y_{ab}=\xi \hot (\div u)_{ab}$ 
and  $W_{a}=\xi \c ( \nab\hot  \eta)_a=\xi_b ( \nab\hot  \eta)_{ab}$  and evaluate components. By simply  manipulating the definitions we find
\beaa
Z_1&=& \xi \c \nab \eta_1+( \xi \c \dual \nab)\dual  \eta_1,\\
Z_2&=& \xi \c \nab \eta_2+( \xi \c \dual \nab)\dual  \eta_2,\\
Y_{11}&=&( \xi \c \nab )u_{11}+( \xi \c \dual \nab) \dual u_{11},\\
Y_{12}&=&( \xi \c \nab) u_{12}+(\xi \c \dual \nab) \dual u_{12},\\
W_1&=&  \xi \c \nab \eta_1 - \xi \c\dual \nab \dual \eta_1,\\
W_2&=&  \xi \c \nab \eta_2 - \xi \c\dual \nab \dual \eta_2,
\eeaa
which implies the stated identities.
\end{proof}

\begin{definition}   
Given  an orthonormal basis of horizontal vectors $e_1, e_2$ we define 
   the Hodge type operators  (recall Definition \ref{definition-SS-real}),  as introduced in    \cite{Ch-Kl}.
\begin{itemize}
\item  $\DDd_1 $ takes $\sk_1$ into\footnote{Recall that $\sk_0$ refers to pairs of scalar functions.}  $\sk_0$:
\beaa
\DDd_1 \xi =(\div\xi, \curl \xi),
\eeaa
\item  $\DDd_2 $ takes $\sk_2$ into $\sk_1$:
\beaa
(\DDd_2 \xi)_a =\nab^b \xi_{ab},
\eeaa
\item $\DDs_1$ takes  $\sk_0 $ into  $\sk_1$:
\beaa
\DDs_1( f, f_*) &=& -\nab_a f+\in_{ab}  \nab_b  f_*,
\eeaa
\item 
$\DDs_2 $ takes  $\sk_1$ into $\sk_2$:
\beaa
\DDs_2 \xi&=& -\frac 1 2 \nab\hot \xi.
\eeaa
\end{itemize}   
\end{definition}

\begin{lemma}
Note the following pointwise identities:
\begin{enumerate}
\item Given $(f, f_*)\in \sk_0$, \, $u\in \sk_1$, we have
\bea\label{pointwise-hodge1}
\DDs_1(f, f_*) \c u = (f, f_*) \c \DDd_1 u-  \nab_a\big(  f  u^a  + f_*   (\dual u)^a\big).
\eea

\item Given $ f\in \sk_1$, \, $u\in \sk_2$,  we have
\bea\label{pointwise-hodge2}
(\DDs_2 f)\c u &=&  f\c( \DDd_2 u)- \nab_a\big( f_b u^{ab} \big).
\eea
\end{enumerate}
\end{lemma}

\begin{proof}
To check \eqref{pointwise-hodge2} we write
 \beaa
 ( \nab\hot   f) \c   u  &=&\big(\nab_a f_b+\nab_b f_a- \de_{ab} \div f \big) u_{ab} =2(\nab_a f_b) u_{ab}= 2\nab_a (u_{ab} f_b)-2(\div u) \c f
 \eeaa
 which immediately yields the second identity.   
\end{proof}

In the particular case  when the horizontal structure is  tangent to $2$-spheres $S$  these operators   are elliptic on $S$ and have  remarkable properties  discussed in Chapter 2 of \cite{Ch-Kl} which we recall below.


\subsubsection{Hodge operators on spheres} 
       

       The following results  were derived in Chapter 2 of \cite{Ch-Kl} in the context of general $2$-dimensional 
       compact   surfaces $S$ with strictly  positive Gauss curvature $K$ which we will refer from now on as a 2-sphere.       
\begin{lemma}
 Given  a $2$-sphere $S$, we have the following:
\begin{itemize}
\item[-]  The kernels of both $\DDd_1$ and $\DDd_2$ in $L^2(S)$ are 
trivial  while  the kernel of $\DDs_1$  consists  of pairs of constants   in $\sk_0$.

\item[-]   The operators   $\DDs_1$, resp. $\DDs_2$ are  the $L^2$ adjoints of
 $\DDd_1$, respectively $\DDd_2$.
 
 \item[-] The kernel of  $\DDs_2$   is  the space of conformal Killing vectorfields
on $S$.
\end{itemize}
Moreover the following identities hold true\footnote{Here $\lap_k:\sk_k \to\sk_k$, $k=0,1,2$,   is  defined   by $(\lap_k U)_A= \nab^a\nab_a U_A$.},  see \cite{Ch-Kl}:
\bea
\begin{split}
\label{eq:dcalident}
\DDs_1\DDd_1&=-\lap_1+K,\qquad\,\,\,\, \DDd_1 \DDs_1=-\lap_0,\\
\DDs_2 \DDd_2 &=-\f12\lap_2+K,\qquad \DDd_2\DDs_2=-\f12(\lap_1+K).
\end{split}
\eea
\end{lemma}

\begin{proof}
The statements about $L^2$-adjoints  follow immediately by integrating formulas \eqref{pointwise-hodge1}-\eqref{pointwise-hodge2} on $S$.
The formulas \eqref{eq:dcalident}  follow easily by  using the  definitions and commuting derivatives. See also the more general Lemma \ref{le:Bochner-nointegrable}. Note also that  for $\xi\in \sk_1$
\beaa
\DDs_2\, \, \xi  =-\frac 1 2 \Lie_\xi \ga
\eeaa
where $\ga$ denotes the induced   horizontal metric as in Definition \ref{def:proxyfirstandsecondfundameentalform}.
\end{proof}

As a simple consequence of \eqref{eq:dcalident} one derives the following $L^2$  estimates.
\begin{proposition} 
\label{prop:2D-hodge}
Let $(S,\ga)$ be
a compact manifold with Gauss curvature $K$.

{\bf i.}\quad The following identity holds for vectorfields  $  f  $ 
 on $S$:
\bea
\label{eq:hodgeident1}
\int_S\big(|\nab   f |^2+K|   f |^2\big)=\int_S\big( |\div  f  |^2+|\curl  f   |^2\big)=\int_S|\DDd_1   f   |^2.
\eea

{\bf ii.}\quad The following identity holds for symmetric, traceless,
 2-tensorfields   $  f $ 
 on $S$:
\bea
\label{eq:hodgeident2}
\int_S\big(|\nab    f  |^2+2K| f  |^2\big)=2\int_S |\div f |^2=2\int_S |\DDd_2   f  |^2.
\eea

{\bf iii.}\quad The following identity holds for pairs of functions $(f, f_*)$
 on $S$:
\bea
\label{eq:hodgeident3}
\int_S\big(|\nab  f |^2+|\nab f_* |^2\big)=\int_S |-\nab f +(\dual \nab f_*) |^2=\int_S|\DDs_1
(f ,  f_*)|^2.
\eea

{\bf iv.}\quad  The following identity holds for vectors $   f   $ on $S$:
\bea
\label{eq:hodgeident3*}
\int_S \big(|\nab  f   |^2-K|  f  |^2\big)=2\int_S|\DDs_2   f   |^2.
\eea
\end{proposition}

\begin{proof}
See Chapter 2 in \cite{Ch-Kl}.
\end{proof}


\subsubsection{Bochner identities in the non-integrable case}


  We extend the identities above to the case of non-integrable horizontal structure.
  \begin{lemma}\label{le:Bochner-nointegrable} 
  Given a general possibly non-integrable horizontal structure, the Hodge operators and the Laplacians are related by the following relations for $\xi\in \sk_1$ and $u\in \sk_2$:
\bea
\begin{split}
\label{eq:dcalident-non-integrable}
\DDs_1\DDd_1\xi&=-\lap_1\xi -\frac  1 2  \in_{ab} [ \nab_a, \nab_b] \dual \xi ,\\
\DDd_2\DDs_2 \xi &=-\f12\lap_1\xi +  \frac 1 4\in_{ab} [\nab_a, \nab_b]\dual \xi,\\
\DDs_2 \DDd_2 u &=-\f12\lap_2u -\frac 1 4\in_{ab} [ \nab_a, \nab_b] \dual u.
\end{split}
\eea
  \end{lemma}
  
  \begin{proof}  
  We check the last relation   by  evaluating the components of the tensor
$Y_{ab}:=( \DDs_2(\DDd_2u))_{ab}$
\beaa
Y_{ab}&=& -\frac 1 2 (\nab \hot (\DDd_2u))_{ab}=-\frac 1 2 (\nab_a(\DDd_2u)_b+\nab_b (\DDd_2u)_a -\de_{ab} (\div (\DDd_2u)))\\
&=&-\frac 1 2 (\nab_a\nab^c u_{bc}+\nab_b \nab^c u_{ac} -\de_{ab} \nab^c\nab^d u_{cd}).
\eeaa
 For $a=b=1$ we derive
\beaa
Y_{11}&=&-\frac 1 2 (\nab_1\nab^c u_{1c}+\nab_1 \nab^cu_{1c} -\de_{11} \nab^c\nab^d u_{cd})\\
&=&-\nab_1\nab_1 u_{11} -\nab_1\nab_2 u_{12} +\frac 1 2 ( \nab_1\nab_1 u_{11}+ \nab_2\nab_1 u_{21}+\nab_1\nab_2 u_{12}+\nab_2\nab_2 u_{22})\\
&=&- \frac 1 2( \nab_1\nab_1+\nab_2\nab_2) u_{11} +\frac 1 2 (  \nab_2\nab_1 u_{12}-\nab_1\nab_2 u_{12})
\eeaa
which gives
\beaa
Y_{11}&=& -\frac 1 2 \lap_2 u_{11}-\frac 1 2 [ \nab_1, \nab_2] u_{12}= -\frac 1 2 \lap_2 u_{11}-\frac 1 2 [ \nab_1, \nab_2] \dual u_{11}.
\eeaa
 For $a=1, b=2$ we derive
\beaa
Y_{12}&=&-\frac 1 2 (\nab_1\nab^c u_{2c}+\nab_2 \nab^c u_{1c} -\de_{12} \nab^c\nab^d u_{cd})\\
&=&-\frac 1 2 (\nab_1\nab_1 u_{21}+\nab_1\nab_2 u_{22}+\nab_2 \nab_1 u_{11}+\nab_2 \nab_2 u_{12} )\\
&=&-\frac 1 2 (\nab_1\nab_1 u_{12}+\nab_2 \nab_2 u_{12})+\frac 1 2 ( \nab_1\nab_2 u_{11}-\nab_2 \nab_1 u_{11})
\eeaa
which gives
\beaa
Y_{12}&=&-\frac 1 2 \lap_2 u_{12}+\frac 1 2 [ \nab_1,\nab_2] u_{11}=-\frac 1 2 \lap_2 u_{12}-\frac 1 2 [ \nab_1,\nab_2] \dual u_{12}.
\eeaa
Hence
\beaa
Y_{ab}&=& -\frac 1 2 \lap_2 u_{ab}-\frac 1 2 [ \nab_1, \nab_2] \dual u_{ab}
\eeaa
as stated. The other relations can be checked in the same manner.
  \end{proof}

Using the pointwise relations \eqref{pointwise-hodge1} and \eqref{pointwise-hodge2} and the above lemma, we can deduce the following pointwise version of the $L^2$ estimates of Proposition \ref{prop:2D-hodge}. 

\begin{proposition} 
\label{prop:2D-hodge-non-integrable}
Given a not necessarily integrable horizontal structure, the following pointwise relations hold:

{\bf i.}\quad The following identity holds for $f\in \sk_1$:
\bea
\label{eq:hodgeident1-nonint}
\bsplit
|\nab   f |^2-\frac{1}{2}\in_{ab}[\nab_a, \nab_b]\dual f \c f & = |\DDd_1   f   |^2\\
&+\nab_a \Big( \nab^a f \c f- (\div f) f^{a}-(\curl f) (\dual f )^a \Big).
\end{split}
\eea

{\bf ii.}\quad The following identity holds for $f\in \sk_2$:
\bea
\label{eq:hodgeident2-nonint}
|\nab    f  |^2- \frac{1}{4}\in_{ab}[\nab_a, \nab_b]\dual f \c f=2 |\DDd_2   f  |^2+\nab_a \Big( \nab^a f \c f-2 (\div f)_b f^{ab} \Big).
\eea

{\bf iii.}\quad  The following identity holds for $f\in\sk_1$:
\bea
\label{eq:hodgeident3-nonint}
 |\nab  f   |^2+\frac{1}{4}\in_{ab}[\nab_a, \nab_b]\dual f \c f=2|\DDs_2   f   |^2+\nab_a \Big( \nab^a f \c f+2 (\DDs_2 f)^{ab} f_{b} \Big).
\eea
\end{proposition}

\begin{proof}
The above relations are obtained by multiplying relations \eqref{eq:dcalident-non-integrable} by $f$ and integrating by parts in the horizontal directions. 
\end{proof}

\begin{remark}
 In      the integrable case the  commutator $\in^{ab}[\nab_a, \nab_b]$ is given by the standard Gauss      formula in terms of  $K$.   
     In the non-integrable case  it      can be computed by using the generalized Gauss equation, see Proposition \ref{prop-nullstr:0}. 
     \end{remark}

Observe that in the relations obtained in Proposition \ref{prop:2D-hodge-non-integrable}, the divergence terms cannot be discarded upon integration because of the absence of an integrable  surface.  There are various ways to  deal with this difficulty, such as   to integrate  \eqref{eq:hodgeident1-nonint}-\eqref{eq:hodgeident3-nonint} on the entire spacetime manifold    $\MM$.
 
\begin{remark}
\lab{remark:divergenceterms=Part1}
 Note that the  divergence terms  in Proposition \ref{prop:2D-hodge-non-integrable} can be re-expressed in terms of  spacetime divergences based on the following lemma.
\end{remark}
   
\begin{lemma}\label{lemma:divergence-spacetime-horizontal-ch1} 
For $f \in \sk_1$, we have\footnote{Here, we extend the horizontal 1-form $f$ as a full 1-form on $\MM$ by setting $f_3=f_4=0$.}
\bea
\D^\a f_\a&=& \nab^a f_a + (\eta +\etab) \c f
\eea
where  $\etab_a:= \frac 1 2 \g(e_a, \D_L\Lb)$ and $\eta_a:= \frac 1 2  \g(e_a, \D_\Lb L)$, see Definition \ref{def:Ricci-coefficients}. 
\end{lemma}

\begin{proof} 
We have, using \eqref{eq:expressions-Riccif-formula},
 \beaa
 \D^\a f_\a-  \nab^a f_a &=&-\frac 1 2\big( \D_3 f_4+\D_4 f_3) =-\frac 1 2\big( e_3( f_4)- f_{\D_3 4}+e_4( f_3) - f_{\D_43}) \\
 &=&\frac 1 2\big( 2\eta_a f_{a} + 2\etab_a f_{a})= (\eta+\etab)\c  f
 \eeaa
as stated. 
\end{proof}


 \subsection{The Gauss equation}
 

 Note that in the case of a non-integrable structure,  we  are missing the traditional Gauss equation which  connects the  Gauss curvature   of a  sphere   to  a Riemann curvature component. In what follows we state a result which is its non-integrable analogue. 
 
 \begin{proposition}
 \lab{prop-nullstr:0}
 The following identity holds true.
 \bea
\lab{Gauss-eq-horizontal}
\bsplit
 \nab_a \nab_b X_c- \nab_b \nab_a X_c&= \R_{c d   ab}X^d+\frac 1 2 \in_{ab}\Big(\atrch\nab_3+\atrchb \nab_4\Big) X_c
\\
 &-\frac 1 2 \Big(\chi_{ac}\chib_{bd} + \chib_{ac}\chi_{bd} - \chi_{bc}\chib_{ad}- \chib_{bc}\chi_{ad}\Big) X^d ,
 \end{split}
\eea
where $\R_{cdab}$ denotes the Riemann curvature of $(\MM, \g)$.
 \end{proposition}  
 
 \begin{proof}
Given the importance of the formula we give below a  direct proof of it.

 For   $X\in \sk_1$ we have,
\beaa
\D_b X_c&=&\nab_b X_c, \qquad 
\D_3X_c=\nab_3 X_c,\\
\D_4 X_c&=&\nab_4 X_c, \qquad 
\D_b X_3=-\chib_{bd} X_d, \qquad 
\D_b X_4=-\chi_{bd} X_d.
\eeaa
Also,
\beaa
\D_a \D_b X_c&=&\nab_a \nab_b X_c- \frac 1 2 \chi_{ab} \D_3 X_c -\frac 1 2 \chib_{ab}\D_4 X_c -\frac 1 2\chi_{ac} \D_b X_3-\frac 1 2 \chib_{ac} \D_b X_4\\
&=&\nab_a \nab_b X_c- \frac 1 2 \chi_{ab} \nab_3 X_c -\frac 1 2 \chib_{ab}\nab_4 X_c+\frac 1 2 \chi_{ac}\chib_{bd} X_d+\frac 1 2 \chib_{ac}\chi_{bd} X_d.
\eeaa
Hence,
\beaa
\D_a \D_b X_c&=&\nab_a \nab_b X_c- \frac 1 2 \chi_{ab} \nab_3 X_c -\frac 1 2 \chib_{ab}\nab_4 X_c+\frac 1 2 \chi_{ac}\chib_{bd} X_d+\frac 1 2 \chib_{ac}\chi_{bd} X_d,\\
\D_b\D_a X_c&=&\nab_b \nab_a X_c- \frac 1 2 \chi_{ba} \nab_3 X_c -\frac 1 2 \chib_{ba}\nab_4 X_c+\frac 1 2 \chi_{bc}\chib_{ad} X_d+\frac 1 2 \chib_{bc}\chi_{ad} X_d.
\eeaa
Subtracting we derive
\beaa
 \R_{c d   ab}X^d&=&\D_a \D_b X_c-\D_b\D_a X_c\\
 &=& \nab_a \nab_b X_c- \nab_b \nab_a X_c-\frac 1 2 (\chi_{ab}-\chi_{ba})\nab_3 X_c-\frac 1 2 (\chib_{ab}-\chib_{ba})\nab_4 X_c\\
&+&\frac 1 2 \Big(\chi_{ac}\chib_{bd}+ \chib_{ac}\chi_{bd} - \chi_{bc}\chib_{ad}- \chib_{bc}\chi_{ad}\Big) X^d. 
\eeaa
Thus,
\beaa
 \nab_a \nab_b X_c- \nab_b \nab_a X_c&=&\frac 1 2 (\chi_{ab}-\chi_{ba})\nab_3 X_c
 +\frac 1 2 (\chib_{ab}-\chib_{ba})\nab_4 X_c\\
 &-&\frac 1 2 \Big(\chi_{ac}\chib_{bd} + \chib_{ac}\chi_{bd} - \chi_{bc}\chib_{ad}- \chib_{bc}\chi_{ad}\Big) X^d + \R_{c d   ab}X^d.
\eeaa
Since
\beaa
\chi_{ab}-\chi_{ba}=\in_{ab} \atrch, \qquad \chib_{ab}-\chib_{ba}=\in_{ab} \atrchb,
\eeaa
this concludes the proof of the proposition.
\end{proof}

\begin{remark}
 We note  that \eqref{Gauss-eq-horizontal}  can be  derived  from 
  Corollary   \ref{Corollory:Commutatornab_anab_b}  according to which,
  relative to an arbitrary frame $e_\mu$,  
  \beaa
\big(\nab_\mu \nab_\nu -\nab_\mu\nab_\nu) X=\nab_{[e_\mu, e_\nu]} X+ 
   \Rdot(e_\mu, e_\nu ) X
\eeaa
with $\Rdot=\R+\frac 12 \B$ and $\B$   defined   in  \eqref{eq:DefineRdot}. The Gauss formula follows then easily   by  evaluating  the  components   $\B_{cd ab}$ of the tensor $\B$ and the term  $\nab_{[e_a, e_b]} X$.
\end{remark}

We now specialize the Gauss equation \eqref{Gauss-eq-horizontal} to tensors.
\begin{proposition}\label{Gauss-equation-2-tensors} 
The following identities   hold true.
\begin{enumerate}
\item  For a scalar $\psi$:
\bea
\lab{Gauss-eq-real-scalars}
[ \nab_a, \nab_b] \psi &=&\left( \frac 1 2 \big(\atrch\nab_3+\atrchb \nab_4\big) \psi\right)\in_{ab}.
\eea
\item The only non-vanishing component of $\B_{abcd}$
is given by 
\bea
\lab{eq:Non-vanishingcomponentsBabcd}
\B_{1212}=-\B_{1221}=\B_{2121}= -\frac 12  \trch \trchb-\frac 1 2 \atrch \atrchb+\chih \c \chibh.
\eea

\item 
For 
$\psi \in \sk_k$ for $k=1,2$,
\bea
\lab{Gauss-eq-real-2-tensors}\lab{gauss-real-Ka}\label{Gauss-eq-real-tensors}
[ \nab_a, \nab_b] \psi &=&\left( \frac 1 2 \big(\atrch\nab_3+\atrchb \nab_4\big) \psi +k  \, \Kh \dual \psi\right)  \in_{ab}
 \eea
 where
 \bea\label{eq:definition-K-in}
\Kh&:=&- \frac 14  \trch \trchb-\frac 1 4 \atrch \atrchb+\frac 1 2 \chih \c \chibh-  \frac 1 4 \R_{3434}.  
\eea
\end{enumerate}
\end{proposition}
      
\begin{proof} 
The case of scalars can be easily checked directly.

We consider  below the case   $\psi \in \sk_2$. 
From Corollary \ref{Corollory:Commutatornab_anab_b} applied to $\psi \in \sk_2$, we have
\beaa
\big( \nab_a \nab_b- \nab_b \nab_a\big)  \psi_{st} &=&\frac 1 2 \in_{ab}(\atrch\nab_3+\atrchb \nab_4) \psi_{st}
 +\frac  12 \B_{sdab} \psi_{dt}    +\frac 1 2 \B_{tdab} \psi_{sd} \\
    &+& \R_{sdab}\psi_{dt}+\R_{tdab} \psi_{sd}
\eeaa
where, by definition of $\B$ given in \eqref{eq:DefineRdot}, 
\bea\label{definition-E}
\B_{cdab}:&=&\chi_{bc}\chib_{ad}+\chib_{bc}\chi_{ad} - \chi_{ac}\chib_{bd} - \chib_{ac}\chi_{bd} . 
\eea
Note that by the symmetries of $\B$, all components  of $\B_{abcd}$ vanish  except for $\B_{1212} $. We have
\beaa
\B_{1212}&=& -\chi_{11}\chib_{22} - \chib_{11}\chi_{22} + \chi_{21}\chib_{12}+\chib_{21}\chi_{12}\\
&=&-\left(\frac 1 2 \trch + \chih_{11} \right)\left(\frac 1 2 \trchb + \chibh_{22} \right) - \left(\frac 1 2 \trchb + \chibh_{11} \right)\left(\frac 1 2 \trch + \chih_{22} \right)\\
&& + \left(-\frac 1 2\atrch +\chih_{21} \right)\left(\frac 1 2\atrchb +\chibh_{12} \right)+\left(-\frac 1 2\atrchb +\chibh_{21} \right)\left(\frac 1 2\atrch +\chih_{12} \right)\\
&=&- \frac 12  \trch \trchb-\frac 1 2 \atrch \atrchb-\chih_{11} \chibh_{22}-\chih_{22} \chibh_{11}+\chih_{21} \chibh_{12}+\chih_{12} \chibh_{21}\\
&=& -\frac 12  \trch \trchb-\frac 1 2 \atrch \atrchb+\chih \c \chibh.
\eeaa
 This implies for $\psi \in \sk_2$:
 \beaa
[ \nab_1, \nab_2] \psi &=&\frac 1 2 (\atrch\nab_3+\atrchb \nab_4) \psi\\
&& -\left( \frac 12  \trch \trchb+\frac 1 2 \atrch \atrchb-\chih\c\chibh+\frac 1 2 \R_{3434}  \right) \dual \psi
\eeaa
as stated. The case $\psi \in \sk_1$  can be treated in the same manner.
\end{proof}

\begin{remark}
 The quantity  $\Kh$  defined by 
 \eqref{eq:definition-K-in} 
  becomes the standard Gauss curvature in the case of an integrable structure.
  We note also  that the  value of $\Kh$ for the standard  non-integrable  structure (induced by the standard  principal null  directions, see  Chapter \ref{SECTION:KERR})
  of Kerr is given by the formula
  \beaa
\Kh= \frac{r^4+a^2r^2\sin^2\th-4ma^2r\cos^2\th-a^4\cos^2\th}{|q|^6}.
\eeaa
\end{remark}

Here is a more general version of Proposition \ref{Gauss-equation-2-tensors}.
\begin{proposition}
\lab{Gauss-equation-k-tensors} 
The following identity holds true for  any horizontal tensor $\psi\in \O_k$  and set of  horizontal indices  $I= i_1\ldots i_k$
\bea
\bsplit
[ \nab_a, \nab_b] \psi_I &=\left( \frac 1 2 \big(\atrch\nab_3+\atrchb \nab_4\big) \psi _I \right)  \in_{ab}\\
&+\Kh \Big[ \big( g_{i_1 a }  g_{tb} - g_{i_1 b}  g_{ta} \big) \psi^t\,_{i_2\ldots i_k}+\cdots  \big( g_{i_k a }  g_{tb} - g_{i_k b}  g_{ta} \big) \psi_{i_1\ldots }\,^ t\Big]
\end{split}
 \eea
 with  $\Kh$ given by \eqref{eq:definition-K-in}.
      \end{proposition}      
      
      \begin{proof}
      The proof  is  a simple extension of the proof of Proposition  \ref{Gauss-equation-2-tensors}, and is  left to the reader.
      \end{proof}

\begin{remark}
 Observe that in the case  when the horizontal structure is tangent to a   $S$-foliation, $\Kh$ reduces to the Gauss curvature of  $S$.  In the integrable case, for $k=1$, we can calculate directly\footnote{One can  check directly that  $g_{sa} \psi_b- g_{sb}\psi_a=\in_{ab}\dual \psi_s $.} on  any surface of integrability  $S$ with Gauss curvature $K$,
 \beaa
 \, [\nab_a, \nab_b] \psi_s&=&K\big( g_{sa} g_{tb}- g_{sb} g_{ta}\big)\psi^t =K\big(g_{sa} \psi_b- g_{sb}\psi_a\big)=K \in_{ab} \dual \psi_s
 \eeaa
 which coincides with formula \eqref{Gauss-eq-real-2-tensors} in this case. 
 \medskip
 Also for  $\psi\in\O_2$ (but not necessarily in $\sk_2$),
 \beaa
  \, [\nab_a, \nab_b] \psi_{s_1 s_2}&=&  K\big( g_{s_1 a}  g_{tb} - g_{s_1 b}  g_{ta} \big) \psi^t\,_{s_2}+ K \big( g_{s_2 a}  g_{tb} - g_{s_2 b}  g_{ta} \big)  \psi_{s_1}\,^t\\\
  &=& K\big(  g_{s_1 a} \psi_{b s_2} -g_{s_1 b} \psi_{as_2}  \big)+  K\big(  g_{s_2 a} \psi_{ s_1 b } -g_{s_2 b} \psi_{s_1 a}  \big).
 \eeaa
\end{remark}

Using \eqref{gauss-real-Ka}  we can rewrite  Proposition \ref{prop:2D-hodge-non-integrable} as follows.

\begin{proposition} 
\lab{prop:2D-hodge-non-integrable-pert-Kerr-div} 
Given a not necessarily integrable horizontal structure, the following pointwise relations hold\footnote{Note that according to   Lemma  \ref{lemma:divergence-spacetime-horizontal-ch1},    the  divergence terms  in the proposition can be re-expressed in terms of  the  spacetime divergences, see       Remark \ref{remark:divergenceterms=Part1}.}:

{\bf i.}\quad The following identity holds for $f\in \sk_1$:
\bea
\label{eq:hodgeident1-nonint-spacet-fin-div}
\begin{split}
|\nab   f |^2+ \Kh |f|^2&=|\DDd_1   f   |^2+\frac 1 2 \bigg(\Big(\atrch\nab_3+\atrchb \nab_4\Big) \dual f \bigg) \c f 
+\div[\DDd_1 f],
\\
 \div[\DDd_1 f]&:= \nab_a \Big( \nab^a f \c f-             (\div f) f^{a}-(\curl f) (\dual f )^a \Big).
\end{split}
\eea

{\bf ii.}\quad The following identity holds for $f\in \sk_2$:
\bea
\label{eq:hodgeident2-nonint-spacet-fin-div}
\begin{split}
|\nab    f  |^2+ 2 \Kh |f|^2 &=2 |\DDd_2   f  |^2+\frac 1 2 \bigg(\Big(\atrch\nab_3+\atrchb \nab_4\Big) \dual f \bigg) \c f+\div[\DDd_2f],\\
\div[\DDd_2 f]&:=\nab_a \big( \nab^a f \c f-2 (\div f)_b  f^{ab} \big).
\end{split}
\eea

{\bf iii.}\quad  The following identity holds for $f\in\sk_1$:
\bea
\label{eq:hodgeident3-nonint-spacet-fin-div}
\begin{split}
 |\nab  f   |^2-\Kh |f|^2&=2|\DDs_2   f   |^2-\frac 1 2 \bigg(\Big(\atrch\nab_3+\atrchb \nab_4\Big) \dual f \bigg) \c f
 +\div[\DDs_2 f],\\
 \div[\DDs_2 f]&:=\nab_a \Big( \nab^a f \c f+2 (\DDs_2 f)^{ab} f_{b} \Big).
 \end{split}
\eea
\end{proposition}

\begin{proof} 
From \eqref{gauss-real-Ka}, we have for $f \in \sk_1$ and $u \in \sk_2$:
\beaa
\frac{1}{2}\in_{ab} [ \nab_a, \nab_b] \dual f &=&\frac 1 2 (\atrch\nab_3+\atrchb \nab_4) \dual f 
 -\Kh f,\\
\frac{1}{2}\in_{ab}[ \nab_a, \nab_b] \dual u &=&\frac 1 2 (\atrch\nab_3+\atrchb \nab_4) \dual u 
 -2\Kh u,
\eeaa
from which we obtain the stated identities. 
\end{proof}


\subsection{Bochner identities  for the horizontal Laplacian}


\begin{proposition}
\lab{PROP:BOCHNERFOR-HORIZONTAL-LAP}
The following identities hold true. 
\begin{enumerate}
\item
Given a scalar function $\psi$ we have
\beaa
\big|\lap \psi|^2 &= \big|\nab^2 \psi\big|^2 +\Kh |\nab\psi|^2 +\err_0[\lap\psi]+\div_0[\lap\psi],
\eeaa
with
\beaa
\err_0[\lap \psi]&:=& -\frac 1 2 \nab \psi\c   \big(\atrch\nab_3+\atrchb \nab_4\big)\dual \nab\psi, \\
\div_0[\lap\psi]&:=& \nab_a\Big(\nab^a\psi\c  \lap\psi  - \frac 1 2 \nab^a|\nab\psi|^2 \Big).
\eeaa

\item For $\psi\in \sk_1$ we have 
\beaa
\big|\lap \psi|^2 &= |\nab^2 \psi|^2   +\Kh \big(|\nab\psi|^2-2\Kh |\psi|^2 \big) +\err_1[\lap \psi]+\div_1[\lap\psi],
\eeaa
with
\beaa
\err_1[\lap\psi]&:=&-\frac 1 2 \nab \psi\c   \big(\atrch\nab_3+\atrchb \nab_4\big)\dual \nab\psi +\frac 1 2 \Big| \big(\atrch\nab_3+\atrchb \nab_4 \big)\psi\Big|^2\\
&&-\frac 3 2  \Kh  
\big(\atrch\nab_3+\atrchb \nab_4\big) \psi\c \dual \psi +\Kh\div[\DDd_1\psi], 
\\
  \div_1[\lap\psi]  &:=&\nab_a\Big(\nab^a\psi\c  \lap\psi  -  \nab_c \psi \c \nab^c\nab^a\psi\Big).
\eeaa

\item   For $\psi\in \sk_2$ we have 
\bea
\lab{eq:BochnerLaplacian}
\big|\lap \psi|^2 &=& |\nab^2 \psi|^2   +\Kh\Big( |\nab\psi|^2- 6 \Kh |\psi|^2\Big) +\err_2[\lap\psi]+\div_2[\lap\psi],
\eea
where
\beaa
\err_2[\lap \psi]&:=&-\frac 1 2 \nab \psi\c   \big(\atrch\nab_3+\atrchb \nab_4\big)\dual \nab\psi +\frac 1 2 \Big| \big(\atrch\nab_3+\atrchb \nab_4 \big)\psi\Big|^2\\
&&- 3 \Kh  
\big(\atrch\nab_3+\atrchb \nab_4\big) \psi\c \dual \psi + 2\Kh\div[\DDd_2\psi], \\
\div_{2}[\lap\psi]&:=&\nab_a\Big(\nab^a\psi\c\lap\psi -\nab_c\psi \c \nab^c\nab^a\psi\Big).
\eeaa
\end{enumerate}
\end{proposition} 

\begin{proof}
See Appendix \ref{subsection:appendix-prop-lap2}.
\end{proof}

\begin{remark}
In the  integrable case,  the horizontal structure  is tangent to   spheres  $S$, and  we derive the following by integration:
\begin{enumerate}
\item  Given a scalar function $\psi$ we have
 \bea
 \lab{eq:remark.Bochnerforsk_0}
\int_S \big|\lap \psi|^2 &=&\int_S \big|\nab^2 \psi\big|^2 +\int_SK |\nab\psi|^2.
\eea
\item For $\psi\in \sk_1$ we have 
\bea
 \lab{eq:remark.Bochnerforsk_1}
\int_S \big|\lap \psi|^2 =\int_S  |\nab^2 \psi|^2   +K \big(|\nab\psi|^2-2K |\psi|^2 \big) +\int_S  K \div[\DDd_1\psi] 
\eea
with
\beaa
 \div[\DDd_1 \psi]&= \nab_a \Big( \nab^a \psi \c \psi-             (\div \psi \psi ^{a}-\curl \psi  \dual \psi^a) \Big).
 \eeaa
 \item    For $\psi\in \sk_2$ we have 
\bea
 \lab{eq:remark.Bochnerforsk_2}
\int_S \big|\lap \psi|^2 &=&\int_S  |\nab^2 \psi|^2   + \int_S K \Big( |\nab\psi|^2- 6 K  |\psi|^2\Big) +\int_S 2 K \div[\DDd_2\psi]
\eea
with
\beaa
\div[\DDd_2 \psi]&=\nab_a \Big( \nab^a \psi \c \psi-2 (\div \psi)_b  \psi^{ab} \Big).
 \eeaa
\end{enumerate}
\end{remark}


\section{Horizontal  structures and Einstein equations}\label{section:einstein-equation}


We apply the general formalism for non-integrable structures to the case of a spacetime solution to the Einstein vacuum equation. For an application of the formalism to the non-vacuum case, such as the   Einstein-Maxwell equation, see \cite{Giorgi} \cite{Giorgi2}.


 \subsection{Ricci coefficients}
 \label{subsection:Ricci-coefficients}
 
 
 \begin{definition}\label{def:Ricci-coefficients}
 We define the following horizontal  $1$-forms
 \beaa
 \etab(X)&:=& \frac 1 2 \g(X, \D_L\Lb),\qquad \eta(X):= \frac 1 2  \g(X, \D_\Lb L),\\
 \xib(X)&:=& \frac 1 2  \g(X, \D_\Lb \Lb),\qquad \xi(X):= \frac 1 2 \g(X, \D_L L).
 \eeaa
 With these definitions we have
 \beaa
 \nab_L X&:=&^{(h)}(\D_L X)=\D_L X-\etab(X)L -\xi(X) \Lb,\\
 \nab_\Lb X&:=&^{(h)}(\D_\Lb X)=\D_\Lb X-\xib(X)L -\eta(X) \Lb. 
 \eeaa
 \end{definition}

 In addition to the horizontal tensor-fields $\chi,\chib,\etab, \eta, \xi,\xib$ introduced above,  
 we also define the scalars
 \beaa
 \omb&:=& \frac 1  4 \g(\D_\Lb\Lb, L),\qquad\quad  \om:=\frac 1 4 \g(\D_L L, \Lb),
 \eeaa
 and the horizontal 1-form
 \beaa
 \ze(X)&=&\frac 1 2 \g(\D_XL,\Lb).
 \eeaa
 We summarize below the definition of the the horizontal $1$-forms $\xi, \xib, \eta, \etab, \ze\in\mathbf{O}_1$:
\begin{equation}\label{fo2}
\begin{cases}
&\xi(X)=\frac 1 2 \g(\D_LL,X),\quad \xib(X) =\frac 1 2 \g(\D_{\Lb}\Lb,X),\\
&\eta(X)=\frac 1 2 \g(\D_{\Lb}L,X),\quad \etab(X)=\frac 1 2 \g(\D_L\Lb,X),\\
&\ze(X)=\frac 1 2 \g(\D_X L,\Lb),
\end{cases}
\end{equation}
and the real scalars
\begin{equation}\label{fo3}
\om=\frac 1 4 \g(\D_LL,\Lb),\qquad\omb=\frac 1 4 \g(\D_{\Lb}\Lb,L).
\end{equation}
 
 \begin{definition}
 \lab{defin:Ricci-coefficients}
 The horizontal tensor-fields $\chi,\chib, \eta, \etab, \ze,  \xi,\xib,\om, \omb$ are  called
 the connection coefficients of the null pair $(L,\Lb)$. Given  an  arbitrary  basis of 
  horizontal vectorfields $e_1,  e_2$, we write using the short hand notation $\D_a=\D_{e_a}, a=1,2$, 
  \beaa
\chib_{ab}&=&\g(\D_a\Lb, e_b),\qquad\,\,\,\, \,\chi_{ab}=\g(\D_aL, e_b),\\
\xib_a&=&\frac 1 2 \g(\D_\Lb\Lb, e_a),\qquad\,\,\, \xi_a=\frac 1 2 \g(\D_L L, e_a),\\
\omb&=&\frac 1 4 \g(\D_\Lb\Lb, L),\qquad\quad  \om=\frac 1 4 \g(\D_L L, \Lb),\qquad \\
\etab_a&=&\frac 1 2 \g (\D_L\Lb, e_a),\qquad \quad \eta_a=\frac 1 2 \g(\D_\Lb L, e_a),\qquad\\
 \ze_a&=&\frac 1 2 \g(\D_{a}L,\Lb).
\eeaa
 \end{definition}
 
We easily  derive the  Ricci formulae,
\bea\label{eq:expressions-Riccif-formula}
\D_a e_b&=&\nab_a e_b+\frac 1 2 \chi_{ab} e_3+\frac 1 2  \chib_{ab}e_4,\nn\\
\D_a e_4&=&\chi_{ab}e_b -\ze_a e_4,\nn\\
\D_a e_3&=&\chib_{ab} e_b +\ze_ae_3,\nn\\
\D_3 e_a&=&\nab_3 e_a +\eta_a e_3+\xib_a e_4,\nn\\
\D_3 e_3&=& -2\omb e_3+ 2 \xib_b e_b,\label{ricci}\\
\D_3 e_4&=&2\omb e_4+2\eta_b e_b,\nn\\
\D_4 e_a&=&\nab_4 e_a +\etab_a e_4 +\xi_a e_3,\nn\\
\D_4 e_4&=&-2 \om e_4 +2\xi_b e_b,\nn\\
\D_4 e_3&=&2 \om e_3+2\etab_b e_b.\nn
\eea


\subsection{Curvature and Weyl fields}\label{subsection:curvature}


Assume that $W\in\mathbf{T}_4^0(\MM)$ is a Weyl field, i.e.
\begin{equation}\label{fo4}
\begin{cases}
&W_{\al\be\mu\nu}=-W_{\be\al\mu\nu}=-W_{\al\be\nu\mu}=W_{\mu\nu\al\be},\\
&W_{\al\be\mu\nu}+W_{\al\mu\nu\be}+W_{\al\nu\be\mu}=0,\\
&\g^{\be\nu}W_{\al\be\mu\nu}=0.
\end{cases}
\end{equation}

We define the null components of the Weyl field $W$, $\al(W),\aa(W),\varrho(W)\in \mathbf{O}_2(\MM)$ and $\be(W),\bb(W)\in\mathbf{O}_1(\MM)$ by the formulas
\begin{equation}\label{fo5}
\begin{cases}
\al(W)(X,Y)=W(L,X,L,Y),\\
\aa(W)(X,Y)=W(\Lb,X,\Lb,Y),\\
\b(W)(X)=\frac 1 2 W(X,L,\Lb,L),\\
\bb(W)(X)=\frac 1 2 W(X,\Lb,\Lb, L),\\
\varrho(W)(X,Y)= W(X,\Lb,Y,L).
\end{cases}
\end{equation}

Recall that if $W$ is a Weyl field its 
Hodge dual $\dual W$, defined by 
${}^{\ast}W_{\al\be\mu\nu}=\frac{1}{2}{\in_{\mu\nu}}^{\rho\si}W_{\al\be\rho\si}$,  is also a Weyl field. We easily
check the formulas,
\begin{equation}
\label{eq:dualW}
\begin{cases}
&\aa(\dual W)=\dual \aa(W),\qquad \a(\dual W)=-
\dual \a(W), \\
&\bb(\dual W)=\dual\bb(W),\qquad \b(\dual W)=-\dual \b(W),\\
&\varrho(\dual W)=\dual \varrho(W). 
\end{cases}
\end{equation}

It is easy to check that $\a,\aa$ are symmetric traceless  horizontal tensor-fields.
 On the other hand    the horizontal 2-tensorfield  $\varrho$ is   neither symmetric nor 
 traceless.  It is convenient to express it in terms
 of  the following  two scalar quantities
 \bea
 \label{fo5'}
 \rho(W)=\frac 1 4 W(L,\Lb,L,\Lb),\qquad \dual\rho(W)=
\frac 1 4  \dual W(L,\Lb,L,\Lb)\label{rho-dualrho}.
 \eea
 Observe also that,
 \beaa
\rho(\dual W)=\rhod(W), \qquad \rhod(\dual W)=-\rho.
\eeaa  
Thus,
\bea
\varrho(X,Y)=\big(-\rho\,\ga(X,Y)+\rhod\, \in(X,Y)\big),\qquad 
\forall\, X,Y\in \O(\MM).
\eea
We have
 \beaa
 W_{a3b4}&=&\varrho_{ab}=(-\rho\de_{ab} +\dual \rho \in_{ab}),\\
 W_{ab34}&=& 2 \in_{ab}\dual\rho,\\
 W_{abcd}&=&-\in_{ab}\in_{cd}\rho,\\
 W_{abc3}&=&\in_{ab}\dual \bb_c,\\
 W_{abc4}&=&-\in_{ab}\dual \b_c.
 \eeaa

\begin{remark}
\label{rem:pairing}
In addition to the Hodge duality we  will need to take into account
 the duality with respect to the interchange of $L, \Lb$, which we call a pairing transformation.  Clearly, under this transformation,   $\a\leftrightarrow \aa$, 
 $\b \leftrightarrow -\bb$,   $\rho   \leftrightarrow \rho$,  $\dual \rho   \leftrightarrow -\dual \rho$,
  $\varrho \leftrightarrow \check{\varrho}$ with 
  $\check{\varrho}_{ab}:=\varrho_{ba}$.    One has to be careful however when combining 
  the Hodge dual  and pairing  transformations. In that case we have,  
   $\dual \aa    \leftrightarrow  -\dual \a$,   $\dual\bb \leftrightarrow \dual \b$. 
   This is due to the  fact that under the pairing transformation 
$\in_{ab}\to-\in_{ab}$ (since $\in_{ab}=\in_{ab34}$).  Indeed, for example,
   \beaa
   \dual\aa_{ab}&=&  \aa(\dual W)_{ab}    = \dual W_{a3b3} =  - \in_{a3c4} W_{c3 b3}=\in_{ac34}W_{c3 b3}=\in_{ac}\aa_{cb},\\
      \dual\a_{ab}&=&  \a(\dual W)_{ab}    = \dual W_{a4b4} =  - \in_{a4c3} W_{c4 b4}=-\in_{cb34}W_{c4 b4}=
   -   \in_{ac}\a_{cb}.
   \eeaa
   \end{remark}

The decomposition above for Weyl fields applies in particular to the Riemann curvature tensor  $\R$ of a vacuum spacetime.

 In the case of a vacuum spacetime, the non-integrable Gauss curvature defined by \eqref{eq:definition-K-in} becomes
\bea\lab{eq:definition-K}
\Kh&=&- \frac 14  \trch \trchb-\frac 1 4 \atrch \atrchb+\frac 1 2 \chih \c \chibh-  \rho .
\eea


\subsection{Horizontal  tensor  $\B$}


We  calculate below the components of the  horizontal   curvature tensor $\B$ defined
 by the formula, see \eqref{eq:DefineRdot},
 \beaa
  \B_{ab   \mu\nu} &:=  (\La_\mu)_{3a} (\La_\nu)_{b4}+  (\La_\mu)_{4a} (\La_\nu)_{b3}- (\La_\nu)_{3a} (\La_\mu)_{b4}-  (\La_\nu)_{4a} (\La_\mu)_{b3}.
 \eeaa
\begin{proposition}
\lab{proposition:componentsofB}
The components of $\B$   are given   by the following formulas:
\bea
\bsplit
\B_{ a   b  c 3} &=       2 \big( -  \chib_{ca} \eta_b +  \chib_{cb} \eta_a -  \chi_{ca} \xib_b+  \chi_{cb} \xib_a\big), \\
\B_{ a   b  c 4} &=2\big( -\chi_{ca} \etab_b +\chi_{cb}  \etab_a -\chib_{ca} \xi_b+ \chib_{cb} \xi_a\big),\\
\B_{ a   b  3 4} &=4\big(-\xib_a \xi_b+\xi_a \xib_b-\eta_a \etab_b+\etab_a\eta_b\big),\\
\B_{ab cd} &= \chi_{bc}\chib_{ad}+ \chib_{bc}\chi_{ad}  -\chi_{ac}\chib_{bd} - \chib_{ac}\chi_{bd} . 
\end{split}
\eea
The above can also be written as 
\bea
\begin{split}
\B_{ a   b  c 3}&=     -  \trchb  \big( \de_{ca}\eta_b-  \de_{cb} \eta_a\big)  -  \atrchb \big( \in_{ca}  \eta_b -  \in_{cb}  \eta_a\big) \\
&+ 2 \big(- \chibh_{ca}  \eta_b + \chibh_{cb} \eta_a-  \chi_{ca} \xib_b+  \chi_{cb} \xib_a\big),\\
\B_{ a   b  c 4}&=     -  \trch  \big( \de_{ca}\etab_b-  \de_{cb} \etab_a\big)  -  \atrch \big( \in_{ca}  \etab_b -  \in_{cb}  \etab_a\big) \\
&+ 2 \big(- \chih_{ca}  \etab_b + \chih_{cb} \etab_a-  \chib_{ca} \xi_b+  \chib_{cb} \xi_a\big).
\end{split}
\eea
Also, $\B_{ab cd}$ is given by 
\beaa
\B_{abcd} = \left(- \frac 12  \trch \trchb-\frac 1 2 \atrch \atrchb+\chih \c \chibh\right)\in_{ab}\in_{cd}.
\eeaa
\end{proposition}

\begin{proof}
We write recalling the definition  $(\La_\mu)_{\a\b}=\g(\D_\mu e_\b, e_\a)$ and   definition of Ricci  coefficients, see Definition         \ref{defin:Ricci-coefficients},
\beaa
\B_{ a   b  c 3}&=&(\La_c)_{3a} (\La_3)_{b4}+  (\La_c)_{4a} (\La_3)_{b3}- (\La_3)_{3a} (\La_c )_{b4}-  (\La_3)_{4a} (\La_ c)_{b3}\\
&=&-2\chib_{ca} \eta_b -2\chi_{ca} \xib_b+2\xib_a \chi_{cb}+2\eta_a\chib_{cb}
\eeaa
and 
\beaa
 \B_{ab  34} &= & (\La_3)_{3a} (\La_4)_{b4}+  (\La_3)_{4a} (\La_4)_{b3}- (\La_4)_{3a} (\La_3)_{b4}-  (\La_4)_{4a} (\La_3)_{b3}\\
 &=&4(-\xib_a)\xi_b+4(- \eta_a) \etab_b -4(-\etab_a)\eta_b- 4(\etab_a) \eta_b-(-\xi_a) \xib_b\\
 &=&4\big(-\xib_a \xi_b+\xi_a \xib_b-\eta_a \etab_b+\etab_a\eta_b\big).
\eeaa
For the remaining formulas see \eqref{definition-E} and \eqref{eq:Non-vanishingcomponentsBabcd}.
\end{proof}


\subsection{Connection to  the Newman-Penrose formalism} 
\lab{section:NPformalism}


In the Newman-Penrose NP  formalism,  one chooses a   specific orthonormal  basis of  horizontal  vectors $(e_1, e_2)$ and defines  all connection coefficients relative   to the  complexified  frame $(n, l, m, \ov{m})$ where $n=\frac 1 2 e_3$, $l=e_4$,  $m=e_1+i e_2$, $\ov{m}= e_1-i e_2.$ Thus, all quantities of interest are complex scalars instead of our horizontal  tensors such as  $\sk_1, \sk_2$.   The NP  formalism works well  for deriving  the basic equations, but has the disadvantage 
of  substantially  increasing  the number of  variables.    Moreover, the calculations become far  more    cumbersome when deriving equations involving higher   derivatives   of the main quantities, in perturbations of Kerr.
 Another  advantage of  the formalism used here is that all important equations look similar to the ones  in \cite{Ch-Kl}. We refer  to \cite{NP}  for the original form of the NP  formalism.
 
The formalism used here is also related to the so-called Geroch-Held-Penrose formalism GHP formalism, which also introduced derivatives with boost weights, which are the scalar equivalent of the conformal derivatives used here, see Lemma \ref{lemma:definition-conformal-derivatives}. Nevertheless, the GHP formalism still involves complex scalars instead of horizontal tensors\footnote{It also introduces  spin weights  which partially  cure   the dependance  of these scalars on the choice of frames.}. We refer to \cite{GHP} for the original form of the GHP formalism.


\subsection{Null  structure equations}


We state below   the  null structure equation in the general setting  discussed above.   We assume given a  vacuum spacetime endowed with  a general null frame $(e_3, e_4, e_1, e_2)$ relative to which we define our  connection and curvature coefficients.
 
\begin{proposition}[Null structure equations] 
The  connection coefficients  verify the following   equations:
\label{prop-nullstr}
\beaa
\nab_3\trchb&=&-|\chibh|^2-\frac 1 2 \big( \trchb^2-\atrchb^2\big)+2\div\xib  - 2\omb \trchb +  2 \xib\c(\eta+\etab-2\ze),\\
\nab_3\atrchb&=&-\trchb\atrchb +2\curl \xib -2\omb\atrchb+ 2 \xib\wedge(-\eta+\etab+2\ze),\\
\nab_3\chibh&=&-\trchb\,  \chibh+  \nab\hot \xib- 2 \omb \chibh+    \xib\hot(\eta+\etab-2\ze)-\aa,
\eeaa
\beaa
\nab_3\trch
&=& -\chibh\c\chih -\frac 1 2 \trchb\trch+\frac 1 2 \atrchb\atrch    +   2   \div \eta+ 2 \omb \trch + 2 \big(\xi\c \xib +|\eta|^2\big)+ 2\rho,\\
\nab_3\atrch
&=&-\chibh\wedge\chih-\frac 1 2(\atrchb \trch+\trchb\atrch)+ 2 \curl \eta + 2 \omb \atrch + 2 \xib\wedge\xi  -  2 \dual \rho,\\
\nab_3\chih
&=&-\frac 1 2 \big( \trch \chibh+\trchb \chih\big)-\frac 1 2 \big(-\dual \chibh \, \atrch+\dual \chih\,\atrchb\big)
+ \nab\hot \eta +2 \omb \chih\\
&+& \xib\hot\xi +\eta\hot\eta,
\eeaa
\beaa
\nab_4\trchb
&=& -\chih\c\chibh -\frac 1 2 \trch\trchb+\frac 1 2 \atrch\atrchb    +  2   \div \etab+ 2 \om \trchb + 2\big( \xi\c \xib +|\etab|^2\big)+2\rho,\\
\nab_4\atrchb
&=&-\chih\wedge\chibh-\frac 1 2(\atrch \trchb+\trch\atrchb)+ 2 \curl \etab + 2 \om \atrchb + 2 \xi\wedge\xib+2 \dual \rho,\\
\nab_4\chibh
&=&-\frac 1 2 \big( \trchb \chih+\trch \chibh\big)-\frac 1 2 \big(-\dual \chih \, \atrchb+\dual \chibh\,\atrch\big)
+\nab\hot \etab +2 \om \chibh\\
&+& \xi\hot\xib + \etab\hot\etab,
\eeaa
\beaa
\nab_4\trch&=&-|\chih|^2-\frac 1 2 \big( \trch^2-\atrch^2\big)+ 2 \div\xi  - 2 \om \trch + 2   \xi\c(\etab+\eta+2\ze),\\
\nab_4\atrch&=&-\trch\atrch + 2 \curl \xi - 2 \om\atrch+ 2 \xi\wedge(-\etab+\eta-2\ze),\\
\nab_4\chih&=&-\trch\,  \chih+ \nab\hot \xi- 2 \om \chih+    \xi\hot(\etab+\eta+2\ze)-\a.
\eeaa
Also,
\beaa
\nab_3 \ze+2\nab\omb&=& -\chibh\c(\ze+\eta)-\frac{1}{2}\trchb(\ze+\eta)-\frac{1}{2}\atrchb(\dual\ze+\dual\eta)+ 2 \omb(\ze-\eta)\\
&&+\hch\c\xib+\frac{1}{2}\trch\,\xib+\frac{1}{2}\atrch\dual\xib +2 \om \xib -\bb,
\\
\nab_4 \ze -2\nab\om&=& \chih\c(-\ze+\etab)+\frac{1}{2}\trch(-\ze+\etab)+\frac{1}{2}\atrch(-\dual\ze+\dual\etab)+2 \om(\ze+\etab)\\
&& -\chibh\c\xi -\frac{1}{2}\trchb\,\xi-\frac{1}{2}\atrchb\dual\xi -2 \omb \xi -\b,
\\
\nab_3 \etab -\nab_4\xib &=& -\chibh\c(\etab-\eta) -\frac{1}{2}\trchb(\etab-\eta)+\frac{1}{2}\atrchb(\dual\etab-\dual\eta) -4 \om \xib  +\bb, \\
\nab_4 \eta    -    \nab_3\xi &=& -\chih\c(\eta-\etab) -\frac{1}{2}\trch(\eta-\etab)+\frac{1}{2}\atrch(\dual\eta-\dual\etab)-4\omb \xi -\b,\\
\eeaa
and
\beaa
\nab_3\om+\nab_4\omb -4\om\omb -\xi\c \xib -(\eta-\etab)\c\ze +\eta\c\etab&=&   \rho.
\eeaa
Also,
\beaa
\div\chih +\ze\c\chih &=& \frac{1}{2}\nab\trch+\frac{1}{2}\trch\ze -\frac{1}{2}\dual\nab\atrch-\frac{1}{2}\atrch\dual\ze -\atrch\dual\eta-\atrchb\dual\xi -\b,\\
\div\chibh -\ze\c\chibh &=& \frac{1}{2}\nab\trchb-\frac{1}{2}\trchb\ze -\frac{1}{2}\dual\nab\atrchb+\frac{1}{2}\atrchb\dual\ze -\atrchb\dual\etab-\atrch\dual\xib +\bb,
\eeaa
and\footnote{Note that this equation follows from  expanding  $\R_{34ab}$.}
\beaa
\curl\ze&=&-\frac 1 2 \chih\wedge\chibh   +\frac 1 4 \big(  \trch\atrchb-\trchb\atrch   \big)+\om \atrchb -\omb\atrch+\dual \rho.
\eeaa
\end{proposition}

\begin{proof}
 Except for  the fact that the order of indices  in $\chi, \chib$ is important, since they are no longer symmetric,    the  derivation is   exactly as in section 7.4 of  \cite{Ch-Kl}.  
 \end{proof}


\subsection{Null Bianchi identities}


We state  below the  equations verified by the null curvature components of an Einstein vacuum space-time.
 
    \begin{proposition}[Null Bianchi identities]\label{prop:bianchi} 
       The  curvature components  verify the following   equations:
    \beaa
    \nab_3\a-  \nab\hot \b&=&-\frac 1 2 \big(\trchb\a+\atrchb\dual \a)+4\omb \a+
  (\ze+4\eta)\hot \b - 3 (\rho\chih +\rhod\dual\chih),\\
\nab_4\beta - \div\a &=&-2(\trch\beta-\atrch \dual \b) - 2  \om\b +\a\c  (2 \ze +\etab) + 3  (\xi\rho+\dual \xi\rhod),\\
     \nab_3\b+\div\varrho&=&-(\trchb \b+\atrchb \dual \b)+2 \omb\,\b+2\bb\c \chih+3 (\rho\eta+\rhod\dual \eta)+    \a\c\xib,\\
 \nab_4 \rho-\div \b&=&-\frac 3 2 (\trch \rho+\atrch \rhod)+(2\etab+\ze)\c\b-2\xi\c\bb-\frac 1 2 \chibh \c\a,\\
   \nab_4 \rhod+\curl\b&=&-\frac 3 2 (\trch \rhod-\atrch \rho)-(2\etab+\ze)\c\dual \b-2\xi\c\dual \bb+\frac 1 2 \chibh \c\dual \a, \\
     \nab_3 \rho+\div\bb&=&-\frac 3 2 (\trchb \rho -\atrchb \rhod) -(2\eta-\ze) \c\bb+2\xib\c\b-\frac{1}{2}\chih\c\aa,
 \\
   \nab_3 \rhod+\curl\bb&=&-\frac 3 2 (\trchb \rhod+\atrchb \rho)- (2\eta-\ze) \c\dual \bb-2\xib\c\dual\b-\frac 1 2 \chih\c\dual \aa,\\
     \nab_4\bb-\div\varoc&=&-(\trch \bb+ \atrch \dual \bb)+ 2\om\,\bb+2\b\c \chibh
    -3 (\rho\etab-\rhod\dual \etab)-    \aa\c\xi,\\
     \nab_3\bb +\div\aa &=&-2(\trchb\,\bb-\atrchb \dual \bb)- 2  \omb\bb-\aa\c(-2\ze+\eta) - 3  (\xib\rho-\dual \xib \rhod),\\
     \nab_4\aa+ \nab\hot \bb&=&-\frac 1 2 \big(\trch\aa-\atrch\dual \aa)+4\om \aa+
 (\ze-4\etab)\hot \bb - 3  (\rho\chibh -\rhod\dual\chibh).
\eeaa
Here,
\beaa
\div\varo&=&- (\nab\rho+\dual\nab\rhod),\\
\div\varoc&=&- (\nab\rho-\dual\nab\rhod).
\eeaa
    \end{proposition} 
    
    \begin{proof}
    The proof  follows line by line from the derivation  in  section 7.3 of \cite{Ch-Kl} except, once more, for keeping track of the lack of symmetry for $\chi, \chib$. Note also  that $\varoc_{ab}=\varo_{ba}$ and  that $ (\div\varo)_b=\nab^a\varo_{ab}$.
    \end{proof}

     
\subsection{Commutation formulas}
\lab{sec:commutationformulasfirsttime}


\begin{lemma}
   \lab{LEMMA:COMM-GEN-B}
Let $U_{A}= U_{a_1\ldots a_k} $ be a general $k$-horizontal  tensorfield.
\begin{enumerate}
\item  We have
\bea
\bsplit
\,[\nab_3, \nab_b] U_A &=- \chib_{bc} \nab_c U_A+( \eta_b-\ze_b) \nab_3 U_A +\xib_b \nab_4 U_A \\
&+\sum_{i=1}^k\Big(-\in_{a_i c} \dual\bb_b +  \frac 1 2 \B_{a_i c 3b} \Big) U_{a_1\ldots }\,^ c \,_{\ldots a_k}. 
\end{split}
\eea

\item We have
\bea
\bsplit
\,[\nab_4, \nab_b] U_A &=- \chi_{bc} \nab_c U_a+( \etab_b+\ze_b) \nab_4 U_a +\xi_b \nab_3 U_a \\
&+\sum_{i=1}^k\Big(\in_{a_i c} \dual\b_b +  \frac 1 2 \B_{a_i c 4b} \Big) U_{a_1\ldots }\,^ c \,_{\ldots a_k}.
\end{split}
\eea

\item We have
\bea
\bsplit
\, [\nab_4, \nab_3] U_A&= 2(\etab_b-\eta_b ) \nab_b U_A + 2 \om \nab_3 U_A -2\omb \nab_4 U_A\\
&+ \sum_{i=1}^k\Big( - \in_{a_i b}\dual \rho+  \frac 1 2 \B_{a_i b 43} \Big) U_{a_1\ldots}\,^b\,_{\ldots a_k}. 
\end{split}
\eea
\end{enumerate}
\end{lemma}

\begin{proof}
It suffices to consider the case $k=1$. Using Proposition \ref{Proposition:commutehorizderivatives}, we have 
\beaa
\Ddot_3 \Ddot_b U_a&=& \nab_3\nab_bU_a -\eta_b \nab_3 U_a -\xib_b \nab_4 U_a, \\
\Ddot_b \Ddot_3 U_a&=& \nab_b \nab_3 U_a -\chib_{bc} \nab_c U_a -\ze_b \nab_3 U_a,\\
\Ddot_3 \Ddot_b U_a-\Ddot_b \Ddot_3 U_a&=&  \Rdot_{ac  3b}U_c= \R_{ac3b} U_c+ \frac 1 2  \B_{ac3b}U_c= -\in_{ac}\dual \bb_b U_c+ \frac 1 2  \B_{ac3b}U_c.
\eeaa
Hence,
\beaa
\,[\nab_3, \nab_b] U_a &=&  \nab_3\nab_bU_a- \nab_b \nab_3 U_a\\
&=&- \chib_{bc} \nab_c U_a+( \eta_b-\ze_b) \nab_3 U_a +\xib_b \nab_4 U_a -\in_{ac}\dual \bb_b U_c+ \frac 1 2  \B_{ac3b}U_c,
\eeaa
as stated. 
The commutator formula  for   $\,[\nab_4, \nab_b] U_a $ is derived easily by symmetry. 
Also,
\beaa
\Ddot_4 \Ddot_3 U_a&=&\nab_4 \nab_3 U_a -2\om \nab_3 U_a-2\etab^b \nab_b U_a , \\
\Ddot_3 \Ddot_4 U_a&=&\nab_3\nab_4 U_a   -2 \omb \nab_4 U_a-  2\eta^b \nab_b U_a,\\
\Ddot_4 \Ddot_3 U_a-\Ddot_3 \Ddot_4 U_a&=&  \Rdot_{ab  43}U_b= \R_{ab43} U_b+ \frac 1 2  \B_{ab43}U_b= -2\rhod \in_{ab} U^b+ \frac 1 2  \B_{ab43}U_b.
\eeaa
Hence
\beaa
\,[\nab_4, \nab_3]U_a  &=&2(\etab_b-\eta_b)\nab_b  U_a+ 2\om \nab_3 U_a-2\omb \nab_4 U  -2\rhod \in_{ab} U^b+ \frac 1 2  \B_{ab43}U_b,
\eeaa
as stated.
\end{proof}

Using the values of $\B$ given by Proposition \ref{proposition:componentsofB}, we obtain Corollary \ref{cor:comm-gen-B}.   In the following Lemma we specialize  to the case of $\sk_0$, $\sk_1$ and $\sk_2$.
  
 \begin{lemma}
   \lab{lemma:comm}
   The following commutation formulas hold true:
   \begin{enumerate}
\item Given   $f \in \sk_0$, we have
       \bea\label{eq:comm-nab3-nab4-naba-f-general}
       \begin{split}
        \,[\nab_3, \nab_a] f &=-\frac 1 2 \left(\trchb \nab_a f+\atrchb \dual \nab_a f\right)+(\eta_a-\ze_a) \nab_3 f-\chibh_{ab}\nab_b f\\
        &  +\xib_a \nab_4 f,\\
         \,[\nab_4, \nab_a] f &=-\frac 1 2 \left(\trch \nab_a f+\atrch \dual \nab_a f\right)+(\etab_a+\ze_a) \nab_4 f-\chih_{ab}\nab_b f \\
         & +\xi_a \nab_3 f, \\
         \, [\nab_4, \nab_3] f&= 2(\etab-\eta ) \c \nab f + 2 \om \nab_3 f -2\omb \nab_4 f. 
         \end{split}
       \eea

  \item   Given  $u\in \sk_1$, we have
    \bea\label{commutator-3-a-u-b}\label{commutator-u-in-SS1}
         \bsplit            
\,  [\nab_3,\nab_a] u_b    &=-\frac 1 2 \trchb \big( \nab_a u_b+\eta_b u_a-\de_{ab} \eta \c u \big)\\
& -\frac 1 2 \atrchb \big( \dual \nab_a u_b+\eta_b \dual u_a-\in_{ab} \eta\c u\big) \\
&+(\eta-\ze)_a \nab_3 u_b+\err_{3ab}[u],\\
  \err_{3ab}[u] &=-\dual \bb_a\dual u_b+\xib_a\nab_4 u_b-\xib_b \chi_{ac} u_c+\chi_{ab} \,\xib\c u-\chibh_{ac}\nab_c u_b-\eta_b\chibh_{ac}u_c\\
  &+\chibh_{ab}\eta\c u,
   \end{split}
   \eea
   \bea\label{commutator-4-a-u-b}
   \bsplit
\,  [\nab_4,\nab_a] u_b    &=-\frac 1 2 \trch \big( \nab_a u_b+\etab_b u_a-\de_{ab} \etab \c u \big) \\
&-\frac 1 2 \atrch \big( \dual \nab_a u_b+\etab_b \dual u_a-\in_{ab} \etab\c u\big)+(\etab +\ze)_a \nab_4 u_b \\
&+\err_{4ab}[u],\\
   \err_{4ab}[u]&=\dual \b_a\dual u_b+\xi_a\nab_3 u_b-\xi_b \chib_{ac} u_c+\chib_{ab} \,\xi\c u-\chih_{ac}\nab_c u_b-\etab_b\chih_{ac}u_c\\
   &+\chih_{ab}\etab\c u, 
      \end{split}
   \eea
   \bea
   \bsplit
 \, [\nab_4, \nab_3] u_a&=2 \om \nab_3 u_a -2\omb \nab_4 u_a+ 2(\etab_b-\eta_b ) \nab_b u_a +2(\etab \c u ) \eta_{a} -2 (\eta \c u )\etab_{a}\\
 &  -2 \dual \rho \dual u_a +\err_{43a}[u],\\
 \err_{43a}[u]&= 2 \big( \xib_{a}  \xi_b- \xi_{a}  \xib_b )u^b.
\end{split}
\eea

\item  Given  $u\in \sk_2$, we have 
    \bea\label{commutator-u-in-SS2}\label{commutator-3-a-u-bc}
         \bsplit            
\,  [\nab_3,\nab_a] u_{bc}    &=-\frac 1 2 \trchb\, (\nab_a u_{bc}+\eta_bu_{ac}+\eta_c u_{ab}-\de_{a b}(\eta \c u)_c-\de_{a c}(\eta \c u)_b )\\
&-\frac 1 2 \atrchb\, (\dual \nab_a u_{bc} +\eta_b\dual u_{ac}+\eta_c\dual u_{ab}- \in_{a b}(\eta \c u)_c- \in_{a c}(\eta \c u)_b )\\
&+(\eta_a-\ze_a)\nab_3 u_{bc}+\err_{3abc}[u],\\
\err_{3abc}[u]&= -2\dual \bb_a \dual u_{bc}+\xib_a \nab_4 u_{bc} -\xib_b\chi_{ad}u_{dc} -\xib_c\chi_{ad}u_{bd}+\chi_{ab}\xib_d u_{dc} \\
&+\chi_{ac}\xib_d u_{bd}-\chibh_{ad} \nab_d u_{bc} -\eta_b\chibh_{ad}u_{dc} - \eta_c\chibh_{ad}u_{bd}+\chibh_{ab}\eta_du_{dc} +\chibh_{ac}\eta_du_{bd},
   \end{split}
   \eea
   \bea
   \bsplit
\,  [\nab_4,\nab_a] u_{bc}    &=-\frac 1 2 \trch\, (\nab_a u_{bc}+\etab_bu_{ac}+\etab_c u_{ab}-\de_{a b}(\etab \c u)_c-\de_{a c}(\etab \c u)_b )\\
&-\frac 1 2 \atrch\, (\dual \nab_a u_{bc} +\etab_b\dual u_{ac}+\etab_c\dual u_{ab}- \in_{a b}(\etab \c u)_c- \in_{a c}(\etab \c u)_b )\\
&+(\etab_a+\ze_a)\nab_4 u_{bc}+\err_{4abc}[u],\\
\err_{4abc}[u]&= 2\dual \b_a \dual u_{bc}+\xi_a \nab_3 u_{bc} -\xi_b\chib_{ad}u_{dc} -\xi_c\chib_{ad}u_{bd}+\chib_{ab}\xi_d u_{dc} +\chib_{ac}\xi_d u_{bd}\\
& -\chih_{ad} \nab_d u_{bc} -\etab_b\chih_{ad}u_{dc} - \etab_c\chih_{ad}u_{bd}+\chih_{ab}\etab_du_{dc} +\chih_{ac}\etab_du_{bd}, 
     \end{split}
   \eea
   \bea\label{commutator-3-a-u-bc}
   \bsplit
   \, [\nab_4, \nab_3] u_{ab}&=  2\om \nab_3 u_{ab} - 2\omb \nab_4 u_{ab}  + 2 (\etab_c-\eta_c) \nab_c u_{ab}\\
&- 2\etab_a\eta_c u_{bc}-2\etab_b \eta_c u_{ac}+ 2\eta_a\etab_c u_{bc}+2\eta_b \etab_c u_{ac}-4\rhod \dual u_{ab}+\err_{43ab}[u]\\
   &=2 \om \nab_3 u_{ab} -2\omb \nab_4 u_{ab} + 2(\etab_c-\eta_c ) \nab_c u_{ab} + 4 \eta \hot (\etab \c u)  \\
   &-4 \etab \hot (\eta \c u)-4 \dual \rho \dual u_{ab}+\err_{43ab}[u],\\
\err_{43ab}[u]&= 2 \big( \xib_{a}  \xi_c- \xi_{a}  \xib_c )u^c\,_{b}+2 \big( \xib_{b}  \xi_c- \xi_{b}  \xib_c )u_{a} \,^c.
   \end{split}
\eea
       \end{enumerate}
 \end{lemma}

We   deduce the following corollary.
   \begin{corollary}\label{corr:comm} 
   The following commutation formulas hold true:
      \begin{enumerate}
      \item    Given  $u\in \sk_1$, we have
    \bea
    \bsplit
  \,    [\nab_3,\div] u&=-\frac 1 2 \trchb \,\big( \div u-\eta\c u\big) +
  \frac 1 2 \atrchb\,\big(  \div\dual u  -\eta \c \dual u\big)      +(\eta-\ze)\c\nab_3 u  \\
&+\err_{3\div}[u],\\
\err_{3\div}[u]&=-\dual \bb \c \dual u +\xib \c \nab_4 u -\xib \c \chih \c  u-\chibh \c \nab u-\eta \c \chibh \c u,
    \\
\,   [\nab_4,\div] u&=-\frac 1 2 \trch \, \big(\div u-\etab\c u\big) +\frac 1 2 \atrch\, \big(\div\dual u -   \etab\c \dual u\big) 
+(\etab+\ze)\c\nab_4 u   \\
& +\err_{4\div}[u], \\
\err_{4\div}[u]&=\dual \b \c \dual u +\xi \c \nab_3 u -\xi\c \chibh \c  u-\chih \c \nab u-\etab \c \chih \c u.
    \end{split}
    \eea
 Also, 
\bea\label{last-statement-item1}
\bsplit
\, [\nab_3, \nab \hot] u&=- \frac 1 2 \trchb \left( \nab\hot u +\eta \hot u\right)- \frac 1 2 \atrchb\, \dual  \left( \nab \hot  u+\etab \hot  u\right)+ (\eta-\ze) \hot \nab_3 u\\
&+\err_{3\hot}[u],\\
\err_{3\hot}[u] &=-\dual \bb \hot \dual u+\xib \hot \nab_4 u-\xib \hot ( \chi \c u)+\chih \,(\xib\c u)-\chibh \c \nab u-\eta \hot (\chibh \c u)\\
&+\chibh (\eta\c u),\\
\\
\, [\nab_4, \nab \hot] u&=- \frac 1 2 \trch \left( \nab\hot u +\etab \hot u\right)- \frac 1 2 \atrch\,  \dual  \left(\nab \hot  u+\etab \hot  u\right)+ (\etab+\ze) \hot \nab_4 u\\
&+\err_{4\hot}[u], \\
\err_{4\hot}[u] &=\dual \b \hot \dual u+\xi \hot \nab_3 u-\xi \hot ( \chib \c u)+\chibh \,(\xi\c u)-\chih \c \nab u-\etab \hot (\chih \c u)\\
&+\chih (\etab\c u).
\end{split}
\eea

\item Given $u\in \sk_2$, we have
\bea\label{eq:comm-nab3-nab4-div-all-errors}
\bsplit
\,[\nab_3,  \div] u &= -\frac 1 2 \trchb \big(  \div u -  2 \eta \c u\big) +\frac 1 2 \atrchb\big(  \div\dual u -2 \eta\c \dual u\big)\\
& +(\eta-\ze)\c\nab_3 u +\err_{3\div}[u],\\
 \err_{3\div}[u] &=-2\dual \bb \c \dual u +\xib\c \nab_4 u -\xib\c  \chi\c  u -(\chi\c u)\xib +\xib\c u\c\chi -\chibh\c \nab u\\
 & -\eta\c \chibh \c u-(\chibh\c u)\eta+\eta\c u\c\chibh, \\
\,[\nab_4,  \div] u &=-\frac 1 2 \trch \big(  \div u - 2\etab \c u\big) +\frac 1 2 \atrch\big(  \div\dual u -2 \etab\c \dual u\big)\\
& +(\etab+\ze)\c\nab_4 u +\err_{4\div}[u],\\
\err_{4\div}[u] &=2\dual \b \c \dual u +\xi\c \nab_3 u -\xi\c  \chib\c  u -(\chib\c u)\xi +\xi\c u\c\chib -\chih\c \nab u\\
& -\etab\c \chih \c u-(\chih\c u)\etab+\etab\c u\c\chih.
\end{split}
\eea
\end{enumerate} 
\end{corollary}
     
     \begin{proof}
     We check \eqref{last-statement-item1}. From \eqref{commutator-4-a-u-b} we have 
 \beaa
2 [\nab_4, \nab \hot] u_{ab}&=&  [\nab_4,\nab_a] u_b   + \,  [\nab_4,\nab_b] u_a - \de_{ab} [\nab_4, \div] u   \\
 &=&-   \trch \left( \nab\hot u +\etab \hot u\right) +(\etab +\ze)_a \nab_4 u_b-\frac 1 2  \atrch  H_{ab}\\
 &&+\err_{4ab}[u] +\err_{4ba}[u]-\de_{ab}\err_{4\div}[u]
 \eeaa
  where $H_{ab}$  denotes
 \beaa
 H_{ab}:&=&\big( \dual \nab_a u_b+\eta_b \dual u_a-\in_{ab} \eta\c u\big)+\big( \dual \nab_b u_a+\eta_b \dual u_a-\in_{ba} \eta\c u\big)\\
 && - \de_{ab}\big( \dual\nab \c u +\eta\c \dual u \big)\\
 &=& 2 (\dual \nab \hot u)_{ab}+ 2 ( \eta\hot\dual  u)_{ab}.
 \eeaa
 Recalling that $\dual\xi \hot  \eta = \xi\hot\dual\eta=\dual\big(\xi \hot  \eta\big)$ we infer that
 $  H= 2 \dual \big(\nab \hot u+ \eta\hot  u\big)$.
  This proves the desired result.

 We check the last statement in item 2. From \eqref{commutator-3-a-u-bc}
 \beaa
 \,  [\nab_4,\nab_a] u_{bc}    &=&-\frac 1 2 \trch\, (\nab_a u_{bc}+\etab_bu_{ac}+\etab_c u_{ab}-\de_{a b}(\etab \c u)_c-\de_{a c}(\etab \c u)_b )\\
&&-\frac 1 2 \atrch\, (\dual \nab_a u_{bc} +\etab_b\dual u_{ac}+\etab_c\dual u_{ab}- \in_{a b}(\etab \c u)_c- \in_{a c}(\etab \c u)_b )\\
&&+(\etab_a+\ze_a)\nab_4 u_{bc}
 \eeaa
 we deduce, recalling that $\delta_{ab} u_{ab}=0$,  
 \beaa
 [\nab_4,  \div] u_c &=& \delta_{ab} [\nab_4, \nab_a] u_{bc}\\
 &=&-\frac 1 2 \trch\, (\div u_{c}+(\etab \c u)_{c}-2(\etab \c u)_c-(\etab \c u)_c )\\
&&-\frac 1 2 \atrch\, (-\div \dual u_{c} +(\etab \c\dual u)_{c}+ \dual(\etab \c u)_c )+((\etab+\ze) \c\nab_4 u)_{c}
 \eeaa
 which proves the desired result.
   \end{proof}


\subsection{Commutation formulas with horizontal Lie derivatives}
\lab{section:horizLieDerivatives}


Recall that the Lie derivative of a $k$-covariant tensor $U$ relative to a vectorfield  $X$ is given by
\beaa
\Lie_X (Y_1, \ldots , Y_k) = X U(Y_1, \ldots, Y_k)-  U(\Lie_XY_1, \ldots,Y_k) - U(Y_1, \ldots, \Lie_X Y_k),
\eeaa
where $\Lie_X Y=[X, Y]$.
In components  relative to an arbitrary frame
\beaa
 \Lie_X U_{\a_1\ldots \a_k}:&=&\D_X U_{\a_1\ldots \a_k}  +\D_{\a_1} X^\b U_{\b \a_1\ldots \a_k} +\cdots  \D_{\a_k} X^\b U_{ \a_1\ldots  \b}.
\eeaa
 Recall also  the general commutation Lemma, see  chapter 7 in \cite{Ch-Kl}.
\begin{lemma}
\lab{Lemma:Generalcomm-Lie_X}
The following formula\footnote{This holds true  for an arbitrary pseudo-riemannian  space $(\MM, \g)$.}  for a vectorfield $X$ and   a $k$-covariant  tensor-field  $U$ holds true:
\begin{equation}\lab{commute.Generalcomm-Lie_X}
\D_\b(\LL_X U_{\a_1\ldots\a_k})-\LL_X(\D_\b U_{\a_1\ldots\a_k})=\sum_{j=1}^k\GaX_{\a_j\b \rho}U_{\a_1\ldots\,\,\,\,\ldots\a_k}^{\,\,\,\,\,\,\,\,\,\,\,\rho},
\end{equation}
where
\bea
\lab{definition:GaX}
\GaX_{\a\b\mu}&=&\frac 1 2 (\D_\a\piX_{\b\mu}+\D_\b\piX_{\a\mu}-\D_\mu\piX_{\a\b}). 
\eea
\end{lemma}

The proof  of the Lemma was  given in  \cite{Ch-Kl}, see Lemma 7.1.3,  based on the  following
\begin{lemma}
\lab{Lemma:R-PiXformula}
Given an arbitrary  vectorfield $X$ we have the identity
\beaa
\D_\mu \D_\nu  X_\b=\R_{\b  \mu\nu  \ga }X^\ga+ \GaX_{\mu\nu \b}.
\eeaa
\end{lemma}

\begin{proof}
Consider the tensor  $A_{\mu\nu\b} =\D_\mu \D_\nu  X_\b- \R_{\b  \mu\nu  \ga }X^\ga- \GaX_{\mu\nu \b}$ and observe that it verifies the symmetries
\beaa
A_{\mu\nu\b}= A_{\nu \mu  \b }=- A_{\mu\b\nu}.
\eeaa 
The proof   of  Lemma \ref{Lemma:R-PiXformula} follows  by observing that any such
 tensor must vanish identically\footnote{Indeed $A_{\mu\nu\b}=-A_{\mu\b \nu} = -  A_{ \b \mu\nu} =A_{ \b \nu\mu} =  A_{ \nu \b \mu} = -  A_{ \nu \mu \b }= -  A_{\mu  \nu  \b }$.}.
\end{proof}

We are now ready to define the horizontal Lie derivative operator $\Lieb$ as follows.

\begin{definition}[Horizontal Lie derivatives]
Given   vectorfields  $X$, $Y$,  the horizontal Lie  derivative  $\Lieb_XY$ is given by
 \beaa
 \Lieb_X Y :=\Lie_X Y+ \frac 1 2 \g(\Lie_XY, e_3) e_4+  \frac 1 2 \g(\Lie_XY, e_4) e_3.
 \eeaa
 Given  a horizontal covariant k-tensor $U$,  the horizontal  Lie derivative $\Lieb_X U $ is defined   to be the projection of $\Lie_X U $
  to the  horizontal space. Thus, for horizontal indices $A=a_1\ldots a_k$,
  \bea
  \lab{eq:projectedLie}
 ( \Lieb_X U) _{A}:&=&\nab_X U_{A}  +\D_{a_1} X^b U_{b \ldots a_k} +\cdots  +\D_{a_k} X^b U_{ a_1\ldots  b}.
  \eea
 \end{definition}

  \begin{lemma}
  \lab{lemma:commutation-Lieb-nab}
  The following commutation formulas hold true for a horizontal covariant k-tensor $U$ and  a vectorfield $X$
  \bea
  \bsplit
  \nab_b(\Lieb_X U_{A})-\Lieb_X(\nab_b U_{A})&=\sum_{j=1}^k  \GabbX_{a_jb  c }U_{a_1\ldots\,\,\,\,\ldots a_k}^{\,\,\,\,\,\,\,\,\,\,\,c},\\
   \nab_4(\Lieb_X U_A)-\Lieb_X(\nab_4 U_A) +\nab_{\Lieb_X e_4} U_{A}&=\sum_{j=1}^k  \GabbX_{a_j4  c }U_{a_1\ldots\,\,\,\,\ldots a_k}^{\,\,\,\,\,\,\,\,\,\,\,c},\\
   \nab_3(\Lieb_X U_A)-\Lieb_X(\nab_3 U_A) +\nab_{\Lieb_X e_3} U_{A}&=\sum_{j=1}^k  \GabbX_{a_j3  c }U_{a_1\ldots\,\,\,\,\ldots a_k}^{\,\,\,\,\,\,\,\,\,\,\,c},
  \end{split}
  \eea
  with\footnote{Here,  $\piX_{ab}$ is treated as a horizontal  symmetric 2-tensor, and $\piX_{a4}$,  $\piX_{a3}$, as   horizontal 1-forms.}
  \bea
  \lab{definition:Gabb_X}
  \bsplit
  \GabbX_{abc }&=\frac 1 2 (\nab_a\piX_{bc }+\nab_b\piX_{a c }-\nab_c \piX_{ab}),\\
  \GabbX_{a4b} &= \frac 1 2 (\nab_a\piX_{4b }+\nab_4\piX_{a b  }-\nab_b  \piX_{a4}),\\
\GabbX_{a3b} &= \frac 1 2 (\nab_a\piX_{3b }+\nab_3\piX_{a b  }-\nab_b  \piX_{a3}).
  \end{split}
  \eea
  \end{lemma}
 
\begin{proof}
Follows   easily by projecting   formula \eqref{commute.Generalcomm-Lie_X}  in Lemma \ref{Lemma:Generalcomm-Lie_X}, see also Lemma 9.1 in \cite{Chr-BH}.
\end{proof}

 We now extend the definition of horizontal Lie derivative to any $ U\in \T_k (\MM)\otimes   \O_l (\MM)$.
\begin{definition}\label{definition:hor-Lie-derivative} 
We define  the general horizontal derivatives as follows.
\begin{enumerate}
\item Given  $X\in \T(\MM)$ and a general,  horizontal  tensor-field $U \in \O_k(\MM)$, we define
\beaa
\Lied_XU &:=& \Lieb_XU.
\eeaa

 \item 
Given a  tensor  in   $ U\in \T_k (\MM)\otimes   \O_l (\MM) $    and $X\in T(\MM)$ 
 we define, for     $Z=  Z_1, \ldots, Z_k \in O(\M) $      and       $ Y= Y_1, \ldots  Y_l   \in  \O_1 (\MM)$ 
 
 \beaa
 \Lied_X U\big( Z,  Y\big) = X  U\big( Z,  Y\big)  &-   U\big( \Lie_X Z_1, \cdots Z_k, Y \big) - \ldots-  U\big( Z_1, \cdots \Lie_X  Z_k, Y \big)\\
& - U\big( Z,  \Lied_X Y_1, \ldots, Y_l  \big) -\ldots-   U\big( Z,  Y_1, \ldots,   \Lied_XY_l  \big). 
 \eeaa

\item We have
\beaa
\Lied_X  (U\otimes V)=  \Lied_X  U\otimes V      + U\otimes \Lied_X V.
\eeaa

\item  The definition can be extended by duality  to any mixed tensors  tensors in  $ \T^{k_1}_{k_2}  (\MM)\otimes   \O^{l_1}_{l_2}  (\MM)  $.
 \end{enumerate}
 \end{definition}

\begin{lemma}\label{lemma:commutator-Lied-Ddot}
The following commutation formulas hold true\footnote{With $\Gabb_X$   defined in  \eqref{definition:Gabb_X}.}  for $U\in \O_k(\MM)$ and $X \in \T(\MM)$,
\beaa
\Ddot_\mu (\Lied_X U_{a_1\ldots a_k})-\Lied_X (\Ddot_\mu U_{a_1\ldots a_k})&=& \sum_{j=1}^k\GabbX_{a_j \mu c}U_{a_1\ldots\,\,\,\,\ldots a_k}^{\,\,\,\,\,\,\,\,\,\,\, c}.
\eeaa

The following commutation formula holds true for $U\in \T(\MM) \otimes  \O_k(\MM)$ and $X \in \T(\MM)$,
\beaa
\Ddot_\mu(\Lied_X U_{\gamma a_1\ldots a_k})-\Lied_X(\Ddot_\mu U_{\gamma  a_1\ldots a_k})=\GaX_{\gamma \mu  \rho}{U^{\rho}}_{a_1\ldots\,\,\,\,\ldots a_k}+\sum_{j=1}^k\GabbX_{a_j\mu c}U_{\gamma a_1\ldots\,\,\,\,\ldots a_k}^{\,\,\,\,\,\,\,\,\,\,\, c }.
\eeaa
\end{lemma}

\begin{proof} 
Follows   easily by projecting   formula \eqref{commute.Generalcomm-Lie_X}  in Lemma \ref{Lemma:Generalcomm-Lie_X}, see also Lemma 9.1 in \cite{Chr-BH}.  
\end{proof}


\subsection{Main equations using conformally invariant derivatives}
\lab{section:conf.inv.Derivatives}


Consider  frame transformations of the form
\beaa
e_3'=\la^{-1} e_3, \qquad  e'_4 = \la e_4 , \qquad e_a'= e_a.
\eeaa
Note that  under   the   above  mentioned  frame transformation we have
\beaa
 \trchb'&=&\la^{-1} \trchb, \quad \atrchb'=\la^{-1} \atrchb,  \quad   \trch'=\la\trch, \quad \atrch'=\la \atrch, \\
  \xi'&=& \la^2 \xi, \quad   \eta'=\eta, \quad \etab'=\etab,  \quad \xib'=\la^{-2}\xib,\\
   \a'&=&\la^2 \a,\quad \b'=\la \b,    \quad \rho'=\rho,    \quad \rhod'=\rhod,  \quad \bb'=\la^{-1} \bb,\quad  \aa'=\la^{-2} \aa,
   \eeaa
   and
   \beaa
 \omb'&=& \la^{-1}\left(\omb +\frac{1}{2} e_3(\log \la)\right), \quad \om'= \la\left(\om -\frac{1}{2} e_4(\log \la)\right), \quad
 \ze'= \ze - \nab (\log \la).
\eeaa

\begin{definition}[$s$-conformally invariants]
\lab{def:sconformalinvariants}
We say that  a horizontal tensor $f$ is $s$-conformally invariant  if, under the  conformal  frame transformation above,  it changes as $f'=\la^s f $. 
\end{definition}

\begin{remark}
If $f$ $s$-conformal invariant, then  $\nab_3 f, \nab_4 f, \nab_a f$ are not conformal invariant.
\end{remark} 

 We correct the lacking of being conformal invariant by making the following  definition.  

\begin{lemma}\label{lemma:definition-conformal-derivatives}
If $f$ is $s$-conformal invariant, then: 
 \begin{enumerate}
 \item $\nabc_3 f:= \nab_3f-2 s \omb f$ is $(s-1)$-conformally invariant.
 
 \item $\nabc_4 f:= \nab_4f+2 s \om f$ is $(s+1)$-conformally invariant.
 
 \item $\nabc_A f:= \nab_Af+ s \ze_A f$ is $s$-conformally invariant. 
  \end{enumerate}
\end{lemma}

\begin{proof}
Immediate verification.
\end{proof}

\begin{remark} 
Note that $s$ is precisely what   in \cite{Ch-Kl} is called       the  signature of the tensor. In GHP formalism \cite{GHP}, the signature is related to the boost weights of the complex scalars.
\end{remark}

Using these definitions we rewrite the main equations as follows.
\begin{proposition}
\label{prop-nullstr-conformal}
We have
\beaa
\nabc_3\trchb&=&-|\chibh|^2-\frac 1 2 \big( \trchb^2-\atrchb^2\big)+2\divc\xib +  2 \xib\c(\eta+\etab),\\
\nabc_3\atrchb&=&-\trchb\atrchb +2\curlc \xib + 2 \xib\wedge(-\eta+\etab),\\
\nabc_3\chibh&=&-\trchb\,  \chibh+  \nabc\hot \xib+     \xib\hot(\eta+\etab)-\aa,
\eeaa
\beaa
\nabc_3\trch
&=& -\chibh\c\chih -\frac 1 2 \trchb\trch+\frac 1 2 \atrchb\atrch    +   2   \divc\eta+ 2 \big(\xi\c \xib +|\eta|^2\big)+ 2\rho,\\
\nabc_3\atrch
&=&-\chibh\wedge\chih-\frac 1 2(\atrchb \trch+\trchb\atrch)+ 2 \curlc \eta  + 2 \xib\wedge\xi  -  2 \dual \rho,\\
\nabc_3\chih
&=&-\frac 1 2 \big( \trch \chibh+\trchb \chih\big)-\frac 1 2 \big(-\dual \chibh \, \atrch+\dual \chih\,\atrchb\big)
+  \nabc\hot \eta + \xib\hot\xi +\eta\hot\eta,
\eeaa
\beaa
\nabc_4\trchb
&=& -\chih\c\chibh -\frac 1 2 \trch\trchb+\frac 1 2 \atrch\atrchb    +  2   \divc \etab+ 2\big( \xi\c \xib +|\etab|^2\big)+2\rho,\\
\nabc_4\atrchb
&=&-\chih\wedge\chibh-\frac 1 2(\atrch \trchb+\trch\atrchb)+ 2 \curlc\etab  + 2 \xi\wedge\xib+2 \dual \rho,\\
\nabc_4\chibh
&=&-\frac 1 2 \big( \trchb \chih+\trch \chibh\big)-\frac 1 2 \big(-\dual \chih \, \atrchb+\dual \chibh\,\atrch\big)
+ \nabc\hot \etab + \xi\hot\xib + \etab\hot\etab,
\eeaa
\beaa
\nabc_4\trch&=&-|\chih|^2-\frac 1 2 \big( \trch^2-\atrch^2\big)+ 2 \divc\xi  + 2   \xi\c(\etab+\eta),\\
\nabc_4\atrch&=&-\trch\atrch + 2 \curlc\xi + 2 \xi\wedge(-\etab+\eta),\\
\nabc_4\chih&=&-\trch\,  \chih+  \nabc\hot \xi+     \xi\hot(\etab+\eta)-\a,
\eeaa
\beaa
\nabc_3 \etab -\nabc_4\xib &=& -\chibh\c(\etab-\eta) -\frac{1}{2}\trchb(\etab-\eta)+\frac{1}{2}\atrchb(\dual\etab-\dual\eta)  +\bb, \\
\nabc_4 \eta    -    \nabc_3\xi &=& -\chih\c(\eta-\etab) -\frac{1}{2}\trch(\eta-\etab)+\frac{1}{2}\atrch(\dual\eta-\dual\etab) -\b.\\
\eeaa
Also,
\beaa
\divc\chih &=& \frac{1}{2}\nabc(\trch) -\frac{1}{2}\dual\nabc (\atrch) -\atrch\dual\eta-\atrchb\dual\xi -\b,\\
\divc\chibh  &=& \frac{1}{2}\nabc(\trchb) -\frac{1}{2}\dual\nabc(\atrchb) -\atrchb\dual\etab-\atrch\dual\xib +\bb.
\eeaa
\end{proposition}

    \begin{proposition}\label{prop:bianchi-conformal} 
    We have
    \beaa
    \nabc_3\a-   \nabc\hot \b&=&-\frac 1 2 \big(\trchb\a+\atrchb\dual \a)+
   4\eta\hot \b - 3 (\rho\chih +\rhod\dual\chih),\\
\nabc_4\beta - \divc\a &=&-2(\trch\beta-\atrch \dual \b)  +\a\c  \etab + 3  (\xi\rho+\dual \xi\rhod),\\
     \nabc_3\b+\divc\varrho&=&-(\trchb \b+\atrchb \dual \b)+2\bb\c \chih+3 (\rho\eta+\rhod\dual \eta)+    \a\c\xib,\\
 \nabc_4 \rho-\divc \b&=&-\frac 3 2 (\trch \rho+\atrch \rhod)+2\etab \c\b-2\xi\c\bb-\frac 1 2 \chibh \c\a,\\
   \nabc_4 \rhod+\curlc \b&=&-\frac 3 2 (\trch \rhod-\atrch \rho)-2\etab\c\dual \b-2\xi\c\dual \bb+\frac 1 2 \chibh \c\dual \a, \\
     \nabc_3 \rho+\divc\bb&=&-\frac 3 2 (\trchb \rho -\atrchb \rhod) -2\eta \c\bb+2\xib\c\b-\frac 12\chih\c\aa,
 \\
   \nabc_3 \rhod+\curlc \bb&=&-\frac 3 2 (\trchb \rhod+\atrchb \rho)- 2\eta\c\dual \bb-2\xib\c\dual\b-\frac 1 2 \chih\c\dual \aa,\\
     \nabc_4\bb-\divc\varoc&=&-(\trch \bb+ \atrch \dual \bb)+2\b\c \chibh
    -3 (\rho\etab-\rhod\dual \etab)-    \aa\c\xi,\\
     \nabc_3\bb +\divc\aa &=&-2(\trchb\,\bb-\atrchb \dual \bb)   -\aa\c \eta - 3  (\xib\rho-\dual \xib \rhod),\\
     \nabc_4\aa+   \nabc\hot \bb&=&-\frac 1 2 \big(\trch\aa-\atrch\dual \aa)
 -   4\etab\hot \bb - 3  (\rho\chibh -\rhod\dual\chibh).
\eeaa
    \end{proposition}


\section{Commutations  with   horizontal  wave operators}\lab{section:Horizwaveoperators}


 Consider a   spacetime $(\MM, \g)$  with a horizontal structure   induced by a  null pair $(e_3, e_4) $.
 \begin{definition}
 We  define  the wave operator for  tensor-fields  $\psi \in \O_k(\MM) $ to be
 \bea\label{eq:def=squared-2}
 \squared_k\psi:= \g^{\mu\nu} \Ddot_\mu\Ddot_ \nu \psi.
\eea
\end{definition}


\subsection{Commutation with  $\protect\Lied_X$}


\begin{proposition}\label{prop:commutator-Lied-squared} 
The following commutation formula\footnote{Recall that $\Lied$ has been introduced in Definition \ref{definition:hor-Lie-derivative}.} holds true for $\psi \in \sk_2$ and $X \in \T(\MM)$,
\beaa
\, [\Lied_X, \squared_2]\psi_{ab}&=& -\piX^{\mu\nu} \Ddot_\mu \Ddot_\nu\psi_{ab}\\
&& - {{\GaX}^\mu}_{\mu  \rho}\Ddot^{\rho} \psi_{ab}-2\GabbX_{a\mu c}{\Ddot^{\mu}}{\psi^c}_{ b}-2\GabbX_{b\mu c}\Ddot^{\mu}{\psi_{ a}}^c \\
&&-  \Ddot^\nu ( \GabbX_{a \nu c}){\psi^c}_b-\Ddot^\nu(\GabbX_{b \nu c}){\psi_a}^c.
\eeaa
\end{proposition}

\begin{proof} 
We have
\beaa
\, [\Lied_X, \squared_2]\psi_{ab} &=& [\Lied_X, \g^{\mu\nu} \Ddot_\mu \Ddot_\nu]\psi_{ab}\\
&=& ( \Lied_X \g^{\mu\nu} )\Ddot_\mu \Ddot_\nu\psi_{ab}+\g^{\mu\nu}  [\Lied_X ,\Ddot_\mu] \Ddot_\nu\psi_{ab}+\g^{\mu\nu}  \Ddot_\mu [ \Lied_X, \Ddot_\nu]\psi_{ab}.
\eeaa
Using Lemma \ref{lemma:commutator-Lied-Ddot} where recall  the definition of $\GaX,\GabbX$ in section \ref{section:horizLieDerivatives}, we obtain
\beaa
\, [\Lied_X, \squared_2]\psi_{ab} &=& -\piX^{\mu\nu} \Ddot_\mu \Ddot_\nu\psi_{ab}-  \Ddot^\nu \big( \GabbX_{a \nu c}{\psi^c}_b+\GabbX_{b \nu c}{\psi_a}^c \big)\\
&&- {{\GaX}^\mu}_{\mu  \rho}\Ddot^{\rho} \psi_{ab}-\GabbX_{a\mu c}{\Ddot^{\mu}}{\psi^c}_{ b}-\GabbX_{b\mu c}\Ddot^{\mu}{\psi_{ a}}^c\\
&=& -\piX^{\mu\nu} \Ddot_\mu \Ddot_\nu\psi_{ab} - {{\GaX}^\mu}_{\mu  \rho}\Ddot^{\rho} \psi_{ab}-2\GabbX_{a\mu c}{\Ddot^{\mu}}{\psi^c}_{ b}-2\GabbX_{b\mu c}\Ddot^{\mu}{\psi_{ a}}^c \\
&&-  \Ddot^\nu ( \GabbX_{a \nu c}){\psi^c}_b-\Ddot^\nu(\GabbX_{b \nu c}){\psi_a}^c,
\eeaa
as stated.
\end{proof}


\subsection{Commutation with $\Ddot_X$}


\begin{lemma}
\lab{LEMMA:COMMUTATIONNAB_XSQUARED}
We have in a vacuum spacetime
\beaa
\squared  (  X^\b \Db_\b  U_a)- X^\b \Db_\b \squared U_a&=& \pi^{\mu\nu}  \Db_\mu \Db_\nu  U_a+
\big( \D^\mu  \pi_\mu \,^\b-\frac 1 2  \D^\b  \tr \pi) \Db_\b U_a\\
&& - 2 X^\b  \R_{ac \b\mu} \Db^\mu U_c+  \D^\b X^\mu \R_{ac \b\mu}  U^c  \\
&&-  X^\b \B_{ac \b\mu}  \Db^\mu U_c+ \frac 12  \D^\b X^\mu \B_{ac \b\mu} U^c +\frac 1 2  X^\b \D^\mu \B_{a  c\mu \b} U^c.
\eeaa
\end{lemma}

\begin{proof} 
Straightforward computation using Lemma \ref{Lemma:R-PiXformula} and Proposition \ref{Proposition:commutehorizderivatives},  see section \ref{sec:proof-lemma-comm-NabXsquared}.
\end{proof}


\subsection{Killing tensor and Carter operator}


Recall that the  deformation tensor    of a vectorfield  $\piX$ is defined as
\beaa
\piX_{\mu\nu}:=\D_{(\mu} X_{\nu)}= \D_{\mu} X_{\nu}+\D_{\nu} X_{\mu}.
\eeaa
The vectorfield is  said to be Killing if  $\piX\equiv 0$.     The Kerr spacetime   has, in addition to  the symmetries  generated by  its two linearly independent Killing vectorfields $\T$ and $\Z$, a  higher order symmetry  defined by a Killing tensor.

\begin{definition}\label{def:killing-tensor}
 A symmetric $2$-tensor $K_{\mu\nu}$  is said to be  a Killing tensor if  its deformation  3-tensor $\Pi$,  defined below,  vanishes identically.
 \bea\label{definition-Pi}
\Pi_{\mu\nu\rho} := \D_{(\mu} K_{\nu \rho)}=\D_{\mu} K_{\nu \rho}+\D_{\nu} K_{\rho\mu}+\D_{\rho} K_{\mu\nu}.
\eea
\end{definition}

\begin{remark}
Observe that if $X, Y$ are Killing vectorfields  then the symmetric $2$-tensor   $K =\frac 1 2 (X \otimes Y+ Y\otimes X) $ is a Killing tensor.
\end{remark}

We define the second order differential operator associated to a tensor-field  $\psi \in \sk_k $. 

\begin{definition}\label{definition:operator-KK-general} 
Given a symmetric tensor $K$ its associated second order differential operator $\KK$  applied to a tensor $\psi\in \sk_k$  is defined as
\bea\label{eq:definition-KK-operators}
\mathcal{K}(\psi)= \Db_{\mu}( K^{\mu\nu} \Db_\nu(\psi)).
\eea
\end{definition}

We now compute the commutators of $\mathcal{K}$ with $\square_\g$ in terms of the symmetric tensor $\Pi$.

\begin{proposition}\label{COMMUTATION-KK-SQUARE} 
In a vacuum spacetime, the commutator between the differential operator $\KK$ and the $\square_\g$ operator applied to a scalar function $\phi$ is given by
\beaa
[\KK, \square_\g] \phi&=&\err[\Pi](\phi)
\eeaa
where $\err[\Pi](\phi)$ denotes terms involving $\Pi$ given by
\beaa
\err[\Pi](\phi)&:=&\D^\mu \left( \big(\D^\a \Pi_{\a\nu\mu}-\frac 1 2 \D_\mu {\Pi^{\a}}_{\a\nu} +\frac 1 2 \D_\nu {\Pi^{\a}}_{\a\mu} \big) \Ddot^\nu \phi -2 \Pi_{\mu\a\nu} \Ddot^{\a}\Ddot^\nu \phi\right)\\
&&- 2 \left(\D^\a \Pi_{\a\nu\mu} \right) \Ddot^\mu \Ddot^\nu \phi.
\eeaa
\end{proposition}

\begin{proof} 
See section \ref{sec:proof-comm-KK}.
\end{proof}

    
\section{Main equations in complex notations}\label{section:complex-notations}


In this section we introduce complex notations for the Ricci coefficients and the curvature components with the objective of simplifying the main equations. From the real scalars, 1-tensors and symmetric traceless 2-tensors already introduced, we define their complexified version which results in anti-self dual tensors.


\subsection{Complex notations}


Recall Definition \ref{definition-SS-real} of the set of real horizontal $k$-tensors $\sk_k=\sk_k(\MM, \mathbb{R})$ on $\MM$. For instance, 
\begin{itemize}
\item $(a, b) \in \sk_0$ is a pair of real scalar function on $\MM$, 
\item $f \in \sk_1$ is a real horizontal 1-tensor on $\MM$,
\item $u \in \sk_2$ is a real horizontal symmetric traceless 2-tensor on $\MM$.
\end{itemize}

By Definition \ref{definition-hodge-duals}, the duals of real horizontal tensors are real horizontal tensors of the same type, i.e. $\dual f \in \sk_1$ and $\dual u \in \sk_2$.

We define the complexified version of horizontal tensors on $\MM$.

\begin{definition} 
We denote by $\sk_k(\mathbb{C})=\sk_k(\MM, \mathbb{C})$ the set of complex anti-self dual $k$-tensors on $\MM$. More precisely, 
\begin{itemize}
\item $a+ i b \in \sk_0(\mathbb{C})$ is a complex scalar function on $\MM$ if $(a, b) \in \sk_0$,
\item $F= f+ i \dual f  \in \sk_1(\mathbb{C})$ is a complex anti-self dual 1-tensor on $\MM$ if $f \in \sk_1$,
\item $U=u + i \dual u \in \sk_2(\mathbb{C})$ is a complex anti-self dual symmetric traceless 2-tensor on $\MM$ if $u \in \sk_2$.
\end{itemize}
\end{definition}

Observe that $F\in \sk_1(\mathbb{C})$ and $ U\in \sk_2(\mathbb{C})$ are indeed anti-self dual tensors, i.e.
\beaa
\dual F=-i F, \qquad \dual U=-i U. 
\eeaa
More precisely 
\beaa
U_{12}=U_{21}= i \dual U_{12}=i \in_{12} U_{22}=-i U_{11}, \qquad U_{11} =i U_{12}.
\eeaa

Recall that the derivatives $\nab_3$, $\nab_4$ and $\nab_a$ are real derivatives. We can use the dual operators to define the complexified version of the $\nab_a$ derivative, which allows to simplify the notations in the main equations.

\begin{definition}
We define the complexified version of the horizontal derivative as 
\beaa
\DD= \nab+ i \dual \nab, \qquad \DDov=\nab- i \dual \nab.
\eeaa
More precisely, we have:
\begin{itemize}
\item for $a+i b \in \sk_0(\mathbb{C})$, 
\beaa
\DD(a+ib) &:=& (\nabla+i\dual\nabla)(a+ib),\qquad \DDov(a+ib) := (\nabla-i\dual\nabla)(a+ib).
\eeaa

\item For $f+ i \dual f \in\sk_1(\mathbb{C}) $, 
\beaa
\DD\c(f+i\dual f) &:=& (\nabla+i\dual\nabla)\c(f+i\dual f)\, =0, \\
\DDov\c(f+i\dual f) &:=& (\nabla-i\dual\nabla)\c(f+i\dual f),\\
\DD\hot(f+i\dual f) &:=& (\nabla+i\dual\nabla)\hot(f+i\dual f).
\eeaa
\item For  $u+ i \dual u \in \sk_2(\mathbb{C}) $, 
\beaa
\DD\c (u+i\dual u) &:=& (\nabla+i\dual\nabla)\c (u+i\dual u)=0,\\
\DDov \c (u+i\dual u) &:=& (\nabla- i\dual\nabla) \c (u+i\dual u).
\eeaa
\end{itemize}
\end{definition}

Note that $\dual \DD=-i\DD$.

For $F= f+i \dual f\in\sk_1(\mathbb{C})$ the operator $-\frac 1 2 \DD\hot $ is formally adjoint to  the operator $\DDov \c U$ applied to $U\in \sk_2(\mathbb{C})$.
For $ h=a+ib \in \sk_0(\mathbb{C})$ the operator $-\DD h $ is formally adjoint to the operator $\DDov \c F$ applied to $F\in\sk_1(\mathbb{C})$.  These notions makes sense  literally only  if  the horizontal structure is integrable.

 \begin{lemma}\label{lemma:spacetime-elliptic-estimates}
 For $F=f+i\dual f \in \sk_1(\CCC)$ and $U=u+i\dual u \in \sk_2(\CCC)$, we have
   \bea\label{eq:adjoint-operators}
 ( \DD \hot   F) \c   \ov{U}  &=&  - 2 F \c (\DD \c \ov{U}) -( (H+\Hb ) \hot F )\c \ov{U} + 2 \DD\c (F \c \ov{U}).
 \eea
 \end{lemma}
 
 \begin{proof}  
 We look at the real parts. Then
 \beaa
 ( \nab\hot   f) \c   u  &=&\big(\nab_a f_b+\nab_b f_a- \de_{ab} \div f \big) u_{ab} =2(\nab_a f_b) u_{ab}=2 \nab_a (u_{ab} f_b)-2(\div u) \c f. 
 \eeaa
Using Lemma \ref{lemma:divergence-spacetime-horizontal-ch1} applied to $\xi= u \c f$ we obtain
  \beaa
 ( \nab\hot   f) \c   u  &=&2 \nab^a (u_{a b} f_b)-2 (\eta+ \etab) \c (u \c f) -2(\div u) \c f \\
 &=&  -2(\div u) \c f -( (\eta+ \etab) \hot f )\c u +2\div(u\c f).
 \eeaa
By complexifying, we obtain the stated identity.
 \end{proof}

\begin{lemma}
The following holds:
\begin{itemize}
\item If $\xi, \eta\in \sk_1$
\beaa
\xi\c \eta+i\dual\xi\c \eta &=& \frac{1}{2}\Big((\xi+i\dual\xi)\c (\ov{\eta+i\dual \eta})\Big),\\
\xi\hot \eta+i\dual(\xi\hot \eta) &=& \frac{1}{2}\Big((\xi+i\dual \xi)\hot (\eta+i\dual \eta)\Big).
\eeaa
\item If $\eta\in \sk_1 $,    $u\in\sk_2$ 
\beaa
u\c \eta+i\dual u\c \eta &=& \frac{1}{2}(u+i\dual u)\c (\ov{\eta+i\dual \eta}).
\eeaa

\item If $u, v\in \sk_2$ 
\beaa
u\c v+i\dual u\c v &=& \frac{1}{2}(u+i\dual u)\c (\ov{v+i\dual v}).
\eeaa

\item If $(a, b) \in \sk_0 $
\beaa
\nabla a-\dual\nabla b+i(\dual\nabla a+\nabla b) &=& \DD(a+ib).
\eeaa

\item If $\xi \in \sk_1 $
\beaa
\div \xi+i\curl \xi &=& \frac{1}{2}\ov{\DD}\c(\xi+i\dual \xi),\\
\nabla\hot \xi+i\dual(\nabla\hot \xi) &=& \frac{1}{2}\mathcal{D}\hot(\xi+i\dual \xi).
\eeaa

\item If $u\in \sk_2 $
\beaa
\div u+i\dual(\div u) &=& \frac{1}{2}\ov{\DD}\c(u+i\dual u).
\eeaa
\end{itemize}
\end{lemma}

\begin{proof}
The first identities rely on Lemma \ref{lemma:usefulidentitiesforcomplexification}. The other rely on the following identities, for  $\xi\in \sk_1 $,  $u\in\sk_2 $,
\beaa
&&\nab\c\dual \xi=\curl \xi, \,\,  \dual\nab\c \xi=-\curl \xi, \,\,  \dual\nab\c\dual \xi=\nab \xi,\,\,   \nabla\hot \dual \xi=\dual\nabla\hot \xi=\dual(\nabla\hot \xi), \,\,\\
&& \dual\nabla\hot \dual \xi = -\nabla\hot \xi, \dual(\div u)=\nab\c \dual u, \,\,  \dual\nab\c u=-\dual(\div u), \,\, \dual\nab\c\dual u=\nab\c u.
\eeaa
\end{proof}

\begin{lemma}\label{dot-hot-complex} 
Let  $E, F \in \sk_1(\mathbb{C})$ and $U\in \sk_2(\mathbb{C}) $.  Then
\bea\label{simil-Leibniz}
E \hot ( \ov{F} \c U)+F \hot ( \ov{E} \c U)&=&2 ( E \c \ov{F}+ \ov{E} \c F) \ U.
\eea
Also, for $E=e+i \dual e$, $F=f+i \dual f$
\bea\label{simil-Leib-half}
E \hot ( \ov{F} \c U) &=& 4 \big(e \c f - i e \wedge f \big) U.
\eea
\end{lemma}

\begin{proof} 
Recall, see Lemma \ref{dot-hot}, 
\beaa
\xi \hot ( \eta  \c u) + \eta  \hot ( \xi \c u)=2 (\xi\c \eta  ) u. 
\eeaa
For $E= \xi +i \dual \xi$, $F=\eta+i\dual \eta$ and $U=u+i \dual u $, we have
\beaa
E \hot ( \ov{F} \c U)&=& (\xi+ i \dual \xi ) \hot \big( \ov{(\eta + i \dual \eta)} \c (u + i \dual u)\big) \\
&=& 2(\xi+ i \dual \xi ) \hot \big( (u \c \eta) + i \dual ( u \c \eta)\big) \\
&=& 4\Big (\xi \hot  (u \c \eta) + i \dual ( \xi \hot  (u \c \eta) ) \Big), 
\\
F \hot ( \ov{E} \c U)&=& 4 \Big(\eta \hot  (u \c \xi) + i \dual ( \eta \hot  (u \c \xi) ) \Big) .
\eeaa
Therefore
\beaa
E \hot ( \ov{F} \c U)+F \hot ( \ov{E} \c U)&=&4 (\xi \hot  (u \c \eta) + i \dual ( \xi \hot  (u \c \eta) ) ) +4 (\eta \hot  (u \c \xi) + i \dual ( \eta \hot  (u \c \xi) ) )\\
 &=&4 (\xi \hot  (u \c \eta) +\eta \hot  (u \c \xi))+ 4i \dual ( \xi \hot  (u \c \eta)  +  \eta \hot  (u \c \xi) )\\
 &=&8 ( (\xi \c \eta) \ u)+ 8i \dual (  (\xi \c \eta)  u)=8 (\xi \c \eta) U
\eeaa
while
\beaa
E \c \ov{F}+ \ov{E} \c F&=& 2(\xi \c \eta + i \dual \xi \c \eta)  + 2(\eta \c \xi + i \dual \eta \c \xi)=4(\xi \c \eta ). 
\eeaa
Hence,
\beaa
E \hot ( \ov{F} \c U)+F \hot ( \ov{E} \c U)&=&2\big(E \c \ov{F}+ \ov{E} \c F\big)U
\eeaa
as stated. The second identity can be derived in the same manner.
\end{proof}


\subsubsection{Leibniz formulas}


We collect here Leibniz formulas involving the derivative operators defined above.

\begin{lemma}\label{SIMPLIFICATION-ANGULAR}
Let $h$ be a scalar function, $F \in \sk_1(\mathbb{C}) $,  $U\in \sk_2(\mathbb{C})$.  Then
\bea
\bsplit
 \ov{\DD} \c (h F)&= h \ov{\DD} \c F+  \ov{\DD}(h) \c F \lab{ov-HH-hF},\\
 \DD\hot (h F)&= h  \DD\hot  F+ \DD( h) \hot F  \lab{DD-hot-hF},\\
\ov{\DD}\c( h U)&=  \ov{\DD}(h) \c U  + h(\ov{\DD} \c U) \lab{ov-DD-hu},\\
 \DD\hot(\ov{F}\c U)&=  2   (\DD\c \ov{F} ) U + 2   (\ov{F}\c \DD) U\label{Leibniz-hot},\\
U \c \ov{\DD} F&= U (\ov{\DD} \c F).\label{rule-1}
\end{split}
\eea
Also,
\bea
\bsplit
F \hot (\DDb \c U)&= 2(F\c \DDb) U= 4F \c \nab U , \lab{Leib-eq-DDb} \\
(F \c \DDb) U + (\ov{F} \c \DD) U&= 4f \c \nab U= 2(F +\ov{F}) \c \nab U. \lab{Leib-eq-DDb-DD-nab}\lab{relation0angular=der}
\end{split}
\eea
\end{lemma}

\begin{proof} 
 Straightforward verifications, see section 
\ref{Appendix:ProofLemma-simplificationangular}.
\end{proof}

\begin{lemma}
\lab{Lemma:DDhot(DDb)=lap}
As a corollary  of \eqref{Leib-eq-DDb} we derive the following  formula for  $U\in \sk_2(\mathbb{C})$
\bea
\DD \hot (\DDb \c U)&=& 2 \lap_2 U- 4 \Kh U -i  \big(\atrch\nab_3+\atrchb \nab_4\big) U
\eea
where
\beaa
\Kh=- \frac 14  \trch \trchb-\frac 1 4 \atrch \atrchb+\frac 1 2 \chih \c \chibh-  \frac 1 4 \rho.
\eeaa
\end{lemma}

\begin{proof}
According to   \eqref{Leib-eq-DDb} we have
\beaa
\frac 1 2  \DD \hot (\DDb \c U)&=( \DD\c \DDb)  U= (\nab^a+i \dual\nab^a)(  (\nab^a-i \dual\nab^a) U= 2\lap_2 U -2 i \in^{ab} \nab_a\nab_b  U.
\eeaa
On the other hand, appealing to  Proposition \ref{Gauss-equation-2-tensors} 
 we have
 \beaa
 2  \in^{ab} \nab_a\nab_b  U&=&\in^{ab} [ \nab_a, \nab_b]  A= 2 \left( \frac 1 2 (\atrch\nab_3+\atrchb \nab_4) \psi +2  \, \Kh \dual A \right)  \\
&=&   (\atrch\nab_3+\atrchb \nab_4) \psi - 4 i \Kh A
 \eeaa
Hence $\frac 1 2  \DD \hot (\DDb \c U)=2\lap_2 U- 4\Kh U-i  (\atrch\nab_3+\atrchb \nab_4) U$ as stated.
\end{proof}


\subsection{Main equations  in complex form}


We now extend the definitions for the Ricci coefficients and curvature components given in Sections \ref{subsection:Ricci-coefficients} and \ref{subsection:curvature}, to the complex case by using the anti-self dual tensors defined above. 
 
\begin{definition}\lab{def:complexRicciandcurvaturecoefficients} 
We define the following complex anti-self dual tensors:
\beaa
A:=\a+i\dual\a, \quad B:=\b+i\dual\b, \quad P:=\rho+i\dual\rho,\quad \Bb:=\bb+i\dual\bb, \quad \Ab:=\aa+i\dual\aa,
\eeaa      
 and   
\beaa
&& X=\chi+i\dual\chi, \quad \Xb=\chib+i\dual\chib, \quad H=\eta+i\dual \eta, \quad \Hb=\etab+i\dual \etab, \quad Z=\ze+i\dual\ze, \\ 
&& \Xi=\xi+i\dual\xi, \quad \Xib=\xib+i\dual\xib.
\eeaa    
In particular, note that 
\beaa
\tr X = \trch-i\atrch, \quad \Xh=\chih+i\dual\chih, \quad \tr\Xb = \trchb -i\atrchb, \quad \Xbh=\chibh+i\dual\chibh.
\eeaa
\end{definition}

The complex notations allow us to rewrite the Ricci equations in  a more compact  form. 
\begin{proposition}
\label{prop-nullstr:complex}
\beaa
\nab_3\tr\Xb +\frac{1}{2}(\tr\Xb)^2+2\omb\,\tr\Xb &=& \DD\c\ov{\Xib}+\Xib\c\ov{\Hb}+\ov{\Xib}\c(H-2Z)-\frac{1}{2}\Xbh\c\ov{\Xbh},\\
\nab_3\Xbh+\Re(\tr\Xb) \Xbh+ 2\omb\,\Xbh&=&\frac 1 2  \DD\hot \Xib+\frac 1 2   \Xib\hot(H+\Hb-2Z)-\Ab,
\eeaa
\beaa
\nab_3\tr X +\frac{1}{2}\tr\Xb\tr X-2\omb\tr X &=& \DD\c\ov{H}+H\c\ov{H}+2P+\Xib\c\ov{\Xi}-\frac{1}{2}\Xbh\c\ov{\Xh},\\
\nab_3\widehat{X} +\frac{1}{2}\tr\Xb\, \widehat{X} -2\omb\widehat{X} &=&\frac 1 2  \DD\hot H  +\frac 1 2 H\hot H -\frac{1}{2}\ov{\tr X} \widehat{\Xb}+\frac{1}{4}\Xib\hot\Xi,
\eeaa
\beaa
\nab_4\tr\Xb +\frac{1}{2}\tr X\tr\Xb -2\om\tr\Xb &=& \DD\c\ov{\Hb}+\Hb\c\ov{\Hb}+2\ov{P}+\Xi\c\ov{\Xib}-\frac{1}{2}\Xh\c\ov{\Xbh},\\
\nab_4\widehat{\Xb} +\frac{1}{2}\tr X\, \widehat{\Xb} -2\om\widehat{\Xb} &=&\frac  12  \DD\hot\Hb  +\frac 1 2 \Hb\hot\Hb -\frac{1}{2}\ov{\tr\Xb} \widehat{X}+\frac{1}{4}\Xi\hot\Xib,
\eeaa
\beaa
\nab_4\tr X +\frac{1}{2}(\tr X)^2+2\om\tr X &=& \DD\c\ov{\Xi}+\Xi\c\ov{H}+\ov{\Xi}\c(H+2Z)-\frac{1}{2}\Xh\c\ov{\Xh},\\
\nab_4\Xh+\Re(\tr X)\Xh+ 2\om\Xh&=&\frac 1 2  \DD\hot \Xi+\frac 12   \Xi\hot(\Hb+H+2Z)-A.
\eeaa
Also,
\beaa
\nab_3Z +\frac{1}{2}\tr\Xb(Z+H)-2\omb(Z-H) &=& -2\DD\omb -\frac{1}{2}\widehat{\Xb}\c(\ov{Z}+\ov{H})\\
&&+\frac{1}{2}\tr X\Xib+2\om\Xib -\Bb+\frac{1}{2}\ov{\Xib}\c\Xh,\\
\nab_4Z +\frac{1}{2}\tr X(Z-\Hb)-2\om(Z+\Hb) &=& 2\DD\om +\frac{1}{2}\widehat{X}\c(-\ov{Z}+\ov{\Hb})\\
&&-\frac{1}{2}\tr\Xb\Xi-2\omb\Xi -B-\frac{1}{2}\ov{\Xi}\c\Xbh,\\
\nab_3\Hb -\nab_4\Xib &=&  -\frac{1}{2}\ov{\tr\Xb}(\Hb-H) -\frac{1}{2}\Xbh\c(\ov{\Hb}-\ov{H}) -4\om\Xib+\Bb,\\
\nab_4H -\nab_3\Xi &=&  -\frac{1}{2}\ov{\tr X}(H-\Hb) -\frac{1}{2}\Xh\c(\ov{H}-\ov{\Hb}) -4\omb\Xi-B,
\eeaa
and
\beaa
\nab_3\om+\nab_4\omb -4\om\omb -\xi\c \xib -(\eta-\etab)\c\ze +\eta\c\etab&=&   \rho.
\eeaa
Also,
\beaa
\frac{1}{2}\ov{\DD}\c\Xh +\frac{1}{2}\Xh\c\ov{Z} &=& \frac{1}{2}\DD\ov{\tr X}+\frac{1}{2}\ov{\tr X}Z-i\Im(\tr X)H-i\Im(\tr \Xb)\Xi-B,\\
\frac{1}{2}\ov{\DD}\c\Xbh -\frac{1}{2}\Xbh\c\ov{Z} &=& \frac{1}{2}\DD\ov{\tr\Xb}-\frac{1}{2}\ov{\tr\Xb}Z-i\Im(\tr\Xb)\Hb-i\Im(\tr X)\Xib+\Bb,
\eeaa
and,
\beaa
\curl\ze&=&-\frac 1 2 \chih\wedge\chibh   +\frac 1 4 \big(  \trch\atrchb-\trchb\atrch   \big)+\om \atrchb -\omb\atrch+\dual \rho.
\eeaa
\end{proposition}

We rewrite the Gauss equation in Proposition \ref{Gauss-equation-2-tensors} for complex tensors. 
\begin{proposition}\label{Gauss-equation-2-tensors-complex}\label{Gauss-equation}  
The following identity holds true for $\Psi \in \sk_k(\CCC)$ for $k=1,2$:
\bea\label{Gauss-eq-first-com}\label{Gauss-eq-2-complex}
[ \nab_a, \nab_b] \Psi &=&\left( \frac 1 2 (\atrch\nab_3+\atrchb \nab_4) \Psi -i k \, \Kh  \Psi\right)  \in_{ab}
 \eea
 where
 \beaa
\Kh &=& - \frac 1 8 \tr X \ov{\tr \Xb}- \frac 1 8 \tr \Xb  \ov{\tr X}+ \frac 1 4 \Xh \c \ov{\Xbh}+\frac 1 4 \ov{\Xh} \c \Xbh -\frac 1 2 P -\frac 1 2  \ov{P}.
\eeaa
      \end{proposition}

The complex notations allow us to rewrite the Bianchi identities as follows.  
  \begin{proposition}\label{prop:bianchi:complex} 
    We have,
 \beaa
 \nab_3A -\frac 1 2 \DD\hot B &=& -\frac{1}{2}\tr\Xb A+4\omb A +\frac 1 2 (Z+4H)\hot B -3\ov{P}\Xh,\\
\nab_4B -\frac{1}{2}\ov{\DD}\c A &=& -2\ov{\tr X} B -2\om B +\frac{1}{2}A\c  (\ov{2Z +\Hb})+3\ov{P} \,\Xi,\\
\nab_3B-\DD\ov{P} &=& -\tr\Xb B+2\omb B+\ov{\Bb}\c \Xh+3\ov{P}H +\frac{1}{2}A\c\ov{\Xib},\\
\nab_4P -\frac{1}{2}\DD\c \ov{B} &=& -\frac{3}{2}\tr X P +\frac{1}{2}(2\Hb+Z)\c\ov{B} -\ov{\Xi}\c\Bb -\frac{1}{4}\Xbh\c \ov{A}, \\
\nab_3P +\frac{1}{2}\ov{\DD}\c\Bb &=& -\frac{3}{2}\ov{\tr\Xb} P -\frac{1}{2}(\ov{2H-Z})\c\Bb +\Xib\c \ov{B} -\frac{1}{4}\ov{\Xh}\c\Ab, \\
\nab_4\Bb+\DD P &=& -\tr X\Bb+2\om\Bb+\ov{B}\c \Xbh-3P\Hb -\frac{1}{2}\Ab\c\ov{\Xi},\\
\nab_3\Bb +\frac{1}{2}\ov{\DD}\c\Ab &=& -2\ov{\tr\Xb}\,\Bb -2\omb\,\Bb -\frac{1}{2}\Ab\c (\ov{-2Z +H})-3P \,\Xib,\\
\nab_4\Ab +\frac{1}{2}\DD\hot\Bb &=& -\frac{1}{2}\tr X \Ab+4\om\Ab +\frac{1}{2}(Z-4\Hb)\hot \Bb -3P\Xbh.
\eeaa
    \end{proposition} 
    
    \begin{proof} 
    We derive the equations for $A$ and $B$. Observe that from the Bianchi identity
  in Proposition   
   \ref{prop:bianchi}, 
    \beaa
    \nab_3\a-  \nab\hot \b&=&-\frac 1 2 \big(\trchb\a+\atrchb\dual \a)+4\omb \a+
  (\ze+4\eta)\hot \b - 3 (\rho\chih +\rhod\dual\chih),
\eeaa
    we obtain
    \beaa
     \dual\nab_3\a&=&\dual(\nab\hot \b)-\frac 1 2 \big(\trchb \dual\a-\atrchb 
     \a)+4\omb \dual \a+
\dual( (\ze+4\eta)\hot \b )- 3 (\rho \dual\chih -\rhod\chih).
    \eeaa
    This implies
    \beaa
    \nab_3 A&=& \nab_3( \a + i \dual \a)\\
    &=& \nab\hot \b+i \dual(\nab\hot \b)-\frac 1 2 \big(\trchb\a+\atrchb\dual \a)-\frac 1 2 i\big(\trchb \dual\a-\atrchb 
     \a)+4\omb (\a+i\dual a)\\
     &&+(\ze+4\eta)\hot \b+ i\dual( (\ze+4\eta)\hot \b ) - 3 (\rho\chih +\rhod\dual\chih)- 3i (\rho \dual\chih -\rhod\chih)\\
      &=& \frac 1 2  \DD \hot (\b + i \dual \b)-\frac 1 2 (\trchb- i \atrchb )
     \a-\frac 1 2 (\trchb- i \atrchb) i\dual\a+4\omb (\a+i\dual a)\\
     &&+ \frac 1 2 ((\ze+4\eta+ i \dual (\ze +4 \eta))\hot (\b+i \dual \b))  - 3 (\rho- i \rhod)\chih - 3 (\rho- i \rhod) i\dual\chih     \eeaa
     which finally gives
     \beaa
      \nab_3 A      &=& \frac 1 2 \DD \hot B-\frac 1 2 \tr\Xb A+4\omb A+ \frac 1 2 (Z+4H)\hot B  - 3 \ov{P}\hat{X},
     \eeaa
     as stated.
     From the equation 
       \beaa
\nab_4\beta - \div\a &=&-2(\trch\beta-\atrch \dual \b) - 2  \om\b +\a\c  (2 \ze +\etab) + 3  (\xi\rho+\dual \xi\rhod),
\eeaa
     we obtain
         \beaa
     \dual\nab_4\beta  &=&\dual\div\a-2(\trch \dual\beta+\atrch \b) - 2  \om\dual\b +\dual(\a\c  (2 \ze +\etab) )+ 3  (\dual\xi\rho-\xi\rhod).
     \eeaa
     This implies
     \beaa
     \nab_4B &=& \nab_4(\b +i \dual \b)\\
     &=&\div\a+i\dual\div\a-2(\trch\beta-\atrch \dual \b)-2i(\trch \dual\beta+\atrch \b)  - 2  \om(\b+i \dual\b) \\
     &&+\a\c  (2 \ze +\etab)+i\dual(\a\c  (2 \ze +\etab) ) + 3  (\xi\rho+\dual \xi\rhod)+ 3i  (\dual\xi\rho-\xi\rhod)\\
      &=&\frac 1 2 \ov{\DD} \c (\a + i \dual \a)-2(\trch+i\atrch )\b-2(\trch  +i\atrch )i\dual \b - 2  \om(\b+i \dual\b) \\
     &&+\frac 1 2 (\a+ i \dual \a)\c  (2 \ze +\etab-i\dual(2 \ze +\etab) ) + 3  (\rho-i\rhod) \xi+ 3 (\rho- i \rhod) i\dual\xi 
     \eeaa
     which finally gives
     \beaa
      \nab_4B       &=&\frac 1 2 \ov{\DD} \c A-2\ov{\tr X}B - 2  \om B+\frac 1 2 A\c  (2 \ov{Z} + \ov{\Hb}  ) + 3  \ov{P} \Xi
     \eeaa
     as stated. The remaining equations are checked in the same fashion.
    \end{proof}


\subsection{Main complex equations using conformal derivatives}


\begin{definition}
 We define the following conformal angular derivatives in the complex notation:
 \begin{itemize}
\item For $a+i b \in \sk_0(\mathbb{C}) $  we define
\beaa
\DDc(a+ib) &:=& \big(\nabc +i\dual \nabc\big)(a+ib).
\eeaa

\item For  $f+i \dual f \in\sk_1(\mathbb{C}) $ we define
\beaa
\DDc(f+i\dual f) &:=& \big(\nabc+i\dual \nabc\big) \c (f+i\dual f),
\\
\DDc \hot(f+i\dual f) &:=& (\nabc+i\dual\nabc )\hot(f+i\dual f).
\eeaa
\item For $u+ i \dual u \in \sk_2(\mathbb{C})$ we define
\beaa
\DDc \c(u+i\dual u) &:=& \big(\nabc +i\dual\nabc \big)\c(u+i\dual u).
\eeaa
\item In all the above cases we set
\beaa
\DDbc&:=&\nabc-i\nabc.
\eeaa
\end{itemize}
\end{definition}

These complex notations allow us to rewrite the null structure equations as follows.  
\begin{proposition}
\label{prop-nullstr:complex-conf}
We have
\beaa
\nabc_3\tr\Xb +\frac{1}{2}(\tr\Xb)^2 &=& \DDc\c\ov{\Xib}+\Xib\c\ov{\Hb}+\ov{\Xib}\c H-\frac{1}{2}\Xbh\c\ov{\Xbh},\\
\nabc_3\Xbh+\Re(\tr\Xb) \Xbh&=&\frac 1 2  \DDc\hot \Xib+  \frac 1 2  \Xib\hot(H+\Hb)-\Ab,
\eeaa
\beaa
\nabc_3\tr X +\frac{1}{2}\tr\Xb\tr X &=& \DDc\c\ov{H}+H\c\ov{H}+2P+\Xib\c\ov{\Xi}-\frac{1}{2}\Xbh\c\ov{\Xh},\\
\nabc_3\widehat{X} +\frac{1}{2}\tr\Xb\, \widehat{X} &=&\frac 1 2 \DDc\hot H  +\frac 1 2 H\hot H -\frac{1}{2}\ov{\tr X} \widehat{\Xb}+\frac 1 4 \Xib\hot\Xi,
\eeaa
\beaa
\nabc_4\tr\Xb +\frac{1}{2}\tr X\tr\Xb &=& \DDc\c\ov{\Hb}+\Hb\c\ov{\Hb}+2\ov{P}+\Xi\c\ov{\Xib}-\frac{1}{2}\Xh\c\ov{\Xbh},\\
\nabc_4\widehat{\Xb} +\frac{1}{2}\tr X\, \widehat{\Xb} &=&\frac 1 2  \DDc\hot\Hb  +\frac 1 2 \Hb\hot\Hb -\frac{1}{2}\ov{\tr\Xb} \widehat{X}+\frac 1 4 \Xi\hot\Xib,
\eeaa
\beaa
\nabc_4\tr X +\frac{1}{2}(\tr X)^2 &=& \DDc\c\ov{\Xi}+\Xi\c\ov{H}+\ov{\Xi}\c H-\frac{1}{2}\Xh\c\ov{\Xh},\\
\nabc_4\Xh+\Re(\tr X)\Xh&=&\frac 1 2  \DDc\hot \Xi+\frac 1 2   \Xi\hot(\Hb+H)-A,
\eeaa
\beaa
\nabc_3\Hb -\nabc_4\Xib &=&  -\frac{1}{2}\ov{\tr\Xb}(\Hb-H) -\frac{1}{2}\Xbh\c(\ov{\Hb}-\ov{H}) +\Bb,\\
\nabc_4H -\nabc_3\Xi &=&  -\frac{1}{2}\ov{\tr X}(H-\Hb) -\frac{1}{2}\Xh\c(\ov{H}-\ov{\Hb}) -B.
\eeaa
Also,
\beaa
\frac{1}{2}\ov{\DDc}\c\Xh &=& \frac{1}{2}\DDc\ov{\tr X}-i\Im(\tr X)H-i\Im(\tr \Xb)\Xi-B,\\
\frac{1}{2}\ov{\DDc}\c\Xbh &=& \frac{1}{2}\DDc\ov{\tr\Xb}-i\Im(\tr\Xb)\Hb-i\Im(\tr X)\Xib+\Bb.
\eeaa
\end{proposition}
    
The complex notations allow us to rewrite the Bianchi identities as follows.  
  \begin{proposition}\label{prop:bianchi:complex-conf} 
    We have
 \beaa
 \nabc_3A -\frac 1 2 \DDc\hot B &=& -\frac{1}{2}\tr\Xb A + 2 H   \hot B -3\ov{P}\Xh,\\
\nabc_4B -\frac{1}{2} \DDbc \c A &=& -2\ov{\tr X} B +\frac{1}{2}A\c \ov{\Hb}+3\ov{P} \,\Xi,\\
\nabc_3B-\DDc\ov{P} &=& -\tr\Xb B+\ov{\Bb}\c \Xh+3\ov{P}H +\frac{1}{2}A\c\ov{\Xib},\\
\nabc_4P -\frac{1}{2}\DDc\c \ov{B} &=& -\frac{3}{2}\tr X P + \Hb \c\ov{B} -\ov{\Xi}\c\Bb -\frac{1}{4}\Xbh\c \ov{A}, 
\eeaa
\beaa
\nabc_3P +\frac{1}{2}\DDbc \c\Bb &=& -\frac{3}{2}\ov{\tr\Xb} P - \ov{H} \c\Bb +\Xib\c \ov{B} -\frac{1}{4}\ov{\Xh}\c\Ab, \\
\nabc_4\Bb+\DDc P &=& -\tr X\Bb+\ov{B}\c \Xbh-3P\Hb -\frac{1}{2}\Ab\c\ov{\Xi},\\
\nabc_3\Bb +\frac{1}{2}\DDbc \c\Ab &=& -2\ov{\tr\Xb}\,\Bb  -\frac 1 2  \Ab\c \ov{H}-3P \,\Xib,\\
\nabc_4\Ab +\frac 1 2 \DDc\hot\Bb &=& -\frac{1}{2}\tr X \Ab - 2 \Hb\hot \Bb -3P\Xbh.
\eeaa
    \end{proposition}

    
\subsection{Renormalized Bianchi identities}\label{section:RenormBianchiIds}
    
    
    We  define renormalized  derivatives  for the curvature  components
     $A, B, P, \Bb, \Ab$  in the spirit of Definition 7.3.2 of \cite{Ch-Kl}.  These
      renormalizations play an important role  in the derivation of the generalized Regge Wheeler equation for  the quantity   $\qfb$ in section \ref{section:gRW-equation-qfb}.

 \begin{definition}\label{definition-Psi_3-Psi_4}   
 Given a conformally invariant curvature component $\Psi$ (i.e. either $A, \Ab, B, \Bb, P$) with signature $s$ we define the operators
    \beaa
    \Psi_3&:=&\nabc_3\Psi+ \frac 1 2 \big( 3-s \big) \ov{\tr \Xb} \Psi,\\
    \Psi_4&:=& \nabc_4\Psi+ \frac 1 2 \big( 3+s \big)\tr X \Psi.
    \eeaa
    \end{definition}

     \begin{proposition}
        Using this  definition the Bianchi identities take the form\footnote{Note that $B$  was  changed to $\ov{B}$  to maintain  the correct definition.}
      \beaa
      \bsplit
\ov{A}_3 -\frac 1 2 \ov{\DDc\hot B} &= 2 \ov{H}   \hot \ov{B} -3P\ov{\Xh},\\
\ov{B}_4 -\frac{1}{2} \DDc \c \ov{A} &=\frac{1}{2}\ov{A}\c \Hb+3P \,\ov{\Xi},\\
\ov{B}_3-\ov{\DDc}P &=\Bb\c \ov{\Xh}+3P\ov{H} +\frac{1}{2}\ov{A}\c\Xib,\\
P_4 -\frac{1}{2}\DDc\c \ov{B} &= \Hb \c\ov{B} -\ov{\Xi}\c\Bb -\frac{1}{4}\Xbh\c \ov{A}, 
\end{split}
\eeaa
\beaa
\bsplit
P_3 +\frac{1}{2}\DDbc \c\Bb &= - \ov{H} \c\Bb +\Xib\c \ov{B} -\frac{1}{4}\ov{\Xh}\c\Ab, \\
\Bb_4+\DDc P &= \ov{B}\c \Xbh-3P\Hb -\frac{1}{2}\Ab\c\ov{\Xi},\\
\Bb_3 +\frac{1}{2}\DDbc \c\Ab &=   -\frac 1 2  \Ab\c \ov{H}-3P \,\Xib,\\
\Ab_4 +\frac 1 2 \DDc\hot\Bb &=  - 2 \Hb\hot \Bb -3P\Xbh.
\end{split}
\eeaa
    \end{proposition}

    
\chapter{The Kerr spacetime}
\lab{SECTION:KERR}


In this chapter, we provide basic facts concerning the Kerr spacetime.


\section{Boyer-Lindquist coordinates}


We consider the Kerr metric in standard Boyer-Lindquist coordinates $(t, r, \th, \phi)$,
$$\g_{a,m}=-\frac{|q|^2\Delta}{\Sigma^2}(dt)^2+\frac{\Sigma^2(\sin\theta)^2}{|q|^2}\left(d\phi-\frac{2amr}{\Sigma^2}dt\right)^2+\frac{|q|^2}{\Delta}(dr)^2+|q|^2(d\theta)^2,$$
where
\bea\label{definition-q}
q=r+ i a \cos\th,
\eea
and
\beaa
\left\{\ba{lll}
\Delta &=& r^2-2mr+a^2,\\
|q|^2 &=& r^2+a^2(\cos\theta)^2,\\
\Sigma^2 &=& (r^2+a^2)|q|^2+2mra^2(\sin\theta)^2=(r^2+a^2)^2-a^2(\sin\theta)^2\Delta.
\ea\right.
\eeaa
Observe that
\beaa
(2mr-|q|^2)\Sigma^2=-|q|^4\Delta+4a^2m^2r^2(\sin\theta)^2.
\eeaa
The metric  $\g=\g_{a,m}$ can also be written in the form
 \beaa
\g &=& \,    -\frac{\left(\Delta-a^2\sin^2\theta\right)}{|q|^2}dt^2-
\frac{ 4 amr }{|q|^2}   \sin^2\theta     dt  d\phi+\frac{|q|^2}{\Delta}dr^2+ |q|^2 d\theta^2+
\frac{ \Si^2}{|q|^2}\sin^2\theta
d\phi^2.
\eeaa
Note that  $\g_{tt} \g_{\phi \phi}- \g_{t\phi}^2 =-\De \sin^2 \th$ and that
the non-vanishing  components of the inverse metric are given by
\bea
\lab{eq:inversemetric-Kerr}
\bsplit
\g^{00}&=-\frac{\Si^2}{|q|^2 \De}, \qquad 
\g^{0\phi}=-\frac{2 a mr }{|q|^2\De},\qquad 
\g^{\phi\phi} =\frac{\De- a^2 \sin^2\th}{|q|^2 \De \sin^2 \th}, \\ 
\g^{rr}&=\frac{\De}{|q|^2}, \qquad \g^{\th\th}=\frac{1}{|q|^2}. 
\end{split}
\eea
The volume element $d\mu$  of $\g$  is given by  
\beaa
d\mu&=&|q|^2 \sin\th dt dr d\th d\phi, \qquad \sqrt{|g|} =|q|^2 \sin\th.
\eeaa
We also   note   that  
\bea
\T=\pr_t , \qquad \Z=\pr_\phi,
\eea
 are both Killing and $\T$  is only time-like in the     complement of the ergoregion,
 i.e.  $|q|^2> 2 Mr $. 
 The    domain of outer communication    of the Kerr metric is given by, 
\beaa
\RR=\{(\theta,r,t,\phi)\in(0,\pi)\times(r_+,\infty)\times\mathbb{R}\times\mathbb{S}^1\},
\eeaa
where     $r_+:=m+\sqrt{m^2-a^2}$,  the larger root of $\De$, corresponds to the      event horizon.


\section{Vectorfields $\That, \Rhat$}


\begin{definition}
We introduce the  vectorfields $\That, \Rhat $ as follows:  
\bea
\That:&=&\pr_t+\frac{a}{r^2+a^2} \pr_\phi \lab{define:That},\\
\Rhat:&=&\frac{\De}{r^2+a^2} \pr_r.   \lab{define:Rhat}
\eea
\end{definition}

Note that  $\Rhat$ is regular at the horizon, as opposed to $\partial_r$.  Unlike $\T$,  which  is  spacelike in the ergoregion, 
 the vectorfield  $\That$, to which we  refer as   the  Hawking vectorfield,
  is  time-like  in the domain of outer communication. 
More precisely we have
 \begin{proposition} 
 \lab{Lemma:Hawking-vf}
 The  vectorfield  $\That$ 
 is  timelike        for $r>r_+$ and null   on the horizon $r=r_+$.  More precisely
   \bea\label{eq:g-That-That}
   \g(\That, \That)&=& -\De\frac{|q|^2}{ (r^2+a^2)^2}. 
   \eea   
   \end{proposition} 
   
 \begin{proof}
 The proof follows from the following more general computation.
 \begin{lemma}\label{lemma:inside-Hawking-vf}
 The vectorfield  $T_\la=\T+\la \Z$, for a scalar function  $\la $,  verifies
 \beaa
 \g(T_\la, T_\la) &=&-\frac{\De}{|q|^2 } \Big( 1+ a^2\la^2(\sin\theta)^4 \Big) + \frac{\sin^2 \th}{|q|^2} E_\la(r),
 \\
   E_\la(r)&=& a^2   - 4 am r\la         +\la ^2(r^2+a^2)^2.
   \eeaa
 \end{lemma}

To check \eqref{eq:g-That-That} we 
    replace  $\la$ by $\frac{a}{r^2+a^2}$ to deduce
 \beaa
 \g(\That,\That) &=&-\frac{\De}{|q|^2 } \Big( 1+ \frac{a^4}{(a^2+r^2)^2} (\sin\theta)^4 \Big)+  \frac{\sin^2 \th}{|q|^2} \Big(  a^2   - 4 am r\frac{a}{a^2+r^2}   +       \frac{a^2}{(a^2+r^2)^2} (r^2+a^2)^2\Big)\\
 &=& -\frac{\De}{|q|^2 } \Big( 1+ \frac{a^4}{(a^2+r^2)^2} (\sin\theta)^4 \Big) +  \frac{ 2 a^2\sin^2 \th}{|q|^2} \Big(  1   -  2mr  \frac{1}{a^2+r^2}   \Big)\\
 &=& -\frac{\De}{|q|^2 } \Big( 1+ \frac{a^4}{(a^2+r^2)^2} (\sin\theta)^4 \Big) +  \frac{ 2 a^2\sin^2 \th}{|q|^2} \frac{\De}{a^2+r^2}\\
 &=& -\frac{\De}{|q|^2 }  \Big( 1+ \frac{a^4}{(a^2+r^2)^2} (\sin\theta)^4- 2 \frac{a^2}{r^2+a^2} \sin^2\th  \Big)\\
  &=& -\frac{\De}{|q|^2 } \Big( 1-  \frac{a^2}{r^2+a^2} \sin^2\th  \Big)^2= -\frac{\De}{|q|^2 } \frac{|q|^4}{(r^2+a^2)^2}\\
  &=&-\De\frac{|q|^2}{ (r^2+a^2)^2} 
 \eeaa
 as stated.  To check Lemma \ref{lemma:inside-Hawking-vf}  we write in  BL coordinates  $\T+\la  \Z=\pr_t +\la \pr_\phi$. Hence
 \beaa
 \g(\T+\la  \Z, \T+\la  \Z) &=&\g(\pr_t +\la \pr_\phi, \pr_t +\la \pr_\phi)=\g_{tt} + 2\la \g_{t\phi}+\la^2 \g_{\phi\phi}\\
 &=&   -\frac{\left(\Delta-a^2\sin^2\theta\right)}{|q|^2}+  2\la \big(\frac{ - 2 amr }{|q|^2}   \sin^2\theta \big)+ \lambda^2 \frac{ \Si^2}{|q|^2}\sin^2\theta\\
 &=&\frac{1}{|q|^2} \Big(-\De +a^2 \sin^2\th  -  4 am r\la \sin^2\th+\la^2\Si^2 \sin^2 \th\Big) \\
 &=&-\frac{\De}{|q|^2 }+ \frac{\sin^2 \th}{|q|^2} \Big(a^2   - 4 am r\la         +\la^2(r^2+a^2)^2-a^2\la^2(\sin\theta)^2\Delta\Big) \\
 &=&-\frac{\De}{|q|^2 } \Big( 1+ a^2\la^2(\sin\theta)^4 \Big) + \frac{\sin^2 \th}{|q|^2} \Big( a^2   - 4 am r\la         +\la^2(r^2+a^2)^2\Big)\\
 &=&-\frac{\De}{|q|^2 } \Big( 1+ a^2\la^2(\sin\theta)^4 \Big) + \frac{\sin^2 \th}{|q|^2} E_\la(r),
 \eeaa
 as stated. 
   \end{proof}
   
\begin{remark}  
As  a consequence  of the lemma, we also   deduce that  the Killing vectorfield $T_\HH=\T+\om_\HH \Z$, with $\om_\HH=\frac{a}{r_{+}^2+a^2}$ the angular velocity  of the horizon, is null on the horizon and timelike in a small neighborhood of it in $r>r_+$. 
\end{remark}


\section{Principal null frames}\label{section:values-Kerr}


The Kerr metric is a spacetime of Petrov Type D, i.e. its Weyl curvature can be diagonalized with two linearly independent eigenvectors, the so-called principal null (PN) directions. We now present the ingoing PN frame $( e_4^{(in)}, e_3^{(in)})$ and the outgoing PN frame $( e_4^{(out)}, e_3^{(out)})$.


\subsection{Ingoing PN frame}
\lab{section:ingoingPNframe}


The ingoing PN  frame (with $\D_3 e_3=0$), regular towards the future  for all $r>0$,  is given by
\bea
\lab{eq:null-pair-in}
\bsplit
e^{(in)}_4&=\frac{r^2+a^2}{|q|^2} \pr_t +\frac{\De}{|q|^2} \pr_r +\frac{a}{|q|^2} \pr_\phi, \\
e^{(in)}_3  & =\frac{r^2+a^2}{\De} \pr_t -\pr_r +\frac{a}{\De} \pr_\phi.
\end{split}
\eea
Note that
\bea
e^{(in)}_4(r)=\frac{\Delta}{|q|^2}, \qquad e^{(in)}_3(r)=-1.
\eea
We complete the PN frame with  the following specific  choice
of horizontal frames $e_1, e_2$,
\bea
\lab{eq:canonicalHorizBasisKerr}
e_1=\frac{1}{|q|}\pr_\th,\qquad e_2=\frac{a\sin\th}{|q|}\pr_t+\frac{1}{|q|\sin\th}\pr_\phi.
\eea
We refer to \eqref{eq:canonicalHorizBasisKerr} as  the \textit{canonical horizontal basis} of Kerr.

Using the Hawking vectorfield     $\That=\pr_t+\frac{a}{r^2+a^2} \pr_\phi$, we  have 
\bea
\lab{formula:e_3e_4-ThatRhat}
e^{(in)}_4=\frac{r^2+a^2}{|q|^2} \left(\That+\frac{\De}{ r^2+a^2} \pr_r \right), \qquad  e^{(in)}_3=\frac{r^2+a^2}{\De}\left( \That -
\frac{\De}{r^2+a^2} \pr_r \right)
\eea
from which we  deduce
\bea\label{eq-incoming-3-4-T-R}
e^{(in)}_4=\frac{r^2+a^2}{|q|^2} \Big(\That+\Rhat\Big), \qquad  e^{(in)}_3=\frac{r^2+a^2}{\De}\Big( \That -\Rhat \Big).
\eea

  \begin{lemma}
  \lab{Lemma:T-frame}
 The following identities hold true.
   \bea
   \lab{eq:ThatRhat-e_3e_4-Kerr}
   \begin{split}
 \That &=\frac 1 2 \left( \frac{|q|^2}{r^2+a^2} e^{(in)}_4+\frac{\De}{r^2+a^2}  e^{(in)}_3\right), \\
  \Rhat &= \frac 1 2 \left( \frac{|q|^2}{r^2+a^2} e^{(in)}_4-\frac{\De}{r^2+a^2}  e^{(in)}_3\right),
  \end{split}
 \eea
 and
 \bea\label{eq:T-Z-That-e_2}
 \begin{split}
\T &=\frac{r^2+a^2}{|q|^2} \That-\frac{a\sin\th}{|q|} e_2, \\
 \Z&=-\frac{a(r^2+a^2)\sin^2\th}{|q|^2 }\That+\frac{(r^2+a^2)\sin\th}{|q|} e_2.
 \end{split}
 \eea
Combining we derive
 \bea\label{eq:T-Z-ingoing-frame-Kerr}
 \bsplit
 \T &= \frac 1 2 \left( e^{(in)}_4+\frac{\De}{|q|^2} e^{(in)}_3\right)- \frac{a\sin\th}{|q|} e_2, \\ 
 \Z &=\frac{(r^2+a^2)\sin\th}{|q|} e_2-\frac 1 2 a  \sin^2\th \left(e^{(in)}_4+\frac{\De}{|q|^2}  e^{(in)}_3\right).
 \end{split}
 \eea
  Note also that $\That, \Rhat$ are perpendicular to  the horizontal  structure and $\g(\That, \Rhat)=0$.
  \end{lemma}
  
  \begin{proof}
Straightforward verification.
  \end{proof}

  
\subsubsection{Ingoing Ricci and curvature coefficients}


The real Ricci coefficients in the ingoing PN frame are given by
\beaa
&&\chih=\chibh=\xi=\xib=\omb=0, \\
&& \trch=\frac{2\De r}{|q|^4},\quad \atrch=\frac{2a \De \cos\th}{|q|^4}, \quad  \trchb=-\frac{2r}{|q|^2}, \quad \atrchb=\frac{2a\cos\th}{|q|^2},\\
&& \eta=\ze ,\quad  \om=-\frac{a^2\cos^2\th (r-m)+mr^2-a^2r}{|q|^4}=- \frac 12\pr_r\left(\frac{\Delta}{|q|^2}\right).
\eeaa
Also, we have 
\beaa
\eta_1&=& -\frac{a^2\sin\th \cos\th}{|q|^3}, \qquad  \eta_2=\frac{a r \sin\th }{|q|^3}, \qquad
\dual \eta_1= \frac{ar\sin\th }{|q|^3}, \qquad\qquad  \eta_2= \frac{a^2\sin\th \cos\th}{|q|^3},\\
\etab_1&=& -\frac{a^2\sin\th \cos\th }{|q|^3}, \qquad \etab_2 =-\frac{ar\sin\th }{|q|^3},\qquad 
\dual \etab_1= -\frac{ar\sin\th }{|q|^3}, \qquad \dual \etab_2 =\frac{a^2\sin\th \cos\th}{|q|^3}.
\eeaa

The complex Ricci coefficients    in the ingoing PN  frame      are given by  
\beaa
&& \Xh=\Xbh=\Xi=\Xib=\omb=0, \quad  \tr X=\frac{2\De \ov{q}}{|q|^4}, \quad \tr\Xb=-\frac{2}{\ov{q}},  \quad H=Z.
\eeaa

The real curvature components in any PN frame\footnote{By definition of principal null frame, $\a=\b=\bb=\aa=0$; on the other hand, $\rho$ and $\dual \rho$ are gauge invariant quantities.} are given by
\beaa
&& \a=\b=\bb=\aa=0, \quad 
\rho=-\frac{2m}{|q|^6}(r^3-3ra^2\cos^2\th), \quad \dual \rho= \frac{2am\cos\th }{|q|^6} (3r^2- a^2 \cos^2\th).
\eeaa
The complex curvature components are given by
\beaa
A=B=\Bb=\Ab=0,\quad P=-\frac{2m}{q^3}.
\eeaa
In any PN frame the components of $\B$  are given   by, see Proposition  \ref{proposition:componentsofB}, 
\bea
\bsplit
\B_{ a   b  c 3}=- \B_{ a   b  3c}&=      -  \trchb  \big( \de_{ca}\eta_b-  \de_{cb} \eta_a\big)  -  \atrchb \big( \in_{ca}  \eta_b -  \in_{cb}  \eta_a\big), \\
\B_{ a   b  c 4}=- \B_{ a   b   4 c} &= -  \trch  \big( \de_{ca}\etab_b-  \de_{cb} \etab_a\big)  -  \atrch \big( \in_{ca}  \etab_b -  \in_{cb}  \etab_a\big), \\
\B_{ a   b  3 4}=- \B_{ a   b  43} &=-4\big(\eta_a \etab_b-\etab_a\eta_b\big),\\
\B_{1212}=-\B_{1221}&=\B_{2121}= \frac 12  \trch \trchb+\frac 1 2 \atrch \atrchb.
\end{split}
\eea

  
\subsubsection{Ingoing Eddington-Finkelstein coordinates}


Let $r_0$ be a constant $r_0>r_+$. We introduce  the  adapted ingoing Eddington-Finkelstein function $\ub$ defined by\footnote{Note that the choice of $u$ and $\ub$ is such that we have $u=\ub=t$ on the timelike hypersurface $r=r_0$.} 
\beaa
\ub &=& t+\int_{r_0}^r\frac{{r'}^2+a^2}{\Delta(r')}dr'.
\eeaa
Note that
\beaa
e^{(in)}_4(\ub)=\frac{2(r^2+a^2)}{|q|^2}, \qquad e^{(in)}_3(\ub)=0, \qquad e_1(\ub)=0, \qquad e_2(\ub)=\frac{a\sin\th}{|q|}.
\eeaa

\begin{remark}
Note that  the non-vanishing of  $e_2(\ub)$ in Kerr is  connected with the lack of  integrability of the null   pair  $( e^{(in)}_3, e^{(in)}_4)$.
\end{remark}

\begin{definition}
The principal null  pair  $( e^{(in)}_3, e^{(in)}_4)$ together with the BL  function $r$, such that $ e^{(in)}_3(r)=1$, is called the canonical, ingoing,   principal geodesic  structure (PG) 
of  Kerr.   The associated, ingoing,  Eddington-Finkelstein coordinates  $(\ub, r, \th, \vphi_+)$   are given by   
\beaa
\ub:= t+ f(r), \quad f'(r)=\frac{r^2+ a^2}{ \De}, \qquad \vphi_+:= \phi +h(r), \quad h'(r)=\frac{a}{\Delta},
\eeaa
such that,
\beaa
e^{(in)}_3(r)=1, \qquad  e^{(in)}_3(\ub)=e^{(in)}_3(\th)=e^{(in)}_3(\vphi_+)=0. 
\eeaa
\end{definition}


\subsubsection{Calculations in the canonical horizontal basis}


Remind that we call \eqref{eq:canonicalHorizBasisKerr} the canonical horizontal basis of Kerr. Note that
\beaa
e_1(r)=0, \qquad e_2(r)=0.
\eeaa
Also
\beaa
\g(\D_4e_1, e_2) &=& \chi_{12}=\frac 1 2 \atrch , \qquad \g(\D_3e_1, e_2) = \chib_{12}=\frac 1 2 \atrchb,
\eeaa
which gives
\bea\label{eq:expressions-nab41}
\nab_4 e_1&=&\frac 1 2 \atrch e_2, \qquad \nab_3 e_1=\frac 1 2 \atrchb e_2.
\eea
Also
\beaa
(\La_1)_{21}:=\g(\D_1e_1, e_2) &=& 0,\\
(\La_2)_{21}:=\g(\D_2e_1, e_2) &=& \frac{r^2+a^2}{|q|^3}\cot\th,\\
(\La_1)_{12}:=\g(\D_1e_2, e_1) &=& 0,\\
(\La_2)_{12}:=\g(\D_2e_2, e_1) &=& - \frac{r^2+a^2}{|q|^3}\cot\th,
\eeaa
or 
\bea\label{eq:nab-1-1-1-2-Kerr}
\nab_{e_1} e_1 = \nab_{e_1} e_2=0, \quad \nab_{e_2} e_1 =\La e_2, \quad \nab_{e_2} e_2=-\La e_1, \quad \La := \frac{r^2+a^2}{|q|^3}\cot\th.
\eea


\subsection{Outgoing PN frame}


The  outgoing PN  frame  (with $ \D_4 e_4=0$),  regular towards the future  for all $r>r_+$,   is given by 
\bea
\lab{eq:Out.PGdirections-Kerr}
\bsplit
 e^{(out)}_4&=\frac{|q|^2}{\De} e_4^{(in)}=\frac{r^2+a^2}{\Delta}\pr_t+\pr_r+\frac{a}{\Delta}\pr_\phi,\\
  e^{(out)}_3&=\frac{\De}{|q|^2} e_3^{(in)}=\frac{r^2+a^2}{|q|^2}\pr_t-\frac{\Delta}{|q|^2}\pr_r+\frac{a}{|q|^2}\pr_\phi,
  \end{split}
\eea
with $e_1, e_2 $  the canonical horizontal basis  \eqref{eq:canonicalHorizBasisKerr}. 
Note that we have
\beaa
e^{(out)}_4(r)=1, \qquad e^{(out)}_3(r)=-\frac{\Delta}{|q|^2}, \qquad e_1(r)=0, \qquad e_2(r)=0.
\eeaa
Using the Hawking vectorfield $\That=\pr_t+\frac{a}{r^2+a^2} \pr_\phi$  we have,
\beaa
e^{(out)}_4&=&\frac{r^2+a^2}{\De} \left( \That+\frac{\De}{r^2+a^2} \pr_r \right), \qquad e^{(out)}_3=\frac{r^2+a^2}{|q|^2} \left(\That- \frac{\De}{r^2+a^2} \pr_r\right),
\eeaa
and, using also the definition of the  vectorfield $\Rhat$,
\bea
\lab{eq:Out.PGdirections-e-3e-4Kerr}
e^{(out)}_4&=&\frac{r^2+a^2}{\De} \big(\That+\Rhat\big), \qquad 
e^{(out)}_3=\frac{r^2+a^2}{|q|^2}  \big(\That-\Rhat\big).
\eea

  \begin{lemma}
  \lab{Lemma:T-frame-outgoing} 
 The following identities hold true. 
   \bea
   \lab{eq:ThatRhat-e_3e_4-outgoing}
   \begin{split}
 \That &=\frac 1 2 \left( \frac{\De}{r^2+a^2} e^{(out)}_4+\frac{|q|^2}{r^2+a^2}  e^{(out)}_3\right), \\
  \Rhat &= \frac 1 2 \left( \frac{\De}{r^2+a^2} e^{(out)}_4-\frac{|q|^2}{r^2+a^2}  e^{(out)}_3\right).
  \end{split}
 \eea
Also,
  \bea\label{eq:T-Z-outgoing-frame-Kerr}
  \begin{split}
\T &= \frac 1 2 \left( \frac{\De}{|q|^2} e^{(out)}_4+  e^{(out)}_3\right)-\frac{a\sin\th}{|q|} e_2, \\ 
\Z &= \frac{(r^2+a^2)\sin\th}{|q|} e_2-\frac 1 2a\sin^2\th \left( \frac{\De}{|q|^2} e^{(out)}_4+ e^{(out)}_3\right).
\end{split}
 \eea
  \end{lemma}
  
  \begin{proof}
Straightforward verification.
  \end{proof}

The real Ricci coefficients in the outgoing PN frame are given by
\beaa
&&\chih=\chibh=\xi=\xib=\om=0, \\
&& \trch=\frac{2r}{|q|^2},\quad \atrch=\frac{2a\cos\th}{|q|^2}, \quad  \trchb=-\frac{2r\Delta}{|q|^4}, \quad \atrchb=\frac{2a\Delta\cos\th}{|q|^4},\\
&&  \etab=-\ze,\quad  \omb=\frac{a^2\cos^2\th (r-m)+mr^2-a^2r}{|q|^4}= \frac 12\pr_r\left(\frac{\Delta}{|q|^2}\right).
\eeaa
Also, we have 
\beaa
\eta_1&=& -\frac{a^2\sin\th \cos\th}{|q|^3}, \qquad  \eta_2=\frac{a r \sin\th }{|q|^3}, \qquad
\dual \eta_1= \frac{ar\sin\th }{|q|^3}, \qquad\qquad  \eta_2= \frac{a^2\sin\th \cos\th}{|q|^3},\\
\etab_1&=& -\frac{a^2\sin\th \cos\th }{|q|^3}, \qquad \etab_2 =-\frac{ar\sin\th }{|q|^3},\qquad 
\dual \etab_1= -\frac{ar\sin\th }{|q|^3}, \qquad \dual \etab_2 =\frac{a^2\sin\th \cos\th}{|q|^3}.
\eeaa

The complex Ricci coefficients    in the outgoing PN  frame      are given by  
\beaa
&& \Xh=\Xbh=\Xi=\Xib=\om=0, \quad  \tr X=\frac{2}{q}, \quad \tr\Xb=-\frac{2\Delta q}{|q|^4},  \quad \Hb=-Z,
\eeaa
and
\beaa
H_1=\frac{a q i\sin\th\, }{|q|^3}, \quad H_2=\frac{a q \sin\th \, }{|q|^3},\quad   Z_1=\frac{a \ov{q} i\sin\th }{|q|^3},\qquad\,\, Z_2=\frac{a \ov{q} \sin\th\,}{|q|^3}.
\eeaa

\begin{remark}
Note the identities 
\beaa
H_1=-\ov{Z_1}, \quad H_2=\ov{Z_2}, \qquad H_1=\ov{\Hb_1}, \quad H_2=-\ov{\Hb_2}.
\eeaa
\end{remark}


\subsubsection{Outgoing Eddington-Finkelstein coordinates}


Let $r_0$ be a constant $r_0>r_+$. We introduce  the  adapted outgoing Eddington-Finkelstein function $u$ defined by
\beaa
u &:=& t-\int_{r_0}^r\frac{{r'}^2+a^2}{\Delta(r')}dr'.
\eeaa
Note that
\beaa
e^{(out)}_4(u)=0, \quad e^{(out)}_3(u)=\frac{2(r^2+a^2)}{|q|^2}, \quad e_1(u)=0, \quad e_2(u)=\frac{a\sin\th}{|q|},
\eeaa
and
\beaa
 \g(\D u, \D u)=\frac{a^2\sin^2\th}{|q|^2}.
\eeaa

\begin{definition}
The principal null  pair  $( e^{(out)}_3, e^{(out)}_4)$ together with the BL  function $r$, such that $ e^{(out)}_4(r)=1$, is called the canonical, outgoing,  PG  structure 
of  Kerr.   The associated, outgoing,  Eddington-Finkelstein coordinates  $(u, r, \th, \vphi_-)$   are given by   
\beaa
u:= t- f(r), \quad f'(r)=\frac{r^2+ a^2}{ \De}, \qquad \vphi_-:= \phi -h(r), \quad h'(r)=\frac{a}{\Delta},
\eeaa
such that,
\beaa
e^{(out)}_4(r)=1, \qquad  e^{(out)}_4(u)=e^{(out)}_4(\th)=e^{(out)}_4(\vphi_-)=0. 
\eeaa
\end{definition}

\begin{lemma}
Relative to  the  outgoing  Eddington-Finkelstein coordinates    $(u, r, \th, \vphi_-)$ we have:
\begin{enumerate}
\item  The action of  the  outgoing PG  frame on the coordinates $(u, r, \th, \vphi_-) $ is given by
 \bea
\lab{eq:OPG-valuesKerr}
\begin{array}{llll}
e_4(r)=1, & e_4(u)=0, & e_4(\th)=0, & e_4(\vphi_-)=0,\\[2mm]
\displaystyle e_3(r)=-\frac{\Delta}{|q|^2}, & \displaystyle e_3(u)=\frac{2(r^2+a^2)}{|q|^2},  & e_3(\th)=0,  & \displaystyle e_3(\vphi_-)=\frac{2a}{|q|^2}, \\[2mm] 
e_1(r)=0, & e_1(u)=0,  & \displaystyle e_1(\th)=\frac{1}{|q|},  & e_1(\vphi_-)=0, \\[2mm] 
e_2(r)=0, & \displaystyle e_2(u)=\frac{a\sin\th}{|q|},  & e_2(\th)=0,  & \displaystyle e_2(\vphi_-)=\frac{1}{|q|\sin\th}.
\end{array}
\eea

\item
In particular
\beaa
\left(\ba{cccc}
e_4\\
e_3\\
e_2\\
e_1
\ea
\right) &=&  \left(\ba{cccc}
1&0&0&0\\
-\frac{\De}{|q|^2} &\frac{2(r^2+a^2)}{|q|^2}  &0& \frac{2a}{|q|^2} \\
0&\frac{a\sin\th}{|q|}& 0 & \frac{1}{|q|\sin\th}\\
0&0&  \frac{1}{|q|} &  0
\ea\right)
  \left(\ba{cccc}
\pr_r \\ \pr_u \\ \pr_\th \\ \pr_{\vphi_-}
\ea\right).
\eeaa

 \item In the outgoing EF coordinates,  the metric takes the form 
 \beaa
\g &=& -\left(1-\frac{2mr}{|q|^2}\right)(du)^2 -2dr du +2a(\sin\th)^2dr d\vphi_-\\
&&-\frac{4mra(\sin\th)^2}{|q|^2}du d\vphi_- +|q|^2(d\th)^2+\frac{\Si^2(\sin\th)^2}{|q|^2}(d\vphi_-)^2.
\eeaa
\end{enumerate}
\end{lemma} 

\begin{proof}
Straightforward verification.
\end{proof}


\section{Additional relations}


Observe that, as a consequence of the null structure equations and Bianchi identities, in Kerr we have
\bea
\bsplit
& \DD P+3P\Hb=0, \qquad \DD\ov{P}+3\ov{P}H=0, \label{eq:DDP-Kerr}\\
&\DD\hot\Hb  +\Hb\hot\Hb=0, \qquad \DD\hot H  +H \hot H=0, \label{eq:DD-hot-H-Kerr}
\end{split}
\eea
which are valid in any frame.

We also have
\bea
\bsplit
\trch  \atrchb+ \trchb\atrch &=0, \label{eq:vanishing-trch-atrch}\\
|\eta|^2-|\etab|^2&= 0, \label{eq:vanishing-eta2-etab2}\\
\div(\eta - \etab) &= 0, \label{eq:vanishing-div-eta-etab} \\
\div(\dual \eta+\dual \etab)&=0, \label{eq:vanishing-div-dual-eta-etab}
\end{split}
\eea
which is immediate from the values given in the previous section.

We also have the following relations, valid in any frame:
\bea
\bsplit
\atrch e_3 + \atrchb e_4 &=  \frac{4a\cos\th (r^2+a^2)}{|q|^4}\That, \label{eq:atrch-e3-atrchb-e4-Kerr}\\
 \atrch e_3+\atrchb e_4+ 2(\eta+\etab) \c \dual \nab&=\frac{4a\cos\th}{|q|^2} \T,\label{eq:atrch-e3-atrch-e4-etaetab-kerr}\\
  \atrch  e_3+ \atrchb  e_4 -4  \La  e_2&= -\frac{4\cos\th}{|q|^2\sin^2\th}\Z.\label{eq:atrch-e3-atrch-e4-Lae2-kerr}
  \end{split}
\eea
Indeed, for example in the outgoing frame we have
\beaa
\atrch e^{(out)}_3 + \atrchb e^{(out)}_4 &=& \frac{2a\cos\th}{|q|^2} e^{(out)}_3 + \frac{2a \De \cos\th}{|q|^4} e^{(out)}_4\\
&=& \frac{4a\cos\th (r^2+a^2)}{|q|^4} \frac 1 2 \left( \frac{\De}{r^2+a^2} e^{(out)}_4+\frac{|q|^2}{r^2+a^2}  e^{(out)}_3\right)\\
&=& \frac{4a\cos\th (r^2+a^2)}{|q|^4}\That.
\eeaa
Also,
  \beaa
  2(\eta+\etab) \c \dual \nab&=& 2(\eta+\etab)_1  \dual e_1+2(\eta+\etab)_2 \dual e_2= 2(\eta+\etab)_1  e_2\\
  &=&-\frac{4a^3\sin^2\th\cos\th}{|q|^4} \T -\frac{4a^2\cos\th}{|q|^4}\Z,
  \eeaa
  which gives
  \beaa
 && \atrch e^{(out)}_3 + \atrchb e^{(out)}_4 +  2(\eta+\etab) \c \dual \nab\\
 &=&  \frac{4a\cos\th (r^2+a^2)}{|q|^4}\big( \T +\frac{a}{r^2+a^2} \Z\big)-\frac{4a^3\sin^2\th\cos\th}{|q|^4} \T -\frac{4a^2\cos\th}{|q|^4}\Z= \frac{4a\cos\th}{|q|^2} \T.
  \eeaa
  We also have, since $\La = \frac{r^2+a^2}{|q|^3}\cot\th$,
  \beaa
    \atrch  e_3+\atrchb  e_4  -4  \La  e_2&=&  \frac{4a\cos\th (r^2+a^2)}{|q|^4}\big( \T +\frac{a}{r^2+a^2} \Z\big)\\
    && - 4  \frac{r^2+a^2}{|q|^3}\cot\th \big(\frac{a\sin\th}{|q|}\T+\frac{1}{|q|\sin\th}\Z \big)\\
    &=&  \frac{4a^2\cos\th }{|q|^4}  \Z - 4  \frac{r^2+a^2}{|q|^4}\frac{\cos\th}{\sin^2\th}\Z =-\frac{4\cos\th}{|q|^2\sin^2\th}\Z.
  \eeaa


\subsection{The scalar quantity $q$}


Recall the definition \eqref{definition-q} of  $q=r+ i a \cos\th$. We have the following equations for $q$. 
\begin{lemma}\label{lemma:equations-q} 
The scalar $q$ satisfies for both the outgoing and incoming PN frames,
\bea\lab{equations:forq}
\bsplit
\nab_4 q&= \frac 1 2 \tr X q,  \quad  \nab_3 q=\frac 1 2 \ov{\tr \Xb }  q, \quad
  \DD q= q  \Hb , \quad   \DDb q=q \ov{H}, \\
   \nab q&=\frac 1 2 (\Hb +\ov{H}) q, \quad \nab (|q|^2)=(\eta+\etab) |q|^2, \quad 
  q \Hb=-\ov{q} H.
  \end{split}
\eea
In particular $ |H|^2=|\Hb|^2 
$.
\end{lemma}

\begin{proof} 
In the outgoing frame, 
from the value of $\tr X=\frac{2}{q}$, and the reduced equation $\nab_4\tr X +\frac{1}{2}(\tr X)^2 =0$
we deduce the equation for $\nab_4 q$. From the value of $P=-\frac{2m}{q^3}$ and the reduced Bianchi identity $\nab_3P  = -\frac{3}{2}\ov{\tr\Xb} P$ we deduce the equation for $\nab_3 q$. 
Similarly, equation \eqref{eq:DDP-Kerr} becomes $\DD\ov{q} =H\ov{q}$
  or $ \DDb q=q \ov{H}$. The last equation in \eqref{equations:forq}  follows in the same manner from \eqref{eq:DDP-Kerr}.
  Similarly in the incoming frame.

We write, recalling that $q=r+ia \cos\th$ and using $e_a(r)=0$, i.e. $\DD r=0$, 
  \beaa
   q  \Hb&=& \DD q=\DD( r+i a \cos\th) =\nab(r+i a \cos \th)+ i \dual \nab(r+i a \cos\th)\\
   &=&ia \nab\cos\th -a \dual\nab \cos\th, \\
   q \ov{H}&=&\DDb q =\DDb( r+i a \cos\th)  =\nab(r+i a \cos \th)- i \dual \nab(r+i a \cos\th)\\
   &=& ia \nab\cos\th +a \dual\nab \cos\th.
  \eeaa
   We deduce $ \ov{ q  \Hb}= - q \ov{H}$ i.e. $  q\Hb=-\ov{q} H$  as stated.
 \end{proof}


\subsection{The canonical complex 1-form $\Jk$}


\begin{definition}
\lab{def:JkandJ}
We define the  following   complex horizontal  $1$-tensor in Kerr,  given  in components relative to $e_1, e_2$ by
\beaa
\Jk_1=\frac{i\sin\th}{|q|}, \qquad \Jk_2=\frac{\sin\th}{|q|}.
\eeaa
Note that $\Jk $ is regular (even at the axis), as well as anti-selfadjoint, i.e. $\dual \Jk=-i \Jk$.
\end{definition}

\begin{remark}
The complex 1-tensor $\Jk$ can also be written as 
\beaa
\Jk&=& j+ i \dual j
\eeaa
with 
\beaa
j_1=-\dual j_2=0, \qquad j_2=\dual j_1=\frac{\sin\th}{|q|}.
\eeaa
\end{remark}

Using the  form $\Jk$ we can rewrite  the expressions for  $H$, $\Hb$ in the form
\bea
\lab{eq:expressionsforHHBZ-usingJk}
\Hb&=&-\frac{a}{q}\Jk= -\frac{a\ov{q}}{|q|^2}\Jk, \qquad\, H=\frac{a}{\ov{q}}\Jk=\frac{aq}{|q|^2}\Jk.
\eea
In particular,
\beaa
q\Hb+\ov{q} H=0,
\eeaa
as obtained in \eqref{equations:forq}.

\begin{lemma}
 The complex $1$-form $\Jk$ verifies
\beaa
\qquad \dual \Jk=-i \Jk, \qquad \Jk\c\ov{\Jk} = \frac{2(\sin\th)^2}{|q|^2}, \quad |\Re(\Jk)|^2 = \frac{(\sin\th)^2}{|q|^2}.
\eeaa
 Also
\bea
\lab{transport-Jk-Kerr}
\nab_4\Jk+\frac 1 2 \tr X \Jk=0, \qquad \nab_3\Jk+\frac 1 2 \tr \Xb  \Jk=0,
\eea
or
\beaa
\nab_4(q \Jk)=\nab_3(\ov{q} \Jk )=0.
\eeaa
Also 
\beaa
\ov{\DD}\c\Jk = \frac{4i(r^2+a^2)\cos\th}{|q|^4},\qquad \DD\hot\Jk = 0,
\eeaa
and
\beaa 
\DD (q)= -a\Jk, \qquad \DD(\ov{q})= a \Jk.
\eeaa
\end{lemma}

\begin{proof}
Straightforward verification. 
  One  can  check \eqref{transport-Jk-Kerr} using the complex structure  equations in Kerr
\beaa
\nab_3 \Hb=-\frac{1}{2} \ov{\tr \Xb}  (\Hb-H), \qquad \nab_4 H=-\frac 1 2 \ov{\tr X} (H-\Hb), 
\eeaa
and \eqref{equations:forq} together with the expressions \eqref{eq:expressionsforHHBZ-usingJk}  of $H,\Hb$ in terms of $\Jk$.
Note that \eqref{transport-Jk-Kerr}   holds true in both the incoming and outgoing frame. In the  outgoing frame it becomes
\beaa
\nab_4\Jk = -\frac{1}{q}\Jk,     \qquad
\nab_3\Jk = \frac{\De q}{|q|^4}\Jk.
\eeaa
From the complex structure equations in Kerr \eqref{eq:DD-hot-H-Kerr}
and \eqref{equations:forq}, we deduce that $\DD \hot \Jk=0$.
We then compute
\beaa
\frac 1 2 \DDb \c \Jk&=& \div j + i \curl j= \nab_1 j_1+\nab_2 j_2 + i (\nab_1 j_2-\nab_2 j_1)\\
&=& - j_{\nab_11}- j_{\nab_22} + i (e_1( j_2)-j_{\nab_1 2}+j_{\nab_2 1})=  i (e_1( j_2)+\La j_{2})\\
&=&  i \left(\frac{1}{|q|}\partial_\th\left( \frac{\sin\th}{|q|}\right)+ \frac{r^2+a^2}{|q|^3}\cot\th \frac{\sin\th}{|q|}\right)\\
&=&  i \left(\frac{\cos\th}{|q|^4}(r^2+a^2)+ \frac{r^2+a^2}{|q|^4}\cos\th\right)=2i \frac{r^2+a^2}{|q|^4}\cos\th,
\eeaa
as stated.
\end{proof}

Using the canonical form $\Jk$, we can also deduce the following identities, which are not immediately implied by the null structure equations.

\begin{lemma}\label{lemma:additional-relations-Kerr}
In Kerr,  relative to any frame, the following relations hold true:
\begin{align}
\begin{split}
\nab_4\Hb+\tr X\Hb&=0, \qquad \nab_3 H + \tr \Xb H=0, \\
 \ov{\tr X} \Hb + \tr X H&=0, \qquad   \ov{\tr \Xb} H + \tr \Xb \Hb=0,\\
  \DDc\tr X +2\tr X\Hb&=0, \qquad  \DDc\tr \Xb +2\tr \Xb H=0. \label{eq:vanishing-relations-Kerr}
 \end{split}
\end{align}
The above can also be written as
\begin{align}\label{eq:vanishing-relations-Kerr-real}
\begin{split}
\nab_4 \etab+\trch \etab +\atrch \dual \etab&=0, \qquad \nab_3 \eta+\trchb \eta +\atrchb \dual \eta=0,\\
\trch (\eta + \etab) + \atrch (\dual \eta -\dual \etab)&=0, \qquad \trchb (\eta + \etab) + \atrchb (\dual \etab -\dual \eta)=0,
\end{split}
\end{align}
and
\begin{align}
\begin{split}
\nabc \trch +\dual \nabc \atrch +  2\trch \etab+  2\atrch \dual \etab&=0, \\
 \nabc \trchb +\dual \nabc \atrchb + 2 \trchb \eta+  2\atrchb \dual \eta&=0.
 \end{split}
\end{align}
We also have
\beaa
\nab \tr\Xb&=& -  \left(\frac 1 2\Hb+ \frac 1 2\ov{H}+\ov{\Hb}+ H\right)\tr\Xb,
\eeaa
which can be written as
\bea\label{eq:vanishing-relations-Kerr-real-2}
\begin{split}
\nab (\trchb)&= -\frac 3 2 \trchb   \left(\etab + \eta\right) -\frac 1 2 \atrchb  \left(  \dual \eta-   \dual \etab\right),\\
\nab(\atrchb)&= -\frac 3 2  \atrchb  (\etab+ \eta) +\frac 1 2\trchb  (\dual \eta-  \dual \etab ).
\end{split}
\eea
\end{lemma}

\begin{proof}
One can check, using \eqref{eq:expressionsforHHBZ-usingJk} and \eqref{transport-Jk-Kerr}, that
\beaa
\nab_4(q \Hb)&=&-a \nab_4 \Jk=\frac 1 2 a \tr X \Jk=-\frac 1 2  \tr X q \Hb.
\eeaa
On the other hand, using \eqref{equations:forq} we have
\beaa
\nab_4(q \Hb)&=&\frac 1 2 \tr X q \Hb + q \nab_4 \Hb,
\eeaa
which implies the first relation. Similarly for $H$.
Also, we have
\beaa
\ov{\tr X} \Hb&=& \frac{2}{\ov{q}} \left(-\frac{a\ov{q}}{|q|^2}\Jk\right)=-\frac{2a}{|q|^2}\Jk, \\
\tr X H&=& \frac{2}{q}  \left(\frac{aq}{|q|^2}\Jk\right)=\frac{2a}{|q|^2}\Jk,
\eeaa
which implies the second relation, and similarly for $\tr\Xb$. 
Finally, using \eqref{equations:forq} we compute
\beaa
\DDc(\tr X) &=& \DD\left(\frac{2}{q}\right)+\tr XZ=-\frac{2}{q^2}\DD(q)+\frac{2a}{q^2}\Jk=\frac{4a}{q^2}\Jk=-2\tr X\Hb,
\eeaa
as desired. Similarly for $\tr \Xb$.
Taking the real parts of the expressions we obtain \eqref{eq:vanishing-relations-Kerr-real}. 

To obtain the last relation, we use that $\nab q=\frac 1 2 q (\Hb+ \ov{H})$, and since $\tr\Xb= -\frac{2\Delta}{q \ov{q}^2}$ we have
\beaa
\nab \tr\Xb &=&  \frac{2\Delta}{q^2 \ov{q}^2}\nab q  +\frac{4\Delta}{q \ov{q}^3}\nab\ov{q}=  \frac{2\Delta}{q \ov{q}^2}\frac 1 2  (\Hb+ \ov{H})  +\frac{2\Delta}{q \ov{q}^2} (\ov{\Hb}+ H)\\
&=&  \frac{2\Delta}{q \ov{q}^2}  (\frac 1 2\Hb+ \frac 1 2\ov{H}+\ov{\Hb}+ H)= -\tr\Xb  \left(\frac 1 2\Hb+ \frac 1 2\ov{H}+\ov{\Hb}+ H\right).
\eeaa
By taking the real and the imaginary part we obtain the stated identities.
\end{proof}


\section{Inverse Kerr metric and  Killing tensors}


We summarize here a computation to write the inverse of the Kerr metric which will be crucial in Chapter \ref{chapter-proof-mor-2}.

\begin{lemma}
\lab{lemma:inversemetricexpressioninKerr}
The inverse Kerr metric can be written in the form
\bea\label{inverse-metric}
|q|^2 \g^{\a\b}=\De \pr_r^\a \pr_r^\b +\frac{1}{\De}\RR^{\a\b}
\eea
with
\bea
\lab{eq:expressionRR-O}
\bsplit
\RR^{\a\b}&= -(r^2+a^2) ^2\pr_t^\a\pr_t^\b- 2a(r^2+a^2) \pr_t^{(\a} \pr_\phi^{\b)}- a^2 \pr_\phi^\a \pr_\phi^\b  +\De O^{\a\b},\\
 O^{\a\b}&= \pr_\th^\a \pr_\th^\b +\frac{1}{\sin^2\th} \pr_\phi^\a \pr_\phi^\b +  2 a \pr_t^{(\a} \pr_\phi^{\b)}+ a^2\sin^2 \th \pr_t^\a \pr_t^\b.
\end{split}
\eea
Note that
\bea
\bsplit
\lab{eq:exxpressionO-e1e2}
\lab{eq:expressionRR-TT}
\RR^{\a\b} = -(r^2+a^2) ^2\That^\a \That^\b  +\De O^{\a\b}, \qquad O^{\a\b}= |q|^2 \big( e_1^\a e_1^\b+ e_2^\a e_2^\b\big),
 \end{split}
\eea
 thus the inverse metric can also be written in the form
 \bea
 \lab{inverse-metric-vfs}
|q|^2 \g^{\a\b}&=&\frac{(r^2+a^2)^2}{\De} \big( -\That^\a\That^\b +      \Rhat^\a \Rhat^\b\big) + O^{\a\b}.
 \eea
\end{lemma}

\begin{proof}
From the expression of the Kerr metric, the inverse metric can be written in the form
\beaa
|q|^2 \g^{\a\b}=\De \pr_r^\a \pr_r^\b +\frac{1}{\De}\RR^{\a\b}
\eeaa
with
\beaa
\RR^{\a\b}&=&-\Si^2\pr_t^\a\pr_t^\b-2 amr \pr_t^\a\pr_\phi^\b  -2 amr \pr_\phi^\a\pr_t^\b+\De \pr_\th^\a \pr_\th^\b+\frac{\De-a^2\sin^2\th}{\sin^2\th} \pr_\phi ^\a\pr_\phi^\b
\eeaa
which establishes \eqref{eq:expressionRR-O}.  
According to the  definition \eqref{define:That} of $\That$, we can write 
\beaa
\RR^{\a\b}&=& -(r^2+a^2) ^2 \left( \pr_t^\a\pr_t^\b+ \frac{2a}{r^2+a^2} \pr_t^{(\a} \pr_\phi^{\b)}+ \frac{a^2}{(r^2+a^2)^2} \pr_\phi^\a \pr_\phi^\b \right) +\De O^{\a\b}\\
&=& -(r^2+a^2) ^2 \That^\a \That^\b +\De O^{\a\b}
\eeaa
which establishes  the first  expression in \eqref{eq:expressionRR-TT}. Finally the second expression \eqref{eq:exxpressionO-e1e2}  can be easily checked  from the expressions   of $e_1, e_2$ in terms of the BL coordinates in \eqref{eq:Out.PGdirections-Kerr}.  
\end{proof}

\begin{remark}\label{remark:d-ff-A-B} 
Observe that in \cite{A-B}, Andersson-Blue use  instead the following expression for $\RR$:
\bea
\lab{eq:expressionRR-Q}
\RR^{\a\b}&= -(r^2+a^2) ^2\pr_t^\a\pr_t^\b- 4 amr \pr_t^{(\a}\pr_\phi^{\b)}+(\De-a^2)   \pr_\phi ^\a\pr_\phi^\b+\De Q^{\a\b}
\eea
where
\bea
\lab{eq:exxpressionQ-O}
Q^{\a\b} &=& O^{\a\b}-  2 a \pr_t^\a \pr_\phi^\b- \pr_\phi^\a \pr_\phi^\b =\pr_\th^\a \pr_\th^\b +\frac{\cos^2\th}{\sin^2\th} \pr_\phi^\a \pr_\phi^\b + a^2\sin^2 \th \pr_t^\a \pr_t^\b.
\eea
\end{remark}

The relevance of the decomposition of the metric in \eqref{inverse-metric} is in the fact that the operator $\RR^{\a\b}$ can be written in terms of $\partial_t^\a \partial_t^\b$, $a\partial_t^{(\a} \partial_\phi^{\b)}$, $a^2\partial_\phi^\a \partial_\phi^\b$ and $O^{\a\b}$.

\begin{definition}\label{definition:tensors-S} 
We define the following symmetric spacetime 2-tensors 
\beaa
S_1^{\a\b}&:=& \T^\a \T^\b=\partial_t^\a \partial_t^\b, \\
S_2^{\a\b}&:=& a \T^{(\a} \Z^{\b)}=a\partial_t^{(\a} \partial_\phi^{\b)}, \\
S_3^{\a\b}&:=& a^2\Z^\a \Z^\b=a^2\partial_\phi^\a \partial_\phi^\b, \\
S_4^{\a\b}&:=& O^{\a\b} =  |q|^2 \big( e_1^\a e_1^\b+ e_2^\a e_2^\b\big).
\eeaa
We denote the set of the above tensors as $S_\aund$, for $\aund=1,2,3,4$. 
\end{definition}

\begin{remark}\label{remark:why-a-SS} 
Observe that the tensor $S_2$ and $S_3$ are defined with a factor of $a$ and $a^2$ respectively. The reason for this choice, which differs from the definition in \cite{A-B}, is our application in Chapter \ref{chapter-proof-mor-2} of the method to 2-tensors as opposed to scalars.  Note also that  $S_1-S_3$ are Killing  tensor  while $S_4$ is    related to the Carter tensor $K$, see section \ref{section:CartertensorKerr}.
\end{remark}

With the above definition, from \eqref{eq:expressionRR-O} we write
\beaa
\RR^{\a\b}&= -(r^2+a^2) ^2S_1^{\a\b}- 2(r^2+a^2)S_2^{\a\b}- S_3^{\a\b}  +\De S_4^{\a\b}.
\eeaa
More compactly, using the repetition in $\aund$ to signify summation over $\aund=1,2,3,4$, we denote
          \bea\label{eq:RR-Sa}
          \RR^{\a\b} =\RR^\aund  \Sa^{\a\b},
          \eea
          with $\RR^\aund$, $\aund=1,2,3,4$, given by
          \bea\label{components-RR-aund}
          \RR^1&=&-(r^2+a^2)^2, \quad \RR^2 = -2(r^2+a^2), \quad \RR^3 =-1, \quad \RR^4=\De.
          \eea


\section{Commutation properties for $\T$ and $\Z$}


We start by collecting the following commutation property between $\nab_\T$ and $\nab_\Z$.

\begin{lemma}\label{lemma:commutation-nabT-nabZ} 
In Kerr spacetime, for $\psi \in \sk_2$ we have
\beaa
[\nab_\T, \nab_\Z]\psi=0.
\eeaa
\end{lemma}

\begin{proof} 
Recall,    see Corollary \ref{Corollory:Commutatornab_anab_b}, that for $X$ and $Y$ vectorfield 
\beaa
\big(\nab_X \nab_Y -\nab_Y\nab_X) \psi_{ab}=\nab_{[X, Y]} \psi_{ab}+    X^\mu Y^\nu  (\Rdot_{ac\mu\nu} {\psi^{c}}_{b}+   \Rdot_{bc\mu\nu} {\psi^{c}}_{a}).
\eeaa
Since $[\T, \Z]=0$, we obtain
\beaa
\big(\nab_\T \nab_\Z -\nab_\Z\nab_\T) \psi_{ab}=    \T^\mu \Z^\nu  (\Rdot_{ac\mu\nu} {\psi^{c}}_{b}+   \Rdot_{bc\mu\nu} {\psi^{c}}_{a}),
\eeaa
with $ \Rdot_{ab   \mu\nu}= \R_{ab    \mu\nu}+ \frac 1 2  \B_{ab   \mu\nu}$.
Using that the only non-vanishing Riemann curvature terms are
\beaa
 \R_{a3b4}&=&-\rho\de_{ab} +\dual \rho \in_{ab},\qquad  \R_{ab34}= 2 \in_{ab}\dual\rho, \qquad \R_{abcd}=-\in_{ab}\in_{cd}\rho,
 \eeaa
 we obtain, using from \eqref{eq:T-Z-ingoing-frame-Kerr} that $\T^3\Z^4=\T^4\Z^3$, 
 \beaa
\T^\mu \Z^\nu  \R_{ac\mu\nu} {\psi^{c}}_{b}&=&    \T^d \Z^e  \R_{acde} {\psi^{c}}_{b}+ \T^3 \Z^4 \R_{ac34} {\psi^{c}}_{b}+ \T^4 \Z^3 \R_{ac43} {\psi^{c}}_{b}\\
&=&  -2\rho \in_{de}  \T^d \Z^e \dual \psi_{ab}+4\dual \rho( \T^3 \Z^4-\T^4 \Z^3) \dual \psi_{ab}=0,
\eeaa
and similarly $\T^\mu \Z^\nu  \R_{bc\mu\nu} {\psi^{c}}_{a}=0$.
We also have 
\beaa
\T^\mu \Z^\nu  \B_{ac\mu\nu} {\psi^{c}}_{b}&=&\T^d \Z^e  \B_{acde} {\psi^{c}}_{b}+\T^d \Z^3  \B_{acd3} {\psi^{c}}_{b}+\T^d \Z^4  \B_{acd4} {\psi^{c}}_{b}\\
&&+\T^3 \Z^d  \B_{ac3d} {\psi^{c}}_{b}+\T^3 \Z^4  \B_{ac34} {\psi^{c}}_{b}+\T^4 \Z^d  \B_{ac4d} {\psi^{c}}_{b}+\T^4 \Z^3  \B_{ac43} {\psi^{c}}_{b}\\
&=&(\T^d \Z^3 -\T^3 \Z^d) \B_{acd3} {\psi^{c}}_{b}+(\T^d \Z^4-\T^4\Z^d)  \B_{acd4} {\psi^{c}}_{b}\\
&=&(\T^d \Z^3 -\T^3 \Z^d)\big(  -  \trchb  \big( \de_{da}\eta_c-  \de_{dc} \eta_a\big)  -  \atrchb \big( \in_{da}  \eta_c -  \in_{dc}  \eta_a\big) \big) {\psi^{c}}_{b}\\
&&+(\T^d \Z^4-\T^4\Z^d) \big(  -  \trch  \big( \de_{da}\etab_c-  \de_{dc} \etab_a\big)  -  \atrch \big( \in_{da}  \etab_c -  \in_{dc}  \etab_a\big) \big) {\psi^{c}}_{b}.
\eeaa
Using that, in the ingoing frame, 
\beaa
&&\T^2\Z^3-\T^3 \Z^2=-\frac{\De \sin\th}{2|q|}, \qquad \T^2\Z^4-\T^4 \Z^2=-\frac{\sin\th}{2} |q|,\\
&& \trch=\frac{2\De r}{|q|^4},\quad \atrch=\frac{2a \De \cos\th}{|q|^4}, \quad  \trchb=-\frac{2r}{|q|^2}, \quad \atrchb=\frac{2a\cos\th}{|q|^2},
\eeaa
we have
\beaa
\T^\mu \Z^\nu  \B_{ac\mu\nu} {\psi^{c}}_{b}&=&-\frac{\De  r \sin\th}{|q|^3} \big( \de_{2a}(\eta_c-\etab_c) -  \de_{2c} (\eta_a-\etab_a)\big) {\psi^{c}}_{b}\\
&&+\frac{\De a\cos\th\sin\th}{|q|^3} \big( \in_{2a} ( \eta_c+\etab_c) -  \in_{2c}  (\eta_a+\etab_a) \big) {\psi^{c}}_{b}.
\eeaa
In components, recalling that $\eta_1-\etab_1=0$, $\eta_2+\etab_2=0$, this gives 
\beaa
\T^\mu \Z^\nu  \B_{1c\mu\nu} {\psi^{c}}_{b}&=&\frac{\De a\cos\th\sin\th}{|q|^3} \big( \in_{21} ( \eta_c+\etab_c){\psi^{c}}_{b} -  \in_{2c}  (\eta_1+\etab_1){\psi^{c}}_{b} \big)  \\
&=&\frac{\De a\cos\th\sin\th}{|q|^3} \big( \in_{21} ( \eta_1+\etab_1){\psi^{1}}_{b} -  \in_{21}  (\eta_1+\etab_1){\psi^{1}}_{b} \big)=0, \\
\T^\mu \Z^\nu  \B_{2c\mu\nu} {\psi^{c}}_{b}&=&-\frac{\De  r \sin\th}{|q|^3} \big( (\eta_c-\etab_c) {\psi^{c}}_{b}-  \de_{2c} (\eta_2-\etab_2){\psi^{c}}_{b}\big) \\
&=&-\frac{\De  r \sin\th}{|q|^3} \big( (\eta_2-\etab_2) {\psi^{2}}_{b}-  (\eta_2-\etab_2){\psi^{2}}_{b}\big) =0.
\eeaa
We therefore infer $\T^\mu \Z^\nu  (\Rdot_{ac\mu\nu} {\psi^{c}}_{b}+   \Rdot_{bc\mu\nu} {\psi^{c}}_{a})=0$, proving the lemma. For a different proof in perturbations of Kerr spacetime, see also Corollary \ref{cor:basicpropertiesLiebTfasdiuhakdisug:chap9}.
\end{proof}

As Killing vectorfields, $\T$ and $\Z$ also have favorable properties regarding commutation with the D'Alembertian operator $\square_\g$ for scalars. Nevertheless, the commutation with the D'Alembertian operator $\squared_2$ for tensors presents in addition terms involving the Riemann curvature. In the case of $\Z$, those term do not vanish even in Schwarzschild.

\begin{proposition}\label{prop:commutators-nabT-nabZ-squared2-Kerr}
The first order differential operators $\nab_\T$ and $\nab_\Z$ satisfy the following commutation relations with $\squared_2$ for $\psi \in \sk_2$:
\beaa
\, [\nab_\T, \squared_2]\psi&=& O(ar^{-4})\dk^{\leq 1}\psi , \\
\, [\nab_\Z, \squared_2]\psi&=&O(mr^{-2})\dk^{\leq 1}\psi.
\eeaa
\end{proposition}

\begin{proof} 
We first derive the following general computation, which specializes Lemma \ref{LEMMA:COMMUTATIONNAB_XSQUARED} to the case of Kerr spacetime. 
\begin{lemma}\label{lemma:squared-nabX-kerr}
In Kerr spacetime we have for $\psi\in\sk_2$:
\beaa
[\squared_2,    \nab_X]  \psi &=& \pi^{\mu\nu}  \Db_\mu \Db_\nu \psi+
\left( \D^\mu  \pi_\mu \,^\b-\frac 1 2  \D^\b  \tr \pi\right) \Db_\b \psi \\
&& +O(am r^{-4})\big(  X^3   \nab_3 - X^4  \nab_4\big) \dual \psi +O(mr^{-3})  X^a \nab \dual  \psi\\
&&+O(am r^{-4})\big(  \D^3 X^4 -  \D^4 X^3\big)\dual   \psi+O(mr^{-3}) \in_{ab} \D^a X^b  \dual  \psi  \\
&&+O(ar^{-3}) (X^3+X^4)\big(   \nab \dual  \psi +r^{-1} \dual \psi \big)+O(r^{-2})  X^a \big(   \nab \dual  \psi +r^{-1} \dual \psi \big)\\
&&+  O(ar^{-3})(  \D^3 X^d +\D^4 X^d) \dual \psi.
\eeaa
\end{lemma}

\begin{proof} 
According to Lemma \ref{LEMMA:COMMUTATIONNAB_XSQUARED} we have for $\psi \in \sk_2$, 
\beaa
[\squared_2,    \nab_X]  \psi_{ab} &=& \pi^{\mu\nu}  \Db_\mu \Db_\nu \psi_{ab}+
\left( \D^\mu  \pi_\mu \,^\b-\frac 1 2  \D^\b  \tr \pi\right) \Db_\b \psi_{ab} \\
&& - 2 X^\b  \R_{ac \b\mu} \Db^\mu \psi_{cb}- 2 X^\b  \R_{bc \b\mu} \Db^\mu \psi_{ac}+  \D^\b X^\mu \R_{ac \b\mu}  \psi^{cb}+  \D^\b X^\mu \R_{bc \b\mu}  \psi^{ac}  \\
&&-  X^\b \B_{ac \b\mu}  \Db^\mu \psi_{cb}-  X^\b \B_{bc \b\mu}  \Db^\mu \psi_{ac}+ \frac 12  \D^\b X^\mu \B_{ac \b\mu} \psi^{cb}+ \frac 12  \D^\b X^\mu \B_{bc \b\mu} \psi^{ac}\\
&& +\frac 1 2  X^\b \D^\mu \B_{a  c\mu \b} \psi^{cb} +\frac 1 2  X^\b \D^\mu \B_{b  c\mu \b} \psi^{ac}.
\eeaa
Using that the only non-vanishing Riemann curvature terms are
\beaa
 \R_{a3b4}&=&O(mr^{-3})\de_{ab} +O(am r^{-4}) \in_{ab},\qquad  \R_{ab34}= O(am r^{-4}) \in_{ab}, \\
  \R_{abcd}&=&O(mr^{-3})\in_{ab}\in_{cd},
 \eeaa
and the only non-vanishing $\B$ terms are
\beaa
\B_{ a   b  c 3}&=&    O(ar^{-3}) , \qquad \B_{ a   b  c 4}=O(ar^{-3}) , \qquad \B_{ a   b  3 4}= O(a^2r^{-4}), \qquad \B_{1212}=O(r^{-2}),
\eeaa
we have
\beaa
[\squared_2,    \nab_X]  \psi&=& \pi^{\mu\nu}  \Db_\mu \Db_\nu \psi+
\left( \D^\mu  \pi_\mu \,^\b-\frac 1 2  \D^\b  \tr \pi\right) \Db_\b \psi \\
&& +O(am r^{-4})\big(  X^3   \nab_3 - X^4  \nab_4\big) \dual \psi +O(mr^{-3})  X^a \nab \dual  \psi\\
&&+O(am r^{-4})\big(  \D^3 X^4 -  \D^4 X^3\big)\dual   \psi+O(mr^{-3}) \in_{ab} \D^a X^b  \dual  \psi  \\
&&+O(ar^{-3}) (X^3+X^4)\big(   \nab \dual  \psi +r^{-1} \dual \psi \big)+O(r^{-2})  X^a \big(   \nab \dual  \psi +r^{-1} \dual \psi \big)\\
&&+  O(ar^{-3})(  \D^3 X^d +\D^4 X^d) \dual \psi,
\eeaa
as stated. This concludes the proof of Lemma \ref{lemma:squared-nabX-kerr}.
\end{proof}
 
Using Lemma \ref{lemma:squared-nabX-kerr}, with vanishing $^{(\T)}\pi$ and $^{(\Z)}\pi$, and writing that $\T^3   \nab_3 - \T^4  \nab_4=\nab_\Rhat$ and $\Z^3   \nab_3 - \Z^4  \nab_4=a \nab_\Rhat$, and $\T^a=O(ar^{-1})$, $\Z^a=O(r)$, we obtain
\beaa
[ \nab_\T, \squared_2]  \psi &=&O(am r^{-4})\nab_\Rhat \dual \psi +O(amr^{-4})  \nab \dual  \psi+O(am r^{-5})\dual   \psi  \\
&&+O(ar^{-3})\big(   \nab \dual  \psi +r^{-1} \dual \psi \big),
\eeaa
and
\beaa
[\nab_\Z, \squared_2]  \psi &=&O(a^2m r^{-4})\nab_\Rhat \dual \psi +O(mr^{-2})  \nab \dual  \psi+O(am r^{-5})\dual   \psi \\
&&+O(r^{-1})  \big(   \nab \dual  \psi +r^{-1} \dual \psi \big),
\eeaa
which can be schematically written as stated. This concludes the proof of Proposition \ref{prop:commutators-nabT-nabZ-squared2-Kerr}.
\end{proof}


\section{Carter tensor and Carter operator}
\lab{section:CartertensorKerr}


Recall, see Definition \ref{def:killing-tensor}, that a Killing 2-tensor $K$ is a symmetric 2-tensor satisfying $\D_{(\mu} K_{\a\b)}=0$.
One of the fundamental properties of the Kerr metric is that it admits a Killing tensor $K$ which cannot be written in terms of the Killing vectorfields $\T$ and $\Z$. This 2-tensor  is called the Carter tensor \cite{Carter}.

\begin{definition}[Carter tensor]\label{definition:Carter-tensor-Kerr} 
In Kerr spacetime, the Carter tensor is defined as
\bea\label{definition-K-kerr}
K^{\a\b}&=& -(a^2\cos^2\th)  \g^{\a\b} +O^{\a\b} 
\eea
where  the tensor $O$ is defined in \eqref{eq:expressionRR-O}.
\end{definition}

Note that the only non-vanishing components of $K$ are:
\beaa
K_{ab}&=& r^2 \de_{ab}, \qquad  K_{34}= 2a^2\cos^2\th. 
\eeaa

\begin{proposition}\label{PROP:COMPUTATIONS-PI} 
The Carter tensor defined in \eqref{definition-K-kerr} is a Killing tensor of the Kerr metric, i.e. $\D_{(\mu} K_{\nu \rho)}=0$.
\end{proposition}

\begin{proof}  
See Section \ref{sec:proof-prop-Pi}.
\end{proof}

Associated to a Killing tensor $K^{\a\b}$, one can construct a corresponding second order operator $\KK$ for horizontal tensors $\psi \in \sk_k$ according to Definition \ref{definition:operator-KK-general}, given by
\beaa
\KK(\psi):=\Ddot_\b(K^{\a\b}\Ddot_\a \psi),
\eeaa
which has the fundamental property of commuting with the D'Alembertian operator for scalars in a vacuum spacetime, see Proposition \ref{COMMUTATION-KK-SQUARE}. Its explicit expression in Kerr for $K$ as in \eqref{definition-K-kerr} is given in the following.

\begin{proposition}\label{PROP:CARTER-OPERATOR} 
In Kerr spacetime, the Carter operator $\KK$ for $\psi \in \sk_k$ is given by
\bea
\KK &=& -( a^2 \cos^2 \th ) \squared_k    + \OO \lab{eq: BL.ExpressionOO}
\eea
where $\OO$ is the following second order angular operator:
\bea
\OO(\psi)&:=& |q|^2\big( \lap_k \psi +(\eta+\etab) \c \nab \psi).\label{def-OO-in-frames}
\eea
\end{proposition}

\begin{proof} 
See Section \ref{sec:proof-carter-operator}.
\end{proof}

\begin{lemma}
In Kerr, the second order operator $\OO$ defined in \eqref{def-OO-in-frames} is equivalent to 
\bea
\OO(\psi) &=& |q|^2 \left(\lap_k \psi - \frac{2a^2\cos\th}{|q|^2} \dual\Re(\Jk)^b \nab_b \psi   \right). \label{eq:OO-Jk-kerr}
\eea
Also, the operator $\OO$ can be written in the following ways, all equivalent to \eqref{def-OO-in-frames}:
\bea
\OO(\psi) &=& |q|^2\left(\lap_k \psi+\frac{2\nab(|q|)}{|q|}\c\nab \psi \right), \label{eq:OO-lap-nab|q|} \\
\OO(\psi)&=& \nab \c \left( |q|^2 \nab \psi \right), \label{eq:OO-nab-nab-|q|2}\\
\OO(\psi)&=& \Ddot_\b(O^{\a\b}\Ddot_\a \psi)- \nab(a^2\cos^2\th) \c  \nab \psi,\label{eq-OO-D-O-}\\
\OO(\psi)&=&  |q|^2\Ddot_\b(|q|^{-2}O^{\a\b}\Ddot_\a \psi). \label{eq:OO-psi-Ddot-S_4}
\eea
\end{lemma}

\begin{proof} 
Using \eqref{eq:expressionsforHHBZ-usingJk}, we see that in Kerr
\beaa
\eta+\etab= \Re \big( H+\Hb \big)=\Re\left(\frac{aq}{|q|^2}\Jk -\frac{a\ov{q}}{|q|^2}\Jk\right)=\frac{a}{|q|^2} \Re\big(2ia \cos\th \Jk\big)=\frac{2a^2\cos\th}{|q|^2} \Re\big(i  \Jk\big),
\eeaa
and since $\dual \Jk = -i\Jk$, we obtain
\beaa
\eta+\etab=-\frac{2a^2\cos\th}{|q|^2} \Re(\dual \Jk)=-\frac{2a^2\cos\th}{|q|^2} \dual\Re(\Jk),
\eeaa
which gives \eqref{eq:OO-Jk-kerr}. 
Using from \eqref{equations:forq} that $\nab(|q|)=\frac 1 2 (\eta+\etab) |q| $ and $\nab(|q|^2)=(\eta+\etab)  |q|^2$, we obtain \eqref{eq:OO-lap-nab|q|} and \eqref{eq:OO-nab-nab-|q|2}.

Using the expression of the Carter tensor $K^{\a\b}$ given by \eqref{definition-K-kerr}, we obtain
\beaa
\KK(\psi)&=&\Ddot_\b( -(a^2\cos^2\th)  \g^{\a\b}\Ddot_\a \psi +O^{\a\b}\Ddot_\a \psi)\\
&=& -(a^2\cos^2\th)  \g^{\a\b}\Ddot_\b\Ddot_\a \psi  -\Ddot_\b(a^2\cos^2\th)  \g^{\a\b}\Ddot_\a \psi +\Ddot_\b(O^{\a\b}\Ddot_\a \psi)\\
&=& -(a^2\cos^2\th)  \squared_k \psi  +\Ddot_\b(O^{\a\b}\Ddot_\a \psi)- \g^{\a\b}\Ddot_\b(a^2\cos^2\th) \Ddot_\a \psi.
\eeaa
By comparing the above with \eqref{eq: BL.ExpressionOO} we deduce \eqref{eq-OO-D-O-}.
Using that $\nab (a^2\cos^2\th)=\nab(r^2+a^2\cos^2\th)=\nab(|q|^2)$, we deduce from the above
\beaa
\OO(\psi)&=& |q|^2\Ddot_\b(|q|^{-2}O^{\a\b}\Ddot_\a \psi)- |q|^2O^{\a\b}\Ddot_\b(|q|^{-2})\Ddot_\a \psi- \nab(a^2\cos^2\th) \c  \nab \psi\\
&=& |q|^2\Ddot_\b(|q|^{-2}O^{\a\b}\Ddot_\a \psi)+ |q|^{-2}O^{\a\b}\Ddot_\b(|q|^{2})\Ddot_\a \psi- \nab(a^2\cos^2\th) \c  \nab \psi\\
&=& |q|^2\Ddot_\b(|q|^{-2}O^{\a\b}\Ddot_\a \psi),
\eeaa
which proves \eqref{eq:OO-psi-Ddot-S_4}.
\end{proof}

\begin{remark}
Even though the above definitions are all equivalent in Kerr, in perturbations of Kerr it is important to define the operator $\OO$ using relation \eqref{eq:OO-Jk-kerr} in order to obtain acceptable error terms, see Definition \ref{definition:OO-perturbations}.
\end{remark}

From Proposition \ref{COMMUTATION-KK-SQUARE}, we deduce that the operator $\KK$ commutes with the D'Alembertian for scalars, as $\Pi =0$. 
In the case of horizontal 2-tensors, the commutator between $\KK$ and $\squared_2$ gives rise to lower order terms involving the Riemann curvature which vanish in Schwarzschild. 
Using the relation between $\KK$ and $\OO$ given by \eqref{eq: BL.ExpressionOO}, one can prove that the modified Laplacian $\OO$ inherits the commutation properties with the conformal D'Alembertian.

\begin{proposition}\label{LEMMA:MOD-LAPLACIAN-KERR}\label{lemma:lot-curvature} 
In Kerr, the following commutation formula holds true for $\psi \in \sk_2$,
\bea\label{eq-commutator-OO-q^2square}
\,[\OO, |q|^2 \squared_2 ]\psi = |q|^2 \left[\nab\left(\frac{8a(r^2+a^2)\cos\th}{|q|^2}\right)\c\nab\nab_\That\dual\psi +O(ar^{-2})\nab^{\leq 1}_{\Rhat}\dk^{\leq 1}\psi \right],
\eea
where $\dk=(\nab_3, r\nab_4, r\nab)$ denotes weighted derivatives as in \cite{KS} and \cite{KS:Kerr}.
\end{proposition}

\begin{proof} 
See Section \ref{sec:proof-mod-lap}. 
\end{proof}

\begin{remark} 
Observe that the terms on the right hand side of \eqref{eq-commutator-OO-q^2square} come from Riemann curvature terms as commutators for a 2-tensor $\psi$. In the proof of \eqref{eq-commutator-OO-q^2square}, we use the following relation, see Lemma \ref{lemma:basicpropertiesLiebTfasdiuhakdisug:chap9},
\beaa
\nab_\T\psi &=& \Lieb_\T\psi + \frac{4amr\cos\th}{|q|^4}\dual\psi,
\eeaa
where the second term is only present for tensors due to curvature. 
\end{remark}


\subsection{The symmetry operators}


We are now ready to define the set of second order symmetry operators which have remarkable commutation properties with the D'Alembertian operator $\squared_2$ for horizontal 2-tensors.

  \begin{definition}\lab{definition:symmetry-tensors} 
  We define the following second order differential operators,  acting on $\mathfrak{s}_2$ tensors, as
\bea
\SS_1= \nab_\T \nab_\T, \quad \SS_2= a \nab_\T \nab_\Z, \quad \SS_3= a^2 \nab_\Z \nab_\Z, \quad \SS_4= \OO.
\eea
\end{definition}

Observe that since $\T$ and $\Z$ are Killing vectorfields which commute, as obtained in Lemma \ref{lemma:commutation-nabT-nabZ}, and because of relation \eqref{eq:OO-psi-Ddot-S_4}, we can also write\footnote{Note that we use the calligraphic $\SS$ for the differential operators, and the normal $S^{\a\b}$ for the symmetric tensors.} as a unique formula
\bea\label{eq:SS-aund-Ddot-kerr}
\SS_\aund\psi&=&|q|^2 \Ddot_\a(|q|^{-2}S_\aund^{\a\b} \Ddot_\b \psi) \qquad \text{ for $\aund=1,2,3, 4$,}
\eea
where $S_{\aund}^{\a\b}$ are given in Definition \ref{definition:tensors-S}. 

 The commutation properties of the symmetry operators in Kerr can be deduced from the more general ones for approximate symmetry operators in perturbations, see Section \ref{symmetry-operators-pert}.


\section{Null geodesics in Kerr}\label{section:geodesics-Kerr}


In this section, we derive estimates for null geodesics in Kerr using the vectorfield $\partial_r$ as multiplier.  The computations appearing in these estimates present similar calculations to the Morawetz estimates for solutions to the wave equations, and for this reason we derive here their precise forms.


\subsection{The constants of motion for geodesics}


Let $\ga(\la)$ be a null geodesic in Kerr. Using the expression for the inverse of the metric given by \eqref{inverse-metric}, along $\ga(\la)$, since $\g(\gadot, \gadot)=0$ we have, with $\gadot_r=\pr_r^\a \gadot_\a$, $\gadot_t=\pr_t^\a \gadot_\a$, $\gadot_\vphi=\pr_\vphi^\a \gadot_\a$
\beaa
0= |q|^2 \g^{\a\b} \gadot_\a \gadot_\b=\left( \De \pr_r^\a \pr_r^\b +\frac{1}{\De}\RR^{\a\b}\right)\gadot_\a \gadot_\b=
\De\gadot_r \gadot_r+ \frac{1}{\De} \RR^{\a\b}\gadot_\a \gadot_\b
\eeaa
with
\beaa
\RR^{\a\b}\gadot_\a \gadot_\b&=& -(r^2+a^2) ^2\gadot_t \gadot_t- 2a(r^2+a^2) \gadot_t \gadot_\vphi- a^2  \gadot_\vphi \gadot_\vphi  +\De O^{\a\b} \gadot_\a \gadot_\b.
\eeaa

Since $\pr_t=T$ and $\pr_\vphi=Z$ are Killing vectorfields, we deduce that  $\gadot_t=\g(\gadot, T)$ and $\gadot_\vphi= \g(\gadot, Z)$ are constants of the motion, i.e. constants along  $\ga$, and respectively called the energy and the azimuthal angular momentum. We write,
\beaa
\e:=-\g(\gadot, T) , \qquad \lz :=-\g(\gadot, Z).
\eeaa
We also define\footnote{Observe that $\k^2$ is a positive constant of motion.}
\beaa
\k^2&:=&K^{\a\b} \gadot_\a\gadot_\b
\eeaa
for the Carter tensor $K$ in Kerr. Since $K$ is Killing, $\k^2$ is also a constant of motion. Indeed, we have
\beaa
\frac{d}{d\lambda} \k^2\big(\ga(\lambda) \big)&=&\D_\la K_{\a\b}\,\gadot^\la \gadot^\a\gadot^\b=0.
\eeaa

Since from \eqref{definition-K-kerr},  $K=   -(a^2\cos^2\th)  \g +O$ and 
$\ga$ is null we deduce, with $\gadot_a=\g(\gadot, e_a)$
\beaa
\k^2&=& O^{\a\b} \gadot_\a \gadot_\b=|q|^2  \big( e_1^\a e_1^\b+ e_2^\a e_2^\b\big) \gadot_\a \gadot_\b =|q|^2 \big(|\gadot_1|^2+|\gadot_2|^2 \big).
\eeaa

We summarize the result in the following.
\begin{proposition}
The quantities 
\beaa
\e=-\g(\gadot, T) , \qquad \lz=-\g(\gadot, Z), \qquad \k^2=K^{\a\b} \gadot_\a\gadot_\b,
\eeaa
are constant along  null geodesics. Moreover, relative to the null frame
\beaa
\k^2&=|q|^2 \big(|\gadot_1|^2+|\gadot_2|^2 \big).
\eeaa
\end{proposition}

With these constants  we have
\beaa
\RR^{\a\b}\gadot_\a \gadot_\b&=& -(r^2+a^2) ^2\gadot_t \gadot_t- 2a(r^2+a^2) \gadot_t \gadot_\vphi- a^2  \gadot_\vphi \gadot_\vphi  +\De O^{\a\b} \gadot_\a \gadot_\b\\
&=& -(r^2+a^2) ^2 \e^2- 2a(r^2+a^2) \e \, \lz - a^2  \lz^2  +\De \k^2
\eeaa
which is only a function of $r$ along  any fixed  null geodesic. 
We introduce the notation
\beaa
\RR(r; a, m, \e,\lz, \k):=  -(r^2+a^2)^2 \e^2- 2a(r^2+a^2) \e \, \lz - a^2  \lz^2  +\De \k^2.
\eeaa
Note that we  have the identity
\bea
\lab{identity:RR(r,a,m,e,lz,q}
-\RR(r; a, m, \e,\lz, \k)&=\big(( r^2+a^2)\e + a \lz\big)^2-\De\k^2.
\eea
In view of the above,  we  infer that
\beaa
0=
\De\gadot_r \gadot_r+ \frac{1}{\De} \RR^{\a\b}\gadot_\a \gadot_\b =\De\gadot_r \gadot_r+ \frac{1}{\De}\RR(r; a, m, \e,\lz, \k).
\eeaa
Since 
\beaa
\frac{dr}{d\la} &=& \gadot^\a\frac{\pr r}{\pr\alpha}=\gadot^r=\g^{rr}\gadot_r=\frac{\De}{|q|^2}\gadot_r,
\eeaa
we finally obtain
\beaa
|q|^4\Big(\frac{dr}{d\la} \Big)^2  =-\RR(r; a, m, \e,\lz, \k)
\eeaa
which is the equation for a null geodesic with constants of motion $\e$, $\lz$, $\k$.


\subsection{Trapped null geodesics}


There exist null geodesics along which  $\RR(r; a, m, \e,\lz, \q)=0$  i.e. $r$ remains constant. These are called 
orbital null geodesics.

\begin{remark}
Trapped null geodesics correspond to null geodesics that stay in a region $[r_1, r_2]$ of $r$ with $r_+<r_1<r_2<+\infty$ for all values of $\la$, and are thus a priori more general than orbital null geodesics. As it turns out, see for example Proposition 2 in \cite{CeJa}, all trapped null geodesics in Kerr are in fact orbital null geodesics. Thus, from now on, we do not distinguish between trapped and orbital null geodesics.
\end{remark}

If $r$ is constant we also have
\beaa
-\partial_r\RR(r; a, m, \e,\lz, \k)=\partial_r(|q|^4)\Big(\frac{dr}{d\la} \Big)^2+2 |q|^4 \frac{dr}{d\la} \partial_r \Big(\frac{dr}{d\la} \Big)=0.
\eeaa
The $r$ values for which such solutions are possible  must then verify the equations
\beaa
\RR(r; a, m, \e,\lz, \k)=\pr_r\RR(r; a, m, \e,\lz, \k)=0.
\eeaa
Thus, introducing
\beaa
\Pi:=( r^2+a^2)\e + a \lz
\eeaa
we write from \eqref{identity:RR(r,a,m,e,lz,q}
\beaa
-\RR(r; a, m, \e,\lz, \k)&=&\Pi^2-\De \k^2=0,\\
-\pr_r \RR(r; a, m, \e,\lz, \k)&=&2\Pi (\pr_r \Pi) -(\pr_r \De)\k^2=0.
\eeaa
From the second equation, we deduce
\bea
\lab{eq:orbittingnullgeod0}
\k^2&=&2\Pi \frac{\pr_r\Pi }{\pr_r \De}.
\eea
Thus, substituting in the first equation,
\beaa
\Pi^2-2\Pi \frac{\De\pr_r\Pi}{\pr_r\De}=0
\eeaa
or, if $\Pi \neq 0$,
\bea
\lab{eq:orbittingnullgeod1}
\Pi(\pr_r\De)-2 (\pr_r \Pi) \De=0.
\eea

We make use of the following calculation.
\begin{lemma}
\lab{lemma:calculationPiDe}
We have the identity
\bea\label{-2T-Pi-De}
\Pi(\pr_r\De)-2 (\pr_r \Pi) \De &=&  - 2 \TT_{\e, \lz}
\eea
where
\beaa
\TT_{\e, \lz}&:=& \big( r^3-3mr^2 + r a^2+ma^2\big)\e-  (r-m) a\lz.
\eeaa
\end{lemma}

\begin{proof}
We have
\beaa
(\pr_r \De) \Pi- 2\De (\pr_r\Pi) &=& 2(r-m) \big( (r^2+a^2) \e + a \lz\big)- 4 r  \big( r^2+ a^2- 2rm\big) \e\\
&=&2\Big(  (r-m) \big(( r^2+a^2)\e+  a \lz\big) - 2 r  \big( r^2+ a^2- 2rm\big)\e\Big)\\
&=& 2\Big( \big(-r^3+3mr^2 - r a^2-ma^2\big) \e+  (r-m) a\lz\Big)\\
&=& - 2\TT_{\e, \lz}
\eeaa
as stated.
\end{proof}

As a consequence of the Lemma we deduce that all orbital null geodesics are given by the equation
\beaa
\TT_{\e, \lz}= \big( r^3-3mr^2 + r a^2+ma^2\big)\e-  (r-m) a\lz=0.
\eeaa

\begin{remark}
The following hold true.
\begin{enumerate}
\item There are no trapped null geodesics perpendicular to $T=\pr_t$ in the exterior of a non-extremal Kerr.

\item  The values of $r$ for which trapped null geodesics exist  depends on  the ratio $\lz/ \e$. More precisely, at trapped null geodesics, we have
\beaa
\frac{r^3-3mr^2 + r a^2+ma^2}{r-m}=\frac{ a\lz}{\e}.
\eeaa
In particular, for $\lz=0$,   the trapped null geodesics  are given by  the equation
\bea\label{definition-TT}
\TT:= r^3- 3m r^2+ a^2r+a^2m=0.
\eea
\end{enumerate}
\end{remark}

\begin{remark}\lab{rmk:rangeofrfortrappednullgeodesics}
Note that one may specify possible values of $r$ for which trapped null geodesics exists. Indeed, let 
\bea
\bsplit
\hat{r}_1 &:=2m\left(1+\cos\left(\frac{2}{3}\arccos\left(-\frac{|a|}{m}\right)\right)\right),\\
\hat{r}_2 &:=2m\left(1+\cos\left(\frac{2}{3}\arccos\left(\frac{|a|}{m}\right)\right)\right).
\end{split}
\eea
Then, a trapped null geodesics satisfies $r\in [\hat{r}_1, \hat{r}_2]$, see for example \cite{Teo}.
\end{remark}

We now show that the trapped null geodesics are unstable, i.e. that $\pr_r^2 \RR(r; a, m, \e,\lz, \k) \leq 0$.
We have
\beaa
\Pi&=&( r^2+a^2) \e + a \lz, \\
\pr_r \Pi &=&2r \e,\\
\pr_r^2 \Pi&=& 2 \e,
\eeaa
and using \eqref{eq:orbittingnullgeod0} to write $\k^2=4r \frac{\Pi }{\pr_r \De} \e$, we have
\beaa
-\pr^2_r \RR(r; a, m, \e,\lz, \k)&=&2(\pr_r \Pi)^2+2\Pi (\pr^2_r \Pi) -2\k^2\\
&=&8r^2 \e^2+4\Pi \e-8r \frac{\Pi }{\pr_r \De} \e\\
&=&\frac{4}{\partial_r \De} \left( 2r^2 \pr_r\De \e^2+(\pr_r\De) \Pi \e-2r \Pi  \e\right)\\
&=&\frac{4\e}{\partial_r \De} \left( 4r^2 (r-m) \e-2m \Pi \right).
\eeaa
Using \eqref{eq:orbittingnullgeod1} to write $\Pi=\frac{4r\De   \e }{\pr_r\De}$
 we deduce
 \beaa
-\pr^2_r \RR(r; a, m, \e,\lz, \k)&=&\frac{4\e}{\partial_r \De} \left( 4r^2 (r-m) \e-2m\frac{4r\De   \e }{\pr_r\De} \right)\\
&=&\frac{32r\e^2}{(\partial_r \De)^2} \left( r (r-m)^2 -m \De  \right)\\
&=& \frac{8r}{(r-m)^2}\e^2\Big( r(r-m)^2 -m (r^2+a^2- 2mr)\Big)\\
&=& \frac{8r}{(r-m)^2}\e^2\Big((r-m)^3 +m( m^2-a^2)\Big),
\eeaa
which is positive since $r\ge m$ and $|a|\le m$.


\subsection{Morawetz estimates for geodesics}\label{section:morawetz-geodesics}


 We now derive estimates for null geodesics in Kerr using the vectorfield $X=\FF(r)\partial_r$, for some function $\FF(r)$ as multiplier. The energy\footnote{Observe that this energy is not positive definite since $X$ is not causal.} relative to $X$ is given by 
\beaa
e_{X} [\ga] (\tau):=-\g^{\a\b}{X}_\a \gadot_\b.
\eeaa

Our goal is to find a function $\FF(r)$ such that the energy $e_{X} [\ga]$ is non-increasing for all $\tau$.
We have,  recalling that $\D_{\gadot}\gadot=0$ and $\g^{\a\b}\gadot_\a\gadot_\b=0$,
\beaa
\frac{d}{d\tau} e_{X} [\ga] (\tau) &=&-\D_\a X_\b \gadot^\a\gadot^\b= -\frac 1 2( \Lie_X\g_{\a\b} )\gadot^\a\gadot^\b=\frac 1 2\Lie_X(\g^{\a\b}) \gadot_\a\gadot_\b=\frac 1 2|q|^{-2}\Lie_X(|q|^2 \g^{\a\b})  \gadot_\a\gadot_\b.
\eeaa

We perform the following computations as done in \cite{A-B}. 
\begin{lemma}
\lab{lemma:Lie_Xgdot}
We have, for $X=\FF(r) \pr_r$ and $\RR^{\a\b}$ defined in \eqref{eq:expressionRR-O}, 
\beaa
\Lie_X(|q|^2 \g^{\a\b})&=&\big( \FF\pr_r \De -2 \De \pr_r \FF\big)\pr_r^\a\pr_r^\b+ \FF\pr_r\left(\frac 1 \De\RR^{\a\b} \right).
\eeaa
\end{lemma}

\begin{proof}
Recall from \eqref{inverse-metric} that we have $|q|^2 \g^{\a\b}=\De \pr_r^\a \pr_r^\b +\frac{1}{\De} \RR^{\a\b}$.
Hence,
\beaa
 \Lie_X(|q|^2 \g^{\a\b})&=&\Lie_X\big(\De \pr_r^\a \pr_r^\b \big)+\Lie_X\left(\frac 1 \De\RR^{\a\b} \right)\\
 &=& X(\De) \pr_r^\a\pr_r^\b+ \De [X, \pr_r]^\a \pr_r^\b +\De \pr_r^\a [X, \pr_r]^\b  +\Lie_X\left(\frac 1 \De\RR^{\a\b} \right)\\
 &=&  \FF\pr_r \De  \pr_r^\a\pr_r^\b- 2 \De \pr_r \FF \pr_r^\a\pr_r^\b+\Lie_{\FF\pr_r}\left(\frac 1 \De\RR^{\a\b} \right).
\eeaa
For  a tensor $U^{\a\b}$ we have,
\beaa
\Lie_{X}   U^{\a\b}=X^\la\pr_\la U^{\a\b}-U^{\la \b}\pr_\la X^\a-U^{\a\la} \pr_\la X^\b.
\eeaa
For $X= \FF(r) \pr_r $   and $U^{r\b}=0=U^{\a r} $, we deduce
\beaa
\Lie_{X}   U^{\a\b}= \FF \pr_r  U^{\a\b}.
\eeaa
Hence
\beaa
 \Lie_X(|q|^2 \g^{\a\b})&=&  \FF\pr_r \De  \pr_r^\a\pr_r^\b- 2 \De \pr_r \FF \pr_r^\a\pr_r^\b+ \FF\pr_r\left(\frac 1 \De\RR^{\a\b} \right)\\
 &=&\big( \FF\pr_r \De -2 \De \pr_r \FF\big)\pr_r^\a\pr_r^\b + \FF\pr_r\left(\frac 1 \De\RR^{\a\b} \right)
 \eeaa
 as stated.
 \end{proof}
 
Using the lemma we deduce,
\beaa
|q|^2\frac{d}{d\tau} e_{\FF\pr_r} [\ga] (\tau) &=&\frac 1 2  \Big( \FF\pr_r \De- 2 \De \pr_r \FF \Big)\gadot_r\gadot_r+\frac 1 2 \FF\pr_r \left(\frac 1 \De\RR^{\a\b} \right)  \gadot_\a\gadot_\b.
\eeaa
We now recall that,
\beaa
\bsplit
\RR&=\RR^{\a\b}\gadot_\a\gadot_\b,\\
\RR^{\a\b}&= -(r^2+a^2) ^2\pr_t^\a\pr_t^\b- 2a(r^2+a^2) \pr_t^{(\a} \pr_\vphi^{\b)}- a^2 \pr_\vphi^\a \pr_\vphi^\b  +\De O^{\a\b}
\end{split}
\eeaa
Thus, we have
\beaa
\pr_r \RR^{\a\b}&=- 4 r(r^2+a^2)  \pr_t^\a\pr_t^\b- 4 ar \pr_t^{(\a}\pr_\vphi^{\b)}+ 2(r-m) O^{\a\b}. 
\eeaa
This yields
\beaa
 \FF\pr_r \left(\frac 1 \De\RR^{\a\b} \right)  \gadot_\a\gadot_\b&=&-\FF \frac{\pr_r\De}{\De^2} \RR + \FF\frac{1}{\De}(\pr_r \RR ^{\a\b} )\gadot_\a\gadot_\b\\
 &=& -\FF \frac{\pr_r\De}{\De^2} \RR  + \FF\frac{1}{\De}\left(- 4 r (r^2+a^2) \e^2- 4 ar  \e \, \lz+ 2(r-m)\k^2\right)\\
 &=& \FF\pr_r \left(\frac 1 \De\RR\right)
\eeaa
 and therefore 
 \bea
 \lab{eq:maingeodidentity1}
 |q|^2\frac{d}{d\tau} e_{\FF\pr_r} [\ga] (\tau) &=&\frac 1 2  \Big( \FF\pr_r \De- 2 \De \pr_r \FF \Big)\gadot_r^2+\frac 1 2  \FF\pr_r \left(\frac 1 \De\RR\right).
 \eea

\begin{remark}
In the particular case of trapped null geodesics, for which $\RR=\partial_r \RR=0$, the non-increasing of the above energy can be obtained by choosing $\FF=-\partial_r\RR$. Indeed, for that choice, evaluating on a trapped null geodesic $\ga$:
\beaa
|q|^2\frac{d}{d\tau} e_{\FF\pr_r} [\ga] (\tau) &=& \frac 1 2  \Big(- \partial_r\RR\pr_r \De+ 2 \De \pr^2_r \RR \Big)\gadot_r^2-\frac 1 2( \partial_r\RR)\pr_r \left(\frac 1 \De\RR\right)\\
&=&   \De( \pr^2_r \RR) \gadot_r^2 \leq 0.
\eeaa
In order to obtain a similar result for all null geodesics, we need to add more freedom in the choice of additional scalar functions.
\end{remark}

 We obtain   the following proposition.
\begin{proposition}
\lab{proposition:geodidentity1}
The following identity holds true for any scalar functions $\FF$,  $ w=w_{red}$,
\bea\lab{d-dtau-e-FF}
\bsplit
|q|^2\frac{d}{d\tau} e_{\FF\pr_r} [\ga] (\tau) &=  \frac 1 2  \Big( \FF\pr_r \De- 2 \De \pr_r \FF +\De w_{red}\Big)\gadot_r^2 +\frac 1 2  \FF\pr_r \left(\frac 1 \De\RR\right)+\frac 1 2   w_{red}\frac 1 \De \RR,
 \end{split}
\eea
We  write the above in the form
\bea
\lab{d-dtau-e-FF2}
|q|^2\frac{d}{d\tau} e_{\FF\pr_r} [\ga] (\tau)  &=&-\AA  \gadot_r^2 - \UU
\eea
where
\bea
\AA&:=& - \frac 1 2\FF\pr_r \De+  \De \pr_r \FF - \frac 1 2 \De w_{red},\label{AA-geodesics} \\
\UU&:=& -\frac 1 2  \FF\pr_r \left(\frac 1 \De\RR\right)-\frac 1 2   w_{red}\frac 1 \De \RR. \label{UU-geodesics}
\eea
\end{proposition}

\begin{proof}
  Using again  \eqref{inverse-metric}, 
and the fact that  $\ga$ is  null we can rewrite   \eqref{eq:maingeodidentity1}, for some scalar $w_{red}$
\beaa
|q|^2\frac{d}{d\tau} e_{\FF\pr_r} [\ga] (\tau) &=& \frac 1 2  \Big( \FF\pr_r \De- 2 \De \pr_r \FF \Big)\gadot_r^2+\frac 1 2  \FF\pr_r \left(\frac 1 \De\RR\right)+\frac 1 2 |q|^2  w_{red}\g^{\a\b} \, \gadot_\a\gadot_\b\\
&=& 
 \frac 1 2  \Big( \FF\pr_r \De- 2 \De \pr_r \FF \Big)\gadot_r^2+\frac 1 2  \FF\pr_r \left(\frac 1 \De\RR\right)+\frac 1 2   w_{red}\left(\De\,\gadot_r^2+\frac{1}{\De} \RR\right)\\
 &=&  \frac 1 2  \Big( \FF\pr_r \De- 2 \De \pr_r \FF + \De w_{red}\Big)\gadot_r^2 +\frac 1 2  \FF\pr_r \left(\frac 1 \De\RR\right)+\frac 1 2   w_{red}\frac 1 \De \RR
\eeaa
as stated. 
\end{proof}

In order to obtain that the energy $e_{X} [\ga]$ is non-increasing for all $\tau$, our goal  is  to choose $\FF$ and $w_{red}$ so that $\AA \geq 0$ and $\UU \geq 0$. We collect here the conditions we need.

\medskip
   {\bf  Calculation of  $\UU$.} 
Given a scalar function $z$ we write
 \beaa
 - \frac{ 1}{2}\FF   \partial_r\left(\frac 1 \De\RR\right)&=& - \frac{ 1}{2}\FF z^{-1}  \partial_r\left( \frac z \De\RR\right)+ \frac{ 1}{2}\FF z^{-1}\partial_r z   \frac{ \RR}{\Delta}.
 \eeaa
  We then have from \eqref{UU-geodesics}
 \beaa
\UU &=&  -\frac 1 2  \FF\pr_r \left(\frac 1 \De\RR\right)-\frac 1 2   w_{red}\frac 1 \De \RR=  - \frac{ 1}{2}\FF z^{-1}  \partial_r\left( \frac z \De\RR\right)+ \frac{ 1}{2}\left(\FF z^{-1}\partial_r z  -w_{red}\right) \frac{ \RR}{\Delta}.
 \eeaa
By choosing
 \bea
 \lab{eq:w_{red}-geodesic}
w_{red}&=&  \FF z^{-1}\partial_r z, 
 \eea
 the coefficient of $\frac{ \RR}{\Delta}$ cancels out.  Thus 
 \beaa
\UU&=&  - \frac{ 1}{2}\FF z^{-1}  \partial_r\left( \frac z \De\RR\right)=  - \frac{ 1}{2}\FF z^{-1}  \RRt',\qquad \RRt':=\partial_r\left( \frac z \De\RR\right).
 \eeaa
Finally, given a scalar function $h$, by choosing\footnote{The relation for $w_{red}$ is implied by the previous dependence on $\FF$.}
\bea\label{def-w-red-in-fun-FF-00-geodesics}\lab{def-w-red-in-fun-FF-00}
\FF&=& - z h  \tilde{\RR}', \qquad  w_{red}=  - (\partial_r z ) h  \tilde{\RR}', 
\eea
we deduce
\beaa
\UU&=&  - \frac{ 1}{2}\FF z^{-1}  \tilde{\RR}'=   \frac{ 1}{2}  h ( \tilde{\RR}' )^2.
\eeaa
In particular,  $\UU$ is positive as long as $h$ is positive.

\medskip

  {\bf  Calculation of  $\AA$.}   
With the choices of $\FF$ and $w_{red}$ given by \eqref{def-w-red-in-fun-FF-00-geodesics}, we compute from \eqref{AA-geodesics}
\beaa
\AA&=&  - \frac 1 2\FF\pr_r \De+  \De \pr_r \FF - \frac 1 2 \De w_{red}=  \partial_r\left(\frac{\FF}{\Delta^{1/2}}  \right) \Delta^{3/2}-\frac 12\Delta   w_{red} \\
&=&\partial_r\left(\frac{- z h \tilde{\RR}' }{\Delta^{1/2}}  \right) \Delta^{3/2}-\frac 12\Delta  (- ( \partial_r z) h \tilde{\RR}' )=-\partial_r\left(\frac{ z^{1/2} z^{1/2}  h \tilde{\RR}' }{\Delta^{1/2}}  \right) \Delta^{3/2}+\frac 12\Delta   ( \partial_r z) h \tilde{\RR}' \\
&=&- \frac 1 2 \partial_rz  \left(\frac{  h \tilde{\RR}' }{\Delta^{1/2}}  \right) \Delta^{3/2}-z^{1/2} \partial_r\left(\frac{ z^{1/2}  h \tilde{\RR}' }{\Delta^{1/2}}  \right) \Delta^{3/2}+\frac 12\Delta   ( \partial_r z) h \tilde{\RR}' \\
&=&-z^{1/2}\Delta^{3/2} \partial_r\left(h \frac{ z^{1/2}  \tilde{\RR}' }{\Delta^{1/2}}  \right).
\eeaa
We obtain
\beaa
\AA&=&-z^{1/2}\Delta^{3/2} \RRt'', \qquad \RRt'':=\pr_r\left(  h \frac{z^{1/2}}{\De^{1/2} } \RRt'\right).
\eeaa

We summarize the result in the following
\begin{lemma}
With the choice of the functions $\FF=- zh \RRt'$, $w_{red}=  - (\partial_r z ) h  \tilde{\RR}' $, for functions $z$ and $h$, the identity \eqref{d-dtau-e-FF2} takes the form
\beaa
|q|^2\frac{d}{d\tau} e_{\FF\pr_r} [\ga] (\tau)  &=&\Big(z^{1/2}\De^{3/2}\RRt'' \Big)\gadot_r^2 -\frac 1 2  h (\RRt')^2
\eeaa
where $\RRt'=\partial_r\left( \frac z \De\RR\right)$ and $\RRt''=\pr_r\Big(  h \frac{z^{1/2}}{\De^{1/2} } \RRt'\Big)$. 
\end{lemma}

To obtain a non-increasing energy $e_{\FF\pr_r} [\ga]$, we are therefore left to choose positive scalar functions $z$ and $h$ for which $\RRt'' \leq 0$.


\subsubsection{Choice of $z$ and $h$}
\lab{subsubsect:chiocezg-geodesic}


From \eqref{identity:RR(r,a,m,e,lz,q}, we write
\beaa
\RR&=-( r^2+a^2)^2\e^2 - 2a (r^2+a^2) \e \,\lz - a^2 \lz^2+\De\k^2,
\eeaa
and therefore we have
\beaa
\RRt'&=& \pr_r\left(\frac{z}{\De} \RR \right)\\
&=& -\pr_r\left(\frac{z}{\De} ( r^2+a^2)^2 \right) \e^2 - 2a \pr_r\left(\frac{z}{\De} (r^2+a^2)\right) \e \,\lz - a^2\pr_r\left(\frac{z}{\De} \right) \lz^2+(\pr_rz)  \k^2 .
\eeaa

We choose $z$ to cancel the coefficient of $\e^2$ in $\RRt'$, i.e.
\bea\label{choice-z-geodesics}
z&=& \frac{\De}{(r^2+a^2)^2}.
\eea
This choice implies
\beaa
\RRt'&=&  - 2a \pr_r\left( \frac{1}{(r^2+a^2)} \right) \e \c\lz - a^2\pr_r\left(\frac{1}{(r^2+a^2)^2} \right) \lz^2+\pr_r \left(\frac{\De}{(r^2+a^2)^2}\right)  \k^2 \\
&=&     \frac{4ar}{(r^2+a^2)^2}  \e \c\lz +\frac{4a^2r}{(r^2+a^2)^3}  \lz^2-2\left(\frac{r^3-3mr^2+a^2r+ma^2}{(r^2+a^2)^3}\right) \k^2 .
\eeaa

\begin{remark} 
Observe that 
\bea\lab{eq:pr_rz}
\pr_r z&=& \pr_r \left(\frac{\De}{(r^2+a^2)^2}\right) =-2\left(\frac{r^3-3mr^2+a^2r+ma^2}{(r^2+a^2)^3}\right)= -\frac{2\TT}{(r^2+a^2)^3}
\eea
where $\TT$ is defined as in \eqref{definition-TT} to be the locus  of trapped geodesics with $\lz=0$, i.e. the trapped set in the axially symmetric case.
\end{remark}

We now compute $\RRt''$:
\beaa
\RRt''&=&\pr_r\left(  h \frac{z^{1/2}}{\De^{1/2} } \RRt'\right)=\pr_r\left(   \frac{h}{r^2+a^2 } \RRt'\right)\\
&=&\pr_r\left(    h  \frac{4ar}{(r^2+a^2)^3} \right) \e \c\lz + \pr_r \left( h\frac{4a^2r}{(r^2+a^2)^4} \right) \lz^2-2 \pr_r \left( h \frac{r^3-3mr^2+a^2r+ma^2}{(r^2+a^2)^4}\right) \k^2. 
\eeaa
We choose\footnote{In the Morawetz estimates for the wave equation we will choose a different $h$, i.e. $h=\frac{(r^2+a^2)^4}{r(r^2-a^2)}$. See Remark \ref{remark:choice-h-subextremal}.} $h$ to cancel the coefficient of $\e \c \lz$ in $\RRt''$, i.e. 
\bea\label{choice-h-geodesics}
h&=& \frac{(r^2+a^2)^3}{r}.
\eea
This choice implies
\beaa
\RRt''&=& \pr_r \left( \frac{4a^2}{(r^2+a^2)} \right) \lz^2-2 \pr_r \left(  \frac{r^3-3mr^2+a^2r+ma^2}{r(r^2+a^2)}\right) \k^2 \\
&=&- \frac{8a^2r}{(r^2+a^2)^2}  \lz^2-  \frac{2m\big( 3 r^4 -6a^2 r^2 -a^4\big)}{r^2( r^2+a^2)^2} \k^2 .
\eeaa
Observe that $\RRt'' \leq 0$ as long as $3 r^4 -6a^2 r^2 -a^4  \ge 0$. 
Note that  $3 r^4 -6a^2 r^2 -a^4  \ge 0$ if
\beaa
r^2\ge \frac{a^2(3+2\sqrt{3} )}{3}, \qquad \mbox{for} \quad r\ge r_+.
\eeaa
It thus suffices to check the sign  on the horizon $r=r_+= m+\sqrt{m^2-a^2} $. We  need therefore,
$\big( m+\sqrt{m^2-a^2} \big) ^2\ge  \frac{a^2(3+2\sqrt{3}) }{3} $
or, taking $\la=\frac{a}{m}$,
\beaa
\big(1+\sqrt{1-\la^2} \big) ^2\ge  \la^2 \frac{(3+2\sqrt{3}) }{3} 
\eeaa
This implies $\sqrt{1-\la^2}  \ge |\la|\sqrt{\frac{(3+2\sqrt{3}) }{3}} -1$, which is verified for
\beaa
\frac{|a|}{m}=|\la|\leq \frac{3}{3+\sqrt{3}}\sqrt{\frac{(3+2\sqrt{3}) }{3}}\simeq 0.9306.
\eeaa
In particular, for Kerr spacetime with $|\frac{a}{m}| \leq  0.93$, we have
\beaa
-\RRt'' \geq 0 \qquad \text{for $r \geq r_{+}.$}
\eeaa

\begin{remark}\label{remark:choice-h-subextremal} 
Observe that if we restrict our analysis to axially symmetric geodesics, for which $\lz=0$, then the term $\RRt''$ reduces to
\beaa
\RRt''&=&-2 \pr_r \left( h \frac{r^3-3mr^2+a^2r+ma^2}{(r^2+a^2)^4}\right) \k^2 
\eeaa
for any positive function $h$. In particular, by choosing 
\bea\label{choice-h-subextremal}
h=\frac{(r^2+a^2)^4}{r(r^2-a^2)}
\eea
we obtain
\beaa
\RRt''&=&-2 \pr_r \left( \frac{r^3-3mr^2+a^2r+ma^2}{r(r^2-a^2)}\right) \k^2 =- 2 \frac{3mr^4-4a^2r^3+ma^4}{r^2(r^2-a^2)^2} \k^2 
\eeaa
which is negative in the exterior region in the full sub-extremal range $|a|<m$. This can be seen by setting $a^2=\gamma m^2$ and $r=(m+\sqrt{m^2-a^2})x=(1+\sqrt{1-\gamma })mx$, where the exterior region is given by $x>1$. We then obtain
\beaa
3mr^4-4a^2r^3+ma^4&=&m^5 \big(3(1+\sqrt{1-\gamma })^4x^4-4\gamma (1+\sqrt{1-\gamma })^3x^3+\gamma^2\big)
\eeaa
which is positive for $x\geq 1$ and $0 \leq \gamma \leq 1$.
\end{remark}

The above choices allow to prove the Morawetz estimates for geodesics in Kerr, as summarized in the following proposition.
\begin{proposition} 
\lab{Prop:Morawetz-nullgeodesics}  
Let $\gamma$ be a null geodesics in a Kerr spacetime with $|\frac{a}{m}| \leq 0.93$.  Then, with the choices
\beaa
\FF=- z h  \tilde{\RR}', \qquad z= \frac{\De}{(r^2+a^2)^2}, \qquad h= \frac{(r^2+a^2)^3}{r},
\eeaa
the energy $e_{\FF \partial_r} [\ga]$ is non-increasing for all $\tau$, i.e.
\beaa
|q|^2\frac{d}{d\tau} e_{\FF\pr_r} [\ga] (\tau)  &=& \Big(z^{1/2}\De^{3/2}\RRt'' \Big)\gadot_r^2 -\frac 1 2  h (\RRt')^2 \\
&=&- \frac{\De^{2}}{(r^2+a^2)}(-\RRt''  ) \gadot_r^2 - \frac{(r^2+a^2)^3}{2r} (\RRt')^2  \leq 0.
\eeaa
\end{proposition}


\chapter{Perturbations of Kerr}\label{PERTURBATIONS-SECTION}



\section{Set-up and linearized quantities}
\lab{sec:setupandlinearizedquantities}


We consider a given vacuum spacetime $(\MM, \g)$ together with a null pair $(e_3, e_4)$ and its corresponding horizontal structure as in section \ref{sec:nullparisandhorizontalstruct}. We will use the complexified Ricci and curvature coefficients  of Definition \ref{def:complexRicciandcurvaturecoefficients}. To be able to  talk about  small perturbations of Kerr we  also  need:
\begin{itemize}
\item $\MM$ is endowed with a pair of constants $(a, m)$.

\item $\MM$ is endowed with a pair of scalar functions $(r, \th)$. 

\item The complex valued scalar function $q$ is defined as $q := r+ia \cos\th. $

\item $\MM$ is endowed with a complex horizontal 1-form $\Jk$. 
\item  Define linearized  quantities, such as  $\Gac=\Ga-\Ga_{Kerr},\,  \Rc=R-R_{Kerr}$, i.e.  subtract   from the Ricci and curvature  coefficients  the  corresponding values  in Kerr, expressed relative to $(a,m, r, \cos \th, \Jk)$.
\end{itemize}


\subsection{Definition of linearized quantities}
\lab{sec:definitionoflinearizedquantities:chap4}


Recall from section \ref{section:values-Kerr} that the following Ricci and curvature coefficients  vanish in Kerr
\beaa
\Xh, \quad \Xbh,\quad \Xi, \quad \Xib,\quad A, \quad B, \quad \Bb, \quad \Ab,
\eeaa
as well as the following additional quantities  
\beaa
\nab(r), \qquad e_4(\th), \qquad e_3(\th), \qquad \DD\hot\Jk.
\eeaa

We renormalize below all other quantities, not vanishing in Kerr. We start with quantities which are 0-conformally invariant in the sense of Definition \ref{def:sconformalinvariants}.

\begin{definition}[Renormalization for 0-conformally invariant quantities] 
\lab{def:renormalizationofallnonsmallquantitiesinPGstructurebyKerrvalue:1}
We  define  the following renormalizations.
\begin{enumerate}
\item Linearization of 0-conformally invariant  Ricci and curvature coefficients:
\bea
 \Hc := H-\frac{aq}{|q|^2}\Jk, \qquad\quad \Hbc:=\Hb+\frac{a\ov{q}}{|q|^2}\Jk,\qquad  \Pc := P+\frac{2m}{q^3}.
\eea

\item Linearization of 0-conformally invariant  derivatives of $r$, $\cos\th$ and $q$:
\bea
\bsplit
\widecheck{\DD q} :=\DD q+a\Jk, \qquad\quad \widecheck{\DD \ov{q}} :=\DD \ov{q}-a\Jk,\qquad \widecheck{\DD(\cos\th)} := \DD(\cos\th) -i\Jk.
\end{split}
\eea

\item Linearization of 0-conformally invariant derivatives of $\Jk$:
\bea
\bsplit
\widecheck{\ov{\DD}\c\Jk}& := \ov{\DD}\c\Jk-\frac{4i(r^2+a^2)\cos\th}{|q|^4}. 
  \end{split}
\eea
 \end{enumerate}
\end{definition}

\begin{definition}[Outgoing renormalization for the remaining quantities] 
\lab{def:renormalizationofallnonsmallquantitiesinPGstructurebyKerrvalue:2} 
We  define  the following renormalizations.\footnote{Note that in Kerr, we have $\om=0$ in the outgoing principal frame so that this quantity does not need to be renormalized. By convention, we thus define $\omc=\om$ if the normalization of the null pair $(e_3, e_4)$ is outgoing.  In addition, note that in Kerr we have $e_3(q)=e_3(r)$ and $e_4(q)=e_4(r)$ so that it suffices to linearize $e_3(r)$ and $e_4(r)$.}
\begin{enumerate}
\item Linearization of the remaining Ricci and curvature coefficients:
\bea
\bsplit
\trXc &:= \tr X-\frac{2}{q}, \qquad\,\qquad \,\,\,\,    \trXbc := \tr\Xb+\frac{2q\Delta}{|q|^4},\\
\Zc &:= Z-\frac{a\ov{q}}{|q|^2}\Jk, \qquad\qquad\quad\,\,  \ombc := \omb  - \frac 1 2 \pr_r\left(\frac{\De}{|q|^2} \right).\\
\end{split}
\eea

\item Linearization of the remaining derivatives of $r$:
\bea
\widecheck{e_3(r)} := e_3(r)+\frac{\Delta}{|q|^2}, \qquad \widecheck{e_4(r)} := e_4(r)-1.
\eea

\item Linearization of the remaining derivatives of $\Jk$:
\bea
\widecheck{\nab_3\Jk}:=\nab_3\Jk -\frac{\De q}{|q|^4}\Jk, \qquad \widecheck{\nab_4\Jk}:=\nab_4\Jk +\frac{1}{q}\Jk. 
\eea
 \end{enumerate}
\end{definition}

\begin{definition}[Ingoing renormalization for the remaining quantities] 
\lab{def:renormalizationofallnonsmallquantitiesinPGstructurebyKerrvalue:3} 
We  define  the following renormalizations.\footnote{Note that in Kerr, we have $\omb=0$ in the ingoing principal frame so that this quantity does not need to be renormalized. By convention, we thus define $\ombc=\omb$ if the normalization of the null pair $(e_3, e_4)$ is ingoing. In addition, note that in Kerr we have $e_3(q)=e_3(r)$ and $e_4(q)=e_4(r)$ so that it suffices to linearize $e_3(r)$ and $e_4(r)$.}
\begin{enumerate}
\item Linearization of the remaining Ricci and curvature coefficients.
\bea
\bsplit
\trXc &:= \tr X-\frac{2\ov{q}\De}{|q|^4}, \qquad\,\qquad     \trXbc := \tr\Xb+\frac{2}{\ov{q}},\\
\Zc &:= Z-\frac{aq}{|q|^2}\Jk,\qquad \qquad \quad\,\,\,\omc  := \om  + \frac{1}{2}\pr_r\left(\frac{\De}{|q|^2} \right).\\
\end{split}
\eea

\item Linearization of remaining derivatives of $r$.
\bea
\bsplit
\widecheck{e_3(r)} := e_3(r)+1, \qquad \widecheck{e_4(r)} := e_4(r)-\frac{\Delta}{|q|^2}.
\end{split}
\eea

\item Linearization of remaining derivatives of $\Jk$.
\bea
\widecheck{\nab_3\Jk}:=\nab_3\Jk -\frac{1}{\ov{q}}\Jk, \qquad \widecheck{\nab_4\Jk}:=\nab_4\Jk +\frac{\De \ov{q}}{|q|^4}\Jk.
\eea
\end{enumerate}
\end{definition}

\begin{remark}
In Part II and Part III, we will always consider normalizations which are either outgoing, i.e. for which $\om$ is small in a suitable sense, or ingoing, i.e.  for which $\omb$ is small in a suitable sense. Also, in the region $r\leq r_++\deh$, we will only consider the ingoing renormalization.
\end{remark}


\subsection{Definition of the notations $\Ga_b$ and $\Ga_g$ for error terms}
\lab{sec:definitionofGabandGagfirsttime}


\begin{definition}
\lab{definition.Ga_gGa_b}
The set of all linearized quantities is of the form $\Ga_g\cup \Ga_b$ with  $\Ga_g,  \Ga_b$
 defined as follows.
 \begin{enumerate}
\item 
 The set   $\Ga_g=\Ga_{g,1}\cup \Ga_{g, 2}$   with
 \bea
 \bsplit
 \Ga_{g,1} &= \Big\{\Xi, \quad \omc, \quad\trXc,\quad  \Xh,\quad \Zc,\quad \Hbc, \quad \trXbc , \quad r\Pc, \quad  rB, \quad  rA\Big\},\\
 \Ga_{g,2} &= \Big\{\widecheck{e_4(r)}, \quad r^{-1}\nab(r), \quad e_4(\cos\th), \quad r\widecheck{\nab_4\Jk}\Big\}. 
 \end{split}
 \eea
 
 \item  The set  $\Ga_b=\Ga_{b,1}\cup \Ga_{b, 2}\cup \Ga_{b,3}$   with
 \bea
 \bsplit
 \Ga_{b,1}&= \Big\{\Hc, \quad \Xbh, \quad \ombc, \quad \Xib,\quad  r\Bb, \quad \Ab\Big\},\\
  \Ga_{b, 2}&= \Big\{r^{-1}\widecheck{e_3(r)}, \quad  \widecheck{\DD(\cos\th)}, \quad e_3(\cos\th)\Big\}, \\
   \Ga_{b,3}&=\Bigg\{ r\,\widecheck{\ov{\DD}\c\Jk}, \quad r\,\DD\hot\Jk, \quad r\,\widecheck{\nab_3\Jk} \Bigg\}. 
   \end{split}
 \eea
\end{enumerate}
\end{definition}

\begin{remark}
The justification for the above decompositions has to do with the expected  decay properties of the linearized  components in perturbations of Kerr. More precisely, we will consider perturbations of Kerr for which $\Ga_g$ and $\Ga_b$ satisfy the following estimates,    see section \ref{sec:controlofGabandGagfromBA} for details, 
\bea\lab{eq:expectedbehaviorGabGag:chap2}
\bsplit
\big|\dk^{\leq s}\Ga_g|&\les \ep \min\Big\{ r^{-2 }\tau^{-1/2-\dec},  \, r^{-1}\tau^{-1-\dec} \Big\}, \\
\big|\nab_3\dk^{\leq s-1}\Ga_g| &\les \ep  r^{-2 }\tau^{-1-\dec},\\
\big|\dk^{\leq s}\Ga_b\big| &\les \ep  r^{-1 }\tau^{-1-\dec},
\end{split}
\eea
for a small constant $\dec>0$, where $\dk=\{\nab_3, r\nab_4, \dkb=r\nab \}$ denotes weighted derivatives, and $\tau$ is a scalar function on $\MM$ whose properties are given in section \ref{sec:basicpropertiestau}.
\end{remark}

In addition to the above, as a consequence of the definition of the linearized quantities and of  $(\Ga_b, \Ga_g)$, as well as of the relations \eqref{eq:DDP-Kerr} and \eqref{eq:vanishing-relations-Kerr}, we also have the following:
\bea\label{eq:vanishing-relations-perturbations}
\begin{split}
 \DD P+3P\Hb &\in r^{-2} \dk^{\leq 1} \Ga_g, \\
  \DD\ov{P}+3\ov{P}H &\in r^{-2} \dk^{\leq 1} \Ga_g, \\
\DD\hot\Hb  +\Hb\hot\Hb &\in r^{-1} \dk^{\leq1} \Ga_g, \\
 \nab_4\Hb+\tr X\Hb &\in r^{-1} \dk^{\leq1} \Ga_g, \\
   \DDc\tr X +2\tr X\Hb &\in r^{-1} \dk^{\leq1} \Ga_g, \\
 \DD\hot H  +H \hot H &\in r^{-1} \dk^{\leq 1} \Ga_b, \\
 \ov{\tr X} \Hb + \tr X H &\in r^{-1} \Ga_b.
 \end{split}
\eea
In view of \eqref{eq:atrch-e3-atrchb-e4-Kerr} we also have 
\bea\label{eq:atrch-e3-atrchb-e4-pert-kerr}
\begin{split}
\atrch\nab_3+\atrchb \nab_4&= \frac{4a\cos\th (r^2+a^2)}{|q|^4}\That+ \Ga_g \c \dk, \\
\atrch e_3+\atrchb e_4+ 2(\eta+\etab) \c \dual \nab&=\frac{4a\cos\th}{|q|^2} \T+ \Ga_g \c \dk.
\end{split}
\eea


\section{Commutation formulas revisited}


We revisit the commutation formulas of section \ref{sec:commutationformulasfirsttime}. In some cases, we write the schematic structure of the error terms by making use of the definition of $\Ga_b$ and $\Ga_g$ introduced in section \ref{sec:definitionofGabandGagfirsttime}, by keeping track of different level of precision for the error terms, as they will be useful in different contexts.

   \begin{lemma}\label{LEMMA:COMMUTATION-FORMULAS-1}
   The following commutation formulas hold true.
   \begin{enumerate}
   \item 
   Let $h \in \sk_0(\CCC)$ $s$-conformally invariant. Then 
   \bea\label{eq:comm-nab4-nab3-DD-h-precise}
   \begin{split}
 \, [\nab_4 , \DD]h  &=  -\frac{1}{2}\tr X\DD h+(\Hb+Z)\nab_4 h -\frac 1 2 \Xh \c\ov{\DD} h+\Xi \nab_3h , \\
 \, [\nab_3 , \DD]h  &=   -\frac{1}{2}\tr \Xb\DD h+(H-Z)\nab_3 h -\frac 1 2 \Xbh\c\ov{\DD} h+\Xib \nab_4h. 
\end{split}
 \eea    
    
   \item
    Let $F\in  \sk_1 (\mathbb{C})$.  Then
    \bea\label{eq:comm:nab4-nab3-DDhot-precise}
    \begin{split}
\, [\nab_4,  \DD \hot] F  &=-\frac 1 2 \tr X\left( \DD \hot  F+\Hb\hot F\right)+(\Hb+Z)\hot\nab_4 F+ \Xi \hot \nab_3 F \\
 &-B \hot F - \frac 1 2 \tr \Xb \Xi \hot  F-\frac 1 2\Xh \c \ov{\DD} F+\frac12\Xh (\ov{\Hb}\c F)+ (\Ga_b \c \Ga_g) F, \\
\, [\nab_3,  \DD \hot] F  &=-\frac 1 2 \tr \Xb\left( \DD \hot  F+H \hot F\right)+(H-Z)\hot\nab_3 F+ \Xib \hot \nab_4 F \\
 &+\Bb \hot F - \frac 1 2 \tr X \Xib \hot  F-\frac 1 2\Xbh \c \ov{\DD} F+\frac12\Xbh (\ov{H}\c F)+ (\Ga_b \c \Ga_g) F.
 \end{split}
\eea
Using the schematic structure of the error terms, the above can be written as
\bea
 \lab{commutator-nab-43-D-hot} 
\bsplit
\, [\nab_4, \mathcal{D}\hot ]F&=- \frac 1 2 \tr X( \mathcal{D}\hot F + \underline{H} \hot F)+ (\underline{H}+Z) \hot \nab_4 F\\
&+\Xi \c \nabc_3 F+ r^{-1} \Ga_g \c  \dk^{\leq 1} F, \\
\, [\nab_3, \mathcal{D}\hot] F&=- \frac 1 2 \tr \Xb( \mathcal{D}\hot F + H \hot F)+ (H-Z) \hot \nab_3 F+r^{-1}\Ga_b \c \dk^{\leq 1} F.
\end{split}
\eea

\item Let $U\in \sk_2(\mathbb{C})$. Then
\bea\label{eq:comm:nab4-nab3-DDc-precise}
\begin{split}
\,[\nab_4,  \ov{\DD} \c ] U &=-\frac 1 2\ov{\tr X} \big(  \ov{\DD}\c U - 2\ov{\Hb} \c U\big) +\ov{(\Hb + Z)}\c\nab_4 U+ \ov{\Xi} \c \nab_3 U \\
&+2\ov{B} \c U  -\frac 1 2 \ov{\tr \Xb} \ov{\Xi}\c  U  -\frac 1 2 \Xh \c \ov{\DD} U-\frac 1 2 (\ov{\Xh}\c U)\ov{\Hb}+ (\Ga_b \c \Ga_g) U, \\
\,[\nab_3,  \ov{\DD} \c ] U &=-\frac 1 2\ov{\tr \Xb} \big(  \ov{\DD}\c U - 2\ov{H} \c U\big) +\ov{(H- Z)}\c\nab_3 U+ \ov{\Xib} \c \nab_4 U \\
&-2\ov{\Bb} \c U  -\frac 1 2 \ov{\tr X} \ov{\Xib}\c  U  -\frac 1 2 \Xbh \c \ov{\DD} U-\frac 1 2 (\ov{\Xbh}\c U)\ov{H}+ (\Ga_b \c \Ga_g) U.
\end{split}
\eea
Using the schematic structure of the error terms, the above can be written as
 \bea\label{commutator-nabc-3-ov-DDc-U}
 \bsplit
\, [\nab_4, \ov{\DD}\c] U&=- \frac 1 2\ov{\tr X}\, ( \ov{\DD} \c U - 2 \ov{\Hb} \c U)+\ov{(\Hb+Z)} \c \nab_4 U\\
&+\Xi \c \nabc_3 U +r^{-1} \Ga_g  \c \dk^{\leq 1} U, \\
\, [\nab_3, \ov{\DD}\c] U&=- \frac 1 2\ov{\tr\Xb}\, ( \ov{\DD} \c U -  2 \ov{H} \c U)+\ov{(H-Z)} \c \nab_3 U + r^{-1}\Ga_b \c  \dk^{\leq 1} U.
\end{split}
\eea
Similarly, for $F\in  \sk_1 (\mathbb{C})$ we have
 \bea\label{commutator-nab4-ov-DDcF}
 \bsplit
\, [\nab_4, \ov{\DD}\c] F&=- \frac 1 2\ov{\tr X}\, ( \ov{\DD} \c F -  \ov{\Hb} \c F)+\ov{(\Hb+Z)} \c \nab_4 F \\
&+\ov{\Xi} \c \nabc_3 F+ r^{-1} \Ga_g  \c \dk^{\leq 1} F, \\
\, [\nab_3, \ov{\DD}\c] F&=- \frac 1 2\ov{\tr\Xb}\, ( \ov{\DD} \c F -   \ov{H} \c F)+\ov{(H-Z)} \c \nab_3 F +r^{-1}\Ga_b \c  \dk^{\leq 1} F.
\end{split}
\eea
\item Let $U\in \sk_2(\mathbb{C})$. Then
\bea\label{correct-commutator-1} \lab{commutator-nab-4-D-c} \label{commutator-nab-3-nab-4-Psi}\label{correct-commutator}
\begin{split}
\,[\nab_3, \nab_4]U &= - 2\om \nab_3 U+ 2\omb \nab_4 U  + 2 (\eta_c-\etab_c) \nab_c U +4i \left(- \rhod+ \eta \wedge \etab  \right) U\\
&+\left(\Ga_b  \c \Ga_g \right) U.
\end{split}
\eea
Also,
\bea\label{commutator-nab-3-nab-a-U}
\bsplit
\, [\nab_3, \nab_a] U_{bc}&=-\frac  1 2   \trchb\, \Big(\nab_a U_{bc}+\eta_bU_{ac}+\eta_c U_{ab}-\de_{a b}(\eta \c U)_c-\de_{a c}(\eta \c U)_b \Big)\\
&-\frac 1 2 \atrchb\, \Big(\dual \nab_a  U_{bc} +\eta_b \dual U_{ac}+\eta_c \dual U_{ab}- \in_{a b}(\eta \c  U)_c- \in_{a c}(\eta \c  U)_b \Big)\\
 &+(\eta_a-\ze_a)\nab_3 U_{bc}+r^{-1}\Ga_b \c \dk^{\leq 1} U,
 \end{split}
\eea
\end{enumerate}
where the above error terms may also contain terms which are 
quadratic   in  the perturbation and  enjoy  better  decay  properties,  or are higher order  and decay at least as good.
\end{lemma}

\begin{proof} 
See section \ref{sec:proof-lemma-comm-form-1}. 
\end{proof}

 We collect here the commutation formulas for the conformal derivatives introduced in Lemma \ref{lemma:definition-conformal-derivatives}.
  
     \begin{lemma}\label{COMMUTATOR-NAB-C-3-DD-C-HOT}\label{commutator-nab-3-nab-4}\label{commtator-3-a}\label{comm:scalar}\label{comm:1form-div}\label{comm:1form}\label{commutator-nabc-4-F-formula}
   The following commutation formulas hold true.
   \begin{enumerate}
   \item Let $h \in \sk_0(\CCC)$ $s$-conformally invariant. Then 
   \bea\label{eq:comm-nabc4-DDc-h-precise}
   \begin{split}
 \, [\nabc_4 , \DDc]h  &= -\frac{1}{2}\tr X\DDc h+\Hb\nabc_4 h -\frac{1}{2}\Xh\c\ov{\DDc}h+\Xi\nabc_3h\\
&  +s\left(\frac{1}{2}\tr X\Hb  +\frac{1}{2}\widehat{X}\c\ov{\Hb}-\frac{1}{2}\tr\Xb\Xi -B\right)h  + (\Ga_b \c \Ga_g) h, \\
 \, [\nabc_3 , \DDc]h  &= -\frac{1}{2}\tr \Xb\DDc h+H \nabc_3 h -\frac{1}{2}\Xbh\c\ov{\DDc}h+\Xib\nabc_4h\\
&  -s\left(\frac{1}{2}\tr \Xb H  +\frac{1}{2}\Xbh\c\ov{H}-\frac{1}{2}\tr X \Xib +\Bb\right)h  + (\Ga_b \c \Ga_g) h.
\end{split}
 \eea
 Using the schematic structure of the error terms, the above can be written as
     \bea\label{eq:comm-nabc4-DDc-h-err}
     \begin{split}
 \, [\nabc_4 , \DDc]h  &= -\frac{1}{2}\tr X\DDc h  +s \frac{1}{2}\tr X\Hb  h+\Hb\nabc_4 h  \\
 &+\Xi  \nabc_3 h+ r^{-1} \Ga_g  \c \dk^{\leq 1} h, \\
  \, [\nabc_3 , \DDc]h  &= -\frac{1}{2}\tr \Xb \DDc h  -s \frac{1}{2}\tr \Xb H   h+H \nabc_3 h  +r^{-1}\Ga_b  \c \dk^{\leq 1} h.
  \end{split}
 \eea
 We also have
 \bea\label{eq:comm-nabc3-nabc4-h}
 \, [\nabc_3, \nabc_4] h  &=&  2(\eta-\etab) \c \nabc  h +2s(\rho   -\eta\c\etab)h +(\Ga_b \cdot \Ga_g) h.
 \eea

\item  Let $F \in \sk_1 (\mathbb{C})$ $s$-conformally invariant. Then
\bea\label{eq:comm-nabc4nabc3DDchot-precise}
\begin{split}
[ \nabc_4, \DDc \hot] F   &= -\frac 1 2 \tr X\left( \DDc \hot  F+(1-s)\Hb\hot F\right)+\Hb\hot\nabc_4 F+ \Xi \hot \nabc_3 F\\
 &-(s+1)B \hot F - (s+1)\frac 1 2 \tr \Xb \Xi \hot  F-\frac 1 2\Xh \c \ov{\DDc} F\\
 & +\frac12\Xh (\ov{\Hb}\c F) +s  \frac{1}{2}(\widehat{X}\c\ov{\Hb})  \hot F+ (\Ga_b \c \Ga_g) F, \\
 [ \nabc_3 , \DDc \hot] F   &=-\frac 1 2 \tr \Xb\left( \DDc \hot  F+(s+1)H \hot F\right)+H\hot\nabc_3 F + \Xib \hot \nabc_4 F \\
 &-(s-1)\Bb \hot F+( s-1) \frac 1 2 \tr X \Xib \hot  F-\frac 1 2\Xbh \c \ov{\DDc} F\\
 &+\frac12\Xbh (\ov{H}\c F) -s\frac{1}{2}(\widehat{\Xb}\c\ov{H}) \hot F + (\Ga_b \c \Ga_g) F.
 \end{split}
\eea
Using the schematic structure of the error terms, the above can be written as
\bea\label{eq:comm-nabc4nabc3-DDchot-err}
\begin{split}
 \, [\nabc_4 , \DDc \hot ]F &=- \frac 1 2 \tr X\left( \DDc\hot F + (1-s)\Hb\hot F\right)+ \underline{H} \hot \nabc_4 F\\
  &+\Xi \hot \nabc_3 F+ r^{-1} \Ga_g  \c \dk^{\leq 1} F,  \\
 \, [\nabc_3, \DDc \hot ]F &=- \frac 1 2 \tr \Xb \left( \DDc \hot F + (1+s)H \hot F \right)  + H \hot \nabc_3 F+ r^{-1}\Ga_b  \c \dk^{\leq 1} F.
 \end{split}
\eea
\item Let $U\in \sk_2(\mathbb{C})$ $s$-conformally invariant. Then
\bea\label{eq:comm-nabc4nabc3-ovDDc-U-precise}
\begin{split}
\,[\nabc_4, \ov{\DDc}\c] U&= -\frac 1 2\ov{\tr X} \big(  \ov{\DDc}\c U - (s+2)\ov{\Hb} \c U\big) +\ov{\Hb}\c\nabc_4 U\\
&+ \ov{\Xi} \c \nabc_3 U -(s-2)\ov{B} \c U  -(s+1)\frac 1 2 \ov{\tr \Xb} \ov{\Xi}\c  U  -\frac 1 2 \Xh \c \ov{\DDc} U\\
&-\frac 1 2 (\ov{\Xh}\c U)\ov{\Hb}+ s \frac{1}{2}(\ov{\widehat{X}}\c\Hb  )\c U+ (\Ga_b \c \Ga_g) U,\\
\,[\nabc_3, \ov{\DDc}\c] U&= -\frac 1 2\ov{\tr \Xb} \big(  \ov{\DDc}\c U +(s-2)\ov{H} \c U\big) +\ov{H}\c\nabc_3 U\\
&+ \ov{\Xib} \c \nabc_4 U +(s+2)\ov{\Bb} \c U  +(s-1)\frac 1 2 \ov{\tr X} \ov{\Xib}\c  U  -\frac 1 2 \Xbh \c \ov{\DDc} U\\
&-\frac 1 2 (\ov{\Xbh}\c U)\ov{H}- s \frac{1}{2}(\ov{\Xbh}\c H  )\c U+ (\Ga_b \c \Ga_g) U.
\end{split}
\eea

Using the schematic structure of the error terms, the above can be written as
 \bea\label{eq:comm-nabc4nabc3-ovDDc-U-err}
 \begin{split}
 \, [\nabc_4, \ov{\DDc}\c] U&=- \frac 1 2\ov{\tr X}\, ( \ov{\DDc} \c U -(s+ 2) \ov{\Hb} \c U)+\ov{\Hb} \c \nabc_4 U\\
 &+\ov{\Xi} \c \nabc_3 U+ r^{-1}\Ga_g \c  \dk^{\leq 1} U,\\
\, [\nabc_3, \ov{\DDc}\c] U&=- \frac 1 2\ov{\tr\Xb}\, ( \ov{\DDc} \c U +(s-  2) \ov{H} \c U)+\ov{H} \c \nabc_3 U+ r^{-1}\Ga_b \c  \dk^{\leq 1} U.
\end{split}
\eea
Similarly, for $F \in \sk_1(\CCC)$,
  \bea\label{eq:comm-nabc4nabc3-ovDDc-F-err}
  \begin{split}
 \, [ \nabc_4 ,\ov{\DDc} \c] F     &=- \frac 1 2\ov{\tr X}\, \left( \ov{\DDc} \c F - (s+1) \ov{\Hb} \c F \right)+\ov{\Hb} \c \nabc_4 F\\
 &+\ov{\Xi} \c \nabc_3 F+ r^{-1} \Ga_g  \c \dk^{\leq 1} F,\\
 \, [\nabc_3, \ov{\DDc}\c] F&=- \frac 1 2\ov{\tr\Xb}\, ( \ov{\DDc} \c F +(s-  1) \ov{H} \c F)+\ov{H} \c \nabc_3 F+ r^{-1}\Ga_b \c  \dk^{\leq 1} F.
 \end{split}
 \eea
 \item Let $U\in \sk_2(\mathbb{C})$ $s$-conformally invariant. Then
  \beaa
\, [\nabc_3, \nabc_4] U  &=&  2(\eta-\etab ) \c  \nabc  U  +\Big( 2s\left(\rho   -\eta\c\etab\right) +4i \left(- \rhod+ \eta \wedge \etab \right) \Big)U\\
 &&+(\Ga_b \cdot \Ga_g) U,\\
 &=&  2(\eta-\etab ) \c  \nabc  U  +\Big((s-2)P+(s+2)\ov{P} -2s\eta\c\etab  +4i \eta \wedge \etab  \Big)U\\
 &&+(\Ga_b \cdot \Ga_g) U.
 \eeaa
Also,
\beaa
 [\nabc_3, \nabc_a] U_{bc}&=&\eta_a\nabc_3 U_{bc}-\frac  1 2   \trchb\, \Big[\nabc_a U_{bc}+s(\eta_a) U_{bc}+\eta_bU_{ac}+\eta_c U_{ab}\\
 &&-\de_{a b}(\eta \c U)_c-\de_{a c}(\eta \c U)_b \Big]\\
&&-\frac 1 2 \atrchb\, \Big[\dual \nabc_a  U_{bc}+s (\dual\eta_a) U_{bc} +\eta_b \dual U_{ac}+\eta_c \dual U_{ab}\\
&&- \in_{a b}(\eta \c  U)_c- \in_{a c}(\eta \c  U)_b \Big]+r^{-1}\Ga_b \c  \dk^{\leq 1} U.
 \eeaa
\end{enumerate}
\end{lemma}

\begin{proof} 
See section \ref{sec:proof-lemma-comm-2}. 
\end{proof}


\section{Approximate Killing vectorfields $\T$ and $\Z$}


\begin{definition}\lab{Definition:vfsTZ} 
In $\MM$, we define $\T$ and $\Z$ as follows:
\begin{itemize}
\item If the normalization of $(e_3, e_4)$ is ingoing, we have
\beaa
\T &:=& \frac{1}{2}\left(e_4+\frac{\Delta}{|q|^2}e_3 -2a\Re(\Jk)^be_b\right),\\
\Z &:=& \frac 1 2 \left(2(r^2+a^2)\Re(\Jk)^be_b -a(\sin\th)^2 e_4 -\frac{a(\sin\th)^2\De}{ |q|^2} e_3\right).
\eeaa

\item If the normalization of $(e_3, e_4)$ is outgoing, we have
\beaa
\T &:=& \frac{1}{2}\left(e_3+\frac{\Delta}{|q|^2}e_4 -2a\Re(\Jk)^be_b\right),\\
\Z &:=& \frac 1 2 \left(2(r^2+a^2)\Re(\Jk)^be_b -a(\sin\th)^2 e_3 -\frac{a(\sin\th)^2\De}{ |q|^2} e_4\right).
\eeaa

\end{itemize}
\end{definition} 

From relations \eqref{eq:T-Z-ingoing-frame-Kerr} and \eqref{eq:T-Z-outgoing-frame-Kerr} and the values of $\Jk$ as in Definition \ref{def:JkandJ}, one can see that the above vectorfields reduce to the Killing vectorfields $\T$ and $\Z$ in Kerr.

The following lemma shows that $\T$ and $\Z$  are almost Killing vectorfields.

\begin{lemma}\label{LEMMA:DEFORMATION-TENSORS-T}\label{lemma:deformation-tensor-T}
Let $\piT$ and $\piZ$ be the deformation tensors of $\T$ and $\Z$ as defined above. We have
\beaa
 \piT_{44}, \piT_{4a} \in \Ga_g, \qquad\,\, \piT_{33}, \piT_{34}, \piT_{3a}, \piT_{ab}  \in \Ga_b,
\eeaa
and
\beaa
\,\piZ_{44}, \, r^{-1} \piZ_{4a} \in \Ga_g,  \qquad   \piZ_{33}, \piZ_{34}, \piZ_{3a}, \piZ_{ab} \in  r\Ga_b.
\eeaa
Moreover, 
\beaa
\tr \piT \in \Ga_g, \qquad  \tr \piZ  \in r\Ga_b,
\eeaa
and 
\begin{align*}
(\Div \piT)_3 &\in r^{-1} \dk^{\leq 1} \Ga_b,  &\qquad (\Div \piT)_4, (\Div \piT)_a &\in \dk^{\leq 1} \Ga_g, \\
(\Div \piZ)_4 &\in \dk^{\leq 1} \Ga_g ,  & \qquad (\Div \piZ)_3 &\in \dk^{\leq 1} \Ga_b, \\
 (\Div \piZ)_a &\in r \dk^{\leq 1} \Ga_g.
\end{align*}
\end{lemma} 

\begin{proof}
See section \ref{appendix-proof-deformation-tensor}.
\end{proof}

We collect here the following   commutator identities for $\T$ and $\Z$ and the D'Alembertian, for scalar functions and for horizontal symmetric 2-tensors. For horizontal symmetric 2-tensors additional terms coming from curvature appear.

\begin{proposition}
\lab{LE:COMMTZSQUARE}
The following  commutation  formulas hold true  for a scalar $\psi$:
\bea
\lab{eq:commTZsquare}
\bsplit
 \,[ \T, \square_\g] \psi&=\dk \big(\Ga_g \c \dk \psi\big)+\Ga_b \c \square_\g\psi ,\\
   \,[ \Z, \square_\g] \psi &= \dk \big(\Ga_g \c \dk \psi\big)+r\Ga_b \c \square_\g\psi.
 \end{split}
 \eea
 The following  commutation  formulas hold true  for $\psi\in \sk_2$:
\bea
\lab{eq:commTZsquare-2tensors}
\bsplit
 \,[ \nab_\T, \squared_2] \psi&=O(ar^{-4})\dk^{\leq 1}\psi +\dk \big(\Ga_g \c \dk \psi\big)+\Ga_b \c \squared_2\psi ,\\
   \,[ \nab_\Z, \squared_2] \psi &= O(r^{-2})\dk^{\leq 1}\psi +\dk \big(\Ga_g \c \dk \psi\big)+r\Ga_b \c \squared_2\psi.
 \end{split}
 \eea
\end{proposition}

\begin{proof} 
See section \ref{proof:comm-T-Z-square} for the proof of \eqref{eq:commTZsquare}.
By putting together Proposition \ref{prop:commutators-nabT-nabZ-squared2-Kerr} and \eqref{eq:commTZsquare} we deduce \eqref{eq:commTZsquare-2tensors}.
\end{proof}

Recall the definition of horizontal Lie derivative $\Lied_X$ for $X \in \T(\MM)$ as given in Definition \ref{definition:hor-Lie-derivative}. We collect here the commutator identities for $\Lied_\T$ and $\Lied_\Z$ and the D'Alembertian operator.

\begin{corollary}\label{cor:commutator-Lied-squared} 
The following commutation formulas holds true for $\psi \in \sk_2$:
\beaa
\, [\Lied_\T, \squared_2]\psi_{ab}&=&  \dk \big(\Ga_g \c \dk \psi\big)+\Ga_b \c \squared_2\psi,\\
\, [\Lied_\Z, \squared_2]\psi_{ab} &=&  \dk \big(\Ga_g \c \dk \psi\big)+r\Ga_b \c \squared_2\psi.
\eeaa
\end{corollary}
\begin{proof}
This is a corollary of Proposition \ref{prop:commutator-Lied-squared}, according to which
\beaa
\, [\Lied_X, \squared_2]\psi_{ab}&=& -\piX^{\mu\nu} \Ddot_\mu \Ddot_\nu\psi_{ab} - {{\GaX}^\mu}_{\mu  \rho}\Ddot^{\rho} \psi_{ab}\\
&&-2\GabbX_{a\mu c}{\Ddot^{\mu}}{\psi^c}_{ b}-2\GabbX_{b\mu c}\Ddot^{\mu}{\psi_{ a}}^c -  \Ddot^\nu ( \GabbX_{a \nu c}){\psi^c}_b-\Ddot^\nu(\GabbX_{b \nu c}){\psi_a}^c
\eeaa
where 
\beaa
\GaX_{\a\b\mu}&=&\frac 1 2 (\D_\a\piX_{\b\mu}+\D_\b\piX_{\a\mu}-\D_\mu\piX_{\a\b}),
\eeaa
and similarly for $\GabbX$. In particular, we can write
\beaa
\, [\Lied_X, \squared_2]\psi_{ab}&=& -\piX^{\mu\nu} \Ddot_\mu \Ddot_\nu\psi_{ab} - \big( \Div \piX_\rho - \frac 1 2 \D_\rho \piX\big)\Ddot^{\rho} \psi_{ab}\\
&&-2\GabbX_{a\mu c}{\Ddot^{\mu}}{\psi^c}_{ b}-2\GabbX_{b\mu c}\Ddot^{\mu}{\psi_{ a}}^c -  \Ddot^\nu ( \GabbX_{a \nu c}){\psi^c}_b-\Ddot^\nu(\GabbX_{b \nu c}){\psi_a}^c.
\eeaa
Observe that the first line gives the same expression as for $[X, \square_\g]\psi$, i.e. \eqref{eq:generalcommutation-scalarwave}, and can then
 be computed as in Proposition \ref{LE:COMMTZSQUARE}.
For the second line we compute
\beaa
 \Ddot^\nu ( \GabbX_{a \nu c})&=&\frac 1 2 (\D^\mu\nab_a\piX_{\mu c}+\D^\mu\nab_\mu\piX_{ac}-\D^\mu\nab_c\piX_{a\mu})\\
 &=&\frac 1 2 (\nab_a \Div\piX_{ c}-\nab_c \Div \piX_{a}+\square_\g\piX_{ac})+O(r^{-2}) \piX,
\eeaa
and
\beaa
\GabbX_{a\mu c}{\Ddot^{\mu}}{\psi^c}_{ b}&=&\Big( r^{-1} \big(\GabbX_{ad c}, \GabbX_{a3 c}\big),\GabbX_{a4 c}\Big) \dk \psi_{cb}.
\eeaa 
Using Lemma \ref{LEMMA:DEFORMATION-TENSORS-T} and the improved decay for $\piZ$ in
\eqref{eq:D4-D3-trpiZ}, we obtain
\beaa
 -  \Ddot^\nu ( {}^{(\T)}\Gabb_{a \nu c}){\psi^c}_b-\Ddot^\nu({}^{(\T)}\Gabb_{b \nu c}){\psi_a}^c&=& r^{-1} \dk^{\leq 2} \Ga_b \c \psi, \\
   -  \Ddot^\nu ( {}^{(\Z)}\Gabb_{a \nu c}){\psi^c}_b-\Ddot^\nu({}^{(\Z)}\Gabb_{b \nu c}){\psi_a}^c&=&  \dk^{\leq 2} \Ga_g  \c \psi,
 \eeaa
 and
 \beaa
 {}^{(\T)}\Gabb_{ad c}&=&\frac 1 2 (\nab_a\piT_{d c}+\nab_d\piT_{ac}-\nab_c\piT_{ad})=r^{-1} \dk\Ga_b,\\
 {}^{(\T)}\Gabb_{a3 c}&=&\frac 1 2 (\nab_a\piT_{3 c}+\nab_3\piT_{ac}-\nab_c\piT_{a3})= \dk \Ga_b, \\
 {}^{(\T)}\Gabb_{a4 c}&=&\frac 1 2 (\nab_a\piT_{4 c}+\nab_4\piT_{ac}-\nab_c\piT_{a4})=r^{-1} \dk \Ga_b,\\
 {}^{(\Z)}\Gabb_{ad c}&=&\frac 1 2 (\nab_a\piZ_{d c}+\nab_d\piZ_{ac}-\nab_c\piZ_{ad})=\dk \Ga_b,\\
{}^{(\Z)}\Gabb_{a3 c}&=&\frac 1 2 (\nab_a\piZ_{3 c}+\nab_3\piZ_{ac}-\nab_c\piZ_{a3})= \Ga_b, \\
{}^{(\Z)}\Gabb_{a4 c}&=&\frac 1 2 (\nab_a\piZ_{4 c}+\nab_4\piZ_{ac}-\nab_c\piZ_{a4})=\Ga_g.
\eeaa
Putting the above together, this concludes the proof.
\end{proof}


\section{Inverse metric and $\RR$, $O$ tensors}


In $\MM$, we define the vectorfields $\That$ and $\Rhat$ as follows:
\begin{itemize}
\item If the normalization of $(e_3, e_4)$ is ingoing, we have
\beaa
\That &:=&\frac 1 2 \left( \frac{|q|^2}{r^2+a^2} e_4+\frac{\De}{r^2+a^2}  e_3\right) ,\\
\Rhat &:=&\frac 1 2 \left( \frac{|q|^2}{r^2+a^2} e_4-\frac{\De}{r^2+a^2}  e_3\right).
\eeaa
\item If the normalization of $(e_3, e_4)$ is outgoing, we have
\beaa
\That &:=& \frac 1 2 \left( \frac{\De}{r^2+a^2} e_4+\frac{|q|^2}{r^2+a^2}  e_3\right),\\
\Rhat &:=&\frac 1 2 \left( \frac{\De}{r^2+a^2} e_4-\frac{|q|^2}{r^2+a^2}  e_3\right).
\eeaa
\end{itemize}

From relations \eqref{eq:ThatRhat-e_3e_4-outgoing} and \eqref{eq:ThatRhat-e_3e_4-Kerr}, one can see that the above vectorfields reduce to $\That$ and $\Rhat$ in Kerr.  Note also that $\That, \Rhat$ are perpendicular to  the horizontal  structure and $\g(\That, \Rhat)=0$.

\begin{lemma}\label{lemma:commutator-Rhat-That} 
The vectorfields $\Rhat$ and  $\That$ satisfy
\beaa
\, [\Rhat, \That]&=&  - \frac{ar \De}{(r^2+a^2)^3}\Z+r^{-1}  \Ga_b \c  \dk.
\eeaa
\end{lemma}

\begin{proof} 
Using the definition in the outgoing normalization, we have 
\beaa
[\Rhat, \That]&=& \frac 1 4 \left[  \frac{\De}{r^2+a^2} e_4-\frac{|q|^2}{r^2+a^2}  e_3,  \frac{\De}{r^2+a^2} e_4+\frac{|q|^2}{r^2+a^2}  e_3\right]\\
&=& \frac 1 2 \left[  \frac{\De}{r^2+a^2} e_4, \frac{|q|^2}{r^2+a^2}  e_3\right]\\
&=& \frac 1 2  \frac{|q|^2 \De}{(r^2+a^2)^2} [e_4, e_3] + \frac 1 2  \frac{\De}{r^2+a^2} e_4\left(\frac{|q|^2}{r^2+a^2} \right) e_3 - \frac 1 2 \frac{|q|^2}{r^2+a^2} e_3\left(\frac{\De}{r^2+a^2} \right) e_4.
\eeaa
Using Lemma \ref{lemma:comm} to write in the outgoing normalization 
\beaa
 [e_4, e_3] &=& 2(\etab-\eta ) \c \nab + 2 \om e_3  -2\omb e_4 \\
 &=&-\frac{4ar}{|q|^2} \Re(\Jk^b)e_b  - \pr_r (\frac{\De}{|q|^2}) e_4+r^{-1}  \Ga_b \c  \dk ,
 \eeaa
 and using \eqref{eq:computations-e3e4-out}
we obtain
\beaa
[\Rhat, \That]&=&- \frac{ar \De}{(r^2+a^2)^3} \Big(  2(r^2+a^2)\Re(\Jk^b)e_b- a \sin^2\th e_3  -  \frac{\De a\sin^2\th}{|q|^2} e_4 \Big)+r^{-1}  \Ga_b \c  \dk\\
&=& - \frac{ar \De}{(r^2+a^2)^3}\Z+r^{-1}  \Ga_b \c  \dk,
\eeaa
as stated.
\end{proof}

Observe that with the above definitions we have 
\bea\label{eq:That-T-Z-perturbations}
\That=\T+\frac{a}{r^2+a^2}\Z.
\eea
 Indeed, for example in the outgoing normalization,
\beaa
\T+\frac{a}{r^2+a^2}\Z&=&\frac{1}{2}\left(e_3+\frac{\Delta}{|q|^2}e_4 -2a\Re(\Jk)^be_b\right)\\
&&+\frac{a}{r^2+a^2} \frac 1 2 \left(2(r^2+a^2)\Re(\Jk)^be_b -a(\sin\th)^2 e_3 -\frac{a(\sin\th)^2\De}{ |q|^2} e_4\right)\\
&=& \frac 1 2 \left( \frac{\De}{r^2+a^2} e_4+\frac{|q|^2}{r^2+a^2}  e_3\right)=\That.
\eeaa
 
 We also define, in terms of the horizontal structure, 
\bea\label{def-O-a-b-pert-Kerr-new}
O^{\a\b}&:=& |q|^2\ga^{ab}e_a^\a e_b^\b,
\eea
where we recall that $\ga_{ab}=\g(e_a, e_b)$, i.e. $\ga$ is the metric induced by $\g$ on the horizontal structure. From \eqref{eq:exxpressionO-e1e2}, one can see that the above tensor reduce to $O^{\a\b}$ in Kerr.

The definition $O^{\a\b}$ allows us to express the inverse metric in perturbations of Kerr as follows, see Lemma \ref{lemma:inversemetricexpressioninKerr} for the case of Kerr. 
\begin{lemma}\label{lemma:inverse-metric-perturbations}
Let $(\MM, \g)$ be a spacetime, with $\Rhat$, $\That$ and $O$ defined as above. Then 
\bea
|q|^2\g^{\a\b}&=&\frac{(r^2+a^2)^2}{\De} \Rhat^\a \Rhat^\b +\frac{1}{\De}\RR^{\a\b}
\eea
where
\bea\label{def-RR-pert-Kerr}
\RR^{\a\b}&:=& -(r^2+a^2)^2 \That^\a \That^\b+ \De O^{\a\b}.
\eea
 Thus the inverse metric can also be written in the form\footnote{Observe that this expression of the metric is not regular at the horizon.}
 \bea
 \lab{inverse-metric-vfs-perturbations}
|q|^2 \g^{\a\b}&=&\frac{(r^2+a^2)^2}{\De} \big( -\That^\a\That^\b +      \Rhat^\a \Rhat^\b\big) + O^{\a\b}.
 \eea
\end{lemma}

\begin{proof} 
Using the definition of frames in $(\MM, \g)$, we can write
\beaa
\g^{\a\b}&=& -\frac 1 2 \big( e_3^\a e_4^\b+ e_3^\b e_4^\a \big)+ \ga^{ab}e_a^\a e_b^\b.
\eeaa
Using \eqref{def-O-a-b-pert-Kerr-new}, we obtain
\beaa
|q|^2\g^{\a\b}&=& -\frac 1 2 |q|^2 \big( e_3^\a e_4^\b+ e_3^\b e_4^\a \big)+ O^{\a\b}.
\eeaa
 We then compute, for example in the ingoing normalization,
 \beaa
 \Rhat^\a \Rhat^\b -  \That^\a \That^\b&=&\frac 1 4  \left( \frac{|q|^2}{r^2+a^2} e^\a_4-\frac{\De}{r^2+a^2}  e^\a_3\right)  \left( \frac{|q|^2}{r^2+a^2} e^\b_4-\frac{\De}{r^2+a^2}  e^\b_3\right)\\
&& - \frac 1 4  \left( \frac{|q|^2}{r^2+a^2} e^\a_4+\frac{\De}{r^2+a^2}  e^\a_3\right) \left( \frac{|q|^2}{r^2+a^2} e^\b_4+\frac{\De}{r^2+a^2}  e^\b_3\right)\\
&=&\frac{\De}{(r^2+a^2)^2} \left(-\frac 1 2 |q|^2  e_3^\a e_4^\b -\frac 1 2 |q|^2 e_3^\b e_4^\a\right),
 \eeaa
which proves the lemma.
\end{proof}

As in Definition \ref{definition:tensors-S}, we define the following approximate symmetric tensors.

\begin{definition}\label{definition:tensors-S-pert:0} 
We define the following symmetric spacetime 2-tensors:
\beaa
S_1^{\a\b}&:=& \T^\a \T^\b, \\
S_2^{\a\b}&:=& a\T^{(\a} \Z^{\b)}, \\
S_3^{\a\b}&:=& a^2\Z^\a \Z^\b, \\
S_4^{\a\b}&:=& O^{\a\b}=|q|^2\ga^{ab}e_a^\a e_b^\b.
\eeaa
We denote the set of the above tensors as $S_\aund$, for $\aund=1,2,3,4$. 
\end{definition}

Using \eqref{eq:That-T-Z-perturbations} to write $\That=\T+\frac{a}{r^2+a^2}\Z$ in \eqref{def-RR-pert-Kerr}, we can write 
 \bea\label{eq:RR-Sa-pert}
          \RR^{\a\b} =\RR^\aund  \Sa^{\a\b},
          \eea
          with $\RR^{\aund}$, $\aund=1, 2, 3, 4$ given by
          \bea\label{components-RR-aund-pert}
          \RR^1&=&-(r^2+a^2)^2, \quad \RR^2 = -2(r^2+a^2), \quad \RR^3 =-1, \quad \RR^4=\De.
          \eea


\section{Approximate Carter tensor and Carter operator}\label{section:modified-laplacian-pert}


We now extend the definition of Carter tensor as given in Definition \ref{definition:Carter-tensor-Kerr} to the case of perturbations of Kerr.

\begin{definition} 
\lab{Definition:CarterOperator}
In $(\MM, \g)$ the Carter tensor  is defined as the following symmetric 2-tensor $K$:
\bea\label{definition-carter}
K^{\a\b}&=& -(a^2\cos^2\th)  \g^{\a\b} +O^{\a\b},
\eea
where the tensor $O$ is defined in \eqref{def-O-a-b-pert-Kerr-new}.
\end{definition}

Recall the definition of Carter operator $\KK$ for horizontal tensors $\psi \in \sk_k$ according to Definition \ref{definition:operator-KK-general}, given by
\beaa
\KK(\psi):=\Ddot_\b(K^{\a\b}\Ddot_\a \psi).
\eeaa

Associated to the Carter operator $\KK$ we define the following.

\begin{definition}\label{definition:OO-perturbations} 
In $(\MM, \g)$ we define the following second order angular operator for $\psi \in \sk_k$:
\bea\label{definition-SS}
\OO(\psi) &:=& |q|^2 \left(\lap_k \psi - \frac{2a^2\cos\th}{|q|^2} \dual\Re(\Jk)^b \nab_b \psi   \right).
\eea
\end{definition}

Note that the above definition reduces to the operator $\OO$ in Kerr,  see  \eqref{eq:OO-Jk-kerr}.

We now show that $\OO$ is an approximate symmetry operator for scalars and, up to Riemann curvature terms, for tensors in perturbations of Kerr. 
\begin{proposition}\label{LEMMA:MOD-LAPLACIAN-PERT-KERR} 
The operator $\OO$ defined in \eqref{definition-SS} for a scalar function $\psi$ satisfies the following commutation formula:
\bea\label{eq:comm-OO-q^2squared-scalars}
\, [\OO, |q|^2 \square_\g]\psi  &=& |q|^2\Big[  \dk^2 \big( \Ga_g \c \dk \psi)+ \Ddot_3 \dk \big(|q|^2 \xi \c \Ddot_a \psi \big)\Big].
\eea

The operator $\OO$ for $\psi \in \sk_2$ satisfies the following commutation formula:  
\bea\label{eq:comm-OO-q^2squared-tensors}
\begin{split}
\, [\OO, |q|^2 \squared_2]\psi  &= |q|^2\Bigg[ \nab\left(\frac{8a(r^2+a^2)\cos\th}{|q|^2}\right)\c\nab\nab_\That\dual\psi +O( a r^{-2}) \nab_{\Rhat}^{\leq 1}\dk^{\leq 1} \psi\\
&+  \dk^2 \big( \Ga_g \c \dk \psi)+ \Ddot_3 \dk \big(|q|^2 \xi \c \Ddot_a \psi \big)\Bigg].
\end{split}
\eea
\end{proposition}

\begin{proof} 
See Section \ref{proof-lemma:commutator-OO-|q|2square-scalar}. Observe that in order to obtain acceptable error terms like the ones on the right hand side of \eqref{eq:comm-OO-q^2squared-scalars} we need to derive the commutator using the decomposition in null frames of $\square_\g$ as in Lemma \ref{lemma:expression-wave-operator}. Similarly for $\squared_2$. The curvature terms on the first line of \eqref{eq:comm-OO-q^2squared-tensors} can be obtained as in Proposition \ref{LEMMA:MOD-LAPLACIAN-KERR}.
\end{proof}

\begin{corollary}\label{lemma:commutator-triangle} 
The following commutation formula holds true for a scalar $\psi$:
\beaa
\, [|q|^2\lap,|q|^2 \square_\g]\psi&=& |q|^2\Big[O(a^2r^{-3}) \dk^{\leq 2}\psi+ \dk^2 \big( \Ga_g \c \dk \psi)+ \Ddot_3 \dk \big(|q|^2 \xi \c \Ddot_a \psi \big)\Big].
\eeaa
\end{corollary}

\begin{proof} 
Using the definition \eqref{definition-SS} to write that $|q|^2\lap =\OO+O(a^2r^{-2}) \dkb+r^{-1} \Ga_b \c \dk$, we deduce from Proposition \ref{LEMMA:MOD-LAPLACIAN-PERT-KERR},
\beaa
\, [|q|^2\lap, |q|^2\square_\g]\psi&=&[\OO, |q|^2\square_\g] +O(a^2r^{-2})[ \dkb, |q|^2\square_\g] \psi+ [r^{-1} \Ga_b \c \dk, |q|^2\square_\g]\psi\\
&=& |q|^2\Big[O(a^2r^{-3}) \dk^{\leq 2}\psi+ \dk^2 \big( \Ga_g \c \dk \psi)+ \Ddot_3 \dk \big(|q|^2 \xi \c \Ddot_a \psi \big)\Big],
\eeaa
as stated.
\end{proof}


\section{Approximate symmetry operators}\label{symmetry-operators-pert}


We can finally define in $\MM$ the following approximate symmetry operators as in Definition \ref{definition:symmetry-tensors}. 

 \begin{definition}\lab{definition:symmetry-tensors-pert} 
 We define the following second order differential operators,  acting on $\mathfrak{s}_2$ tensors,
\beaa
\SS_1= \nab_\T \nab_\T, \quad \SS_2= a \nab_\T \nab_\Z,\quad \SS_3= a^2 \nab_\Z \nab_\Z,\quad \SS_4=\OO.
\eeaa
\end{definition}

Recall that in Kerr we had also the alternative definition of $\SS_\aund$ give by \eqref{eq:SS-aund-Ddot-kerr}. In perturbations, they differ by error terms in the following way. 

\begin{lemma}\label{LEMMA:SYMM-OPERATORS}
Let
\beaa
\widetilde{\SS}_\aund\psi&=&|q|^2 \Ddot_\a(|q|^{-2}S_\aund^{\a\b} \Ddot_\b \psi), \qquad \textrm{ for }\aund=1,2,3, 4,
\eeaa
where $S_{\aund}^{\a\b}$ are given in Definition \ref{definition:tensors-S-pert:0}. Then, we have the following comparison with the approximate symmetry operators of Definition \ref{definition:symmetry-tensors-pert}:
\beaa
\widetilde{\SS}_1 &=& \SS_1  + \Ga_b\c \dk, \\
\widetilde{\SS}_2 &=& \SS_2 +r\Ga_b \c \dk ,\\
\widetilde{\SS}_3 &=& \SS_3 + r\Ga_b \c \dk,\\
\widetilde{\SS}_4 &=& \OO +r\Ga_b\c\dk.
 \eeaa
\end{lemma}

\begin{proof} 
See Section \ref{appendix-proof-symmetry-operators}.
\end{proof}

We now collect the formula for the commutators of the approximate symmetry operators $\SS_\aund$ for $\aund=1,2,3, 4$ with $\squared_2$.

\begin{proposition}\label{prop:commutators-SSaund-squared2-pert}
 The following  commutation  formulas hold true  for $\psi\in \sk_2$:
\beaa
\, [\SS_1, \squared_2] \psi &=& O(ar^{-2})\dk^{\leq 2}\psi+\dk^2 \big(\Ga_g \c \dk \psi\big)+\dk\big( \Ga_b \c \squared_2\psi  \big), \\
\, [\SS_2, \squared_2] \psi&=&O(ar^{-2})\dk^{\leq 2}\psi+\dk^2 \big(\Ga_g \c \dk \psi\big)+r\dk \big(\Ga_b \c \squared_2\psi\big),\\
\, [\SS_3, \squared_2] \psi&=& O(ar^{-2})\dk^{\leq 2}\psi+\dk^2 \big(\Ga_g \c \dk \psi\big)+r\dk \big(\Ga_b \c \squared_2\psi\big),
\eeaa
and
\beaa
\,  [\SS_4, |q|^2\squared_2] \psi&=& |q|^2\Big[ O( a r^{-2}) \  \dk^{\leq 2} \psi+  \dk^2 \big( \Ga_g \c \dk \psi)+ \Ddot_3 \dk \big(|q|^2 \xi \c \Ddot_a \psi \big)\Big].
\eeaa
\end{proposition}

\begin{proof} 
The first three relations are straightforward from \eqref{eq:commTZsquare-2tensors} in 
Proposition \ref{LE:COMMTZSQUARE}. 
 The commutator for $\SS_4=\OO$ is obtained from Proposition \ref{LEMMA:MOD-LAPLACIAN-PERT-KERR}.
\end{proof}


\section{Wave operator and Energy Momentum tensor}\label{section:energy-momentum}


Consider  variational  wave equations  for real-valued tensors  $\psi\in \sk_k$, of the form 
\bea\lab{eq:Gen.RW-general}
\squared_k \psi-V\psi=N,
\eea
where $V$ is a real potential. The variational wave equation \eqref{eq:Gen.RW-general} has Lagrangian
\beaa
 \LL[\psi]&=& \g^{\mu\nu}\Db_\mu  \psi\c\Db_\nu  \psi  +V\psi\c\psi,
 \eeaa
  where  the dot product   here denotes full contraction with respect to the  horizontal indices.

The  corresponding   energy-momentum tensor associated to \eqref{eq:Gen.RW-general} is given by 
 \bea\label{eq:definition-QQ-mu-nu}
 \QQ_{\mu\nu}:=\Db_\mu  \psi \c \Db _\nu \psi 
          -\frac 12 \g_{\mu\nu} \left(\Db_\la \psi\c\Db^\la \psi + V\psi \c \psi\right)= \Db_\mu  \psi \c \Db _\nu \psi -\frac 1 2\g_{\mu\nu} \LL[\psi].
 \eea

\begin{lemma}\label{lemma-divergence-QQ}
Given a solution $\psi \in \sk_k$ of equation  \eqref{eq:Gen.RW-general} we have
 \beaa
 \D^\nu\QQ_{\mu\nu}
  &=& \Db_\mu  \psi \c  \left(\squared_k \psi- V\psi\right)+ \Db^\nu  \psi ^A\Rdot_{ A   B   \nu\mu}\psi^B-\frac 1 2 \D_\mu V |\psi|^2.
 \eeaa
\end{lemma}

\begin{proof}
We have, making us of Proposition \ref{Proposition:commutehorizderivatives}
 \beaa
 \D^\nu\QQ_{\mu\nu}&=&\Db^\nu \Db_\nu  \psi \c \Db _\mu \psi
 +  \Db^\nu  \psi\c  \left( \Db_\nu \Db _\mu -\Db_\mu \Db _\nu \right)\psi -V\D_\mu \psi \c \psi      -\frac 1 2 \D_\mu V \psi\c\psi\\
  &=& \Db_\mu  \psi \c\Db^\nu \Db _\nu \psi+ \Db^\nu  \psi ^a\Rdot_{ ab   \nu\mu}\psi^b -V\D_\mu \psi \c \psi   -\frac 1 2 \D_\mu V \psi\c\psi\\
  &=& \Db_\mu \psi \c \left(\squared_k \psi- V\psi\right)+ \Db^\nu  \psi ^a\Rdot_{ ab   \nu\mu}\psi^b -\frac 1 2 \D_\mu V |\psi|^2
 \eeaa
with $|\psi|^2:=\psi\c\psi$.
\end{proof}


\subsubsection{Standard calculation for generalized currents}


We collect here some general calculations for generalized currents associated to equation \eqref{eq:Gen.RW-general}.

 \begin{proposition}\lab{prop-app:stadard-comp-Psi}
 Let   $\psi\in \mathfrak{s}_k$ be a solution of \eqref{eq:Gen.RW-general} and   $X$ be  a vectorfield. Then,
    \begin{enumerate}
\item 
 The $1$-form   $\PP_\mu=\QQ_{\mu\nu} X^\nu$   verifies
\beaa
\D^\mu \PP_\mu&=& \frac 1 2 \QQ  \c\piX+ X( \psi) \c  \left(\squared_k \psi- V\psi\right)-\frac 1 2 X(V) |\psi|^2 + X^\mu \Db^\nu  \psi ^a\Rdot_{ ab   \nu\mu}\psi^b.
\eeaa

\item
   Let $X$ as above,   $w$ a scalar   and $M$  a one form. Define
 \beaa
\PP_\mu[X, w, M]&:=&\QQ_{\mu\nu} X^\nu +\frac 1 2  w \psi \c \Db_\mu \psi -\frac 1 4|\psi|^2   \pr_\mu w +\frac 1 4 |\psi|^2 M_\mu.
  \eeaa
 Then,  
  \beaa
  \D^\mu  \PP_\mu[X, w, M] &=& \frac 1 2 \QQ  \c\piX - \frac 1 2 X( V ) |\psi|^2+\frac 12  w \LL[\psi] -\frac 1 4|\psi|^2   \square_\g  w \\
  &+&  X^\mu \Db^\nu  \psi ^a\Rdot_{ ab   \nu\mu}\psi^b+ \frac 1 4  \Div(|\psi|^2 M\big)\\
  &&+  \left(X( \psi )+\frac 1 2   w \psi\right)\c \left(\squared_k \psi- V\psi\right).
 \eeaa
 \end{enumerate}
\end{proposition}

\begin{proof}
 Let     $\PP_\mu[X, 0, 0]=\QQ_{\mu\nu} X^\nu$. Then,
\beaa
\D^\mu \PP_\mu[X, 0, 0]&=&\QQ_{\mu\nu} \D^\mu  X^\nu+X^\nu \D^\mu \QQ_{\mu\nu} \\
&=& \frac 1 2 \QQ  \c\piX + X^\mu \Db_\mu  \psi \c \left(\squared_k \psi- V\psi\right)+ X^\mu \Db^\nu  \psi ^a\Rdot_{ ab   \nu\mu}\psi^b-\frac 1 2 X^\mu \D_\mu V |\psi|^2
\eeaa
where we used Lemma \ref{lemma-divergence-QQ}. 
This proves the first part of the proposition. 
To prove the second part we write  
  \beaa
 \D^\mu  \PP_\mu[X, w, M] &=& \frac 1 2 \QQ  \c\piX + X( \psi )\c \left(\squared \psi- V\psi\right) 
-\frac 1 2  X( V ) |\psi|^2+  X^\mu \Db^\nu  \psi ^a\Rdot_{ ab   \nu\mu}\psi^b\\
&+&\frac 1 2 \D^\mu w  \,   \psi\c  \Db_\mu \psi  + \frac 1 2  w\, \Db^\mu  \psi \c \Db_\mu \psi\\
&+& \frac 1 2  w \psi \squared_k \psi              -\frac 1 2\psi \c \Db^\mu \psi    \pr_\mu w -\frac 1 4|\psi|^2   \square_\g  w+ \frac 1 4  \Div(|\psi|^2 M\big)\\
 &=& \frac 1 2 \QQ  \c\piX - \frac 1 2 X( V ) |\psi|^2 + \frac 1 2  w\, \Db^\mu  \psi \c \Db_\mu \psi+  \frac 1 2  w \psi\left( \squared \psi \right) \\
  &+&  X^\mu \Db^\nu  \psi ^a\R_{ ab   \nu\mu}\psi^b
- \frac 1 4|\psi|^2   \square_\g  w+ \frac 1 4  \Div(|\psi|^2 M\big)+ X( \psi )\c \left(\squared_k \psi- V\psi\right)
 \eeaa
which gives the desired result.
\end{proof}

We now specialize Proposition \ref{prop-app:stadard-comp-Psi} to the case of equation \eqref{eq:Gen.RW-general} in perturbations of Kerr. 

 \begin{proposition}\label{prop-app:stadard-comp-Psi-perturbations-Kerr}
 Let   $\psi\in \mathfrak{s}_k(\MM)$ be a solution of \eqref{eq:Gen.RW-general} and   $X$ be  a vectorfield of the form
 \beaa
    X=  X^3  e_3 +X^4 e_4.
    \eeaa
 Then, 
    \begin{enumerate}
\item 
 The $1$-form   $\PP_\mu=\QQ_{\mu\nu} X^\nu$   verifies
\beaa
\D^\mu \PP_\mu&=& \frac 1 2 \QQ  \c\piX+ X( \psi) \c N -\frac 1 2 X(V) |\psi|^2 - \big(\rhod +\etab\wedge\eta\big)\nab_{X^4e_4-X^3e_3}  \psi\c\dual\psi \\
&-& \frac{1}{2}\Im\Big(\tr\Xb H X^3 +\tr X\Hb X^4\Big)\c\nab\psi\c\dual\psi +r^{-2} \big(X^3\Ga_b+ X^4\Ga_g\big) \dk \psi \c \psi.
\eeaa
\item
   Let $X$ as above,   $w$ a scalar   and $M$  a one form. Define
 \beaa
\PP_\mu[X, w, M]&:=&\QQ_{\mu\nu} X^\nu +\frac 1 2  w \psi \c \Db_\mu \psi -\frac 1 4|\psi|^2   \pr_\mu w +\frac 1 4 |\psi|^2 M_\mu.
  \eeaa
 Then,  
  \beaa
  \D^\mu  \PP_\mu[X, w, M] &=& \frac 1 2 \QQ  \c\piX - \frac 1 2 X( V ) |\psi|^2+\frac 12  w \LL[\psi] -\frac 1 4|\psi|^2   \square_\g  w \\
  &-& \big(\rhod +\etab\wedge\eta\big)\nab_{X^4e_4-X^3e_3}  \psi\c\dual\psi \\
&-& \frac{1}{2}\Im\Big(\tr\Xb H X^3 +\tr X\Hb X^4\Big)\c\nab\psi\c\dual\psi \\
  &+& \frac 1 4  \Div(|\psi|^2 M\big)+  \left(X( \psi )+\frac 1 2   w \psi\right)\c N\\
  &+&r^{-2} \big(X^3\Ga_b+ X^4\Ga_g\big) \dk \psi \c \psi .
 \eeaa
 \end{enumerate}
\end{proposition}

\begin{proof} 
By Proposition \ref{prop-app:stadard-comp-Psi}, we only need to specialize the computation of the term $X^\mu  \Db^\nu  \psi ^A\Rdot_{ A   B   \nu \mu }\psi^B$ to the case of $X= X^3 e_3 +X^4 e_4$ and perturbations of Kerr. 
Since $\Rdot_{ ab   \nu\mu}$ is antisymmetric with respect to $(a, b)$, we have
\beaa
X^\mu \Db^\nu  \psi ^a\Rdot_{ ab   \nu\mu}\psi^b &=& \frac{1}{2}\in^{bc}X^\mu\Rdot_{bc   \nu\mu} \Db^\nu  \psi\c\dual\psi.
\eeaa
Introducing the spacetime 1-form
\bea
F_\nu &:=& \in^{bc}\Rdot_{bc   \nu\mu}X^\mu,
\eea 
we infer
\beaa
X^\mu \Db^\nu  \psi ^a\Rdot_{ ab   \nu\mu}\psi^b &=& \frac{1}{2} F_\mu\D^\mu\psi\c\dual\psi.
\eeaa
Next, we rewrite $F_\mu$ as
\beaa
F_\mu &=& \in^{bc}\Rdot_{bc\mu 4}X^4 +\in^{bc}\Rdot_{bc\mu 3}X^3
\eeaa
and we compute the various components of $F_\mu$. To this end, recall that we have the following decomposition of the curvature: 
 \beaa
 \R_{ab34}&=& 2 \in_{ab}\dual\rho,\\
 \R_{abc3}&=&\in_{ab}\dual \bb_c=r^{-1}\Ga_b,\\
 \R_{abc4}&=&-\in_{ab}\dual \b_c =r^{-1}\Ga_g.
 \eeaa
 Also, recalling Proposition \ref{proposition:componentsofB},
the components of $\B$   are given   by the formula
\beaa
\begin{split}
\B_{ a   b  c 3}&=     -  \trchb  \big( \de_{ca}\eta_b-  \de_{cb} \eta_a\big)  -  \atrchb \big( \in_{ca}  \eta_b -  \in_{cb}  \eta_a\big) \\
&+ 2 \big(- \chibh_{ca}  \eta_b + \chibh_{cb} \eta_a-  \chi_{ca} \xib_b+  \chi_{cb} \xib_a\big)\\
&=     -  \trchb  \big( \de_{ca}\eta_b-  \de_{cb} \eta_a\big)  -  \atrchb \big( \in_{ca}  \eta_b -  \in_{cb}  \eta_a\big) +r^{-1}\Ga_b,\\
\B_{ a   b  c 4}&=     -  \trch  \big( \de_{ca}\etab_b-  \de_{cb} \etab_a\big)  -  \atrch \big( \in_{ca}  \etab_b -  \in_{cb}  \etab_a\big) \\
&+ 2 \big(- \chih_{ca}  \etab_b + \chih_{cb} \etab_a-  \chib_{ca} \xi_b+  \chib_{cb} \xi_a\big)\\
&=     -  \trch  \big( \de_{ca}\etab_b-  \de_{cb} \etab_a\big)  -  \atrch \big( \in_{ca}  \etab_b -  \in_{cb}  \etab_a\big)+r^{-1}\Ga_g, \\
\B_{ a   b  3 4} &=4\big(-\xib_a \xi_b+\xi_a \xib_b-\eta_a \etab_b+\etab_a\eta_b\big)\\
&=4\big(-\eta_a \etab_b+\etab_a\eta_b\big)+\Ga_b\c\Ga_g.
\end{split}
\eeaa
Using the definition  \eqref{eq:DefineRdot} of $\Rdot$, we infer
 \beaa
 \Rdot_{ab34}&=& 2 \Big(\in_{ab}\dual\rho+\big(-\eta_a \etab_b+\etab_a\eta_b\big)\Big),\\
 \Rdot_{abc3}&=&    -  \frac{1}{2}\trchb  \big( \de_{ca}\eta_b-  \de_{cb} \eta_a\big)  -  \frac{1}{2}\atrchb \big( \in_{ca}  \eta_b -  \in_{cb}  \eta_a\big) +r^{-1}\Ga_b,\\
 \Rdot_{abc4}&=&   -  \frac{1}{2}\trch  \big( \de_{ca}\etab_b-  \de_{cb} \etab_a\big)  -  \frac{1}{2}\atrch \big( \in_{ca}  \etab_b -  \in_{cb}  \etab_a\big)+r^{-1}\Ga_g.
 \eeaa

We deduce from the definition for $F_\mu$ and the above identities for the components of $B_{ab\mu\nu}$
 \beaa
F_4 &=& \in^{bc}\Rdot_{bc43}X^3 \\
&=& \in^{bc}\Big(-2\in_{bc}\dual\rho -2(\etab_b\eta_c-\eta_b\etab_c)\Big)X^3\\
&=& -4\rhod X^3 -4(\etab\wedge\eta)X^3,
\eeaa
\beaa
F_3 &=& \in^{bc}\Rdot_{bc34}X^4 = 4\rhod X^4 +4(\etab\wedge\eta)X^4,
\eeaa
and
\beaa
F_e &=& \in^{bc}\Rdot_{bce3}X^3 +\in^{bc}\Rdot_{bce4}X^4\\
&=& \frac{1}{2}\in^{bc}\Big( -  \trchb  \big( \de_{eb}\eta_c -  \de_{ec} \eta_b\big)  -  \atrchb \big( \in_{eb}  \eta_c -  \in_{ec}  \eta_b\big)+r^{-1}\Ga_b\Big)X^3\\
&&+ \frac{1}{2}\in^{bc}\Big( -  \trch  \big( \de_{eb}\etab_c -  \de_{ec} \etab_b\big)  -  \atrch \big( \in_{eb}  \etab_c -  \in_{ec}  \etab_b\big)+r^{-1}\Ga_g\Big)X^4\\
&=& \Big( -  \trchb\dual\eta_e   + \atrchb\eta_e\Big)X^3 +\Big( -  \trch\dual\etab_e   + \atrch\etab_e\Big)X^4+r^{-1}(X^3\Ga_b+X^4\Ga_g),
\eeaa
i.e.
\beaa
F_4 &=& -4\rhod X^3 -4(\etab\wedge\eta)X^3,\\
F_3 &=& 4\rhod X^4 +4(\etab\wedge\eta)X^4,\\
F_e &=& \Big( -  \trchb\dual\eta_e   + \atrchb\eta_e\Big)X^3 +\Big( -  \trch\dual\etab_e   + \atrch\etab_e\Big)X^4 +r^{-1}(X^3\Ga_b+X^4\Ga_g).
\eeaa
Since we have 
\beaa
X^\mu \Db^\nu  \psi ^a\Rdot_{ ab   \nu\mu}\psi^b &=& \frac{1}{2}F_\mu\Db^\mu  \psi\c\dual\psi\\
&=& \frac{1}{2}F_4\Db^4\psi\c\dual\psi+\frac{1}{2}F_3\Db^3\psi\c\dual\psi+\frac{1}{2}F_b\Db^b  \psi\c\dual\psi\\
&=& -\frac{1}{4}F_4\nab_3\psi\c\dual\psi -\frac{1}{4}F_3\nab_4\psi\c\dual\psi +\frac{1}{2}F_b\nab^b\psi\c\dual\psi,
\eeaa
we infer
\beaa
&& X^\mu \Db^\nu  \psi ^a\Rdot_{ ab   \nu\mu}\psi^b\\
 &=& -\big(\rhod +\etab\wedge\eta\big)\nab_{X^4e_4-X^3e_3}  \psi\c\dual\psi \\
&&+\frac{1}{2}\Big(\big( -  \trchb\dual\eta   + \atrchb\eta\big)X^3 +\big( -  \trch\dual\etab   + \atrch\etab\big)X^4\Big)\c\nab\psi\c\dual\psi\\
&& +r^{-2} \big(X^3\Ga_b+ X^4\Ga_g\big) \dk \psi \c \psi
\eeaa
and hence
\beaa
X^\mu \Db^\nu  \psi ^a\Rdot_{ ab   \nu\mu}\psi^b &=& -\big(\rhod +\etab\wedge\eta\big)\nab_{X^4e_4-X^3e_3}  \psi\c\dual\psi \\
&&-\frac{1}{2}\Im\Big(\tr\Xb H X^3 +\tr X\Hb X^4\Big)\c\nab\psi\c\dual\psi +r^{-2} \big(X^3\Ga_b+ X^4\Ga_g\big) \dk \psi \c \psi
\eeaa
which concludes the proof of Proposition \ref{prop-app:stadard-comp-Psi}.
\end{proof}


\subsection{Decomposition of the wave operator in null frames}


\begin{lemma}\label{lemma:expression-wave-operator}
The wave operator for $\psi\in \sk_k$ is given by
\bea
\bsplit
\squared_k \psi&=-\frac 1 2 \big(\nab_3\nab_4\psi+\nab_4 \nab_3 \psi\big)+\left(\omb -\frac 1 2 \trchb\right) \nab_4\psi+\left(\om -\frac 1 2 \trch\right) \nab_3\psi \\
&+\lap_k \psi + (\eta+\etab) \c\nab \psi,
\end{split}
\eea
where $\lap=\nab^a \nab_a$ denotes the horizontal Laplacian for $k$-tensors. 
\end{lemma}

\begin{proof}
By definition
\beaa
\squared_k \psi&=& \g^{34} \Ddot_3\Ddot_4 \psi + \g^{43} \Ddot_4\Ddot_3 \psi+ \g^{cd}\Ddot_c\Ddot_d \psi.
\eeaa
We write, using \eqref{eq:expressions-Riccif-formula}, 
\beaa
\Ddot_4 \psi&=& \nab_4 \psi,\\
\Ddot_3 \Ddot_4 \psi&=& \nab_3 \nab_4 \psi - 2\omb \nab_4 \psi- 2\eta\c  \nab \psi, \\
\Ddot_4 \Ddot_3 \psi&=& \nab_4 \nab_3 \psi - 2\om \nab_3 \psi- 2\etab\c  \nab \psi, \\
\Ddot_d\psi&=&\nab_d\psi,\\
\Ddot_c \Ddot_d \psi&=&\nab_c\nab_d \psi - \frac 1 2 \chi_{cd}\nab_3 \psi -\frac 1 2 \chib_{cd}\nab_4 \psi. 
\eeaa
Hence
\beaa
\squared_k \psi&=&-\frac 1 2  \Ddot_3\Ddot_4 \psi    -\frac 1 2  \Ddot_4\Ddot_3 \psi     + \g^{cd}\Ddot_c\Ddot_d \psi\\
&=&-\frac 1 2 \big(\nab_3\nab_4\psi+\nab_4 \nab_3 \psi\big)+\g^{cd}\left(\nab_c\nab_d \psi - \frac 1 2 \chi_{cd}\nab_3\psi -\frac 1 2 \chib_{cd}\nab_4\psi \right)\\
& + &\omb \nab_4 \psi+\eta\c  \nab \psi +\omb\nab_3 \psi+\etab\c  \nab\psi\\
&=&-\frac 1 2 \big(\nab_3\nab_4\psi+\nab_4 \nab_3\psi\big)+\lap_2\psi- \frac 1 2 \trch\nab_3\psi -\frac1  2 \trchb \nab_4\psi\\
& + &\omb \nab_4 \psi+\eta\c  \nab\psi +\omb\nab_3 \psi+\etab\c  \nab \psi.
\eeaa
Hence
\beaa
\squared_k \psi&=& -\frac 1 2 \big(\nab_3\nab_4\psi+\nab_4 \nab_3 \psi\big)+\lap_k \psi +\left(\omb -\frac 1 2 \trchb\right) \nab_4\psi+\left(\om -\frac 1 2 \trch\right) \nab_3\psi \\
&&+(\eta+\etab) \c\nab \psi,
\eeaa
as stated.
\end{proof}

\begin{lemma}\label{lemma:expression-wave-operator-pert}
The wave operator for $\psi \in \sk_2(\CCC)$ is given by
\bea\label{first-equation-square}
\begin{split}
\squared_2 \psi&=-\nab_4 \nab_3 \psi  -\frac 1 2 \trchb \nab_4\psi+\left(2\om -\frac 1 2 \trch\right) \nab_3\psi+\lap_2 \psi+2\etab \c\nab \psi \\
&+ 2i \left( \rhod- \eta \wedge \etab \right) \psi+(\Ga_b \c \Ga_g) \c \psi.
\end{split}
\eea

Moreover, if $\psi \in \sk_2(\CCC)$ is $0$-conformally invariant, the above can be written as
\bea\label{first-equation-square-conformal-sign0}
\begin{split}
\squared_2 \psi&=-\nabc_4 \nabc_3 \psi  -\frac 1 2 \trchb \nabc_4\psi -\frac 1 2 \trch \nabc_3\psi+\lapc_2 \psi+2\etab \c\nabc \psi \\
&+ 2i \left( \rhod- \eta \wedge \etab \right) \psi+(\Ga_b \c \Ga_g) \c \psi.
\end{split}
\eea
\end{lemma}

\begin{proof} 
Using Lemma \ref{lemma:expression-wave-operator} for $\psi \in \sk_2(\CCC)$ and using  \eqref{correct-commutator-1} ,
we obtain
\beaa
\squared_2 \psi&=&-\nab_4 \nab_3 \psi -\frac 1 2 [\nab_3, \nab_4]\psi+\lap_2\psi +\left(\omb -\frac 1 2 \trchb\right) \nab_4\psi+\left(\om -\frac 1 2 \trch\right) \nab_3\psi \\
&&+(\eta+\etab) \c\nab \psi\\
&=&-\nab_4 \nab_3 \psi +\lap_2 \psi -\frac 1 2 \trchb \nab_4\psi+\left(2\om -\frac 1 2 \trch\right) \nab_3\psi+2\etab \c\nab \psi \\
&&+ 2i \left( \rhod-\eta \wedge \etab  \right) \psi+(\Ga_b \c \Ga_g) \c \psi,
\eeaa
as stated. The second relation for $\psi$ of signature $0$ is straightforward. 
\end{proof}


\subsection{Representation of the wave operator using $\That,$ $\Rhat$}


\begin{lemma}\label{LEMMA:SQUARED-K-THAT-RHAT-KERR}\label{lemma:squaredkpsi-perturbations}
We have\footnote{Observe that the expression
 in \eqref{eq:q2squaredkpsi-perturbations} is not regular at the horizon,  i.e. for $\De=0$.}  for $\psi\in \sk_k$, 
\bea\label{eq:q2squaredkpsi}\label{eq:q2squaredkpsi-perturbations}
\begin{split}
|q|^2 \squared_k \psi &=\frac{(r^2+a^2)^2}{\De} \big( -  \nab_\That \nab_\That \psi+   \nab_\Rhat \nab_\Rhat \psi \big) +2r \nab_\Rhat \psi\\
&+  |q|^2 \lap_k \psi   + |q|^2  (\eta+\etab) \c \nab \psi  + r^2 \Ga_g \c \dk \psi,
\end{split}
\eea
where $\lap_k$ denotes the horizontal Laplacian for $k$-tensors.
\end{lemma}

\begin{proof} 
See Section \ref{proof:alternative-repr-square}.
\end{proof}


\subsection{The wave operator using complex derivatives}


We now express the laplacian in terms of complex derivatives. We summarize the result in the following.
\begin{lemma}
\lab{COMPLEXWAVE-DECOMP}
We have for $\psi \in \sk_2(\CCC)$,
\bea\label{expression-DD-laplacian-0-K}
   \DD \hot (\DDb \c \psi)&=& 4\lap_2  \psi -2 i \left(\atrch\nab_3+\atrchb \nab_4\right) \psi  -8\Kh  \psi
\eea
where $\Kh$ is defined in \eqref{eq:definition-K}. In particular, in perturbations of Kerr we have
\bea
\DD\hot( \DDb \c \psi)&=& 4\lap_2  \psi -2 i \left(\atrch\nab_3+\atrchb \nab_4\right) \psi  \nn\\
&& + 2  \left( \trch\trchb+ \atrch\atrchb+4\rho\right) \psi +( \Ga_g \c \Ga_b) \c \psi,\label{eq:expression-DD-laplacian-real}\\
   \DD \hot (\DDb \c \psi)&=& 4\lap_2  \psi - 2i \left(\atrch\nab_3+\atrchb \nab_4\right) \psi \nn \\
   &&+ 2 \left(  \frac 1 2 \tr X \ov{\tr \Xb}+ \frac 1 2 \tr \Xb  \ov{\tr X} + 2P + 2\ov{P}\right) \psi  +( \Ga_g \c \Ga_b) \c \psi.\label{expression-DD-laplacian-0}
\eea
\end{lemma}

\begin{proof} 
See section \ref{sec:proof-complex-wave}.  
\end{proof}

We rewrite the  above using the conformal derivatives introduced in Lemma \ref{lemma:definition-conformal-derivatives}.

\begin{lemma}\label{LEMMA:CONFORMAL-DD-LAP} 
We have for $\psi \in \sk_2(\CCC)$ $s$-conformally invariant,
\bea
\begin{split}
\DDc\hot( \DDbc \c \psi)&=4\lapc_2  \psi -2 i \left(\atrch \nabc_3+\atrchb \nabc_4\right)\psi  \\
&+2\Big[    \left( \trch\trchb+ \atrch\atrchb+4\rho\right)  \\
& - i s \left(  \frac 1 2 \big(  \trch\atrchb-\trchb\atrch   \big)+2\dual \rho \right) \Big]\psi\\
& +( \Ga_g \c \Ga_b) \c \psi
   \end{split}
\eea
where $\lapc_2:=\ga^{ab}\nabc_a \nabc_b$ is the conformal Laplacian operator for horizontal $2$-tensors. 
\end{lemma}

\begin{proof} 
See section \ref{sec:proof-lemma-conformal-DD}. 
\end{proof}

By putting together the canonical expression for the wave operator given in Lemma \ref{lemma:expression-wave-operator} and the expression for the Laplacian given in Lemma \ref{COMPLEXWAVE-DECOMP}, we obtain the following. 

\begin{corollary}\label{corollary-wave-complex} 
We have, for $\psi\in \sk_2(\CCC)$,
\bea\label{wave-equation-Psi}
\begin{split}
\squared_2 \psi&=-\nab_4 \nab_3 \psi +\frac 1 4  \DD\hot( \DDb \c \psi)+\left(2\om -\frac 1 2 \tr X\right) \nab_3\psi- \frac 1 2 \tr\Xb \nab_4\psi+2\etab \c\nab \psi  \\
& +  \left( -\frac 1 2\trch\trchb- \frac 1 2 \atrch\atrchb-2\rho\right) \psi+ 2i \left( \rhod-\eta \wedge \etab  \right) \psi+(\Ga_b \c \Ga_g) \c \psi,
\end{split}
\eea
which can be rewritten as
\bea
\begin{split}
\squared_2 \psi&=-\nab_4 \nab_3 \psi +\frac 1 4  \DD\hot( \DDb \c \psi)+\left(2\om -\frac 1 2 \tr X\right) \nab_3\psi- \frac 1 2 \tr\Xb \nab_4\psi+2\etab \c\nab \psi \\
& +  \left( - \frac 1 4 \tr X \ov{\tr \Xb}- \frac 1 4 \tr \Xb  \ov{\tr X}  - 2\ov{P}\right) \psi- 2i \left(\eta \wedge \etab\right)  \psi+(\Ga_b \c \Ga_g) \c \psi.
\end{split}
\eea
\end{corollary}


\subsection{Commutators with the D'Alembertian}


We collect here some additional commutators with the horizontal laplacian, the operator $\OO$ and the D'Alembertian.

\begin{lemma}\label{LEMMA:COMMUTATOR-NAB3-NAB4-LAP} 
The following commutation formulas hold true for a 2-tensor $\psi \in \sk_2$:
\beaa
\, [\nab_3, |q|^2 \lap]\psi&=&(\eta-\ze) \c  |q|^2\nab_3 \nab \psi+(\eta-\ze) \c  |q|^2\nab\nab_3 \psi +\div(\eta-\ze)|q|^2 \nab_3 \psi \\
&&-\frac 1 2|q|^2 \left( \nab \trchb  \c \nab \psi+\nab \atrchb \c  \dual \nab \psi\right)+O(ar^{-4}) \dk^{\leq 1} \psi+\dk \big( \Ga_b \c \dk \psi), \\
\, [\nab_4, |q|^2 \lap]\psi&=&(\etab+\ze) \c |q|^2\nab_4 \nab\psi+(\etab+\ze) \c |q|^2\nab \nab_4 \psi +\div(\etab+\ze) |q|^2\nab_4 \psi \\
&&-\frac 1 2|q|^2 \left( \nab \trch  \c \nab \psi+\nab \atrch \c  \dual \nab \psi\right)+O(ar^{-4}) \dk^{\leq 1} \psi\\
&&+r^2\Ddot_3 \big( \xi \c \Ddot_a \psi \big)+ \dk \big( \Ga_g \c \dk \psi).
\eeaa
Similarly we have for a 2-tensor $\psi \in \sk_2$:
\beaa
\, [\nab_3, \OO]\psi&=& (\eta-\ze) \c  |q|^2\nab_3 \nab \psi+(\eta-\ze) \c  |q|^2\nab\nab_3 \psi \\
&&+\big( \div(\eta-\ze)+(\eta+\etab) \c (\eta - \ze) \big) |q|^2 \nab_3 \psi +O(ar^{-2}) \dk^{\leq 1} \psi +\dk \big( \Ga_b \c \dk \psi), \\
\, [\nab_4, \OO ]\psi&=& (\etab+\ze) \c  |q|^2\nab_4 \nab \psi+(\etab+\ze) \c  |q|^2\nab\nab_4 \psi \\
&&+\big( \div(\etab+\ze)+(\eta+\etab) \c (\etab + \ze) \big) |q|^2 \nab_4 \psi+O(ar^{-3}) \dk^{\leq 1} \psi \\
&&+r^2\Ddot_3 \big( \xi \c \Ddot_a \psi \big)+ \dk \big( \Ga_g \c \dk \psi).
\eeaa
\end{lemma}

\begin{proof} 
See section \ref{proof-lemma-comm-nab3-nab4-lap}. 
\end{proof}

\begin{lemma}\label{LEMMA:COMMUTATOR-NAB3-NAB4-SQUARE} 
The following commutation formulas hold true for a scalar $\psi$:
\beaa
\, [\nab_3, \square_\g]\psi&=&2 \om \nab_3\nab_3 \psi -(\trchb +2\omb) \nab_3\nab_4 \psi -\trchb\square_\g \psi \\
&&+r^{-2}  \dk^{\leq1} \psi +O(ar^{-2})  \nab_3\nab \psi   +r^{-2} \dk \big( \Ga_b \c \dk \psi), \\
\, [\nab_4, \square_\g]\psi&=& 2 \omb \nab_4\nab_4 \psi -(\trch +2\om) \nab_4\nab_3 \psi -\trch\square_\g \psi \\
&&+r^{-2}  \dk^{\leq1} \psi +O(ar^{-2})  \nab_4\nab \psi   +r^{-1} \dk \big( \Ga_g \c \dk \psi).
\eeaa
The following commutation formulas hold true for $\psi \in \sk_2$:
\beaa
[\nab_4, \squared_2]\psi&=&2\omb \nab_4 \nab_4 \psi -(\trch +2\om) \nab_4 \nab_3 \psi + 2 (\eta+\zeta) \c \nab_4 \nab \psi  -\trch \squared_2 \psi  \\
&& -\frac 1 4( \trch)^2  \nab_3\psi +O(r^{-3}) \dk^{\leq 1} \psi  +\Ddot_3 \big( \xi \c \Ddot_a \psi \big)+ r^{-2}\dk \big( \Ga_g \c \dk \psi).
\eeaa
In particular, the above can be written as
\beaa
[\nab_4, \squared_2]\psi&=& -(\trch +2\om) \nab_4 \nab_3 \psi  -\trch \squared_2 \psi -\frac 1 4( \trch)^2  \nab_3\psi \\
&&+O(r^{-3}) \dk^{\leq 1} \psi +O(r^{-4}) \dk^{\leq 2} \psi +\Ddot_3 \big( \xi \c \Ddot_a \psi \big)+ r^{-2}\dk \big( \Ga_g \c \dk \psi).
\eeaa
We can also deduce for $\psi \in \sk_2$:
\bea\label{eq:comm-rnab4-squared2}
\begin{split}
[r \nab_4, \squared_2]\psi &=-\nab_4\nab_4\psi -r\big(\frac 1 2 \trch -2\om\big) \squared_2 \psi -r\big(\frac 1 2 \trch +2\om\big) \lap_2 \psi  \\
&+O(r^{-2}) \dk^{\leq 1} \psi+O(r^{-3}) \dk^{\leq 2} \psi+r\Ddot_3 \big( \xi \c \Ddot_a \psi \big)+ r^{-1}\dk \big( \Ga_g \c \dk \psi).
\end{split}
\eea
\end{lemma}

\begin{proof} 
See section \ref{proof:lemma-comm-nab3-nab4-square}.
\end{proof}

\begin{lemma}\label{LEMMA:COMMUTATOR-NAB-RHAT-SQUARE} 
The following commutation formula holds true for a scalar $\psi$:
  \beaa
   \, [\nab_\Rhat, |q|^2\square_\g]\psi&=&O(r)\square_\g \psi+O(r)\lap\psi +O(ar^{-1})\dk^{\le 2}  \psi+ O(1) \nab_\Rhat \psi +O(ar^{-1}) \nab \psi\\
 &&+ r \dk \big( \Ga_b \c \dk \psi\big).
 \eeaa

We also have for a scalar $\psi$:
\bea\label{eq:commutator-nab-Rhat-lap}
[\nab_\Rhat, |q|^2 \lap] \psi&=&O(ar^{-3}) \dkb \psi+ r \dk \big( \Ga_g \c \dk \psi \big).
\eea
In particular,
\beaa
[\nab_\Rhat, \lap]\psi&=&O(r^{-3} \De) \lap \psi +O(ar^{-5}) \dkb \psi+ r^{-1} \dk \big( \Ga_g \c \dk \psi \big).
\eeaa
\end{lemma} 

\begin{proof} 
See Section \ref{proof:lemma-comm-nab-Rhat-square}. 
\end{proof}

\begin{lemma}\lab{LEMMA:COMMUTATIONOFHODGEELLIPTICORDER1WITHSQAURED2FDILUHS} 
The following commutation formulas hold true for $\psi \in \sk_2$:
\beaa
|q|\DDd_2\squared_2 - \squared_1|q|\DDd_2\psi&=& 3 \Kh |q|\DDd_2\psi -\frac{2a\cos\th}{|q|}\dual \DDd_2 \Lieb_\T\psi  -|q|( \eta+\etab) \c \squared_2\psi   \\
&&+ O(ar^{-2})\dk^{\leq 1}  \psi +O(ar^{-3})\dk^{\leq 2}\psi\\
&&+ \dk^{\leq 2} (\Ga_g \c \psi)+r \Hc \c \squared_2 \psi + \Ddot_3 (r\xi \c \nab_3 \psi),
\eeaa
and similarly for $|q|\DDs_1$. In particular,
\beaa
|q|\DDd_2\squared_2\psi - \squared_1|q|\DDd_2\psi &=& \frac{3}{r^2} |q|\DDd_2\psi +O(ar^{-2} )\dk^{\leq 2}\psi+ \dk^{\leq 2} (\Ga_g \c \psi)\\
&&+r \Hc \c \squared_2 \psi + \Ddot_3 (r\xi \c \nab_3 \psi ), \\
|q|\DDs_1\,\squared_0 \psi - \squared_1|q|\DDs_1\,\psi &=&- \frac{1}{r^2}|q|\DDs_1\,\psi+O(ar^{-2})\dk^{\leq 2}\psi+\dk^{\leq 2}(\Ga_g\c\psi)\\
&&+r \Hc \c \squared_2 \psi + \Ddot_3 (r\xi \c \nab_3 \psi ).
\eeaa

The following commutation formulas hold true for $\psi \in \sk_1$:
\beaa
|q|\DDs_2 \, \squared_1\psi - \squared_2|q|\DDs_2\,\psi &=&- \frac{3}{r^2}|q|\DDs_2\,\psi+O(ar^{-2})\dk^{\leq 2}\psi+\dk^{\leq 2}(\Ga_g\c\psi)\\
&&+r \Hc \c \squared_1 \psi + \Ddot_3 (r\xi \c \nab_3 \psi ), \\
|q|\DD_1\squared_1\psi - \squared_0|q|\DD_1 \psi &=& \frac{1}{r^2}|q|\DD_1\psi+O(ar^{-2})\dk^{\leq 2}\psi+\dk^{\leq 2}(\Ga_g\c\psi)\\
&&+r \Hc \c \squared_1 \psi + \Ddot_3 (r\xi \c \nab_3 \psi ).
\eeaa
\end{lemma}

\begin{proof} 
See section \ref{proof:lemma:commofHodgewithsquare}.
\end{proof}


\section{Spacetime elliptic identities}
\label{Section:SpacetimeellipticI}


We adapt the spacetime elliptic identities in Proposition \ref{prop:2D-hodge-non-integrable-pert-Kerr-div}  to the case of perturbations of Kerr.

\begin{lemma}\label{prop:2D-hodge-non-integrable-pert-Kerr} 
Given a not necessarily integrable horizontal structure, the following pointwise relations hold:

{\bf i.)}\quad The following identity holds for $f\in \sk_1$:
\bea
\label{eq:hodgeident1-nonint-pert-Kerr}
\begin{split}
|\nab   f |^2+ \Kh |f|^2&=|\DDd_1   f   |^2+  \frac{2a\cos\th}{|q|^2}\dual  \nab_\T  f  \c f \\
&+\D_\a \big( \nab^\a f \c f- (\div f) f^{\a}-(\curl f) (\dual f )^\a \big)+\Ga_g \cdot \dk f \cdot f,
\end{split}
\eea
and
\beaa
\begin{split}
|\nab   f |^2+ \Kh |f|^2&=|\DDd_1   f   |^2+\frac{2a\cos\th (r^2+a^2)}{|q|^4}\dual \nab_\That f  \c f \\
&+ \nab_a \Big( \nab^a f \c f-             (\div f) f^{a}-(\curl f) (\dual f )^a \Big)+\Ga_g \cdot \dk f \cdot f.
\end{split}
\eeaa

{\bf ii.)}\quad The following identity holds for $f\in \sk_2$:
\bea
\label{eq:hodgeident2-nonint-pert-Kerr}
\begin{split}
|\nab    f  |^2+ 2 \Kh |f|^2 &=2 |\DDd_2   f  |^2+  \frac{2a\cos\th}{|q|^2}\dual  \nab_\T  f  \c f\\
&+\D_\a \big( \nab^\a f \c f-2 (\div f)_\b f^{\a\b} \big)+\Ga_g \cdot \dk f \cdot f,
\end{split}
\eea
and
\beaa
\begin{split}
|\nab    f  |^2+ 2 \Kh |f|^2 &=2 |\DDd_2   f  |^2+\frac{2a\cos\th (r^2+a^2)}{|q|^4}\dual \nab_\That f  \c f\\
&+\nab_a \big( \nab^a f \c f-2 (\div f)_b  f^{ab} \big)+\Ga_g \cdot \dk f \cdot f.
\end{split}
\eeaa

{\bf iii.)}\quad  The following identity holds for $f\in\sk_1$:
\bea
\label{eq:hodgeident3-nonint-pert-Kerr}
\begin{split}
 |\nab  f   |^2-\Kh |f|^2&=2|\DDs_2   f   |^2-  \frac{2a\cos\th}{|q|^2}\dual  \nab_\T  f  \c f\\
 &+\D_\a \big( \nab^\a f \c f+2 (\DDs_2 f)^{\a\b} f_{\b} \big)+\Ga_g \cdot \dk f \cdot f,
 \end{split}
\eea
and
\beaa
\begin{split}
 |\nab  f   |^2-\Kh |f|^2&=2|\DDs_2   f   |^2-\frac{2a\cos\th (r^2+a^2)}{|q|^4}\dual \nab_\That f  \c f
 \\
 &+\nab_a \Big( \nab^a f \c f+2 (\DDs_2 f)^{ab} f_{b} \Big)+\Ga_g \cdot \dk f \cdot f.
 \end{split}
\eeaa
\end{lemma}

\begin{proof} 
Straightforward application of Proposition \ref{prop:2D-hodge-non-integrable-pert-Kerr-div} and \eqref{eq:atrch-e3-atrchb-e4-pert-kerr}.
\end{proof}

We collect here some elliptic identities involving the operator $\OO$ and $\lap_2$.
\begin{lemma} 
We have for $\psi \in \sk_2$, 
\bea
\bsplit
\lap_2 \psi \c \psi &= - |\nab \psi|^2- ((\eta+ \etab) \c  \nab  \psi) \c \psi+ \D^\a( \nab_\a \psi \c \psi), \\
\OO(\psi) \c \psi&=  -  |q|^2|\nab \psi|^2- |q|^2((\eta+ \etab) \c  \nab \psi) \c \psi+  \D^\a(|q|^2 \nab_\a \psi \c \psi)+  \Ga_b \c \dk\psi \c \psi . \label{eq:Ou-u}
\end{split}
\eea
We have for $\psi \in \sk_2(\CCC)$,
 \bea
 \bsplit
 \lap_2 \psi \c \ov{\psi} &= - |\nab \psi |^2 - ((H+\Hb) \c  \nab \psi) \c \ov{\psi}+ \D^\a( \nab_\a \psi \c \ov{\psi}),\label{eq:lap-nab}\\
  \OO( \psi )\c \ov{\psi} &= - |q|^2|\nab \psi |^2 - |q|^2((H+\Hb) \c  \nab\psi ) \c \ov{\psi}+ \D^\a(|q|^2 \nab_\a \psi \c \ov{\psi})+  \Ga_b \c \dk\psi \c \psi .
 \end{split}
 \eea
\end{lemma}

\begin{proof} 
Using Lemma \ref{lemma:divergence-spacetime-horizontal-ch1}, we obtain
\beaa
\lap_2 \psi \c \psi&=& \nab^a \nab_a \psi \c \psi= \nab^a( \nab_a \psi \c \psi) - |\nab \psi |^2 \\
&=& \D^\a( \nab_\a \psi \c \psi)- ((\eta+ \etab) \c  \nab  \psi) \c \psi - |\nab \psi|^2 
\eeaa
which gives the first identity. We then deduce, using $\nab(|q|^2)=(\eta+\etab) |q|^2+ r \Ga_b$,
\beaa
\OO(\psi) \c \psi&=& |q|^2 \left( \lap_2 \psi+(\eta+ \etab) \c  \nab \psi \right)  \c \psi \\
&=& -  |q|^2|\nab \psi|^2+ |q|^2 \D^\a( \nab_\a \psi \c \psi)\\
&=& -  |q|^2|\nab \psi|^2- |q|^2((\eta+ \etab) \c  \nab  \psi) \c \psi+  \D^\a(|q|^2 \nab_\a \psi \c \psi )+  \Ga_b \c \dk\psi \c \psi 
\eeaa
as stated.
\end{proof}


\chapter{Derivation of the main equations}\label{CHAPTER-DERIVATION-MAIN-EQS}



\section{Teukolsky equation for $A$}\lab{section:teukolsky-equations}


It is known that the curvature components $A$ and $\Ab$ satisfy  wave equations  which   decouple from  all other components at the linear level, the celebrated Teukolsky equations. 
In this section we derive, using our formalism,  the corresponding   Teukolsky equation  for $A$ while keeping track of  the error terms  generated  by the perturbation from Kerr expressed in terms of $(\Ga_b, \Ga_g)$.


\subsection{The Teukolsky equation for $A$}\label{section:teukolsky-equation}


 \begin{proposition}\label{TEUKOLSKY-PROPOSITION} 
 The complex tensor $A \in \sk_2(\CCC)$ satisfies the following equation:
 \bea\label{Teukolsky-equation-tens}
 \LL(A)&=& \err[\LL(A)]
 \eea
 where
\bea\label{Teukolsky-operator-ch5}
\begin{split}
\LL(A) &=-\nabc_4\nabc_3A+ \frac{1}{4}\DDc\hot (\DDbc \c A)+\left(- \frac 1 2 \tr X -2\ov{\tr X} \right)\nabc_3A\\
&-\frac{1}{2}\tr\Xb \nabc_4A+\left( 4H+\Hb +\ov{\Hb} \right)\c \nabc A+ \left(-\ov{\tr X} \tr \Xb +2\ov{P}\right) A+  H   \hot (\ov{\Hb} \c A),
\end{split}
\eea
with error term  expressed schematically
\bea
\err[\LL(A)]&=&r^{-1}  \dk^{\leq 1}\big( \Ga_g \c  B\big) + \nabc_3\Xi  \c B +  \Ga_b \c \Ga_g \c A.
\eea
\end{proposition}

\begin{proof} 
See section \ref{section:proof-teukolsky-eq}.
\end{proof}

For completeness, we collect here the real and imaginary part of the Teukolsky operator $\LL(A)$ defined in \eqref{Teukolsky-operator-ch5}. 

\begin{corollary}\label{COROLLARY-REAL-TEUK} 
The Teukolsky operator in \eqref{Teukolsky-operator-ch5} can also be written as
\beaa
\LL(A) &=&\Re(\LL)(A) + i \Im(\LL)(A)
\eeaa
where
\beaa
\Re(\LL)(A)&=& -\nabc_4\nabc_3A+\lapc_2  A +(4\eta+2\etab) \c \nabc A -\frac{1}{2}\trchb  \nabc_4A-\frac 5 2 \trch \nabc_3 A\\
&& +   \left( -\frac 1 2\trch\trchb-\frac 1 2 \atrch\atrchb+4\rho+ 4 \eta \c \etab \right) A, \\
\Im(\LL)(A)&=& -2\atrch \nabc_3A+  4  \dual \eta \c \nabc A\\
&&+ \left( \frac 1 2  \trch \atrchb- \frac 1 2 \trchb\atrch-4\dual \rho- 4 \eta \wedge \etab \right) A.
\eeaa
\end{corollary}

\begin{proof} 
See section \ref{proof-corollary-teukolsky}.
\end{proof}


\subsection{Connection to the classical Teukolsky equation}
\lab{section:classicalTeuk-eq}


The Teukolsky equation \eqref{Teukolsky-equation-tens} is a tensorial equation for  $A\in \sk_2(\mathbb{C})$, as defined in our formalism.   The standard derivation  of the equation, in linear theory,  is done instead with respect to  the Newman-Penrose  formalism, see section \ref{section:NPformalism}.  To  relate the Teukolsky equation in our formalism to the classical one  in NP formalism we have to project it with respect to the standard  horizontal frame $e_1, e_2$ of Kerr, see
 \eqref{eq:canonicalHorizBasisKerr}, 
\beaa
e_1=\frac{1}{|q|}\pr_\th,  \qquad e_2=\frac{a\sin\th}{|q|}\pr_t+\frac{1}{|q|\sin\th}\pr_\phi,
\eeaa
 for which  the relations  \eqref{eq:expressions-nab41} are  verified.

 One can check in fact  that  the standard  Teukolsky variable, which we denote by 
$\a^{[+2]}$, is related  to  our $A$ via the formula 
\bea
\lab{relation-a-A}
\a^{[+2]}&:=& -\frac{\ov{q}}{ q} A_{11}, \qquad  A_{11}:=A(e_1, e_1).
\eea
 Indeed, in the physics literature, see \cite{Chand1}, the curvature component $\a^{[+2]}$ is a complex scalar defined in Newman-Penrose formalism as 
\beaa
\a^{[+2]}&=& -W(l, m, l, m)
\eeaa
where $l$ and $m$ are related to  our  (ingoing) frame  $e_4$, $e_3$, $e_1, e_3$   by 
\beaa
l = e_4, \qquad m= \frac{|q|}{\sqrt{2} q}\left(e_1 + ie_2\right), 
\eeaa
 and $W$ coincides with  the Riemann curvature tensor $\R$   for vacuum spacetimes.       We therefore deduce
\beaa
\a^{[+2]}&=& -\frac{|q|^2}{2 q^2}W(e_4, \left(e_1 + ie_2\right), e_4, \left(e_1 + ie_2\right))= -\frac{\ov{q}}{2 q}\left(W_{4141}+i W_{4142}+i W_{4241}-W_{4242}\right)\\
&=& -\frac{\ov{q}}{2 q}\left(W_{4141}-W_{4242}+2i W_{4142}\right)= -\frac{\ov{q}}{ q}\left(W_{4141}+i W_{4142}\right).
\eeaa
On the other hand
\beaa
A_{11}&=& W_{4141}+ i W_{4142}
\eeaa
and therefore  $\a^{[+2]}= -\frac{\ov{q}}{ q}A_{11}$.

  One can   check that
in the particular case of the Kerr metric  we have   
\bea\label{eq:standard-Teukolsky-equation-literature}
\begin{split}
|q|^2\square_{{m,a}} \a^{[+2]}&=-4(r-m) \pr_r \a^{[+2]}-4 \left(\frac{m(r^2-a^2)}{\Delta} - r -i a \cos\th \right)\pr_t \a^{[+2]} \\
&- 4\left(\frac{a(r-m)}{\Delta} + i \frac{\cos\th}{\sin^2\th} \right) \pr_\vphi \a^{[+2]}+ (4 \cot^2\th-2)\a^{[+2]}, 
\end{split}
\eea
 where $\square_{{m,a}}  $ is the D'Alembertian  relative to the  Kerr metric.  This  is  the standard form of the Teukolsky equation in Boyer-Lindquist coordinates, see \cite{Teuk}.


\section{Generalized Regge-Wheeler equation for $\qf$}
\label{section:gRW-equation}


In this section we derive the generalized Regge-Wheeler-type equation.


\subsection{The invariant quantities $Q$ and $\qf$}


We start with the following lemma.
\begin{lemma}\label{lemma:basicinvarianceandconformalinvarianceofQofA} 
Let $C_1$ and $C_2$ be scalar functions. The expression 
\bea\lab{definition-Q(A)}
 Q(A)&=& \nabc_3\nabc_3 A + C_1  \nabc_3A + C_2  A \in \sk_2(\CCC)
 \eea
is  $0$-conformally invariant provided $C_1$ is $-1$-conformally invariant and $C_2$ is $-2$-conformally invariant. 
\end{lemma}

\begin{proof} 
Direct verification in view of the definition of  the conformal derivative $\nabc_3$.  
\end{proof}

\begin{definition} 
\lab{Definition:Define-qf}
Given a  fixed  null  pair  $(e_3, e_4)$  and  scalar functions $r$ and $\th$ as in Section \ref{sec:setupandlinearizedquantities}, we define our main quantity $\qf \in \sk_2(\CCC)$ as 
\bea
\lab{eq:definition-qf}
\qf&=& q \ov{q}^{3} Q(A)=q \ov{q}^{3} \left( \nabc_3\nabc_3 A + C_1  \nabc_3A + C_2   A\right)
\eea
where $q=r+ i a \cos\th$, and  the scalar function $C_1$, $C_2$ are given by
\bea\lab{eq:C1-C2-comparison-Ma}
\begin{split}
C_1&=2\trchb - 2\frac {\atrchb^2}{ \trchb}  -4 i \atrchb, \\
C_2  &= \frac 1 2 \trchb^2- 4\atrchb^2+\frac 3 2 \frac{\atrchb^4}{\trchb^2} +  i \left(-2\trchb\atrchb +4\frac{\atrchb^3}{\trchb}\right).
\end{split}
\eea
\end{definition}

\begin{remark} 
Note that   $\qf$  is  independent of the particular normalization. 
More precisely if $ e_3' =\la^{-1} e_3, e_4'=\la e_4 $ and $A'=\la^2 A$ then $\qf'=\qf$.
\end{remark}

In perturbations of Kerr, the quantity $\qf$ defined above can be factorized as follows\footnote{In practice, we will rely on a more precise factorization, see Lemma \ref{lemma:factorizationofqfusefulfortransportequations}.}. 
\begin{proposition}\lab{PROP:FACTORIZATION-QF}
In perturbations of Kerr,  the quantity $\qf$  defined in \eqref{eq:definition-qf} with $C_1, C_2$ given by \eqref{eq:C1-C2-comparison-Ma} can be factorized as
\bea
r \nabc_3\left( r^2   \left(   \nabc_3  \left( r \frac{\ov{q}^4}{r^4}  A\right)\right)\right)&=& \frac{\ov{q} }{q } \qf+  r^4 \Ga_b \c \nab_3^{\leq 1} A,
\eea
or also as 
\bea
 r^2 \nabc_3\nabc_3\left(\frac{\ov{q}^4}{r^2} A \right) &=& \frac{\ov{q} }{q } \qf+  r^4 \Ga_b \c \nab_3^{\leq 1} A.
\eea
In particular, in Kerr we have 
\bea
r \nabc_3\left( r^2   \left(   \nabc_3  \left( r \frac{\ov{q}^4}{r^4}  A\right)\right)\right)&=& \frac{\ov{q} }{q } \qf, \lab{eq:Prop-factorization-qf} \\
r^2\nabc_3\nabc_3\left(\frac{\ov{q}^4}{r^2} A \right) &=&  \frac{\ov{q} }{q }  \qf. \lab{eq:secondfactorization-qf}
\eea
\end{proposition}

\begin{proof} 
See section \ref{proof:prop-factorization-qf}. 
\end{proof}


\subsection{Comparison of  $\qf$ with Ma's quantity}
\label{section:projection-gRW}


The quantity in    Ma,  \cite{Ma},     which corresponds to our  $ \qf $,   is chosen to be
 a complex scalar   $\a^{[+2]}$,  obtained by a    Chadrasekhar  type transformation from  what Ma denotes as  $\a^{[+2]}$, i.e. the  complex scalar 
  verifying  the standard Teukolsky equation \eqref{eq:standard-Teukolsky-equation-literature} as derived in \cite{Teuk} by using the NP formalism. Thus
  $\a^{[+2]}$  relates to  our $A$ according to \eqref{relation-a-A}, i.e. 
  \beaa
\a^{[+2]}= -\frac{\ov{q}}{ q}A_{11}, \qquad A_{11}= A(e_1, e_2).
\eeaa
\begin{proposition}
\lab{Prop:Projectionpqf-Ma}
We have
\bea
\lab{eq:Projectionpqf-Ma}
\qf_{11} &=&-r \nab_3 \left( r^2 \nab_3 \left( r \frac{|q|^4}{r^4}  \a^{[+2]}      \right)\right).
\eea
\end{proposition}

\begin{proof}
The proof is based on the following lemma.
\begin{lemma}
Let $U$ be an anti-selfadjoint  complex $2$ tensor $U\in \sk_2(\CCC)$ in  Kerr. Then,  setting $p= q \ov{q}^{-1}$ we   have
\beaa
(\nab_3 U)_{11}=  p^{-1}  e_3( p U_{11}).
\eeaa
\end{lemma}

\begin{proof}
We  calculate, using the relations \eqref{eq:expressions-nab41} and $U_{12}=U_{21}=( \dual U)_{11}=-i U_{11}$,
\beaa
\nab_3 U_{11}&=& e_3 (U_{11} )  -     2U_{\nab_3 1 1}=  e_3 U_{11} -\atrchb U_{12}= 
 e_3( U_{11})+ i \atrchb U_{11}.
 \eeaa
 Note that\footnote{We have $e_3(p)= e_3( q\ov{q}^{-1} )= -\ov{q}^{-1}+ q \ov{q}^{-2}=  \ov{q}^{-2}\big(-\ov{q}+q \big)=    2 ai \cos\th \ov{q}^{-2}    = p i \atrchb $.}
 $ e_3 p =  i\atrchb p$. Hence
 \beaa
 \nab_3 U_{11}&=& e_3( U_{11})+ e_3(p)  U_{11} =  p^{-1}  e_3(p U_{11} )
 \eeaa
 as stated.
  \end{proof}
  
 According to Proposition \ref{PROP:FACTORIZATION-QF}  we have, in the ingoing normalization, 
 \beaa
 \frac{\ov{q} }{q } \qf= r \nab_3\left( r^2   \left(   \nab_3  \left( r \frac{\ov{q}^4}{r^4}  A\right)\right)\right).
 \eeaa
 Therefore, using the above lemma, 
 \beaa
  \frac{\ov{q} }{q } \qf_{11} &=&  r ( p^{-1} \nab_3  p )   \left( r^2   \left(    \big( p^{-1}  \nabc_3  p\big)   \left( r \frac{\ov{q}^4}{r^4}  A_{11}\right)\right)\right)= p^{-1}  r \nab_3 \left( r^2 \nab_3 \left( r p  \frac{\ov{q}^4}{r^4}  A_{11}\right)\right).
 \eeaa
We deduce, using  $A_{11}=-p  \a^{[+2]} $,
\beaa
\qf_{11} &=& r \nab_3 \left( r^2 \nab_3 \left( r \frac{q}{\ov{q}}  \frac{\ov{q}^4}{r^4}   A_{11} \right)\right)=-r \nab_3 \left( r^2 \nab_3 \left( r \frac{q}{\ov{q}}  \frac{\ov{q}^4}{r^4}  \frac{q}{\ov{q} } \a^{[+2]}      \right)\right)\\
&=& -r \nab_3 \left( r^2 \nab_3 \left( r \frac{|q|^4}{r^4}  \a^{[+2]}      \right)\right)
\eeaa
as stated.
\end{proof}

\begin{remark}
Note that all quantities  in Proposition \ref{Prop:Projectionpqf-Ma} are defined with respect to the  ingoing normalization.  If we replace  $ \a^{[+2]} $ with
  $ \a_{(out)}^{[+2]}= \frac{\De^2}{|q|^4}  \a^{[+2]}  $,   corresponding to  the outgoing  normalization   $e_4^{(out)} =\frac{|q|^2}{\De} e_4^{(in)} $,  then
  \eqref{eq:Projectionpqf-Ma} becomes 
  \beaa
\lab{eq:Ma'smainquantity}
\qf_{11} &=&-r \nab_3 \left( r^2 \nab_3 \left( r \frac{\De^2}{r^4}   \a_{(out)}^{[+2]}       \right)\right).
\eeaa
Ma denotes    $   \phi_{+2}^0 = \frac{\De^2}{r^4}   \a_{(out)}^{[+2]}  $  and defines
  \bea\label{eq:definition-phi2-Ma}
\phi_{+2}^2&=& (r Yr)(r Y r)(\phi_{+2}^0)
\eea
 where $Y$ is precisely $e_3=\frac{r^2+a^2}{\De} \partial_t + \frac{a}{\Delta} \partial_\phi - \partial_r$ in BL coordinates. Thus $\phi_{+2}^2=-\qf_{11}$.
\end{remark} 

In \cite{Ma}, Ma   states the following proposition.
\begin{proposition}[Equation (24.c) in \cite{Ma}]
The quantity  
\beaa
\phi^2_{+2}:=r \nab_3 \Big( r^2 \nab_3 \big( r \phi^0_{+2}     \big)\Big)
\eeaa  
verifies  the following equation in Kerr
\bea\label{eq:Ma-equation-1}
\nn&&|q|^2 \square_{a, m} \phi_{+2}^2+ 4 i  \left(\frac{\cos\th}{\sin^2\th} \partial_\phi - a \cos\th \partial_t \right)\phi_{+2}^2-4 \left( \cot^2\theta +\frac{r^2-2mr+2a^2}{r^2} \right)\phi_{+2}^2\\
&=& -8 (a^2 \partial_t + a \partial_\phi) \phi_{+2}^1 -12 a^2 \phi_{+2}^0,
\eea
 where $\square_{{a, m}}  $ is the  standard D'Alembertian  relative to the  Kerr metric and 
 \beaa
 \phi_{+2}^0=\frac{\De^2 \a_{(out)}^{[+2]}}{r^4}, \qquad \phi_{+2}^1= r e_3 r  (\phi_{+2}^0).
 \eeaa
\end{proposition}


\subsection{The derivation of the gRW  equation  for $\qf$}
\lab{subsection:DerivationgRW-eq}


We now state the first main result of  Part I   concerning  the wave equation satisfied by $\qf$.  To start with we note that that  we have decided to  use  as definition of $\qf$ the more complicated expression  in \eqref{eq:definition-qf} rather than  the more direct formula  appearing on the left hand side  of \eqref{eq:Prop-factorization-qf}, which is directly comparable  (by projection) with the  quantity  $\phi^2_{+2}$ of Ma,
    due to the fact that the second formula generates   un-acceptable  error terms.  
  
 \begin{theorem}\label{MAIN-THEOREM-PART1} 
 The invariant symmetric traceless $2$-tensor $\qf \in \sk_2(\CCC)$ in Definition \ref{Definition:Define-qf} satisfies the equation
 \bea\label{wave-equation-qf}
 \squared_2 \qf   -i \frac{4 a\cos\th}{|q|^2} \nab_\T \qf   - V  \qf &=&   L_{\qf}[A] + \err[\squared_2 \qf],
 \eea
 where:
 \begin{itemize}
  \item $\T$ is the vectorfield given by Definition \ref{Definition:vfsTZ}, see also Remark \ref{rmk:inabT-term} below.

 \item The potential $V$ is the \textbf{real} scalar function given by
 \bea
 V&=& \frac{4}{|q|^2}\frac{r^2-2mr+2a^2}{r^2}-\frac{4a^2\cos^2\th}{|q|^6}( r^2+6mr+a^2\cos^2\th ),
 \eea
  which for $a=0$ coincides with the potential of the Regge-Wheeler equation  in Schwarzschild, i.e. $V=-\trch\trchb+ O(\frac{|a|}{r^4})$, see also Remark \ref{rmk:propertiesofthepotentialV} below.
 
 \item $ L_{\qf}[A ]$ is  a  linear second order   operator  in $A$,   given in the outgoing frame by
  \beaa
 L_\qf[A]&=& q\ov{q}^3    \Bigg(- \frac{8a^2 \De}{r^2|q|^4}\nab_\T \nab_3 A -\frac{8a  \De }{r^2|q|^4} \nab_\Z \nab_3A\\
&&+ W_4 \  \nab_4A+ W_3\nab_3A+W\c\nab A +W_0  A\Bigg),
 \eeaa
where $W_4$, $W_3$, $W_0$ are complex functions  of $(r, \th)$ and $W$ is the product of a complex function  of $(r, \th)$ with $\dual\Re(\Jk)$, with the following fall-off in $r$ 
\beaa
q\ov{q}^3 W_4= q\ov{q}^3 W_3=q\ov{q}^3 W=O\left(a\right), \qquad  q\ov{q}^3 W_0=O\left(\frac{a}{r}\right).
\eeaa

 \item $\err[\squared_2 \qf]$ is  the nonlinear correction term, which under the additional condition\footnote{This additional condition makes the structure of $\err[\squared_2 \qf]$ in \eqref{eq:MaiThmParq-err} possible. This structure is essential in the control of the nonlinear term in Chapter \ref{chapter-full-RWforqf}, see also Chapter 5 in \cite{KS} in the particular case of perturbations of Schwarzschild.}  
\beaa
\Hc \in \Ga_g
\eeaa 
is given schematically by  the expression
   \bea
   \lab{eq:MaiThmParq-err}
   \bsplit
 \err[\squared_2 \qf]&= r^2 \frak{d}^{\leq 3} (\Ga_g \c (A, B))+ \nab_3 (r^2 \frak{d}^{\leq 2}( \Ga_b \c (A, B)))\\
&+\frak{d}^{\leq 1} (\Ga_g \c \qf) + r^{3}\dk^{\leq 2} \big( \Ga_b \c \Ga_g \c \Ga_g\big).
\end{split}
 \eea
 \end{itemize}
 \end{theorem}
 
 We now collect some remarks on Theorem \ref{MAIN-THEOREM-PART1}.

 \begin{remark}\label{rmk:inabT-term} 
 The first order term in $\nab_\T$ on the LHS of \eqref{wave-equation-qf} presents good divergence properties because of its structure. More precisely, since it is given by an imaginary function multiplied by $\nab_\T\qf$, this term cancels out in the derivation of the energy estimates in the trapping region.
 
Notice that there is a conformally invariant definition of $\T$ (where in Definition \ref{Definition:vfsTZ} the vectors are replaced by their conformally invariant counterparts), but since $\qf$ is 0-conformally invariant the two definitions coincide in this case, and equation \eqref{wave-equation-qf} is fully conformally invariant.
 \end{remark}

\begin{remark}\label{rmk:propertiesofthepotentialV}
The potential term $V$ is a real function which coincides, for zero angular momentum, with the potential  in Schwarzschild in \cite{KS}, given by
\beaa
V_0=-\trch\trchb=\frac{4\De}{(r^2+a^2) |q|^2}.
\eeaa
 Observe that the fact that the potential is \textbf{real} is crucial in the derivation of the estimates for the Regge-Wheeler equation. This is obtained through the choice of the imaginary part of the scalar function $C_1$ which gives $\Im(V)=0$. This is obtained in Proposition \ref{prop:potential-real}. 
\end{remark}

\begin{remark}\label{rmk:choice-C1-C2} 
The scalar functions $C_1$ and $C_2$ are chosen in order to obtain cancellation of terms in the derivation of the gRW equation. The real parts of $C_1$ and $C_2$ are chosen to obtain the cancellation of the highest order terms in the commutator $[Q, \LL]$, see Proposition \ref{first-intermediate-step-main-theorem}, and they coincide with the values in Schwarzschild in \cite{KS}. As mentioned above, the imaginary part of $C_1$ is chosen so to cancel the imaginary part of the potential $\Im(V)$ and the imaginary part of the functions appearing as coefficients of the highest order terms in  $ L_\qf[A]$. Finally, the imaginary part of  $C_2$ is chosen  in order to have a good transport relation between $\qf$ and $A$, see Lemma \ref{lemma:factorizationofqfusefulfortransportequations}, a fact  used in the derivation of the estimates.
\end{remark}

\begin{remark}\label{rmk:properties-lot} 
In the derivation of the estimates for the Regge-Wheeler equation, the linear second order operator $L_{\qf}[A]$, which contains at most two derivatives of $A$, will be treated as a lower order term through transport estimates. The form of the highest order terms in $L_{\qf}[A]$, i.e. $\nab_\T (\nab_3 A)$ and $\nab_\Z (\nab_3 A)$, is crucial for the estimates in the trapping and relies on an important cancellation obtained in Proposition \ref{prop:lot-real}.
\end{remark}

We now describe the steps of the proof of Theorem \ref{MAIN-THEOREM-PART1}, relying on the computations collected in Appendix \ref{proof:theorem-main-Part1}.

\begin{enumerate}
\item In Step 1, obtained in section \ref{proof:step1-main-thm-part1}, we apply to the Teukolsky equation $ \LL(A)= \err[\LL(A)]$ the operator $Q=\nabc_3\nabc_3  + C_1  \nabc_3 + C_2  $. 
We then compute the commutator $[Q, \LL]$, and find conditions on the real part of $C_1$ and $C_2$ in order to have lower order terms in the commutator, denoted $L_Q(A)$, which are $O(|a|)$, see Proposition \ref{first-intermediate-step-main-theorem}.

\item In Step 2, obtained in section \ref{proof:step2-main-thm-Part1}, we derive the wave equation for $Q(A)$ and $\qf$. We first obtain, see Proposition \ref{prop:wave-eq-Q}, the wave equation for $Q(A)$, and we then rescale $Q(A)$ through the function $f=q\ov{q}^3$ by defining $\qf= f Q(A)$, which satisfies an equation of the form, see Proposition \ref{prop:rescaling-f},
\beaa
 \squared_2 \qf  - i   \frac{4a\cos\th}{|q|^2} \nab_\T  \qf- V_1 \qf   &=&  \widetilde{L_{\qf}[A]} + f \Big( \err[\squared_2 Q]+ \Ga_g \c \dk^{\leq 1} Q \Big).
\eeaa

\item In Step 3, obtained in section \ref{proof:step-3-main-thm-Part1}, we derive additional conditions on the scalar functions $C_1$ and $C_2$ by imposing the reality of the potential of the equation for $\qf$ as well as a specific structure for the lower order terms. More precisely, this is obtained by first imposing that the potential and the terms in $L_\qf[A]$ involving two derivatives of $A$ should be real, see Proposition \ref{prop:potential-real}, and then by combining the potential term with the lower order term to have only $\T$ and $\Z$ derivatives of $\nab_3A$. 
\end{enumerate}


\subsection{The real part of the gRW equation}
\lab{sec:realpartofgRWequationforqf:partI}


 Since $\qf\in\sk_2(\CCC)$ is  a complex anti-self dual tensor, we can decompose it as
 \bea
 \qf=\psi + i \dual \psi
 \eea
for some $\psi=\Re(\qf) \in \sk_2(\RRR)$. Taking the real part of \eqref{wave-equation-qf}, since $V$ is real, we then obtain an equation for $\psi$, which is given by
\beaa
 \squared_2 \psi  +\frac{4 a\cos\th}{|q|^2} \dual \nab_\T\psi   -V  \psi&=&\Re(   L_{\qf}[A]) + \Re(\err[\squared_2 \qf]).
\eeaa
 We summarize in the following.
\begin{proposition}\label{prop:eq-real-gen-RW} 
The tensor $\psi \in \sk_2(\RRR)$ satisfies
 \bea\label{eq:Gen.RW-pert}
\squared_2 \psi -V_0\psi=- \frac{4 a\cos\th}{|q|^2}\dual \nab_T  \psi+N, \qquad  V_0= \frac{4\De}{ (r^2+a^2) |q|^2}, 
\eea
 with the right hand side $N$ being given by
 \bea\label{eq:definition-N-psi}
 \begin{split}
 N&:= \big( V-V_0\big) \psi + \Re( L_{\qf}[A])+\Re(\err[\squared_2 \qf])\\
 &= N_0+N_L+N_{\err}
 \end{split}
 \eea
 where:
 \begin{itemize}
 \item[-] $N_0$ denotes the zero-th order term in $\psi$, i.e.
 \bea\lab{eq:definition-N-0-psi}\lab{eq:definition-N-0-psi-ch7}
 N_0:= \left(V- \frac{4\De}{ (r^2+a^2) |q|^2}\right)\psi=O\left(\frac{a}{r^4}\right) \psi.
 \eea
 
 \item[-] $N_L$ denotes the lower order terms in $\psi$, i.e.
 \bea\lab{eq:definition-N-L-psi}\lab{eq:definition-N-L-psi-ch7}
 \begin{split}
 N_L&:= \Re\Bigg(q\ov{q}^3    \Bigg[- \frac{8a^2 \De}{r^2|q|^4}\nab_\T \nab_3 A -\frac{8a  \De }{r^2|q|^4} \nab_\Z \nab_3A\\
&+ W_4\nab_4A+ W_3\nab_3A+W\c\nab A +W_0  A\Bigg] \Bigg)
\end{split}
 \eea
where $W_4$, $W_3$, $W_0$ are complex functions  of $(r, \th)$, and $W$ is the product of a complex function of $(r, \th)$ with $\dual\Re(\Jk)$, having the following fall-off in $r$ 
\beaa
q\ov{q}^3 W_4= q\ov{q}^3 W_3=q\ov{q}^3 W=O\left(a\right), \qquad  q\ov{q}^3 W_0=O\left(\frac{a}{r}\right).
\eeaa

 \item[-] $N_{\err}[\psi]$ denotes the error terms, i.e.
 \bea
 N_{\err}[\psi] := \Re(\err[\squared_2 \qf])
 \eea
 which are schematically given by
  \beaa
 N_{\err}[\psi]&=& r^2 \frak{d}^{\leq 2} (\Ga_g \c (\a, \b)) + \nab_3 (r^2 \frak{d}^{\leq 2}( \Ga_b \c (\a, \b)))\\
 &&+\frak{d}^{\leq 1} (\Ga_g \psi)  + r^{3}\dk^{\leq 2} \big( \Ga_b \c \Ga_g \c \Ga_g\big).
 \eeaa
 \end{itemize}
 
   Also, recall that $\psi$ and $A$ are related by the differential relation:
  \beaa
 \psi= \Re\Big(q \ov{q}^{3} \big( \nabc_3\nabc_3 A + C_1  \nabc_3A +C_2   A\big)\Big), 
 \eeaa
 with
\beaa
\begin{split}
C_1&=2\trchb - 2\frac {\atrchb^2}{ \trchb}  -4 i \atrchb, \\
C_2  &= \frac 1 2 \trchb^2- 4\atrchb^2+\frac 3 2 \frac{\atrchb^4}{\trchb^2} +  i \left(-2\trchb\atrchb +4\frac{\atrchb^3}{\trchb}\right).
\end{split}
\eeaa
 \end{proposition}


\section{Generalized Regge-Wheeler equation  for $\protect\underline{\qf}$}
\label{section:gRW-equation-qfb}


In this section, we derive the generalized Regge-Wheeler equation for $\qfb$.


\subsection{The Teukolsky equation for $\Ab$}


Here we derive the Teukolsky equation for $\Ab$.  In order to capture correctly the non linear terms in the equation, we express the Bianchi identity for $\Ab$ in terms of 
\beaa
\Ab_4=\nabc_4\Ab +\frac 1 2 \tr X \Ab,
\eeaa 
see Definition \ref{definition-Psi_3-Psi_4},  which has an improved decay rate as compared to $\nabc_4 \Ab$. In the derivation of the Teukolsky equation below, we express explicitly the error terms which decay less than $r^{-2} \dk^{\leq 1}(\Ga_g \c \Ga_b)$.

\begin{proposition}\label{PROP:TEUK-AB} 
We have 
     \bea\label{eq:Teuk-Ab-Ab4}
     \begin{split}
\left(\nabc_3+2 \ov{\tr \Xb}+\frac 1 2 \tr \Xb\right)\Ab_4&=    \frac 1 4 \big( \DDc +H +4\Hb\big) \hot \big(\DDbc \c\Ab +  \ov{H}  \c  \Ab\big)\\
&+3P\Ab+ \err_{TE}
 \end{split}
\eea
where $ \err_{TE}$ is given schematically by
\beaa
\err_{TE}&=& \tr X \Xib \hot  \Bb+( \DDbc \c\Bb) \Xbh+ (\widehat{\Xb}\c\ov{\Hc}) \Bb +(\Ga_b \c \Ga_b)\c (A, B)+ r^{-2} \dk^{\le 1} (\Ga_g \c \Ga_b).
\eeaa
\end{proposition}

\begin{proof} 
See section \ref{proof-teukolsky-Ab}.
\end{proof}

\begin{remark} 
Observe that, if we were to obtain the Teukolsky equation for $\Ab$ in the symmetric way as the Teukolsky equation for $A$ in Proposition \ref{TEUKOLSKY-PROPOSITION}, we would obtain
 \beaa
 \underline{\LL}(\Ab)&=& \err[\underline{\LL}(\Ab)]
 \eeaa
 where
\beaa
\begin{split}
\underline{\LL}(\Ab) &=-\nabc_3\nabc_4\Ab+ \frac{1}{4}\DDc\hot (\DDbc \c \Ab)+\left(- \frac 1 2 \tr \Xb -2\ov{\tr \Xb} \right)\nabc_4\Ab\\
&-\frac{1}{2}\tr X \nabc_3\Ab+\left( 4\Hb+H +\ov{H} \right)\c \nabc \Ab+ \left(-\ov{\tr \Xb} \tr X +2P\right) \Ab+  \Hb   \hot (\ov{H} \c \Ab)
\end{split}
\eeaa
with error term  expressed schematically as
\beaa
\err[\underline{\LL}(\Ab)]&=& r^{-1}\dk^{\leq 1}\big(  \Ga_b \c \Bb \big) + \Ga_b \c \Ga_b \c \Ga_g.
\eeaa
The above error terms are not acceptable in the forthcoming derivation of the generalized Regge-Wheeler equation for $\qfb$, and we therefore rely instead on Proposition \ref{PROP:TEUK-AB}, i.e. we express the Teukolsky equation in terms of $\Ab_4=\nabc_4\Ab +\frac 1 2 \tr X \Ab$, which has an improved decay rate as compared to $\nabc_4 \Ab$.
\end{remark}


\subsection{The invariant quantities $\protect\underline{Q}$ and $\protect\underline{\qf}$}


In  this section  we consider the analog $\qfb$ of $\qf$ and derive its corresponding gRW equation. 

\begin{definition}
\lab{Definition:Define-qfb}
Given a  fixed  null  pair  $(e_3, e_4)$  and  scalar functions $r$ and $\th$ as in Section \ref{sec:setupandlinearizedquantities}, we define our second  main quantity $\qfb \in \sk_2(\CCC)$ as 
   \bea\label{eq:definition-qfb}
\qfb&=&  \ov{q} q^3  \underline{Q}(\Ab)=\ov{q} q^3 \left(  \nabc_4\nabc_4 \Ab + \underline{C}_1  \nabc_4\Ab + \underline{C}_2  \Ab\right),
\eea
with complex scalars
\bea\label{eq:Cb1-Cb2-comparison-Ma}
\begin{split}
\und{C}_1&=2\trch - 2\frac {\atrch^2}{ \trch}  -4 i \atrch, \\
\und{C}_2  &= \frac 1 2 \trch^2- 4\atrch^2+\frac 3 2 \frac{\atrch^4}{\trch^2} +  i \left(-2\trch\atrch +4\frac{\atrch^3}{\trch}\right).
\end{split}
\eea
\end{definition}
 
 In the particular case of Kerr,  the quantity $\qfb$  defined above can be factorized as follows.
 
\begin{proposition}
\lab{Prop:factorization-qfb}
In Kerr,  the quantity $\qfb$  defined in \eqref{eq:definition-qfb} with $\und{C_1}, \und{C_2}$ given by \eqref{eq:Cb1-Cb2-comparison-Ma} can be factorized as
\bea
\lab{eq:Prop-factorization-qfb}
r \nabc_4\left( r^2   \left(   \nabc_4  \left( r \frac{q^4}{r^4}  \Ab\right)\right)\right)= \frac{q }{\ov{q} } \qfb.
\eea
Alternatively,
\bea
\lab{eq:secondfactorization-qfb}
\nabc_4\nabc_4\left(\frac{q^4}{r^2} \Ab \right) &=&  \frac{q }{\ov{q} } r^{-2} \qfb.
\eea
\end{proposition}  

  \begin{proof}
  We  proceed as in the proof of Proposition \ref{PROP:FACTORIZATION-QF}. Since  the formulas are manifestly scale invariant    we chose the  outgoing normalization  
  of  $e_4$ such that  $\om=0$ and  $e_4(r)=1$ and 
  $ \trch=\frac{2r}{|q|^2},\,\, \atrch=\frac{2a\cos\th}{|q|^2}$. Thus  to check 
 \eqref{eq:Prop-factorization-qfb} we have to  show that 
 \beaa
 I:=r \nab_4\Big( r^2   \big(   \nab_4   ( r f  \Ab)\big)\Big)= r^4  f\Big( \nab_4 \nab_4 \Ab+ I_1    \nab_4 \Ab+ I_2 \Ab\Big)
 \eeaa
 with  $f=\frac{q^4}{r^4}= f_1 f^2 _2, \,  f_1=\frac{|q|^4}{r^4}, \,  f_2=\frac{q}{\ov {q}}$
 and 
 \beaa
  I_1&=& 2 f^{-1}( e_4 f +2 r^{-1}  f),  \\
 I_2&=& f^{-1} \Big(   \nab_4 e_4 (f)+  4 r^{-1}   e_4 f + 2 r^{-2} f\Big).
 \eeaa
 Note that $e_4 f_1= -\frac{ 4a^2 \cos^2 \th}{ r|q|^2 }  f_1, \, 
 e_4 f_2 =-i\atrch f_2$ and therefore
\beaa
  e_4f&=&  \left( -\frac{ 4a^2 \cos^2 \th}{ r|q|^2 } - 2 i \atrch\right) f.
 \eeaa
 Thus,
 \beaa
 I_1&=&    2 \big( f^{-1}  e_4 f +2 r^{-1}  \big) =
  2  \left( -\frac{ 4a^2 \cos^2 \th}{ r|q|^2 }  +   \frac{2r}{|q|^2}- 2 i \atrch     \right)\\
  &=& 2  \left( -\frac{ 2a^2 \cos^2 \th}{ r|q|^2 }  +   \frac{2r}{|q|^2}- 2 i \atrch     \right)\\
  &=&2\trch - 2\frac {\atrch^2}{ \trch}  -4 i \atrch=\und{C_1}. 
 \eeaa
  Similarly, see the proof  of  Proposition \ref{PROP:FACTORIZATION-QF},  $I_2=\und{C}_2$ which establishes \eqref{eq:Prop-factorization-qfb}.
  \end{proof}
  
In the derivation of the gRW equation for $\qfb$, we will rely on the following more involved factorization of $\qfb$.
\begin{lemma}
Let
\bea\label{eq:widetilde-Qb-Ab}
\widetilde{\underline{Q}(\Ab)}:=\left(\nabc_4 +\frac{5}{2}\tr X-\ov{\tr X} - 2\frac {\atrch^2}{ \trch} \right)\Ab_4.
\eea
Then, we the following holds 
\bea\label{relation-qfb-widetilde-underline-Q-Ab}
\qfb=\ov{q}q^3\widetilde{\underline{Q}(\Ab)}+ O(a^2)\Ab+r^{-2}\dk^{\leq 1}(\Ga_g\c\Ga_b).
\eea
\end{lemma}

\begin{proof}
We have, see Lemma \ref{lemma:almostfactorizationfoqfbneededforcontralangularderivativesPsib:chap12},
\beaa
\bsplit
 &\ov{q}q^3\left(\nabc_4 +\frac{3}{2}\tr X - 2\frac {\atrch^2}{ \trch}  -2i \atrch\right)\left(\nabc_4\Ab +\frac{1}{2}\tr X\Ab\right)\\ 
 &= \qfb + O(a^2)\Ab+r^{-2}\dk^{\leq 1}(\Ga_g\c\Ga_b),
 \end{split}
\eeaa
which together with the fact that $\ov{\tr X}=\trch+i\atrch$ and $\Ab_4=\nabc_4\Ab +\frac 1 2 \tr X \Ab$ proves the lemma.
\end{proof}


\subsection{The derivation of the gRW equation for $\protect\underline{\qf}$}
\lab{subsection:DerivationgRW-eqbar}


We state below   the  gRW equation satisfied by $\qfb$.

 \begin{theorem}\lab{THEOREM:EQ-QFB}
 The invariant symmetric traceless $2$-tensor $\qfb \in \sk_2(\CCC)$  in Definition \ref{Definition:Define-qfb} satisfies the equation
 \bea\label{wave-equation-qfb}\lab{eq:gRW-qfb}
 \squared_2 \qfb + i \frac{4 a\cos\th}{|q|^2} \nab_\T \qfb   - V  \qfb &=&   L_{\qfb}[\Ab] + \err[\squared_2 \qfb]
 \eea
 where:
 \begin{itemize}
 \item The potential $V$ is the \textbf{real} scalar function given by
 \bea
 V&=& \frac{4}{|q|^2}\frac{r^2-2mr+2a^2}{r^2}-\frac{4a^2\cos^2\th}{|q|^6}( r^2+6mr+a^2\cos^2\th ),
 \eea
  which for $a=0$ coincides with the potential of the Regge-Wheeler equation  in Schwarzschild, i.e. $V=-\trch\trchb+ O(\frac{|a|}{r^4})$.
 
 \item $ L_{\qfb}[\Ab] $ is  a  linear second order   operator  in $\Ab$   given in the ingoing frame by
  \beaa
L_{\qfb}[\Ab] = q\ov{q}^3    \left( \frac{8a^2 \De}{r^2|q|^4}\nab_\T \Ab_4 +\frac{8a  \De }{r^2|q|^4} \nab_\Z  \Ab_4+ \underline{W}_4\Ab_4+ \underline{W}_3 \nab_3 \Ab+\underline{W}\c\nab\Ab +\underline{W}_0  \Ab\right),
 \eeaa
where $\underline{W}_4$, $\underline{W}_3$, $\underline{W}_0$ are complex functions  of $(r, \th)$ and $\underline{W}$ is the product of a complex function of $(r, \th)$ with $\dual\Re(\Jk)$, with the following fall-off\footnote{Note that the fall-off provided here  for  $\underline{W}_4$, $\underline{W}_3$, $\underline{W}_0$ and $\underline{W}$  is stronger than the one  for $W_4$, $W_3$, $W_0$ and $W$ in $L_{\qf}[A]$, see section \ref{subsection:DerivationgRW-eq}. In fact, provided we replace $\nab_3A$ with $A_3$ in the definition of $L_{\qf}[A]$, $W_4$, $W_3$, $W_0$ and $W$ satisfy the same fall-off is $r$, but this stronger fall-off is unnecessary for the  weighted estimates derived for $\qf$ in Chapter \ref{chapter-full-RWforqf}.} in $r$
\beaa
 q\ov{q}^3 \underline{W}_4, q\ov{q}^3 \underline{W}=O\left(\frac{a^2}{r}\right), \qquad  \qquad  q\ov{q}^3 \underline{W}_3, q\ov{q}^3 \underline{W}_0=O\left(\frac{a^2}{r^2}\right).
\eeaa

 \item $\err[\squared_2 \qfb]$ is  the nonlinear correction term, which under the additional conditions\footnote{In fact, it suffices to assume that $\Xi\in r^{-2}\Ga_g$ and $\Hbc\in r^{-1}\Ga_g$. These additional conditions make the structure of $\err[\squared_2 \qfb]$ in \eqref{eq:MaiThmParq-err-bar} possible. This structure is essential in the control of the nonlinear term in Chapter \ref{chapter-full-RWforqfb}.}  
\beaa
\Xi=0, \qquad \Hbc=0, \quad\textrm{for}\quad r\geq r_0,
\eeaa
is given schematically by  the expression
 \bea\lab{eq:MaiThmParq-err-bar}
\nn \err[\square_2 \qfb]&=&  r^2 \dk^{\leq 2}(\Ga_b \c (A, B))+ \dk^{\leq 3} (\Ga_g \c \Ga_b).
 \eea
 \end{itemize}
 \end{theorem}

We now describe  the steps of the proof of Theorem \ref{THEOREM:EQ-QFB}, relying on the computations collected in Appendix \ref{proof:thm-eq-qfb}.
\begin{enumerate}
\item In Step 1,  we take the first derivative in the $\nabc_4$ direction of the Teukolsky equation for $\Ab$. We express explicitly the error terms which decay less than $r^{-3} \dk^{\leq 2}(\Ga_g \c \Ga_b)$. We obtain, see Proposition \ref{prop:appendix-nabc43Ab4},
\beaa
&&\left( \nabc_4+  \tr X+ \frac 1 2\ov{\tr X}\right) \left( \nabc_3+ 2 \ov{\tr \Xb}+\frac 1 2 \tr \Xb\right) \Ab_4\\
&=&    \frac 1 4\big(  \DDc+H+5\Hb \big) \hot \Big( \DDbc \c \Ab_4+(\ov{H}+\ov{\Hb}) \c \Ab_4\Big)+3P\Ab_4\\
&& +3\left( \frac 1 2\ov{\tr X}-\tr X\right) P\Ab+\mathcal{J}_{434} \hot \big(\DDbc \c\Ab +  \ov{H} \c \Ab \big)+  \err_{434},
\eeaa
for a one-form $\mathcal{J}_{434}$ and error terms given by 
   \beaa
  \err_{434}&=& \nabc_4\err_{TE} +\left( \frac 1 2 \tr X+ \ov{\tr X}\right)\err_{TE}\\
  &&+\DDc\hot( \Xh \c \ov{\DDc} \Ab)+ \DDc \hot ((  \DDc\c\ov{\Xh})\c \Ab)\\
 &&+ r^{-1} \dk^{\le 1 } ( (A,B)\c \Ga_b) + r^{-3}\dk^{\leq1}( \Ga_g  \c \Ga_b).
  \eeaa

\item In Step 2, we take the second derivative in the $\nabc_4$ direction of the Teukolsky equation for $\Ab$. Again, we express explicitly the error terms which decay less than $r^{-4} \dk^{\leq 3}(\Ga_g \c \Ga_b)$. We obtain, see Proposition \ref{prop:appendix-nabc4nabc4nabc3Ab},
\beaa
&&\left(\nabc_4 + 3  \tr X-\frac 1 2 \ov{\tr X}- 2\frac {\atrch^2}{ \trch}\right)\left( \nabc_4+  \tr X+ \frac 1 2\ov{\tr X}\right) \left( \nabc_3+ 2 \ov{\tr \Xb}+\frac 1 2 \tr \Xb\right) \Ab_4\\
         &&+ \big(2\tr X- \ov{\tr X}  \big)3P\Ab_4 \\
&= &  \frac 1 4  \big( \DDc+H+6\Hb\big)\hot \Big( \DDbc \c \widetilde{\underline{Q}(\Ab)}+(\ov{H}+2\ov{\Hb}) \c \widetilde{\underline{Q}(\Ab)} \Big)+3P \ \widetilde{\underline{Q}(\Ab)}\\
   &&+\mathscr{L}[\Ab, \Ab_4, \DD \Ab_4, \DD \Ab]+\err_{4434}
   \eeaa
where $\widetilde{\underline{Q}(\Ab)}$ is defined in \eqref{eq:widetilde-Qb-Ab}, and is such that 
\beaa
 \qfb=\ov{q}q^3\widetilde{\underline{Q}(\Ab)}+ O(a^2)\Ab+r^{-2}\dk^{\leq 1}(\Ga_g\c\Ga_b) .
 \eeaa

We denote by $\mathscr{L}[\Ab, \Ab_4, \DD \Ab_4, \DD \Ab]$   linear order terms in $\Ab, \Ab_4, \DD \Ab, \DD \Ab_4$ schematically given by \eqref{eq:mathscr-L-Ab}, i.e.
\beaa
\mathscr{L}[\Ab, \Ab_4,  \DD \Ab, \DD \Ab_4]&=&O(ar^{-3}) \big( \DDbc \c \Ab_4+(\ov{H} +\ov{\Hb}) \c \Ab_4\big)\\
&&+O(a^2r^{-4})\big(\DDbc \c\Ab +  \ov{H} \c \Ab \big) +O(a^2r^{-7}) \Ab.
\eeaa
The error terms are given by 
\beaa
\err_{4434}&=& \nabc_4  \err_{434}+ \left(   \tr X+\frac 3 2 \ov{\tr X}+2\frac {\atrch^2}{ \trch} \right)  \err_{434}\\
&&+ r^{-2} \dk^{\leq 1}( (A,B) \c \Ga_b)+ r^{-4} \dk^{\leq 2}(\Ga_g \c  \Ga_b)\\
&&+ r^{-1}\left(\nabc_4(\Xh \c \Hc)+\frac 3 2 \tr X (\Xh \c \Hc)\right)\c  \dk^{\leq 1} \Ab.
\eeaa
By commuting the operators on the left hand side of \eqref{eq:nabc4-nabc4-nabc4-Ab4} and using the expression of $\qfb$ in terms of $\widetilde{\underline{Q}(\Ab)}$ given in \eqref{relation-qfb-widetilde-underline-Q-Ab}, we can deduce the form of the equation for $\qfb$. 

\item In Step 4, we show that the error terms $\err_{4434}$ can be simplified to obtain, see Proposition \ref{prop:appendix-error-terms-qfb}, 
\beaa
\err_{4434}= r^{-2} \dk^{\leq 2}((A, B)\c \Ga_b)+ r^{-4} \dk^{\leq 3}(\Ga_g \c  \Ga_b) .
\eeaa
We show that the terms which behave worse that $r^{-4} \dk^{\le  3 } (\Ga_g \c \Ga_b)$ get improved once applied the differential operators $( \nabc_4+  \tr X+ \frac 1 2\ov{\tr X})$ and $\nabc_4 + 2  \tr X+\frac 1 2 \ov{\tr X}$ by making use of the renormalized Bianchi identities and improved decay in the null structure equations. 
 This ends the proof of Theorem \ref{THEOREM:EQ-QFB}.
\end{enumerate}

As in section \ref{sec:realpartofgRWequationforqf:partI} in the case of $\qf$, we infer from Theorem \ref{THEOREM:EQ-QFB}, for the real part of $\qfb$ denoted  $\underline{\psi}=\Re(\qfb)$, the following real equation:
 \bea\label{eq:Gen.RW-pert-qfb}
\squared_2 \underline{\psi} -V_0\underline{\psi}= \frac{4 a\cos\th}{|q|^2}\dual \nab_T  \underline{\psi}+N, \qquad  V_0= \frac{4\De}{ (r^2+a^2) |q|^2}, 
\eea
 where the right hand side $N$ is given by $ N = N_0+N_L+N_{\err}$, with
 \beaa
 N_0=O\left(\frac{a}{r^4}\right)\underline{\psi}, \qquad N_L=\Re( L_{\qfb}[\Ab]), \qquad N_{\err}=\Re(\err[\square_2 \qfb]).
  \eeaa


\section{Teukolsky-Starobinski identity}


We state here, in the context of perturbations of Kerr, one of the Teukolsky-Starobinski identities, which relate the complex curvature components $A$ and $\Ab$ through fourth-order differential operators. 

\begin{proposition}\label{THEOREM:TEUK-STAR} 
Assume that $\Xi=0$ in $r \leq r_0$. The complex tensors $A$, $\Ab \in \sk_2(\CCC)$ satisfy the following relation in the region $r \leq r_0$
\bea\label{eq:TS-id}
\Big(\nabc_4+2\tr X\Big)^4\Ab&=& r^{-4} \dk^{\leq 4} A+ \dk^{\leq 3}\big(\Ga_b\c\Ga_g\big).
\eea
\end{proposition}

\begin{proof} 
See Appendix \ref{proof:teukolsky-starobinski}.
\end{proof}

\begin{remark} 
Proposition \ref{THEOREM:TEUK-STAR} is stated without proof  in Chapter 7 of \cite{KS:Kerr}. Both the assumption $\Xi=0$ and the restriction to $r\leq r_0$ are unnecessary and assumed only for convenience, as they hold when applying Proposition \ref{THEOREM:TEUK-STAR}  in Chapter 7 of \cite{KS:Kerr}. In particular, the restriction to $r\leq r_0$ allows us to avoid having to track the precise powers of $r$ in the nonlinear terms.
\end{remark}

\begin{remark}
Choosing a normalization such that $\om\in\Ga_g$, and hence $\tr X=\frac{2}{q}+\Ga_g$, we infer from Proposition \ref{THEOREM:TEUK-STAR}, for $\Xi=0$ in the region $r\leq r_0$, 
\beaa
 \frac{1}{q^7}\nab_4\left(q^2\nab_4\left(q^2\nab_4\left(q^2\nab_4(q\Ab)\right)\right)\right)&=& r^{-4} \dk^{\leq 4} A+ \dk^{\leq 3}\big(\Ga_b\c\Ga_g\big).
\eeaa
\end{remark}


\section{The wave equation for $P$}
\lab{section:waveEqforP}


Here we derive the wave equation satisfied by the curvature component $P$.

\begin{lemma}
\lab{LEMMA:WAVEEQP}
The curvature component $P$ satisfies the following scalar wave equation:
\bea
\lab{eq:WaveEq-forP}
\bsplit
\square_\g P&=\tr X \nab_3P+ \ov{\tr\Xb} \, \nab_4P-\ov{H} \c \DD P- \Hb\c \ov{\DD} P   \\
&+\frac{3}{2}\Big[ \ov{\tr\Xb} \tr X+2P  -2 \Hb \c  \ov{H}  \Big] P + \err[\square_\g P],
\end{split}
\eea
with error terms given by
\bea
\label{eq:err-square-P-precise}
\begin{split}
\err[\square_\g P]&=  -\nabc_3(\ov{\Xi}\c\Bb) -\frac{1}{4}\nabc_3(\Xbh\c \ov{A})+\frac 1 2 \DDc \c \left( \Bb\c \ov{\Xh}+\frac{1}{2}\ov{A}\c \Xib\right)\\
& +\Hb \c \left( \Bb\c \ov{\Xh}+\frac{1}{2}\ov{A}\c \Xib \right)+\frac 12 \Big(-\ov{\Bb} \c  B +\Xib \c \nab_4 \ov{B} -\ov{\Xbh} \c \DDc \ov{B}-\ov{H} \c \ov{\Xbh} \c B \Big)\\
&- \frac 1 2\big(2\ov{\tr\Xb}+ \tr\Xb \big)\left(\ov{\Xi}\c\Bb+\frac{1}{4}\Xbh\c \ov{A}\right)+\frac 12 H \c \left( \Bb\c \ov{\Xh}+\frac{1}{2}\ov{A}\c \Xib \right)\\
&+\left(-\frac{1}{2}\ov{\DDc}\c\Xbh-i\Im(\tr X)\Xib+\Bb+\nabc_4\Xib -\frac{1}{2}\Xbh\c(\ov{\Hb}-\ov{H}) +\Bb\right)\c\ov{B}\\
&-\frac 32 \big(\Xib\c\ov{\Xi}-\frac{1}{2}\Xbh\c\ov{\Xh}\big)P.
\end{split}
\eea
\end{lemma}

\begin{proof} 
See Appendix \ref{proof-wave-P}. 
\end{proof}

\begin{remark}\label{rem:error-terms-square-P} 
From the expression in \eqref{eq:err-square-P-precise} and using the null structure equation for $\nabc_3 \Xbh$, observe that the error terms $\err[\square_\g P]$ can be schematically written as 
\beaa
\err[\square_\g P]&=&r^{-1}   \dk^{\leq 1} (\Ga_g \c\Bb) +\nab_3 ( \Xi \c \Bb)+r^{-1}\dk^{\leq 1} (\Ga_b \c (A, B))+ r^{-3} (\Ga_g \c \Ga_b)-\Ab\c \ov{A}.
\eeaa
The above structure will be used in Chapter \ref{Chapter:EN-MorforPc}, see \eqref{eq:P-WaveEq-M8}.
\end{remark}


\subsection{A renormalized  wave equation} 


\begin{lemma}
\lab{Le:squareq^2Psi}
Let $\Psi$ a scalar function solution to the following wave equation 
\bea
\lab{eq:Le-squareq^2Psi}
\square_\g\Psi &=&  \tr X\nab_3\Psi +\ov{\tr\Xb}\nab_4\Psi  -\ov{H}\c\DD\Psi  - \Hb\c\ov{\DD}\Psi +V\Psi+F, 
\eea
where $V$ is a potential and $F$ a scalar function. Then, we have
\beaa
\square_\g(q^2\Psi) &=& \Big[V +q^{-2}\square_\g(q^2)\Big]q^2\Psi +r\Ga_b\c \dk\Psi+q^2F.
\eeaa
\end{lemma}

\begin{proof}
We have
\beaa
\square_\g(q^2\Psi) &=&  q^2\square_\g(\Psi)+\square_\g(q^2)\Psi+2\g^{\a\b}\pr_\a(q^2)\pr_\b(\Psi)\\
 &=& q^2\square_\g(\Psi)+\square_\g(q^2)\Psi - e_3(q^2)e_4\Psi- e_4(q^2)e_3\Psi+2\nab(q^2)\nab\Psi\\
&=& q^2\left[\square_\g(\Psi)+q^{-2}\square_\g(q^2)\Psi -\frac{2e_3(q)}{q}e_4\Psi-\frac{2e_4(q)}{q}e_3\Psi+4\frac{\nab(q)}{q}\c\nab\Psi\right].
\eeaa
Since
\beaa
 \frac{\ov{\DD}(q)}{q}\c\DD \Psi  + \frac{\DD(q)}{q}\c\ov{\DD}\Psi &=& \frac{(\nab-i\dual\nab)(q)}{q}\c(\nab+i\dual\nab) \Psi  + \frac{(\nab+i\dual\nab)(q)}{q}\c(\nab-i\dual\nab)\Psi\\
&=& 4\frac{\nab q}{q}\c \nab\Psi,
\eeaa
we infer
\beaa
\square_\g(q^2\Psi) = q^2\left[\square_\g(\Psi)+q^{-2}\square_\g(q^2)\Psi -\frac{2e_3(q)}{q}e_4\Psi-\frac{2e_4(q)}{q}e_3\Psi +\frac{\ov{\DD}(q)}{q}\c\DD \Psi  + \frac{\DD(q)}{q}\c\ov{\DD}\Psi\right].
\eeaa
Using
\beaa
\frac{e_4(q)}{q}=\frac{1}{2}\tr X+r^{-1}\Ga_b, \quad \frac{e_3(q)}{q}=\frac{1}{2}\ov{\tr\Xb}+\Ga_b, \quad H=\frac{\DD(\ov{q})}{\ov{q}}+r^{-1} \Ga_b, \quad \Hb=\frac{\DD(q)}{q}+\Ga_b,
\eeaa
we deduce
\beaa
\square_\g(q^2\Psi) &=& q^2\Bigg[\square_\g(\Psi)+q^{-2}\square_\g(q^2)\Psi - \ov{\tr\Xb}e_4\Psi-\tr X e_3\Psi +\ov{H}\c\DD \Psi  + \Hb\c\ov{\DD}\Psi\Bigg]\\
&&+ r\Ga_b \c \dk \Psi+ q^2 F.
\eeaa
Since $\Psi$ satisfies \eqref{eq:Le-squareq^2Psi} we infer
\beaa
\square_\g(q^2\Psi) &=& \Big[V +q^{-2}\square_\g(q^2)\Big]q^2\Psi + r\Ga_b\c  \dk \Psi +q^2F
\eeaa
as stated. 
\end{proof}


\subsection{Wave equations for $ q^2(\T, \Z)P$}


\begin{lemma}
\lab{Le:waveeq.TZ(P)}
The linearized  quantities $q^2 \T P,\, q^2  \Z P$ verify the following wave equations
\bea
\lab{waveeq:TZ(P)}
\bsplit
\square_\g(q^2\T P) &=W q^2\T P+ r\Ga_b\c  \dk^2 P +r^2\dk \err[\square_\g P]+r^2 \dk^{\leq 1} \big(\Ga_g \c \dk^{\leq 1} P\big)+r^2\Ga_b \c \square_\g P,\\
\square_\g(q^2\Z P) &= W q^2 \Z P+ r\Ga_b\c  \dk^2 P +r^2\dk \err[\square_\g P]+r^2 \dk^{\leq 1} \big(\Ga_g \c \dk^{\leq 1} P\big)+r^3\Ga_b \c \square_\g P,
\end{split}
\eea
where the potential  $W$, given by
\beaa
W=\frac{3}{2}\Big[ \ov{\tr\Xb} \tr X+2P  -2 \Hb \c  \ov{H}  \Big]+q^{-2}\square_\g(q^2),
\eeaa
 is complex, and satisfies
\beaa
\Re(W)=W_{Sch}+O( a r^{-4} ), \qquad  \Im(W)=O(a r^{-3}).
\eeaa
where $W_{Sch}=O(mr^{-3})$ is the respective real potential in Schwarzschild. 
\end{lemma}

\begin{proof}
We commute the  wave equation for $P$ \eqref{eq:WaveEq-forP}
\beaa
\square_\g P &=  \tr X\nab_3P +\ov{\tr\Xb}\nab_4P  -\ov{H}\c\DD P  - \Hb\c\ov{\DD}P + V P  +F,
\eeaa
with $\T$. Using \eqref{eq:commTZsquare} and, as a consequence of Lemma \ref{LEMMA:DEFORMATION-TENSORS-T},
\beaa
[\T, e_3]=[\T, e_4]=[\T, e_a]=\Ga_b \c \dk
\eeaa
we deduce
\beaa
\square_\g(\T P)&=&   \T \square_\g P + \dk \big(\Ga_g \c \dk P\big)+\Ga_b \c \square_\g P\\
&=&  \tr X\nab_3(\T P) +\ov{\tr\Xb}\nab_4(\T P)  -\ov{H}\c\DD( \T P)  - \Hb\c\ov{\DD}(\T P) +V \T P +\T F\\
&+&  \T (\tr X ) \nab_3P + \T( \ov{\tr\Xb}) \nab_4P  - \T (\ov{H})\c\DD P  -\T( \Hb) \c\ov{\DD}P +\T(V)   P  \\
&+& \tr X [\T, \nab_3] P +\ov{\tr\Xb} [\T, \nab_4] P+ \dk \big(\Ga_g \c \dk P\big)+\Ga_b \c \square_\g P\\
&=&  \tr X\nab_3(\T P) +\ov{\tr\Xb}\nab_4(\T P)  -\ov{H}\c\DD( \T P)  - \Hb\c\ov{\DD}(\T P) +V \T P\\
&& +\dk \err[\square_\g P]+ \dk^{\leq 1} \big(\Ga_g \c \dk^{\leq 1} P\big)+\Ga_b \c \square_\g P.
\eeaa

To the above we can apply Lemma \ref{Le:squareq^2Psi} and deduce
  \beaa
   \square_\g ( q^2 \T P  )&=&W q^2 \T P+ r\Ga_b\c  \dk^2 P +r^2\dk \err[\square_\g P]+r^2 \dk^{\leq 1} \big(\Ga_g \c \dk^{\leq 1} P\big)+r^2\Ga_b \c \square_\g P
  \eeaa
  where 
  \beaa
  W=V+q^{-2}\square_\g(q^2).
  \eeaa
The   second   equation in \eqref{waveeq:TZ(P)}    can be derived in the same manner.
\end{proof}


\section{An identity for $\DDc\hot \DDc P$}
\lab{section:Otheridentites.Chap5}


The following identity will be used in Chapter \ref{Chapter:EN-MorforPc}.
 \begin{proposition}\label{PROP:RELATION-QF-P}
  \lab{prop:Formula-qf-DDhotDDPc-ch5}
  The following relation holds true:
  \bea
\begin{split}
& \nabc_4 \nabc_4\Ab +2\tr X  \nabc_4\Ab  +\frac 12(\tr X)^2 \Ab\\
&=\frac 1 2 \DDc\hot \DDc P+\big( 4 \DDc P   + 6P  \Hb + \tr X \Bb  \big)\hot \Hb   +\frac 3 2 P \big(\tr X\, \widehat{\Xb}  +\ov{\tr\Xb} \widehat{X} \big) \\
&+\Xi \hot \nabc_3 \Bb+r^{-1} \dk^{\leq 1} ( \Xi \c \Ab) + r^{-1} \dk^{\leq 1} (\Ga_g  \c \Bb)+ r^{-1} \dk^{\leq 1} (\Ga_b \c B ).
\end{split}
\eea
  From the above, we deduce the following relation between $\qfb$ and $\Pc$
  \bea
  \begin{split}
\qfb &= \frac 1 2\ov{q} q^3 \DDc\hot \DDc \Pc + \dkb^{\leq 1}( \Ga_b, r \Ga_g)+ O(ar) \dk^{\leq 1} \Bb + O(a^2) \Ab+O(ar) \dk^{\leq 1} \Pc  \\
&+ r^{3} \dk^{\leq 1}( \Ga_b \c (\Pc, B)) + r^{2} \dk^{\leq 1} (\Ga_g  \c (r\Bb, \Ab))  +r^4\Xi \hot \nabc_3 \Bb+r^{3} \dk^{\leq 1} ( \Xi \c \Ab)
\end{split}
\eea
 where the linear term $\dkb^{\leq 1}( \Ga_b, r \Ga_g)$ does not contain $\Xib$ or $\widecheck{\tr\Xb}$. 
  \end{proposition}
  
  \begin{proof} 
  See Appendix \ref{proof-relation-qf-P}.  
  \end{proof}


\part{Analysis of the wave equations}



\chapter{Estimates  for the model gRW equation in perturbations of Kerr}
\label{chapter-estimates-RW-model}



\section{Preliminaries}
\lab{sec:pfdoisdvhdifuhgiwhdniwbvoubwuyf}


In this chapter we introduce the model problem for the full generalized Regge-Wheeler equation in perturbations of Kerr obtained in Proposition \ref{prop:eq-real-gen-RW}.  The model problem consists in the following gRW equation  for a real tensor  $\psi\in \mathfrak{s}_2$ in the spacetime $\MM$ of section \ref{sec:setupandlinearizedquantities}:
\bea\lab{eq:Gen.RW}
\squared_2 \psi -V\psi=- \frac{4 a\cos\th}{|q|^2}\dual \nab_T  \psi+N, \qquad V= \frac{4\De}{ (r^2+a^2) |q|^2},
\eea
for some right-hand side $N$. Here $\squared_2$ denotes the D'Alembertian operator for horizontal 2-tensors in $\MM$.

\begin{remark}
We note that in applications   to the gRW equation for $\qf$,   $N$ contains   linear terms, in particular   
$\Re(L_{\qf}[A])$, as well as quadratic  $\Re(\err[\square_2 \qf])$, see  Proposition \ref{prop:eq-real-gen-RW}.
\end{remark}

   
     \subsection{The spacetime $\MM$}
     \lab{section:SpacetimeMM-chap6}


We consider a given vacuum spacetime $\MM$ satisfying the properties in section \ref{sec:setupandlinearizedquantities}:   
\begin{itemize}
\item $\MM$ comes together with a null pair $(e_4, e_3)$ and its corresponding horizontal structure as in section \ref{sec:nullparisandhorizontalstruct}. 

\item $\MM$ is endowed with a pair of constants $(a, m)$.

\item $\MM$ is endowed with a pair of scalar functions $(r, \th)$. 

\item The complex valued scalar function $q$ is defined as
\beaa
q := r+i\cos\th. 
\eeaa

\item $\MM$ is endowed with a complex horizontal 1-form $\Jk$. 
\end{itemize}

In addition, we assume:
\begin{enumerate}
\item $\MM$ is also endowed with a scalar function $\tau$ whose level sets $\Si(\tau)$ are spacelike. $\tau\in[1,\tau_*]$ on $\MM$ for some arbitrary large constant $\tau_*$. 

\item The boundary of $\MM$ is given by 
\bea
\pr\MM=\AA\cup\Si_*\cup\Si(1)\cup\Si(\tau_*)
\eea
where 
\bea
\AA:=\Big\{r=r_+-\deh, \, 1\leq\tau\leq\tau_*\Big\}, 
\eea
and $\Si_*$ is a spacelike hypersurface on which $\tau$ takes the values $[1,\tau_*]$ and $r\geq r_*$ with $r_*\gg \tau_*$. 

\item Let $r_0$ a large enough fixed constant. We decompose $\MM$ as follows
\bea
\Mint:=\MM\cap\{r\leq r_0\}, \qquad \Mext:=\MM\cap\{r\geq r_0\}. 
\eea
\end{enumerate}


   \subsection{Admissible perturbations of Kerr}
   \lab{sec:controlofGabandGagfromBA}


Recall that $\MM$ comes together three scalar functions $(r, \th, \tau)$, and with a null pair $(e_4, e_3)$ and its corresponding horizontal structure as in section \ref{sec:nullparisandhorizontalstruct}. Then:
\begin{itemize}
\item We use the complexified Ricci and curvature coefficients  of Definition \ref{def:complexRicciandcurvaturecoefficients}. 

\item We define  the linearized quantities corresponding to these complexified coefficients as in Definition \ref{def:renormalizationofallnonsmallquantitiesinPGstructurebyKerrvalue:1} and 
\ref{def:renormalizationofallnonsmallquantitiesinPGstructurebyKerrvalue:3}, i.e. we consider that the   normalization of $(e_3, e_4)$ is ingoing.

\item With respect to these linearized quantities, we the notations $\Ga_g$ and $\Ga_b$ for error terms are given by Definition \ref{definition.Ga_gGa_b}.
\end{itemize}

\begin{remark}
We define the following linearized quantity 
\beaa
\widecheck{|\Re(\Jk)|^2}:=|\Re(\Jk)|^2 -\frac{(\sin\th)^2}{|q|^2}.
\eeaa
Note that $\widecheck{|\Re(\Jk)|^2}=0$ in Kerr. We assume in addition  that 
$\dk^{\leq 1}(\widecheck{|\Re(\Jk)|^2})\in r^{-2}\Ga_b$.
\end{remark}

We this definition of $\Ga_g$ and $\Ga_b$, we can now state our main assumptions on $\MM$. Let $\kl$ a large enough integer. We make assumptions on decay and on boundedness on $(\Ga_b, \Ga_g)$.


   \subsubsection{Assumptions on $\MM$}
  

On $\MM$, we will prove energy Morawetz estimates. To this end, we introduce the scalar function 
$\tau_{trap}$ defined by
\bea
\tau_{trap} := \left\{\ba{lll}
1+\tau & \textrm{on} & \MM_{trap},\\
1& \textrm{on} & \Mntrap.
\ea\right.
\eea
Then, we assume that the linearized quantities satisfy the following estimates on $\MM$ 
\bea\lab{eq:assumptionsonMMforpartII}
\bsplit
r^3|\dk^{\leq k}\xi|+ r^2|\dk^{\leq k}\Ga_g|+r|\dk^{\leq k}\Ga_b| &\leq \ep, \qquad\qquad\,\,\, k\leq k_L,\\
r^3|\dk^{\leq k}\xi|+ r^2|\dk^{\leq k}\Ga_g|+r|\dk^{\leq k}\Ga_b| &\leq \frac{\ep}{\tau_{trap}^{1+\dec}}, \qquad k\leq \frac{k_L}{2}.
\end{split}
\eea

\begin{remark}
In this section $k_L$  is an unspecified  large positive  integer.
The  bounds  \eqref{eq:assumptionsonMMforpartII} will be assumed in all results and proofs of Chapters\footnote{Chapters 7 and 8 concern proofs in Kerr and do thus not require assumptions on $(\Ga_g, \Ga_b)$.} 6, 9 and 10. 
  In applications   to chapters  11 and 12 we will specify $k_L$ and make additional assumptions.
\end{remark}
 
\begin{remark}\lab{rmk:xiisactuallybetterthanGag}
Note that the assumptions for $\xi$ in \eqref{eq:assumptionsonMMforpartII} are consistent with $\xi\in r^{-1}\Ga_g$, while $\xi$ is a priori only in $\Ga_g$ according to Definition \ref{definition.Ga_gGa_b}. These  stronger assumptions \eqref{eq:assumptionsonMMforpartII} for $\xi$ will always hold whenever we apply the results  of Chapter \ref{chapter-estimates-RW-model}:
\begin{itemize}
\item In Chapter \ref{chapter-full-RWforqf}, this follows from the assumptions \eqref{eq:GlobalFrame-HcinGa_g2} and the fact that $\xi\in\Ga_g$. 
\item In Chapter \ref{chapter-full-RWforqfb}, this follows from the assumptions \eqref{eq:Xi-Hb-chapter12}. 
\item In Part III, this follows from the assumptions \eqref{eq:specialidentityforthegloablframeofMMinpartIII}.
\end{itemize}
\end{remark}

   
     \subsection{Basic properties of the $\tau$ function}
     \lab{sec:basicpropertiestau}


   
     \subsubsection{Choice of $\tau$}
     

Recall that $\MM$ is also endowed with a scalar function $\tau$ whose level sets $\Si(\tau)$ are spacelike. We provide in this section the basic properties of the $\tau$ function that will be used later.  Recall that, given a time function $\tau$, the vectorfield  $=-\g^{\a\b}\pr_\b\tau \pr_\a$ is timelike future oriented.  Given a level hypersurface  $\Si=\Si(\tau)$, we denote  
  \beaa
  N_\Si :=-\g^{\a\b}\pr_\b\tau \pr_\a.
  \eeaa

\begin{definition}[Choice of $\tau$]
\lab{definition:definition-oftau}
Let  $\de_\HH>0$ small enough. We choose the smooth scalar function $\tau$ on $r\geq r_+(1-\deh)$  such that  we have on $r\geq r_+(1-\deh)$
\beaa
\g(N_{\Si}, N_{\Si}) \leq -\frac{m^2}{8r^2}, \qquad e_4(\tau)>0, \qquad e_3(\tau)>0, \qquad |\nab\tau|^2  \leq \frac{8}{9}e_4(\tau)e_3(\tau).
\eeaa
In addition, we have the following asymptotic behavior  for $r$ large
\beaa
\frac{m^2}{r^2}\les e_4(\tau)\les \frac{m^2}{r^2}, \qquad 1\les e_3(\tau)\les 1.
\eeaa
Finally, we assume on $\MM$ 
\beaa
\T(\tau)=1+ r\Ga_b, \qquad \nab(\tau)=a \Re( \Jk)+\Ga_b,
\eeaa
where $\T$ is the vectorfield introduced below in Definition \ref{def:TandZiningoing:chap6asoihafeoigha}. 
\end{definition}

\begin{remark}
We refer to  the statement and proof of Proposition 9.3.5 in \cite{KS:Kerr} for an explicit example of function $\tau$ verifying the above properties. 
\end{remark}

   
 \subsubsection{Coordinates systems on $S(\tau, r)$}


We assume that the spheres $S(\tau, r)$ are covered by three coordinates systems.
\begin{definition}[Coordinates systems on $S(\tau, r)$]\lab{def:coordinatessystemSoftauandr:chap6}
On each $S(\tau, r)$, let 
\begin{itemize}
\item $(x^1_S, x^2_S)$ a coordinates system defined on $\frac{2\pi}{3}<\theta\leq\pi$, 

\item $(x^1_E, x^2_E)$ a coordinates system defined on $\frac{\pi}{4}<\theta\leq\frac{3\pi}{4}$, 

\item $(x^1_N, x^2_N)$ a coordinates system defined on $0\leq\theta<\pi$, 
\end{itemize}
so that we have the following control on each corresponding coordinate chart 
\beaa
\max_{b,c=1,2}|\pr^{\leq 2}(g_{bc} - (g_{a,m})_{bc})| &\les& r^2\ep,
\eeaa
where $g_{bc}$ denotes the induced metric coefficients in these coordinates systems, $(g_{a,m})_{bc}$ the corresponding expression in Kerr, and $\pr^{\leq 2}$ at most two coordinates derivatives. 
\end{definition}

\begin{remark}\lab{rmk:whatthesecoordinatessystemsonSoftaurareinKerrandrefmainKerr:chap6}
In \cite{KS:Kerr}, these coordinates systems are constructed using the scalar function $\th$ and an auxiliary scalar function $\vphi$ as follows
\beaa
(x^1_S, x^2_S)=(x^1_N, x^2_N)=(\sin\th\cos\vphi, \sin\th\sin\vphi), \qquad  (x^1_E, x^2_E)=(\th, \vphi).
\eeaa
We refer to Propositions 4.1 and 4.2 of \cite{KS:Kerr} for the control of such coordinates systems. Also, see Lemma 2.4.10 in \cite{KS:Kerr} for the form of $(g_{a,m})_{bc}$ in the $(x^1_E, x^2_E)$ coordinates system, and see Lemma 2.4.24 in \cite{KS:Kerr} for the form of $(g_{a,m})_{bc}$ in the $(x^1_S, x^2_S)$ and $(x^1_N, x^2_N)$ coordinates systems. 
\end{remark}

   
\subsection{Regions of integration and basic vectorfields}


   
\subsubsection{Regions of integration}

  
  Recall the time function $\tau$ introduced in Definition \ref{definition:definition-oftau}. 
  We denote by $\Sigma_\tau$  the level sets of the  function $\tau$.

     \begin{definition}\label{def:causalregions}
    We define the following regions of $\MM$.
    \begin{enumerate}
      \item   We define the trapping region of $\MM$ to be the set 
    \bea\lab{eq:def-MM-trap}
     \MM_{trap}(\de_{trap})=\MM\cap\left\{     \frac{|\TT|}{r^3}  \le \de_{trap} \right\}, \qquad \delta_{trap} = \frac{1}{10},
    \eea
    where  $\TT$ is  the   polynomial in $r$ defined in \eqref{definition-TT}, i.e.
    \beaa
    \TT=r^3-3mr^2+ a^2r+ma^2.
    \eeaa
  
  \item We denote $\Mntrap$ the complement to the trapping region $\MM_{trap}$. 
 
 \item We denote $\MM_{red}$ the be the region 
 \bea
 \MM_{red}:=\MM\cap\big\{r\leq r_+(1+2\de_{red})\big\}.
 \eea
 where the small enough constant $\de_{red}>0$ depends only on $m-|a|$ and is such that $\deh\ll \de_{red}^{20}$.
    
   \item We define the  domain $\MM(\tau_1,\tau_2)$ to be the  region of $\MM$ where  $\tau_1\le \tau\le \tau_2$, where $\tau$ is the time function defined in Definition \ref{definition:definition-oftau}.
      \end{enumerate}
     \end{definition}

   
\subsubsection{Basic vectorfields}


We start with the definition of $\T$ and $\Z$ as in section \ref{symmetry-operators-pert}, taking into account  that the normalization of $(e_3, e_4)$ is ingoing. 
\begin{definition}\lab{def:TandZiningoing:chap6asoihafeoigha}
In $\MM$, we define $\T$ and $\Z$ as follows
\beaa
\T &:=& \frac{1}{2}\left(e_4+\frac{\Delta}{|q|^2}e_3 -2a\Re(\Jk)^be_b\right),\\
\Z &:=& \frac 1 2 \left(2(r^2+a^2)\Re(\Jk)^be_b -a(\sin\th)^2 e_4 -\frac{a(\sin\th)^2\De}{ |q|^2} e_3\right).
\eeaa
\end{definition} 

\begin{remark}\lab{rmk:computationofgTTfortypeofvectorfieldT}
Note that we have
\beaa
\g(\T, \T) &=& -\frac{\De}{|q|^2}+a^2|\Re(\Jk)|^2 = -\frac{r^2+a^2(\cos\th)^2-2mr}{|q|^2}+a^2\left(|\Re(\Jk)|^2-\frac{(\sin\th)^2}{|q|^2}\right)\\
&=&  -\frac{r^2+a^2(\cos\th)^2-2mr}{|q|^2}+r^{-1}\Ga_b=-\frac{r^2+a^2(\cos\th)^2-2mr}{|q|^2}+O(r^{-2}\ep)
\eeaa
where we have used the definition of $\Ga_b$ and its control by \eqref{eq:assumptionsonMMforpartII}. 
In particular, we have the following non sharp estimate for any $|a|\leq m$
\beaa
\g(\T, \T)<0\quad\textrm{on}\quad r\geq \frac{5m}{2}
\eeaa
so that $\T$ is  timelike on the region $r\geq \frac{5m}{2}$ of $\MM$.  
\end{remark}

 \begin{lemma}\label{rem:de-trap}
 For  $|a|/m$ sufficiently small  and  $\delta_{trap} = \frac{1}{10}$,   the vectorfield  $\T$ is strictly timelike in $\MM_{trap}$.
    \end{lemma}
    
\begin{proof}
Observe that $\MM_{trap}=[\tilde{r}_-, \tilde{r}_+]$ where $\tilde{r}_\pm$ is the unique root to $\frac{\TT}{r^3}=\pm \de_{trap}$. By writing $a^2=\gamma m^2$, for $0 \leq \gamma \ll 1$, and defining $\tilde{x}_\pm=\frac{\tilde{r}_\pm}{m}$, $\tilde{x}_\pm$ is  the  unique root of the following equation
\beaa
(1\mp\de_{trap})(\tilde{x}_\pm)^3 -3(\tilde{x}_\pm)^2+\gamma \tilde{x}_\pm+\gamma =0.
\eeaa
Since  $0 \leq \gamma \ll 1$, we easily infer
\beaa
\tilde{x}_\pm = \frac{3}{1 \mp\de_{trap}}+O(\ga), \qquad \tilde{r}_\pm=\left(\frac{3}{1 \mp\de_{trap}}+O(\ga)\right)m.
\eeaa
In particular, fixing $\de_{trap}=\frac{1}{10}$, we have clearly $\tilde{r}_->\frac{5m}{2}$ for $\ga$ small enough and hence $\T$ is strictly timelike in $\MM_{trap}(\de_{trap})$ in view of Remark \ref{rmk:computationofgTTfortypeofvectorfieldT}. 
\end{proof} 

Also, we define the following vectorfields
 \bea
 \That &=&\frac 1 2 \left( \frac{|q|^2}{r^2+a^2} e_4+\frac{\De}{r^2+a^2}  e_3\right), \qquad \Rhat= \frac 1 2 \left( \frac{|q|^2}{r^2+a^2} e_4-\frac{\De}{r^2+a^2}  e_3\right).
 \eea

Finally, we introduce the vectorfield $\That_\de$ that will be used for energy estimates. 
 \begin{definition}
\lab{definition:vfppppp}\lab{definitionThat_de}
We define the  vectorfield 
\bea
\That_\de:= \T+\frac{a}{r^2+a^2} \chi_0\left( \de^{-1} \frac{\TT}{r^3} \right) \Z
\eea
with $\de=\de_{trap}$ and with $\chi_0$ the smooth  bump function
 \bea
 \lab{definition:chi-T_chippppp}
 \chi_0(x)=\begin{cases}
 &0  \qquad \qquad  \mbox{if}   \quad  |x| \le  1, \\
 &1  \qquad  \qquad \mbox{if}   \quad   |x|\ge 2. 
 \end{cases} 
 \eea
 We also write
 \beaa
 \That_\de:= \T+\chi_\de \Z, \qquad \chi_\de:= \frac{a}{r^2+a^2} \chi_0\left( \de^{-1} \frac{\TT}{r^3} \right).
 \eeaa
\end{definition}

    
    \subsection{Main norms}
        \lab{subsection:basicnormsforpsi}
        

   We introduce in this section the main norms needed to state the main results  of this chapter   concerning combined Energy-Morawetz     and $r^p$-weighted estimates for solutions to \eqref{eq:Gen.RW}:

   {\bf 1. Reduced  basic Morawetz  norms.} 
       \bea
\bsplit
 \Mor[\psi](\tau_1, \tau_2)&:=\int_{\MM(\tau_1, \tau_2) } 
      r^{-2}  | \nab_\Rhat  \psi|^2 +r^{-3}|\psi|^2\\
      &+ \int_{\Mntrap(\tau_1, \tau_2)} \left(  r^{-2}|\nab_3\psi|^2 + r^{-1}  |\nab  \psi|^2\right),\\
     \Morr[\psi](\tau_1, \tau_2)&:= \Mor[\psi](\tau_1, \tau_2)  + \int_{\MM_{r\geq 4m}(\tau_1,\tau_2)}  r^{-1-\de} |\nab_3 \psi|^2.
   \end{split}
\eea

{\bf 2. Basic Energy norm.} 
\bea
 E[\psi](\tau) :=\int_{\Si (\tau)} \Big( |\nab_4\psi|^2 +   r^{-2}|\nab_3\psi|^2 +|\nab\psi|^2 + r^{-2}|\psi|^2\Big).
\eea

{\bf 3. Basic Flux  norm.} 
\bea
\bsplit
F[\psi](\tau_1,\tau_2) : =& F_{\AA}[\psi](\tau_1,\tau_2)+F_{\Si_*}[\psi](\tau_1,\tau_2),\\
F_{\AA}[\psi](\tau_1,\tau_2) := & \int_{\AA(\tau_1, \tau_2)}\Big( |\nab_4\psi |^2+|\nab_3\psi|^2+|\nab\psi|^2+r^{-2} | \psi |^2\Big), \\
F_{\Si_*}[\psi](\tau_1,\tau_2) := & \int_{\Si_*(\tau_1, \tau_2)}\Big( |\nab_4\psi |^2+|\nab_3\psi|^2+|\nab\psi|^2+ r^{-2} | \psi |^2\Big).
\end{split}
\eea

{\bf 4. Basic  $N$- norm.} 
\bea
\bsplit
\NN[\psi,  N](\tau_1, \tau_2) :=&\int_{\MM(\tau_1, \tau_2)}\big(|\nab_{\Rhat} \psi|+r^{-1}|\psi|\big) |N|+\left|\int_{\Mtrap} \nab_{\That_\de} \psi \c N\right|+\int_{\Mntrap} |D \psi||N| \\
&+\int_{\MM(\tau_1, \tau_2)}|N|^2+\sup_{\tau\in[\tau_1, \tau_2]}\int_{\Si(\tau)}|N|^2+\int_{\Si_*(\tau_1, \tau_2)}|N|^2.
\end{split}
\eea

{\bf 5. Weighted  bulk norm.} For $0<p<2$, we define
\bea
B_{p}[\psi](\tau_1, \tau_2)&:=& \Morr[\psi](\tau_1, \tau_2) +\int_{\MM_{r\geq 4m}(\tau_1,\tau_2)} r^{p-3}\Big(|\dk \psi|^2 +|\psi|^2  \Big).
\eea

{\bf 6. Weighted  energy  norm.} For $0<p<2$, we define
\bea
E_p[\psi](\tau) := \left\{\ba{l} 
E[\psi]+\displaystyle\int_{\Si_{r\geq 4m}(\tau)} r^{p}\Big(|\nab_4\psi|^2  + r^{-2}|\psi|^2  \Big) \quad\textrm{for }p\leq 1-\de,\\[3mm]
E[\psi]+\displaystyle\int_{\Si_{r\geq 4m}(\tau)} r^{p}\Big(|r^{-1}\nab_4(r\psi)|^2  + r^{-p-1-\de}|\psi|^2  \Big) \quad\textrm{for }p> 1-\de.
\ea\right.
\eea

\begin{remark}
By a slight abuse of notation,  we will often identify  $E_p[\psi](\tau_2)$ with 
$\sup_{\tau\in[\tau_1,\tau_2] } E_p[\psi](\tau)$. 
\end{remark}

{\bf 7. Weighted flux   norm.} For $0<p<2$, we define
\bea
\nn F_p[\psi](\tau_1,\tau_2):&=& F[\psi](\tau_1,\tau_2)\\
&& +\int_{\Si_*(\tau_1, \tau_2)} r^p\Big(|\nab_4 \psi|^2 +|\nab\psi|^2+ r^{-2}|\psi|^2 \Big).
\eea

{\bf 8. Combined norms.}  We denote the combined norm
\bea
\BEF_p[\psi](\tau_1,\tau_2):= \sup_{\tau\in[\tau_1, \tau_2]} E_p[\psi](\tau) +B_p[\psi](\tau_1, \tau_2) +F_p[\psi](\tau_1, \tau_2).
\eea
and 
\bea
\EF_p[\psi](\tau_1,\tau_2):= \sup_{\tau\in[\tau_1, \tau_2]} E_p[\psi](\tau)  +F_p[\psi](\tau_1, \tau_2).
\eea

{\bf 9. Weighted   $N$- norm.} For $0<p<2$, we define
\bea
\NN_p[\psi,  N](\tau_1, \tau_2) &=& \NN[\psi,  N](\tau_1, \tau_2)+\left| \int_{\MM_{r\geq 4m}}  r^{p-1}  \, \nab_4 (r\psi ) \c  N\right|.
\eea

{\bf 10. $\Mext$ norms.} We denote by $\Bext_p, \Eext_p, \Next_p$ the restrictions of the norms 
$B_p, E_p, \NN_p$  to $\Mext$, i.e. the region in $\MM$ where $r\ge r_0$.

\bigskip

\medskip

{\bf 11. Higher order norms.}  We define the higher  derivative norms $\Mor^s[\psi]$, $E^s[\psi]$, $F^s[\psi]$,    $\NN^s[\psi,  N]$, $B_p^s[\psi]$, $ E_p^s[\psi]$,  $F_p^s[\psi]$, $\NN_p^s[\psi]$, 
 by the  general procedure for a norm $Q[\psi]$, i.e. 
 \beaa
 Q^s[\psi]=\sum_{k\le s} Q[\dk^k\psi].
 \eeaa

    \bigskip

\begin{remark}
The $\int_{\MM}\big(|\nab_{\Rhat} \psi|+r^{-1}|\psi|\big) |N|$ part in the definition of  of $\NN[\psi, N]$ is obtained in the  proof of the Morawetz estimate while  the $\Big|\int_{\MM} \nab_{\That_\de} \psi \c N\Big|$ part is needed for the energy estimate. We stress here the presence of the term involving $\nab_{\That_\de}$, as its specific form will be needed in Chapter \ref{chapter-full-RW} to treat the full Regge-Wheeler equation. 
\end{remark}

   
 \section{Main theorems for the model gRW}
   
   
    We start by stating the main results  of this chapter   concerning combined Energy-Morawetz   
     and $r^p$-weighted estimates for solutions to \eqref{eq:Gen.RW} on a spacetime $\MM$ which is 
     an admissible perturbation of Kerr in the sense that \eqref{eq:assumptionsonMMforpartII} holds.

       \begin{theorem}[Basic $r^p$-weighted estimates]
       \lab{THEOREM:GENRW1-P}
       The following estimates hold true  for solutions $\psi\in\sk_2$  of  \eqref{eq:Gen.RW} on $\MM$, for all $\de\le p\le  2-\de$ and  $2\leq s\le \kl$,
       \bea
       \lab{eqtheorem:GenRW1-p}
       \sup_{\tau\in[\tau_1, \tau_2]} E_p^s[\psi](\tau) +B_p^s[\psi](\tau_1, \tau_2) +F_p^s[\psi](\tau_1, \tau_2) \les
       E_p^s[\psi](\tau_1)+\NN_p^s[\psi, N](\tau_1, \tau_2).
       \eea
       \end{theorem}
       
       To state the  second theorem we  need to introduce  
        the quantity 
        \bea
        \psiwc :=  f_2 \left(e_4\psi+\frac{r}{|q|^2}\psi\right).
        \eea
        with $f_2= r^2$  for $r\ge R$ and $f_2= 0$ for $r\le R/2$.
      
       \begin{theorem}[Improved $r^p$-weighted estimates]
       \lab{THEOREM:GENRW2-Q}
       The following  estimates hold true   for  the quantity $\psiwc$ 
       for solutions $\psi\in\sk_2$  of  \eqref{eq:Gen.RW} on $\MM$,
        for all $-1+\de<q\le 1-\de$, $s\le\kl -1$,
         \bea
         \lab{eqtheorem:GenRW2-q}
         \BEF^s_q[\psiwc](\tau_1, \tau_2)
       \les
       \Et_q^s[\psiwc](\tau_1)+\NNt_q^s[\psiwc, N](\tau_1, \tau_2) + \NN_{\max\{q,\de\}}^{s+1}[\psi, N](\tau_1, \tau_2).
       \eea
         where\footnote{Recall that
          $\BEF_q^s[\psiwc](\tau_1,\tau_2)= \sup_{\tau\in[\tau_1, \tau_2]} E_q^s[\psiwc](\tau) +B_q^s[\psiwc](\tau_1, \tau_2) +F_q^s[\psiwc](\tau_1, \tau_2)$.} the norms on the right are given by
               \bea
       \lab{eq:TildeNormsN-chap6}
        \Et_q^s[\psiwc](\tau)=E_q^s[\psiwc](\tau)+ E^{s+1}_{\max\{q,\de\} }[\psi](\tau)
         \eea
         and
         \bea
          \NNt_q^s[\psiwc, N](\tau_1, \tau_2)&=& \left|\int_{\MM_{\ge R}(\tau_1, \tau_2)} r^{q+2} \dk^{\le s} \psiwc \c \left(\nab_4 \dk^{\le s }N+ \frac 3 r \dk^{\le s} N\right)\right|. 
         \eea        
        \end{theorem}

\begin{remark}
These results are the analog in perturbations of Kerr to Theorem 5.17 and  Theorem 5.18.  in \cite{KS} for perturbations of Schwarzschild.
\end{remark}
      
Theorems \ref{THEOREM:GENRW1-P} and \ref{THEOREM:GENRW2-Q} will be proved in Chapter \ref{chapter-rp-estimates} by relying on $r^p$-weighted estimates, and Morawetz-Energy estimates derived in Chapter \ref{chapter-perturbations-Kerr}. In the next section, we discuss these Morawetz-Energy estimates.

    
\section{Main Morawetz-Energy results}
\lab{section:Main-results-forRWmodelproblem}


The following theorem is our main Morawetz-Energy result  for solutions to \eqref{eq:Gen.RW} on a spacetime $\MM$ which is   an admissible perturbation of Kerr in the sense that \eqref{eq:assumptionsonMMforpartII} holds. 

\begin{theorem}[Morawetz-Energy]
\lab{THM:HIGHERDERIVS-MORAWETZ-CHP3}
Let $\psi$  an $\sk_2$  solution of \eqref{eq:Gen.RW} in $\MM$. For $|a|/ m \ll 1$ sufficiently small, we have, for all   $2\leq s\le \kl$, and for any $\de>0$, 
\bea
\nn&& \Mor[(\nab_3, \nab_4, \dkb)^{\leq s}\psi](\tau_1, \tau_2) + E[(\nab_3, \nab_4, \dkb)^{\leq s}\psi](\tau_2) +F[(\nab_3, \nab_4, \dkb)^{\leq s}\psi](\tau_1,\tau_2)  \\
\nn&\les & E^s[\psi](\tau_1)  +\NN^s [\psi, N] (\tau_1, \tau_2)\\
&&+\left(\frac{|a|}{m}+\ep\right)\left(\sup_{\tau\in [\tau_1, \tau_2]}E^s[\psi]+B_\de^s[\psi](\tau_1, \tau_2)+F^s[\psi](\tau_1, \tau_2)\right).
\eea
\end{theorem}

\begin{remark}
This result is the analog in perturbations of Kerr to Theorem 10.1 of \cite{KS} for perturbations of Schwarzschild.
\end{remark}

Theorem \ref{THM:HIGHERDERIVS-MORAWETZ-CHP3} will be proved in section \ref{sec:proofofThm:HigherDerivs-Morawetz-chp3:chap9}. In the rest of this section, we introduce norms needed to state intermediary Morawetz-Energy estimates. Then we state these results while providing  the outline of the proof of Theorem \ref{THM:HIGHERDERIVS-MORAWETZ-CHP3}.


 \subsection{Additional energy flux and bulk quantities}\label{section:energy-flux-ch3}
 \lab{sec:defintionadditionalfrluxandbulkforaxandSSmorawetz:chap6}


In this section, we introduce additional energy flux and bulk quantities that are needed for the proof of the Energy-Morawetz estimates. 

Together with the above gRW equation \eqref{eq:Gen.RW}, we consider the commuted gRW equations, given by 
  \bea
  \lab{eq:waveeqfor-psia-chp3}
  \square_2\psia -V\psia&=&- \frac{ 4 a\cos\th}{|q|^2} \dual \nab_T  \psia+N_\aund, \qquad V= \frac{4\De}{ (r^2+a^2) |q|^2},
  \eea
where  $\psia$ is defined by 
\bea
  \psia:= \SS_\aund \psi \quad \textrm{for}\quad\aund=1,2,3,4,
  \eea
  with $\SS_\aund$ denoting  the set of second order differential operators,   see  Definition \ref{definition:symmetry-tensors-pert},  
\beaa
\SS_1\psi&=& \nab_T \nab_T \psi, \\
\SS_2\psi&=& a \nab_T \nab_Z \psi,\\
\SS_3\psi&=& a^2 \nab_Z \nab_Z \psi,\\
\SS_4\psi&=& \OO (\psi),
\eeaa  
and where the right-hand sides $N_\aund$ can be explicitly computed from $N$ and $\psi$, see Lemma \ref{lemma-commuted-equ-Naund:perturbation}.


 \subsubsection{Pointwise notation}
 
 
\begin{definition}
  We introduce the following pointwise  notation for $\psi\in \sk_2$.
  \begin{enumerate}
 \item 
  We denote 
\beaa
 |\psi|^2_{\SS}&:=&\sum_{\aund=1}^4\big|\psia|^2.
\eeaa
\item 
Given a vectorfield $Y$ we denote
\beaa |\nab_Y\psi|^2_{\SS}&:=&\sum_{\aund=1}^4\big|\nab_Y\psia|^2.
\eeaa
  \end{enumerate}
  \end{definition}

    
\subsubsection{Degenerate energy norm} 


\begin{definition}[Degenerate energy norm]
\lab{def:degenerate-energy}\lab{def:SSenergy-flux}
We define the following degenerate energy   for $\psi\in \sk_2$ along $\Si(\tau)$:
\beaa
 E_{deg}[\psi](\tau) :=\int_{\Si(\tau)} \left( |\nab_4\psi|^2 +  \frac{|\De|}{r^4} |\nab_3\psi|^2 +|\nab\psi|^2 + r^{-2}|\psi|^2\right).
\eeaa
\end{definition}


\subsubsection{Refined Morawetz norms} 


    \begin{definition}[Refined Morawetz norms]\lab{definition:SS0Morawetznorms} We define the following Morawetz  norms for $\psi \in \sk_2$.
    \begin{enumerate}
    \item 
  The degenerate  axially symmetric Morawetz norms  in $\MM=\MM(\tau_1, \tau_2) $:
  \beaa
     \Mordot^{ax}_{deg}[\psi](\tau_1, \tau_2)&:=&\int_{\MM(\tau_1, \tau_2)} \frac{m}{r^2} |\nab_{\Rhat} \psi|^2 +\frac{\TT^2}{r^6} \left(\frac{m}{r^2} |\nab_\That \psi|^2 + r^{-1}|\nab\psi|^2\right), \\
     \Mor^{ax}_{deg}[\psi](\tau_1, \tau_2)&:=& \int_{\MM(\tau_1, \tau_2) } 
      \frac{m}{r^2} | \nab_\Rhat \psi|^2 +r^{-3}|\psi|^2+\frac{\TT^2}{r^6} \left(\frac{m}{r^2} |\nab_\That \psi|^2 + r^{-1}|\nab\psi|^2\right).
\eeaa

\item We also define the higher degenerate and non-degenerate Morawetz norms in $\MM=\MM(\tau_1, \tau_2) $, for a scalar function $z$:
 \beaa
 \Mordot_{\SSz, deg}[\psi](\tau_1, \tau_2)&:=& \int_{\MM} 
      \frac{m}{r^2}  | \nab_\Rhat \psi|_\SS^2+ r^{3} \Big(\big|\nab_\That \Psi_z \big |^2+r^2\big |\nab \Psi_z\big|^2\Big),\\
    \Mor_{\SSz, deg}[\psi](\tau_1, \tau_2)&:=& \int_{\MM} 
      \frac{m}{r^2}  | \nab_\Rhat \psi|_\SS^2+ r^{-3}|\psi|_{\SS}^2+ r^{3}\Big(\big|\nab_\That  \Psi_z \big |^2+r^2\big |\nab  \Psi_z\big|^2  \Big),
 \eeaa
 where\footnote{Note that $z$, $\RRa$ and $\De$ are functions depending only on $r$ so that $\pr_r$ in the formula for $\RRtp^\aund[z]$  simply  denotes  differentiation w.r.t. $r$.} 
 \beaa
 \Psi_z=\Psi_z[\psi]&:=& \RRtp^\aund[z] \psia, \qquad   \RRtp^\aund[z] := \pr_r\left( \frac{z}{\De} \RRa\right),
 \eeaa
with $z$ a suitable function of $r$ to be chosen later, and  with the scalar functions $\RRa$ given by \eqref{components-RR-aund}, i.e. 
         \beaa
          \RR^1&=&-(r^2+a^2)^2, \qquad \RR^2 = -2(r^2+a^2), \qquad \RR^3 =-1, \qquad \RR^4=\De.
          \eeaa
\end{enumerate}
  \end{definition}

 \begin{remark} 
 Observe the following:
\begin{enumerate} 
\item The axially symmetric norms $\Mordot^{ax}$ and $\Mor^{ax}$ have trapped $\nab_\That$ and $\nab$ derivatives at $\TT=0$, which describes the trapping region for axially symmetric solutions. For general solutions, those norms cannot be bounded by the initial energy. Nevertheless those Morawetz energies will be used in Part 1 of our proof.

\item The higher norms $ \Mordot_{\SSz}$ and $\Mor_{\SSz}$ are positive definite norms where the trapping properties are encoded in the term $\Psi_z$ defined above. For our choice $z=z_0-\de_0 z_0^2$, for $z_0=\frac{\De}{(r^2+a^2)^2}$, the term $\Psi_z$ is given by, see \eqref{definition-Psiz},
\beaa
\Psi_z = -\frac{2\TT}{(r^2+a^2)^3}   \big(\de_0  \SS_1\psi+ (1+O(r^{-2} \de_0)) \OO\psi\big)+ \frac{4ar}{(r^2+a^2)^2}  \nab_\That \nab_Z \psi   \big(1+O(r^{-2} \de_0) \big).
\eeaa
The overall expression of $\Psi_z$ describes the trapping structure of general solutions.

\item The norms $\Mordot_{\SSz}$ and $\Mor_{\SSz}$  involve a sum of positive terms, each of which is given in terms of a linear combination of derivatives of $\psi$. 
In order to bound them  by a sum of positive terms involving $\psi$ directly, one needs to have a full degeneracy in the trapping region $\MM_{trap}$, i.e. for small $|a|/m$,
\beaa
&&\int_{\MM(\tau_1, \tau_2) } 
      \frac{m}{r^2} \big( | \nab_\Rhat \psi|_{\SS}^2 +r^{-2}|\psi|_{\SS}^2\big)+\int_{\Mntrap(\tau_1, \tau_2)}
     r^{-1}\Big( \frac{m}{r}|\nab_T\psi|_{\SS}^2+ |\nab \psi|_{\SS}^2 \Big)\\
     &\les& \sup_{\tau\in[\tau_1, \tau_2]}E_{deg}[(\nab_T, \dkb)^{\leq 2}\psi](\tau) + \Mor_{\SSz, deg}[\psi](\tau_1, \tau_2) +\deh F^2_\AA[\psi](\tau_1, \tau_2)\\
     &&+F_{\Si_*}[(\nab_T, \dkb)^{\leq 2}\psi](\tau_1, \tau_2)+\frac{|a|}{m}B^2_\de[\psi](\tau_1, \tau_2).
\eeaa
This will be achieved thanks to Lemma \ref{LEMMA:LOWERBOUNDPHIZOUTSIDEMTRAP}.
 \end{enumerate}
\end{remark}


\subsection{Outline of the proof of Theorem 
\ref{THM:HIGHERDERIVS-MORAWETZ-CHP3}}\label{section:outline-proofs}


In what follows, we give a description of the main steps in  the proof of Theorem 
\ref{THM:HIGHERDERIVS-MORAWETZ-CHP3}.

   
\subsubsection{Part 1: conditional Morawetz and Energy estimates in Kerr} 
\lab{section:Proofs-part1}


In the first part of the proof of the Energy-Morawetz estimates, we prove bounds for the first derivatives of the solution in what we call  \textit{conditional} Morawetz and energy estimates.  Such estimates are conditional as they  depend on the control of a derivative of the solution with respect to $Z$ on the right hand side of the estimate. In the case of axially symmetric solutions, those would become unconditional.

To ease the exposition, the proof of the conditional Morawetz and energy estimates are first derived  in the case of Kerr in Chapter \ref{chapter-proof-part1}, and will then be extended to perturbations of Kerr in section \ref{sec:proofofresultschapter7inperturbationofKerr}. 

\medskip

First, we derive the following basic degenerate, conditional, Morawetz  estimate, see Proposition \ref{proposition:Morawetz1-step1:perturbation} for the extension to perturbations of Kerr.  
\begin{proposition}
\lab{proposition:Morawetz1-step1}
The following  estimates  hold true in Kerr:
\begin{enumerate}
\item   For  all $|a|/m < 1$,  we have
\bea\lab{eq:conditional-mor-par1-1-I}
\bsplit
\Mordot^{ax}_{deg}[\psi](\tau_1, \tau_2) &\les  \int_{\pr\MM(\tau_1, \tau_2)}|M(\psi)| +\int_{\MM(\tau_1, \tau_2)}      \big(   a^2 r^{-4 } |\nab_Z\psi|^2 +r^{-3} |\psi|^2\big)\\
& +\int_{\MM(\tau_1, \tau_2)}\big(|\nab_{\Rhat} \psi|+r^{-1}|\psi|\big) |N|.
 \end{split}
\eea

\item For   $|a|/m \ll 1$ sufficiently small, we have
\bea\lab{eq:conditional-mor-par1-1-II}
\bsplit
\Mor^{ax}_{deg}[\psi](\tau_1, \tau_2) &\les \int_{\pr\MM(\tau_1, \tau_2)}|M(\psi)|+  \int_{\MM(\tau_1, \tau_2)}ar^{-1}\big(|\nab\psi|^2+r^{-2}|\nab_t\psi|^2\big) \\
&   +\int_{\MM(\tau_1, \tau_2)}\big(|\nab_{\Rhat} \psi|+r^{-1}|\psi|\big) |N|.
 \end{split}
\eea
\end{enumerate}
In both cases  $M(\psi)$   denotes a quadratic  expression in  $\psi$  and its first derivatives for which we  have 
\beaa
\int_{\pr\MM(\tau_1, \tau_2)}|M(\psi)| &\les & \sup_{[\tau_1, \tau_2]}E_{deg}[\psi](\tau)+\deh F_{\AA}[\psi](\tau_1, \tau_2)+F_{\Si_*}[\psi](\tau_1, \tau_2).
\eeaa
\end{proposition}

\begin{remark}
Observe that the above conditional estimates contain respectively the integrals
\beaa
 \int_{\MM(\tau_1, \tau_2)}  a^2   r^{-4}    |\nab_Z\psi|^2, \qquad  \int_{\MM(\tau_1, \tau_2)}ar^{-1}\big(|\nab\psi|^2+r^{-2}|\nab_t\psi|^2\big),
\eeaa
on the right hand side. These term can be absorbed only by considering higher derivative estimates, as described in Part 2. Also,  the difference in the interval for the parameter $|a|/m$ in the two estimates of Proposition \ref{proposition:Morawetz1-step1} is related to the control of the zero-th order term.
To be able to show that  this  term  comes with the right sign, and thus is only on  the left hand side,
 we make use of a Poincar\'e and a Hardy type  inequality, whose validity is restricted to small angular momentum.  
\end{remark}

We summarize here the main ideas of the proof of Proposition \ref{proposition:Morawetz1-step1}, which is obtained in Section \ref{section:cond:mor}. 

\begin{enumerate}
\item The Morawetz estimate is obtained by applying the vector field method with the vectorfield $X=\FF(r) \partial_r$ as multiplier. The choice of the function $\FF$ is inspired by \cite{A-B}, and the current relative to the vectorfield $X$ has the following form (see \eqref{first-low-bound-ref}):
 \beaa
|q|^2 \D^\mu\PP_\mu[X, w]&=& \AA |\nab_r \psi|^2+ \UU^{\a\b} (\Db_\a \psi) (\Db_\b \psi) + \VV |\psi|^2.
\eeaa
 The function $\FF$ is chosen so that the term $\UU^{\a\b} (\Db_\a \psi) (\Db_\b \psi)$ vanishes at $\TT=0$ and the coefficient $\AA$ of $\nab_r \psi$ is positive for all $|a|/m<1$. This is done in Section \ref{subsection:first-lower-bound}. 
\item By making use of the Lagrangian of the wave equation, the trapped term $\UU^{\a\b} (\Db_\a \psi) (\Db_\b \psi)$ can be upgraded to contain also the time derivative of $\psi$. This requires to absorb a term by $\AA$, which is still positive in the full sub-extremal range $|a|/m < 1$. This is done in Section \ref{section:lowerboundscontainigThat}. These steps imply estimate \eqref{eq:conditional-mor-par1-1-I}, as shown in Section \ref{section:Proof-proposition:Morawetz-Energy1}. 

\item To prove \eqref{eq:conditional-mor-par1-1-II}, we need to obtain positivity of the coefficient $\VV$ of $|\psi|^2$. This is obtained for sufficiently small $|a|/m$ through a combined Poincar\'e inequality (which extracts extra positivity from the trapped term) and a Hardy inequality (which extract extra positivity from the $\AA$ term). This is done in Section \ref{section:proof-hardy-poincare}, and completes the proof of Proposition \ref{proposition:Morawetz1-step1}. 
\end{enumerate}

Next, we derive a  conditional degenerate  energy estimate. 

We define the vectorfield $\That_\de$ such that $\That_\de=T$  at the trapped set and  $\That_\de=\That$ away from it, with $\de_{trap}=\frac{1}{10}$ so that, in view of Lemma    \ref{rem:de-trap}, $T$ is strictly timelike in $\MM_{trap}$. 
We then prove the following, see Proposition \ref{proposition:Energy1:perturbation} for the extension to perturbations of Kerr.
\begin{proposition}
\lab{proposition:Energy1}
The following  estimate  holds true in Kerr for solutions of    \eqref{eq:Gen.RW} in $\MM$, 
 for $|a|/m \ll 1$ sufficiently small,  
 \bea
\nn E_{deg}[\psi](\tau_2) +F_{\Si_*}[\psi](\tau_1, \tau_2) &\les&  E_{deg}[\psi](\tau_1) +\deh\big(E_{r\leq r_+(1+\deh)}(\tau_2)[\psi]+F_\AA[\psi](\tau_1, \tau_2)\big)\\
\nn&&+ \frac{|a|}{m}\Mor^{ax}_{deg}[\psi](\tau_1, \tau_2)   +\left|\int_{\MM(\tau_1, \tau_2)}  \nab_{\That_\de } \psi  \c N\right|\\
&& +\int_{\MM(\tau_1, \tau_2)}|N|^2.
\eea
\end{proposition} 

The proof of this proposition is obtained in Section \ref{subsection:energy-identity}, and is an application of the vectorfield method applied with the multiplier $\That_\de$.


\subsubsection{Part  2: $\SS$-derivatives Morawetz estimates in Kerr}
\lab{section:Higher derivative Morawetz  Estimates}


The  main limitation of the results of Proposition \ref{proposition:Morawetz1-step1} is the presence of $\nab_Z\psi$ on the right hand side of  the estimates.  To correct for this, we follow the approach of   Andersson and Blue  in \cite{A-B}    based on a remarkable extension of the classical  vectorfield method to   include   commutation with the second order  Carter operator. To ease the exposition, the proof of the $\SS$-derivatives Morawetz  estimates are first derived  in the case of Kerr in Chapter \ref{chapter-proof-mor-2}, and will then be extended to perturbations of Kerr in section \ref{sec:proofofresultschapter8inperturbationofKerr}. 

\medskip
 
 We have the following, see Proposition \ref{prop:morawetz-higher-order:perturbation} for the extension to perturbations of Kerr.
 \begin{proposition}[$\SS$-derivatives Morawetz estimates]
 \label{prop:morawetz-higher-order}
 Let the scalar function $z$ given by
\beaa
z=z_0-\de_0 z_0^2, \qquad z_0=\frac{\De}{(r^2+a^2)^2}.
\eeaa
Then, for $|a|/m \ll 1$ sufficiently small, the following  estimate  holds true for solutions of    \eqref{eq:Gen.RW} in Kerr:
 \bea\label{eq:mor-higher-unconditional-}
 \bsplit
  \Mor_{\SSz, deg}[\psi](\tau_1, \tau_2) \les& \int_{\pr\MM(\tau_1, \tau_2) }|M_\SS(\psi)| +\frac{|a|}{m}B^2_\de[\psi](\tau_1, \tau_2)\\
  &+\sum_{\aund=1}^4\int_{\MM(\tau_1, \tau_2)}\big(|\nab_{\Rhat} \psia|+r^{-1}|\psia|\big) |N_{\aund}|
  \end{split}
 \eea
where $M_\SS(\psi)$  denotes  an expression in  $\psi$ for which we  have a bound of the form
\beaa
&&\int_{\pr\MM(\tau_1, \tau_2) }|M_\SS(\psi)| \\
&\les &   \sum_{\aund=1}^4\left(\sup_{[\tau_1, \tau_2]}E_{deg}[\psi_\aund](\tau) +\deh F_{\AA}[\psi_{\aund}](\tau_1, \tau_2)+F_{\Si_*}[\psi_{\aund}](\tau_1, \tau_2)\right) \\
&+&\left(\sup_{[\tau_1, \tau_2]}E_{deg}[(\nab_T, \dkb)^{\leq 1}\psi](\tau)+\deh F_{\AA}[(\nab_T, \dkb)^{\leq 1}\psi](\tau_1, \tau_2)+F_{\Si_*}[(\nab_T, \dkb)^{\leq 1}\psi](\tau_1, \tau_2)\right)^{\frac{1}{2}}\\
&\times&\left(\sup_{[\tau_1, \tau_2]}E_{deg}[(\nab_T, \dkb)^{\leq 2}\psi](\tau)+\deh F_{\AA}[(\nab_T, \dkb)^{\leq 2}\psi](\tau_1, \tau_2)+F_{\Si_*}[(\nab_T, \dkb)^{\leq 2}\psi](\tau_1, \tau_2)\right)^{\frac{1}{2}}.
\eeaa 
 \end{proposition}

 We summarize here the main ideas in the proof of Proposition \ref{prop:morawetz-higher-order}, which is obtained in Section \ref{section:conditional-SS-values-Mor} and Section \ref{section:unconditional-SS-mOr}.
 
 \begin{enumerate}
\item The Morawetz estimate is obtained by applying the generalized vector field method with the double-indexed vectorfield $X^{\aund\bund}=\FF^{\aund\bund}(r) \partial_r$ as multiplier, to obtain the current, see \eqref{generalized-current-operator}, 
\beaa
|q|^2 \D^\mu\PP_\mu[X, w]&=&  \AA^{\aund\bund}  \nab_r\psia \nab_r\psib + \UU^{\a\b\aund\bund} \, \D_\a \psia \, \D_\b \psib  +\VV^{\aund\bund} \psia\psib.
  \eeaa
 The double-indexed function $\FF^{\aund\bund}$ is chosen so that the term $\UU^{\a\b\aund\bund} \, \D_\a \psia \, \D_\b \psib $ presents a quadratic expressions in higher derivatives and the coefficients $\AA^{\aund\bund}$ are positive for all $|a|/m<1$. 
 More precisely, by performing a crucial integration by parts (see Lemma \ref{lemma:IntegrationbypartsP}), we can write
\beaa
\UU^{\a\b\aund\bund} \, \D_\a \psia \, \D_\b \psib &=&\frac 1 2 h L^{\a\b}  \D_\a \Psi    \  \D_\b \Psi +\text{boundary terms,} \qquad \Psi:=\RRtp^{\aund}\psia,
\eeaa
for some constant coefficients $L^{\a\b}$. 
 This is done in Section \ref{section:preliminaries-chapter5}. 
 
\item By using another integration by parts lemma, see Lemma \ref{Lemma:integrationbypartsSS_3SS_4}, we prove that the term $ \AA^{\aund\bund}  \nab_r\psia \nab_r\psib $ is positive. This is done in Section \ref{section:conditional-SS-values-Mor}.

\item Finally, to prove \eqref{eq:mor-higher-unconditional-}, we apply a combined Poincar\'e inequality and Hardy inequality to obtain positivity for the term $\VV^{\aund\bund} \psia\psib$, for sufficiently small $|a|/m$. This is done in Section \ref{section:unconditional-SS-mOr}, and completes the proof of Proposition \ref{prop:morawetz-higher-order}. 
\end{enumerate}

Finally, to show that  $\Mor_{\SSz}$ controls $\psi$ in $\Mntrap$, we will rely on the following lemma proved in section \ref{section:lowerboundPhizoutsideMtrap}, see Lemma \ref{lemma:lowerboundPhizoutsideMtrap:perturbation} for the extension to perturbations of Kerr.
\begin{lemma}\lab{LEMMA:LOWERBOUNDPHIZOUTSIDEMTRAP}
For $\de_0>0$ small enough\footnote{Recall that the constant $\de_0>0$ is involved in the definition of $z=z_0-\de_0 z_0^2$.} and $|a|/m\ll \de_0$, there exists a universal constant $c_0>0$ such that the following holds on $\Mntrap$ in Kerr:
\beaa
r^3\Big(|\nab_T\Psi_z|^2+r^2|\nab\Psi_z|^2\Big)+r^{-3}|\psi|_{\SS}^2 &\geq& c_0r^{-3}\Big(|\nab_T\psi|^2_{\SS}+|\nab_Z\psi|^2_{\SS}+r^2|\nab\psi|^2_{\SS}\Big)\\
&& -O(ar^{-3})\big|(\nab_T, \dkb)^{\leq 1}\dk^{\leq 2}\psi\big|^2 +\Ddot_\a F^\a
\eeaa
where the 1-form $F$ denotes  an expression in  $\psi$ for which we  have a bound of the form
\beaa
&&\left|\int_{\pr\MM(\tau_1, \tau_2) }F^\mu N_\mu\right| \\
&\les &   \left(\sup_{[\tau_1, \tau_2]}E_{deg}[(\nab_T, \dkb)^{\leq 1}\psi](\tau)+\deh F_{\AA}[(\nab_T, \dkb)^{\leq 1}\psi](\tau_1, \tau_2)+F_{\Si_*}[(\nab_T, \dkb)^{\leq 1}\psi](\tau_1, \tau_2)\right)^{\frac{1}{2}}\\
&\times&\left(\sup_{[\tau_1, \tau_2]}E_{deg}[(\nab_T, \dkb)^{\leq 2}\psi](\tau)+\deh F_{\AA}[(\nab_T, \dkb)^{\leq 2}\psi](\tau_1, \tau_2)+F_{\Si_*}[(\nab_T, \dkb)^{\leq 2}\psi](\tau_1, \tau_2)\right)^{\frac{1}{2}}.
\eeaa 
\end{lemma}


\subsubsection{Part  3: Outline of the proof of Theorem \ref{THM:HIGHERDERIVS-MORAWETZ-CHP3}}
\lab{section:HigherDerivativesEnergy-Morawetz}


Theorem \ref{THM:HIGHERDERIVS-MORAWETZ-CHP3} is proved in Chapter \ref{chapter-perturbations-Kerr} according to the following steps:
\begin{enumerate}
\item First, we revisit the proof of Propositions \ref{proposition:Morawetz1-step1}, \ref{proposition:Energy1}  and  \ref{prop:morawetz-higher-order}, and of Lemma \ref{LEMMA:LOWERBOUNDPHIZOUTSIDEMTRAP}, by exhibiting the extra terms in perturbations of Kerr, and prove that the conclusions of Propositions \ref{proposition:Morawetz1-step1}, \ref{proposition:Energy1}  and  \ref{prop:morawetz-higher-order}, and of Lemma \ref{LEMMA:LOWERBOUNDPHIZOUTSIDEMTRAP},   also hold in perturbations of Kerr up to the addition of suitable error terms see sections \ref{sec:proofofresultschapter7inperturbationofKerr} and \ref{sec:proofofresultschapter8inperturbationofKerr}.

\item Next, we prove redshift estimates to remove the degeneracy on the horizon, see section \ref{section:Redshift-estimates-chp3}.

\item Then, we derive the conclusions of Theorem \ref{THM:HIGHERDERIVS-MORAWETZ-CHP3} in the particular case  $s=2$, see section \ref{sec:proofofThm:HigherDerivs-Morawetz-chp3:cases=2}.

\item Finally, we argue by iteration from $s=2$ to recover higher order derivatives which concludes the proof of Theorem \ref{THM:HIGHERDERIVS-MORAWETZ-CHP3}, see section \ref{sec:proofofThm:HigherDerivs-Morawetz-chp3:generalcase}.
\end{enumerate}


\chapter{Proof of conditional Morawetz and Energy estimates in Kerr}\label{chapter-proof-part1}


In this chapter we  prove the conditional Morawetz and energy estimates of Propositions \ref{proposition:Morawetz1-step1} and \ref{proposition:Energy1}. Recall that we are in Kerr throughout this chapter and that the results of this chapter will be extended to perturbations of Kerr in section \ref{sec:proofofresultschapter7inperturbationofKerr}.


\section{Preliminaries}\label{section-preliminaries-proof-part1}
\lab{sec:appendixonenergymomentumtensorwaveeqpsi}


In this section we collect preliminary results to apply the vector field method to obtain the desired energy-Morawetz estimates.


\subsection{Deformation tensors of   basic vectorfields}


Recall the Hawking timelike vectorfield $\That=\pr_t+\frac{a}{r^2+a^2} \pr_\phi$ defined in \eqref{define:That}. 

\begin{definition}
\lab{definition:vf-Tmod}
We define the  vectorfield 
\beaa
\That_\de:= \pr_t +\frac{a}{r^2+a^2} \chi_0\left( \de^{-1} \frac{\TT}{r^3} \right) \pr_\phi
\eeaa
with $\de=\de_{trap}$ and with $\chi_0$ the smooth  bump function
 \bea
 \lab{definition:chi-T_chi}
 \chi_0(x)=\begin{cases}
 &0  \qquad \qquad  \mbox{if}   \quad  |x| \le  1, \\
 &1  \qquad  \qquad \mbox{if}   \quad   |x|\ge 2. 
 \end{cases} 
 \eea
 We also write
 \beaa
 \That_\de:= \pr_t+\chi_\de \pr_\phi, \qquad \chi_\de:= \frac{a}{r^2+a^2} \chi_0\left( \de^{-1} \frac{\TT}{r^3} \right).
 \eeaa
\end{definition}

According to the definition \eqref{eq:def-MM-trap} of the trapped set, we have that  $\That_\de=T$  at the trapped set and  $\That_\de=\That$ away from it.
  In view of Lemma    \ref{rem:de-trap} we can  fix $\de=\frac{1}{10}$ such that $T=\partial_t$ is strictly timelike in $\MM_{trap}$.

In the derivation of  the Energy-Morawetz  inequalities we make use of  the following vectorfields:
\begin{enumerate}
\item  The radial vectorfield  $X=\FF(r) \pr_r$, for a well chosen function $\FF$.

\item The modified timelike  vectorfield  $\That_\de = \pr_t +\chi_\de(r) \pr_\phi $  as  defined  in Definition  \ref{definition:chi-T_chi}
 with $\de=\frac{1}{10}$.
\end{enumerate}

\begin{lemma}
\lab{Lemma:dualityThatRhat}
The vectorfields $\That$ and $\Rhat $ are dual to each other in the following sense. 
\begin{enumerate}
\item If $X= \Rhat$ then
\beaa
X^4 e_4-X^3  e_3=\That.
\eeaa
\item   If $X= \That$ then
\beaa
X^4 e_4-X^3  e_3=\Rhat.
\eeaa
\item  If $X=\pr_\phi$ then
\beaa
X^4 e_4-X^3  e_3=   -\frac{a\sin^2\th(r^2+a^2)}{|q|^2} \, \Rhat.
\eeaa
\item If $ X=\That_\de$ then, with $ \widetilde{\chi}_\de=( \chi_\de-\frac{a}{r^2+a^2}) =\frac{a}{r^2+a^2}\Big(  \chi_0\big( \frac{\TT}{r^3} \de^{-1}  \big)  -1\Big) $,
\beaa
 \That_\de^4 e_4 - \That_\de ^3 e_3&=& \left(1 -   \widetilde{\chi}_\de \frac{a\sin^2\th(r^2+a^2)}{|q|^2}\right)\Rhat.
\eeaa
\end{enumerate} 
\end{lemma}

\begin{proof}
The first two identities  follow  from,   see  \eqref{eq:ThatRhat-e_3e_4-Kerr},
  \beaa
 \That &=&\frac 1 2 \left( \frac{|q|^2}{r^2+a^2} e_4+\frac{\De}{r^2+a^2}  e_3\right), \qquad \Rhat= \frac 1 2 \left( \frac{|q|^2}{r^2+a^2} e_4-\frac{\De}{r^2+a^2}  e_3\right).
 \eeaa
To check the third identity we  note in view of the formula for $Z=\pr_\phi$ in  \eqref{eq:T-Z-That-e_2} that 
\beaa
\pr_\phi^4=-\frac{a(\sin\th)^2}{2}, \qquad  \pr_\phi^3=-\frac{\De}{|q|^2}\frac{a(\sin\th)^2}{2}.
\eeaa
Thus,
\beaa
\pr_\phi^4 e_4- \pr_\phi^3 e_3&=& -\frac{a\sin^2\th(r^2+a^2)}{2|q|^2}\left( \frac{|q|^2}{r^2+a^2} e_4-\frac{\De}{r^2+a^2}  e_3\right)=-\frac{a\sin^2\th(r^2+a^2)}{|q|^2} \, \Rhat
\eeaa
as stated. Finally, with $ \widetilde{\chi}_\de=( \chi_\de-\frac{a}{r^2+a^2})$,   
   \beaa
   \That_\de^4 e_4 - \That_\de ^3 e_3&=&  \That^4 e_4 - \That ^3 e_3 +  \widetilde{\chi}_\de   \big(\pr^4_\phi e_4- \pr^3_\phi e_3\big)
=\Rhat -  \widetilde{\chi}_\de   \frac{a\sin^2\th(r^2+a^2)}{|q|^2} \, \Rhat \\
&=&\left(1 -   \widetilde{\chi}_\de\frac{a\sin^2\th(r^2+a^2)}{|q|^2}\right)\Rhat
   \eeaa
as stated. 
\end{proof}

In the lemma below we calculate the   deformation tensors of these vectorfields.
\begin{lemma}
\lab{lemma:deformationtensors}
The following  identities hold true.
\begin{enumerate}
\item  For $X=\FF\pr_r $ we have
\beaa
\Lie_X(|q|^2 \g^{\a\b})&=&\big( \FF\pr_r \De -2 \De \pr_r \FF\big)\pr_r^\a\pr_r^\b+ \FF\pr_r\left(\frac 1 \De\RR^{\a\b} \right),
\eeaa
and
\beaa
\piX^{\a\b}=-|q|^{-2} \left( \big(\FF\pr_r \De -2 \De \pr_r \FF\big)\pr_r^\a\pr_r^\b+ \FF\pr_r\left(\frac 1 \De\RR^{\a\b}\right) \right) +|q|^{-2} X\big(|q|^2\big) \g^{\a\b}.
\eeaa

\item For  $\That_\de=  \pr_t+\chi_\de  \pr_\phi$  we have
\beaa
\Lie_{\That_\de}(|q|^2 \g^{\a\b})&=&- 2 \De (\pr_r \chi_\de )\,\pr_\phi^\a \pr_r^\b
\eeaa
and
\bea
\, ^{(\That_\de)} \pi^{\a\b}&=& \frac{2 \De( \pr_r \chi_\de )}{|q|^2} \,   \pr_\phi^\a \pr_r^\b, 
\eea
with
\bea
\lab{eq:pr_rchi_de}
 \pr_r \chi_\de =-\frac{2ar}{(r^2+a^2)^2}\chi_0 + O(a \de^{-1} ) \chi'_0.
\eea

\item  In particular, for $\That= \pr_t+\frac{a}{r^2+a^2} \pr_\phi$,  we have 
\beaa
\Lie_{\That}(|q|^2 \g^{\a\b})&=&\frac{4ar\De}{(r^2+a^2)^2}\pr_\phi^{(\a} \pr_r^{\b)}
\eeaa
and
\bea
\, ^{(\That)} \pi^{\a\b}&=&-\frac{4a r\De}{(r^2+a^2)^2|q|^2}\pr_\phi^{(\a} \pr_r^{\b)}.
\eea
\end{enumerate}  
\end{lemma}

\begin{proof}
The  first part  of the lemma has already been established  in Lemma \ref{lemma:Lie_Xgdot}. Since $\pr_t, \pr_\phi$ are Killing, since $[\That_\de, \pr_\th]=0$,  and since $ [\That_\de, \pr_r]=-\pr_ r \chi_\de \pr_\phi$, we have, using the formula \eqref{inverse-metric} for $|q|^2\g^{\a\b}$, 
\beaa
\Lie_{\That_\de}(|q|^2 \g^{\a\b})&=& \Lie_{\That_\de}\left(\De \pr_r^\a \pr_r^\b +\frac{1}{\De}\RR^{\a\b}\right)\\
&=& \Lie_{\That_\de}\Bigg(\De \pr_r^\a \pr_r^\b + \pr_\th^\a \pr_\th^\b\\
&&+\frac{1}{\De}\left(-\Si^2\pr_t^\a\pr_t^\b-2 amr \pr_t^\a\pr_\phi^\b  -2 amr \pr_\phi^\a\pr_t^\b+\frac{\De-a^2\sin^2\th}{\sin^2\th} \pr_\phi ^\a\pr_\phi^\b\right)\Bigg)\\
&=& \De [\That_\de, \pr_r]^\a \pr_r^\b +\De \pr_r^\a [\That_\de, \pr_r^\b]=-2\De (\pr_r\chi_\de)  \pr_\phi^\a \pr_r^\b 
\eeaa
as stated.

To calculate the deformation tensors  we simply remark that, for any vectorfield $X$,
\beaa
 \piX^{\a\b}&=&-\Lie_X\big( |q|^{-2}  |q|^2 \g^{\a\b}\big)=-|q|^{-2} \Lie_X\big(  |q|^2 \g^{\a\b}\big)- |q|^2\Lie_X\big(|q|^{-2}\big)\g^{\a\b}\\
 &=& -|q|^{-2} \Lie_X\big(  |q|^2 \g^{\a\b}\big)+ |q|^{-2}X\big(|q|^{2}\big)\g^{\a\b}.
 \eeaa
 Note that  the second term vanishes for $X=\That_\de$. The identities  for $\That$ correspond to the ones for $\That_\de$ in the particular case $\chi_0=1$. 
\end{proof} 

Recall the   energy-momentum tensor as defined in Section \ref{section:energy-momentum}, i.e.
 \beaa
 \QQ_{\mu\nu}:=\Db_\mu  \psi \c \Db _\nu \psi 
          -\frac 12 \g_{\mu\nu} \left(\Db_\la \psi\c\Db^\la \psi + V\psi \c \psi\right)= \Db_\mu  \psi \c \Db _\nu \psi -\frac 1 2\g_{\mu\nu} \LL[\psi].
 \eeaa

As a corollary to the above lemma we derive the following.

\begin{lemma}
\lab{lemma:QQcdotpi}
The following identities hold true.
\begin{enumerate}
 \item For 
  $X=\FF\pr_r$ we have the identity
 \bea\lab{eq:PPpi(X)}
 \bsplit
   |q|^2  \QQ  \c\piX  
   &=  \big(2 \De \pr_r \FF- \FF\pr_r \De\big)|\nab_r\psi|^2 
   - \FF\pr_r\left(\frac 1 \De\RR^{\a\b}\right)  \Db_\a \psi\c \Db_\b \psi \\
   &+X\big( |q|^2\big)\big(\LL[\psi]-V|\psi|^2\big)  -\LL[\psi]  |q|^2\Div_\g X.
   \end{split}
   \eea
    \item For $\That_\de$  we have the identity, 
    \bea
    \lab{eq:QQpi(That_de)}
 \bsplit
   |q|^2  \QQ  \c \,^{(\That_\de)} \pi&= 2 (r^2+a^2)    (\pr_r \chi_\de) \nab_\phi\psi\c\nab_\Rhat \psi\\
   &=\left(-\frac{4ar}{r^2+a^2}\chi_0+  O(a \de^{-1} ) \chi'_0 \right)\nab_\phi\psi \c \nab_\Rhat \psi.
   \end{split}
  \eea
  
\item  In particular
 \bea
   |q|^2  \QQ  \c \,^{(\That)} \pi&=&- \frac{4ar }{r^2+a^2}   \nab_\phi\psi\c\nab_\Rhat \psi.
 \eea
 \end{enumerate}
 \end{lemma}
 
 \begin{proof} 
 We compute
  \beaa
     \QQ  \c\piX&=&\piX^{\a\b}\left(\Db_\a \psi \Db_\b \psi - \frac 12 \g_{\a\b}\LL[\psi] \right)=
      \piX^{\a\b} \Db_\a \psi\c \Db_\b \psi -\frac 1 2\LL[\psi] \g^{\a\b} \piX _{\a\b}\\
      &=&   \piX^{\a\b} \Db_\a\psi\c \Db_\b \psi -\LL[\psi] \Div_\g X.
    \eeaa
    According to Lemma \ref{lemma:deformationtensors}, we have\footnote{Here $\nab_r =\nab_{\pr_r}$}  
    \beaa
   |q|^2    \piX^{\a\b} \Db_\a \psi \Db_\b \psi
   &=& \big(2 \De \pr_r \FF- \FF\pr_r \De\big)|\nab_r\psi|^2 
   - \FF\pr_r\left(\frac 1 \De\RR^{\a\b}\right)  \Db_\a \psi\c \Db_\b \psi \\
   &&+X\big( |q|^2\big) \g^{\a\b}\Db_\a\psi\c\Db_\b \psi\\
   &=& \big(2 \De \pr_r \FF- \FF\pr_r \De\big)|\nab_r\psi|^2 
   - \FF\pr_r\left(\frac 1 \De\RR^{\a\b}\right)  \Db_\a \psi\c \Db_\b \psi \\
   &&+X\big( |q|^2\big)\big(\LL[\psi]-V|\psi|^2\big)
      \eeaa
    and hence
    \beaa
   |q|^2  \QQ  \c\piX  
   &=&  \big(2 \De \pr_r \FF- \FF\pr_r \De\big)|\nab_r\psi|^2 
   - \FF\pr_r\left(\frac 1 \De\RR^{\a\b}\right)  \Db_\a \psi\c \Db_\b \psi \\
   &&+X\big( |q|^2\big)\big(\LL[\psi]-V\psi^2\big)  -\LL[\psi]  |q|^2\Div_\g X
    \eeaa
    as stated. 
    
    Similarly, we have
    \beaa
|q|^2   \QQ\c  \, ^{(\That_\de)} \pi&=&|q|^2\,       ^{(\That_\de)} \pi^{\a\b} \Db_\a \psi\c \Db_\b \psi - \frac 1 2|q|^2 \g_{\a\b}\, ^{(\That_\de)} \pi^{\a\b}  \LL[\psi]\\
   &=&\Big(  2\De  (\pr_r \chi_\de) \pr_\phi^\a \pr_r^\b \Big) \Db_\a \psi\c \Db_\b \psi   - \frac 1 2\Big(  2\De  (\pr_r \chi_\de) \pr_\phi^\a \pr_r^\b \Big)\g_{\a\b}\LL[\psi]\\
   &=& 2\De ( \pr_r \chi_\de)  \nab_\phi\psi\c\nab_r \psi =2 (r^2+a^2)    (\pr_r \chi_\de) \nab_\phi\psi\c\nab_\Rhat \psi
    \eeaa
 since $\nab_\Rhat=\frac{\De}{r^2+a^2} \nab_r$.   Finally, the expression for  $\That$ corresponds to the ones for $\That_\de$ in the particular case $\chi_0=1$.
    \end{proof}


\subsection{Basic spacetime  identity for $X=\FF \pr_r$}


In this section we prove a fundamental spacetime identity for $X=\FF(r)\partial_r$, as summarized in Proposition \ref{proposition:Morawetz1}. The computations in this section follow closely the corresponding ones in \cite{A-B}.

Recall, see   Proposition \ref{prop-app:stadard-comp-Psi-perturbations-Kerr}, that for $X$ a linear combination of $e_3, e_4$ we can write in Kerr
  \beaa
  \D^\mu  \PP_\mu[X, w, M] &=& \frac 1 2 \QQ  \c\piX - \frac 1 2 X( V ) |\psi|^2+\frac 12  w \LL[\psi] -\frac 1 4|\psi|^2   \square_\g  w + \frac 1 4  \Div(|\psi|^2 M\big)\\
  &-& \big(\rhod +\etab\wedge\eta\big)\nab_{X^4e_4-X^3e_3}  \psi\c\dual\psi \\
&-& \frac{1}{2}\Im\Big(\tr\Xb H X^3 +\tr X\Hb X^4\Big)\c\nab\psi\c\dual\psi\\
  &+&  \left(\nab_X\psi +\frac 1 2   w \psi\right)\c \big(\squared_k\psi - V \psi \big).
 \eeaa
We introduce     the expression
      \bea\lab{definition-EE-gen}
      \begin{split}
\EE[X, w, M] &:= \D^\mu  \PP_\mu[X, w, M] -  \left(\nab_X\psi+\frac 1 2   w \psi\right)\c  \big(\squared_k\psi - V \psi \big)\\
&  +\big(\rhod +\etab\wedge\eta\big)\nab_{X^4e_4-X^3e_3}  \psi\c\dual\psi +\frac{1}{2}\Im\Big(\tr\Xb H X^3 +\tr X\Hb X^4\Big)\c\nab\psi\c\dual\psi
\end{split}
\eea
which represents the current for the Morawetz bulk.

In what follows, we assume,  in BL coordinates\footnote{Or,  relative to the ingoing null  frame, $X=\frac 1 2\FF \Big(\frac{|q|^2}{\De}  e_4- e_3\Big)$.}, 
  $X=\FF\pr_r=\FF \frac{r^2+a^2}{\De}\Rhat$.   In that case,  in view of Lemma \ref{Lemma:dualityThatRhat},    
  \beaa
  \nab_{X^4e_4-X^3e_3} = \FF \frac{r^2+a^2}{\De}\nab_{\That} \psi.
  \eeaa
Also, since $X^4=\FF\frac{|q|^2}{2\De}$ and $X^3=-\frac{1}{2}\FF$, we have
\beaa
\Im\Big(\tr\Xb H X^3 +\tr X\Hb X^4\Big) &=& \Im\left(\frac{1}{\ov{q}}\frac{aq}{|q|^2}\FF\Jk -\frac{1}{q}\frac{a\ov{q}}{|q|^2}\FF\Jk\right) = \frac{a\FF(r)}{|q|^4}\Im((q^2-\ov{q}^2)\Jk)\\
&=& \frac{4a^2r\cos\th\FF(r)}{|q|^4}\Re(\Jk)\\
&=& \frac{4a^2r\cos\th\FF(r)}{(r^2+a^2)|q|^4}\left(\pr_\phi +a(\sin\th)^2\frac{r^2+a^2}{|q|^2}\That\right)\\
&=& \frac{4a^2r\cos\th\FF(r)}{(r^2+a^2)|q|^4}\pr_\phi +\frac{4a^3r\cos\th(\sin\th)^2\FF(r)}{|q|^6}\That.
\eeaa
 Thus \eqref{definition-EE-gen} takes the form  
      \bea\lab{definition-EE-X=FFR}
      \begin{split}
\EE[X, w, M] &= \D^\mu  \PP_\mu[X, w, M]\\
&-  \left(\nab_X\psi+\frac 1 2   w \psi\right)\c  \big(\squared_k\psi - V \psi \big) +\frac{2a^2r\cos\th\FF(r)}{(r^2+a^2)|q|^4}\nab_\phi\psi\c\dual\psi\\
& +\left(\big(\rhod +\etab\wedge\eta\big)\frac{r^2+a^2}{\De}+\frac{2a^3r\cos\th(\sin\th)^2}{|q|^6}\right)\FF\nab_\That\psi\c\dual\psi.
\end{split}
\eea
With this definition of $\EE$ we write, using \eqref{eq:PPpi(X)}, 
\beaa
|q|^2\EE[X, w, M]   &=& \frac 1 2 |q|^2\QQ  \c\piX +|q|^2\left( - \frac 1 2 X( V )  |\psi|^2+\frac 12  w \LL[\psi] -\frac 1 4|\psi|^2   \square_\g  w\right)\\
&& + \frac 1 4 |q|^2 \Div(|\psi|^2 M\big)\\
   &=&  \left( \De \pr_r \FF- \frac 1 2 \FF\pr_r \De\right)|\nab_r\psi|^2 
   -\frac 1 2  \FF\pr_r\left(\frac 1 \De\RR^{\a\b}\right)  \Db_\a \psi\c \Db_\b \psi \\
     &&+\frac 1 2 \Big(  X\big( |q|^2\big)- |q|^2 \Div_\g X+ |q|^2  w\Big)\LL[\psi]\\
   &&-\frac 1 2 \left( X\big(|q|^2\big) V  +|q|^2  X(V)+\frac 1 2|q|^2 \square_\g  w \right)|\psi|^2+ \frac 1 4 |q|^2 \Div(|\psi|^2 M\big).
 \eeaa
     To simplify the coefficient of $\LL[\psi]$ we introduce a reduced function $w_{red}$, given by
       \beaa
       w_{red}&:=& |q|^2 \D_\a\big( |q|^{-2}  X^\a\big)- w=  \Div_\g X-|q|^{-2}  X(|q|^2) -w
       \eeaa
     and therefore we write the coefficient of $\LL[\psi]$ as
       \beaa
   \frac 1 2 \Big(    X\big( |q|^2\big)- |q|^2 \Div_\g X+ |q|^2  w\Big)=\frac 1 2 |q|^2 \Big( |q|^{-2} X\big( |q|^2\big)-\Div_\g X +w\Big)=-\frac 1 2 |q|^2w_{red}.
       \eeaa
  Hence
  \beaa
|q|^2\EE[X, w, M]      &=&  \left( \De \pr_r \FF- \frac 1 2 \FF\pr_r \De\right)|\nab_r\psi|^2 
   -\frac 1 2  \FF\pr_r\left(\frac 1 \De\RR^{\a\b}\right)  \Db_\a \psi\c \Db_\b \psi \\
   && -\frac 1 2 |q|^2w_{red} \LL[\psi] -\frac 1 2 \Big( X\big(|q|^2\big) V  +|q|^2  X(V)+\frac 1 2|q|^2 \square_\g  w \Big)|\psi|^2\\
   &&+ \frac 1 4 |q|^2 \Div(|\psi|^2 M\big).
 \eeaa
    Finally writing
    \beaa
   |q|^2  \LL[\psi]&=&|q|^2 \g^{\a\b}\Db_\a\psi\c\Db_\b\psi+|q|^2 V|\psi|^2=\left( \De \pr_r^\a \pr_r^\b +\frac{1}{\De}\RR^{\a\b}\right) \Db_\a\psi\c\Db_\b\psi+ |q|^2 V|\psi|^2\\
   &=&\De | \nab_r\psi|^2+\frac 1 \De  \RR^{\a\b}\Db_\a\psi\c\Db_\b\psi+|q|^2 V|\psi|^2
    \eeaa
    we obtain\footnote{Observe that the first line is identical (with opposite sign) to the computations obtained in \eqref{d-dtau-e-FF} in Proposition \ref{proposition:geodidentity1} in the case of geodesics.}
      \beaa
|q|^2\EE[X, w, M]      &=&  \left( \De \pr_r \FF- \frac 1 2 \FF\pr_r \De\right)|\nab_r\psi|^2 
   -\frac 1 2  \FF\pr_r\left(\frac 1 \De\RR^{\a\b}\right)  \Db_\a \psi\c \Db_\b \psi \\
   &&-\frac 1 2 w_{red}  \left(\De | \nab_r\psi|^2+\frac 1 \De  \RR^{\a\b}\Db_\a\psi\c\Db_\b\psi+|q|^2 V|\psi|^2 \right) \\
   &&-\frac 1 2 \Big( X\big(|q|^2\big) V  +|q|^2  X(V)+\frac 1 2|q|^2 \square_\g  w \Big)|\psi|^2+ \frac 1 4 |q|^2 \Div(|\psi|^2 M\big)\\
   &=&    \left( \De \pr_r \FF- \frac 1 2 \FF\pr_r \De-\frac 1 2 \De w_{red}\right)|\nab_r\psi|^2\\
   && -\frac 1 2 \left( \FF \pr_r\left(\frac 1 \De \RR^{\a\b}\right) + w_{red} \frac{1}{\De} \RR^{\a\b} \right)  \Db_\a \psi\c \Db_\b \psi\\
   &&-\frac 1 2 \Big( X\big(|q|^2\big) V  +|q|^2  X(V)+\frac 1 2|q|^2 \square_\g  w  + |q|^2  w_{red} V \Big)|\psi|^2\\
   &&+ \frac 1 4 |q|^2 \Div(|\psi|^2 M\big).
 \eeaa

        We summarize the result in the first part of the  following.
   \begin{proposition}   
        \lab{proposition:Morawetz1}
        The following statements hold true.
        \begin{enumerate}
        \item  Let  $\FF, w_{red}$   given  functions of  $r$.
         With the choice  of vectorfield $X=\FF \pr_r$ and scalar function $w=  |q|^2 \D_\a\big( |q|^{-2}  X^\a\big)-w_{red}$,  the generalized current $\EE[X, w, M]$ defined in \eqref{definition-EE-gen}
  verifies  
  \bea
  \lab{identity:prop.Morawetz1}
  \bsplit
   |q|^2\EE[X, w, M]  &=\AA |\nab_r\psi|^2 + \UU^{\a\b}(\Db_\a \psi )\c(\Db_\b \psi )+\VV |\psi|^2 \\
   &+\frac 1 4 |q|^2  \D^\mu (|\psi|^2 M_\mu)  
   \end{split}
   \eea
   where\footnote{Observe that the expressions for $\UU^{\a\b}$ and for $\AA$ are the same as the ones for geodesics, see Proposition \ref{proposition:geodidentity1}.}
   \beaa
   \bsplit
   \AA&= \De \pr_r \FF- \frac 1 2 \FF\pr_r \De-\frac 1 2 \De w_{red},\\
   \UU^{\a\b}&=  -\frac 1 2  \FF\pr_r \left(\frac 1 \De\RR^{\a\b}\right)-\frac 1 2   w_{red}\frac 1 \De \RR^{\a\b},\\
   \VV&= -\frac 1 2 \Big( X\big(|q|^2\big) V  +|q|^2  X(V)+\frac 1 2|q|^2 \square_\g  w  + |q|^2  w_{red} V \Big).
   \end{split}
   \eeaa
   \item  If in addition we choose, for  fixed  functions  $z, f, h$ depending on $r$,
\bea
\lab{def-w-red-in-fun-FF-00-wave}
\FF=-z h f, \qquad               w_{red}=  \FF  z^{-1}\partial_r z = - (\partial_r z ) h  f , \qquad w =- z \pr_r \big( h  f  \big), \lab{Equation:w}
\eea
then
\bea
\lab{eq:coeeficientsUUAAVV}
\bsplit
 \AA&=-z^{1/2}\Delta^{3/2} \partial_r\left(h \frac{ z^{1/2}  f }{\Delta^{1/2}}  \right),   \\
  \UU^{\a\b}&=   \frac{ 1}{2}  h f \pr_r\left( \frac z \De\RR^{\a\b}\right),\\
\VV&=   \frac 1 4 \left( \pr_r\Bigg(\De \pr_r \Big(
 z \pr_r \big( h f \big)  \Big)  \Bigg)+8  h    \pr_r \left(\frac{z\De}{r^2+a^2} \right)   f \right).
 \end{split}
\eea

\item If  $M = v(r) \pr_r$, for some function $v=v(r)$, we have
\bea\label{expression-Div-M-I}
\frac 1 4 |q|^2 \Div(|\psi|^2 M\big)&=& \frac 1 4 |q|^2\left( 2 v(r)\psi\c \nab_r \psi + \left(\pr_r v+ \frac{2r}{|q|^2} v\right) |\psi|^2 \right).
\eea  
    \end{enumerate}
        \end{proposition}
        
        \begin{proof}
        It  remains to  check the last two  parts of the proposition.  
        \medskip
        
        {\bf  Calculation of  $\UU^{\a\b}$.} As in the computations after Proposition \ref{proposition:geodidentity1}, we have 
 \beaa
\UU^{\a\b} &=&  -\frac 1 2  \FF\pr_r \left(\frac 1 \De\RR^{\a\b}\right)-\frac 1 2   w_{red}\frac 1 \De \RR^{\a\b}=  - \frac{ 1}{2}\FF z^{-1}  \partial_r\left( \frac z \De\RR^{\a\b}\right)+ \frac{ 1}{2}\left(\FF z^{-1}\partial_r z  -w_{red}\right) \frac{ \RR^{\a\b}}{\Delta}.
 \eeaa
Choosing  $w_{red}=  \FF z^{-1}\partial_r z $,  the coefficient of $\frac{ \RR^{\a\b}}{\Delta}$ cancels out, and setting $\FF=-zhf$, we  deduce the stated expression for $\UU^{\a\b}$. 
\medskip

  {\bf  Calculation of  $\AA$.}   
As in the computations after Proposition \ref{proposition:geodidentity1}, with the choices of $\FF$ and $w_{red}$ given by \eqref{def-w-red-in-fun-FF-00-wave}, we compute 
\beaa
\AA&=&  - \frac 1 2\FF\pr_r \De+  \De \pr_r \FF - \frac 1 2 \De w_{red}=  \partial_r\left(\frac{\FF}{\Delta^{1/2}}  \right) \Delta^{3/2}-\frac 12\Delta   w_{red} \\
&=&\partial_r\left(\frac{- z h f }{\Delta^{1/2}}  \right) \Delta^{3/2}-\frac 12\Delta  (- ( \partial_r z) hf )=- \frac 1 2 \partial_rz  \left(\frac{  h f}{\Delta^{1/2}}  \right) \Delta^{3/2}\\
&& -z^{1/2} \partial_r\left(\frac{ z^{1/2}  h f}{\Delta^{1/2}}  \right) \Delta^{3/2}+\frac 12\Delta   ( \partial_r z) hf \\
&=&-z^{1/2}\Delta^{3/2} \partial_r\left(h \frac{ z^{1/2}  f}{\Delta^{1/2}}  \right) 
\eeaa
as stated.

\medskip

{\bf Calculation of $\VV$.}
We calculate $\VV$ with the choices we have made so far. We have
\beaa
 \VV&=& -\frac 1 2 \Big( X\big(|q|^2\big) V  +|q|^2  X(V)+\frac 1 2|q|^2 \square_\g  w  + |q|^2  w_{red} V \Big)=\VV_0+\VV_1
\eeaa
with
\beaa
 \VV_0:= -\frac 1 4 |q|^2 \square_\g  w, \qquad  \VV_1:= -\frac 1 2 \Big( X\big(|q|^2\big) V  +|q|^2  X(V) + |q|^2  w_{red} V \Big).
\eeaa
We first calculate $\VV_0$.
Recalling  the definition of $w$ and $w_{red}=  \FF z^{-1}\partial_r z $,
\beaa
w &=&  |q|^2 \D_\a\big( |q|^{-2}  \big(X \big)^\a\big)-w _{red}=
  |q|^2 \D_\a\big( |q|^{-2} \FF\pr_r ^\a\big)-\FF z^{-1}\partial_r z \\
  &=&   |q|^2\pr_r\big(|q|^{-2}  \FF\big)+\FF( \D_\a\pr_r^\a) -\FF z^{-1}\partial_r z.
\eeaa
We write, for $Y=\pr_r$,
\beaa
\D_\a Y^\a=\frac{1}{\sqrt{|g|}} \pr_\a\big( \sqrt{|\g|} Y^\a \big)= \frac{1}{\sqrt{|g|}} \pr_r \big( \sqrt{|\g|} \big)=\frac{1}{|q|^2} \pr_r \big(|q|^2\big).
\eeaa
Hence
\beaa
w &=&   |q|^2\pr_r\big(|q|^{-2}  \FF\big)+ \frac{1}{|q|^2} \pr_r \big(|q|^2\big)\FF -\FF z^{-1}\partial_r z=\pr_r \FF -\FF z^{-1}\partial_r z= z \pr_r\left( \frac{\FF}{z} \right).
\eeaa
Thus, in view of our choice for $\FF = - z h f$ in \eqref{def-w-red-in-fun-FF-00-wave}
\bea\label{definition-w-}
w &=& z \pr_r\left( \frac{\FF}{z} \right)=- z \pr_r \big( h  f  \big).
\eea
Now, for a function $H=H(r)$, 
\bea
\lab{eq:waveH(r)}
 \square_\g H  =\frac{1}{\sqrt{|\g|} } \pr_\a \big(\sqrt{|\g|} \g^{\a\b} \pr_\b\big)  H= \frac{1}{\sqrt{|\g|} } \pr_r  \big(\sqrt{|\g|} \g^{rr} \pr_r \big)  H=\frac{1}{|q|^2} \pr_r \big( \De\pr_r H\big).
\eea
Thus
\beaa
 |q|^2 \square_\g  w=\pr_r \Big(\De \pr_r \big(w  \big)\Big)= \pr_r\Bigg(\De \pr_r \Big(-
 z \pr_r \big( h f\big)  \Big)  \Bigg).
\eeaa
We deduce,
\beaa
 \VV_0:= -\frac 1 4 |q|^2 \square_\g  w=\frac 1 4  \pr_r\Bigg(\De \pr_r \Big(
 z \pr_r \big( h  f\big)  \Big)  \Bigg).
\eeaa
It remains to calculate
\beaa
\VV_1&=& -\frac 1 2 \Big( X\big(|q|^2\big) V  +|q|^2  X(V) + |q|^2  w_{red} V \Big)=-\frac 1 2  \Big( X\big(|q|^2 V\big)    + |q|^2  w_{red} V \Big).
\eeaa
Recalling that  $|q|^2V=\frac{4\De}{r^2+a^2} $,  $w_{red}=  \FF z^{-1}\partial_r z$  and $\FF=-zhf$    we deduce
\beaa
\VV_1&=&- 2  \left( X\left(\frac{\De}{r^2+a^2} \right)    +   w_{red} \frac{\De}{r^2+a^2}  \right)=- 2  \left( \FF \pr_r \left(\frac{\De}{r^2+a^2} \right)    +    \FF z^{-1}\partial_r z\frac{\De}{r^2+a^2}  \right)\\
&=& -2  z^{-1}  \pr_r \left(\frac{z\De}{r^2+a^2} \right)  \FF  = -2 z^{-1}   \pr_r \left(\frac{z\De}{r^2+a^2} \right)  \big(-z  h f \big)= 2 hf  \pr_r  \left(\frac{z\De}{r^2+a^2} \right).
\eeaa
Thus, as stated,
\beaa
\VV&=& \VV_0+    \VV_1= \frac 1 4  \pr_r\Bigg(\De \pr_r \Big(
 z \pr_r \big( h  f \big)  \Big)  \Bigg)+2  h    \pr_r \left(\frac{z\De}{r^2+a^2} \right)   f.
\eeaa

\medskip
   {\bf Calculation of   $ \frac 1 4 |q|^2 \Div(|\psi|^2 M\big)$.}
We choose   $M = v(r) \pr_r$ for some function $v=v(r)$ and write
\beaa
  \D^\mu (|\psi|^2 M_\mu)&=&2  v(r)\psi\c \nab_r \psi +|\psi|^2 \Div M.
  \eeaa
On the other hand
\beaa
\Div M &=&\frac{1}{\sqrt{ |\g|} }\pr_\a \big(\g^{\a\b}  \sqrt{ |\g|}M_\b\big) = \frac{1}{\sqrt{ |\g|} }\pr_r \Big(\sqrt{ |\g|} v\Big)= \frac{1}{|q|^ 2} \pr_r \big( |q|^2  v)= \pr_r v+ \frac{2r}{|q|^2} v.
\eeaa
Hence
\bea
\lab{eq:expression-DivM}
\frac 1 4 |q|^2 \Div(|\psi|^2 M\big)&=& \frac 1 4 |q|^2\left(2  v(r)\psi\c \nab_r \psi + \left(\pr_r v+ \frac{2r}{|q|^2} v\right) |\psi|^2 \right)
\eea   
     as stated.
        \end{proof}  
        
                We collect here a lemma which will be used in the derivation of the estimates to analyze the right hand side of the main equation \eqref{eq:Gen.RW}.
                
\begin{lemma}
\lab{identity:X+frac12wpsiN_1}
The following identity holds true with $X$ and $w$ as in Proposition \ref{proposition:Morawetz1}.
\begin{enumerate}
\item We have, with $T=\pr_t$ in BL coordinates
\bea
\bsplit
\left(\nab_X\psi  +\frac 12 w \psi\right) \c \nab_T (\dual \psi)&= \frac 1 2 (\pr_r z) hf \psi\c \nab_T\dual \psi  +zhf  \rhod\frac{|q|^2}{\De}|\psi|^2\\
   & -\frac 1 2 \nab_r\Big(  zh f \psi\c\nab_T\dual \psi\Big)   +\frac 1 2\nab_T\Big(  zhf \psi\c  \nab_r\dual \psi\Big).
   \end{split}
\eea

\item  We have with $Z=\pr_\phi$ in BL coordinates
\bea
\bsplit
\left(\nab_X\psi  +\frac 12 w \psi \right) \c \nab_Z (\dual \psi)&= \frac 1 2 (\pr_r z) hf \psi\c \nab_Z \dual \psi -zhf  \rhod\frac{a(\sin\th)^2|q|^2}{\De}|\psi|^2\\
   & -\frac 1 2 \nab_r\Big(  zh f \psi\c\nab_Z\dual \psi\Big)   +\frac 1 2\nab_Z \Big(  zhf \psi\c  \nab_r\dual \psi\Big).
   \end{split}
\eea
\end{enumerate}
\end{lemma}

\begin{proof}
 We  set recalling \eqref{def-w-red-in-fun-FF-00-wave}
 \beaa
  J:=- \left(\nab_X\psi  +\frac 12 w \psi \right) \c \nab_T (\dual \psi)= \left(  z  h  f  \nab_r\psi+\frac 1 2  z \pr_r \big( h  f  \big)  \psi\right)\c\nab_T\dual \psi
  \eeaa
 and proceed   as follows
    \beaa
  J&=&
     z h  f \nab_r\psi \c\nab_T\dual \psi+ \frac 1 2  z \pr_r \big( h  f  \big)  \psi \c\nab_T\dual \psi\\
     &=&  z hf  \nab_r\psi \c\nab_T\dual \psi+\frac 1 2 \nab_r\Big(  zh f \psi\c\nab_T\dual \psi\Big)-\frac 1 2 
     zh f\nab_r \psi\c\nab_T\dual \psi\\
     &&-\frac 1 2 (\pr_r z) hf \psi\c \nab_T\dual \psi -\frac 1 2 zhf \psi\c \nab_r\nab_t \dual \psi\\
     &=&\frac 1 2   zh f\nab_r \psi\c\nab_T\dual \psi -\frac 1 2 (\pr_r z) hf \psi\c \nab_T\dual \psi  +\frac 1 2 \pr_r\Big(  zh f \psi\c\nab_T\dual \psi\Big)\\
     && -\frac 1 2 zhf \psi\c \nab_r\nab_t \dual \psi.
    \eeaa
    For the last term we write
    \beaa
    -\frac 1 2 zhf \psi\c \nab_r\nab_t \dual \psi &=&   -\frac 1 2 zhf \psi\c \nab_t \nab_r\dual \psi -\frac 1 2 zhf \psi\c [\nab_r, \nab_t]\dual \psi\\
    &=& -\frac 1 2 \pr_t\Big(  zhf \psi\c \nab_r \dual \psi\Big)+\frac 1 2 zhf \nab_T\psi\c\nab_r\dual \psi\\
    && -\frac 1 2 zhf \psi\c \left(-2\rhod\frac{|q|^2}{\De}\dual(\dual\psi)\right)\\
    &=& -\frac 1 2 \pr_t\Big(  zhf \psi\c \nab_r \dual \psi\Big)+\frac 1 2 zhf \nab_T\psi\c\nab_r\dual \psi
    - zhf  \rhod\frac{|q|^2}{\De}|\psi|^2.
    \eeaa
    Hence
    \beaa
    J&=& \frac 1 2   zh f\nab_r \psi\c\nab_T\dual \psi + \frac 1 2 zhf \nab_T\psi\c\nab_r\dual \psi -\frac 1 2 (\pr_r z) hf \psi\c \nab_T\dual \psi  - zhf  \rhod\frac{|q|^2}{\De}|\psi|^2\\
    &&+\frac 1 2 \pr_r\Big(  zh f \psi\c\nab_T\dual \psi\Big)   -\frac 1 2\pr_t\Big( zhf \psi\c  \nab_r\dual \psi\Big).
    \eeaa
    Note that $ \nab_r \psi\c\nab_T\dual \psi +\nab_T\psi\c\nab_r\dual \psi =0$.
    Therefore,
    \beaa
    J&=& -\frac 1 2 (\pr_r z) hf \psi\c \nab_T\dual \psi   - zhf  \rhod\frac{|q|^2}{\De}|\psi|^2
    +\frac 1 2 \nab_r\Big(  zh f \psi\c\nab_T\dual \psi\Big)\\
    &&    -\frac 1 2 \pr_t \Big(zhf \psi\c  \nab_r\dual \psi\Big).
    \eeaa
     The second statement is proved in the same way.
\end{proof}

        
 \subsection{Choice of $z$, $f$ and $h$}
 \lab{section:choicezhf}
 
      
        In this section, we present a choice for the functions $z$, $h$, $f$.
            Our goal is to choose such functions  so that the generalized current $\EE[X, w, M]$ as given in Proposition \ref{proposition:Morawetz1} is positive definite. 
      According to   \eqref{identity:prop.Morawetz1}, we have 
 \beaa
  \bsplit
   |q|^2 \EE[X, w, M] &=\AA |\nab_r\psi|^2 + \UU^{\a\b}(\Db_\a \psi )\c(\Db_\b \psi )+\VV |\psi|^2 +\frac 1 4 |q|^2  \D^\mu (|\psi|^2 M_\mu).
   \end{split}
   \eeaa
   We start by looking at what we call \textit{principal term} in derivative of $\psi$, i.e.
   \bea
    P:=\UU^{\a\b}(\Db_\a\psi )\c( \Db_\b \psi)=      \frac{ 1}{2}  h f\RRtp^{\a\b}
 \Db_\a\psi\c  \Db_\b \psi, \quad  \RRtp^{\a\b}:= \pr_r\left( \frac z \De\RR^{\a\b}\right).
   \eea
    From \eqref{eq:expressionRR-O}, we write
\beaa
\RR^{\a\b}&= -(r^2+a^2) ^2\pr_t^\a\pr_t^\b- 2a(r^2+a^2) \pr_t^{(\a} \pr_\phi^{\b)}- a^2 \pr_\phi^\a \pr_\phi^\b  +\De O^{\a\b}.
     \eeaa   
          Hence
\beaa
\RRtp^{\a\b}&=& -\pr_r\left(\frac{z}{\De} ( r^2+a^2)^2 \right) \pr_t^\a\pr_t^\b - 2a \pr_r\left(\frac{z}{\De} (r^2+a^2)\right) \pr_t^{(\a} \pr_\phi^{\b)} - a^2\pr_r\left(\frac{z}{\De} \right) \pr_\phi^\a \pr_\phi^\b\\
&& +(\pr_rz)  O^{\a\b}.
\eeaa
Just as in \eqref{choice-z-geodesics} in the case of geodesics, we similarly choose $z$ to cancel the coefficient of $ \pr_t^\a\pr_t^\b$ in $\RRt'^{\a\b}$, i.e.
\bea\label{choice-z-wave}
z&=& \frac{\De}{(r^2+a^2)^2}.
\eea
Recall \eqref{eq:pr_rz}, i.e.
\beaa
\pr_r z=-\frac{2\TT}{(r^2+a^2)^3} ,
\eeaa
and thus
\bea\label{RRt'-general-case}
\RRt'^{\a\b}&=&     \frac{4ar}{(r^2+a^2)^2}  \pr_t^{(\a} \pr_\phi^{\b)} +\frac{4a^2r}{(r^2+a^2)^3} \pr_\phi^\a \pr_\phi^\b-\frac{2\TT}{(r^2+a^2)^3}  O^{\a\b}.
\eea
The principal term $P$ then becomes
 \beaa
  P&=&\frac 1 2 h f  \left( -\frac{2\TT}{ (r^2+a^2)^3} O^{\a\b}+  \frac{4ar}{(r^2+a^2)^2} \partial_t^{(\a} \partial_\phi^{\b)}+ \frac{4a^2r}{(r^2+a^2)^3}  \partial_\phi^\a \partial_\phi^\b\right)\Db_\a \psi\c\Db_\b\psi\\
 &=& \frac 1 2 h f  \left(  -\frac{2\TT}{ (r^2+a^2)^3} O^{\a\b} \Db_\a \psi\c\Db_\b\psi\ \right)\\
 &&+\frac 1 2 h f \frac{4ar}{(r^2+a^2)^2} \left( \nab_t\psi \c \nab_\phi \psi+\frac{a}{(r^2+a^2)}\nab_\phi \psi \c\nab_\phi\psi\right)\\
&=& \frac 1 2 h f  \left(  -\frac{2\TT}{ (r^2+a^2)^3} O^{\a\b} \Db_\a \psi\c\Db_\b\psi+ \frac{4ar}{(r^2+a^2)^2} \left( \nab_t \psi+ \frac{a}{r^2+a^2} \nab_\phi \psi \right) \c \nab_\phi \psi\right),
 \eeaa
i.e. using the Hawking vector field $\That=\pr_t+\frac{a}{r^2+a^2} \pr_\phi$, 
\bea
 \lab{eq:P-firstMorawetz}
P &=& \frac 1 2 h f  \left(  -\frac{2\TT}{ (r^2+a^2)^3} O^{\a\b} \Db_\a \psi\c\Db_\b\psi+ \frac{4ar}{(r^2+a^2)^2} \nab_\That\psi \c\nab_\phi \psi\right).
 \eea
 
\begin{remark}
\lab{remark:P-axisymmetric}
In the particular case of axial symmetry, i.e. $\nab_\phi \psi=0$, \eqref{eq:P-firstMorawetz} becomes
\beaa
 P&=&-    \frac 1 2 h f \frac{2\TT}{ (r^2+a^2)^3} O^{\a\b} \Db_\a \psi\c\Db_\b\psi
\eeaa
and  thus, to  have $P$ non-negative, it makes sense to choose, for a  non-negative  $h$, 
\bea
\lab{choice-f-And-Mor}\lab{choice-f-wave}
f=\pr_r z=-\frac{2\TT}{ (r^2+a^2)^3}.
\eea
\end{remark}

In the general case, we keep the choice of $f$ given by \eqref{choice-f-wave}. 
With these choices, we compute according to formula \eqref{eq:coeeficientsUUAAVV}  for $\AA$
\beaa
\AA&=& -z^{1/2}\Delta^{3/2} \partial_r\left(\frac{   h f }{r^2+a^2}  \right)= \frac{2\Delta^{2} }{r^2+a^2}\partial_r\left(  h   \frac{\TT}{(r^2+a^2)^4} \right).
\eeaa

For the choice of $h$, in order to obtain a bound for $\AA$ which is valid in the full sub-extremal range $|a|<m$, we choose $h$ to be (as in \eqref{choice-h-subextremal} in the case of axially symmetric geodesics) 
\bea\lab{choice-h-wave-subextremal}
h=\frac{(r^2+a^2)^4}{r(r^2-a^2)}.
\eea
  Thus, with this choice,
 \beaa
 \AA&=& \frac{2\Delta^{2} }{r^2+a^2}\pr_r\left(\frac{\TT}{r(r^2-a^2)}\right)= \frac{2\Delta^{2} }{r^2(r^2-a^2)^2(r^2+a^2)} \big(3mr^4-4a^2r^3+ma^4\big),
 \eeaa
 which is positive in the exterior region in the sub-extremal range.

  It remains to calculate the  coefficient $\VV$  of  the lower order term. 
According to formula \eqref{eq:coeeficientsUUAAVV}  for $\VV$, in  view of our choices of $z, f, h$ we derive
\beaa
\VV_1&=&2  hf    \pr_r \left(z \frac{\De}{r^2+a^2} \right)  = -4  \frac{(r^2+a^2)\TT}{  r(r^2-a^2)}  \pr_r \left(\frac{\De^2}{(r^2+a^2)^3} \right)\\
&=&
 -4   \frac{(r^2+a^2)\TT}{  r(r^2-a^2)} \De \frac{-2r^3+8mr^2-2a^2r-4ma^2 }{(r^2+a^2)^4}\\
 &=&8  \De \frac{r^3-4mr^2+a^2r+2ma^2 }{r(r^2+a^2)^4} \frac{r^2+a^2}{r^2-a^2}   \TT,
\eeaa
and
\beaa
\VV_0&=&  \frac 1 4  \pr_r\left(\De \pr_r \left(
 z \pr_r \left( h   \frac{-2\TT}{(r^2+a^2)^3}  \right)  \right)  \right)= \frac 1 4  \pr_r\left(\De \pr_r \left(
 z \pr_r \left(    \frac{-2\TT}{r} \frac{r^2+a^2}{r^2-a^2}  \right)  \right)  \right)\\
 &=&\frac 1 2  \pr_r\left(\De \pr_r \left(
 z   \left(-2r+3m-\frac{5ma^2}{r^2}+\frac{8a^4}{r(r^2-a^2)}-\frac{12ma^4}{r^2(r^2-a^2)}+8 a^6\frac{r -m}{r^2(r^2-a^2)^2} \right) \right)\right).
\eeaa
 Continuing to differentiate  we  find
 \beaa
\bsplit
\VV_0&=  \frac{9mr^6-46m^2r^5+54m^3r^4}{(r^2+a^2)^4}+ O(a^2 r^{-3}).
\end{split}
\eeaa

We summarize the results above in the following.
\begin{proposition}
\lab{prop:Choice-zhf}
The generalized current 
induced by 
\beaa
X=-zhf \pr_r,  \qquad  w =- z \pr_r \big( h  f  \big),
\eeaa
with the choices 
\beaa
z=\frac{\De}{(r^2+a^2)^2}, \qquad f=-\frac{2\TT}{ (r^2+a^2)^3},  \qquad  h=\frac{(r^2+a^2)^4}{r(r^2-a^2)},
\eeaa
satisfies the following:
\begin{enumerate}
\item We have
\bea
\lab{eq:X=FFpr_r}\lab{eq:X-nabRh}
\FF=\frac{2\De \TT}{r(r^2-a^2)( r^2+a^2)}, \quad X= \frac{2\De\TT}{r(r^2-a^2)( r^2+a^2)}\pr_r=\frac{2\TT}{r(r^2-a^2)}\Rhat,
\eea
and
\beaa
w_{red}=  - h f \partial_r z= -  \frac{4\TT^2}{r(r^2-a^2)(r^2+a^2)^2}.
\eeaa

\item  The principal term  $   P=\UU^{\a\b}(\Db_\a\psi )( \Db_\b \psi)$ is given by
\bea\label{eq:principal-term}
P=\frac{\TT}{r}\frac{r^2+a^2}{r^2-a^2}\left(  \frac{2\TT}{ (r^2+a^2)^3}   O^{\a\b}\nab_\a\psi\c\nab_\b \psi - \frac{4ar}{(r^2+a^2)^2} \nab_\That \psi \c \nab_\phi \psi\right).
\eea

\item The  coefficients $ \AA, \VV$ in equation  \eqref{eq:coeeficientsUUAAVV}   of Proposition \ref{proposition:Morawetz1}
take the form
\bea\label{eq:prop-Choice-AA}
\bsplit
\AA&=\frac{2\Delta^{2} }{r^2(r^2-a^2)^2(r^2+a^2)} \big(3mr^4-4a^2r^3+ma^4\big) ,
\end{split}
\eea
and
\bea
\lab{eq:prop:Choice-zhf}
\bsplit
\VV&=\VV_0+\VV_1,\\
\VV_1&=8  \De \frac{r^3-4mr^2+a^2r+2ma^2 }{r(r^2+a^2)^4} \frac{r^2+a^2}{r^2-a^2}   \TT,\\
\VV_0&=   \frac{9mr^6-46m^2r^5+54m^3r^4}{(r^2+a^2)^4}+ O(a^2 r^{-3}).
\end{split}
\eea

\item For small $a$, we  have
\beaa
\VV&=&  \frac{8r^3-63mr^2+162m^2r-138m^3}{r^4} +O(a^2 r^{-3} ).
\eeaa
\end{enumerate}
\end{proposition}


\section{Conditional Morawetz estimates}\label{section:cond:mor}


In this section we derive the first conditional Morawetz estimates, proving Proposition \ref{proposition:Morawetz1-step1}.


\subsection{A first lower  bound}\label{subsection:first-lower-bound}


    According to  Proposition \ref{prop:Choice-zhf},  for  $M=0$ (i.e. $v=0$), the generalized current is given by
 \bea\label{first-low-bound-ref}
  \bsplit
   |q|^2 \EE[X, w, M=0] &=\AA |\nab_r\psi|^2 + P+\VV |\psi|^2,\quad   \\
    \AA&=\frac{2\Delta^{2} }{r^2(r^2-a^2)^2(r^2+a^2)} \big(3mr^4-4a^2r^3+ma^4\big),\\
   P&=\frac{\TT}{r}\frac{r^2+a^2}{r^2-a^2} \left(  \frac{2\TT}{ (r^2+a^2)^3}   O^{\a\b}\nab_\a\psi\c\nab_\b \psi - \frac{4ar}{(r^2+a^2)^2} \nab_\That\psi \c \nab_\phi \psi\right), \\
   \end{split}
   \eea
   with $\VV$ given by \eqref{eq:prop:Choice-zhf}, for which one easily  checks that $\VV= O( r^{-1}) $.

 In Remark \ref{remark:choice-h-subextremal} we have shown that  the polynomial  $3mr^4-4a^2r^3+ma^4  \ge 0$  for all values of $r\ge r_+$  in the range $|a|/m \leq 1$.
  More precisely, there exists a small constant $\de_0>0$  such that\footnote{In particular, note the following computation at $r=r_+$ with $\ga=|a|/m$, 
\beaa
\frac{3mr^4-4a^2r^3+ma^4}{(r^2-a^2)^2} &=& \frac{m(2(3-2\ga)+(6-\ga)\sqrt{1-\ga})}{(1+\sqrt{1-\ga})^2}>0, \qquad 0\leq \ga\leq 1.
\eeaa}  for $|a|/m \leq 1$  
 \beaa
 \AA\ge \de_0\frac{ m \De^2}{r^4}.
 \eeaa
We therefore obtain for $|a|/m < 1$
  \beaa
|q|^2 \EE[X, w, M=0]&\ge  &\de_0\frac{ m \De^2}{r^4}|\nab_r\psi|^2  + P -  O( r^{-1}) |\psi|^2 
\eeaa
which implies
 \bea\label{eq:first-lower-bound-conditional}
 \EE[X, w, M=0]&\ge &\de_0\frac{ m \De^2}{r^6}|\nab_r\psi|^2  + r^{-2}  P -  O( r^{-3}) |\psi|^2.
  \eea

 
\subsection{A  second  lower bound containing $\nab_\That$}
\lab{section:lowerboundscontainigThat}


We now want to incorporate the identity \eqref{eq:first-lower-bound-conditional} with a (trapped) control for the $\nab_\That$ derivative as well. 
We achieve  this with the help of a new 
 current of the form
 \beaa
 \PP'_\mu&=&\frac 1 2  w' \psi \c \Db_\mu \psi -\frac 1 4\psi^2   \pr_\mu w' , \qquad  w'=-w_{red}',
 \eeaa
 which corresponds to a current associated to a zero vector field $X=0$ and a scalar function $w'$ to be chosen. 
 In view  of Proposition \ref{proposition:Morawetz1} we  
derive
\beaa
|q|^2 \EE'[0, w', 0] &=&\frac 1 2 \De  \, w'     |\nab_r\psi|^2+\frac 1 2   w' \frac 1 \De \RR^{\a\b} \, \Db_\a \psi \c \Db_\b \psi  -\frac 1 2 \left( \frac 1 2|q|^2 \square_\g  w ' -  |q|^2  w'  V \right) |\psi|^2.
\eeaa
Recall that,  in view of  \eqref{eq:expressionRR-TT}
\beaa
 \RR^{\a\b} \, \Db_\a \psi \c \Db_\b \psi  &=&  -(r^2+a^2) ^2 |\nab_\That\psi|^2 + \De O^{\a\b}\Db_\a \psi\c\Db_\b \psi.
\eeaa
Thus, we have
\beaa
|q|^2 \EE'[0, w', 0] &=&\frac 1 2 \De  \, w'     |\nab_r\psi|^2-\frac{w' (r^2+a^2)^2}{ 2 \De}|\nab_\That \psi|^2 +\frac 1 2  w' O^{\a\b}\Db_\a\psi\c\Db_\b \psi \\
&& -\frac 1 2 \Big( \frac 1 2|q|^2 \square_\g  w ' -  |q|^2  w'  V \Big) |\psi|^2.
\eeaa
By summing the above to \eqref{first-low-bound-ref} we obtain
\beaa
 |q|^2\Big( \EE+  \EE' \Big)&=& -\frac{w' (r^2+a^2)^2}{ 2 \De}|\nab_\That \psi|^2 +\left(\AA+\frac 1 2\De w'\right)|\nab_r \psi|^2\\
 &&+\left(\frac{2\TT^2}{r (r^2+a^2)^2(r^2-a^2)} +\frac 1 2 w'\right)O^{\a\b} \nab_\a\psi\c\nab_\b \psi  \\
 &&- \frac{\TT}{r} \frac{4ar}{(r^2+a^2)(r^2-a^2)} \nab_\That\psi \c \nab_\phi \psi\\
 &&+\left(\VV - \frac 1 2 \left( \frac 1 2|q|^2 \square_\g  w ' -  |q|^2  w'  V \right)\right)|\psi|^2.
 \eeaa
We choose for some $\de_1>0$
\beaa
w'&=&- \de_1 \frac{4 m \De \TT^2}{r^2 (r^2+a^2)^4}
\eeaa
so that
\beaa
|q|^2\Big( \EE +  \EE' \Big)&=& \de_1 \frac{2 m\TT^2}{r^2 (r^2+a^2)^2} |\nab_\That \psi|^2\\
&&+\left(1- \de_1\frac{m\De(r^2-a^2) }{r(r^2+a^2)^2}\right)  \frac{2\TT^2}{r (r^2+a^2)^2(r^2-a^2)}O^{\a\b} \nab_\a\psi\c\nab_\b \psi \\
&&+\left(\AA+\frac 1 2\De w'\right)|\nab_r \psi|^2- \frac{\TT}{r} \frac{4ar}{(r^2+a^2)(r^2-a^2)} \nab_\That\psi \c \nab_\phi \psi\\
&&+\left(\VV - \frac 1 2 \left( \frac 1 2|q|^2 \square_\g  w ' -  |q|^2  w'  V \right)\right)|\psi|^2.
\eeaa
Observe that for $\de_1<1$, the coefficient $1- \de_1\frac{m\De(r^2-a^2)}{r(r^2+a^2)^2}$ is positive in the exterior in the full subextremal range. 
Therefore, we can write
\beaa
&& |q|^2\Big( \EE +  \EE' \Big)\\
&=&\left(\AA-  \de_1 \frac{2 m \De^2 \TT^2}{r^2 (r^2+a^2)^4}\right)|\nab_r \psi|^2\\
&&+  \frac{2 \TT^2}{r (r^2+a^2)^2} \left(\de_1 \frac{m}{r} |\nab_\That \psi|^2+\left(1- \de_1\frac{m\De(r^2-a^2)}{r(r^2+a^2)^2}\right)  \frac{1}{ (r^2 - a^2)}O^{\a\b} \nab_\a\psi\c\nab_\b \psi\right) \\
&&- \frac{\TT}{r} \frac{4ar}{(r^2+a^2)(r^2-a^2)} \nab_\That\psi \c \nab_\phi \psi+\left(\VV- \frac 1 2 \left( \frac 1 2|q|^2 \square_\g  w ' -  |q|^2  w'  V \right)\right)|\psi|^2.
\eeaa
Note also that $\EE+\EE'$ is in fact the   generalized current
associated to  $X=\FF\pr_r$ and $w$  the sum between the  old   $w =- z \pr_r \big( h  f  \big)$ (see  Proposition   \ref{prop:Choice-zhf})
 and $w'=- \de_1 \frac{4 m \De \TT^2}{r^2 (r^2+a^2)^4}$,  i.e. 
 \bea
 \lab{eq:new-w}
 w=  w_X- \de_1 \frac{4 m \De \TT^2}{r^2 (r^2+a^2)^4}, \qquad w_X= - z \pr_r \big( h  f  \big).
 \eea
   We  replace the vectorfield $\pr_r$ with $\Rhat=\frac{\De}{r^2+a^2}   \pr_r $,  see \eqref{define:Rhat},  to deduce
\beaa
|q|^2(\EE+\EE')=  |q|^2\EE[X, w]&\ge &\AA_{\de_1} \frac{(r^2+a^2)^2} {\De^2}    | \nab_\Rhat \psi|^2+P'_{\de_1} - \frac{\TT}{r} \frac{4ar}{(r^2+a^2)^2} \nab_\That\psi \c \nab_\phi \psi\\
&& +\VV_{\de_1} |\psi|^2.
\eeaa
Hence
\beaa
\EE[X, w]&\ge& \AA_{\de_1} \frac{ (r^2+a^2)^2} {|q|^2\De^2}   | \nab_\Rhat \psi|^2+\frac{1}{|q|^2} P'_{\de_1} - \frac{\TT}{r|q|^2} \frac{4ar}{(r^2+a^2)(r^2-a^2)} \nab_\That\psi \c \nab_\phi \psi\\
&& +\frac{1}{|q|^2}\VV_{\de_1} |\psi|^2,
\eeaa
with
\beaa
\bsplit
\AA_{\de_1}&= \AA-  \de_1 \frac{2 m \De^2 \TT^2}{r^2 (r^2+a^2)^4},  \qquad 
\VV_{\de_1}=\VV - \frac 1 2 \left( \frac 1 2|q|^2 \square_\g  w ' -  |q|^2  w'  V \right),\\
P'_{\de_1}&=  \frac{2 \TT^2}{r (r^2+a^2)^2} \left(\de_1 \frac{m}{r} |\nab_\That \psi|^2+\left(1- \de_1\frac{m\De }{r(r^2+a^2)}\right)  \frac{1}{ (r^2 - a^2)}O^{\a\b} \nab_\a\psi\c\nab_\b \psi\right).
\end{split}
\eeaa

Observe that
\beaa
\AA_{\de_1}&=&\frac{2\Delta^{2} }{r^2(r^2-a^2)^2(r^2+a^2)} \big(3mr^4-4a^2r^3+ma^4\big)-  \de_1 \frac{2 m \De^2 \TT^2}{r^2 (r^2+a^2)^4}\\
&=&\frac{2\Delta^{2} }{r^2(r^2+a^2)} \left(\frac{ \big(3mr^4-4a^2r^3+ma^4\big)}{(r^2-a^2)^2}-  \de_1 \frac{ m  \TT^2}{(r^2+a^2)^3}\right).
\eeaa
Observe that
\beaa
r\to \frac{ \big(3mr^4-4a^2r^3+ma^4\big)}{(r^2-a^2)^2}-  \de_1 \frac{ m  \TT^2}{(r^2+a^2)^3}
\eeaa
is increasing on $r\geq r+$ for $\de_1>0$ small enough. In particular, we have
\beaa
\AA_{\de_1} &\geq&\frac{2\Delta^{2} }{r^2(r^2+a^2)}\left[ \left(\frac{ \big(3mr^4-4a^2r^3+ma^4\big)}{(r^2-a^2)^2}-  \de_1 \frac{ m  \TT^2}{(r^2+a^2)^3}\right)\right]_{\Big|_{r=r_+}}.
\eeaa 
To compute the expression in bracket on $r=r_+$, we used $r_+^2=2mr_+-a^2$, and therefore, for $r=r_+$, we have
\beaa
\frac{\big(3mr^4-4a^2r^3+ma^4\big)}{(r^2-a^2)^2}&=&\frac{\big(3m(2mr-a^2)^2-4a^2r(2mr-a^2)+ma^4\big)}{(2mr-2a^2)^2}\\
&=&\frac{ \big(3m(4m^2r^2-4ma^2r+a^4)-4a^2r(2mr-a^2)+ma^4\big)}{4(mr-a^2)^2}\\
&=&\frac{ (3m^3-2ma^2)r^2-(3m^2a^2-a^4)r+ma^4}{(mr-a^2)^2}.
\eeaa
Using that at the horizon $(mr-a^2)^2=(2mr-a^2)(m^2-a^2)$, we have for $r=r_+$
\beaa
&=&\frac{ (6m^4-7m^2a^2+a^4)r-3a^2m^3+3ma^4}{(mr-a^2)^2}=\frac{ (6m^2-a^2)(m^2-a^2)r-3a^2m(m^2-a^2)}{(2mr-a^2)(m^2-a^2)}\\
&=& \frac{ (6m^2-a^2)r-3a^2m}{(2mr-a^2)}.
\eeaa
On the other hand, we have on $r=r_+$
\beaa
\frac{ m  \TT^2}{(r^2+a^2)^3} &=& \frac{(mr-a^2)^2}{2r^3}=\frac{r^2(m^2-a^2)}{2r^3}=\frac{m^2-a^2}{2r}
\eeaa
and hence, with the notation $\ga=a^2/m^2$, we infer
\beaa
&&\left[ \left(\frac{ \big(3mr^4-4a^2r^3+ma^4\big)}{(r^2-a^2)^2}-  \de_1 \frac{ m  \TT^2}{(r^2+a^2)^3}\right)\right]_{\Big|_{r=r_+}}\\ 
&=& \frac{ (6m^2-a^2)r_+-3a^2m}{(2mr_+-a^2)}-\de_1\frac{m^2-a^2}{2r_+}\\
&=& m\left(\frac{(6-\ga)(1+\sqrt{1-\ga})-3\ga}{2(1+\sqrt{1-\ga})-\ga}-\de_1\frac{1-\ga}{2(1+\sqrt{1-\ga})}\right)\\
&=& m\left(\frac{(6-4\ga)+(6-\ga)\sqrt{1-\ga}}{(2-\ga)+2\sqrt{1-\ga})}-\de_1\frac{1-\ga}{2(1+\sqrt{1-\ga})}\right).
\eeaa
Since $0\leq\ga\leq 1$, we infer the existence of a constant $\de_*>0$ such that for $|a|<m$
\beaa
\left[ \left(\frac{ \big(3mr^4-4a^2r^3+ma^4\big)}{(r^2-a^2)^2}-  \de_1 \frac{ m  \TT^2}{(r^2+a^2)^3}\right)\right]_{\Big|_{r=r_+}} &\geq& \de_*
\eeaa
and hence, for $|a|<m$, we obtain 
\beaa
\AA_{\de_1}&\geq & \de_* \frac{m \De^2}{(r^2+a^2)^4} r^4.
\eeaa
Also, as before we have $\frac{1}{|q|^2}\VV_{\de_1} =O(r^{-3})$.  
Finally,  the term 
 $  -  \frac{4a \TT}{|q|^2(r^2+a^2)^2} \nab_\That\psi \nab_\phi \psi$ can be bounded by Cauchy-Schwarz by a term containing $\nab_\That \psi$ which can be absorbed by $P'_{\delta_1}$ and a term  of the form 
 $O(a^2 r^{-4})  |\nab_\phi \psi|^2$.
 
 Using \eqref{eq:exxpressionO-e1e2} to write $ O^{\a\b} \nab_\a\psi\c\nab_\b\psi  = |q|^2 |\nab\psi|^2$, we summarize the final estimate in the following.

\begin{proposition}
\lab{prop:positivity-bulk} 
The generalized current induced by the vectorfield $X=\FF \pr_r$ and the scalar function $w=w_X- \de_1\frac{4 m \De \TT^2}{r^2 (r^2+a^2)^4}$ with $\FF$ and $w_X$ defined as in Proposition \ref{prop:Choice-zhf} satisfies the following estimate, for all $|a|/m <1$ and for a constant $\de_*>0$ depending\footnote{The term $|\nab_\phi \psi|^2$ is multiplied in particular by $r^4(r^2-a^2)^{-2}\les \left(1-\frac{|a|}{m}\right)^{-1}$ so that we need $\de_*\ll 1-\frac{|a|}{m}$ for $|a|$ close to $m$.}  on $1-\frac{|a|}{m}$,
\bea\lab{eq:second-bound-including-T}
\begin{split}
\lab{lemma:lower-boundEE+EE'}
\EE[X, w]&\ge\de_* \frac{m}{r^2}   | \nab_\Rhat \psi|^2+   \de_*\frac{ \TT^2}{r^6} \left( \frac{m}{ r^2} |\nab_\That \psi|^2+ r^{-1} |\nab \psi|^2\right)\\
& -O(a^2 r^{-4})  |\nab_\phi \psi|^2        -O(r^{-3})  |\psi|^2.
\end{split}
\eea
\end{proposition}

Observe that the first line on the right hand side of the above is precisely the density in $\Mordot^{ax}_{deg}[\psi](\tau_1, \tau_2)$ as defined in Definition \ref{definition:SS0Morawetznorms}.


\subsection{Proof of Proposition  \ref{proposition:Morawetz1-step1}, estimate \eqref{eq:conditional-mor-par1-1-I}}
\lab{section:Proof-proposition:Morawetz-Energy1}


We are now ready to prove the first part of Proposition \ref{proposition:Morawetz1-step1}, i.e. estimate \eqref{eq:conditional-mor-par1-1-I}.

  For $X=\FF \pr_r$ and $w=w_X-\de_1 \frac{ 4m \De \TT^2}{r^2 (r^2+a^2)^4}$ with $\FF$ and $w_X$ defined as in Proposition \ref{prop:Choice-zhf}, using \eqref {definition-EE-X=FFR}, we have
\beaa
\EE[X, w] &=& \D^\mu  \PP_\mu[X, w] -  \left(\nab_X\psi+\frac 1 2   w \psi\right)\c  \big(\squared_2 \psi - V \psi \big)  +\frac{2a^2r\cos\th\FF(r)}{(r^2+a^2)|q|^4}\nab_\phi\psi\c\dual\psi\\
&& +\left(\big(\rhod +\etab\wedge\eta\big)\frac{r^2+a^2}{\De}+\frac{2a^3r\cos\th(\sin\th)^2}{|q|^6}\right)\FF\nab_\That\psi\c\dual\psi.
\eeaa
    
We now show how to absorb the last three terms on the right hand side in the following.
\begin{lemma}\lab{lemma-control-rhs} 
We have, for arbitrarily small positive constants $\de_2$, $\de_3$, to be fixed later:
\bea
\begin{split}
 \left(\nab_X\psi+\frac 1 2   w \psi\right)\c  \left(\squared_2 \psi - V \psi \right)  & \geq - \de_2\frac{m}{r^2} |\nab_\Rhat \psi|^2 - \de_2\frac{\TT^2}{ r^6}  \frac{m}{r^2} |\nab_\That \psi|^2\\
 & +O(1)(|\nab_\Rhat\psi|+r^{-1}|\psi|)|N|\\
 & +O(a^3r^{-6}) |\nab_\phi \psi|^2+O(a r^{-4}) |\psi|^2 \\
&+ \D_\mu\left(\frac{2a\cos\th}{|q|^2}(\pr_r)^\mu  zh f \psi\c\nab_T\dual \psi\right)\\
& -\pr_t \left(  \frac{2a\cos\th}{|q|^2} zhf \psi\c  \nab_r\dual \psi\right)
\end{split}
\eea
and
\bea
\nn&&\left|\left(\big(\rhod +\etab\wedge\eta\big)\frac{r^2+a^2}{\De}+\frac{2a^3r\cos\th(\sin\th)^2}{|q|^6}\right)\FF\nab_\That\psi\c\dual\psi\right|\\
\nn&&+\left|\frac{2a^2r\cos\th\FF(r)}{(r^2+a^2)|q|^4}\nab_\phi\psi\c\dual\psi\right|\\
& \leq &\de_3  \frac{ \TT^2}{ r^6}  \left(\frac{m}{r^2} |\nab_\That \psi|^2+\frac{1}{r^4} |\nab_\phi \psi|^2\right) + O(a^2r^{-6}) |\psi|^2.
\eea
\end{lemma}

\begin{proof}
According to equation \eqref{eq:Gen.RW}, we have
\beaa
\left(\nab_X\psi+\frac 1 2   w \psi\right)\c  \big(\squared_2 \psi - V \psi \big) &=&\left(\nab_X\psi+\frac 1 2   w \psi\right)\c  \left(- \frac{4 a\cos\th}{|q|^2}\dual \nab_T  \psi+N \right).
\eeaa
We consider the first order term. Using \eqref{eq:new-w}, we write
\beaa
\left(\nab_X\psi+\frac 1 2   w \psi\right)\c  \left(- \frac{4 a\cos\th}{|q|^2}\dual \nab_T  \psi \right)&=&- \frac{4 a\cos\th}{|q|^2}\Big(\nab_X\psi +\frac 1 2 w_X \psi \Big) \c  \big(\dual \nab_T  \psi \big)\\
&&- \frac{4 a\cos\th}{|q|^2}\frac 1 2 w' \psi  \c  \big(\dual \nab_T  \psi \big).
\eeaa
According to Lemma \ref{identity:X+frac12wpsiN_1} and recalling that $\pr_r z=f=-\frac{2\TT}{ (r^2+a^2)^3}$,
we deduce
\beaa
&&- \frac{4 a\cos\th}{|q|^2}\left(\nab_X\psi+\frac 1 2   w_X \psi\right)\c \big(\dual \nab_T  \psi \big)\\
&=& - \frac{2a\cos\th}{|q|^2} \Bigg(  h f^2 \psi\c \nab_T\dual \psi +2zhf\dual\rho\frac{|q|^2}{\De}|\psi|^2 -\pr_r\Big(  zh f \psi\c\nab_T\dual \psi\Big)\\
 && +\pr_t\Big(  zhf \psi\c  \nab_r\dual \psi\Big)\Bigg).
\eeaa
Hence, 
\beaa
&&- \frac{4 a\cos\th}{|q|^2}\left(\nab_X\psi+\frac 1 2   w \psi\right)\c \big(\dual \nab_T  \psi \big)\\
&=& -\frac{2a\cos\th}{|q|^2}   h f^2 \psi\c \nab_T\dual \psi + \de_1\frac{4 a\cos\th}{|q|^2}  \frac{2 m \De \TT^2}{ r^2 (r^2+a^2)^4} \psi  \c  \nab_T\dual \psi -\frac{4a\cos\th}{|q|^2}zhf\dual\rho\frac{|q|^2}{\De}|\psi|^2\\
&&+ \frac{2a\cos\th}{|q|^2} \Bigg(\pr_r\Big(  zh f \psi\c\nab_T\dual \psi\Big) -\pr_t\Big(  zhf \psi\c  \nab_r\dual \psi\Big)\Bigg),
\eeaa
i.e.
\beaa
&& - \frac{4 a\cos\th}{|q|^2}\left(\nab_X\psi+\frac 1 2   w \psi\right) \c \big(\dual \nab_T  \psi \big)\\
&=&- \frac{8a\cos\th}{|q|^2} \frac{\TT^2}{r (r^2+a^2)^2(r^2-a^2)} \left(1 - \de_1 \frac{m \De}{r(r^2+a^2)} \right) \psi\c \nab_T\dual \psi\\
&&  -\frac{4a\cos\th}{|q|^2}zhf\dual\rho\frac{|q|^2}{\De}|\psi|^2\\
&&+ \frac{2a\cos\th}{|q|^2} \Bigg(\pr_r\Big(  zh f \psi\c\nab_T\dual \psi\Big) -\pr_t\Big(  zhf \psi\c  \nab_r\dual \psi\Big)\Bigg).
\eeaa
Also, notice that\footnote{Indeed, we have $\sqrt{|\g|}=\sin\th|q|^2$ and hence
    \beaa
    \D_\mu\Big(\cos\th |q|^{-2}(\pr_r)^\mu\Big) = \frac{1}{\sqrt{|\g|}}\pr_\mu(\sqrt{|\g|}\cos\th |q|^{-2}(\pr_r)^\mu)=\frac{1}{\sin\th|q|^2}\pr_r(\sin\th|q|^2\cos\th |q|^{-2})=0.
    \eeaa}  
    $\D_\mu(\cos\th|q|^{-2}(\pr_r)^\mu)=0$ which implies
    \beaa
     \frac{2a\cos\th}{|q|^2}\pr_r\Big(  zh f \psi\c\nab_T\dual \psi\Big) &=& \D_\mu\left(\frac{2a\cos\th}{|q|^2}(\pr_r)^\mu  zh f \psi\c\nab_T\dual \psi\right)
    \eeaa
    and hence
    \beaa
&& - \frac{4 a\cos\th}{|q|^2}\left(\nab_X\psi+\frac 1 2   w \psi\right) \c \big(\dual \nab_T  \psi \big)\\
&=&- \frac{8a\cos\th}{|q|^2} \frac{\TT^2}{r (r^2+a^2)^2(r^2-a^2)} \left(1 -\de_1 \frac{m \De}{r(r^2+a^2)} \right) \psi\c \nab_T\dual \psi\\
&&  -\frac{4a\cos\th}{|q|^2}zhf\dual\rho\frac{|q|^2}{\De}|\psi|^2\\
&&+ \D_\mu\left(\frac{2a\cos\th}{|q|^2}(\pr_r)^\mu  zh f \psi\c\nab_T\dual \psi\right) -\pr_t\left(\frac{2a\cos\th}{|q|^2}   zhf \psi\c  \nab_r\dual \psi\right).
\eeaa

The first two lines on the right hand side can be bounded as follows:
\beaa
&& \left|\frac{8a\cos\th}{|q|^2} \frac{\TT^2}{r (r^2+a^2)^3} \left(1 +\de_1 \frac{m \De}{r(r^2+a^2)} \right) \psi\c \nab_T\dual \psi -\frac{4a\cos\th}{|q|^2}zhf\dual\rho\frac{|q|^2}{\De}|\psi|^2\right|\\
&\leq & \frac{a}{|q|^2} \frac{\TT^2}{r (r^2+a^2)^3}  \left(\de_2 r |\nab_T \psi|^2 + \de_2^{-1} r^{-1} |\psi|^2\right) +O(ar^{-6})|\psi|^2\\
&\leq & \frac{a}{|q|^2} \frac{\TT^2}{r (r^2+a^2)^3}  \left(\de_2 r |\nab_\That \psi|^2-\de_2\frac{2ar}{r^2+a^2} \nab_\That \psi \c \nab_\phi  \psi +\de_2\frac{a^2r}{(r^2+a^2)^2} |\nab_\phi \psi|^2 \right)\\
&&+O(a r^{-4}) |\psi|^2
\eeaa
which finally gives
\beaa
\left(\nab_X\psi+\frac 1 2   w \psi\right) \c \left(- \frac{4 a\cos\th}{|q|^2}\dual \nab_T  \psi \right) &\geq& -\de_2 \frac{\TT^2}{ r^6}  \frac{m}{r^2} |\nab_\That \psi|^2 +O(a^3r^{-6}) |\nab_\phi \psi|^2\\
&&+O(ar^{-4}) |\psi|^2 +\D_\mu\left(\frac{2a\cos\th}{|q|^2}(\pr_r)^\mu  zh f \psi\c\nab_T\dual \psi\right)\\
&& -\pr_t\left( \frac{2a\cos\th}{|q|^2}  zhf \psi\c  \nab_r\dual \psi\right).
\eeaa
Since $X=\FF\pr_r=\frac{2\TT}{r(r^2 -a^2)}\Rhat=O(1) \Rhat $ and $w=O(r^{-1})$, we can bound the second product by
\beaa
\Big| \left(\nab_X\psi+\frac 1 2   w \psi\right)\c N \Big|&\les&  \Big( | \nab_\Rhat \psi | + r^{-1}|\psi| \Big)    |  N|. 
\eeaa
By putting together with the previous bound we obtain the first desired estimate. 

  Next, notice that 
\beaa
|\rhod\,|\les \frac{am}{r^4},\qquad |\eta|+|\etab|\les \frac{a}{r^2}, \qquad |\FF| \les \frac{|\De|}{r^2}\frac{|\TT|}{r^3}.
\eeaa  
We therefore deduce
\beaa
\nn&&\left|\left(\big(\rhod +\etab\wedge\eta\big)\frac{r^2+a^2}{\De}+\frac{2a^3r\cos\th(\sin\th)^2}{|q|^6}\right)\FF\nab_\That\psi\c\dual\psi\right|\\
\nn&&+\left|\frac{2a^2r\cos\th\FF(r)}{(r^2+a^2)|q|^4}\nab_\phi\psi\c\dual\psi\right|\\
& \leq & \frac{am}{r^7} | \TT| |\nab_\That\psi|\,|\psi| +\frac{am}{r^8} | \TT| |\nab_\phi\psi|\,|\psi|
\eeaa
and hence 
\beaa
\nn&&\left|\left(\big(\rhod +\etab\wedge\eta\big)\frac{r^2+a^2}{\De}+\frac{2a^3r\cos\th(\sin\th)^2}{|q|^6}\right)\FF\nab_\That\psi\c\dual\psi\right|\\
\nn&&+\left|\frac{2a^2r\cos\th\FF(r)}{(r^2+a^2)|q|^4}\nab_\phi\psi\c\dual\psi\right|\\
& \leq &\de_3  \frac{ \TT^2}{ r^6}  \left(\frac{m}{r^2} |\nab_\That \psi|^2+\frac{1}{r^4} |\nab_\phi \psi|^2\right) + O(a^2r^{-6}) |\psi|^2
\eeaa
as stated. 
\end{proof}

By putting together Proposition \ref{prop:positivity-bulk} and Lemma \ref{lemma-control-rhs}, 
we obtain from \eqref {definition-EE-X=FFR} 
\beaa
\D^\mu  \PP_\mu[X, w]&=& \EE[X, w] +\left(\nab_X\psi+\frac 1 2   w \psi\right)\c  \big(\squared_2 \psi - V \psi \big) -\frac{2a^2r\cos\th\FF(r)}{(r^2+a^2)|q|^4}\nab_\phi\psi\c\dual\psi\\
&& -\left(\big(\rhod +\etab\wedge\eta\big)\frac{r^2+a^2}{\De}+\frac{2a^3r\cos\th(\sin\th)^2}{|q|^6}\right)\FF\nab_\That\psi\c\dual\psi\\
&\geq & \de_*\frac{m}{r^2}   | \nab_\Rhat \psi|^2+  \left( \de_*-\de_2-\de_3\right)\frac{ \TT^2}{r^6} \left( \frac{m}{ r^2} |\nab_\That \psi|^2+ r^{-1}|\nab\psi|^2\right)\\
&& -O(a^2 r^{-4})  |\nab_\phi \psi|^2        -O(r^{-3})  |\psi|^2-O(1)\Big( | \nab_\Rhat \psi | + r^{-1}|\psi| \Big)    |  N|\\
&& +\D_\mu\left(\frac{2a\cos\th}{|q|^2}(\pr_r)^\mu  zh f \psi\c\nab_T\dual \psi\right) -\pr_t \left(  \frac{2a\cos\th}{|q|^2} zhf \psi\c  \nab_r\dual \psi\right).
\eeaa
By choosing $\de_2$ and $\de_3$ sufficiently small, integrating the above inequality on the region $\MM(\tau_1, \tau_2) $ and applying the divergence theorem we deduce
\beaa
&&\int_{\MM(\tau_1, \tau_2)}  \frac{m}{r^2} |\nab_{\Rhat} \psi|^2 +\frac{\TT^2}{r^6} \left(\frac{m}{r^2} |\nab_\That \psi|^2 + r^{-1}|\nab\psi|^2\right)   \\
\les &&\int_{\pr\MM(\tau_1, \tau_2)}|M(\psi)|+\int_{\MM(\tau_1, \tau_2)}      \big(  r^{-4}  a^2 |\nab_Z\psi|^2 + r^{-3}|\psi|^2\big)
 +\int_{\MM(\tau_1, \tau_2)}\Big( | \nab_\Rhat \psi | + r^{-1}|\psi| \Big)    |  N|
\eeaa
where
\beaa
M(\psi):= \PP\c N+O(a r^{-2}) zhf \psi\c  \nab_T\dual \psi +O(a r^{-2}) zhf \psi\c  \nab_r\dual \psi.
\eeaa
Recalling that $X=\frac{2\TT}{r(r^2+a^2)}\Rhat =O(1) \Rhat$, $w=O(r^{-1}) $, the properties of $N_\Si$ in Definition \ref{definition:definition-oftau} and
 \beaa
 \PP_\mu &=&\PP_\mu[X, w, M]=\QQ_{\mu\nu} X^\nu +\frac 1 2  w \psi \c \Db_\mu \psi -\frac 1 4|\psi|^2   \pr_\mu w +\frac 1 4 |\psi|^2 M_\mu,
  \eeaa
 we easily deduce 
 \beaa
 \int_{\pr\MM(\tau_1, \tau_2)}|M(\psi)| &\les&  \sup_{[\tau_1, \tau_2]}E_{deg}[\psi](\tau)+\deh F_{\AA}[\psi](\tau_1, \tau_2)+F_{\Si_*}[\psi](\tau_1, \tau_2).
\eeaa
This ends the proof of estimate \eqref{eq:conditional-mor-par1-1-I} of Proposition  \ref{proposition:Morawetz1-step1}.


\subsection{Proof of Proposition  \ref{proposition:Morawetz1-step1}, estimate \eqref{eq:conditional-mor-par1-1-II}}\label{section:proof-hardy-poincare}


In this section we provide the proof for the   estimate in the second part of Proposition  \ref{proposition:Morawetz1-step1}.

The main difference between estimate \eqref{eq:conditional-mor-par1-1-I} and estimate \eqref{eq:conditional-mor-par1-1-II} in Proposition \ref{proposition:Morawetz1-step1} is in the control on the left hand side of the zero-th order term $|\psi|^2$. To extend the estimate to control such lower order term we will make use of a Poincar\'e inequality and the one-form $M$ in the definition of the current through a Hardy estimate.


\subsubsection{Poincar\'e inequality}


Recall that   according to Proposition \ref{proposition:Morawetz1}, for  choices   of functions $z, f, h$ and $v$ in the definition of $M= v\pr_r$,  the generalized current associated to the Morawetz vectorfield in \eqref{identity:prop.Morawetz1} is given by
 \beaa
|q|^2 \EE[X, w, M]&=&\AA|\nab_r\psi|^2+ P+\VV |\psi|^2 +\frac 1 4 |q|^2  \D^\mu (|\psi|^2 M_\mu).
 \eeaa

 Here we consider the choices for $z, f, h$ made in Proposition \ref{prop:Choice-zhf}   and show that there exists a choice of  $v$ 
 such that the above expression is positive definite. 
\bigskip

With the choices of $z, f, h$ mentioned above,  the coefficients $\AA$ and $\VV$ are given by \eqref{eq:prop-Choice-AA} and \eqref{eq:prop:Choice-zhf} respectively, and  the principal term $P$ given by \eqref{eq:principal-term}, i.e.
 \beaa
 P&=& \UU^{\a\b}(\Db_\a \psi )\c(\Db_\b \psi )=\frac{\TT}{r}\frac{r^2+a^2}{r^2-a^2} \left(  \frac{2\TT}{ (r^2+a^2)^3}   O^{\a\b}\nab_\a\psi\c\nab_\b \psi - \frac{4ar}{(r^2+a^2)^2} \nab_\That\psi\c\nab_\phi \psi\right).
 \eeaa
We rewrite $P$ as
\beaa
P = r^2H|\nab\psi|^2 -O(a)\big(r|\nab\psi|^2+r^{-1}|\nab_t\psi|^2\big),\qquad H := \frac{2\TT^2}{ r(r^2+a^2)^2(r^2-a^2)}.
\eeaa
To obtain a lower bound for $\int_SP$, we will rely on the following Poincar\'e inequality.
\begin{lemma}\lab{lemma:poincareinequalityfornabonSasoidfh:chap6}
For $\psi\in\sk_2$, we have 
\beaa
\int_S|\nab\psi|^2 &\geq& \frac{2}{r^2}\int_S|\psi|^2  -O(a)\int_S\big(|\nab\psi|^2+r^{-2}|\nab_t\psi|^2+r^{-4}|\psi|^2\big).
\eeaa
\end{lemma}

\begin{proof}
Denoting by $\nab^S$ the covariant derivative for the induced metric on $S$, one easily checks 
\beaa
\nab^S &=& \big(1+O(a^2r^{-2})\big)\nab -\frac{a\sin\th}{r}\big(1+O(a^2r^{-2})\big)\nab_t
\eeaa
and hence 
\beaa
|\nab^S\psi|^2=|\nab\psi|^2-O(a)\big(|\nab\psi|^2+r^{-2}|\nab_t\psi|^2\big)
\eeaa
so that 
\beaa
\int_S|\nab^S\psi|^2 &=& \int_S|\nab\psi|^2 -O(a)\int_S\big(|\nab\psi|^2+r^{-2}|\nab_t\psi|^2\big).
\eeaa
Also, note that in the three coordinates systems of Remark \ref{rmk:whatthesecoordinatessystemsonSoftaurareinKerrandrefmainKerr:chap6}, we have, in view of Lemma 2.4.10 and Lemma 2.4.24 in \cite{KS:Kerr}
\beaa
\max_{b,c=1,2}\big|\pr^{\leq 2}((g_S)_{x^ax^b} - r^2(\ga_{\SSS^2})_{x^ax^b})\big| &\les& a^2
\eeaa
where $\pr^{\leq 2}$ denotes at most 2 coordinates derivatives, and $(\ga_{\SSS^2})_{x^ax^b}$ the metric coefficients on $\SSS^2$ in the corresponding coordinates system. We deduce in particular 
\beaa
K_S &=& \frac{1}{r^2}\big(1+O(a^2r^{-2})\big), \qquad r_S=r\big(1+O(a^2r^{-2})\big).
\eeaa
Applying the effective uniformization result of Corollary 3.8 in \cite{KS-GCM2}, we obtain a map $\Phi:\SSS^2\to S$ and a scalar function $u$ on $\SSS^2$ such that 
\beaa
\Phi^\#(g_S) &=& (r_S)^2e^{2u}\ga_{\SSS^2}, \qquad \|\dkb^{\leq 2}(u\circ\Phi^{-1})\|_{L^2(S)}\les (a^2r^{-2})r_S.
\eeaa
We infer, relying on the well known Poincar\'e inequality for $\psi\in\sk_2(\SSS^2)$, see for example \cite{GKS1}, 
\beaa
\int_S|\nab^S\psi|^2 &=& \frac{1}{r^2}\big(1+O(a^2r^{-2})\big)\int_{\SSS^2}|\nab_{\SSS^2}\Phi^{\#}\psi|^2+ \frac{1}{r^2}O(a^2r^{-2})\int_{\SSS^2}|\Phi^{\#}\psi|^2\\
&\geq& \frac{2}{r^2}\big(1+O(a^2r^{-2})\big)\int_{\SSS^2}|\Phi^{\#}\psi|^2\\
&\geq&  \frac{2}{r^2}\int_S\big(1+O(a^2r^{-2})\big)|\psi|^2  
\eeaa
and hence
\beaa
\int_S|\nab\psi|^2 &\geq& \frac{2}{r^2}\int_S|\psi|^2  -O(a)\int_S\big(|\nab\psi|^2+r^{-2}|\nab_t\psi|^2+r^{-4}|\psi|^2\big)
\eeaa
as stated.
\end{proof}

Next, we define, for $\de>0$ sufficiently small chosen later, the following quadratic form
  \bea
  \lab{def:Qr_red:new}
  \bsplit
 \mbox{Qr}_\de[\psi]  = & (1-\de)\AA |\nab_r\psi|^2   + \left(\VV+  (1-\de)\frac{2\TT^2}{(r^2+a^2)^2(r^2-a^2)}\right) |\psi|^2 \\
 & +\frac 1 4 |q|^2  \D^\mu (|\psi|^2 M_\mu)
 \end{split}
 \eea
 so that, in view of the above, we have
  \beaa
|q|^2 \EE[X, w, M]&=& \AA|\nab_r\psi|^2+ P+\VV |\psi|^2 +\frac 1 4 |q|^2  \D^\mu (|\psi|^2 M_\mu)\\
&=& \de\AA|\nab_r\psi|^2+\de P + (1-\de)\left(P -\frac{2\TT^2}{r(r^2+a^2)^2(r^2-a^2)}|\psi|^2\right) +\mbox{Qr}_\de[\psi].
 \eeaa 
Now, since
\beaa
P = \frac{2r\TT^2}{(r^2+a^2)^2(r^2-a^2)}|\nab\psi|^2 -O(a)\big(r|\nab\psi|^2+r^{-1}|\nab_t\psi|^2\big),
\eeaa
we infer, in view of Lemma \ref{lemma:poincareinequalityfornabonSasoidfh:chap6}, 
\beaa
&&\int_S\frac{1}{|q|^2}\left(\de P + (1-\de)\left(P -\frac{2\TT^2}{ r(r^2+a^2)^3}|\psi|^2\right)\right)\\
 &\gtrsim& \de\frac{\TT^2}{r^7}\int_S|\nab\psi|^2 -O(ar^{-1})\int_S\big(|\nab\psi|^2+r^{-2}|\nab_t\psi|^2+r^{-2}|\psi|^2\big)
\eeaa
and hence
 \beaa
\int_S\EE[X, w, M] &\gtrsim& \de\int_S\left(\frac{\AA}{|q|^2}|\nab_r\psi|^2+\frac{\TT^2}{r^7}|\nab\psi|^2\right)+\int_S\frac{1}{|q|^2}\mbox{Qr}_\de[\psi] \\
&& -O(ar^{-1})\int_S\big(|\nab\psi|^2+r^{-2}|\nab_t\psi|^2+r^{-2}|\psi|^2\big),
 \eeaa 
or
 \bea\lab{eq:intermediarylowerboundEEafterPoincaresfasuighalsuidh}
\nn\int_S\EE[X, w, M] &\geq& \de\int_S\left(\frac{m}{r^2}|\nab_{\Rhat}\psi|^2+\frac{\TT^2}{r^7}|\nab\psi|^2\right)+\int_S\frac{1}{|q|^2}\mbox{Qr}_\de[\psi] \\
&& -O(ar^{-1})\int_S\big(|\nab\psi|^2+r^{-2}|\nab_t\psi|^2+r^{-2}|\psi|^2\big).
 \eea

  
 \subsubsection{The Hardy inequality}
 
  
Note that, using \eqref{expression-Div-M-I} and \eqref{def:Qr_red:new}, we have
  \beaa
  \bsplit
 \mbox{Qr}_\de[\psi]  = & (1-\de)\AA |\nab_r\psi|^2   + \left(\VV+  (1-\de)\frac{2\TT^2}{(r^2+a^2)^2(r^2-a^2)}\right) |\psi|^2 \\
 & + \frac 1 4 |q|^2\left( 2 v(r)\psi\c \nab_r \psi + \left(\pr_r v+ \frac{2r}{|q|^2} v\right) |\psi|^2 \right).
 \end{split}
 \eeaa
In view of \eqref{eq:intermediarylowerboundEEafterPoincaresfasuighalsuidh}, it remains to derive a lower bound for  the term  $\mbox{Qr}_\de[\psi]$ which is done in the following lemma.

 \begin{lemma}[Half Poincar\'e + Hardy]\lab{Le:Crucialpositivity-a=0} 
 For $|a|/m\ll 1$, there exists a function $v(r)$ with $v(r)= O( m^{1/2}\Delta r^{-9/2})$, and $\de>0$ sufficiently small,  such that  for all $r\ge r_+(1-\deh)$,
 \beaa
\mbox{Qr}_\de[\Phi]  &\geq& O(\de)  \left( \frac{m \De^2}{r^4} \big|\nab_r\Phi|^2 + r^{-1} |\Phi|^2 \right).
 \eeaa
 \end{lemma}
 
 \begin{proof}
 We have
  \beaa
  \bsplit
 \mbox{Qr}_\de[\Phi]  = & (1-\de)\AA \left|\nab_r\Phi+\frac{|q|^2}{4(1-\de)\AA}v(r)\Phi\right|^2 -\frac{|q|^4}{16(1-\de)\AA}v^2|\Phi|^2\\
 &  + \left(\VV+  (1-\de)\frac{2\TT^2}{(r^2+a^2)^2(r^2-a^2)}\right) |\Phi|^2  + \frac 1 4 |q|^2 \left(\pr_r v+ \frac{2r}{|q|^2} v\right) |\Phi|^2 
 \end{split}
 \eeaa
 and hence
 \beaa
  \bsplit
 \mbox{Qr}_\de[\Phi]  \geq & \left(\VV+  (1-\de)\frac{2\TT^2}{(r^2+a^2)^2(r^2-a^2)} + \frac 1 4 |q|^2 \left(\pr_r v+ \frac{2r}{|q|^2} v\right) -\frac{|q|^4}{16(1-\de)\AA}v^2\right)|\Phi|^2 
 \end{split}
 \eeaa
 Thus, it suffices to prove that  for $|a|/m \ll 1$ and $\de>0$ small enough, there exists a function $v$ such that 
  \bea
 \lab{condition-crucial-2}
E:=  \VV+  \frac{2(1-\de)\TT^2}{r(r^2+a^2)^2(r^2-a^2)} + \frac 1 4 |q|^2\left(\pr_r v+ \frac{2r}{|q|^2} v\right)-\frac{1}{16(1-\de)\AA}|q|^4  v^2>0.
 \eea

By continuity it remains to prove that there exists a function $v(r)$ for which the condition \eqref{condition-crucial-2} is valid  for $a=0$ and $\de=0$. 
 For $a=0$, i.e. for Schwarzschild spacetime, we have, with $\Up:=1-\frac{2m}{r}$, 
 \beaa
 \VV&=&\VV_0+\VV_1,\\
\VV_1&=&8  \De \frac{r^3-4mr^2}{r^9}  (r^3-3 mr^2) = 8 r^{-1} \Up \left(1-\frac{4m}{r} \right)    \left(1-\frac{3m}{r} \right),    \\
\VV_0&=&  \frac{9mr^6-46m^2r^5+54m^3r^4}{r^8}= \frac{9mr^2-46m^2r+54m^3}{r^4},\\
\VV&=&  \frac{8r^3-63mr^2+162m^2r-138m^3}{r^4}, \\
\AA&=& \frac{2m \De^2}{r^8}  3 r^4 = 6m \Up^2.
 \eeaa
 The expression $ E$ in \eqref{condition-crucial-2} becomes, setting $\widetilde{v}= r^2 v$,
 \beaa
 E&=&\VV+ \frac{2}{r} \left( 1-\frac{3m}{r} \right)^2 + \frac 1 4  r^2\left(\pr_r v+ \frac{2r}{ r ^2} v\right)-\frac{1}{16\c 6 m }    \Up^{-2}  r^4  v^2\\
&=&  \VV+ \frac{2}{r}\left( 1-\frac{3m}{r} \right)^2 + \frac 1 4 \pr_r(\widetilde{v} ) -\frac{1}{96 m }    \Up^{-2}    \widetilde{v}^2\\
&=&\frac{8r^3-63mr^2+162m^2r-138m^3}{r^4} + \frac 2 r\left( 1-\frac{3m}{r} \right)^2 + \frac 1 4 \pr_r \widetilde{v} -\frac{1}{96 m }    \Up^{-2}     \widetilde{v}^2\\
&=& \frac{10r^3-75mr^2+180m^2r-138m^3}{r^4}   + \frac 1 4 \pr_r \widetilde{v}  -\frac{1}{96 m }    \Up^{-2}    \widetilde{v}^2.
 \eeaa
 Introducing $x=\frac{r}{2m} $   and assuming that  $\widetilde{v}= \widetilde{v}\big(\frac{r}{2m} \big)= \widetilde{v}_0(x)$, we derive
  \beaa
 E&=& \frac{10(2mx)^3-75m(2mx)^2+180m^2(2mx)-138m^3}{(2mx)^4}  +  \frac{1}{8m} \widetilde{v}_0' -\frac{1}{96 m }  \frac{x^2}{(x-1)^2}   \widetilde{v}_0^2\\
 &=& \frac{40 x^3-150x^2+180x-69}{8m x^4}  +  \frac{1}{8m} \widetilde{v}_0' -\frac{1}{96 m }  \frac{x^2}{(x-1)^2}   \widetilde{v}_0^2.
 \eeaa
 By setting $\widetilde{v}_0(x)= (x-1) k_0(x)$, and $\widetilde{v}_0'= k_0+ (x-1) k_0'$, then 
  \beaa
8m  E &=&40 x^{-1}-150x^{-2}+180x^{-3}-69 x^{-4}  + k_0+ (x-1) k_0' -\frac{1}{12 } x^2  k_0^2.
 \eeaa
 In order to have that $ x^4 k^2_0 $ with same degree in $x$ than $x^{-1}$ we assume
\beaa
k_0(x)= A x^{-3/2} , \qquad k_0'(x)=- \frac 3 2 A x^{-5/2}.
\eeaa
This yields
   \beaa
8m  E &=&40 x^{-1}-150x^{-2}+180x^{-3}-69 x^{-4}  + A x^{-3/2}- \frac 3 2 A  (x-1) x^{-5/2}-\frac{1}{12 } x^2  A^2 x^{-3}\\
&=&\left( 40 -\frac{1}{12 }   A^2\right)  x^{-1} -\frac 1 2  A x^{-3/2}-150x^{-2}+ \frac 3 2 A   x^{-5/2}+180x^{-3}-69 x^{-4} .
 \eeaa
For example for $A=2$, the above is positive for $x \geq 1$. This results in 
\beaa
 E &=& \frac{119}{12 }  \frac{1}{r} - \frac{1}{8m}\left(\frac{2m}{r}\right)^{3/2}-150\frac{1}{8m}\left(\frac{2m}{r}\right)^{2}+  3\frac{1}{8m}    \left(\frac{2m}{r}\right)^{5/2}+180\frac{1}{8m}\left(\frac{2m}{r}\right)^{3}\\
 && -69\frac{1}{8m} \left(\frac{2m}{r}\right)^{4}.
\eeaa
 Note also that in this case
 \beaa
v=r^{-2} \widetilde{v}= 2\frac{(2m)^{3/2}}{r^{7/2}}\left(\frac{r}{2m}-1\right),
 \eeaa
where $v$ satisfies $v(r)= O( m^{1/2}\Delta r^{-9/2})$. This ends the proof of Lemma \ref{Le:Crucialpositivity-a=0}. 
 \end{proof}

 \begin{proposition}
 \lab{Prop.StongversionMorawetz1}
  There exists a  choice of $v(r)= O\big( m^{1/2}\Delta r^{-9/2} \big) $  and a small universal constant $c_0>0$ such that, for $|a|/m\ll 1$ and  $r\ge r_+(1-\deh)$,
  \beaa
\nn\int_S\EE[X, w, M] &\geq& c_0\int_S\left(\frac{m}{r^2}|\nab_{\Rhat}\psi|^2+\frac{\TT^2}{r^7}|\nab\psi|^2  + r^{-3} |\psi|^2\right)\\
&&  -O(ar^{-1})\int_S\big(|\nab\psi|^2+r^{-2}|\nab_t\psi|^2\big).
 \eeaa
  \end{proposition}
 
 \begin{proof}
 Immediate consequence of \eqref{eq:intermediarylowerboundEEafterPoincaresfasuighalsuidh} and Lemma \ref{Le:Crucialpositivity-a=0}. 
 \end{proof}

 We can proceed as in Section \ref{section:lowerboundscontainigThat} to derive a version of Proposition \ref{Prop.StongversionMorawetz1} which also contains  $\nab_\That \psi$. We obtain the following proposition. 
   \begin{proposition}
 \lab{Prop.StongversionMorawetz1-withT}
  There exists a  choice of $v(r)= O\big( m^{1/2}\Delta r^{-9/2} \big) $, a redefinition of $w$   and a small universal  constant $c_0>0$ such that, for all $a$ verifying $|a|/m\ll 1$ and  all $r\ge r_+(1-\deh)$,
   \beaa
  \int_S\EE[X, w, M]&\ge &  c_0\int_S\left(   \frac{m}{r^2}  |\nab_\Rhat\psi|^2+   \frac{\TT^2 }{ r^6} \left(\frac{m}{r^2} |\nab_\That\psi|^2+        r^{-1}  |\nab \psi|^2\right) + r^{-3} |\psi|^2 \right)\\  
 && -O(ar^{-1})\int_S\big(|\nab\psi|^2+r^{-2}|\nab_t\psi|^2\big).
\eeaa
 \end{proposition}

  
 \subsubsection{End of the proof}
 
 
 We now proceed as in Section \ref{section:Proof-proposition:Morawetz-Energy1} and we apply Lemma \ref{lemma-control-rhs} to control the right hand side in the divergence of $\PP[X, w, M]$. Since the zero-th order terms on the right hand side appear with higher decay in $r$ and are multiplied by the angular momentum $a$ they can be absorbed by the positive Morawetz bulk.

 By applying the divergence theorem we then obtain
  the following estimate 
 \beaa
&&\int_{\MM(\tau_1, \tau_2)}  \frac{m}{r^2} |\nab_{\Rhat} \psi|^2+r^{-3}|\psi|^2 +\frac{\TT^2}{r^6} \left(\frac{m}{r^2} |\nab_\That \psi|^2 + r^{-1}|\nab\psi|^2\right)    \\
\les &&\int_{\pr\MM(\tau_1, \tau_2)}|\PP\c N_\Si|+\int_{\MM(\tau_1, \tau_2)}ar^{-1}\big(|\nab\psi|^2+r^{-2}|\nab_t\psi|^2\big)\\
&& +\int_{\MM(\tau_1, \tau_2)}\Big(|\nab_\Rhat\psi|+r^{-1}|\psi|\Big)| N|
\eeaa
with the same bounds as before for $M(\psi):= \PP\c N_{\Si}$, which is precisely  \eqref{eq:conditional-mor-par1-1-II}, and hence concludes the proof of Proposition  \ref{proposition:Morawetz1-step1}.

     
\section{Energy  estimates} 
\lab{section:generalized.energy}\label{subsection:energy-identity}


We start with the following lemma.
\begin{lemma}
\lab{Lemma:QQ(That,NSi)}
The following hold true with a sufficiently small $c_0>0$, for any $|a| \ll m$,
\bea
\lab{eq:Lemma-forQQ3}
\bsplit
\int_{\Si(\tau)}\QQ(\That, N_\Si)&\ge   c_0E_{deg}[\psi](\tau)- O(\deh) E_{r\leq r_+}[\psi](\tau),\\
\int_{\Si_*(\tau_1, \tau_2)}\QQ(\That, N_{\Si_*}) &\ge c_0F_{\Si_*}[\psi](\tau_1, \tau_2),\\
\int_{\AA(\tau_1, \tau_2)}\QQ(\That, N_{\AA}) &\gtrsim - \deh F_{\AA}[\psi](\tau_1, \tau_2).
 \end{split}
\eea
\end{lemma}

\begin{proof}
We have
\beaa
\QQ(e_4, N_\Si)&=&\QQ\left( e_4, \frac 1 2 e_3(\tau) e_4 +\frac 1 2 e_4(\tau) e_3- \nab(\tau)\right)\\
&=&  \frac 1 2e_3(\tau)\QQ_{44} +\frac 1 2 e_4(\tau)\QQ_{34} -\nab^a(\tau)\QQ_{4a}, \\
\QQ(e_3, N_\Si)
&=&  \frac 1 2e_3(\tau)\QQ_{34} +\frac 1 2 e_4(\tau)\QQ_{33} -\nab^a(\tau)\QQ_{3a}.
\eeaa
In view of \eqref{eq:definition-QQ-mu-nu}, we have
\beaa
\bsplit
\QQ_{33}&=|\nab_3\psi|^2, \qquad  \QQ_{44}=|\nab_4\psi|^2, \qquad \QQ_{34}=|\nab\psi|^2 +V|\psi|^2,\\
\QQ_{4a}&=\nab_4\Psi\c\nab_a\Psi, \qquad \QQ_{3a}=\nab_3\Psi\c\nab_a\Psi.
\end{split}
\eeaa
We infer
\beaa
\QQ(e_4, N_\Si)&\geq&  \frac 1 2e_3(\tau)|\nab_4\psi|^2 +\frac 1 2 e_4(\tau)\Big(|\nab\psi|^2 +V|\psi|^2\Big) -|\nab \tau||\nab_4\psi||\nab\psi|,\\
\QQ(e_3, N_\Si)
&\geq &  \frac 1 2e_3(\tau)\Big(|\nab\psi|^2 +V|\psi|^2\Big) +\frac 1 2 e_4(\tau)|\nab_3\psi|^2 -|\nab \tau||\nab_3\psi||\nab\psi|.
\eeaa

Next, since we have in view of the choice of $\tau$
\beaa
e_4(\tau)>0, \qquad e_3(\tau)>0, \qquad |\nab\tau|^2  \leq \frac{8}{9}e_4(\tau)e_3(\tau),
\eeaa
we infer
\beaa
|\nab \tau||\nab_4\psi||\nab\psi| &\leq& \sqrt{\frac{8}{9}e_4(\tau)e_3(\tau)}|\nab_4\psi||\nab\psi|\leq \sqrt{\frac{8}{9}}\left( \frac 1 2e_3(\tau)|\nab_4\psi|^2 +\frac 1 2 e_4(\tau)|\nab\psi|^2\right),\\
|\nab \tau||\nab_3\psi||\nab\psi| &\leq&\sqrt{\frac{8}{9}e_4(\tau)e_3(\tau)}|\nab_3\psi||\nab\psi|\leq\sqrt{\frac{8}{9}}\left(\frac 1 2e_3(\tau)|\nab\psi|^2  +\frac 1 2 e_4(\tau)|\nab_3\psi|^2\right),
\eeaa
and thus
\beaa
\QQ(e_4, N_\Si)&\geq& \left(1-\sqrt{\frac{8}{9}}\right)\left( \frac 1 2e_3(\tau)|\nab_4\psi|^2 +\frac 1 2 e_4(\tau)|\nab\psi|^2\right)+ \frac 1 2 e_4(\tau)V|\psi|^2,\\
\QQ(e_3, N_\Si)
&\geq &  \left(1-\sqrt{\frac{8}{9}}\right)\left(\frac 1 2e_3(\tau)|\nab\psi|^2 +\frac 1 2 e_4(\tau)|\nab_3\psi|^2\right)+ \frac 1 2 e_3(\tau)V|\psi|^2.
\eeaa
Using the fact that we have\footnote{Recall that we consider the explicit potential $V=\frac{4\De}{ (r^2+a^2) |q|^2}$, see \eqref{eq:Gen.RW}, which indeed satisfies $V\geq -O(\deh)\mathbb{1}_{r\leq r_+}$.}  on $\MM$
\beaa
e_4(\tau)\gtrsim \frac{m^2}{r^2}, \qquad e_3(\tau)\gtrsim 1, \qquad V\geq -O(\deh)\mathbb{1}_{r\leq r_+},
\eeaa
we infer the existence of a constant $c_0>0$ such that, on $\MM$, 
\beaa
\QQ( e_4, N_\Si)  &\ge & c_0\Big( |\nab_4\psi|^2+ r^{-2}|\nab\psi|^2\Big) -O(\deh)|\psi|^2\mathbb{1}_{r\leq r_+},\\
\QQ( e_3, N_\Si)  &\ge & c_0\Big(r^{-2}  |\nab_3\psi|^2+  |\nab\psi|^2\Big)  -O(\deh)|\psi|^2\mathbb{1}_{r\leq r_+}.
\eeaa
We infer on $\MM$
\beaa
\QQ(\That, N_\Si) &=&\frac 1 2  \frac{|q|^2}{r^2+a^2} \QQ( e_4, N_\Si)+\frac 1 2 
 \frac{\De}{r^2+a^2}  \QQ( e_3, N_\Si)\\
 &\ge & c_0\left( |\nab_4\psi|^2 + \frac{|\De|}{r^4} |\nab_3\psi|^2 +|\nab\psi|^2\right) - O(\deh)\Big(|\nab_3\psi|^2+|\psi|^2\Big)\mathbb{1}_{r\leq r_+}.
\eeaa
Similarly, using $\De\gtrsim r^{-2}$ and $V\gtrsim r^{-2}$ on $\Si_*$,  
\beaa
\QQ(\That, N_{\Si_*})  &\ge & c_0\Big( |\nab_4\psi|^2 +  |\nab_3\psi|^2 +|\nab\psi|^2 +r^{-2}|\psi|^2\Big).
\eeaa
Also, we have on $\AA$
\beaa
\QQ(\That, N_{\AA})  \ge  -O(\deh)\QQ_{34}-O(\deh^2)\QQ_{33}\ge -O(\deh)|\nab\psi|^2-O(\deh^2)(|\nab_3\psi|^2+|\psi|^2).
\eeaa
This yields 
\beaa
\bsplit
\int_{\Si(\tau)}\QQ(\That, N_\Si)&\ge   c_0\int_{\Si(\tau)}\left( |\nab_4\psi|^2 + \frac{|\De|}{r^4} |\nab_3\psi|^2 +|\nab\psi|^2\right) - O(\deh) E_{r\leq r_+}[\psi](\tau),\\
\int_{\Si_*(\tau_1, \tau_2)}\QQ(\That, N_{\Si_*}) &\ge c_0F_{\Si_*}[\psi](\tau_1, \tau_2),\\
\int_{\AA(\tau_1, \tau_2)}\QQ(\That, N_{\AA}) &\gtrsim - \deh F_{\AA}[\psi](\tau_1, \tau_2).
 \end{split}
\eeaa
In particular, we have obtained the desired estimates on $\Si_*$ and $\AA$.

Also, we have in view of Lemma \ref{lemma:poincareinequalityfornabonSasoidfh:chap6}, for $\psi\in\sk_2$ and $|a|\ll m$, 
\beaa
\frac{1}{r^2}\int_S|\psi|^2 &\les& \int_S\big(|\nab\psi|^2+r^{-2}|\nab_t\psi|^2)\\
&\les& \int_S\left( |\nab_4\psi|^2 + \frac{|\De|}{r^4} |\nab_3\psi|^2 +|\nab\psi|^2\right). 
\eeaa
Together with the above, we infer
\beaa
\int_{\Si(\tau)}\QQ(\That, N_\Si)&\ge&   c_0\int_{\Si(\tau)}\left( |\nab_4\psi|^2 + \frac{|\De|}{r^4} |\nab_3\psi|^2 +|\nab\psi|^2+r^{-2}|\psi|^2\right) \\
&& - O(\deh) E_{r\leq r_+}[\psi](\tau),\\
\eeaa
and hence
\beaa
\int_{\Si(\tau)}\QQ(\That, N_\Si)&\ge   c_0E_{deg}[\psi](\tau)- O(\deh) E_{r\leq r_+}[\psi](\tau)
\eeaa
as stated. This concludes the proof of Lemma \ref{Lemma:QQ(That,NSi)}.
\end{proof}

We consider  the energy current associated  to the modified timelike vectorfield $\That_\de=\pr_t +\chi_\de(r) \pr_\phi$, as defined in Definition \ref{definition:vf-Tmod}, with $\de=\frac{1}{10}$ and $|a|/m\ll 1$ small enough. Recall that $\chi_\de= \frac{a}{r^2+a^2} \chi_0\Big( \de^{-1} \frac{\TT}{r^3} \Big)$, with $\chi_0=0$ in $\MM_{trap}$.

From   Proposition \ref{prop-app:stadard-comp-Psi}, we have for the current associated to $\That_\de$:
  \beaa
  \begin{split}
  \D^\mu  \PP_\mu[\That_\de, 0, 0] &= \frac 1 2 \QQ  \c \,^{(\That_\de)} \pi +\That_\de^\mu \Db^\nu  \psi ^a\Rdot_{ ab   \nu\mu}\psi^b+ \nab_{\That_\de} \psi  \c \big(\squared_2 \psi - V \psi \big)
  \end{split}
 \eeaa
 and hence
  \beaa
  \D^\mu  \PP_\mu[\That_\de, 0, 0] &=& \frac 1 2 \QQ  \c \,^{(\That_\de)} \pi   +\pr_t^\mu \Db^\nu  \psi ^a\Rdot_{ ab   \nu\mu}\psi^b\\
&&   +(\That_\de-\pr_t)^\mu \Db^\nu  \psi ^a\Rdot_{ ab   \nu\mu}\psi^b+ \nab_{\That_\de} \psi  \c \big(\squared_2 \psi - V \psi \big).
 \eeaa
 Since $\Rdot_{ ab   \nu\mu}$ is antisymmetric with respect to $(a, b)$, we rewrite 
  \beaa
  \D^\mu  \PP_\mu[\That_\de, 0, 0] &=& \frac 1 2 \QQ  \c \,^{(\That_\de)} \pi   +\frac{1}{2}\pr_t^\mu \in^{ab}\Rdot_{ ab   \nu\mu} \dual\psi\c\Db^\nu\psi\\
&&   +(\That_\de-\pr_t)^\mu \Db^\nu  \psi ^a\Rdot_{ ab   \nu\mu}\psi^b+ \nab_{\That_\de} \psi  \c \big(\squared_2 \psi - V \psi \big).
 \eeaa 
 Introducing the following spacetime 1-form
\bea\lab{eq:thespacetime1formAappearingintheenergyestimate}
A_\mu &:=& \in^{bc}\Rdot_{bc \mu\nu}\pr_t^\nu,
\eea 
 we infer
 \bea\lab{eq:basiccomputationenergyestimatesKerrwithThatdeltaflsdiuh}
 \nn \D^\mu  \PP_\mu[\That_\de, 0, 0] &=& \frac 1 2 \QQ  \c \,^{(\That_\de)} \pi   +\frac{1}{2}A_\nu \dual\psi\c\Db^\nu\psi  +(\That_\de-\pr_t)^\mu \Db^\nu  \psi ^a\Rdot_{ ab   \nu\mu}\psi^b\\
  &&+ \nab_{\That_\de} \psi  \c \big(\squared_2 \psi - V \psi \big).
 \eea 
 
Next, we compute the components of $A$.
\begin{lemma}\lab{lemma:computationofthecomponentsofthetensorAforexactenergyconservationKerr}
Let $A$ the spacetime 1-form given by \eqref{eq:thespacetime1formAappearingintheenergyestimate}. Then, 
we have
\beaa
A_4 &=& -4\rhod\pr_t^3 -4(\etab\wedge\eta)\pr_t^3 +  \trch\big({}^{(h)}\pr_t\wedge\etab)  - \atrch\big(\etab\c {}^{(h)}\pr_t),\\
A_3 &=& 4\rhod\pr_t^4 +4(\etab\wedge\eta)\pr_t^4 +  \trchb\big({}^{(h)}\pr_t\wedge\eta)  -\atrchb\big(\eta\c {}^{(h)}\pr_t\big),\\
A_e &=& \Big( -  \trchb\dual\eta_e   + \atrchb\eta_e\Big)\pr_t^3 +\Big( -  \trch\dual\etab_e   + \atrch\etab_e\Big)\pr_t^4\\
&& -\frac{1}{2}\Big(4\rho  + \trch\trchb+\atrch\atrchb\Big)\dual ({}^{(h)}\pr_t)_e.
\eeaa
\end{lemma}

\begin{proof}
We rewrite $A_\mu$ as 
\beaa
A_\mu &=& \in^{bc}\Rdot_{bc \mu 3}\pr_t^3 +\in^{bc}\Rdot_{bc \mu 4}\pr_t^4 +\in^{bc}\Rdot_{bc \mu d}\pr_t^d.
\eeaa
Next, we compute the various components of $A_\mu$. We have in Kerr, using the horizontal tensor ${}^{(h)}\pr_t$ defined by $({}^{(h)}\pr_t)_b=(\pr_t)_b$, the definition  \eqref{eq:DefineRdot} of $\Rdot$,  and Proposition \ref{proposition:componentsofB}, 
 \beaa
A_4 &=& \in^{bc}\Rdot_{bc43}\pr_t^3  +\in^{bc}\Rdot_{bc4d}\pr_t^d,\\
&=& \in^{bc}\Big(-2\in_{bc}\dual\rho -2(\etab_b\eta_c-\eta_b\etab_c)\Big)\pr_t^3\\
&& +  \frac{1}{2}\in^{bc}\Big(\trch(\de_{db}\etab_c-\de_{dc}\etab_b) + \atrch(\in_{db}\etab_c - \in_{dc}\etab_b)\Big)\pr_t^d\\
&=& -4\rhod\pr_t^3 -4(\etab\wedge\eta)\pr_t^3 +  \frac{1}{2}\in^{bc}\Big(\trch(\etab_c\T_b - \etab_b\T_c) + \atrch(-\etab_c\dual ({}^{(h)}\pr_t)_b + \dual ({}^{(h)}\pr_t)_c \etab_b)\Big)\\
&=& -4\rhod\pr_t^3 -4(\etab\wedge\eta)\pr_t^3 +  \trch\big({}^{(h)}\pr_t\wedge\etab)  - \atrch\big(\etab\c {}^{(h)}\pr_t),
\eeaa
\beaa
A_3 &=& \in^{bc}\Rdot_{bc34}\pr_t^4 +\in^{bc}\Rdot_{bc3d}\pr_t^d\\
&=& 4\rhod\pr_t^4 +4(\etab\wedge\eta)\pr_t^4 +  \trchb\big({}^{(h)}\pr_t\wedge\eta)  - \atrchb\big(\eta\c {}^{(h)}\pr_t\big),
\eeaa
and
\beaa
A_e &=& \in^{bc}\Rdot_{bce3}\pr_t^3 +\in^{bc}\Rdot_{bce4}\pr_t^4 +\in^{bc}\Rdot_{bced}\pr_t^d\\
&=& \frac{1}{2}\in^{bc}\Big( -  \trchb  \big( \de_{eb}\eta_c -  \de_{ec} \eta_b\big)  -  \atrchb \big( \in_{eb}  \eta_c -  \in_{ec}  \eta_b\big)\Big)\pr_t^3\\
&&+ \frac{1}{2}\in^{bc}\Big( -  \trch  \big( \de_{eb}\etab_c -  \de_{ec} \etab_b\big)  -  \atrch \big( \in_{eb}  \etab_c -  \in_{ec}  \etab_b\big)\Big)\pr_t^4\\
&& +\frac{1}{2}\in^{bc}\left(-2\in_{bc}\in_{ed}\rho  -\frac{1}{2}\left(\trch\trchb+\atrch\atrchb\right)\in_{bc}\in_{ed}\right)\pr_t^d\\
&=& \Big( -  \trchb\dual\eta_e   + \atrchb\eta_e\Big)\pr_t^3 +\Big( -  \trch\dual\etab_e   + \atrch\etab_e\Big)\pr_t^4\\
&& -\frac{1}{2}\Big(4\rho  + \trch\trchb+\atrch\atrchb\Big)\dual ({}^{(h)}\pr_t)_e
\eeaa
as stated. This concludes the proof of Lemma \ref{lemma:computationofthecomponentsofthetensorAforexactenergyconservationKerr}.
\end{proof}

We infer the following corollary.
\begin{corollary}\lab{cor:computationofthecomponentsofthetensorAforexactenergyconservationKerr}
We have
\bea
A_\mu=-\D_\mu\left(\Im\left(\frac{2m}{q^2}\right)\right).
\eea
\end{corollary}

\begin{proof}
We have in Kerr
\beaa
\pr_t^4=\frac{1}{2}, \qquad \pr_t^3=\frac{1}{2}\frac{\De}{|q|^2}, \qquad (\pr_t)_b=({}^{(h)}\pr_t)_b=-a\Re(\Jk)_b.
\eeaa
Plugging in the identities of Lemma \ref{lemma:computationofthecomponentsofthetensorAforexactenergyconservationKerr}, we infer
\beaa
A_4 &=& -2\rhod\frac{\De}{|q|^2} -2(\etab\wedge\eta)\frac{\De}{|q|^2} -  a\trch\big(\Re(\Jk)\wedge\etab)  +a\atrch\big(\etab\c\Re(\Jk)),\\
A_3 &=& 2\rhod +2(\etab\wedge\eta) -  a\trchb\big(\Re(\Jk)\wedge\eta)  +a\atrchb\big(\eta\c\Re(\Jk)\big),\\
A_e &=& \frac{1}{2}\Big( -  \trchb\dual\eta_e   + \atrchb\eta_e\Big)\frac{\De}{|q|^2} +\frac{1}{2}\Big( -  \trch\dual\etab_e   + \atrch\etab_e\Big)\\
&& +\frac{a}{2}\Big(4\rho  + \trch\trchb+\atrch\atrchb\Big)\dual\Re(\Jk)_e.
\eeaa

Next, we rewrite $2A_e$ as
\beaa
2A_e &=& \Big( -  \Re(\tr\Xb)\Im(H_e)   -\Im(\tr\Xb)\Re(H_e)\Big)\frac{\De}{|q|^2} +\Big( -   \Re(\tr X)\Im(\Hb_e)   -\Im(\tr X)\Re(\Hb_e)\Big)\\
&& +a\Big(4\Re(P)  + \Re(\tr X)\Re(\tr\Xb)+\Im(\tr X)\Im(\tr\Xb)\Big)\dual\Re(\Jk)_e\\
&=& -\Im(\tr\Xb H_e)\frac{\De}{|q|^2} -\Im(\tr X\Hb_e)+a\Re\big(4P+\tr X \ov{\tr\Xb}\big)\dual\Re(\Jk)_e\\
&=& -\Im\left(-\frac{2}{\ov{q}}\frac{aq}{|q|^2}\Jk_e\right)\frac{\De}{|q|^2} -\Im\left(\frac{2}{q}\left(-\frac{\De}{|q|^2}\frac{a\ov{q}}{|q|^2}\Jk_e\right)\right)+a\Re\left(-\frac{8m}{q^3} -\frac{\De}{|q|^2}\frac{4}{q^2} \right)\dual\Re(\Jk)_e\\
&=& \frac{2a\De}{|q|^2}\left(\Im\left(\left(\frac{1}{\ov{q}^2}+\frac{1}{q^2}\right)\Jk_e\right) -\Re\left(\frac{2}{q^2}\right)\dual\Re(\Jk)_e\right) -a\Re\left(\frac{8m}{q^3}\right)\dual\Re(\Jk)_e
\eeaa
and hence
\beaa
2A_e &=& \frac{2a\De}{|q|^2}\left(\Re\left(\frac{1}{\ov{q}^2}+\frac{1}{q^2}\right)\Im(\Jk_e) -\Re\left(\frac{2}{q^2}\right)\dual\Re(\Jk)_e\right) -a\Re\left(\frac{8m}{q^3}\right)\dual\Re(\Jk_e)\\
&=& \frac{4a\De}{|q|^2}\left(\Re\left(\frac{1}{q^2}\right)\Im(\Jk_e) -\Re\left(\frac{1}{q^2}\right)\dual\Re(\Jk_e)\right) -a\Re\left(\frac{8m}{q^3}\right)\dual\Re(\Jk)_e\\
&=& -a\Re\left(\frac{8m}{q^3}\right)\dual\Re(\Jk)_e.
\eeaa
Since 
\beaa
\nab(\cos\th) = \Re(i\Jk)=-\dual\Re(\Jk), 
\eeaa
we infer
\beaa
A_e &=& a\Re\left(\frac{4m}{q^3}\right)\nab_e(\cos\th) =\Re\left(\frac{4m}{q^3}\nab_e(a\cos\th)\right)\\
&=& \Im\left(\frac{4m}{q^3}\nab_e(ia\cos\th)\right)=\Im\left(\frac{4m}{q^3}\nab_e(q)\right)
\eeaa
and thus
\beaa
A_e &=& -\nab_e\left(\Im\left(\frac{2m}{q^2}\right)\right).
\eeaa

Next, we rewrite $A_4$ as 
\beaa
A_4 &=& -2\rhod\frac{\De}{|q|^2} -2(\etab\c\dual\eta)\frac{\De}{|q|^2} -  a\trch\big(\Re(\Jk)\c\dual\etab)  +a\atrch\big(\etab\c\Re(\Jk))\\
&=& -2\rhod\frac{\De}{|q|^2} -2(\Re(\Hb)\c\Im(H))\frac{\De}{|q|^2} -  a\Re(\tr X)\big(\Re(\Jk)\c\Im(\Hb))  -a\Im(\tr X)\big(\Re(\Hb)\c\Re(\Jk))\\
&=& -2\rhod\frac{\De}{|q|^2} -2(\Re(\Hb)\c\Im(H))\frac{\De}{|q|^2} -  a\Big(\Re(\tr X)\Im(\Hb)+\Im(\tr X)\Re(\Hb)\Big)\c\Re(\Jk)\\
&=& -2\rhod\frac{\De}{|q|^2} -2(\Re(\Hb)\c\Im(H))\frac{\De}{|q|^2} -  a\Im(\tr X \Hb)\c\Re(\Jk)\\
\eeaa
and hence, since $\Im(\Jk)=\dual\Re(\Jk)$, 
\beaa
A_4 &=& -2\rhod\frac{\De}{|q|^2} +2a^2\frac{\De}{|q|^6}\Re(\ov{q}\Jk)\c\Im(q\Jk) -  a\Im\left(\frac{2}{q}\frac{\De}{|q|^2}\left(- \frac{a\ov{q}}{|q|^2}\Jk\right)\right)\c\Re(\Jk)\\
&=& -2\rhod\frac{\De}{|q|^2} +2a^2\frac{\De}{|q|^6}\Big(\Re(q)\Re(\Jk)+\Im(q)\Im(\Jk)\Big)\c\Big(\Im(q)\Re(\Jk)+\Re(q)\Im(\Jk)\Big)\\
&& +  2a^2\frac{\De}{|q|^6}\Im(\ov{q}^2\Jk)\c\Re(\Jk)\\
&=& -2\rhod\frac{\De}{|q|^2} +4a^2\frac{\De}{|q|^6}\Re(q)\Im(q)|\Re(\Jk)|^2 +  2a^2\frac{\De}{|q|^6}\Im(\ov{q}^2)|\Re(\Jk)|^2\\
&=& -2\rhod\frac{\De}{|q|^2} +4a^2\frac{\De}{|q|^6}\Re(q)\Im(q)|\Re(\Jk)|^2 +  4a^2\frac{\De}{|q|^6}\Re(\ov{q})\Im(\ov{q})|\Re(\Jk)|^2\\
&=& -2\rhod\frac{\De}{|q|^2}.
\eeaa
Since $e_4(q)=\frac{\De}{|q|^2}$, this yields
\beaa
A_4 &=& 4\Im\left(\frac{m}{q^3}\right)\frac{\De}{|q|^2}=4\Im\left(\frac{m}{q^3}\right)e_4(q)\\
&=& -e_4\left(\Im\left(\frac{2m}{q^2}\right)\right)
\eeaa
and similarly 
\beaa
A_3 &=& -e_3\left(\Im\left(\frac{2m}{q^2}\right)\right).
\eeaa
Thus, we have obtained 
\beaa
A_e = \nab_e\left(\Im\left(\frac{2m}{q^2}\right)\right), \qquad A_4 = e_4\left(\Im\left(\frac{2m}{q^2}\right)\right), \qquad A_3=e_3\left(\Im\left(\frac{2m}{q^2}\right)\right),
\eeaa
and hence
\beaa
A_\mu=-\D_\mu\left(\Im\left(\frac{2m}{q^2}\right)\right)
\eeaa
as stated. This concludes the proof of Corollary  \ref{cor:computationofthecomponentsofthetensorAforexactenergyconservationKerr}.
\end{proof}

\eqref{eq:basiccomputationenergyestimatesKerrwithThatdeltaflsdiuh} and Corollary  \ref{cor:computationofthecomponentsofthetensorAforexactenergyconservationKerr} imply 
 \bea\label{eq:modifieddivergenceTmod}
 \nn \D^\mu  \PP_\mu[\That_\de, 0, 0] &=& \frac 1 2 \QQ  \c \,^{(\That_\de)} \pi   -\D_\nu\left(\Im\left(\frac{m}{q^2}\right)\right) \dual\psi\c\Db^\nu\psi  +(\That_\de-\pr_t)^\mu \Db^\nu  \psi ^a\Rdot_{ ab   \nu\mu}\psi^b\\
  &&+ \nab_{\That_\de} \psi  \c \big(\squared_2 \psi - V \psi \big).
 \eea
  
Next, we modify the identity \eqref{eq:modifieddivergenceTmod} to cancel the second term on the RHS. To this end, we consider the following modified current
\bea
\widetilde{\PP}_\mu &:=&  \PP_\mu[\That_\de, 0, 0] + \tilde{w}\dual\psi\c\Ddot_\mu\psi,
\eea
for a scalar function $\tilde{w}=\tilde{w}(r, \cos\th)$ to be chosen below. Since 
\beaa
\bsplit
\D^\mu\Big[\tilde{w}\dual\psi\c\Ddot_\mu\psi\Big] &= \tilde{w}\dual\psi\c\Ddot^\mu\Ddot_\mu\psi + \tilde{w}\dual\Ddot^\mu\psi\c\Ddot_\mu\psi+\D^\mu(\tilde{w})\dual\psi\c\Ddot_\mu\psi\\
&= \tilde{w}\dual\psi\c\squared_2\psi +\D^\mu(\tilde{w})\dual\psi\c\Ddot_\mu\psi\\
&= \tilde{w}\dual\psi\c\left(V\psi - \frac{4 a\cos\th}{|q|^2}\dual \nab_t  \psi+N\right)+\D^\mu(\tilde{w})\dual\psi\c\Ddot_\mu\psi\\
&= - \tilde{w}\frac{4 a\cos\th}{|q|^2}\nab_t(|\dual\psi|^2)+\tilde{w}\dual\psi\c N+\D^\mu(\tilde{w})\dual\psi\c\Ddot_\mu\psi\\
&= \D^\mu(\tilde{w})\dual\psi\c\Ddot_\mu\psi - \D_\mu\left(\pr_t^\mu \tilde{w}\frac{4 a\cos\th}{|q|^2}|\psi|^2\right) + \tilde{w}\dual\psi\c N
\end{split}
\eeaa
we infer
\beaa
\bsplit
\DD^\mu\left(\widetilde{\PP}_\mu +\pr_t^\mu \tilde{w}\frac{4 a\cos\th}{|q|^2}|\psi|^2\right) =&  \frac 1 2 \QQ  \c \,^{(\That_\de)} \pi +\D_\nu\left(-\Im\left(\frac{m}{q^2}\right)+\tilde{w}\right) \dual\psi\c\Db^\nu\psi  \\
&+(\That_\de-\pr_t)^\mu \Db^\nu  \psi ^a\Rdot_{ ab   \nu\mu}\psi^b   \\
& + \nab_{\That_\de} \psi  \c \big(\squared_2 \psi - V \psi \big)+ \tilde{w}\dual\psi\c N.
\end{split}
\eeaa
 Next, we make the following choice for $\tilde{w}$
\bea
\tilde{w} := \Im\left(\frac{m}{q^2}\right)= -\frac{2amr\cos\th}{|q|^4}
\eea
which yields
\bea\lab{eq:modifieddivergenceTmod:ter}
\bsplit
\DD^\mu\left(\widetilde{\PP}_\mu + \pr_t^\mu \tilde{w}\frac{4 a\cos\th}{|q|^2}|\psi|^2\right) =&  \frac 1 2 \QQ  \c \,^{(\That_\de)} \pi +(\That_\de-\pr_t)^\mu \Db^\nu  \psi ^a\Rdot_{ ab   \nu\mu}\psi^b   \\
& + \nab_{\That_\de} \psi  \c \big(\squared_2 \psi - V \psi \big)+ \tilde{w}\dual\psi\c N.
\end{split}
\eea

 Using  \eqref{eq:QQpi(That_de)}, we have 
 \beaa
   \big|\QQ  \c \,^{(\That_\de)} \pi\big|&\les&\left( \frac{4 |a|  r} {|q|^2 (r^2+a^2)}|\chi_0|+ O(|a| \de^{-1} ) |\chi'_0| \right) | \nab_\phi\psi | | \nab_\Rhat \psi|\\
   &\les&\mathbb{1}_{\Mntrap}\de^{-1}   \frac{ |a| }{ r^3}   | \nab_\phi\psi | | \nab_\Rhat \psi|.
 \eeaa
 Also, we have $\That_\de-\pr_t=\chi_\de(r)\pr_\phi$ and hence
 \beaa
 && (\That_\de-\pr_t)^\mu \Db^\nu  \psi ^a\Rdot_{ ab   \nu\mu}\psi^b\\
 &=& \frac{1}{2}\chi_\de(r)\pr_\phi^\mu \in^{ab}\Rdot_{ ab   \nu\mu} \dual\psi\c\Db^\nu  \psi\\
 &=& \frac{1}{2}\chi_\de(r)\left(-\frac{1}{2}\pr_\phi^\mu \in^{ab}\Rdot_{ ab3\mu} \dual\psi\c\nab_4\psi
  -\frac{1}{2}\pr_\phi^\mu \in^{ab}\Rdot_{ ab4\mu} \dual\psi\c\nab_3\psi + \pr_\phi^\mu \in^{ab}\Rdot_{ abc\mu} \dual\psi\c\nab^c\psi \right).
 \eeaa
 Since we have
 \beaa
 \pr_\phi^3=O(ar^{-2}\De), \qquad  \pr_\phi^4=O(a), \qquad \pr_\phi^c=O(r^2)\Re(\Jk)^c,
 \eeaa
we infer, together with  the definition  \eqref{eq:DefineRdot} of $\Rdot$,  and Proposition \ref{proposition:componentsofB}, 
\beaa
\pr_\phi^\mu \in^{ab}\Rdot_{ ab3\mu} &=& \pr_\phi^4 \in^{ab}\Rdot_{ ab34}+\pr_\phi^c \in^{ab}\Rdot_{ ab3c}= O(a)\big(\rhod, \eta\etab\big)+O(r)\chib\eta\\
&=& O(ar^{-2}),\\
\pr_\phi^\mu \in^{ab}\Rdot_{ ab4\mu} &=& \pr_\phi^3 \in^{ab}\Rdot_{ ab43}+\pr_\phi^c \in^{ab}\Rdot_{ ab4c}=O(ar^{-2}\De)\big(\rhod, \eta\etab\big)+O(r)\chi\etab\\
&=& O(ar^{-4}\De),\\
\pr_\phi^\mu \in^{ab}\Rdot_{ abc\mu} &=& \pr_\phi^4 \in^{ab}\Rdot_{ abc4}+\pr_\phi^3 \in^{ab}\Rdot_{ abc3}+\pr_\phi^d \in^{ab}\Rdot_{ abcd}\\
&=& O(a)\chi\etab+O(ar^{-2}\De)\chib\eta +O(r)(\rho, \chi\chib)\\
&=& O(r^{-1}),
\eeaa 
 which implies, since $\chi_\de= \frac{a}{r^2+a^2} \chi_0\Big( \de^{-1} \frac{\TT}{r^3} \Big)$, 
  \beaa
 && \big|(\That_\de-\pr_t)^\mu \Db^\nu  \psi ^a\Rdot_{ ab   \nu\mu}\psi^b\big|\\
 &\les& \frac{|a|}{r^2}\mathbb{1}_{\Mntrap}\Big(|\pr_\phi^\mu \in^{ab}\Rdot_{ ab3\mu}| |\dual\psi\c\nab_4\psi|+|\pr_\phi^\mu \in^{ab}\Rdot_{ ab4\mu}||\dual\psi\c\nab_3\psi| \\
 && + |\pr_\phi^\mu \in^{ab}\Rdot_{ abc\mu}| |\dual\psi\c\nab^c\psi| \Big)\\
 &\les& \frac{|a|}{r^3}\mathbb{1}_{\Mntrap}\Big(|\nab_4\psi|+r^{-2}|\De||\nab_3\psi|+|\nab\psi|\Big)|\psi|\\
 &\les& \frac{|a|}{r^3}\mathbb{1}_{\Mntrap}\Big(|\nab_{\Rhat}\psi|+|\nab_T\psi|+|\nab\psi|\Big)|\psi|.
 \eeaa
    Finally, using equation \eqref{eq:Gen.RW}, we have
   \beaa
   \nab_{\That_\de} \psi  \c \big(\squared_2 \psi - V \psi \big)&=&- \frac{4 a\cos\th}{|q|^2}  \nab_{\That_\de} \psi  \c  \dual \nab_T  \psi+ \nab_{\That_\de} \psi  \c  N \\
   &=&- \frac{4 a\cos\th}{|q|^2}(  \nab_{T} \psi +\chi_\de \nab_\phi \psi ) \c  \dual \nab_T  \psi+ \nab_{\That_\de} \psi  \c  N \\
   &=&- \frac{4 a\cos\th}{|q|^2}\chi_\de \nab_\phi \psi   \c  \dual \nab_T  \psi+ \nab_{\That_\de} \psi  \c  N 
   \eeaa
      where we have the crucial cancellation $\nab_T \psi \c \dual \nab_T \psi=0$. 
     
      We summarize the result in the following.
      \begin{lemma}
    \lab{lemma:CurrentTmod}
    Consider the modified current 
    \beaa
    \widetilde{\PP}_\mu &=&  \PP_\mu[\That_\de, 0, 0] + \tilde{w}\dual\psi\c\Ddot_\mu\psi, \qquad \tilde{w} = \Im\left(\frac{m}{q^2}\right),
    \eeaa
    where  the vectorfield $\That_\de$ is given by $\That_\de= \pr_t +\chi_\de(r) \pr_\phi$ for $\de=\frac{1}{10}$ and $|a|/m\ll 1$ small enough.   Then,  for all $r\ge r_+$,
    \beaa
  &&\left|   \D^\mu\left(\widetilde{\PP}_\mu + \pr_t^\mu \tilde{w}\frac{4 a\cos\th}{|q|^2}|\psi|^2\right) -\Big(\nab_{\That_\de} \psi + \tilde{w}\dual\psi\Big) \c  N \right|\\ 
  &\les&  \mathbb{1}_{\Mntrap}  \left(\de^{-1}   \frac{ |a| }{ r^3}  |\nab_\Rhat \psi||\nab_\phi\psi|+\frac{|a|}{r^4} | \nab_T  \psi|  | \nab_\phi \psi | + \frac{|a|}{r^3}\Big[|\nab_\Rhat\psi|+|\nab_T\psi|+|\nab\psi|\Big]|\psi| \right).
    \eeaa
    \end{lemma} 
 
 Integrating the above inequality on $\MM=\MM(\tau_1, \tau_2) $ and applying the divergence theorem we deduce, in view of the definition of $\Mor^{ax}_{deg}[\psi](\tau_1, \tau_2)$,  
 \beaa
 &&\int_{\Si(\tau_2)}\left(\widetilde{\PP}_\mu + \pr_t^\mu \tilde{w}\frac{4 a\cos\th}{|q|^2}|\psi|^2\right)\c N_\Si +\int_{\Si_*(\tau_1, \tau_2)}\left(\widetilde{\PP}_\mu + \pr_t^\mu \tilde{w}\frac{4 a\cos\th}{|q|^2}|\psi|^2\right)\c N_{\Si_*}\\ 
 &&+\int_{\AA(\tau_1, \tau_2)}\left(\widetilde{\PP}_\mu + \pr_t^\mu \tilde{w}\frac{4 a\cos\th}{|q|^2}|\psi|^2\right)\c N_{\AA} \\
 &\les&  \int_{\Si(\tau_1)}\left(\widetilde{\PP}_\mu + \pr_t^\mu \tilde{w}\frac{4 a\cos\th}{|q|^2}|\psi|^2\right)\c N_\Si+ \frac{|a|}{m}\Mor^{ax}_{deg}[\psi](\tau_1, \tau_2)   +\left|\int_{\MM(\tau_1, \tau_2)}  \nab_{\That_\de } \psi  \c N\right|\\
 &&+\int_{\MM(\tau_1, \tau_2)}|N|^2
 \eeaa
where we used the fact that $|\tilde{w}|\les amr^{-3}$ to control the term $\tilde{w}\dual\psi\c N$.   

Finally, since $|\tilde{w}|\les amr^{-3}$, we have for a vectorfield $N$
\bea\lab{eq:theextraterminthemodifycurrentissmallontheboundaryforsmalla}
\nn\left(\widetilde{\PP}_\mu + \pr_t^\mu \tilde{w}\frac{4 a\cos\th}{|q|^2}|\psi|^2\right)\c N &=&  \QQ(\That_\de, N) + \tilde{w}\dual\psi\c \nab_N\psi +\g( \pr_t, N)\tilde{w}\frac{4 a\cos\th}{|q|^2}|\psi|^2\\
\nn&=& \QQ(\That_\de, N)  -O(amr^{-3})|\psi||\nab_N\psi|\\
&& -O(a^2mr^{-5})|\g( \pr_t, N)||\psi|^2.
\eea
Also, we have, in view of  Lemma \ref{Lemma:QQ(That,NSi)} and the definition of $E_{deg}[\psi]$, for some constant $c_0>0$
\beaa
\bsplit
\int_{\Si(\tau)}\QQ(\That, N_\Si)&\ge   c_0E_{deg}[\psi](\tau)- O(\deh) E_{r\leq r_+}[\psi](\tau),\\
\int_{\Si_*(\tau_1, \tau_2)}\QQ(\That, N_{\Si_*}) &\ge c_0F_{\Si_*}[\psi](\tau_1, \tau_2),\\
\int_{\AA(\tau_1, \tau_2)}\QQ(\That, N_{\AA}) &\gtrsim - \deh F_{\AA}[\psi](\tau_1, \tau_2).
 \end{split}
\eeaa
Since we have, in view of the definition of $\That_\de$,
\beaa
\That_\de &=& \That +\frac{a}{r^2+a^2}(\chi_0-1)\left(\de^{-1}\frac{\TT}{r^3}\right)\pr_\phi,
\eeaa
we have $\That_\de= \That $ on $\Si_*$ and $\AA$, and hence
\beaa
\bsplit
\int_{\Si_*(\tau_1, \tau_2)}\QQ(\That_\de, N_{\Si_*}) &\ge c_0F_{\Si_*}[\psi](\tau_1, \tau_2),\\
\int_{\AA(\tau_1, \tau_2)}\QQ(\That_\de, N_{\AA}) &\gtrsim - \deh F_{\AA}[\psi](\tau_1, \tau_2).
 \end{split}
\eeaa
In view of \eqref{eq:theextraterminthemodifycurrentissmallontheboundaryforsmalla}, we deduce for $|a|\ll m, m\deh$,  
\beaa
\bsplit
\int_{\Si_*(\tau_1, \tau_2)}\left(\widetilde{\PP}_\mu + \pr_t^\mu \tilde{w}\frac{4 a\cos\th}{|q|^2}|\psi|^2\right)\c N_{\Si_*} &\ge c_0F_{\Si_*}[\psi](\tau_1, \tau_2),\\
\int_{\AA(\tau_1, \tau_2)}\left(\widetilde{\PP}_\mu + \pr_t^\mu \tilde{w}\frac{4 a\cos\th}{|q|^2}|\psi|^2\right)\c N_{\AA} &\gtrsim - \deh F_{\AA}[\psi](\tau_1, \tau_2).
 \end{split}
\eeaa

Also, we have
\beaa
\PP[\That_\de]\c N_\Si &=& \QQ(\That, N_\Si) +\frac{a}{r^2+a^2}(\chi_0-1)\left(\de^{-1}\frac{\TT}{r^3}\right)\QQ(\pr_\phi, N_\Si).
\eeaa
In view of the above and the properties of $\chi_0$, this yields 
\beaa
\int_{\Si(\tau)}\PP[\That_\de]\c N_\Si &\ge &  c_0E_{deg}[\psi](\tau) - O(a)\int_{\Si_{+}(\tau)}|\QQ(\pr_\phi, N_\Si)|\mathbb{1}_{\frac{|\TT|}{r^3}\leq 2\de} - O(\deh) E_{r\leq r_+}[\psi](\tau).
\eeaa
Since the set $\frac{|\TT|}{r^3}\leq 2\de$ is localized near $r=3m$ for $|a|\ll m$ and $\de>0$ small enough, we infer, using also \eqref{eq:theextraterminthemodifycurrentissmallontheboundaryforsmalla}, 
\beaa
\int_{\Si(\tau)}\left(\widetilde{\PP}_\mu + \pr_t^\mu \tilde{w}\frac{4 a\cos\th}{|q|^2}|\psi|^2\right)\c N_\Si &\ge &  \frac{c_0}{2} E_{deg}[\psi](\tau)  - O(\deh) E_{r\leq r_+}[\psi](\tau),
\eeaa
and hence 
 \beaa
\nn E_{deg}[\psi](\tau_2) +F_{\Si_*}[\psi](\tau_1, \tau_2) &\les&  E_{deg}[\psi](\tau_1) +\deh\big(E_{r\leq r_+(1+\deh)}(\tau_2)[\psi]+F_\AA[\psi](\tau_1, \tau_2)\big)\\
&&+ \frac{|a|}{m}\Mor^{ax}_{deg}[\psi](\tau_1, \tau_2)   +\left|\int_{\MM(\tau_1, \tau_2)}  \nab_{\That_\de } \psi  \c N\right|\\
&& +\int_{\MM(\tau_1, \tau_2)}|N|^2.
\eeaa
 This concludes the proof of  Proposition \ref{proposition:Energy1}.

     
\section{Conditional control of Energy--Morawetz estimate for $r$ large enough} 


In this section, we prove the following conditional control of Energy--Morawetz estimate in $\MM(r\geq r_1)$ for $r_1$ large enough.
\begin{proposition}\lab{prop:recoverEnergyMorawetzwithrweightfromnoweight}
Let $\psi$ a solution to the gRW equation \eqref{eq:Gen.RW}. Also, recall the norms $E[\psi]$ and $\textrm{Mor}[\psi]$ introduced in section \ref{subsection:basicnormsforpsi}. Then, for $r_1=r_1(m)$ large enough, we have
\bea
\bsplit
\textrm{Mor}_{r\geq 2r_1}[\psi](\tau_1, \tau_2) \les&  \sup_{\tau\in[\tau_1, \tau_2]}E_{r\geq r_1}[\psi](\tau) +r_1\textrm{Mor}_{r_1\leq r\leq 2r_1}[\psi](\tau_1, \tau_2)\\
& +\mathcal{N}_{r\geq r_1}[\psi, N](\tau_1, \tau_2),
\end{split}
\eea
and
\bea
\bsplit
\sup_{\tau\in[\tau_1, \tau_2]}E_{r\geq 2r_1}[\psi](\tau) \les& 
E_{r\geq r_1}[\psi](\tau_1) +\mathcal{N}_{r\geq r_1}[\psi, N](\tau_1, \tau_2)+r_1\textrm{Mor}_{r_1\leq r\leq 2r_1}[\psi](\tau_1, \tau_2)\\
&+\frac{|a|}{r_1}\textrm{Mor}_{r\geq 2r_1}[\psi](\tau_1, \tau_2).
\end{split}
\eea
\end{proposition}

\begin{proof}
Let a smooth cut-off function $\chi$ such that $\chi=0$ for $r\leq 1$ and $\chi=1$ for $r\geq 2$, and let $\chi_{r_1}$ and $\psi_{r_1}$ defined by 
\beaa
\chi_{r_1}(r):=\left(\frac{r}{r_1}\right), \qquad \psi_{r_1}:=\chi_{r_1}\psi,
\eeaa
so that the support of $\psi_{r_1}$ is included in $r\geq r_1$ and $\psi_{r_1}=\psi$ for $r\geq 2r_1$. 
Then, since $\psi$ satisfies \eqref{eq:Gen.RW}, $\psi_{r_1}$ satisfies the following gRW equation
\bea
\bsplit
\squared_2 \psi_{r_1} -V\psi_{r_1}&=- \frac{4 a\cos\th}{|q|^2}\dual \nab_T  \psi_{r_1}+N_{r_1}, \\ 
N_{r_1}&:=\g^{rr}\chi_{r_1}'(r)\nab_r\psi+\square(\chi_{r_1})\psi+\chi\left(\frac{r}{r_1}\right)N.
\end{split}
\eea

We start by deriving a Morawetz estimate for $\psi_{r_1}$. According to Proposition \ref{proposition:Morawetz1}, with the choice  $M= 0$,  the generalized current associated to the Morawetz vectorfield in \eqref{identity:prop.Morawetz1} is given by
 \beaa
|q|^2 \EE_{r_1}[X, w, M=0]&=&\AA|\nab_r\psi_{r_1}|^2+ P[\psi_{r_1}]+\VV |\psi_{r_1}|^2.
 \eeaa
Here we consider the choices for $z, f, h$ made in Proposition \ref{prop:Choice-zhf}  so that the coefficients $\AA$ and $\VV$ are given by \eqref{eq:prop-Choice-AA} and \eqref{eq:prop:Choice-zhf} respectively, and  the principal term $P[\psi_{r_1}]$ is given by  
 \beaa
 P&=& \frac{\TT}{r}\frac{r^2+a^2}{r^2-a^2} \left(  \frac{2\TT}{ (r^2+a^2)^3}|q|^2|\nab\psi|^2 - \frac{4ar}{(r^2+a^2)^2} \nab_\That\psi\c\nab_\phi \psi\right).
 \eeaa
Note that we have for $r\geq r_1$
\beaa
\bsplit
\AA &=6m(1+O(mr^{-1})), \qquad \VV=\frac{8}{r}(1+O(mr^{-1})), \\
P &=   2r(1+O(mr^{-1})) |\nab\psi_{r_1}|^2+O(m)|\nab_{\That}\psi_{r_1}||\nab\psi_{r_1}|+O(m^2r^{-1})|\nab_\That\psi_{r_1}|^2,
\end{split}
\eeaa
where we used in particular \eqref{eq:prop-Choice-AA} and \eqref{eq:prop:Choice-zhf}. This implies, for $r\geq r_1$ with $r_1=r_1(m)$ large enough, 
 \beaa
\EE_{r_1}[X, w, M=0]&\geq & \frac{m}{r^2}|\nab_r\psi_{r_1}|^2 +\frac{1}{r}|\nab\psi_{r_1}|^2
+\frac{1}{r^3}|\psi_{r_1}|^2 -O(m^2r^{-3})|\nab_\That\psi_{r_1}|^2.
 \eeaa

Also, since we have for $r\geq r_1$
\beaa
z=\frac{1}{r^2}(1+O(mr^{-1})), \qquad f=-\frac{2}{r^3}(1+O(mr^{-1})), \qquad h=r^5(1+O(a^2r^{-2})),
\eeaa
we easily obtain the following simpler analogs of the estimates of Lemma \ref{lemma-control-rhs}
\beaa
&&\left(\nab_X\psi_{r_1}+\frac 1 2   w \psi_{r_1}\right)\c  \left(\squared_2 \psi_{r_1} - V \psi_{r_1} \right) \\
  &=& O(ar^{-3})|\psi_{r_1}||\nab_T\psi_{r_1}|+O(a^2mr^{-6})|\psi_{r_1}|^2+O(1)(|\nab_\Rhat\psi_{r_1}|+r^{-1}|\psi_{r_1}|)|N_{r_1}|\\
&&+ \D_\mu\left(\frac{2a\cos\th}{|q|^2}(\pr_r)^\mu  zh f \psi_{r_1}\c\nab_T\dual \psi_{r_1}\right) -\pr_t \left(  \frac{2a\cos\th}{|q|^2} zhf \psi_{r_1}\c  \nab_r\dual \psi_{r_1}\right)
\eeaa
and
\beaa
\nn&&\left|\left(\big(\rhod +\etab\wedge\eta\big)\frac{r^2+a^2}{\De}+\frac{2a^3r\cos\th(\sin\th)^2}{|q|^6}\right)\FF\nab_\That\psi_{r_1}\c\dual\psi_{r_1}\right|\\
\nn&&+\left|\frac{2a^2r\cos\th\FF(r)}{(r^2+a^2)|q|^4}\nab_\phi\psi_{r_1}\c\dual\psi_{r_1}\right|\\
& \leq & \frac{|a|m}{r^4}|\nab_{\That}\psi_{r_1}||\psi_{r_1}| +\frac{a^2}{r^5}|\nab_\phi\psi_{r_1}||\psi_{r_1}|.
\eeaa

In view of the above, and since $\psi_{r_1}$ is supported in $r\geq r_1$, we obtain from \eqref {definition-EE-X=FFR}, for $r_1=r_1(m)$ large enough, 
\beaa
&&\D^\mu  (\PP_{r_1})_\mu[X, w]\\
&=& \EE_{r_1}[X, w] +\left(\nab_X\psi_{r_1}+\frac 1 2   w \psi_{r_1}\right)\c  \big(\squared_2 \psi_{r_1} - V \psi_{r_1} \big) -\frac{2a^2r\cos\th\FF(r)}{(r^2+a^2)|q|^4}\nab_\phi\psi_{r_1}\c\dual\psi_{r_1}\\
&& -\left(\big(\rhod +\etab\wedge\eta\big)\frac{r^2+a^2}{\De}+\frac{2a^3r\cos\th(\sin\th)^2}{|q|^6}\right)\FF\nab_\That\psi_{r_1}\c\dual\psi_{r_1}\\
&\geq & \frac{m}{r^2}|\nab_r\psi_{r_1}|^2 +\frac{1}{r}|\nab\psi_{r_1}|^2
+\frac{1}{2r^3}|\psi_{r_1}|^2 -O(m^2r^{-3})|\nab_{\That}\psi_{r_1}|^2 -O(1)\Big( | \nab_\Rhat \psi_{r_1} | + r^{-1}|\psi_{r_1}| \Big)    |  N_{r_1}|\\
&& +\D_\mu\left(\frac{2a\cos\th}{|q|^2}(\pr_r)^\mu  zh f \psi_{r_1}\c\nab_T\dual \psi_{r_1}\right) -\pr_t \left(  \frac{2a\cos\th}{|q|^2} zhf \psi_{r_1}\c  \nab_r\dual \psi_{r_1}\right).
\eeaa
By applying the divergence theorem, and using that $\psi_{r_1}$ is supported in $r\geq r_1$, we then easily infer
 \beaa
&&\int_{\MM_+}  \frac{m}{r^2} |\nab_{\Rhat} \psi_{r_1}|^2+r^{-3}|\psi_{r_1}|^2 + r^{-1}|\nab\psi_{r_1}|^2   \\
\les && \sup_{\tau\in[\tau_1, \tau_2]}E_{r\geq r_1}[\psi](\tau)+\int_{\MM_{+} }\frac{m^2}{r^3}|\nab_{\That}\psi_{r_1}|^2 
 +\int_{\MM_+} \Big(|\nab_\Rhat\psi_{r_1}|+r^{-1}|\psi_{r_1}|\Big)| N_{r_1}|.
\eeaa
Also, given the form of $N_{r_1}$ and $\psi_{r_1}$, we deduce
 \beaa
&&\int_{\MM_+}  \frac{m}{r^2} |\nab_{\Rhat} \psi_{r_1}|^2+r^{-3}|\psi_{r_1}|^2 + r^{-1}|\nab\psi_{r_1}|^2   \\
\les && \sup_{\tau\in[\tau_1, \tau_2]}E_{r\geq r_1}[\psi](\tau)+\int_{\MM_{+} }\frac{m^2}{r^3}|\nab_{\That}\psi_{r_1}|^2 
+\frac{r_1}{m}\textrm{Mor}_{r_1\leq r\leq 2r_1}[\psi](\tau_1, \tau_2) +\mathcal{N}_{r\geq r_1}[\psi, N](\tau_1, \tau_2).
\eeaa
  
To complete the desired Morawetz estimate, we still need to recover the control of $\nab_{\That}\psi_{r_1}$. 
We compute 
\beaa
\Ddot^\a\left(\frac{m}{r^2}\psi_{r_1}\c\Ddot_\a\psi_{r_1}\right) &=& \frac{m}{r^2}\psi_{r_1}\c\squared_2\psi_{r_1}+\frac{m}{r^2}\Ddot^\a\psi_{r_1}\c\Ddot_\a\psi_{r_1} -\frac{2m}{r^3}\g^{rr}\psi_{r_1}\c\Ddot_r\psi_{r_1}\\
&=& \frac{m}{r^2}\psi_{r_1}\c\left(V\psi_{r_1} - \frac{4 a\cos\th}{|q|^2}\dual \nab_T  \psi_{r_1}+N_{r_1}\right)+\frac{m}{r^2}\Ddot^\a\psi_{r_1}\c\Ddot_\a\psi_{r_1}\\
&& -\frac{2m}{r^3}\g^{rr}\psi_{r_1}\c\Ddot_r\psi_{r_1}.
\eeaa
Since $\psi_{r_1}$ is supported in $r\geq r_1$, we deduce, for $r_1=r_1(m)$ large enough,
\beaa
\Ddot^\a\left(\frac{m}{r^2}\psi_{r_1}\c\Ddot_\a\psi_{r_1}\right) &=& -\frac{m}{r^2}(1+O(mr^{-1}))|\nab_{\That}\psi|^2\\
&&+O(1)\left(\frac{m}{r^2} |\nab_{\Rhat} \psi_{r_1}|^2+r^{-3}|\psi_{r_1}|^2 + r^{-1}|\nab\psi_{r_1}|^2\right)\\
&&+O(mr^{-2})|\psi_{r_1}||N_{r_1}|.
\eeaa
Together with the above, we infer
 \beaa
&&\int_{\MM_+}  \frac{m}{r^2} |\nab_{\Rhat} \psi_{r_1}|^2+\frac{m}{r^2} |\nab_{\That} \psi_{r_1}|^2+r^{-3}|\psi_{r_1}|^2 + r^{-1}|\nab\psi_{r_1}|^2   \\
\les && \sup_{\tau\in[\tau_1, \tau_2]}E_{r\geq r_1}[\psi](\tau)+\int_{\MM_{+} }\frac{m^2}{r^3}|\nab_{\That}\psi_{r_1}|^2 
+\frac{r_1}{m}\textrm{Mor}_{r_1\leq r\leq 2r_1}[\psi](\tau_1, \tau_2) +\mathcal{N}_{r\geq r_1}[\psi, N](\tau_1, \tau_2).
\eeaa
Since $\psi_{r_1}$ is supported in $r\geq r_1$, and since $\psi_{r_1}=\psi$ for $r\geq 2r_1$, we deduce, for $r_1=r_1(m)$ large enough,
\beaa
\textrm{Mor}_{r\geq 2r_1}[\psi](\tau_1, \tau_2) &\les&  \sup_{\tau\in[\tau_1, \tau_2]}E_{r\geq r_1}[\psi](\tau) 
+\frac{r_1}{m}\textrm{Mor}_{r_1\leq r\leq 2r_1}[\psi](\tau_1, \tau_2) +\mathcal{N}_{r\geq r_1}[\psi, N](\tau_1, \tau_2)
\eeaa
as stated.

Finally, it remains to control the energy. In view of Lemma \ref{lemma:CurrentTmod} applied to $\psi_{r_1}$, we have
     \beaa
     &&\left|   \D^\mu\left((\widetilde{\PP_{r_1}})_\mu + \pr_t^\mu \tilde{w}\frac{4 a\cos\th}{|q|^2}|\psi_{r_1}|^2\right) -\Big(\nab_{\That_\de} \psi_{r_1} + \tilde{w}\dual\psi_{r_1}\Big) \c  N_{r_1} \right|\\ 
  &\les&     \frac{ |a| }{ r^3}  |\nab_\Rhat \psi_{r_1}||\nab_\phi\psi_{r_1}|+\frac{|a|}{r^4} | \nab_T  \psi_{r_1}|  | \nab_\phi \psi_{r_1}|  +  \frac{|a|m}{r^4}\Big[|\nab_\Rhat\psi_{r_1}|+|\nab_T\psi_{r_1}|+|\nab\psi_{r_1}|\Big] \,|\psi_{r_1}|
    \eeaa
    where
      \beaa
    (\widetilde{\PP_{r_1}})_\mu &=&  (\PP_\mu)_{r_1}[\That_\de, 0, 0] + \tilde{w}\dual\psi_{r_1}\c\Ddot_\mu\psi_{r_1}, \qquad \tilde{w} = -\Im\left(\frac{m}{q^2}\right).
    \eeaa
By applying the divergence theorem, since $\psi_{r_1}$ is supported in $r\geq r_1$, and since $\psi_{r_1}=\psi$ for $r\geq 2r_1$, we easily infer
\beaa
\nn\sup_{\tau\in[\tau_1, \tau_2]}E_{r\geq 2r_1}[\psi](\tau) &\les& 
E_{r\geq r_1}[\psi](\tau_1) +\mathcal{N}_{r\geq r_1}[\psi, N](\tau_1, \tau_2)\\
&&+\frac{r_1}{m}\textrm{Mor}_{r_1\leq r\leq 2r_1}[\psi](\tau_1, \tau_2)+\frac{|a|}{r_1}\textrm{Mor}_{r\geq 2r_1}[\psi](\tau_1, \tau_2)
\eeaa
as stated. This concludes the proof of Proposition \ref{prop:recoverEnergyMorawetzwithrweightfromnoweight}.
\end{proof}

 
 \chapter{Proof of  $\SS$-derivative Morawetz estimates in Kerr}\label{chapter-proof-mor-2}
 
 
 In this chapter we provide a complete proof of the $\SS$-derivative  Morawetz estimates as stated in Section \ref{section:Higher derivative Morawetz  Estimates}. Recall that we are in Kerr throughout this chapter and that the results of this chapter will be extended to perturbations of Kerr in section \ref{sec:proofofresultschapter8inperturbationofKerr}.

 
 \section{Preliminaries}\label{section:preliminaries-chapter5}
 

In this section we collect preliminary results to extend the vector field method to include commutation with  the following second order differential operators,   see  Definition \ref{definition:symmetry-tensors}, 
\beaa
\SS_1(\psi)&=& \nab_T \nab_T \psi, \qquad \SS_2(\psi)= a \nab_{T} \nab_{Z} \psi, \qquad \SS_3(\psi)= a^2 \nab_Z \nab_Z \psi, \qquad \SS_4\psi= \OO (\psi),
\eeaa
which we denote $\SS_\aund$, for $\aund=1,2,3,4$.

\begin{lemma}\label{lemma-commuted-equ-Naund}
Given a $\mathfrak{s}_2$ tensor  $\psi$ solution  of the equation  \eqref{eq:Gen.RW}.  Then the commuted $\mathfrak{s}_2$ tensor  $\psia:=\SS_\aund \psi$ satisfies
  \bea
  \lab{eq:waveeqfor-psia}
  \squared_2\psia -V\psia&=&- \frac{ 4 a\cos\th}{|q|^2} \dual \nab_T  \psia+N_\aund, \qquad \aund=1, 2, 3, 4,
  \eea
  where
\beaa
|N_\aund| &\les& |\dk^{\leq 2}N| + ar^{-2}|\dk^{\leq 2}\psi|, \qquad \aund=1, 2, 3, 4.
\eeaa
\end{lemma}

\begin{proof} 
Since $\psi$ satisfies \eqref{eq:Gen.RW}, i.e. 
  \beaa
  \squared_2\psi -V\psi &=&- \frac{4 a\cos\th}{|q|^2} \dual \nab_T  \psi+N,
  \eeaa
we infer
 \beaa
  \squared_2\psia -V\psia &=&- \frac{4 a\cos\th}{|q|^2} \dual \nab_T  \psia+N_\aund,\\
  N_\aund &:=& -[\SS_\aund, \squared_2]\psi +[\SS_\aund, V]\psi -\dual\left[\SS_\aund, \frac{4 a\cos\th}{|q|^2} \nab_T\right]\psi+\dk^{\leq 2}N, \qquad \aund=1, 2, 3,\\
  N_4 &:=& -\frac{1}{|q|^2}[\SS_4, |q|^2\squared_2]\psi +\frac{1}{|q|^2}[\SS_4, |q|^2V]\psi -\frac{1}{|q|^2}\dual[\SS_4, 4a\cos\th\nab_T]\psi+\dk^{\leq 2}N,
  \eeaa
and hence
\beaa
|N_\aund| &\les& |\dk^{\leq 2}N|+|[\SS_\aund, \squared_2]\psi| +|[\SS_\aund, V]\psi| +a\left|\left[\SS_\aund, \frac{\cos\th}{|q|^2}\nab_T\right]\psi\right|, \qquad \aund=1, 2, 3,\\
  |N_4| &\les& |\dk^{\leq 2}N|+ r^{-2}|[\SS_4, |q|^2\squared_2]\psi| +r^{-2}|[\SS_4, |q|^2V]\psi| +ar^{-2}|[\SS_4, \cos\th\nab_T]\psi|.
\eeaa

Next, since 
\beaa
V= \frac{4\De}{ (r^2+a^2) |q|^2}, \qquad |q|^2V= \frac{4\De}{ (r^2+a^2)},
\eeaa
see \eqref{eq:Gen.RW}, we infer $[\SS_\aund, V]\psi=0$, $\aund=1, 2, 3$, and $[\SS_4, |q|^2V]\psi=0$, and hence
\beaa
|N_\aund| &\les& |\dk^{\leq 2}N|+|[\SS_\aund, \squared_2]\psi|  +a\left|\left[\SS_\aund, \frac{\cos\th}{|q|^2}\nab_T\right]\psi\right|, \qquad \aund=1, 2, 3,\\
  |N_4| &\les& |\dk^{\leq 2}N|+ r^{-2}|[\SS_4, |q|^2\squared_2]\psi| +ar^{-2}|[\SS_4, \cos\th\nab_T]\psi|.
\eeaa

Also, using Lemma \ref{lemma:commutation-nabT-nabZ}, we have
\beaa
\left[\SS_\aund, \frac{\cos\th}{|q|^2}\nab_T\right]\psi =0,\quad  \aund=1, 2, 3.
\eeaa
We also have
\beaa
r^{-2}|[\SS_4, \cos\th\nab_T]\psi| &\les&  r^{-1}|\dk^{\leq 1}\nab(\cos\th)||\dk^{\leq 1}\nab_T\psi|+r^{-2}|[\SS_4, \nab_T]\psi| \\
&\les&  r^{-2}|\dk^{\leq 2}\psi|.
\eeaa
We infer
\beaa
|N_\aund| &\les& |\dk^{\leq 2}N|+|[\SS_\aund, \squared_2]\psi|, \qquad \aund=1, 2, 3,\\
  |N_4| &\les& |\dk^{\leq 2}N|+ r^{-2}|[\SS_4, |q|^2\squared_2]\psi| +ar^{-2}|\dk^{\leq 2}\psi|.
\eeaa

Finally, we have 
\beaa
|[\SS_\aund, \squared_2]\psi| &\les& ar^{-2}|\dk^{\leq 2}\psi|, \qquad \aund=1, 2, 3.
\eeaa
Also, in view of Proposition \ref{LEMMA:MOD-LAPLACIAN-KERR}, we have
\beaa
r^{-2}|[\SS_4, |q|^2 \squared_2]\psi|  &\les& ar^{-2}|\dk^{\leq 2}\psi|.
\eeaa
We infer
\beaa
|N_\aund| &\les& |\dk^{\leq 2}N| + ar^{-2}|\dk^{\leq 2}\psi| \qquad \aund=1, 2, 3, 4,
\eeaa
as stated. This concludes the proof of  Lemma \ref{lemma-commuted-equ-Naund}.
 \end{proof}


\subsection{Basic spacetime  $\SS$-valued  identity}


We extend the definition of generalized current as given in Proposition \ref{prop-app:stadard-comp-Psi} to the case of double-indexed vector fields, functions and 1-forms.

\begin{definition}[Generalized Current] 
\lab{def:generalizedcurrent}  Let $\X$ be a double-indexed collection of vector fields $\X=\{ X^{\underline{a} \underline{b}}\}$, $\bold{w}$ be a double-indexed collection of functions $\bold{w}=\{ w^{\underline{ab}} \}$, and $\M=\{M^{\aund\bund} \}$   a double-indexed collection of  $1$-forms, all symmetric in the indices $\aund, \bund$.

Consider a solution $\psi\in \mathfrak{s}_2$ of equation  \eqref{eq:Gen.RW}  and let $\psi_\aund=\SSa\psi$ verifying  \eqref{eq:waveeqfor-psia}, i.e.
\beaa
 \squared_2\psia -V\psia&=&- \frac{ 4 a\cos\th}{|q|^2} \dual \nab_T  \psia+N_\aund.
 \eeaa
The generalized  current 
$\PP_\mu =  \PP_\mu^{(\bold{X}, \bold{w}, \M)} [\psi] $  associated to  $\psia, \psib$   is given by
 \bea
 \lab{eq:generalizedcurrent}
 \PP_\mu&=& \QQ[\psi]_{\underline{ab} \mu \nu} X^{\underline{ab} \nu}+\frac 1 2  w^{\underline{ab}} \, \Db_\mu  \psia\c  \psib-\frac 1 4 (\partial_\mu w^{\underline{ab}})\psia\c\psib  +\frac 1 4 M^{\aund\bund}_\mu\psi_\aund\c \psi_\bund.
  \eea
  We  also define
\beaa
\QQ_{\mu\nu}(\psi_\aund, \psi_\bund)&=& \Db_\mu  \psi_\aund \c \Db_\nu \psi _\bund
          -\frac 12 \g_{\mu\nu} \left(\g^{\a\b} \Db_\a \psi_\aund \c\Db_\b  \psi_\bund+ V\psi_\aund \c \psi_\bund\right)\\
          &=&
           \Db_\mu  \psi _\aund\c \Db _\nu \psi_\bund    -\frac 12 \g_{\mu\nu}  \LL[\psi_\aund, \psi_\bund], \\
          \LL[\psi_\aund, \psi_\bund]&=&\g^{\a\b} \Db_\a \psi_\aund\c\Db_\b  \psi _\bund+ V\psi_\aund \c \psi_\bund.
          \eeaa
   \end{definition}
   
   We also define, in analogy with \eqref{definition-EE-gen}, the  modified divergence   
     \bea\lab{definition-EE-gen-SSvalued}
      \begin{split}
\EE[\X, \w, \M] &:= \D^\mu  \PP_\mu[\X, \w, \M]- \NN[\X, \w], \\
 \NN[\X, \w]&:=\left(\nab_{X^{\aund\bund}}\psia+\frac 1 2   w^{\aund\bund} \psia\right)\c \big( \squared_2\psib -V\psib \big)\\
 & -\big(\rhod +\etab\wedge\eta\big)\nab_{(X^{\aund\bund})^4 e_4- (X^{\aund\bund})^3e_3}  \psi_\aund\c\dual\psi_\bund\\
 & -\frac{1}{2}\Im\Big(\tr\Xb H (X^{\aund\bund})^3 +\tr X\Hb (X^{\aund\bund})^4\Big)\c\nab\psi_\aund\c\dual\psi_\bund.
\end{split}
\eea

Using the above  notation we  derive the following analogue of Proposition \ref{proposition:Morawetz1}.
\begin{proposition}
\lab{proposition:Morawetz3}
The following  hold true
\begin{enumerate}
\item If  we choose 
\beaa
X^{\aund\bund} =\FF^{\aund\bund} \pr_r, \qquad w^{\aund\bund}=  |q|^2 \Div \big( |q|^{-2}  X ^{\aund\bund} \big)-w_{red}  ^{\aund\bund},
\eeaa
then the generalized current   defined in \eqref{definition-EE-gen-SSvalued}  verifies the identity
\bea\label{generalized-current-operator}
\begin{split}
 |q|^2\EE[\bold{X}, \bold{w}, \bold{M}]   &=\AA^{\aund\bund}  \nab_r\psia\c \nab_r\psib + \UU^{\a\b\aund\bund} \, \Db_\a \psia \c \Db_\b \psib  +\VV^{\aund\bund} \psia\c\psib \\
 & +\frac 1 4 |q|^2 \D^\mu\Big(M^{\aund\bund}_\mu\psi_\aund\c \psi_\bund \Big)
   \end{split}
   \eea
   where
   \beaa
   \bsplit
   \AA^{\aund\bund}&= \De \pr_r \FF^{\aund\bund}- \frac 1 2 \FF^{\aund\bund}\pr_r \De-\frac 1 2 \De w^{\aund\bund}_{red},\\
   \UU^{\a\b\aund\bund}&=  -\frac 1 2  \FF^{\aund\bund}\pr_r \left(\frac 1 \De\RR^{\a\b}\right)-\frac 1 2   w^{\aund\bund}_{red}\frac 1 \De \RR^{\a\b},\\
   \VV^{\aund\bund}&= -\frac 1 2 \left( X^{\aund\bund}\big(|q|^2\big) V  +|q|^2  X^{\aund\bund}(V)+\frac 1 2|q|^2 \square_\g  w^{\aund\bund}  + |q|^2  w^{\aund\bund}_{red} V \right).
   \end{split}
   \eeaa 
   
\item If in addition we choose, for  functions  $z$ and $h$, and a double-indexed function $f^{\aund\bund}$
\bea\lab{choice-FF-w-operator}
\FF^{\aund\bund}=- z h f^{\aund\bund}, \qquad   w^{\aund\bund} =- z \pr_r \big( h  f^{\aund\bund}  \big), \qquad            w^{\aund\bund} _{red}=  \FF^{\aund\bund}  z^{-1}\partial_r z, 
\eea
then
\beaa
\bsplit
  \UU^{\a\b\aund\bund}&=   \frac{ 1}{2}  h f^{\aund\bund} \pr_r\left( \frac z \De\RR^{\a\b}\right),\\
  \AA^{\aund\bund}&=-z^{1/2}\Delta^{3/2} \partial_r\left(h \frac{ z^{1/2}  f^{\aund\bund} }{\Delta^{1/2}}  \right),   
 \\
   \VV^{\aund\bund}&=\frac 1 4 \pr_r \left(\De \pr_r \Big( z \pr_r \big( h  f^{\aund\bund} \big) \Big)\right)+2  hf^{\aund\bund}  \pr_r \left(z \frac{\De}{r^2+a^2} \right). 
 \end{split}
\eeaa

\item If  $M^{\aund\bund} = v^{\aund\bund}(r) \pr_r$, for a double-indexed function $v=v^{\aund\bund}(r)$, we have
\bea\label{expression-Div-M}
\nn&&\frac 1 4 |q|^2 \Div\big( (\psia\c\psib)M^{\aund\bund} \big)\\ 
&=& \frac 1 4 |q|^2\left( 2 v^{\aund\bund}(r)\psia\c \nab_r \psib + \left(\pr_r v^{\aund\bund}+ \frac{2r}{|q|^2} v^{\aund\bund}\right) \psia\c\psib \right).
\eea  
\end{enumerate}
\end{proposition}

\begin{proof}
The proof follows exactly the same steps as in the proof of Proposition \ref{proposition:Morawetz1}.  For the convenience of the reader we    derive below the formula for $\VV$.  We write
\beaa
 \VV^{\aund\bund}&=& \VV_0^{\aund\bund}+ \VV_1^{\aund\bund},\\
 \VV_0^{\aund\bund}&=& -\frac 1 4 |q|^2 \square_\g  w^{\aund\bund}, \\
   \VV_1^{\aund\bund}&=& -\frac 1 2 \Big( X^{\aund\bund}\big(|q|^2\big) V  +|q|^2  X^{\aund\bund}(V) + |q|^2  w^{\aund\bund}_{red} V \Big).
 \eeaa
We first calculate $\VV_0$. Recall that
\beaa
w^{\aund\bund} &=&  |q|^2 \D_\a\big( |q|^{-2}  \big(X^{\aund\bund} \big)^\a\big)-w^{\aund\bund} _{red}=
  |q|^2 \D_\a\big( |q|^{-2} \FF^{\aund\bund}\pr_r ^\a\big)-\FF^{\aund\bund} z^{-1}\partial_r z \\
  &=&   |q|^2\pr_r\big(|q|^{-2}  \FFab\big)+\FFab( \D_\a\pr_r^\a) -\FF^{\aund\bund} z^{-1}\partial_r z.
\eeaa
Using that $\D_\a \pr_r^\a=\frac{1}{|q|^2} \pr_r \big(|q|^2\big)$, we have
\beaa
w^{\aund\bund} =   |q|^2\pr_r\big(|q|^{-2}  \FFab\big)+ \frac{1}{|q|^2} \pr_r \big(|q|^2\big)\FFab -\FF^{\aund\bund} z^{-1}\partial_r z=\pr_r \FFab -\FF^{\aund\bund} z^{-1}\partial_r z
= z \pr_r\left( \frac{\FFab}{z} \right).
\eeaa

Thus, in view of our choice for $\FF^{\aund\bund}=- z h f^{\aund\bund}$ in \eqref{choice-FF-w-operator}, 
\beaa
w^{\aund\bund} &=& z \pr_r\left( \frac{\FF^{\aund\bund}}{z} \right)=- z \pr_r \big( h  f^{\aund\bund} \big).
\eeaa

Recalling  formula   \eqref{eq:waveH(r)}   for  $\square  H(r) $,  we deduce
\beaa
\VV_0^{\aund\bund}&=& -\frac 1 4 |q|^2 \square_\g  w^{\aund\bund}=-\frac 1 4 \pr_r \Big(\De \pr_r \big(w^{\aund\bund}  \big)\Big)=\frac 1 4 \pr_r \left(\De \pr_r \Big( z \pr_r \big( h  f^{\aund\bund} \big) \Big)\right).
\eeaa
Remains to calculate
\beaa
\VV_1^{\aund\bund}&=& -\frac 1 2 \Big( X^{\aund\bund}\big(|q|^2\big) V  +|q|^2  X^{\aund\bund}(V) + |q|^2  w^{\aund\bund}_{red} V \Big)=-\frac 1 2  \Big( X^{\aund\bund}\big(|q|^2 V\big)    + |q|^2  w^{\aund\bund}_{red} V \Big).
\eeaa
Recalling that $|q|^2V=4\frac{\De}{r^2+a^2} $, $w^{\aund\bund}_{red}=  \FF^{\aund\bund} z^{-1}\partial_r z$  and $\FF^{\aund\bund}=-zhf^{\aund\bund}$    we deduce
\beaa
\VV_1^{\aund\bund}&=&- 2  \left( X^{\aund\bund}\left(\frac{\De}{r^2+a^2} \right)    +   w^{\aund\bund}_{red} \frac{\De}{r^2+a^2}  \right)\\
&=& - 2  \left( \FFab \pr_r \left(\frac{\De}{r^2+a^2} \right)    +    \FFab z^{-1}\partial_r z\frac{\De}{r^2+a^2}  \right)\\
&=& -2  z^{-1}  \pr_r \left(z \frac{\De}{r^2+a^2} \right)  \FFab = 2  hf^{\aund\bund}  \pr_r \left(z \frac{\De}{r^2+a^2} \right). 
\eeaa
This proves the proposition. 
\end{proof}

We simplify the notations by writing \eqref{generalized-current-operator} as
\bea\label{eq:separation-EE-I-J-K}
\begin{split}
 |q|^2\EE[\bold{X}, \bold{w}, \bold{M}]   &=P+I  +J +K
   \end{split}
   \eea
where
\beaa
P&:=& \UU^{\a\b\aund\bund} \, \Db_\a \psia \c \Db_\b \psib,\\
I&:=& \AA^{\aund\bund}  \nab_r\psia\c \nab_r\psib, \\
J&:=& \VV^{\aund\bund} \psia\c\psib,\\
K&:=& \frac 1 4 |q|^2 \D^\mu\Big(M^{\aund\bund}_\mu\psi_\aund \c\psi_\bund \Big).
\eeaa

In what follows, we will choose functions $z$, $h$ and the double-indexed function $f^{\aund\bund}$ to obtain positivity of the generalized current. Each of the above term will produce lower order terms in derivatives and in angular momentum $a$.


\subsection{Principal term}


We consider now the principal  term $P$, i.e. from Proposition \ref{proposition:Morawetz3}
\beaa
P&=& \UU^{\a\b\aund\bund} \, \Db_\a \psia \c \Db_\b \psib= \frac{ 1}{2}  h f^{\aund\bund} \pr_r\left( \frac z \De\RR^{\a\b}\right) \, \Db_\a \psia \c \Db_\b \psib.
\eeaa
Recall that, see \eqref{eq:RR-Sa}, 
          \beaa
          \RR^{\a\b} =\RR^\aund  \Sa^{\a\b}.
          \eeaa
We similarly set
\beaa
\RRtp^{\aund}:&=&  \pr_r\left( \frac z \De\RR^{\aund}\right),
\eeaa
and thus write
\beaa
 \pr_r\left( \frac z \De\RR^{\a\b}\right)&=& \pr_r\left( \frac z \De \RR^\aund  \Sa^{\a\b}\right)=  \pr_r \Big(\frac{z}{\De} \RR^\aund \Big) \Sa^{\a\b} =\RRtp^{\aund} \Sa^{\a\b},
\eeaa
which gives
\beaa
P&=& \frac{ 1}{2}  h f^{\aund\bund} \RRtp^{\underline{c}} S_{\underline{c}}^{\a\b} \, \Db_\a \psia \c \Db_\b \psib.
\eeaa
To create a quadratic  expression in $\RRtp^{\aund}$ in $P$, we  choose
 \bea\lab{choice-f-operator}
  f^{\underline{a}\underline{b}}&=&   \tilde{\RR}'^{(\underline{a}} \LL^{\underline{b})}= \frac 1 2 \big(\RRtp^{\aund}\LL^{\bund}+\RRtp^{\bund}\LL^{\aund} \big),
 \eea
where $\LL^{\aund} $ are constant coefficients to be chosen later of  a given 2-tensor of the form
\bea
\lab{eq:operato-LL2}\lab{eq:operato-LL}
 L^{\a\b} &=& \LL^{\underline{a}} S^{\a\b}_{\underline{a}}= \LL^{1} T^\a T^\b+ \LL^{2} aT^{(\a}  Z^{\b)} +\LL^{3}a^2Z^\a Z^\b+\LL^4 O^{\a\b}.
\eea

 With the choice \eqref{choice-f-operator} for $f^{\aund\bund}$ we deduce
 \beaa
P&=& \frac{ 1}{4}  h  \big(\RRtp^{\aund}\LL^{\bund}+\RRtp^{\bund}\LL^{\aund} \big) \RRtp^{\underline{c}} \Sc^{\a\b} \, \Db_\a \psia \c \Db_\b \psib\\
&=& \frac{ 1}{4}  h  \big(\RRtp^{\aund}\RRtp^{\underline{c}}\LL^{\bund}+\RRtp^{\bund}\RRtp^{\underline{c}}\LL^{\aund} \big)  \Sc^{\a\b} \, \Db_\a \psia \c \Db_\b \psib\\
&=& \UU^{\aund\bund\,\a\b}  \, \Db_\a \psia \c \Db_\b \psib,
\eeaa
where
\bea
\lab{definition:UUab}
\UU^{\aund\bund}:=\frac 1 2  h\RRtp^{\aund}  \RRtp^{\bund} , \qquad   \UU^{\aund\bund\,\a\b} := \frac 1  2\big(\UU^{\aund\cund} \LL^{\bund}+\UU^{\bund\cund} \LL^\aund \big)\Sc^{\a\b}= \UU^{\cund (\aund} \LL^{\bund)} \Sc^{\a\b}.
\eea

We now isolate a positive part in the principal term $P$ from a divergence component. This is the crucial step which allows to obtain a positive trapped term in the case of higher derivative Morawetz estimates.

   \begin{lemma}
   \lab{lemma:IntegrationbypartsP}
Let  $P$ the principal term defined as above. 
We then have the identity
\beaa
 P&=&\frac 1 2  h  L^{\a\b}  \Db_\a \Psi\c    \Db_\b \Psi -\frac 1 2  h\Psi\c ( \RRtp^{\cund}  \LL^{\bund}[\SS_{\cund}, \SS_{\bund}]\psi) +|q|^2\Db_\a \BB^\a,
\eeaa
where $\Psi$ is defined as 
\bea
\Psi:= \RRtp^{\aund} \psia ,
\eea
and the boundary term $\BB$ is given by
\bea\label{eq:definition-BB}
\BB^\a&:=& |q|^{-2} \frac 1 2  h \Psi  \RRtp^{\cund}  \LL^{\bund}\c\left(\Sc^{\a\b}   \, \Db_\b \psib-  \Sb^{\a\b}   \Db_\b \psic \right).
\eea
\end{lemma}

\begin{proof}
We consider
\beaa
|q|^{-2} P&=& \frac 1  2|q|^{-2}\big(\UU^{\aund\cund} \LL^{\bund}+\UU^{\bund\cund} \LL^\aund \big)\Sc^{\a\b}  \, \Db_\a \psia \c \Db_\b \psib=|q|^{-2} \UU^{\aund\cund} \LL^{\bund}\Sc^{\a\b}  \, \Db_\a \psia \c \Db_\b \psib\\
&=& - \UU^{\aund\cund} \LL^{\bund} \psia\c \Db_\a(|q|^{-2} \Sc^{\a\b}   \, \Db_\b \psib)+\Db_\a( |q|^{-2}\UU^{\aund\cund} \LL^{\bund}\Sc^{\a\b}   \psia \c \Db_\b \psib)\\
&& - \Db_\a( \UU^{\aund\cund} \LL^{\bund}) |q|^{-2}\Sc^{\a\b}   \psia \c \Db_\b \psib.
\eeaa
Note  that        $\LL^{\aund}$  and  $ \UU^{\aund\bund}$ only depend on $r$. Therefore, since  for all $\cund$  the 2-tensors $ \Sc^{\a\b}$   do not contain $r$ derivatives we deduce
\beaa
\Db_\a\big( \UU^{\aund\cund} \LL^{\bund} \big)|q|^{-2} \Sc^{\a\b}  \psia \c \Db_\b \psib=   \pr_r\big( \UU^{\aund\cund}  \LL^{\bund}  \big)  |q|^{-2}\Sc^{r \b} \psia \c \Db_\b \psib=0.
\eeaa
Now consider $\Db_\a(|q|^{-2} \Sc^{\a\b}   \, \Db_\b \psib)$. Recall Definition \ref{definition:symmetry-tensors}, i.e. for $\psi \in \sk_2$,
\beaa
\SS_\aund\psi&=&|q|^2 \Ddot_\a(|q|^{-2}S_\aund^{\a\b} \Ddot_\b \psi) \qquad \text{ for $a=1,2,3, 4$.}
\eeaa
Therefore we write
\beaa
 \Db_\a(|q|^{-2} \Sc^{\a\b}   \, \Db_\b \psib)&=& \Db_\a(|q|^{-2} \Sc^{\a\b}   \, \Db_\b \SS_\bund\psi)=|q|^{-2}\SS_\cund \SS_\bund\psi=|q|^{-2}\SS_\bund \SS_\cund\psi+|q|^{-2}[\SS_{\cund}, \SS_{\bund}]\psi\\
 &=& \Db_\a(|q|^{-2} \Sb^{\a\b}   \, \Db_\b \psic)+|q|^{-2}[\SS_{\cund}, \SS_{\bund}]\psi.
\eeaa

 Thus, repeating the integration by parts procedure and noting, as above, that the last term vanishes, we obtain 
\beaa
  &&\UU^{\aund\cund} \LL^{\bund} \psia\c    \Db_\a \big(|q|^{-2}\Sc^{\a\b}  \Db_\b  \psib\big)\\
  &=&  \UU^{\aund\cund} \LL^{\bund} \psia\c   \Db_\a(|q|^{-2} \Sb^{\a\b}   \, \Db_\b \psic) +|q|^{-2}\UU^{\aund\cund} \LL^{\bund} \psia\c [\SS_{\cund}, \SS_{\bund}]\psi  \\
&=& -|q|^{-2}\UU^{\aund\cund} \LL^{\bund} \Sb^{\a\b}  \Db_\a \psia\c      \Db_\b \psic+\Db_\a \Big(|q|^{-2} \UU^{\aund\cund} \LL^{\bund}  \Sb^{\a\b}   \psia\c   \Db_\b \psic\Big)\\
&&+|q|^{-2}\UU^{\aund\cund} \LL^{\bund} \psia\c [\SS_{\cund}, \SS_{\bund}]\psi.
\eeaa
Therefore, recalling that $  \LL^{\bund}   \Sb^{\a\b}=L^{\a\b}  $,
\beaa
 P&=& \UU^{\aund\cund}L^{\a\b}  \Db_\a \psia\c     \Db_\b \psic+|q|^2\Db_\a \left(|q|^{-2} \UU^{\aund\cund} \LL^{\bund}  \psia\c \left(\Sc^{\a\b}   \, \Db_\b \psib-  \Sb^{\a\b}   \Db_\b \psic \right) \right ) \\
&& -\UU^{\aund\cund} \LL^{\bund} \psia\c [\SS_{\cund}, \SS_{\bund}]\psi.
\eeaa

Finally, using \eqref{definition:UUab}, and using that $\RRtp^{\aund}$ only depend on $r$, we obtain
\beaa
 P&=&\frac 1 2  h\RRtp^{\aund}  \RRtp^{\cund} L^{\a\b}  \Db_\a \psia\c     \Db_\b \psic -\frac 1 2  h\RRtp^{\aund}  \RRtp^{\cund}  \LL^{\bund} \psia\c [\SS_{\cund}, \SS_{\bund}]\psi\\
&& +|q|^2\Db_\a \left(|q|^{-2} \frac 1 2  h\RRtp^{\aund}  \RRtp^{\cund}  \LL^{\bund}  \psia\c \left(\Sc^{\a\b}   \, \Db_\b \psib-  \Sb^{\a\b}   \Db_\b \psic \right) \right ) \\
&=&\frac 1 2  h  L^{\a\b}  \Db_\a (\RRtp^{\aund} \psia )\c    \Db_\b(\RRtp^{\cund} \psic) -\frac 1 2  h \RRtp^{\cund}  \LL^{\bund}(\RRtp^{\aund}  \psia )\c[\SS_{\cund}, \SS_{\bund}]\psi\\
&& +|q|^2\Db_\a \left(|q|^{-2} \frac 1 2  h  \RRtp^{\cund}  \LL^{\bund} (\RRtp^{\aund} \psia )\c\left(\Sc^{\a\b}   \, \Db_\b \psib-  \Sb^{\a\b}   \Db_\b \psic \right) \right ) .
\eeaa
By denoting $\Psi=\RRtp^{\aund} \psia$ we obtain the stated expression.
\end{proof}

We summarize the result in the following  proposition.
\begin{proposition}
\lab{proposition:Morawetz4}\lab{proposition:expression-for-Psi}
Consider the generalized bilinear current  \eqref{eq:generalizedcurrent}  with
 \beaa
   X^{\aund\bund} &=&\FF^{\aund\bund}  \pr_r, \qquad  w^{\aund\bund} =  |q|^2 \D_\a\big( |q|^{-2}  \big(X^{\aund\bund} \big)^\a\big)-w^{\aund\bund} _{red},\qquad  w^{\aund\bund}_{red}=  \FF^{\aund\bund} z^{-1}\partial_r z, 
  \eeaa
  and 
 \beaa
  \FF^{\underline{a}\underline{b}}&=& - z h  \tilde{\RR}'^{(\underline{a}} \LL^{\underline{b})}=-\frac 1 2 z h  \big(\RRtp^{\aund}\LL^{\bund}+\RRtp^{\bund}\LL^{\aund} \big), \qquad  \RRtp^{\aund}=      \pr_r\Big( \frac{z}{\De} \RRa\Big),
 \eeaa
for constant $\LL^\aund$.
 Then, 
  \bea\lab{Morawetz-effective-principal-term}
\begin{split}
 |q|^2\EE[\bold{X}, \bold{w}, \bold{M}] -|q|^2 \D^\mu \BB_\mu   &=\widetilde{P}+P_{lot}+I  +J +K
   \end{split}
   \eea
 where
 \begin{itemize}
 \item $\BB$ is the boundary term defined in \eqref{eq:definition-BB},
 
 \item The principal term $   \widetilde{P}$ is positive definite and given by
 \bea\label{eq:definition-widetilde}\label{eq:widetilde-P}
   \widetilde{P}&:=&\frac 1 2  h  L^{\a\b}  \Db_\a \Psi\c    \Db_\b \Psi , \qquad \Psi=\RRtp^{\aund} \psia,
 \eea
 with $L^{\a\b} $  as in  \eqref{eq:operato-LL2},
 
 \item The lower order term $P_{lot}$ is given by 
 \bea\label{eq:Plot}
 P_{lot}&=&-\frac 1 2  h\Psi\c ( \RRtp^{\cund}  \LL^{\bund}[\SS_{\cund}, \SS_{\bund}]\psi).
 \eea
 
 \item The quadratic form $I$ is given by 
 \bea
I&=& \big(\AAa[z]\nab_r \psia\big)\c \big( \LL^{\aund}\nab_r \psia\big) \label{eq:expr-I}, 
\eea
where
\beaa
 \AAa[z]=  - z^{1/2}\De^{3/2}   \RRtpp^{\aund},  \qquad  \RRtpp^\aund:=  \pr_r\Big( \frac{ h  z^{1/2}  \RRtp^{\aund}  }{\De^{1/2} } \Big).
\eeaa

\item The quadratic form $J$ is given by 
\bea
J&=& \big(\VVa[z]  \psia\big)\c \big( \LL^{\aund} \psia\big),\label{eq:expr-J}
\eea
where
\beaa
   \VV^{\aund}= \frac 1 4  \pr_r\left(\De \pr_r \Big(
 z \pr_r \big( h  \tilde{\RR}'^{\underline{a}}  \big)  \Big)  \right)+2   h   \tilde{\RR}'^{\underline{a}}  \pr_r \left(\frac{z\De}{r^2+a^2} \right) .
\eeaa

\item The quadratic form $K$ is given by
\bea
K&=&\frac 1 4 |q|^2 \D^\mu\Big((M^{\aund}_\mu\psi_\aund )\c (\LL^\bund   \psi_\bund) \Big).\label{eq:expr-K}
   \eea 
 \end{itemize}
 \end{proposition}
 
\begin{proof} 
The above follows from \eqref{eq:separation-EE-I-J-K} and Lemma \ref{lemma:IntegrationbypartsP}. Also from Proposition \ref{proposition:Morawetz3} and the choice \eqref{choice-f-operator}
 $ f^{\underline{a}\underline{b}}=   \tilde{\RR}'^{(\underline{a}} \LL^{\underline{b})}$, we deduce 
    \beaa
    \AA^{\aund\bund}=- z^{1/2}\De^{3/2} \RRtpp^{(\aund}\LL^{\bund)}, \qquad  \RRtpp^\aund:=  \pr_r\left( \frac{ h  z^{1/2}  \RRtp^{\aund}  }{\De^{1/2} } \right),
\eeaa
and
 \beaa
    \VV^{\aund\bund}&=& \frac 1 4  \pr_r\left(\De \pr_r \Big(
 z \pr_r \big( h  \tilde{\RR}'^{(\underline{a}} \LL^{\underline{b})}  \big)  \Big)  \right)+2   h   \tilde{\RR}'^{(\underline{a}} \LL^{\underline{b})}  \pr_r \left(\frac{z\De}{r^2+a^2} \right).  
 \eeaa
The expressions for $I$, $J$ and $K$ are implied by writing $\AA^{\aund\bund}=\AA^{(\aund} \LL^{\bund)}$, $\VV^{\aund\bund}=\VV^{(\aund} \LL^{\bund)}$ and $M^{\aund\bund}=M^{(\aund} \LL^{\bund)}$.
\end{proof}

We denote the effective generalized current 
\bea\lab{widetilde-EE}
 |q|^2 \widetilde{\EE}[\bold{X}, \bold{w}, \bold{M}]:= |q|^2 \EE[\bold{X}, \bold{w}, \bold{M}]  -|q|^2 \D^\mu \BB_\mu 
\eea
which differs from the generalized bilinear current by the boundary terms given by $\BB$. Our goal in deriving the higher order Morawetz estimates is to show positivity of the effective quadratic form $ |q|^2 \widetilde{\EE}[\bold{X}, \bold{w}, \bold{M}]$.


\subsection{Choice of $z$}


In this section, we present two choices for the function $z$. 

 \begin{enumerate}
 \item The first is the same choice as in Section   \ref{section:choicezhf}, i.e.
 \beaa
 z=z_0=\frac{\De}{(r^2+a^2)^2}.
 \eeaa
 This choice, which was good enough to derive conditional Morawetz estimates for first derivatives, has its limitations in the case of higher derivatives estimates. In particular, the contribution in $\SS_1=\nab_T \nab_T$ vanishes with this choice. We therefore have to modify it, but we use it as reference for the next choice.
 \item The crucial choice in the next section is the following (see also \cite{A-B}):
 \bea
 \lab{choice-z-operator-est} 
 z=z_0-\de_0 z_0^2, \qquad \de_0>0.
 \eea
The parameter $\de_0$ will be chosen later to be a small positive constant. This allows to obtain a contribution in $\SS_1$ with this choice, and for $\de_0$ small enough we can transpose most information from the first choice of $z$ to the second one.
 \end{enumerate}

We now collect the main coefficients for the two choices.


\subsubsection{Main coefficients for  the choice $z=z_0$}


\begin{lemma}\label{lemma:coeff-z=z0}
With the choice  $z= \frac{\De}{(r^2+a^2)^2} $        and $h=\frac{(r^2+a^2)^4}{r(r^2-a^2)} $
 (as   in Section \ref{section:choicezhf})
 we have the following coefficients.
\begin{enumerate}
\item The coefficients $\RRtp^\aund= \pr_r\Big( \frac{z}{\De} \RRa\Big)$ are given by
 \bea
  \lab{eq:valuesforRRtp}
  \bsplit
\RRtp^1&=0, \qquad \RRtp^2=  \frac{4r}{(r^2+a^2)^2} , \\ 
\RRtp^3&=\frac{4r}{(r^2+a^2)^3}, \qquad \RRtp^4=f=  -\frac{2\TT}{ (r^2+a^2)^3}.
\end{split}
\eea

\item The coefficients  $\RRtpp^\aund=\pr_r \Big(  \frac{ h  z^{1/2}  \RRtp^{\aund}  }{\De^{1/2} } \Big)$
  are given by 
 \beaa
 \bsplit
 \RRtpp^1&=0, \qquad  \RRtpp^2 =-\frac{16a^2r}{(r^2-a^2)^2},\qquad \RRtpp^3= -\frac{8r}{(r^2-a^2)^2}, \\
 \RRtpp^4&=- 2 \frac{3mr^4-4a^2r^3+ma^4}{r^2(r^2-a^2)^2}.
 \end{split}
 \eeaa
 
 \item The coefficients $\AAa=  - z^{1/2}\De^{3/2}   \RRtpp^{\aund}$
   are given by
   \bea
   \lab{eq:AAcoefficients-z=z_0} 
   \begin{split}
   \AA^1&=0, \quad \AA^2=\De^2\frac{16a^2r}{(r^2+a^2)(r^2-a^2)^2}, \quad \AA^3= \De^2 \frac{8r}{(r^2+a^2)(r^2-a^2)^2},\\
    \AA^4&=  2\De^2 \frac{3mr^4-4a^2r^3+ma^4}{r^2(r^2+a^2)(r^2-a^2)^2}.
    \end{split}
   \eea
   
\item The coefficients $\VVa$ verify
\beaa
\VV^1 = 0, \qquad \VV^2 = O(r^{-1}), \qquad \VV^3 = O(r^{-3}), \qquad \VV^4=\VV, 
\eeaa
where $\VV$ is the scalar function given by \eqref{eq:prop:Choice-zhf}.
\end{enumerate}
\end{lemma}

\begin{proof}
The values of $\RRtp$ are the same as the one computed in Section \ref{section:morawetz-geodesics} in the case of geodesics. For $\RRtpp$ we compute
\beaa
 \RRtpp^2&=&\pr_r \left( \frac{(r^2+a^2)^3}{r(r^2-a^2)}   \RRtp^2 \right)=4\pr_r \left( \frac{(r^2+a^2)}{(r^2-a^2)}    \right)=-\frac{16a^2r}{(r^2-a^2)^2},\\
 \RRtpp^3&=&\pr_r \left( \frac{(r^2+a^2)^3}{r(r^2-a^2)}   \RRtp^3 \right)=4\pr_r \left( \frac{1}{(r^2-a^2)} \right)=-\frac{8r}{(r^2-a^2)^2},\\
  \RRtpp^4&=&\pr_r \left( \frac{(r^2+a^2)^3}{r(r^2-a^2)}   \RRtp^4 \right)=-2\pr_r \left( \frac{\TT}{r(r^2-a^2)} \right)=- 2 \frac{3mr^4-4a^2r^3+ma^4}{r^2(r^2-a^2)^2},
\eeaa
as desired. The computations of $\AA^{\aund}$ are immediate. Finally, $\VV^1=0$ follows from $\RRtp^1=0$, $\VV^2 = O(r^{-1})$ from $\RRtp^2=O(r^{-3})$, $\VV^3 = O(r^{-3})$ from $\RRtp^3=O(r^{-5})$, and $\VV^4=\VV$ from $\RRtp^4=f$ and Proposition \ref{prop:Choice-zhf}.
\end{proof}


\subsubsection{Main coefficients for  the choice $z= z_0-\de_0z_0^2$}


\begin{lemma}\lab{lemma:mainchoicesofzandhandcorrecpondingquantities:faposdifah}
With the choice  $z =z_0- \de_0 z_0^2$ as in \eqref{choice-z-operator-est}  with $z_0= \frac{\De}{(r^2+a^2)^2} $        and $h=\frac{(r^2+a^2)^4}{r(r^2-a^2)} $
 we have the following coefficients.
\begin{enumerate}
\item The coefficients $\RRtp^a[z] =  \pr_r\Big( \frac{z}{\De} \RRa\Big)$ are given by
\bea
\label{asymptotics-RRtp-new}
\bsplit
\RRtp^1[z]&=\de_0 f, \\
\RRtp^2[z]& =\frac{4r}{(r^2+a^2)^2} \big(1+O(r^{-2} \de_0) \big), \\
\RRtp^3[z]  &=\frac{4r}{(r^2+a^2)^3} \big(1+O(r^{-2} \de_0) \big), \\
\RRtp^4[z] &= f( 1-2\de_0 z_0), 
\end{split}
\eea
where,  as in \eqref{choice-f-wave}, $f=\pr_r z_0=-\frac{2\TT}{ (r^2+a^2)^3}$.

\item The coefficients  \, $  \RRtpp^\aund[z]=\pr_r \Big(  \frac{ h  z^{1/2}  \RRtp^{\aund}  }{\De^{1/2} } \Big)$
  are given by 
 \bea
 \label{asymptotics-RRtpp-new}
 \bsplit
 \RRtpp^1[z]&= -\frac{2(3mr^4-4a^2r^3+ma^4)}{r^2(r^2-a^2)^2}\de_0+O(\de_0^2 r^{-3} ),\\
 \RRtpp^2[z]&=-\frac{16a^2r}{(r^2-a^2)^2} +O( \de_0 r^{-3} ), \\
 \RRtpp^3[z]&= -\frac{8r}{(r^2-a^2)^2}\Big( 1+O(\de_0 r^{-2})\Big),\\
 \RRtpp^4[z]&=\frac{-2(3mr^4-4a^2r^3+ma^4)}{r^2(r^2-a^2)^2}+O(\de_0 r^{-3} ).
 \end{split}
 \eea 
 
 \item The coefficients $ \AAa[z]=  - z^{1/2}\De^{3/2}   \RRtpp^{\aund}$
   are given by
   \bea\lab{eq:explicit-AA-}
\bsplit
 \AA^1[z]&= \de_0\De^2\left( \frac{2(3mr^4-4a^2r^3+ma^4)}{r^2(r^2+a^2)(r^2-a^2)^2}+ O(\de_0 r^{-5} )\right),\\
 \AA^2[z]&= \frac{2\De^2}{r^2+a^2} \left(\frac{8a^2r}{(r^2-a^2)^2}+O(\de_0 r^{-3})\right), \\
 \AA^3[z]&= \De^2\left( \frac{8r}{(r^2+a^2)(r^2-a^2)^2}+ O(\de_0 r^{-7} )\right), \\
  \AA^4[z]&= \De^2 \left(  \frac{2(3mr^4-4a^2r^3+ma^4)}{r^2(r^2+a^2)(r^2-a^2)^2}+ O(\de_0 r^{-5} )\right).
 \end{split}
\eea 

\item The coefficients $\VVa$ verify
\bea\label{expressions-VV-asymp}
\VV^1 = \de_0\VV, \qquad \VV^2 = O(r^{-1}), \qquad \VV^3 = O(r^{-3}), \qquad \VV^4=\VV+O(\de_0 r^{-3}), 
\eea
where $\VV$ is the scalar function given by \eqref{eq:prop:Choice-zhf}.
\end{enumerate}
\end{lemma}

\begin{proof}
With this choice  for $z$ we deduce
\beaa
\RRtp^a[z] &=&  \pr_r\Big( \frac{z}{\De} \RRa\Big)=  \pr_r\Big( \frac{z_0(1-\de z_0)}{\De} \RRa\Big)= (1-\de z_0)\RRtp^a[z_0]-\de    ( \pr_r z_0 ) \frac{z_0}{\De} \RR^\aund.
\eeaa
Hence,
\bea
\lab{eq:NewRRtpz}
\RRtp^a[z] = (1-\de_0 z_0)\RRtp^a[z_0]- \de_0 f \frac{1}{(r^2+a^2)^2} \RR^\aund.
\eea
In particular, using Lemma \ref{lemma:coeff-z=z0} for the coefficients $\RRtp^a[z_0]$, and \eqref{components-RR-aund} for $\RR^\aund$, we obtain for $\de_0>0$ sufficiently small and all $r\ge r_+$,
\beaa
\RRtp^1[z]&=&\de_0 f, \\
\RRtp^2[z]&=&\frac{4r}{(r^2+a^2)^2} -\de_0 \frac{4r\De}{(r^2+a^2)^4}   + \de_0 f \frac{2}{(r^2+a^2)} =  \frac{4r}{(r^2+a^2)^2} \big(1+O(r^{-2} \de_0) \big), \\
\RRtp^3[z]  &=& \frac{4r}{(r^2+a^2)^3} -\de_0 \frac{4 r\De}{(r^2+a^2)^5}+ \de_0 f \frac{1}{(r^2+a^2)^2} =\frac{4r}{(r^2+a^2)^3} \big(1+O(r^{-2} \de_0) \big), \\
\RRtp^4[z] &=& f -\de_0 \frac{2\De}{(r^2+a^2)^2} f= f( 1-2\de_0 z_0),
\eeaa
which ends the proof of part 1. To check part 2,  we write
 \beaa
z^{1/2}=z_0^{1/2}(1-\de_0 z_0)^{1/2}=\frac{\De^{1/2}}{r^2+a^2}(1-\de_0 z_0)^{1/2}.
\eeaa
 Relying on   the   same function $h$ as before  i.e.  $h=\frac{(r^2+a^2)^4}{ r(r^2+a^2)} $, we deduce
\beaa
\RRtpp^\aund&=&  \pr_r\left( \frac{ h  z^{1/2}  \RRtp^{\aund}[z]  }{\De^{1/2} } \right)= \pr_r\left(\frac{(r^2+a^2)^3}{ r(r^2-a^2)}(1-\de_0 z_0)^{1/2} \RRtp^{\aund}[z]  \right) \\
&=&  \pr_r\left(\frac{(r^2+a^2)^3}{ r(r^2 - a^2)}    \RRtp^\aund[z]  \right)(1-\de_0 z_0)^{1/2} +\frac{(r^2+a^2)^3}{ r(r^2 - a^2)}    \RRtp^\aund[z] \left(-\frac{\de_0}{2}\right)\frac{\pr_r z_0}{(1-\de_0 z_0)^{1/2}}\\
&=&  \pr_r\left(\frac{(r^2+a^2)^3}{ r(r^2 - a^2)}    \RRtp^\aund[z]  \right)(1-\de_0 z_0)^{1/2} +\frac{\de_0\TT}{ r(r^2 - a^2)(1-\de_0 z_0)^{1/2}}    \RRtp^\aund[z].
\eeaa
We therefore compute, relying on   \eqref{asymptotics-RRtp-new},
\beaa
\RRtpp^1&=& \pr_r\left(\frac{(r^2+a^2)^3}{ r(r^2 - a^2)}    \RRtp^1[z]  \right)(1-\de_0 z_0)^{1/2}+\frac{\de_0\TT}{ r(r^2 - a^2)(1-\de_0 z_0)^{1/2}}    \RRtp^1[z]\\
&=& \de_0\pr_r\left(\frac{(r^2+a^2)^3f}{ r(r^2 - a^2)}    \right)(1-\de_0 z_0)^{1/2}+\frac{\de_0^2f\TT}{ r(r^2 - a^2)(1-\de_0 z_0)^{1/2}}\\
&=& -\de_0  \pr_r\left(\frac{2\TT}{r (r^2-a^2)} \right) + O(\de_0^2 r^{-3} )=- 2 \frac{3mr^4-4a^2r^3+ma^4}{r^2(r^2-a^2)^2}\de_0  + O(\de_0^2r^{-3}), \\
\RRtpp^3&=& \pr_r\left(\frac{(r^2+a^2)^3}{ r(r^2 - a^2)}    \RRtp^3[z]  \right)(1-\de_0 z_0)^{1/2}+\frac{\de_0\TT}{ r(r^2 - a^2)(1-\de_0 z_0)^{1/2}}    \RRtp^3[z]\\
&=&\pr_r \left( \frac{4}{r^2-a^2}\right) + O(\de_0 r^{-5} )=-\frac{8r}{(r^2-a^2)^2} + O(\de_0 r^{-5} ),\\
\RRtpp^4&=& \pr_r\left(\frac{(r^2+a^2)^3}{ r(r^2 - a^2)}    \RRtp^4[z]  \right)(1-\de_0 z_0)^{1/2}+\frac{\de_0\TT}{ r(r^2 - a^2)(1-\de_0 z_0)^{1/2}}    \RRtp^4[z]\\
&=& - \pr_r\Big(\frac{2\TT}{r (r^2-a^2)} \Big)  + O(\de_0 r^{-3})  =\frac{-2(3mr^4-4a^2r^3+ma^4)}{r^2(r^2-a^2)^2} + O(\de_0 r^{-3}).
\eeaa
We also have
\beaa
\RRtpp^2&=& \pr_r\left(\frac{(r^2+a^2)^3}{ r(r^2 - a^2)}    \RRtp^2[z]  \right)(1-\de_0 z_0)^{1/2}+\frac{\de_0\TT}{ r(r^2 - a^2)(1-\de_0 z_0)^{1/2}}    \RRtp^2[z]\\
&=& \pr_r\left(\frac{(r^2+a^2)^3}{ r(r^2 - a^2)}\left(\frac{4r}{(r^2+a^2)^2} -\de_0 \frac{4r\De}{(r^2+a^2)^4}+ \de_0 f \frac{2}{(r^2+a^2)}\right)\right)(1-\de_0 z_0)^{1/2}\\
&& +\frac{4\de_0r\TT}{ r(r^2 - a^2)(r^2+a^2)^2(1-\de_0 z_0)^{1/2}}= -\frac{16a^2r}{(r^2-a^2)^2}  +O( \de_0 r^{-3} ).
\eeaa

The second part of the lemma is immediate from  the definition
 \beaa
  \AA^{\aund}=-  z^{1/2}\De^{3/2} \RRtpp^{\aund}=-\frac{\De^2}{r^2+a^2} \Big( 1+O(\de_0 r^{-2}) \Big)  \RRtpp^{\aund}.
  \eeaa
  The last part follows easily from the definition of $\VVa$ and the case $\de_0=0$ in Lemma \ref{lemma:coeff-z=z0}.
\end{proof}

\begin{remark}
Note that    we can write 
\bea
\lab{simplifications-AA}
\bsplit
\AA^1&= \de_0 \AA\big(1+O(r^{-1} \de_0)\big), \quad \,\AA^4=\AA\big(1+O(r^{-1}\de_0)\big), \\
 \AA^2&=\AAt \big(2a^2+O(\de_0 )\big), \qquad\,\,\,\,   \AA^3=  \tilde{\AA}\big(1+O(r^{-2}\de_0 )\big),
 \end{split}
\eea
where 
\beaa
\AA&=&   \frac{2\De^2}{r^2(r^2+a^2)(r^2-a^2)^2}(3mr^4-4a^2r^3+ma^4), \qquad \widetilde{\AA}= \frac{8\De^2r}{(r^2+a^2)(r^2-a^2)^2},
\eeaa
are positive coefficients for $|a|/m <1$.  We also note the term  $\AA$ is precisely the one appearing  in  Proposition \ref{prop:Choice-zhf}.
\end{remark}

 
 \section{Computation of the effective generalized current}
 \label{section:conditional-SS-values-Mor}
 
 
 From now on, we consider the choices for $z$ and $h$ made in Lemma \ref{lemma:mainchoicesofzandhandcorrecpondingquantities:faposdifah}, i.e.
\beaa
z =z_0- \de_0 z_0^2, \qquad z_0= \frac{\De}{(r^2+a^2)^2}, \qquad \de_0>0, \qquad h=\frac{(r^2+a^2)^4}{r(r^2-a^2)}.
\eeaa
Recall that according to Proposition \ref{proposition:Morawetz4}, the effective generalized current $\widetilde{\EE}[\bold{X}, \bold{w}, \bold{M}]$ introduced in \eqref{widetilde-EE}  is given by, see \eqref{Morawetz-effective-principal-term},
 \beaa
 |q|^2 \widetilde{\EE}[\bold{X}, \bold{w}, \bold{M}] &=&\widetilde{P}+P_{lot}+I  +J +K.
   \eeaa 
In this section, we compute the terms $\widetilde{P}$, $P_{lot}$, $I$, $J$ and $K$.


\subsection{Principal trapping term}


In this section, we analyze the principal trapping term $\widetilde{P}$ in \eqref{eq:definition-widetilde}, i.e.
\beaa
   \widetilde{P}&=&\frac 1 2  h  L^{\a\b}  \Db_\a \Psi_z    \Db_\b \Psi_z , \qquad \Psi_z=\RRtp^{\aund}[z] \psia.
   \eeaa
Observe that with the choice $z=z_0-\de_0 z_0^2 $, using  \eqref{asymptotics-RRtp-new}, we have
\beaa
\Psi_z&=& f \big(\de_0  \SS_1\psi+ (1+O(r^{-2} \de_0)) \OO\psi\big)+ \frac{4r}{(r^2+a^2)^2}  \Big( \SS_2\psi+\frac{1}{r^2+a^2} \SS_3\psi\Big)     \big(1+O(r^{-2} \de_0) \big).
\eeaa
Observe that 
\beaa
\SS_2 \psi + \frac{1}{r^2+a^2} \SS_3\psi&=& a \nab_T \nab_Z \psi + \frac{a^2}{r^2+a^2} \nab_Z \nab_Z \psi= a \nab_\That \nab_Z \psi
\eeaa
so we can write
\bea\label{eq:Psiz-Explicit}
\Psi_z&=& f \big(\de_0  \SS_1\psi+ (1+O(r^{-2} \de_0)) \OO\psi\big)+ \frac{4ar}{(r^2+a^2)^2}  \nab_\That \nab_Z \psi   \big(1+O(r^{-2} \de_0) \big).
\eea

By defining\footnote{The relevance of this particular choice will become evident in Section \ref{subsection-Poincare-higher}, where it will be used to derive a Poincar\'e inequality.}
\bea
\lab{eq:remark:choiceofLL-weak}
 \LL^{1} =  \de_0, \qquad    \LL^{2} =0, \qquad \LL^3=\LL^4=1,
\eea
 from \eqref{eq:operato-LL2}\lab{eq:operato-LL} and \eqref{eq:widetilde-P}, we obtain
\bea\label{eq:expression-widetilde-Psiz}
\widetilde{P} = \frac 1 2 h  L^{\a\b} \Db_\a \Psi_z\c \Db_\b  \Psi_z=\frac 1 2  h \Big(\de_0 \big|   \nab_T  \Psi_z \big|^2 + a^2\big|  \nab_Z \Psi_z \big|^2 + O^{\a\b} \Db_\a \Psi_z \Db_\b   \Psi_z\Big).
\eea

We can also simplify the term $P_{lot}$. From \eqref{eq:Plot} and \eqref{asymptotics-RRtp-new}, we can write
 \beaa
 P_{lot}&=&-\frac 1 2  h\Psi\c ( \RRtp^{\cund}  \LL^{\bund}[\SS_{\cund}, \SS_{\bund}]\psi)\\
 &=&-\frac 1 2  h\Psi\c \Big( \big(\RRtp^{1}  \LL^{4}- \RRtp^{4}  \LL^{1}\big)[\SS_{1}, \OO]\psi+\RRtp^{2}  \LL^{4}[\SS_{2}, \OO]\psi+ \big( \RRtp^{3}  \LL^{4}-  \RRtp^{4}  \LL^{3}\big) [\SS_{3}, \OO]\psi\Big)\\
  &=&-\frac 1 2  h\Psi\c \Big(  O(\de^2_0 r^{-5}) [\SS_{1}, \OO]\psi+O(r^{-3}) [\SS_{2}, \OO]\psi+ O(r^{-3}) [\SS_{3}, \OO]\psi\Big)\\
    &=&-\frac 1 2  h\Psi\c \Big( O(a r^{-2} )\nab \nab_T \psi + O( a^2 r^{-2})\nab\nab_Z \psi + O(a r^{-3})\dk^{\leq 1} \psi \Big),
 \eeaa
 which, since $h=O(r^5)$ and $\Psi=O(r^{-3})\dk^{\leq 2}\psi$,  can be summarized as 
 \bea\label{eq:expression-Plot}
  P_{lot}    &=&  O(a r^{-1})(\dk^{\leq 2}\psi)^2.
 \eea

 
 \subsection{Integration by parts   identities}
 \lab{sec:integrationbypartsidentitiesSS}
 

In the next section, we compute the quadratic forms $I$, $J$, $K$ as given in \eqref{eq:expr-I}, \eqref{eq:expr-J}, \eqref{eq:expr-K}. In order to do that, we will make use of some integration by parts identities, which we collect in Lemma \ref{Lemma:integrationbypartsSS_3SS_4}  below. 
 Such computations are important to observe that the mixed products between the symmetry operators $\SS_1$  and $\OO$ contain positive definite norms, modulo lower order or boundary terms.
 
We start with the following definition that will be useful to take care of boundary terms in the integrations by parts.  
\begin{definition}\label{def:boundarytermsinintegrationbypartsareok}
We denote by $M(\psi)$   quadratic terms of the following type
\beaa
\nab_T \psi \c \SS_{\aund}\psi, \qquad |q|^2\nab  \psi \c \nab \nab_T\psi, \qquad |q|^2\nab\psi\c\nab\nab_Z\psi,
\eeaa
Also, we denote by $M(\nab_r\psi)$ denote quadratic terms of the following type
\beaa
\nab_T \nab_r\psi \c \SS_{\aund}\nab_r\psi, \qquad |q|^2\nab  \nab_r\psi \c \nab \nab_T\nab_r\psi, \qquad |q|^2\nab  \nab_r\psi \c \nab \nab_Z\nab_r\psi.
\eeaa
\end{definition}

\begin{remark}
Notice in view of Definition \ref{def:boundarytermsinintegrationbypartsareok} the following pointwise  estimates for $M(\psi)$ and $M(\nab_r\psi)$:  
\bea
\bsplit
|M(\psi)| &\les |(\nab_T, \dkb)^{\leq 1}\psi||(\nab_T, \dkb)^{\leq 2}\psi|, \\
|M(\nab_r\psi)| &\les |\nab_r^{\leq 1}(\nab_T, \dkb)^{\leq 1}\psi||\nab_r^{\leq 1}(\nab_T, \dkb)^{\leq 2}\psi|.
\end{split}
\eea
\end{remark}

\begin{lemma}\label{Lemma:integrationbypartsSS_3SS_4} 
For any function $H=H(r)$, the following identities hold true:
\beaa
H \OO(\psi) \c \SS_1 \psi  &=&H|q|^2| \nab \nab_T\psi|^2  -O(ar^{-2})H(\dk^{\leq 2}\psi)^2 +|q|^2 \Ddot_\b(H|q|^{-2}O^{\a\b}\Ddot_\a \psi \c \SS_1 \psi )\\
&& +\partial_t (HM(\psi)), \\
H \nab_r \OO\psi \c \nab_r \SS_1\psi &=& H|q|^2|\nab \nab_T\nab_r\psi|^2 -O(a r^{-2})H(\nab_r\dk^{\leq 2}\psi)^2 -O(a r^{-2})H(\dk^{\leq 2}\psi)^2\\
&&+|q|^2 \Ddot_\b(H|q|^{-2}O^{\a\b}\Ddot_\a \nab_r\psi \c \SS_1 \nab_r\psi ) +\partial_t (HM(\nab_r\psi)), \\
  H \nab_r\OO(\psi) \c \SS_1 \psi &=&H |q|^2 \nab\nab_T\nab_r \psi \c \nab \nab_T \psi -O(a r^{-2})H(\nab_r\dk^{\leq 2}\psi)^2 -O(a r^{-2})H(\dk^{\leq 2}\psi)^2\\
  &&+|q|^2 \Ddot_\b(H|q|^{-2}O^{\a\b}\Ddot_\a \nab_r\psi \c \SS_1 \psi ) +\partial_t (HM(\nab_r\psi)),\\
  H\OO(\psi) \c \nab_r \SS_1\psi &=&H |q|^2 \nab\nab_T\nab_r \psi \c \nab \nab_T \psi -O(a r^{-2})H(\nab_r\dk^{\leq 2}\psi)^2 -O(a r^{-2})H(\dk^{\leq 2}\psi)^2\\
  &&+|q|^2 \Ddot_\b(H|q|^{-2}O^{\a\b}\Ddot_\a\psi \c \SS_1\nab_r \psi ) +\partial_t (HM(\nab_r\psi)).
 \eeaa
  In all the above,  $M(\psi)$ and $M(\nab_r\psi)$ denote the quadratic expressions in $\psi$ and its derivatives of Definition \ref{def:boundarytermsinintegrationbypartsareok}.
 \end{lemma}
 
 \begin{proof} 
Using \eqref{eq:OO-psi-Ddot-S_4}, i.e. $\OO(\psi)= |q|^2\Ddot_\b(|q|^{-2}O^{\a\b}\Ddot_\a \psi)$, we obtain for any $\Phi$,
\bea\label{eq:intermediate-integration-by-parts-lemma}
\nn |q|^{-2} \OO(\psi) \c \Phi &=& \Ddot_\b(|q|^{-2}O^{\a\b}\Ddot_\a \psi) \c \Phi\\ 
&=& \Ddot_\b(|q|^{-2}O^{\a\b}\Ddot_\a \psi \c \Phi )-|q|^{-2}O^{\a\b}\Ddot_\a \psi \c \Ddot_\b\Phi.
\eea
Applying the above to $\Phi=\SS_1 \psi=\nab_T\nab_T \psi$, and using that $[\nab, \nab_T]\psi=O(ar^{-4}) \dkb^{\leq 1}\psi$ we obtain
\beaa
 H\OO(\psi) \c \SS_1 \psi&=&-H|q|^{2}\nab  \psi \c \nab\nab_T \nab_T\psi+|q|^2 \Ddot_\b(H|q|^{-2}O^{\a\b}\Ddot_\a \psi \c \SS_1 \psi )\\
 &=&-H|q|^2\nab  \psi \c \nab_T\nab \nab_T\psi  -O(ar^{-2})H(\dk^{\leq 2}\psi)^2 \\
 &&+|q|^2 \Ddot_\b(H|q|^{-2}O^{\a\b}\Ddot_\a \psi \c \SS_1 \psi )\\
  &=&|q|^2H|\nab \nab_T\psi|^2  -O(ar^{-2})H(\dk^{\leq 2}\psi)^2+|q|^2 \Ddot_\b(H|q|^{-2}O^{\a\b}\Ddot_\a \psi \c \SS_1 \psi )\\
  && -\partial_t \big(H|q|^2\nab  \psi \c \nab \nab_T\psi \big) \\
    &=&H|q|^2|\nab \nab_T\psi|^2  -O(ar^{-2})H(\dk^{\leq 2}\psi)^2+|q|^2 \Ddot_\b(H|q|^{-2}O^{\a\b}\Ddot_\a \psi \c \SS_1 \psi )\\
    && +\partial_t (HM(\psi))
\eeaa
as stated.

Using $[\nab_r, \OO]\psi= O(a r^{-2}) \dkb^{\leq1}\psi$, we obtain
\beaa
&& H \nab_r \OO\psi \c \nab_r \SS_1\psi\\ 
&=& H \OO\nab_r \psi \c  \SS_1\nab_r\psi +HO(a r^{-2}) \dkb^{\leq1}\psi\c \nab_r \SS_1\psi+H \OO\nab_r \psi \c O(a r^{-4}) \dk^{\leq1}\psi\\
&=& H|q|^2|\nab \nab_T\nab_r\psi|^2-O(ar^{-2})H\nab  \nab_r\psi \c\dkb^{\leq 1}\nab_T\nab_r\psi+|q|^2 \Ddot_\b(H|q|^{-2}O^{\a\b}\Ddot_\a \nab_r\psi \c \SS_1 \nab_r\psi )\\
&&-\partial_t (|q|^2H\nab  \nab_r\psi \c \nab \nab_T\nab_r\psi) +HO(a r^{-2}) \dkb^{\leq1}\psi\c \nab_r \SS_1\psi +H \OO\nab_r \psi \c O(a r^{-4}) \dk^{\leq1}\psi
\eeaa
which can be written as 
\beaa
H \nab_r \OO\psi \c \nab_r \SS_1\psi &=& H|q|^2|\nab \nab_T\nab_r\psi|^2+|q|^2 \Ddot_\b(H|q|^{-2}O^{\a\b}\Ddot_\a \nab_r\psi \c \SS_1 \nab_r\psi )\\
&&+\partial_t (HM(\nab_r\psi)) -O(ar^{-2})H(\nab_r\dk^{\leq 2}\psi)^2 -O(a r^{-2})H(\dk^{\leq 2}\psi)^2,
\eeaa
as stated. 

The mixed products $H  \nab_r\OO\psi \c  \SS_1\psi$ and $H \OO\psi \c  \nab_r\SS_1\psi$ are treated similarly. 
This concludes the proof of Lemma \ref{Lemma:integrationbypartsSS_3SS_4}.
  \end{proof}

We use the above lemma to derive the following general computation. 

\begin{lemma}\label{general-computations-for-BB} 
Let $\Phi_\aund=\SS_\aund \Phi$ for some\footnote{We will apply it to $\Phi=\psi$ and $\Phi=\nab_r\psi$.} $\Phi$, and let $\mathcal{Y}^\aund$ be some coefficients only depending on $r$, such that
\beaa
\mathcal{Y}^{1} =\de_0 \mathcal{Y}, \qquad \mathcal{Y}^{4}=\mathcal{Y}.
\eeaa
 Then for $\LL^{\aund}$ given by \eqref{eq:remark:choiceofLL-weak}, i.e.
\beaa
 \LL^{1} =  \de_0, \qquad    \LL^{2} =0, \qquad \LL^3=\LL^4=1,
\eeaa
 we have
 \beaa
\big(\mathcal{Y}^{\aund}\Phi_\aund\big)\c \big( \LL^{\aund}\Phi_\aund\big)&=& \mathcal{Y}\Big(\de_0^2 |\SS_1\Phi|^2+ |\OO\Phi|^2 +2\de_0|q|^2| \nab \nab_T\Phi|^2\Big)\\
&& -O(a)(|\mathcal{Y}|+|\mathcal{Y}^2|+|\mathcal{Y}^3|)(\dk^{\leq 2}\psi)^2+\text{Bdr}[\Phi],
\eeaa
where the boundary term is given by
\beaa
\text{Bdr}[\Phi]&=&\partial_t \Big(\de_0\mathcal{Y}  r^2 M(\Phi)\Big)+|q|^2 \Ddot_\b\Big(2\de_0 |q|^{-2}O^{\a\b}\Ddot_\a \Phi \c \mathcal{Y} \SS_1 \Phi \Big).
\eeaa
\end{lemma}

\begin{proof}
By writing
\beaa
 \big(\mathcal{Y}^{\aund}\Phi_\aund\big)&=&\de_0 \mathcal{Y}( \SS_1\Phi)+ \mathcal{Y}( \OO\Phi)+ \mathcal{Y}^2 (\SS_2\Phi)+ \mathcal{Y}^3 (\SS_3\Phi ) \\
 \big( \LL^{\aund}\Phi_\aund\big)&=& \de_0 \SS_1\Phi+ \OO\Phi +  \SS_3\Phi 
\eeaa
we obtain
\beaa
\big(\mathcal{Y}^{\aund}\Phi_\aund\big)\c \big( \LL^{\aund}\Phi_\aund\big)&=&\de_0^2 \mathcal{Y} |\SS_1\Phi|^2+\mathcal{Y} |\OO\Phi|^2+\mathcal{Y}^3 |\SS_3\Phi|^2+(\mathcal{Y}+\mathcal{Y}^3) \de_0\SS_1\Phi \c \SS_3\Phi \\
&&+2\de_0 \mathcal{Y}(\SS_1\Phi \c \OO \Phi )+(\mathcal{Y}+\mathcal{Y}^3) \SS_3 \Phi \c \OO \Phi \\
&& +( \mathcal{Y}^2 \SS_2\Phi )\c ( \de_0 \SS_1\Phi+ \OO\Phi +  \SS_3\Phi ).
\eeaa
Since $\SS_1\psi, \OO\psi=O(1)\dk^{\leq 2}\psi$, $\SS_2\psi=O(a)\dk^{\leq 2}\psi$ and $\SS_3\psi=O(a^2)\dk^{\leq 2}\psi$, we infer
\beaa
\big(\mathcal{Y}^{\aund}\Phi_\aund\big)\c \big( \LL^{\aund}\Phi_\aund\big)&=&\de_0^2 \mathcal{Y} |\SS_1\Phi|^2+\mathcal{Y} |\OO\Phi|^2 +2\de_0 \mathcal{Y}(\SS_1\Phi \c \OO \Phi )\\
&& +O(a)(|\mathcal{Y}|+|\mathcal{Y}^2|+|\mathcal{Y}^3|)(\dk^{\leq 2}\psi)^2.
\eeaa
Using Lemma \ref{Lemma:integrationbypartsSS_3SS_4}, we obtain the desired identity. 
\end{proof}


   \subsection{The quadratic  forms  $I$, $J$ and $K$}
   \lab{sec:thequadraticformsIJKarenowcomputed}
  
   
The goal of this section is to compute the quadratic forms $I$, $J$ and $K$ appearing in the computation \eqref{Morawetz-effective-principal-term} of the effective generalized current $\widetilde{\EE}[\bold{X}, \bold{w}, \bold{M}]$ by making use of the integration by parts identities of section  \ref{sec:integrationbypartsidentitiesSS}.  Recall that $I$, $J$ and $K$ are given by \eqref{eq:expr-I}, \eqref{eq:expr-J}, \eqref{eq:expr-K} and that we make use of the choice of $z$ and $h$ in Lemma \ref{lemma:mainchoicesofzandhandcorrecpondingquantities:faposdifah}, i.e.
\beaa
z =z_0- \de_0 z_0^2, \qquad z_0= \frac{\De}{(r^2+a^2)^2}, \qquad \de_0>0, \qquad h=\frac{(r^2+a^2)^4}{r(r^2-a^2)}.
\eeaa


   \subsubsection{The quadratic form  $I$}
  
   
Recall that the quadratic form $I$ is given by \eqref{eq:expr-I}, i.e.
 \beaa
I&=& \big(\AAa[z]\nab_r \psia\big)\c \big( \LL^{\aund}\nab_r \psia\big),
\eeaa   
where  $\LL^{\aund}$ is given by \eqref{eq:remark:choiceofLL-weak}, and where $\AAa[z]$ is given by \eqref{eq:explicit-AA-}, which may be rewritten under the form \eqref{simplifications-AA},i.e.
\beaa
\bsplit
\AA^1&= \de_0 \AA\big(1+O(r^{-1} \de_0)\big), \quad \,\AA^4=\AA\big(1+O(r^{-1}\de_0)\big), \\
 \AA^2&=\AAt \big(2a^2+O(\de_0 )\big), \qquad\,\,\,\,   \AA^3=  \tilde{\AA}\big(1+O(r^{-2}\de_0 )\big),
 \end{split}
\eeaa
where 
\beaa
\AA&=&   \frac{2\De^2}{r^2(r^2+a^2)(r^2-a^2)^2}(3mr^4-4a^2r^3+ma^4), \qquad \widetilde{\AA}= \frac{8\De^2r}{(r^2+a^2)(r^2-a^2)^2},
\eeaa
are positive coefficients for $|a|/m <1$.  We may thus apply Lemma \ref{general-computations-for-BB}  to $I$ with the choice $\mathcal{Y}^\aund=\AA^\aund$. Using the structure or the commutators $[\nab_r, \SS_{\aund}]$, we obtain
 \bea\lab{eq:computationofthequadraticformIforchap8:fspoighs}
 \bsplit
 I =& {\AA}(1+O(r^{-1} \de_0))\Big(\de_0^2 |\nab_r\SS_1\psi|^2+ |\nab_r\OO\psi|^2 +2\de_0|q|^2| \nab \nab_T\nab_r\psi|^2\Big)\\
& -O(a)(|\AA|+|\AAt|)  (1+O(r^{-1} \de_0))\Big((\nab_r\dk^{\leq 2}\psi)^2+r^{-2}(\dk^{\leq 2}\psi)^2\Big) +\text{Bdr}[\psi]_I,
\end{split}
\eea
where the boundary term is given by
\beaa
\text{Bdr}[\psi]_I&=&\partial_t \Big(\de_0\AA M(\nab_r\psi)\Big) +|q|^2 \Ddot_\b\Big(2\de_0|q|^{-2}O^{\a\b}\Ddot_\a \nab_r \psi  \c {\AA} \SS_1\nab_r\psi \Big).
\eeaa
   
 \begin{remark}
  Observe that the  estimate  for the quadratic form $I$ for the choice $z=z_0$ can be deduced by  plugging  $\de_0=0$ in \eqref{eq:computationofthequadraticformIforchap8:fspoighs}.  In that case one can see that one does not control the term $|\nab_r\SS_1\psi|^2$. This is the main reason one is led to use the choice $z=z_0-\de_0z_0^2$.  
\end{remark}

     
\subsubsection{The quadratic forms $J$ and $K$}


As a corollary of Lemma \ref{general-computations-for-BB}, we obtain the expressions for $J$ and $K$ in Proposition \ref{proposition:expression-for-Psi}, i.e. 
\beaa
J= \big(\VVa[z]  \psia\big)\c \big( \LL^{\aund} \psia\big), \qquad K=\frac 1 4 |q|^2 \D^\mu\Big((M^{\aund}_\mu\psi_\aund )\c (\LL^\bund   \psi_\bund) \Big).
\eeaa

\begin{lemma}\lab{lemma:computationofthequadraticformJandKforchap8:fspoighs} 
The term $J$ is given by 
 \bea
 \begin{split}
J =& \big( \VV +O( \de_0   r^{-3}) \big)\Big(\de_0^2 |\SS_1\psi|^2+ |\OO\psi|^2 +2\de_0|q|^2| \nab \nab_T\psi|^2\Big)  - O(ar^{-1})(\dk^{\leq 2}\psi)^2\\
&+\text{Bdr}[\psi]_J,
\end{split}
\eea
 with boundary term
\beaa
\text{Bdr}_J[\psi] =\partial_t \Big(\de_0(\VV +O(   r^{-3}))  r^2   M(\psi)\Big) +|q|^2 \Ddot_\b\Big(2\de_0 
|q|^{-2}O^{\a\b}\Ddot_\a \psi \c (\VV +O( a^2   r^{-3})) \SS_1\psi \Big).
\eeaa
For a one-form of the type  $M^\aund :=v^{\aund} \pr_ r$, with $v^2=v^3=0$, $v^1=\de_0 v$ and $v^4=v$ for some given function $v=v(r)$, the term $K$ is given by
\beaa
 K&=&\frac{|q|^2}{2} v \Big( \de_0^2 \nab_r \SS_1\psi \c \SS_1\psi+ \nab_r \OO\psi \c \OO\psi+2\de_0   |q|^2 \nab\nab_T\nab_r \psi \c \nab \nab_T \psi\Big)\\
&&+\frac{|q|^2}{4}v'\Big( \de_0^2  |\SS_1\psi|^2+ |\OO\psi|^2 +2\de_0 |q|^2| \nab \nab_T\psi|^2\Big) \\
&&- vO(a r^{\frac{5}{2}})(\nab_r\dk^{\leq 2}\psi)^2 - O(a r^{\frac{3}{2}}) v(\dk^{\leq 2}\psi)^2 - O(a r^2)v' (\dk^{\leq 2}\psi)^2+\text{Bdr}[\psi]_K, 
\eeaa
where we denoted
\beaa
v'^{\aund}:=\pr_r v^\aund+ \frac{2r}{|q|^2} v^\aund,
\eeaa 
and with boundary term 
\beaa
\text{Bdr}[\psi]_K &=& \partial_t \big(vr^4M(\nab_r\psi)\big) +\frac{|q|^4}{4} \Ddot_\b(2\de_0 v|q|^{-2}O^{\a\b}\Ddot_\a\psi \c \SS_1\nab_r \psi )\\
&& +\frac{|q|^4}{4} \Ddot_\b(2\de_0 v|q|^{-2}O^{\a\b}\Ddot_\a\nab_r\psi \c \SS_1\psi ) +\partial_t \Big(\de_0v' r^4 M(\psi)\Big)  \\
&&+\frac{|q|^4}{4}\Ddot_\b\Big(2\de_0 v' |q|^{-2}O^{\a\b}\Ddot_\a \psi \c\SS_1\psi \Big).
\eeaa
\end{lemma}

  \begin{proof} 
  The expressions for $J$ is a straightforward application of Lemma \ref{general-computations-for-BB}, using that, see \eqref{expressions-VV-asymp},
  \beaa
  \VV^1&=& \de_0\VV, \qquad \VV^4=  \VV +O( \de_0   r^{-3}), \qquad \VV^2=O( r^{-1}), \qquad  \quad \VV^3= O( r^{-3}).
  \eeaa
    
  To write the term $K$, we first compute
\beaa
\frac{4}{|q|^2} K&=&\D^\mu\Big((M^{\aund}_\mu\psi_\aund )\c (\LL^\bund   \psi_\bund) \Big)\\
&=&(\D^\mu (M^{\aund}_\mu\psi_\aund ))\c (\LL^\bund   \psi_\bund) +(M^{\aund}_\mu\psi_\aund )\c \D^\mu(\LL^\bund   \psi_\bund)  \\
&=&(\D^\mu M^{\aund}_\mu) \psi_\aund \c (\LL^\bund   \psi_\bund) + (M^{\aund}_\mu \D^\mu\psi_\aund )\c (\LL^\bund   \psi_\bund) +(M^{\aund}_\mu\psi_\aund )\c \D^\mu(\LL^\bund   \psi_\bund).
\eeaa
For $M^\aund =v^{\aund} \pr_ r$, we obtain
\beaa
\frac{4}{|q|^2} K&=&\left(\pr_r v^\aund+ \frac{2r}{|q|^2} v^\aund\right) \psi_\aund \c (\LL^\bund   \psi_\bund) + v^{\aund} \nab_r \psi_\aund \c (\LL^\bund   \psi_\bund) +(v^{\aund}\psi_\aund )\c \nab_r (\LL^\bund   \psi_\bund)  \\
&=&\left(\pr_r v^\aund+ \frac{2r}{|q|^2} v^\aund\right) \psi_\aund \c (\LL^\bund   \psi_\bund) + v^{\aund} \nab_r \psi_\aund \c (\LL^\bund   \psi_\bund) +(v^{\aund}\psi_\aund )\c  \LL^\bund   \nab_r\psi_\bund\\
&=& (v'^{\aund} \psi_\aund) \c (\LL^\bund   \psi_\bund) + v^{\aund} \nab_r \psi_\aund \c (\LL^\bund   \psi_\bund) +(v^{\aund}\psi_\aund )\c  \LL^\bund   \nab_r\psi_\bund
\eeaa
where we wrote $v'^{\aund}:=\pr_r v^\aund+ \frac{2r}{|q|^2} v^\aund$. By defining for some $v$
\beaa
v^1=\de_0 v, \qquad v'^1=\de_0 v', \qquad v^4=v, \qquad v'^4=v',\qquad v^2=v^3=0,
\eeaa
we can apply Lemma \ref{general-computations-for-BB} to the first term, and obtain
 \beaa
\big(v'^{\aund}\psi_\aund\big)\c \big( \LL^{\aund}\psi_\aund\big)&=&v'\big( \de_0^2  |\SS_1\psi|^2+ |\OO\psi|^2 +2\de_0 |q|^2| \nab \nab_T\psi|^2\big) - O(a)v'(\dk^{\leq 2}\psi)^2+\text{Bdr}[\psi],
\eeaa
where the boundary term is given by
\beaa
\text{Bdr}[\psi]&=&\partial_t \Big(\de_0v' M(\psi)\Big)  +|q|^2 \Ddot_\b\Big(2\de_0|q|^{-2}O^{\a\b}\Ddot_\a \psi\c v' \SS_1\psi \Big).
\eeaa

For the other two terms we write, using Lemma \ref{Lemma:integrationbypartsSS_3SS_4},
\beaa
&&v^{\aund} \nab_r \psi_\aund \c (\LL^\bund   \psi_\bund) +(v^{\aund}\psi_\aund )\c  \LL^\bund   \nab_r\psi_\bund \\
&=&\big( \de_0 v \nab_r\SS_1\psi+  v \nab_r\OO\psi\big) \c ( \de_0 \SS_1\psi+ \SS_3\psi+\OO\psi)\\
&& +\big( \de_0 v \SS_1\psi+  v\OO\psi\big) \c ( \de_0 \nab_r\SS_1\psi+ \nab_r\SS_3\psi+\nab_r\OO\psi)\\
&=& 2\de_0^2v \nab_r \SS_1\psi \c \SS_1\psi+2v \nab_r \OO\psi \c \OO\psi+2\de_0 v \big( \nab_r \SS_1\psi \c \OO \psi+ \SS_1\psi \c \nab_r \OO\psi\big)\\
&&+O(a)v|\nab_r\dk^{\leq 2}\psi||\dk^{\leq 2}\psi|\\
&=& 2\de_0^2v \nab_r \SS_1\psi \c \SS_1\psi+2v \nab_r \OO\psi \c \OO\psi+4\de_0 v  |q|^2 \nab\nab_T\nab_r \psi \c \nab \nab_T \psi\\
&&+2\de_0 v \big(    \nab_r \SS_2 \psi \c  \SS_2 \psi\big)+ 2v |q|^2a^2 \nab\nab_Z\nab_r \psi \c \nab \nab_Z \psi\\
&&- vO(ar^{\frac{1}{2}})(\nab_r\dk^{\leq 2}\psi)^2 - vO(ar^{-\frac{1}{2}})(\dk^{\leq 2}\psi)^2 +\partial_t \big(vM(\nab_r\psi)\big)\\
&&+|q|^2 \Ddot_\b(2\de_0 v|q|^{-2}O^{\a\b}\Ddot_\a\psi \c\SS_1\nab_r \psi ) +|q|^2 \Ddot_\b( 2\de_0v|q|^{-2}O^{\a\b}\Ddot_\a\nab_r\psi \c\SS_1\psi ).
\eeaa
Summing the two above, we obtain the lemma. 
\end{proof}


   \subsection{Conclusion}
  

We summarize the results obtained so far in the following proposition.

\begin{proposition}
\label{proposition-with-quadratic-form}
The  effective generalized current is given by
 \bea\label{generalized-current-operator-2}
 |q|^2 \widetilde{\EE}[\bold{X}, \bold{w}, \bold{M}] &=&\widetilde{P}+\Qr_{\SS_1, \OO, \nab \nab_T}+   \EE_{lot}+\text{Bdr}
   \eea
   with the following terms.
\begin{enumerate}
\item The principal trapping term $\widetilde{P}$ is given by, see \eqref{eq:expression-widetilde-Psiz},
\bea
\widetilde{P}&=&\frac 1 2  h \Big( \de_0 \big|   \nab_T  \Psi_z \big|^2+a^2|\nab_Z \Psi_z|^2 + O^{\a\b} \Db_\a \Psi_z \Db_\b   \Psi_z\Big),
\eea
where
\bea\label{definition-Psiz}
\begin{split}
\Psi_z&= -\frac{2\TT}{(r^2+a^2)^3}   \big(\de_0  \SS_1\psi+ (1+O(r^{-2} \de_0)) \OO\psi\big)\\
&+ \frac{4ar}{(r^2+a^2)^2}  \nab_\That \nab_Z \psi   \big(1+O(r^{-2} \de_0) \big).
\end{split}
\eea

\item The quadratic   form $\Qr_{\SS_1, \OO, \nab \nab_T}$ is given by
   \beaa
   \Qr_{\SS_1, \OO, \nab \nab_T}&:=& {\AA}(1+O(r^{-1} \de_0))\Big(\de_0^2 |\nab_r\SS_1\psi|^2+ |\nab_r\OO\psi|^2 +2\de_0|q|^2| \nab \nab_T\nab_r\psi|^2\Big)\\
   && +\frac{|q|^2}{2} v \Big( \de_0^2 \nab_r \SS_1\psi \c \SS_1\psi+ \nab_r \OO\psi \c \OO\psi+2\de_0   |q|^2 \nab\nab_T\nab_r \psi \c \nab \nab_T \psi\Big)\\
 &&  +\left( \VV+\frac{|q|^2}{4}v' +O( \de_0   r^{-3}) \right)\Big(\de_0^2 |\SS_1\psi|^2+ |\OO\psi|^2 +2\de_0|q|^2| \nab \nab_T\psi|^2\Big).
   \eeaa

   \item The terms $\EE_{lot}$ are lower order terms in $a$, given by
   \bea\label{definition-E-lot}
   \begin{split}
   \EE_{lot}&= -O(a)({\AA}+{\AAt})  (1+O(r^{-1} \de_0))((\nab_r\dk^{\leq 2}\psi)^2+r^{-2}(\dk^{\leq 2}\psi)^2)\\
   &  -O(ar^{-1})(\dk^{\leq 2}\psi)^2 - vO(a r^{\frac{5}{2}})(\nab_r\dk^{\leq 2}\psi)^2 - O(a r^{\frac{3}{2}}) v(\dk^{\leq 2}\psi)^2\\
   & - O(a r^2)v' (\dk^{\leq 2}\psi)^2.
   \end{split}
   \eea
   In particular, since $\AA, \AAt=O(\De^2 r^{-4})$, and  for $v=O(m^{1/2}\Delta r^{-9/2})$, we can bound the above  as
   \bea\label{eq:bound-EE-lot}
      \EE_{lot} \geq -O(a)\big( |\nab_{\Rhat}\dk^{\leq 2}\psi|^2+r^{-1} |\dk^{\leq 2}\psi|^2\big).
   \eea
   
\item The boundary terms are given by
   \beaa
   \text{Bdr}&=& \text{Bdr}[\psi]_I+\text{Bdr}[\psi]_J+\text{Bdr}[\psi]_K\\
   &=& \partial_t \Big(  M(\nab_{\Rhat}\psi)+M(\psi)\Big)    +|q|^2 \D_\b  \widehat{\BB}^\b
   \eeaa
   with 
   \beaa
   \widehat{\BB}^\b&:=& 2\de_0 |q|^{-2}O^{\a\b}\Ddot_\a \nab_r \psi  \c {\AA} \SS_1\nab_r\psi +2\de_0|q|^{-2}O^{\a\b}\Ddot_\a \psi \c(\VV +O( a^2   r^{-3})) \SS_1\psi\\
   &&+\frac{1}{4}\de_0 vO^{\a\b}\Ddot_\a\psi \c\SS_1\nab_r \psi +\frac{1}{4}\de_0 vO^{\a\b}\Ddot_\a\nab_r\psi \c \SS_1\psi +\frac{1}{2}\de_0 O^{\a\b}\Ddot_\a \psi \c v' \SS_1\psi,
   \eeaa
  where $M(\psi)$ denotes the  quadratic expressions in $\psi$ and  $M(\nab_{\Rhat}\psi)$  denotes the  quadratic expressions in $\psi$ and its derivatives of Definition \ref{def:boundarytermsinintegrationbypartsareok}.
\end{enumerate}
\end{proposition}
   
\begin{proof}
Recall that according to Proposition \ref{proposition:Morawetz4}, the effective generalized current $\widetilde{\EE}[\bold{X}, \bold{w}, \bold{M}]$ introduced in \eqref{widetilde-EE}  is given by, see \eqref{Morawetz-effective-principal-term},
 \beaa
 |q|^2 \widetilde{\EE}[\bold{X}, \bold{w}, \bold{M}] &=&\widetilde{P}+P_{lot}+I  +J +K.
   \eeaa 
The proof then follows immediately from the control for $\widetilde{P}$ in \eqref{eq:expression-widetilde-Psiz}, for 
$P_{lot}$ in \eqref{eq:expression-Plot}, for $I$ in \eqref{eq:computationofthequadraticformIforchap8:fspoighs}, and for $J$ and $K$ in Lemma \ref{lemma:computationofthequadraticformJandKforchap8:fspoighs}. 
\end{proof}


\section{Proof of $\SS$-derivative Morawetz estimate}\label{section:unconditional-SS-mOr}


 In this section  we prove the $\SS$-valued Morawetz estimates, i.e. we prove Proposition \ref{prop:morawetz-higher-order}.  Recall that we make use of the choices for $z$ and $h$ made in Lemma \ref{lemma:mainchoicesofzandhandcorrecpondingquantities:faposdifah}, i.e.
\beaa
z =z_0- \de_0 z_0^2, \qquad z_0= \frac{\De}{(r^2+a^2)^2}, \qquad \de_0>0, \qquad h=\frac{(r^2+a^2)^4}{r(r^2-a^2)}.
\eeaa


\subsection{Control of the quadratic form $ \Qr_{\SS_1, \OO, \nab \nab_T}$}
\label{subsection-Poincare-higher}


In this section, we prove that  the quadratic form $ \Qr_{\SS_1, \OO, \nab \nab_T}$ appearing in Proposition \ref{proposition-with-quadratic-form}  is positive definite. In order to do so, as in Section \ref{section:proof-hardy-poincare}, we will apply Poincar\'e and Hardy inequalities.

Recalling that 
\beaa
\frac{2}{h} \widetilde{P}&=&\de_0 \big|   \nab_T  \Psi_z \big|^2+a^2|\nab_Z\Psi_z|^2 + O^{\a\b} \Db_\a \Psi_z \Db_\b   \Psi_z\\
&=& \de_0 \big|   \nab_T  \Psi_z \big|^2+a^2|\nab_Z\Psi_z|^2 + |q|^2|\nab\Psi_z|^2
\eeaa
we obtain, using the Poincar\'e inequality of Lemma \ref{lemma:poincareinequalityfornabonSasoidfh:chap6}, for $|a|\ll m$ and $|a|\ll \de_0m$, 
 \bea\label{eq:poincare-inequality-widetilde-P}
  \int_S\widetilde{P}   &\ge &  \int_S  h   |\Psi_z|^2 -O(ar^7)\int_S|\nab\Psi_z|^2.
 \eea
 Using the integration by parts identities as in Lemma \ref{Lemma:integrationbypartsSS_3SS_4}, we can write
 \beaa
 |\Psi_z|^2&=& \big( \RRtp^\aund[z] \psi_\aund)\c\big( \RRtp^\bund[z] \psi_\bund) \\
 &=& ( \RRtp^{1})^2| \SS_1\psi|^2+ ( \RRtp^{2})^2| \SS_2\psi|^2+(\RRtp^{3})^2| \SS_3\psi|^2+(\RRtp^{4} )^2|\OO\psi|^2\\
 &&+ 2 \RRtp^{1}\RRtp^{3} \SS_1 \psi \c \SS_3\psi + 2 \RRtp^{1}\RRtp^{4} \SS_1 \psi \c \OO\psi+ 2 \RRtp^{3}\RRtp^{4} \SS_3 \psi \c \OO\psi  \\
 &&+2\RRtp^{2} \SS_2\psi \c ( \RRtp^{1} \SS_1\psi+\RRtp^{3} \SS_3\psi+\RRtp^{4} \OO\psi)\\
 &=& ( \RRtp^{1})^2| \SS_1\psi|^2+(\RRtp^{4} )^2|\OO\psi|^2+ 2 \RRtp^{1}\RRtp^{4} \SS_1 \psi \c \OO\psi  \\
 &&+ O(a)\Big(( \RRtp^{1})^2+( \RRtp^{2})^2+( \RRtp^{3})^2+( \RRtp^{4})^2\Big)(\dk^{\leq 2}\psi)^2\\
  &=& ( \RRtp^{1})^2| \SS_1\psi|^2+(\RRtp^{4} )^2|\OO\psi|^2+ 2 \RRtp^{1}\RRtp^{4}|q|^2| \nab \nab_T\psi|^2\\
 &&+ O(a)\Big(( \RRtp^{1})^2+( \RRtp^{2})^2+( \RRtp^{3})^2+( \RRtp^{4})^2\Big)(\dk^{\leq 2}\psi)^2\\
&& +\partial_t \big(\RRtp^{1}\RRtp^{4}M(\psi)\big)  +|q|^2 \Ddot_\b(|q|^{-2}\RRtp^{1}\RRtp^{4}O^{\a\b}\Ddot_\a \psi \c\SS_1\psi ).
 \eeaa
 Writing from \eqref{asymptotics-RRtp-new},
 \beaa
 \RRtp^1=\de_0 f, \qquad \RRtp^4=f(1+O(\de_0r^{-2})), \qquad \RRtp^2=O(r^{-3}), \qquad \RRtp^3=O(r^{-5}),
 \eeaa
 we obtain 
  \beaa
 |\Psi_z|^2 &=& f^2\big(\de_0^2| \SS_1\psi|^2+|\OO\psi|^2+ 2 \de_0| \nab \nab_T\psi|^2\big) -O(ar^{-6})(\dk^{\leq 2}\psi)^2\\
 &&+\partial_t \big(r^{-6}M(\psi)\big)  +|q|^2 \Ddot_\b(O(r^{-6})|q|^{-2}O^{\a\b}\Ddot_\a \psi \c (\SS_1+\SS_3) \psi ).
 \eeaa

  By combining Proposition \ref{proposition-with-quadratic-form} and the above for the term $(1-\de)\widetilde{P}$ in $\widetilde{P}=\de\widetilde{P}+(1-\de) \widetilde{P}$, we obtain the following

\begin{proposition}\label{proposition-with-quadratic-form-poincare} 
Upon applying the Poincar\'e inequality to the principal term, the generalized current  verifies the identity, for  any sphere $S=S(t, r)$,
\beaa
\int_S |q|^2 \widetilde{\EE}[\bold{X}, \bold{w}, \bold{M}] &\geq& \int_S\Big(\de\widetilde{P}+\widetilde{\Qr}_{\SS_1, \OO, \nab \nab_T}+   \widetilde{\EE}_{lot}+\text{Bdr}\Big)
 \eeaa
 where the quadratic form $\widetilde{\Qr}_{\SS_1, \OO, \nab \nab_T}$ is given by
\beaa
\widetilde{\Qr}_{\SS_1, \OO, \nab \nab_T}&:=& \Qr_{\SS_1, \OO, \nab \nab_T}+(1-\de)h  f^2\big(\de_0^2| \SS_1\psi|^2+|\OO\psi|^2+ 2 \de_0| \nab \nab_T\psi|^2\big)\\
&=&  \AA(1+O(r^{-1} \de_0))\Big(\de_0^2 |\nab_r\SS_1\psi|^2+ |\nab_r\OO\psi|^2 +2\de_0|q|^2| \nab \nab_T\nab_r\psi|^2\Big)\\
   && +\frac{|q|^2}{2} v \big( \de_0^2 \nab_r \SS_1\psi \c \SS_1\psi+ \nab_r \OO\psi \c \OO\psi+2\de_0   |q|^2 \nab\nab_T\nab_r \psi \c \nab \nab_T \psi\big)\\
 &&  +\left( \VV+\frac{|q|^2}{4}v' +(1-\de)h  f^2+O( \de_0   r^{-3}) \right)\Big(\de_0^2 |\SS_1\psi|^2+ |\OO\psi|^2 +2\de_0|q|^2| \nab \nab_T\psi|^2\Big)
\eeaa
and the lower order term $\widetilde{\EE}_{lot}$ is given by
\beaa
\widetilde{\EE}_{lot}&:=&\EE_{lot}- O(ar^7)|\nab\Psi_z|^2 -O(ar^{-1})|\dk^{\leq 2}\psi|^2.
\eeaa
\end{proposition}

We can now apply the Hardy inequality where $v$ is given by Proposition \ref{Prop.StongversionMorawetz1}, 
to show that for $|a|/m \ll 1$ and a universal constant $c_1$ we have
\beaa
\widetilde{\Qr}_{\SS_1, \OO, \nab \nab_T} &\geq& c_1  \Big( m\big(\big|\nab_{\Rhat}\SS_1\psi|^2 + \big|\nab_{\Rhat}\OO\psi|^2+|q|^2\big|\nab\nab_T\nab_r\psi|^2\big)\\
&& + r^{-1}\big( |\SS_1\psi|^2+|\OO\psi|^2+|q|^2 |\nab\nab_T\psi|^2\big) \Big)\\
&& -O(ar^{-2})\big(|\nab_{\Rhat}\dk^{\leq 2}\psi|^2+|\dk^{\leq 2}\psi|^2\big).
\eeaa

Indeed, we can separate the quadratic form $\widetilde{\Qr}_{\SS_1, \OO, \nab \nab_T} $ in the terms involving $\SS_1$, $\OO$ and $\nab\nab_T$, and obtain by Lemma \ref{Le:Crucialpositivity-a=0} that
\beaa
\widetilde{\Qr}_{\SS_1}&:=& \de_0^2 \bigg(\AA \big| \nab_r \SS_1 \psi \big|^2 + \frac {|q|^2}{ 2}  v  \nab_r \SS_1\psi \c \SS_1\psi+\left(  \VV+ \frac {|q|^2}{ 4} v' +(1-\de) f^2h+O( \de_0  r^{-3})\right) |\SS_1\psi|^2\bigg)\\
&\geq& c_0  \de_0^2 \Big( \frac{m \De^2}{r^4} \big|\nab_r\SS_1\psi|^2 + r^{-1} |\SS_1\psi|^2 \Big).
\eeaa
Similarly, by applying again Lemma \ref{Le:Crucialpositivity-a=0} we have
\beaa
\widetilde{\Qr}_{\OO}&:=& \AA \big| \nab_r \OO \psi \big|^2 + \frac {|q|^2}{ 2}  v  \nab_r \OO\psi \c \OO\psi+\left(  \VV+ \frac {|q|^2}{ 4} v' +(1-\de) f^2h+O( \de_0   r^{-3})\right) |\OO\psi|^2\\
&\geq& c_0 \Big( \frac{m \De^2}{r^4} \big|\nab_r\OO\psi|^2 + r^{-1} |\OO\psi|^2 \Big)
\eeaa
and, recalling $[\nab_r, |q|\nab]=O(ar^{-2})\psi$ and $[\nab_r, \nab_T]=O(ma r^{-4})\dk^{\leq 1}$, 
\beaa
\widetilde{\Qr}_{\nab \nab_T}&:=& \de_0|q|^2 \Bigg(\AA \big|  \nab\nab_T\nab_r \psi \big|^2 + \frac {|q|^2}{ 2}  v  \nab \nab_T\nab_r\psi \c \nab \nab_T\psi\\
&&+\left(  \VV+ \frac {|q|^2}{ 4} v' +(1-\de) f^2h+O( \de_0   r^{-3})\right) |\nab\nab_T\psi|^2\Bigg)\\
&=& \de_0 \bigg(\AA \big| \nab_r( |q|\nab\nab_T \psi) \big|^2 + \frac {|q|^2}{ 2}  v  \nab_r (|q|\nab \nab_T\psi) \c (|q|\nab \nab_T\psi)\\
&&+\left(  \VV+ \frac {|q|^2}{ 4} v' +(1-\de) f^2h+O( (\de_0+a)   r^{-3})\right) ||q|\nab\nab_T\psi|^2\bigg)\\
&& -O(ar^{-2})\big(|\nab_{\Rhat}\dk^{\leq 2}\psi|^2+|\dk^{\leq 2}\psi|^2\big)\\
&\geq& c_0  \de_0|q|^2 \Big( \frac{m \De^2}{r^4}|q|^{-2} \big|\nab_r(|q|\nab\nab_T\psi)|^2 + r^{-1} |\nab\nab_T\psi|^2 \Big)\\
&& -O(ar^{-2})\big(|\nab_{\Rhat}\dk^{\leq 2}\psi|^2+|\dk^{\leq 2}\psi|^2\big)\\
&\geq& c_0  \de_0|q|^2 \Big( \frac{m \De^2}{r^4} \big|\nab\nab_T\nab_r\psi)|^2 + r^{-1} |\nab\nab_T\psi|^2 \Big) -O(ar^{-2})\big(|\nab_{\Rhat}\dk^{\leq 2}\psi|^2+|\dk^{\leq 2}\psi|^2\big).
\eeaa

This gives for the generalized current, for  any sphere $S=S(t, r)$,
\bea\label{eq:positivity-first-quadratic-form}
\begin{split}
\int_S |q|^2 \widetilde{\EE}[\bold{X}, \bold{w}, \bold{M}] \geq& \de\int_S\widetilde{P}\\
 &+c_1  \int_S\Big( m\big(\big|\nab_{\Rhat}\SS_1\psi|^2 + \big|\nab_{\Rhat}\OO\psi|^2+|q|^2\big|\nab\nab_T\nab_{\Rhat}\psi|^2\big)\\
 & + r^{-1}\big( |\SS_1\psi|^2+|\OO\psi|^2+|q|^2 |\nab\nab_T\psi|^2\big) \Big)\\
 &+ \int_S\Big(  \widetilde{\EE}_{lot} -O(ar^{-2})\big(|\nab_{\Rhat}\dk^{\leq 2}\psi|^2+|\dk^{\leq 2}\psi|^2\big) +\text{Bdr}\Big).
 \end{split}
\eea


\subsection{Control of the effective generalized current}
\lab{sec:finallowerboundforeffectivegeneralizedcurrentaosidhfauidh}


We have, recall \eqref{eq:bound-EE-lot}, 
\beaa
\widetilde{\EE}_{lot}&=&\EE_{lot}- O(ar^7)|\nab\Psi_z|^2 -O(ar^{-1})|\dk^{\leq 2}\psi|\\
&\geq& -O(a)\big( |\nab_{\Rhat}\dk^{\leq 2}\psi|^2+r^{-1} |\dk^{\leq 2}\psi|^2\big) - O(ar^7)|\nab\Psi_z|^2.
\eeaa
Plugging in \eqref{eq:positivity-first-quadratic-form}, we infer
\beaa
\begin{split}
\int_S |q|^2 \widetilde{\EE}[\bold{X}, \bold{w}, \bold{M}] \geq& \de\int_S\widetilde{P} +c_1  \int_S\Big( m\big(\big|\nab_{\Rhat}\SS_1\psi|^2 + \big|\nab_{\Rhat}\OO\psi|^2\big)+ r^{-1}\big( |\SS_1\psi|^2+|\OO\psi|^2\big) \Big)\\
 & -O(a)\int_S\big( |\nab_{\Rhat}\dk^{\leq 2}\psi|^2+r^{-1} |\dk^{\leq 2}\psi|^2\big) +\int_S\text{Bdr}.
 \end{split}
\eeaa
Since $\SS_2=O(a)\dk^{\leq 2}$ and $\SS_3=O(a^2)\dk^{\leq 2}$, we infer
\beaa
\begin{split}
\int_S |q|^2 \widetilde{\EE}[\bold{X}, \bold{w}, \bold{M}] \geq& \de\int_S\widetilde{P} +c_1  \int_S\Big( m\big(\big|\nab_{\Rhat}\SS_1\psi|^2+\big|\nab_{\Rhat}\SS_2\psi|^2+\big|\nab_{\Rhat}\SS_3\psi|^2 + \big|\nab_{\Rhat}\OO\psi|^2\big)\\
&+ r^{-1}\big( |\SS_1\psi|^2+|\SS_2\psi|^2+|\SS_3\psi|^2+|\OO\psi|^2\big) \Big)\\
 & -O(a)\int_S\big( |\nab_{\Rhat}\dk^{\leq 2}\psi|^2+r^{-1} |\dk^{\leq 2}\psi|^2\big) +\int_S\text{Bdr}
 \end{split}
\eeaa
and hence
\beaa
\begin{split}
\int_S |q|^2 \widetilde{\EE}[\bold{X}, \bold{w}, \bold{M}] \geq& \de\int_S\widetilde{P} +c_1  \int_S\Big( m|\nab_{\Rhat}\psi|^2_{\SS}+r^{-1}|\psi|^2_{\SS} \Big)\\
 & -O(a)\int_S\big( |\nab_{\Rhat}\dk^{\leq 2}\psi|^2+r^{-1} |\dk^{\leq 2}\psi|^2\big) +\int_S\text{Bdr}.
 \end{split}
\eeaa
Since 
\beaa
\widetilde{P} &=&\frac 1 2  h \Big( \de_0 \big|   \nab_T  \Psi_z \big|^2+a^2|\nab_Z \Psi_z|^2 + O^{\a\b} \Db_\a \Psi_z \Db_\b   \Psi_z\Big)  \\
&\geq &  c_0r^5\Big(  \big|   \nab_T  \Psi_z \big|^2 + |q|^2|\nab\Psi_z|^2\Big)
\eeaa
we infer
\beaa
\begin{split}
\int_S |q|^2 \widetilde{\EE}[\bold{X}, \bold{w}, \bold{M}] \geq& \de c_0\int_S r^5\Big(  \big|   \nab_T  \Psi_z \big|^2 + |q|^2|\nab\Psi_z|^2\Big) +c_1  \int_S\Big( m|\nab_{\Rhat}\psi|^2_{\SS}+r^{-1}|\psi|^2_{\SS} \Big)\\
 & -O(a)\int_S\big( |\nab_{\Rhat}\dk^{\leq 2}\psi|^2+r^{-1} |\dk^{\leq 2}\psi|^2\big) +\int_S\text{Bdr}.
 \end{split}
\eeaa
We finally obtain, for $|a|/m \ll1$,
\bea\lab{eq:veryfinalconclusionforthelowerboundgeneralizedcurrent:fdsoiugah}
\bsplit
\int_S \widetilde{\EE}[\bold{X}, \bold{w}, \bold{M}] \ges& \int_S  \frac{m}{r^2}  | \nab_\Rhat \psi|_\SS^2 + r^{-3}|\psi|_{\SS}^2+ r^{3} \Big(\big|\nab_\That \Psi_z \big |^2+r^2\big |\nab \Psi_z\big|^2\Big)+\text{Bdr}\\
&  -O(a)\int_S\big( r^{-2}|\nab_{\Rhat}\dk^{\leq 2}\psi|^2+r^{-3} |\dk^{\leq 2}\psi|^2\big)
\end{split}
\eea
on any sphere $S=S(t, r)$.


\subsection{$\SS$-derivative version of Lemma \ref{lemma-control-rhs}}
\lab{sec:SSderivativeversionofLemmacontrolrhs:aodihas}


In view of \eqref{definition-EE-gen-SSvalued} and \eqref{widetilde-EE}, we have
\bea\lab{definition-EE-gen-SSvalued:fsidufsdofgus}
\begin{split}
\D^\mu  \PP_\mu[\X, \w, \M] =&  \widetilde{\EE}[\bold{X}, \bold{w}, \bold{M}]  +  \D^\mu \BB_\mu + \left(\nab_{X^{\aund\bund}}\psia+\frac 1 2   w^{\aund\bund} \psia\right)\c \big( \squared_2\psib -V\psib \big)\\
& -\big(\rhod +\etab\wedge\eta\big)\nab_{(X^{\aund\bund})^4 e_4- (X^{\aund\bund})^3e_3}  \psi_\aund\c\dual\psi_\bund\\
 & -\frac{1}{2}\Im\Big(\tr\Xb H (X^{\aund\bund})^3 +\tr X\Hb (X^{\aund\bund})^4\Big)\c\nab\psi_\aund\c\dual\psi_\bund.
 \end{split}
\eea
Also, arguing as for the proof of \eqref{definition-EE-X=FFR}, we have
   \beaa
  \nab_{(X^{\aund\bund})^4 e_4- (X^{\aund\bund})^3e_3}\psi_\aund = \FF^{\aund\bund}\frac{r^2+a^2}{\De}\nab_{\That} \psi_\aund
  \eeaa
and
\beaa
\Im\Big(\tr\Xb H (X^{\aund\bund})^3 +\tr X\Hb (X^{\aund\bund})^4\Big) 
&=& \frac{4a^2r\cos\th\FF^{\aund\bund}(r)}{(r^2+a^2)|q|^4}\pr_\phi +\frac{4a^3r\cos\th(\sin\th)^2\FF^{\aund\bund}(r)}{|q|^6}\That.
\eeaa
 Thus we infer
  \bea\lab{definition-EE-gen-SSvalued:fsidufsdofgus}
\begin{split}
\D^\mu  \PP_\mu[\X, \w, \M] =&  \widetilde{\EE}[\bold{X}, \bold{w}, \bold{M}]  +  \D^\mu \BB_\mu + \left(\nab_{X^{\aund\bund}}\psia+\frac 1 2   w^{\aund\bund} \psia\right)\c \big( \squared_2\psib -V\psib \big)\\
& -\frac{2a^2r\cos\th\FF^{\aund\bund}(r)}{(r^2+a^2)|q|^4}\nab_\phi\psi_\aund\c\dual\psi_\bund\\
& -\left(\big(\rhod +\etab\wedge\eta\big)\FF^{\aund\bund}\frac{r^2+a^2}{\De}+\frac{2a^3r\cos\th(\sin\th)^2\FF^{\aund\bund}(r)}{|q|^6}\right)\nab_{\That}\psi_\aund\c\dual\psi_\bund.
 \end{split}
\eea

The goal of this section is to control the last three terms on the RHS of \eqref{definition-EE-gen-SSvalued:fsidufsdofgus}. To this end, we derive the $\SS$-derivative version of Lemma \ref{lemma-control-rhs}, being careful to obtain terms involving $\Psi_z$.

\begin{lemma}\lab{lemma-control-rhs-higher} 
We have, for sufficiently small positive constants $\de_2$, $\de_3$:
\bea
\begin{split}
&\left(\nab_{X^{\aund\bund}}\psia+\frac 1 2   w^{\aund\bund} \psia\right)\c \big( \squared_2\psib -V\psib \big)\\
\geq& -\de_2r^{-2} h|\nab_{\That} \Psi_z|^2 -\de_2 a^2r^{-6}h|\nab_Z \Psi_z|^2+O(ar^{-3}) |\psi|_\SS^2\\
&+O(1)\Big( | \nab_\Rhat \psi |_{\SS} + r^{-1}|\psi|_{\SS} \Big)\sum_{\aund=1}^4| N_{\aund}| + \D_\mu\left( \frac{ 2 a\cos\th}{|q|^2}(\pr_r)^\mu zh f^{\aund\bund} \psia\c\nab_T\dual \psib\right)\\
&   - \nab_T\left( \frac{ 2 a\cos\th}{|q|^2}\left( zhf^{\aund\bund} \psia\c  \nab_r\dual \psib +\frac 1 2 (\pr_r z) h \Psi_z\c\dual\LL^{\aund}\psia\right)\right),
\end{split}
\eea
and
\bea
\bsplit
& \left(\big(\rhod +\etab\wedge\eta\big)\FF^{\aund\bund}\frac{r^2+a^2}{\De}+\frac{2a^3r\cos\th(\sin\th)^2\FF^{\aund\bund}(r)}{|q|^6}\right)\nab_{\That}\psi_\aund\c\dual\psi_\bund\\
& +\frac{2a^2r\cos\th\FF^{\aund\bund}(r)}{(r^2+a^2)|q|^4}\nab_\phi\psi_\aund\c\dual\psi_\bund\\
\leq& \de_3 r^3\Big(|\nab_\That \Psi_z|^2+r^2|\nab\Psi_z|^2\Big)  +O(a r^{-3}) |\psi|_\SS^2 -\frac{1}{2}\nab_\phi\left(z h\frac{2a^2r\cos\th}{(r^2+a^2)|q|^4}\Psi_z\c\dual(\LL^{\underline{b}}\psi_\bund)\right)\\
 & -\frac{1}{2}\nab_{\That}\left(z h  \left(\big(\rhod +\etab\wedge\eta\big)\frac{r^2+a^2}{\De}+\frac{2a^3r\cos\th(\sin\th)^2}{|q|^6}\right)\Psi_z\c\dual(\LL^{\underline{b}}\psi_\bund)\right).
\end{split}
\eea
\end{lemma}

\begin{proof}
According to equation   \eqref{eq:waveeqfor-psia-chp3}, we have
\beaa
&&\left(\nab_{X^{\aund\bund}}( \psia)+\frac 1 2   w^{\aund\bund} \psia\right)\c \big( \squared_2\psib -V\psib \big)\\ 
&=& \left(\nab_{X^{\aund\bund}}( \psia)+\frac 1 2   w^{\aund\bund} \psia\right)\c \left(- \frac{ 4 a\cos\th}{|q|^2} \dual \nab_T  \psib+N_\bund \right).
\eeaa 
We consider the first order term, and we write, as in Lemma \ref{identity:X+frac12wpsiN_1},
\beaa
 \left(\nab_{X^{\aund\bund}}( \psia)+\frac 1 2   w^{\aund\bund} \psia\right)\c  \dual \nab_T  \psib&=& \left(-zhf^{\aund\bund}\nab_r \psia-\frac 1 2    z \pr_r \big( h  f^{\aund\bund}  \big)\psia\right)\c  \dual \nab_T  \psib\\
 &=& \frac 1 2 (\pr_r z) hf^{\aund\bund} \psia\c \nab_T\dual \psib +zhf^{\aund\bund}\rhod\frac{|q|^2}{\De} \psia\psib \\
   && -\frac 1 2 \nab_r\Big(  zh f^{\aund\bund} \psia\c\nab_T\dual \psib\Big)   +\frac 1 2\nab_T\Big(  zhf^{\aund\bund} \psia\c  \nab_r\dual \psib\Big).
\eeaa
Using that $f^{\aund\bund}= \tilde{\RR}'^{(\underline{a}} \LL^{\underline{b})}$ and $\Psi_z=\tilde{\RR}'^{\underline{a}}\psia$, we obtain
\beaa
&& \left(\nab_{X^{\aund\bund}}( \psia)+\frac 1 2   w^{\aund\bund} \psia\right)\c  \dual \nab_T  \psib\\ 
&=& \frac 1 4 (\pr_r z) h \Big(\LL^{\aund}\psia\c \nab_T\dual\Psi_z+\Psi_z\c \nab_T\dual\LL^{\aund}\psia\Big)  +zh \rhod\frac{|q|^2}{\De} \Psi_z\LL^{\underline{a}}\psia\\
   && -\frac 1 2 \nab_r\Big(  zh \LL^\aund \psia\c\nab_T\dual \Psi_z\Big)   +\frac 1 2\nab_T\Big(  zh\Psi_z\c  \LL^\bund\nab_r\dual \psib\Big)\\
&=& \frac 1 2 (\pr_r z) h \LL^{\aund}\psia\c \nab_T\dual\Psi_z  +zh \rhod\frac{|q|^2}{\De} \Psi_z\LL^{\underline{a}}\psia\\
   && -\frac 1 2 \nab_r\Big(  zh f^{\aund\bund} \psia\c\nab_T\dual \psib\Big)   +\frac 1 2\nab_T\left(  zhf^{\aund\bund} \psia\c  \nab_r\dual \psib   +\frac 1 2 (\pr_r z) h \Psi_z\c\dual\LL^{\aund}\psia\right).
\eeaa

Also, since $\X=O(1) \Rhat $ and $w^{\aund\bund}=O(r^{-1})$, we can bound the second product by
\beaa
\left| \left(\nab_{X^{\aund\bund}}( \psia)+\frac 1 2   w^{\aund\bund} \psia\right)\c N_{\bund} \right|&\les&  \Big( | \nab_\Rhat \psi |_{\SS} + r^{-1}|\psi|_{\SS} \Big)\sum_{\aund=1}^4| N_{\aund}|.
\eeaa
By putting together with the previous bounds we obtain
\beaa
&&\left(\nab_{X^{\aund\bund}}( \psia)+\frac 1 2   w^{\aund\bund} \psia\right)\c \big( \squared_2\psib -V\psib \big)\\
&=& \left(\nab_{X^{\aund\bund}}( \psia)+\frac 1 2   w^{\aund\bund} \psia\right)\c \left(- \frac{ 4 a\cos\th}{|q|^2} \dual \nab_T  \psib+N_\bund \right)\\
&\geq& -\de_2r^{-2} h|\nab_{T} \Psi_z|^2+O(ar^{-3}) |\psi|_\SS^2+O(1)\Big( | \nab_\Rhat \psi |_{\SS} + r^{-1}|\psi|_{\SS} \Big)\sum_{\aund=1}^4| N_{\aund}|\\
&& + \frac{ 2 a\cos\th}{|q|^2}\nab_r\left(  zh f^{\aund\bund} \psia\c\nab_T\dual \psib\right) \\
&&  - \nab_T\left( \frac{ 2 a\cos\th}{|q|^2}\left( zhf^{\aund\bund} \psia\c  \nab_r\dual \psib +\frac 1 2 (\pr_r z) h \Psi_z\c\dual\LL^{\aund}\psia\right)\right)\\
&\geq& -\de_2r^{-2} h|\nab_{\That} \Psi_z|^2 -\de_2 a^2r^{-6}h|\nab_Z \Psi_z|^2+O(ar^{-3}) |\psi|_\SS^2\\
&&+O(1)\Big( | \nab_\Rhat \psi |_{\SS} + r^{-1}|\psi|_{\SS} \Big)\sum_{\aund=1}^4| N_{\aund}| + \frac{ 2 a\cos\th}{|q|^2}\nab_r\left( zh f^{\aund\bund} \psia\c\nab_T\dual \psib\right) \\
&&  - \nab_T\left( \frac{ 2 a\cos\th}{|q|^2}\left( zhf^{\aund\bund} \psia\c  \nab_r\dual \psib +\frac 1 2 (\pr_r z) h \Psi_z\c\dual\LL^{\aund}\psia\right)\right).
\eeaa
Since we have $\D_\mu(\cos\th |q|^{-2}(\pr_r)^\mu)=0$, we infer
\beaa
&&\left(\nab_{X^{\aund\bund}}( \psia)+\frac 1 2   w^{\aund\bund} \psia\right)\c \big( \squared_2\psib -V\psib \big)\\
&\geq& -\de_2r^{-2} h|\nab_{\That} \Psi_z|^2 -\de_2 a^2r^{-6}h|\nab_Z \Psi_z|^2+O(ar^{-3}) |\psi|_\SS^2\\
&&+O(1)\Big( | \nab_\Rhat \psi |_{\SS} + r^{-1}|\psi|_{\SS} \Big)\sum_{\aund=1}^4| N_{\aund}| + \D_\mu\left(\frac{ 2 a\cos\th}{|q|^2}(\pr_r)^\mu zh f^{\aund\bund} \psia\c\nab_T\dual \psib\right)\\ 
&&  - \nab_T\left( \frac{ 2 a\cos\th}{|q|^2}\left( zhf^{\aund\bund} \psia\c  \nab_r\dual \psib +\frac 1 2 (\pr_r z) h \Psi_z\c\dual\LL^{\aund}\psia\right)\right).
\eeaa
as desired. 

For the second estimate, we write
\beaa
&& \left(\big(\rhod +\etab\wedge\eta\big)\FF^{\aund\bund}\frac{r^2+a^2}{\De}+\frac{2a^3r\cos\th(\sin\th)^2\FF^{\aund\bund}(r)}{|q|^6}\right)\nab_{\That}\psi_\aund\c\dual\psi_\bund\\
&& +\frac{2a^2r\cos\th\FF^{\aund\bund}(r)}{(r^2+a^2)|q|^4}\nab_\phi\psi_\aund\c\dual\psi_\bund\\
&=& z h  \tilde{\RR}'^{(\underline{a}} \LL^{\underline{b})}\left(\big(\rhod +\etab\wedge\eta\big)\frac{r^2+a^2}{\De}+\frac{2a^3r\cos\th(\sin\th)^2}{|q|^6}\right)\nab_{\That}\psi_\aund\c\dual\psi_\bund\\
&& +z h  \tilde{\RR}'^{(\underline{a}} \LL^{\underline{b})}\frac{2a^2r\cos\th}{(r^2+a^2)|q|^4}\nab_\phi\psi_\aund\c\dual\psi_\bund\\
&=&  z h  \tilde{\RR}'^{(\underline{a}} \LL^{\underline{b})}\left(\big(\rhod +\etab\wedge\eta\big)\frac{r^2+a^2}{\De}+\frac{2a^3r\cos\th(\sin\th)^2}{|q|^6}\right)\nab_{\That}\psi_\aund\c\dual\psi_\bund\\
&& +z h  \tilde{\RR}'^{(\underline{a}} \LL^{\underline{b})}\frac{2a^2r\cos\th}{(r^2+a^2)|q|^4}\nab_\phi\psi\c\dual\psi
\eeaa
and hence
\beaa
&& \left(\big(\rhod +\etab\wedge\eta\big)\FF^{\aund\bund}\frac{r^2+a^2}{\De}+\frac{2a^3r\cos\th(\sin\th)^2\FF^{\aund\bund}(r)}{|q|^6}\right)\nab_{\That}\psi_\aund\c\dual\psi_\bund\\
&& +\frac{2a^2r\cos\th\FF^{\aund\bund}(r)}{(r^2+a^2)|q|^4}\nab_\phi\psi_\aund\c\dual\psi_\bund\\ 
&=& z h  \left(\big(\rhod +\etab\wedge\eta\big)\frac{r^2+a^2}{\De}+\frac{2a^3r\cos\th(\sin\th)^2}{|q|^6}\right)\nab_{\That}\Psi_z\c\dual(\LL^{\underline{b}}\psi_\bund)\\
&& -\frac{1}{2}\nab_{\That}\left(z h  \left(\big(\rhod +\etab\wedge\eta\big)\frac{r^2+a^2}{\De}+\frac{2a^3r\cos\th(\sin\th)^2}{|q|^6}\right)\Psi_z\c\dual(\LL^{\underline{b}}\psi_\bund)\right)\\
&& +z h\frac{2a^2r\cos\th}{(r^2+a^2)|q|^4}\nab_\phi\Psi_z\c\dual(\LL^{\underline{b}}\psi_\bund) -\frac{1}{2}\nab_\phi\left(z h\frac{2a^2r\cos\th}{(r^2+a^2)|q|^4}\Psi_z\c\dual(\LL^{\underline{b}}\psi_\bund)\right)\\
&=& O(ar^{-1})\nab_{\That}\Psi_z\c\dual(\LL^{\underline{b}}\psi_\bund)+ O(a^2r^{-2})\nab_\phi\Psi_z\c\dual(\LL^{\underline{b}}\psi_\bund)\\
&& -\frac{1}{2}\nab_{\That}\left(z h  \left(\big(\rhod +\etab\wedge\eta\big)\frac{r^2+a^2}{\De}+\frac{2a^3r\cos\th(\sin\th)^2}{|q|^6}\right)\Psi_z\c\dual(\LL^{\underline{b}}\psi_\bund)\right)\\
&& -\frac{1}{2}\nab_\phi\left(z h\frac{2a^2r\cos\th}{(r^2+a^2)|q|^4}\Psi_z\c\dual(\LL^{\underline{b}}\psi_\bund)\right)\\
 &\leq & \de_3 r^3\Big(|\nab_\That \Psi_z|^2+r^2|\nab\Psi_z|^2\Big)  +O(a r^{-3}) |\psi|_\SS^2 -\frac{1}{2}\nab_\phi\left(z h\frac{2a^2r\cos\th}{(r^2+a^2)|q|^4}\Psi_z\c\dual(\LL^{\underline{b}}\psi_\bund)\right)\\
 && -\frac{1}{2}\nab_{\That}\left(z h  \left(\big(\rhod +\etab\wedge\eta\big)\frac{r^2+a^2}{\De}+\frac{2a^3r\cos\th(\sin\th)^2}{|q|^6}\right)\Psi_z\c\dual(\LL^{\underline{b}}\psi_\bund)\right)
\eeaa
as stated.
\end{proof}


\subsection{Proof of Proposition \ref{prop:morawetz-higher-order}}
\lab{section-SS-valued.unconditionalMorawetz}


We are now ready to prove Proposition \ref{prop:morawetz-higher-order}. Recall \eqref{definition-EE-gen-SSvalued:fsidufsdofgus}, i.e. 
\beaa
\begin{split}
\D^\mu  \PP_\mu[\X, \w, \M] =&  \widetilde{\EE}[\bold{X}, \bold{w}, \bold{M}]  +  \D^\mu \BB_\mu + \left(\nab_{X^{\aund\bund}}\psia+\frac 1 2   w^{\aund\bund} \psia\right)\c \big( \squared_2\psib -V\psib \big)\\
& -\frac{2a^2r\cos\th\FF^{\aund\bund}(r)}{(r^2+a^2)|q|^4}\nab_\phi\psi_\aund\c\dual\psi_\bund\\
& -\left(\big(\rhod +\etab\wedge\eta\big)\FF^{\aund\bund}\frac{r^2+a^2}{\De}+\frac{2a^3r\cos\th(\sin\th)^2\FF^{\aund\bund}(r)}{|q|^6}\right)\nab_{\That}\psi_\aund\c\dual\psi_\bund.
 \end{split}
\eeaa
We apply  the divergence theorem to the above on $\MM(\tau_1, \tau_2)$, which yields
\beaa
&&\int_{\MM(\tau_1, \tau_2)}\Bigg[\widetilde{\EE}[\bold{X}, \bold{w}, \bold{M}]  + \left(\nab_{X^{\aund\bund}}\psia+\frac 1 2   w^{\aund\bund} \psia\right)\c \big( \squared_2\psib -V\psib\big)\\
&& -\left(\big(\rhod +\etab\wedge\eta\big)\FF^{\aund\bund}\frac{r^2+a^2}{\De}+\frac{2a^3r\cos\th(\sin\th)^2\FF^{\aund\bund}(r)}{|q|^6}\right)\nab_{\That}\psi_\aund\c\dual\psi_\bund\\
&& -\frac{2a^2r\cos\th\FF^{\aund\bund}(r)}{(r^2+a^2)|q|^4}\nab_\phi\psi_\aund\c\dual\psi_\bund\Bigg]\\
&\leq& \int_{\pr\MM(\tau_1, \tau_2)}\big(|\PP_\mu[\X, \w, \M] N^\mu|+|\BB_\mu N^\mu|\big).
\eeaa

Using the lower bound \eqref{eq:veryfinalconclusionforthelowerboundgeneralizedcurrent:fdsoiugah} for $\widetilde{\EE}[\bold{X}, \bold{w}, \bold{M}]$,  and estimating the three other terms on the LHS thanks to Lemma \ref{lemma-control-rhs-higher}, we obtain, for sufficiently small $\de_2$ and $\de_3$, 
 \beaa
 \bsplit
  \Mor_{\SSz, deg}[\psi](\tau_1, \tau_2) \les& \int_{\pr\MM(\tau_1, \tau_2) }|M_\SS(\psi)| + |a|\int_{\MM(\tau_1, \tau_2)}\big( r^{-2}|\nab_{\Rhat}\dk^{\leq 2}\psi|^2+r^{-3} |\dk^{\leq 2}\psi|^2\big)\\
  &+\sum_{\aund=1}^4\int_{\MM(\tau_1, \tau_2)}\big(|\nab_{\Rhat} \psia|+r^{-1}|\psia|\big) |N_{\aund}|
  \end{split}
 \eeaa
where we recall that
\beaa
  \Mor_{\SSz, deg}[\psi](\tau_1, \tau_2)&=& \int_{\MM(\tau_1, \tau_2)} 
      \frac{m}{r^2}  | \nab_\Rhat \psi|_\SS^2 + r^{-3}|\psi|_{\SS}^2+ r^{3} \Big(\big|\nab_\That \Psi_z \big |^2+r^2\big |\nab \Psi_z\big|^2\Big),
\eeaa
and where $M_\SS(\psi)$  denotes  an expression in  $\psi$ for which we  have a bound of the form
\beaa
&&\int_{\pr\MM(\tau_1, \tau_2) }|M_\SS(\psi)| \\
&\les &  \sum_{\aund=1}^4\left(\sup_{[\tau_1, \tau_2]}E_{deg}[\psi_\aund](\tau) +\deh F_{\AA}[\psi_{\aund}](\tau_1, \tau_2)+F_{\Si_*}[\psi_{\aund}](\tau_1, \tau_2)\right)\\
&+&\left(\sup_{[\tau_1, \tau_2]}E_{deg}[(\nab_T, \dkb)^{\leq 1}\psi](\tau)+\deh F_{\AA}[(\nab_T, \dkb)^{\leq 1}\psi](\tau_1, \tau_2)+F_{\Si_*}[(\nab_T, \dkb)^{\leq 1}\psi](\tau_1, \tau_2)\right)^{\frac{1}{2}}\\
&\times&\left(\sup_{[\tau_1, \tau_2]}E_{deg}[(\nab_T, \dkb)^{\leq 2}\psi](\tau)+\deh F_{\AA}[(\nab_T, \dkb)^{\leq 2}\psi](\tau_1, \tau_2)+F_{\Si_*}[(\nab_T, \dkb)^{\leq 2}\psi](\tau_1, \tau_2)\right)^{\frac{1}{2}}.
\eeaa 
Finally, since we have 
\beaa
&&\int_{\MM(\tau_1, \tau_2)}\big( r^{-2}|\nab_{\Rhat}\dk^{\leq 2}\psi|^2+r^{-3} |\dk^{\leq 2}\psi|^2\big)\\
&\les& \int_{\MM_{trap}(\tau_1, \tau_2)}\big( r^{-2}|\nab_{\Rhat}\dk^{\leq 2}\psi|^2+r^{-3} |\dk^{\leq 2}\psi|^2\big) +\int_{\Mntrap(\tau_1, \tau_2)}\big( r^{-2}|\nab_3\dk^{\leq 2}\psi|^2+r^{-3} |\dk^{\leq 3}\psi|^2\big)\\
&\les& B^2_\de[\psi](\tau_1, \tau_2)
\eeaa
for any $\de>0$ in view of the definition of $B^2_\de[\psi](\tau_1, \tau_2)$, we infer 
 \beaa
 \bsplit
  \Mor_{\SSz, deg}[\psi](\tau_1, \tau_2) \les& \int_{\pr\MM(\tau_1, \tau_2) }|M_\SS(\psi)| +\frac{|a|}{m}B^2_\de[\psi](\tau_1, \tau_2)\\
  &+\sum_{\aund=1}^4\int_{\MM(\tau_1, \tau_2)}\big(|\nab_{\Rhat} \psia|+r^{-1}|\psia|\big) |N_{\aund}|
  \end{split}
 \eeaa
which concludes the proof of  Proposition \ref{prop:morawetz-higher-order}.


\section{Proof of Lemma \ref{LEMMA:LOWERBOUNDPHIZOUTSIDEMTRAP}}\label{section:lowerboundPhizoutsideMtrap}


In this section, we prove Lemma \ref{LEMMA:LOWERBOUNDPHIZOUTSIDEMTRAP} on the lower bound for $\Psi_z$ on $\Mntrap$. Recall the definition \eqref{definition-Psiz} of $\Psi_z$
\beaa
\begin{split}
\Psi_z&= -\frac{2\TT}{(r^2+a^2)^3}   \big(\de_0  \SS_1\psi+ (1+O(r^{-2} \de_0)) \OO\psi\big)+ \frac{4ar}{(r^2+a^2)^2}  \nab_\That \nab_Z \psi   \big(1+O(r^{-2} \de_0) \big),
\end{split}
\eeaa
as well as the definition of $\Mntrap$, see \eqref{eq:def-MM-trap},
\beaa
\Mntrap=\left\{     \frac{|\TT|}{r^3}  \geq \frac{1}{10} \right\}.
\eeaa
In particular, we have on $\Mntrap$
\beaa
|\nab_T\Psi_z|+r|\nab\Psi_z| &\geq& \frac{1}{5}\frac{r^3}{(r^2+a^2)^3}\Big(\left|\nab_T\big(\de_0  \SS_1\psi+  \OO\psi\big)\right|+|q|\left|\nab\big(\de_0  \SS_1\psi+  \OO\psi\big)\right|\Big)\\
&& -O(\de_0r^{-3})\Big(|\nab_T\OO\psi|+r|\nab\OO \psi|\Big)\\
&& -O(ar^{-3})\Big(|\nab_T\nab_\That \nab_Z \psi|+r|\nab\nab_\That \nab_Z \psi|\Big) -O((\de_0+|a|)r^{-3})|\psi|_{\SS}
\eeaa
and hence for a universal constant $c_0>0$
\beaa
\bsplit
r^3\Big(|\nab_T\Psi_z|^2+r^2|\nab\Psi_z|^2\Big) \geq& c_0r^{-3}\Big(\left|\nab_T\big(\de_0  \SS_1\psi+  \OO\psi\big)\right|^2+|q|^2\left|\nab\big(\de_0  \SS_1\psi+  \OO\psi\big)\right|^2\Big)\\
& -O(a^2r^{-3})\Big(|\nab_T\nab_\That \nab_Z \psi|^2+r^2|\nab\nab_\That \nab_Z \psi|^2\Big)\\
& -O(\de_0^2r^{-3})\Big(|\nab_T\OO\psi|^2+r^2|\nab\OO \psi|^2\Big)  - O(r^{-3})|\psi|_{\SS}^2.
\end{split}
\eeaa
which we rewrite 
\bea\lab{eq:lowerboundPhizoutsideMtrap:1}
\bsplit
r^3\Big(|\nab_T\Psi_z|^2+r^2|\nab\Psi_z|^2\Big) \geq& c_0r^{-3}\Big(\left|\nab_T\big(\de_0  \SS_1\psi+  \OO\psi\big)\right|^2+|q|^2\left|\nab\big(\de_0  \SS_1\psi+  \OO\psi\big)\right|^2\Big)\\
& -O(\de_0^2r^{-3})\Big(|\nab_T\OO\psi|^2+r^2|\nab\OO \psi|^2\Big)\\
& -O(ar^{-3})|(\nab_T, \dkb)\dk^{\leq 2}\psi|^2 - O(r^{-3})|\psi|_{\SS}^2.
\end{split}
\eea

Next, we focus on the first term on the RHS of \eqref{eq:lowerboundPhizoutsideMtrap:1}. We start with the following computation 
\beaa
&&\left|\nab_T\big(\de_0  \SS_1\psi+  \OO\psi\big)\right|^2+|q|^2\left|\nab\big(\de_0  \SS_1\psi+  \OO\psi\big)\right|^2\\ 
&=& \de_0^2\Big(|\nab_T\SS_1\psi|^2+|q|^2|\nab\SS_1\psi|^2\Big)+2\de_0\Big(\nab_T\SS_1\psi\c\nab_T\OO\psi+|q|^2\nab\SS_1\psi\c\nab\OO\psi\Big)\\
&& + |\nab_T\OO\psi|^2+|q|^2|\nab\OO\psi|^2.
\eeaa
Using the structure of commutators\footnote{Note in particular that $[\nab, \OO] = -4\nab\psi+O(ar^{-1})(\nab_\T, \dkb)^{\leq 2}\psi$.}, we write the second term on the RHS as follows 
\beaa
&&\nab_T\SS_1\psi\c\nab_T\OO\psi+|q|^2\nab\SS_1\psi\c\nab\OO\psi\\
 &=& \SS_1\nab_T\psi\c\OO\nab_T\psi+|q|^2\SS_1\nab\psi\c\OO\nab\psi\\
 &&+ \SS_1\nab_T\psi\c[\nab_T,\OO]\psi+|q|^2\nab\SS_1\psi\c[\nab, \OO]\psi+|q|^2[\nab,\SS_1]\psi\c\OO\nab\psi\\
 &=& \SS_1\nab_T\psi\c\OO\nab_T\psi+|q|^2\SS_1\nab\psi\c\OO\nab\psi -|q|^2\nab\SS_1\psi\c\nab\psi +O(a)\big((\nab_T, \dkb)\dk^{\leq 2}\psi\big)^2\\
&=& \SS_1\nab_T\psi\c\OO\nab_T\psi+|q|^2\SS_1\nab\psi\c\OO\nab\psi +|q|^2|\nab\nab_T\psi|^2
+O(a)\big((\nab_T, \dkb)\dk^{\leq 2}\psi\big)^2\\
&& -\pr_t(|q|^2\nab\nab_T\psi\c\nab\psi).
\eeaa
Next, we rely on Lemma \ref{Lemma:integrationbypartsSS_3SS_4} and obtain 
\beaa
&&\nab_T\SS_1\psi\c\nab_T\OO\psi+|q|^2\nab\SS_1\psi\c\nab\OO\psi\\
 &=& |q|^2| \nab \nab_T^2\psi|^2+|q|^2 \Ddot_\b(|q|^{-2}O^{\a\b}\Ddot_\a \nab_T\psi \c \SS_1 \nab_T\psi ) +\partial_t (M(\nab_T\psi))\\
&&+|q|^4| \nab^2\nab_T\psi|^2+|q|^4 \Ddot_\b(|q|^{-2}O^{\a\b}\Ddot_\a \nab\psi \c \SS_1\nab\psi ) +\partial_t (r^2M(\nab\psi))\\
&&+|q|^2|\nab\nab_T\psi|^2+O(a)\big((\nab_T, \dkb)^{\leq 1}\dk^{\leq 2}\psi\big)^2 -\pr_t(|q|^2\nab\nab_T\psi\c\nab\psi).
\eeaa
Rearranging, this yields
\beaa
&&\nab_T\SS_1\psi\c\nab_T\OO\psi+|q|^2\nab\SS_1\psi\c\nab\OO\psi\\ 
&=& |q|^2| \nab \nab_T^2\psi|^2 +|q|^4| \nab^2\nab_T\psi|^2 +|q|^2|\nab\nab_T\psi|^2 -O(a)\big|(\nab_T, \dkb)^{\leq 1}\dk^{\leq 2}\psi\big|^2\\
 &&+|q|^2 \Ddot_\b(|q|^{-2}O^{\a\b}\Ddot_\a \nab_T\psi \c \SS_1 \nab_T\psi ) +|q|^4 \Ddot_\b(|q|^{-2}O^{\a\b}\Ddot_\a \nab\psi \c \SS_1\nab\psi )\\
 && +\partial_t (M(\nab_T\psi)) +\partial_t (M(|q|\nab\psi)) -\pr_t(|q|^2\nab\nab_T\psi\c\nab\psi)
\eeaa
and hence
\bea\lab{eq:lowerboundPhizoutsideMtrap:2}
\nn&&\left|\nab_T\big(\de_0  \SS_1\psi+  \OO\psi\big)\right|^2+|q|^2\left|\nab\big(\de_0  \SS_1\psi+  \OO\psi\big)\right|^2\\ 
\nn&=& \de_0^2\Big(|\nab_T\SS_1\psi|^2+|q|^2|\nab\SS_1\psi|^2\Big)+2\de_0\Big( |q|^2| \nab \nab_T^2\psi|^2 +|q|^4| \nab^2\nab_T\psi|^2 +|q|^2|\nab\nab_T\psi|^2 \Big)\\
\nn&& + |\nab_T\OO\psi|^2+|q|^2|\nab\OO\psi|^2 -O(a)\big|(\nab_T, \dkb)^{\leq 1}\dk^{\leq 2}\psi\big|^2\\
 \nn&&+\de_0|q|^2 \Ddot_\b(|q|^{-2}O^{\a\b}\Ddot_\a \nab_T\psi \c \SS_1 \nab_T\psi ) +\de_0|q|^4 \Ddot_\b(|q|^{-2}O^{\a\b}\Ddot_\a \nab\psi \c \SS_1\nab\psi )\\
 && +\partial_t (M(\nab_T\psi)) +\partial_t (M(|q|\nab\psi)) -\pr_t(2\de_0|q|^2\nab\nab_T\psi\c\nab\psi).
\eea

Next, notice that $\SS_2=O(a)\dk^{\leq 2}$ and $\SS_3=O(a^2)\dk^{\leq 2}$ so that 
\beaa
|\nab_T\SS_1\psi|^2+|q|^2|\nab\SS_1\psi|^2  + |\nab_T\OO\psi|^2+|q|^2|\nab\OO\psi|^2 &\geq & |\nab_T\psi|_{\SS}+r^2|\nab\psi|_{\SS}\\
&& -O(a^2)\big|(\nab_T, \dkb)^{\leq 1}\dk^{\leq 2}\psi\big|^2.
\eeaa
Together with \eqref{eq:lowerboundPhizoutsideMtrap:2}, we infer
\beaa
\nn&&\left|\nab_T\big(\de_0  \SS_1\psi+  \OO\psi\big)\right|^2+|q|^2\left|\nab\big(\de_0  \SS_1\psi+  \OO\psi\big)\right|^2\\ 
\nn&=& \de_0^2\Big(|\nab_T\psi|_{\SS}+r^2|\nab\psi|_{\SS}\Big) +(1-\de_0^2)\Big( |\nab_T\OO\psi|^2+|q|^2|\nab\OO\psi|^2\Big)-O(a)\big|(\nab_T, \dkb)^{\leq 1}\dk^{\leq 2}\psi\big|^2\\
 \nn&&+\de_0|q|^2 \Ddot_\b(|q|^{-2}O^{\a\b}\Ddot_\a \nab_T\psi \c \SS_1 \nab_T\psi ) +\de_0|q|^4 \Ddot_\b(|q|^{-2}O^{\a\b}\Ddot_\a \nab\psi \c \SS_1\nab\psi )\\
 && +\partial_t (M(\nab_T\psi)) +\partial_t (M(|q|\nab\psi)) -\pr_t(2\de_0|q|^2\nab\nab_T\psi\c\nab\psi).
\eeaa
Next, plugging in \eqref{eq:lowerboundPhizoutsideMtrap:1}, we deduce
\beaa
\bsplit
& r^3\Big(|\nab_T\Psi_z|^2+r^2|\nab\Psi_z|^2\Big)\\
 \geq& c_0r^{-3}\Big(\left|\nab_T\big(\de_0  \SS_1\psi+  \OO\psi\big)\right|^2+|q|^2\left|\nab\big(\de_0  \SS_1\psi+  \OO\psi\big)\right|^2\Big)\\
& -O(\de_0^2r^{-3})\Big(|\nab_T\OO\psi|^2+r^2|\nab\OO \psi|^2\Big) -O(ar^{-3})|(\nab_T, \dkb)\dk^{\leq 2}\psi|^2 - O(r^{-3})|\psi|_{\SS}^2\\
\geq& c_0\de_0^2r^{-3}\Big(|\nab_T\psi|_{\SS}+r^2|\nab\psi|_{\SS}\Big) +c_0r^{-3}(1-O(\de_0^2))\Big( |\nab_T\OO\psi|^2+|q|^2|\nab\OO\psi|^2\Big)\\
& -O(ar^{-3})\big|(\nab_T, \dkb)^{\leq 1}\dk^{\leq 2}\psi\big|^2 - O(r^{-3})|\psi|_{\SS}^2\\
& +\de_0r^{-3}|q|^2 \Ddot_\b(|q|^{-2}O^{\a\b}\Ddot_\a \nab_T\psi \c \SS_1 \nab_T\psi ) +\de_0r^{-3}|q|^4 \Ddot_\b(|q|^{-2}O^{\a\b}\Ddot_\a \nab\psi \c \SS_1\nab\psi )\\
 & +\partial_t (r^{-3}M(\nab_T\psi)) +\partial_t (r^{-3}M(|q|\nab\psi)) -\pr_t(2\de_0r^{-3}|q|^2\nab\nab_T\psi\c\nab\psi).
\end{split}
\eeaa
Hence, fixing $\de_0>0$ small enough, we infer
\beaa
\bsplit
& r^3\Big(|\nab_T\Psi_z|^2+r^2|\nab\Psi_z|^2\Big)\\
\geq& c_0\de_0^2r^{-3}\Big(|\nab_T\psi|_{\SS}+r^2|\nab\psi|_{\SS}\Big) -O(ar^{-3})\big|(\nab_T, \dkb)^{\leq 1}\dk^{\leq 2}\psi\big|^2 - O(r^{-3})|\psi|_{\SS}^2\\
& + \Ddot_\b(\de_0 O(r^{-3})O^{\a\b}\Ddot_\a \nab_T\psi \c \SS_1 \nab_T\psi ) +\Ddot_\b(\de_0 O(r^{-1})O^{\a\b}\Ddot_\a \nab\psi \c \SS_1\nab\psi )\\
 & +\partial_t (r^{-3}M(\nab_T\psi)) +\partial_t (r^{-3}M(|q|\nab\psi)) -\pr_t(2\de_0r^{-3}|q|^2\nab\nab_T\psi\c\nab\psi).
\end{split}
\eeaa
In particular, we deduce the existence of a universal constant $c_0>0$ such that the following holds on $\Mntrap$
\beaa
r^3\Big(|\nab_T\Psi_z|^2+r^2|\nab\Psi_z|^2\Big)+r^{-3}|\psi|_{\SS}^2 &\geq& c_0r^{-3}\Big(|\nab_T\psi|^2_{\SS}+|\nab_Z\psi|^2_{\SS}+r^2|\nab\psi|^2_{\SS}\Big) \\
&& -O(ar^{-3})\big|(\nab_T, \dkb)^{\leq 1}\dk^{\leq 2}\psi\big|^2+\textrm{Bdr}
\eeaa
where
\beaa
\textrm{Bdr} &=&  \Ddot_\b(O(r^{-3})O^{\a\b}\Ddot_\a \nab_T\psi \c \SS_1 \nab_T\psi ) + \Ddot_\b(O(r^{-1})O^{\a\b}\Ddot_\a \nab\psi \c \SS_1\nab\psi )\\
 && +\partial_t (r^{-3}M(\nab_T\psi)) +\partial_t (r^{-3}M(|q|\nab\psi)) -\pr_t(O(r^{-1})\nab\nab_T\psi\c\nab\psi).
\eeaa
Finally, note that the boundary terms verify 
\beaa
\textrm{Bdr} &=& \Ddot_\a F^\a, \qquad |F^\mu N_\mu| \les r^{-2}|(\nab_{\Rhat}, \nab_T, \dkb)^{\leq 1}(\nab_T, \dkb)^{\leq 1}\psi||(\nab_{\Rhat}, \nab_T, \dkb)^{\leq 1}(\nab_T, \dkb)^{\leq 2}\psi|,
\eeaa
where $N$ denotes the normal to either $\Si(\tau)$, $\AA$ or $\Si_*$. This concludes the proof of Lemma \ref{LEMMA:LOWERBOUNDPHIZOUTSIDEMTRAP}.


\chapter{Energy-Morawetz in perturbations of Kerr}
\label{chapter-perturbations-Kerr}


In this chapter, we prove Theorem \ref{THM:HIGHERDERIVS-MORAWETZ-CHP3}, i.e. we establish Energy-Morawetz estimates for solutions to the model gRW equation \eqref{eq:Gen.RW}
\beaa
\squared_2 \psi -V\psi=- \frac{4 a\cos\th}{|q|^2}\dual \nab_T  \psi+N, \qquad V= \frac{4\De}{ (r^2+a^2) |q|^2},
\eeaa
in perturbations of Kerr. To this end, we proceed as follows:
Theorem \ref{THM:HIGHERDERIVS-MORAWETZ-CHP3} is proved in Chapter \ref{chapter-perturbations-Kerr} according to the following steps:
\begin{enumerate}
\item First, we revisit the proof of Propositions \ref{proposition:Morawetz1-step1}, \ref{proposition:Energy1}  and  \ref{prop:morawetz-higher-order} by exhibiting the extra terms in perturbations of Kerr, and prove that the conclusions of Propositions \ref{proposition:Morawetz1-step1}, \ref{proposition:Energy1}  and  \ref{prop:morawetz-higher-order}   also hold in perturbations of Kerr up to the addition of suitable error terms see section \ref{section:adaptiation-errors}.

\item Next, we prove redshift estimates to remove the degeneracy on the horizon, see section \ref{section:Redshift-estimates-chp3}.

\item Then, we derive the conclusion of Theorem \ref{THM:HIGHERDERIVS-MORAWETZ-CHP3} for $s=2$, see section \ref{sec:proofofThm:HigherDerivs-Morawetz-chp3:cases=2}.

\item Finally, we argue by iteration from $s=2$ to recover higher order derivatives which and conclude the proof of Theorem \ref{THM:HIGHERDERIVS-MORAWETZ-CHP3}, see section \ref{sec:proofofThm:HigherDerivs-Morawetz-chp3:generalcase}.
\end{enumerate}


\section{Preliminaries}\label{section:perturbations-Kerr}


In this section, we recall the basic set up and state the analog of Propositions \ref{proposition:Morawetz1-step1}, \ref{proposition:Energy1}  and  \ref{prop:morawetz-higher-order} in perturbations of Kerr.


\subsection{Admissible perturbations of Kerr}
\lab{sec:setupchap9}


In this chapter, we prove Theorem \ref{THM:HIGHERDERIVS-MORAWETZ-CHP3} for solution of the model gRW equation \eqref{eq:Gen.RW} on the spacetime $\MM$, where $\MM$ is an admissible perturbation of Kerr which satisfies the assumptions of section \ref{sec:pfdoisdvhdifuhgiwhdniwbvoubwuyf}. We briefly recall them below for the convenience of the reader.

   
     \subsubsection{The spacetime $\MM$}


We consider a given vacuum spacetime $\MM$ satisfying the following properties:   
\begin{itemize}
\item $\MM$ comes together with a null pair $(e_4, e_3)$ and its corresponding horizontal structure as in section \ref{sec:nullparisandhorizontalstruct}. 

\item $\MM$ is endowed with a pair of constants $(a, m)$.

\item $\MM$ is endowed with a pair of scalar functions $(r, \th)$. 

\item The complex valued scalar function $q$ is defined as
\beaa
q := r+i\cos\th. 
\eeaa

\item $\MM$ is endowed with a complex horizontal 1-form $\Jk$. 

\item $\MM$ is endowed with a scalar function $\tau$ whose level sets $\Si(\tau)$ are spacelike. Also:
\begin{itemize}
\item $\tau\in[1,\tau_*]$ on $\MM$ for some arbitrary large constant $\tau_*$. 

\item  Given a level hypersurface  $\Si=\Si(\tau)$, we denote  
  \beaa
  N_\Si :=-\g^{\a\b}\pr_\b\tau \pr_\a.
  \eeaa

\item $\tau$ satisfies the properties of Definition \ref{definition:definition-oftau}.
\end{itemize}

\item The boundary of $\MM$ is given by 
\beaa
\pr\MM=\AA\cup\Si_*\cup\Si(1)\cup\Si(\tau_*)
\eeaa
where 
\beaa
\AA:=\Big\{r=r_+-\deh, \, 1\leq\tau\leq\tau_*\Big\}, 
\eeaa
and $\Si_*$ is a spacelike hypersurface on which $\tau$ takes the values $[1,\tau_*]$ and $r\geq r_*$ with $r_*\gg \tau_*$. 

\item Let $r_0$ a large enough fixed constant. We decompose $\MM$ as follows
\beaa
\Mint:=\MM\cap\{r\leq r_0\}, \qquad \Mext:=\MM\cap\{r\geq r_0\}. 
\eeaa
\end{itemize}


   \subsubsection{Admissible perturbations of Kerr}


Recall that $\MM$ comes together three scalar functions $(r, \th, \tau)$, and with a null pair $(e_4, e_3)$ and its corresponding horizontal structure as in section \ref{sec:nullparisandhorizontalstruct}. Then:
\begin{itemize}
\item We use the complexified Ricci and curvature coefficients  of Definition \ref{def:complexRicciandcurvaturecoefficients}. 

\item We define  the linearized quantities corresponding to these complexified coefficients as in Definition \ref{def:renormalizationofallnonsmallquantitiesinPGstructurebyKerrvalue:1} and 
\ref{def:renormalizationofallnonsmallquantitiesinPGstructurebyKerrvalue:3}, i.e. we consider that the   normalization of $(e_3, e_4)$ is ingoing.

\item With respect to these linearized quantities, we the notations $\Ga_g$ and $\Ga_b$ for error terms are given by Definition \ref{definition.Ga_gGa_b}.
\end{itemize}

We this definition of $\Ga_g$ and $\Ga_b$, we can now state our main assumptions on $\MM$ allowing to prove the Energy-Morawetz estimates of Theorem \ref{THM:HIGHERDERIVS-MORAWETZ-CHP3}. Let $\kl$ a large enough integer, and let the scalar function 
$\tau_{trap}$ defined by
\beaa
\tau_{trap} := \left\{\ba{lll}
1+\tau & \textrm{on} & \MM_{trap},\\
1& \textrm{on} & \Mntrap.
\ea\right.
\eeaa
Then, we assume that $(\Ga_g, \Ga_b)$ satisfy the following estimates on $\MM$
\bea\lab{eq:assumptionsonMMforpartII:again}
\bsplit
r^3|\dk^{\leq k}\xi|+ r^2|\dk^{\leq k}\Ga_g|+r|\dk^{\leq k}\Ga_b| &\leq \frac{\ep}{\tau_{trap}^{1+\dec}}, \quad 0\leq k\leq \frac{\kl}{2},\\
r^3|\dk^{\leq k}\xi|+ r^2|\dk^{\leq k}\Ga_g|+r|\dk^{\leq k}\Ga_b| &\leq \ep, \quad k\leq \kl.
\end{split}
\eea

\begin{remark}\lab{rmk:strongerassumptionxiinrm1Gagforchapter9}
Note that the assumptions for $\xi$ in \eqref{eq:assumptionsonMMforpartII:again} are consistent with $\xi\in r^{-1}\Ga_g$, see Remark \ref{rmk:xiisactuallybetterthanGag} for the justification of these stronger assumptions. We may thus assume in this chapter that $\xi\in r^{-1}\Ga_g$. This will be needed to deal with various commutators with the wave operator. 
\end{remark}


\subsection{Regions of integration and vectorfields}
\lab{sec:setupchap9:bis}


   
\subsubsection{Regions of integration}

  
  Recall the time function $\tau$ introduced in Definition \ref{definition:definition-oftau}. 
  We denote by $\Sigma_\tau$  the level sets of the  function $\tau$.

     \begin{definition}\label{def:causalregions}
    We define the following regions of $\MM$.
    \begin{enumerate}
      \item   We define the trapping region of $\MM$ to be the set 
    \bea\lab{eq:def-MM-trap}
     \MM_{trap}(\de_{trap})=\left\{     \frac{|\TT|}{r^3}  \le \de_{trap} \right\}, \qquad \delta_{trap} = \frac{1}{10},
    \eea
    where  $\TT$ is  the   polynomial in $r$ defined in \eqref{definition-TT}, i.e.
    \beaa
    \TT=r^3-3mr^2+ a^2r+ma^2.
    \eeaa
  
  \item We denote $\Mntrap$ the complement to the trapping region $\MM_{trap}$. 
    
   \item We define the  domain $\MM(\tau_1,\tau_2)$ to be the  region of $\MM$ where  $\tau_1\le \tau\le \tau_2$, where $\tau$ is the time function defined in Definition \ref{definition:definition-oftau}.
      \end{enumerate}
     \end{definition}

   
\subsubsection{Basic vectorfields}


$\T$ and $\Z$ are defined in $\MM$  as follows
\beaa
\T &=& \frac{1}{2}\left(e_4+\frac{\Delta}{|q|^2}e_3 -2a\Re(\Jk)^be_b\right),\\
\Z &=& \frac 1 2 \left(2(r^2+a^2)\Re(\Jk)^be_b -a(\sin\th)^2 e_4 -\frac{a(\sin\th)^2\De}{ |q|^2} e_3\right).
\eeaa

\begin{remark}
Recall from Lemma \ref{rem:de-trap}  that for  $|a|/m$ sufficiently small  and  $\delta_{trap} = \frac{1}{10}$,   the vectorfield  $\T$ is strictly timelike in $\MM_{trap}$.
\end{remark}

Also, we define the following vectorfields
 \beaa
 \That &=&\frac 1 2 \left( \frac{|q|^2}{r^2+a^2} e_4+\frac{\De}{r^2+a^2}  e_3\right), \qquad \Rhat= \frac 1 2 \left( \frac{|q|^2}{r^2+a^2} e_4-\frac{\De}{r^2+a^2}  e_3\right).
 \eeaa

Finally,  the vectorfield $\That_\de$ that will be used for energy estimates is given by
\beaa
\That_\de= \T+\frac{a}{r^2+a^2} \chi_0\left( \de^{-1} \frac{\TT}{r^3} \right) \Z
\eeaa
with $\de=\de_{trap}$ and with $\chi_0$ the smooth  bump function
 \beaa
 \chi_0(x)=\begin{cases}
 &0  \qquad \qquad  \mbox{if}   \quad  |x| \le  1, \\
 &1  \qquad  \qquad \mbox{if}   \quad   |x|\ge 2. 
 \end{cases} 
 \eeaa
 We also write
 \beaa
 \That_\de= \T+\chi_\de \Z, \qquad \chi_\de:= \frac{a}{r^2+a^2} \chi_0\left( \de^{-1} \frac{\TT}{r^3} \right).
 \eeaa

    
    \subsection{Main norms}
        
    
We recall below the relevant main norms  introduced in section \ref{subsection:basicnormsforpsi}:

   {\bf 1. Reduced  basic Morawetz  norms.} 
       \bea
\bsplit
 \Mor[\psi](\tau_1, \tau_2)&:=\int_{\MM(\tau_1, \tau_2) } 
      r^{-2}  | \nab_\Rhat  \psi|^2 +r^{-3}|\psi|^2\\
      &+ \int_{\Mntrap(\tau_1, \tau_2)} \left(  r^{-2}|\nab_3\psi|^2 + r^{-1}  |\nab  \psi|^2\right),\\
     \Morr[\psi](\tau_1, \tau_2)&:= \Mor[\psi](\tau_1, \tau_2)  + \int_{\MM_{r\geq 4m}(\tau_1,\tau_2)}  r^{-1-\de} |\nab_3 \psi|^2.
   \end{split}
\eea

{\bf 2. Basic Energy norm.} 
\bea
 E[\psi](\tau) :=\int_{\Si (\tau)} \Big( |\nab_4\psi|^2 +   r^{-2}|\nab_3\psi|^2 +|\nab\psi|^2 + r^{-2}|\psi|^2\Big).
\eea

{\bf 3. Basic Flux  norm.} 
\bea
\bsplit
F[\psi](\tau_1,\tau_2) : =& F_{\AA}[\psi](\tau_1,\tau_2)+F_{\Si_*}[\psi](\tau_1,\tau_2),\\
F_{\AA}[\psi](\tau_1,\tau_2) := & \int_{\AA(\tau_1, \tau_2)}\Big( |\nab_4\psi |^2+|\nab_3\psi|^2+|\nab\psi|^2+r^{-2} | \psi |^2\Big), \\
F_{\Si_*}[\psi](\tau_1,\tau_2) := & \int_{\Si_*(\tau_1, \tau_2)}\Big( |\nab_4\psi |^2+|\nab_3\psi|^2+|\nab\psi|^2+ r^{-2} | \psi |^2\Big).
\end{split}
\eea

{\bf 4. Basic  $N$- norm.} 
\bea
\bsplit
\NN[\psi,  N](\tau_1, \tau_2) :=&\int_{\MM(\tau_1, \tau_2)}\big(|\nab_{\Rhat} \psi|+r^{-1}|\psi|\big) |N|+\left|\int_{\MM} \nab_{\That_\de} \psi \c N\right| \\
& +\int_{\MM(\tau_1, \tau_2)\cap\left\{r\leq r_+(1+2\de_{red})\right\}}|\dk\psi||N|\\
&+\int_{\MM(\tau_1, \tau_2)}|N|^2+\sup_{\tau\in[\tau_1, \tau_2]}\int_{\Si(\tau)}|N|^2+\int_{\Si_*(\tau_1, \tau_2)}|N|^2.
\end{split}
\eea

{\bf 5. Weighted  bulk norm.} For $0<p<2$, we define
\bea
B_{p}[\psi](\tau_1, \tau_2)&:=& \Morr[\psi](\tau_1, \tau_2) +\int_{\MM_{r\geq 4m}(\tau_1,\tau_2)} r^{p-3}\Big(|\dk \psi|^2 +|\psi|^2  \Big).
\eea

\medskip

\bigskip

{\bf 6. Higher order norms.}  We define the higher  derivative norms $\Mor^s[\psi]$, $E^s[\psi]$, $F^s[\psi]$,    $\NN^s[\psi,  N]$, $B_p^s[\psi]$,  by the  general procedure for a norm $Q[\psi]$, i.e. 
 \beaa
 Q^s[\psi]=\sum_{k\le s} Q[\dk^k\psi].
 \eeaa


\section{Basic energy-Morawetz in perturbations of Kerr}
\label{section:adaptiation-errors}


In this section, we revisit the proofs of Chapters \ref{chapter-proof-part1} and \ref{chapter-proof-mor-2} in Kerr and show that  the conclusions of Chapters \ref{chapter-proof-part1} and \ref{chapter-proof-mor-2}  also hold in perturbations of Kerr up to the addition suitable error terms, see the statements in section \ref{sec:statementofresultschapter7and8in Kerralsoholdinperturbations}.


 \subsection{Some basic commutations with $\T$ and $\Z$}


In the following lemma, we establish basic properties of $\Lieb_\T$ and $\Lieb_\Z$ that will be used repeatedly.
 \begin{lemma}\lab{lemma:basicpropertiesLiebTfasdiuhakdisug:chap9}
For a horizontal covariant k-tensor $U$, we have
\beaa
\nab_\T U_{b_1\cdots b_k} &=& \Lieb_\T U_{b_1\cdots b_k} +\frac{2amr\cos\th}{|q|^4}\sum_{j=1}^k\in_{b_jc} U_{b_1\cdots c\cdots b_k}+\Ga_b U,
\eeaa
\beaa
\nab_\Z U_{b_1\cdots b_k} &=& \Lieb_\Z U_{b_1\cdots b_k} -\frac{\cos\th((r^2+a^2)^2-a^2(\sin\th)^2\De)}{|q|^4}\sum_{j=1}^k\in_{b_jc} U_{b_1\cdots c\cdots b_k}+r\Ga_b U,
\eeaa
and
\beaa
[\Lieb_\T, \dk]U = \dk(\Ga_b U), \qquad [\Lieb_\Z, \dk]U = r\dk(\Ga_b U).
\eeaa
\end{lemma}

\begin{proof}
First,  we compute 
\beaa
2\g(\nab_b\T, e_c) &=& \g\left(\nab_b\left(e_4+\frac{\De}{|q|^2}e_3 -2a\Re(\Jk)_de_d\right), e_c\right)\\
&=& \chi_{bc}+\frac{\De}{|q|^2}\chib_{bc} -2a\nab_b\Re(\Jk)_c\\
&=& \frac{1}{2}\left(\trch +\frac{\De}{|q|^2}\trchb\right)\de_{bc}+\frac{1}{2}\left(\atrch +\frac{\De}{|q|^2}\atrchb\right)\in_{bc}\\
&& - a\div(\Re(\Jk))\de_{bc} - a\curl(\Re(\Jk))\in_{bc}+\Ga_b\\
&=& \left(\frac{2a\cos\th \De}{|q|^4}-\frac{2a(r^2+a^2)\cos\th}{|q|^4}\right)\in_{bc}+\Ga_b\\
&=& -\frac{4amr\cos\th}{|q|^4}\in_{bc}+\Ga_b.
\eeaa
Since we have
\beaa
\Lieb_\T U_{b_1\cdots b_k} &=& \nab_\T U_{b_1\cdots b_k} +\g(\D_{b_1}\T, e_c)U_{cb_2\cdots b_k}+\cdots,
\eeaa
we infer
\beaa
\nab_\T U_{b_1\cdots b_k} &=& \Lieb_\T U_{b_1\cdots b_k} +\frac{2amr\cos\th}{|q|^4}\sum_{j=1}^k\in_{b_jc} U_{b_1\cdots c\cdots b_k}+\Ga_b U
\eeaa
as stated.

Next, we compute 
\beaa
2\g(\nab_b\Z, e_c) &=& \g\left(\nab_b\left(2(r^2+a^2)\Re(\Jk)_de_d -a(\sin\th)^2e_4-\frac{a(\sin\th)^2\De}{|q|^2}e_3\right), e_c\right)\\
&=& 2(r^2+a^2)\nab_b\Re(\Jk)_c +4re_b(r)\Re(\Jk)_c -a(\sin\th)^2\left(\chi_{bc}+\frac{\De}{|q|^2}\chib_{bc}\right) \\
&=& (r^2+a^2)\div(\Re(\Jk))\de_{bc} +(r^2+a^2)\curl(\Re(\Jk))\in_{bc}\\
&&-\frac{a(\sin\th)^2}{2}\left(\trch +\frac{\De}{|q|^2}\trchb\right)\de_{bc}-\frac{a(\sin\th)^2}{2}\left(\atrch +\frac{\De}{|q|^2}\atrchb\right)\in_{bc}\\
&& +r\Ga_b\\
&=& \left(-\frac{2a^2(\sin\th)^2\cos\th \De}{|q|^4}+\frac{2(r^2+a^2)^2\cos\th}{|q|^4}\right)\in_{bc}+r\Ga_b\\
&=& \frac{2\cos\th((r^2+a^2)^2-a^2(\sin\th)^2\De)}{|q|^4}\in_{bc}+r\Ga_b.
\eeaa
Since we have
\beaa
\Lieb_\Z U_{b_1\cdots b_k} &=& \nab_\Z U_{b_1\cdots b_k} +\g(\D_{b_1}\Z, e_c)U_{cb_2\cdots b_k}+\cdots,
\eeaa
we infer
\beaa
\nab_\Z U_{b_1\cdots b_k} &=& \Lieb_\Z U_{b_1\cdots b_k} -\frac{\cos\th((r^2+a^2)^2-a^2(\sin\th)^2\De)}{|q|^4}\sum_{j=1}^k\in_{b_jc} U_{b_1\cdots c\cdots b_k}+r\Ga_b U
\eeaa
as stated.

Next, in view of Lemma \ref{lemma:commutation-Lieb-nab}, together with the fact that $\piT\in \Ga_b$, see Lemma \ref{LEMMA:DEFORMATION-TENSORS-T}, we have for a horizontal covariant k-tensor $U$
  \beaa
  \bsplit
  \nab_b(\Lieb_\T U_{A})-\Lieb_\T(\nab_b U_{A})&= r^{-1}\dk\Ga_b U,\\
   \nab_4(\Lieb_\T U_A)-\Lieb_\T(\nab_4 U_A) +\nab_{\Lieb_\T e_4} U_{A}&=r^{-1}\dk\Ga_b U,\\
   \nab_3(\Lieb_\T U_A)-\Lieb_\T(\nab_3 U_A) +\nab_{\Lieb_\T e_3} U_{A}&=\dk\Ga_b U.
  \end{split}
  \eeaa
Also, we have 
\beaa
2\Lie_\T e_4 &=& 2[\T, e_4]=\left[e_4+\frac{\De}{|q|^2}e_3 -2a\Re(\Jk)^de_d, e_4\right]\\
&=&  \frac{\De}{|q|^2}[e_3, e_4] -e_4\left(\frac{\De}{|q|^2}\right)e_3 -2a[\Re(\Jk)^de_d, e_4]
\eeaa
 and hence
\beaa
2\Lieb_\T e_4 &=& \frac{2\De}{|q|^2}(\widecheck{\eta}_b- \widecheck{\etab}_b)e_b -2a\Re(\widecheck{\nab_4\Jk})_be_b 
\eeaa 
which yields
\beaa
\nab_{\Lieb_\T e_4} &=& r^{-1}\Ga_b\dkb.
\eeaa
Similarly, we have
\beaa
\nab_{\Lieb_\T e_3} &=& r^{-1}\Ga_b\dkb
\eeaa 
and hence
  \beaa
  \bsplit
  \nab_b(\Lieb_\T U_{A})-\Lieb_\T(\nab_b U_{A})&= r^{-1}\dk\Ga_b U,\\
   \nab_4(\Lieb_\T U_A)-\Lieb_\T(\nab_4 U_A) &=r^{-1}\dk(\Ga_b U),\\
   \nab_3(\Lieb_\T U_A)-\Lieb_\T(\nab_3 U_A)  &=\dk(\Ga_b U).
  \end{split}
  \eeaa
Since $\T(r)\in r\Ga_b$, we deduce 
\beaa
[\Lieb_\T, \dk]U &=& \dk(\Ga_b U)
\eeaa
as stated.

Finally, in view of Lemma \ref{lemma:commutation-Lieb-nab}, together with the fact that $\piZ\in r\Ga_b$, see Lemma \ref{LEMMA:DEFORMATION-TENSORS-T}, we have for a horizontal covariant k-tensor $U$
  \beaa
  \bsplit
  \nab_b(\Lieb_\Z U_{A})-\Lieb_\Z(\nab_b U_{A})&= \dk\Ga_b U,\\
   \nab_4(\Lieb_\Z U_A)-\Lieb_\Z(\nab_4 U_A) +\nab_{\Lieb_\Z e_4} U_{A}&=\dk\Ga_b U,\\
   \nab_3(\Lieb_\Z U_A)-\Lieb_\Z(\nab_3 U_A) +\nab_{\Lieb_\Z e_3} U_{A}&=r\dk\Ga_b U.
  \end{split}
  \eeaa
Also, we have 
\beaa
2\Lie_\Z e_4 &=& 2[\Z, e_4]=\left[2(r^2+a^2)\Re(\Jk)_de_d -a(\sin\th)^2e_4-\frac{a(\sin\th)^2\De}{|q|^2}e_3, e_4\right]\\
&=& 2(r^2+a^2)[\Re(\Jk)_de_d, e_4] -4re_4(r)\Re(\Jk)_de_d +ae_4((\sin\th)^2)e_4\\
&& -\frac{a(\sin\th)^2\De}{|q|^2}[e_3, e_4]+e_4\left(\frac{a(\sin\th)^2\De}{|q|^2}\right)e_3
\eeaa
 and hence
\beaa
2\Lieb_\Z e_4 &=& 2(r^2+a^2)\Re(\widecheck{\nab_4\Jk})_be_b -4r\widecheck{e_4(r)}\Re(\Jk)_de_d -\frac{2a(\sin\th)^2\De}{|q|^2}(\widecheck{\eta}_b- \widecheck{\etab}_b)e_b  
\eeaa 
which yields
\beaa
\nab_{\Lieb_\Z e_4} &=& r^{-1}\Ga_b\dkb.
\eeaa
Similarly, we have
\beaa
\nab_{\Lieb_\T e_3} &=& \Ga_b\dkb
\eeaa 
and hence
  \beaa
  \bsplit
  \nab_b(\Lieb_\Z U_{A})-\Lieb_\Z(\nab_b U_{A})&= \dk\Ga_b U,\\
   \nab_4(\Lieb_\Z U_A)-\Lieb_\Z(\nab_4 U_A) &=\dk(\Ga_b U),\\
   \nab_3(\Lieb_\Z U_A)-\Lieb_\Z(\nab_3 U_A)  &=r\dk(\Ga_b U).
  \end{split}
  \eeaa
Since $\Z(r)\in r^2\Ga_g$, we deduce 
\beaa
[\Lieb_\Z, \dk]U &=& r\dk(\Ga_b U)
\eeaa
as stated. This concludes the proof of Lemma \ref{lemma:basicpropertiesLiebTfasdiuhakdisug:chap9}.
\end{proof}

\begin{corollary}\lab{cor:basicpropertiesLiebTfasdiuhakdisug:chap9}
We have the following commutator identities 
\beaa
\,[\nab_\T, \nab_\Z]\psi &=& \dk(\Ga_b\c\psi),\\
\,[\nab_\T, \dkb]\psi &=& \frac{2amr\cos\th}{|q|^4}r\dual\nab_b\psi+O(ar^{-3})\psi+\dk(\Ga_b\c\psi),\\
\,[\nab_\Z, \dkb]\psi &=& -\frac{2\cos\th((r^2+a^2)^2-a^2(\sin\th)^2\De)}{|q|^4}r\dual\nab_b\psi+O(1)\psi+r\dk(\Ga_b\c\psi).
\eeaa
\end{corollary}

\begin{proof}
In view of Lemma \ref{lemma:basicpropertiesLiebTfasdiuhakdisug:chap9}, we have
\beaa
[\nab_\T, \nab_\Z]\psi &=& [\Lieb_\T, \dk]\psi +\Z\left(\frac{4amr\cos\th}{|q|^4}\right)\dual\psi+\dk(\Ga_b\c\psi)\\
&=&  O(r^{-3})(r^{-1}\Z(r), \Z(\cos\th))\psi+\dk(\Ga_b\c\psi)
\eeaa
and hence
\beaa
[\nab_\T, \nab_\Z]\psi &=& \dk(\Ga_b\c\psi).
\eeaa

Using again Lemma \ref{lemma:basicpropertiesLiebTfasdiuhakdisug:chap9}, we have 
\beaa
[\nab_\T, \dkb]\psi &=& [\Lieb_\T, \dk]\psi +\frac{2amr\cos\th}{|q|^4}r\dual\nab\psi +r\nab\left(\frac{4amr\cos\th}{|q|^4}\right)\dual\psi+\dk(\Ga_b\c\psi)\\
&=& \frac{2amr\cos\th}{|q|^4}r\dual\nab\psi +r\nab\left(\frac{4amr\cos\th}{|q|^4}\right)\dual\psi +\dk(\Ga_b\c\psi)
\eeaa
and hence, since $\nab(r)\in r\Ga_g$ and $\widecheck{\nab(\cos\th)}\in\Ga_b$, we infer
\beaa
[\nab_\T, \dkb]\psi &=& \frac{2amr\cos\th}{|q|^4}r\dual\nab\psi +O(ar^{-3})\psi+\dk(\Ga_b\c\psi).
\eeaa

Finally, using again Lemma \ref{lemma:basicpropertiesLiebTfasdiuhakdisug:chap9}, we have 
\beaa
[\nab_\Z, \dkb]\psi &=& [\Lieb_\Z, \dk]\psi -\frac{2\cos\th((r^2+a^2)^2-a^2(\sin\th)^2\De)}{|q|^4}r\dual\nab\psi\\
&& +r\nab\left(-\frac{2\cos\th((r^2+a^2)^2-a^2(\sin\th)^2\De)}{|q|^4}\right)\dual\psi+r\dk(\Ga_b\c\psi)\\
&=& -\frac{2\cos\th((r^2+a^2)^2-a^2(\sin\th)^2\De)}{|q|^4}r\dual\nab\psi \\
&& -r\nab\left(\frac{2\cos\th((r^2+a^2)^2-a^2(\sin\th)^2\De)}{|q|^4}\right)\dual\psi +r\dk(\Ga_b\c\psi)
\eeaa
and hence, since $\nab(r)\in r\Ga_g$ and $\widecheck{\nab(\cos\th)}\in\Ga_b$, we infer
\beaa
[\nab_\Z, \dkb]\psi &=& -\frac{2\cos\th((r^2+a^2)^2-a^2(\sin\th)^2\De)}{|q|^4}r\dual\nab\psi+O(1)\psi+r\dk(\Ga_b\c\psi).
\eeaa
This concludes the proof of Corollary \ref{cor:basicpropertiesLiebTfasdiuhakdisug:chap9}.
\end{proof}


 \subsection{Approximate symmetry operators}


Recall   the set of second order differential operators  $\SS_\aund$,   see  Definition \ref{definition:symmetry-tensors-pert},  
\beaa
\SS_1\psi&=& \nab_T \nab_T \psi, \\
\SS_2\psi&=& a \nab_T \nab_Z \psi,\\
\SS_3\psi&=& a^2 \nab_Z \nab_Z \psi,\\
\SS_4\psi&=& \OO (\psi),
\eeaa  
where $\OO$ is given by \eqref{definition-SS}, i.e.
\beaa
\OO(\psi) &=& |q|^2 \left(\lap_k \psi + \frac{2a^2\cos\th}{|q|^2} \dual\Re(\Jk)^b \nab_b \psi   \right).
\eeaa

These symmetry operators satisfy the following commutation properties.
\begin{lemma}\lab{lemma:commutationpropertiesofthesymmetryoperatorsafsoiudf:chap9}
We have
\beaa
\,[\SS_1, \SS_2], \,[\SS_1, \SS_3],\,[\SS_2, \SS_3]=\dk^2(\Ga_b\c\psi).
\eeaa
Also, we have
\beaa
\,[\nab_\T, \OO] &=& O(ar^{-3})\dkb^{\leq 1}\psi+ \dk^{\leq 2}(\Ga_b\psi),\\
\,[\nab_\Z, \OO] &=& O(1)\dkb^{\leq 1}\psi+ r\dk^{\leq 2}(\Ga_b\psi),
\eeaa
and 
\beaa
\,[\SS_1, \OO] &=& O(ar^{-3})\dk^{\leq 2}\psi+ \dk^{\leq 3}(\Ga_b\psi),\\
\,[\SS_2, \OO] &=& O(a)\dkb^{\leq 2}\psi+ r\dk^{\leq 3}(\Ga_b\psi),\\
\,[\SS_3, \OO] &=& O(a^2)\dkb^{\leq 2}\psi+ r\dk^{\leq 3}(\Ga_b\psi).
\eeaa
\end{lemma}

\begin{proof}
First, since $[\nab_\T, \nab_\Z]\psi = \dk(\Ga_b\c\psi)$ in view of Corollary \ref{cor:basicpropertiesLiebTfasdiuhakdisug:chap9}, we immediately have in view of the definition of $\SS_1$, $\SS_2$ and $\SS_3$,   
\beaa
\,[\SS_1, \SS_2], \,[\SS_1, \SS_3],\,[\SS_2, \SS_3]=\dk^2(\Ga_b\c\psi)
\eeaa
as stated.

Next, using Lemma \ref{lemma:basicpropertiesLiebTfasdiuhakdisug:chap9}  and the definition of $\OO$, we have
\beaa
[\nab_\T, \OO] &=& [\Lieb_\T, \OO] +\left[\frac{4amr\cos\th}{|q|^4}, \OO\right]\dual\psi +\dk^{\leq 2}(\Ga_b\psi)\\
&=& \left[\frac{4amr\cos\th}{|q|^4}, \OO\right]\dual\psi +\dk^{\leq 2}(\Ga_b\psi)\\
&=& O(ar^{-3})\dkb^{\leq 1}\psi+ \dk^{\leq 2}(\Ga_b\psi)
\eeaa
and
\beaa
[\nab_\Z, \OO] &=& [\Lieb_\Z, \OO] -\left[\frac{\cos\th((r^2+a^2)^2-a^2(\sin\th)^2\De)}{|q|^4}, \OO\right]\dual\psi +r\dk^{\leq 2}(\Ga_b\psi)\\
&=& -\left[\frac{\cos\th((r^2+a^2)^2-a^2(\sin\th)^2\De)}{|q|^4}, \OO\right]\dual\psi +r\dk^{\leq 2}(\Ga_b\psi)\\
&=& O(1)\dkb^{\leq 1}\psi+ r\dk^{\leq 2}(\Ga_b\psi)
\eeaa
as stated. 

Finally, the identities for $[\SS_1, \OO]$, $[\SS_2, \OO]$ and $[\SS_3, \OO]$ following immediately from the definition of $\SS_1$, $\SS_2$ and $\SS_3$, and the above commutator identities for $[\nab_\T, \OO]$ and $[\nab_\Z, \OO]$. This concludes the proof of the lemma.
\end{proof}

We have the following  analog of Lemma \ref{lemma-commuted-equ-Naund} in perturbations of Kerr.

 \begin{lemma}\label{lemma-commuted-equ-Naund:perturbation}
Given a $\mathfrak{s}_2$ tensor  $\psi$ solution  of the equation  \eqref{eq:Gen.RW}.  Then the commuted $\mathfrak{s}_2$ tensor  $\psia:=\SS_\aund \psi$ satisfies
  \bea
  \lab{eq:waveeqfor-psia:perturbation}
  \squared_2\psia -V\psia&=&- \frac{ 4 a\cos\th}{|q|^2} \dual \nab_\T  \psia+N_\aund
  \eea
  where $N_\aund$ satisfies
\bea\lab{eq:estimaeforMundinperturbationsofKerr}
|N_{\aund}| &\les& |\dk^{\leq 2}N|+|a|r^{-2}|\dk^{\leq 1}\nab_3\psi|+|a|r^{-3}|\dk^{\leq 2}\psi|+|\dk^{\leq 3}(\Ga_g\c\psi)|.
\eea
\end{lemma}

\begin{proof}
Since $\psi$ satisfies \eqref{eq:Gen.RW}, i.e. 
  \beaa
  \squared_2\psi -V\psi &=&- \frac{4 a\cos\th}{|q|^2} \dual \nab_T  \psi+N,
  \eeaa
we infer
 \beaa
  \squared_2\psia -V\psia &=&- \frac{4 a\cos\th}{|q|^2} \dual \nab_T  \psia+N_\aund,\\
  N_\aund &:=& -[\SS_\aund, \squared_2]\psi +[\SS_\aund, V]\psi -\dual\left[\SS_\aund, \frac{4 a\cos\th}{|q|^2} \nab_T\right]\psi+\dk^{\leq 2}N, \qquad \aund=1, 2, 3,\\
  N_4 &:=& -\frac{1}{|q|^2}[\SS_4, |q|^2\squared_2]\psi +\frac{1}{|q|^2}[\SS_4, |q|^2V]\psi -\frac{1}{|q|^2}\dual[\SS_4, 4a\cos\th\nab_T]\psi+\dk^{\leq 2}N,
  \eeaa
and hence
\beaa
|N_\aund| &\les& |\dk^{\leq 2}N|+|[\SS_\aund, \squared_2]\psi| +|[\SS_\aund, V]\psi| +a\left|\left[\SS_\aund, \frac{\cos\th}{|q|^2}\nab_T\right]\psi\right|, \qquad \aund=1, 2, 3,\\
  |N_4| &\les& |\dk^{\leq 2}N|+ r^{-2}|[\SS_4, |q|^2\squared_2]\psi| +r^{-2}|[\SS_4, |q|^2V]\psi| +ar^{-2}|[\SS_4, \cos\th\nab_T]\psi|.
\eeaa

Next, since 
\beaa
V= \frac{4\De}{ (r^2+a^2) |q|^2}, \qquad |q|^2V= \frac{4\De}{ (r^2+a^2)},
\eeaa
see \eqref{eq:Gen.RW}, and since $\T(r)\in r\Ga_b$, $\T(\cos\th)\in \Ga_b$, $\Z(r)\in r^2\Ga_g$, $\Z(\cos\th)\in r\Ga_b$, and $\nab(r)\in r\Ga_g$, we infer
\beaa
|[\SS_\aund, V]\psi| &\les& r^{-3}|\dk^{\leq 1}((\T ,\Z)(r, \cos\th))||\dk^{\leq 1}\psi|\les r^{-1}|\dk^{\leq 1}(\Ga_g\c\psi)|,  \qquad \aund=1, 2, 3,\\
r^{-2}|[\SS_4, |q|^2V]\psi| &\les& r^{-2}|\dk^{\leq 1}(\nab(r))||\dk^{\leq 1}\psi|\les r^{-1}|\dk^{\leq 1}(\Ga_g\c\psi)|,
\eeaa
and hence
\beaa
|N_\aund| &\les& |\dk^{\leq 2}N|+r^{-1}|\dk^{\leq 1}(\Ga_g\c\psi)|+|[\SS_\aund, \squared_2]\psi|  +a\left|\left[\SS_\aund, \frac{\cos\th}{|q|^2}\nab_T\right]\psi\right|, \qquad \aund=1, 2, 3,\\
  |N_4| &\les& |\dk^{\leq 2}N|+ r^{-1}|\dk^{\leq 1}(\Ga_g\c\psi)|+ r^{-2}|[\SS_4, |q|^2\squared_2]\psi| +ar^{-2}|[\SS_4, \cos\th\nab_T]\psi|.
\eeaa

Also, we have
\beaa
\left|\left[\SS_\aund, \frac{\cos\th}{|q|^2}\nab_T\right]\psi\right| &\les& r^{-2}|\dk^{\leq 1}((\T ,\Z)(r, \cos\th))||\dk^{\leq 2}\psi|+r^{-2}|[\SS_\aund, \nab_T]\psi|\\
&\les& |\dk^{\leq 2}(\Ga_g\c\psi)|+r^{-2}|[\SS_\aund, \nab_T]\psi|\qquad \aund=1, 2, 3,
\eeaa
and 
\beaa
r^{-2}|[\SS_4, \cos\th\nab_T]\psi| &\les&  r^{-1}|\dk^{\leq 1}\nab(\cos\th)||\dk^{\leq 1}\nab_T\psi|+r^{-2}|[\SS_4, \nab_T]\psi| \\
&\les&  r^{-1}(r^{-1}+|\dk^{\leq 1}\Ga_b|)|\dk^{\leq 1}\nab_T\psi|+r^{-2}|[\SS_4, \nab_T]\psi| \\
&\les&  r^{-2}|\dk^{\leq 1}\nab_3\psi|+r^{-3}|\dk^{\leq 2}\psi|+|\dk^{\leq 2}(\Ga_g\c\psi)|+r^{-2}|[\SS_4, \nab_T]\psi| 
\eeaa
and hence
\beaa
|N_\aund| &\les& |\dk^{\leq 2}N|+|\dk^{\leq 2}(\Ga_g\c\psi)|+|[\SS_\aund, \squared_2]\psi|  +ar^{-2}\left|\left[\SS_\aund, \nab_T\right]\psi\right|, \qquad \aund=1, 2, 3,\\
  |N_4| &\les& |\dk^{\leq 2}N|+ ar^{-2}|\dk^{\leq 1}\nab_3\psi|+ar^{-3}|\dk^{\leq 2}\psi|+ |\dk^{\leq 2}(\Ga_g\c\psi)|\\
  &&+ r^{-2}|[\SS_4, |q|^2\squared_2]\psi| +ar^{-2}|[\SS_4, \nab_T]\psi|.
\eeaa

Next, we use Lemma \ref{lemma:basicpropertiesLiebTfasdiuhakdisug:chap9} which implies
\beaa
|[\SS_\aund, \nab_T]\psi| &\les& |[\SS_\aund, \Lieb_T]\psi|+ar^{-2}|\dk^{\leq 2}\psi|+|\dk^{\leq 2}(\Ga_b\c\psi)|\\
&\les& ar^{-2}|\dk^{\leq 2}\psi|+|\dk^{\leq 2}(\Ga_b\c\psi)|, \aund=1, 2, 3,4,
\eeaa
and hence 
\beaa
|N_\aund| &\les& |\dk^{\leq 2}N|+a^2r^{-4}|\dk^{\leq 2}\psi|+|\dk^{\leq 2}(\Ga_g\c\psi)|+|[\SS_\aund, \squared_2]\psi|, \qquad \aund=1, 2, 3,\\
  |N_4| &\les& |\dk^{\leq 2}N|+ ar^{-2}|\dk^{\leq 1}\nab_3\psi|+ar^{-3}|\dk^{\leq 2}\psi|+ |\dk^{\leq 2}(\Ga_g\c\psi)|+ r^{-2}|[\SS_4, |q|^2\squared_2]\psi|.
\eeaa

Finally, recall from Proposition \ref{prop:commutators-SSaund-squared2-pert} that the following  commutation  formulas hold true  for $\psi\in \sk_2$:
\beaa
\, [\SS_1, \squared_2] \psi &=& O(ar^{-2})\dk^{\leq 2}\psi +\dk^2 \big(\Ga_g \c \dk \psi\big)+\dk\big( \Ga_b \c \squared_2\psi  \big), \\
\, [\SS_2, \squared_2] \psi&=& O(ar^{-2})\dk^{\leq 2}\psi +\dk^2 \big(\Ga_g \c \dk \psi\big)+r\dk \big(\Ga_b \c \squared_2\psi\big),\\
\, [\SS_3, \squared_2] \psi&=& O(ar^{-2})\dk^{\leq 2}\psi +\dk^2 \big(\Ga_g \c \dk \psi\big)+r\dk \big(\Ga_b \c \squared_2\psi\big),
\eeaa
and
\beaa
\, [\SS_4, |q|^2\squared_2] \psi&=& |q|^2\Big[ O( a r^{-2}) \  \dk^{\leq 2} \psi+  \dk^2 \big( \Ga_g \c \dk \psi)+ \Ddot_3 \dk \big(|q|^2 \xi \c \Ddot_a \psi \big)\Big].
\eeaa
We deduce 
\beaa
|N_{\aund}| &\les& |\dk^{\leq 2}N|+|a|r^{-2}|\dk^{\leq 1}\nab_3\psi|+|a|r^{-3}|\dk^{\leq 2}\psi|+|\dk^{\leq 3}(\Ga_g\c\psi)|+r|\dk \big(\Ga_b \c \squared_2\psi\big)|\\
&&+r|\dk^{\leq 3}(\xi\c\psi)|.
\eeaa
Next, recall from Remark \ref{rmk:strongerassumptionxiinrm1Gagforchapter9} that we assume $\xi\in r^{-1}\Ga_g$ in this chapter which implies
\beaa
|N_{\aund}| &\les& |\dk^{\leq 2}N|+|a|r^{-2}|\dk^{\leq 1}\nab_3\psi|+|a|r^{-3}|\dk^{\leq 2}\psi|+|\dk^{\leq 3}(\Ga_g\c\psi)|+r|\dk \big(\Ga_b \c \squared_2\psi\big)|.
\eeaa
Plugging the model RW equation \eqref{eq:Gen.RW} for $\psi$ in the last term, we have
\beaa
r|\dk \big(\Ga_b \c \squared_2\psi\big)| &\les& r\left|\dk \left(\Ga_b \c\left(-V\psi - \frac{4 a\cos\th}{|q|^2} \dual \nab_T  \psi+N\right)\right)\right| \\
&\les& r^{-1}|\dk^{\leq 2}(\Ga_b\c\psi)|+|\dk^{\leq 1}N|
\eeaa
and hence
\beaa
|N_{\aund}| &\les& |\dk^{\leq 2}N|+|a|r^{-2}|\dk^{\leq 1}\nab_3\psi|+|a|r^{-3}|\dk^{\leq 2}\psi|+|\dk^{\leq 3}(\Ga_g\c\psi)|
\eeaa
as stated. This concludes the proof of Lemma \ref{lemma-commuted-equ-Naund:perturbation}.
\end{proof}


 \subsection{Commutation with $\nab_{\Rhat}$}


Recall that $\Rhat$ is given by 
\beaa
\Rhat &=& \frac 1 2 \left( \frac{|q|^2}{r^2+a^2} e_4-\frac{\De}{r^2+a^2}  e_3\right)= \frac{1}{2}\frac{|q|^2}{r^2+a^2}X, \qquad X=e_4-\frac{\De}{|q|^2}e_3.
\eeaa
We use the following commutation lemma for vectorfields $X$ spanned by $e_3$ and $e_4$.
\begin{lemma}\lab{lemma:commutationofXwithnabinKerrwhenXspannedbye3e4only}
Let $X$ such that $X=X^4e_4+X^3e_3$. Then, we have
\beaa
\,[\nab_b, \nab_X]\psi &=&  \frac{1}{2}(X^4\trch+X^3\trchb)\nab_b\psi+\frac{1}{2}(X^4\atrch+X^3\atrchb)\dual\nab_b\psi \\
&&+(e_b(X^4)-X^4(\etab_b+\eta_b))\nab_4\psi +e_b(X^3)\nab_3\psi \\
&&+O(ar^{-3})(X^4, X^3)\psi +(|X^3|+|X^4|)\big(\Ga_b\nab_3\psi+r^{-1}\Ga_b\dk^{\leq 1}\psi\big).
\eeaa
\end{lemma}

\begin{proof}
In view of Corollary \ref{cor:comm-gen-B} and our definition of $\Ga_b$, $\Ga_g$, we have 
\beaa
\,[\nab_3, \nab_b]\psi &=& -\frac 1 2 \big( \trchb \nab_b\psi +\atrchb \dual \nab_b\psi\big)  + O(ar^{-3})\psi+\Ga_b\nab_3\psi+r^{-1}\Ga_b\dk^{\leq 1}\psi
\eeaa
and
\beaa
\,[\nab_4, \nab_b]\psi &=& -\frac 1 2 \big( \trch \nab_b\psi +\atrch \dual \nab_b\psi\big)+( \etab_b+\ze_b) \nab_4\psi\\
&&  + O(ar^{-3})\psi+ \Ga_g\nab_3\psi+r^{-1}\Ga_g\dk^{\leq 1}\psi.
\eeaa
We infer
\beaa
(\nab_b\nab_X -\nab_X\nab_b)\psi &=& e_b(X^4)\nab_4\psi+e_b(X^3)\nab_3\psi+X^4[\nab_b, \nab_4]\psi+X^3[\nab_b, \nab_3]\psi\\
&=& X^4\left(\frac{1}{2}\trch\nab_b+\frac{1}{2}\atrch\dual\nab_b -(\etab_b+\eta_b)\nab_4+O(ar^{-3})\right)\psi\\
&&+X^3\left(\frac{1}{2}\trchb\nab_b+\frac{1}{2}\atrchb\dual\nab_b +O(ar^{-3})\right)\psi \\
&&+e_b(X^4)\nab_4\psi+e_b(X^3)\nab_3\psi +(|X^3|+|X^4|)\big(\Ga_b\nab_3\psi+r^{-1}\Ga_b\dk^{\leq 1}\psi\big)\\
&=& \frac{1}{2}(X^4\trch+X^3\trchb)\nab_b\psi+\frac{1}{2}(X^4\atrch+X^3\atrchb)\dual\nab_b\psi \\
&&+(e_b(X^4)-X^4(\etab_b+\eta_b))\nab_4\psi +e_b(X^3)\nab_3\psi \\
&&+O(ar^{-3})(X^4, X^3)\psi +(|X^3|+|X^4|)\big(\Ga_b\nab_3\psi+r^{-1}\Ga_b\dk^{\leq 1}\psi\big)
\eeaa
as stated.
\end{proof}

\begin{corollary}\lab{cor:commutatornabRhatwith|q|nab}
We have 
\beaa
\,[|q|\nab_b, \nab_\Rhat]\psi &=&  O(ar^{-2})\psi +r\Ga_b\dk^{\leq 1}\psi.
\eeaa
In particular, we have
\beaa
[\OO, \nab_\Rhat]\psi &=& O(ar^{-2})\dk^{\leq 1}\psi +r\Ga_b\dk^{\leq 2}\psi.
\eeaa
\end{corollary}

\begin{proof}
We apply Lemma \ref{lemma:commutationofXwithnabinKerrwhenXspannedbye3e4only} with $X^4=1$ and $X^3=-\frac{\De}{|q|^2}$. This yields
\beaa
(\nab_b\nab_X -\nab_X\nab_b)\psi &=& \frac{1}{2}\left(\trch -\frac{\De}{|q|^2}\trchb\right)\nab_b\psi+\frac{1}{2}\left(\atrch -\frac{\De}{|q|^2}\atrchb\right)\dual\nab_b\psi \\
&&-(\etab_b+\eta_b)\nab_4\psi -e_b\left(\frac{\De}{|q|^2}\right)\nab_3\psi +O(ar^{-3})\psi \\
&&+\Ga_b\nab_3\psi+r^{-1}\Ga_b\dk^{\leq 1}\psi.
\eeaa
Since 
\beaa
\tr X -\frac{\De}{|q|^2}\tr\Xb &=& \frac{\De}{|q|^2}\left(\frac{2}{q}+\frac{2}{\ov{q}}\right)+\Ga_g=\frac{4r\De}{|q|^4}+\Ga_g,
\eeaa
\beaa
e_b\left(\frac{\De}{|q|^2}\right) &=& \pr_r\left(\frac{\De}{|q|^2}\right)e_b(r)+\pr_{\cos\th}\left(\frac{\De}{|q|^2}\right)e_b(\cos\th)\\
&=& O(r^{-2})e_b(r) -\frac{2a^2\cos\th\De}{|q|^4}e_b(\cos\th)=-\frac{2a^2\cos\th\De}{|q|^4}\Re(i\Jk)_b+r^{-1}\Ga_g,
\eeaa
and
\beaa
H+\Hb &=& \frac{a(q-\ov{q})}{|q|^2}\Jk+\Ga_b=\frac{2a^2\cos\th}{|q|^2}i\Jk+\Ga_b,
\eeaa
we infer, using also $\dual\Jk=-i\Jk$, 
\beaa
(\nab_b\nab_X -\nab_X\nab_b)\psi = \frac{2r\De}{|q|^4}\nab_b\psi +\frac{2a^2\cos\th}{|q|^2}\dual\Re(\Jk)_b\nab_X +O(ar^{-3})\psi +\Ga_b\nab_3\psi+r^{-1}\Ga_b\dk^{\leq 1}\psi.
\eeaa
Also, note that 
\beaa
X(|q|) &=& \frac{1}{2}\frac{X(|q|^2)}{|q|}=\frac{rX(r)}{|q|}+O(r^{-1})X(\cos\th)=\frac{r}{|q|}\left(e_4(r)-\frac{\De}{|q|^2}e_3(r)\right)+r^{-1}\Ga_b\\
&=& \frac{2r\De}{|q|^3}+r\Ga_b
\eeaa
and hence 
\beaa
\,[|q|\nab_b, \nab_X]\psi &=& |q|[\nab_b, \nab_X]\psi -\nab_X(|q|)\nab_b\\
&=& \frac{2a^2\cos\th}{|q|}\dual\Re(\Jk)_b\nab_X +O(ar^{-2})\psi +r\Ga_b\nab_3\psi+\Ga_b\dk^{\leq 1}\psi.
\eeaa
Now, notice that 
\beaa
\Rhat &=& \frac{|q|^2}{r^2+a^2}X
\eeaa
and 
\beaa
\nab_b\left(\frac{|q|^2}{r^2+a^2}\right) &=& \pr_r\left(\frac{|q|^2}{r^2+a^2}\right)e_b(r)+\pr_{\cos\th}\left(\frac{|q|^2}{r^2+a^2}\right)e_b(\cos\th)\\
&=& -\frac{2a^2\cos\th}{r^2+a^2}\dual\Re(\Jk)_b+r^{-2}\Ga_b.
\eeaa
We infer
\beaa
\,[|q|\nab_b, \nab_\Rhat]\psi &=& \left[|q|\nab_b, \frac{|q|^2}{r^2+a^2}X\right]\psi \\
&=& \frac{|q|^2}{r^2+a^2}[|q|\nab_b, \nab_X]\psi+|q|\nab_b\left(\frac{|q|^2}{r^2+a^2}\right)\nab_X\psi\\
&=& O(ar^{-2})\psi +r\Ga_b\nab_3\psi+\Ga_b\dk^{\leq 1}\psi
\eeaa
as stated. 
\end{proof}

Next, we compute the commutators $[\nab_{\Rhat}, \nab_\T]$ and $[\nab_{\Rhat}, \nab_\Z]$. 
\begin{lemma}\lab{lemma:commutationofnabRhatwithnabTandnabZ}
We have
\beaa
\,[\nab_\T, \nab_\Rhat]\psi &=& O(amr^{-4})\psi +\dk(\Ga_b\c\psi),\\
\,[\nab_\Z, \nab_\Rhat]\psi &=&  O(a^2r^{-3})\psi +r\dk(\Ga_b\c\psi),
\eeaa
and 
\beaa
[\nab_\That, \nab_\Rhat]\psi &=& O(ar^{-3})\nab_\Z\psi+O(amr^{-4})\psi +\dk(\Ga_b\c\psi).
\eeaa
In particular, we have
\beaa
\,[\nab_\Rhat, \SS_1] &=& O(amr^{-4})\dk^{\leq 1}\psi +\dk^{\leq 2}(\Ga_b\c\psi),\\
\,[\nab_\Rhat, \SS_2] &=& O(a^2r^{-3})\dk^{\leq 1}\psi +r\dk^{\leq 2}(\Ga_b\c\psi),\\
\,[\nab_\Rhat, \SS_3] &=&  O(a^4r^{-3})\dk^{\leq 1}\psi +r\dk^{\leq 2}(\Ga_b\c\psi).
\eeaa
\end{lemma}

\begin{proof}
In view of Lemma \ref{lemma:basicpropertiesLiebTfasdiuhakdisug:chap9}, we have
\beaa
[\nab_\T, \nab_\Rhat]\psi &=& [\Lieb_\T, \dk]\psi +\Rhat\left(\frac{4amr\cos\th}{|q|^4}\right)\dual\psi+\dk(\Ga_b\c\psi)\\
&=&  \Rhat\left(\frac{4amr\cos\th}{|q|^4}\right)\dual\psi +\dk(\Ga_b\c\psi)\\
&=& O(amr^{-4})\psi +\dk(\Ga_b\c\psi)
\eeaa
and 
\beaa
[\nab_\Z, \nab_\Rhat]\psi &=& [\Lieb_\Z, \dk]\psi +\Rhat\left(-\frac{2\cos\th((r^2+a^2)^2-a^2(\sin\th)^2\De)}{|q|^4}\right)\dual\psi+r\dk(\Ga_b\c\psi)\\
&=&  \Rhat\left(-\frac{2\cos\th((r^2+a^2)^2-a^2(\sin\th)^2\De)}{|q|^4}\right)\dual\psi +r\dk(\Ga_b\c\psi)\\
&=& O(1)\Rhat(\cos\th)\psi+O(r^{-3})\Rhat(r)\psi +r\dk(\Ga_b\c\psi)\\
&=& O(a^2r^{-3})\psi +r\dk(\Ga_b\c\psi)
\eeaa
as stated.

Also, we have
\beaa
[\nab_\That, \nab_\Rhat]\psi &=& \left[\nab_\T+\frac{a}{r^2+a^2}\nab_\Z, \nab_\Rhat\right]\psi\\
&=& [\nab_\T, \nab_\Rhat]\psi+O(ar^{-2})[\nab_\Z, \nab_\Rhat]\psi -\Rhat\left(\frac{a}{r^2+a^2}\right)\nab_\Z\psi\\
&=& O(ar^{-3})\nab_\Z\psi+O(amr^{-4})\psi +\dk(\Ga_b\c\psi)
\eeaa
as stated. 
\end{proof}


 \subsection{The modified $\widetilde{\OO}$ operator}


The commutation properties of the operator $\OO$ with $\squared_2$, see Proposition \ref{LEMMA:MOD-LAPLACIAN-PERT-KERR}, are not good enough to derive energy estimates for $\OO\psi$ with $\psi \in \sk_2$. In the lemma below, we introduce a modified operator $\widetilde{\OO}$ which enjoys better properties. 
\begin{lemma}\lab{lemma:theoperqatorwidetildeOOcommutingwellwtihRWmodel}
Let 
\bea\lab{eq:defintionoftheoperqatorwidetildeOOcommutingwellwtihRWmodel}
\widetilde{\OO}\psi := \OO\psi  +\frac{4a(r^2+a^2+|q|^2)\cos\th}{|q|^2}\nab_\T\dual\psi +\frac{4a^2\cos\th}{|q|^2}\nab_\Z\dual\psi.
\eea
Then, we have
\bea
\nn\frac{1}{|q|^2}\left[|q|^2\left(\squared_2-V+\frac{4a\cos\th}{|q|^2}\dual\nab_\T\right), \widetilde{\OO}\right]\psi &=& O(ar^{-2})\nab_{\Rhat}\dk^{\leq 1}\psi +O(ar^{-2})\dk^{\leq 1}\psi\\
&& +\dk^{\leq 3}(\Ga_g\c\psi) +\Ga_b \c \squared_2\psi.
\eea
\end{lemma}

\begin{proof}
Recall from Proposition \ref{LEMMA:MOD-LAPLACIAN-PERT-KERR} that we have  
\beaa
\frac{1}{|q|^2}[|q|^2 \squared_2, \OO]\psi  &=&   -\nab\left(\frac{8a(r^2+a^2)\cos\th}{|q|^2}\right)\c\nab\nab_\That\dual\psi +O( a r^{-2}) \nab_{\Rhat}^{\leq 1}\dk^{\leq 1} \psi\\
&&+  \dk^2 \big( \Ga_g \c \dk \psi)+ \Ddot_3 \dk \big(|q|^2 \xi \c \Ddot_a \psi \big).
\eeaa 
We infer, using the fact that $|q|^2V$ only depends on $r$ and satisfies $|q|^2V=4+O(r^{-1})$, as well as the computation of the commutator $[\nab_\T, \OO]\psi$ in Lemma  \ref{lemma:commutationpropertiesofthesymmetryoperatorsafsoiudf:chap9},
\beaa
&&\frac{1}{|q|^2}\left[|q|^2\left(\squared_2-V+\frac{4a\cos\th}{|q|^2}\dual\nab_\T\right), \OO\right]\psi \\
&=& \frac{1}{|q|^2}[|q|^2\squared_2, \OO]\psi+\frac{4}{|q|^2}[a\cos\th\dual\nab_\T, \OO]\psi +O(r^{-1})\dkb^{\leq 1}\nab(r)\c\dkb^{\leq 1}\psi\\
&=& -\nab\left(\frac{8a(r^2+a^2)\cos\th}{|q|^2}\right)\c\nab\nab_\That\dual\psi - 8a\nab(\cos\th)\c\nab\nab_\T\dual\psi\\
&& +O(a^2r^{-2})\nab_{\Rhat}\dk^{\leq 1}\psi +O(ar^{-2})\dk^{\leq 1}\psi\\
&&+\dk^2 \big( \Ga_g \c \dk \psi)+ \Ddot_3 \dk \big(|q|^2 \xi \c \Ddot_a \psi \big)+r^{-2}\dk^{\leq 2}(\Ga_b\psi)+\dk^{\leq 1}\Ga_g\dk^{\leq 1}\psi.
\eeaa
Since $\xi\in r^{-1}\Ga_g$, see Remark \ref{rmk:strongerassumptionxiinrm1Gagforchapter9}, we deduce
\beaa
&&\frac{1}{|q|^2}\left[|q|^2\left(\squared_2-V+\frac{4a\cos\th}{|q|^2}\dual\nab_\T\right), \OO\right]\psi \\
&=& -\nab\left(\frac{8a(r^2+a^2)\cos\th}{|q|^2}\right)\c\nab\nab_\That\dual\psi - 8a\nab(\cos\th)\c\nab\nab_\T\dual\psi\\
&& +O(a^2r^{-2})\nab_{\Rhat}\dk^{\leq 1}\psi +O(ar^{-2})\dk^{\leq 1}\psi+\dk^{\leq 3}(\Ga_g\c\psi)
\eeaa
and hence, since $\That=\T+\frac{a}{r^2+a^2}\Z$, 
\beaa
&&\frac{1}{|q|^2}\left[|q|^2\left(\squared_2-V+\frac{4a\cos\th}{|q|^2}\dual\nab_\T\right), \OO\right]\psi \\
&=& -\nab\left(\frac{8a(r^2+a^2+|q|^2)\cos\th}{|q|^2}\right)\c\nab\nab_\T\dual\psi  -\nab\left(\frac{8a^2\cos\th}{|q|^2}\right)\c\nab\nab_\Z\dual\psi\\ 
&& +O(ar^{-2})\nab_{\Rhat}\dk^{\leq 1}\psi +O(ar^{-2})\dk^{\leq 1}\psi +\dk^{\leq 3}(\Ga_g\c\psi).
\eeaa

Next, recall from Proposition \ref{LE:COMMTZSQUARE} that we have
\beaa
\,[ \nab_\T, \squared_2] \psi &=& O(ar^{-4})\dk^{\leq 1}\psi +\dk \big(\Ga_g \c \dk \psi\big)+\Ga_b \c \squared_2\psi.
\eeaa
Given a smooth function $f_1(r, \cos\th)$ such that $f_1=O(a)$, $\pr_rf_1=O(ar^{-3})$ and $\pr_{\cos\th}f_1=O(a)$, we deduce
\beaa
\,[f_1(r, \cos\th)\nab_\T, \squared_2]\psi &=& f_1(r, \cos\th)\,[\nab_\T, \squared_2]\psi+[f_1(r, \cos\th), \squared_2]\nab_\T\psi\\
&=& -2\g^{\a\b}\nab_\a(f_1)\nab_\b\nab_\T\psi +O\Big(|\square_\g(f_1)|+ar^{-4}\Big)\dk^{\leq 1}\psi\\
&&+\dk \big(\Ga_g \c \dk \psi\big)+\Ga_b \c \squared_2\psi\\
&=& -2\pr_rf_1(r, \cos\th)\g^{\a\b}\nab_\a(r)\nab_\b\nab_\T\psi \\
&& -2\pr_{\cos\th}f_1(r, \cos\th)\g^{\a\b}\nab_\a(\cos\th)\nab_\b\nab_\T\psi \\
&&+O(ar^{-2})\dk^{\leq 1}\psi +\dk \big(\Ga_g \c \dk \psi\big)+\Ga_b \c \squared_2\psi\\
&=&  -2\pr_{\cos\th}f_1(r, \cos\th)\nab(\cos\th)\c\nab\nab_\T\psi \\
&&+O(\pr_rf_1)\nab_{\Rhat}\nab_\T\psi+O(ar^{-2})\dk^{\leq 1}\psi +\dk \big(\Ga_g \c \dk \psi\big)+\Ga_b \c \squared_2\psi\\
&=&  -2\nab(f_1(r, \cos\th))\c\nab\nab_\T\psi \\
&&+O(ar^{-3})\nab_{\Rhat}\dk^{\leq 1}\psi +O(ar^{-2})\dk^{\leq 1}\psi +\dk \big(\Ga_g \c \dk \psi\big)+\Ga_b \c \squared_2\psi.
\eeaa

Also, recall from Proposition \ref{LE:COMMTZSQUARE} that we have
\beaa
   \,[ \nab_\Z, \squared_2] \psi &=& O(r^{-2})\dk^{\leq 1}\psi +\dk \big(\Ga_g \c \dk \psi\big)+r\Ga_b \c \squared_2\psi.
 \eeaa
Given a smooth function $f_2(r, \cos\th)$ such that $f_2=O(ar^{-2})$, $\pr_rf_2=O(ar^{-3})$ and $\pr_{\cos\th}f_2=O(ar^{-2})$,  we deduce similarly 
\beaa
\,[f_2(r, \cos\th)\nab_\Z, \squared_2]\psi &=& f_2(r, \cos\th)\,[\nab_\Z, \squared_2]\psi+[f_2(r, \cos\th), \squared_2]\nab_\Z\psi\\
&=&  -2\nab(f_2(r, \cos\th))\c\nab\nab_\Z\psi \\
&&+O(ar^{-3})\nab_{\Rhat}\dk^{\leq 1}\psi +O(ar^{-2})\dk^{\leq 1}\psi +\dk \big(\Ga_g \c \dk \psi\big)+\Ga_b \c \squared_2\psi.
\eeaa
Setting
\beaa
f_1(r, \cos\th)=\frac{4a(r^2+a^2+|q|^2)\cos\th}{|q|^2}, \qquad f_2(r, \cos\th)=\frac{4a^2\cos\th}{|q|^2},
\eeaa
we infer
\beaa
&& \,\left[\squared_2, \frac{4a(r^2+a^2+|q|^2)\cos\th}{|q|^2}\nab_\T +\frac{4a^2\cos\th}{|q|^2}\nab_\Z\right]\psi\\
&=&  \nab\left(\frac{8a(r^2+a^2+|q|^2)\cos\th}{|q|^2}\right)\c\nab\nab_\T\psi  +\nab\left(\frac{8a^2\cos\th}{|q|^2}\right)\c\nab\nab_\Z\psi\\
&&+O(ar^{-3})\nab_\Rhat\dk^{\leq 1}\psi+O(ar^{-2})\dk^{\leq 1}\psi +\dk \big(\Ga_g \c \dk \psi\big)+\Ga_b \c \squared_2\psi
\eeaa
and hence
\beaa
&& \,\frac{1}{|q|^2}\left[|q|^2\left(\squared_2-V+\frac{4a\cos\th}{|q|^2}\dual\nab_\T\right), \frac{4a(r^2+a^2+|q|^2)\cos\th}{|q|^2}\dual\nab_\T + \frac{4a^2\cos\th}{|q|^2}\dual\nab_\Z\right]\psi\\
&=&  -\nab\left(\frac{8a^3(\sin\th)^2\cos\th}{|q|^2}\right)\c\nab\nab_\T\dual\psi  -\nab\left(\frac{8a^2\cos\th}{|q|^2}\right)\c\nab\nab_\Z\dual\psi\\
&&+O(ar^{-3})\nab_\Rhat\dk^{\leq 1}\psi+O(ar^{-2})\dk^{\leq 1}\psi +\dk \big(\Ga_g \c \dk \psi\big)+\Ga_b \c \squared_2\psi.
\eeaa
Setting 
\beaa
\widetilde{\OO}\psi = \OO\psi  +\frac{4a(r^2+a^2+|q|^2)\cos\th}{|q|^2}\nab_\T\dual\psi +\frac{4a^2\cos\th}{|q|^2}\nab_\Z\dual\psi,
\eeaa
we infer
\beaa
\frac{1}{|q|^2}\left[|q|^2\left(\squared_2-V+\frac{4a\cos\th}{|q|^2}\dual\nab_\T\right), \widetilde{\OO}\right]\psi &=& O(ar^{-2})\nab_{\Rhat}\dk^{\leq 1}\psi +O(ar^{-2})\dk^{\leq 1}\psi\\
&& +\dk^{\leq 3}(\Ga_g\c\psi) +\Ga_b \c \squared_2\psi
\eeaa
as stated. This concludes the proof of Lemma \ref{eq:defintionoftheoperqatorwidetildeOOcommutingwellwtihRWmodel}.
\end{proof}


 \subsection{Additional energy flux and bulk quantities}


We recall below the relevant main norms  introduced in section \ref{sec:defintionadditionalfrluxandbulkforaxandSSmorawetz:chap6}:


 \subsubsection{Pointwise notation}
 
 
\begin{definition}
  We introduce the following pointwise  notation for $\psi\in \sk_2$.
  \begin{enumerate}
 \item 
  We denote 
\beaa
 |\psi|^2_{\SS}&:=&\sum_{\aund=1}^4\big|\psia|^2.
\eeaa
\item 
Given a vectorfield $Y$ we denote
\beaa |\nab_Y\psi|^2_{\SS}&:=&\sum_{\aund=1}^4\big|\nab_Y\psia|^2.
\eeaa
  \end{enumerate}
  \end{definition}

    
\subsubsection{Degenerate energy norm} 


\begin{definition}[Degenerate energy norm]
We define the following degenerate energy   for $\psi\in \sk_2$ along $\Si(\tau)$:
\beaa
 E_{deg}[\psi](\tau) :=\int_{\Si(\tau)} \left( |\nab_4\psi|^2 +  \frac{|\De|}{r^4} |\nab_3\psi|^2 +|\nab\psi|^2 + r^{-2}|\psi|^2\right).
\eeaa
\end{definition}


\subsubsection{Refined Morawetz norms} 


    \begin{definition}[Refined Morawetz norms]
    We define the following Morawetz  norms for $\psi \in \sk_2$.
    \begin{enumerate}
    \item 
  The degenerate  axially symmetric Morawetz norms  in $\MM=\MM(\tau_1, \tau_2) $:
  \beaa
     \Mordot^{ax}_{deg}[\psi](\tau_1, \tau_2)&:=&\int_{\MM(\tau_1, \tau_2)} \frac{m}{r^2} |\nab_{\Rhat} \psi|^2 +\frac{\TT^2}{r^6} \left(\frac{m}{r^2} |\nab_\That \psi|^2 + r^{-1}|\nab\psi|^2\right), \\
     \Mor^{ax}_{deg}[\psi](\tau_1, \tau_2)&:=& \int_{\MM(\tau_1, \tau_2) } 
      \frac{m}{r^2} | \nab_\Rhat \psi|^2 +r^{-3}|\psi|^2+\frac{\TT^2}{r^6} \left(\frac{m}{r^2} |\nab_\That \psi|^2 + r^{-1}|\nab\psi|^2\right).
\eeaa

\item We also define the higher degenerate and non-degenerate Morawetz norms in $\MM=\MM(\tau_1, \tau_2) $, for a scalar function $z$:
 \beaa
 \Mordot_{\SSz, deg}[\psi](\tau_1, \tau_2)&:=& \int_{\MM} 
      \frac{m}{r^2}  | \nab_\Rhat \psi|_\SS^2+ r^{3} \Big(\big|\nab_\That \Psi_z \big |^2+r^2\big |\nab \Psi_z\big|^2\Big),\\
    \Mor_{\SSz, deg}[\psi](\tau_1, \tau_2)&:=& \int_{\MM} 
      \frac{m}{r^2}  | \nab_\Rhat \psi|_\SS^2+ r^{-3}|\psi|_{\SS}^2+ r^{3}\Big(\big|\nab_\That  \Psi_z \big |^2+r^2\big |\nab  \Psi_z\big|^2  \Big),
 \eeaa
 where\footnote{Note that $z$, $\RRa$ and $\De$ are functions depending only on $r$ so that $\pr_r$ in the formula for $\RRtp^\aund[z]$  simply  denotes  differentiation w.r.t. $r$.} 
 \beaa
 \Psi_z=\Psi_z[\psi]&:=& \RRtp^\aund[z] \psia, \qquad   \RRtp^\aund[z] := \pr_r\left( \frac{z}{\De} \RRa\right),
 \eeaa
with $z$ a suitable function of $r$ to be chosen later, and  with the scalar functions $\RRa$ given by \eqref{components-RR-aund}, i.e. 
         \beaa
          \RR^1&=&-(r^2+a^2)^2, \qquad \RR^2 = -2(r^2+a^2), \qquad \RR^3 =-1, \qquad \RR^4=\De.
          \eeaa
\end{enumerate}
  \end{definition}


 \subsection{Basic energy and Morawetz in perturbations of Kerr}
 \lab{sec:statementofresultschapter7and8in Kerralsoholdinperturbations}


Let $\MM$ satisfying the properties in section \ref{sec:setupchap9}. In particular, we assume the estimates \eqref{eq:assumptionsonMMforpartII:again}  for $(\Ga_b, \Ga_b)$. We now state the main results of section \ref{section:adaptiation-errors} whose goal is to show that the results proved in Kerr in Chapters \ref{chapter-proof-part1} and \ref{chapter-proof-mor-2} also hold on $\MM$.

First, we have the following analog of Proposition \ref{proposition:Morawetz1-step1}.
\begin{proposition}
\lab{proposition:Morawetz1-step1:perturbation}
For   $|a|/m \ll 1$ sufficiently small, the solution $\psi$ to the model gRW equation \eqref{eq:Gen.RW} satisfies in $\MM$
\bea\lab{eq:conditional-mor-par1-1-II:perturbation}
\bsplit
\Mor^{ax}_{deg}[\psi](\tau_1, \tau_2) \les& \sup_{[\tau_1, \tau_2]}E_{deg}[\psi](\tau)+\deh F_{\AA}[\psi](\tau_1, \tau_2)+F_{\Si_*}[\psi](\tau_1, \tau_2)\\
&+   \int_{\MM(\tau_1, \tau_2)}ar^{-2}\big(|\nab\psi|^2+|\nab_{\T}\psi|^2\big) \\
&+\int_{\MM(\tau_1, \tau_2)}\big(|\nab_{\Rhat} \psi|+r^{-1}|\psi|\big) |N|\\
& +\ep \left(\sup_{[\tau_1, \tau_2]}E[\psi](\tau)+ B^1_\de[\psi](\tau_1, \tau_2)\right).
 \end{split}
\eea
\end{proposition}

Also, we have the following analog of Proposition \ref{proposition:Energy1} for perturbations of Kerr. 
\begin{proposition}
\lab{proposition:Energy1:perturbation}
For   $|a|/m \ll 1$ sufficiently small, the solution $\psi$ to the model gRW equation \eqref{eq:Gen.RW} satisfies in $\MM$
 \bea
\bsplit
E_{deg}[\psi](\tau_2) +F_{\Si_*}[\psi](\tau_1, \tau_2) \les&  E_{deg}[\psi](\tau_1)+\deh\Big(E_{r\leq r_+(1+\deh)}[\psi](\tau_2)+F_\AA[\psi](\tau_1, \tau_2)\Big)\\
&+ \frac{|a|}{m}\Mor^{ax}_{deg}[\psi](\tau_1, \tau_2)     +\left|\int_{\MM(\tau_1, \tau_2)}  \nab_{\That_\de } \psi  \c N\right| \\
& +\int_{\MM(\tau_1, \tau_2)} |N|^2+\ep\left(\sup_{[\tau_1, \tau_2]}E[\psi](\tau)+ B_\de[\psi](\tau_1, \tau_2)\right).
\end{split}
\eea
\end{proposition} 

Next, we have the following analog of Proposition \ref{prop:recoverEnergyMorawetzwithrweightfromnoweight} for perturbations of Kerr. 
\begin{proposition}\lab{prop:recoverEnergyMorawetzwithrweightfromnoweight:perturbation}
Let $\psi$ a solution to the gRW equation \eqref{eq:Gen.RW}. Also, recall the norms $E[\psi]$ and $\textrm{Mor}[\psi]$ introduced in section \ref{subsection:basicnormsforpsi}. 
\bea
\bsplit
\textrm{Mor}_{r\geq 2r_1}[\psi](\tau_1, \tau_2) \les&  \sup_{\tau\in[\tau_1, \tau_2]}E_{r\geq r_1}[\psi](\tau) +r_1\textrm{Mor}_{r_1\leq r\leq 2r_1}[\psi](\tau_1, \tau_2)\\
& +\mathcal{N}_{r\geq r_1}[\psi, N](\tau_1, \tau_2) +\ep B_\de[\psi](\tau_1,\tau_2),
\end{split}
\eea
and
\bea
\bsplit
\sup_{\tau\in[\tau_1, \tau_2]}E_{r\geq 2r_1}[\psi](\tau) \les& 
E_{r\geq r_1}[\psi](\tau_1) +\mathcal{N}_{r\geq r_1}[\psi, N](\tau_1, \tau_2)+r_1\textrm{Mor}_{r_1\leq r\leq 2r_1}[\psi](\tau_1, \tau_2)\\
&+\frac{|a|}{r_1}\textrm{Mor}_{r\geq 2r_1}[\psi](\tau_1, \tau_2) +\ep B_\de[\psi](\tau_1,\tau_2).
\end{split}
\eea
\end{proposition}

Next, we have the following analog of Proposition \ref{prop:morawetz-higher-order} for perturbations of Kerr. 
\begin{proposition}[$\SS$-derivatives Morawetz estimates]
 \label{prop:morawetz-higher-order:perturbation}
 Let the scalar function $z$ given by
\beaa
z=z_0-\de_0 z_0^2, \qquad z_0=\frac{\De}{(r^2+a^2)^2}.
\eeaa
Then, for   $|a|/m \ll 1$ sufficiently small, the solution $\psi$ to the model gRW equation \eqref{eq:Gen.RW} satisfies in $\MM$
\bea\label{eq:mor-higher-unconditional-pertubation}
 \nn && \Mor_{\SSz, deg}[\psi](\tau_1, \tau_2) \\
\nn  &\les& \sum_{\aund=1}^4\left(\sup_{[\tau_1, \tau_2]}E_{deg}[\psi_\aund](\tau) +\deh F_{\AA}[\psi_{\aund}](\tau_1, \tau_2)+F_{\Si_*}[\psi_{\aund}](\tau_1, \tau_2)\right)\\
\nn&&+\left(\sup_{[\tau_1, \tau_2]}E_{deg}[(\nab_T, \dkb)^{\leq 1}\psi](\tau)+\deh F_{\AA}[(\nab_T, \dkb)^{\leq 1}\psi](\tau_1, \tau_2)+F_{\Si_*}[(\nab_T, \dkb)^{\leq 1}\psi](\tau_1, \tau_2)\right)^{\frac{1}{2}}\\
\nn &&\times\left(\sup_{[\tau_1, \tau_2]}E_{deg}[(\nab_T, \dkb)^{\leq 2}\psi](\tau)+\deh F_{\AA}[(\nab_T, \dkb)^{\leq 2}\psi](\tau_1, \tau_2)+F_{\Si_*}[(\nab_T, \dkb)^{\leq 2}\psi](\tau_1, \tau_2)\right)^{\frac{1}{2}} \\
  &&+\sum_{\aund=1}^4\int_{\MM(\tau_1, \tau_2)}\big(|\nab_{\Rhat} \psia|+r^{-1}|\psia|\big) |N_{\aund}| +\ep\left(\sup_{[\tau_1, \tau_2]}E^2[\psi](\tau)+ B^2_\de[\psi](\tau_1, \tau_2)\right).
 \eea 
 \end{proposition}

Finally, to show that $\Mor_{\SSz}$ controls $\psi$ in $\Mntrap$, we will rely on the following lemma which is the analog of Lemma \ref{LEMMA:LOWERBOUNDPHIZOUTSIDEMTRAP} for perturbations of Kerr. 
\begin{lemma}\lab{lemma:lowerboundPhizoutsideMtrap:perturbation}
For $\de_0>0$ small enough\footnote{Recall that the constant $\de_0>0$ is involved in the definition of $z=z_0-\de_0 z_0^2$.} and $|a|/m\ll \de_0$, there exists a universal constant $c_0>0$ such that the following holds on $\Mntrap$:
\beaa
r^3\Big(|\nab_T\Psi_z|^2+r^2|\nab\Psi_z|^2\Big)+r^{-3}|\psi|_{\SS}^2 &\geq& c_0r^{-3}\Big(|\nab_T\psi|^2_{\SS}+|\nab_Z\psi|^2_{\SS}+r^2|\nab\psi|^2_{\SS}\Big)\\
&& -O(ar^{-3})\big|(\nab_T, \dkb)^{\leq 1}\dk^{\leq 2}\psi\big|^2 +\Ddot_\a F^\a+\err_\ep
\eeaa
where the 1-form $F$ denotes  an expression in  $\psi$ for which we  have a bound of the form
\beaa
&&\left|\int_{\pr\MM(\tau_1, \tau_2) }F^\mu N_\mu\right| \\
&\les &   \left(\sup_{[\tau_1, \tau_2]}E_{deg}[(\nab_T, \dkb)^{\leq 1}\psi](\tau)+\deh F_{\AA}[(\nab_T, \dkb)^{\leq 1}\psi](\tau_1, \tau_2)+F_{\Si_*}[(\nab_T, \dkb)^{\leq 1}\psi](\tau_1, \tau_2)\right)^{\frac{1}{2}}\\
&\times&\left(\sup_{[\tau_1, \tau_2]}E_{deg}[(\nab_T, \dkb)^{\leq 2}\psi](\tau)+\deh F_{\AA}[(\nab_T, \dkb)^{\leq 2}\psi](\tau_1, \tau_2)+F_{\Si_*}[(\nab_T, \dkb)^{\leq 2}\psi](\tau_1, \tau_2)\right)^{\frac{1}{2}},
\eeaa
and the scalar function $\err_\ep$ satisfies 
\beaa
\int_{\MM}|\err_\ep| &\les& \ep\left(\sup_{[\tau_1, \tau_2]}E^2[\psi](\tau)+ B^2_\de[\psi](\tau_1, \tau_2)\right).
\eeaa
\end{lemma}

The rest of section \ref{section:adaptiation-errors} is devoted to the proof of the above results. Propositions \ref{proposition:Morawetz1-step1:perturbation}, \ref{proposition:Energy1:perturbation} and 
\ref{prop:recoverEnergyMorawetzwithrweightfromnoweight:perturbation} are proved in section \ref{sec:proofofresultschapter7inperturbationofKerr} and Proposition \ref{prop:morawetz-higher-order:perturbation} and 
Lemma \ref{lemma:lowerboundPhizoutsideMtrap:perturbation} are proved in section \ref{sec:proofofresultschapter8inperturbationofKerr}.

    
    \subsection{Acceptable error terms}
        

Recall the definition of the  energy-momentum tensor associated to the model gRW equation \eqref{eq:Gen.RW},  see \eqref{eq:definition-QQ-mu-nu}, 
 \bea\label{eq:definition-QQ-mu-nu:perturbation}
 \QQ_{\mu\nu}:=\Db_\mu  \psi \c \Db _\nu \psi 
          -\frac 12 \g_{\mu\nu} \left(\Db_\la \psi\c\Db^\la \psi + V\psi \c \psi\right)= \Db_\mu  \psi \c \Db _\nu \psi -\frac 1 2\g_{\mu\nu} \LL[\psi].
 \eea
 
When revisiting the proofs of Chapters \ref{chapter-proof-part1} and \ref{chapter-proof-mor-2} in the context of admissible perturbations $\MM$ of Kerr, we generate additional terms. We introduce below acceptable terms.
 
  \begin{definition}[Acceptable error terms]
 \lab{Def:acceptable-errors-ch9}
 The following quantity 
 \beaa
 F^{\mu\nu}\QQ_{\mu\nu}[\psi]+G^\mu \psi\c\Ddot_\mu\psi +H|\psi|^2
 \eeaa
 is said to be of the acceptable type, and denoted $\textrm{Good}$, if: 
 
\begin{itemize}
\item $F_{44}\in \Ga_g$, and all other components of $F_{\mu\nu}$ belong to $\Ga_b$.

\item All components of $G_\mu$ belong to $\Ga_g$.

\item $H\in r^{-1}\Ga_g$. 
\end{itemize} 
 \end{definition} 

The justification for Definition \ref{Def:acceptable-errors-ch9} is provided by the following lemma.
\begin{lemma}\lab{lemma:accetableerrortermsenergymorawetzundercontrol}
Assume that the quantity 
 \beaa
 F^{\mu\nu}\QQ_{\mu\nu}[\psi]+G^\mu \psi\c\Ddot_\mu\psi +H|\psi|^2
 \eeaa
 is of the acceptable type in the sense of Definition \ref{Def:acceptable-errors-ch9}. Then, it satisfies the following estimate
\beaa
\int_{\MM_{trap}(\tau_1, \tau_2)}\Big|F^{\mu\nu}\QQ_{\mu\nu}[\psi]+G^\mu \psi\c\Ddot_\mu\psi +H|\psi|^2\Big|^2 \les \ep\sup_{[\tau_1, \tau_2]}E[\psi](\tau)
\eeaa
and 
\beaa
\int_{\Mntrap(\tau_1, \tau_2)}\Big|F^{\mu\nu}\QQ_{\mu\nu}[\psi]+G^\mu \psi\c\Ddot_\mu\psi +H|\psi|^2\Big|^2 \les \ep B_\de[\psi](\tau_1, \tau_2).
\eeaa
\end{lemma}

\begin{remark}\lab{rmk::accetableerrortermsenergymorawetzundercontrolafsodsudfh}
Recall from Definition \ref{Def:acceptable-errors-ch9} that terms of the acceptable type are denoted $\textrm{Good}$. In view of Lemma \ref{lemma:accetableerrortermsenergymorawetzundercontrol}, we infer that such term satisfy the following estimate
\beaa
\int_{\MM}|\textrm{Good}| &\les& \ep\left(\sup_{[\tau_1, \tau_2]}E[\psi](\tau)+ B_\de[\psi](\tau_1, \tau_2)\right).
\eeaa
\end{remark}

\begin{proof}
We start with the control on $\MM_{trap}$. We have, using the control of $\Ga_g$ and $\Ga_b$, 
\beaa
&&\int_{\MM_{trap}(\tau_1, \tau_2)}\Big|F^{\mu\nu}\QQ_{\mu\nu}[\psi]+G^\mu \psi\c\Ddot_\mu\psi +H|\psi|^2\Big|^2\\
&\les&  \ep\int_{\MM_{trap}(\tau_1, \tau_2)}\frac{1}{\tau_{trap}^{1+\dec}}|\dk^{\leq 1}\psi|^2\\
&\les& \ep\left(\int_{\tau_1}^{\tau_2}\frac{1}{\tau_{trap}^{1+\dec}}\right)E[\psi]\\
&\les& \ep\sup_{\tau\in [\tau_1, \tau_2]}E[\psi]
\eeaa 
as stated. 

Concerning the control on $\Mntrap$, we have
\beaa
&&\int_{\Mntrap(\tau_1, \tau_2)}\Big|F^{\mu\nu}\QQ_{\mu\nu}[\psi]+G^\mu \psi\c\Ddot_\mu\psi +H|\psi|^2\Big|^2\\
&\les& \ep\int_{\Mntrap(\tau_1, \tau_2)}\Big(|\Ga_g||\nab_3\psi||\dk^{\leq 1}\psi|+r^{-1}|\Ga_g||\dk^{\leq 1}\psi|^2\Big)\\
&\les& \ep\int_{\Mntrap(\tau_1, \tau_2)}r^{-2}\Big(|\nab_3\psi||\dk^{\leq 1}\psi|+r^{-1}|\dk^{\leq 1}\psi|^2\Big)\\
&\les& \ep\int_{\Mntrap(\tau_1, \tau_2)}\Big(r^{-1-\de}|\nab_3\psi|^2+r^{\de-3}|\dk^{\leq 1}\psi|^2\Big)\\
&\les& \ep B^s_\de[\psi](\tau_1, \tau_2)
\eeaa
as stated. This concludes the proof of the lemma.
\end{proof}

Next, we introduce the linearization of the quantities $F_{\mu\nu}$, $G_\mu$ and $H$. 
\begin{definition}\lab{def:linearizationofMorawetzquantitiesinperturbationofKerr:chap9}
Let $F_{\mu\nu}$, $G_\mu$ and $H$. We define their linearization as follows
\beaa
F_{\mu\nu}=(F_{\mu\nu})_K+\widecheck{F_{\mu\nu}}, \qquad G_{\mu}=(G_{\mu})_K+\widecheck{G_{\mu}}, \qquad H=H_K+\widecheck{H},
\eeaa
where\footnote{This reflects the fact that in the principal null frame of Kerr, the only nontrivial 1-form is $\Re(\Jk)$, and that there are no nontrivial symmetric traceless 2-tensors.}: 
\begin{enumerate}
\item the quantities 
\beaa
(F_{44})_K, \qquad (F_{34})_K, \qquad (F_{33})_K, \qquad (G_4)_K, \qquad (G_3)_K, \qquad H_K,
\eeaa
are given as explicit functions of $(r, \cos\th)$ coinciding with the corresponding expressions in Kerr,

\item the quantities 
\beaa
(F_{4a})_K, \qquad (F_{3a})_K, \qquad (G_a)_K, 
\eeaa
are given as the 1-form $\Re(\Jk)_a$ multiplied by explicit functions of $(r, \cos\th)$ coinciding with the corresponding expressions in Kerr,

\item the quantity $(F_{ab})_K$ is given by the the symmetric 2-tensor $\ga_{ab}$ multiplied by explicit functions of $(r, \cos\th)$ coinciding with the corresponding expressions in Kerr.
\end{enumerate}
\end{definition}

In view of Definition \ref{def:linearizationofMorawetzquantitiesinperturbationofKerr:chap9}, we can decompose the above expressions in their main part and error terms as follows
\bea\lab{eq:deompostionofMorawetzquantitiesinperturvbationfadsfaosdiuasdfoiahdsflaiuhf}
\bsplit
 F^{\mu\nu}\QQ_{\mu\nu}[\psi]+G^\mu \psi\c\Ddot_\mu\psi +H|\psi|^2 &= \Big[F^{\mu\nu}\QQ_{\mu\nu}[\psi]+G^\mu \psi\c\Ddot_\mu\psi +H|\psi|^2\Big]_K+\err,\\
 \Big[F^{\mu\nu}\QQ_{\mu\nu}[\psi]+G^\mu \psi\c\Ddot_\mu\psi +H|\psi|^2\Big]_K &=(F^{\mu\nu})_K\QQ_{\mu\nu}[\psi]+(G^\mu)_K \psi\c\Ddot_\mu\psi +H_K|\psi|^2,\\
 \err &= \widecheck{F^{\mu\nu}}\QQ_{\mu\nu}[\psi]+\widecheck{G^\mu} \psi\c\Ddot_\mu\psi +\widecheck{H}|\psi|^2.
 \end{split}
\eea 
The proof of the results of section \ref{sec:statementofresultschapter7and8in Kerralsoholdinperturbations} will rely in particular on showing that the extra terms in perturbations of Kerr appearing in the various divergence identities involved in energy and Morawetz estimates are of the  acceptable type in the sense of Definition \ref{Def:acceptable-errors-ch9}, i.e, that $\err\in\textrm{Good}$  in \eqref{eq:deompostionofMorawetzquantitiesinperturvbationfadsfaosdiuasdfoiahdsflaiuhf}.


\subsection{Some deformation tensors}


We star with computing the deformation tensor of $e_3$ and $e_4$.
\begin{lemma}\lab{lemma:deformationtensore3}
We have
\bea
\bsplit
{}^{(e_3)}\pi_{44} &= 4\pr_r\left(\frac{\De}{|q|^2}\right)+\Ga_g,\quad {}^{(e_3)}\pi_{34} = \Ga_b,\quad {}^{(e_3)}\pi_{33} = 0,\\
{}^{(e_3)}\pi_{4a} &= -\frac{2ar}{|q|^2}\Re(\Jk)_a+\Ga_g,\quad {}^{(e_3)}\pi_{3a} = \Ga_b,\quad {}^{(e_3)}\pi_{ab} =  -\frac{2r}{|q|^2}\de_{ab}+\Ga_b,
\end{split}
\eea
and
\bea
\bsplit
{}^{(e_4)}\pi_{33} &=\Ga_b,\quad {}^{(e_4)}\pi_{34} =-2\pr_r\left(\frac{\De}{|q|^2}\right)+\Ga_g,\quad {}^{(e_4)}\pi_{44} =0,\\
{}^{(e_4)}\pi_{3a} &= \Ga_b\quad {}^{(e_4)}\pi_{4a} =\Ga_g,\quad {}^{(e_4)}\pi_{ab} = \frac{2r\De}{|q|^4}\de_{ab}+\Ga_g.
\end{split}
\eea
\end{lemma}

\begin{proof}
We have
\beaa
{}^{(e_3)}\pi_{44} &=& 2\g(\D_4e_3, e_4)=-8\om=4\pr_r\left(\frac{\De}{|q|^2}\right)+\Ga_g,\\
{}^{(e_3)}\pi_{34} &=& \g(\D_3e_3, e_4)=4\omb=\Ga_b,\\
{}^{(e_3)}\pi_{33} &=& 2\g(\D_3e_3, e_3)=0,\\
{}^{(e_3)}\pi_{4a} &=& \g(\D_4e_3, e_a)+\g(\D_ae_3, e_4)=2(\etab_a-\ze_a)\\
&=&-\Re\left(\left(\frac{aq}{|q|^2}+\frac{a\ov{q}}{|q|^2}\right)\Jk_a\right)+\Ga_g=-\frac{2ar}{|q|^2}\Re(\Jk)_a+\Ga_g,\\
{}^{(e_3)}\pi_{3a} &=& \g(\D_3e_3, e_a)+\g(\D_ae_3, e_3)=2\xib_a=\Ga_b,\\
{}^{(e_3)}\pi_{ab} &=& \g(\D_ae_3, e_b)+\g(\D_be_3, e_a)=\chib_{ab}+\chib_{ab}=\trchb\de_{ab}+2\chibh_{ab}\\
&=& -\frac{2r}{|q|^2}\de_{ab}+\Ga_b,
\eeaa
and
\beaa
{}^{(e_4)}\pi_{33} &=& 2\g(\D_3e_4, e_3)=-8\omb=\Ga_b,\\
{}^{(e_4)}\pi_{34} &=& \g(\D_4e_4, e_3)=4\om=-2\pr_r\left(\frac{\De}{|q|^2}\right)+\Ga_g,\\
{}^{(e_4)}\pi_{44} &=& 2\g(\D_4e_4, e_4)=0,\\
{}^{(e_4)}\pi_{3a} &=& \g(\D_3e_4, e_a)+\g(\D_ae_4, e_3)=2(\eta_a+\ze_a)=\Ga_b\\
{}^{(e_4)}\pi_{4a} &=& \g(\D_4e_4, e_a)+\g(\D_ae_4, e_4)=2\xi_a=\Ga_g,\\
{}^{(e_4)}\pi_{ab} &=& \g(\D_ae_4, e_b)+\g(\D_be_4, e_a)=\chi_{ab}+\chi_{ab}=\trch\de_{ab}+2\chih_{ab}\\
&=& \frac{2r\De}{|q|^4}\de_{ab}+\Ga_g,
\eeaa
as desired.
\end{proof}

Recall that we have
\beaa
\That &=&\frac 1 2 \left( \frac{|q|^2}{r^2+a^2} e_4+\frac{\De}{r^2+a^2}  e_3\right),\\
\Rhat &=&\frac 1 2 \left( \frac{|q|^2}{r^2+a^2} e_4-\frac{\De}{r^2+a^2}  e_3\right).
\eeaa
Recall also that $\That_\de=\That$ on $\Mntrap$ by construction. We now show that particular quantities $\widecheck{F^{\mu\nu}}\QQ_{\mu\nu}[\psi]$ are acceptable error terms.
\begin{lemma}\lab{lemma:deformationtensorRhatThatThatdeletagenerateGoodasldf}
We have
\beaa
\widecheck{{}^{(\Rhat)}\pi^{\mu\nu}}\QQ_{\mu\nu}[\psi]=\textrm{Good}, \qquad \widecheck{{}^{(\That)}\pi^{\mu\nu}}\QQ_{\mu\nu}[\psi]=\textrm{Good}, \qquad \widecheck{{}^{(\That_\de)}\pi^{\mu\nu}}\QQ_{\mu\nu}[\psi]=\textrm{Good}.
\eeaa
\end{lemma}

\begin{proof}
Since 
\beaa
\That &=&\frac 1 2 \left( \frac{|q|^2}{r^2+a^2} e_4+\frac{\De}{r^2+a^2}  e_3\right),\\
\Rhat &=&\frac 1 2 \left( \frac{|q|^2}{r^2+a^2} e_4-\frac{\De}{r^2+a^2}  e_3\right),
\eeaa
we infer
\beaa
\widecheck{{}^{(\Rhat)}\pi_{\mu\nu}}, \,\widecheck{{}^{(\That)}\pi_{\mu\nu}}  &=& O(1)\widecheck{{}^{(e_4)}\pi_{\mu\nu}}+O(1)\widecheck{{}^{(e_3)}\pi_{\mu\nu}}+O(r^{-1})\widecheck{e_4(r)}+O(r^{-1})\widecheck{e_3(r)}\\
&&+O(r^{-1})e_4(\cos\th)+O(r^{-1})e_3(\cos\th).
\eeaa
In view of Lemma \ref{lemma:deformationtensore3} and the fact that $\widecheck{e_3(r)}\in r\Ga_b$, while the other components behave better, we infer $\widecheck{{}^{(\Rhat)}\pi_{\mu\nu}}, \widecheck{{}^{(\That)}\pi_{\mu\nu}}\in\Ga_b$.  

Also, we have for the particular case $\mu=\nu=4$
\beaa
\widecheck{{}^{(\Rhat)}\pi_{44}}, \,\widecheck{{}^{(\That)}\pi_{44}}  &=& O(1)\widecheck{{}^{(e_4)}\pi_{44}}+O(1)\widecheck{{}^{(e_3)}\pi_{44}}+O(r^{-1})\widecheck{e_4(r)}+O(r^{-1})e_4(\cos\th).
\eeaa
In view of Lemma \ref{lemma:deformationtensore3} and the fact that $\widecheck{e_4(r)}\in r\Ga_g$, while $e_4(\cos\th)$ behaves better, we infer $\widecheck{{}^{(\Rhat)}\pi_{44}}, \widecheck{{}^{(\That)}\pi_{44}}\in\Ga_g$. 

Since we have shown that  $\widecheck{{}^{(\Rhat)}\pi_{\mu\nu}}, \widecheck{{}^{(\That)}\pi_{\mu\nu}}\in\Ga_b$ and $\widecheck{{}^{(\Rhat)}\pi_{44}}, \widecheck{{}^{(\That)}\pi_{44}}\in\Ga_g$, we have, in view of Definition \ref{Def:acceptable-errors-ch9}
\beaa
\widecheck{{}^{(\Rhat)}\pi^{\mu\nu}}\QQ_{\mu\nu}[\psi]=\textrm{Good}, \qquad \widecheck{{}^{(\That)}\pi^{\mu\nu}}\QQ_{\mu\nu}[\psi]=\textrm{Good}.
\eeaa
Also, since $\That_\de=\That$ on $\Mntrap$, we also have $\widecheck{{}^{(\That_\de)}\pi_{\mu\nu}}\in\Ga_b$ and $\widecheck{{}^{(\That_\de)}\pi_{44}}\in\Ga_g$, and thus
\beaa
\widecheck{{}^{(\That_\de)}\pi^{\mu\nu}}\QQ_{\mu\nu}[\psi]=\textrm{Good}.
\eeaa
This concludes the proof of the Lemma \ref{lemma:deformationtensorRhatThatThatdeletagenerateGoodasldf}.
\end{proof}


\subsection{Poincar\'e inequality}


Recall from Definition \ref{def:coordinatessystemSoftauandr:chap6} that the spheres $S(\tau, r)$ are covered by three coordinates systems $(x^1_S, x^2_S)$, $(x^1_E, x^2_E)$ and $(x^1_N, x^2_N)$ so that we have the following control on each corresponding coordinate chart
\beaa
\max_{b,c=1,2}|g_{bc} - (g_{a,m})_{bc}| &\les& r^2\ep,
\eeaa
where $g_{bc}$ denotes the induced metric coefficients in these coordinates systems, and $(g_{a,m})_{bc}$ the corresponding expression in Kerr. 

The following Poincar\'e inequality is the analog of Lemma \ref{lemma:poincareinequalityfornabonSasoidfh:chap6} in perturbations of Kerr. 
\begin{lemma}\lab{lemma:poincareinequalityfornabonSasoidfh:chap9}
For $\psi\in\sk_2$, we have 
\beaa
\int_S|\nab\psi|^2 &\geq& \frac{2\big(1+O(\ep+a^2r^{-2})\big)}{r^2}\int_S|\psi|^2 -O(a+\ep)\int_S\big(|\nab\psi|^2+r^{-2}|\nab_\T\psi|^2\big)\\
&& -O(\ep r^{-2})\int_S|\nab_3\psi|^2.
\eeaa
\end{lemma}

\begin{proof}
We denote by $e_b$, $b=1,2$, an orthonormal basis for the horizontal structure associated to $(e_3, e_4)$, and we look for vectorfields $X^S_b$, $b=1,2$, tangent to $S$ as 
\beaa
X^S_b &=& e_b -a\Re(\Jk)_b\T+\la_b e_3+\underline{\la}_be_4,\quad b=1,2,
\eeaa   
for some 1-forms $\la$ and $\underline{\la}$. The vectorfields $X^S_b$, $b=1,2$, are tangent to $S(\tau,r)$ if and only if $X^S_b(r)=X^S_b(\tau)=0$, and hence, if and only if 
\beaa
0 &=& \nab(r) - a\Re(\Jk)\T(r) +\la e_3(r)+\underline{\la}e_4(r), \\
0 &=& \nab(\tau) - a\Re(\Jk)\T(\tau) +\la e_3(\tau)+\underline{\la}e_4(\tau),
\eeaa
or
\beaa
 \la e_3(r)+\underline{\la}e_4(r) &=& -\nab(r) + a\Re(\Jk)\T(r), \\
\la e_3(\tau)+\underline{\la}e_4(\tau) &=& -\nab(\tau) + a\Re(\Jk)\T(\tau).
\eeaa
Since we have
\beaa
e_4(r)e_3(\tau)-e_3(r)e_4(\tau)=\frac{\De}{|q|^2}e_3(\tau)+e_4(\tau)+O(\ep)\gtrsim 1+O(\ep+\deh)>0,
\eeaa
we infer the existence and uniqueness of $(\la, \underline{\la})$ with 
\beaa
|\la|+|\underline{\la}| &\les& \big|\nab(r) - a\Re(\Jk)\T(r)|+\big|\nab(\tau) - a\Re(\Jk)\T(\tau)|.
\eeaa 
Together with 
\beaa
|\nab(r)|\les r^{-1}\ep, \qquad |\T(r)|\les \ep,\qquad  |\T(\tau)-1|\les \ep, \qquad |\nab(\tau)-a\Re(\Jk)|\les r^{-1}\ep,
\eeaa
we deduce
\beaa
|\la|+|\underline{\la}| &\les& \frac{\ep}{r}.
\eeaa 

Coming back to the tangent vectorfields 
\beaa
X^S_b &=& e_b -a\Re(\Jk)_b\T+\la_b e_3+\underline{\la}_be_4,\quad b=1,2,
\eeaa 
we infer
\beaa
\g(X^S_b, X^S_c) &=& \de_{bc} -2a^2\Re(\Jk)_b\Re(\Jk)_c+a^4|\Re(\Jk)|^2\Re(\Jk)_b\Re(\Jk)_c +a\la_b\Re(\Jk)_c  +a\la_c\Re(\Jk)_b\\
&& +a\frac{\De}{|q|^2}\underline{\la}_b\Re(\Jk)_c  +a\frac{\De}{|q|^2}\underline{\la}_c\Re(\Jk)_b -\la_b\underline{\la}_c -\la_c\underline{\la}_b,\quad b=1,2,
\eeaa
and hence, since $|\la|=O(r^{-1}\ep)$, $|\underline{\la}|=O(r^{-1}\ep)$, $|\Re(\Jk)|=O(r^{-1})$,  
\beaa
\g(X^S_b, X^S_c) &=& \de_{bc} +O((a^2+\ep^2)r^{-2}),\quad b=1,2.
\eeaa
We deduce
\beaa
\nab^S &=& \big(1+O((a^2+\ep^2)r^{-2})\big)\Big(\nab+O(ar^{-1})\nab_T+O(r^{-1}\ep) e_3+O(r^{-1}\ep)e_4\Big)\\
&=& \big(1+O((a^2+\ep^2)r^{-2})\big)\Big(\nab+O((a+\ep)r^{-1})\nab_T+O(r^{-1}\ep) e_3\Big).
\eeaa
This yields
\beaa
|\nab\psi|^2=|\nab^S\psi|^2-O(a+\ep)\big(|\nab\psi|^2+r^{-2}|\nab_\T\psi|^2\big) -O(\ep r^{-2})|\nab_3\psi|^2
\eeaa
so that 
\beaa
\int_S|\nab\psi|^2 &=& \int_S|\nab^S\psi|^2 -O(a+\ep)\int_S\big(|\nab\psi|^2+r^{-2}|\nab_\T\psi|^2\big)  -O(\ep r^{-2})\int_S|\nab_3\psi|^2.
\eeaa

Next, recall from the proof of Lemma \ref{lemma:poincareinequalityfornabonSasoidfh:chap6} that in the three coordinates systems of Remark \ref{rmk:whatthesecoordinatessystemsonSoftaurareinKerrandrefmainKerr:chap6}, we have
\beaa
\max_{b,c=1,2}\big|\pr^{\leq 2}((g_{a,m})_{x^ax^b} - r^2(\ga_{\SSS^2})_{x^ax^b})\big| &\les& a^2
\eeaa
where $\pr^{\leq 2}$ denotes at most 2 coordinates derivatives, and $(\ga_{\SSS^2})_{x^ax^b}$ the metric coefficients on $\SSS^2$ in the corresponding coordinates system. Together with Definition \ref{def:coordinatessystemSoftauandr:chap6}, we infer, for the induced metric on $S(\tau,r)$, 
\beaa
\max_{b,c=1,2}\big|\pr^{\leq 2}((g_{x^ax^b} - r^2(\ga_{\SSS^2})_{x^ax^b})\big| &\les& a^2+r^2\ep.
\eeaa
We deduce in particular, for the Gauss curvature $K_S$ of $S$ and the area radius $r_S$ of $S$, 
\beaa
K_S &=& \frac{1}{r^2}\big(1+O(\ep+a^2r^{-2})\big), \qquad r_S=r\big(1+O(\ep+a^2r^{-2})\big).
\eeaa
Applying the effective uniformization result of Corollary 3.8 in \cite{KS-GCM2}, we obtain a map $\Phi:\SSS^2\to S$ and a scalar function $u$ on $\SSS^2$ such that 
\beaa
\Phi^\# g &=& (r_S)^2e^{2u}\ga_{\SSS^2}, \qquad \|\dkb^{\leq 2}(u\circ\Phi^{-1})\|_{L^2(S)}\les (a^2r^{-2})r_S.
\eeaa
We infer, relying on the well known Poincar\'e inequality for $\psi\in\sk_2(\SSS^2)$, see for example \cite{GKS1}, 
\beaa
\int_S|\nab^S\psi|^2 &=& \frac{1}{r^2}\big(1+O(\ep+a^2r^{-2})\big)\int_{\SSS^2}|\nab_{\SSS^2}\Phi^{\#}\psi|^2+ \frac{1}{r^2}O(\ep+a^2r^{-2})\int_{\SSS^2}|\Phi^{\#}\psi|^2\\
&\geq& \frac{2}{r^2}\big(1+O(\ep+a^2r^{-2})\big)\int_{\SSS^2}|\Phi^{\#}\psi|^2\\
&\geq&  \frac{2}{r^2}\int_S\big(1+O(\ep+a^2r^{-2})\big)|\psi|^2  
\eeaa
and hence
\beaa
\int_S|\nab\psi|^2 &=& \int_S|\nab^S\psi|^2  -O(a+\ep)\int_S\big(|\nab\psi|^2+r^{-2}|\nab_\T\psi|^2\big)  -O(\ep r^{-2})\int_S|\nab_3\psi|^2\\
&\geq& \frac{2\big(1+O(\ep+a^2r^{-2})\big)}{r^2}\int_S|\psi|^2\\
&&  -O(a+\ep)\int_S\big(|\nab\psi|^2+r^{-2}|\nab_\T\psi|^2\big)  -O(\ep r^{-2})\int_S|\nab_3\psi|^2
\eeaa
as stated. This concludes the proof of Lemma \ref{lemma:poincareinequalityfornabonSasoidfh:chap9}.
\end{proof}


\subsection{Proof of Propositions \ref{proposition:Morawetz1-step1:perturbation}, \ref{proposition:Energy1:perturbation} and 
\ref{prop:recoverEnergyMorawetzwithrweightfromnoweight:perturbation}}
\lab{sec:proofofresultschapter7inperturbationofKerr}



\subsubsection{Preliminaries for the Morawetz estimate}


We start with some preliminaries for the Morawetz estimate. 
Let   $\psi\in \mathfrak{s}_k(\MM)$ be a solution of \eqref{eq:Gen.RW},   $X$  a vectorfield of the form
 \beaa
    X=  X^3  e_3 +X^4 e_4,
    \eeaa
 $w$ a scalar   and $M$  a one form. Define
 \beaa
\PP_\mu[X, w, M]&:=&\QQ_{\mu\nu} X^\nu +\frac 1 2  w \psi \c \Db_\mu \psi -\frac 1 4|\psi|^2   \pr_\mu w +\frac 1 4 |\psi|^2 M_\mu.
  \eeaa
Then, we have in view of  Proposition \ref{prop-app:stadard-comp-Psi-perturbations-Kerr}
  \beaa
  \D^\mu  \PP_\mu[X, w, M] &=& \frac 1 2 \QQ  \c\piX - \frac 1 2 X( V ) |\psi|^2+\frac 12  w \LL[\psi] -\frac 1 4|\psi|^2   \square_\g  w \\
& -& \big(\rhod +\etab\wedge\eta\big)\nab_{X^4e_4-X^3e_3}  \psi\c\dual\psi \\
&-& \frac{1}{2}\Im\Big(\tr\Xb H X^3 +\tr X\Hb X^4\Big)\c\nab\psi\c\dual\psi\\
  &+& \frac 1 4  \Div(|\psi|^2 M\big)+  \left(\nab_X\psi +\frac 1 2   w \psi\right)\c\big(\squared_2-V\psi\big)\\
  &+&r^{-2} \big(X^3\Ga_b+ X^4\Ga_g\big) \dk \psi \c \psi.
 \eeaa

We make the following choices for $(X, w, M)$, consistent with the ones in Chapter \ref{chapter-proof-part1}:
\begin{enumerate}
\item  $X$ of the type\footnote{Recall that we have in Kerr
\beaa
\Rhat=\frac{\De}{r^2+a^2}\pr_r.
\eeaa}
\beaa
X=\FF\frac{(r^2+a^2)}{\De}\Rhat, \qquad \Rhat=\frac{|q|^2}{2(r^2+a^2)}e_4-\frac{\De}{2(r^2+a^2)}e_3,
\eeaa
where $\FF=\FF(r)$ is such that $\FF/\De$ is a smooth function of $r$, and $(r\pr_r)^k\FF(r)$ is uniformly bounded on $\MM$. In particular, there holds $X= X^3  e_3 +X^4 e_4$ as anticipated, with $X^3$ and $X^4$ smooth functions of $r$ such that $(r\pr_r)^kX^3$ and $(r\pr_r)^kX^4$ are uniformly bounded on $\MM$. 

\item The scalar function $w(r)$ is such that $w/\De$ is a smooth function of $r$ and $(r\pr_r)^k(rw(r))$ is uniformly bounded on $\MM$.

\item  $M$ is of the type
\beaa
M=v(r)\frac{(r^2+a^2)}{\De}\Rhat
\eeaa 
where the scalar function $v(r)$ is such that $v/\De$ is a smooth function of $r$ and $(r\pr_r)^k(r^{\frac{5}{2}}v(r))$ is uniformly bounded on $\MM$.
\end{enumerate}

We have the following lemma.
\begin{lemma}\lab{lemma:decompositionofthedivergenceinerrorandtermKersoiduhf}
Using the decomposition \eqref{eq:deompostionofMorawetzquantitiesinperturvbationfadsfaosdiuasdfoiahdsflaiuhf}, we have
\beaa
  \D^\mu  \PP_\mu[X, w, M] &=& \Bigg[\frac 1 2 \QQ  \c\piX - \frac 1 2 X( V ) |\psi|^2+\frac 12  w \LL[\psi] -\frac 1 4|\psi|^2   \square_\g  w + \frac 1 4  \Div(|\psi|^2 M\big)\Bigg]_K\\
 && - \big(\rhod +\etab\wedge\eta\big)\nab_{X^4e_4-X^3e_3}  \psi\c\dual\psi \\
&& - \frac{1}{2}\Im\Big(\tr\Xb H X^3 +\tr X\Hb X^4\Big)\c\nab\psi\c\dual\psi\\
  &&+  \left(\nab_X\psi+\frac 1 2   w \psi\right)\c\big(\squared_2-V\psi\big) +\textrm{Good}
 \eeaa
where $\textrm{Good}$ is given by Definition \ref{Def:acceptable-errors-ch9}.
\end{lemma}

\begin{proof}
First, since $(r\pr_r)^kX^3$ and $(r\pr_r)^kX^4$ are uniformly bounded on $\MM$, we easily check that 
\beaa
r^{-2} \big(X^3\Ga_b+ X^4\Ga_g\big) \dk \psi \c \psi &=& \textrm{Good}
\eeaa
in the sense of Definition \ref{Def:acceptable-errors-ch9}. Together with the above, we deduce
 \beaa
  \D^\mu  \PP_\mu[X, w, M] &=& \frac 1 2 \QQ  \c\piX - \frac 1 2 X( V ) |\psi|^2+\frac 12  w \LL[\psi] -\frac 1 4|\psi|^2   \square_\g  w \\
&& - \big(\rhod +\etab\wedge\eta\big)\nab_{X^4e_4-X^3e_3}  \psi\c\dual\psi \\
&& - \frac{1}{2}\Im\Big(\tr\Xb H X^3 +\tr X\Hb X^4\Big)\c\nab\psi\c\dual\psi\\
  && +\frac 1 4  \Div(|\psi|^2 M\big)+  \left(\nab_X\psi+\frac 1 2   w \psi\right)\c\big(\squared_2-V\psi\big) +\textrm{Good}.
 \eeaa

Next, applying the decomposition \eqref{eq:deompostionofMorawetzquantitiesinperturvbationfadsfaosdiuasdfoiahdsflaiuhf}, and using the fact that $X^4$, $X^3$, $v$ and $w$ are functions of $r$, we obtain 
\beaa
&& \frac 1 2 \QQ  \c\piX - \frac 1 2 X( V ) |\psi|^2+\frac 12  w \LL[\psi] -\frac 1 4|\psi|^2   \square_\g  w + \frac 1 4  \Div(|\psi|^2 M\big)\\
  &=&\Bigg[\frac 1 2 \QQ  \c\piX - \frac 1 2 X( V ) |\psi|^2+\frac 12  w \LL[\psi] -\frac 1 4|\psi|^2   \square_\g  w + \frac 1 4  \Div(|\psi|^2 M\big)\Bigg]_K+\err
  \eeaa
  where 
  \beaa
  \err &=& \frac 1 2 \QQ  \c\widecheck{\piX} - \frac 1 2 \widecheck{X( V )} |\psi|^2 -\frac 1 4|\psi|^2   \widecheck{\square_\g  w} + \frac 1 4  \widecheck{\Div( M)}|\psi|^2.
\eeaa
Since 
\beaa
\widecheck{X( V )} &=& O(r^{-3})\Big(\widecheck{e_3(r)}, \widecheck{e_4(r)}, e_3(\cos\th), e_4(\cos\th)\Big)=r^{-2}\Ga_b=r^{-1}\Ga_g,\\
\widecheck{\square_\g  w} &=& r^{-3}\dk^{\leq 1}\big(\widecheck{e_3(r)}, \nab(r)\big)+r^{-2}\dk^{\leq 1}\widecheck{e_4(r)}+ r^{-2}\Ga_g=r^{-2}\Ga_b=r^{-1}\Ga_g,\\
\widecheck{\Div( M)} &=& O\left(r^{-7/2}\right)\big(\widecheck{e_3(r)}, \widecheck{e_4(r)}, \nab(r)\big)+O\left(r^{-5/2}\right)\widecheck{\Div(\Rhat)}\\
&=& r^{-\frac{5}{2}}\Ga_b+O\left(r^{-5/2}\right)\widecheck{\Div(\Rhat)}=r^{-\frac{3}{2}}\Ga_g+O\left(r^{-5/2}\right)\widecheck{\Div(\Rhat)}
\eeaa
we infer in view of Definition \ref{Def:acceptable-errors-ch9}
  \beaa
  \err &=& \frac 1 2 \QQ  \c\widecheck{\piX} +O\left(r^{-5/2}\right)\widecheck{\Div(\Rhat)}|\psi|^2 +\textrm{Good}.
\eeaa

Also, recall from the proof of Lemma \ref{lemma:deformationtensorRhatThatThatdeletagenerateGoodasldf} that $\widecheck{{}^{(\Rhat)}\pi}\in \Ga_b$ so that $\widecheck{\Div(\Rhat)}\in \Ga_b$ and hence
\beaa
O\left(r^{-5/2}\right)\widecheck{\Div(\Rhat)}|\psi|^2=O\left(r^{-5/2}\right)\Ga_b|\psi|^2=O\left(r^{-3/2}\right)\Ga_g|\psi|^2=\textrm{Good}
\eeaa
so that 
 \beaa
  \err &=& \frac 1 2 \QQ  \c\widecheck{\piX}  +\textrm{Good}.
\eeaa
Also, in view of the definition of $X$, we have 
\beaa
\widecheck{\piX_{\mu\nu}} &=& O(1)\widecheck{{}^{(\Rhat)}\pi_{\mu\nu}}+O(r^{-1})\big(\widecheck{e_3(r)}, \widecheck{e_4(r)}, \nab r\big)=O(1)\widecheck{{}^{(\Rhat)}\pi_{\mu\nu}}+\Ga_b,\\
\widecheck{\piX_{44}} &=& O(1)\widecheck{{}^{(\Rhat)}\pi_{44}}+O(r^{-1})\widecheck{e_4(r)}=O(1)\widecheck{{}^{(\Rhat)}\pi_{44}}+\Ga_g,
\eeaa
which together with Lemma \ref{lemma:deformationtensorRhatThatThatdeletagenerateGoodasldf} implies
\beaa
\QQ  \c\widecheck{\piX} =\textrm{Good}
\eeaa
and hence
 \beaa
  \err &=& \textrm{Good}.
\eeaa
This yields 
\beaa
&& \frac 1 2 \QQ  \c\piX - \frac 1 2 X( V ) |\psi|^2+\frac 12  w \LL[\psi] -\frac 1 4|\psi|^2   \square_\g  w + \frac 1 4  \Div(|\psi|^2 M\big)\\
  &=&\Bigg[\frac 1 2 \QQ  \c\piX - \frac 1 2 X( V ) |\psi|^2+\frac 12  w \LL[\psi] -\frac 1 4|\psi|^2   \square_\g  w + \frac 1 4  \Div(|\psi|^2 M\big)\Bigg]_K+\textrm{Good}
  \eeaa
  and thus 
  \beaa
  \D^\mu  \PP_\mu[X, w, M] &=& \Bigg[\frac 1 2 \QQ  \c\piX - \frac 1 2 X( V ) |\psi|^2+\frac 12  w \LL[\psi] -\frac 1 4|\psi|^2   \square_\g  w + \frac 1 4  \Div(|\psi|^2 M\big)\Bigg]_K\\
&& - \big(\rhod +\etab\wedge\eta\big)\nab_{X^4e_4-X^3e_3}  \psi\c\dual\psi \\
&& - \frac{1}{2}\Im\Big(\tr\Xb H X^3 +\tr X\Hb X^4\Big)\c\nab\psi\c\dual\psi\\
  &&+  \left(\nab_X\psi+\frac 1 2   w \psi\right)\c\big(\squared_2-V\psi\big) +\textrm{Good}
 \eeaa
as stated. This concludes the proof of Lemma \ref{lemma:decompositionofthedivergenceinerrorandtermKersoiduhf}.
\end{proof}


\subsubsection{Proof of Proposition \ref{proposition:Morawetz1-step1:perturbation}}


We choose $(X, w, M)$ as in Chapter \ref{chapter-proof-part1}:
\begin{enumerate}
\item  $X$ is given by
\beaa
\bsplit
X &=\FF\frac{(r^2+a^2)}{\De}\Rhat, \qquad \FF=-zhf, \\ 
z &=\frac{\De}{(r^2+a^2)^2}, \qquad f=-\frac{2\TT}{(r^2+a^2)^3}, \qquad h=\frac{(r^2+a^2)^4}{r(r^2-a^2)},
\end{split}
\eeaa
i.e.  $z, f, h$ correspond to the choices made in Proposition \ref{prop:Choice-zhf} and $X$ is given by 
\beaa
X &=& \frac{2\TT}{r(r^2-a^2)}\Rhat.
\eeaa
 
\item The scalar function $w(r)$ is given, as in Proposition \ref{prop:Choice-zhf}, by
\beaa
w=-z\pr_r(hf).
\eeaa

\item  $M$ is given by
\beaa
M=v(r)\frac{(r^2+a^2)}{\De}\Rhat
\eeaa 
where the scalar function $v(r)$ is the one of Lemma \ref{Le:Crucialpositivity-a=0} and satisfies in particular  $v(r)= O( m^{1/2}\Delta r^{-9/2})$.
\end{enumerate}

Next, we introduce the expression
    \bea\lab{definition-EE-gen:chap9}
      \begin{split}
\EE_K[X, w, M] &:=  \Bigg[\frac 1 2 \QQ  \c\piX - \frac 1 2 X( V ) |\psi|^2+\frac 12  w \LL[\psi] -\frac 1 4|\psi|^2   \square_\g  w \\
&+ \frac 1 4  \Div(|\psi|^2 M\big)\Bigg]_K.
\end{split}
\eea
Notice that $\EE_K[X, w, M]$ coincides in fact with the quantity $\EE[X, w, M]$ in \eqref{definition-EE-gen}. 
Thus, we have according to Proposition \ref{proposition:Morawetz1}
 \beaa
|q|^2 \EE_K[X, w, M]&=&\AA|\nab_r\psi|^2+ P+\VV |\psi|^2 +\frac 1 4 |q|^2  \Big[\D^\mu (|\psi|^2 M_\mu)\Big]_K,
\eeaa
where  the coefficients $\AA$ and $\VV$ are given by \eqref{eq:prop-Choice-AA} and \eqref{eq:prop:Choice-zhf} respectively, and  the principal term $P$ given by \eqref{eq:principal-term}, i.e.
 \beaa
 P&=& \frac{\TT}{r}\frac{r^2+a^2}{r^2-a^2} \left(  \frac{2\TT}{ (r^2+a^2)^3} |q|^2|\nab\psi|^2 - \frac{4ar}{(r^2+a^2)^2} \nab_\That\psi\c\nab_{\Z} \psi\right).
 \eeaa

Next, as in \eqref{def:Qr_red:new}, we define, for $\de>0$ sufficiently small chosen later, the following quadratic form
  \bea
  \lab{def:Qr_red:new:chap9}
  \bsplit
 \mbox{Qr}_\de[\psi]  = & (1-\de)\AA |\nab_r\psi|^2   + \left(\VV+  (1-\de)\frac{2\TT^2}{(r^2+a^2)^2(r^2-a^2)}\right) |\psi|^2 \\
 & +\frac 1 4 |q|^2  \Big[\D^\mu (|\psi|^2 M_\mu)\Big]_K
 \end{split}
 \eea
 so that, in view of the above, we have
  \beaa
|q|^2 \EE_K[X, w, M] &=& \de\AA|\nab_r\psi|^2+\de P + (1-\de)\left(P -\frac{2\TT^2}{r(r^2+a^2)^2(r^2-a^2)}|\psi|^2\right) +\mbox{Qr}_\de[\psi].
 \eeaa 
 Together with  Lemma \ref{lemma:poincareinequalityfornabonSasoidfh:chap9}, we obtain the following 
 analog of \eqref{eq:intermediarylowerboundEEafterPoincaresfasuighalsuidh} 
 \bea\lab{eq:intermediarylowerboundEEafterPoincaresfasuighalsuidh:chap9}
\nn\int_S\EE_K[X, w, M] &\geq& \de\int_S\left(\frac{m}{r^2}|\nab_{\Rhat}\psi|^2+\frac{\TT^2}{r^7}|\nab\psi|^2\right)+\int_S\frac{1}{|q|^2}\mbox{Qr}_\de[\psi] \\
\nn&& -O((a+\ep)r^{-1})\int_S\big(|\nab\psi|^2+r^{-2}|\nab_\T\psi|^2\big) -O(\ep r^{-3})\int_S|\nab_3\psi|^2\\
&& -O(\ep r^{-1}+ar^{-2})\int_S\big(r^{-2}|\psi|^2\big).
 \eea 
 Also, in view of the definition of $\mbox{Qr}_\de[\psi]$  in \eqref{def:Qr_red:new:chap9}, Lemma \ref{Le:Crucialpositivity-a=0} applies to $\mbox{Qr}_\de[\psi]$, and hence, for $|a|/m\ll 1$ and $\de>0$ sufficiently small,  we have
 \beaa
\mbox{Qr}_\de[\psi] &\geq& O(\de)  \left( m \big|\nab_{\Rhat}\Phi|^2 + r^{-1} |\Phi|^2 \right).
 \eeaa
Together with \eqref{eq:intermediarylowerboundEEafterPoincaresfasuighalsuidh:chap9}, we deduce the following analog of Proposition \ref{Prop.StongversionMorawetz1}. 

\begin{proposition}
 \lab{Prop.StongversionMorawetz1:chap9}
  There exists a   small universal constant $c_0>0$ such that, for $|a|/m\ll 1$ and $\ep>0$ small enough, we have on $\MM$
  \beaa
\nn\int_S\EE_K[X, w, M] &\geq& c_0\int_S\left(\frac{m}{r^2}|\nab_{\Rhat}\psi|^2+\frac{\TT^2}{r^7}|\nab\psi|^2  + r^{-3} |\psi|^2\right)\\
&& -O((a+\ep)r^{-1})\int_S\big(|\nab\psi|^2+r^{-2}|\nab_\T\psi|^2\big) -O(\ep r^{-3})\int_S|\nab_3\psi|^2.
 \eeaa
  \end{proposition}
 
In view of Lemma \ref{lemma:decompositionofthedivergenceinerrorandtermKersoiduhf} and \eqref{definition-EE-gen:chap9}, we have
\beaa
  \D^\mu  \PP_\mu[X, w, M] &=& \EE_K[X, w, M]  - \big(\rhod +\etab\wedge\eta\big)\nab_{X^4e_4-X^3e_3}  \psi\c\dual\psi \\
&& - \frac{1}{2}\Im\Big(\tr\Xb H X^3 +\tr X\Hb X^4\Big)\c\nab\psi\c\dual\psi\\
  &&+  \left(\nab_X\psi +\frac 1 2   w \psi\right)\c\big(\squared_2-V\psi\big) +\textrm{Good}
 \eeaa
with  $\EE_K[X, w, M]$ satisfying the estimate of Proposition \ref{Prop.StongversionMorawetz1:chap9}. Also, following the computation leading to \eqref{definition-EE-X=FFR}, we have  
  \beaa
  \nab_{X^4e_4-X^3e_3} = \FF \frac{r^2+a^2}{\De}\nab_{\That} \psi
  \eeaa
and
\beaa
\Im\Big(\tr\Xb H X^3 +\tr X\Hb X^4\Big)\c\nab &=&  \frac{4a^2r\cos\th\FF(r)}{|q|^4}\Re(\Jk) \c\nab+r^{-1}\Ga_b\nab\\
&=& \frac{4a^2r\cos\th\FF(r)}{(r^2+a^2)|q|^4}\left(\nab_\Z+a(\sin\th)^2\frac{r^2+a^2}{|q|^2}\nab_\That\right)+r^{-1}\Ga_b\nab\\
&=& \frac{4a^2r\cos\th\FF(r)}{(r^2+a^2)|q|^4}\nab_Z +\frac{4a^3r\cos\th(\sin\th)^2\FF(r)}{|q|^6}\nab_\That +r^{-1}\Ga_b\nab
\eeaa
which yields
\beaa
  \D^\mu  \PP_\mu[X, w, M] &=& \EE_K[X, w, M]  - \big(\rhod +\etab\wedge\eta\big)\nab_{X^4e_4-X^3e_3}  \psi\c\dual\psi \\
&&  -\left(\big(\rhod +\etab\wedge\eta\big)\frac{r^2+a^2}{\De}+\frac{2a^3r\cos\th(\sin\th)^2}{|q|^6}\right)\FF\nab_\That\psi\c\dual\psi\\
  &&  -\frac{2a^2r\cos\th\FF(r)}{(r^2+a^2)|q|^4}\nab_\Z\psi\c\dual\psi +  \left(\nab_X\psi +\frac 1 2   w \psi\right)\c\big(\squared_2-V\psi\big) +\textrm{Good}
 \eeaa
where we used the fact that $r^{-1}\Ga_b\psi\nab\psi$ is of type $\textrm{Good}$. 

Next, following Section \ref{section:lowerboundscontainigThat}, we consider
\beaa
\PP_\mu[0, w', 0] &=& \frac 1 2  w' \psi \c \Db_\mu \psi -\frac 1 4|\psi|^2   \pr_\mu w',\\
w'&=&- \de_1 \frac{4 m \De \TT^2}{r^2 (r^2+a^2)^4},
\eeaa
for some $\de_1>0$ small enouh. Lemma \ref{lemma:decompositionofthedivergenceinerrorandtermKersoiduhf} applies, and yields
\beaa
  \D^\mu  \PP_\mu[0, w', 0] &=& \left[\frac 12  w' \LL[\psi] -\frac 1 4|\psi|^2   \square_\g  w' \right]_K + \frac 1 2   w' \psi\c\big(\squared_2-V\psi\big) +\textrm{Good}.
 \eeaa
We introduce the expression
    \bea\lab{definition-EE-gen:chap9:bis}
      \begin{split}
\EE_K[0, w', 0] &:=  \left[\frac 12  w' \LL[\psi] -\frac 1 4|\psi|^2   \square_\g  w' \right]_K
\end{split}
\eea
which yields
\bea\lab{eq:intermediaryalmostfinaladaptaionbasicMorawetzperturbationaofsdihf:chap9}
\nn  \D^\mu\PP_\mu[X, w+w', M] &=& \EE_K[X, w, M] +\EE_K[0, w', 0] \\
\nn&&  -\left(\big(\rhod +\etab\wedge\eta\big)\frac{r^2+a^2}{\De}+\frac{2a^3r\cos\th(\sin\th)^2}{|q|^6}\right)\FF\nab_\That\psi\c\dual\psi\\
 \nn &&  -\frac{2a^2r\cos\th\FF(r)}{(r^2+a^2)|q|^4}\nab_\Z\psi\c\dual\psi\\
 &&+  \left(\nab_X\psi+\frac 1 2   (w+w') \psi\right)\c\big(\squared_2-V\psi\big) +\textrm{Good}.
 \eea
Also, in view of Section \ref{section:lowerboundscontainigThat}, we have
\beaa
|q|^2 \EE_K[0, w', 0] &=&\frac 1 2 \frac{w'}{\De} (r^2+a^2)^2    |\nab_{\Rhat}\psi|^2-\frac{w' (r^2+a^2)^2}{ 2 \De}|\nab_\That \psi|^2 +\frac 1 2  w' |q|^2|\nab\psi|^2\\
&& -\frac 1 2 \Big( \frac 1 2|q|^2 \big[\square_\g  w '\big]_K -  |q|^2  w'  V \Big) |\psi|^2.
\eeaa
Proposition \ref{Prop.StongversionMorawetz1:chap9} together with the choice of $w'$ easily yields, for $\de_1>0$ small enough, 
  \bea\lab{eq:intermediaryalmostfinaladaptaionbasicMorawetzperturbationaofsdihf:chap9:bis}
\nn&&\int_S\Big(\EE_K[X, w, M]+\EE_K[0, w', 0]\Big)\\
\nn  &\geq& c_0\int_S\left(\frac{m}{r^2}|\nab_{\Rhat}\psi|^2+\frac{\TT^2}{r^6}\left(\frac{m}{r^2}|\nab_{\That}\psi|^2+r^{-1}|\nab\psi|^2\right)  + r^{-3} |\psi|^2\right)\\
&& -O((a+\ep)r^{-1})\int_S\big(|\nab\psi|^2+r^{-2}|\nab_\T\psi|^2\big) -O(\ep r^{-3})\int_S|\nab_3\psi|^2.
 \eea
  
Finally, we derive the following analog of Lemma \ref{lemma-control-rhs}.
\begin{lemma}\lab{lemma-control-rhs:chap9} 
We have, for arbitrarily small positive constants $\de_2$, $\de_3$, to be fixed later:
\bea
\begin{split}
 \left(\nab_X\psi+\frac 1 2   w \psi\right)\c  \left(\squared_2 \psi - V \psi \right)  & \geq - \de_2\frac{m}{r^2} |\nab_\Rhat \psi|^2 - \de_2\frac{\TT^2}{ r^6}  \frac{m}{r^2} |\nab_\That \psi|^2\\
 & +O(1)(|\nab_\Rhat\psi|+r^{-1}|\psi|)|N|\\
 & +O(a^3r^{-6}) |\nab_{\Z} \psi|^2+O(a r^{-4}) |\psi|^2 +\textrm{Good}\\
& + \D_\mu\left(\frac{2a\cos\th}{|q|^2}(\Rhat)^\mu (r^2+a^2)\frac{z}{\De}h f \psi\c\nab_T\dual \psi\right)\\
& -\D_\mu\left(\T^\mu\frac{2a\cos\th}{|q|^2}   \frac{z}{\De}(r^2+a^2)hf \psi\c  \nab_{\Rhat}\dual \psi\right)
\end{split}
\eea
and
\bea
\nn&& \left|\left(\big(\rhod +\etab\wedge\eta\big)\frac{r^2+a^2}{\De}+\frac{2a^3r\cos\th(\sin\th)^2}{|q|^6}\right)\FF\nab_\That\psi\c\dual\psi\right|\\
\nn  &&  +\left|\frac{2a^2r\cos\th\FF(r)}{(r^2+a^2)|q|^4}\nab_\Z\psi\c\dual\psi\right|\\
& \leq & \de_3 \frac{ \TT^2}{r^6 }\left(\frac{m}{r^2}|\nab_\That \psi|^2 +\frac{1}{r^4}|\nab_\Z\psi|^2\right) + O(a^2r^{-6}) |\psi|^2+\textrm{Good}.
\eea
\end{lemma}

\begin{proof}
According to equation \eqref{eq:Gen.RW}, we have
\beaa
\left(\nab_X\psi+\frac 1 2   w \psi\right)\c  \big(\squared_2 \psi - V \psi \big) &=&\left(\nab_X\psi+\frac 1 2   w \psi\right)\c  \left(- \frac{4 a\cos\th}{|q|^2}\dual \nab_T  \psi+N \right).
\eeaa
We have
   \beaa
&& - \frac{4 a\cos\th}{|q|^2}\left(\nab_X\psi+\frac 1 2   w \psi\right) \c \big(\dual \nab_T  \psi \big)\\
&=& - \left[\frac{4 a\cos\th}{|q|^2}\left(\nab_X\psi+\frac 1 2   w \psi\right) \c \big(\dual \nab_T  \psi \big)\right]_K\\
&&+O(r^{-2})\widecheck{{}^{(\Rhat)}\pi_{\mu\nu}}\psi\c\nab_T\dual \psi + O(r^{-2}){}^{(\T)}\pi_{\mu\nu}\psi\c\nab_{\Rhat}\dual \psi +O(r^{-2})\psi\c\dual \widecheck{[\nab_{\Rhat}, \nab_\T]\psi}\\
&&+O(r^{-3})\Big(\widecheck{\Rhat(r)}, \Rhat(\cos\th), \T(r), \T(\cos\th)\Big)|\psi|^2
\eeaa
where the first term in the RHS is the one computed in the proof of Lemma \ref{lemma-control-rhs}, i.e it verifies 
\beaa
  \left[\left(\nab_X\psi+\frac 1 2   w \psi\right) \c \left(- \frac{4 a\cos\th}{|q|^2}\dual \nab_T  \psi \right)\right]_K &\geq& -\de_2 \frac{\TT^2}{ r^6}  \frac{m}{r^2} |\nab_\That \psi|^2 +O(a^3r^{-6}) |\nab_{\Z} \psi|^2+O(ar^{-4}) |\psi|^2\\
  && + \D_\mu\left(\frac{2a\cos\th}{|q|^2}(\Rhat)^\mu (r^2+a^2)\frac{z}{\De}h f \psi\c\nab_T\dual \psi\right)\\
&& -\D_\mu\left(\T^\mu\frac{2a\cos\th}{|q|^2}   \frac{z}{\De}(r^2+a^2)hf \psi\c  \nab_{\Rhat}\dual \psi\right).
\eeaa

Also, since $\widecheck{{}^{(\Rhat)}\pi_{\mu\nu}}, {}^{(\T)}\pi_{\mu\nu}\in\Ga_b$, $\widecheck{\Rhat(r)}, \T(r)\in r\Ga_b$, and $\Rhat(\cos\th), \T(\cos\th)\in\Ga_b$, we have
 \beaa
&& - \frac{4 a\cos\th}{|q|^2}\left(\nab_X\psi+\frac 1 2   w \psi\right) \c \big(\dual \nab_T  \psi \big)\\
&=& - \Big[\frac{4 a\cos\th}{|q|^2}\left(\nab_X\psi+\frac 1 2   w \psi\right) \c \big(\dual \nab_T  \psi \big)\Big]_K +r^{-2}\Ga_b\c\psi\c\nab_T\dual \psi + r^{-2}\Ga_b\c\psi\c\nab_{\Rhat}\dual \psi\\
&& +O(r^{-2})\psi\c\dual \widecheck{[\nab_{\Rhat}, \nab_\T]\psi} +r^{-2}\Ga_b|\psi|^2.
\eeaa
Moreover, we have
\beaa
[\nab_{\Rhat}, \nab_\T]\psi &=& \nab_{[\Rhat, \T]}\psi + \R^{bc}\,_{\mu\nu}\psi_{cd}{\Rhat}^\mu{\T}^\nu\\
&=& \Ga_b\c(\nab_3, \nab_4, \nab)\psi+\Rhat^3\R^{bc}\,_{3\nu}\psi_{cd}{\T}^\nu+\Rhat^4\R^{bc}\,_{4\nu}\psi_{cd}{\T}^\nu
\eeaa
and hence
\beaa
\widecheck{[\nab_{\Rhat}, \nab_\T]\psi} &=& \Ga_b\c(\nab_3, \nab_4, \nab)\psi+\Big(O(1)\widecheck{\rhod}+O(ar^{-1})(\bb,\b)\Big)\psi\\
&=& \Ga_b\c(\nab_3, \nab_4, \nab)\psi+r^{-1}\Ga_g\c\psi
\eeaa
which yields
 \beaa
&& - \frac{4 a\cos\th}{|q|^2}\left(\nab_X\psi+\frac 1 2   w \psi\right) \c \big(\dual \nab_T  \psi \big)\\
&=& - \Big[\frac{4 a\cos\th}{|q|^2}\left(\nab_X\psi+\frac 1 2   w \psi\right) \c \big(\dual \nab_T  \psi \big)\Big]_K +r^{-2}\Ga_b\c\psi\c(\nab_3, \nab_4, \nab)\psi+r^{-2}\Ga_b|\psi|^2.
\eeaa
In view of Definition \ref{Def:acceptable-errors-ch9} for $\textrm{Good}$, we infer
 \beaa
 - \frac{4 a\cos\th}{|q|^2}\left(\nab_X\psi+\frac 1 2   w \psi\right) \c \big(\dual \nab_T  \psi \big) &=& - \Big[\frac{4 a\cos\th}{|q|^2}\left(\nab_X\psi+\frac 1 2   w \psi\right) \c \big(\dual \nab_T  \psi \big)\Big]_K +\textrm{Good}.
\eeaa
Together with the above, we infer
\beaa
 \left(\nab_X\psi+\frac 1 2   w \psi\right) \c \left(- \frac{4 a\cos\th}{|q|^2}\dual \nab_T  \psi \right) &\geq& -\de_2 \frac{\TT^2}{ r^6}  \frac{m}{r^2} |\nab_\That \psi|^2 +O(a^3r^{-6}) |\nab_{\Z} \psi|^2+O(ar^{-4}) |\psi|^2\\
  && + \D_\mu\left(\frac{2a\cos\th}{|q|^2}(\Rhat)^\mu (r^2+a^2)\frac{z}{\De}h f \psi\c\nab_T\dual \psi\right)\\
&& -\D_\mu\left(\T^\mu\frac{2a\cos\th}{|q|^2}   \frac{z}{\De}(r^2+a^2)hf \psi\c  \nab_{\Rhat}\dual \psi\right)+\textrm{Good}.
\eeaa
As in the proof of Lemma \ref{lemma-control-rhs}, we  bound the second product by
\beaa
\Big| \left(\nab_X\psi+\frac 1 2   w \psi\right)\c N \Big|&\les&  \Big( | \nab_\Rhat \psi | + r^{-1}|\psi| \Big)    |  N|. 
\eeaa
By putting together with the previous bound we obtain the first desired estimate. 

 Next, since $\widecheck{\rhod}\in r^{-1}\Ga_g$, $\widecheck{\eta}\in \Ga_b$, and $\widecheck{\etab}\in \Ga_g$, and since $\frac{r^2+a^2}{\De} \,\FF$ is bounded, we have
  \beaa
\big(\rhod +\etab\wedge\eta\big)\frac{r^2+a^2}{\De}\FF\nab_\That\psi\c\dual\psi &=& \left[\big(\rhod +\etab\wedge\eta\big)\frac{r^2+a^2}{\De}\FF\nab_\That\psi\c\dual\psi\right]_K+\textrm{Good}.
\eeaa
 Also, proceeding as in Lemma \ref{lemma-control-rhs}, we have
 \beaa
&&\left|\left[\big(\rhod +\etab\wedge\eta\big)\frac{r^2+a^2}{\De}\FF\nab_\That\psi\c\dual\psi\right]_K\right|+\left|\frac{2a^3r\cos\th(\sin\th)^2}{|q|^6}\FF\nab_\That\psi\c\dual\psi\right|\\
&&+\left|\frac{2a^2r\cos\th\FF(r)}{(r^2+a^2)|q|^4}\nab_\Z\psi\c\dual\psi\right|\\
&\leq& \de_3 \frac{ \TT^2}{r^6 }\left(\frac{m}{r^2}|\nab_\That \psi|^2 +\frac{1}{r^4}|\nab_\Z\psi|^2\right)+O(a^2r^{-6}) |\psi|^2
\eeaa
and hence
\beaa
\nn&& \left|\left(\big(\rhod +\etab\wedge\eta\big)\frac{r^2+a^2}{\De}+\frac{2a^3r\cos\th(\sin\th)^2}{|q|^6}\right)\FF\nab_\That\psi\c\dual\psi\right|\\
\nn  &&  +\left|\frac{2a^2r\cos\th\FF(r)}{(r^2+a^2)|q|^4}\nab_\Z\psi\c\dual\psi\right|\\
& \leq &  \de_3 \frac{ \TT^2}{r^6 }\left(\frac{m}{r^2}|\nab_\That \psi|^2 +\frac{1}{r^4}|\nab_\Z\psi|^2\right) + O(a^2r^{-6}) |\psi|^2+\textrm{Good}
\eeaa
as stated. 
\end{proof}

 By applying the divergence theorem to \eqref{eq:intermediaryalmostfinaladaptaionbasicMorawetzperturbationaofsdihf:chap9}, relying on the lower bound \eqref{eq:intermediaryalmostfinaladaptaionbasicMorawetzperturbationaofsdihf:chap9:bis}, and using Lemma \ref{lemma-control-rhs}, we then obtain  the following estimate 
 \beaa
&&\int_{\MM(\tau_1, \tau_2)}  \frac{m}{r^2} |\nab_{\Rhat} \psi|^2+r^{-3}|\psi|^2 +\frac{\TT^2}{r^6} \left(\frac{m}{r^2} |\nab_\That \psi|^2 + r^{-1}|\nab\psi|^2\right)    \\
\les &&\int_{\pr\MM(\tau_1, \tau_2)}|\PP\c N_\Si|+\int_{\MM(\tau_1, \tau_2)}\Big((a+\ep)r^{-1}\big(|\nab\psi|^2+r^{-2}|\nab_\T\psi|^2\big)+\ep r^{-3}|\nab_3\psi|^2+\textrm{Good}\Big)\\
&& +\int_{\MM(\tau_1, \tau_2)}\Big(|\nab_\Rhat\psi|+r^{-1}|\psi|\Big)| N|.
\eeaa
Since the extra terms in perturbations of Kerr satisfy 
\beaa
&&\int_{\MM(\tau_1, \tau_2)}\Big(\ep r^{-1}\big(|\nab\psi|^2+r^{-2}|\nab_\T\psi|^2+r^{-2}|\nab_3\psi|^2\big)+\textrm{Good}\Big)\\
&\les& \ep\int_{\MM(\tau_1, \tau_2)}\Big(r^{-3}|\nab_3\psi|^2+r^{-4}|\dk\psi|^2\Big)+\ep\left(\sup_{[\tau_1, \tau_2]}E[\psi](\tau)+ B_\de[\psi](\tau_1, \tau_2)\right)\\
&\les& \ep\left(\sup_{[\tau_1, \tau_2]}E[\psi](\tau)+ B^1_\de[\psi](\tau_1, \tau_2)\right),
\eeaa
where we used in particular Remark \ref{rmk::accetableerrortermsenergymorawetzundercontrolafsodsudfh}, this concludes the proof of Proposition  \ref{proposition:Morawetz1-step1:perturbation}.


\subsubsection{Proof of Proposition \ref{proposition:Energy1:perturbation}}


We start with the following analog of Lemma \ref{Lemma:QQ(That,NSi)}.
\begin{lemma}
\lab{Lemma:QQ(That,NSi):chap9}
The following hold true with a sufficiently small $c_0>0$, for any $|a| \ll m$,
\bea
\lab{eq:Lemma-forQQ3:chap9}
\bsplit
\int_{\Si(\tau)}\QQ(\That, N_\Si)&\ge   c_0E_{deg}[\psi](\tau)- O(\deh) E_{r\leq r_+(1+\deh)}[\psi](\tau),\\
\int_{\Si_*(\tau_1, \tau_2)}\QQ(\That, N_{\Si_*}) &\ge c_0F_{\Si_*}[\psi](\tau_1, \tau_2),\\
\int_{\AA(\tau_1, \tau_2)}\QQ(\That, N_{\AA}) &\gtrsim - \deh F_{\AA}[\psi](\tau_1, \tau_2).
 \end{split}
\eea
\end{lemma}

\begin{proof}
For $\ep>0$ sufficiently small, satisfying in particular $\ep\ll\deh$, proceeding as in the proof of  Lemma \ref{Lemma:QQ(That,NSi)}, we obtain 
\beaa
\bsplit
\int_{\Si(\tau)}\QQ(\That, N_\Si)&\ge   c_0\int_{\Si(\tau)}\left( |\nab_4\psi|^2 + \frac{|\De|}{r^4} |\nab_3\psi|^2 +|\nab\psi|^2\right) - O(\deh) E_{r\leq r_+}[\psi](\tau),\\
\int_{\Si_*(\tau_1, \tau_2)}\QQ(\That, N_{\Si_*}) &\ge c_0F_{\Si_*}[\psi](\tau_1, \tau_2),\\
\int_{\AA(\tau_1, \tau_2)}\QQ(\That, N_{\AA}) &\gtrsim - \deh F_{\AA}[\psi](\tau_1, \tau_2).
 \end{split}
\eeaa
In particular, we have obtained the desired estimates on $\Si_*$ and $\AA$.

Also, we have in view of Lemma \ref{lemma:poincareinequalityfornabonSasoidfh:chap9}, for $\psi\in\sk_2$ and $|a|\ll m$, 
\beaa
\frac{1}{r^2}\int_S|\psi|^2 &\les& \int_S\big(|\nab\psi|^2+r^{-2}|\nab_{\T}\psi|^2+\ep r^{-2}|\nab_3\psi|^2)\\
&\les& \int_S\left( |\nab_4\psi|^2 + \frac{|\De|}{r^4} |\nab_3\psi|^2 +|\nab\psi|^2+\ep r^{-2}|\nab_3\psi|^2\right). 
\eeaa
Together with the above, we infer, for $\ep\ll\deh$, 
\beaa
\int_{\Si(\tau)}\QQ(\That, N_\Si)&\ge&   c_0\int_{\Si(\tau)}\left( |\nab_4\psi|^2 + \frac{|\De|}{r^4} |\nab_3\psi|^2 +|\nab\psi|^2+r^{-2}|\psi|^2\right) \\
&& - O(\deh) E_{r\leq r_+(1+\deh)}[\psi](\tau),\\
\eeaa
and hence
\beaa
\int_{\Si(\tau)}\QQ(\That, N_\Si)&\ge   c_0E_{deg}[\psi](\tau)- O(\deh) E_{r\leq r_+(1+\deh)}[\psi](\tau)
\eeaa
as stated. This concludes the proof of Lemma \ref{Lemma:QQ(That,NSi):chap9}.
\end{proof}

We consider  the energy current associated  to the modified timelike vectorfield $\That_\de$ introduced in Definition \ref{definitionThat_de} by
 \beaa
 \That_\de= \T+\chi_\de \Z, \qquad \chi_\de= \frac{a}{r^2+a^2} \chi_0\left( \de^{-1} \frac{\TT}{r^3} \right),
 \eeaa
with $\de=\frac{1}{10}$ and $|a|/m\ll 1$ small enough, where $\chi_\de= \frac{a}{r^2+a^2} \chi_0\Big( \de^{-1} \frac{\TT}{r^3} \Big)$ with $\chi_0$ given by \eqref{definition:chi-T_chippppp} satisfying in particular $\chi_0=0$ in $\MM_{trap}$.

From   Proposition \ref{prop-app:stadard-comp-Psi}, we have for the current associated to $\That_\de$:
  \beaa
  \begin{split}
  \D^\mu  \PP_\mu[\That_\de, 0, 0] &= \frac 1 2 \QQ  \c \,^{(\That_\de)} \pi -\frac 1 2 \That_\de(V) |\psi|^2 + \That_\de^\mu \Db^\nu  \psi ^a\R_{ ab   \nu\mu}\psi^b + \nab_{\That_\de} \psi  \c \big(\squared_2 \psi - V \psi \big).
  \end{split}
 \eeaa
 Since
 \beaa
 \That_\de(V) &=& O(r^{-3})\Big(\That_\de(r),\, \That_\de(\cos\th)\Big) =  O(r^{-3})\Big(\T(r), \,\T(\cos\th), \,\nab(r), \,\widecheck{\nab(\cos\th)}\Big)\\
 &=& r^{-2}\Ga_b
 \eeaa
 and hence $\That_\de(V) |\psi|^2=\textrm{Good}$ which yields
  \beaa
  \begin{split}
  \D^\mu  \PP_\mu[\That_\de, 0, 0] &= \frac 1 2 \QQ  \c \,^{(\That_\de)} \pi  + \That_\de^\mu \Db^\nu  \psi ^a\R_{ ab   \nu\mu}\psi^b + \nab_{\That_\de} \psi  \c \big(\squared_2 \psi - V \psi \big) +\textrm{Good}
  \end{split}
 \eeaa 
or
 \beaa
  \D^\mu  \PP_\mu[\That_\de, 0, 0] &=& \frac 1 2 \QQ  \c \,^{(\That_\de)} \pi   +\T^\mu \Db^\nu  \psi ^a\Rdot_{ ab   \nu\mu}\psi^b +(\That_\de-\T)^\mu \Db^\nu  \psi ^a\Rdot_{ ab   \nu\mu}\psi^b\\
  &&+ \nab_{\That_\de} \psi  \c \big(\squared_2 \psi - V \psi \big) +\textrm{Good}.
 \eeaa
Since $\Rdot_{ ab   \nu\mu}$ is antisymmetric with respect to $(a, b)$, we rewrite 
  \beaa
  \D^\mu  \PP_\mu[\That_\de, 0, 0] &=& \frac 1 2 \QQ  \c \,^{(\That_\de)} \pi   +\frac{1}{2}\T^\mu \in^{ab}\Rdot_{ ab   \nu\mu} \dual\psi\c\Db^\nu\psi  +(\That_\de-\T)^\mu \Db^\nu  \psi ^a\Rdot_{ ab   \nu\mu}\psi^b\\
  &&+ \nab_{\That_\de} \psi  \c \big(\squared_2 \psi - V \psi \big) +\textrm{Good}.
 \eeaa 
 Introducing the following spacetime 1-form
\bea\lab{eq:thespacetime1formAappearingintheenergyestimate:chap9}
A_\mu &:=& \in^{bc}\Rdot_{bc \mu\nu}\T^\nu,
\eea 
 we infer the following analog of \eqref{eq:basiccomputationenergyestimatesKerrwithThatdeltaflsdiuh}
 \bea\lab{eq:basiccomputationenergyestimatesKerrwithThatdeltaflsdiuh:chap9}
 \nn \D^\mu  \PP_\mu[\That_\de, 0, 0] &=& \frac 1 2 \QQ  \c \,^{(\That_\de)} \pi   +\frac{1}{2}A_\nu \dual\psi\c\Db^\nu\psi  +(\That_\de-\T)^\mu \Db^\nu  \psi ^a\Rdot_{ ab   \nu\mu}\psi^b\\
  &&+ \nab_{\That_\de} \psi  \c \big(\squared_2 \psi - V \psi \big)+\textrm{Good}.
 \eea 
 
Next, we compute the components of $A$.
\begin{lemma}\lab{lemma:computationofthecomponentsofthetensorAforexactenergyconservationKerr:chap9}
Let $A$ the spacetime 1-form given by \eqref{eq:thespacetime1formAappearingintheenergyestimate:chap9}. Then, 
we have
\beaa
A_4 &=& -4\rhod\T^3 -4(\etab\wedge\eta)\T^3 +  \trch\big({}^{(h)}\T\wedge\etab)  -\atrch\big(\etab\c {}^{(h)}\T)+r^{-1}\Ga_g,\\
A_3 &=& 4\rhod\T^4 +4(\etab\wedge\eta)\T^4 +  \trchb\big({}^{(h)}\T\wedge\eta)  -\atrchb\big(\eta\c {}^{(h)}\T\big)+r^{-2}\Ga_b,\\
A_e &=& \Big( -  \trchb\dual\eta_e   + \atrchb\eta_e\Big)\T^3 +\Big( -  \trch\dual\etab_e   + \atrch\etab_e\Big)\T^4\\
&& -\frac{1}{2}\Big(4\rho  + \trch\trchb+\atrch\atrchb\Big)\dual ({}^{(h)}\T)_e+r^{-1}\Ga_b.
\eeaa
\end{lemma}

\begin{proof}
First, note that we have in view of Proposition \ref{proposition:componentsofB} and the definition of $\Ga_b$ and $\Ga_g$
\bea\lab{eq:consequenceproposition:componentsofB}
\bsplit
\B_{ a   b  c 3}&=- \B_{ a   b  3c}=      -  \trchb  \big( \de_{ca}\eta_b-  \de_{cb} \eta_a\big)  -  \atrchb \big( \in_{ca}  \eta_b -  \in_{cb}  \eta_a\big)  +r^{-1}\Ga_b,\\
\B_{ a   b  c 4}&=- \B_{ a   b  c 4 c} =  -  \trch  \big( \de_{ca}\etab_b-  \de_{cb} \etab_a\big)  -  \atrch \big( \in_{ca}  \etab_b -  \in_{cb}  \etab_a\big)  +r^{-1}\Ga_g,\\
\B_{ a   b  3 4}&=- \B_{ a   b  43} =4(\etab_a\eta_b - \eta_a\etab_b) +\Ga_b\c\Ga_g,\\
\B_{ab cd}&= -\B_{ab dc} =  -\frac{1}{2}\left(\trch\trchb+\atrch\atrchb\right)\in_{ab}\in_{cd} +r^{-1}\Ga_b. 
\end{split}
\eea

Next, we rewrite $A_\mu$ as 
\beaa
A_\mu &=& \in^{bc}\Rdot_{bc \mu 3}\T^3 +\in^{bc}\Rdot_{bc \mu 4}\T^4 +\in^{bc}\Rdot_{bc \mu d}\T^d.
\eeaa
and compute the various components of $A_\mu$. We have, using the horizontal tensor ${}^{(h)}\T$ defined by ${}^{(h)}\T_b=\T_b$, the definition  \eqref{eq:DefineRdot} of $\Rdot$, \eqref{eq:consequenceproposition:componentsofB}, and the computations of the main terms in the proof of Lemma \ref{lemma:computationofthecomponentsofthetensorAforexactenergyconservationKerr},  
\beaa
A_4 &=& \in^{bc}\Rdot_{bc43}\T^3  +\in^{bc}\Rdot_{bc4d}\T^d,\\
&=& -4\rhod\T^3 -4(\etab\wedge\eta)\T^3 +  \trch\big({}^{(h)}\T\wedge\etab)  -\atrch\big(\etab\c {}^{(h)}\T)\\
&&+O(\T^3)\Ga_b\c\Ga_g+O(|{}^{(h)}\T|)\big(\b, r^{-1}\Ga_g\big)\\
&=& -4\rhod\T^3 -4(\etab\wedge\eta)\T^3 +  \trch\big({}^{(h)}\T\wedge\etab)  -\atrch\big(\etab\c {}^{(h)}\T)+r^{-1}\Ga_g,
\eeaa
\beaa
A_3 &=& \in^{bc}\Rdot_{bc34}\T^4 +\in^{bc}\Rdot_{bc3d}\T^d\\
&=& 4\rhod\T^4 +4(\etab\wedge\eta)\T^4 +  \trchb\big({}^{(h)}\T\wedge\eta)  -\atrchb\big(\eta\c {}^{(h)}\T\big)\\
&&+O(\T^4)\Ga_b\c\Ga_g+O(|{}^{(h)}\T|)\big(\bb, r^{-1}\Ga_b\big)\\
&=& 4\rhod\T^4 +4(\etab\wedge\eta)\T^4 +  \trchb\big({}^{(h)}\T\wedge\eta)  -\atrchb\big(\eta\c {}^{(h)}\T\big)+r^{-2}\Ga_b,
\eeaa
and
\beaa
A_e &=& \in^{bc}\Rdot_{bce3}\T^3 +\in^{bc}\Rdot_{bce4}\T^4 +\in^{bc}\Rdot_{bced}\T^d\\
&=& \Big( -  \trchb\dual\eta_e   + \atrchb\eta_e\Big)\T^3 +\Big( -  \trch\dual\etab_e   + \atrch\etab_e\Big)\T^4\\
&& -\frac{1}{2}\Big(4\rho  + \trch\trchb+\atrch\atrchb\Big)\dual ({}^{(h)}\T)_e\\
&& +O(\T^3)\big(\bb, r^{-1}\Ga_b\big)+O(\T^4)\big(\b, r^{-1}\Ga_g\big)+O(|{}^{(h)}\T|)r^{-1}\Ga_b\\
&=& \Big( -  \trchb\dual\eta_e   + \atrchb\eta_e\Big)\T^3 +\Big( -  \trch\dual\etab_e   + \atrch\etab_e\Big)\T^4\\
&& -\frac{1}{2}\Big(4\rho  + \trch\trchb+\atrch\atrchb\Big)\dual ({}^{(h)}\T)_e+r^{-1}\Ga_b
\eeaa
as stated. This concludes the proof of Lemma \ref{lemma:computationofthecomponentsofthetensorAforexactenergyconservationKerr:chap9}.
\end{proof}

We infer the following corollary.
\begin{corollary}\lab{cor:computationofthecomponentsofthetensorAforexactenergyconservationKerr:chap9}
We have
\bea
A_\mu=-\D_\mu\left(\Im\left(\frac{2m}{q^2}\right)\right)+r^{-1}\Ga_b.
\eea
\end{corollary}

\begin{proof}
We have
\beaa
\T^4=\frac{1}{2}, \qquad \T^3=\frac{1}{2}\frac{\De}{|q|^2}, \qquad \T_b={}^{(h)}\T_b=-a\Re(\Jk)_b.
\eeaa
Plugging in the identities of Lemma \ref{lemma:computationofthecomponentsofthetensorAforexactenergyconservationKerr:chap9}, we infer
\beaa
A_4 &=& -2\rhod\frac{\De}{|q|^2} -2(\etab\wedge\eta)\frac{\De}{|q|^2} -  a\trch\big(\Re(\Jk)\wedge\etab)  +a\atrch\big(\etab\c\Re(\Jk))+r^{-1}\Ga_g,\\
A_3 &=& 2\rhod +2(\etab\wedge\eta) -  a\trchb\big(\Re(\Jk)\wedge\eta)  +a\atrchb\big(\eta\c\Re(\Jk)\big)+r^{-2}\Ga_b,\\
A_e &=& \frac{1}{2}\Big( -  \trchb\dual\eta_e   + \atrchb\eta_e\Big)\frac{\De}{|q|^2} +\frac{1}{2}\Big( -  \trch\dual\etab_e   + \atrch\etab_e\Big)\\
&& +\frac{a}{2}\Big(4\rho  + \trch\trchb+\atrch\atrchb\Big)\dual\Re(\Jk)_e+r^{-1}\Ga_b.
\eeaa

Next, since $\trXc, \trXbc, \, \Hbc\in \Ga_g$, $\rhoc, \rhodc\in r^{-1}\Ga_g$ and $\Hc\in\Ga_b$, and using the proof of the main terms in the proof of Corollary \ref{cor:computationofthecomponentsofthetensorAforexactenergyconservationKerr}, we infer
\beaa
A_e &=& -a\Re\left(\frac{4m}{q^3}\right)\dual\Re(\Jk)_e +r^{-1}\Ga_b,\\
A_4 &=& 4\Im\left(\frac{m}{q^3}\right)\frac{\De}{|q|^2} +r^{-1}\Ga_g,\\
A_3 &=& -4\Im\left(\frac{m}{q^3}\right)+r^{-1}\Ga_g.
\eeaa
Using again the proof of the main terms in the proof of Corollary \ref{cor:computationofthecomponentsofthetensorAforexactenergyconservationKerr}, and since $\widecheck{\nab(q)}\in r\Ga_g$, $\widecheck{e_3(q})\in r\Ga_b$ and $\widecheck{e_4(q)}\in r^{-1}\Ga_g$, we infer
\beaa
A_e &=& -\nab_e\left(\Im\left(\frac{2m}{q^2}\right)\right)+r^{-1}\Ga_b,\\
A_4 &=& -e_4\left(\Im\left(\frac{2m}{q^2}\right)\right)+r^{-1}\Ga_g,\\
A_3 &=& -e_3\left(\Im\left(\frac{2m}{q^2}\right)\right)+r^{-1}\Ga_g,
\eeaa
and hence
\beaa
A_\mu=-\D_\mu\left(\Im\left(\frac{2m}{q^2}\right)\right)+r^{-1}\Ga_b
\eeaa
as stated. This concludes the proof of Corollary  \ref{cor:computationofthecomponentsofthetensorAforexactenergyconservationKerr:chap9}.
\end{proof}

Recall from Definition \ref{Def:acceptable-errors-ch9} that the terms $G_\mu\psi\c\Db^\mu\psi$ with $G_\mu\in \Ga_g$ belong to $\textrm{Good}$. Given that the terms $G_\mu \dual\psi\c\Db^\mu\psi$ with $G_\mu\in \Ga_g$ satisfy the same estimates, i.e the ones of Lemma \ref{lemma:accetableerrortermsenergymorawetzundercontrol}, we will from now on also assume that they are part of $\textrm{Good}$. In particular, in view of Corollary  \ref{cor:computationofthecomponentsofthetensorAforexactenergyconservationKerr:chap9}, this implies 
\beaa
\frac{1}{2}A_\nu \dual\psi\c\Db^\nu\psi &=& -\D_\nu\left(\Im\left(\frac{m}{q^2}\right)\right) \dual\psi\c\Db^\nu\psi +\textrm{Good}
\eeaa
which together with \eqref{eq:basiccomputationenergyestimatesKerrwithThatdeltaflsdiuh:chap9} implies
 \bea\label{eq:modifieddivergenceTmod:chap9}
 \nn \D^\mu  \PP_\mu[\That_\de, 0, 0] &=& \frac 1 2 \QQ  \c \,^{(\That_\de)} \pi   - \D_\nu\left(\Im\left(\frac{m}{q^2}\right)\right) \dual\psi\c\Db^\nu\psi  +(\That_\de-\T)^\mu \Db^\nu  \psi ^a\Rdot_{ ab   \nu\mu}\psi^b\\
  &&+ \nab_{\That_\de} \psi  \c \big(\squared_2 \psi - V \psi \big)+\textrm{Good}.
 \eea

Next, we modify the identity \eqref{eq:modifieddivergenceTmod:chap9} to cancel the second term on the RHS. To this end, we consider the following modified current
\bea
\widetilde{\PP}_\mu &:=&  \PP_\mu[\That_\de, 0, 0] + \tilde{w}\dual\psi\c\Ddot_\mu\psi,
\eea
for a scalar function $\tilde{w}=\tilde{w}(r, \cos\th)$ to be chosen below and satisfying $w=O(mar^{-3})$. We have 
\beaa
\bsplit
\D^\mu\Big[\tilde{w}\dual\psi\c\Ddot_\mu\psi\Big] &= \tilde{w}\dual\psi\c\Ddot^\mu\Ddot_\mu\psi + \tilde{w}\dual\Ddot^\mu\psi\c\Ddot_\mu\psi+\D^\mu(\tilde{w})\dual\psi\c\Ddot_\mu\psi\\
&= \tilde{w}\dual\psi\c\squared_2\psi +\D^\mu(\tilde{w})\dual\psi\c\Ddot_\mu\psi\\
&= \tilde{w}\dual\psi\c\left(V\psi - \frac{4 a\cos\th}{|q|^2}\dual \nab_\T  \psi+N\right)+\D^\mu(\tilde{w})\dual\psi\c\Ddot_\mu\psi\\
&= - \tilde{w}\frac{4 a\cos\th}{|q|^2}\nab_\T(|\dual\psi|^2)+\tilde{w}\dual\psi\c N+\D^\mu(\tilde{w})\dual\psi\c\Ddot_\mu\psi.
\end{split}
\eeaa
Since $\tilde{w}=\tilde{w}(r, \cos\th)$ and $w=O(mar^{-3})$, we have
\beaa
- \tilde{w}\frac{4 a\cos\th}{|q|^2}\nab_\T(|\dual\psi|^2) &=& - \D_\mu\left(\T^\mu \tilde{w}\frac{4 a\cos\th}{|q|^2}|\psi|^2\right)+ O(a^2 r^{-6})\Big(\T(r), \T(\cos\th)\Big)|\psi|^2\\
&&+O(a^2 r^{-5})\Div(\T)|\psi|^2\\
&=& - \D_\mu\left(\T^\mu \tilde{w}\frac{4 a\cos\th}{|q|^2}|\psi|^2\right) +r^{-5}\Ga_b|\psi|^2\\
&=& - \D_\mu\left(\T^\mu \tilde{w}\frac{4 a\cos\th}{|q|^2}|\psi|^2\right)+\textrm{Good}
\eeaa
and hence
\beaa
\bsplit
\DD^\mu\left(\widetilde{\PP}_\mu +\T^\mu \tilde{w}\frac{4 a\cos\th}{|q|^2}|\psi|^2\right) =&  \frac 1 2 \QQ  \c \,^{(\That_\de)} \pi +\D_\nu\left(-\Im\left(\frac{m}{q^2}\right)+\tilde{w}\right) \dual\psi\c\Db^\nu\psi  \\
&+(\That_\de-\T)^\mu \Db^\nu  \psi ^a\Rdot_{ ab   \nu\mu}\psi^b \\
&+ \nab_{\That_\de} \psi  \c \big(\squared_2 \psi - V \psi \big)+ \tilde{w}\dual\psi\c N +\textrm{Good}.
\end{split}
\eeaa
 Next, we make the following choice for $\tilde{w}$
\bea
\tilde{w} := \Im\left(\frac{m}{q^2}\right)= -\frac{2amr\cos\th}{|q|^4}
\eea
which yields
\bea\lab{eq:modifieddivergenceTmod:ter}
\bsplit
\DD^\mu\left(\widetilde{\PP}_\mu +\T^\mu \tilde{w}\frac{4 a\cos\th}{|q|^2}|\psi|^2\right) =&  \frac 1 2 \QQ  \c \,^{(\That_\de)} \pi +(\That_\de-\T)^\mu \Db^\nu  \psi ^a\Rdot_{ ab   \nu\mu}\psi^b   \\
& + \nab_{\That_\de} \psi  \c \big(\squared_2 \psi - V \psi \big)+ \tilde{w}\dual\psi\c N  +\textrm{Good}.
\end{split}
\eea

Next, recalling that  $\That_\de= \T+\chi_\de \Z$, we have
 \beaa
 \,^{(\That_\de)} \pi_{\mu\nu} &=& \,^{(\T)}\pi_{\mu\nu} +\chi_\de\,^{(\Z)}\pi_{\mu\nu} +\chi_\de'(r)\D_\mu(r)\Z_\nu+\chi_\de'(r)\D_\nu(r)\Z_\mu\\
 &=& \,^{(\T)}\pi_{\mu\nu} +O(r^{-2})\,^{(\Z)}\pi_{\mu\nu} +\chi_\de'(r)\D_\mu(r)\Z_\nu+\chi_\de'(r)\D_\nu(r)\Z_\mu. \eeaa
 Together with the control of $\piT$ and $\piZ$ provided by  Lemma \ref{LEMMA:DEFORMATION-TENSORS-T}, and the fact that $\widecheck{e_3(r)}\in r\Ga_b$, $\widecheck{e_4(r)}\in r\Ga_g$, $\nab(r)\in r\Ga_g$, and $\chi_\de'(r)=O(r^{-3})$, we infer
 \beaa
 \QQ  \c \,^{(\That_\de)} \pi &=& \Big[ \QQ  \c \,^{(\That_\de)} \pi \Big]_K +\textrm{Good}.
 \eeaa
 In view of the estimate for $[ \QQ  \c \,^{(\That_\de)} \pi ]_K$ in section \ref{section:generalized.energy}, we infer
 \beaa
   \big|\QQ  \c \,^{(\That_\de)} \pi\big|&\les& \mathbb{1}_{\Mntrap}\de^{-1}   \frac{ |a| }{ r^3}   | \nab_{\Z}\psi | | \nab_\Rhat \psi| +\textrm{Good}.
 \eeaa
 Also, we have $\That_\de-\T=\chi_\de(r)\Z$ and hence
 \beaa
 && (\That_\de-\T)^\mu \Db^\nu  \psi ^a\Rdot_{ ab   \nu\mu}\psi^b\\
 &=& \frac{1}{2}\chi_\de(r)\Z^\mu \in^{ab}\Rdot_{ ab   \nu\mu} \dual\psi\c\Db^\nu  \psi\\
 &=& \frac{1}{2}\chi_\de(r)\left(-\frac{1}{2}\Z^\mu \in^{ab}\Rdot_{ ab3\mu} \dual\psi\c\nab_4\psi
  -\frac{1}{2}\Z^\mu \in^{ab}\Rdot_{ ab4\mu} \dual\psi\c\nab_3\psi + \Z^\mu \in^{ab}\Rdot_{ abc\mu} \dual\psi\c\nab^c\psi \right).
 \eeaa
 Since we have
 \beaa
 \Z^3=O(ar^{-2}\De), \qquad  \Z^4=O(a), \qquad \Z^c=O(r^2)\Re(\Jk)^c,
 \eeaa
we infer, together with  the definition  \eqref{eq:DefineRdot} of $\Rdot$,  and \eqref{eq:consequenceproposition:componentsofB}, 
\beaa
\Z^\mu \in^{ab}\Rdot_{ ab3\mu} &=& \Z^4 \in^{ab}\Rdot_{ ab34}+\Z^c \in^{ab}\Rdot_{ ab3c}= O(ar^{-2})+\Ga_b,\\
\Z^\mu \in^{ab}\Rdot_{ ab4\mu} &=& \Z^3 \in^{ab}\Rdot_{ ab43}+\Z^c \in^{ab}\Rdot_{ ab4c}=O(ar^{-4}\De)+\Ga_g,\\
\Z^\mu \in^{ab}\Rdot_{ abc\mu} &=& \Z^4 \in^{ab}\Rdot_{ abc4}+\Z^3 \in^{ab}\Rdot_{ abc3}+\Z^d \in^{ab}\Rdot_{ abcd}\\
&=& O(r^{-1})+\Ga_b,
\eeaa 
 which implies, since $\chi_\de= \frac{a}{r^2+a^2} \chi_0\Big( \de^{-1} \frac{\TT}{r^3} \Big)$, 
  \beaa
  \big|(\That_\de-\T)^\mu \Db^\nu  \psi ^a\Rdot_{ ab   \nu\mu}\psi^b\big| &\les& \frac{|a|}{r^3}\mathbb{1}_{\Mntrap}\Big(|\nab_{\Rhat}\psi|+|\nab_\T\psi|+|\nab\psi|\Big)|\psi|+|\textrm{Good}|.
 \eeaa
   Finally, using equation \eqref{eq:Gen.RW}, we have
   \beaa
   \nab_{\That_\de} \psi  \c \big(\squared_2 \psi - V \psi \big)   &=&- \frac{4 a\cos\th}{|q|^2}\chi_\de \nab_\Z\psi   \c  \dual \nab_\T  \psi+ \nab_{\That_\de} \psi  \c  N 
   \eeaa
      where we have the crucial cancellation $\nab_\T \psi \c \dual \nab_\T \psi=0$. 
    
      We summarize the result in the following.
      \begin{lemma}
    \lab{lemma:CurrentTmod:chap9}
    Consider the modified current 
    \beaa
    \widetilde{\PP}_\mu &=&  \PP_\mu[\That_\de, 0, 0] + \tilde{w}\dual\psi\c\Ddot_\mu\psi, \qquad \tilde{w} = \Im\left(\frac{m}{q^2}\right),
    \eeaa
    where  the vectorfield $\That_\de$ is given by $\That_\de= \T+\chi_\de(r) \Z$ for $\de=\frac{1}{10}$ and $|a|/m\ll 1$ small enough.   Then,  we have on $\MM$,
    \beaa
  &&\left|   \D^\mu\left(\widetilde{\PP}_\mu +\T^\mu \tilde{w}\frac{4 a\cos\th}{|q|^2}|\psi|^2\right) -\Big(\nab_{\That_\de} \psi + \tilde{w}\dual\psi\Big) \c  N \right|\\ 
  &\les&  \mathbb{1}_{\Mntrap}  \left(\de^{-1}   \frac{ |a| }{ r^3}  |\nab_\Rhat \psi||\nab_\Z\psi|+\frac{|a|}{r^4} | \nab_\T  \psi|  | \nab_\Z\psi | + \frac{|a|m}{r^4}\Big[|\nab_{\Rhat}\psi|+|\nab_\T\psi|+|\nab\psi|\Big]|\psi| \right) \\
  && +|\textrm{Good}|.
    \eeaa
    \end{lemma} 
 
 Integrating the above inequality on $\MM=\MM(\tau_1, \tau_2) $ and applying the divergence theorem we deduce, in view of the definition of $\Mor^{ax}_{deg}[\psi](\tau_1, \tau_2)$,  
 \beaa
 &&\int_{\Si(\tau_2)}\left(\widetilde{\PP}_\mu +\T^\mu \tilde{w}\frac{4 a\cos\th}{|q|^2}|\psi|^2\right)\c N_\Si +\int_{\Si_*(\tau_1, \tau_2)}\left(\widetilde{\PP}_\mu +\T^\mu \tilde{w}\frac{4 a\cos\th}{|q|^2}|\psi|^2\right)\c N_{\Si_*}\\ 
 &&+\int_{\AA(\tau_1, \tau_2)}\left(\widetilde{\PP}_\mu +\T^\mu \tilde{w}\frac{4 a\cos\th}{|q|^2}|\psi|^2\right)\c N_{\AA} \\
 &\les&  \int_{\Si(\tau_1)}\left(\widetilde{\PP}_\mu +\T^\mu \tilde{w}\frac{4 a\cos\th}{|q|^2}|\psi|^2\right)\c N_\Si+ \frac{|a|}{m}\Mor^{ax}_{deg}[\psi](\tau_1, \tau_2)   +\left|\int_{\MM(\tau_1, \tau_2)}  \nab_{\That_\de } \psi  \c N\right|\\
 &&+\int_{\MM(\tau_1, \tau_2)}|N|^2+\int_{\MM(\tau_1, \tau_2)}|\textrm{Good}|.
 \eeaa
 The rest of the proof is a simple adaptation of the one of  Proposition \ref{proposition:Energy1}, using Lemma \ref{Lemma:QQ(That,NSi):chap9} for the sign of boundary terms, and Remark \ref{rmk::accetableerrortermsenergymorawetzundercontrolafsodsudfh} to control the term $\textrm{Good}$. This concludes the proof of  Proposition \ref{proposition:Energy1:perturbation}.


\subsubsection{Proof of Proposition 
\ref{prop:recoverEnergyMorawetzwithrweightfromnoweight:perturbation}}


The proof of Proposition 
\ref{prop:recoverEnergyMorawetzwithrweightfromnoweight:perturbation} is a simple adaptation of the one of Proposition 
\ref{prop:recoverEnergyMorawetzwithrweightfromnoweight}, where the extra error terms are tracked as in the proof of Propositions \ref{proposition:Morawetz1-step1:perturbation} and \ref{proposition:Energy1:perturbation}.


\subsection{Proof of Proposition \ref{prop:morawetz-higher-order:perturbation} and 
Lemma \ref{lemma:lowerboundPhizoutsideMtrap:perturbation}}
\lab{sec:proofofresultschapter8inperturbationofKerr}


Let   the  second order differential operators  $\SS_\aund$, $\aund=1,2,3,4$,   introduced in  Definition \ref{definition:symmetry-tensors-pert}. Given a $\mathfrak{s}_2$ tensor  $\psi$ solution  of the equation  \eqref{eq:Gen.RW}, we consider in this section the  commuted $\mathfrak{s}_2$ tensors
\bea
\psia:=\SS_\aund \psi, \qquad \aund=1,2,3,4, 
\eea
that satisfy according to Lemma \ref{lemma-commuted-equ-Naund:perturbation} the equation \eqref{eq:waveeqfor-psia:perturbation}, i.e. 
  \bea
  \squared_2\psia -V\psia&=&- \frac{ 4 a\cos\th}{|q|^2} \dual \nab_\T  \psia+N_\aund
  \eea
  where $N_\aund$ satisfies \eqref{eq:estimaeforMundinperturbationsofKerr}.

    
    \subsubsection{Acceptable error terms for $\SS$-Morawetz}
        

To \eqref{eq:waveeqfor-psia:perturbation}, we associate the following generalized energy-momentum tensor
\beaa
\QQ_{\mu\nu}(\psi_\aund, \psi_\bund)&=& \Db_\mu  \psi_\aund \c \Db_\nu \psi _\bund
          -\frac 12 \g_{\mu\nu} \left(\g^{\a\b} \Db_\a \psi_\aund \c\Db_\b  \psi_\bund+ V\psi_\aund \c \psi_\bund\right)\\
          &=&
           \Db_\mu  \psi _\aund\c \Db _\nu \psi_\bund    -\frac 12 \g_{\mu\nu}  \LL[\psi_\aund, \psi_\bund], \\
          \LL[\psi_\aund, \psi_\bund]&=&\g^{\a\b} \Db_\a \psi_\aund\c\Db_\b  \psi _\bund+ V\psi_\aund \c \psi_\bund.
          \eeaa
 
When revisiting the proofs of Chapter  \ref{chapter-proof-mor-2} in the context of admissible perturbations $\MM$ of Kerr, we generate additional terms. We introduce below acceptable terms.
 
  \begin{definition}[Acceptable error terms]
 \lab{Def:acceptable-errors-ch9:SSversion}
 The following quantity 
 \beaa
 F^{\aund\bund\,\mu\nu}\QQ_{\mu\nu}(\psi_\aund, \psi_\bund)+G^{\aund\bund\,\mu} \psi_\aund\c\Ddot_\mu\psi_\bund +H^{\aund\bund}\psia\c\psib +I[\psi]
 \eeaa
 is said to be of the acceptable type, and denoted $\textrm{Good}_\SS$, if: 
 
\begin{itemize}
\item $F^{\aund\bund}_{44}\in \Ga_g$, and all other components of $F^{\aund\bund}_{\mu\nu}$ belong to $\Ga_b$.

\item All components of $G^{\aund\bund}_\mu$ belong to $\Ga_g$.

\item $H^{\aund\bund}\in r^{-1}\Ga_g$.

\item $I[\psi]=\dk^{\leq 3}(\Ga_g\c\psi)\c\nab_3(\dk^{\leq 2}\psi)+r^{-1}\dk^{\leq 3}(\Ga_g\c\psi)\c\dk^{\leq 3}\psi$. 
\end{itemize} 
 \end{definition} 

The justification for Definition \ref{Def:acceptable-errors-ch9} is provided by the following lemma.
\begin{lemma}\lab{lemma:accetableerrortermsenergymorawetzundercontrol:SSversion}
Assume that the quantity 
 \beaa
 F^{\aund\bund\,\mu\nu}\QQ_{\mu\nu}(\psi_\aund, \psi_\bund)+G^{\aund\bund\,\mu} \psi_\aund\c\Ddot_\mu\psi_\bund +H^{\aund\bund}\psia\c\psib +I[\psi]
 \eeaa
 is of the acceptable type in the sense of Definition \ref{Def:acceptable-errors-ch9:SSversion}. Then, it satisfies the following estimate
\beaa
\int_{\MM_{trap}(\tau_1, \tau_2)}\Big|F^{\aund\bund\,\mu\nu}\QQ_{\mu\nu}(\psi_\aund, \psi_\bund)+G^{\aund\bund\,\mu} \psi_\aund\c\Ddot_\mu\psi_\bund +H^{\aund\bund}\psia\c\psib +I[\psi]\Big|^2 \les \ep\sup_{[\tau_1, \tau_2]}E^2[\psi](\tau)
\eeaa
and 
\beaa
\int_{\Mntrap(\tau_1, \tau_2)}\Big|F^{\aund\bund\,\mu\nu}\QQ_{\mu\nu}(\psi_\aund, \psi_\bund)+G^{\aund\bund\,\mu} \psi_\aund\c\Ddot_\mu\psi_\bund +H^{\aund\bund}\psia\c\psib +I[\psi]\Big|^2 \les \ep B^2_\de[\psi](\tau_1, \tau_2).
\eeaa
\end{lemma}

\begin{remark}\lab{rmk::accetableerrortermsenergymorawetzundercontrolafsodsudfh:SSversion}
Recall from Definition \ref{Def:acceptable-errors-ch9:SSversion} that terms of the acceptable type are denoted $\textrm{Good}_\SS$. In view of Lemma \ref{lemma:accetableerrortermsenergymorawetzundercontrol:SSversion}, we infer that such term satisfy the following estimate
\beaa
\int_{\MM}|\textrm{Good}_\SS| &\les& \ep\left(\sup_{[\tau_1, \tau_2]}E^2[\psi](\tau)+ B^2_\de[\psi](\tau_1, \tau_2)\right).
\eeaa
\end{remark}

\begin{proof}
The proof is analogous to the one of Lemma \ref{lemma:accetableerrortermsenergymorawetzundercontrol} noticing that $\psi_\aund$ is schematically of the type $\dk^{\leq 2}\psi$.
\end{proof}

Next, we introduce the linearization of the quantities $F_{\mu\nu}^{\aund\bund}$, $G_\mu^{\aund\bund}$ and $H^{\aund\bund}$. 
\begin{definition}\lab{def:linearizationofMorawetzquantitiesinperturbationofKerr:chap9:SSversion}
Let $F_{\mu\nu}^{\aund\bund}$, $G_\mu^{\aund\bund}$ and $H^{\aund\bund}$. We define their linearization as follows
\beaa
F_{\mu\nu}^{\aund\bund}=(F_{\mu\nu}^{\aund\bund})_K+\widecheck{F_{\mu\nu}^{\aund\bund}}, \qquad G_{\mu}^{\aund\bund}=(G_{\mu}^{\aund\bund})_K+\widecheck{G_{\mu}^{\aund\bund}}, \qquad H^{\aund\bund}=H^{\aund\bund}_K+\widecheck{H^{\aund\bund}},
\eeaa
where: 
\begin{enumerate}
\item the quantities 
\beaa
(F_{44}^{\aund\bund})_K, \qquad (F_{34}^{\aund\bund})_K, \qquad (F_{33}^{\aund\bund})_K, \qquad (G_4^{\aund\bund})_K, \qquad (G_3^{\aund\bund})_K, \qquad H_K^{\aund\bund},
\eeaa
are given as explicit functions of $(r, \cos\th)$ coinciding with the corresponding expressions in Kerr,

\item the quantities 
\beaa
(F_{4a}^{\aund\bund})_K, \qquad (F_{3a}^{\aund\bund})_K, \qquad (G_a^{\aund\bund})_K, 
\eeaa
are given as the 1-form $\Re(\Jk)_a$ multiplied by explicit functions of $(r, \cos\th)$ coinciding with the corresponding expressions in Kerr,

\item the quantity $(F^{\aund\bund}_{ab})_K$ is given by the the symmetric 2-tensor $\ga_{ab}$ multiplied by explicit functions of $(r, \cos\th)$ coinciding with the corresponding expressions in Kerr.
\end{enumerate}
\end{definition}

In view of Definition \ref{def:linearizationofMorawetzquantitiesinperturbationofKerr:chap9:SSversion}, we can decompose the above expressions in their main part and error terms as follows
\bea\lab{eq:deompostionofMorawetzquantitiesinperturvbationfadsfaosdiuasdfoiahdsflaiuhf:SSversion}
\bsplit
 & F^{\aund\bund\,\mu\nu}\QQ_{\mu\nu}(\psi_\aund, \psi_\bund)+G^{\aund\bund\,\mu} \psi_\aund\c\Ddot_\mu\psi_\bund +H^{\aund\bund}\psia\c\psib +I[\psi]\\
 =& \Big[F^{\aund\bund\,\mu\nu}\QQ_{\mu\nu}(\psi_\aund, \psi_\bund)+G^{\aund\bund\,\mu} \psi_\aund\c\Ddot_\mu\psi_\bund +H^{\aund\bund}\psia\c\psib\Big]_K+\err,\\
& \Big[F^{\aund\bund\,\mu\nu}\QQ_{\mu\nu}(\psi_\aund, \psi_\bund)+G^{\aund\bund\,\mu} \psi_\aund\c\Ddot_\mu\psi_\bund +H^{\aund\bund}\psia\c\psib\Big]_K\\
  =& (F^{\aund\bund\,\mu\nu})_K\QQ_{\mu\nu}[\psi]+(G^{\aund\bund\,\mu})_K \psi\c\Ddot_\mu\psi +H^{\aund\bund}_K|\psi|^2,\\
 \err =& \widecheck{F^{\aund\bund\,\mu\nu}}\QQ_{\mu\nu}[\psi]+\widecheck{G^{\aund\bund\,\mu}} \psi\c\Ddot_\mu\psi +\widecheck{H^{\aund\bund}}|\psi|^2 +I[\psi].
 \end{split}
\eea 
The proof of section \ref{sec:proofofresultschapter8inperturbationofKerr} will rely in particular on showing that the extra terms in perturbations of Kerr appearing in the various divergence identities involved in energy and Morawetz estimates are of the  acceptable type in the sense of Definition \ref{Def:acceptable-errors-ch9:SSversion}, i.e, that $\err\in\textrm{Good}_\SS$  in \eqref{eq:deompostionofMorawetzquantitiesinperturvbationfadsfaosdiuasdfoiahdsflaiuhf:SSversion}.


\subsubsection{Preliminaries for $\SS$-Morawetz}


We start with some preliminaries for the $\SS$-Morawetz estimate. Let $\X$ be a double-indexed collection of vector fields $\X=\{ X^{\underline{a} \underline{b}}\}$, $\bold{w}$ be a double-indexed collection of functions $\bold{w}=\{ w^{\underline{ab}} \}$, and $\M=\{M^{\aund\bund} \}$   a double-indexed collection of  $1$-forms, all symmetric in the indices $\aund, \bund$. Assume furthermore that  $\X$  are of the form
 \beaa
    X^{\aund\bund}=  (X^{\aund\bund})^3  e_3 +(X^{\aund\bund})^4 e_4.
    \eeaa
  Define
   \beaa
 \PP_\mu[\bold{X}, \bold{w}, \M] &:=& \QQ_{\mu\nu}(\psi_\aund, \psi_\bund) X^{\underline{ab} \nu}+\frac 1 2  w^{\underline{ab}} \, \Db_\mu  \psia\c  \psib-\frac 1 4 (\partial_\mu w^{\underline{ab}})\psia\c\psib  +\frac 1 4 M^{\aund\bund}_\mu\psi_\aund\c \psi_\bund.
  \eeaa

We make the following choices for $(\bold{X}, \bold{w}, \M)$, consistent with the ones in Chapter \ref{chapter-proof-mor-2}:
\begin{enumerate}
\item  $\X$ of the type
\beaa
X^{\aund\bund}=\FF^{\aund\bund}\frac{(r^2+a^2)}{\De}\Rhat,
\eeaa
where $\FF^{\aund\bund}=\FF^{\aund\bund}(r)$ are such that $\FF^{\aund\bund}/\De$ is a smooth function of $r$, and $(r\pr_r)^k\FF^{\aund\bund}(r)$ is uniformly bounded on $\MM$. In particular, there holds $X^{\aund\bund}= (X^{\aund\bund})^3  e_3 +(X^{\aund\bund})^4 e_4$ as anticipated, with $(X^{\aund\bund})^3$ and $(X^{\aund\bund})^4$ smooth functions of $r$ such that $(r\pr_r)^k(X^{\aund\bund})^3$ and $(r\pr_r)^k(X^{\aund\bund})^4$ are uniformly bounded on $\MM$. 

\item The scalar functions $w^{\aund\bund}(r)$ are such that $w^{\aund\bund}/\De$ is a smooth function of $r$ and $(r\pr_r)^k(rw^{\aund\bund}(r))$ is uniformly bounded on $\MM$.

\item  $M^{\aund\bund}$ are of the type
\beaa
M^{\aund\bund}=v^{\aund\bund}(r)\frac{(r^2+a^2)}{\De}\Rhat
\eeaa 
where the scalar functions $v^{\aund\bund}(r)$ are such that $v^{\aund\bund}/\De$ is a smooth function of $r$ and $(r\pr_r)^k(r^{\frac{5}{2}}v^{\aund\bund}(r))$ is uniformly bounded on $\MM$.
\end{enumerate}

We have the following lemma.
\begin{lemma}\lab{lemma:decompositionofthedivergenceinerrorandtermKersoiduhf:SSversion}
Using the decomposition \eqref{eq:deompostionofMorawetzquantitiesinperturvbationfadsfaosdiuasdfoiahdsflaiuhf:SSversion}, we have
\beaa
  \D^\mu  \PP_\mu[\bold{X}, \bold{w}, \M] &=& \Bigg[\frac 1 2 \QQ_{\aund\bund}\c{}^{(X^{\aund\bund})}\pi - \frac 1 2 X^{\aund\bund}( V )\psia\c\psib+\frac 12  w^{\aund\bund} \LL[\psia, \psib] \\
  &&-\frac 1 4\psia\c\psib   \square_\g  w^{\aund\bund} + \frac 1 4  \Div(\psia\c\psib M^{\aund\bund}\big)\Bigg]_K\\
     && - \big(\rhod +\etab\wedge\eta\big)\nab_{(X^{\aund\bund})^4e_4-(X^{\aund\bund})^3e_3}  \psi_\aund\c\dual\psi_\bund \\
&& - \frac{1}{2}\Im\Big(\tr\Xb H(X^{\aund\bund})^3 +\tr X\Hb(X^{\aund\bund})^4\Big)\c\nab\psi_\aund\c\dual\psi_\bund\\
  &&+  \left(\nab_{X^{\aund\bund}}\psia+\frac 1 2   w^{\aund\bund} \psia\right)\c\big(\squared_2\psib-V\psib\big) +\textrm{Good}_\SS
 \eeaa
where $\textrm{Good}_\SS$ is given by Definition \ref{Def:acceptable-errors-ch9:SSversion}.
\end{lemma}

\begin{proof}
The proof is analogous to the one of Lemma \ref{lemma:decompositionofthedivergenceinerrorandtermKersoiduhf}.
\end{proof}

 
 \subsubsection{Integration by parts   identities}
 

We now generalize the results of section \ref{sec:integrationbypartsidentitiesSS} to admissible perturbations of Kerr. We start with the following analog of Definition \ref{def:boundarytermsinintegrationbypartsareok} that will be useful to take care of boundary terms in the integrations by parts.  
\begin{definition}\label{def:boundarytermsinintegrationbypartsareok:SSversion:chap9}
We denote by $M(\psi)$   quadratic terms of the following type
\beaa
\nab_T \psi \c \SS_{\aund}\psi, \qquad |q|^2\nab  \psi \c \nab \nab_T\psi, \qquad |q|^2\nab\psi\c\nab\nab_Z\psi,
\eeaa
Also, we denote by $M(\nab_{\Rhat}\psi)$ denote quadratic terms of the following type
\beaa
\nab_T \nab_{\Rhat}\psi \c \SS_{\aund}\nab_{\Rhat}\psi, \qquad |q|^2\nab  \nab_{\Rhat}\psi \c \nab \nab_T\nab_{\Rhat}\psi, \qquad |q|^2\nab  \nab_{\Rhat}\psi \c \nab \nab_Z\nab_{\Rhat}\psi.
\eeaa
\end{definition}

\begin{remark}
Notice the following pointwise  estimates for $M(\psi)$ and $M(\nab_{\Rhat}\psi)$:  
\bea
\bsplit
|M(\psi)| &\les |(\nab_T, \dkb)^{\leq 1}\psi||(\nab_T, \dkb)^{\leq 2}\psi|, \\
|M(\nab_{\Rhat}\psi)| &\les |\nab_{\Rhat}^{\leq 1}(\nab_T, \dkb)^{\leq 1}\psi||\nab_{\Rhat}^{\leq 1}(\nab_T, \dkb)^{\leq 2}\psi|.
\end{split}
\eea
\end{remark}

\begin{lemma}\label{Lemma:integrationbypartsSS_3SS_4:SSversion:chap9} 
For any function $H=H(r)$, the following identities hold true:
\beaa
H \OO(\psi) \c \SS_1 \psi  &=&H|q|^2| \nab \nab_\T\psi|^2  -O(ar^{-2})H(\dk^{\leq 2}\psi)^2\\
&& +|q|^2 \Ddot_\b(H|q|^{-2}O^{\a\b}\Ddot_\a \psi \c \SS_1 \psi ) +\Ddot_\mu(HM(\psi)\T^\mu)\\
&& +r^2(|H|+|H'|)\textrm{Good}_\SS, \\
H \nab_\Rhat \OO\psi \c \nab_\Rhat \SS_1\psi &=& H|q|^2|\nab \nab_T\nab_\Rhat\psi|^2 -O(a r^{-2})H(\nab_\Rhat\dk^{\leq 2}\psi)^2 -O(a r^{-2})H(\dk^{\leq 2}\psi)^2\\
&&+|q|^2 \Ddot_\b(H|q|^{-2}O^{\a\b}\Ddot_\a \nab_\Rhat\psi \c \SS_1 \nab_\Rhat\psi ) +\Ddot_\mu(HM(\nab_\Rhat\psi)\T^\mu)\\
&&+r^2(|H|+|H'|)\textrm{Good}_\SS, \\
  H \nab_\Rhat\OO(\psi) \c \SS_1 \psi &=&H |q|^2 \nab\nab_\T\nab_\Rhat \psi \c \nab \nab_\T \psi -O(a r^{-2})H(\nab_\Rhat\dk^{\leq 2}\psi)^2\\
  && -O(a r^{-2})H(\dk^{\leq 2}\psi)^2 +|q|^2 \Ddot_\b(H|q|^{-2}O^{\a\b}\Ddot_\a \nab_\Rhat\psi \c \SS_1 \psi )\\
  && +\Ddot_\mu(HM(\nab_\Rhat\psi)\T^\mu) +r^2(|H|+|H'|)\textrm{Good}_\SS,\\
  H\OO(\psi) \c \nab_\Rhat \SS_1\psi &=&H |q|^2 \nab\nab_\T\nab_\Rhat \psi \c \nab \nab_\T \psi -O(a r^{-2})H(\nab_\Rhat\dk^{\leq 2}\psi)^2\\
  && -O(a r^{-2})H(\dk^{\leq 2}\psi)^2 +|q|^2 \Ddot_\b(H|q|^{-2}O^{\a\b}\Ddot_\a\psi \c \SS_1\nab_\Rhat\psi )\\ 
  && +\Ddot_\mu(HM(\nab_\Rhat\psi)\T^\mu) +r^2(|H|+|H'|)\textrm{Good}_\SS.
 \eeaa
In all the above,  $M(\psi)$ and $M(\nab_\Rhat\psi)$ denote the quadratic expressions in $\psi$ and its derivatives of Definition \ref{def:boundarytermsinintegrationbypartsareok}.
 \end{lemma}
 
\begin{proof}
First, we have in view of Lemma \ref{LEMMA:SYMM-OPERATORS}
\beaa
 |q|^2\Ddot_\b(|q|^{-2}O^{\a\b}\Ddot_\a \psi) &=& \OO(\psi)+r\Ga_b\dk\psi.
 \eeaa
 We obtain for any $\Phi$,
\bea\label{eq:intermediate-integration-by-parts-lemma:chap9}
\nn |q|^{-2} \OO(\psi) \c \Phi &=& \Ddot_\b(|q|^{-2}O^{\a\b}\Ddot_\a \psi) \c \Phi +r^{-1}\Ga_b\Phi\c\dk\psi\\ 
\nn &=& \Ddot_\b(|q|^{-2}O^{\a\b}\Ddot_\a \psi \c \Phi )-|q|^{-2}O^{\a\b}\Ddot_\a \psi \c \Ddot_\b\Phi \\&& +r^{-1}\Ga_b\Phi\c\dk\psi.
\eea
Applying the above to $\Phi=\SS_1 \psi=\nab_T\nab_T \psi$, and using that 
\beaa
[\nab, \nab_T]\psi=O(ar^{-4}) \dkb^{\leq 1}\psi+r^{-1}\dk(\Ga_b\c\psi)
\eeaa 
from Corollary \ref{cor:basicpropertiesLiebTfasdiuhakdisug:chap9},  $\tr\piT\in \Ga_g$ from Lemma \ref{LEMMA:DEFORMATION-TENSORS-T}, and $\nab(r)\in r\Ga_g$ and $\T(r)\in r\Ga_b$, we infer, proceeding as in Lemma \ref{Lemma:integrationbypartsSS_3SS_4},
\beaa
 H\OO(\psi) \c \SS_1 \psi&=&-H|q|^{2}\nab  \psi \c \nab\nab_\T \nab_\T\psi+H|q|^2 \Ddot_\b(|q|^{-2}O^{\a\b}\Ddot_\a \psi \c \SS_1 \psi )\\
 &&+Hr\dk(\Ga_b\c\psi)\nab_\T\dk^{\leq 1}\psi\\
    &=&H|q|^2|\nab \nab_\T\psi|^2  -O(ar^{-2})H(\dk^{\leq 2}\psi)^2+|q|^2 \Ddot_\b(H|q|^{-2}O^{\a\b}\Ddot_\a \psi \c \SS_1 \psi )\\
    && +\Ddot_\mu(HM(\psi)\T^\mu) +H|q|^2\nab\psi\c[\nab_\T, \nab]\nab_\T\psi +H|q|^2[\nab_\T, \nab]\psi\c\nab\nab_\T\psi\\
    && + H|q|^{2}\nab  \psi \c \nab\nab_\T\psi\D_\mu\T^\mu +H'(r)\big(\T(r), r\nab(r)\big)\dk^{\leq 1}\psi\nab_\T\dk^{\leq 1}\psi\\
    && +Hr\dk(\Ga_b\c\psi)\nab_\T\dk^{\leq 1}\psi\\
    &=& |q|^2|\nab \nab_\T\psi|^2  -O(ar^{-2})(\dk^{\leq 2}\psi)^2+|q|^2 \Ddot_\b(|q|^{-2}O^{\a\b}\Ddot_\a \psi \c \SS_1 \psi )\\
    && +\Ddot_\mu(HM(\psi)\T^\mu) +(|H|+|H'|)r^2\dk^{\leq 2}(\Ga_g\c\psi)\nab_\T\dk^{\leq 1}\psi\\
    &=& H|q|^2|\nab \nab_\T\psi|^2  -O(ar^{-2})H(\dk^{\leq 2}\psi)^2+|q|^2 \Ddot_\b(H|q|^{-2}O^{\a\b}\Ddot_\a \psi \c \SS_1 \psi )\\
    && +\Ddot_\mu(HM(\psi)\T^\mu) +r^2(|H|+|H'|)\textrm{Good}_\SS
\eeaa
as stated, where we have used the fact that 
\beaa
\dk^{\leq 2}(\Ga_g\c\psi)\nab_\T\dk^{\leq 1}\psi=\dk^{\leq 2}(\Ga_g\c\psi)\nab_3\dk^{\leq 1}\psi+r^{-1}\dk^{\leq 2}(\Ga_g\c\psi)\dk^{\leq 2}\psi=\textrm{Good}_\SS.
\eeaa

Similarly, using in addition $[\OO, \nab_\Rhat]\psi = O(ar^{-2})\dk^{\leq 1}\psi +r\Ga_b\dk^{\leq 2}\psi$ from Corollary \ref{cor:commutatornabRhatwith|q|nab}, and $[\nab_\Rhat, \SS_1] = O(amr^{-4})\dk^{\leq 1}\psi +\dk^{\leq 2}(\Ga_b\c\psi)$ from Lemma \ref{lemma:commutationofnabRhatwithnabTandnabZ}, we have
\beaa
&& H \nab_\Rhat \OO\psi \c \nab_\Rhat \SS_1\psi\\ 
&=& H \OO\nab_\Rhat \psi \c  \SS_1\nab_\Rhat\psi +HO(a r^{-2}) \dkb^{\leq1}\psi\c \nab_\Rhat\SS_1\psi+H \OO\nab_\Rhat\psi \c O(a r^{-4}) \dk^{\leq1}\psi\\
&& +Hr\Ga_b\dk^{\leq 2}\psi \dk^{\leq 3}\psi\\
&=& H|q|^2|\nab \nab_\T\nab_\Rhat\psi|^2-O(ar^{-2})H\nab  \nab_\Rhat\psi \c\dkb^{\leq 1}\nab_\T\nab_\Rhat\psi\\
&& +|q|^2 \Ddot_\b(H|q|^{-2}O^{\a\b}\Ddot_\a \nab_\Rhat\psi \c \SS_1 \nab_\Rhat\psi ) -\Ddot_\mu(|q|^2H\nab  \nab_\Rhat\psi \c \nab \nab_\T\nab_\Rhat\psi \T^\mu)\\
&& +HO(a r^{-2}) \dkb^{\leq1}\psi\c \nab_\Rhat\SS_1\psi +H \OO\nab_\Rhat\psi \c O(a r^{-4}) \dk^{\leq1}\psi +(|H|+|H'|)r^2\Ga_g\dk^{\leq 2}\psi \dk^{\leq 3}\psi
\eeaa
which can be written as 
\beaa
H \nab_\Rhat\OO\psi \c \nab_\Rhat\SS_1\psi &=& H|q|^2|\nab \nab_\T\nab_\Rhat\psi|^2+|q|^2 \Ddot_\b(H|q|^{-2}O^{\a\b}\Ddot_\a \nab_\Rhat\psi \c \SS_1 \nab_\Rhat\psi )\\
&&+\Ddot_\mu(HM(\nab_\Rhat\psi)\T^\mu) -O(ar^{-2})H(\nab_\Rhat\dk^{\leq 2}\psi)^2 -O(a r^{-2})H(\dk^{\leq 2}\psi)^2\\
&& +r^2\textrm{Good}_\SS,
\eeaa
as stated. 

The mixed products $H  \nab_\Rhat\OO\psi \c  \SS_1\psi$ and $H \OO\psi \c  \nab_\Rhat\SS_1\psi$ are treated similarly. 
This concludes the proof of Lemma \ref{Lemma:integrationbypartsSS_3SS_4:SSversion:chap9}. 
\end{proof}

We use the above lemma to derive the following analog of Lemma \ref{general-computations-for-BB}. 

\begin{lemma}\label{general-computations-for-BB:chap9} 
Let $\Phi_\aund=\SS_\aund \Phi$ for some\footnote{We will apply it to $\Phi=\psi$ and $\Phi=\nab_r\psi$.} $\Phi$, and let $\mathcal{Y}^\aund$ be some coefficients only depending on $r$, such that
\beaa
\mathcal{Y}^{1} =\de_0 \mathcal{Y}, \qquad \mathcal{Y}^{4}=\mathcal{Y}.
\eeaa
 Then for $\LL^{\aund}$ given by \eqref{eq:remark:choiceofLL-weak}, i.e.
\beaa
 \LL^{1} =  \de_0, \qquad    \LL^{2} =0, \qquad \LL^3=\LL^4=1,
\eeaa
 we have
 \beaa
\big(\mathcal{Y}^{\aund}\Phi_\aund\big)\c \big( \LL^{\aund}\Phi_\aund\big)&=& \mathcal{Y}\Big(\de_0^2 |\SS_1\Phi|^2+ |\OO\Phi|^2 +2\de_0|q|^2| \nab \nab_T\Phi|^2\Big)\\
&& -O(a)(|\mathcal{Y}|+|\mathcal{Y}^2|+|\mathcal{Y}^3|)(\dk^{\leq 2}\psi)^2+\text{Bdr}[\Phi] +r^2(|\mathcal{Y}|+|\mathcal{Y}'|)\textrm{Good}_\SS,
\eeaa
where the boundary term is given by
\beaa
\text{Bdr}[\Phi]&=&\Ddot_\mu\Big(\de_0\mathcal{Y}  r^2 M(\Phi)\T^\mu\Big)+|q|^2 \Ddot_\b\Big(2\de_0 |q|^{-2}O^{\a\b}\Ddot_\a \Phi \c \mathcal{Y} \SS_1 \Phi \Big).
\eeaa
\end{lemma}

\begin{proof}
As in the proof of Lemma \ref{general-computations-for-BB}, we have
\beaa
\big(\mathcal{Y}^{\aund}\Phi_\aund\big)\c \big( \LL^{\aund}\Phi_\aund\big)&=&\de_0^2 \mathcal{Y} |\SS_1\Phi|^2+\mathcal{Y} |\OO\Phi|^2 +2\de_0 \mathcal{Y}(\SS_1\Phi \c \OO \Phi )\\
&& +O(a)(|\mathcal{Y}|+|\mathcal{Y}^2|+|\mathcal{Y}^3|)(\dk^{\leq 2}\psi)^2.
\eeaa
Using Lemma \ref{Lemma:integrationbypartsSS_3SS_4:SSversion:chap9}, we obtain the desired identity. 
\end{proof}


\subsubsection{Proof of Proposition \ref{prop:morawetz-higher-order:perturbation}}


We choose $(\X, \w, \M)$ as in Chapter \ref{chapter-proof-mor-2}:
\begin{enumerate}
\item  $\X$ is given by
\beaa
\bsplit
X^{\aund\bund} &=\FF^{\aund\bund}\frac{(r^2+a^2)}{\De}\Rhat, \qquad \FF^{\aund\bund}=-zhf^{\aund\bund}, \\ 
z &=z_0-\de_0z_0^2. \qquad z_0=\frac{\De}{(r^2+a^2)^2}, \qquad \de_0>0, \qquad h=\frac{(r^2+a^2)^4}{r(r^2-a^2)},\\
f^{\underline{a}\underline{b}}&=\frac 1 2 \big(\RRtp^{\aund}\LL^{\bund}+\RRtp^{\bund}\LL^{\aund} \big), \qquad \RRtp^{\aund} =      \pr_r\Big( \frac{z}{\De} \RRa\Big),\\
\LL^{1} &=  \de_0, \qquad    \LL^{2} =0, \qquad \LL^3=\LL^4=1,\\
\RR^1 &=-(r^2+a^2)^2, \qquad \RR^2 = -2(r^2+a^2), \qquad \RR^3 =-1, \qquad \RR^4=\De.
\end{split}
\eeaa
 
\item $\w$ is given  by
\beaa
w^{\aund\bund}=-z\pr_r(hf^{\aund\bund}).
\eeaa

\item  $\M$ is given by
\beaa
M^{\aund\bund}=v^{\aund}\LL^{\bund}\frac{(r^2+a^2)}{\De}\Rhat, \qquad v^1=\de_0v, \qquad v^2=v^3=0, \qquad v^4=v,
\eeaa 
where the scalar function $v(r)$ is the one of Lemma \ref{Le:Crucialpositivity-a=0} and satisfies in particular  $v(r)= O( m^{1/2}\Delta r^{-9/2})$.
\end{enumerate}

Next, we introduce the expression
    \bea\lab{definition-EE-gen:chap9:bisrepetita}
      \begin{split}
\EE_K[\X, \w, \M] &:=  \Bigg[\frac 1 2 \QQ_{\aund\bund}\c{}^{(X^{\aund\bund})}\pi - \frac 1 2 X^{\aund\bund}( V )\psia\c\psib+\frac 12  w^{\aund\bund} \LL[\psia, \psib] \\
  &-\frac 1 4\psia\c\psib   \square_\g  w^{\aund\bund} + \frac 1 4  \Div(\psia\c\psib M^{\aund\bund}\big)\Bigg]_K.
\end{split}
\eea
Notice that $\EE_K[\X, \w, \M]$ coincides in fact with the quantity $\EE[\X, \w, \M]$ in \eqref{definition-EE-gen-SSvalued}. 
Thus, we have according to \eqref{eq:separation-EE-I-J-K}
\bea\label{eq:separation-EE-I-J-K:chap9}
\begin{split}
 |q|^2\EE_K[\bold{X}, \bold{w}, \bold{M}]   &=P+I  +J +K
   \end{split}
   \eea
where
\beaa
P&:=& \UU^{\a\b\aund\bund} \, \Db_\a \psia \c \Db_\b \psib,\\
I&:=& \AA^{\aund\bund}  \nab_r\psia\c \nab_r\psib, \\
J&:=& \VV^{\aund\bund} \psia\c\psib,\\
K&:=& \frac 1 4 |q|^2 \D^\mu\Big(M^{\aund\bund}_\mu\psi_\aund \c\psi_\bund \Big),
\eeaa
and where the coefficients $\UU^{\a\b\aund\bund}$, $\AA^{\aund\bund}$ and $\VV^{\aund\bund}$ are provided by Proposition \ref{proposition:Morawetz3}.

We now derive the analog of Lemma \ref{lemma:IntegrationbypartsP}.
   \begin{lemma}
   \lab{lemma:IntegrationbypartsP:chap9}
Let  $P$ the principal term defined as above. 
We then have the identity
\beaa
 P&=&\frac 1 2  h  L^{\a\b}  \Db_\a \Psi\c    \Db_\b \Psi -\frac 1 2  h\Psi\c ( \RRtp^{\cund}  \LL^{\bund}[\SS_{\cund}, \SS_{\bund}]\psi) +|q|^2\Db_\a \BB^\a +r^2\textrm{Good}_\SS,
\eeaa
where $\Psi$ is defined as 
\bea
\Psi:= \RRtp^{\aund} \psia ,
\eea
and the boundary term $\BB$ is given by
\bea\label{eq:definition-BB:chap9}
\BB^\a&:=& |q|^{-2} \frac 1 2  h \Psi  \RRtp^{\cund}  \LL^{\bund}\c\left(\Sc^{\a\b}   \, \Db_\b \psib-  \Sb^{\a\b}   \Db_\b \psic \right).
\eea
\end{lemma}

\begin{proof}
As in the proof of Lemma \ref{lemma:IntegrationbypartsP}, we have
\beaa
|q|^{-2} P &=& - \UU^{\aund\cund} \LL^{\bund} \psia\c \Db_\a(|q|^{-2} \Sc^{\a\b}   \, \Db_\b \psib)+\Db_\a( |q|^{-2}\UU^{\aund\cund} \LL^{\bund}\Sc^{\a\b}   \psia \c \Db_\b \psib)\\
&& - \Db_\a( \UU^{\aund\cund} \LL^{\bund}) |q|^{-2}\Sc^{\a\b}   \psia \c \Db_\b \psib.
\eeaa
Since $\LL^{\aund}$  and  $ \UU^{\aund\bund}$ only depend on $r$, and in view of the form of the 2-tensors $ \Sc^{\a\b}$, we have 
\beaa
\Db_\a\big( \UU^{\aund\cund} \LL^{\bund} \big)|q|^{-2} \Sc^{\a\b}  \psia \c \Db_\b \psib &=&   \pr_r\big( \UU^{\aund\cund}  \LL^{\bund}  \big) e_\a(r) |q|^{-2}\Sc^{\a\b} \psia \c \Db_\b \psib\\
&=& O(r^{-2})\pr_r\big( \UU^{\aund\cund}  \LL^{\bund}  \big)\Big(\T(r), \Z(r), r\nab(r)\Big)\psi\dk^{\leq 1}\psi\\
&=& \pr_r\big( \UU^{\aund\cund}  \LL^{\bund}  \big)\Ga_g\psi\dk^{\leq 1}\psi\\
&=& \pr_r\big( \UU^{\aund\cund}  \LL^{\bund}  \big)r\textrm{Good}_\SS.
\eeaa
Now, note that the choices of Chapter \ref{chapter-proof-mor-2} imply $\UU^{\aund\cund}  \LL^{\bund}=O(r^{-1})$ and $\pr_r( \UU^{\aund\cund}  \LL^{\bund})=O(r^{-2})$ and hence
\beaa
\Db_\a\big( \UU^{\aund\cund} \LL^{\bund} \big)|q|^{-2} \Sc^{\a\b}  \psia \c \Db_\b \psib &=& r^{-1}\textrm{Good}_\SS.
\eeaa
Also, we have in view of Lemma \ref{LEMMA:SYMM-OPERATORS}
\beaa
\widetilde{\SS}_1 = \SS_1  + \Ga_b\c \dk, \qquad \widetilde{\SS}_\aund &=& \SS_\aund +r\Ga_b \c \dk, \quad \aund=2,3,4,
 \eeaa
where 
\beaa
\widetilde{\SS}_\aund\psi&=&|q|^2 \Ddot_\a(|q|^{-2}S_\aund^{\a\b} \Ddot_\b \psi), \qquad \textrm{ for }\aund=1,2,3, 4.
\eeaa
Therefore we write
\beaa
 \Db_\a(|q|^{-2} \Sc^{\a\b}   \, \Db_\b \psib)&=&  |q|^{-2}\widetilde{\SS}_\cund(\SS_\bund\psi)=|q|^{-2}\SS_\cund \SS_\bund\psi +r^{-1}\Ga_b\dk\SS_\bund\psi\\
&=& |q|^{-2}\SS_\bund \SS_\cund\psi+|q|^{-2}[\SS_{\cund}, \SS_{\bund}]\psi  +r^{-1}\Ga_b\dk^{\leq 3}\psi\\
 &=& \Db_\a(|q|^{-2} \Sb^{\a\b}   \, \Db_\b \psic)+|q|^{-2}[\SS_{\cund}, \SS_{\bund}]\psi +r^{-1}\Ga_b\dk^{\leq 3}\psi.
\eeaa

 Thus, repeating the integration by parts procedure and noting, as above, that $\UU^{\aund\cund}  \LL^{\bund}=O(r^{-1})$ and that the last term is of the type $r^{-1}\textrm{Good}_\SS$, we obtain 
\beaa
  &&\UU^{\aund\cund} \LL^{\bund} \psia\c    \Db_\a \big(|q|^{-2}\Sc^{\a\b}  \Db_\b  \psib\big)\\
  &=&  \UU^{\aund\cund} \LL^{\bund} \psia\c   \Db_\a(|q|^{-2} \Sb^{\a\b}   \, \Db_\b \psic) +|q|^{-2}\UU^{\aund\cund} \LL^{\bund} \psia\c [\SS_{\cund}, \SS_{\bund}]\psi   + r^{-1}\UU^{\aund\cund} \LL^{\bund} \psia\c\Ga_b\dk^{\leq 3}\psi\\
&=& -|q|^{-2}\UU^{\aund\cund} \LL^{\bund} \Sb^{\a\b}  \Db_\a \psia\c      \Db_\b \psic+\Db_\a \Big(|q|^{-2} \UU^{\aund\cund} \LL^{\bund}  \Sb^{\a\b}   \psia\c   \Db_\b \psic\Big)\\
&&+|q|^{-2}\UU^{\aund\cund} \LL^{\bund} \psia\c [\SS_{\cund}, \SS_{\bund}]\psi+r^{-1}\textrm{Good}_\SS+ r^{-2}\Ga_b\dk^{\leq 2}\psi\c\dk^{\leq 3}\psi\\
&=& -|q|^{-2}\UU^{\aund\cund} \LL^{\bund} \Sb^{\a\b}  \Db_\a \psia\c      \Db_\b \psic+\Db_\a \Big(|q|^{-2} \UU^{\aund\cund} \LL^{\bund}  \Sb^{\a\b}   \psia\c   \Db_\b \psic\Big)\\
&&+|q|^{-2}\UU^{\aund\cund} \LL^{\bund} \psia\c [\SS_{\cund}, \SS_{\bund}]\psi+\textrm{Good}_\SS.
\eeaa
Therefore, recalling that $  \LL^{\bund}   \Sb^{\a\b}=L^{\a\b}  $,
\beaa
 P&=& \UU^{\aund\cund}L^{\a\b}  \Db_\a \psia\c     \Db_\b \psic+|q|^2\Db_\a \left(|q|^{-2} \UU^{\aund\cund} \LL^{\bund}  \psia\c \left(\Sc^{\a\b}   \, \Db_\b \psib-  \Sb^{\a\b}   \Db_\b \psic \right) \right ) \\
&& -\UU^{\aund\cund} \LL^{\bund} \psia\c [\SS_{\cund}, \SS_{\bund}]\psi +r^2\textrm{Good}_\SS.
\eeaa

Finally, since we have $\RRtp^{\aund}=O(r^{-3})$ and $\pr_r\RRtp^{\aund}=O(r^{-4})$, and since $\RRtp^{\aund}$ only depend on $r$, we obtain
\beaa
&& L^{\a\b}\Db_\a (\RRtp^{\aund} \psia )\c    \Db_\b(\RRtp^{\cund} \psic)\\ 
&=& L^{\a\b}\RRtp^{\aund}  \RRtp^{\cund} \Db_\a \psia\c     \Db_\b \psic+\big(\T(r), \Z(r), r\nab(r)\big)O(r^{-7})\dk^{\leq 2}\psi\dk^{\leq 3}\psi\\
&=& L^{\a\b}\RRtp^{\aund}  \RRtp^{\cund} \Db_\a \psia\c     \Db_\b \psic+r^{-5}\Ga_g\dk^{\leq 2}\psi\dk^{\leq 3}\psi\\
&=& L^{\a\b}\RRtp^{\aund}  \RRtp^{\cund} \Db_\a \psia\c     \Db_\b \psic+r^{-4}\textrm{Good}_\SS
\eeaa
which together with  \eqref{definition:UUab} implies
\beaa
 P&=&\frac 1 2  h\RRtp^{\aund}  \RRtp^{\cund} L^{\a\b}  \Db_\a \psia\c     \Db_\b \psic -\frac 1 2  h\RRtp^{\aund}  \RRtp^{\cund}  \LL^{\bund} \psia\c [\SS_{\cund}, \SS_{\bund}]\psi\\
&& +|q|^2\Db_\a \left(|q|^{-2} \frac 1 2  h\RRtp^{\aund}  \RRtp^{\cund}  \LL^{\bund}  \psia\c \left(\Sc^{\a\b}   \, \Db_\b \psib-  \Sb^{\a\b}   \Db_\b \psic \right) \right ) +r^2\textrm{Good}_\SS\\
&=&\frac 1 2  h  L^{\a\b}  \Db_\a (\RRtp^{\aund} \psia )\c    \Db_\b(\RRtp^{\cund} \psic) -\frac 1 2  h \RRtp^{\cund}  \LL^{\bund}(\RRtp^{\aund}  \psia )\c[\SS_{\cund}, \SS_{\bund}]\psi\\
&& +|q|^2\Db_\a \left(|q|^{-2} \frac 1 2  h  \RRtp^{\cund}  \LL^{\bund} (\RRtp^{\aund} \psia )\c\left(\Sc^{\a\b}   \, \Db_\b \psib-  \Sb^{\a\b}   \Db_\b \psic \right) \right ) +r^2\textrm{Good}_\SS.
\eeaa
By denoting $\Psi=\RRtp^{\aund} \psia$ we obtain the stated expression.
\end{proof}

In view of \eqref{eq:separation-EE-I-J-K:chap9} and Lemma \ref{lemma:IntegrationbypartsP:chap9}, we obtain the following analog of \eqref{Morawetz-effective-principal-term}
 \bea\lab{Morawetz-effective-principal-term:chap9}
\begin{split}
 |q|^2\EE_K[\bold{X}, \bold{w}, \bold{M}] -|q|^2 \D^\mu \BB_\mu   &=\widetilde{P}+P_{lot}+I  +J +K +r^2\textrm{Good}_\SS
   \end{split}
   \eea
 where the quadratic forms $I$, $J$ and $K$ are defined below \eqref{eq:separation-EE-I-J-K:chap9} and 
 \begin{itemize}
 \item $\BB$ is the boundary term defined in \eqref{eq:definition-BB:chap9},
 
 \item The principal term $\widetilde{P}$ is positive definite and given by
 \bea\label{eq:definition-widetilde:chap9}\label{eq:widetilde-P:chap9}
   \widetilde{P}&:=&\frac 1 2  h  L^{\a\b}  \Db_\a \Psi\c    \Db_\b \Psi , \qquad \Psi=\RRtp^{\aund} \psia,
 \eea
 with $L^{\a\b} $  as in  \eqref{eq:operato-LL2},
 
 \item The lower order term $P_{lot}$ is given by 
 \bea\label{eq:Plot:chap9}
 P_{lot}&=&-\frac 1 2  h\Psi\c ( \RRtp^{\cund}  \LL^{\bund}[\SS_{\cund}, \SS_{\bund}]\psi).
 \eea
  \end{itemize}

As in \eqref{widetilde-EE}, we denote the effective generalized current 
\bea\lab{widetilde-EE:chap9}
 |q|^2 \widetilde{\EE}_K[\bold{X}, \bold{w}, \bold{M}]:= |q|^2 \EE_K[\bold{X}, \bold{w}, \bold{M}]  -|q|^2 \D^\mu \BB_\mu. 
\eea

Next, we evaluate each term on the RHS of \eqref{Morawetz-effective-principal-term:chap9}:
\begin{enumerate}
\item For $\widetilde{P}$, we immediately have the analog of \eqref{eq:expression-widetilde-Psiz} of, i.e.    
\bea\label{eq:expression-widetilde-Psiz:chap9}
\widetilde{P} =\frac 1 2  h \Big(\de_0 \big|   \nab_T  \Psi_z \big|^2 + a^2\big|  \nab_Z \Psi_z \big|^2 + O^{\a\b} \Db_\a \Psi_z \Db_\b   \Psi_z\Big).
\eea

\item Concerning $ P_{lot}$, we make use of Lemma \ref{lemma:commutationpropertiesofthesymmetryoperatorsafsoiudf:chap9} according to which we have
\beaa
\,[\SS_1, \SS_2], \,[\SS_1, \SS_3],\,[\SS_2, \SS_3]=\dk^2(\Ga_b\c\psi)
\eeaa
and 
\beaa
\,[\SS_1, \OO] &=& O(ar^{-3})\dk^{\leq 2}\psi+ \dk^{\leq 3}(\Ga_b\psi),\\
\,[\SS_2, \OO] &=& O(a)\dkb^{\leq 2}\psi+ r\dk^{\leq 3}(\Ga_b\psi),\\
\,[\SS_3, \OO] &=& O(a^2)\dkb^{\leq 2}\psi+ r\dk^{\leq 3}(\Ga_b\psi).
\eeaa
In view of the definition of $P_{lot}$ in \eqref{eq:Plot:chap9}, we infer, since $h=O(r^5)$, $\RRtp^{\aund}=O(r^{-3})$, and $\Psi=O(r^{-3})\dk^{\leq 2}\psi$, 
\beaa
P_{lot} &=& -\frac 1 2  h\Psi\c ( \RRtp^{\cund}  \LL^{\bund}[\SS_{\cund}, \SS_{\bund}]\psi) \\
&=& [P_{lot}]_K+h\Psi\c \RRtp^{\cund}  \LL^{\bund} r\dk^{\leq 3}(\Ga_b\psi)= [P_{lot}]_K+\dk^{\leq 2}\psi\c\dk^{\leq 3}(\Ga_b\psi)\\
&=& [P_{lot}]_K+r^2\textrm{Good}_\SS
\eeaa
where $[P_{lot}]_K$ denotes the corresponding expression in Kerr, i.e the one of \eqref{eq:Plot}. In particular, $[P_{lot}]_K$ satisfies \eqref{eq:expression-Plot} which implies the following analog of \eqref{eq:Plot} for $P_{lot}$
\bea\label{eq:expression-Plot:chap9}
P_{lot}    &=&  O(a r^{-1})(\dk^{\leq 2}\psi)^2+r^2\textrm{Good}_\SS.
\eea

\item $I$ is estimated as in section \ref{sec:thequadraticformsIJKarenowcomputed}, where Lemma \ref{general-computations-for-BB} is replaced with Lemma \ref{general-computations-for-BB:chap9} which yields the following analog of \eqref{eq:computationofthequadraticformIforchap8:fspoighs}
\bea\lab{eq:computationofthequadraticformIforchap8:fspoighs:chap9}
 \bsplit
 I =& {\AA}(1+O(r^{-1} \de_0))\Big(\de_0^2 |\nab_r\SS_1\psi|^2+ |\nab_r\OO\psi|^2 +2\de_0|q|^2| \nab \nab_T\nab_r\psi|^2\Big)\\
& -O(a)(|\AA|+|\AAt|)  (1+O(r^{-1} \de_0))\Big((\nab_r\dk^{\leq 2}\psi)^2+r^{-2}(\dk^{\leq 2}\psi)^2\Big)\\
& +\text{Bdr}[\psi]_I +r^2\textrm{Good}_\SS,
\end{split}
\eea
where the boundary term is given by
\beaa
\text{Bdr}[\psi]_I&=&\Ddot_\mu\Big(\de_0\AA M(\nab_r\psi)\T^\mu\Big) +|q|^2 \Ddot_\b\Big(2\de_0|q|^{-2}O^{\a\b}\Ddot_\a \nab_r \psi  \c {\AA} \SS_1\nab_r\psi \Big).
\eeaa

\item $J$ is estimated as in Lemma \ref{lemma:computationofthequadraticformJandKforchap8:fspoighs} where Lemma \ref{general-computations-for-BB} is replaced with Lemma \ref{general-computations-for-BB:chap9} which yields the following analog of Lemma \ref{lemma:computationofthequadraticformJandKforchap8:fspoighs}
\bea
 \begin{split}
J =& \big( \VV +O( \de_0   r^{-3}) \big)\Big(\de_0^2 |\SS_1\psi|^2+ |\OO\psi|^2 +2\de_0|q|^2| \nab \nab_T\psi|^2\Big)\\
&  - O(ar^{-1})(\dk^{\leq 2}\psi)^2+\text{Bdr}[\psi]_J +r\textrm{Good}_\SS,
\end{split}
\eea
 with boundary term
\beaa
\text{Bdr}_J[\psi] &=& \Ddot_\mu\Big(\de_0(\VV +O(   r^{-3}))  r^2   M(\psi)\T^\mu\Big) \\
&&+|q|^2 \Ddot_\b\Big(2\de_0 
|q|^{-2}O^{\a\b}\Ddot_\a \psi \c (\VV +O( a^2   r^{-3})) \SS_1\psi \Big).
\eeaa

\item For $K$, we assume as in Lemma \ref{lemma:computationofthequadraticformJandKforchap8:fspoighs} that \beaa
M^\aund :=v^{\aund}\frac{r^2+a^2}{\De}\Rhat,\qquad v^1=\de_0 v, \qquad v^2=v^3=0, \qquad v^4=v,
\eeaa
for some given function $v=v(r)$. 
Furthermore, under the assumption $v=O(m^{1/2}\Delta r^{-9/2})$, we have
 \beaa
 \D_\mu((M^\aund)^\mu) &=& \D_\mu\left(\frac{v^\aund}{\De}(r^2+a^2)\Rhat^\mu\right)\\
 &=& \left[\D_\mu\left(\frac{v^\aund}{\De}(r^2+a^2)\Rhat^\mu\right)\right]_K+O(r^{-\frac{7}{2}})\widecheck{\Rhat(r)}+O(r^{-\frac{5}{2}})\widecheck{\D_\mu\Rhat^\mu}\\
  &=& \left[\D_\mu\left(\frac{v^\aund}{\De}(r^2+a^2)\Rhat^\mu\right)\right]_K+r^{-\frac{5}{2}}\Ga_b
 \eeaa
 which together with the proof of Lemma \ref{lemma:computationofthequadraticformJandKforchap8:fspoighs} yields
 \beaa
\frac{4}{|q|^2} K &=&(\D^\mu M^{\aund}_\mu) \psi_\aund \c (\LL^\bund   \psi_\bund) + (M^{\aund}_\mu \D^\mu\psi_\aund )\c (\LL^\bund   \psi_\bund) +(M^{\aund}_\mu\psi_\aund )\c \D^\mu(\LL^\bund   \psi_\bund)\\
&=&  (v'^{\aund} \psi_\aund) \c (\LL^\bund   \psi_\bund) + v^{\aund} \nab_r \psi_\aund \c (\LL^\bund   \psi_\bund) +(v^{\aund}\psi_\aund )\c  \LL^\bund   \nab_r\psi_\bund + r^{-\frac{5}{2}}\Ga_b \psi_\aund \c (\LL^\bund   \psi_\bund)\\
&=&  (v'^{\aund} \psi_\aund) \c (\LL^\bund   \psi_\bund) + v^{\aund} \nab_r \psi_\aund \c (\LL^\bund   \psi_\bund) +(v^{\aund}\psi_\aund )\c  \LL^\bund   \nab_r\psi_\bund + r^{-\frac{5}{2}}\Ga_b(\dk^{\leq 2}\psi)^2\\
&=&  (v'^{\aund} \psi_\aund) \c (\LL^\bund   \psi_\bund) + v^{\aund} \nab_r \psi_\aund \c (\LL^\bund   \psi_\bund) +(v^{\aund}\psi_\aund )\c  \LL^\bund   \nab_r\psi_\bund + r^{-\frac{1}{2}}\textrm{Good}_\SS
\eeaa
where we wrote $v'^{\aund}:=\pr_r v^\aund+ \frac{2r}{|q|^2} v^\aund$. Then, we follow the rest of the proof of 
 Lemma \ref{lemma:computationofthequadraticformJandKforchap8:fspoighs} where Lemma \ref{general-computations-for-BB} is replaced with Lemma \ref{general-computations-for-BB:chap9} and Lemma \ref{Lemma:integrationbypartsSS_3SS_4} is replaced by Lemma \ref{Lemma:integrationbypartsSS_3SS_4:SSversion:chap9} 
  which yields the following analog of Lemma \ref{lemma:computationofthequadraticformJandKforchap8:fspoighs}
\bea
\nn K&=&\frac{|q|^2}{2} v \Big( \de_0^2 \nab_r \SS_1\psi \c \SS_1\psi+ \nab_r \OO\psi \c \OO\psi+2\de_0   |q|^2 \nab\nab_T\nab_r \psi \c \nab \nab_T \psi\Big)\\
\nn&&+\frac{|q|^2}{4}v'\Big( \de_0^2  |\SS_1\psi|^2+ |\OO\psi|^2 +2\de_0 |q|^2| \nab \nab_T\psi|^2\Big) - vO(a r^{\frac{5}{2}})(\nab_r\dk^{\leq 2}\psi)^2\\
&& - O(a r^{\frac{3}{2}}) v(\dk^{\leq 2}\psi)^2 - O(a r^2)v' (\dk^{\leq 2}\psi)^2 +\text{Bdr}[\psi]_K+r^{\frac{3}{2}}\textrm{Good}_\SS, 
\eea
where we denoted
\beaa
v'^{\aund}:=\pr_r v^\aund+ \frac{2r}{|q|^2} v^\aund,
\eeaa 
and with boundary term 
\beaa
\text{Bdr}[\psi]_K &=& \Ddot_\mu\Big(\big(vr^4M(\nab_r\psi)\big)\T^\mu\Big) +\frac{|q|^4}{4} \Ddot_\b(2\de_0 v|q|^{-2}O^{\a\b}\Ddot_\a\psi \c \SS_1\nab_r \psi )\\
&& +\frac{|q|^4}{4} \Ddot_\b(2\de_0 v|q|^{-2}O^{\a\b}\Ddot_\a\nab_r\psi \c \SS_1\psi ) +\Ddot_\mu\Big(\de_0v' r^4 M(\psi)\T^\mu\Big)  \\
&&+\frac{|q|^4}{4}\Ddot_\b\Big(2\de_0 v' |q|^{-2}O^{\a\b}\Ddot_\a \psi \c\SS_1\psi \Big).
\eeaa
\end{enumerate}

The above computations of $\widetilde{P}$, $P_{lot}$, $I$, $J$ and $K$ immediately yields the following analog of   Proposition \ref{proposition-with-quadratic-form}.
\begin{proposition}
\label{proposition-with-quadratic-form:chap9}
The  effective generalized current is given by
 \bea\label{generalized-current-operator-2:chap9}
 |q|^2 \widetilde{\EE}_K[\bold{X}, \bold{w}, \bold{M}] &=&\widetilde{P}+\Qr_{\SS_1, \OO, \nab \nab_T}+   \EE_{lot}+\text{Bdr}+r^2\textrm{Good}_\SS
   \eea
   with the following terms.
\begin{enumerate}
\item The principal trapping term $\widetilde{P}$ is given by, see \eqref{eq:expression-widetilde-Psiz:chap9},
\bea
\widetilde{P}&=&\frac 1 2  h \Big( \de_0 \big|   \nab_T  \Psi_z \big|^2+a^2|\nab_Z \Psi_z|^2 + O^{\a\b} \Db_\a \Psi_z \Db_\b   \Psi_z\Big),
\eea
where
\bea\label{definition-Psiz:chap9}
\begin{split}
\Psi_z&= -\frac{2\TT}{(r^2+a^2)^3}   \big(\de_0  \SS_1\psi+ (1+O(r^{-2} \de_0)) \OO\psi\big)\\
&+ \frac{4ar}{(r^2+a^2)^2}  \nab_\That \nab_Z \psi   \big(1+O(r^{-2} \de_0) \big).
\end{split}
\eea

\item The quadratic   form $\Qr_{\SS_1, \OO, \nab \nab_T}$ is given by
   \beaa
   \Qr_{\SS_1, \OO, \nab \nab_T}&:=& {\AA}(1+O(r^{-1} \de_0))\Big(\de_0^2 |\nab_r\SS_1\psi|^2+ |\nab_r\OO\psi|^2 +2\de_0|q|^2| \nab \nab_T\nab_r\psi|^2\Big)\\
   && +\frac{|q|^2}{2} v \Big( \de_0^2 \nab_r \SS_1\psi \c \SS_1\psi+ \nab_r \OO\psi \c \OO\psi+2\de_0   |q|^2 \nab\nab_T\nab_r \psi \c \nab \nab_T \psi\Big)\\
 &&  +\left( \VV+\frac{|q|^2}{4}v' +O( \de_0   r^{-3}) \right)\Big(\de_0^2 |\SS_1\psi|^2+ |\OO\psi|^2 +2\de_0|q|^2| \nab \nab_T\psi|^2\Big).
   \eeaa

   \item The terms $\EE_{lot}$ are lower order terms in $a$ satisfying 
   \bea\label{eq:bound-EE-lot:chap9}
      \EE_{lot} \geq -O(a)\big( |\nab_{\Rhat}\dk^{\leq 2}\psi|^2+r^{-1} |\dk^{\leq 2}\psi|^2\big).
   \eea
   
\item The boundary terms are given by
   \beaa
   \text{Bdr}&=& \Ddot_\mu\Big(  \big(M(\nab_{\Rhat}\psi)+M(\psi)\big)\T^\mu\Big)    +|q|^2 \D_\b  \widehat{\BB}^\b
   \eeaa
   with 
   \beaa
   \widehat{\BB}^\b&:=& 2\de_0 |q|^{-2}O^{\a\b}\Ddot_\a \nab_r \psi  \c {\AA} \SS_1\nab_r\psi +2\de_0|q|^{-2}O^{\a\b}\Ddot_\a \psi \c(\VV +O( a^2   r^{-3})) \SS_1\psi\\
   &&+\frac{1}{4}\de_0 vO^{\a\b}\Ddot_\a\psi \c\SS_1\nab_r \psi +\frac{1}{4}\de_0 vO^{\a\b}\Ddot_\a\nab_r\psi \c \SS_1\psi +\frac{1}{2}\de_0 O^{\a\b}\Ddot_\a \psi \c v' \SS_1\psi,
   \eeaa
  where $M(\psi)$ denotes the  quadratic expressions in $\psi$ and  $M(\nab_{\Rhat}\psi)$  denotes the  quadratic expressions in $\psi$ and its derivatives of Definition \ref{def:boundarytermsinintegrationbypartsareok:SSversion:chap9}.
\end{enumerate}
\end{proposition}

Next, we proceed as in section \ref{subsection-Poincare-higher} in order to control  the quadratic form $ \Qr_{\SS_1, \OO, \nab \nab_T}$ appearing in Proposition \ref{proposition-with-quadratic-form:chap9}. 
The proof is analogous provided we use: 
\begin{itemize}
\item the Poincar\'e inequality of Lemma \ref{lemma:poincareinequalityfornabonSasoidfh:chap9} which is the analog of the Poincar\'e inequality of Lemma \ref{lemma:poincareinequalityfornabonSasoidfh:chap6},

\item the integration by parts identities of Lemma \ref{Lemma:integrationbypartsSS_3SS_4:SSversion:chap9} which is the analog of the ones of Lemma \ref{Lemma:integrationbypartsSS_3SS_4}.
\end{itemize}

This leads to the following analog of \eqref{eq:positivity-first-quadratic-form}, for  any sphere $S=S(\tau, r)$,
\bea\label{eq:positivity-first-quadratic-form:chap9}
\begin{split}
\int_S |q|^2 \widetilde{\EE}_K[\bold{X}, \bold{w}, \bold{M}] \geq& \de\int_S\widetilde{P}\\
 &+c_1  \int_S\Big( m\big(\big|\nab_{\Rhat}\SS_1\psi|^2 + \big|\nab_{\Rhat}\OO\psi|^2+|q|^2\big|\nab\nab_T\nab_{\Rhat}\psi|^2\big)\\
 & + r^{-1}\big( |\SS_1\psi|^2+|\OO\psi|^2+|q|^2 |\nab\nab_T\psi|^2\big) \Big)\\
 &+ \int_S\Big(  \widetilde{\EE}_{lot} -O(ar^{-2})\big(|\nab_{\Rhat}\dk^{\leq 2}\psi|^2+|\dk^{\leq 2}\psi|^2\big) +\text{Bdr}+r^2\textrm{Good}_\SS\Big).
 \end{split}
\eea

Next, we proceed exactly as in section \ref{sec:finallowerboundforeffectivegeneralizedcurrentaosidhfauidh} and obtain the following analog of \eqref{eq:veryfinalconclusionforthelowerboundgeneralizedcurrent:fdsoiugah}
\bea\lab{eq:veryfinalconclusionforthelowerboundgeneralizedcurrent:fdsoiugah:chap9}
\bsplit
\int_S \widetilde{\EE}_K[\bold{X}, \bold{w}, \bold{M}] \ges& \int_S  \frac{m}{r^2}  | \nab_\Rhat \psi|_\SS^2 + r^{-3}|\psi|_{\SS}^2+ r^{3} \Big(\big|\nab_\That \Psi_z \big |^2+r^2\big |\nab \Psi_z\big|^2\Big)+\text{Bdr}\\
&  -O(a)\int_S\big( r^{-2}|\nab_{\Rhat}\dk^{\leq 2}\psi|^2+r^{-3} |\dk^{\leq 2}\psi|^2\big)+\int_S\textrm{Good}_\SS
\end{split}
\eea
on any sphere $S=S(\tau, r)$. 

Next, proceeding as in the beginning of section \ref{sec:SSderivativeversionofLemmacontrolrhs:aodihas}, relying on Lemma \ref{lemma:decompositionofthedivergenceinerrorandtermKersoiduhf:SSversion}, \eqref{definition-EE-gen:chap9:bisrepetita} and \eqref{widetilde-EE:chap9}, we derive the following analog of \eqref{definition-EE-gen-SSvalued:fsidufsdofgus}
\beaa
\begin{split}
\D^\mu  \PP_\mu[\X, \w, \M] =&  \widetilde{\EE}_K[\bold{X}, \bold{w}, \bold{M}]  +  \D^\mu \BB_\mu + \left(\nab_{X^{\aund\bund}}\psia+\frac 1 2   w^{\aund\bund} \psia\right)\c \big( \squared_2\psib -V\psib \big)\\
& -\frac{2a^2r\cos\th\FF^{\aund\bund}(r)}{(r^2+a^2)|q|^4}\nab_\Z\psi_\aund\c\dual\psi_\bund\\
& -\left(\big(\rhod +\etab\wedge\eta\big)\FF^{\aund\bund}\frac{r^2+a^2}{\De}+\frac{2a^3r\cos\th(\sin\th)^2\FF^{\aund\bund}(r)}{|q|^6}\right)\nab_{\That}\psi_\aund\c\dual\psi_\bund\\
&+r^{-1}\Ga_b\Big(|(X^{\aund\bund})^3|+(X^{\aund\bund})^4|\Big)\nab\psi_{\aund}\c\dual\psi_{\bund}
 \end{split}
\eeaa
where we used the fact that $\Hc\in\Ga_b$ and $\trXbc\in \Ga_g$. Since $|(X^{\aund\bund})^3|+(X^{\aund\bund})^4|\les 1$, we have
\beaa
r^{-1}\Ga_b\Big(|(X^{\aund\bund})^3|+(X^{\aund\bund})^4|\Big)\nab\psi_{\aund}\c\dual\psi_{\bund} &=& r^{-2}\Ga_b\dk^{\leq 1}\psi\c\dual\psi =\textrm{Good}_\SS
\eeaa
and we deduce
\bea\lab{definition-EE-gen-SSvalued:fsidufsdofgus:chap9}
\begin{split}
& \D^\mu  \PP_\mu[\X, \w, \M] \\
=&  \widetilde{\EE}_K[\bold{X}, \bold{w}, \bold{M}]  +  \D^\mu \BB_\mu + \left(\nab_{X^{\aund\bund}}\psia+\frac 1 2   w^{\aund\bund} \psia\right)\c \big( \squared_2\psib -V\psib \big)\\
& -\left(\big(\rhod +\etab\wedge\eta\big)\FF^{\aund\bund}\frac{r^2+a^2}{\De}+\frac{2a^3r\cos\th(\sin\th)^2\FF^{\aund\bund}(r)}{|q|^6}\right)\nab_{\That}\psi_\aund\c\dual\psi_\bund\\
& -\frac{2a^2r\cos\th\FF^{\aund\bund}(r)}{(r^2+a^2)|q|^4}\nab_\Z\psi_\aund\c\dual\psi_\bund+\textrm{Good}_\SS.
 \end{split}
\eea

We now derive the following analog of Lemma \ref{lemma-control-rhs-higher}. 
\begin{lemma}\lab{lemma-control-rhs-higher:chap9} 
We have, for sufficiently small positive constants $\de_2$, $\de_3$:
\bea
\begin{split}
&\left(\nab_{X^{\aund\bund}}\psia +\frac 1 2   w^{\aund\bund} \psia\right)\c \big( \squared_2\psib -V\psib \big)\\
\geq& -\de_2r^{-2} h|\nab_{\That} \Psi_z|^2 -\de_2 a^2r^{-6}h|\nab_\Z \Psi_z|^2+O(ar^{-3}) |\psi|_\SS^2\\
&+O(1)\Big( | \nab_\Rhat \psi |_{\SS} + r^{-1}|\psi|_{\SS} \Big)\sum_{\aund=1}^4| N_{\aund}|\\
& + \Ddot_\mu\left( \frac{ 2 a\cos\th}{|q|^2}\frac{r^2+a^2}{\De}\Rhat^\mu zh f^{\aund\bund} \psia\c\nab_\T\dual \psib\right)\\
&   - \Ddot_\mu\left( \frac{ 2 a\cos\th}{|q|^2}\left( zhf^{\aund\bund} \psia\c  \frac{r^2+a^2}{\De}\nab_\Rhat\dual \psib +\frac 1 2 (\pr_r z) h \Psi_z\c\dual\LL^{\aund}\psia\right)\T^\mu\right)+\textrm{Good}_\SS,
\end{split}
\eea
and
\bea
\bsplit
& \left(\big(\rhod +\etab\wedge\eta\big)\FF^{\aund\bund}\frac{r^2+a^2}{\De}+\frac{2a^3r\cos\th(\sin\th)^2\FF^{\aund\bund}(r)}{|q|^6}\right)\nab_{\That}\psi_\aund\c\dual\psi_\bund\\
& +\frac{2a^2r\cos\th\FF^{\aund\bund}(r)}{(r^2+a^2)|q|^4}\nab_\Z\psi_\aund\c\dual\psi_\bund\\
\leq& \de_3 r^3\Big(|\nab_\That \Psi_z|^2+r^2|\nab\Psi_z|^2\Big)  +O(a r^{-3}) |\psi|_\SS^2 -\frac{1}{2}\Ddot_\mu\left(z h\frac{2a^2r\cos\th}{(r^2+a^2)|q|^4}\Psi_z\c\dual(\LL^{\underline{b}}\psi_\bund)\Z^\mu\right)\\
 & -\frac{1}{2}\Ddot_\mu\left(z h  \left(\big(\rhod +\etab\wedge\eta\big)\frac{r^2+a^2}{\De}+\frac{2a^3r\cos\th(\sin\th)^2}{|q|^6}\right)\Psi_z\c\dual(\LL^{\underline{b}}\psi_\bund)\That^\mu\right)+\textrm{Good}_\SS.
\end{split}
\eea
\end{lemma}

\begin{proof}
We proceed as in the proof of Lemma \ref{lemma-control-rhs-higher}. It then suffices to check that the extra error terms are all of the type $\textrm{Good}_\SS$. For the first estimate, the extra error terms are of the type
\beaa
&& O(r^{-2})\Big(r^{-1}\widecheck{\Rhat(r)}, r^{-1}\T(r), (\D_\mu\T^\mu), \widecheck{\D_\mu\Rhat^\mu}, \T(\cos\th), \Rhat(\cos\th)\Big)\dk^{\leq 2}\psi\c\dk^{\leq 3}\psi\\
&=& r^{-2}\Ga_b\dk^{\leq 2}\psi\c\dk^{\leq 3}\psi=\textrm{Good}_\SS
\eeaa
as desired. 

For the second estimate,  the extra error terms are of the type
\beaa
O(1)\Big(\widecheck{\rhod},O(r^{-2})\Hc, O(r^{-2})\Hbc\Big)\dk^{\leq 2}\psi\c\dk^{\leq 3}\psi = r^{-1}\Ga_g\dk^{\leq 2}\psi\c\dk^{\leq 3}\psi=\textrm{Good}_\SS
\eeaa
and 
\beaa
&& O(r^{-4})\Big(r^{-1}\T(r),  r^{-1}\Z(r), (\D_\mu\T^\mu), (\D_\mu\Z^\mu), \T(\cos\th), \Z(\cos\th)\Big)\dk^{\leq 2}\psi\c\dk^{\leq 3}\psi\\
&=& r^{-2}\Ga_b\dk^{\leq 2}\psi\c\dk^{\leq 3}\psi=\textrm{Good}_\SS
\eeaa
as desired. This concludes the proof of Lemma \ref{lemma-control-rhs-higher:chap9}.
\end{proof}

We are now ready to prove Proposition \ref{prop:morawetz-higher-order:perturbation}. Recall \eqref{definition-EE-gen-SSvalued:fsidufsdofgus:chap9}, i.e. 
\beaa
\begin{split}
& \D^\mu  \PP_\mu[\X, \w, \M] \\
=&  \widetilde{\EE}_K[\bold{X}, \bold{w}, \bold{M}]  +  \D^\mu \BB_\mu + \left(\nab_{X^{\aund\bund}}\psia+\frac 1 2   w^{\aund\bund} \psia\right)\c \big( \squared_2\psib -V\psib \big)\\
& -\left(\big(\rhod +\etab\wedge\eta\big)\FF^{\aund\bund}\frac{r^2+a^2}{\De}+\frac{2a^3r\cos\th(\sin\th)^2\FF^{\aund\bund}(r)}{|q|^6}\right)\nab_{\That}\psi_\aund\c\dual\psi_\bund\\
& -\frac{2a^2r\cos\th\FF^{\aund\bund}(r)}{(r^2+a^2)|q|^4}\nab_\Z\psi_\aund\c\dual\psi_\bund+\textrm{Good}_\SS.
 \end{split}
\eeaa
We apply  the divergence theorem to the above on $\MM(\tau_1, \tau_2)$, which yields
\beaa
&&\int_{\MM(\tau_1, \tau_2)}\Bigg[\widetilde{\EE}_K[\bold{X}, \bold{w}, \bold{M}]  + \left(\nab_{X^{\aund\bund}}\psia+\frac 1 2   w^{\aund\bund} \psia\right)\c \big( \squared_2\psib -V\psib\big)\\
&& -\left(\big(\rhod +\etab\wedge\eta\big)\FF^{\aund\bund}\frac{r^2+a^2}{\De}+\frac{2a^3r\cos\th(\sin\th)^2\FF^{\aund\bund}(r)}{|q|^6}\right)\nab_{\That}\psi_\aund\c\dual\psi_\bund\\
&& -\frac{2a^2r\cos\th\FF^{\aund\bund}(r)}{(r^2+a^2)|q|^4}\nab_\Z\psi_\aund\c\dual\psi_\bund\Bigg]\\
&\leq& \int_{\pr\MM(\tau_1, \tau_2)}\big(|\PP_\mu[\X, \w, \M] N^\mu|+|\BB_\mu N^\mu|+\textrm{Good}_\SS\big).
\eeaa
Then, we proceed as in section \ref{section-SS-valued.unconditionalMorawetz} relying on  the lower bound \eqref{eq:veryfinalconclusionforthelowerboundgeneralizedcurrent:fdsoiugah:chap9} for $\widetilde{\EE}_K[\bold{X}, \bold{w}, \bold{M}]$ and Lemma \ref{lemma-control-rhs-higher:chap9} which leads to 
 \beaa
 \bsplit
  \Mor_{\SSz, deg}[\psi](\tau_1, \tau_2) \les& \int_{\pr\MM(\tau_1, \tau_2) }|M_\SS(\psi)| +\frac{|a|}{m}B^2_\de[\psi](\tau_1, \tau_2)\\
  &+\sum_{\aund=1}^4\int_{\MM(\tau_1, \tau_2)}\big(|\nab_{\Rhat} \psia|+r^{-1}|\psia|\big) |N_{\aund}|+\int_{\MM(\tau_1, \tau_2)}|\textrm{Good}_\SS|.
  \end{split}
 \eeaa
 where $M_\SS(\psi)$  denotes  an expression in  $\psi$ for which we  have a bound of the form
\beaa
&&\int_{\pr\MM(\tau_1, \tau_2) }|M_\SS(\psi)| \\
&\les &  \sum_{\aund=1}^4\left(\sup_{[\tau_1, \tau_2]}E_{deg}[\psi_\aund](\tau) +\deh F_{\AA}[\psi_{\aund}](\tau_1, \tau_2)+F_{\Si_*}[\psi_{\aund}](\tau_1, \tau_2)\right)\\
&+&\left(\sup_{[\tau_1, \tau_2]}E_{deg}[(\nab_T, \dkb)^{\leq 1}\psi](\tau)+\deh F_{\AA}[(\nab_T, \dkb)^{\leq 1}\psi](\tau_1, \tau_2)+F_{\Si_*}[(\nab_T, \dkb)^{\leq 1}\psi](\tau_1, \tau_2)\right)^{\frac{1}{2}}\\
&\times&\left(\sup_{[\tau_1, \tau_2]}E_{deg}[(\nab_T, \dkb)^{\leq 2}\psi](\tau)+\deh F_{\AA}[(\nab_T, \dkb)^{\leq 2}\psi](\tau_1, \tau_2)+F_{\Si_*}[(\nab_T, \dkb)^{\leq 2}\psi](\tau_1, \tau_2)\right)^{\frac{1}{2}}.
\eeaa 
 Together with the control of $\textrm{Good}_\SS$ provided by Remark \ref{rmk::accetableerrortermsenergymorawetzundercontrolafsodsudfh:SSversion}, this concludes the proof of  Proposition \ref{prop:morawetz-higher-order:perturbation}.


\subsubsection{Proof of Lemma \ref{lemma:lowerboundPhizoutsideMtrap:perturbation}}


Proceeding exactly as for the proof of Lemma \ref{LEMMA:LOWERBOUNDPHIZOUTSIDEMTRAP} in section \ref{section:lowerboundPhizoutsideMtrap}, we obtain on $\Mntrap$ the following analog of \eqref{eq:lowerboundPhizoutsideMtrap:1}
\bea\lab{eq:lowerboundPhizoutsideMtrap:1:chap9}
\bsplit
r^3\Big(|\nab_\T\Psi_z|^2+r^2|\nab\Psi_z|^2\Big) \geq& c_0r^{-3}\Big(\left|\nab_\T\big(\de_0  \SS_1\psi+  \OO\psi\big)\right|^2+|q|^2\left|\nab\big(\de_0  \SS_1\psi+  \OO\psi\big)\right|^2\Big)\\
& -O(\de_0^2r^{-3})\Big(|\nab_\T\OO\psi|^2+r^2|\nab\OO \psi|^2\Big)\\
& -O(ar^{-3})|(\nab_\T, \dkb)\dk^{\leq 2}\psi|^2 - O(r^{-3})|\psi|_{\SS}^2.
\end{split}
\eea
as well a the following computation 
\beaa
&&\left|\nab_\T\big(\de_0  \SS_1\psi+  \OO\psi\big)\right|^2+|q|^2\left|\nab\big(\de_0  \SS_1\psi+  \OO\psi\big)\right|^2\\ 
&=& \de_0^2\Big(|\nab_\T\SS_1\psi|^2+|q|^2|\nab\SS_1\psi|^2\Big)+2\de_0\Big(\nab_\T\SS_1\psi\c\nab_\T\OO\psi+|q|^2\nab\SS_1\psi\c\nab\OO\psi\Big)\\
&& + |\nab_\T\OO\psi|^2+|q|^2|\nab\OO\psi|^2.
\eeaa

Next, recalling the definition of $\OO$ in \eqref{definition-SS}, we have
\beaa
[\nab, \OO] &=& \left[\nab,  |q|^2 \left(\lap_2 \psi + \frac{2a^2\cos\th}{|q|^2}\dual\Re(\Jk)^b \nab_b \psi   \right)\right]\\
&=& |q|^2[\nab, \De_2]\psi+O(ar^{-3})\dkb^{\leq 2}\psi+\Big(r^{-1}\nab(r), r^{-2}\widecheck{\nab(\cos\th)}, r^{-3}\widecheck{\nab\Re(\Jk)}\Big)\dkb^{\leq 2}\psi\\
&=& |q|^2[\nab, \De_2]\psi+O(ar^{-3})\dkb^{\leq 2}\psi+\Ga_g\dkb^{\leq 2}\psi.
\eeaa
Now, in view of Proposition \ref{Gauss-equation-2-tensors}, we have
\beaa
[ \nab_a, \nab_b] \psi &=&\Big( \frac 1 2 (\atrch\nab_3+\atrchb \nab_4) \psi +2\, \Kh \dual \psi\Big)  \in_{ab}\\
&=& \frac{2}{r^2}\dual\psi\in_{ab} +O(ar^{-2})(\nab_\T, \dkb)^{\leq 1}\psi +\Ga_g\dk^{\leq 1}\psi
 \eeaa
and hence
\beaa
[\nab_c, \De_2]\psi &=& \g^{ab}[\nab_c, \nab_a]\nab_b+\g^{ab}\nab_a[\nab_c,\nab_b]\\
&=& -\frac{4}{r^2}\nab_c\psi+O(ar^{-3})(\nab_\T, \dkb)^{\leq 2}\psi+r^{-1}\dk^{\leq 1}(\Ga_g\dk^{\leq 1}\psi)
\eeaa
 which implies
 \beaa
[\nab, \OO] &=& -4\nab\psi+O(ar^{-1})(\nab_\T, \dkb)^{\leq 2}\psi+\dk^{\leq 2}(\Ga_g\c\psi).
\eeaa
 Together with the commutators in Lemma \ref{lemma:commutationpropertiesofthesymmetryoperatorsafsoiudf:chap9} and Corollary \ref{cor:basicpropertiesLiebTfasdiuhakdisug:chap9}, we infer
\beaa
&&\nab_\T\SS_1\psi\c\nab_\T\OO\psi+|q|^2\nab\SS_1\psi\c\nab\OO\psi\\
 &=& \SS_1\nab_\T\psi\c\OO\nab_\T\psi+|q|^2\SS_1\nab\psi\c\OO\nab\psi\\
 &&+ \SS_1\nab_\T\psi\c[\nab_\T,\OO]\psi+|q|^2\nab\SS_1\psi\c[\nab, \OO]\psi+|q|^2[\nab,\SS_1]\psi\c\OO\nab\psi\\
 &=& \SS_1\nab_\T\psi\c\OO\nab_\T\psi+|q|^2\SS_1\nab\psi\c\OO\nab\psi -4|q|^2\nab\SS_1\psi\c\nab\psi +O(a)\big((\nab_\T, \dkb)\dk^{\leq 2}\psi\big)^2\\
 &&+r\dk^{\leq 3}\psi\c\dk^{\leq 2}(\Ga_g\c\psi)\\
&=& \SS_1\nab_\T\psi\c\OO\nab_\T\psi+|q|^2\SS_1\nab\psi\c\OO\nab\psi +|q|^2|\nab\nab_\T\psi|^2
+O(ar^{-3})\big((\nab_\T, \dkb)\dk^{\leq 2}\psi\big)^2\\
&& -\D_\mu(|q|^2\nab\nab_\T\psi\c\nab\psi \T^\mu) +r\dk^{\leq 3}\psi\c\dk^{\leq 2}(\Ga_g\c\psi).
\eeaa
Next, we rely on Lemma \ref{Lemma:integrationbypartsSS_3SS_4:SSversion:chap9} and obtain 
\beaa
&&\nab_\T\SS_1\psi\c\nab_\T\OO\psi+|q|^2\nab\SS_1\psi\c\nab\OO\psi\\
 &=& |q|^2| \nab \nab_\T^2\psi|^2+|q|^2 \Ddot_\b(|q|^{-2}O^{\a\b}\Ddot_\a \nab_\T\psi \c \SS_1 \nab_\T\psi ) +\Ddot_\mu(M(\nab_\T\psi)\T^\mu)\\
&&+|q|^4| \nab^2\nab_\T\psi|^2+|q|^4 \Ddot_\b(|q|^{-2}O^{\a\b}\Ddot_\a \nab\psi \c \SS_1\nab\psi ) +\Ddot_\mu (r^2M(\nab\psi)\T^\mu)\\
&&+|q|^2|\nab\nab_\T\psi|^2+O(a)\big((\nab_\T, \dkb)^{\leq 1}\dk^{\leq 2}\psi\big)^2\\ 
&&  -\D_\mu(|q|^2\nab\nab_\T\psi\c\nab\psi \T^\mu) +r\dk^{\leq 3}\psi\c\dk^{\leq 2}(\Ga_g\c\psi)+r^2\textrm{Good}_\SS.
\eeaa
Rearranging, and using the definition of $\textrm{Good}_\SS$, this yields
\beaa
&&\nab_\T\SS_1\psi\c\nab_\T\OO\psi+|q|^2\nab\SS_1\psi\c\nab\OO\psi\\ 
&=& |q|^2| \nab \nab_\T^2\psi|^2 +|q|^4| \nab^2\nab_\T\psi|^2 +|q|^2|\nab\nab_\T\psi|^2 -O(a)\big|(\nab_\T, \dkb)^{\leq 1}\dk^{\leq 2}\psi\big|^2\\
 &&+|q|^2 \Ddot_\b(|q|^{-2}O^{\a\b}\Ddot_\a \nab_\T\psi \c \SS_1 \nab_\T\psi ) +|q|^4 \Ddot_\b(|q|^{-2}O^{\a\b}\Ddot_\a \nab\psi \c \SS_1\nab\psi )\\
 && +\Ddot_\mu(M(\nab_\T\psi)\T^\mu) +\Ddot_\mu(M(|q|\nab\psi)\T^\mu) -\Ddot_\mu(|q|^2\nab\nab_\T\psi\c\nab\psi\T^\mu)+r^2\textrm{Good}_\SS.
\eeaa
The remainder of the proceeds exactly as for the proof of Lemma \ref{LEMMA:LOWERBOUNDPHIZOUTSIDEMTRAP} in section \ref{section:lowerboundPhizoutsideMtrap}, and leads to  a universal constant $c_0>0$ such that the following holds on $\Mntrap$
\beaa
r^3\Big(|\nab_\T\Psi_z|^2+r^2|\nab\Psi_z|^2\Big)+r^{-3}|\psi|_{\SS}^2 &\geq& c_0r^{-3}\Big(|\nab_\T\psi|^2_{\SS}+|\nab_Z\psi|^2_{\SS}+r^2|\nab\psi|^2_{\SS}\Big) \\
&& -O(ar^{-3})\big|(\nab_\T, \dkb)^{\leq 1}\dk^{\leq 2}\psi\big|^2+\Ddot_\a F^\a+\err_\ep,
\eeaa
where the additional error term $\err_\ep$ is given in view of the above by 
\beaa
\err_\ep = \textrm{Good}_\SS
\eeaa
and thus satisfies in view of  the control of $\textrm{Good}_\SS$ provided by Remark \ref{rmk::accetableerrortermsenergymorawetzundercontrolafsodsudfh:SSversion}
\beaa
\int_{\MM}|\err_\ep| &\les& \ep\left(\sup_{[\tau_1, \tau_2]}E^2[\psi](\tau)+ B^2_\de[\psi](\tau_1, \tau_2)\right)
\eeaa
as stated. This concludes the proof of Lemma \ref{lemma:lowerboundPhizoutsideMtrap:perturbation}.


\section{Non-integrable Hodge estimates}
\lab{section:nonintegrableHodge}


 In this section,  we  derive elliptic type estimates by projecting the basic horizontal Hodge operators on  surfaces $S$  of constant   $(\tau, r)$ in $\MM$  and by relying on the  standard elliptic estimates   on spheres,  see  Proposition \ref{prop:2D-hodge}. We start with the following lemma.

  \begin{lemma}\lab{Lemma:projectionS}
  Let $\eS_b$  denote  the projections of  $e_b$ to the spheres $S$ of  fixed  $(\tau, r)$. 
  Then
  \bea
\eS_b&=&\big(1+ O( a r^{-2} )\big) e_b+O(ar^{-1} )\T+r\Ga_g( e_3, e_4). 
\eea  
  \end{lemma}
 
 \begin{proof}
 Note that
\beaa
\T(r)= r\Ga_b,\qquad \T(\tau)=1+ r\Ga_b, \qquad \nab(r)=r\Ga_g, \qquad \nab(\tau)=a \Re( \Jk)+\Ga_b.
\eeaa
  Thus, setting  $e_b'=  e_b+ \La_b\T$ for a 1-form $\La$, and assuming $\La=O(r^{-1})$, we have
  \beaa
  e_b'(r)=r\Ga_g, \qquad e_b'(\tau)=  a \Re( \Jk)_b+\La_b+\Ga_b.
  \eeaa
  We may thus  choose $\La= - a \Re( \Jk)$, which indeed satisfies $\La=O(r^{-1})$, and deduce $ e_b'(r),  e_b'(\tau)=r\Ga_g$.
     
We then look for a pair of vectorfields  $\widetilde{e}_b^S$, $b=1,2$, under the form $\widetilde{e}_b^S= e_b'+\la e_3 +\underline{\la}  e_4$ and choose $\la$ and $\underline{\la}$ to enforce $\widetilde{e}_b^S(\tau)=0$ and $\widetilde{e}_b^S(r)=0$ so that $\widetilde{e}_b^S$, $b=1,2$ are tangent to $S$. We  infer
\beaa
0 &=& \widetilde{e}_b^S(r)=e_b'(r)+\la e_3(r) +\underline{\la}  e_4(r)=\left(-1+r\Ga_b\right)\la +\left(\frac{\De}{|q|^2}+r\Ga_b\right)\underline{\la}+ r\Ga_g,\\
0 &=& \widetilde{e}_b^S(\tau)=e_b'(\tau)+\la e_3(\tau) +\underline{\la}  e_4(\tau)=\la e_3(\tau) +\underline{\la}  e_4(\tau)+ r\Ga_g.
\eeaa  
Since $e_3(\tau)>0$ and $e_4(\tau)>0$, see Definition \ref{definition:definition-oftau}, we infer that $\la, \underline{\la}$ exist and satisfy $\la, \underline{\la}\in r\Ga_g $. Hence
 \beaa
\widetilde{e}_b^S &=& e_b'+\la e_3 +\underline{\la}  e_4=e_b - a \Re( \Jk)_b\T +r\Ga_g(e_3, e_4)\\
&=& e_b +O(ar^{-1})\T +r\Ga_g(e_3, e_4).
\eeaa
We then apply Gram-Schmidt to $\widetilde{e}_b^S$, $b=1,2$ so obtain an orthonormal frame $e_b^S$, $b=1,2$ of $S$ satisfying the desired properties.   This concludes the proof of Lemma \ref{Lemma:projectionS}.
    \end{proof}

    \begin{proposition}
    \lab{Prop:HodgeThmM8}
       Consider a sphere   $S \subset \MM$ of the type $S(\tau, r)$.
The following    estimates    hold true for $a$ small enough:
\begin{enumerate}
\item For $f\in \sk_p$, $p=1,2$, we have
\bea
\lab{eq:Prop-HodgeMint12}
\bsplit
\int_S \big(r^2|\nab f|^2 + |f|^2\big) \les&  \int_S r^2 |\DDd_p f|^2 +O(a^2)\int_S   |\nab_\T f|^2 \\
& + O(\ep^2)  \int_S  | (\nab_3, \nab_4)  f|^2. 
\end{split}
\eea
Moreover for higher derivatives, in a simplified form, we have
\bea
\lab{eq:Prop-HodgeMint12-higherderiv}
\bsplit
\int_S \big(r^2| \nab \dk^{ \le k}   f|^2 +|\dk^{ \le k}   f|^2\big)  &\les  \int_S r^2 |\dk^{\le k} \DDd_p f|^2 +
O(a^2, \ep^2 )\int_S    |\dk^{\le k+1}  f|^2.
\end{split}
\eea

\item For $f\in \sk_p$, $p=0,1$, we have
\bea
\lab{eq:Prop-HodgeMint01star}
\bsplit
\int_S r^2  |\nab f|^2 \les&  \int_S\big(r^2 |\DDs_p f|^2+ |f|^2\big)  +O(a^2)\int_S   |\nab_\T f|^2  \\
&+ O(\ep^2)  \int_S  | (\nab_3, \nab_4)  f|^2. 
\end{split}
\eea
Moreover for higher derivatives, in a simplified form, we have 
\bea
\lab{eq:Prop-HodgeMint01star-higherderiv}
\bsplit
\int_S  r^2| \nab \dk^{ \le k}   f|^2  &\les  \int_S \big(r^2 |\dk^{\le k} \DDd_p f|^2 +|\dk^{\le k} \DDd_p f|^2\big) +
O(a^2, \ep^2 )\int_S    |\dk^{\le k+1}  f|^2.
\end{split}
\eea
\end{enumerate}
\end{proposition} 

     \begin{proof}
     We start with the proof of \eqref{eq:Prop-HodgeMint12}, and consider the case $p=2$, i.e. $f\in \sk_2$.  Recalling from Lemma \ref{Lemma:projectionS} that $\eS_b=\big(1+ O( a r^{-2} )\big) e_b+O(ar^{-1} )\T+r\Ga_g( e_3, e_4) $, we infer
        \beaa
     \nab^b f_{ab} &=& \nabS ^b f_{ab} +O(ar^{-2})\nab f+ O(ar^{-1} ) \nab_\T f +r\Ga_g \nab_3 f. 
     \eeaa
     Hence
     \beaa
\DDd_2 f =   \DDdS_2 f +O(ar^{-2})\nab f +O(ar^{-1} ) \nab_\T f  +r\Ga_g( \nab_3, \nab_4) f.
     \eeaa
         Integrating over $S$ and using  the standard   elliptic estimates of Proposition \ref{prop:2D-hodge} we have
     \beaa
     \int_S|\nabS f|^2 + 2 \, K^S |f|^2  =2 \int_S|\DDdS_2 f|^2
     \eeaa
      where $K^S$ is the Gauss curvature of  $S$.
      We deduce
      \beaa
        \int_S\big(|\nabS f|^2 + 2 \, K^S |f|^2\big)  &\les& 2 \int_S|\DDd_2 f|^2+a^2\int_Sr^{-4}|\nab f|^2+ a^2\int_S   r^{-2} |\nab_\T f|^2 \\
        && + \int_S|r\Ga_g( \nab_3, \nab_4) f|^2. 
      \eeaa
      Also, recall from the proof of Lemma \ref{lemma:poincareinequalityfornabonSasoidfh:chap9} that the Gauss curvature $K_S$ verifies 
\beaa
K_S &=& \frac{1}{r^2}\big(1+O(\ep+a^2r^{-2})\big).
\eeaa
Hence
\beaa
        \int_S\left(|\nabS f|^2 + \frac{2}{r^2} |f|^2\right)  &\les& 2 \int_S|\DDd_2 f|^2+a^2\int_Sr^{-4}|\nab f|^2+a^2\int_Sr^{-4}|f|^2\\
        &&+ a^2\int_S   r^{-2} |\nab_\T f|^2  + \int_S |r\Ga_g( \nab_3, \nab_4) f|^2. 
      \eeaa      
            Also
      \beaa
     \int_S |\nabS f|^2 &=&  \int_S \big(1+ O( a^2 r^{-4} )\big) |\nab f|^2  + O(a^2)\int_S   r^{-2} |\nab_\T f|^2  + \int_S |r\Ga_g (\nab_3, \nab_4) f|^2
      \eeaa
      and thus, absorbing terms on the RHS for $a$ small enough, we infer
      \beaa
       \int_S\left(|\nab f|^2 + \frac{2}{r^2}|f|^2\right)  &\les & \int_S|\DDd_2 f|^2+a^2\int_S   r^{-2} |\nab_\T f|^2  + \ep^2  \int_S r^{-2} |(\nab_3, \nab_4) f|^2. 
      \eeaa
      This ends the proof of  \eqref{eq:Prop-HodgeMint12} for $p=2$.  The case $p=1$ is derived in the same manner.  
      
Next, we consider the higher derivative estimates \eqref{eq:Prop-HodgeMint12-higherderiv}. For $p=1, 2$ and $k\geq 1$, we have the following non sharp commutator estimate 
\beaa
\DD_p\dk^kf &=& \dk^k\DD_pf+[\DD_p, \dk^k]f\\
&=& \dk^k\DD_pf+O(r^{-1})\nab\dk^{\leq k-1}f +O(ar^{-2})\dk^{\leq k}f +O(r^{-1})\dk^{\leq k-1}f+O(\ep r^{-1})\dk^{\leq k}f.
\eeaa
Hence, applying \eqref{eq:Prop-HodgeMint12}, we deduce
\beaa
\int_S r^2|\nab\dk^k f|^2+|\dk^k f|^2 &\les& \int_S|\dk^k\DD_pf|^2 +\int_S r^2|\nab\dk^{\leq k-1} f|^2+|\dk^{\leq k-1} f|^2\\
&& +O(a^2, \ep^2 )\int_S    |\dk^{\le k+1}  f|^2
\eeaa
and \eqref{eq:Prop-HodgeMint12-higherderiv} follows immediately by iteration starting from the case $k=0$ provided by  \eqref{eq:Prop-HodgeMint12}.
      
The proof of \eqref{eq:Prop-HodgeMint01star} is similar to the one of  \eqref{eq:Prop-HodgeMint12}. Then, the proof of \eqref{eq:Prop-HodgeMint01star-higherderiv} is derived by  iteration starting from the case $k=0$ provided by \eqref{eq:Prop-HodgeMint01star} as above. This concludes the proof of  Proposition \ref{Prop:HodgeThmM8}.    
\end{proof}


\section{Redshift estimates}
\lab{section:Redshift-estimates-chp3}


To remove degeneracies of the energy in the neighborhood of the horizon, we make use of the Dafermos-Rodnianski redshift vectorfield. The goal of this section is to prove the following proposition. 

\begin{proposition}[Redshift estimates]
\lab{Prop:Redshift-estimates-chp3}
Let $\psi$ a solution to \eqref{eq:Gen.RW}. Then, for $|a|<m$, there exists a small enough constant $\de_{red}>0$ such that $\de_{red}=\de_{red}(m-|a|)$ with $\de_{red}\geq \deh$, and a small constant $c_0>0$ with $c_0=c_0(m-|a|)$, such that the following estimate  holds true in $\MM(\tau_1, \tau_2) $:
\bea
\bsplit
&c_0E_{r\leq r_+(1+\de_{red})}[\psi](\tau_2)+c_0\textrm{Mor}_{r\leq r_+(1+\de_{red})}[\psi](\tau_1, \tau_2)+c_0F_{\AA}[\psi](\tau_1, \tau_2)\\
  \leq& E_{r\leq r_+(1+2\de_{red})}[\psi](\tau_1) +\de_{red}^{-1}\textrm{Mor}_{r_+(1+\de_{red})\leq r\leq r_+(1+2\de_{red})}[\psi](\tau_1, \tau_2)\\
&+   \int_{\MM(\tau_1, \tau_2)\cap\left\{\frac{r}{r_+}\le 1+2\de_{red}\right\}}|N|^2.
\end{split}
\eea 
\end{proposition}

The proof of Proposition \ref{Prop:Redshift-estimates-chp3} is postponed to the end of this section. We start with the following lemma.
\begin{lemma}\lab{lemma:basiccomputationforredshiftvectorfield}
Given the vectorfield 
\beaa
Y=\underline{d}(r)e_3+d(r)e_4,
\eeaa
and assuming
\beaa
\sup_{r\leq 4m}\Big(|d(r)|+|d'(r)|+|\underline{d}(r)|+|\underline{d}'(r)|\Big) &\les& 1,
\eeaa 
we have for $r\leq 4m$
\beaa
\QQ\c{}^{(Y)}\pi &=& \left(\pr_r\left(\frac{\De}{|q|^2}\right)\underline{d}(r) -\frac{\De}{|q|^2}\underline{d}'(r)\right)|\nab_3\psi|^2 +d'(r)|\nab_4\psi|^2\\
&&+\left(\underline{d}'(r) -\pr_r\left(\frac{\De}{|q|^2}\right)d(r) -\frac{\De}{|q|^2}d'(r)\right)\Big(|\nab\psi|^2 +V|\psi|^2\Big)\\
&& +\frac{2ar}{|q|^2}\Re(\Jk)^b\underline{d}(r)\nab_b\Psi\c\nab_3\Psi -\frac{2r}{|q|^2}\left(\underline{d}(r)-\frac{\De}{|q|^2}d(r)\right)\Big(\nab_3\psi\c\nab_4\psi -V|\psi|^2\Big)\\
&&+O(\ep)\Big(|\nab_3\psi|^2+|\nab_4\psi|^2+|\nab\psi|^2+r^{-2}|\psi|^2\Big)
\eeaa
and 
 \beaa
 \D^\mu  \PP_\mu[Y, 0, 0]  &=& \frac 1 2 \QQ  \c {}^{(Y)}\pi  +2\left(\underline{d}(r) -\frac{\De}{|q|^2}d(r)\right)\pr_r\left(\frac{\De}{ (r^2+a^2) |q|^2}\right) |\psi|^2\\
 && - \rhod  \in_{AB} \Big(d(r)  \nab_4  \psi ^A \psi^B - \underline{d}(r)  \nab_3   \psi ^A \psi^B\Big) +  Y( \psi )\c\left(\squared_2\psi- V\psi\right)\\
 && +O(\ep)\Big(|\nab_3\psi|^2+|\nab_4\psi|^2+|\nab\psi|^2+r^{-2}|\psi|^2\Big).
 \eeaa

Also, we have for $r\leq 4m$
 \beaa
 \D^\mu  \PP_\mu[\T, 0, 0]  &=&    -\frac{1}{2}\rhod\in_{AB}\left(\frac{\De}{|q|^2}\nab_4\psi^A\psi^B -\nab_3\psi^A\psi^B\right)  -a\dual\Re(\Jk)^d\nab_d\psi_A\rho\dual\psi^A\\
 && +  \T( \psi )\c \left(\squared_2 \psi- V\psi\right)+ O(\ep)\Big(|\nab_3\psi|^2+|\nab_4\psi|^2+|\nab\psi|^2+r^{-2}|\psi|^2\Big).
 \eeaa
\end{lemma}

\begin{proof}
We have 
\beaa
{}^{(Y)}\pi_{\mu\nu} &=& \underline{d}(r){}^{(e_3)}\pi_{\mu\nu}+e_\mu(r)\underline{d}'(r)\g(e_3, e_\nu)+e_\nu(r)\underline{d}'(r)\g(e_3, e_\mu)\\
&&+d(r){}^{(e_4)}\pi_{\mu\nu}+e_\mu(r)d'(r)\g(e_4, e_\nu)+e_\nu(r)d'(r)\g(e_4, e_\mu).
\eeaa
Together with the fact, in the region $r\leq 4m$,
\beaa
e_3(r)=-1+O(\ep), \qquad e_4(r)=\frac{\De}{|q|^2}+O(\ep), \qquad \nab(r)=O(\ep),
\eeaa
and using Lemma \ref{lemma:deformationtensore3}, we infer
\beaa
{}^{(Y)}\pi_{44} &=& \underline{d}(r){}^{(e_3)}\pi_{44} +2e_4(r)\underline{d}'(r)\g(e_3, e_4) +d(r){}^{(e_4)}\pi_{44}\\
&=& 4\pr_r\left(\frac{\De}{|q|^2}\right)\underline{d}(r) -\frac{4\De}{|q|^2}\underline{d}'(r)+O(\ep),\\
{}^{(Y)}\pi_{34} &=& \underline{d}(r){}^{(e_3)}\pi_{34}+e_3(r)\underline{d}'(r)\g(e_3, e_4)+d(r){}^{(e_4)}\pi_{34}+e_4(r)d'(r)\g(e_4, e_3) ,\\
&=&  2\underline{d}'(r) -2\pr_r\left(\frac{\De}{|q|^2}\right)d(r) -\frac{2\De}{|q|^2}d'(r) +O(\ep),\\
{}^{(Y)}\pi_{33} &=& \underline{d}(r){}^{(e_3)}\pi_{33}+d(r){}^{(e_4)}\pi_{33}+2e_3(r)d'(r)\g(e_4, e_3) = 4d'(r)+O(\ep),\\
{}^{(Y)}\pi_{4b} &=& \underline{d}(r){}^{(e_3)}\pi_{4b}+d(r){}^{(e_4)}\pi_{4b} = -\frac{2ar}{|q|^2}\Re(\Jk)_b\underline{d}(r)+O(\ep),\\
{}^{(Y)}\pi_{3b} &=& \underline{d}(r){}^{(e_3)}\pi_{3b}+d(r){}^{(e_4)}\pi_{3b} = O(\ep),\\
{}^{(Y)}\pi_{bc} &=& \underline{d}(r){}^{(e_3)}\pi_{bc}+d(r){}^{(e_4)}\pi_{bc} = -\frac{2r}{|q|^2}\left(\underline{d}(r)-\frac{\De}{|q|^2}d(r)\right)\de_{bc}+O(\ep).
\eeaa
This yields, in the region $r\leq 4m$,
\beaa
\QQ\c{}^{(Y)}\pi &=& \left(\pr_r\left(\frac{\De}{|q|^2}\right)\underline{d}(r) -\frac{\De}{|q|^2}\underline{d}'(r)\right)\QQ_{33} +d'(r)\QQ_{44}\\
&&+\left(\underline{d}'(r) -\pr_r\left(\frac{\De}{|q|^2}\right)d(r) -\frac{\De}{|q|^2}d'(r)\right)\QQ_{34}\\
&& +\frac{2ar}{|q|^2}\Re(\Jk)^b\underline{d}(r)\QQ_{3b} -\frac{2r}{|q|^2}\left(\underline{d}(r)-\frac{\De}{|q|^2}d(r)\right)\de^{bc}\QQ_{bc}\\
&&+O(\ep)\Big(|\nab_3\psi|^2+|\nab_4\psi|^2+|\nab\psi|^2+r^{-2}|\psi|^2\Big).
\eeaa
Since we have in view of \eqref{eq:definition-QQ-mu-nu}
\beaa
\bsplit
\QQ_{33}&=|\nab_3\psi|^2, \qquad  \QQ_{44}=|\nab_4\psi|^2, \qquad \QQ_{34}=|\nab\psi|^2 +V|\psi|^2,\\
\QQ_{4a}&=\nab_4\Psi\c\nab_a\Psi, \qquad \QQ_{3a}=\nab_3\Psi\c\nab_a\Psi, \qquad \de^{bc}\QQ_{bc} = \nab_3\psi\c\nab_4\psi -V|\psi|^2,
\end{split}
\eeaa
we deduce
\beaa
\QQ\c{}^{(Y)}\pi &=& \left(\pr_r\left(\frac{\De}{|q|^2}\right)\underline{d}(r) -\frac{\De}{|q|^2}\underline{d}'(r)\right)|\nab_3\psi|^2 +d'(r)|\nab_4\psi|^2\\
&&+\left(\underline{d}'(r) -\pr_r\left(\frac{\De}{|q|^2}\right)d(r) -\frac{\De}{|q|^2}d'(r)\right)\Big(|\nab\psi|^2 +V|\psi|^2\Big)\\
&& +\frac{2ar}{|q|^2}\Re(\Jk)^b\underline{d}(r)\nab_b\Psi\c\nab_3\Psi -\frac{2r}{|q|^2}\left(\underline{d}(r)-\frac{\De}{|q|^2}d(r)\right)\Big(\nab_3\psi\c\nab_4\psi -V|\psi|^2\Big)\\
&&+O(\ep)\Big(|\nab_3\psi|^2+|\nab_4\psi|^2+|\nab\psi|^2+r^{-2}|\psi|^2\Big).
\eeaa
as stated. 

Also, in view of Proposition \ref{prop-app:stadard-comp-Psi-perturbations-Kerr}, we have, in the region $r\leq 4m$,
 \beaa
 \D^\mu  \PP_\mu[Y, 0, 0]  &=& \frac 1 2 \QQ  \c {}^{(Y)}\pi - \frac 1 2 Y( V ) |\psi|^2 - \rhod  \in_{AB} \Big(d(r)  \nab_4  \psi ^A \psi^B - \underline{d}(r)  \nab_3   \psi ^A \psi^B\Big)\\
 && +  Y( \psi )\c\left(\squared_2\psi- V\psi\right) +O(\ep)\Big(|\nab_3\psi|^2+|\nab_4\psi|^2+|\nab\psi|^2+r^{-2}|\psi|^2\Big).
 \eeaa
Since $V= \frac{4\De}{ (r^2+a^2) |q|^2}$ in view of \eqref{eq:Gen.RW}, and since, in the region $r\leq 4m$,
\beaa
e_3(r)=-1+O(\ep), \qquad e_4(r)=\frac{\De}{|q|^2}+O(\ep), \qquad e_4(\cos\th)=O(\ep), \quad e_3(\cos\th)=O(\ep),
\eeaa
  we have
  \beaa
  Y( V ) &=&  -4\left(\underline{d}(r) -\frac{\De}{|q|^2}d(r)+O(\ep)\right)\pr_r\left(\frac{\De}{ (r^2+a^2) |q|^2}\right). \eeaa
We deduce
 \beaa
 \D^\mu  \PP_\mu[Y, 0, 0]  &=& \frac 1 2 \QQ  \c {}^{(Y)}\pi  +2\left(\underline{d}(r) -\frac{\De}{|q|^2}d(r)\right)\pr_r\left(\frac{\De}{ (r^2+a^2) |q|^2}\right) |\psi|^2\\
 && - \rhod  \in_{AB} \Big(d(r)  \nab_4  \psi ^A \psi^B - \underline{d}(r)  \nab_3   \psi ^A \psi^B\Big) +  Y( \psi )\c\left(\squared_2\psi- V\psi\right)\\
 && +O(\ep)\Big(|\nab_3\psi|^2+|\nab_4\psi|^2+|\nab\psi|^2+r^{-2}|\psi|^2\Big)
 \eeaa
 as stated.

Finally, we have $\piT_{\mu\nu}=O(\ep)$ in the region $r\leq 4m$ so that 
\beaa
\QQ\c\piT &=& O(\ep)\Big(|\nab_3\psi|^2+|\nab_4\psi|^2+|\nab\psi|^2+r^{-2}|\psi|^2\Big).
\eeaa
Together with Proposition \ref{prop-app:stadard-comp-Psi}, and using that $\T(V)=O(\ep)$  in the region $r\leq 4m$, this implies
 \beaa
 \D^\mu  \PP_\mu[\T, 0, 0]  &=&    \T^\mu \Db^\nu  \psi ^a\R_{ ab   \nu\mu}\psi^b+  \T( \psi )\c \left(\squared_2 \psi- V\psi\right)\\
  &&+ O(\ep)\Big(|\nab_3\psi|^2+|\nab_4\psi|^2+|\nab\psi|^2+r^{-2}|\psi|^2\Big).
 \eeaa
 Also, using the computations in the proof of Proposition \ref{prop-app:stadard-comp-Psi-perturbations-Kerr}, we have,  in the region $r\leq 4m$,
 \beaa
 &&\T^\mu \Db^\nu  \psi^A\R_{AB\nu\mu}\psi^B\\
  &=& -\rhod\in_{AB}\Big(\T^4\Db_4\psi^A\psi^B -\T^3\Db_3\psi^A\psi^B\Big)\\
 &&+\T^c\left(-\frac{1}{2}\Db_3\psi^A\R_{AB4c}\psi^B -\frac{1}{2}\Db_4\psi^A\R_{AB3c}\psi^B + \Db^d\psi^A\R_{ABdc}\psi^B\right)\\
 && + O(\ep)\Big(|\nab_3\psi|^2+|\nab_4\psi|^2+|\nab\psi|^2+r^{-2}|\psi|^2\Big)\\
 &=& -\frac{1}{2}\rhod\in_{AB}\left(\frac{\De}{|q|^2}\nab_4\psi^A\psi^B -\nab_3\psi^A\psi^B\right)\\
 && -a\Re(\Jk)^c\left(-\frac{1}{2}\nab_3\psi^A\in_{AB}\dual\b_c\psi^B +\frac{1}{2}\nab_4\psi^A\in_{AB}\dual\bb_c\psi^B + \nab^d\psi^A\in_{AB}\in_{cd}\rho\psi^B\right)\\
 && + O(\ep)\Big(|\nab_3\psi|^2+|\nab_4\psi|^2+|\nab\psi|^2+r^{-2}|\psi|^2\Big)\\
&=& -\frac{1}{2}\rhod\in_{AB}\left(\frac{\De}{|q|^2}\nab_4\psi^A\psi^B -\nab_3\psi^A\psi^B\right)  -a\dual\Re(\Jk)^d\nab_d\psi_A\rho\dual\psi^A\\
 && + O(\ep)\Big(|\nab_3\psi|^2+|\nab_4\psi|^2+|\nab\psi|^2+r^{-2}|\psi|^2\Big)\\
 \eeaa
 and hence
 \beaa
 \D^\mu  \PP_\mu[\T, 0, 0]  &=&    -\frac{1}{2}\rhod\in_{AB}\left(\frac{\De}{|q|^2}\nab_4\psi^A\psi^B -\nab_3\psi^A\psi^B\right)  -a\dual\Re(\Jk)^d\nab_d\psi_A\rho\dual\psi^A\\
 && +  \T( \psi )\c \left(\squared_2 \psi- V\psi\right)+ O(\ep)\Big(|\nab_3\psi|^2+|\nab_4\psi|^2+|\nab\psi|^2+r^{-2}|\psi|^2\Big)
 \eeaa
as stated. This concludes the proof of Lemma \ref{lemma:basiccomputationforredshiftvectorfield}.
\end{proof}

\begin{corollary}\lab{cor:basiccomputationforredshiftvectorfieldatr=r+}
Let 
\beaa
Y_{(0)} :=Y+2T, \qquad Y=\underline{d}(r)e_3+d(r)e_4.
\eeaa
If we choose 
\beaa
d(r_+)=0, \qquad \underline{d}(r_+)=1, \qquad d'(r_+)\geq \frac{c_1}{r_+-m}, \qquad \underline{d}'(r_+)\geq \frac{c_1}{r_+-m},
\eeaa
for some large enough universal constant $c_1$, then at $r=r_+$, we have,  for any sphere $S=S(\tau, r_+)$,
\bea
\nn\int_S\D^\mu  \PP_\mu[Y_{(0)}, 0, 0]  &\geq& \frac{(r_+-m)}{32r_+^2}\int_S\Big(|\nab_3\psi|^2+|\nab_4\psi|^2+|\nab\psi|^2+r^{-2}|\psi|^2\Big)  \\
&& +\int_S Y_{(0)}( \psi )\c \left(\squared_2 \psi- V\psi +\frac{4a\cos\th}{|q|^2}\dual\nab_T\psi\right).
\eea
\end{corollary}

\begin{proof}
At $r=r+$, assuming that $\underline{d}(r_+)=1$ and $d(r_+)=0$, and using $\De(r_+)=0$, $\pr_r\left(\frac{\De}{|q|^2}\right)_{|_{r=r_+}}=\frac{2(r_+-m)}{|q|^2}$ and $V(r_+)=0$, we obtain from Lemma \ref{lemma:basiccomputationforredshiftvectorfield}
\beaa
\QQ\c{}^{(Y)}\pi &=& \frac{2(r_+-m)}{|q|^2}|\nab_3\psi|^2 +d'(r_+)|\nab_4\psi|^2+\underline{d}'(r_+)|\nab\psi|^2\\
&& +\frac{2ar}{|q|^2}\Re(\Jk)^b\nab_b\Psi\c\nab_3\Psi -\frac{2r}{|q|^2}\nab_3\psi\c\nab_4\psi\\
&&+O(\ep)\Big(|\nab_3\psi|^2+|\nab_4\psi|^2+|\nab\psi|^2+r^{-2}|\psi|^2\Big).
\eeaa
Also, since $\pr_r\left(\frac{\De}{ (r^2+a^2) |q|^2}\right)_{|_{r=r_+}}=\frac{2(r_+-m)}{(r^2+a^2)|q|^2}$, we obtain from Lemma \ref{lemma:basiccomputationforredshiftvectorfield}, at $r=r+$,
 \beaa
 \D^\mu  \PP_\mu[Y, 0, 0]  &=& \frac 1 2 \QQ  \c {}^{(Y)}\pi  +\frac{4(r_+-m)}{(r^2+a^2)|q|^2} |\psi|^2 + \rhod  \in_{AB}\nab_3   \psi ^A \psi^B \\
 &&+  Y( \psi )\c\left(\squared_2\psi- V\psi\right) +O(\ep)\Big(|\nab_3\psi|^2+|\nab_4\psi|^2+|\nab\psi|^2+r^{-2}|\psi|^2\Big).
 \eeaa
Together with Lemma \ref{lemma:basiccomputationforredshiftvectorfield}, and since $Y_{(0)}=Y+2T$, we infer, at $r=r+$, 
\beaa
\D^\mu  \PP_\mu[Y_{(0)}, 0, 0]  &=& \D^\mu  \PP_\mu[Y, 0, 0] +2\D^\mu  \PP_\mu[\T, 0, 0]\\
&=& \frac{(r_+-m)}{|q|^2}|\nab_3\psi|^2 +\frac{d'(r_+)}{2}|\nab_4\psi|^2+\frac{\underline{d}'(r_+)}{2}|\nab\psi|^2  +\frac{4(r_+-m)}{(r^2+a^2)|q|^2} |\psi|^2\\
&& +\frac{ar}{|q|^2}\Re(\Jk)^b\nab_b\Psi\c\nab_3\Psi -\frac{r}{|q|^2}\nab_3\psi\c\nab_4\psi + 2\rhod  \in_{AB}\nab_3   \psi ^A \psi^B \\
 && -2a\dual\Re(\Jk)^d\nab_d\psi_A\rho\dual\psi^A + Y_{(0)}( \psi )\c \left(\squared_2 \psi- V\psi\right)\\
 && + O(\ep)\Big(|\nab_3\psi|^2+|\nab_4\psi|^2+|\nab\psi|^2+r^{-2}|\psi|^2\Big).
\eeaa
We deduce, at $r=r_+$, 
\beaa
\D^\mu  \PP_\mu[Y_{(0)}, 0, 0]  &\geq& \frac{(r_+-m)}{4|q|^2}|\nab_3\psi|^2 +\left(\frac{d'(r_+)}{2} - \frac{r_+^2}{(r_+-m)|q|^2}\right)|\nab_4\psi|^2\\
&&+\left(\frac{\underline{d}'(r_+)}{2} - \frac{a^2r_+^2}{(r_+-m)|q|^2}|\Re(\Jk)|^2 -a|\Re(\Jk)||\rho|r_+\right)|\nab\psi|^2 \\
&& +\left(\frac{4(r_+-m)}{(r^2+a^2)|q|^2} - \frac{4|q|^2\rhod^2}{(r_+-m)} -\frac{a|\Re(\Jk)||\rho|}{r_+}\right)|\psi|^2\\
 && + Y_{(0)}( \psi )\c \left(\squared_2 \psi- V\psi\right)+ O(\ep)\Big(|\nab_3\psi|^2+|\nab_4\psi|^2+|\nab\psi|^2+r^{-2}|\psi|^2\Big).
\eeaa

Also, since $Y_{(0)}=Y+2T$, we have
\beaa
-\frac{4a\cos\th}{|q|^2}Y_{(0)}( \psi )\c\dual\nab_T\psi &=&  -\frac{4a\cos\th}{|q|^2}\nab_Y\psi\c\dual\nab_T\psi -\frac{8a\cos\th}{|q|^2}\nab_T\psi\c\dual\nab_T\psi\\
&=&  -\frac{4a\cos\th}{|q|^2}\nab_Y\psi\c\dual\nab_T\psi.
\eeaa
At $r=r_+$, we obtain 
\beaa
-\frac{4a\cos\th}{|q|^2}Y_{(0)}( \psi )\c\dual\nab_T\psi &=&  -\frac{2a\cos\th}{|q|^2}\nab_3\psi\c\nab_{e_4 -2a\Re(\Jk)^be_b}\dual\psi\\
&\geq& -\frac{2|a|}{|q|^2}|\nab_3\psi||\nab_4\psi| -\frac{4a^2}{|q|^2}|\Re(\Jk)||\nab_3\psi||\nab\psi|\\
&\geq& -\frac{(r_+-m)}{8|q|^2}|\nab_3\psi|^2 -\frac{16a^2}{(r_+-m)|q|^2}|\nab_4\psi|^2\\
&& -\frac{64a^4}{(r_+-m)|q|^2}|\Re(\Jk)|^2|\nab\psi|^2
\eeaa
and hence
\beaa
\D^\mu  \PP_\mu[Y_{(0)}, 0, 0]  &\geq& \frac{(r_+-m)}{8|q|^2}|\nab_3\psi|^2 +\left(\frac{d'(r_+)}{2} - \frac{r_+^2}{(r_+-m)|q|^2} -\frac{16a^2}{(r_+-m)|q|^2}\right)|\nab_4\psi|^2\\
&&+\left(\frac{\underline{d}'(r_+)}{2} - \frac{a^2r_+^2}{(r_+-m)|q|^2}|\Re(\Jk)|^2 -a|\Re(\Jk)||\rho|r_+  -\frac{64a^4}{(r_+-m)|q|^2}|\Re(\Jk)|^2\right)|\nab\psi|^2 \\
&& +\left(\frac{4(r_+-m)}{(r^2+a^2)|q|^2} - \frac{4|q|^2\rhod^2}{(r_+-m)} -\frac{a|\Re(\Jk)||\rho|}{r_+}\right)|\psi|^2\\
 && + Y_{(0)}( \psi )\c \left(\squared_2 \psi- V\psi +\frac{4a\cos\th}{|q|^2}\dual\nab_T\psi\right)\\
 &&+ O(\ep)\Big(|\nab_3\psi|^2+|\nab_4\psi|^2+|\nab\psi|^2+r^{-2}|\psi|^2\Big).
\eeaa

Since we have, for $r\leq 4m$,
\beaa
r_+>m, \quad r\leq |q|\leq 2r, \quad |\rho|\leq \frac{2m}{r^3}+O(\ep), \quad |\rhod|\leq \frac{6|a|m}{r^4}+O(\ep), \quad |\Re(\Jk)|\leq \frac{1}{r}+O(\ep),
\eeaa
we infer, at $r=r_+$, 
\beaa
\D^\mu  \PP_\mu[Y_{(0)}, 0, 0]  &\geq& \frac{(r_+-m)}{16r_+^2}|\nab_3\psi|^2 +\left(\frac{d'(r_+)}{2} - \frac{17}{r_+-m}\right)|\nab_4\psi|^2\\
&&+\left(\frac{\underline{d}'(r_+)}{2} - \frac{67m|a|}{r_+^2(r_+-m)}\right)|\nab\psi|^2  +\left(\frac{(r_+-m)}{4r_+^4} - \frac{578|a|m}{(r_+-m)r_+^4}\right)|\psi|^2\\
 && + Y_{(0)}( \psi )\c \left(\squared_2 \psi- V\psi +\frac{4a\cos\th}{|q|^2}\dual\nab_T\psi \right)\\
 &&+ O(\ep)\Big(|\nab_3\psi|^2+|\nab_4\psi|^2+|\nab\psi|^2+r^{-2}|\psi|^2\Big).
\eeaa

Next, using the Poincar\'e inequality of Lemma \ref{lemma:poincareinequalityfornabonSasoidfh:chap9}, we have for some universal constant $c_0>0$, and for any sphere $S=S(\tau, r)$, 
\beaa
\int_S\left( |\nab_4\psi|^2 + \frac{\De}{r^4} |\nab_3\psi|^2 +|\nab\psi|^2\right) \geq \frac{c_0}{r^2}\int_S|\psi|^2.
\eeaa
In particular, we deduce at $r=r_+$
\beaa
\int_S\left( |\nab_4\psi|^2  +|\nab\psi|^2\right) \geq \frac{c_0}{r^2}\int_S|\psi|^2,
\eeaa
and hence, for any sphere $S=S(\tau, r_+)$
\beaa
\int_S\D^\mu  \PP_\mu[Y_{(0)}, 0, 0]  &\geq& \frac{(r_+-m)}{16r_+^2}\int_S|\nab_3\psi|^2 +\frac{(r_+-m)}{4r_+^4}\int_S|\psi|^2\\
&&+\left(\frac{d'(r_+)}{2} - \frac{17}{r_+-m} -\frac{1}{c_0}\frac{578m|a|}{(r_+-m)r_+^2}\right)\int_S|\nab_4\psi|^2\\
&&+\left(\frac{\underline{d}'(r_+)}{2} - \frac{67m|a|}{r_+^2(r_+-m)} -\frac{1}{c_0}\frac{578m|a|}{(r_+-m)r_+^2}\right)\int_S|\nab\psi|^2 \\
 && + \int_SY_{(0)}( \psi )\c \left(\squared_2 \psi- V\psi +\frac{4a\cos\th}{|q|^2}\dual\nab_T\psi\right)\\
 &&+ O(\ep)\int_S\Big(|\nab_3\psi|^2+|\nab_4\psi|^2+|\nab\psi|^2+r^{-2}|\psi|^2\Big).
\eeaa
Thus, if we choose 
\beaa
d(r_+)=0, \qquad \underline{d}(r_+)=1, \qquad d'(r_+)\geq \frac{c_1}{r_+-m}, \qquad \underline{d}'(r_+)\geq \frac{c_1}{r_+-m},
\eeaa
for some large enough universal constant $c_1$, then at $r=r_+$, we have,  for any sphere $S=S(\tau, r_+)$,
\beaa
\nn\int_S\D^\mu  \PP_\mu[Y_{(0)}, 0, 0]  &\geq& \frac{(r_+-m)}{32r_+^2}\int_S\Big(|\nab_3\psi|^2+|\nab_4\psi|^2+|\nab\psi|^2+r^{-2}|\psi|^2\Big)  \\
&& +\int_S Y_{(0)}( \psi )\c \left(\squared_2 \psi- V\psi +\frac{4a\cos\th}{|q|^2}\dual\nab_T\psi\right)
\eeaa
as stated. This concludes the proof of Corollary \ref{cor:basiccomputationforredshiftvectorfieldatr=r+}.
 \end{proof}

\begin{proposition}\lab{prop:basiccomputationforredshiftvectorfieldnearr=r+fordivergenceandboundary}
Let $\ka(r)$ a positive bump function supported in $[2,2]$ and equal to 1 on $[-1,1]$. Also, for $|a|<m$, let a small enough constant $\de_{red}>0$ such that $\de_{red}=\de_{red}(m-|a|)$ with $\de_{red}\geq \deh$. Let $Y_\HH$ the vectorfield given by 
\bea
Y_\HH := \ka_\HH Y_{(0)}, \qquad \ka_\HH=\ka\left(\frac{\frac{r}{r_+}-1}{\de_{red}}\right),
\eea
where $Y_{(0)}$ is the vectorfield of Corollary \ref{cor:basiccomputationforredshiftvectorfieldatr=r+}. 
Then, the following estimate holds, for any sphere $S=S(\tau, r)$,
\bea
\nn\int_S\D^\mu  \PP_\mu[Y_\HH, 0, 0]  &\geq& \frac{(r_+-m)}{64r_+^2}\int_S\Big(|\nab_3\psi|^2+|\nab_4\psi|^2+|\nab\psi|^2+r^{-2}|\psi|^2\Big)\mathbb{1}_{\left|\frac{r}{r_+}-1\right|\leq \de_{red}}  \\
\nn&& -O(r_+^{-1}\de_{red}^{-1})\int_S\Big(|\nab_3\psi|^2+|\nab_4\psi|^2+|\nab\psi|^2+r^{-2}|\psi|^2\Big)\mathbb{1}_{\de_{red}\leq\left|\frac{r}{r_+}-1\right|\leq 2\de_{red}}\\
&& +\int_S Y_\HH( \psi )\c \left(\squared_2 \psi- V\psi +\frac{4a\cos\th}{|q|^2}\dual\nab_T\psi\right).
\eea
Also, recallling that the vectorfield $N_\Si$ is given by $N_\Si=-\g^{\a\b}\pr_\a(\tau)\pr_\b$,  there holds on $\MM$
\bea
\PP_\mu[Y_\HH, 0, 0]\c N_\Si &\geq& 0,
\eea
and, there exists, for $|a|<m$, a constant $c_0>0$, with $c_0=c_0(m-|a|)$, such that we have for any sphere $S=S(\tau, r)$ with $|\frac{r}{r_+}-1|\leq \de_{red}$,
\bea
\int_S\PP_\mu[Y_\HH, 0, 0]\c N_\Si &\geq& c_0\int_S\Big(|\nab_3\psi|^2+|\nab_4\psi|^2+|\nab\psi|^2+r^{-2}|\psi|^2\Big).
\eea
Finally,  there exists, for $|a|<m$, a constant $c_0>0$, with $c_0=c_0(m-|a|)$, such that we have for any sphere $S=S(\tau, r)$ with $r=(1-\deh)r_+$, i.e. any sphere in $\AA$, 
\bea
\int_S\PP_\mu[Y_\HH, 0, 0]\c N_\AA &\geq &  c_0\int_S\Big(\deh|\nab_3\psi|^2+|\nab_4\psi|^2+|\nab\psi|^2+r^{-2}|\psi|^2\Big).
\eea
\end{proposition}

\begin{proof}
In view of the definition of $Y_\HH$ and the properties of $\ka_\HH$, we have
\beaa
\D^\mu  \PP_\mu[Y_\HH, 0, 0]  &=& \frac 1 2 \QQ  \c {}^{(Y_\HH)}\pi+ Y_\HH( \psi) \c  \left(\squared_2 \psi- V\psi\right)-\frac 1 2 Y_\HH(V) |\psi|^2 \\
&&+ Y_\HH^\mu \Db^\nu  \psi ^a\R_{ ab   \nu\mu}\psi^b\\
&=& \ka_\HH \D^\mu  \PP_\mu[Y_{(0)}, 0, 0] +\QQ(Y_{(0)}, d\ka_\HH)\\
&=& \ka_\HH \D^\mu  \PP_\mu[Y_{(0)}, 0, 0]\\
&& - O(r_+^{-1}\de_{red}^{-1})\mathbb{1}_{\de_{red}\leq\left|\frac{r}{r_+}-1\right|\leq 2\de_{red}}\Big(|\nab_3\psi|^2+|\nab_4\psi|^2+|\nab\psi|^2+r^{-2}|\psi|^2\Big).
\eeaa
Together with the lower bound for $\D^\mu  \PP_\mu[Y_{(0)}, 0, 0]$ of Corollary \ref{cor:basiccomputationforredshiftvectorfieldatr=r+} and the properties of $\ka_\HH$, we infer, provided $\de_{red}$ is chosen sufficiently small, for any sphere $S=S(\tau, r)$,
\beaa
\nn\int_S\D^\mu  \PP_\mu[Y_\HH, 0, 0]  &\geq& \frac{(r_+-m)}{64r_+^2}\int_S\Big(|\nab_3\psi|^2+|\nab_4\psi|^2+|\nab\psi|^2+r^{-2}|\psi|^2\Big)\mathbb{1}_{\left|\frac{r}{r_+}-1\right|\leq \de_{red}}  \\
\nn&& -O(r_+^{-1}\de_{red}^{-1})\int_S\Big(|\nab_3\psi|^2+|\nab_4\psi|^2+|\nab\psi|^2+r^{-2}|\psi|^2\Big)\mathbb{1}_{\de_{red}\leq\left|\frac{r}{r_+}-1\right|\leq 2\de_{red}}\\
&& +\int_S Y_\HH( \psi )\c \left(\squared_2 \psi- V\psi +\frac{4a\cos\th}{|q|^2}\dual\nab_T\psi\right)
\eeaa
as stated.

Next, we have
\beaa
\PP_\mu[Y_\HH, 0, 0]\c N_\Si &=& \QQ(Y_\HH, N_\Si)=\ka_\HH\QQ(Y_{(0)}, N_\Si)\\
&=& \ka_\HH\QQ\left(\underline{d}(r)e_3+d(r)e_4+e_4+\frac{\De}{|q|^2}e_3 -2a\Re(\Jk)^be_b, N_\Si\right)\\
&=& \ka_\HH\QQ\Big(e_3+e_4+O(r-r_+)e_3+O(r-r_+)e_4 -2a\Re(\Jk)^be_b, N_\Si\Big).
\eeaa
In view of the support of $\ka_\HH$, we infer
\beaa
\PP_\mu[Y_\HH, 0, 0]\c N_\Si &=& \ka_\HH\Bigg[\QQ\Big(e_3+e_4 -2a\Re(\Jk)^be_b, N_\Si\Big)\\
&&+O(\de_{red})\Big(|\nab_3\psi|^2+|\nab_4\psi|^2+|\nab\psi|^2+r^{-2}|\psi|^2\Big)\Bigg].
\eeaa

Next, in view of \eqref{eq:definition-QQ-mu-nu}, and since $V(r_+)=0$ in view of its explicit formula in \eqref{eq:Gen.RW}, we have on the support of $\ka_\HH$
\beaa
\bsplit
\QQ_{33}&=|\nab_3\psi|^2, \qquad  \QQ_{44}=|\nab_4\psi|^2, \qquad \QQ_{34}=|\nab\psi|^2 +O(\de_{red})|\psi|^2,\\
\QQ_{4a}&=\nab_4\Psi\c\nab_a\Psi, \qquad \QQ_{3a}=\nab_3\Psi\c\nab_a\Psi, \\ 
\tr\QQ &= \nab_4\psi\c\nab_3\psi+O(\de_{red})|\psi|^2,\qquad \widehat{\QQ}_{ab} =\frac{1}{2}(\nab\psi\hot\nab\psi)_{ab} +O(\de_{red})|\psi|^2,
\end{split}
\eeaa
and hence
\beaa
\QQ\Big(e_3+e_4 -2a\Re(\Jk)^be_b, N_\Si\Big)=\widetilde{\QQ}\Big(e_3+e_4 -2a\Re(\Jk)^be_b, N_\Si\Big)+O(\de_{red})|\psi|^2
\eeaa
where the symmetric 2-tensor $\widetilde{\QQ}$ is given by
\beaa
 \widetilde{\QQ}_{\mu\nu}:=\Db_\mu  \psi \c \Db _\nu \psi  -\frac 12 \g_{\mu\nu}\Db_\la \psi\c\Db^\la \psi.
 \eeaa
We deduce, on the support of $\ka_\HH$, 
 \beaa
\PP_\mu[Y_\HH, 0, 0]\c N_\Si &=& \ka_\HH\Bigg[\widetilde{\QQ}\Big(e_3+e_4 -2a\Re(\Jk)^be_b, N_\Si\Big)\\
&&+O(\de_{red})\Big(|\nab_3\psi|^2+|\nab_4\psi|^2+|\nab\psi|^2+r^{-2}|\psi|^2\Big)\Bigg].
\eeaa

Next, recall that the choice of $\tau$ in Definition \ref{definition:definition-oftau} is such that
\beaa
\g(N_\Si, N_\Si) \leq -\frac{m^2}{8r^2}<0.  
\eeaa
Also, we have
\beaa
\g\Big(e_3+e_4 -2a\Re(\Jk)^be_b, e_3+e_4 -2a\Re(\Jk)^be_b\Big) &=& -4\Big(1-a^2|\Re(\Jk)|^2\Big)
\eeaa
and hence,  on the support of $\ka_\HH$, 
\beaa
\g\Big(e_3+e_4 -2a\Re(\Jk)^be_b, e_3+e_4 -2a\Re(\Jk)^be_b\Big) &=& -4\left(1-\frac{a^2(\sin\th)^2}{|q|^2}+O(\ep)\right)\\
&\leq& -4\left(1-\frac{a^2}{r_+^2}+O(\de_{red}+\ep)\right)<0
\eeaa
for $|a|<m$, and $\de_{red}$ and $\ep$ small enough. Since both vectorfields are uniformly timelike on the support of $\ka_\HH$, we deduce from the above,  together with the Poincar\'e inequality of Lemma \ref{lemma:poincareinequalityfornabonSasoidfh:chap9}, the existence of a constant $c_0>0$ such that, for $\de_{red}$  small enough, for any sphere $S=S(\tau, r)$, 
 \beaa
\int_S\PP_\mu[Y_\HH, 0, 0]\c N_\Si &\geq& c_0\ka_\HH\int_S\Big(|\nab_3\psi|^2+|\nab_4\psi|^2+|\nab\psi|^2+r^{-2}|\psi|^2\Big).
\eeaa
In particular, since $\ka_\HH\geq 0$, we have on $\MM$
\beaa
\int_S\PP_\mu[Y_\HH, 0, 0]\c N_\Si &\geq& 0,
\eeaa
and moreover, in view of the definition of $\ka_\HH$, we have, for any sphere $S(\tau, r)$ with $|\frac{r}{r_+}-1|\leq \de_{red}$
\beaa
\int_S\PP_\mu[Y_\HH, 0, 0]\c N_\Si &\geq& c_0\int_S\Big(|\nab_3\psi|^2+|\nab_4\psi|^2+|\nab\psi|^2+r^{-2}|\psi|^2\Big)
\eeaa
as stated. 

Finally, we have on $\AA=\{r=(1-\deh)r_+\}$
\beaa
N_\AA &=& \g^{\a\b}\pr_\a(r)\pr_\b = -\frac{1}{2}e_3(r)e_4  -\frac{1}{2}e_4(r)e_3+\nab(r)\\
&=& \frac{1}{2}\big(1+O(\ep)\big)e_4 -\frac{1}{2}\left(\frac{\De}{|q|^2}+O(\ep)\right)e_3+O(\ep)\nab\\
&=& \frac{1}{2}\big(1+O(\ep)\big)e_4  +\frac{1}{2}\left(\frac{|\De|}{|q|^2}+O(\ep)\right)e_3+O(\ep)\nab
\eeaa
and 
\beaa
Y_{(0)} &=& \underline{d}(r)e_3+d(r)e_4+e_4+\frac{\De}{|q|^2}e_3 -2a\Re(\Jk)^be_b\\
&=& \big(1+O(\deh)\big)e_3+\big(1+O(\deh)\big)e_3 -2a\Re(\Jk)^be_b.
\eeaa
In view of the definition of $\ka_\HH$ and the fact that $\de_{red}\geq\deh$, we have $\ka_\HH=1$ on $\AA$. We infer
\beaa
\PP_\mu[Y_\HH, 0, 0]\c N_\AA &=& \QQ(Y_\HH, N_\AA)=\ka_\HH\QQ(Y_{(0)}, N_\AA)=\QQ(Y_{(0)}, N_\AA)\\
&=& \frac{1}{2}\big(1+O(\deh+\ep)\big)\QQ_{43} +\frac{1}{2}\big(1+O(\deh+\ep)\big)\QQ_{44} -a\big(1+O(\ep)\big)\Re(\Jk)^b\QQ_{4b}\\
&&+\frac{1}{2}\left(\frac{|\De|}{|q|^2}+O(\ep)\right)\Big(\big(1+O(\deh+\ep)\big)\QQ_{33}-2a\Re(\Jk)^b\QQ_{3b}\Big)\\
&&+O(\ep)\Big(\QQ_{b4}+\QQ_{b3}-2a\Re(\Jk)^c\QQ_{bc}\Big).
\eeaa
Next, in view of \eqref{eq:definition-QQ-mu-nu}, and since $V(r_+)=0$ in view of its explicit formula in \eqref{eq:Gen.RW}, we have on $\AA=\{r=(1-\deh)r_+\}$
\beaa
\bsplit
\QQ_{33}&=|\nab_3\psi|^2, \qquad  \QQ_{44}=|\nab_4\psi|^2, \qquad \QQ_{34}=|\nab\psi|^2 +O(\deh)|\psi|^2,\\
\QQ_{4a}&=\nab_4\Psi\c\nab_a\Psi, \qquad \QQ_{3a}=\nab_3\Psi\c\nab_a\Psi, \\ 
\tr\QQ &= \nab_4\psi\c\nab_3\psi+O(\deh)|\psi|^2,\qquad \widehat{\QQ}_{ab} =\frac{1}{2}(\nab\psi\hot\nab\psi)_{ab} +O(\deh)|\psi|^2.
\end{split}
\eeaa
We infer on $\AA$
\beaa
\PP_\mu[Y_\HH, 0, 0]\c N_\AA &=&  \frac{1}{2}|\nab\psi|^2 +\frac{1}{2}|\nab_4\psi|^2  -a\Re(\Jk)^b\nab_4\Psi\c\nab_b\Psi\\
&&+\frac{1}{2}\frac{|\De|}{|q|^2}\Big(|\nab_3\psi|^2 -2a\Re(\Jk)^b\nab_3\Psi\c\nab_b\Psi\Big)\\
&&+O(\deh^2+\ep)|\nab_3\psi|^2+O(\deh+\ep)\Big(|\nab_4\psi|^2+|\nab\psi|^2+r^{-2}|\psi|^2\Big).
\eeaa
Arguing as above, for $|a|<m$, and $\deh$ and $\ep$ small enough, this yields the existence of a constant $c_0>0$ depending on $a$ and $m$ such that 
\beaa
\PP_\mu[Y_\HH, 0, 0]\c N_\AA &\geq &  c_0\Big(\deh|\nab_3\psi|^2+|\nab_4\psi|^2+|\nab\psi|^2\Big)+O(\deh)|\psi|^2.
\eeaa 
Using again the above Poincar\'e inequality, we infer, for a possibly smaller $c_0>0$,
\beaa
\int_S\PP_\mu[Y_\HH, 0, 0]\c N_\AA &\geq &  c_0\int_S\Big(\deh|\nab_3\psi|^2+|\nab_4\psi|^2+|\nab\psi|^2+r^{-2}|\psi|^2\Big)
\eeaa 
for any sphere $S(\tau, r)$ with $r=(1-\deh)r_+$ as stated. This concludes the proof of Proposition \ref{prop:basiccomputationforredshiftvectorfieldnearr=r+fordivergenceandboundary}
\end{proof}

We are now ready to prove the redshift estimates of Proposition \ref{Prop:Redshift-estimates-chp3}. 

\begin{proof}[Proof of Proposition \ref{Prop:Redshift-estimates-chp3}]
Let $Y_\HH$ the vectorfield of Proposition \ref{prop:basiccomputationforredshiftvectorfieldnearr=r+fordivergenceandboundary}. We integrate $\D^\mu  \PP_\mu[Y_\HH, 0, 0]$ on $\MM(\tau_1, \tau_2)$ and apply the divergence theorem. The proof of Proposition \ref{Prop:Redshift-estimates-chp3} follows then immediately from the lower bounds for $\D^\mu  \PP_\mu[Y_\HH, 0, 0]$, $\PP_\mu[Y_\HH, 0, 0]\c N_\Si$ and $\PP_\mu[Y_\HH, 0, 0]\c N_\AA$ derived in Proposition \ref{prop:basiccomputationforredshiftvectorfieldnearr=r+fordivergenceandboundary}. 
\end{proof}

The following lemma shows that $e_3$ commute well with $\squared_k$ in the redshift region.
\begin{lemma}\lab{lemma:commutationwithe3forredshift}
We have,  for $r\leq 4m$, 
\beaa
[\nab_3, \squared_k] &=& -\pr_r\left(\frac{\De}{|q|^2}\right)\nab_3^2\psi+O(1)\nab\nab_3\psi +O(1)\nab_4 \nab_3\psi\\
&& +O(1)\squared_k \psi +O(1)\dk^{\leq 1}\psi+O(\ep)\dk^{\leq 2}\psi
\eeaa
\end{lemma}

\begin{proof}
Recall the following decomposition, see \eqref{first-equation-square}, 
\beaa
\begin{split}
\squared_k \psi&=-\nab_4 \nab_3 \psi  -\frac 1 2 \trchb \nab_4\psi+\left(2\om -\frac 1 2 \trch\right) \nab_3\psi+\lap_k \psi+2\etab \c\nab \psi \\
&+ 2i \left( \rhod- \eta \wedge \etab \right) \psi+(\Ga_b \c \Ga_g) \c \psi.
\end{split}
\eeaa
We infer, for $r\leq 4m$, 
\beaa
[\nab_3, \squared_k] &=& \Bigg[\nab_3, -\nab_4 \nab_3 \psi  -\frac 1 2 \trchb \nab_4+\left(2\om -\frac 1 2 \trch\right) \nab_3+\lap_k+2\etab \c\nab \\
&&+ 2i \left( \rhod- \eta \wedge \etab \right) \psi+(\Ga_b \c \Ga_g) \c\Bigg]\psi\\
&=& -[\nab_3, \nab_4]\nab_3\psi +[\nab_3, \De_k]\psi+O(1)\dk^{\leq 1}\psi\\
&=& 2\om\nab_3^2\psi+2(\eta-\etab)\c\nab\nab_3\psi+2\omb\nab_4\nab\psi -2\chib_{bc}\nab_c\nab_b\psi\\
&&+2(\eta-\ze)\c\nab\nab_3\psi+2\xib\c\nab_4\nab\psi+O(1)\dk^{\leq 1}\psi\\
&=& 2\om\nab_3^2\psi+2(\eta-\etab)\c\nab\nab_3\psi -\trchb\De_k\psi+O(1)\dk^{\leq 1}\psi+O(\ep)\dk^{\leq 2}\psi.
\eeaa
Using again \eqref{first-equation-square}, we infer
\beaa
[\nab_3, \squared_k] &=& 2\om\nab_3^2\psi+2(\eta-\etab)\c\nab\nab_3\psi -\trchb\nab_4 \nab_3\psi\\
&&-\trchb\squared_k \psi +O(1)\dk^{\leq 1}\psi+O(\ep)\dk^{\leq 2}\psi
\eeaa
and hence
\beaa
[\nab_3, \squared_k] &=& -\pr_r\left(\frac{\De}{|q|^2}\right)\nab_3^2\psi+O(1)\nab\nab_3\psi +O(1)\nab_4 \nab_3\psi\\
&& +O(1)\squared_k \psi +O(1)\dk^{\leq 1}\psi+O(\ep)\dk^{\leq 2}\psi
\eeaa
as stated.
\end{proof}

In view of Lemma \ref{lemma:commutationwithe3forredshift}, the following corollary of Proposition \ref{Prop:Redshift-estimates-chp3} will be useful when commuting the wave equation with $\nab_3$. 
\begin{corollary}
\lab{cor:Redshift-estimates-chp3}
Let $\psi$ a solution, in $\MM\cap\{r\leq 4m\}$, to 
\beaa
\squared_2 \psi -V\psi=- \frac{4 a\cos\th}{|q|^2}\dual \nab_T  \psi+N+\pr_r\left(\frac{\De}{|q|^2}\right)\nab_3\psi+O(1)\nab\psi +O(1)\nab_4\psi.
\eeaa
Then, for $|a|<m$, there exists a small enough constant $\de_{red}>0$ such that $\de_{red}=\de_{red}(m-|a|)$ with $\de_{red}\geq \deh$, and a small constant $c_0>0$ with $c_0=c_0(m-|a|)$, such that the following estimate  holds true in $\MM(\tau_1, \tau_2) $:
\bea
\bsplit
&c_0E_{r\leq r_+(1+\de_{red})}[\psi](\tau_2)+c_0\textrm{Mor}_{r\leq r_+(1+\de_{red})}[\psi](\tau_1, \tau_2)+c_0F_{\AA}[\psi](\tau_1, \tau_2)\\
  \leq& E_{r\leq r_+(1+2\de_{red})}[\psi](\tau_1) +\de_{red}^{-1}\textrm{Mor}_{r_+(1+\de_{red})\leq r\leq r_+(1+2\de_{red})}[\psi](\tau_1, \tau_2)\\
&+   \int_{\MM(\tau_1, \tau_2)\cap\left\{\frac{r}{r_+}\le 1+2\de_{red}\right\}}|N|^2.
\end{split}
\eea 
\end{corollary}

\begin{proof}
Let $Y_\HH$ the redshift vectorfield of Proposition \ref{prop:basiccomputationforredshiftvectorfieldnearr=r+fordivergenceandboundary}. As in the proof of that proposition, we have
\beaa
\D^\mu  \PP_\mu[Y_\HH, 0, 0]  &=& \ka_\HH \D^\mu  \PP_\mu[Y_{(0)}, 0, 0]\\
&& - O(r_+^{-1}\de_{red}^{-1})\mathbb{1}_{\de_{red}\leq\left|\frac{r}{r_+}-1\right|\leq 2\de_{red}}\Big(|\nab_3\psi|^2+|\nab_4\psi|^2+|\nab\psi|^2+r^{-2}|\psi|^2\Big).
\eeaa
Together with the lower bound for $\D^\mu  \PP_\mu[Y_{(0)}, 0, 0]$ of Corollary \ref{cor:basiccomputationforredshiftvectorfieldatr=r+} and the properties of $\ka_\HH$, we infer, provided $\de_{red}$ is chosen sufficiently small, for any sphere $S=S(\tau, r)$, choosing $d'(r_+)$ and $\underline{d}'(r_+)$ even larger than in Corollary \ref{cor:basiccomputationforredshiftvectorfieldatr=r+}, 
\beaa
\nn\int_S\D^\mu  \PP_\mu[Y_\HH, 0, 0]  &\geq& \frac{(r_+-m)}{64r_+^2}\int_S\Big(|\nab_3\psi|^2+|\nab_4\psi|^2+|\nab\psi|^2+r^{-2}|\psi|^2\Big)\mathbb{1}_{\left|\frac{r}{r_+}-1\right|\leq \de_{red}}  \\
\nn&& -O(r_+^{-1}\de_{red}^{-1})\int_S\Big(|\nab_3\psi|^2+|\nab_4\psi|^2+|\nab\psi|^2+r^{-2}|\psi|^2\Big)\mathbb{1}_{\de_{red}\leq\left|\frac{r}{r_+}-1\right|\leq 2\de_{red}}\\
&& +\int_S Y_\HH( \psi )\c \left(\squared_2 \psi- V\psi +\frac{4a\cos\th}{|q|^2}\dual\nab_T\psi\right)\\
&&+\frac{1}{2}\left(d'(r_+) - \frac{c_1}{r_+-m}\right)\int_S|\nab_4\psi|^2\mathbb{1}_{\left|\frac{r}{r_+}-1\right|\leq \de_{red}}\\
&&+\frac{1}{2}\left(\underline{d}'(r_+) - \frac{c_1}{r_+-m}\right)\int_S|\nab\psi|^2\mathbb{1}_{\left|\frac{r}{r_+}-1\right|\leq \de_{red}}.
\eeaa
Plugging the equation for $\psi$, we infer
\beaa
\nn\int_S\D^\mu  \PP_\mu[Y_\HH, 0, 0]  &\geq& \frac{(r_+-m)}{64r_+^2}\int_S\Big(|\nab_3\psi|^2+|\nab_4\psi|^2+|\nab\psi|^2+r^{-2}|\psi|^2\Big)\mathbb{1}_{\left|\frac{r}{r_+}-1\right|\leq \de_{red}}  \\
\nn&& -O(r_+^{-1}\de_{red}^{-1})\int_S\Big(|\nab_3\psi|^2+|\nab_4\psi|^2+|\nab\psi|^2+r^{-2}|\psi|^2\Big)\mathbb{1}_{\de_{red}\leq\left|\frac{r}{r_+}-1\right|\leq 2\de_{red}}\\
&& +\int_S Y_\HH( \psi )\c \left(N+\pr_r\left(\frac{\De}{|q|^2}\right)\nab_3\psi+O(1)\nab\psi +O(1)\nab_4\psi\right)\\
&&+\frac{1}{2}\left(d'(r_+) - \frac{c_1}{r_+-m}\right)\int_S|\nab_4\psi|^2\mathbb{1}_{\left|\frac{r}{r_+}-1\right|\leq \de_{red}}\\
&&+\frac{1}{2}\left(\underline{d}'(r_+) - \frac{c_1}{r_+-m}\right)\int_S|\nab\psi|^2\mathbb{1}_{\left|\frac{r}{r_+}-1\right|\leq \de_{red}}.
\eeaa
Now, in view of the choice of $Y_\HH$, we have
\beaa
&&\int_S Y_\HH( \psi )\c \left(\pr_r\left(\frac{\De}{|q|^2}\right)\nab_3\psi+O(1)\nab\psi +O(1)\nab_4\psi\right)\\
&=& \int_S \ka_\HH \Big(\big(1+O(\de_{red})\big)\nab_3+O(1)\nab_4+O(1)\nab\Big)\psi\c \left(\pr_r\left(\frac{\De}{|q|^2}\right)\nab_3\psi+O(1)\nab\psi +O(1)\nab_4\psi\right). 
\eeaa
Since $\pr_r\left(\frac{\De}{|q|^2}\right)\geq 0$ on the support of $\ka_\HH$, we infer
\beaa
&&\int_S Y_\HH( \psi )\c \left(\pr_r\left(\frac{\De}{|q|^2}\right)\nab_3\psi+O(1)\nab\psi +O(1)\nab_4\psi\right)\\
&\geq & -O(1)\int_S \ka_\HH\Big(|\nab_4\psi|^2+|\nab\psi|^2+|\nab_3\psi|\big(|\nab_4\psi|+|\nab\psi|\big)\Big) 
\eeaa
and hence, for $d'(r_+)$ and $\underline{d}'(r_+)$ large enough, 
\beaa
\nn\int_S\D^\mu  \PP_\mu[Y_\HH, 0, 0]  &\geq& \frac{(r_+-m)}{128r_+^2}\int_S\Big(|\nab_3\psi|^2+|\nab_4\psi|^2+|\nab\psi|^2+r^{-2}|\psi|^2\Big)\mathbb{1}_{\left|\frac{r}{r_+}-1\right|\leq \de_{red}}  \\
\nn&& -O(r_+^{-1}\de_{red}^{-1})\int_S\Big(|\nab_3\psi|^2+|\nab_4\psi|^2+|\nab\psi|^2+r^{-2}|\psi|^2\Big)\mathbb{1}_{\de_{red}\leq\left|\frac{r}{r_+}-1\right|\leq 2\de_{red}}\\
&& -O(1)\int_S|N|^2\mathbb{1}_{\left|\frac{r}{r_+}-1\right|\leq 2\de_{red}}.
\eeaa
We integrate $\D^\mu  \PP_\mu[Y_\HH, 0, 0]$ on $\MM(\tau_1, \tau_2)$ and apply the divergence theorem. The proof of Corollary \ref{cor:Redshift-estimates-chp3} follows then immediately from the above lower bounds for $\D^\mu  \PP_\mu[Y_\HH, 0, 0]$, and from the lower bound for $\PP_\mu[Y_\HH, 0, 0]\c N_\Si$ and $\PP_\mu[Y_\HH, 0, 0]\c N_\AA$ derived in Proposition \ref{prop:basiccomputationforredshiftvectorfieldnearr=r+fordivergenceandboundary}. 
\end{proof}


\section{Proof of Theorem \ref{THM:HIGHERDERIVS-MORAWETZ-CHP3}}
\lab{sec:proofofThm:HigherDerivs-Morawetz-chp3:chap9}


In this section, we finally prove our main result for Energy-Morawetz estimates in perturbations of Kerr. 
We first start with the particular case $s=2$, and then treat the general case.


\subsection{Proof of Theorem \ref{THM:HIGHERDERIVS-MORAWETZ-CHP3} in the case $s=2$}
\lab{sec:proofofThm:HigherDerivs-Morawetz-chp3:cases=2}


Commutating the model RW equation \eqref{eq:Gen.RW} satisfied by $\psi$ with various vectorfields will generate in particular error terms of the type $\dk^{\leq 1}(\Ga_g\c\psi)$. We start with the following simple lemma which will allow us to control the contribution of these error terms.

\begin{lemma}\lab{lemma:controloftheerrortermsGagpsiinwaveeq}
We have
\beaa
\int_{\MM_{trap}(\tau_1, \tau_2)}|\dk^{\leq s+1}(\Ga_g\c\psi)||\dk^{\leq s+1}\psi| \les \ep\sup_{\tau\in [\tau_1, \tau_2]}E^s[\psi]
\eeaa
and
\beaa
\int_{\Mntrap(\tau_1, \tau_2)}|\dk^{\leq s+1}(\Ga_g\c\psi)|\Big(|\nab_{\That_\de}\dk^{\leq s}\psi|+|\nab_{\Rhat}\dk^{\leq s}\psi|+r^{-1}|\dk^{\leq s}\psi|\Big) \les \ep B^s_\de[\psi](\tau_1, \tau_2).
\eeaa
\end{lemma}

\begin{proof}
We start with the control on $\MM_{trap}$. We have
\beaa
&&\int_{\MM_{trap}(\tau_1, \tau_2)}|\dk^{\leq s+1}(\Ga_g\c\psi)||\nab_{\That_\de}\dk^{\leq s}\psi|\\
&\les&  \ep\int_{\MM_{trap}(\tau_1, \tau_2)}\frac{1}{\tau_{trap}^{1+\dec}}|\dk^{\leq s+1}\psi|^2\\
&\les& \ep\left(\int_{\tau_1}^{\tau_2}\frac{1}{\tau_{trap}^{1+\dec}}\right)E^s[\psi]\\
&\les& \ep\sup_{\tau\in [\tau_1, \tau_2]}E^s[\psi]
\eeaa 
as stated. 

Concerning the control on $\Mntrap$, we have
\beaa
&&\int_{\Mntrap(\tau_1, \tau_2)}|\dk^{\leq s+1}(\Ga_g\c\psi)|\Big(|\nab_{\That_\de}\dk^{\leq s}\psi|+|\nab_{\Rhat}\dk^{\leq s}\psi|+r^{-1}|\dk^{\leq s}\psi|\Big)\\
&\les& \ep\int_{\Mntrap(\tau_1, \tau_2)}r^{-2}|\dk^{\leq s+1}\psi|\Big(|\nab_{\That_\de}\dk^{\leq s}\psi|+|\nab_{\Rhat}\dk^{\leq s}\psi|+r^{-1}|\dk^{\leq s}\psi|\Big)\\
&\les& \ep\int_{\Mntrap(\tau_1, \tau_2)}r^{-2}|\dk^{\leq s+1}\psi|\Big(|\nab_3\dk^{\leq s}\psi|+r^{-1}|\dk^{\leq s+1}\psi|\Big)\\
&\les& \ep\int_{\Mntrap(\tau_1, \tau_2)}\Big(r^{-1-\de}|\nab_3\dk^{\leq s}\psi|^2+r^{\de-3}|\dk^{\leq s+1}\psi|^2\Big)\\
&\les& \ep B^s_\de[\psi](\tau_1, \tau_2)
\eeaa
as stated. This concludes the proof of the lemma.
\end{proof}

First, we derive the control of energy for at most one derivative of $\psi$ and Morawetz for at most one $(\nab_\T, \nab_\Z)$ derivative of $\psi$. 
\begin{proposition}\lab{prop:energyforzeroandfirstderivativeinchap9}
The solution $\psi$ of of the model RW equation \eqref{eq:Gen.RW} satisfies the following energy estimate 
\bea
\nn&& E_{deg}[(\dkb, \nab_T, \nab_{\Rhat})^{\leq 1}\psi](\tau_2)+F_{\Si_*}[(\dkb, \nab_T, \nab_{\Rhat})^{\leq 1}\psi](\tau_1, \tau_2)+\textrm{Mor}[(\nab_\T, \nab_\Z)^{\leq 1}\psi](\tau_1, \tau_2)\\
\nn&\les& \deh\Big(E_{r\leq r_+}[\dkb\psi](\tau_2)+ F_{\AA}[\dkb\psi](\tau_1, \tau_2)\Big)+E^1[\psi](\tau_1)+\NN^1[\psi, N] (\tau_1, \tau_2)\\
&&+\left(\frac{|a|}{m}+\ep\right)\left(\sup_{\tau\in [\tau_1, \tau_2]}E^2[\psi]+B_\de^2[\psi](\tau_1, \tau_2)+F^2[\psi](\tau_1, \tau_2)\right).
\eea
\end{proposition}

\begin{proof}
We proceed in several steps.

{\bf Step 1.} From Proposition \ref{proposition:Morawetz1-step1:perturbation}, we have 
\beaa
\Mor^{ax}_{deg}[\psi](\tau_1, \tau_2) &\les& \sup_{\tau\in[\tau_1, \tau_2]}E_{deg}[\psi]+   \int_{\MM(\tau_1, \tau_2)}ar^{-2}\big(|\nab\psi|^2+|\nab_{\T}\psi|^2\big)\\
&&+\int_{\MM(\tau_1, \tau_2)}(|\nab_{\Rhat}\psi|+r^{-1}|\psi|)|N|+\ep\left(\sup_{[\tau_1, \tau_2]}E[\psi](\tau)+ B^1_\de[\psi]\right).
\eeaa
Also, from Proposition \ref{proposition:Energy1:perturbation}, we have
 \beaa
&& E_{deg}[\psi](\tau_2)+F_{\Si_*}[\psi](\tau_1, \tau_2) \\ 
&\les&  \de_{\HH}\left(E_{r\leq r_+(1+\deh)}(\tau_2)+F_{\AA}[\psi](\tau_1, \tau_2)\right)+ E_{deg}[\psi](\tau_1)+ \frac{|a|}{m}\Mor^{ax}_{deg}[\psi](\tau_1, \tau_2) \\
&&  +\left|\int_{\MM}  \nab_{\That_\de } \psi  \c N\right|+\ep\left(\sup_{[\tau_1, \tau_2]}E[\psi](\tau)+ B_\de[\psi]\right).
\eeaa
Also, we have the redshift estimate of Proposition \ref{Prop:Redshift-estimates-chp3}
\beaa
&& \textrm{Mor}_{r\leq r_+(1+\de_{red})}[\psi](\tau_1, \tau_2)+E_{r\leq r_+(1+\de_{red})}[\psi](\tau_2)+F_{\AA}[\psi](\tau_1, \tau_2)\\
  &\les& E_{r\leq r_+(1+2\de_{red})}[\psi](\tau_1) +\de_{red}^{-3}\textrm{Mor}_{deg}[\psi](\tau_1, \tau_2)+   \int_{\MM(\tau_1, \tau_2)\cap\left\{\frac{r}{r_+}\le 1+2\de_{red}\right\}}|N|^2.
\eeaa 
Combining the three estimates, where the second is multiplied by the large constant $\La$ and the third by $\de_{red}^4$, we infer, for a constant $c_0$ only depending on $m-|a|$, 
\beaa
&& c_0\Big(\de_{red}^4\big(E_{red}[\psi](\tau_2)+\textrm{Mor}_{red}[\psi](\tau_1, \tau_2)+F_{\AA}[\psi](\tau_1, \tau_2)\big)+\Lambda \big(E_{deg}(\tau_2)+F_{\Si_*}[\psi](\tau_1, \tau_2)\big)\\
&&+\Mor^{ax}_{deg}[\psi](\tau_1, \tau_2)\Big)\\
&\leq& \sup_{\tau\in[\tau_1, \tau_2]}E_{deg}[\psi]+   \int_{\MM(\tau_1, \tau_2)}ar^{-2}\big(|\nab\psi|^2+|\nab_{\T}\psi|^2\big)+\int_{\MM(\tau_1, \tau_2)}(|\nab_{\Rhat}\psi|+r^{-1}|\psi|)|N|\\
&&+\int_{\MM(\tau_1, \tau_2)}|\dk^{\leq 1}(\Ga_g\c\psi)||\dk^{\leq 1}\psi|\\
&&+\La\Bigg( \de_{\HH}\left(E_{r\leq r_+(1+\deh)}(\tau_2)+F_{\AA}[\psi](\tau_1, \tau_2)\right)+ E_{deg}[\psi](\tau_1)+ \frac{|a|}{m}\Mor^{ax}_{deg}[\psi](\tau_1, \tau_2) \\
&&  +\left|\int_{\MM}  \nab_{\That_\de } \psi  \c N\right|+\ep\left(\sup_{[\tau_1, \tau_2]}E[\psi](\tau)+ B^1_\de[\psi]\right)\Bigg)\\
&& +\de_{red}^4E_{r\leq r_+(1+2\de_{red})}[\psi](\tau_1) +\de_{red}\textrm{Mor}_{deg}[\psi](\tau_1, \tau_2)+   \de_{red}^4\int_{\MM(\tau_1, \tau_2)\cap\left\{\frac{r}{r_+}\le 1+2\de_{red}\right\}}|N|^2.
\eeaa
We now choose $\La$ and $\de_{red}$, depending on $m-|a|$, such that 
\beaa
c_0\La\geq 2, \qquad  c_0\de_{red}^4\geq 2\La\de_\HH, \qquad c_0\geq 2\de_{red}
\eeaa
which is possible provided $\de_\HH$ is sufficiently small. We infer
\beaa
&& \frac{c_0}{2}\Big(\de_{red}^4\big(E_{red}[\psi](\tau_2)+\textrm{Mor}_{red}[\psi](\tau_1, \tau_2)+F_{\AA}[\psi](\tau_1, \tau_2)\big)+\Lambda \big(E_{deg}(\tau_2)+F_{\Si_*}[\psi](\tau_1, \tau_2)\big)\\
&&+\Mor^{ax}_{deg}[\psi](\tau_1, \tau_2)\Big)\\
&\leq&  \int_{\MM(\tau_1, \tau_2)}ar^{-2}\big(|\nab\psi|^2+|\nab_{\T}\psi|^2\big)+\int_{\MM(\tau_1, \tau_2)}(|\nab_{\Rhat}\psi|+r^{-1}|\psi|)|N|+\int_{\MM(\tau_1, \tau_2)}|\dk^{\leq 1}(\Ga_g\c\psi)||\dk^{\leq 1}\psi|\\
&&+\La\Bigg(  E_{deg}[\psi](\tau_1)+ \frac{|a|}{m}\Mor^{ax}_{deg}[\psi](\tau_1, \tau_2)   +\left|\int_{\MM}  \nab_{\That_\de } \psi  \c N\right|+\ep\left(\sup_{[\tau_1, \tau_2]}E[\psi](\tau)+ B^1_\de[\psi]\right)\Bigg)\\
&& +\de_{red}^4E_{r\leq r_+(1+2\de_{red})}[\psi](\tau_1) +   \de_{red}^4\int_{\MM(\tau_1, \tau_2)\cap\left\{\frac{r}{r_+}\le 1+2\de_{red}\right\}}|N|^2.
\eeaa
We deduce
\beaa
&& E[\psi](\tau_2)+F[\psi](\tau_1, \tau_2)+\textrm{Mor}[\psi](\tau_1, \tau_2)\\
&\les& E[\psi](\tau_1)+\NN[\psi, N] (\tau_1, \tau_2)+\left(\frac{|a|}{m}+\ep\right)\left(\sup_{\tau\in [\tau_1, \tau_2]}E^1[\psi]+B_\de^1[\psi](\tau_1, \tau_2)+F^1[\psi](\tau_1, \tau_2)\right).
\eeaa

{\bf Step 2.} We commute with $(\Lied_\T, \Lied_\Z)$ the RW model equation satisfied by $\psi$ and obtain
\beaa
\squared_2(\Lieb_\T\psi) &=& -[\Lied_\T, \squared_2]\psi+\Lieb_\T\left(V\psi -\frac{4a\cos\th}{|q|^2}\nab_\T\dual\psi+N\right)\\
&=& -[\Lied_\T, \squared_2]\psi+V\Lieb_\T\psi -\frac{4a\cos\th}{|q|^2}\nab_\T\dual\Lieb_\T\psi+\dk^{\leq 1}N +\Ga_g\dk^{\leq 1}\psi
\eeaa
and
\beaa
\squared_2(\Lieb_\Z\psi) &=& -[\Lied_\Z, \squared_2]\psi+\Lieb_\Z\left(V\psi -\frac{4a\cos\th}{|q|^2}\nab_\T\dual\psi+N\right)\\
&=& -[\Lied_\Z, \squared_2]\psi+V\Lieb_\Z\psi -\frac{4a\cos\th}{|q|^2}\nab_\T\dual\Lieb_\Z\psi+\dk^{\leq 1}N +\Ga_g\dk^{\leq 1}\psi.
\eeaa
Recalling Corollary \ref{cor:commutator-Lied-squared}, i.e. 
\beaa
\, [\Lied_\T, \squared_2]\psi &=&  \dk \big(\Ga_g \c \dk \psi\big)+\Ga_b \c \square_\g\psi,\\
\, [\Lied_\Z, \squared_2]\psi &=& \dk \big(\Ga_g \c \dk \psi\big)+r\Ga_b \c \square_\g\psi,
\eeaa
we infer, using again the RW model equation for $\psi$ to simplify the RHS, 
\beaa
\squared_2(\Lieb_\T\psi) - V\Lieb_\T\psi &=&  -\frac{4a\cos\th}{|q|^2}\dual\nab_\T\Lieb_\T\psi +N_{\Lieb_\T\psi},\\
N_{\Lieb_\T\psi} &=& \dk^{\leq 1}N +\dk^{\leq 1}(\Ga_g\dk^{\leq 1}\psi)  +\Ga_b\squared_2\psi\\
&=&\dk^{\leq 1}N +\dk^{\leq 1}(\Ga_g\dk^{\leq 1}\psi)
\eeaa
and
\beaa
\squared_2(\Lieb_\Z\psi) - V\Lieb_\Z\psi &=&  -\frac{4a\cos\th}{|q|^2}\dual\nab_\T\Lieb_\Z\psi +N_{\Lieb_\Z\psi},\\
N_{\Lieb_\Z\psi} &=& \dk^{\leq 1}N +\dk^{\leq 1}(\Ga_g\dk^{\leq 1}\psi)  +r\Ga_b\squared_2\psi\\
&=&\dk^{\leq 1}N +\dk^{\leq 1}(\Ga_g\dk^{\leq 1}\psi).
\eeaa
We may apply the control derived in Step 1 to these RW model equations which yields, using Lemma \ref{lemma:controloftheerrortermsGagpsiinwaveeq} to control the $\dk^{\leq 1}(\Ga_g\dk^{\leq 1}\psi)$ error terms, 
\beaa
&& E[(\Lieb_\T, \Lieb_\Z)\psi](\tau_2)+F[(\Lieb_\T, \Lieb_\Z)\psi](\tau_1, \tau_2)+\textrm{Mor}[(\Lieb_\T, \Lieb_\Z)\psi](\tau_1, \tau_2)\\
&\les& E^1[\psi](\tau_1)+\NN^1[\psi, N] (\tau_1, \tau_2)+\left(\frac{|a|}{m}+\ep\right)\left(\sup_{\tau\in [\tau_1, \tau_2]}E^2[\psi]+B_\de^2[\psi](\tau_1, \tau_2)+F^2[\psi](\tau_1, \tau_2)\right).
\eeaa
Given the relation between $\Lieb_\T$ and $\nab_\T$, and the one between $\Lieb_\Z$ and $\nab_\Z$, see Lemma \ref{lemma:basicpropertiesLiebTfasdiuhakdisug:chap9}, we infer, together with the estimate of Step 1, 
\beaa
&& E[(\nab_\T, \nab_\Z)^{\leq 1}\psi](\tau_2)+F[(\nab_\T, \nab_\Z)^{\leq 1}\psi](\tau_1, \tau_2)+\textrm{Mor}[(\nab_\T, \nab_\Z)^{\leq 1}\psi](\tau_1, \tau_2)\\
&\les&  E^1[\psi](\tau_1)+\NN^1[\psi, N] (\tau_1, \tau_2)+\left(\frac{|a|}{m}+\ep\right)\left(\sup_{\tau\in [\tau_1, \tau_2]}E^2[\psi]+B_\de^2[\psi](\tau_1, \tau_2)+F^2[\psi](\tau_1, \tau_2)\right).
\eeaa 

{\bf Step 3.} Next, we  commute the wave equation \eqref{eq:Gen.RW} satisfied by $\psi$ by $|q|\DD_2$ using the commutation formula of Lemma \ref{LEMMA:COMMUTATIONOFHODGEELLIPTICORDER1WITHSQAURED2FDILUHS} 
\beaa
\squared_1(|q|\DD_2\psi) -V|q|\DD_2\psi &=& - \frac{4 a\cos\th}{|q|^2}\dual \nab_T|q|\DD_2\psi+N_{|q|\DD_2},\\
N_{\dkb_2} &:=& \dk^{\leq 1}N -\frac{3}{r^2}|q|\DD_2\psi +O(ar^{-2})\dk^{\leq 2}\psi+\dk^{\leq 2}(\Ga_g\c\psi).
\eeaa
From Proposition \ref{proposition:Energy1:perturbation}, we infer
\beaa
&&E_{deg}[|q|\DD_2\psi](\tau_2) +F_{\Si_*}[|q|\DD_2\psi](\tau_1, \tau_2)\\
&\les& \deh\Big( E_{r\leq r_+}[\dkb\psi](\tau_2)+ F_{\AA}[\dkb\psi](\tau_1, \tau_2)\Big)+E^1[\psi](\tau_1) +\left|\int_{\MM(\tau_1, \tau_2)}  \nab_{\That_\de }(|q|\DD_2\psi)  \c N_{|q|\DD_2}\right|\\
&&+\left(\frac{|a|}{m}+\ep\right)\left(\sup_{\tau\in [\tau_1, \tau_2]}E^2[\psi]+B_\de^2[\psi](\tau_1, \tau_2)+F^2[\psi](\tau_1, \tau_2)\right).
\eeaa
Now, from the definition of $N_{|q|\DD_2}$, we have, integrating by parts  the second term in $N_{|q|\DD_2}$,  
\beaa
&&\left|\int_{\MM(\tau_1, \tau_2)}  \nab_{\That_\de }(|q|\DD_2\psi)  \c N_{|q|\DD_2}\right|\\
 &\les& \left|\int_{\MM(\tau_1, \tau_2)} \frac{1}{r^2} \nab_{\That_\de }(|q|\DD_2\psi)  \c |q|\DD_2\psi\right|+\NN^1[\psi, N] (\tau_1, \tau_2)\\
 &&+\left(\frac{|a|}{m}+\ep\right)\left(\sup_{\tau\in [\tau_1, \tau_2]}E^2[\psi]+B_\de^2[\psi](\tau_1, \tau_2)+F^2[\psi](\tau_1, \tau_2)\right)\\
&\les&  \sup_{\tau\in[\tau_1, \tau_2]}E[\psi]+\NN^1[\psi, N] (\tau_1, \tau_2)\\
 &&+\left(\frac{|a|}{m}+\ep\right)\left(\sup_{\tau\in [\tau_1, \tau_2]}E^2[\psi]+B_\de^2[\psi](\tau_1, \tau_2)+F^2[\psi](\tau_1, \tau_2)\right),
\eeaa
where we used the fact that $\That_\de(r)\in r\Ga_b$, the control of $\D_\mu \That^\mu$ induced from the one of ${}^{(\That_\de)}\pi$ in the proof of Proposition \ref{proposition:Energy1:perturbation} in section \ref{sec:proofofresultschapter7inperturbationofKerr}, and Lemma \ref{lemma:controloftheerrortermsGagpsiinwaveeq} to control the various error terms. 

In view of the above, and in particular the estimate for $E[\psi]$ of Step 1, we infer
\beaa
&&E_{deg}[|q|\DD_2\psi](\tau_2) +F_{\Si_*}[|q|\DD_2\psi](\tau_1, \tau_2)\\
&\les& \deh\Big( E_{r\leq r_+}[\dkb\psi](\tau_2)+ F_{\AA}[\dkb\psi](\tau_1, \tau_2)\Big)+E^1[\psi](\tau_1)\\
&&+\NN^1[\psi, N] (\tau_1, \tau_2)+\left(\frac{|a|}{m}+\ep\right)\left(\sup_{\tau\in [\tau_1, \tau_2]}E^2[\psi]+B_\de^2[\psi](\tau_1, \tau_2)+F^2[\psi](\tau_1, \tau_2)\right).
\eeaa
Together with Step 2, and using the Hodge estimates of Proposition \ref{Prop:HodgeThmM8} to control $\dkb$ from $|q|\DD_2$ and $\nab_\T$, we deduce 
\beaa
&& E_{deg}[(\dkb, \nab_T)^{\leq 1}\psi](\tau_2)+F_{\Si_*}[(\dkb, \nab_T)^{\leq 1}\psi](\tau_1, \tau_2) +\textrm{Mor}[(\nab_\T, \nab_\Z)^{\leq 1}\psi](\tau_1, \tau_2)\\
&\les& \deh\Big(E_{r\leq r_+}[\dkb\psi](\tau_2)+ F_{\AA}[\dkb\psi](\tau_1, \tau_2)\Big)+E^1[\psi](\tau_1)\\
&&+\NN^1[\psi, N] (\tau_1, \tau_2)+\left(\frac{|a|}{m}+\ep\right)\left(\sup_{\tau\in [\tau_1, \tau_2]}E^2[\psi]+B_\de^2[\psi](\tau_1, \tau_2)+F^2[\psi](\tau_1, \tau_2)\right).
\eeaa 

{\bf Step 4.} Using   the representation of the wave operator provided by \eqref{eq:q2squaredkpsi}, i.e. 
\beaa
\begin{split}
|q|^2 \squared_2 \psi &=\frac{(r^2+a^2)^2}{\De} \big( -  \nab_\That \nab_\That \psi+   \nab_\Rhat \nab_\Rhat \psi \big) +2r \nab_\Rhat \psi\\
&+  |q|^2 \lap_2 \psi   + |q|^2  (\eta+\etab) \c \nab \psi  + r^2 \Ga_g \c \dk \psi,
\end{split}
\eeaa
together with the wave equation \eqref{eq:Gen.RW} satisfied by $\psi$,  we have
\beaa
\int_{\Si(\tau_2)}|\nab_{\Rhat}^2\psi|^2+\int_{\Si_*(\tau_1, \tau_2)}|\nab_{\Rhat}^2\psi|^2 &\les& E_{deg}[(\nab_T, \nab_Z)^{\leq 1}\psi](\tau)+F_{\Si_*}[(\nab_T, \nab_Z)^{\leq 1}\psi]\\
&&+\int_{\Si(\tau_2)}|N|^2+\int_{\Si_*(\tau_1, \tau_2)}|N|^2.
\eeaa
Since we have 
\beaa
&& E_{deg}[(\dkb, \nab_T, \nab_{\Rhat})^{\leq 1}\psi](\tau_2)+F_{\Si_*}[(\dkb, \nab_T, \nab_{\Rhat})^{\leq 1}\psi](\tau_1, \tau_2)\\
&\les& E_{deg}[(\dkb, \nab_T)^{\leq 1}\psi](\tau_2)+F_{\Si_*}[(\dkb, \nab_T)^{\leq 1}\psi](\tau_1, \tau_2)+\int_{\Si(\tau_2)}|\nab_{\Rhat}^2\psi|^2+\int_{\Si_*(\tau_1, \tau_2)}|\nab_{\Rhat}^2\psi|^2,
\eeaa
we infer
\beaa
&& E_{deg}[(\dkb, \nab_T, \nab_{\Rhat})^{\leq 1}\psi](\tau_2)+F_{\Si_*}[(\dkb, \nab_T, \nab_{\Rhat})^{\leq 1}\psi](\tau_1, \tau_2)\\
&\les& E_{deg}[(\dkb, \nab_T)^{\leq 1}\psi](\tau_2)+F_{\Si_*}[(\dkb, \nab_T)^{\leq 1}\psi](\tau_1, \tau_2)+\int_{\Si(\tau_2)}|N|^2+\int_{\Si_*(\tau_1, \tau_2)}|N|^2.
\eeaa
Together with Step 3, we deduce 
\beaa
&& E_{deg}[(\dkb, \nab_T, \nab_{\Rhat})^{\leq 1}\psi](\tau_2)+F_{\Si_*}[(\dkb, \nab_T, \nab_{\Rhat})^{\leq 1}\psi](\tau_1, \tau_2) +\textrm{Mor}[(\nab_\T, \nab_\Z)^{\leq 1}\psi](\tau_1, \tau_2)\\
&\les& \deh\Big(E_{r\leq r_+}[\dkb\psi](\tau_2)+ F_{\AA}[\dkb\psi](\tau_1, \tau_2)\Big)+E^1[\psi](\tau_1)\\
&&+\NN^1[\psi, N] (\tau_1, \tau_2)+\left(\frac{|a|}{m}+\ep\right)\left(\sup_{\tau\in [\tau_1, \tau_2]}E^2[\psi]+B_\de^2[\psi](\tau_1, \tau_2)+F^2[\psi](\tau_1, \tau_2)\right)
\eeaa 
as stated. This concludes the proof of Proposition \ref{prop:energyforzeroandfirstderivativeinchap9}.
\end{proof}

Next, we control the energy of $\psi_\aund=\SS_\aund\psi$ for $\aund=1,2,3,4$. 
\begin{lemma}\lab{lemma:contorloftheenergyofpsiaund}
We have
\bea
\nn&& \sum_{\aund=1}^4\Big(E_{deg}[\psi_\aund](\tau_2)+F_{\Si_*}[\psi_\aund](\tau_1, \tau_2)\Big)\\
\nn&\les& \deh\Big(E^2_{r\leq r_+}[\psi](\tau_2)+ F^2_{\AA}[\psi](\tau_1, \tau_2)\Big)+E^2[\psi](\tau_1)+\NN^2[\psi, N] (\tau_1, \tau_2)\\
&&+\left(\frac{|a|}{m}+\ep\right)\left(\sup_{\tau\in [\tau_1, \tau_2]}E^2[\psi]+B_\de^2[\psi](\tau_1, \tau_2)+F^2[\psi](\tau_1, \tau_2)\right).
\eea 
\end{lemma}

\begin{proof}
We first consider $\psi_\aund$ for $\aund=1,2,3$. We commute with $(\Lied_\T, \Lied_\Z)^2$ the RW model equation satisfied by $\psi$ and obtain, proceeding as in Step 2 of the proof of Proposition \ref{prop:energyforzeroandfirstderivativeinchap9},
\beaa
\squared_2((\Lieb_\T, \Lieb_\Z)^2\psi) - V(\Lieb_\T, \Lieb_\Z)^2\psi &=&  -\frac{4a\cos\th}{|q|^2}\dual\nab_\T(\Lieb_\T, \Lieb_\Z)^2\psi +N_{(\Lieb_\T, \Lieb_\Z)^2},\\
N_{(\Lieb_\T, \Lieb_\Z)^2} &=& \dk^{\leq 2}N +\dk^{\leq 2}(\Ga_g\dk^{\leq 1}\psi).
\eeaa
From Proposition \ref{proposition:Energy1:perturbation}, we infer
\beaa
&&E_{deg}[(\Lieb_\T, \Lieb_\Z)^2\psi](\tau_2) +F_{\Si_*}[(\Lieb_\T, \Lieb_\Z)^2\psi](\tau_1, \tau_2)\\
&\les& \deh\Big( E^2_{r\leq r_+}[\psi](\tau_2)+ F^2_{\AA}[\psi](\tau_1, \tau_2)\Big)+E^2[\psi](\tau_1)\\
&& +\left|\int_{\MM(\tau_1, \tau_2)}  \nab_{\That_\de }((\Lieb_\T, \Lieb_\Z)^2\psi)  \c N_{(\Lieb_\T, \Lieb_\Z)^2}\right|\\
&&+\left(\frac{|a|}{m}+\ep\right)\left(\sup_{\tau\in [\tau_1, \tau_2]}E^2[\psi]+B_\de^2[\psi](\tau_1, \tau_2)+F^2[\psi](\tau_1, \tau_2)\right).
\eeaa
Together with Lemma \ref{lemma:controloftheerrortermsGagpsiinwaveeq} and the structure of  $N_{(\Lieb_\T, \Lieb_\Z)^2}$, we deduce 
\beaa
&&E_{deg}[(\Lieb_\T, \Lieb_\Z)^2\psi](\tau_2) +F_{\Si_*}[(\Lieb_\T, \Lieb_\Z)^2\psi](\tau_1, \tau_2)\\
&\les& \deh\Big( E^2_{r\leq r_+}[\psi](\tau_2)+ F^2_{\AA}[\psi](\tau_1, \tau_2)\Big)+E^2[\psi](\tau_1) +\NN^2[\psi, N] (\tau_1, \tau_2)\\
&&+\left(\frac{|a|}{m}+\ep\right)\left(\sup_{\tau\in [\tau_1, \tau_2]}E^2[\psi]+B_\de^2[\psi](\tau_1, \tau_2)+F^2[\psi](\tau_1, \tau_2)\right).
\eeaa
Given the relation between $\Lieb_\T$ and $\nab_\T$, and the one between $\Lieb_\Z$ and $\nab_\Z$, see Lemma \ref{lemma:basicpropertiesLiebTfasdiuhakdisug:chap9}, and using the control of $E_{deg}[(\nab_\T, \nab_\Z)^{\leq 1}\psi]$ and $F_{\Si_*}[(\nab_\T, \nab_\Z)^{\leq 1}\psi]$ provided by Proposition \ref{prop:energyforzeroandfirstderivativeinchap9}, we obtain 
\beaa
&&E_{deg}[(\nab_\T, \nab_\Z)^2\psi](\tau_2) +F_{\Si_*}[(\nab_\T, \nab_\Z)^2\psi](\tau_1, \tau_2)\\
&\les& \deh\Big( E^2_{r\leq r_+}[\psi](\tau_2)+ F^2_{\AA}[\psi](\tau_1, \tau_2)\Big)+E^2[\psi](\tau_1) +\NN^2[\psi, N] (\tau_1, \tau_2)\\
&&+\left(\frac{|a|}{m}+\ep\right)\left(\sup_{\tau\in [\tau_1, \tau_2]}E^2[\psi]+B_\de^2[\psi](\tau_1, \tau_2)+F^2[\psi](\tau_1, \tau_2)\right),
\eeaa
and hence, in view of the definition of $\psi_\aund$ for $\aund=1,2,3$, 
\beaa
\nn&& \sum_{\aund=1}^3\Big(E_{deg}[\psi_\aund](\tau_2)+F_{\Si_*}[\psi_\aund](\tau_1, \tau_2)\Big)\\
\nn&\les& \deh\Big(E^2_{r\leq r_+}[\psi](\tau_2)+ F^2_{\AA}[\psi](\tau_1, \tau_2)\Big)+E^2[\psi](\tau_1)+\NN^2[\psi, N] (\tau_1, \tau_2)\\
&&+\left(\frac{|a|}{m}+\ep\right)\left(\sup_{\tau\in [\tau_1, \tau_2]}E^2[\psi]+B_\de^2[\psi](\tau_1, \tau_2)+F^2[\psi](\tau_1, \tau_2)\right).
\eeaa 

It remains to control the energy of $\psi_4=\OO\psi$. To this end, we rely on the modified operator $\widetilde{\OO}$ defined by \eqref{eq:defintionoftheoperqatorwidetildeOOcommutingwellwtihRWmodel}, i.e.
\beaa
\widetilde{\OO}\psi = \OO\psi  +\frac{4a(r^2+a^2+|q|^2)\cos\th}{|q|^2}\nab_\T\dual\psi +\frac{4a^2\cos\th}{|q|^2}\nab_\Z\dual\psi
\eeaa
and derive the model RW equation satisfied by $\widetilde{\OO}\psi$. In view of the model RW equation for $\psi$, we have
\beaa
\squared_2(\widetilde{\OO}\psi) - V\widetilde{\OO}\psi &=& -\frac{4a\cos\th}{|q|^2}\dual\nab_\T\widetilde{\OO}\psi +\frac{1}{|q|^2}\widetilde{\OO}(|q|^2N)\\
&&+\frac{1}{|q|^2}\left[|q|^2\left(\squared_2-V+\frac{4a\cos\th}{|q|^2}\dual\nab_\T\right), \widetilde{\OO}\right]\psi. 
\eeaa
Using Lemma \ref{lemma:theoperqatorwidetildeOOcommutingwellwtihRWmodel} to estimate the last term, we infer
\beaa
\squared_2(\widetilde{\OO}\psi) - V\widetilde{\OO}\psi &=& -\frac{4a\cos\th}{|q|^2}\dual\nab_\T\widetilde{\OO}\psi +N_{\widetilde{\OO}},
\eeaa
with
\beaa
N_{\widetilde{\OO}} &=& \dk^{\leq 2}N+O(ar^{-2})\nab_{\Rhat}\dk^{\leq 1}\psi +O(ar^{-2})\dk^{\leq 1}\psi +\dk^{\leq 3}(\Ga_g\c\psi) +\Ga_b \c \squared_2\psi\\
&=& \dk^{\leq 2}N+O(ar^{-2})\nab_{\Rhat}\dk^{\leq 1}\psi +O(ar^{-2})\dk^{\leq 1}\psi +\dk^{\leq 3}(\Ga_g\c\psi)
\eeaa
where we used again  the model RW equation for $\psi$. In view of Proposition \ref{proposition:Energy1:perturbation}, we deduce
\beaa
&&E_{deg}[\widetilde{\OO}\psi](\tau_2) +F_{\Si_*}[\widetilde{\OO}\psi](\tau_1, \tau_2)\\
&\les& \deh\Big( E^2_{r\leq r_+}[\psi](\tau_2)+ F^2_{\AA}[\psi](\tau_1, \tau_2)\Big)+E^2[\psi](\tau_1) +\left|\int_{\MM(\tau_1, \tau_2)}  \nab_{\That_\de }(\widetilde{\OO}\psi)  \c N_{\widetilde{\OO}}\right|\\
&&+\left(\frac{|a|}{m}+\ep\right)\left(\sup_{\tau\in [\tau_1, \tau_2]}E^2[\psi]+B_\de^2[\psi](\tau_1, \tau_2)+F^2[\psi](\tau_1, \tau_2)\right).
\eeaa
Together with Lemma \ref{lemma:controloftheerrortermsGagpsiinwaveeq} and the structure of  $N_{(\Lieb_\T, \Lieb_\Z)^2}$, we deduce 
\beaa
&&E_{deg}[\widetilde{\OO}\psi](\tau_2) +F_{\Si_*}[\widetilde{\OO}\psi](\tau_1, \tau_2)\\
&\les& \deh\Big( E^2_{r\leq r_+}[\psi](\tau_2)+ F^2_{\AA}[\psi](\tau_1, \tau_2)\Big)+E^2[\psi](\tau_1) +\NN^2[\psi, N] (\tau_1, \tau_2)\\
&&+\left|\int_{\MM(\tau_1, \tau_2)}  \nab_{\That_\de }(\widetilde{\OO}\psi)  \c\Big(O(ar^{-2})\nab_{\Rhat}\dk^{\leq 1}\psi +O(ar^{-2})\dk^{\leq 1}\psi\Big)\right|\\
&&+\left(\frac{|a|}{m}+\ep\right)\left(\sup_{\tau\in [\tau_1, \tau_2]}E^2[\psi]+B_\de^2[\psi](\tau_1, \tau_2)+F^2[\psi](\tau_1, \tau_2)\right).
\eeaa
Integrating by parts the $\That_\de$ derivative, we obtain 
\beaa
&&E_{deg}[\widetilde{\OO}\psi](\tau_2) +F_{\Si_*}[\widetilde{\OO}\psi](\tau_1, \tau_2)\\
&\les& \deh\Big( E^2_{r\leq r_+}[\psi](\tau_2)+ F^2_{\AA}[\psi](\tau_1, \tau_2)\Big)+E^2[\psi](\tau_1) +\NN^2[\psi, N] (\tau_1, \tau_2)\\
&&+|a|\int_{\MM(\tau_1, \tau_2)}r^{-2}|\dk^{\leq 2}\psi|\Big(|\nab_{\Rhat}\dk^{\leq 2}\psi| +|\nab_{\That_\de}\dk^{\leq 1}\psi|\Big)\\
&&+\left(\frac{|a|}{m}+\ep\right)\left(\sup_{\tau\in [\tau_1, \tau_2]}E^2[\psi]+B_\de^2[\psi](\tau_1, \tau_2)+F^2[\psi](\tau_1, \tau_2)\right).
\eeaa
where we used the fact that $\That_\de(r)\in r\Ga_b$ and $\That_\de(\cos\th)\in \Ga_b$, the structure of the commutator $[\nab_\T, \nab_\Rhat]$ provided by Lemma \ref{lemma:commutationofnabRhatwithnabTandnabZ}, 
the control of $\D_\mu \That^\mu$ induced from the one of ${}^{(\That_\de)}\pi$ in the proof of Proposition \ref{proposition:Energy1:perturbation} in section \ref{sec:proofofresultschapter7inperturbationofKerr}, and Lemma \ref{lemma:controloftheerrortermsGagpsiinwaveeq} to control the various error terms. Since 
\beaa
&&\int_{\MM(\tau_1, \tau_2)}r^{-2}|\dk^{\leq 2}\psi|\Big(|\nab_{\Rhat}\dk^{\leq 2}\psi| +|\nab_{\That_\de}\dk^{\leq 1}\psi|\Big)\\
&\les& \int_{\MM_{trap}\tau_1, \tau_2)}\Big(|\nab_{\Rhat}\dk^{\leq 2}\psi|^2+|\dk^{\leq 2}\psi|^2\Big)+\int_{\Mntrap\tau_1, \tau_2)}r^{-2}\Big(|\nab_3\dk^{\leq 2}\psi||\dk^{\leq 3}\psi|+r^{-1}|\dk^{\leq 3}\psi|^2\Big)\\
&\les& B_\de^2[\psi](\tau_1, \tau_2),
\eeaa
we infer
\beaa
&&E_{deg}[\widetilde{\OO}\psi](\tau_2) +F_{\Si_*}[\widetilde{\OO}\psi](\tau_1, \tau_2)\\
&\les& \deh\Big( E^2_{r\leq r_+}[\psi](\tau_2)+ F^2_{\AA}[\psi](\tau_1, \tau_2)\Big)+E^2[\psi](\tau_1) +\NN^2[\psi, N] (\tau_1, \tau_2)\\
&&+\left(\frac{|a|}{m}+\ep\right)\left(\sup_{\tau\in [\tau_1, \tau_2]}E^2[\psi]+B_\de^2[\psi](\tau_1, \tau_2)+F^2[\psi](\tau_1, \tau_2)\right).
\eeaa
Together with the definition of $\widetilde{\OO}$ and  the control of $E_{deg}[(\nab_\T, \nab_\Z)^{\leq 1}\psi]$ and $F_{\Si_*}[(\nab_\T, \nab_\Z)^{\leq 1}\psi]$ provided by Proposition \ref{prop:energyforzeroandfirstderivativeinchap9}, we obtain 
\beaa
&&E_{deg}[\OO\psi](\tau_2) +F_{\Si_*}[\OO\psi](\tau_1, \tau_2)\\
&\les& \deh\Big( E^2_{r\leq r_+}[\psi](\tau_2)+ F^2_{\AA}[\psi](\tau_1, \tau_2)\Big)+E^2[\psi](\tau_1) +\NN^2[\psi, N] (\tau_1, \tau_2)\\
&&+\left(\frac{|a|}{m}+\ep\right)\left(\sup_{\tau\in [\tau_1, \tau_2]}E^2[\psi]+B_\de^2[\psi](\tau_1, \tau_2)+F^2[\psi](\tau_1, \tau_2)\right).
\eeaa
The is the stated energy estimate for $\psi_4=\OO\psi$. Together with the above energy estimates for $\psi_\aund$, $\aund=1,2,3$, this concludes the proof of Lemma \ref{lemma:contorloftheenergyofpsiaund}.
\end{proof}

The following proposition provides the control of Morawetz from the energy.
\begin{proposition}\lab{prop:Morawetzfromenergyinchap9}
The solution $\psi$ of of the model RW equation \eqref{eq:Gen.RW} satisfies the following Morawetz estimate 
\bea
\nn&&\de_{red}^7\Big(E^2_{r\leq r_+(1+\de_{red})}[\psi](\tau_2)+\textrm{Mor}[(\dkb, \nab_3, \nab_4)^{\leq 2}\psi](\tau_1, \tau_2)+F^2_{\AA}[\psi](\tau_1, \tau_2)\Big)\\
\nn&\les& \deh\left(\sup_{\tau\in[\tau_1, \tau_2]}E_{deg}[(\dkb, \nab_T, \nab_{\Rhat})^{\leq 2}\psi]+F_{\Si_*}[(\dkb, \nab_T, \nab_{\Rhat})^{\leq 2}\psi]\right)\\
\nn&& +\Bigg(E^1[\psi](\tau_1)+\NN^1[\psi, N] (\tau_1, \tau_2)\\
\nn&&+\left(\frac{|a|}{m}+\ep\right)\left(\sup_{\tau\in [\tau_1, \tau_2]}E^2[\psi]+B_\de^2[\psi](\tau_1, \tau_2)+F^2[\psi](\tau_1, \tau_2)\right)\Bigg)^{\frac{1}{2}}\\
\nn&&\times\left(\sup_{\tau\in[\tau_1, \tau_2]}E_{deg}[(\dkb, \nab_T, \nab_{\Rhat})^{\leq 2}\psi]+F_{\Si_*}[(\dkb, \nab_T, \nab_{\Rhat})^{\leq 2}\psi]\right)^{\frac{1}{2}}\\
\nn&&+\NN^2[\psi, N] (\tau_1, \tau_2)+\left(\frac{|a|}{m}+\ep\right)\left(\sup_{\tau\in [\tau_1, \tau_2]}E^2[\psi]+B_\de^2[\psi](\tau_1, \tau_2)+F^2[\psi](\tau_1, \tau_2)\right)\\
&&+E^2[\psi](\tau_1)
\eea
\end{proposition}

\begin{proof}
We proceed in several steps.

{\bf Step 1.} From Proposition \ref{prop:morawetz-higher-order:perturbation}, together with Lemma \ref{lemma:lowerboundPhizoutsideMtrap:perturbation}, we infer
\beaa
&& \textrm{Mor}_{deg}[\nab_T^2\psi](\tau_1, \tau_2)+\textrm{Mor}_{deg}[\OO\psi](\tau_1, \tau_2)+\textrm{Mor}_{deg}[\dkb\nab_T\psi](\tau_1, \tau_2)\\ 
&\les& \sum_{\aund=1}^4\left(\sup_{[\tau_1, \tau_2]}E_{deg}[\psi_\aund](\tau) +\deh F_{\AA}[\psi_{\aund}](\tau_1, \tau_2)+F_{\Si_*}[\psi_{\aund}](\tau_1, \tau_2)\right)\\
&&+\left(\sup_{\tau\in[\tau_1, \tau_2]}E_{deg}[(\dkb, \nab_T)^{\leq 2}\psi]+F_{\Si_*}[(\dkb, \nab_T)^{\leq 2}\psi]+\deh F_{\AA}[(\dkb, \nab_T)^{\leq 2}\psi]\right)^{\frac{1}{2}}\\
&&\times \left(\sup_{\tau\in[\tau_1, \tau_2]}E_{deg}[(\dkb, \nab_T)^{\leq 1}\psi]+F_{\Si_*}[(\dkb, \nab_T)^{\leq 1}\psi]+\deh F_{\AA}[(\dkb, \nab_T)^{\leq 1}\psi]\right)^{\frac{1}{2}}\\
&&+\int_{\MM(\tau_1, \tau_2)}(|\nab_{\Rhat}\psi_{\aund}|+r^{-1}|\psi_{\aund}|)|N_{\aund}|+\ep\left(\sup_{[\tau_1, \tau_2]}E^2[\psi](\tau)+ B^2_\de[\psi]\right).
\eeaa
Plugging the control of the energy provided by Proposition \ref{prop:energyforzeroandfirstderivativeinchap9} 
and Lemma \ref{lemma:contorloftheenergyofpsiaund} in the RHS, we infer 
\beaa
&& \textrm{Mor}_{deg}[\nab_T^2\psi](\tau_1, \tau_2)+\textrm{Mor}_{deg}[\OO\psi](\tau_1, \tau_2)+\textrm{Mor}_{deg}[\dkb\nab_T\psi](\tau_1, \tau_2)\\ 
&\les& \deh\Big( E^2_{r\leq r_+}[\psi](\tau_2)+ F^2_{\AA}[\psi](\tau_1, \tau_2)\Big)+\mathcal{I}_{RHS}+\int_{\MM(\tau_1, \tau_2)}(|\nab_{\Rhat}\psi_{\aund}|+r^{-1}|\psi_{\aund}|)|N_{\aund}|
\eeaa
where we have introduce the following notation 
\beaa
\mathcal{I}_{RHS} &:=& \deh\left(\sup_{\tau\in[\tau_1, \tau_2]}E_{deg}[(\dkb, \nab_T, \nab_{\Rhat})^{\leq 2}\psi]+F_{\Si_*}[(\dkb, \nab_T, \nab_{\Rhat})^{\leq 2}\psi]\right)\\
\nn&& +\Bigg(E^1[\psi](\tau_1)+\NN^1[\psi, N] (\tau_1, \tau_2)\\
\nn&&+\left(\frac{|a|}{m}+\ep\right)\left(\sup_{\tau\in [\tau_1, \tau_2]}E^2[\psi]+B_\de^2[\psi](\tau_1, \tau_2)+F^2[\psi](\tau_1, \tau_2)\right)\Bigg)^{\frac{1}{2}}\\
\nn&&\times\left(\sup_{\tau\in[\tau_1, \tau_2]}E_{deg}[(\dkb, \nab_T, \nab_{\Rhat})^{\leq 2}\psi]+F_{\Si_*}[(\dkb, \nab_T, \nab_{\Rhat})^{\leq 2}\psi]\right)^{\frac{1}{2}}\\
\nn&&+\NN^2[\psi, N] (\tau_1, \tau_2)+\left(\frac{|a|}{m}+\ep\right)\left(\sup_{\tau\in [\tau_1, \tau_2]}E^2[\psi]+B_\de^2[\psi](\tau_1, \tau_2)+F^2[\psi](\tau_1, \tau_2)\right)\\
&&+E^2[\psi](\tau_1).
\eeaa

Also, recall from \eqref{eq:estimaeforMundinperturbationsofKerr} that we have
\beaa
|N_{\aund}| &\les& |\dk^{\leq 2}N|+\frac{|a|}{m}r^{-2}|\dk^{\leq 1}\nab_3\psi|+\frac{|a|}{m}r^{-3}|\dk^{\leq 2}\psi|+|\dk^{\leq 3}(\Ga_g\c\psi)|\\
\eeaa
so that, using Lemma \ref{lemma:controloftheerrortermsGagpsiinwaveeq},  
\beaa
&&\int_{\MM(\tau_1, \tau_2)}(|\nab_{\Rhat}\psi_{\aund}|+r^{-1}|\psi_{\aund}|)|N_{\aund}|\\
 &\les& \NN^2[\psi, N] (\tau_1, \tau_2)+\left(\frac{|a|}{m}+\ep\right)\left(\sup_{\tau\in [\tau_1, \tau_2]}E^2[\psi]+B_\de^2[\psi](\tau_1, \tau_2)+F^2[\psi](\tau_1, \tau_2)\right)\\
 &\les& \mathcal{I}_{RHS}.
\eeaa
We infer
\beaa
 &&\textrm{Mor}_{deg}[\nab_T^2\psi](\tau_1, \tau_2)+\textrm{Mor}_{deg}[\OO\psi](\tau_1, \tau_2)+\textrm{Mor}_{deg}[\dkb\nab_T\psi](\tau_1, \tau_2)\\ 
 &\les& \mathcal{I}_{RHS}+ \deh\Big( E^2_{r\leq r_+}[\psi](\tau_2)+ F^2_{\AA}[\psi](\tau_1, \tau_2)\Big).
\eeaa
Moreover, using the definition  of $\OO$ and the Hodge estimates of Proposition \ref{Prop:HodgeThmM8}, we have
\beaa
\textrm{Mor}_{deg}[\dkb^{\leq 2}\psi](\tau_1, \tau_2) &\les& \textrm{Mor}_{deg}[r^2\Delta\psi](\tau_1, \tau_2)+\textrm{Mor}_{deg}[\nab_T^{\leq 2}\psi](\tau_1, \tau_2)\\
&\les& \textrm{Mor}_{deg}[\OO\psi](\tau_1, \tau_2)+\textrm{Mor}_{deg}[\nab_T^{\leq 2}\psi](\tau_1, \tau_2)+\frac{|a|}{m} B_\de^1[\psi](\tau_1, \tau_2)
\eeaa
and hence
\beaa
 &&\textrm{Mor}_{deg}[\nab_T^2\psi](\tau_1, \tau_2)+\textrm{Mor}_{deg}[\dkb^{\leq 2}\psi](\tau_1, \tau_2)+\textrm{Mor}_{deg}[\dkb\nab_T\psi](\tau_1, \tau_2)\\
 &\les& \deh\Big( E^2_{r\leq r_+}[\psi](\tau_2)+ F^2_{\AA}[\psi](\tau_1, \tau_2)\Big)+\mathcal{I}_{RHS} +\textrm{Mor}_{deg}[\nab_T^{\leq 1}\psi](\tau_1, \tau_2)+\frac{|a|}{m} B_\de^1[\psi](\tau_1, \tau_2)\\
 &\les& \deh\Big( E^2_{r\leq r_+}[\psi](\tau_2)+ F^2_{\AA}[\psi](\tau_1, \tau_2)\Big)+\mathcal{I}_{RHS} +\textrm{Mor}_{deg}[\nab_T^{\leq 1}\psi](\tau_1, \tau_2).
\eeaa
Together with the control of $\textrm{Mor}_{deg}[\nab_T^{\leq 1}\psi]$ provided by Proposition \ref{prop:energyforzeroandfirstderivativeinchap9}, we deduce
\beaa
\textrm{Mor}_{deg}[(\nab_T, \dkb)^{\leq 2}\psi](\tau_1, \tau_2) &\les& \deh\Big( E^2_{r\leq r_+}[\psi](\tau_2)+ F^2_{\AA}[\psi](\tau_1, \tau_2)\Big)+\mathcal{I}_{RHS}.
\eeaa

{\bf Step 2.} We now introduce the operators $\dkb_2^j$ acting on $\sk_2$ as follows, for $j\geq 0$, 
\bea
\dkb_2^{2j}:=(|q|^2\De_2)^j, \qquad \dkb_2^{2j+1}:=|q|\DD_2(|q|^2\De_2)^j.
\eea
Then, we commute the equation for $\psi$ by $\dkb_2^{j_1}\Lieb_\T^{j_2}$ with $j_1+j_2\leq 2$, and obtain in $r\leq r_+(1+2\de_{red})$, using the commutation formulas of Corollary \ref{cor:commutator-Lied-squared}, Lemma \ref{lemma:commutator-triangle} and Lemma \ref{LEMMA:COMMUTATIONOFHODGEELLIPTICORDER1WITHSQAURED2FDILUHS} 
 \beaa
\squared_{k_{j_1}}(\dkb_2^{j_1}\Lieb_\T^{j_2}\psi) -V\dkb_2^{j_1}\Lieb_\T^{j_2}\psi &=& - \frac{4 a\cos\th}{|q|^2}\dual \nab_T(\dkb_2^{j_1}\Lieb_\T^{j_2}\psi)+N_{\dkb_2^{j_1}\Lieb_\T^{j_2}},
\eeaa
where $k_0=k_2=2$ and $k_1=1$, and, we have, for $r\leq 4m$,
\beaa
N_{\dkb_2^{j_1}\Lieb_\T^{j_2}} := \dk^{\leq 2}N+\de_{j_1 1}O(1)\dkb_2^{j_1}\Lieb_\T^{j_2}\psi+O\left(\frac{|a|}{m}+\ep\right)\dk^{\leq 3}\psi.
\eeaa
Using the redshift estimates of Proposition \ref{Prop:Redshift-estimates-chp3}, we deduce, for $|j|=1,2$,
\beaa
\bsplit
&c_0E_{r\leq r_+(1+\de_{red})}[\dkb_2^{j_1}\Lieb_\T^{j_2}\psi](\tau_2)+c_0\textrm{Mor}_{r\leq r_+(1+\de_{red})}[\dkb_2^{j_1}\Lieb_\T^{j_2}\psi](\tau_1, \tau_2)+c_0F_{\AA}[\dkb_2^{j_1}\Lieb_\T^{j_2}\psi](\tau_1, \tau_2)\\
  \leq& E^2_{r\leq r_+(1+2\de_{red})}[\psi](\tau_1) +\de_{red}^{-1}\textrm{Mor}_{r_+(1+\de_{red})\leq r\leq r_+(1+2\de_{red})}[(\dkb, \Lieb_\T)^{\leq 2}\psi](\tau_1, \tau_2)\\
&+  \int_{\MM(\tau_1, \tau_2)\cap\left\{\frac{r}{r_+}\le 1+2\de_{red}\right\}}|N_{\dkb_2^{j_1}\Lieb_\T^{j_2}}|^2.
\end{split}
\eeaa
In view of the form of $N_{\dkb_2^{j_1}\Lieb_\T^{j_2}}$, we infer
\beaa
\bsplit
&c_0E_{r\leq r_+(1+\de_{red})}[\dkb_2^{j_1}\Lieb_\T^{j_2}\psi](\tau_2)+c_0\textrm{Mor}_{r\leq r_+(1+\de_{red})}[\dkb_2^{j_1}\Lieb_\T^{j_2}\psi](\tau_1, \tau_2)+c_0F_{\AA}[\dkb_2^{j_1}\Lieb_\T^{j_2}\psi](\tau_1, \tau_2)\\
  \leq& E^2_{r\leq r_+(1+2\de_{red})}[\psi](\tau_1)+(\de_{j_1 1}+\de_{j_1 2})\textrm{Mor}_{r\leq r_+(1+\de_{red})}[\dkb_2^{j_1-1}\Lieb_\T^{j_2}\psi](\tau_1, \tau_2) \\
&+\de_{red}^{-1}\textrm{Mor}_{r_+(1+\de_{red})\leq r\leq r_+(1+2\de_{red})}[(\dkb, \Lieb_\T)^{\leq 2}\psi](\tau_1, \tau_2)\\
&+\NN^2[\psi, N] (\tau_1, \tau_2)+\left(\frac{|a|}{m}+\ep\right)\left(\sup_{\tau\in [\tau_1, \tau_2]}E^2[\psi]+B_\de^2[\psi](\tau_1, \tau_2)+F^2[\psi](\tau_1, \tau_2)\right).
\end{split}
\eeaa
Arguing by iteration on $j_1$, using also the comparison of $\Lieb_\T$ and $\nab_\T$ of Lemma \ref{lemma:basicpropertiesLiebTfasdiuhakdisug:chap9}, and  the Hodge estimates of Proposition \ref{Prop:HodgeThmM8}, we obtain
\beaa
\bsplit
&c_0E_{r\leq r_+(1+\de_{red})}[(\dkb, \nab_\T)^{\leq 2}\psi](\tau_2)+c_0\textrm{Mor}_{r\leq r_+(1+\de_{red})}[(\dkb, \nab_\T)^{\leq 2}\psi](\tau_1, \tau_2)\\
&+c_0F_{\AA}[(\dkb, \nab_\T)^{\leq 2}\psi](\tau_1, \tau_2)\\
  \leq& E^2_{r\leq r_+(1+2\de_{red})}[\psi](\tau_1)+\de_{red}^{-1}\textrm{Mor}_{r_+(1+\de_{red})\leq r\leq r_+(1+2\de_{red})}[(\dkb, \nab_\T)^{\leq 2}\psi](\tau_1, \tau_2)\\
&+\NN^2[\psi, N] (\tau_1, \tau_2)+\left(\frac{|a|}{m}+\ep\right)\left(\sup_{\tau\in [\tau_1, \tau_2]}E^2[\psi]+B_\de^2[\psi](\tau_1, \tau_2)+F^2[\psi](\tau_1, \tau_2)\right).
\end{split}
\eeaa
Multiplying this estimate by $\de_{red}^6$, and summing it with the above control obtained in Step 1 for $\textrm{Mor}_{deg}[(\dkb, \nab_T)^{\leq 2}\psi](\tau_1, \tau_2)$, we infer
\beaa
&& \de_{red}^6\Big(E_{r\leq r_+(1+\de_{red})}[(\dkb, \nab_\T)^{\leq 2}\psi](\tau_2)+\textrm{Mor}_{r\leq r_+(1+\de_{red})}[(\dkb, \nab_\T)^{\leq 2}\psi](\tau_1, \tau_2)\\
&&+F_{\AA}[(\dkb, \nab_\T)^{\leq 2}\psi](\tau_1, \tau_2)\Big)+\textrm{Mor}_{deg}[(\dkb, \nab_T)^{\leq 2}\psi](\tau_1, \tau_2)\\ 
&\les& \deh\Big( E^2_{r\leq r_+}[\psi](\tau_2)+ F^2_{\AA}[\psi](\tau_1, \tau_2)\Big)+\mathcal{I}_{RHS}\\
&&+\de_{red}^5\textrm{Mor}_{r_+(1+\de_{red})\leq r\leq r_+(1+2\de_{red})}[(\dkb, \nab_\T)^{\leq 2}\psi](\tau_1, \tau_2).
\eeaa
Since 
\beaa
\de_{red}^5\textrm{Mor}_{r_+(1+\de_{red})\leq r\leq r_+(1+2\de_{red})}[(\dkb, \nab_\T)^{\leq 2}\psi](\tau_1, \tau_2) \les \de_{red}\textrm{Mor}_{deg}[(\dkb, \nab_T)^{\leq 2}\psi](\tau_1, \tau_2)
\eeaa
we may absorb the last term of the RHS by the LHS for $\de_{red}$ small enough, and hence
\beaa
&& \de_{red}^6\Big(E_{r\leq r_+(1+\de_{red})}[(\dkb, \nab_\T)^{\leq 2}\psi](\tau_2)+\textrm{Mor}_{r\leq r_+(1+\de_{red})}[(\dkb, \nab_\T)^{\leq 2}\psi](\tau_1, \tau_2)\\
&&+F_{\AA}[(\dkb, \nab_\T)^{\leq 2}\psi](\tau_1, \tau_2)\Big)+\textrm{Mor}_{deg}[(\dkb, \nab_T)^{\leq 2}\psi](\tau_1, \tau_2)\\ 
&\les& \deh\Big( E^2_{r\leq r_+}[\psi](\tau_2)+ F^2_{\AA}[\psi](\tau_1, \tau_2)\Big)+\mathcal{I}_{RHS}.
\eeaa

{\bf Step 3.} Next, recalling the representation of the wave operator provided by \eqref{eq:q2squaredkpsi}, i.e. 
\beaa
\begin{split}
|q|^2 \squared_2 \psi &=\frac{(r^2+a^2)^2}{\De} \big( -  \nab_\That \nab_\That \psi+   \nab_\Rhat \nab_\Rhat \psi \big) +2r \nab_\Rhat \psi\\
&+  |q|^2 \lap_2 \psi   + |q|^2  (\eta+\etab) \c \nab \psi  + r^2 \Ga_g \c \dk \psi,
\end{split}
\eeaa
we infer
\beaa
\int_{\MM}\Big(|\nab_{\Rhat}^3\psi|^2+|\nab_{\Rhat}^2\psi|^2\Big) &\les& \textrm{Mor}_{deg}[(\dkb, \nab_T)^{\leq 2}\psi](\tau_1, \tau_2)\\
&&+\int_{\MM}|\dk^{\leq 1}N|^2+\int_{\MM}|\dk^{\leq 1}(\Ga_g\c\psi)|^2.
\eeaa
Then,  we have by integration by parts
\beaa
&&\int_{\MM}\Big(|(\dkb, \nab_T)\nab_{\Rhat}^2\psi|^2+|(\dkb, \nab_T)\nab_{\Rhat}\psi|^2\Big)\\ 
&\les& \int_{\MM}\Big(|\nab_{\Rhat}(\dkb, \nab_T)^2\psi||\nab_{\Rhat}^3\psi|+|(\dkb, \nab_T)^2\psi||\nab_{\Rhat}^2\psi|\Big)\\
&&+\left(\sup_{\tau\in[\tau_1, \tau_2]}E_{deg}[(\dkb, \nab_T, \nab_{\Rhat})^{\leq 1}\psi]+F_{\Si_*}[(\dkb, \nab_T, \nab_{\Rhat})^{\leq 1}\psi]+\de_\HH F_\AA[(\dkb, \nab_T, \nab_{\Rhat})^{\leq 1}\psi]\right)^{\frac{1}{2}}\\&&\times\left(\sup_{\tau\in[\tau_1, \tau_2]}E_{deg}[(\dkb, \nab_T, \nab_{\Rhat})^{\leq 2}\psi]+F_{\Si_*}[(\dkb, \nab_T, \nab_{\Rhat})^{\leq 2}\psi]+\de_\HH F_\AA[(\dkb, \nab_T, \nab_{\Rhat})^{\leq 2}\psi]\right)^{\frac{1}{2}},
\eeaa
where the last term is used to  control the boundary terms. Since we have
\beaa
M_{deg}[(\nab_T, \nab_{\Rhat}, \dkb)^{\leq 2}\psi] &\les& M_{deg}[(\nab_T, \dkb)^{\leq 2}\psi] +\int_{\MM}\Big(|\nab_{\Rhat}^3\psi|^2+|\nab_{\Rhat}^2\psi|^2\Big)\\
&&+\int_{\MM}\Big(|(\dkb, \nab_T)\nab_{\Rhat}^2\psi|^2+|(\dkb, \nab_T)\nab_{\Rhat}\psi|^2\Big),
\eeaa
we infer from the above 
\beaa
&& M_{deg}[(\nab_T, \nab_{\Rhat}, \dkb)^{\leq 2}\psi]\\ 
&\les& M_{deg}[(\nab_T, \dkb)^{\leq 2}\psi] +\int_{\MM}|\dk^{\leq 1}N|^2+\int_{\MM}|\dk^{\leq 1}(\Ga_g\c\psi)|^2\\
&&+\left(\sup_{\tau\in[\tau_1, \tau_2]}E_{deg}[(\dkb, \nab_T, \nab_{\Rhat})^{\leq 1}\psi]+F_{\Si_*}[(\dkb, \nab_T, \nab_{\Rhat})^{\leq 1}\psi]+\de_\HH F_\AA[(\dkb, \nab_T, \nab_{\Rhat})^{\leq 1}\psi]\right)^{\frac{1}{2}}\\&&\times\left(\sup_{\tau\in[\tau_1, \tau_2]}E_{deg}[(\dkb, \nab_T, \nab_{\Rhat})^{\leq 2}\psi]+F_{\Si_*}[(\dkb, \nab_T, \nab_{\Rhat})^{\leq 2}\psi]+\de_\HH F_\AA[(\dkb, \nab_T, \nab_{\Rhat})^{\leq 2}\psi]\right)^{\frac{1}{2}}
\eeaa
and hence
\beaa
M_{deg}[(\nab_T, \nab_{\Rhat}, \dkb)^{\leq 2}\psi] &\les& M_{deg}[(\nab_T, \dkb)^{\leq 2}\psi] + \mathcal{I}_{RHS}
\eeaa
where we used the control of the energy in Proposition \ref{prop:energyforzeroandfirstderivativeinchap9} and the definition of $\mathcal{I}_{RHS}$. Together with the above estimate of Step 2, we infer
\beaa
 &&\de_{red}^6\Big(E_{r\leq r_+(1+\de_{red})}[(\dkb, \nab_\T)^{\leq 2}\psi](\tau_2)+\textrm{Mor}_{r\leq r_+(1+\de_{red})}[(\dkb, \nab_\T)^{\leq 2}\psi](\tau_1, \tau_2)\\
&&+F_{\AA}[(\dkb, \nab_\T)^{\leq 2}\psi](\tau_1, \tau_2)\Big)+\textrm{Mor}_{deg}[(\dkb, \nab_T, \nab_{\Rhat})^{\leq 2}\psi](\tau_1, \tau_2)\\ 
&\les& \deh\Big( E^2_{r\leq r_+}[\psi](\tau_2)+ F^2_{\AA}[\psi](\tau_1, \tau_2)\Big)+\mathcal{I}_{RHS}.
\eeaa

{\bf Step 4.} Next, we commute  the wave equation \eqref{eq:Gen.RW} for $\psi$ in $r\leq 4m$ by $(\nab_3, \dkb_2, \Lieb_\T)^{\leq 2}$, using the commutation formulas of Corollary \ref{cor:commutator-Lied-squared}, Lemma \ref{lemma:commutator-triangle} and Lemma \ref{LEMMA:COMMUTATIONOFHODGEELLIPTICORDER1WITHSQAURED2FDILUHS}, as well as Lemma \ref{lemma:commutationwithe3forredshift} for the commutations w.r.t to $\nab_3$. Together with Corollary \ref{cor:Redshift-estimates-chp3}, we infer
\beaa
\bsplit
&c_0E_{r\leq r_+(1+\de_{red})}[(\nab_3, \dkb_2, \Lieb_\T)^{\leq 2}\psi](\tau_2)+c_0\textrm{Mor}_{r\leq r_+(1+\de_{red})}[(\nab_3, \dkb_2, \Lieb_\T)^{\leq 2}\psi](\tau_1, \tau_2)\\
&+c_0F_{\AA}[(\nab_3, \dkb_2, \Lieb_\T)^{\leq 2}\psi](\tau_1, \tau_2)\\
  \leq& E_{r\leq r_+(1+2\de_{red})}^2[\psi](\tau_1) +\de_{red}^{-1}\textrm{Mor}_{r_+(1+\de_{red})\leq r\leq r_+(1+2\de_{red})}[(\nab_3, \dkb, \Lieb_\T)^{\leq 2}\psi](\tau_1, \tau_2)\\
&+\NN^2[\psi, N] (\tau_1, \tau_2)+\left(\frac{|a|}{m}+\ep\right)\left(\sup_{\tau\in [\tau_1, \tau_2]}E^2[\psi]+B_\de^2[\psi](\tau_1, \tau_2)+F^2[\psi](\tau_1, \tau_2)\right).
\end{split}
\eeaa 
Using the comparison of $\Lieb_\T$ and $\nab_\T$ of Lemma \ref{lemma:basicpropertiesLiebTfasdiuhakdisug:chap9}, and  the Hodge estimates of Proposition \ref{Prop:HodgeThmM8}, and using 
\beaa
&&\de_{red}^{-1}\textrm{Mor}_{r_+(1+\de_{red})\leq r\leq r_+(1+2\de_{red})}[(\nab_3, \dkb, \Lieb_\T)^{\leq 2}\psi](\tau_1, \tau_2)\\
&\les& \de_{red}^{-7}\textrm{Mor}_{r_+(1+\de_{red})\leq r\leq r_+(1+2\de_{red})}[(\dkb, \nab_\T, \nab_{\Rhat})^{\leq 2}\psi](\tau_1, \tau_2),
\eeaa
we infer
\beaa
\bsplit
&c_0E_{r\leq r_+(1+\de_{red})}[(\nab_3, \dkb, \nab_4)^{\leq 2}\psi](\tau_2)+c_0\textrm{Mor}_{r\leq r_+(1+\de_{red})}[(\nab_3, \dkb, \nab_4)^{\leq 2}\psi](\tau_1, \tau_2)\\
&+c_0F_{\AA}[(\nab_3, \dkb, \nab_4)^{\leq 2}\psi](\tau_1, \tau_2)\\
  \leq& E^2_{r\leq r_+(1+2\de_{red})}[\psi](\tau_1) \\
  &+\de_{red}^{-7}\textrm{Mor}_{r_+(1+\de_{red})\leq r\leq r_+(1+2\de_{red})}[(\dkb, \nab_\T, \nab_{\Rhat})^{\leq 2}\psi](\tau_1, \tau_2)\\
&+\NN^2[\psi, N] (\tau_1, \tau_2)+\left(\frac{|a|}{m}+\ep\right)\left(\sup_{\tau\in [\tau_1, \tau_2]}E^2[\psi]+B_\de^2[\psi](\tau_1, \tau_2)+F^2[\psi](\tau_1, \tau_2)\right).
\end{split}
\eeaa 
Together with the control of Step 3, we deduce
\beaa
&& \de_{red}^7\Big(E_{r\leq r_+(1+\de_{red})}[\dk^{\leq 2}\psi](\tau_2)+\textrm{Mor}[(\dkb, \nab_3, \nab_4)^{\leq 2}\psi](\tau_1, \tau_2)+F_{\AA}[\dk^{\leq 2}\psi](\tau_1, \tau_2)\Big)\\ 
&\les& \deh\Big( E^2_{r\leq r_+}[\psi](\tau_2)+ F^2_{\AA}[\psi](\tau_1, \tau_2)\Big)+\mathcal{I}_{RHS}.
\eeaa
Choosing $\de_{red}$ such that $\deh\ll \de_{red}^7$, we finally obtain 
\beaa
\de_{red}^7\Big(E^2_{r\leq r_+(1+\de_{red})}[\psi](\tau_2)+\textrm{Mor}[(\dkb, \nab_3, \nab_4)^{\leq 2}\psi](\tau_1, \tau_2)+F^2_{\AA}[\psi](\tau_1, \tau_2)\Big) \les \mathcal{I}_{RHS}.
\eeaa
In view of the definition of $\mathcal{I}_{RHS}$, this yields 
\beaa
&&\de_{red}^7\Big(E^2_{r\leq r_+(1+\de_{red})}[\psi](\tau_2)+\textrm{Mor}[(\dkb, \nab_3, \nab_4)^{\leq 2}\psi](\tau_1, \tau_2)+F^2_{\AA}[\psi](\tau_1, \tau_2)\Big)\\
&\les& \deh\left(\sup_{\tau\in[\tau_1, \tau_2]}E_{deg}[(\dkb, \nab_T, \nab_{\Rhat})^{\leq 2}\psi]+F_{\Si_*}[(\dkb, \nab_T, \nab_{\Rhat})^{\leq 2}\psi]\right)\\
\nn&& +\Bigg(E^1[\psi](\tau_1)+\NN^1[\psi, N] (\tau_1, \tau_2)\\
\nn&&+\left(\frac{|a|}{m}+\ep\right)\left(\sup_{\tau\in [\tau_1, \tau_2]}E^2[\psi]+B_\de^2[\psi](\tau_1, \tau_2)+F^2[\psi](\tau_1, \tau_2)\right)\Bigg)^{\frac{1}{2}}\\
\nn&&\times\left(\sup_{\tau\in[\tau_1, \tau_2]}E_{deg}[(\dkb, \nab_T, \nab_{\Rhat})^{\leq 2}\psi]+F_{\Si_*}[(\dkb, \nab_T, \nab_{\Rhat})^{\leq 2}\psi]\right)^{\frac{1}{2}}\\
\nn&&+\NN^2[\psi, N] (\tau_1, \tau_2)+\left(\frac{|a|}{m}+\ep\right)\left(\sup_{\tau\in [\tau_1, \tau_2]}E^2[\psi]+B_\de^2[\psi](\tau_1, \tau_2)+F^2[\psi](\tau_1, \tau_2)\right)\\
&&+E^2[\psi](\tau_1)
\eeaa
as stated. This concludes the proof of Proposition \ref{prop:Morawetzfromenergyinchap9}.
\end{proof}

The following proposition provides the control of the energy from Morawetz.
\begin{proposition}\lab{prop:energyfromMorawetzfrominchap9}
The solution $\psi$ of of the model RW equation \eqref{eq:Gen.RW} satisfies the following energy estimate 
\bea
\nn&& \de_{red}^8\Big(F^2_{\AA}[\psi](\tau_1, \tau_2)+E[(\nab_3, \nab_4, \dkb)^{\leq 2}\psi](\tau_2)+F_{\Si_*}[(\nab_3,  \nab_4, \dkb)^{\leq 2}\psi](\tau_1, \tau_2)\Big)\\
\nn&\les& E^2[\psi](\tau_1)+\textrm{Mor}[(\nab_3, \nab_4, \dkb)^{\leq 2}\psi](\tau_1, \tau_2) +\NN^2[\psi, N] (\tau_1, \tau_2)\\
&&+\left(\frac{|a|}{m}+\ep\right)\left(\sup_{\tau\in [\tau_1, \tau_2]}E^2[\psi]+B_\de^2[\psi](\tau_1, \tau_2)+F^2[\psi](\tau_1, \tau_2)\right).
\eea
\end{proposition}

\begin{proof}
The proof proceeds in several steps.

{\bf Step 1.} First, note from the control of Proposition \ref{prop:energyforzeroandfirstderivativeinchap9} and Lemma \ref{lemma:contorloftheenergyofpsiaund} that we have in particular 
\beaa
&& E_{deg}[(\dkb, \nab_T, \nab_{\Rhat})^{\leq 1}\psi](\tau_2)+F_{\Si_*}[(\dkb, \nab_T, \nab_{\Rhat})^{\leq 1}\psi](\tau_1, \tau_2)\\
&&+\sup_{\tau\in[\tau_1, \tau_2]}E_{deg}[(\nab_T, \nab_Z)^2\psi] +F_{\Si_*}[(\nab_T, \nab_Z)^2\psi](\tau_1, \tau_2)\\
&\les& \deh\Big(E^2_{r\leq r_+}[\psi](\tau_2)+ F^2_{\AA}[\psi](\tau_1, \tau_2)\Big)+E^2[\psi](\tau_1)+\NN^2[\psi, N] (\tau_1, \tau_2)\\
&&+\left(\frac{|a|}{m}+\ep\right)\left(\sup_{\tau\in [\tau_1, \tau_2]}E^2[\psi]+B_\de^2[\psi](\tau_1, \tau_2)+F^2[\psi](\tau_1, \tau_2)\right).
\eeaa

{\bf Step 2.} Commuting the wave equation \eqref{eq:Gen.RW} satisfied by $\psi$ with $\dkb_2^{j_1}\Lieb_\T^{j_2}$  with $j_1+j_2\leq 2$, and obtain, using the commutation formulas of using the commutation formulas of Corollary \ref{cor:commutator-Lied-squared}, Lemma \ref{lemma:commutator-triangle} and Lemma \ref{LEMMA:COMMUTATIONOFHODGEELLIPTICORDER1WITHSQAURED2FDILUHS}
 \beaa
\squared_{k_{j_1}}(\dkb_2^{j_1}\Lieb_\T^{j_2}\psi) -V\dkb_2^{j_1}\Lieb_\T^{j_2}\psi &=& - \frac{4 a\cos\th}{|q|^2}\dual \nab_T(\dkb_2^{j_1}\Lieb_\T^{j_2}\psi)+N_{\dkb_2^{j_1}\Lieb_\T^{j_2}},
\eeaa
where $k_0=k_2=2$ and $k_1=1$, and
\beaa
N_{\dkb_2^{j_1}\Lieb_\T^{j_2}} := \dk^{\leq 2}N -\de_{j_1 1}\frac{3}{r^2}\dkb_2^{j_1}\Lieb_\T^{j_2}\psi+O(ar^{-2})\dk^{\leq 3}\psi+\dk^{\leq 3}(\Ga_g\c\psi).
\eeaa
Next, we rely on  Proposition \ref{prop:recoverEnergyMorawetzwithrweightfromnoweight:perturbation}. We infer, for $r\geq r_1$ with $r_1=r_1(m)$ sufficiently large, 
together with the Hodge estimates of Proposition \ref{Prop:HodgeThmM8}
\beaa
&&E_{r\geq 2r_1}[(\nab_T, \dkb)^{\leq 2}\psi](\tau_2)+F_{\Si_*}[(\nab_T, \dkb)^{\leq 2}\psi](\tau_1, \tau_2)\\
&\les&  E^2[\psi](\tau_1) +\frac{r_1}{m}\textrm{Mor}_{r_1\leq r\leq 2r_1}[(\nab_T, \dkb)^{\leq 2}\psi](\tau_1, \tau_2) \\
&&+\NN^2[\psi, N] (\tau_1, \tau_2)+\left(\frac{|a|}{m}+\ep\right)\left(\sup_{\tau\in [\tau_1, \tau_2]}E^2[\psi]+B_\de^2[\psi](\tau_1, \tau_2)+F^2[\psi](\tau_1, \tau_2)\right).
\eeaa
Also, using again  the representation of the wave operator provided by \eqref{eq:q2squaredkpsi}, we have
\beaa
&& E_{r\geq 2r_1}[\nab_{\Rhat}^2\psi]+F_{\Si_*}[\nab_{\Rhat}^2\psi]\\
&\les& E_{r\geq 2r_1}[(\nab_T, \dkb)^{\leq 2}\psi](\tau_2)+F_{\Si_*}[(\nab_T, \dkb)^{\leq 2}\psi](\tau_1, \tau_2)+\int_{\Si(\tau_2)}|\dk^{\leq 1}N|^2\\
&&+\int_{\Si(\tau_2)}|\dk^{\leq 1}(\Ga_g\c\dk\psi)|^2+\int_{\Si_*(\tau_1, \tau_2)}|\dk^{\leq 1}N|^2+\int_{\Si_*(\tau_1, \tau_2)}|\dk^{\leq 1}(\Ga_g\c\dk\psi)|^2.
\eeaa
Since 
\beaa
&& E_{r\geq 2r_1}[(\nab_T, \nab_{\Rhat}, \dkb)^{\leq 2}\psi](\tau_2)+F_{\Si_*}[(\nab_T, \nab_{\Rhat}, \dkb)^{\leq 2}\psi](\tau_2) \\
 &\les& E_{r\geq 2r_1}[(\nab_T, \dkb)^{\leq 2}\psi](\tau_2)+F_{\Si_*}[(\nab_T, \dkb)^{\leq 2}\psi](\tau_1, \tau_2)+E_{r\geq 2r_1}[\nab_{\Rhat}^2\psi]+F_{\Si_*}[\nab_{\Rhat}^2\psi],
\eeaa
we infer from the above
\beaa
&&E_{r\geq 2r_1}[(\nab_T, \dkb, \nab_{\Rhat})^{\leq 2}\psi](\tau_2)+F_{\Si_*}[(\nab_T, \dkb, \nab_{\Rhat})^{\leq 2}\psi](\tau_1, \tau_2)\\
&\les&  E^2[\psi](\tau_1) +\frac{r_1}{m}\textrm{Mor}_{r_1\leq r\leq 2r_1}[(\nab_T, \dkb)^{\leq 2}\psi](\tau_1, \tau_2) \\
&&+\NN^2[\psi, N] (\tau_1, \tau_2)+\left(\frac{|a|}{m}+\ep\right)\left(\sup_{\tau\in [\tau_1, \tau_2]}E^2[\psi]+B_\de^2[\psi](\tau_1, \tau_2)+F^2[\psi](\tau_1, \tau_2)\right).
\eeaa

{\bf Step 3.} Next, we commute  the wave equation \eqref{eq:Gen.RW} for $\psi$ in $\leq 4m$ by $(\nab_3, |q|\nab, \Lieb_\T)^{\leq 2}$, using the commutation formulas of Corollary \ref{cor:commutator-Lied-squared}, Lemma \ref{lemma:commutator-triangle} and Lemma \ref{LEMMA:COMMUTATIONOFHODGEELLIPTICORDER1WITHSQAURED2FDILUHS}, as well as Lemma \ref{lemma:commutationwithe3forredshift} for the commutations w.r.t to $\nab_3$. 
Together with Corollary \ref{cor:Redshift-estimates-chp3}, we infer, as in Step 4 of the proof of Proposition \ref{prop:Morawetzfromenergyinchap9}, 
\beaa
\bsplit
&c_0E_{r\leq r_+(1+\de_{red})}[(\nab_3, \dkb_2, \Lieb_\T)^{\leq 2}\psi](\tau_2)+c_0\textrm{Mor}_{r\leq r_+(1+\de_{red})}[(\nab_3, \dkb_2, \Lieb_\T)^{\leq 2}\psi](\tau_1, \tau_2)\\
&+c_0F_{\AA}[(\nab_3, \dkb_2, \Lieb_\T)^{\leq 2}\psi](\tau_1, \tau_2)\\
  \leq& E_{r\leq r_+(1+2\de_{red})}^2[\psi](\tau_1) +\de_{red}^{-1}\textrm{Mor}_{r_+(1+\de_{red})\leq r\leq r_+(1+2\de_{red})}[(\nab_3, \dkb, \Lieb_\T)^{\leq 2}\psi](\tau_1, \tau_2)\\
&+\NN^2[\psi, N] (\tau_1, \tau_2)+\left(\frac{|a|}{m}+\ep\right)\left(\sup_{\tau\in [\tau_1, \tau_2]}E^2[\psi]+B_\de^2[\psi](\tau_1, \tau_2)+F^2[\psi](\tau_1, \tau_2)\right).
\end{split}
\eeaa 
from which we deduce, dropping the Morawetz term on the LHS, and using the Hodge estimates of Proposition \ref{Prop:HodgeThmM8}
\beaa
\bsplit
&c_0E_{r\leq r_+(1+\de_{red})}[\dk^{\leq 2}\psi](\tau_2)+c_0F_{\AA}[\dk^{\leq 2}\psi](\tau_1, \tau_2)\\
  \leq& E^2[\psi](\tau_1) +\de_{red}^{-1}\textrm{Mor}_{r_+(1+\de_{red})\leq r\leq r_+(1+2\de_{red})}[(\nab_3, \dkb, \nab_\T)^{\leq 2}\psi](\tau_1, \tau_2)\\
&+\NN^2[\psi, N] (\tau_1, \tau_2)+\left(\frac{|a|}{m}+\ep\right)\left(\sup_{\tau\in [\tau_1, \tau_2]}E^2[\psi]+B_\de^2[\psi](\tau_1, \tau_2)+F^2[\psi](\tau_1, \tau_2)\right).
\end{split}
\eeaa 

{\bf Step 4.} Using again  the representation of the wave operator provided by \eqref{eq:q2squaredkpsi}, i.e. 
\beaa
\begin{split}
|q|^2 \squared_2 \psi &=\frac{(r^2+a^2)^2}{\De} \big( -  \nab_\That \nab_\That \psi+   \nab_\Rhat \nab_\Rhat \psi \big) +2r \nab_\Rhat \psi\\
&+  |q|^2 \lap_2 \psi   + |q|^2  (\eta+\etab) \c \nab \psi  + r^2 \Ga_g \c \dk \psi,
\end{split}
\eeaa
together with the wave equation \eqref{eq:Gen.RW} satisfied by $\psi$,  we infer
\beaa
E_{deg}[U] \les E_{deg}[(\nab_T, \nab_Z)^{\leq 2}\psi](\tau)+E_{deg}[(\nab_T, \nab, \nab_{\Rhat})^{\leq 1}\psi](\tau)+\int_{\Si(\tau)}|\dk^{\leq 1}N|^2,
\eeaa
where the constant in $\les$ is independent of $\tau$, and where the tensor $U$ is given by
\beaa
U :=  \nab_\Rhat \nab_\Rhat \psi +  \frac{|q|^2\De}{(r^2+a^2)^2} \lap_2 \psi. 
\eeaa
Next, we proceed as follows:

{\bf Step 4a.} We consider a cut-off $\ka_{\deh}(r)$ such that $0\leq\ka_{\deh}\leq 1$, $\ka_{\deh}=1$ on $r\geq r_+$, and $\ka_{\deh}$ vanishes on $\AA$. Also, for any integer $2\leq n\leq \tau_*-2$, we consider a cut-off $\ka_n(\tau)$ such that $0\leq \ka_n\leq 1$, $\ka_n(\tau)=1$ on $\{n\leq \tau\leq n+1\}$, and $\ka_n$ vanishes on $\tau\leq n-1$ and $\tau\geq n+2$. We have, by  integrating by parts, and using  the support properties of $\ka_n(\tau)$ and the above control of $U$, 
\beaa
&&\int_{\MM}\ka_{\deh}(r)\ka_n(\tau)\left(|\nab_4\nab_{\Rhat}^2\psi|^2+\frac{|\De|}{r^2}|\nab_4\nab\nab_{\Rhat}\psi|^2+\frac{|\De|^2}{r^4}|\nab_4\nab^2\psi|^2\right)\\
&&+\int_{\MM}\ka_{\deh}(r)\ka_n(\tau)r^{-2}\left(|\nab_{\That}\nab_{\Rhat}^2\psi|^2+\frac{|\De|}{r^2}|\nab_{\That}\nab\nab_{\Rhat}\psi|^2+\frac{|\De|^2}{r^4}|\nab_{\That}\nab^2\psi|^2\right)\\
&&+\int_{\MM}\ka_{\deh}(r)\ka_n(\tau)\left(|\nab\nab_{\Rhat}^2\psi|^2+\frac{|\De|}{r^2}|\nab\nab_{\Rhat}\nab\psi|^2+\frac{|\De|^2}{r^4}|\nab^3\psi|^2\right)\\
&\les& \sup_{\tau\in[n-1,n+2]}E_{deg}[U]\\
&&+\left(\sup_{\tau\in[n-1,n+2]}E_{deg}[(\nab_T, \nab, \nab_{\Rhat})^{\leq 1}\psi](\tau)+F_{\Si_*}[(\nab_T, \nab_{\Rhat}, \nab)^{\leq 1}\psi](n-1, n+2)\right)^{\frac{1}{2}}\\
&&\times\left(\sup_{\tau\in[n-1,n+2]}E_{deg}[(\nab_T, \nab, \nab_{\Rhat})^{\leq 2}\psi](\tau)+F_{\Si_*}[(\nab_T, \nab_{\Rhat}, \nab)^{\leq 2}\psi](n-1, n+2)\right)^{\frac{1}{2}}.
\eeaa
Together with the above estimates, we deduce 
\beaa
&&\int_{\MM}\ka_{\deh}(r)\ka_n(\tau)\left(|\nab_4\nab_{\Rhat}^2\psi|^2+\frac{|\De|}{r^2}|\nab_4\nab\nab_{\Rhat}\psi|^2+\frac{|\De|^2}{r^4}|\nab_4\nab^2\psi|^2\right)\\
&&+\int_{\MM}\ka_{\deh}(r)\ka_n(\tau)r^{-2}\left(|\nab_{\That}\nab_{\Rhat}^2\psi|^2+\frac{|\De|}{r^2}|\nab_{\That}\nab\nab_{\Rhat}\psi|^2+\frac{|\De|^2}{r^4}|\nab_{\That}\nab^2\psi|^2\right)\\
&&+\int_{\MM}\ka_{\deh}(r)\ka_n(\tau)\left(|\nab\nab_{\Rhat}^2\psi|^2+\frac{|\De|}{r^2}|\nab\nab_{\Rhat}\nab\psi|^2+\frac{|\De|^2}{r^4}|\nab^3\psi|^2\right)\\
&\les& \sup_{\tau\in[n-1,n+2]}E_{deg}[(\nab_T, \nab_Z)^{\leq 2}\psi](\tau)+\sup_{\tau\in[n-1,n+2]}E_{deg}[(\nab_T, \nab, \nab_{\Rhat})^{\leq 1}\psi](\tau)\\
&&+\left(\sup_{\tau\in[n-1,n+2]}E_{deg}[(\nab_T, \nab, \nab_{\Rhat})^{\leq 1}\psi](\tau)+F_{\Si_*}[(\nab_T, \nab_{\Rhat}, \nab)^{\leq 1}\psi](n-1, n+2)\right)^{\frac{1}{2}}\\
&&\times\left(\sup_{\tau\in[n-1,n+2]}E_{deg}[(\nab_T, \nab, \nab_{\Rhat})^{\leq 2}\psi](\tau)+F_{\Si_*}[(\nab_T, \nab_{\Rhat}, \nab)^{\leq 2}\psi](n-1, n+2)\right)^{\frac{1}{2}}\\
&&+\sup_{\tau\in[n-1,n+2]}\int_{\Si(\tau)}|\dk^{\leq 1}N|^2.
\eeaa

{\bf Step 4b.} In view of the properties of the cut-offs, the above estimates imply 
\beaa
&&\int_{\MM\cap\{r\geq r_+\}\cap\{n\leq\tau\leq n+1\}}\Bigg[|\nab_4\nab_{\Rhat}^2\psi|^2+\frac{|\De|}{r^2}|\nab_4\nab\nab_{\Rhat}\psi|^2+\frac{|\De|^2}{r^4}|\nab_4\nab^2\psi|^2\\
&&+r^{-2}\left(|\nab_{\That}\nab_{\Rhat}^2\psi|^2+\frac{|\De|}{r^2}|\nab_{\That}\nab\nab_{\Rhat}\psi|^2+\frac{|\De|^2}{r^4}|\nab_{\That}\nab^2\psi|^2\right)+|\nab\nab_{\Rhat}^2\psi|^2\\
&&+\frac{|\De|}{r^2}|\nab\nab_{\Rhat}\nab\psi|^2+\frac{|\De|^2}{r^4}|\nab^3\psi|^2\Bigg]\\
&\les& \sup_{\tau\in[n-1,n+2]}E_{deg}[(\nab_T, \nab_Z)^{\leq 2}\psi](\tau)+\sup_{\tau\in[n-1,n+2]}E_{deg}[(\nab_T, \nab, \nab_{\Rhat})^{\leq 1}\psi](\tau)\\
&&+\left(\sup_{\tau\in[n-1,n+2]}E_{deg}[(\nab_T, \nab, \nab_{\Rhat})^{\leq 1}\psi](\tau)+F_{\Si_*}[(\nab_T, \nab_{\Rhat}, \nab)^{\leq 1}\psi](n-1, n+2)\right)^{\frac{1}{2}}\\
&&\times\left(\sup_{\tau\in[n-1,n+2]}E_{deg}[(\nab_T, \nab, \nab_{\Rhat})^{\leq 2}\psi](\tau)+F_{\Si_*}[(\nab_T, \nab_{\Rhat}, \nab)^{\leq 2}\psi](n-1, n+2)\right)^{\frac{1}{2}}\\
&&+\sup_{\tau\in[n-1,n+2]}\int_{\Si(\tau)}|\dk^{\leq 1}N|^2.
\eeaa
By the mean value theorem, we infer the existence of $\tau_{(n)}\in[n,n+1]$ such that 
\beaa
&&\int_{\Si(\tau_{(n)})\cap\{r\geq r_+\}}\Bigg[|\nab_4\nab_{\Rhat}^2\psi|^2+\frac{|\De|}{r^2}|\nab_4\nab\nab_{\Rhat}\psi|^2+\frac{|\De|^2}{r^4}|\nab_4\nab^2\psi|^2\\
&&+r^{-2}\left(|\nab_{\That}\nab_{\Rhat}^2\psi|^2+\frac{|\De|}{r^2}|\nab_{\That}\nab\nab_{\Rhat}\psi|^2+\frac{|\De|^2}{r^4}|\nab_{\That}\nab^2\psi|^2\right)+|\nab\nab_{\Rhat}^2\psi|^2\\
&&+\frac{|\De|}{r^2}|\nab\nab_{\Rhat}\nab\psi|^2+\frac{|\De|^2}{r^4}|\nab^3\psi|^2\Bigg]\\
&\les&  \sup_{\tau\in[n-1,n+2]}E_{deg}[(\nab_T, \nab_Z)^{\leq 2}\psi](\tau)+\sup_{\tau\in[n-1,n+2]}E_{deg}[(\nab_T, \nab, \nab_{\Rhat})^{\leq 1}\psi](\tau)\\
&&+\left(\sup_{\tau\in[n-1,n+2]}E_{deg}[(\nab_T, \nab, \nab_{\Rhat})^{\leq 1}\psi](\tau)+F_{\Si_*}[(\nab_T, \nab_{\Rhat}, \nab)^{\leq 1}\psi](n-1, n+2)\right)^{\frac{1}{2}}\\
&&\times\left(\sup_{\tau\in[n-1,n+2]}E_{deg}[(\nab_T, \nab, \nab_{\Rhat})^{\leq 2}\psi](\tau)+F_{\Si_*}[(\nab_T, \nab_{\Rhat}, \nab)^{\leq 2}\psi](n-1, n+2)\right)^{\frac{1}{2}}\\
&&+\sup_{\tau\in[n-1,n+2]}\int_{\Si(\tau)}|\dk^{\leq 1}N|^2.
\eeaa

{\bf Step 4c.} Together with the above control of $E_{deg}[(\nab_T, \nab_Z)^{\leq 2}\psi](\tau)$, the above control of $E_{r\geq 2r_1}[(\nab_T, \dkb, \nab_{\Rhat})^{\leq 2}\psi](\tau_2)$ and $F_{\Si_*}[(\nab_T, \dkb, \nab_{\Rhat})^{\leq 2}\psi](\tau_1, \tau_2)$, and the above control of $E_{r\leq r_+(1+\de_{red})}[\dk^{\leq 2}\psi](\tau_2)$ and $F_{\AA}[\dk^{\leq 2}\psi](\tau_1, \tau_2)$, and fixing the value of $r_1=r_1(m)$, we infer, the following non sharp estimate, for any $n$ such that $\tau_{(n)}\geq \tau_1$,  
\beaa
&& \de_{red}^4\Big(F_{\AA}[\dk^{\leq 2}\psi](\tau_1, \tau_{(n)})+E[(\nab_3, \nab_4, \dkb)^{\leq 2}\psi](\tau_{(n)})+F_{\Si_*}[(\nab_3,  \nab_4, \dkb)^{\leq 2}\psi](\tau_1, \tau_{(n)})\Big)\\
&\les&  \sup_{\tau\in[n-1,n+2]}E_{deg}[(\nab_T, \nab, \nab_{\Rhat})^{\leq 1}\psi](\tau)+\deh\sup_{\tau\in[\tau_1, n+2]}E_{r\leq r_+}[(\nab_T, \nab_Z)^{\leq 2}\psi]\\
&&+\deh F_{\AA}[(\nab_T, \nab_Z)^{\leq 2}\psi](\tau_1, n+2)\\
&&+\left(\sup_{\tau\in[n-1,n+2]}E_{deg}[(\nab_T, \nab, \nab_{\Rhat})^{\leq 1}\psi](\tau)+F_{\Si_*}[(\nab_T, \nab_{\Rhat}, \nab)^{\leq 1}\psi](n-1, n+2)\right)^{\frac{1}{2}}\\
&&\times\left(\sup_{\tau\in[n-1,n+2]}E_{deg}[(\nab_T, \nab, \nab_{\Rhat})^{\leq 2}\psi](\tau)+F_{\Si_*}[(\nab_T, \nab_{\Rhat}, \nab)^{\leq 2}\psi](n-1, n+2)\right)^{\frac{1}{2}}\\
&&+E^2[\psi](\tau_1)+\textrm{Mor}[(\nab_3, \nab_4, \dkb)^{\leq 2}\psi](\tau_1, n+2)\\
&&+\NN^2[\psi, N] (\tau_1, n+2)+\left(\frac{|a|}{m}+\ep\right)\left(\sup_{\tau\in [\tau_1, n+2]}E^2[\psi]+B_\de^2[\psi](\tau_1, n+2)+F^2[\psi](\tau_1, n+2)\right).
\eeaa

{\bf Step 4d.} Let $\tau_2\geq \tau_1$. By local energy estimates, it suffices to consider the case $\tau_2\geq \tau_1+5$. We then choose $n$ such that $\tau_1\leq n-1\leq n+2\leq\tau_2<n+3$. In particular, we have $\tau_{(n)}+1\leq \tau_2\leq \tau_{(n)}+3$, and hence, using local energy estimates between $\tau_{(n)}$ and $\tau_2$, we infer from the previous estimate, choosing also $\de_\HH\ll \de_{red}^4$ to absorb some terms on the RHS from the LHS,  
\beaa
&& \de_{red}^4\Big(F_{\AA}[\dk^{\leq 2}\psi](\tau_1, \tau_2)+E[(\nab_3, \nab_4, \dkb)^{\leq 2}\psi](\tau_2)+F_{\Si_*}[(\nab_3,  \nab_4, \dkb)^{\leq 2}\psi](\tau_1, \tau_2)\Big)\\
&\les&  \sup_{\tau\in[\tau_1,\tau_2]}E_{deg}[(\nab_T, \nab, \nab_{\Rhat})^{\leq 1}\psi](\tau)\\
&&+\left(\sup_{\tau\in[\tau_1,\tau_2]}E_{deg}[(\nab_T, \nab, \nab_{\Rhat})^{\leq 1}\psi](\tau)+F_{\Si_*}[(\nab_T, \nab_{\Rhat}, \nab)^{\leq 1}\psi](\tau_1, \tau_2)\right)^{\frac{1}{2}}\\
&&\times\left(\sup_{\tau\in[\tau_1,\tau_2]}E_{deg}[(\nab_T, \nab, \nab_{\Rhat})^{\leq 2}\psi](\tau)+F_{\Si_*}[(\nab_T, \nab_{\Rhat}, \nab)^{\leq 2}\psi](\tau_1, \tau_2)\right)^{\frac{1}{2}}\\
&&+E^2[\psi](\tau_1)+\textrm{Mor}[(\nab_3, \nab_4, \dkb)^{\leq 2}\psi](\tau_1, \tau_2)\\
&&+\NN^2[\psi, N] (\tau_1, \tau_2)+\left(\frac{|a|}{m}+\ep\right)\left(\sup_{\tau\in [\tau_1, \tau_2]}E^2[\psi]+B_\de^2[\psi](\tau_1, \tau_2)+F^2[\psi](\tau_1, \tau_2)\right)
\eeaa
and thus
\beaa
&& \de_{red}^8\Big(F_{\AA}[\dk^{\leq 2}\psi](\tau_1, \tau_2)+E[(\nab_3, \nab_4, \dkb)^{\leq 2}\psi](\tau_2)+F_{\Si_*}[(\nab_3,  \nab_4, \dkb)^{\leq 2}\psi](\tau_1, \tau_2)\Big)\\
&\les&  \sup_{\tau\in[\tau_1,\tau_2]}E_{deg}[(\nab_T, \nab, \nab_{\Rhat})^{\leq 1}\psi](\tau)+F_{\Si_*}[(\nab_T, \nab_{\Rhat}, \nab)^{\leq 1}\psi](\tau_1, \tau_2)\\
&&+E^2[\psi](\tau_1)+\textrm{Mor}[(\nab_3, \nab_4, \dkb)^{\leq 2}\psi](\tau_1, \tau_2)\\
&&+\NN^2[\psi, N] (\tau_1, \tau_2)+\left(\frac{|a|}{m}+\ep\right)\left(\sup_{\tau\in [\tau_1, \tau_2]}E^2[\psi]+B_\de^2[\psi](\tau_1, \tau_2)+F^2[\psi](\tau_1, \tau_2)\right).
\eeaa

{\bf Step 5.} Plugging the control of $E_{deg}[(\nab_T, \nab, \nab_{\Rhat})^{\leq 1}\psi](\tau_2)$ and $F_{\Si_*}[(\dkb, \nab_T, \nab_{\Rhat})^{\leq 1}\psi](\tau_1, \tau_2)$ derived in Proposition \ref{prop:energyforzeroandfirstderivativeinchap9}, we infer, choosing also $\deh$ such that $\deh\ll \de_{red}^8$ to absorb terms from the RHS by the LHS, 
\beaa
&& \de_{red}^8\Big(F_{\AA}[\dk^{\leq 2}\psi](\tau_1, \tau_2)+E[(\nab_3, \nab_4, \dkb)^{\leq 2}\psi](\tau_2)+F_{\Si_*}[(\nab_3,  \nab_4, \dkb)^{\leq 2}\psi](\tau_1, \tau_2)\Big)\\
&\les& E^2[\psi](\tau_1)+\textrm{Mor}[(\nab_3, \nab_4, \dkb)^{\leq 2}\psi](\tau_1, \tau_2)\\ 
&&+\NN^2[\psi, N] (\tau_1, \tau_2)+\left(\frac{|a|}{m}+\ep\right)\left(\sup_{\tau\in [\tau_1, \tau_2]}E^2[\psi]+B_\de^2[\psi](\tau_1, \tau_2)+F^2[\psi](\tau_1, \tau_2)\right)
\eeaa
as stated. This concludes the proof of Proposition \ref{prop:energyfromMorawetzfrominchap9}.
\end{proof}

We are now ready to prove Theorem \ref{THM:HIGHERDERIVS-MORAWETZ-CHP3} in the particular case $s=2$.
\begin{proof}[Proof of Theorem \ref{THM:HIGHERDERIVS-MORAWETZ-CHP3} in the case $s=2$]
Plugging the energy estimate of Proposition \ref{prop:energyfromMorawetzfrominchap9} in the RHS of the Morawetz estimate of Proposition \ref{prop:Morawetzfromenergyinchap9}, we infer
\beaa
\nn&& \de_{red}^7\Big(E^2_{r\leq r_+(1+\de_{red})}[\psi](\tau_2)+\textrm{Mor}[(\dkb, \nab_3, \nab_4)^{\leq 2}\psi](\tau_1, \tau_2)+F^2_{\AA}[\psi](\tau_1, \tau_2)\Big)\\ 
\nn&\les&  \de_{red}^{-8}\deh\Bigg(E^2[\psi](\tau_1)+\textrm{Mor}[(\nab_3, \nab_4, \dkb)^{\leq 2}\psi](\tau_1, \tau_2)\\ 
&&+\NN^2[\psi, N] (\tau_1, \tau_2)+\left(\frac{|a|}{m}+\ep\right)\left(\sup_{\tau\in [\tau_1, \tau_2]}E^2[\psi]+B_\de^2[\psi](\tau_1, \tau_2)+F^2[\psi](\tau_1, \tau_2)\right)\Bigg)\\
\nn&& +\de_{red}^{-4}\Bigg(E^1[\psi](\tau_1)+\NN^1[\psi, N] (\tau_1, \tau_2)\\
\nn&&+\left(\frac{|a|}{m}+\ep\right)\left(\sup_{\tau\in [\tau_1, \tau_2]}E^2[\psi]+B_\de^2[\psi](\tau_1, \tau_2)+F^2[\psi](\tau_1, \tau_2)\right)\Bigg)^{\frac{1}{2}}\\
\nn&&\times\Bigg(E^2[\psi](\tau_1)+\textrm{Mor}[(\nab_3, \nab_4, \dkb)^{\leq 2}\psi](\tau_1, \tau_2)\\ 
&&+\NN^2[\psi, N] (\tau_1, \tau_2)+\left(\frac{|a|}{m}+\ep\right)\left(\sup_{\tau\in [\tau_1, \tau_2]}E^2[\psi]+B_\de^2[\psi](\tau_1, \tau_2)+F^2[\psi](\tau_1, \tau_2)\right)\Bigg)^{\frac{1}{2}}\\
\nn&&+\NN^2[\psi, N] (\tau_1, \tau_2)+\left(\frac{|a|}{m}+\ep\right)\left(\sup_{\tau\in [\tau_1, \tau_2]}E^2[\psi]+B_\de^2[\psi](\tau_1, \tau_2)+F^2[\psi](\tau_1, \tau_2)\right)\\
&&+E^2[\psi](\tau_1).
\eeaa
Choosing $\deh$  such that $\deh\ll \de_{red}^{15}$, we may absorb the Morawetz terms on the LHS from the one on the RHS. We deduce
\beaa
\nn&& \de_{red}^{22}\Big(E_{r\leq r_+(1+\de_{red})}[\dk^{\leq 2}\psi](\tau_2)+\textrm{Mor}[(\dkb, \nab_3, \nab_4)^{\leq 2}\psi](\tau_1, \tau_2)+F_{\AA}[\dk^{\leq 2}\psi](\tau_1, \tau_2)\Big)\\ 
\nn&\les& E^2[\psi](\tau_1)+\NN^2[\psi, N] (\tau_1, \tau_2)\\
&&+\left(\frac{|a|}{m}+\ep\right)\left(\sup_{\tau\in [\tau_1, \tau_2]}E^2[\psi]+B_\de^2[\psi](\tau_1, \tau_2)+F^2[\psi](\tau_1, \tau_2)\right).
\eeaa
Together with  the energy estimate of Proposition \ref{prop:energyfromMorawetzfrominchap9}, and fixing the value of $\de_{red}$, we obtain 
\bea\lab{eq:iterationassumptiontruefors=2froproofThm:HigherDerivs-Morawetz-chp3}
\nn&& E[(\nab_3, \nab_4, \dkb)^{\leq 2}\psi](\tau_2)+\textrm{Mor}[(\dkb, \nab_3, \nab_4)^{\leq 2}\psi](\tau_1, \tau_2)+F[(\nab_3,  \nab_4, \dkb)^{\leq 2}\psi](\tau_1, \tau_2)\\ 
\nn&\les& E^2[\psi](\tau_1)+\NN^2[\psi, N] (\tau_1, \tau_2)\\
&&+\left(\frac{|a|}{m}+\ep\right)\left(\sup_{\tau\in [\tau_1, \tau_2]}E^2[\psi]+B_\de^2[\psi](\tau_1, \tau_2)+F^2[\psi](\tau_1, \tau_2)\right)
\eea
which concludes the proof of  Theorem \ref{THM:HIGHERDERIVS-MORAWETZ-CHP3} in the particular case $s=2$.
\end{proof}


\subsection{Proof of Theorem \ref{THM:HIGHERDERIVS-MORAWETZ-CHP3}}
\lab{sec:proofofThm:HigherDerivs-Morawetz-chp3:generalcase}


We are in position to concludes the proof of Theorem \ref{THM:HIGHERDERIVS-MORAWETZ-CHP3}. We consider the following iteration assumption for $2\leq j\leq s-1$:
\bea\lab{eq:iterationassumptiononnumberderivativesfroproofThm:HigherDerivs-Morawetz-chp3}
\nn&& E[(\nab_3, \nab_4, \dkb)^{\leq j}\psi](\tau_2)+\textrm{Mor}[(\dkb, \nab_3, \nab_4)^{\leq j}\psi](\tau_1, \tau_2)+F[(\nab_3,  \nab_4, \dkb)^{\leq j}\psi](\tau_1, \tau_2)\\ 
\nn&\les& E^j[\psi](\tau_1)+\NN^j[\psi, N] (\tau_1, \tau_2)\\
&&+\left(\frac{|a|}{m}+\ep\right)\left(\sup_{\tau\in [\tau_1, \tau_2]}E^j[\psi]+B_\de^j[\psi](\tau_1, \tau_2)+F^j[\psi](\tau_1, \tau_2)\right).
\eea
Note that \eqref{eq:iterationassumptiononnumberderivativesfroproofThm:HigherDerivs-Morawetz-chp3} holds for $j=2$ in view of \eqref{eq:iterationassumptiontruefors=2froproofThm:HigherDerivs-Morawetz-chp3}. We may thus assume that  \eqref{eq:iterationassumptiononnumberderivativesfroproofThm:HigherDerivs-Morawetz-chp3} holds  for some $2\leq j\leq s-1$, and our goal is to prove that this estimate  also holds for $j$ replaced by $j+1$. 

We start with the following Morawetz estimate.
\begin{proposition}\lab{prop:Morawetzfromenergyinchap9:higherderivatives}
Let $2\leq j\leq s-1$. Assume that the solution $\psi$ of of the model RW equation \eqref{eq:Gen.RW} satisfies the iteration assumption 
\eqref{eq:iterationassumptiononnumberderivativesfroproofThm:HigherDerivs-Morawetz-chp3}. Then, the following Morawetz estimate holds:
\bea
\nn&& E_{r\leq r_+(1+\de_{red})}[\dk^{\leq j+1}\psi](\tau_2)+\textrm{Mor}[(\dkb, \nab_3, \nab_4)^{\leq j+1}\psi](\tau_1, \tau_2)+F_{\AA}[\dk^{\leq j+1}\psi](\tau_1, \tau_2)\\ 
\nn&\les& \Bigg( E^j[\psi](\tau_1)+\NN^j[\psi, N] (\tau_1, \tau_2)\\
\nn&&+\left(\frac{|a|}{m}+\ep\right)\left(\sup_{\tau\in [\tau_1, \tau_2]}E^j[\psi]+B_\de^j[\psi](\tau_1, \tau_2)+F^j[\psi](\tau_1, \tau_2)\right)\Bigg)^{\frac{1}{2}}\\
\nn&&\times\Bigg(\sup_{\tau\in[\tau_1, \tau_2]}E_{deg}[(\dkb, \nab_T, \nab_{\Rhat})^{\leq j+1}\psi]\\
\nn&&+F_{\Si_*}[(\dkb, \nab_T, \nab_{\Rhat})^{\leq j+1}\psi](\tau_1, \tau_2)+\de_\HH F_\AA[(\dkb, \nab_T, \nab_{\Rhat})^{\leq j+1}\psi](\tau_1, \tau_2)\Bigg)^{\frac{1}{2}}\\
\nn&&+E^{j+1}[\psi](\tau_1)+\NN^{j+1}[\psi, N] (\tau_1, \tau_2)\\
&&+\left(\frac{|a|}{m}+\ep\right)\left(\sup_{\tau\in [\tau_1, \tau_2]}E^{j+1}[\psi]+B_\de^{j+1}[\psi](\tau_1, \tau_2)+F^{j+1}[\psi](\tau_1, \tau_2)\right).
\eea
\end{proposition}

\begin{proof}
We proceed in several steps. 

{\bf Step 1.} This step is analogous to  Step 2 of the proof of Proposition \ref{prop:energyforzeroandfirstderivativeinchap9}, so we only sketch it. We commute the wave equation \eqref{eq:Gen.RW} satisfied by $\psi$ with $(\Lied_\T, \Lied_\Z)$ and apply the iteration assumption \eqref{eq:iterationassumptiononnumberderivativesfroproofThm:HigherDerivs-Morawetz-chp3} to the commuted wave equation. This yields, comparing also $(\nab_T, \nab_Z)$ and $(\Lieb_\T, \Lieb_\Z)$, and using the iteration assumption \eqref{eq:iterationassumptiononnumberderivativesfroproofThm:HigherDerivs-Morawetz-chp3} to absorb lower order terms, 
\beaa
\nn&& E[(\nab_3, \nab_4, \dkb)^{\leq j}(\nab_T, \nab_Z)\psi](\tau_2)+\textrm{Mor}[(\dkb, \nab_3, \nab_4)^{\leq j}(\nab_T, \nab_Z)\psi](\tau_1, \tau_2)\\
&&+F[(\nab_3,  \nab_4, \dkb)^{\leq j}(\nab_T, \nab_Z)\psi](\tau_1, \tau_2)\\
\nn&\les& E^{j+1}[\psi](\tau_1)+\NN^{j+1}[\psi, N] (\tau_1, \tau_2)\\
&&+\left(\frac{|a|}{m}+\ep\right)\left(\sup_{\tau\in [\tau_1, \tau_2]}E^{j+1}[\psi]+B_\de^{j+1}[\psi](\tau_1, \tau_2)+F^{j+1}[\psi](\tau_1, \tau_2)\right).
\eeaa

{\bf Step 2.} Using again  the representation of the wave operator provided by \eqref{eq:q2squaredkpsi}, i.e. 
\beaa
\begin{split}
|q|^2 \squared_2 \psi &=\frac{(r^2+a^2)^2}{\De} \big( -  \nab_\That \nab_\That \psi+   \nab_\Rhat \nab_\Rhat \psi \big) +2r \nab_\Rhat \psi\\
&+  |q|^2 \lap_2 \psi   + |q|^2  (\eta+\etab) \c \nab \psi  + r^2 \Ga_g \c \dk \psi,
\end{split}
\eeaa
together with the wave equation \eqref{eq:Gen.RW} satisfied by $\psi$,  we infer 
\beaa
\textrm{Mor}[(\nab_3, \nab_4, \dkb)^{\leq j-1}U] &\les& \textrm{Mor}[(\nab_3, \nab_4, \dkb)^{\leq j-1}(\nab_T, \nab_Z)^{\leq 2}\psi](\tau)\\
&&+\textrm{Mor}[(\nab_3, \nab_4, \dkb)^{\leq j}\psi](\tau)+\int_{\Si(\tau)}|\dk^{\leq j-1}N|^2,
\eeaa
where the constant in $\les$ is independent of $\tau$, and where the tensor $U$ is given by
\beaa
U :=  \nab_\Rhat \nab_\Rhat \psi +  \frac{|q|^2\De}{(r^2+a^2)^2} \lap_2 \psi. 
\eeaa
Together with Step 1 and  the iteration assumption \eqref{eq:iterationassumptiononnumberderivativesfroproofThm:HigherDerivs-Morawetz-chp3}, we deduce
\beaa
&&\textrm{Mor}[(\nab_3, \nab_4, \dkb)^{\leq j-1}U] \\
&\les& E^{j+1}[\psi](\tau_1)+\NN^{j+1}[\psi, N] (\tau_1, \tau_2)\\
&&+\left(\frac{|a|}{m}+\ep\right)\left(\sup_{\tau\in [\tau_1, \tau_2]}E^{j+1}[\psi]+B_\de^{j+1}[\psi](\tau_1, \tau_2)+F^{j+1}[\psi](\tau_1, \tau_2)\right).
\eeaa
Next, using integration by parts and the iteration assumption \eqref{eq:iterationassumptiononnumberderivativesfroproofThm:HigherDerivs-Morawetz-chp3} to absorb lower order terms, we have
\beaa
&&\int_{\MM(\tau_1, \tau_2)}r^{-2}\Big(|\nab_\Rhat^3(\nab_3, \nab_4, \dkb)^{j-1}\psi|^2+\frac{|\De|}{r^2}|\nab_\Rhat^2\nab(\nab_3, \nab_4, \dkb)^{j-1}\psi|^2\\
&&+\frac{\De^2}{r^4}|\De_2\nab_\Rhat(\nab_3, \nab_4, \dkb)^{j-1}\psi|^2 \Big)\\
&\les& \textrm{Mor}[(\nab_3, \nab_4, \dkb)^{\leq j-1}U]\\
&&+\Bigg(\sup_{\tau\in[\tau_1, \tau_2]}E_{deg}[(\dkb, \nab_T, \nab_{\Rhat})^{\leq j}\psi]+F_{\Si_*}[(\dkb, \nab_T, \nab_{\Rhat})^{\leq j}\psi](\tau_1, \tau_2)\\
\nn&&+\de_\HH F_\AA[(\dkb, \nab_T, \nab_{\Rhat})^{\leq j}\psi](\tau_1, \tau_2)\Bigg)^{\frac{1}{2}}\Bigg(\sup_{\tau\in[\tau_1, \tau_2]}E_{deg}[(\dkb, \nab_T, \nab_{\Rhat})^{\leq j+1}\psi]\\
\nn&&+F_{\Si_*}[(\dkb, \nab_T, \nab_{\Rhat})^{\leq j+1}\psi](\tau_1, \tau_2)+\de_\HH F_\AA[(\dkb, \nab_T, \nab_{\Rhat})^{\leq j+1}\psi](\tau_1, \tau_2)\Bigg)^{\frac{1}{2}}.
\eeaa
Using again the iteration assumption \eqref{eq:iterationassumptiononnumberderivativesfroproofThm:HigherDerivs-Morawetz-chp3}, we infer
\beaa
&&\int_{\MM(\tau_1, \tau_2)}r^{-2}\Big(|\nab_\Rhat^3(\nab_3, \nab_4, \dkb)^{j-1}\psi|^2+\frac{|\De|}{r^2}|\nab_\Rhat^2\nab(\nab_3, \nab_4, \dkb)^{j-1}\psi|^2\\
&&+\frac{\De^2}{r^4}|\De_2\nab_\Rhat(\nab_3, \nab_4, \dkb)^{j-1}\psi|^2 \Big)\\
&\les& \textrm{Mor}[(\nab_3, \nab_4, \dkb)^{\leq j-1}U]\\
&&+\Bigg( E^j[\psi](\tau_1)+\NN^j[\psi, N] (\tau_1, \tau_2)\\
&&+\left(\frac{|a|}{m}+\ep\right)\left(\sup_{\tau\in [\tau_1, \tau_2]}E^j[\psi]+B_\de^j[\psi](\tau_1, \tau_2)+F^j[\psi](\tau_1, \tau_2)\right)\Bigg)^{\frac{1}{2}}\\
&&\times\Bigg(\sup_{\tau\in[\tau_1, \tau_2]}E_{deg}[(\dkb, \nab_T, \nab_{\Rhat})^{\leq j+1}\psi]\\
\nn&&+F_{\Si_*}[(\dkb, \nab_T, \nab_{\Rhat})^{\leq j+1}\psi](\tau_1, \tau_2)+\de_\HH F_\AA[(\dkb, \nab_T, \nab_{\Rhat})^{\leq j+1}\psi](\tau_1, \tau_2)\Bigg)^{\frac{1}{2}}.
\eeaa
Together with the above control of $\textrm{Mor}[(\nab_3, \nab_4, \dkb)^{\leq j-1}U]$, we infer
\beaa
\int_{\MM(\tau_1, \tau_2)}r^{-2}\Big(|\nab_\Rhat^3(\nab_3, \nab_4, \dkb)^{j-1}\psi|^2+\frac{|\De|}{r^2}|\nab_\Rhat^2\nab(\nab_3, \nab_4, \dkb)^{j-1}\psi|^2\\
+\frac{\De^2}{r^4}|\De_2\nab_\Rhat(\nab_3, \nab_4, \dkb)^{j-1}\psi|^2 \Big) &\les& \mathcal{J}_{RHS}
\eeaa
where we have introduced the following notation 
\beaa
\mathcal{J}_{RHS} &:=& \Bigg( E^j[\psi](\tau_1)+\NN^j[\psi, N] (\tau_1, \tau_2)\\
&&+\left(\frac{|a|}{m}+\ep\right)\left(\sup_{\tau\in [\tau_1, \tau_2]}E^j[\psi]+B_\de^j[\psi](\tau_1, \tau_2)+F^j[\psi](\tau_1, \tau_2)\right)\Bigg)^{\frac{1}{2}}\\
&&\times\Bigg(\sup_{\tau\in[\tau_1, \tau_2]}E_{deg}[(\dkb, \nab_T, \nab_{\Rhat})^{\leq j+1}\psi]\\
\nn&&+F_{\Si_*}[(\dkb, \nab_T, \nab_{\Rhat})^{\leq j+1}\psi](\tau_1, \tau_2)+\de_\HH F_\AA[(\dkb, \nab_T, \nab_{\Rhat})^{\leq j+1}\psi](\tau_1, \tau_2)\Bigg)^{\frac{1}{2}}\\
&&+E^{j+1}[\psi](\tau_1)+\NN^{j+1}[\psi, N] (\tau_1, \tau_2)\\
&&+\left(\frac{|a|}{m}+\ep\right)\left(\sup_{\tau\in [\tau_1, \tau_2]}E^{j+1}[\psi]+B_\de^{j+1}[\psi](\tau_1, \tau_2)+F^{j+1}[\psi](\tau_1, \tau_2)\right).
\eeaa
Using the Hodge estimates of Proposition \ref{Prop:HodgeThmM8}, as well as the estimate of Step 1 and the iteration assumption \eqref{eq:iterationassumptiononnumberderivativesfroproofThm:HigherDerivs-Morawetz-chp3},  
we deduce
\beaa
\int_{\MM(\tau_1, \tau_2)}r^{-2}\Big(|\nab_\Rhat^3(\nab_3, \nab_4, \dkb)^{j-1}\psi|^2+\frac{|\De|}{r^2}|\nab_\Rhat^2\nab(\nab_3, \nab_4, \dkb)^{j-1}\psi|^2\\
+\frac{\De^2}{r^4}|\nab_\Rhat\nab^2(\nab_3, \nab_4, \dkb)^{j-1}\psi|^2 \Big) &\les& \mathcal{J}_{RHS}.
\eeaa
Together with Step 1, we finally obtain 
\beaa
\int_{\MM(\tau_1, \tau_2)}r^{-6}\frac{\De^2}{r^4}|\nab_\Rhat(\nab_3, \nab_4, \dkb)^{j+1}\psi|^2  &\les& \mathcal{J}_{RHS}.
\eeaa

{\bf Step 3.} Using again  the representation of the wave operator provided by \eqref{eq:q2squaredkpsi}, we have, together with the wave equation \eqref{eq:Gen.RW} satisfied by $\psi$, 
\beaa
&&\int_{\Mntrap(\tau_1, \tau_2)}r^{-6}\frac{\De^4}{r^8}|\De_2(\nab_3, \nab_4, \dkb)^j\psi|^2+\int_{\MM(\tau_1, \tau_2)}r^{-6}\frac{\De^4}{r^8}|\De_2(\nab_3, \nab_4, \dkb)^{j-1}\psi|^2\\
&\les& \int_{\Mntrap(\tau_1, \tau_2)}r^{-6}\frac{\De^2}{r^4}|\nab_\Rhat^2(\nab_3, \nab_4, \dkb)^j\psi|^2 +\textrm{Mor}[(\nab_3, \nab_4, \dkb)^{\leq j-1}(\nab_T, \nab_Z)^{\leq 2}\psi](\tau)\\
&&+\textrm{Mor}[(\nab_3, \nab_4, \dkb)^{\leq j}\psi](\tau)+\int_{\Si(\tau)}|\dk^{\leq j-1}N|^2.
\eeaa
In view of Step 1 and Step 2, we infer
\beaa
\int_{\Mntrap(\tau_1, \tau_2)}r^{-6}\frac{\De^4}{r^8}|\De_2(\nab_3, \nab_4, \dkb)^j\psi|^2+\int_{\MM(\tau_1, \tau_2)}r^{-6}\frac{\De^4}{r^8}|\De_2(\nab_3, \nab_4, \dkb)^{j-1}\psi|^2 &\les& \mathcal{J}_{RHS}.
\eeaa
Using the Hodge estimates of Proposition \ref{Prop:HodgeThmM8}, as well as the estimate of Step 1 and the iteration assumption \eqref{eq:iterationassumptiononnumberderivativesfroproofThm:HigherDerivs-Morawetz-chp3},  
we deduce
\beaa
\int_{\Mntrap(\tau_1, \tau_2)}r^{-6}\frac{\De^4}{r^8}|\nab^2(\nab_3, \nab_4, \dkb)^j\psi|^2+\int_{\MM(\tau_1, \tau_2)}r^{-6}\frac{\De^4}{r^8}|\nab^2(\nab_3, \nab_4, \dkb)^{j-1}\psi|^2 &\les& \mathcal{J}_{RHS}.
\eeaa
Using again Step 1 and Step 2, we obtain 
\beaa
\int_{\MM(\tau_1, \tau_2)}r^{-6}\frac{\De^2}{r^4}|\nab_\Rhat(\nab_3, \nab_4, \dkb)^{j+1}\psi|^2\\
+\int_{\Mntrap(\tau_1, \tau_2)}r^{-8}\frac{\De^4}{r^8}|(\nab, \nab_T)(\nab_3, \nab_4, \dkb)^{j+1}\psi|^2\\
+\int_{\MM(\tau_1, \tau_2)}r^{-10}\frac{\De^4}{r^8}|(\nab_3, \nab_4, \dkb)^{j+1}\psi|^2 &\les& \mathcal{J}_{RHS}.
\eeaa

{\bf Step 4.} This step is analogous to  Step 2 of the proof of Proposition \ref{prop:energyfromMorawetzfrominchap9}, so we only sketch it. Commuting the wave equation \eqref{eq:Gen.RW} satisfied by $\psi$ with $\dkb_2^{j_1}\Lieb_\T^{j_2}$  with $j_1+j_2\leq j+1$, and relying on  Proposition \ref{prop:recoverEnergyMorawetzwithrweightfromnoweight:perturbation}, we infer, for $r\geq r_1$ with $r_1=r_1(m)$ sufficiently large, 
together with the Hodge estimates of Proposition \ref{Prop:HodgeThmM8}
\beaa
\textrm{Mor}_{r\geq 2r_1}[(\nab_T, \dkb)^{\leq j+1}\psi](\tau_1, \tau_2) &\les&  \frac{r_1}{m}\textrm{Mor}_{r_1\leq r\leq 2r_1}[(\nab_T, \dkb)^{\leq 2}\psi](\tau_1, \tau_2) +\mathcal{J}_{RHS},
\eeaa
where we used also the iteration assumption \eqref{eq:iterationassumptiononnumberderivativesfroproofThm:HigherDerivs-Morawetz-chp3} to absorb lower order terms coming from commutations with $\dkb_2$. Noticing that $\textrm{Mor}_{r_1\leq r\leq 2r_1}[(\nab_T, \dkb)^{\leq 2}\psi](\tau_1, \tau_2)$ is controlled by the LHS of the final estimate of Step 3, we infer
\beaa
\textrm{Mor}_{r\geq 2r_1}[(\nab_T, \dkb)^{\leq j+1}\psi](\tau_1, \tau_2) &\les&  \mathcal{J}_{RHS}.
\eeaa

Next, using again  the representation of the wave operator provided by \eqref{eq:q2squaredkpsi}, and arguing by iteration, we infer from the above estimate, for any $0\leq 2l\leq j+1$,  
\beaa
\int_{\MM_{r\geq 2r_1}(\tau_1, \tau_2)}r^{-2}|\nab_{\Rhat}^{1+2l}(\nab_T, \dkb)^{j+1-2l}\psi|^2 &\les&  \mathcal{J}_{RHS},
\eeaa
and for any $2\leq 2l\leq j+2$,  
\beaa
\int_{\MM_{r\geq 2r_1}(\tau_1, \tau_2)}r^{-2}|\nab_{\Rhat}^{2l}(\nab_T, \dkb)^{j+2-2l}\psi|^2 &\les&  \mathcal{J}_{RHS},
\eeaa
and hence, together with the above estimate
\beaa
\textrm{Mor}_{r\geq 2r_1}[(\nab_T, \dkb, \nab_{\Rhat})^{\leq j+1}\psi](\tau_1, \tau_2) &\les&  \mathcal{J}_{RHS}.
\eeaa
Together with Step 3, we deduce
\beaa
\int_{\Mntrap(\tau_1, \tau_2)}\Big(r^{-1}\frac{\De^4}{r^8}|\nab(\nab_3, \nab_4, \dkb)^{j+1}\psi|^2+r^{-2}\frac{\De^4}{r^8}|\nab_T(\nab_3, \nab_4, \dkb)^{j+1}\psi|^2\Big)\\
+\int_{\MM(\tau_1, \tau_2)}r^{-2}\frac{\De^2}{r^4}|\nab_\Rhat(\nab_3, \nab_4, \dkb)^{j+1}\psi|^2+\int_{\MM(\tau_1, \tau_2)}r^{-3}\frac{\De^4}{r^8}|(\nab_3, \nab_4, \dkb)^{j+1}\psi|^2 &\les& \mathcal{J}_{RHS}.
\eeaa

{\bf Step 5.} Next, we commute  the wave equation \eqref{eq:Gen.RW} for $\psi$ in $r\leq 4m$ by $(\nab_3, \dkb_2, \Lieb_\T)^{\leq j+1}$, using the commutation formulas of Corollary \ref{cor:commutator-Lied-squared}, Lemma \ref{lemma:commutator-triangle} and Lemma \ref{LEMMA:COMMUTATIONOFHODGEELLIPTICORDER1WITHSQAURED2FDILUHS}, as well as Lemma \ref{lemma:commutationwithe3forredshift} for the commutations w.r.t to $\nab_3$. Together with Corollary \ref{cor:Redshift-estimates-chp3}, we infer
\beaa
\bsplit
&c_0E_{r\leq r_+(1+\de_{red})}[(\nab_3, \dkb_2, \Lieb_\T)^{\leq j+1}\psi](\tau_2)+c_0\textrm{Mor}_{r\leq r_+(1+\de_{red})}[(\nab_3, \dkb_2, \Lieb_\T)^{\leq j+1}\psi](\tau_1, \tau_2)\\
&+c_0F_{\AA}[(\nab_3, \dkb_2, \Lieb_\T)^{\leq j+1}\psi](\tau_1, \tau_2)\\
  \leq& \de_{red}^{-1}\textrm{Mor}_{r_+(1+\de_{red})\leq r\leq r_+(1+2\de_{red})}[(\nab_3, \dkb, \Lieb_\T)^{\leq j+1}\psi](\tau_1, \tau_2)+\mathcal{J}_{RHS}.
\end{split}
\eeaa 
Using the comparison of $\Lieb_\T$ and $\nab_\T$ of Lemma \ref{lemma:basicpropertiesLiebTfasdiuhakdisug:chap9}, and  the Hodge estimates of Proposition \ref{Prop:HodgeThmM8}, and using 
\beaa
&&\de_{red}^{-1}\textrm{Mor}_{r_+(1+\de_{red})\leq r\leq r_+(1+2\de_{red})}[(\nab_3, \dkb, \Lieb_\T)^{\leq j+1}\psi](\tau_1, \tau_2)\\
&\les& \de_{red}^{-7}\Bigg(\int_{\Mntrap(\tau_1, \tau_2)}\Big(r^{-1}\frac{\De^4}{r^8}|\nab(\nab_3, \nab_4, \dkb)^{j+1}\psi|^2+r^{-2}\frac{\De^4}{r^8}|\nab_T(\nab_3, \nab_4, \dkb)^{j+1}\psi|^2\Big)\\
&&+\int_{\MM(\tau_1, \tau_2)}r^{-2}\frac{\De^2}{r^4}|\nab_\Rhat(\nab_3, \nab_4, \dkb)^{j+1}\psi|^2+\int_{\MM(\tau_1, \tau_2)}r^{-3}\frac{\De^4}{r^8}|(\nab_3, \nab_4, \dkb)^{j+1}\psi|^2\Bigg),
\eeaa
we infer from Step 4
\beaa
\bsplit
E_{r\leq r_+(1+\de_{red})}[(\nab_3, \dkb, \nab_4)^{\leq j+1}\psi](\tau_2)+\textrm{Mor}_{r\leq r_+(1+\de_{red})}[(\nab_3, \dkb, \nab_4)^{\leq j+1}\psi](\tau_1, \tau_2)\\
+F_{\AA}[(\nab_3, \dkb, \nab_4)^{\leq j+1}\psi](\tau_1, \tau_2)  \les& \mathcal{J}_{RHS}.
\end{split}
\eeaa 
Together with the control of Step 4, we deduce
\beaa
\nn&& E_{r\leq r_+(1+\de_{red})}[\dk^{\leq j+1}\psi](\tau_2)+\textrm{Mor}[(\dkb, \nab_3, \nab_4)^{\leq j+1}\psi](\tau_1, \tau_2)+F_{\AA}[\dk^{\leq j+1}\psi](\tau_1, \tau_2)\\ 
\nn&\les& \mathcal{J}_{RHS}
\eeaa
which in view of the definition of $\mathcal{J}_{RHS}$ yields 
\beaa
\nn&& E_{r\leq r_+(1+\de_{red})}[\dk^{\leq j+1}\psi](\tau_2)+\textrm{Mor}[(\dkb, \nab_3, \nab_4)^{\leq j+1}\psi](\tau_1, \tau_2)+F_{\AA}[\dk^{\leq j+1}\psi](\tau_1, \tau_2)\\ 
\nn&\les& \Bigg( E^j[\psi](\tau_1)+\NN^j[\psi, N] (\tau_1, \tau_2)\\
&&+\left(\frac{|a|}{m}+\ep\right)\left(\sup_{\tau\in [\tau_1, \tau_2]}E^j[\psi]+B_\de^j[\psi](\tau_1, \tau_2)+F^j[\psi](\tau_1, \tau_2)\right)\Bigg)^{\frac{1}{2}}\\
&&\times\Bigg(\sup_{\tau\in[\tau_1, \tau_2]}E_{deg}[(\dkb, \nab_T, \nab_{\Rhat})^{\leq j+1}\psi]\\
\nn&&+F_{\Si_*}[(\dkb, \nab_T, \nab_{\Rhat})^{\leq j+1}\psi](\tau_1, \tau_2)+\de_\HH F_\AA[(\dkb, \nab_T, \nab_{\Rhat})^{\leq j+1}\psi](\tau_1, \tau_2)\Bigg)^{\frac{1}{2}}\\
&&+E^{j+1}[\psi](\tau_1)+\NN^{j+1}[\psi, N] (\tau_1, \tau_2)\\
&&+\left(\frac{|a|}{m}+\ep\right)\left(\sup_{\tau\in [\tau_1, \tau_2]}E^{j+1}[\psi]+B_\de^{j+1}[\psi](\tau_1, \tau_2)+F^{j+1}[\psi](\tau_1, \tau_2)\right)
\eeaa
as stated. This concludes the proof of Proposition \ref{prop:Morawetzfromenergyinchap9:higherderivatives}. 
\end{proof}

Next, we consider the following combined energy-Morawetz estimate.
\begin{proposition}\lab{prop:combinedenergyMorawetzfromenergyinchap9:higherderivatives}
Let $2\leq j\leq s-1$. Assume that the solution $\psi$ of of the model RW equation \eqref{eq:Gen.RW} satisfies the iteration assumption 
\eqref{eq:iterationassumptiononnumberderivativesfroproofThm:HigherDerivs-Morawetz-chp3}. Then, the following combined energy-Morawetz estimate holds:
\bea
\nn&& E[(\nab_3, \nab_4, \dkb)^{\leq j+1}\psi](\tau_2)+\textrm{Mor}[(\dkb, \nab_3, \nab_4)^{\leq j+1}\psi](\tau_1, \tau_2)+F[(\nab_3,  \nab_4, \dkb)^{\leq j+1}\psi](\tau_1, \tau_2)\\ 
\nn&\les& E^{j+1}[\psi](\tau_1)+\NN^{j+1}[\psi, N] (\tau_1, \tau_2)\\
&&+\left(\frac{|a|}{m}+\ep\right)\left(\sup_{\tau\in [\tau_1, \tau_2]}E^{j+1}[\psi]+B_\de^{j+1}[\psi](\tau_1, \tau_2)+F^{j+1}[\psi](\tau_1, \tau_2)\right).
\eea
\end{proposition}

\begin{proof}
Since the proof is similar, and in fact easier, than the one of \eqref{eq:iterationassumptiontruefors=2froproofThm:HigherDerivs-Morawetz-chp3}, we only sketch it:
\begin{enumerate}
\item First, we commute the wave equation \eqref{eq:Gen.RW} satisfied by $\psi$ with $(\Lieb_\T, \Lieb_\Z)$ and apply the iteration assumption  \ref{eq:iterationassumptiononnumberderivativesfroproofThm:HigherDerivs-Morawetz-chp3}  to the commuted wave equation. This yields, comparing also $(\nab_T, \nab_Z)$ and $(\Lieb_\T, \Lieb_\Z)$, and using the iteration assumption \eqref{eq:iterationassumptiononnumberderivativesfroproofThm:HigherDerivs-Morawetz-chp3} to absorb lower order terms, 
\beaa
\nn&& E[(\nab_3, \nab_4, \dkb)^{\leq j}(\nab_T, \nab_Z)\psi](\tau_2)+\textrm{Mor}[(\dkb, \nab_3, \nab_4)^{\leq j}(\nab_T, \nab_Z)\psi](\tau_1, \tau_2)\\
&&+F[(\nab_3,  \nab_4, \dkb)^{\leq j}(\nab_T, \nab_Z)\psi](\tau_1, \tau_2)\\ 
\nn&\les& E^{j+1}[\psi](\tau_1)+\NN^{j+1}[\psi, N] (\tau_1, \tau_2)\\
&&+\left(\frac{|a|}{m}+\ep\right)\left(\sup_{\tau\in [\tau_1, \tau_2]}E^{j+1}[\psi]+B_\de^{j+1}[\psi](\tau_1, \tau_2)+F^{j+1}[\psi](\tau_1, \tau_2)\right).
\eeaa

\item Next, we commute the wave equation \eqref{eq:Gen.RW} satisfied by $\psi$ with $(\Lieb_\T, \Lieb_\Z)^{j+1-j}\dkb_2^k$ for $1\leq k\leq j+1$, where we recall that 
\beaa
\dkb_2^{2k}=(|q|^2\De_2)^k, \qquad \dkb_2^{2k+1}=|q|\DD_2(|q|^2\De_2)^k.
\eeaa
Denoting 
\beaa
\psi_{k,j}=\Lieb_\T^{j+1-k}\dkb_2^k\psi,
\eeaa
and using the commutation formulas of using the commutation formulas of Corollary \ref{cor:commutator-Lied-squared}, Lemma \ref{lemma:commutator-triangle} and Lemma \ref{LEMMA:COMMUTATIONOFHODGEELLIPTICORDER1WITHSQAURED2FDILUHS}, we obtain
\beaa 
\bsplit
\square_2\psi_{2k,j} &= N_{2k,j}, \qquad\, N_{2k,j}=\dk^{\leq j+1}N+O(ar^{-2})\dk^{j+2}\psi+\dk^{\leq j+2}(\Ga_g\c\psi),\\
\square_1\psi_{2k+1,j} &= N_{2k+1,j}, \,\, N_{2k+1,j}=\dk^{\leq j+1}N-\frac{3}{r^2}\psi_{2k+1,j}+O(ar^{-2})\dk^{j+2}\psi+\dk^{\leq j+2}(\Ga_g\c\psi).
\end{split}
\eeaa
Relying on Proposition \ref{prop:recoverEnergyMorawetzwithrweightfromnoweight:perturbation}, we infer, for $r\geq r_1$ with $r_1=r_1(m)$ sufficiently large, 
\beaa
&&E_{r\geq 2r_1}[\psi_{k,j}](\tau_2)+F_{\Si_*}[\psi_{k,j}](\tau_1, \tau_2)\\
&\les& \textrm{Mor}_{r_1\leq r\leq 2r_1}[(\nab_\T, \dkb)^{\leq j+1}\psi](\tau_1, \tau_2) + E^{j+1}[\psi](\tau_1)+\NN^{j+1}[\psi, N] (\tau_1, \tau_2)\\
&&+\left(\frac{|a|}{m}+\ep\right)\left(\sup_{\tau\in [\tau_1, \tau_2]}E^{j+1}[\psi]+B_\de^{j+1}[\psi](\tau_1, \tau_2)+F^{j+1}[\psi](\tau_1, \tau_2)\right).
\eeaa
In view of the definition of $\psi_{k,j}$, together with the Hodge estimates of Proposition \ref{Prop:HodgeThmM8} and the iteration assumption \ref{eq:iterationassumptiononnumberderivativesfroproofThm:HigherDerivs-Morawetz-chp3}, we obtain
\beaa
&&E_{r\geq 2r_1}[(\nab_T, \dkb)^{\leq j+1}](\tau_2)+F_{\Si_*}[(\nab_T, \dkb)^{\leq j+1}](\tau_1, \tau_2)\\
&\les& \textrm{Mor}_{r_1\leq r\leq 2r_1}[(\nab_\T, \dkb)^{\leq j+1}\psi](\tau_1, \tau_2) + E^{j+1}[\psi](\tau_1)+\NN^{j+1}[\psi, N] (\tau_1, \tau_2)\\
&&+\left(\frac{|a|}{m}+\ep\right)\left(\sup_{\tau\in [\tau_1, \tau_2]}E^{j+1}[\psi]+B_\de^{j+1}[\psi](\tau_1, \tau_2)+F^{j+1}[\psi](\tau_1, \tau_2)\right).
\eeaa
Also, using again  the representation of the wave operator provided by \eqref{eq:q2squaredkpsi}, we have, for $2k\leq j+1$,
\beaa
&& E_{r\geq 2r_1}[(\nab_T, \dkb)^{\leq j+1-2k}\nab_{\Rhat}^{2k}\psi]+F_{\Si_*}[(\nab_T, \dkb)^{\leq j+1-2k}\nab_{\Rhat}^{2k}\psi]\\
&\les& E_{r\geq 2r_1}[(\nab_T, \dkb)^{\leq j+1}\psi](\tau_2)+F_{\Si_*}[(\nab_T, \dkb)^{\leq j+1}\psi](\tau_1, \tau_2)\\
&&+\NN^{j+1}[\psi, N] (\tau_1, \tau_2)+\ep\left(\sup_{\tau\in [\tau_1, \tau_2]}E^{j+1}[\psi]+F^{j+1}[\psi](\tau_1, \tau_2)\right).
\eeaa
and  for $2k\leq j+2$,
\beaa
&& \int_{\Si_{r\geq 2r_1(\tau_2)}}|\nab_{\Rhat}^{2k}\psi|^2+\int_{\Si_*(\tau_1, \tau_2)}|\nab_{\Rhat}^{2k}\psi|^2\\
&\les& E_{r\geq 2r_1}[(\nab_T, \dkb)^{\leq j+1}\psi](\tau_2)+F_{\Si_*}[(\nab_T, \dkb)^{\leq j+1}\psi](\tau_1, \tau_2)\\
&&+\NN^{j+1}[\psi, N] (\tau_1, \tau_2)+\ep\left(\sup_{\tau\in [\tau_1, \tau_2]}E^{j+1}[\psi]+F^{j+1}[\psi](\tau_1, \tau_2)\right).
\eeaa
We infer from the above
\beaa
&&E_{r\geq 2r_1}[(\nab_T, \dkb, \nab_{\Rhat})^{\leq j+1}\psi](\tau_2)+F_{\Si_*}[(\nab_T, \dkb, \nab_{\Rhat})^{\leq j+1}\psi](\tau_1, \tau_2)\\
&\les&   \textrm{Mor}_{r_1\leq r\leq 2r_1}[(\nab_\T, \dkb)^{\leq j+1}\psi](\tau_1, \tau_2) + E^{j+1}[\psi](\tau_1)+\NN^{j+1}[\psi, N] (\tau_1, \tau_2)\\
&&+\left(\frac{|a|}{m}+\ep\right)\left(\sup_{\tau\in [\tau_1, \tau_2]}E^{j+1}[\psi]+B_\de^{j+1}[\psi](\tau_1, \tau_2)+F^{j+1}[\psi](\tau_1, \tau_2)\right).
\eeaa
\item Using again  the representation of the wave operator provided by \eqref{eq:q2squaredkpsi}, i.e. 
\beaa
\begin{split}
|q|^2 \squared_2 \psi &=\frac{(r^2+a^2)^2}{\De} \big( -  \nab_\That \nab_\That \psi+   \nab_\Rhat \nab_\Rhat \psi \big) +2r \nab_\Rhat \psi\\
&+  |q|^2 \lap_2 \psi   + |q|^2  (\eta+\etab) \c \nab \psi  + r^2 \Ga_g \c \dk \psi,
\end{split}
\eeaa
together with the wave equation \eqref{eq:Gen.RW} satisfied by $\psi$,  we infer
\beaa
E[(\nab_3, \nab_4, \dkb)^{\leq j-1}U] &\les& E[(\nab_3, \nab_4, \dkb)^{\leq j-1}(\nab_T, \nab_Z)^{\leq 2}\psi](\tau)\\
&&+E[(\nab_3, \nab_4, \dkb)^{\leq j}\psi](\tau)+\int_{\Si(\tau)}|\dk^{\leq j-1}N|^2,
\eeaa
where the constant in $\les$ is independent of $\tau$, and where the tensor $U$ is given by
\beaa
U :=  \nab_\Rhat \nab_\Rhat \psi +  \frac{|q|^2\De}{(r^2+a^2)^2} \lap_2 \psi. 
\eeaa
Next, we proceed as follows:
\begin{enumerate}
\item We consider a cut-off $\ka_{\deh}(r)$ such that $0\leq\ka_{\deh}\leq 1$, $\ka_{\deh}=1$ on $r\geq r_+$, and $\ka_{\deh}$ vanishes on $\AA$. Also, for any integer $2\leq n\leq \tau_*-2$, we consider a cut-off $\ka_n(\tau)$ such that $0\leq \ka_n\leq 1$, $\ka_n(\tau)=1$ on $\{n\leq \tau\leq n+1\}$, and $\ka_n$ vanishes on $\tau\leq n-1$ and $\tau\geq n+2$. We have, by  integrating by parts, and using  the support properties of $\ka_n(\tau)$ and the above control of $U$, 
\beaa
&&\int_{\MM}\ka_{\deh}(r)\ka_n(\tau)|(\nab_T, \nab, \nab_{\Rhat})(\nab_3, \nab_4, \dkb)^{\leq j-1}\nab_{\Rhat}^2\psi|^2\\
&&+\int_{\MM}\ka_{\deh}(r)\ka_n(\tau)\frac{|\De|}{r^2}|(\nab_T, \nab, \nab_{\Rhat})(\nab_3, \nab_4, \dkb)^{\leq j-1}\nab\nab_{\Rhat}\psi|^2\\
&&+\int_{\MM}\ka_{\deh}(r)\ka_n(\tau)\frac{|\De|^2}{r^4}|(\nab_T, \nab, \nab_{\Rhat})(\nab_3, \nab_4, \dkb)^{\leq j-1}\nab^2\psi|^2\\
&\les& \sup_{\tau\in[n-1,n+2]}E[(\nab_3, \nab_4, \dkb)^{\leq j-1}U]\\
&&+\left(\sup_{\tau\in[n-1,n+2]}E[(\nab_3, \nab_4, \dkb)^{\leq j}\psi](\tau)+F_{\Si_*}[(\nab_3, \nab_4, \dkb)^{\leq j}\psi](n-1, n+2)\right)^{\frac{1}{2}}\\
&&\times\left(\sup_{\tau\in[n-1,n+2]}E[(\nab_T, \nab, \nab_{\Rhat})^{\leq j+1}\psi](\tau)+F_{\Si_*}[(\nab_T, \nab_{\Rhat}, \nab)^{\leq j+1}\psi](n-1, n+2)\right)^{\frac{1}{2}}.
\eeaa
Together with the above estimate for $U$, we deduce 
\beaa
&&\int_{\MM}\ka_{\deh}(r)\ka_n(\tau)|(\nab_T, \nab, \nab_{\Rhat})(\nab_3, \nab_4, \dkb)^{\leq j-1}\nab_{\Rhat}^2\psi|^2\\
&&+\int_{\MM}\ka_{\deh}(r)\ka_n(\tau)\frac{|\De|}{r^2}|(\nab_T, \nab, \nab_{\Rhat})(\nab_3, \nab_4, \dkb)^{\leq j-1}\nab\nab_{\Rhat}\psi|^2\\
&&+\int_{\MM}\ka_{\deh}(r)\ka_n(\tau)\frac{|\De|^2}{r^4}|(\nab_T, \nab, \nab_{\Rhat})(\nab_3, \nab_4, \dkb)^{\leq j-1}\nab^2\psi|^2\\
&\les& \sup_{\tau\in[n-1,n+2]}\Big[E[(\nab_3, \nab_4, \dkb)^{\leq j-1}(\nab_T, \nab_Z)^{\leq 2}\psi](\tau)+E[(\nab_3, \nab_4, \dkb)^{\leq j}\psi](\tau)\Big]\\
&& +\left(\sup_{\tau\in[n-1,n+2]}E[(\nab_3, \nab_4, \dkb)^{\leq j}\psi](\tau)+F_{\Si_*}[(\nab_3, \nab_4, \dkb)^{\leq j}\psi](n-1, n+2)\right)^{\frac{1}{2}}\\
&&\times\left(\sup_{\tau\in[n-1,n+2]}E[(\nab_T, \nab, \nab_{\Rhat})^{\leq j+1}\psi](\tau)+F_{\Si_*}[(\nab_T, \nab_{\Rhat}, \nab)^{\leq j+1}\psi](n-1, n+2)\right)^{\frac{1}{2}}\\
&&+\NN^{j+1}[\psi, N] (n-1, n+2).
\eeaa
\item In view of the properties of the cut-offs, the above estimates imply 
\beaa
&&\int_{\MM\cap\{r\geq r_+\}\cap\{n\leq\tau\leq n+1\}}\Bigg[|(\nab_T, \nab, \nab_{\Rhat})(\nab_3, \nab_4, \dkb)^{\leq j-1}\nab_{\Rhat}^2\psi|^2\\
&&+\frac{|\De|}{r^2}|(\nab_T, \nab, \nab_{\Rhat})(\nab_3, \nab_4, \dkb)^{\leq j-1}\nab\nab_{\Rhat}\psi|^2\\
&&+\frac{|\De|^2}{r^4}|(\nab_T, \nab, \nab_{\Rhat})(\nab_3, \nab_4, \dkb)^{\leq j-1}\nab^2\psi|^2\Bigg]\\
&\les& \sup_{\tau\in[n-1,n+2]}\Big[E[(\nab_3, \nab_4, \dkb)^{\leq j-1}(\nab_T, \nab_Z)^{\leq 2}\psi](\tau)+E[(\nab_3, \nab_4, \dkb)^{\leq j}\psi](\tau)\Big]\\
&& +\left(\sup_{\tau\in[n-1,n+2]}E[(\nab_3, \nab_4, \dkb)^{\leq j}\psi](\tau)+F_{\Si_*}[(\nab_3, \nab_4, \dkb)^{\leq j}\psi](n-1, n+2)\right)^{\frac{1}{2}}\\
&&\times\left(\sup_{\tau\in[n-1,n+2]}E[(\nab_T, \nab, \nab_{\Rhat})^{\leq j+1}\psi](\tau)+F_{\Si_*}[(\nab_T, \nab_{\Rhat}, \nab)^{\leq j+1}\psi](n-1, n+2)\right)^{\frac{1}{2}}\\
&&+\NN^{j+1}[\psi, N] (n-1, n+2).
\eeaa
By the mean value theorem, we infer the existence of $\tau_{(n)}\in[n,n+1]$ such that 
\beaa
&&\int_{\Si(\tau_{(n)})\cap\{r\geq r_+\}}\Bigg[|\nab_4\nab_{\Rhat}^2\psi|^2+\frac{|\De|}{r^2}|\nab_4\nab\nab_{\Rhat}\psi|^2+\frac{|\De|^2}{r^4}|\nab_4\nab^2\psi|^2\\
&&+r^{-2}\left(|\nab_{\That}\nab_{\Rhat}^2\psi|^2+\frac{|\De|}{r^2}|\nab_{\That}\nab\nab_{\Rhat}\psi|^2+\frac{|\De|^2}{r^4}|\nab_{\That}\nab^2\psi|^2\right)+|\nab\nab_{\Rhat}^2\psi|^2\\
&&+\frac{|\De|}{r^2}|\nab\nab_{\Rhat}\nab\psi|^2+\frac{|\De|^2}{r^4}|\nab^3\psi|^2\Bigg]\\
&\les&  \sup_{\tau\in[n-1,n+2]}\Big[E[(\nab_3, \nab_4, \dkb)^{\leq j-1}(\nab_T, \nab_Z)^{\leq 2}\psi](\tau)+E[(\nab_3, \nab_4, \dkb)^{\leq j}\psi](\tau)\Big]\\
&& +\left(\sup_{\tau\in[n-1,n+2]}E[(\nab_3, \nab_4, \dkb)^{\leq j}\psi](\tau)+F_{\Si_*}[(\nab_3, \nab_4, \dkb)^{\leq j}\psi](n-1, n+2)\right)^{\frac{1}{2}}\\
&&\times\left(\sup_{\tau\in[n-1,n+2]}E[(\nab_T, \nab, \nab_{\Rhat})^{\leq j+1}\psi](\tau)+F_{\Si_*}[(\nab_T, \nab_{\Rhat}, \nab)^{\leq j+1}\psi](n-1, n+2)\right)^{\frac{1}{2}}\\
&&+\NN^{j+1}[\psi, N] (n-1, n+2).
\eeaa
\item Together with the above control of $E[(\nab_3, \nab_4, \dkb)^{\leq j}(\nab_T, \nab_Z)\psi](\tau)$, the above control of $E_{r\geq 2r_1}[(\nab_T, \dkb, \nab_{\Rhat})^{\leq j+1}\psi](\tau_2)$ and $F_{\Si_*}[(\nab_T, \dkb, \nab_{\Rhat})^{\leq j+1}\psi](\tau_1, \tau_2)$, the control of $E_{r\leq r_+(1+\de_{red})}[\dk^{\leq j+1}\psi](\tau_2)$ and $F_{\AA}[\dk^{\leq j+1}\psi](\tau_1, \tau_2)$ provided by Proposition \ref{prop:Morawetzfromenergyinchap9:higherderivatives}, and  the iteration assumption \ref{eq:iterationassumptiononnumberderivativesfroproofThm:HigherDerivs-Morawetz-chp3}, and fixing the value of $r_1=r_1(m)$, we infer, the following estimate, for any $n$ such that $\tau_{(n)}\geq \tau_1$,  
\beaa
&& F_{\AA}[\dk^{\leq j+1}\psi](\tau_1, \tau_{(n)})+E[(\nab_3, \nab_4, \dkb)^{\leq j+1}\psi](\tau_{(n)})+F_{\Si_*}[(\nab_3,  \nab_4, \dkb)^{\leq j+1}\psi](\tau_1, \tau_{(n)})\\
&\les&  \Bigg(E^{j}[\psi](\tau_1)+\NN^{j}[\psi, N] (\tau_1, n+2)\\
&&+\left(\frac{|a|}{m}+\ep\right)\left(\sup_{\tau\in [\tau_1, n+2]}E^{j}[\psi]+B_\de^{j}[\psi](\tau_1, n+2)+F^{j}[\psi](\tau_1, n+2)\right)\Bigg)^{\frac{1}{2}}\\
&&\times\left(\sup_{\tau\in[n-1,n+2]}E[(\nab_T, \nab, \nab_{\Rhat})^{\leq j+1}\psi](\tau)+F_{\Si_*}[(\nab_T, \nab_{\Rhat}, \nab)^{\leq j+1}\psi](n-1, n+2)\right)^{\frac{1}{2}}\\
&&+\textrm{Mor}[(\nab_\T, \dkb)^{\leq j+1}\psi](\tau_1, n+2)+E^{j+1}[\psi](\tau_1)+\NN^{j+1}[\psi, N] (\tau_1, n+2)\\
&&+\left(\frac{|a|}{m}+\ep\right)\left(\sup_{\tau\in [\tau_1, n+2]}E^{j+1}[\psi]+B_\de^{j+1}[\psi](\tau_1, n+2)+F^{j+1}[\psi](\tau_1, n+2)\right).
\eeaa
\item Let $\tau_2\geq \tau_1$. By local energy estimates, it suffices to consider the case $\tau_2\geq \tau_1+5$. We then choose $n$ such that $\tau_1\leq n-1\leq n+2\leq\tau_2<n+3$. In particular, we have $\tau_{(n)}+1\leq \tau_2\leq \tau_{(n)}+3$, and hence, using local energy estimates between $\tau_{(n)}$ and $\tau_2$, we infer from the previous estimate
\beaa
&& F_{\AA}[\dk^{\leq j+1}\psi](\tau_1, \tau_2)+E[(\nab_3, \nab_4, \dkb)^{\leq j+1}\psi](\tau_2)+F_{\Si_*}[(\nab_3,  \nab_4, \dkb)^{\leq j+1}\psi](\tau_1, \tau_2)\\
&\les&   \Bigg(E^{j}[\psi](\tau_1)+\NN^{j}[\psi, N] (\tau_1, \tau_2)\\
&&+\left(\frac{|a|}{m}+\ep\right)\left(\sup_{\tau\in [\tau_1, \tau_2]}E^{j}[\psi]+B_\de^{j}[\psi](\tau_1, \tau_2)+F^{j}[\psi](\tau_1, \tau_2)\right)\Bigg)^{\frac{1}{2}}\\
&&\times\left(\sup_{\tau\in[\tau_1,\tau_2]}E[(\nab_T, \nab, \nab_{\Rhat})^{\leq j+1}\psi](\tau)+F_{\Si_*}[(\nab_T, \nab_{\Rhat}, \nab)^{\leq j+1}\psi](\tau_1, \tau_2)\right)^{\frac{1}{2}}\\
&&+\textrm{Mor}[(\nab_\T, \dkb)^{\leq j+1}\psi](\tau_1, \tau_2)+E^{j+1}[\psi](\tau_1)+\NN^{j+1}[\psi, N] (\tau_1, \tau_2)\\
&&+\left(\frac{|a|}{m}+\ep\right)\left(\sup_{\tau\in [\tau_1, \tau_2]}E^{j+1}[\psi]+B_\de^{j+1}[\psi](\tau_1, \tau_2)+F^{j+1}[\psi](\tau_1, \tau_2)\right).
\eeaa
and thus
\beaa
&& F_{\AA}[\dk^{\leq j+1}\psi](\tau_1, \tau_2)+E[(\nab_3, \nab_4, \dkb)^{\leq j+1}\psi](\tau_2)+F_{\Si_*}[(\nab_3,  \nab_4, \dkb)^{\leq j+1}\psi](\tau_1, \tau_2)\\
&\les&   \textrm{Mor}[(\nab_\T, \dkb)^{\leq j+1}\psi](\tau_1, \tau_2)+E^{j+1}[\psi](\tau_1)+\NN^{j+1}[\psi, N] (\tau_1, \tau_2)\\
&&+\left(\frac{|a|}{m}+\ep\right)\left(\sup_{\tau\in [\tau_1, \tau_2]}E^{j+1}[\psi]+B_\de^{j+1}[\psi](\tau_1, \tau_2)+F^{j+1}[\psi](\tau_1, \tau_2)\right).
\eeaa
\end{enumerate}
\item Together with the control for $\textrm{Mor}[(\dkb, \nab_3, \nab_4)^{\leq j+1}\psi]$ obtained in Proposition \ref{prop:Morawetzfromenergyinchap9:higherderivatives}, we infer
\beaa
\nn&& E[(\nab_3, \nab_4, \dkb)^{\leq j+1}\psi](\tau_2)+\textrm{Mor}[(\dkb, \nab_3, \nab_4)^{\leq j+1}\psi](\tau_1, \tau_2)+F[(\nab_3,  \nab_4, \dkb)^{\leq j+1}\psi](\tau_1, \tau_2)\\ 
\nn&\les& \left(E^j[\psi](\tau_1)+\NN^j[\psi, N] (\tau_1, \tau_2)+\left(\frac{|a|}{m}+\ep\right)\left(\sup_{\tau\in [\tau_1, \tau_2]}E^j[\psi]+B_\de^j[\psi](\tau_1, \tau_2)\right)\right)^{\frac{1}{2}}\\
\nn&&\times\Bigg(\sup_{\tau\in[\tau_1, \tau_2]}E_{deg}[(\dkb, \nab_T, \nab_{\Rhat})^{\leq j+1}\psi]+F_{\Si_*}[(\dkb, \nab_T, \nab_{\Rhat})^{\leq j+1}\psi](\tau_1, \tau_2)\\
\nn&&+\de_\HH F_\AA[(\dkb, \nab_T, \nab_{\Rhat})^{\leq j+1}\psi](\tau_1, \tau_2)\Bigg)^{\frac{1}{2}}+E^{j+1}[\psi](\tau_1)+\NN^{j+1}[\psi, N] (\tau_1, \tau_2)\\
&&+\left(\frac{|a|}{m}+\ep\right)\left(\sup_{\tau\in [\tau_1, \tau_2]}E^{j+1}[\psi]+B_\de^{j+1}[\psi](\tau_1, \tau_2)+F^{j+1}[\psi](\tau_1, \tau_2)\right).
\eeaa
We deduce
\beaa
\nn&& E[(\nab_3, \nab_4, \dkb)^{\leq j+1}\psi](\tau_2)+\textrm{Mor}[(\dkb, \nab_3, \nab_4)^{\leq j+1}\psi](\tau_1, \tau_2)+F[(\nab_3,  \nab_4, \dkb)^{\leq j+1}\psi](\tau_1, \tau_2)\\ 
\nn&\les& E^{j+1}[\psi](\tau_1)+\NN^{j+1}[\psi, N] (\tau_1, \tau_2)\\
&&+\left(\frac{|a|}{m}+\ep\right)\left(\sup_{\tau\in [\tau_1, \tau_2]}E^{j+1}[\psi]+B_\de^j[\psi](\tau_1, \tau_2)+F^{j+1}[\psi](\tau_1, \tau_2)\right)
\eeaa
as stated.
\end{enumerate}
This concludes the proof of Proposition \ref{prop:combinedenergyMorawetzfromenergyinchap9:higherderivatives}.
\end{proof}

We are now ready to prove Theorem \ref{THM:HIGHERDERIVS-MORAWETZ-CHP3}.
\begin{proof}[Proof of Theorem \ref{THM:HIGHERDERIVS-MORAWETZ-CHP3}]
Recall that \eqref{eq:iterationassumptiononnumberderivativesfroproofThm:HigherDerivs-Morawetz-chp3} holds for $j=2$ in view of \eqref{eq:iterationassumptiontruefors=2froproofThm:HigherDerivs-Morawetz-chp3}. Also, if the iteration assumption \eqref{eq:iterationassumptiononnumberderivativesfroproofThm:HigherDerivs-Morawetz-chp3} holds for some $2\leq j\leq s-1$, then, it also holds for $j$ replaced by $j+1$ in view of  Proposition \ref{prop:combinedenergyMorawetzfromenergyinchap9:higherderivatives}. Thus, we infer that  \eqref{eq:iterationassumptiononnumberderivativesfroproofThm:HigherDerivs-Morawetz-chp3} holds for all $2\leq j\leq s$. Thus, we have, for all   $2\leq s\le \kl$, and for any $\de>0$, 
\beaa
\nn&& \Mor[(\nab_3, \nab_4, \dkb)^{\leq s}\psi](\tau_1, \tau_2) + E[(\nab_3, \nab_4, \dkb)^{\leq s}\psi](\tau_2) +F[(\nab_3, \nab_4, \dkb)^{\leq s}\psi](\tau_1,\tau_2)  \\
&\les & E^s[\psi](\tau_1)  +\NN^s [\psi, N] (\tau_1, \tau_2)+\left(\frac{|a|}{m}+\ep\right)\left(\sup_{\tau\in [\tau_1, \tau_2]}E^s[\psi]+B_\de^s[\psi](\tau_1, \tau_2)\right).
\eeaa
This concludes the proof of Theorem \ref{THM:HIGHERDERIVS-MORAWETZ-CHP3}.
\end{proof}

 
\section{Conditional estimate for the scalar wave}
\lab{sec:conditionalMorawetzscalar}
 
  
 \begin{proposition}
\lab{proposition:GeneralWaveComplexPotential:Morawetz}
Let $\psi$ be a solution to the following scalar wave equation 
\bea
\square\psi+V\psi &=& N
\eea 
where $V$ is real and satisfies $V\sim r^{-3}$ for $r$ large. Then the   following estimates hold true  for all $2\leq s\le \kl$ and  some small $ \de>0$,
 \bea
\lab{eq:Estimatesforpsi-M8-1:Morawetz}
 \bsplit 
   \Mor[(\nab_3, \nab_4, \dkb)^{\leq s}\psi](\tau_1, \tau_2) 
\les& \int_{\MM}r^{-3}|(\nab_3, \nab_4, \dkb)^{\leq s}\psi|^2+ \Mor^{s-1}[\psi](\tau_1, \tau_2)\\
&+  F^{s-1}[\psi](\tau_1, \tau_2)+ \sup_{\tau\in[\tau_1,\tau_2] }  E^s[\psi](\tau)\\
&+\NNmor^s[\psi,  N](\tau_1, \tau_2) +\NNmor^s[\psi,  N](\tau_1, \tau_2)+\ep B^s_{\de}[\psi]
\end{split}
\eea
and 
\bea
\lab{eq:Estimatesforpsi-M8-1:energy-Morawetz}
 \bsplit 
& E [(\nab_3, \nab_4, \dkb)^{\leq s}\psi](\tau_2)+ \Mor[(\nab_3, \nab_4, \dkb)^{\leq s}\psi](\tau_1, \tau_2) +F[(\nab_3, \nab_4, \dkb)^{\leq s}\psi](\tau_1, \tau_2)\\
&\les \int_{\MM}r^{-3}|(\nab_3, \nab_4, \dkb)^{\leq s}\psi|^2+ \Mor^{s-1}[\psi](\tau_1, \tau_2)\\
&+  F^{s-1}[\psi](\tau_1, \tau_2)+ \sup_{\tau\in[\tau_1,\tau_2] }  E^{s-1}[\psi](\tau)+\NN^s[\psi,  N](\tau_1, \tau_2)+\ep B^s_{\de}[\psi]
\end{split}
\eea
 where
 \beaa
 \NN^s[\psi,  N](\tau_1, \tau_2):=\int_{\MM} \big(|\nab_{\Rhat}\dk^{\le s} \psi|+r^{-1}|\dk^{\le s}\psi|\big) | \dk^{\le s}N|.
 \eeaa
\end{proposition} 
 
 \begin{remark}
Note that both estimates are conditional on the control of the quantities $\int_{\MM}r^{-3}|(\nab_3, \nab_4, \dkb)^{\leq s}\psi|^2$, $\Mor^{s-1}[\psi](\tau_1, \tau_2)$ and $F^{s-1}[\psi](\tau_1, \tau_2)$. Proposition 
\ref{proposition:GeneralWaveComplexPotential:Morawetz} will be extended to a conditional $r^p$-weighted version in Proposition \ref{Prop:scalarwavePsi-M8-chap10}, see section \ref{sec:conditionalweigthedMorawetzscalar}. These estimates will be used  to control  $\Pc$  in Chapter \ref{Chapter:EN-MorforPc}. 
\end{remark}

The proof of Proposition \ref{proposition:GeneralWaveComplexPotential:Morawetz} is similar to the one of Theorem \ref{THM:HIGHERDERIVS-MORAWETZ-CHP3}.  Given that the estimates are only conditional, and in view of the the strong decay in $r$ for  the potential   $V$, and the fact that extending the Andersson-Blue method to perturbations of Kerr is more straightforward for scalar waves, the proof of Proposition \ref{proposition:GeneralWaveComplexPotential:Morawetz} is in fact simpler than one of Theorem \ref{THM:HIGHERDERIVS-MORAWETZ-CHP3}.

 
 \chapter{Proof of Theorems \ref{THEOREM:GENRW1-P} and \ref{THEOREM:GENRW2-Q}}
 \label{chapter-rp-estimates}
 

In this chapter, we derive the $r^p$-weighted estimates for  the  reduced  gRW equation  \eqref{eq:Gen.RW}
\bea\lab{eq:Gen.RW-chap10}
\squared_2 \psi -V\psi=- \frac{4 a\cos\th}{|q|^2}\dual \nab_T  \psi+N, \qquad V= \frac{4\De}{ (r^2+a^2) |q|^2},
\eea
on a spacetime $\MM$ which is   an admissible perturbation of Kerr in the sense that \eqref{eq:assumptionsonMMforpartII} holds.
Together with the Energy-Morawetz estimates of Theorem \ref{THM:HIGHERDERIVS-MORAWETZ-CHP3}, this will  conclude the proof of Theorems \ref{THEOREM:GENRW1-P} 
 and \ref{THEOREM:GENRW2-Q}.  These $r^p$ weighted  estimates concern only the region $r\ge R$ for a sufficiently large $R$. In such a region  the    equation 
 \eqref{eq:Gen.RW} closely resembles  the     equation 
 \beaa
 \squared_2 \psi -V\psi=N, \qquad V=\left(1-\frac{2m}{r} \right) \frac{4}{r^2}, 
 \eeaa
 which was studied in Chapter 10 of \cite{KS}.  The estimates in this  chapter  are thus  similar to the   $r^p$ estimates in sections 10.2-10.5  of \cite{KS}.   There is however  an important difference in that the  hypersurfaces $\Si(\tau)$ are  not  null, as in 
 \cite{KS}, but spacelike asymptotically null. This   leads to  some significant differences in the   proof of Theorem  \ref{THEOREM:GENRW1-P}.

 
\section{Proof of Theorem \ref{THEOREM:GENRW1-P}}
 

In order to prove Theorem  \ref{THEOREM:GENRW1-P}, we need to show that the following estimates  hold true  for solutions $\psi\in\sk_2$  of  \eqref{eq:Gen.RW-chap10} on an admissible  $\MM$, for all $\de\le p\le  2-\de$ and  $2\leq s\le \kl$,
       \bea
       \lab{eqtheorem:GenRW1-p:chap10}
       \BEF_p^s[\psi](\tau_1,\tau_2) \les
       E_p^s[\psi](\tau_1)+\NN_p^s[\psi, N](\tau_1, \tau_2),
       \eea
    where the   $\BEF_p^s$  norms have been introduced in   section \ref{subsection:basicnormsforpsi}.

We proceed in  steps as follows.  

{\bf Step 0.}  Recall that we have proved in section 
\ref{sec:proofofThm:HigherDerivs-Morawetz-chp3:chap9} the  global   energy-Morawetz estimates of Theorem   \ref{THM:HIGHERDERIVS-MORAWETZ-CHP3}, i.e. the fact that the following estimates  hold true for solutions $\psi\in\sk_2$  of  \eqref{eq:Gen.RW-chap10} on an admissible  $\MM$, for $|a|/ m \ll 1$ sufficiently small, for   $2\leq s\le \kl$, and for any $\de>0$,
\bea\lab{eq:conclusionmainenergymorawetz:Chap10}
\bsplit
& \Mor[(\nab_3, \nab_4, \dkb)^{\leq s}\psi](\tau_1, \tau_2) + E[(\nab_3, \nab_4, \dkb)^{\leq s}\psi](\tau_2) +F[(\nab_3, \nab_4, \dkb)^{\leq s}\psi](\tau_1,\tau_2)  \\
&\les  E^s[\psi](\tau_1)  +\NN^s [\psi, N] (\tau_1, \tau_2) +\left(\frac{|a|}{m}+\ep\right)\BEF^s_{\de}[\psi](\tau_1, \tau_2).
\end{split}
\eea
  In view of \eqref{eq:conclusionmainenergymorawetz:Chap10},  to complete the proof of Theorem \ref{THEOREM:GENRW1-P}, i.e. to prove \eqref{eqtheorem:GenRW1-p:chap10}, it remains     to derive $r$-weighted estimates in the region $r\geq R$ for $R$ large enough.

{\bf Step 1.} We first derive a basic $r$-weighted estimate for $\nab_3\psi$ in the region $r\geq R$.
\begin{proposition}
\lab{Prop:Step1-Chap10}
Let $R\gg m$ large enough. We have
\bea
\bsplit
\int_{\MM_{\geq R}(\tau_1, \tau_2)}r^{-1-\de}|\nab_3\psi|^2 \les& \int_{\MM_{\geq R}(\tau_1, \tau_2)}r^{\de-1}\Big(|\nab_4\psi|^2+|\nab\psi|^2+r^{-2}|\psi|^2\Big)\\
& +E_{\geq \frac{R}{2}}[\psi](\tau_1)+Mor_{\frac{R}{2}\leq r\leq R}[\psi](\tau_1, \tau_2)+\NN[\psi, N](\tau_1, \tau_2). 
\end{split}
\eea
\end{proposition}

Proposition \ref{Prop:Step1-Chap10}  is  proved  in section  \ref{section:Prop:Step1-Chap10}.

{\bf Step 2.} Next, we derive $r$-weighted estimates in the region $r\geq R$.
\begin{proposition}
\lab{Proposition:Step3-Chap10}
Let $R\gg m$ large enough. We have, for $\de\leq p\leq 2-\de$,  $0\le s\le k_L$,
\bea
\lab{eq:Proposition-Step3-Chap10}
\bsplit
\BEF^s_{p, \geq R}[\psi](\tau_1, \tau_2) \les& E^s_{p,\geq \frac{R}{2}}[\psi](\tau_1)+\NN^s_{p, \ge R/2}[\psi, N] +Mor^s_{\frac{R}{2}\leq r\leq R}[\psi](\tau_1, \tau_2).
\end{split}
\eea
\end{proposition}

The proof of Proposition \ref{Proposition:Step3-Chap10} is obtained in section \ref{sec:proof-Step3-Chap10}.

{\bf Step 3.} Combining the energy-Morawetz estimates \eqref{eq:conclusionmainenergymorawetz:Chap10} with the $r$-weighted estimates of Proposition \ref{Proposition:Step3-Chap10}, we infer,  for $\de\leq p\leq 2-\de$,  $0\le s\le k_L$,
\beaa
       \BEF_p^s[\psi](\tau_1,\tau_2) \les
       E_p^s[\psi](\tau_1)+\NN_p^s[\psi, N](\tau_1, \tau_2)+\left(\frac{|a|}{m}+\ep\right)\BEF^s_{\de}[\psi](\tau_1, \tau_2).
 \eeaa
For $a$ and $\ep$ small enough, may absorb the last term on the RHS from the LHS and obtain \eqref{eqtheorem:GenRW1-p:chap10} which concludes the proof of Theorem \ref{THEOREM:GENRW1-P}.

 
\section{Basic setup and control of $\nab_3\psi$}
 

 
\subsection{Renormalization of the horizontal structure}


As in sections 10.2-10.5  of \cite{KS} (see  in particular formula (10.2.6) in \cite{KS}),   it is convenient to work  with   the renormalized  frame 
\bea\lab{eq:renormalizedframe}
e_4' =\la e_4, \qquad e_3'=\la^{-1}e_3, \qquad e_a'=e_a, \qquad \la:=\frac{|q|^2}{\De}.  
\eea
 Note that  for  $r$ large  we have  $\la =\Up^{-1}\big(1+O(a^2r^{-2} )\big)$, with   $\Up= (1-\frac{2m}{r})$ as in Schwarzschild.

The corresponding renormalized quantities  verify the following.
 \begin{lemma}\lab{lemma:renormalizationoftheglobalframeforrgeqr1rpweigthedproof}
Let an unprimed frame and a primed frame related by \eqref{eq:renormalizedframe}. Using the ingoing normalization for the linearized quantities associated to the unprimed horizontal structure, and the  outgoing normalization for the linearized quantities associated to the primed horizontal structure, see section \ref{sec:definitionoflinearizedquantities:chap4} for the definition of the linearized quantities, we have
\beaa
&& \trXc '=\la\trXc, \quad \Xh'=\la\Xh, \quad \trXbc'=\la^{-1}\trXbc, \quad \Xbh'=\la^{-1}\Xbh,\\
&& \Hc'=\Hc,\quad \Hbc'=\Hbc, \quad \Xi' = \la^2\Xi, \quad  \Xib' = \la^{-2}\Xib,
\eeaa
\beaa
 \Zc' &=& \Zc  -\frac{2a^2\cos\th}{|q|^2}\widecheck{\DD(\cos\th)}-\left(\frac{2r}{|q|^2} -\frac{2r-2m}{\De}\right)\DD(r),\\
 \om' &=& \la\left(\omc  +\frac{1}{2}\la\pr_r\left(\frac{\De}{|q|^2}\right)\widecheck{e_4(r)}- \frac{a^2\cos\th}{|q|^2}e_4(\cos\th)\right), \\ 
 \ombc' &=& \la^{-1}\left(\omb  -\frac{1}{2}\la\pr_r\left(\frac{\De}{|q|^2}\right)\widecheck{e_3(r)}+ \frac{a^2\cos\th}{|q|^2}e_3(\cos\th)\right),
 \eeaa
\beaa
A'=\la^2A, \quad B'=\la B, \quad \Pc'=\Pc, \quad \Bb'=\la^{-1}\Bb, \quad \Ab'=\la^{-2}\Ab,
\eeaa
and 
\beaa
&& \widecheck{e_4'(r)}=\la\widecheck{e_4(r)}, \quad e_4'(\cos\th)=\la e_4(\cos\th), \quad \widecheck{\nab_4'\Jk}=\la\widecheck{\nab_4\Jk},\\
&& \widecheck{e_3'(r)}=\la^{-1}\widecheck{e_3(r)}, \quad \widecheck{e_3'(\cos\th)}=\la^{-1}\widecheck{e_3(\cos\th)}, \quad \widecheck{\nab_3'\Jk}=\la^{-1}\widecheck{\nab_3\Jk},\\
&& \DD'(r)=\DD(r), \quad \widecheck{\DD'(\cos\th)}=\widecheck{\DD(\cos\th)}, \quad  \DD'\hot\Jk = \DD\hot\Jk,\quad \widecheck{\ov{\DD'}\c\Jk} = \widecheck{\ov{\DD}\c\Jk}.
\eeaa 
 \end{lemma}
 
 \begin{proof}
 Under the conformal transformation \eqref{eq:renormalizedframe}, the complexified Ricci coefficients  transform as follows:
\beaa
&& \tr X'=\la\tr X, \quad \Xh'=\la\Xh, \quad \tr\Xb'=\la^{-1}\tr\Xb, \quad \Xbh'=\la^{-1}\Xbh,\\
&& Z' = Z -{\DD}'(\log\la), \quad H'=H,\quad \Hb'=\Hb, \quad \Xi' = \la^2\Xi, \quad  \Xib' = \la^{-2}\Xib,\\
&& \om'=\la\left(\om -\frac{1}{2}\la^{-1}e_4'(\log\la)\right), \quad \omb' = \la^{-1}\left(\omb+\frac{1}{2}\la e_3'(\log\la)\right).
\eeaa
Also, the curvature components transform as follows
\beaa
A'=\la^2A, \quad B'=\la B, \quad P'=P, \quad \Bb'=\la^{-1}\Bb, \quad \Ab'=\la^{-2}\Ab.
\eeaa
Using the choice $\la=\frac{|q|^2}{\De}$ in \eqref{eq:renormalizedframe}, the ingoing normalization for the linearized quantities associated to the unprimed horizontal structure, and the  outgoing normalization for the linearized quantities associated to the primed horizontal structure,   we easily infer the stated identities for the linearized Ricci and curvature coefficients. 
 
Also, we have
\beaa
&& e_4'(r)=\la e_4(r), \quad e_4'(\cos\th)=\la e_4(\cos\th), \quad \nab_4'\Jk=\la\nab_4\Jk,\\
&& e_3'(r)=\la^{-1}e_3(r), \quad e_3'(\cos\th)=\la^{-1}e_3(\cos\th), \quad \nab_3'\Jk=\la^{-1}\nab_3\Jk,\\
&& \DD'(r)=\DD(r), \quad \DD'(\cos\th)=\DD(\cos\th), \quad  \DD'\Jk = \DD\Jk,
\eeaa 
 which immediately yields the remaining statements. This concludes the proof of Lemma \ref{lemma:renormalizationoftheglobalframeforrgeqr1rpweigthedproof}.
 \end{proof}

 \begin{lemma}\lab{lemma:assumptionsonMMforpartII:chap10outgoingrenorm}
 Assume that the primed horizontal structure is related to the unprimed one by the conformal transformation \eqref{eq:renormalizedframe}. Also, assume that $(\Ga_b, \Ga_g)$ associated to the unprimed horizontal structure satisfies \eqref{eq:assumptionsonMMforpartII}. Then, $(\Ga_b', \Ga_g')$ associated to the primed horizontal structure satisfies  in the region  $r\ge  4m_0$
 \bea\lab{eq:assumptionsonMMforpartII:chap10outgoingrenorm}
\bsplit
r^3|\dk^{\leq k}\xi'|+ r^2|\dk^{\leq k}\Ga_g'|+r|\dk^{\leq k}\Ga_b'| &\les \ep, \qquad\qquad\,\,\, k\leq k_L,\\
r^3|\dk^{\leq k}\xi'|+ r^2|\dk^{\leq k}\Ga_g'|+r|\dk^{\leq k}\Ga_b'| &\les \frac{\ep}{\tau_{trap}^{1+\dec}}, \qquad k\leq \frac{k_L}{2}.
\end{split}
\eea
\end{lemma}

 \begin{proof}
 The proof is immediate in view of Lemma \ref{lemma:renormalizationoftheglobalframeforrgeqr1rpweigthedproof} and the fact that $\la=\frac{|q|^2}{\De}$ is smooth in the region $r\geq 4m_0$. 
 \end{proof}

\begin{remark}\lab{rmk:fromnowonoutgoingrenormalization}
In view of Lemma \ref{lemma:assumptionsonMMforpartII:chap10outgoingrenorm}, in the remainder of the chapter, we make all the calculations in the renormalized frame and, since there is no danger of confusion, we drop the primes.  In particular, in the frame we  shall use throughout the chapter,  we   thus have
 \beaa
 \bsplit
 e_4(r)&=1+\Ga_g, \quad e_3(r)=-\frac{\De}{|q|^2} +r\Ga_b, \quad \trch= \frac{2r}{|q|^2}+\Ga_g, \quad \trchb=-\frac{2r\Delta}{|q|^4}+\Ga_g,\\
 \om&\in \Ga_g,\quad \omb=\frac 1 2 \pr_r\left(\frac{\De}{|q|^2}\right)+\Ga_b, \quad Z= \frac{a\ov{q}}{|q|^2}\Jk+\Ga_g.
 \end{split}
 \eeaa
 \end{remark}

 
\subsection{Boundaries and integral identities}


 
 \subsubsection{Spacelike, asymptotically null hypersurface   $\Si=\Si(\tau)$} 
 

 In view of Definition \ref{definition:definition-oftau}, and normalizing the normal $N_\Si$ such that $\g(N_\Si, e_3)=-2$, we have
  \bea
 \lab{eq:normalizedN_Si}
 \bsplit
N_\Si&= e_4 +\frac 1 2  r^{-2}  \la e_3 + Y^b e_b, \qquad         |Y|=O(ar^{-1} ), \\
 g(N_\Si, N_\Si) &=-2 r^{-2} \la +|Y|^2  \les -\frac{m^2}{r^2},
 \end{split}
\eea
with $\la$ satisfying 
\bea
\lab{eq:choiceofla-chap10}
2 m^2\les \la\les 2 m^2  , \qquad D\la=O( \ep +R^{-1} ).
\eea
Also, define the  vectorfield   orthogonal to $N_\Si$,
\bea\label{eq:def-nu-Si-chap10}
\nu_\Si:=  e_4 - \frac 1 2   r^{-2} \la e_3 
\eea
and  note that  $\nu_\Si$ is tangent to $\Si$.

 
\subsubsection{Spacelike hypersurface   $\Si_*$} 


  Let $N_*$ be the vectorfield normal to $\Si_*$ of the form\footnote{See for comparison     Lemma 10.44 of \cite{KS}.} 
\bea
\lab{eq:NormalSi_*chap10}
N_*&=&e_4+U e_3 +Y_*, 
\eea
with  $U$ a scalar function    and  $Y_*$ 
 horizontal vectorfield verifying
 \beaa
   |U| = 1 +O(\ep), \qquad  |Y_*|=O(ar^{-1} ).
\eeaa
    Note also  that  the vectorfield  
    \bea
    \lab{eq:NormalTangentSi_*}
    \nu_*=\nu_{\Si_*}=   e_4 - U e_3
    \eea
     is perpendicular  to $N_*$ and thus    tangent to  $\Si_*$.


\subsubsection{Basic pointwise  identity of  Proposition  \ref{prop-app:stadard-comp-Psi-perturbations-Kerr}}
 
  
According  to Proposition  \ref{prop-app:stadard-comp-Psi-perturbations-Kerr},  if   $\psi\in \mathfrak{s}_2(\MM)$ is  a solution of \, $\squared_2\psi -V\psi= N'$, with
$N':=- \frac{4 a\cos\th}{|q|^2}\dual \nab_T  \psi+N$,  and   $X$ be  a vectorfield of the form $ X=  X^3  e_3 +X^4 e_4$, $w$ a scalar, $M$ a one form and $\PP_\mu[X, w, M]$ the current defined by
 \bea
 \lab{eq:defPP-chapter10}
\PP_\mu[X, w, M]&:=&\QQ_{\mu\nu} X^\nu +\frac 1 2  w \psi \c \Db_\mu \psi -\frac 1 4|\psi|^2   \pr_\mu w +\frac 1 4 |\psi|^2 M_\mu,
  \eea
 then,  
  \bea
  \lab{eq:DivergencePP-chapter10}
  \bsplit
  \D^\mu  \PP_\mu[X, w, M] &= \frac 1 2 \QQ  \c\piX - \frac 1 2 X( V ) |\psi|^2+\frac 12  w \LL[\psi] -\frac 1 4|\psi|^2   \square_\g  w + \frac 1 4  \Div(|\psi|^2 M\big) \\
  &+  \left(X( \psi )+\frac 1 2   w \psi\right)\c N'- \big(\rhod +\etab\wedge\eta\big)\nab_{X^4e_4-X^3e_3}  \psi\c\dual\psi \\
&- \frac{1}{2}\Im\Big(\tr\Xb H X^3 +\tr X\Hb X^4\Big)\c\nab\psi\c\dual\psi +r^{-2} \big(X^3\Ga_b+ X^4\Ga_g\big) \dk \psi \c \psi.
  \end{split}
 \eea
   Recalling  the  definition of the  expression $\EE[X, w, M] $ introduced in \eqref{definition-EE-gen},
   i.e.
   \bea
   \lab{eq:EE-chapter10}
   \bsplit
   \EE[X, w, M] &=\frac 1 2 \QQ  \c\piX - \frac 1 2 X( V ) |\psi|^2+\frac 12  w \LL[\psi] -\frac 1 4|\psi|^2   \square_\g  w\\
   & + \frac 1 4  \Div(|\psi|^2 M\big),
   \end{split}
   \eea
    we deduce,  with $ \PP_\mu= \PP_\mu[X, w, M]$, $ \EE=\EE[X, w, M]$,
      \bea
      \lab{eq:DivergencePP-chapter10-2}
      \begin{split}
  \D^\mu  \PP_\mu&= \EE  +  \big(\nab_X\psi+\frac 1 2   w \psi\big)\c  \big(\squared_k\psi - V \psi \big)   +\big(\rhod +\etab\wedge\eta\big)\nab_{X^4e_4-X^3e_3}  \psi\c\dual\psi \\
  &+\frac{1}{2}\Im\Big(\tr\Xb H X^3 +\tr X\Hb X^4\Big)\c\nab\psi\c\dual\psi+r^{-2} \big(X^3\Ga_b+ X^4\Ga_g\big) \dk \psi \c \psi.
\end{split}
\eea


\subsection{Proof of Proposition \ref{Prop:Step1-Chap10}}
\lab{section:Prop:Step1-Chap10}


The proof of Proposition \ref{Prop:Step1-Chap10} follows  easily  by integration  from the  following  lemma\footnote{This is the precise analogue of Proposition 10.36 in \cite{KS}. Note that the identity for $\EE[f_{-\de}T, 0, 0]$ is used together with Proposition \ref{prop-app:stadard-comp-Psi} which generalizes \eqref{eq:DivergencePP-chapter10} to the case of vectorfields which are not spanned by $(e_3, e_4)$.}.
\begin{lemma}
With the notation in \eqref{eq:EE-chapter10} the following   identity holds true in the region $r\ge R$
\beaa
\EE[f_{-\de}T, 0, 0]  = \frac 1 4  \de r^{-1-\de}  |\nab_3\psi|^2- \frac 1 4 \frac{\De^2}{|q|^4} \de r^{-1-\de}  |\nab_4 \psi|^2 + O(\ep+R^{-1})  r^{-1-\de} \big( |D \psi|^2+r^{-2} |\psi|^2\big) .
\eeaa
Moreover
\beaa
\PP_\mu[f_{-\de}T, 0, 0]\c e_4 =f_{-\de} \QQ(T, e_4) \geq 0, \qquad \PP_\mu[f_{-\de}T, 0, 0]\c e_3 =f_{-\de} \QQ(T, e_3) \geq 0.
\eeaa
\end{lemma}

\begin{proof}
The first part of the  lemma can be derived  by using   the       identity  \eqref{eq:EE-chapter10} with the vectorfield   $X= f_{-\de} \T$, $w=0$, $M=0$,   with $f_{-\de} = r^{-\de} $ for $r\ge R$ and supported for $r\ge R/2$ with $R$ sufficiently large.   We refer the reader to the proof of  Proposition 10.36 in \cite{KS}.
\end{proof}

 
 \section{$r^p$ weighted estimates}
 

 
 \subsection{Basic pointwise identities}
 
 
 In what follows  we apply  formulas \eqref{eq:DivergencePP-chapter10}, \eqref{eq:EE-chapter10} with the choice
 \bea
\lab{eqXwM-chapter10}
X=f(r)\left(e_4+ \frac 1 2  r^{-2} \la  e_3 \right), \qquad w=\frac{2r}{|q|^2}f, \qquad M=\frac{2r}{|q|^2} f 'e_4. 
\eea
 In view of \eqref{eq:normalizedN_Si},  we have $X= f(N_\Si - Y^b e_b)$.
We choose  $f=f_p$   non-negative  defined as $f_p=r^p$ for $r \geq R$ and $f_p=0$ for $r \leq R/2$, where $R$ is a fixed sufficiently large constant.
\begin{remark}
\lab{remark:errorterms-Chapt10} 
In view of the  intended choice  $f=f_p$  we  can   write schematically
\beaa
f'= O(R^{-1} ) f, \qquad  f''= O(R^{-2} ) f, \qquad ( r\pr_r)^{\le 2}(f) = O(1+R^{-1} ) f.
\eeaa
We will use this  to simplify     various  error terms in the identities that follow.
\end{remark} 

Based on the remark above  we  rewrite   formula  \eqref{eq:DivergencePP-chapter10-2} for $X$ as in \eqref{eqXwM-chapter10} in the form
 \bea
 \lab{eq:DivergencePP-chapter10-3}
      \begin{split}
& \D^\mu  \PP_\mu[X, w, M]  =\EE[X, w, M]   + \left(\nab_X\psi+\frac 1 2   w \psi\right)\c  \big(\squared_k\psi - V \psi \big) \\
& +\big(\ep+R^{-1}  \big)    r^{-1}  f \Big(|\nab_4\psi|^2 +|\nab\psi|^2 + r^{-2}  |\nab_3\psi|^2 + r^{-2}|\psi|^2 \Big).
\end{split}
\eea
We write $X=\Xone+\Xtwo$ with $\Xone= f e_4$ and  $\Xtwo=  \frac 1 2r^{-2}  \la  f  e_3$.  The following lemma  is similar to  Lemma 10.40 in \cite{KS}.
 \begin{lemma}\label{lemma:piXtilde-chap10}
 The following  hold true.
 \begin{enumerate}
 \item The deformation  tensor of the  vectorfield $\Xone=f(r)e_4$ is given by
 \beaa
\pione=\,  \frac {2 r}{|q|^2} f  \g + \pionet, 
 \eeaa
 with symmetric tensor $\pionet$ which verifies
   \bea
   \lab{eq:expression-piXt}
   \bsplit
 \pionet_{43}&= -2f'+\frac {4 }{r} f+   O(\ep+ R^{-1} ) r^{-1}  f,   \\
 \pionet_{33}&=4f'+   O(\ep+R^{-1})  r^{-1}  f, \\
  \pionet_{44}&=0, \quad  \pionet_{3a}=  O(\ep+R^{-1}) r^{-1} f, \quad 
 \pionet_{4a}, \,\pionet_{ab} =   O(\ep) r^{-2}  f.
 \end{split}
 \eea
 
 \item The deformation tensor  of the vectorfield $\Xtwo =\frac 1 2 r^{-2}  \la f  e_3 $ is given by 
\bea
\lab{eq:expression-piXt2}
\bsplit
\pitwo_{33} &=0,
\\ \pitwo_{44} &= O(\ep+ R^{-1}) r^{-2} f, \, \qquad 
\pitwo_{34}=O(\ep+ R^{-1}) r^{-2} f,\\
\pitwo_{ab}&= O(\ep+ R^{-1}) r^{-2} f, \,\,\qquad 
\pitwo_{3a}=  O(\ep  ) r^{-3} f, \,\,\\
\pitwo_{4a} &= O(\ep+ R^{-1})  r^{-3} f.
\end{split}
\eea
 
  \item For $w= \frac {2 r}{|q|^2} f$
  \bea
   \lab{eq:expression-LaX}
   \bsplit
 \square_\g w    &=\frac{2r}{|q|^2} f''+ O(\ep+ R^{-1} ) r^{-2} f .
 \end{split}
\eea

\item 
For  $M= \frac{2r}{|q|^2}  f' e_4    $  we have
\bea
 \lab{eq:expression-Div..M}
 \bsplit
 \Div(|\psi|^2 M\big)&=  4  r^{-1} f'\nab_4\psi\c\psi \\
 &+\left(\frac{2 f'}{r^2} +\frac{2 f''}{r} + O(\ep+R^{-1}) r^{-3} f  \right) |\psi|^2. 
 \end{split}
\eea
\end{enumerate}
  \end{lemma}
  
  \begin{proof}
   Since 
   $\piX_{\mu\nu}=\g(\D_\mu X, e_\nu) +\g(\D_\nu X, e_\mu) $ we  deduce
 \beaa
 \bsplit
 \pione_{44}&=0, \qquad \quad \,\, \pione_{43}=(e_4f)\g_{34}+4f\om, \qquad \pione_{33}=-8f\omb-4e_3f,\\
 \pione_{4a}&=2f\xi_a, \qquad \pione_{3a}=2f(\eta+\ze)_a, \qquad \qquad \pione_{ab}=2f\left(\chih_{ab}+ \frac 1 2 \frac{2r}{|q|^2} \g_{ab}\right),
 \end{split}
 \eeaa
 from which \eqref{eq:expression-piXt} easily follows.
 
Similarly
\beaa
\pitwo_{33} &=&0, \qquad\qquad \qquad 
\pitwo_{44} = -4 r^{-2}  \la  f \om - 2  e_4(r^{-2}\la  f)\\ 
 \pitwo_{ab}&=&  \frac 1 2 r^{-2}\la  f\chib_{ab}, \qquad \pitwo_{34}= - 2  r^{-2}\la  f \omb-2 e_3(r^{-2}\la  f),\qquad 
\\
\pitwo_{3a}&=&  r^{-2} f\la   \xib_a, \qquad \quad 
\pitwo_{4a}= r^{-2}  f\la ( \etab_a -\ze_a), 
\eeaa
from which, since $\la=(O(1)$  and $D\la=O(\ep +R^{-1}) $, \eqref{eq:expression-piXt2} follows.

 Using  formula \eqref{eq:waveH(r)}   for    $H=w= \frac {2 r}{|q|^2} f $, 
 \beaa
 \square_\g H &=&\frac{1}{\sqrt{|\g|} } \pr_\a \big(\sqrt{|\g|} \g^{\a\b} \pr_\b\big)  H=\frac{1}{|q|^2} \pr_r \big( \De\pr_r H\big)+\dkb^{\le 1 } ( \Ga_g\c  \dk^{ \le 1} H)\\
 &=&\frac{\De}{|q|^2}  H ''+\frac{2(r-m)}{|q|^2}  H'+r^{-2} \Big( (r\pr_r )^{\le 2}  f\c   \Ga_b+ (r\pr_r )^{\le 1}  f\c \dk^{\le 1 }\Ga_b \Big)\\
 &=& \frac{\De}{|q|^2}  H ''+\frac{2(r-m)}{|q|^2}  H'+ O(\ep, R^{-1} ) r^{-3} f. 
 \eeaa
Hence, modulo error terms of the form  $O(\ep, R^{-1} ) r^{-3} f$, 
 \beaa
 \square_\g w  &=&\frac{\De}{|q|^2}  \left(\frac {2 r}{|q|^2} f\right)''+\frac{2(r-m)}{|q|^2}  \left(\frac {2 r}{|q|^2} f\right)'\\
 &=&\frac{\De}{|q|^2}  \left(\frac {2 r}{|q|^2} f''+2\left(\frac {2 r}{|q|^2}\right)' f'+\left(\frac {2 r}{|q|^2}\right)'' f\right)+\frac{2(r-m)}{|q|^2}  \left(\frac {2 r}{|q|^2} f'+\left(\frac {2 r}{|q|^2}\right)' f\right).
\eeaa
Since  $\partial_r(\frac {2 r}{|q|^2})=\frac{-2r^2+2a^2\cos^2\th}{|q|^4}=-\frac{2}{|q|^2}+\frac{4a^2\cos^2\th}{|q|^4}$,
 \beaa
 \square_\g w    &=&\frac{\De}{|q|^2}  \left(\frac {2 r}{|q|^2} f''+\left(-\frac{4}{|q|^2}+\frac{8a^2\cos^2\th}{|q|^4}\right) f'+\left(\frac {4 r}{|q|^4}-\frac{16a^2\cos^2\th r}{|q|^6}\right) f\right)\\
 &&+\frac{2(r-m)}{|q|^2}  \left(\frac {2 r}{|q|^2} f'+\left(-\frac{2}{|q|^2}+\frac{4a^2\cos^2\th}{|q|^4}\right) f \right)+O( R^{-1} ) r^{-3} f\\
  &=&\frac{2r\De}{|q|^4} f''-\frac{4(\De-r^2+mr)}{|q|^4}  f'+\frac{4(\De-r^2+mr) r }{|q|^6}  f+O\left(\frac{a^2}{r^5}\right) (r\pr_r)^{\le 1 } f \\
  &&+O( R^{-1} ) r^{-3} f\\
 &=&\frac{2r}{|q|^2} f''+O( R^{-1} ) r^{-3} f.
\eeaa
Thus $\square_\g w  =\frac{2r}{|q|^2} f''+O(\ep+ R^{-1} ) r^{-3} f$
as stated.

To prove the last part of the lemma we write
\beaa
\Div(|\psi|^2 M\big) &=&4 r^{-1}  f'\nab_4\psi\c\psi +\Div M |\psi|^2, 
\\
\Div M&=& 2\Div(r^{-1} f'  e_4) =\frac{2 f'}{r^2} +\frac{2 f''}{r}  +O(\ep+R^{-1}) r^{-3} f, 
\eeaa
as stated.
  \end{proof}
  
 We are now ready to prove the following.
\begin{proposition}
\lab{proposition:MainIdentity-chap10}
The following hold true:
\begin{enumerate}
 \item We have, for $\EE=\EE[X, w, M]$,   $X= f(e_4 +\frac 1 2 r^{-2}\la  e_3)$ and 
 $w=\frac{2r}{|q|^2 } f$,
\bea
\lab{eq:MainIdentity-chap10-1}
\bsplit
\EE&= f'|\nab_4 \psi|^2 +\frac 1 2 \left( - f'+\frac{2}{r} f \right)\big(|\nab\psi|^2+V|\psi|^2 \big)  -\frac 1 2  \frac{r}{|q|^2} f'' |\psi|^2 + \frac 1 4  \Div(|\psi|^2 M\big)\\
&+O(\ep+R^{-1} )r^{-1} f\Big( |\nab_4 \psi|^2+|\nab\psi|^2 
+r^{-2}|\psi|^2 + r^{-2} |\nab_3\psi|^2 \Big).
\end{split} 
\eea

\item If in addition we choose  
$M= 2\frac{r}{|q|^2} f'  e_4$ we  deduce, with 
$\nabcheck_4\psi= \nab_4\psi+ r^{-1} \psi$,
\bea
\lab{eq:MainIdentity-chap10-2}
\bsplit
\EE&=\frac 1 2  f' |\nabcheck_4  \psi|^2  +\frac 1 2 \left( - f'+\frac{2}{r} f \right)\big(|\nab\psi|^2+V|\psi|^2 \big) \\
&+O(\ep+R^{-1} )r^{-1} f\Big( |\nab_4 \psi|^2+|\nab\psi|^2 
+r^{-2}|\psi|^2 + r^{-2} |\nab_3\psi|^2 \Big).
\end{split}
 \eea
\end{enumerate}
\end{proposition}

\begin{proof}
We   calculate the expression
\beaa
\EE&=&\frac 1 2 \QQ  \c\piX +\frac 12  w \LL[\psi]  -\frac 1 4|\psi|^2   \square_\g  w  - \frac 1 2 X( V ) |\psi|^2 + \frac 1 4  \Div\big(|\psi|^2 M\big)\\
&=&\EE'+  \frac 1 4  \Div\big(|\psi|^2 M\big)
\eeaa
with $X,w$ as in \eqref{eqXwM-chapter10}.
Recalling that 
\beaa
\QQ_{\mu\nu}=\Db_\mu\psi\c \Db_\nu\psi- \frac 1 2 \g_{\mu\nu} \LL, \qquad \LL=\Db ^\mu\psi\c \Db_\mu \psi+V|\psi|^2,  
\eeaa
we calculate, using Lemma \ref{lemma:piXtilde-chap10},
\beaa
\EE'&=&\frac 1 2 \QQ  \c\piX +\frac 12  w \LL[\psi]  - \frac 1 2 X( V ) |\psi|^2-\frac 1 4|\psi|^2   \square_\g  w\\
&=&\frac 1 2 \QQ  \c\pione+\frac 1 2 \QQ  \c\pitwo +\frac 12  w \LL[\psi]  - \frac 1 2 X( V ) |\psi|^2-\frac 1 4|\psi|^2   \square_\g  w\\
&=& \frac 1 2 \QQ\c \left(  \frac {2 r}{|q|^2} f  \g + \pionet\right) +\frac 12  \frac{2r}{|q|^2}f  \Db ^\mu\psi\c \Db_\mu \psi +\frac 12  \frac{2r}{|q|^2}f V |\psi|^2- \frac 1 2 X( V ) |\psi|^2\\
&&+\frac 1 2 \QQ  \c\pitwo  -\frac 1 4|\psi|^2   \square_\g  w\\
&=&  -\frac {2r}{|q|^2} f \LL[\psi] +  \frac{2r}{|q|^2}f  \Db ^\mu\psi\c \Db_\mu \psi+\frac 12  \frac{2r}{|q|^2}f V |\psi|^2 - \frac 1 2 X( V ) |\psi|^2 \\
&&+\frac 1 2 \QQ\c  \pionet +\frac 1 2 \QQ  \c\pitwo  -\frac 1 4|\psi|^2   \square_\g  w\\
&=&-\left( \frac { r}{|q|^2} f  V  + \frac 1 2 X( V ) \right)|\psi|^2 +\frac 1 2 \QQ\c  \pionet +\frac 1 2 \QQ  \c\pitwo \\
&& -\frac 1 4|\psi|^2 \left( \frac{2r}{|q|^2} f'' +O(\ep+ R^{-1} ) r^{-3} f\right).
\eeaa
We deduce
\beaa
\EE' &=& \frac 1 2 \QQ\c  \pionet+\frac 1 2 \QQ  \c\pitwo  -\frac 1 2\left( \frac {2 r}{|q|^2} f  V  +  X( V ) \right)|\psi|^2\\
&& -\frac 1 4|\psi|^2 \left( \frac{2r}{|q|^2} f'' +O(\ep+ R^{-1} ) r^{-3} f\right).
\eeaa
Recall that 
\beaa
 \QQ_{33}= |\nab_3 \psi|^2, \qquad \QQ_{44}=|\nab_4 \psi|^2, \qquad \QQ_{34}=|\nab \psi|^2 + V |\psi|^2, \qquad \QQ_{4a}=\nab_4 \psi \c \nab_a \psi. 
 \eeaa
Therefore in  view of Lemma \ref{lemma:piXtilde-chap10} we have
\beaa
 \QQ\c  \pionet&=& f'|\nab_4 \psi|^2 +\left( - f'+\frac{2}{r} f \right)\big(|\nab\psi|^2+V|\psi|^2 \big)\\
&& +O(\ep+R^{-1} )r^{-1} f\Big( |\nab_4 \psi|^2+|\nab\psi|^2 +V|\psi|^2 \Big)\\
&& O(\ep+R^{-1}) r^{-1} f  \QQ_{4a} +  O(\ep+R^{-1}) r^{-2} f  \QQ_{3a} + O(\ep) r^{-2}  f \QQ_{ab}.
\eeaa
Since 
 \beaa
 |\QQ_{3a}| \leq |\nab_3 \psi ||\nab \psi|, \quad |\QQ_{ab}| \leq |\nab_3 \psi||\nab_4 \psi|+|\nab \psi|^2 + |V| |\psi|^2, \quad  |\QQ_{4a}| \leq |\nab_4 \psi ||\nab \psi|,
 \eeaa 
we easily deduce
\beaa
\QQ\c  \pionet&=& f'|\nab_4 \psi|^2 +\left( - f'+\frac{2}{r} f \right)\big(|\nab\psi|^2+V|\psi|^2 \big)\\
&&+O(\ep+R^{-1} )r^{-1} f\Big( |\nab_4 \psi|^2+|\nab\psi|^2 
+r^{-2}|\psi|^2 + r^{-2} |\nab_3\psi|^2 \Big).
\eeaa
Similarly
\beaa
\QQ\c  \pitwo&=&O(\ep+R^{-1} )r^{-1} f\Big( |\nab_4 \psi|^2+|\nab\psi|^2 
+r^{-2}|\psi|^2 + r^{-2} |\nab_3\psi|^2 \Big).
\eeaa
Recall also that $V= \frac{4\De}{ (r^2+a^2) |q|^2}$. Therefore, since
\footnote{See for example  the proof of Proposition 10.2.5 in \cite{KS}. } 
\beaa
\frac 1 2  e_4(V)+\frac{f}{r} V=O(m r^{-4} )+ O(a^2 r^{-6}) , \qquad \frac { r}{|q|^2} f  V =\frac{1}{r} V+ O( ar^{-4}),
\eeaa
we deduce
\beaa
\frac { r}{|q|^2} f  V  + \frac 1 2 X( V )= \frac { r}{|q|^2} f  V+ \frac 1 2  f\left(
 e_4 ( V)+\frac 1 2 r^{-2}\la  e_3(V) \right)= O(R^{-1} ) r^{-3} f.
\eeaa
Consequently
\beaa
\EE'&=&\frac 1 2  f'|\nab_4 \psi|^2 +\frac 1 2 \left( - f'+\frac{2}{r} f \right)\big(|\nab\psi|^2+V|\psi|^2 \big)  -\frac 1 2  \frac{r}{|q|^2} f'' |\psi|^2  \\
&&+O(\ep+R^{-1} )r^{-1} f\Big( |\nab_4 \psi|^2+|\nab\psi|^2 
+r^{-2}|\psi|^2 + r^{-2} |\nab_3\psi|^2 \Big),
\eeaa
and thus
\beaa
\EE&=&\EE'  + \frac 1 4  \Div(|\psi|^2 M\big)\\
&=&\frac 1 2  f'|\nab_4 \psi|^2 +\frac 1 2 \left( - f'+\frac{2}{r} f \right)\big(|\nab\psi|^2+V|\psi|^2 \big)  -\frac 1 2  \frac{r}{|q|^2} f'' |\psi|^2 + \frac 1 4  \Div(|\psi|^2 M\big)\\
&&+O(\ep+R^{-1} )r^{-1} f\Big( |\nab_4 \psi|^2+|\nab\psi|^2 
+r^{-2}|\psi|^2 + r^{-2} |\nab_3\psi|^2 \Big) 
\eeaa
as stated.

To derive the second part of the Proposition we   choose 
$M=  2\frac{r}{|q|^2}f' e_4$  and make use of the formula \eqref{eq:expression-Div..M}
\beaa
\Div(|\psi|^2 M\big)&=&  4  r^{-1} f'\nab_4\psi\c\psi +\left(\frac{2 f'}{r^2} +\frac{2 f''}{r} + O(\ep+R^{-1}) r^{-3} f  \right) |\psi|^2. 
\eeaa
We deduce
\beaa
\EE&=&\frac 1 2  f'|\nab_4 \psi|^2 +\frac 1 2 \left( - f'+\frac{2}{r} f \right)\big(|\nab\psi|^2+V|\psi|^2 \big)  -\frac 1 2  \frac{r}{|q|^2} f'' |\psi|^2 \\
&&+  \left(   r^{-1} f'\nab_4\psi\c\psi +\frac 1 2 \left(\frac{ f'}{r^2} +\frac{ f''}{r}\right)  |\psi|^2\right)\\
&&+O(\ep+R^{-1} )r^{-1} f\Big( |\nab_4 \psi|^2+|\nab\psi|^2 
+r^{-2}|\psi|^2 + r^{-2} |\nab_3\psi|^2 \Big).
\eeaa

Consider the term
\beaa
J&=&\frac 1 2  f'|\nab_4 \psi|^2   -\frac 1 2  \frac{r}{|q|^2} f'' |\psi|^2  +   r^{-1} f'\nab_4\psi\c\psi +\frac 1 4 \left(\frac{2 f'}{r^2} +\frac{2 f''}{r} + O(\ep+R^{-1}) r^{-3} f  \right) |\psi|^2\\
&=&\frac 1 2 f'\Big( |\nab_4 \psi|^2 +2   r^{-1} \nab_4\psi\c\psi + r^{-2}|\psi|^2 \Big)  + O(\ep+R^{-1}) r^{-3} f |\psi|^2 \\
&=& \frac 1 2 f'\big| \nab_4  \psi+ r^{-1} \psi\big|^2  + O(\ep+R^{-1}) r^{-3} f |\psi|^2.
\eeaa
We  thus infer that, recalling that $\nabcheck_4\psi= \nab_4\psi+ r^{-1} \psi$,
\beaa
\EE&=&\frac 1 2  f' |\nabcheck_4  \psi|^2  +\frac 1 2 \left( - f'+\frac{2}{r} f \right)\big(|\nab\psi|^2+V|\psi|^2 \big) \\
&&+O(\ep+R^{-1} )r^{-1} f\Big( |\nab_4 \psi|^2+|\nab\psi|^2 
+r^{-2}|\psi|^2 + r^{-2} |\nab_3\psi|^2 \Big)
 \eeaa
 as stated.
 \end{proof}

 
 \subsubsection{Boundary terms}
 
 
 When we integrate   formula \eqref{eq:DivergencePP-chapter10-3} we  get, in addition
  to the  bulk terms  expressed in Proposition \ref{proposition:MainIdentity-chap10},
   boundary terms  on  $\Si_{\ge R/2} (\tau)\cup \Si_* $.     In what follows we   deal with these boundary terms. 
   \begin{lemma}
   \lab{lemma:pointwiseBoun-Chap10}
  Given   $ \PP=\PP[X, w, M]$  as in \eqref{eq:defPP-chapter10}  with   $(X, w, M)$ as in \eqref{eqXwM-chapter10},  we have
  \bea
  \bsplit
  \PP\c  e_4&= f\big|\nabcheck_4\psi\big|^2  - \frac 1 2 r^{-2}  \nab_4( r  f |\psi|^2) +\frac 1 2 r^{-2}\la  f   \QQ_{34} +O(\ep+R^{-1} ) r^{-3}  f|\psi|^2, \\
  \PP\c e_3 &=  f  \QQ_{34}+    f\frac 1 2 r^{-2}\la \QQ_{33}  + \frac 1 2 r^{-2}\nab_3\big ( r    f |\psi|^2) + (\ep+R^{-1}) f r^{-3} |\psi|^2.
  \end{split}
  \eea
  Also,
  \bea
  \bsplit
   \big|\PP\c  Y\big|     &\les  f   \big|\nabcheck_4 \psi \big|  |Y| |   \nab \psi|+ O(R^{-\de}) r^{-2} |\nab\psi|^2 +
 r^{-4+\de}|\nab_3\psi|^2\\
 &+ O(\ep +R^{-1} ) r^{-3}  f |\psi|^2 .
 \end{split}
   \eea
  \end{lemma}
     
   \begin{proof}  
   We write, since $X= f( e_4+\frac 1 2 r^{-2} \la   e_3)$
    \beaa
\PP\c  e_4&=&\left(\QQ_{\mu\nu} X^\nu +\frac 1 2  w \psi \c \Db_\mu \psi -\frac 1 4|\psi|^2   \pr_\mu w +\frac 1 4 |\psi|^2 M_\mu\right) e_4^\mu\\
          &=&\QQ(X,  e_4 ) + f \frac{r}{|q|^2} \psi \c \nab_4 \psi -\frac 1  2  e_4 \left( f \frac{r}{|q|^2}\right) |\psi|^2  \\
 &=&f\QQ(e_4 , e_4 )      +\frac 1 2 r^{-2}  f   \la   \QQ(e_4 , e_3 )      + f \frac{r}{|q|^2} \psi\c \nab_4 \psi -\frac 1  2  e_4 \left( f \frac{r}{|q|^2}\right) |\psi|^2\\
 &=&f\big|\nab_4 \psi|^2  + f r^{-1}  \psi\c\nab_4 \psi  -\frac 1  2  e_4 ( f r^{-1}) |\psi|^2+\frac 1 2 r^{-2}  f   \la  \QQ_{34} +O(\ep+R^{-1} ) r^{-3}  f|\psi|^2.
\eeaa

We  rewrite  the expression 
\beaa
I=f\big|\nab_4 \psi|^2  + f r^{-1}  \psi \c \nab_4 \psi  -\frac 1  2  e_4 ( f r^{-1}) |\psi|^2
\eeaa
in the form
\beaa
I&=&f\left(  |\nab_4 \psi|^2 +\frac 1 r  \psi \cdot \nab_4 \psi\right) - \frac 1 2 e_4(r^{-1} f) |\psi|^2\\
&=&f\left| \nab_4\psi+ \frac 1 r   \psi\right|^2 - \frac 1 r    f  \psi \cdot \nab_4\psi- r^{-2}    f  |\psi|^2 - \frac 1 2 e_4(r^{-1} f) |\psi|^2\\
&=&f\left|\nabcheck_4\psi\right|^2 - \frac 1 2 r^{-2}  \nab_4( r  f |\psi|^2)+ \frac 1 2 r^{-2}  e_4( r  f )|\psi|^2   - r^{-2}    f  |\psi|^2      - \frac 1 2 e_4(r^{-1} f) |\psi|^2\\
&=&f\left|\nabcheck_4\psi\right|^2  - \frac 1 2 r^{-2}  \nab_4( r  f |\psi|^2)+r^{-2}   (e_4(r)-1)       f   |\psi|^2.
\eeaa
Hence
\beaa
\PP\c  e_4&=& f\big|\nabcheck_4\psi\big|^2  - \frac 1 2 r^{-2}  \nab_4( r  f |\psi|^2) +\frac 1 2 r^{-2}  f   \la \QQ_{34} +O(\ep+R^{-1} ) r^{-3}  f|\psi|^2 
\eeaa
as stated.

 Also,  since  $M= 2\frac{r}{|q|^2} f ' e_4 $,  
\beaa
\PP\c e_3 &=&\left(\QQ_{\mu\nu} X^\nu +\frac 1 2  w \psi \c \Db_\mu \psi -\frac 1 4|\psi|^2   \pr_\mu w +\frac 1 4 |\psi|^2 M_\mu\right) e_3^\mu\\
&=& \QQ(X, e_3)+f\frac{r}{|q|^2} \psi \c \nab_3 \psi-\frac 1 2   e_3 \left(f\frac{r}{|q|^2} \right) |\psi|^2 -\frac{r}{|q|^2} f'  |\psi|^2 \\
&=& f  \QQ_{34}+  \frac 1 2 r^{-2}  f   \la \QQ_{33}+  \frac 1 2   f\frac{r}{|q|^2} \nab_3(|\psi|^2) -\frac 1  2  e_3  \left(  f \frac{r}{|q|^2} \right) |\psi|^2 -\frac{r}{|q|^2} f'  |\psi|^2\\
&=&  J+  f  \QQ_{34}+ \frac 1 2 r^{-2}  f   \la \QQ_{33}+O(\ep+R^{-1} )r^{-2} |\psi|^2
\eeaa
with   
\beaa
J= \frac 1 2 r^{-1}   f \nab_3(|\psi|^2) -\frac 1  2  e_3  ( r^{-1}   f  ) |\psi|^2 -\frac{1}{r} f'  |\psi|^2.
\eeaa
 We rewrite  $J$  in the form
\beaa
 2 J&=&  r^{-2}\nab_3\big ( r    f |\psi|^2)-   r^{-2} e_3(rf)|\psi|^2  -  e_3  ( r^{-1}   f  ) |\psi|^2 - 2r^{-1} f'  |\psi|^2\\
 &=&r^{-2}\nab_3\big ( r    f |\psi|^2) - \Big(2 r^{-1}  e_3(f) + r^{-2}  f e_3(r) +e_3(r^{-1} ) f  + 2r^{-1} f'\Big)|\psi|^2\\
 &=&r^{-2}\nab_3\big ( r    f |\psi|^2) -\Big( 2 r^{-1}  e_3(r) f' + r^{-2}  f e_3(r) - r^{-2}  e_3(r) f + 2r^{-1} f'\Big)|\psi|^2\\
 &=&r^{-2}\nab_3\big ( r    f |\psi|^2) - 2f' \Big(  r^{-1}  e_3(r)  + r^{-1} \Big)|\psi|^2.
 \eeaa
 Since
$ e_3(r)=-\frac{\De}{|q|^2}+r\Ga_b=-1 +(\ep+R^{-1}) $ we deduce
\beaa
 2 J&=&r^{-2}\nab_3\big ( r    f |\psi|^2) + (\ep+R^{-1}) f' r^{-2} |\psi|^2.
\eeaa
Therefore,
\beaa
\PP\c e_3 &=&  f  \QQ_{34}+   \frac 1 2 r^{-2}  f   \la \QQ_{33} + \frac 1 2 r^{-2}\nab_3\big ( r    f |\psi|^2) + (\ep+R^{-1}) f r^{-3} |\psi|^2
\eeaa
as stated.

 Finally,
 \beaa
\PP\c  Y&=&\left(\QQ_{\mu\nu} X^\nu +\frac 1 2  w \psi \c \Db_\mu \psi -\frac 1 4|\psi|^2   \pr_\mu w +\frac 1 4 |\psi|^2 M_\mu\right) Y^\mu\\
&=&\QQ(X, Y)+   \frac{r}{|q|^2}   f  \psi \c \nab_Y \psi      -|\psi|^2  Y\left(\frac{r}{2|q|^2}  f \right)\\
&=&f\QQ(e_4 , Y)+ \frac 1 2 r^{-2}\la  f \QQ(e_3  , Y)  +  \frac{r}{|q|^2}   f  \psi \c \nab_Y \psi     -|\psi|^2  Y\left(\frac{r}{2|q|^2}  f \right)\\
&=&f\nab_4 \psi \c\nab_Y \psi +\frac 1 2 r^{-2}  f   \la \nab_3 \psi\c \nab_Y\psi   +  \frac{r}{|q|^2}   f  \psi \c \nab_Y \psi    -|\psi|^2  Y\left(\frac{r}{2|q|^2}  f \right)\\
&=&f \big(\nabcheck_4 \psi-  r^{-1} \psi \big)\c  \nab_Y\psi +\frac 1 2 r^{-2}  f   \la \nab_3 \psi\c \nab_Y\psi   +
 \frac{r}{|q|^2}   f  \psi \c \nab_Y \psi    -|\psi|^2  Y\left(\frac{r}{2|q|^2}  f \right).
\eeaa
We deduce
\beaa
\PP\c  Y&=& f\nabcheck_4 \psi\c  \nab_Y\psi +\left( \frac{r}{|q|^2} -\frac 1 r \right) f\psi\c \nab_Y \psi +\frac 1 2 r^{-2}  f   \la\nab_3 \psi\c \nab_Y\psi
 -|\psi|^2  Y\left(\frac{r}{2|q|^2}  f \right)
\eeaa
and hence
\beaa
\big|\PP\c  Y\big|&\les & f   \big|\nabcheck_4 \psi \big|  |Y| |   \nab \psi|
+ O(r^{-4})  f |\psi | |\nab\psi|  +O(R^{-\de/2})|Y||\nab\psi| r^{-2+\de/2} |\nab_3\psi |\\
&&+ O(\ep +R^{-1} ) r^{-3}  f |\psi|^2 \\
&\les&  f   \big|\nabcheck_4 \psi \big|  |Y| |   \nab \psi|+ O(R^{-\de}) r^{-2} |\nab\psi|^2 +
 r^{-4+\de} |\nab_3\psi|^2+ O(\ep +R^{-1} ) r^{-3}  f |\psi|^2 
\eeaa
as stated. This concludes the proof of Lemma \ref{lemma:pointwiseBoun-Chap10}.
   \end{proof}

   \begin{proposition}
   \lab{Prop:pointwiseBoun-Chap10}
   The following   bounds hold true, for   $r\ge R $ sufficiently large:
   \begin{enumerate}
   \item On $\Si=\Si(\tau)$, for $0\le p\le 2-\de$,
   \bea
    \lab{eq:Prop.pointwiseBoun-Chap10-1}
    \bsplit
   \PP\c N_\Si &\ges   f\left|\nabcheck_4\psi\right|^2 +r^{-2}    f |\nab\psi|^2  - \frac 1 2  \div_\Si\big(r^{-1} f|\psi|^2 \nu_\Si\big)\\
 &  - O(\ep +R^{-1} ) r^{-3}  f |\psi|^2 -r^{-4+\de} f\Big(|\nab_3\psi|^2 +|\nab\psi|^2 + r^{-2} |\psi|^2 \Big).
\end{split}
\eea

\item  On $\Si_*$, with\footnote{Recall   \eqref{eq:NormalSi_*chap10}, \eqref{eq:NormalTangentSi_*}.} 
$N_*=e_4+U e_3 +Y_*$,  $\nu_*=  e_4 - U e_3  $,  for $0\le p\le 2-\de$,
   \bea
    \lab{eq:Prop.pointwiseBoun-Chap10-2}
    \bsplit
   \PP\c N_{*} &\ges  f\big(|\nabcheck_4\psi|^2 + |\nab\psi|^2 +r^{-2}  |\nab_3\psi|^2\big)    - \frac 1 2 \div_{\Si_*}\big(r^{-1} f|\psi|^2 \nu_*\big) \\
&   - O(R^{-1}+\ep)r^{-3}f|\psi|^2 .
\end{split}
   \eea    
   
   \item  On $\Si$  we have, with $f=r^p$ and for  $0\le p\le 1-\de$,
    \bea
    \lab{eq:Prop.pointwiseBoun-Chap10-3}
    \bsplit
     \PP\c N_\Si&\ge      \frac{\de^2}{8}    r^{p-2} |\psi|^2 - \frac p 2 r^{-2}  \nu_\Si(  r^{p+1} |\psi|^2)
     +\frac{m}{r^2}r^p\QQ_{34}\\
     &-O(r^{p-3} )\Big( |\nab_3 \psi|^2  +|\psi|^2\Big).
     \end{split}
     \eea
     
    \item  On $\Si_*$  we have, with $f=r^p$ and for $0\le p\le 1-\de$,
    \bea
    \lab{eq:Prop.pointwiseBoun-Chap10-4}
    \bsplit
     \PP\c N_*&\ge      \frac{\de^2}{8}    r^{p-2}|\psi|^2 - \frac p 2 r^{-2}  \nu_*(  r^{p+1} |\psi|^2)
     +  r^p  |\nab\psi|^2 \\
     &-O(r^{p-3} )\Big( |\nab_3 \psi|^2  +|\psi|^2\Big).
     \end{split}
     \eea
\end{enumerate}
   \end{proposition}
   
   \begin{proof}
    We calculate using the definition $N_\Si= e_4 + \frac 1 2 r^{-2}\la e_3 + Y^b e_b+ O(r^{-3}) e_3$, see \eqref{eq:normalizedN_Si}, and the lemma above 
 \beaa
 &&\PP\c N_\Si   -\PP\c Y=   \PP\c \left(e_4+\frac 1 2 r^{-2}\la e_3+O(r^{-3}) e_3\right) \\
&=&  f\left|\nabcheck_4\psi\right|^2  - \frac 1 2 r^{-2}  \nab_4( r  f |\psi|^2) +\frac 1 2 r^{-2}\la   f  \QQ_{34} +O(\ep+R^{-1} ) r^{-3}  f|\psi|^2\\
&&+\frac 1 2 r^{-2}\la\left( f  \QQ_{34}+    f\frac 1 2 r^{-2}\la \QQ_{33}+ \frac 1 2 r^{-2}\nab_3\big ( r    f |\psi|^2) + (\ep+R^{-1}) f r^{-3} |\psi|^2\right)\\
&=&  f\left|\nabcheck_4\psi\right|^2 +\frac 1 2 r^{-2}\la    f  \QQ_{34} +   f\frac 1 4 r^{-4} \la |\nab_3\psi|^2   - \frac 1 2 r^{-2} \left( \nab_4  - \frac 1 2 r^{-2}\la \nab_3\right)  (rf |\psi|^2)\\
&& +O(\ep+R^{-1} ) r^{-3}  f|\psi|^2\\
&=&  f\left|\nabcheck_4\psi\right|^2 +\frac 1 2 r^{-2} \la    f |\nab\psi|^2  +    f\frac 1 4 r^{-4} \la|\nab_3\psi|^2   - \frac 1 2 r^{-2} \nu_\Si (rf |\psi|^2)  +O(\ep+R^{-1} ) r^{-3}  f|\psi|^2,
    \eeaa
    where  we recall the definition of $\nu_\Si=e_4-\frac 1 2  r^{-2} \la  e_3 $, see \eqref{eq:def-nu-Si-chap10}.

     We now write
     \beaa
     \div_\Si\big(r^{-1} f|\psi|^2 \nu_\Si\big)&=&\nu_\Si\big( r^{-1} f|\psi|^2 \big)
      + r^{-1} f|\psi|^2 \div_\Si(\nu_\Si)\\
      &=& r^{-2}  \nu_\Si\big( r f|\psi|^2 \big)+\nu_\Si(r^{-2})  r f |\psi|^2  + r^{-1} f|\psi|^2 \div_\Si(\nu_\Si)\\
      &=& r^{-2}  \nu_\Si\big( r f|\psi|^2 \big)  -2r^{-2}  f |\psi|^2\nu_\Si(r)  + r^{-1} f|\psi|^2 \div_\Si(\nu_\Si).
     \eeaa
     Therefore,
     \beaa
      r^{-2}  \nu_\Si\big( r f|\psi|^2 \big)&=&    \div_\Si\big(r^{-1} f|\psi|^2 \nu_\Si\big)- r^{-1} f|\psi|^2 \div_\Si(\nu_\Si) +2 fr^{-2} |\psi|^2\nu_\Si(r)\\
      &=&  \div_\Si\big(r^{-1} f|\psi|^2 \nu_\Si\big)+fr^{-1} |\psi|^2\big(2 r^{-1}\nu_\Si(r)-  \div_\Si(\nu_\Si) \big).
      \eeaa
    Note that $ \nu_\Si(r)= e_4(r)- \frac 1 2 r^{-2} \la  e_3(r) =1+O(r^{-2})  + r\Ga_g$ and
    \beaa
    \div_\Si (\nu_\Si)&=& \trch +O(r^{-1}+\ep )  r^{-1} = \frac 2 r +O(r^{-1}+\ep )  r^{-2},\\
     2\nu_\Si(r)- r^{-1} \div_\Si(\nu_\Si) &=& O(r^{-1}+\ep )  r^{-2}.
    \eeaa
    Hence
    \beaa
      r^{-2}  \nu_\Si\big( r f|\psi|^2 \big)&=&    \div_\Si\big(r^{-1} f|\psi|^2 \nu_\Si\big) +O(r^{-1}+\ep ) fr^{-3} |\psi|^2.
    \eeaa 
    We deduce 
    \beaa
    \PP\c N_\Si  -\PP\c Y&\ge & f\left|\nabcheck_4\psi\right|^2 +\frac 1 2 r^{-2} \la    f |\nab\psi|^2    - \frac 1 2  \div_\Si\big(r^{-1} f|\psi|^2 \nu_\Si\big)  \\
    &&+O(\ep+R^{-1} ) r^{-3}  f|\psi|^2.
    \eeaa
    In view of the estimate for $  \big|\PP\c  Y\big|  $  of Lemma \ref{lemma:pointwiseBoun-Chap10}, we obtain
    \beaa
     \PP\c N_\Si &\ge & f\left|\nabcheck_4\psi\right|^2 +\frac 1 2 r^{-2} \la   f |\nab\psi|^2   - \frac 1 2  \div_\Si\big(r^{-1} f|\psi|^2 \nu_\Si\big)  \\
&&+O(\ep+R^{-1} ) r^{-3}  f|\psi|^2- f   \big|\nabcheck_4 \psi \big|  |Y| |   \nab \psi|\\
&&+ O(R^{-\de}) r^{-2} f|\nab\psi|^2 +
 r^{-4+\de} f\Big(|\nab_3\psi|^2 +|\nab\psi|^2 +  |\psi|^2 \Big)\\
 &&+ O(\ep +R^{-1} ) r^{-3}  f |\psi|^2. 
    \eeaa
    Note that
    \beaa
     f\left|\nabcheck_4\psi\right|^2 +\frac 1 2 r^{-2} \la    f |\nab\psi|^2 - f   \big|\nabcheck_4 \psi \big|  |Y| |   \nab \psi|\ges  f\left|\nabcheck_4\psi\right|^2+ f r^{-2}|\nab\psi|^2 .
     \eeaa
    Hence, for $R$ large,
    \beaa
     \PP\c N_\Si &\ges &  f\left|\nabcheck_4\psi\right|^2 +r^{-2}    f |\nab\psi|^2 + O(\ep +R^{-1})  r^{-3}  f |\psi|^2  - \frac 1 2  \div_\Si\big(r^{-1} f|\psi|^2 \nu_\Si\big) \\
 &&-fO(r^{-4+\de}) |\nab_3\psi|^2,
    \eeaa
    as stated.

   To prove \eqref{eq:Prop.pointwiseBoun-Chap10-2}  we write, recalling that
    $N_*= e_4+Ue_3 +Y_*$  and $\nu_*=e_4- U e_3$,    
       \beaa
     \PP\c  N_*&=&  U\PP\c  e_3 +\PP\c  e_4 +\PP \c Y_*\\
     &=&U\left(  f  \QQ_{34}+    f\frac 1 2 r^{-2}\la \QQ_{33}  + \frac 1 2 r^{-2}\nab_3\big ( r    f |\psi|^2) + (\ep+R^{-1}) f r^{-3} |\psi|^2\right)\\
     &&+ \left( f\left|\nabcheck_4\psi\right|^2  - \frac 1 2 r^{-2}  \nab_4( r  f |\psi|^2) +\frac 1 2 r^{-2}\la  f   \QQ_{34} +O(\ep+R^{-1} ) r^{-3}  f|\psi|^2 \right)\\
   &=&   f  \left( \big(U    +\frac 1 2 r^{-2}\la \big)  \QQ_{34} + \left|\nabcheck_4\psi\right|^2 + \frac 1 2 r^{-2}\la|\nab_3 \psi|^2 \right)+\frac 12 U f'  r^{-2} |\psi|^2 \\
   && +\frac 1 2 fr^{-2} \Big( \nab_{\nu_*}\big ( r    f |\psi|^2)\Big).
    \eeaa
      As  in the proof of  part 1,   see also the proof of Lemma 10.44 in \cite{KS}, 
    we have
    \beaa
    \frac 1 2 r^{-2} \Big( \nab_{\nu_*}\big ( r    f |\psi|^2)\Big)=\frac 1 2 \div_{\Si_*} \big(r^{-1} f|\psi|^2\nu_* \big)+O(\ep+R^{-1} ) r^{-3} f |\psi|^2 
    \eeaa
     and therefore,  
    \beaa
     \PP\c N_* &\ges&  f\big(|\nabcheck_4\psi|^2 + |\nab\psi|^2 +r^{-2}  |\nab_3\psi|^2\big)    - \frac 1 2 \div_\Si\big(r^{-1} f|\psi|^2 \nu_\Si\big) \\
&&   + O(R^{-1}+\ep) r^{-3}|\psi|^2 
\eeaa
      as stated.

    It remains to derive  the last two  parts. Starting with the identities in Lemma \ref{lemma:pointwiseBoun-Chap10}
    \beaa
\PP\c  e_4 &=&f\big|\nab_4 \psi|^2  + f r^{-1}  \psi\nab_4 \psi  -\frac 1  2  e_4 ( f r^{-1}) |\psi|^2+\frac 1 2 r^{-2}\la   f  \QQ_{34} +O(\ep+R^{-1} ) r^{-3}  f|\psi|^2, \\
\PP\c e_3 &=& f  \QQ_{34}+  f\frac 1 2 r^{-2}\la \QQ_{33}+  r^{-1}   f \psi\c\nab_3 \psi  -\frac 1  2  e_3  ( r^{-1}   f  ) |\psi|^2 -\frac{r}{|q|^2} f'  |\psi|^2\\
&&+O(\ep+R^{-1} ) r^{-3}  f|\psi|^2,
\eeaa
 and combining them  along  $\Si$ we derive
\beaa
\PP\c N_\Si&=& f\big|\nab_4 \psi\big|^2+ f r^{-1}\psi \c  \nu_{\Si} \psi -\frac 1 2 \nu_{\Si} ( f r^{-1}) |\psi|^2 +\frac 1 2 r^{-2} \la    f  \QQ_{34}+ f  \frac 1 4 r^{-4} \la^2 |\nab_3 \psi|^2 \\
&&+ O(R^{-1}+\ep)  r^{-2}f|\psi|^2       +\PP\c Y.
\eeaa
Note that, since $\nu_\Si= e_4-\frac 1 2 r^{-2} \la$,  
\beaa
\big|\nab_4 \psi\big|^2 +   \frac 1 4 r^{-4} \la^2  |\nab_3 \psi|^2        =   \big|\nab_\Si  \psi \big|^2  +  r^{-2}\la   \nab_4 \psi  \c  \nab_3 \psi.
\eeaa
  Hence
  \beaa
  \PP\c N_\Si&=& f    \big|\nab_\Si  \psi \big|^2  + f r^{-1}\psi \c  \nu_{\Si} \psi -\frac 1 2 \nu_{\Si} ( f r^{-1}) |\psi|^2\\
  &&+\frac 1 2 r^{-2}\la   f  \QQ_{34}+ f  r^{-2} \la    \nab_4 \psi  \c  \nab_3 \psi + O(R^{-1}+\ep)  r^{-2}|\psi|^2       +\PP\c Y.
  \eeaa
  Since   
  \beaa
     r^{-2} \la    \nab_4 \psi  \c  \nab_3 \psi &=&    r^{-2} \la \nab_\nu\psi \nab_3 \psi -O(r^{-4})  |\nab_3\psi|^2
    \les  r^{-1}|\nab_\nu\psi|^2 + O(r^{-3})|\nab_3\psi|^2 
   \eeaa
     we deduce\footnote{Note that we  also  estimate  $\PP\c Y$ slightly differently than  in  Lemma \ref{lemma:pointwiseBoun-Chap10}.}   
   \beaa
    \PP\c N_\Si&\ge &  f (1- r^{-1}  )    \big|\nab_\Si  \psi \big|^2  + f r^{-1}\psi \c  \nu_{\Si} \psi -\frac 1 2 \nu_{\Si} ( f r^{-1}) |\psi|^2\\
    &&+\frac 1 2 r^{-2} \la    f  \QQ_{34}+O(r^{-3} )  f\Big( |\nab_3 \psi|^2  +|\nab\psi|^2 +|\psi|^2\Big).
   \eeaa

   \begin{lemma}
       We re-express 
       \beaa
       J=  f (1- r^{-1}  )  \big|\nab_\Si  \psi \big|^2  + f r^{-1}\psi \c  \nu_{\Si} \psi -\frac 1 2 \nu_{\Si} ( f r^{-1}) |\psi|^2,
       \eeaa 
       for $f= r^p$, $0\le p\le 1-\de$,  $r\ge R$  sufficiently large,  as follows
        \beaa
        J+\frac p 2 r^{-2}  \nu_\Si(  r f |\psi|^2)\ge    \frac{\de^2}{8}    r^{p-2} |\psi|^2 .
      \eeaa
     \end{lemma}
     
 \begin{proof}
      We calculate      
      \beaa
     && J+\frac \la 2 r^{-2}  \nu_\Si(  r f |\psi|^2)\\
     &=& f  \big(1- r^{-1}\big) |\nu_\Si\psi|^2 +  r^{-1} f   \psi\c  \nu_\Si(\psi) 
        -\frac 1  2 \nu_\Si(r^{-1} f ) |\psi|^2+ \frac \la 2 r^{-2}  \nu_\Si(  r f |\psi|^2)\\
     &=&  f  \big(1- r^{-1}\big) |\nu_\Si\psi|^2+(1+\la) r^{-1} f \psi\c \nu_\Si\psi +\left(\frac \la 2 r^{-2} \nu_\Si(r f)-\frac 1 2   \nu_\Si(r^{-1} f )\right) |\psi|^2\\
     &=&f\left| \big(1- r^{-1}\big)^{1/2}\nu_\Si \psi  + \frac{1+\la}{ 2 \big(1- r^{-1}\big)^{1/2}  } r^{-1}  \psi\right|^2-\frac{(1+\la)^2}{  4 \big(1- O(r^{-1})\big) } r^{-2} f  |\psi|^2 \\
      && -\frac 1  2 \nu_\Si(r^{-1} f ) |\psi|^2+ \frac \la 2 r^{-2}  \nu_\Si(  r f) |\psi|^2\\
      &\ge & L|\psi|^2,
     \\
     L&:=&  \left( -\frac{(1+\la)^2}{  4 \big(1- r^{-1}\big) } r^{-2} f 
       -\frac 1  2 \nu_\Si(r^{-1} f ) + \frac \la 2 r^{-2}  \nu_\Si(  r f) \right) |\psi|^2 .
      \eeaa
           Using $\nu_{\Si}(r)= 1+ r\Ga_g+ r^{-1} \Ga_b = 1 +r^{-1} \ep $   we deduce for    $f= r^p$,         $r\ge R$,
      \beaa
      L&=&   r^{p-2}\left(- \frac{(1+\la)^2}{  4 } \big(1+ O(r^{-1})\big) +\frac{\la (p+1)}{2} -\frac{p-1}{2} \right)\\
      &=& r^{p-2}\left(- \frac{(1+\la)^2}{  4 }  +\frac{\la (p+1)}{2} -\frac{p-1}{2} \right) + O(r^{p-3})  \frac{(1+\la)^2}{4}.
      \eeaa
     For  $\la=p$, $0\le p\le 1-\de$ we derive
     \beaa
     L&=& r^{p-2}\left(\frac 1 4 (p-1)^2 +O(R^{-2} )\frac{(p+1)^2}{4} \right)\ge r^{-1-\de} \left(\frac{\de^2}{ 4} - O(R^{-2} )  \right).
     \eeaa
     We need $R^{-1} \les  \frac{\de}{ 2}$ to deduce, for $p\le 1-\de$,
     \beaa
        J+\frac p 2 r^{-2}  \nu_\Si(  r f |\psi|^2)\ge    \frac{\de^2}{8}    r^{p-2} |\psi|^2 
      \eeaa
       as stated.
 \end{proof}
 
  Using the above lemma we deduce
  \beaa
     \PP\c N_\Si&\ge &     \frac{\de^2}{8}    r^{p-2}f  |\psi|^2 
     +   f r^{-2}  |\nab \psi|^2  - \frac p 2 r^{-2}  \nu_\Si(  r f |\psi|^2)+O(r^{p-3} )f \Big( |\nab_3 \psi|^2   +|\psi|^2\Big),
     \eeaa
     as stated. The last inequality can be derived in the same manner.  This concludes the proof of Proposition \ref{Prop:pointwiseBoun-Chap10}.       
     \end{proof}


\subsection{Proof of Proposition \ref{Proposition:Step3-Chap10} in the case $s=0$}
\label{sec:proof-Step3-Chap10}


The goal of this section is to prove Proposition \ref{Proposition:Step3-Chap10} in the case $s=0$\footnote{The case of higher order derivatives $0\leq s\leq \kl$ is postponed to section \ref{section:HigherDerivEsti-chap10}.} for solutions  of the wave equation  \eqref{eq:Gen.RW-chap10}
  \beaa
  \squared_2\psi - V \psi = -\frac{4 a\cos\th}{|q|^2}\dual \nab_T  \psi  +N.
  \eeaa
   Using the estimates already derived in Propositions  \ref{proposition:MainIdentity-chap10}
  and \ref{Prop:pointwiseBoun-Chap10}, the proof  is very similar  to that of  Theorem 10.37 in  section 10.2.3  of \cite{KS}.
 
The proof of Proposition \ref{Proposition:Step3-Chap10} in the case $s=0$ proceeds in the following steps.

{\bf Step 1.} We start by  integrating    the expression \eqref{eq:DivergencePP-chapter10-3}, i.e.
\beaa
      \begin{split}
& \D^\mu  \PP_\mu[X, w, M]  =\EE[X, w, M]   + \left(\nab_X\psi+\frac 1 2   w \psi\right)\c  \big(\squared_2\psi - V \psi \big) \\
& +O(\ep+R^{-1} )  r^{-1}  f \Big(|\nab_4\psi|^2 +|\nab\psi|^2 + r^{-2}  |\nab_3\psi|^2 + r^{-2}|\psi|^2 \Big)
\end{split}
\eeaa
with, see \eqref{eq:MainIdentity-chap10-2},
\beaa
\bsplit
\EE&=\frac 1 2  f' |\nabcheck_4  \psi|^2  +\frac 1 2 \left( - f'+\frac{2}{r} f \right)\big(|\nab\psi|^2+V|\psi|^2 \big) \\
&+O(\ep+R^{-1} )r^{-1} f\Big( |\nab_4 \psi|^2+|\nab\psi|^2 
+r^{-2}|\psi|^2 + r^{-2} |\nab_3\psi|^2 \Big)
\end{split}
 \eeaa
 for  the specific  choice  $f= f_p=f_{p, R}$,  defined as $f_p=r^p$ for $r \geq R$ and $f_p=0$ for $r \leq R/2$, where $R$ is a fixed sufficiently large constant.  

   We derive
   \bea
     \lab{eq:Diveregencefor-Chatpter10}
     \bsplit
  & \int_{\Si (\tau_2)}   \PP\c  N_\Si  +\int_{\Sigma_*(\tau_1,\tau_2)}\PP \c N_{\Si_*} +\int_{\MM(\tau_1, \tau_2)} \EE=\int_{\Si (\tau_1)}   \PP\c   N_\Si + \err,\\
  & \err(\tau_1,\tau_2):=- \int_{\MM(\tau_1,\tau_2)} f_p \nabcheck_4 \psi\c N
   -  \int_{\MM(\tau_1,\tau_2)} f_p \nabcheck_4\psi \c  \frac{4 a\cos\th}{|q|^2}\dual \nab_T  \psi.
   \end{split}
   \eea

 Denoting  the boundary terms,   
     \beaa
     K_{\ge R}(\tau_1, \tau_2):&=&\int_{\Si_{\ge R}  (\tau_2)}   \PP\c  N_\Si        +\int_{\Sigma_*(\tau_1,\tau_2)}\PP \c N_{\Si_*}-\int_{\Si_{\ge R} (\tau_1)}   \PP\c   N_\Si,\\
        K_{\le R}(\tau_1, \tau_2):&=&\int_{\Si_{\le R} (\tau_1)}   \PP\c   N_\Si- \int_{\Si_{\le R} (\tau_2)}   \PP\c   N_\Si,
     \eeaa
     we   write,
     \beaa
     \bsplit
      K_{\ge R} (\tau_1, \tau_2)+\int_{\MM_{\ge R} (\tau_1, \tau_2)}\EE         &=   K_{\le R}(\tau_1, \tau_2)    -\int_{\MM_{\le R} (\tau_1, \tau_2)}\EE +\err(\tau_1, \tau_2).
           \end{split}
     \eeaa
     We estimate  the term $ K_{\le R}(\tau_1, \tau_2)    -\int_{\MM_{\le R} (\tau_1, \tau_2)}\EE$ on the right hand side    and deduce (see also  Lemma 10.45  in \cite{KS})
    \bea
    \lab{eq:PropStep2-Chap10-1}
     \nn K_{\ge R} (\tau_1, \tau_2)+\int_{\MM_{\ge R} (\tau_1, \tau_2)}\EE & \les &  R^{p+2}\Big(E_{R/2 \le r\le R}[\psi](\tau_1)     
      +\Mor_{R/2 \le r\le R}[\psi](\tau_1, \tau_2)   \Big)\\
      &&+\err(\tau_1, \tau_2).
\eea

{\bf Step 2.} We make use of the  second  identity in Proposition  \ref{proposition:MainIdentity-chap10}
 with $f=f_p$ to deduce, exactly as 
 in the Lemma 10.46 in \cite{KS}, for $\de\le p\le 2-\de $ and $R$ sufficiently large, 
 \bea
 \lab{eq:PropStep2-Chap10-2}
\nn \int_{\MM_{\ge R}(\tau_1,\tau_2)  }\EE&\geq &\frac 1  4  \int_{\MM_{\ge R}(\tau_1,\tau_2)  }r^{p-1}\Big(  p   |\nabcheck_4(\psi)|^2+ (2-p) (|\nab \psi|^2 + r^{-2}|\psi|^2)\Big)\\
 &-&  O(\ep+R^{-1} )\int_{\MM_{\ge R}(\tau_1,\tau_2)  }  r^{p-3}|\nab_3 \psi|^2.
 \eea
 Indeed, according to \eqref{eq:MainIdentity-chap10-2}, since $f=r^p$ for $r\ge R$,
  \beaa
\bsplit
\EE&= \frac 1 2  r^{p-1}\Big(  p  |\nabcheck_4  \psi|^2  +(2-p)\big(|\nab\psi|^2+V|\psi|^2 \big) \Big)   \\
&  +O(\ep+R^{-1} )  r^{p-1}\Big( |\nabcheck_4 \psi|^2+|\nab\psi|^2 
+r^{-2}|\psi|^2 + r^{-2} |\nab_3\psi|^2 \Big)\Big)\\
&\ge   \frac 1 4 r^{p-1}\Big(  p  |\nabcheck_4  \psi|^2  +(2-p)\big(|\nab\psi|^2+r^{-2} |\psi|^2 \big) \Big)+ O(\ep+R^{-1} )  r^{p-3}|\nab_3 \psi|^2,
\end{split}
 \eeaa
and  integrating  on $\MM(\tau_1, \tau_2)$ immediately yields \eqref{eq:PropStep2-Chap10-2}.

  {\bf Step 3.} According to  Proposition \ref{Prop:Step1-Chap10},  we have
 \beaa
\bsplit
\int_{\MM_{\geq R}(\tau_1, \tau_2)}r^{-1-\de }|\nab_3 \psi|^2 \les&   \int_{\MM_{\geq R}(\tau_1, \tau_2)}  r^{\de-1}     \Big(|\nab_4 \psi|^2+|\nab \psi|^2+r^{-2}|\psi|^2\Big)\\
& +E_{\geq \frac{R}{2}}[\psi](\tau_1)+Mor_{\frac{R}{2}\leq r\leq R}[\psi](\tau_1, \tau_2)+\NN_{\ge R/2}[\psi, N](\tau_1, \tau_2). 
\end{split}
\eeaa
Together with \eqref{eq:PropStep2-Chap10-2},  we deduce, for $\de\le p\le 2-\de$ and for $\ep$ and $R^{-1}$ small enough,
  \bea
  \lab{eq:PropStep2-Chap10-3}
  \bsplit
  B_{p, \ge R}[\psi](\tau_1, \tau_2)\les&  \int_{\MM_{\ge R}(\tau_1,\tau_2)  }\EE+Mor_{\frac{R}{2}\leq r\leq R}[\psi](\tau_1, \tau_2)\\
  &+\NN_{\ge R/2}[\psi, N](\tau_1, \tau_2) +E_{\geq \frac{R}{2}}[\psi](\tau_1)
  \end{split}
  \eea
  where,  see  section \ref{subsection:basicnormsforpsi},
  \beaa
   B_{p, \ge R}[\psi](\tau_1, \tau_2)&=&\int_{\MM_{\geq R}(\tau_1, \tau_2)}\Big(r^{p-1}\big(     |\nab_4 \psi|^2+ |\nab \psi|^2   + r^{-2}|\psi|^2\big) + r^{-1-\de} |\nab_3 \psi|^2\Big).
  \eeaa

{\bf Step 4.}  We treat the boundary terms 
\beaa
 K_{\ge R} (\tau_1, \tau_2)= \int_{\Si_{\ge R}  (\tau_2)}   \PP\c  N_\Si        +\int_{\Sigma_*(\tau_1,\tau_2)}\PP \c N_{\Si_*}-\int_{\Si_{\ge R} (\tau_1)}   \PP\c   N_\Si,
 \eeaa by making use 
of Proposition \ref{Prop:pointwiseBoun-Chap10} (see also    Step 2  in  the proof of 
Theorem 10.37 in \cite{KS}).   Integrating  the inequality \eqref{eq:Prop.pointwiseBoun-Chap10-1}  we deduce
\beaa
 \int_{\Si_{\ge R}  (\tau_2)}   \PP\c N_\Si&\ges  &\int_{\Si_{\ge R}  (\tau_2)}   r^p\Big(    |\nabcheck_4\psi|^2 +r^{-2}    |\nab\psi|^2 \Big)  - \frac 1 2  \int_{\Si_{\ge R}  (\tau_2)} \div_\Si\big(r^{-1}f|\psi|^2 \nu_\Si\big)\\
 &&-O(1) \int_{\Si_{\ge R}  (\tau_2)}   r^{p-4+\de}|\nab_3\psi|^2  -O(\ep +R^{-1} )  \int_{\Si_{\ge R}  (\tau_2)} r^{p-3}|\psi|^2. 
\eeaa
Since $p-4+\de\le -2 $ and, see section \ref{subsection:basicnormsforpsi},  
$
E_{ \ge R}[\psi](\tau) \ge \int_{\Si_{\ge R}(\tau) } r^{-2} |\nab_3\psi|^2
$ 
 we deduce
\beaa
 \int_{\Si_{\ge R}  (\tau_2)}   \PP\c N_\Si&\ges  & \int_{\Si_{\ge R}  (\tau_2)}   r^p\Big(    |\nabcheck_4\psi|^2 +r^{-2}     |\nab\psi|^2 \Big) + O(R^{-1}+\ep) \int_{\Si_{\ge R}  (\tau_2)}   r^{p-3} |\psi|^2\\
 &&   - O(1)  E_{\ge R}[\psi](\tau_2)- \frac 1 2  \int_{\Si  (\tau_2)} \div_\Si\big(r^{-1}f|\psi|^2 \nu_\Si\big).
\eeaa
 Similarly, integrating \eqref{eq:Prop.pointwiseBoun-Chap10-2} on $\Si_*(\tau_1, \tau_2)$,  we  deduce  
   \beaa
          \bsplit
  \int_{\Si_*(\tau_1, \tau_2)} \PP\c N_{\Si_*} &\ge \frac 1 2  \int_{\Si_*(\tau_1, \tau_2)}   r^p \left(|\nabcheck_4\psi|^2 +\frac{m}{r^2}    |\nab\psi|^2  +   \frac{m^2}{2r^4}|\nab_3\psi|^2\right)     \\
&  - \frac 1 2    \int_{\Si_*(\tau_1, \tau_2)} \div_\Si\big(r^{-1} f|\psi|^2 \nu_\Si\big) + O(R^{-1}+\ep)r^{p-3}|\psi|^2.
\end{split}
   \eeaa   
    Applying the divergence theorem on $\Si(\tau_2)\cup\Si_*(\tau_1, \tau_2)$ to $\div_\Si(r^{-1} f|\psi|^2 \nu_\Si)$, and noticing that the corresponding boundary terms from $\Si(\tau)$ and $\Si_*(\tau_1, \tau_2)$ cancel each other at $\Si(\tau)\cap\Si_*$,  we deduce, for $R$ sufficiently large, $\de\le p\le 2-\de$,
    \bea
    \lab{eq:PropStep2-Chap10-4}
    \bsplit
     K_{\ge R} (\tau_1, \tau_2)
     &\ges 
   \int_{\Si_{\ge R}  (\tau_2)}   r^p\Big(   |\nabcheck_4\psi|^2 + r^{-2}   |\nab\psi|^2 \Big) +  \int_{\Si_*(\tau_1,\tau_2)}   r^p \Big(|\nabcheck_4\psi|^2 + r^{-2}   |\nab\psi|^2 \Big) \\
    &- O(R^{-1}+\ep) \int_{\Si_{\ge R}  (\tau_2)}   r^{p-3} |\psi|^2 - O(R^{-1}+\ep)\int_{\Si_*(\tau_1,\tau_2)} r^{p-3}|\psi|^2\\
    & -  E_{\ge R}[\psi](\tau_2) -  E_{p, \ge R}[\psi](\tau_1).
    \end{split}
    \eea
    
    {\bf Step 5.} Combining \eqref{eq:PropStep2-Chap10-4} with \eqref{eq:PropStep2-Chap10-3} and \eqref{eq:PropStep2-Chap10-1} we derive
    \bea
    \lab{eq:PropStep2-Chap10-5}
    \bsplit
    \BEF_{p,\ge R}[\psi](\tau_1, \tau_2)&\les   \int_{\Si_{\ge R}  (\tau_2)}   r^{p-3} |\psi|^2 +\int_{\Si_*} r^{p-3}|\psi|^2 +E_{p, \geq R}[\psi](\tau_1)\\
    &+ R^{p+2}\Big(E_{R/2 \le r\le R}[\psi](\tau_1)     
      +\Mor_{R/2 \le r\le R}[\psi](\tau_1, \tau_2)   \Big)\\
      &+E_{\ge R}[\psi](\tau_2)+\err(\tau_1, \tau_2) +\NN_{p, \ge R/2}[\psi, N](\tau_1, \tau_2).
      \end{split} 
    \eea     
        
    {\bf Step 6.}  We now estimate the error term $\err(\tau_1,\tau_2)$ which we decompose as
\beaa
\bsplit   
    \err(\tau_1,\tau_2)=&\err_1(\tau_1,\tau_2)+\err_2(\tau_1,\tau_2),\\
\err_1(\tau_1,\tau_2) :=&   - \int_{\MM(\tau_1,\tau_2)} f_p \nabcheck_4 \psi\c N,\\
   \err_2(\tau_1,\tau_2) :=&   -  \int_{\MM(\tau_1,\tau_2)} f_p \nabcheck_4\psi \c  \frac{4 a\cos\th}{|q|^2}\dual \nab_T  \psi.
\end{split}
\eeaa    
First, in view of the definition  of the  $\NN_p$ norms in section
   \ref{subsection:basicnormsforpsi}, we have
   \beaa
   \big|\err_1(\tau_1,\tau_2)\big|&\les&\NN_{p, \ge R/2}[\psi, N](\tau_1,\tau_2).
   \eeaa    
   Next, we focus on the control of $\err_2(\tau_1, \tau_2)$.  Recalling  that the vectorfield $\T$ is given by\footnote{The formula for $\T$ takes into account the renormalization of the frame, see Remark \ref{rmk:fromnowonoutgoingrenormalization}.} $\T=  \frac 1 2( \frac{\De}{|q|^2} e_4+  e_3 -2a\Re(\Jk)^be_b)$, 
 we deduce
  \beaa
  - \frac{4 a\cos\th}{|q|^2} f_p\widecheck{\nab}_4\psi\c  \dual \nab_\T  \psi
  &=& - \frac{2 a\cos\th}{|q|^2} f_p\widecheck{\nab}_4\psi\c  \dual   \nab_3 \psi\\
  && +  O(r^{-2})  |f_p|  | \widecheck{\nab}_4\psi |    \c \big(|\nab_4 \psi |+ r^{-1} |\nab\psi| \big).
  \eeaa
  Hence
  \beaa
    \err_2(\tau_1, \tau_2) &=& - \int_{\MM(\tau_1,\tau_2) } f_p \nabcheck_4 \psi\c  \left(\frac{2 a\cos\th}{|q|^2}\dual \nab_3  \psi \right)\\ 
    &&+ R^{-1} \int_{\MM_{r\ge R}(\tau_1,\tau_2)} r^{p-1} \big( | \widecheck{\nab}_4\psi |^2  +   |\nab\psi|^2 \big).
  \eeaa
  To estimate the integral 
  \beaa
  I:= \left|  \int_{\MM(\tau_1,\tau_2) } f_p \nabcheck_4 \psi\c  \left(\frac{2 a\cos\th}{|q|^2}\dual \nab_3  \psi \right)\right|, 
  \eeaa
   we write, with $p\le 2-\de$,
  \beaa
I &\les&  \int_{\MM_{\ge R/ 2}(\tau_1,\tau_2) } r^{p-2}  | \widecheck{\nab}_4\psi | \,|\nab_3 \psi| \\
&\les&\left( \int_{\MM_{\ge R/ 2}(\tau_1,\tau_2) } r^{p-1}  | \widecheck{\nab}_4\psi | ^2 \right)^{\frac{1}{2}}  \left( \int_{\MM_{\ge R/ 2}(\tau_1,\tau_2) } r^{p-3}  | \nab_3 \psi | ^2 \right)^{\frac{1}{2}} \\ &\les &
 r^{-\frac{2-\de-p}{2} }  \left( \int_{\MM_{\ge R/ 2}(\tau_1,\tau_2) } r^{p-1}  | \widecheck{\nab}_4\psi | ^2 \right)^{\frac{1}{2}} \left( \int_{\MM_{\ge R/ 2}(\tau_1,\tau_2) } r^{-1-\de}  | \nab_3 \psi | ^2 \right)^{\frac{1}{2}}. 
  \eeaa
Plugging the above estimates for $\err_1(\tau_1, \tau_2)$ and $\err_2(\tau_1, \tau_2)$ in \eqref{eq:PropStep2-Chap10-5}, and making use once more of Proposition \ref{Prop:Step1-Chap10}  we deduce
 \bea
    \lab{eq:PropStep2-Chap10-6}
    \bsplit
    \BEF_{p,\ge R}[\psi](\tau_1, \tau_2)&\les   \int_{\Si_{\ge R}  (\tau_2)}   r^{p-3} |\psi|^2 +\int_{\Si_*(\tau_1, \tau_2)} r^{p-3}|\psi|^2 +E_{p, \geq R}[\psi](\tau_1)\\
    &+ \NN_{p, \ge R/2}[\psi, N](\tau_1, \tau_2) + E_{\ge R}[\psi](\tau_2) \\
    &+ R^{p+2}\Big(E_{R/2 \le r\le R}[\psi](\tau_1) +    
      \Mor_{R/2 \le r\le R}[\psi](\tau_1, \tau_2)   \Big).
        \end{split} 
    \eea

    {\bf Step 7.}    Next, we eliminate  the term  
    \beaa
    I_p(\tau_1, \tau_2):= \int_{\Si_{\ge R}  (\tau_2)}   r^{p-3} |\psi|^2 +\int_{\Si_*(\tau_1, \tau_2)} r^{p-3}|\psi|^2
    \eeaa 
    on the right hand side of \eqref{eq:PropStep2-Chap10-6}. Note that  for $p\le 1$   we have
    \beaa
    I_p(\tau_1, \tau_2)&\les&  E_{\ge R}[\psi](\tau_2) +F[\psi](\tau_1, \tau_2). 
    \eeaa
    Hence,  for $p \le 1$,
  \bea
    \lab{eq:PropStep2-Chap10-7}
    \bsplit
    \BEF_{p,\ge R}[\psi](\tau_1, \tau_2)&\les  
    \EF_{ \ge R}[\psi](\tau_1, \tau_2) +E_{p, \geq \frac{R}{2}}[\psi](\tau_1)+ \NN_{p, \ge R/2}[\psi, N](\tau_1, \tau_2) \\
    &+ R^{p+2}\Big(E_{R/2 \le r\le R}[\psi](\tau_1) +    
      \Mor_{R/2 \le r\le R}[\psi](\tau_1, \tau_2)   \Big).
        \end{split} 
    \eea      
    For  the remaining range  $1\le  p\le 2-\de$   we have
    \beaa
      I_p(\tau_1, \tau_2)&\les&  I_{2-\de}(\tau_1, \tau_2) =\int_{\Si_{\ge R}  (\tau_2)}   r^{-1-\de} |\psi|^2 +\int_{\Si_*} r^{-1-\de}|\psi|^2\\
      &\les& \EF_{1-\de, \ge R}[\psi](\tau_1, \tau_2)
    \eeaa 
    which together with \eqref{eq:PropStep2-Chap10-7} implies, for $1\le  p\le 2-\de$, 
    \beaa
    \bsplit
    I_p(\tau_1, \tau_2)\les&  \EF_{ \ge R}[\psi](\tau_1, \tau_2) +E_{p, \geq R}[\psi](\tau_1)+ \NN_{1-\de, \ge R/2}[\psi, N](\tau_1, \tau_2) \\
    &+ R^{3-\de}\Big(E_{R/2 \le r\le R}[\psi](\tau_1) +    
      \Mor_{R/2 \le r\le R}[\psi](\tau_1, \tau_2)   \Big).
       \end{split} 
    \eeaa
     Combining with \eqref{eq:PropStep2-Chap10-6}  we  deduce,  for $\de\le p\le 2-\de$,
       \bea
    \lab{eq:PropStep2-Chap10-8}
    \bsplit
    \BEF_{p,\ge R}[\psi](\tau_1, \tau_2)&\les \EF_{ \ge R}[\psi](\tau_1, \tau_2) +E_{p, \geq R}[\psi](\tau_1)+ \NN_{p, \ge R/2}[\psi, N](\tau_1, \tau_2) \\
    &+ R^{p+2}\Big(E_{R/2 \le r\le R}[\psi](\tau_1) +    
      \Mor_{R/2 \le r\le R}[\psi](\tau_1, \tau_2)   \Big).
        \end{split} 
    \eea  
    
     {\bf Step 8.} It remains to eliminate the term $\EF_{ \ge R}[\psi](\tau_1, \tau_2)$ on the RHS of \eqref{eq:PropStep2-Chap10-8}. We rely on Proposition \ref{prop:recoverEnergyMorawetzwithrweightfromnoweight:perturbation} which yields
   \beaa
    \bsplit
   \EF_{ \ge R}[\psi](\tau_1, \tau_2) \les& E_{\ge R}[\psi](\tau_1)+ \NN_{p, \ge R/2}[\psi, N](\tau_1, \tau_2) +\ep B_\de[\psi](\tau_1, \tau_2)\\
    &+ R\Big(E_{R/2 \le r\le R}[\psi](\tau_1) +    
      \Mor_{R/2 \le r\le R}[\psi](\tau_1, \tau_2)   \Big).
       \end{split} 
   \eeaa  
     Together with \eqref{eq:PropStep2-Chap10-8}, we infer
        \bea\lab{eq:PropStep2-Chap10-9}
    \bsplit
    \BEF_{p,\ge R}[\psi](\tau_1, \tau_2)&\les E_{p, \ge R}[\psi](\tau_1)+ \NN_{p, \ge R/2}[\psi, N](\tau_1, \tau_2) \\
    &+ R^{p+2}\Big(E_{R/2 \le r\le R}[\psi](\tau_1) +    
      \Mor_{R/2 \le r\le R}[\psi](\tau_1, \tau_2)   \Big)
      \end{split} 
    \eea
      which concludes the proof of Proposition  \ref{Proposition:Step3-Chap10} in the case $s=0$.

                
\subsection{Proof  of Proposition  \ref{Proposition:Step3-Chap10}} 
\lab{section:HigherDerivEsti-chap10}
 
    
In order to prove Proposition  \ref{Proposition:Step3-Chap10}, we proceed as in section 10.4 in \cite{KS}, i.e. we 
 extend the estimates derived in the previous  section for $s=0$ to $\dk^s$ derivatives of $\psi$ with $0\leq s\leq \kl$ by recovering the derivatives one by one. We indicate below how to go from $s=0$ to $s=1$ and note that the procedure to recover the estimate for $s+1$ from the one for $s$ is completely analogous.
     
To derive the  estimate for    $s=1$ we proceed in the following steps, see also section 10.4 in \cite{KS}.

{\bf Step 1.} We start with the following result, in the spirit of Lemma \ref{lemma:controloftheerrortermsGagpsiinwaveeq}, which will be used to deal with terms generated by  commutators with $\squared_2$.
\begin{lemma}\lab{lemma:controloftheerrortermsGagpsiinwaveeq:rpcasechap10}
Let $\widetilde{N}$ a tensor with the following schematic structure 
\beaa
\widetilde{N} &=& O(r^{-2})\dk F+O(r^{-2})\dk^{\leq 1}\psi +O(r^{-3})\dk^{\leq 2}\psi +\dk^{\leq 2}(\Ga_g\c\psi),
\eeaa
where $F$ is a given tensor. Also, let $\psi^{(1)}$ a tensor satisfying $|\psi^{(1)}|\les |\dk^{\leq 1}\psi|$. Then, we have, for $\de\leq p\leq 2-\de$, 
\beaa
&&\NN_{p, \ge R/2}[\psi^{(1)}, \widetilde{N}](\tau_1, \tau_2)\\
&\les& \sqrt{B_{p, \ge R/2}[\psi^{(1)}](\tau_1, \tau_2)}\left(\sqrt{B_{p, \ge R/2}[F](\tau_1, \tau_2)}+\sqrt{B_{p, \ge R/2}[\psi](\tau_1, \tau_2)}\right)\\
&& +(\ep+R^{-1})B^1_{p, \ge R/2}[\psi](\tau_1, \tau_2).
\eeaa
\end{lemma}

\begin{proof}
In view of the definition of $\NN_p$, we have, for $\de\leq p\leq 2-\de$, 
\beaa
\NN_{p, \ge R/2}[\psi^{(1)}, \widetilde{N}](\tau_1, \tau_2) &\les& \int_{\MM(r\geq R/2)}|\widetilde{N}|\Big(|\nab_3\psi^{(1)}|+r^{p-1}|\dk^{\leq 1}\psi^{(1)}|\Big)\\
&\les& \sqrt{B_{p, \ge R/2}[\psi^{(1)}](\tau_1, \tau_2)}\left(\int_{\MM(r\geq R/2)}r^{1+p}|\widetilde{N}|^2\right)^{\frac{1}{2}}.
\eeaa
Now, in view of the structure of $\widetilde{N}$ and the assumptions on $\Ga_g$, we have
\beaa
\int_{\MM(r\geq R/2)}r^{1+p}|\widetilde{N}|^2 &\les& \int_{\MM(r\geq R/2)}r^{p-3}\Big(|\dk F|^2+|\dk^{\leq 1}\psi|^2+O(\ep^2+R^{-2})|\dk^{\leq 2}\psi|^2\Big)\\
&\les& B_{p, \ge R/2}[F](\tau_1, \tau_2)+B_{p, \ge R/2}[\psi](\tau_1, \tau_2)\\
&&+O(\ep^2+R^{-2})B^1_{p, \ge R/2}[\psi](\tau_1, \tau_2).
\eeaa
Plugging the second estimate in the first one implies  
\beaa
&&\NN_{p, \ge R/2}[\psi^{(1)}, \widetilde{N}](\tau_1, \tau_2)\\
&\les& \sqrt{B_{p, \ge R/2}[\psi^{(1)}](\tau_1, \tau_2)}\left(\sqrt{B_{p, \ge R/2}[F](\tau_1, \tau_2)}+\sqrt{B_{p, \ge R/2}[\psi](\tau_1, \tau_2)}\right)\\
&& +(\ep+R^{-1})B^1_{p, \ge R/2}[\psi](\tau_1, \tau_2)
\eeaa
where we used also  the fact that $B_{p, \ge R/2}[\psi^{(1)}](\tau_1, \tau_2)\les B_{p, \ge R/2}[\psi^{(1)}](\tau_1, \tau_2)$. This concludes the proof of Lemma \ref{lemma:controloftheerrortermsGagpsiinwaveeq:rpcasechap10}.
\end{proof}

{\bf Step 2.} Next, we derive and estimate for $\Lieb_\T\psi$. To this end, we commute  the wave equation  \eqref{eq:Gen.RW-chap10}  with $\Lieb_\T$ and obtain, using Corollary \ref{cor:commutator-Lied-squared},
\beaa
 \squared_2\Lieb_\T \psi-V\Lieb_\T \psi&=& - \frac{4 a\cos\th}{|q|^2}\dual \nab_T\Lieb_\T \psi+N_{\Lieb_\T},
\eeaa
where 
\beaa
N_{\Lieb_\T} &=& \Lieb_\T N+\dk(\Ga_g\c\dk\psi)+\Ga_b\c\squared_2\psi.
\eeaa
In view of \eqref{eq:Gen.RW-chap10}, we may rewrite $N_{\Lieb_\T}$ in the following form
\beaa
N_{\Lieb_\T} &=& \Lieb_\T N+\dk^{\leq 2}(\Ga_g\c\psi)+\Ga_b\c N.
\eeaa
Applying \eqref{eq:PropStep2-Chap10-9} to the above wave equation for $\Lieb_\T\psi$, we obtain 
  \beaa
    \bsplit
    \BEF_{p,\ge R}[\Lieb_\T\psi](\tau_1, \tau_2)&\les E^1_{p, \ge R}[\psi](\tau_1)+ \NN_{p, \ge R/2}[\Lieb_\T\psi, N_{\Lieb_\T}](\tau_1, \tau_2) \\
    &+ R^{p+2}\Big(E^1_{R/2 \le r\le R}[\psi](\tau_1) +    
      \Mor^1_{R/2 \le r\le R}[\psi](\tau_1, \tau_2)   \Big).
      \end{split} 
    \eeaa
Next, we have
\beaa
\NN_{p, \ge R/2}[\Lieb_\T\psi, N_{\Lieb_\T}](\tau_1, \tau_2) &\les& \NN_{p, \ge R/2}[\Lieb_\T\psi, \dk^{\leq 2}(\Ga_g\c\psi)](\tau_1, \tau_2)+\NN^1_{p, \ge R/2}[\psi, N](\tau_1, \tau_2)\\
&\les&  \ep B^1_{p, \ge R/2}[\psi](\tau_1, \tau_2)+\NN^1_{p, \ge R/2}[\psi, N](\tau_1, \tau_2)
\eeaa
where we have applied Lemma \ref{lemma:controloftheerrortermsGagpsiinwaveeq:rpcasechap10} in the particular case where $\widetilde{N}=\dk^{\leq 2}(\Ga_g\c\psi)$. We deduce
 \bea\lab{eq:recoverderivativeonebyonerpweightedest:chap10:LiebT}
    \bsplit
    \BEF_{p,\ge R}[\Lieb_\T\psi](\tau_1, \tau_2)&\les E^1_{p, \ge R}[\psi](\tau_1)+ \NN^1_{p, \ge R/2}[\psi, N](\tau_1, \tau_2) + \ep B^1_{p, \ge R/2}[\psi](\tau_1, \tau_2)\\
    &+ R^{p+2}\Big(E^1_{R/2 \le r\le R}[\psi](\tau_1) +    
      \Mor^1_{R/2 \le r\le R}[\psi](\tau_1, \tau_2)   \Big).
      \end{split} 
    \eea

{\bf Step 3.} Next, we derive and estimate for $\dkb\psi$. To this end, we commute  the wave equation  \eqref{eq:Gen.RW-chap10}  with $|q|\DDd_2$ and obtain, using Lemma \ref{LEMMA:COMMUTATIONOFHODGEELLIPTICORDER1WITHSQAURED2FDILUHS}, 
\beaa
 \squared_1(|q|\DDd_2\psi) -V|q|\DDd_2\psi &=& - \frac{4 a\cos\th}{|q|^2}\dual \nab_T(|q|\DDd_2\psi)+N_{|q|\DDd_2},
\eeaa
where 
\beaa
N_{|q|\DDd_2} &=& O(1)\dk^{\leq 1}N +(O(r^{-1})+r\Ga_b)\c\squared_2\psi+O(r^{-2})\dk\Lieb_\T\psi+O(r^{-2})\dk^{\leq 1}\psi \\
&&+O(r^{-3})\dk^{\leq 2}\psi+ \dk^{\leq 2} (\Ga_g \c \psi).
\eeaa
In view of \eqref{eq:Gen.RW-chap10}, we may rewrite $N_{|q|\DDd_2}$ in the following form
\beaa
N_{|q|\DDd_2} &=& O(1)\dk^{\leq 1}N +O(r^{-2})\dk\Lieb_\T\psi+O(r^{-2})\dk^{\leq 1}\psi +O(r^{-3})\dk^{\leq 2}\psi+ \dk^{\leq 2} (\Ga_g \c \psi).
\eeaa
Applying \eqref{eq:PropStep2-Chap10-9} to the above wave equation for $|q|\DDd_2\psi$, we obtain\footnote{In fact, we apply a variant for tensors in $\sk_1$. Adapting  \eqref{eq:PropStep2-Chap10-9} to this case can be done along the same lines, and is in fact easier as an estimate conditional on $\int_{\MM(r\geq R/2)}r^{p-3}|\psi|^2$ suffices for this step. See Theorem 10.61 in \cite{KS} for a proof of this variant in the case of perturbations of Schwarzschild.}
  \beaa
    \bsplit
    \BEF_{p,\ge R}[|q|\DDd_2\psi](\tau_1, \tau_2)&\les E^1_{p, \ge R}[\psi](\tau_1)+ \NN_{p, \ge R/2}[|q|\DDd_2\psi, N_{|q|\DDd_2}](\tau_1, \tau_2) \\
    &+ R^{p+2}\Big(E^1_{R/2 \le r\le R}[\psi](\tau_1) +    
      \Mor^1_{R/2 \le r\le R}[\psi](\tau_1, \tau_2)   \Big).
      \end{split} 
    \eeaa
Next, we have
\beaa
&&\NN_{p, \ge R/2}[|q|\DDd_2\psi, N_{|q|\DDd_2}](\tau_1, \tau_2)\\ 
&\les& \sqrt{B_{p, \ge R/2}[|q|\DDd_2\psi](\tau_1, \tau_2)}\left(\sqrt{B_{p, \ge R/2}[\Lieb_\T\psi](\tau_1, \tau_2)}+\sqrt{B_{p, \ge R/2}[\psi](\tau_1, \tau_2)}\right)\\
&& +(\ep+R^{-1})B^1_{p, \ge R/2}[\psi](\tau_1, \tau_2)+\NN^1_{p, \ge R/2}[\psi, N](\tau_1, \tau_2)
\eeaa
where we have applied Lemma \ref{lemma:controloftheerrortermsGagpsiinwaveeq:rpcasechap10} in the particular case where $F=\Lieb_\T\psi$. We deduce
\beaa
    \bsplit
    &\BEF_{p,\ge R}[|q|\DDd_2\psi](\tau_1, \tau_2)\\
    \les& E^1_{p, \ge R}[\psi](\tau_1)+ \NN^1_{p, \ge R/2}[\psi, N](\tau_1, \tau_2) +(\ep+R^{-1})B^1_{p, \ge R/2}[\psi](\tau_1, \tau_2)\\
    &+ R^{p+2}\Big(E^1_{R/2 \le r\le R}[\psi](\tau_1) +    
      \Mor^1_{R/2 \le r\le R}[\psi](\tau_1, \tau_2)   \Big)\\
      & \sqrt{B_{p, \ge R/2}[|q|\DDd_2\psi](\tau_1, \tau_2)}\left(\sqrt{B_{p, \ge R/2}[\Lieb_\T\psi](\tau_1, \tau_2)}+\sqrt{B_{p, \ge R/2}[\psi](\tau_1, \tau_2)}\right).
      \end{split} 
    \eeaa
Together with \eqref{eq:PropStep2-Chap10-9} and \eqref{eq:recoverderivativeonebyonerpweightedest:chap10:LiebT}, we infer
\beaa
    \bsplit
    &\BEF_{p,\ge R}[|q|\DDd_2\psi](\tau_1, \tau_2)+\BEF_{p,\ge R}[\Lieb_\T\psi](\tau_1, \tau_2)+\BEF_{p,\ge R}[\psi](\tau_1, \tau_2)\\
    \les& E^1_{p, \ge R}[\psi](\tau_1)+ \NN^1_{p, \ge R/2}[\psi, N](\tau_1, \tau_2) +(\ep+R^{-1})B^1_{p, \ge R/2}[\psi](\tau_1, \tau_2)\\
    &+ R^{p+2}\Big(E^1_{R/2 \le r\le R}[\psi](\tau_1) +    
      \Mor^1_{R/2 \le r\le R}[\psi](\tau_1, \tau_2)   \Big).
      \end{split} 
    \eeaa
Using the Hodge estimates of Proposition \ref{Prop:HodgeThmM8}, and comparing $\nab_\T$ and $\Lieb_\T$, we deduce
\bea\lab{eq:recoverderivativeonebyonerpweightedest:chap10:LiebTanddkb}
    \bsplit
    &\BEF_{p,\ge R}[(\nab_\T, \dkb)\psi](\tau_1, \tau_2)+\BEF_{p,\ge R}[\psi](\tau_1, \tau_2)\\
    \les& E^1_{p, \ge R}[\psi](\tau_1)+ \NN^1_{p, \ge R/2}[\psi, N](\tau_1, \tau_2) +(\ep+R^{-1})B^1_{p, \ge R/2}[\psi](\tau_1, \tau_2)\\
    &+ R^{p+2}\Big(E^1_{R/2 \le r\le R}[\psi](\tau_1) +    
      \Mor^1_{R/2 \le r\le R}[\psi](\tau_1, \tau_2)   \Big).
      \end{split} 
    \eea

{\bf Step 4.} Next, we derive and estimate for $r\nab_4\psi$. To this end, we commute  the wave equation  \eqref{eq:Gen.RW-chap10}  with $r\nab_4$ and obtain, using  Lemma \ref{LEMMA:COMMUTATOR-NAB3-NAB4-SQUARE}, 
\beaa
 \squared_2(r\nab_4\psi) -Vr\nab_4\psi &=& - \frac{4 a\cos\th}{|q|^2}\dual \nab_T(r\nab_4\psi)+N_{r\nab_4},
\eeaa
where 
\beaa
N_{r\nab_4} &=& \frac{1}{r}\nab_4(r\nab_4\psi)+O(1)\dk^{\leq 1}N +O(1)\c\squared_2\psi+O(r^{-2})\dkb^2\psi+O(r^{-2})\dk^{\leq 1}\psi \\
&&+O(r^{-3})\dk^{\leq 2}\psi+ \dk^{\leq 2} (\Ga_g \c \psi).
\eeaa
In view of \eqref{eq:Gen.RW-chap10}, we may rewrite $N_{r\nab_4}$ in the following form
\beaa
N_{r\nab_4} &=&  \frac{1}{r}\nab_4(r\nab_4\psi)+\widetilde{N}_{r\nab_4},\\
\widetilde{N}_{r\nab_4} &=& +O(1)\dk^{\leq 1}N +O(r^{-2})\dkb^2\psi+O(r^{-2})\dk^{\leq 1}\psi +O(r^{-3})\dk^{\leq 2}\psi + \dk^{\leq 2} (\Ga_g \c \psi).
\eeaa    
Next, we apply \eqref{eq:PropStep2-Chap10-9} to the above wave equation for $r\nab_4\psi$, and notice that the $\frac{1}{r}\nab_4(r\nab_4\psi)$ has a favorable sign in the estimate so that we may drop it and obtain 
 \beaa
    \bsplit
    \BEF_{p,\ge R}[r\nab_4\psi](\tau_1, \tau_2)&\les E^1_{p, \ge R}[\psi](\tau_1)+ \NN_{p, \ge R/2}[r\nab_4\psi, \widetilde{N}_{r\nab_4}](\tau_1, \tau_2) \\
    &+ R^{p+2}\Big(E^1_{R/2 \le r\le R}[\psi](\tau_1) +    
      \Mor^1_{R/2 \le r\le R}[\psi](\tau_1, \tau_2)   \Big).
      \end{split} 
    \eeaa
Next, we have
\beaa
&&\NN_{p, \ge R/2}[r\nab_4, \widetilde{N}_{r\nab_4}](\tau_1, \tau_2)\\ 
&\les& \sqrt{B_{p, \ge R/2}[r\nab_4\psi](\tau_1, \tau_2)}\left(\sqrt{B_{p, \ge R/2}[\dkb\psi](\tau_1, \tau_2)}+\sqrt{B_{p, \ge R/2}[\psi](\tau_1, \tau_2)}\right)\\
&& +(\ep+R^{-1})B^1_{p, \ge R/2}[\psi](\tau_1, \tau_2)+\NN^1_{p, \ge R/2}[\psi, N](\tau_1, \tau_2)
\eeaa
where we have applied Lemma \ref{lemma:controloftheerrortermsGagpsiinwaveeq:rpcasechap10} in the particular case where $F=\dkb\psi$. We deduce
 \beaa
    \bsplit
    &\BEF_{p,\ge R}[r\nab_4\psi](\tau_1, \tau_2)\\
    \les& E^1_{p, \ge R}[\psi](\tau_1)+ \NN^1_{p, \ge R/2}[\psi, N](\tau_1, \tau_2) +(\ep+R^{-1})B^1_{p, \ge R/2}[\psi](\tau_1, \tau_2)\\
    &+ R^{p+2}\Big(E^1_{R/2 \le r\le R}[\psi](\tau_1) +    
      \Mor^1_{R/2 \le r\le R}[\psi](\tau_1, \tau_2)   \Big)\\
      &+\sqrt{B_{p, \ge R/2}[r\nab_4\psi](\tau_1, \tau_2)}\left(\sqrt{B_{p, \ge R/2}[\dkb\psi](\tau_1, \tau_2)}+\sqrt{B_{p, \ge R/2}[\psi](\tau_1, \tau_2)}\right).
      \end{split} 
    \eeaa
Together with \eqref{eq:PropStep2-Chap10-9} and \eqref{eq:recoverderivativeonebyonerpweightedest:chap10:LiebTanddkb}, we infer
\beaa
    \bsplit
    &\BEF_{p,\ge R}[(\nab_\T, \dkb, r\nab_4)\psi](\tau_1, \tau_2)+\BEF_{p,\ge R}[\psi](\tau_1, \tau_2)\\
    \les& E^1_{p, \ge R}[\psi](\tau_1)+ \NN^1_{p, \ge R/2}[\psi, N](\tau_1, \tau_2) +(\ep+R^{-1})B^1_{p, \ge R/2}[\psi](\tau_1, \tau_2)\\
    &+ R^{p+2}\Big(E^1_{R/2 \le r\le R}[\psi](\tau_1) +    
      \Mor^1_{R/2 \le r\le R}[\psi](\tau_1, \tau_2)   \Big).
      \end{split} 
    \eeaa
Since $(\nab_\T, \dkb, r\nab_4)$ spans $\dk$ away from the horizon, we infer
\beaa
    \bsplit
    \BEF^1_{p,\ge R}[\psi](\tau_1, \tau_2) \les& E^1_{p, \ge R}[\psi](\tau_1)+ \NN^1_{p, \ge R/2}[\psi, N](\tau_1, \tau_2) +(\ep+R^{-1})B^1_{p, \ge R/2}[\psi](\tau_1, \tau_2)\\
    &+ R^{p+2}\Big(E^1_{R/2 \le r\le R}[\psi](\tau_1) +    
      \Mor^1_{R/2 \le r\le R}[\psi](\tau_1, \tau_2)   \Big).
      \end{split} 
    \eeaa
For $\ep$ and $R^{-1}$ small enough, we may absorb the third term on the RHS from the LHS and we obtain 
\beaa
    \bsplit
    \BEF^1_{p,\ge R}[\psi](\tau_1, \tau_2) \les& E^1_{p, \ge R}[\psi](\tau_1)+ \NN^1_{p, \ge R/2}[\psi, N](\tau_1, \tau_2) \\
    &+ R^{p+2}\Big(E^1_{R/2 \le r\le R}[\psi](\tau_1) +    
      \Mor^1_{R/2 \le r\le R}[\psi](\tau_1, \tau_2)   \Big).
      \end{split} 
    \eeaa
    which establishes the desired estimate \eqref{eq:Proposition-Step3-Chap10} for $s=1$. 
    
 We have shown above how to go from $s=0$ to $s=1$. The procedure to recover the estimate for $s+1$ from the one for $s$ is completely analogous. This concludes the proof of Proposition  \ref{Proposition:Step3-Chap10}.


\section{Proof of Theorem  \ref{THEOREM:GENRW2-Q}}


For the convenience of the reader we restate  Theorem \ref{THEOREM:GENRW2-Q}.

\begin{theorem}[Improved $r^p$-weighted estimates]
       \lab{theorem:GenRW2-q-Chap10}
       The following  estimates hold true   for  the quantity $\psiwc= f_2\left(e_4\psi+\frac{r}{|q|^2}\psi\right)$ 
        corresponding to  solutions $\psi\in\sk_2$  of  \eqref{eq:Gen.RW} on $\MM$,
        for all $-1+\de<q\le 1-\de$, $s\le\kl -1$,
         \bea
         \lab{eqtheorem:GenRW2-q-Chap10}
         \BEF^s_q[\psiwc](\tau_1, \tau_2)
       \les
       \Et_q^s[\psiwc](\tau_1)+\NNt_q^s[\psiwc, N](\tau_1, \tau_2) + \NN_{\max\{q,\de\}}^{s+1}[\psi, N](\tau_1, \tau_2).
       \eea
         where the norms on the right are given by
               \bea
       \lab{eq:TildeNormsN-chap6-chap10}
        \Et_q^s[\psiwc](\tau)=E_q^s[\psiwc](\tau)+ E^{s+1}_{\max\{q,\de\} }[\psi](\tau)
         \eea
         and
         \bea
         \bsplit
          \NNt_q^s[\psiwc, N](\tau_1, \tau_2)& =   \left|\int_{\MM_{\ge R}(\tau_1,\tau_2) } r^{q+2} \dk^{\le s} \psiwc \c \left(\nab_4 \dk^{\le s }N+ \frac 3 r \dk^{\le s} N\right)\right|.
         \end{split}
         \eea
                  \end{theorem}
                  
\begin{proof}    
We sketch the proof below  in the case $s=0$. Like the proof of  the corresponding result in \cite{KS}  (see   Theorem  5.18       in \cite{KS})  the proof  of the result rests on  a commutation formula according to which,  if $\psi$ verifies  equation \eqref{eq:Gen.RW-chap10}, then $\psiwc$ verifies an equation of the form
\bea\lab{eq:Gen.RW-chap10-psiwc}
\squared_2 \psiwc -V\psiwc=- \frac{4 a\cos\th}{|q|^2}\dual \nab_T  \psiwc+\widecheck{N} +f_2\left(\nab_4+\frac 3 r N\right), 
\eea
where
\bea
\lab{eq:checkN-chapter10}
\widecheck{N}=
\begin{cases}
&  \frac{2}{r}\left(1-\frac{3m}{r} +O(r^{-2} )\right) \nab_4\psiwc   +O(r^{-2})\dk^{\leq 1}\psi\\
 & +r\Ga_b \c \nab_4\dk\psi +\dk^{\leq 1}(\Ga_b)\c \dk^{\leq 1}\psi  +r\dk^{\leq 1}(\Ga_g)\c \nab_3\psi +\dk^{\leq 1}(\Ga_g)\c \dk^2\psi, \quad  r\ge R/2,\\
 \\
& O(1)\dk^{\leq 2}\psi, \qquad\qquad \quad \quad\,\, r\le  R.
\end{cases}
\eea
Formula \eqref{eq:checkN-chapter10}, which is the precise analogue of  formula   (10.3.2)   of \cite{KS},  can be verified by a straightforward   calculation.    We refer the reader  to  section  10.3.2  and appendix  D.4  of \cite{KS}  for the details.

We apply the results of Theorem 
\ref{THEOREM:GENRW1-P}   to equation \eqref{eq:Gen.RW-chap10-psiwc} and derive\footnote{We refer the reader to  section 10.3.2  of \cite{KS}  for the  same   type of calculation in the proof of the analogous result, i.e. Theorem  5.18  in \cite{KS}.}  more details
\beaa
\BEF_q[\psiwc] &\les&  E_q[\psiwc] (\tau_1)+\NN_q[\psiwc, \widecheck{N}](\tau_1,\tau_2). 
\eeaa
Note that  the main term\footnote{The sign   can be easily tracked down   from the divergence formula, see for example  \eqref{eq:Diveregencefor-Chatpter10}.} 
 in  $\NN_q[\psiwc, \widecheck{N}]$ is given by
\beaa
-\int_{\MM_{\ge R}(\tau_1,\tau_2)}   \frac{2}{r}\left(1-\frac{3m}{r} +O(r^{-2} )\right)   r^q|\nabcheck_4 \psiwc|^2
\eeaa
Denoting $\EE_q=\EE[\psicheck][f_q,  2 r^{-1} f_q, 2r^{-1} f_q' e_4 ]$  with
$\EE[\psicheck]$ as in  \eqref{eq:EE-chapter10}  (for $\psi$ replaced by $\psicheck$),
 we  derive  the following analogue of Proposition 10.48 in \cite{KS}.

     \begin{proposition}\label{proposition:enhanced Dafermos-Rodn}
     The following estimate holds true,
     \bea
     \bsplit
     \int_{\MM_{\ge R}(\tau_1, \tau_2)} (\EE_q+ r^q \ec_4(\psicheck) \check{N})&\ge  \frac 1 8   \int_{\MM_{\ge R}(\tau_1, \tau_2)} r^{q-1}\left(   (2+q) | \ec_4 \psicheck |^2 + (2-q)|\nabb\psicheck|^2+ 2 r^{-2}|\psicheck|^2\right)\\   
 &- O(\ep)  \sup_{\tau_1\le \tau\le \tau_2} \Edot_{q,R}[\psicheck](\tau)    \\
     &- O(1)\left(E_{\max(q,\de)}^1[\psicheck](\tau_1)+\NN_{\max(q,\de)}^1[\psi,N]\right).
     \end{split}
     \eea
     \end{proposition}
     
The  remaining part of the proof, based on  choosing $R$ large   and making use of the result of  Theorem \ref{THEOREM:GENRW1-P},   is   exactly as in   section 10.3.2 of \cite{KS}.
\end{proof}


\section{Conditional weighted estimate for  scalar wave}
\lab{sec:conditionalweigthedMorawetzscalar}


\begin{proposition}
\lab{Prop:scalarwavePsi-M8-chap10}
Let $\psi$ be a solution to the following scalar wave equation 
\bea\lab{eq:saclarwwavemodelthatwillbeusedcontrolPc}
\square_\g\psi+V\psi &=& N,
\eea 
where $V$ is real and satisfies  $V=O(r^{-3})$ for $r$ large, in a spacetime  $\MM=\MM(1, \tau_*)$   verifying  the assumptions \eqref{eq:assumptionsonMMforpartII}. Then:
  \begin{enumerate}
\item The following conditional Morawetz estimates hold true in $\MM=\MM(1,\tau_*)$
\bea
\lab{eq:Estimatesforpsi-M8-1-chap10}
\bsplit
   B^k_{\de}[\psi](1,\tau_*)  &\les  \EF_\de^{k}[\psi](1,\tau_*)  +B ^{k-1}_{\de}[\psi](1,\tau_*)  +\int_{\MM_{trap}(1,\tau_*)} |\dk^{\le k}\psi|^2 \\
   &+ \NNmor^k[\psi,  N](1,\tau_*)+ \NNred^k [\psi,  N](1,\tau_*)\\
   & +   \int_{\Mext(1,\tau_*)} r^\de \Big(\big|\nab_4\dk^{\le k} \psi \big|+r^{-1}\big|\dk^{\le k} \psi \big|\Big) \big| \dk^{\le k} N \big|.
\end{split}
\eea

 \item The following conditional  Energy-Morawetz estimates hold true
 \bea
\lab{eq:Estimatesforpsi-M8-2chap10}
\bsplit
\BEF_{\de}[\psi] (1,\tau_*)
& \les   E^k _{\de}[\psi](0)       + \BEF^{k-1}_\de[\psi](1,\tau_*) \\
& +\int_{\MM_{trap}(1,\tau_*)} |\dk^{k}\psi|^2  +\NN^k_\de[\psi,  N](1,\tau_*). 
    \end{split}
\eea
 \end{enumerate}
\end{proposition}

\begin{remark}
 Note  that    both  estimates    are conditional   on the control of  the terms $  \BEF^{k-1}_\de[\psi](1,\tau_*)$ and  $\int_{\MM_{trap}} |\dk^{k}\psi|^2 $.  This result will be used  to control  $\Pc$  in Chapter \ref{Chapter:EN-MorforPc}. 
\end{remark}

The proof of the Proposition \ref{Prop:scalarwavePsi-M8-chap10} relies on  the conditional energy-Morawetz estimates of Proposition \ref{proposition:GeneralWaveComplexPotential:Morawetz} and an analog of Proposition \ref{Proposition:Step3-Chap10} for solutions to \eqref{eq:saclarwwavemodelthatwillbeusedcontrolPc}. Given that the estimate is only conditional, and in view of the the strong decay in $r$ for  the potential   $V$, the proof for solutions to \eqref{eq:saclarwwavemodelthatwillbeusedcontrolPc} is similar and in fact simpler that the one of Proposition \ref{Proposition:Step3-Chap10}.

 
 \chapter{Estimates for the full  Regge Wheeler equation}
 \lab{chapter-full-RW}\lab{chapter-full-RWforqf}
 

 The goal of this chapter is to  provide a complete  proof for   Theorem M1. 
 To this end we  proceed as follows:
 \begin{enumerate}
 \item  We use the  gRW equation for $\psi= \Re(\qf)$  coupled with the transport equation provided  by the definition of $\qf$ in terms  of $A$  to derive   combined   $r$-weighted  estimates  for $(\psi, A)$, see Theorem \ref{theorem:unconditional-result-final}.
 \item We  extend  the results of  Theorem  \ref{THEOREM:GENRW2-Q} to the full gRW equations
 to derive  improved   $r$-weighted  estimates based on the  quantity 
  $\psicheck= r^2 (\nab_4 \psi + \frac{r}{|q|^2}  \psi)$.  The result, stated in Theorem \ref{theorem:gRW1-p-strong-psiwc}, is the precise analogue of  of Theorem 5.15 in \cite{KS}.

 \item We use  the $r$ weighted  estimates  of Theorems  \ref{theorem:unconditional-result-final} and  \ref{theorem:gRW1-p-strong-psiwc}
  to prove  Theorem M1 by  relying on the  decay of flux arguments, based on mean value arguments,  following the procedure  detailed in   section 5.4 of \cite{KS}.  
 \end{enumerate}


\section{Preliminaries}
\lab{section:Preliminaries:Transport}
 
 
 The spacetime $\MM$ we are dealing with here is  precisely  that described in 
 section \ref{sec:pfdoisdvhdifuhgiwhdniwbvoubwuyf}. We do however make stronger assumptions  on $(\Ga_g, \Ga_b)$. We assume in fact  that
  for all $k\le \kl$, with\footnote{This is consistent with the  value of $k_L$  used in       the bootstrap assumption needed in the proof of Theorem M1 (see  section \ref{sec:bootstrapassumptions:intro}).}  $ k_L=k_{small}+120 $
\bea\lab{eq:assumptionsonMextforpartII-1}
\bsplit
\Big(r^2\tau^{\frac{1}{2}+\dec}+r\tau^{1+\dec}\Big)|\dk^{\leq k}\Ga_g| & \leq\ep, \\ 
r\tau^{1+\dec}|\dk^{\leq k}\Ga_b| & \leq\ep.
\end{split}
\eea
 We also assume that the curvature components   $A, B$ verify,  for  $ k\leq k_L$,
 \bea
 \lab{eq:assumptionsonMextforpartII-2}
  r^{7/2+\dec} |\dk^{\leq k} (A, B)| &\leq \ep.
 \eea
 It is  important in this chapter  that the frame is such that 
 \bea
 \lab{eq:GlobalFrame-HcinGa_g1}
 \Hc \in \Ga_g
 \eea
 and that $\Xi$ verifies
 \bea
 \lab{eq:GlobalFrame-HcinGa_g2}
 r^3|\dk^{\le k} \Xi|\le \ep, \qquad \nab_3\Xi\in r^{-1}\dk^{\leq 1}\Ga_g.
 \eea
 \begin{remark}
 The additional  conditions \eqref{eq:GlobalFrame-HcinGa_g1}- \eqref{eq:GlobalFrame-HcinGa_g2}   are verified  by the  global frame   constructed in section  3.6.3 of \cite{KS:Kerr}. These are crucial in
  deriving the  correct structure of the nonlinear terms  of the gRW  equation for $\qf$.  
\end{remark}
 \begin{remark}
 \lab{remark:interpolation-chap11}
 We note that in  reality the estimates \eqref{eq:assumptionsonMextforpartII-1}-\eqref{eq:assumptionsonMextforpartII-2} for  $ k_L-120 \le s \le k_{L}$ should be relaxed  by  replacing $\dec$ with $\frac 3 4 \dec$. This  loss is due to the interpolation between the   bootstrap estimates  for $k_{large}$ and  those for $k_{small}$, see for example Lemma 5.1  in \cite{KS}. The loss is  more than compensated  by the fact that   the  resulting  gain  $ \frac 3 4 \dec$   is doubled      in  nonlinear estimates. The remark also applies to Chapter \ref{chapter-full-RWforqfb}.
 \end{remark}


\subsection{Definition of $\qf$}

 
 Recall, see Definition \ref{Definition:Define-qf},  
 \bea
 \lab{eq:definition-qf-again}
 \bsplit
 \qf&=q \ov{q}^{3} \left( \nabc_3\nabc_3 A + C_1  \nabc_3A + C_2   A\right),\\
 C_1&=2\trchb - 2\frac {\atrchb^2}{ \trchb}  -4 i \atrchb, \\
C_2  &= \frac 1 2 \trchb^2- 4\atrchb^2+\frac 3 2 \frac{\atrchb^4}{\trchb^2} +  i \left(-2\trchb\atrchb +4\frac{\atrchb^3}{\trchb}\right).
 \end{split}
\eea

  
  \subsection{Full Regge Wheeler equation for $\qf$}
  \lab{section:FullgRWchap11}

The real part of $\qf$, denoted $\psi=\Re(\qf)$, verifies    according to Proposition \ref{prop:eq-real-gen-RW}   the  equation
 \bea\label{eq:Gen.RW-pert-again}
\squared_2 \psi -V_0\psi=- \frac{4 a\cos\th}{|q|^2}\dual \nab_T  \psi+N, \qquad  V_0= \frac{4\De}{ (r^2+a^2) |q|^2}, 
\eea
 with the right hand side $N$ being given by
 \bea\lab{eq:definition-N-psi-again}
 N&=& N_0+N_L+N_{\err}
 \eea
 where
 \begin{itemize}
 \item[-] $N_0$ denotes the zero-th order term in $\psi$, i.e.
 \bea\lab{eq:definition-N-0-psi-again}
 N_0:= \big(V- V_0\big)  \psi=O\left(\frac{a}{r^4}\right) \psi.
 \eea
 
 \item[-] $N_L$ denotes the linear term in $A$ given by 
 \bea\lab{eq:definition-N-L-psi-again}
 \bsplit
N_L =& \Re\Bigg(q\ov{q}^3    \Bigg[- \frac{8a^2 \De}{r^2|q|^4}\nab_\T \nab_3 A -\frac{8a  \De }{r^2|q|^4} \nab_\Z \nab_3A\\
&+ W_4\nab_4A+ W_3\nab_3A+W\c\nab A +W_0  A\Bigg] \Bigg)
\end{split}
 \eea
where $W_4$, $W_3$, $W_0$ are complex functions  of $(r, \th)$, and $W$ is the product of a complex function of $(r, \th)$ with $\dual\Re(\Jk)$, having the following fall-off in $r$ 
\beaa
q\ov{q}^3 W_4= q\ov{q}^3 W_3=q\ov{q}^3 W=O\left(a\right), \qquad  q\ov{q}^3 W_0=O\left(\frac{a}{r}\right).
\eeaa
Away from the trapping, the following schematic structure will suffice
 \bea\lab{eq:definition-N-L-psi-again:schematic}
N_L&=& O(a)\dk^{\leq 1}\nab_3 A+ O(ar^{-1} )\dk^{\leq1}A.
\eea
 
 \item[-] $N_{\err}$ denotes the  quadratic error terms, given  schematically by the expression 
 \bea
 \lab{eq:N_err=N_g+}
 \bsplit
 N_{\err}
 &=N_g+ \nab_3(rN_g)+ N_m[\qf],\\
  N_g&= r^2 \dk^{\le 2}\big(\Ga_g\c(A,B) \big), \quad N_m[\qf]=\dk^{\le 1}\big(\Ga_g \c \qf\big).
  \end{split}
 \eea
 \end{itemize}
 We refer  to solutions to  \eqref{eq:Gen.RW-pert-again}, with $N$ as above,   as  solutions to  the  real   gRW.


\subsection{Factorization of $\qf$}


\begin{lemma}\lab{lemma:factorizationofqfusefulfortransportequations}
We have
\bea
\nn&&\nabc_3\left(\nabc_3\left(\frac{(\ov{\tr\Xb})^2}{(\Re(\tr\Xb))^2(\tr\Xb)^2}A\right) -\frac{r^2}{2}F_2A\right) \\
&=& \Big(O(r^{-2})+r^{-1}\Ga_g\Big)\qf  +r^2\dk^{\leq 1}\Ga_b\nabc_3A +r\dk^{\leq 1}\Ga_b A,
\eea
where $F_2$ is given by 
\bea\lab{eq:defintionofF2appearingtransportequationA:chap11}
\bsplit
F_1 &:= \DDc\c\ov{\Xib}+\Xib\c\ov{\Hb}+\ov{\Xib}\c H-\frac{1}{2}\Xbh\c\ov{\Xbh},\\
F_2 &:= -\frac{F_1}{\tr\Xb} +\frac{\ov{F_1}}{\ov{\tr\Xb}}  - \frac{\Re(F1)}{\Re(\tr\Xb)}.
\end{split}
\eea
\end{lemma}

\begin{proof}
We compute
\beaa
&&\nabc_3\nabc_3\left(\frac{(\ov{\tr\Xb})^2}{(\Re(\tr\Xb))^2(\tr\Xb)^2}A\right) \\
&=& \nabc_3\left(\frac{(\ov{\tr\Xb})^2}{(\Re(\tr\Xb))^2(\tr\Xb)^2}\nabc_3A+\nabc_3\left(\frac{(\ov{\tr\Xb})^2}{(\Re(\tr\Xb))^2(\tr\Xb)^2}\right)A\right)\\
&=& \frac{(\ov{\tr\Xb})^2}{(\Re(\tr\Xb))^2(\tr\Xb)^2}\Big(\nabc_3\nabc_3A +\widetilde{C}_1\nabc_3A +\widetilde{C}_2A\Big)
\eeaa
where 
\beaa
\widetilde{C}_1 &:=& 2\frac{\nabc_3\left(\frac{(\ov{\tr\Xb})^2}{(\Re(\tr\Xb))^2(\tr\Xb)^2}\right)}{\frac{(\ov{\tr\Xb})^2}{(\Re(\tr\Xb))^2(\tr\Xb)^2}},\\
\widetilde{C}_2 &:=& \frac{\nabc_3\nabc_3\left(\frac{(\ov{\tr\Xb})^2}{(\Re(\tr\Xb))^2(\tr\Xb)^2}\right)}{\frac{(\ov{\tr\Xb})^2}{(\Re(\tr\Xb))^2(\tr\Xb)^2}}.
\eeaa
Together with the definition of $\qf$, we infer
\beaa
&&\nabc_3\nabc_3\left(\frac{(\ov{\tr\Xb})^2}{(\Re(\tr\Xb))^2(\tr\Xb)^2}A\right) \\
&=& \frac{(\ov{\tr\Xb})^2}{(\Re(\tr\Xb))^2(\tr\Xb)^2}\left(\frac{1}{q \ov{q}^{3}}\qf +(\widetilde{C}_1-C_1)\nabc_3A +(\widetilde{C}_2-C_2)A\right).
\eeaa

Next, recall the following null structure equation 
\beaa
\nabc_3\tr\Xb +\frac{1}{2}(\tr\Xb)^2 &=& \DDc\c\ov{\Xib}+\Xib\c\ov{\Hb}+\ov{\Xib}\c H-\frac{1}{2}\Xbh\c\ov{\Xbh}\\
&=:& F_1.
\eeaa
We have, in view of the definition of $\widetilde{C}_1$,  
\beaa
\frac{1}{4}\widetilde{C}_1 &=& -\frac{\nabc_3\tr\Xb}{\tr\Xb} + \frac{\nabc_3\ov{\tr\Xb}}{\ov{\tr\Xb}}  - \frac{\nabc_3\Re(\tr\Xb)}{\Re(\tr\Xb)}\\
&=& \frac{1}{2}\tr\Xb - \frac{1}{2}\ov{\tr\Xb} +\frac{1}{4\Re(\tr\Xb)}((\tr\Xb)^2+(\ov{\tr\Xb})^2) + F_2
\eeaa
where 
\beaa
F_2 &:=& -\frac{F_1}{\tr\Xb} +\frac{\ov{F_1}}{\ov{\tr\Xb}}  - \frac{\Re(F1)}{\Re(\tr\Xb)}.
\eeaa

Also, we have in view of the definition of $\widetilde{C}_2$
\beaa
\widetilde{C}_2 &=& \frac{1}{2}\frac{\nabc_3\left(\frac{(\ov{\tr\Xb})^2}{(\Re(\tr\Xb))^2(\tr\Xb)^2}\widetilde{C}_1\right)}{\frac{(\ov{\tr\Xb})^2}{(\Re(\tr\Xb))^2(\tr\Xb)^2}}\\
&=& \frac{1}{2}\nabc_3\widetilde{C}_1 +\frac{1}{4}(\widetilde{C}_1)^2\\
&=& \frac{1}{2}\nabc_3\left(2\tr\Xb-2\ov{\tr\Xb} +\frac{1}{\Re(\tr\Xb)}((\tr\Xb)^2+(\ov{\tr\Xb})^2) + 4F_2\right)\\
&& +\frac{1}{4}\left(2\tr\Xb-2\ov{\tr\Xb} +\frac{1}{\Re(\tr\Xb)}((\tr\Xb)^2+(\ov{\tr\Xb})^2) + 4F_2\right)^2.
\eeaa
Note that $F_1\in r^{-1}\dk^{\leq 1}\Ga_b$ and $F_2\in \dk^{\leq 1}\Ga_b$. In view of the above, and using the definition of $C_1$ and $C_2$ in \eqref{eq:definition-qf-again}, we infer
\beaa
\widetilde{C}_1 &=& 2\tr\Xb-2\ov{\tr\Xb} +\frac{1}{\Re(\tr\Xb)}((\tr\Xb)^2+(\ov{\tr\Xb})^2) + \dk^{\leq 1}\Ga_b\\
&=& -4i\atrchb  +\frac{2}{\trchb}((\trchb)^2 -(\atrchb)^2) + \dk^{\leq 1}\Ga_b\\
&=& -4i\atrchb  +2\trchb -\frac{2(\atrchb)^2}{\trchb} + \dk^{\leq 1}\Ga_b\\
&=& C_1 + \dk^{\leq 1}\Ga_b
\eeaa
and 
\beaa
\widetilde{C}_2 &=& 2\nabc_3(F_2)+\frac{1}{2}\Bigg(-(\tr\Xb)^2+(\ov{\tr\Xb})^2 -\frac{1}{\Re(\tr\Xb)}((\tr\Xb)^3+(\ov{\tr\Xb})^3)\\
&& +\frac{\Re((\tr\Xb)^2)}{2(\Re(\tr\Xb))^2}((\tr\Xb)^2+(\ov{\tr\Xb})^2) \Bigg)\\
&& +\frac{1}{4}\left(2\tr\Xb-2\ov{\tr\Xb} +\frac{1}{\Re(\tr\Xb)}((\tr\Xb)^2+(\ov{\tr\Xb})^2)\right)^2+r^{-1}\dk^{\leq 1}\Ga_b\\
&=& 2\nabc_3(F_2)+\frac{1}{2}\Bigg(4i\trchb\atrchb -\frac{1}{\trchb}(2\trchb^3 -6\trchb\atrchb^2)\\
&& +\frac{\trchb^2-\atrchb^2}{\trchb^2}(\trchb^2 -\atrchb^2) \Bigg) +\frac{1}{4}\left(-4i\atrchb +\frac{2}{\trchb}(\trchb^2 -\atrchb^2)\right)^2\\
&& +r^{-1}\dk^{\leq 1}\Ga_b\\
&=& C_2+2\nabc_3(F_2)+r^{-1}\dk^{\leq 1}\Ga_b.
\eeaa

We deduce
\beaa
&&\nabc_3\nabc_3\left(\frac{(\ov{\tr\Xb})^2}{(\Re(\tr\Xb))^2(\tr\Xb)^2}A\right) \\
&=& \frac{(\ov{\tr\Xb})^2}{(\Re(\tr\Xb))^2(\tr\Xb)^2}\left(\frac{1}{q \ov{q}^{3}}\qf +\dk^{\leq 1}\Ga_b\nabc_3A +\left(2\nabc_3(F_2)+r^{-1}\dk^{\leq 1}\Ga_b\right)A\right)\\
&=& \frac{(\ov{\tr\Xb})^2}{(\Re(\tr\Xb))^2(\tr\Xb)^2}\frac{1}{q \ov{q}^{3}}\qf +\frac{2(\ov{\tr\Xb})^2}{(\Re(\tr\Xb))^2(\tr\Xb)^2}\nabc_3(F_2)A +r^2\dk^{\leq 1}\Ga_b\nabc_3A +r\dk^{\leq 1}\Ga_b A\\
&=& \Big(O(r^{-2})+r^{-1}\Ga_g\Big)\qf +\frac{r^2}{2}\nabc_3(F_2)A +r^2\dk^{\leq 1}\Ga_b\nabc_3A +r\dk^{\leq 1}\Ga_b A\\
&=& \Big(O(r^{-2})+r^{-1}\Ga_g\Big)\qf +\nabc_3\left(\frac{r^2}{2}F_2A\right) +r^2\dk^{\leq 1}\Ga_b\nabc_3A +r\dk^{\leq 1}\Ga_b A
\eeaa
and hence
\beaa
&&\nabc_3\left(\nabc_3\left(\frac{(\ov{\tr\Xb})^2}{(\Re(\tr\Xb))^2(\tr\Xb)^2}A\right) -\frac{r^2}{2}F_2A\right) \\
&=& \Big(O(r^{-2})+r^{-1}\Ga_g\Big)\qf  +r^2\dk^{\leq 1}\Ga_b\nabc_3A +r\dk^{\leq 1}\Ga_b A
\eeaa
as stated. This concludes the proof of Lemma \ref{lemma:factorizationofqfusefulfortransportequations}.
\end{proof}

Next, we introduce the tensor $\Psi$. 
\begin{definition}\lab{def:PsiforintegrationofAfromqf:chap11}
Let $\Psi\in \sk_2(\CCC)$ given by 
\beaa
\Psi &:=& \nabc_3\left(\frac{(\ov{\tr\Xb})^2}{(\Re(\tr\Xb))^2(\tr\Xb)^2}A\right) -\frac{r^2}{2}F_2A,
\eeaa
where $F_2$ is given by \eqref{eq:defintionofF2appearingtransportequationA:chap11}.
\end{definition}

We have the following corollary of Lemma \ref{lemma:factorizationofqfusefulfortransportequations}.
\begin{corollary}\lab{cor:systemoftransportequationsforPsiandAfromqf:chap11}
Let $\Psi$ as in Definition \ref{def:PsiforintegrationofAfromqf:chap11}. Then, $(\Psi, A)$ satisfies the following system of transport equations 
\beaa
\bsplit
\nabc_3\Psi &=\Big(O(r^{-2})+r^{-1}\Ga_g\Big)\qf  +r^2\dk^{\leq 1}\Ga_b\nabc_3A +r\dk^{\leq 1}\Ga_b A, \\ \nabc_3\left(\frac{(\ov{\tr\Xb})^2}{(\Re(\tr\Xb))^2(\tr\Xb)^2}A\right) & = \Psi+r^2\dk^{\leq 1}(\Ga_b)\c A.
\end{split}
\eeaa
\end{corollary}

\begin{proof}
This is an immediate consequence of Lemma \ref{lemma:factorizationofqfusefulfortransportequations}, Definition  \ref{def:PsiforintegrationofAfromqf:chap11}, and the fact that $F_2\in \dk^{\leq 1}\Ga_b$.
\end{proof}


\subsection{Norms for $\psi$}


     We recall that the  norms $B^s_p[\psi]$,  $E^s_p[\psi]$ and $F^s_p[\psi]$, respectively for the bulk, energy and flux of $\psi$ were defined in section \ref{subsection:basicnormsforpsi}.     To simplify notations, we make use of   the combined norms
      \bea
      \BEF^s_p(\tau_1,\tau_2)&:=&B^s_p[\psi](\tau_1, \tau_2)+ 
      \sup_{\tau\in[\tau_1,\tau_2]} E^s_p[\psi](\tau)+ F^s_p[\psi](\tau_1, \tau_2).
      \eea
  We also  recall    
   \beaa
       \NN_p[\psi,  N](\tau_1, \tau_2) &=&\int_{\MM(\tau_1,\tau_2)} \big(|\nab_{\Rhat} \psi|+r^{-1}|\psi|\big) |N|+\left| \int_{\MM_{r\geq 4m}(\tau_1,\tau_2)}  r^{p-1}  \, \nab_4 (r\psi ) \c  N\right|\\
       &&+\left|\int_{\MM(\tau_1,\tau_2)} \nab_{\That_\de} \psi \c N\right|,
       \eeaa 
       and the corresponding  higher  derivatives  $ \NN^s_p[\psi,  N]$ norms.

 
 \section{Control of the full gRW equation for $\qf$}
 \label{section:main-results-chapter-7}
 


\subsection{Combined norms for $\psi$, $A$}
  

  \begin{definition}
  \lab{Definition:NormsBEF-A}
   We define the following norms  modified bulk and energy-flux norms  for A and 
    combined norms for $(A, \psi)$ 
 \beaa
B_p[A](\tau_1, \tau_2)&=& \int_{\MM(\tau_1, \tau_2)}   r^{p+1} \Big(r^4|\nab_3\nabc_3A|^2+r^4|\nab_4\nab_3A|^2+r^4|\nab\nab_3A|^2+r^2|\nab_3A|^2\\
&&+r^2|\nab_4A|^2+r^2|\nab A|^2+|A|^2\Big),\\
E_p[A](\tau)&=& \int_{\Si(\tau)}   r^{p+2} \Big(r^4\chi_{nt}^2|\nab_4\nab_3A|^2+r^2|\nab_{\Rhat}\nab_3A|^2+r^2|\nab_3A|^2+r^2|\nab_4A|^2\\
&&+|A|^2 \Big),\\
F_p[A](\tau_1, \tau_2)&=&\int_{\AA\cup \Si_*(\tau_1, \tau_2)  }  r^{p+2} \Big(r^4\chi_{nt}^2|\nab_4\nab_3A|^2+r^2|\nab_{\Rhat}\nab_3A|^2+r^2|\nab_3A|^2+r^2|\nab_4A|^2\\
&&+|A|^2 \Big),
\eeaa
where $\chi_{nt}=\chi_{nt}(r)$ denotes a smooth cut-off function equal to 0 on $\MM_{trap}$ and equal to 1 on $r\geq 4m$.
   
    The higher derivative  norms are defined by the usual procedure
    \beaa
    B^s_p[A]=  B_p[\dk^{\le s} A],  \qquad   E^s_p[A]=  E^s_p[\dk^{\le s} A],  \qquad   F^s_p[A]=  F^s_p[\dk^{\le s} A]. 
    \eeaa
        We  also  define the  combined norms:
  \beaa
  E^s_p[\psi, A](\tau)&=&  E^s_p[\psi](\tau) +E^s_p[A](\tau),\\
  B^s_p[\psi, A](\tau_1, \tau_2) &=&  B^s_p[\psi](\tau_1, \tau_2) +B^s_p[A](\tau_1, \tau_2),\\
  F^s_p[\psi, A](\tau_1, \tau_2) &=& F^s_p[\psi](\tau_1, \tau_2)+ F^s_p[A](\tau_1, \tau_2).
  \eeaa
  We  use the short hand notation
  \beaa
  \BEF_p^s[A](\tau_1, \tau_2)  &=& B_p^s[A](\tau_1, \tau_2)  +\sup_{\tau\in[\tau_1,\tau_2]} E_p^s[A](\tau)+
   F_p^s[A](\tau_1, \tau_2), \\
  \BEF_p^s[\psi, A] &=&\BEF_p^s[\psi] +\BEF_p^s[A]. 
  \eeaa
     \end{definition}
     
    \begin{remark}
    \lab{Remark:BEF[A]-norms}
    Note  that the derivatives  $\nab_4^2 A$, $\nab\nab_4 A$ and $\nab^2A$   are missing  in the combined norm $\BEF_p[A](\tau_1, \tau_2)$, as they  cannot be  derived by the  transport  equation methods used here.  Fortunately,
     in view of the structure of  the $N_L$ term, they are not needed to close the estimates for $\psi$ and thus for $\qf$.  Additional derivatives in $E_p[A](\tau)$ and $F_p[A](\tau_1, \tau_2)$ are missing as well and are also  not needed to close the estimates for $\qf$, with the exception of the ones recovered in \eqref{eq:transportA:additionalestimateformissingderivativesenergy}. We  also  remark  that
    \beaa
     \int_{\MM(\tau_1, \tau_2)}   r^{p+1} \Big(r^2 |\dk^{\le s+1}  \nab_3 A|^2+  |\dk^{\le s+1} A|^2 \Big)&\les& B^s_p[A](\tau_1, \tau_2),  
    \eeaa
     but the norm $B_p^s[A]$ is in fact stronger in powers of $r$ for  the $\nab_3$  derivative. 
    \end{remark}

  
 \subsection{Weighted  estimates for  the full gRW system  \ref{eq:Gen.RW-pert-again} }


The main technical  results of Chapter \ref{chapter-full-RWforqf}, which extend Theorems  \ref{THEOREM:GENRW1-P}  and  \ref{THEOREM:GENRW2-Q}
 to the full  real gRW  system \ref{eq:Gen.RW-pert-again},   are   as follows.
  
 \begin{theorem}
\lab{Thm:Nondegenerate-Morawetz}\lab{theorem:unconditional-result-final}
The following holds true, for   all   $\de\le p\le  2-\de$ and  $2\leq s\le \kl$,
 \bea
       \lab{eqtheorem:gRW1-p-chap11}
        \BEF_p^s[\psi, A](\tau_1, \tau_2) \les  E_p^s[\psi, A](\tau_1) +\NN_p^s[\psi, N_\err](\tau_1, \tau_2) +\ep_0^2\tau_1^{-2-3\dec},
       \eea
       where $N_\err$ is given by \eqref{eq:N_err=N_g+}.
\end{theorem}

\begin{theorem}
\lab{THEOREM:GRW1-P-WEAK-PSIWC} Under the same assumptions as before
we have for $2\le s\le \kl-1$, for all $-1+\de\leq q \leq 1-\de$,
 \bea
       \lab{eq:theorem-gRW1-p-weak-psiwc-strong'}
       \begin{split}
       \BEF_q^s[\psiwc](\tau_1, \tau_2)  &\les
       E_q^s[\psiwc](\tau_1)+E_{\max(q, \de)}^{s+1}[\psi, A](\tau_1)+\widecheck{\NN}_{q}^{s}[\psiwc, N_\err](\tau_1, \tau_2)\\
       &+\NN_{\max(q, \de)}^{s+1}[\psi, N_\err](\tau_1, \tau_2)  +\ep_0^2\tau_1^{-2-3\dec},
       \end{split}
       \eea
       where $\psiwc= r^2(\nab_4 \psi+\frac{r}{|q|^2} \psi)$ and 
       \beaa
       \widecheck{\NN}_{q}^{s}[\psiwc, N_\err](\tau_1, \tau_2)&=& \int_{\Mext} r^{q+2} \widecheck{\nab}_4 \dk^{\leq s}\psiwc \c \left( \nab_4 \dk^{\leq s}N_\err + \frac 3 r \dk^{\leq s} N_\err \right).
       \eeaa
\end{theorem}

We postpone the proof of Theorem \ref{THEOREM:GRW1-P-WEAK-PSIWC} to section \ref{section:ProofThm-gRW1-p-weak-psiwc} and concentrate  our attention 
 to the proof of  Theorem  \ref{theorem:unconditional-result-final}.


\subsection{Proof of Theorem \ref{Thm:Nondegenerate-Morawetz}}
\lab{section:SketchProofThmNondegenerate-Morawetz}


       We restate for convenience   the result of Theorem \ref{THEOREM:GENRW1-P}, whose proof was completed in  chapter \ref{chapter-rp-estimates}. 
       \begin{theorem}[Basic $r^p$-weighted estimates]
       The following estimates hold true  for solutions $\psi\in\sk_2$  of   the model gRW  equation, see \eqref{eq:Gen.RW},  on $\MM$, for all $\de\le p\le  2-\de$ and  $2\leq s\le \kl$,
       \bea
       \lab{eqtheorem:GenRW1-p'}
       \BEF^s_p[\psi](\tau_1, \tau_2) \les
       E_p^s[\psi](\tau_1)+\NN_p^s[\psi, N](\tau_1, \tau_2).
       \eea
       \end{theorem}

 The proof of Theorem \ref{Thm:Nondegenerate-Morawetz}  is done in steps as follows.
 
 {\bf Step 1.}  Recall that $N=N_0+N_L+N_\err$, see \eqref{eq:definition-N-psi-again}.  
 We  first eliminate $N_0+N_L$ from  the right hand side  of \eqref{eqtheorem:GenRW1-p'}. 
  \begin{proposition}
  \lab{prop:Step1-psi}
  The following  estimate for solutions  $\psi$ of the  full gRW  equation hold true, for all $s\le k_L$ and all $\de\le p\le 2-\de$,
  \bea
  \lab{eq:Step1-psi}
  \BEF_p^s[\psi](\tau_1, \tau_2) \les
       E_p^s[\psi](\tau_1) +   |a| \BEF^s_\de[\psi, A](\tau_1,\tau_2) +\NN_p^s[\psi, N_\err](\tau_1, \tau_2).
  \eea
  \end{proposition}
  
The proof of Proposition \ref{prop:Step1-psi} is an immediate consequence of \eqref{eqtheorem:GenRW1-p'} and the following lemma.
\begin{lemma}
\lab{LEMMA:DECAYFORTHENTERM-PSI}
 For $\de\leq p\le 2-\de$, $N$ given by \eqref{eq:definition-N-psi-again} satisfies 
 \bea
 \NN^s _p[\psi,  N](\tau_1, \tau_2) 
  \les |a| \BEF_\de^s[\psi, A](\tau_1, \tau_2)+  \NN_p^s[\psi, N_\err](\tau_1, \tau_2) +\ep_0^2\tau_1^{-2-3\dec}.
  \eea   
\end{lemma}

The proof of Lemma \ref{LEMMA:DECAYFORTHENTERM-PSI} is given in section \ref{section:Them:gRW1-p-weak}.

\bigskip

  {\bf Step 2.} We eliminate the terms  $ \BEF^s_p[A]$  from the right hand side 
  of \eqref{eq:Step1-psi} with the help of the proposition  below. 
   \begin{proposition}
 \lab{prop:MaiTransportA-steps} 
 The following  estimates hold true, for all $ s\le k_L$ and for all $\de\leq p \leq 2-\de$, 
 \bea
 \lab{eq:transportA}
\BEF_p^s[A] (\tau_1, \tau_2)&\les&   B^s_p[\psi](\tau_1,\tau_2) 
 + E^s_p[A](\tau_1) +\ep_0^2\tau_1^{-2-3\dec}.
 \eea
Also, we have the following additional control  on $\Si(\tau)$ with $\tau\in[\tau_1, \tau_2]$, for $ s\le k_L$, for all $\de\leq p \leq 2-\de$, 
\bea\lab{eq:transportA:additionalestimateformissingderivativesenergy}
\bsplit
&\int_{\Si(\tau)}r^{p+2}\Big(r^{\min(4, 5-\de-p)}|\nab_3\nabc_3\dk^{\leq s}A|^2+r^4|\nab_4\nab_3\dk^{\leq s}A|^2\\
&+r^4|\nab\nab_3\dk^{\leq s}\Ab|^2+r^2|\nab\dk^{\leq s}\Ab|^2\Big)\\
 \les& EB^s_p[\psi](\tau_1,\tau_2)   + E^s_p[A](\tau_1) +\ep_0^2 \tau_1^{-2-3\dec}.
\end{split}
\eea
 \end{proposition}

The proof of Proposition \ref{prop:MaiTransportA-steps} is given in section  \ref{section:transport-estimates}.

{\bf Step 3.} As a consequence    of Proposition \ref{prop:MaiTransportA-steps}  and  Proposition \ref{prop:Step1-psi}, as well as the smallness of $|a|/m$,  we deduce, for all $s\le k_L$ and for all $\de\le p\le 2-\de$,
\beaa
\BEF_p^s[\psi, A] (\tau_1, \tau_2)&\les&  E^s_p[\psi, A](\tau_1) +\NN_p^s[\psi, N_\err](\tau_1, \tau_2) +\ep_0^2\tau_1^{-2-3\dec}
\eeaa
as stated. This  ends the proof of Theorem  \ref{Thm:Nondegenerate-Morawetz}.

  
\section{Proof of Lemma \ref{LEMMA:DECAYFORTHENTERM-PSI}}
\lab{section:Them:gRW1-p-weak}


Recall that $N=N_0+N_L+N_\err$, see \eqref{eq:definition-N-psi-again}.  Therefore, we have
\beaa
\NN_p^s[\psi, N](\tau_1, \tau_2)&\les& \NN_p^s[\psi, N_0](\tau_1, \tau_2)+
  \NN_p^s[\psi, N_L](\tau_1, \tau_2)+   \NN_p^s[\psi, N_\err](\tau_1, \tau_2).
\eeaa
In order to prove Lemma \ref{LEMMA:DECAYFORTHENTERM-PSI}, we need to eliminate  the terms  $\NN_p^s[\psi, N_0]$ and $\NN_p^s[\psi, N_L]$.

 Recall that by definition of $\NN^s_p[\psi,  N](\tau_1, \tau_2)$ we  have:
\beaa
\NN_p^s[\psi,  N](\tau_1, \tau_2) &=&\, ^{(Mor)} \NN^s[\psi,  N](\tau_1, \tau_2) +\Next^s_p[\psi,  N](\tau_1, \tau_2)+\, ^{(En)} \NN^s[\psi,  N](\tau_1, \tau_2),\\
\, ^{(Mor)} \NN^s[\psi,  N](\tau_1, \tau_2) &=&\sum_{k\le s} \int_{\MM(\tau_1, \tau_2)}\Big(|\nab_{\Rhat}(\dk^k \psi)|+r^{-1}|\dk^{k}\psi|\Big) |\dk^k  N|,\\
\Next^s_p[\psi,  N](\tau_1, \tau_2)&=&\sum_{k\le s}\left| \int_{\Mext(\tau_1,\tau_2)}  r^{p-1}  \, \nab_4 (r \dk^k \psi ) \c  \dk^k N\right|, \\
\, ^{(En)} \NN^s[\psi,  N](\tau_1, \tau_2)&=&\sum_{k\le s}  \left|\int_{\MM(\tau_1, \tau_2)} \nab_{\That_\de} ( \dk^k\psi )\c  \dk^kN\right|.
\eeaa
We then  estimate separately the above  terms as follows. 

\bigskip

{\bf Step 1.}  First, we estimate  $^{(Mor)} \NN^s[\psi,  N_0+N_L] +\Next^s_p[\psi,  N_0+N_L]$. We obtain, for $\de\leq p\leq 2-\de$, 
\bea\label{eq:step1-}
 ^{(Mor)} \NN^s[\psi,  N_0+N_L](\tau_1, \tau_2) +\Next^s_p[\psi,  N_0+N_L](\tau_1, \tau_2) \les   |a|B_\de^s[\psi, A](\tau_1, \tau_2).
\eea
This is done in section \ref{subsection:proof-step1}.

\bigskip

{\bf Step  2.}  Then, we estimate  $\, ^{(En)} \NN^s[\psi,  N_0]$ and obtain
\bea\label{eq:step2-}
 ^{(En)} \NN^s[\psi,  N_0](\tau_1, \tau_2) &\les& |a|\left(\sup_{\tau\in [\tau_1, \tau_2]} E^s[\psi](\tau)+ B^s_\de[\psi](\tau_1, \tau_2)\right).
\eea
This is done in section \ref{subsection:proof-step2}.

\bigskip

{\bf Step 3.}  Next, we estimate  $\, ^{(En)} \NN^s[\psi,  N_L]$ and obtain
\bea\label{eq:step3-}
\nn\, ^{(En)} \NN^s[\psi,  N_L](\tau_1, \tau_2)&\les& |a|\left(\sup_{\tau \in (\tau_1, \tau_2)} E_\de^s[\psi, A](\tau)+B_\de^s [\psi, A](\tau_1, \tau_2) \right)\\
&&+\ep_0^2\tau_1^{-2-3\dec}.
\eea
This is done in section \ref{subsection:proof-step3}.

\bigskip

{\bf Step 4.}  Combining the results of Steps 1-3, we obtain the following, for $\de\leq p\le 2-\de$, 
 \beaa
 \NN^s _p[\psi,  N](\tau_1, \tau_2)  &\les&  ^{(Mor)} \NN^s[\psi,  N_0+N_L](\tau_1, \tau_2) +\Next^s_p[\psi,  N_0+N_L](\tau_1, \tau_2)\\
 &&+ ^{(En)} \NN^s[\psi,  N_0](\tau_1, \tau_2)+ ^{(En)} \NN^s[\psi,  N_L](\tau_1, \tau_2)+\NN^s _p[\psi,  N_\err](\tau_1, \tau_2)\\
  &\les& |a|  \BEF_\de^s[\psi, A](\tau_1, \tau_2)+  \NN_p^s[\psi, N_\err](\tau_1, \tau_2) +\ep_0^2\tau_1^{-2-3\dec}
  \eeaa  
as stated. This concludes the proof of Lemma \ref{LEMMA:DECAYFORTHENTERM-PSI}.

It thus remains  to prove  estimates \eqref{eq:step1-}, \eqref{eq:step2-} and \eqref{eq:step3-}. This is done in sections \ref{subsection:proof-step1}, \ref{subsection:proof-step2} and \ref{subsection:proof-step3} respectively.

  
\subsection{Proof of the estimate \eqref{eq:step1-}}\label{subsection:proof-step1}


In this section, we estimate $^{(Mor)} \NN^s_p[\psi,  N_0+N_L](\tau_1, \tau_2) +\Next^s_p[\psi,  N_0+N_L](\tau_1, \tau_2)$.

We have
\beaa
&& ^{(Mor)}\NN^s[\psi,  N_0+N_L](\tau_1, \tau_2) \\
&=& \int_{\MM(\tau_1, \tau_2)}  \Big(|\nab_\Rhat\dk^{\le s} \psi|+ r^{-1} |\dk^{\le s}\psi|\Big) |\dk^{\le s }N| \\
&\les& \left(  \int_{\MM(\tau_1, \tau_2)}  r^{-1-\de}  \Big(|\nab_\Rhat\dk^{\le s} \psi|^2+ r^{-2} |\psi|^2\Big) \right)^{1/2} \left( \int_{\MM(\tau_1, \tau_2)} r^{1+\de}  |\dk^{\le s }N| ^2 \right)^{1/2}\\
&\les&\left(  \int_{\MM(\tau_1, \tau_2)}  r^{-1-\de}  \Big(|\nab_\Rhat\dk^{\le s} \psi|^2\mathbb{1}_{r\leq 4m}+\big(|\nab_3\psi|^2+r^{-2}|\dk\psi|^2\big)\mathbb{1}_{r\geq 4m}+ r^{-2} |\psi|^2\Big) \right)^{1/2}\\
&&\times \left( \int_{\MM(\tau_1, \tau_2)} r^{1+\de}  |\dk^{\le s }N| ^2 \right)^{1/2}\\
&\les& \Big(B_\de^s[\psi](\tau_1, \tau_2)\Big)^{1/2}\left( \int_{\MM(\tau_1, \tau_2)} r^{1+\de}|\dk^{\le s }N| ^2\right)^{1/2}.
\eeaa
Since   $N_0= O(ar^{-4} ) \psi$, we have
\beaa
 \int_{\MM(\tau_1, \tau_2)} r^{1+\de}|\dk^{\le s }N_0| ^2 &\les&  |a|\Mor^s[\psi](\tau_1, \tau_2).
\eeaa
To estimate the term in  $N_L$ we  make use of the schematic structure  given in 
 \eqref{eq:definition-N-L-psi-again:schematic}
 \beaa
N_L&=&  O(a)\dk^{\leq 1}\nab_3 A+ O(ar^{-1} )\dk^{\leq1}A.
\eeaa
Thus,
\beaa
\int_{\MM(\tau_1, \tau_2)} r^{1+\de}|\dk^{\leq s} N_L|^2 &\les & |a| \int_{\MM(\tau_1, \tau_2)} r^{1+\de}\Big( | \nab_3 \dk^{\leq s+1} A|^2+ r^{-2} | \dk^{\leq s+1} A|^2 \Big).
\eeaa
Therefore,
\beaa
^{(Mor)}\NN^s[\psi,  N_0+N_L] &\les& |a|B_\de^s[\psi](\tau_1, \tau_2) + |a|\Big(B_\de^s[\psi](\tau_1, \tau_2)\Big)^{1/2} \\
&&\times \left(\int_{\MM(\tau_1, \tau_2)}  r^{1+\de}\Big( | \nab_3 \dk^{\leq s+1} A|^2+ r^{-2} | \dk^{\leq s+1} A|^2 \Big)  \right)^{1/2}\\
&\les& |a|\left( B_\de^s[\psi](\tau_1, \tau_2) + \int_{\MM(\tau_1, \tau_2)}  r^{1+\de}\Big( | \nab_3 \dk^{\leq s+1} A|^2+ r^{-2} | \dk^{\leq s+1} A|^2 \Big) \right) \\
&\les&  |a|\Big(  B_\de^s[\psi](\tau_1, \tau_2) + B^s_\de[A](\tau_1, \tau_2)\Big)\\
&\les&  |a| B_\de^s[\psi, A](\tau_1, \tau_2).
\eeaa

Similarly, we obtain
 \beaa
 \Next^s_p[\psi,  N_0+N_L]&=&\left| \int_{\Mext}  r^{p-1}  (r\nab_4) \dk^{\leq s} \psi  \c  \dk^{\leq s} N\right|\\
 &\les& \left(  \int_{\Mext}   r^{p-3} |\dk^{\le s+1} \psi|^2\right)^{1/2} \left( \int_{\Mext} r^{p+1}   |\dk^{\le s }N| ^2 \right)^{1/2}\\
&\les&\Big( \Bext_p^s[\psi] \Big)^{1/2}  \left( \int_{\Mext} r^{p+1}   |\dk^{\le s }N| ^2 \right)^{1/2}.
 \eeaa
Now\footnote{Recall the definition of the $B_p[A]$ norms  and Remark \ref{Remark:BEF[A]-norms}.}, for $\de\leq p\leq 2-\de$,  
\beaa
 \int_{\Mext} r^{p+1}   |\dk^{\le s }N_0| ^2 &\les & |a|\int_{\Mext} r^{p-7 }   |\dk^{\le s }\psi| ^2\les  |a| B_\de^s[\psi],\\
  \int_{\Mext} r^{p+1}   |\dk^{\le s }N_L| ^2 &\les&  |a|\int_{\Mext} r^{p+1} \big(  |\dk^{\le s+1 }\nab_3 A|^2+r^{-2} | \dk^{\leq s+1} A|^2\big) \les |a| B_\de^s[A],
\eeaa
 from which we   deduce, for $\de\leq p\leq 2-\de$, 
 \beaa
 \Next_p^s[\psi,  N_0+N_L] &\les |a|  \Big( B_\de^s[\psi]+ B_\de^s[A]\Big) \les  |a| B_\de^s[\psi, A](\tau_1, \tau_2).
 \eeaa
We conclude, for $\de\leq p\leq 2-\de$, 
\beaa
^{(Mor)} \NN^s[\psi,  N_0+N_L] +\Next^s_p[\psi,  N_0+N_L]&\les&  |a| B_\de^s[\psi, A],
\eeaa
which proves \eqref{eq:step1-}.


\subsection{Proof of the estimate \eqref{eq:step2-}}
\label{subsection:proof-step2}


Here we estimate  $\, ^{(En)} \NN^s[\psi,  N_0]$. We first observe that
\beaa
 ^{(En)} \NN^s[\psi,  N_0](\tau_1, \tau_2) &=& \left|\int_{\MM(\tau_1, \tau_2)} \nab_{\That_\de} ( \dk^{\leq s}\psi )\c  \dk^{\leq s}N\right|\\
 &\les&  \left|\int_{\MM(\tau_1, \tau_2)} \chi \nab_{\That_\de} ( \dk^{\leq s}\psi )\c  \dk^{\leq s}N_0\right|\\
 &&+ \int_{\MM(\tau_1, \tau_2)}(1-\chi) \big|\nab_{\That_\de} ( \dk^{\leq s}\psi ) \big| \big| \dk^{\leq s}N_0\big|,
\eeaa
where $\chi=\chi(r)$ is a smooth cut-off function which is $1$ on $\MM_{trap}$ and vanishes for $r\geq 4m$ and $r\leq r_+(1+\de_{red})$.

Observe that the second integral above can be bounded just like $^{(Mor)}\NN^s[\psi,  N_0](\tau_1, \tau_2) $, and gives
\beaa
\int_{\MM(\tau_1, \tau_2)}(1-\chi) \big|\nab_{\That_\de} ( \dk^{\leq s}\psi ) \big| \big| \dk^{\leq s}N_0\big| \les O(a)B^s_\de[ \psi](\tau_1, \tau_2).
\eeaa 

We are therefore left to estimate $|\int_{\MM(\tau_1, \tau_2)} \chi \nab_{\That_\de} ( \dk^{\leq s}\psi )\c  \dk^{\leq s}N_0|$. Using that  $N_0= O(ar^{-4} ) \psi$, we have 
\beaa
\chi \nab_{\That_\de} ( \dk^{\leq s}\psi )\c  \dk^{\leq s}N_0&=&O(ar^{-4}) \chi \nab_{\That_\de} ( \dk^{\leq s}\psi )\c  \dk^{\leq s} \psi\\
&=&  \D_\a \left(O(ar^{-4}) \chi  |\dk^{\leq s}\psi |^2 \That_\de^\a \right) -O(ar^{-4}) \chi  |\dk^{\leq s}\psi |^2 \D_\a\That_\de^\a\\
&& - O(ar^{-4})\big(\That_\de(r), \That_\de(\cos\th)\big)\big(|\chi'|+r^{-1}|\chi|\big) |\dk^{\leq s}\psi |^2.
\eeaa
Integrating by parts the first term, we therefore conclude\footnote{Notice that weights in $r$ are irrelevant in this estimate since $r$ is bounded on the support of $\chi$.}
\beaa
 ^{(En)} \NN^s[\psi,  N_0](\tau_1, \tau_2) &\les&  \left|\int_{\MM(\tau_1, \tau_2)} \chi \nab_{\That_\de} ( \dk^{\leq s}\psi )\c  \dk^{\leq s}N_0\right|+|a|B^s_\de[\psi](\tau_1, \tau_2)\\
 &\les& |a| \left(\sup_{\tau\in[\tau_1, \tau_2]} E^s[\psi](\tau_1, \tau_2)+ B^s_\de[\psi](\tau_1, \tau_2)\right),
\eeaa
which proves \eqref{eq:step2-}.


\subsection{Proof of the estimate \eqref{eq:step3-}}
\label{subsection:proof-step3}


Here we estimate $\, ^{(En)} \NN^s[\psi,  N_L]$. As above, we write
\beaa
 ^{(En)} \NN^s[\psi,  N_L](\tau_1, \tau_2) &\les&  \left|\int_{\MM(\tau_1, \tau_2)} \chi \nab_{\That_\de} ( \dk^{\leq s}\psi )\c  \dk^{\leq s}N_L\right|\\
 &&+ \int_{\MM(\tau_1, \tau_2)}(1-\chi)\big|\nab_{\That_\de} ( \dk^{\leq s}\psi ) \big| \big| \dk^{\leq s}N_L\big|\\
 &\les&  \left|\int_{\MM(\tau_1, \tau_2)} \chi \nab_{\That_\de} ( \dk^{\leq s}\psi )\c  \dk^{\leq s}N_L\right|+ |a| B^s_\de[ \psi, A](\tau_1, \tau_2),
\eeaa
where $\chi=\chi(r)$ is a smooth cut-off function which is $1$ on $\MM_{trap}$ and vanishes for $r\geq 4m$ and $r\leq r_+(1+\de_{red})$. Therefore, in what follows, we consider the first integral, and notice that $r$ is bounded on the support of $\chi$ so that we can  neglect the powers of $r$.

Also, recall that we have
\beaa
\That_\de&=& \T +\frac{a}{r^2+a^2} \chi_0\left( \de^{-1} \frac{\TT}{r^3} \right) \Z
\eeaa
where $\chi_0( \de^{-1} \frac{\TT}{r^3})=0$ on $\MM_{trap}$. In particular, we can estimate the term involving $\chi_0( \de^{-1} \frac{\TT}{r^3})$ as above, and we obtain 
\beaa
 ^{(En)} \NN^s[\psi,  N_L](\tau_1, \tau_2)  &\les&  \left|\int_{\MM(\tau_1, \tau_2)} \chi \nab_{\T} ( \dk^{\leq s}\psi )\c  \dk^{\leq s}N_L\right|+ |a| B^s_\de[ \psi, A](\tau_1, \tau_2).
\eeaa
Furthermore, since $\psi=\Re(\qf)$, and in view of the definition of $\qf$ in \eqref{eq:definition-qf-again}, we have 
\beaa
\psi &=& \Re\left(q \ov{q}^{3}\left( \nabc_3\nabc_3 A + C_1  \nabc_3A + C_2   A\right)\right)\\
&=& \Re\left(q \ov{q}^{3}\nab_3\nab_3 A\right)+O(a)\dk^{\leq 1}A.
\eeaa
We introduce the notation 
\beaa
\psi_0 &:=& \Re\left(q \ov{q}^{3}\nab_3\nab_3 A\right)
\eeaa
and obtain 
\beaa
\psi &=& \psi_0+O(a)\dk^{\leq 1}A.
\eeaa
Proceeding as above, we infer
\beaa
 ^{(En)} \NN^s[\psi,  N_L](\tau_1, \tau_2)  &\les&  \left|\int_{\MM(\tau_1, \tau_2)} \chi \nab_{\T}(\dk^{\leq s}\psi_0)\c  \dk^{\leq s}N_L\right|+ |a| B^s_\de[ \psi, A](\tau_1, \tau_2).
\eeaa

Next, we focus on the term $\chi \nab_{\T} ( \dk^{\leq s}\psi_0 )\c  \dk^{\leq s}N_L$. We have
\beaa
 \chi \nab_{\T} ( \dk^{\leq s}\psi_0 )\c  \dk^{\leq s}N_L&=& - \chi ( \dk^{\leq s}\psi_0 )\c \nab_{\T} \dk^{\leq s}N_L + \D_\a \big(\chi  ( \dk^{\leq s}\psi_0 )\c  (\dk^{\leq s} N_L)  \T^\a \big)\\
 && -\chi  ( \dk^{\leq s}\psi_0 )\c  (\dk^{\leq s} N_L)  \D_\a\T^\a  -\chi'(r)\T(r)( \dk^{\leq s}\psi_0 )\c  (\dk^{\leq s} N_L)\\
 &=&  - \chi ( \dk^{\leq s}\psi_0 )\c  \dk^{\leq s} \nab_{\T} N_L  +\chi  ( \dk^{\leq s}\psi_0 )\c  \dk^{\leq s}N_L\\
&&  + \D_\a \big(\chi  ( \dk^{\leq s}\psi_0 )\c  (\dk^{\leq s} N_L)  \T^\a \big)\\
 && -\chi  ( \dk^{\leq s}\psi_0 )\c  (\dk^{\leq s} N_L)  \D_\a\T^\a  -\chi'(r)\T(r)( \dk^{\leq s}\psi_0 )\c  (\dk^{\leq s} N_L)
\eeaa
where we used the fact that $[\nab_{\T}, \dk^k]=O(1)\dk^k$ in view of the control of $\Ga_b$ and $\Ga_g$ on the support of $\chi$.

Observe that, since $N_L=O(a) \dkb^{\leq 1} \nab_3^{\leq 1} A$, we have, using integration by parts\footnote{Notice that for the boundary terms, we use the control of $A$ provided by $E_\de[A]$ as well as the one  provided by \eqref{eq:transportA:additionalestimateformissingderivativesenergy}.}, and the fact that $\psi_0 = \psi+O(a)\dk^{\leq 1}A$, 
\beaa
&&\Bigg| \int_{\MM(\tau_1, \tau_2)}\Big[\chi  ( \dk^{\leq s}\psi_0 )\c  \dk^{\leq s}N_L + \D_\a \big(\chi  ( \dk^{\leq s}\psi_0 )\c  (\dk^{\leq s} N_L)  \T^\a \big)\\
 && -\chi  ( \dk^{\leq s}\psi_0 )\c  (\dk^{\leq s} N_L)  \D_\a\T^\a  -\chi'(r)\T(r)( \dk^{\leq s}\psi_0 )\c  (\dk^{\leq s} N_L)\Big]\Bigg| \\
&\les& |a| \left(\sup_{\tau\in[\tau_1, \tau_2]} E_\de[\psi, A]+ B_\de^s[\psi, A](\tau_1, \tau_2)\right)+\ep_0^2\tau_1^{-2-3\dec}.
\eeaa
We deduce from the above that 
\beaa
 ^{(En)} \NN^s[\psi_0,  N_L](\tau_1, \tau_2)  &\les&  \left|\int_{\MM(\tau_1, \tau_2)} \chi ( \dk^{\leq s}\psi_0 )\c  \dk^{\leq s} \nab_{\T} N_L\right|\\
 && +  |a| \left(\sup_{\tau\in[\tau_1, \tau_2]} E_\de[\psi, A]+ B_\de^s[\psi, A](\tau_1, \tau_2)\right)+\ep_0^2\tau_1^{-2-3\dec}.
\eeaa

We now consider the term $ \chi ( \dk^{\leq s}\psi_0 )\c  \dk^{\leq s} \nab_{\T} N_L$. Recall that, see \eqref{eq:definition-N-L-psi-again},
\beaa
N_L  &=& -\Re\left(\frac{8q \ov{q}^{3}r^2}{|q|^6}(a^2\nab_\T+a\nab_\Z)\nab_3A\right) + O(a)\dk^{\leq 1}\a
 \eeaa
 and hence
 \beaa
 \nab_{\T} N_L &=&  -\Re\left(\frac{8q \ov{q}^{3}r^2}{|q|^6}(a^2\nab_\T^2+a\nab_\T\nab_\Z)\nab_3A\right) + O(a)\dk^{\leq 2}\a.
 \eeaa
 Since we have
\beaa
 \T &=& \That -\frac{a}{r^2+a^2}\Z\\
 &=& \Rhat +\frac{\De}{r^2+a^2}e_3 - \frac{a}{r^2+a^2}\Z,
\eeaa 
 we infer
 \beaa
 \nab_{\T} N_L &=& -\Re\left(\frac{8q \ov{q}^{3}r^2}{|q|^6}\left(a^2\nab_\T^2+a\nab_\Rhat\nab_\Z +\frac{a\De}{r^2+a^2}\nab_3\nab_\Z
 - \frac{a^2}{r^2+a^2}\nab_\Z^2\right)\nab_3A\right)\\
 && + O(a)\dk^{\leq 2}\a\\
 &=& -\Re\left(\frac{8q \ov{q}^{3}r^2}{|q|^6}\left(a^2\nab_\T^2 +\frac{a\De}{r^2+a^2}\nab_\Z\nab_3
 - \frac{a^2}{r^2+a^2}\nab_\Z^2\right)\nab_3A\right)\\
 && + O(a)\nab_\Rhat\dk^{\leq 2}\a+ O(a)\dk^{\leq 2}\a
 \eeaa 
and hence 
\beaa
&&\chi ( \dk^{\leq s}\psi_0 )\c  \dk^{\leq s} \nab_{\T} N_L \\
&=& -\chi ( \dk^{\leq s}\psi_0 )\c\dk^{\leq s}\Re\left(\frac{8q \ov{q}^{3}r^2}{|q|^6}\left(a^2\nab_\T^2 +\frac{a\De}{r^2+a^2}\nab_\Z\nab_3
 - \frac{a^2}{r^2+a^2}\nab_\Z^2\right)\nab_3A\right)\\
&&+O(a)\chi\dk^{\leq s}\psi_0\c\nab_\Rhat\dk^{\leq s+2}\a+ O(a)\chi(\dk^{\leq s+2}\a)\c(\dk^{\leq s+2}\a)\\
&=& \chi ( \nab_\T\dk^{\leq s}\psi_0 )\c\dk^{\leq s}\Re\left(\frac{8q \ov{q}^{3}r^2}{|q|^6}a^2\nab_\T\nab_3A\right)\\
&&-\chi ( \dk^{\leq s}\psi_0 )\c\dk^{\leq s}\Re\left(\frac{8q \ov{q}^{3}r^2}{|q|^6}\frac{a\De}{r^2+a^2}\nab_\Z\nab_3^2A\right)\\
&& -\chi ( \nab_\Z\dk^{\leq s}\psi_0 )\c\dk^{\leq s}\Re\left(\frac{8q \ov{q}^{3}r^2}{|q|^6}\frac{a^2}{r^2+a^2}\nab_\Z\nab_3A\right)\\
&&+\D_\mu\Big(O(a)\chi\dk^{\leq s}\psi_0\c\dk^{\leq s+2}\a\big(\Rhat^\mu, \T^\mu, \Z^\mu\big)\Big)+O(a)\chi\nab_\Rhat\dk^{\leq s}\psi_0\dk^{\leq s+2}\a\\
&&+ O(a)(|\chi'|+|\chi|)(\dk^{\leq s+2}\a)\c(\dk^{\leq s+2}\a)
\eeaa 
 Recalling 
 \beaa
\psi_0 = \Re\left(q \ov{q}^{3}\nab_3\nab_3 A\right), \qquad \psi = \psi_0+O(a)\dk^{\leq 1}A,
\eeaa
 we deduce
\beaa
&&\chi ( \dk^{\leq s}\psi_0 )\c  \dk^{\leq s} \nab_{\T} N_L \\
&=& \chi ( \nab_\T\dk^{\leq s}\Re\left(q \ov{q}^{3}\nab_3\nab_3 A\right) )\c\dk^{\leq s}\Re\left(\frac{8q \ov{q}^{3}r^2}{|q|^6}a^2\nab_\T\nab_3A\right)\\
&&-\chi ( \dk^{\leq s}\Re\left(q \ov{q}^{3}\nab_3\nab_3 A\right) )\c\dk^{\leq s}\Re\left(\frac{8q \ov{q}^{3}r^2}{|q|^6}\frac{a\De}{r^2+a^2}\nab_\Z\nab_3^2A\right)\\
&& -\chi ( \nab_\Z\dk^{\leq s}\Re\left(q \ov{q}^{3}\nab_3\nab_3 A\right) )\c\dk^{\leq s}\Re\left(\frac{8q \ov{q}^{3}r^2}{|q|^6}\frac{a^2}{r^2+a^2}\nab_\Z\nab_3A\right)\\
&&+\D_\mu\Big(O(a)\chi\dk^{\leq s}\psi_0\c\dk^{\leq s+2}\a\big(\Rhat^\mu, \T^\mu, \Z^\mu\big)\Big)+O(a)\chi\nab_\Rhat\dk^{\leq s}\psi\dk^{\leq s+2}\a\\
&&+ O(a)(|\chi'|+|\chi|)(\dk^{\leq s+2}\a)\c(\dk^{\leq s+2}\a),
\eeaa  
or
\beaa
&&\chi ( \dk^{\leq s}\psi_0 )\c  \dk^{\leq s} \nab_{\T} N_L \\
&=& \frac{8a^2r^2}{|q|^6}\chi \nab_3\left(\dk^{\leq s}\Re\left(q \ov{q}^{3}\nab_\T\nab_3 A\right) \right)\c\dk^{\leq s}\Re\left(q \ov{q}^{3}\nab_\T\nab_3A\right)\\
&&-\frac{a\De}{r^2+a^2}\frac{8r^2}{|q|^6}\chi ( \dk^{\leq s}\Re\left(q \ov{q}^{3}\nab_3^2A\right) )\c\nab_\Z\left(\dk^{\leq s}\Re\left(q \ov{q}^{3}\nab_3^2A\right)\right)\\
&& -\frac{8r^2}{|q|^6}\frac{a^2}{r^2+a^2}\chi\nab_3\left(\dk^{\leq s}\Re\left(q \ov{q}^{3}\nab_\Z\nab_3 A\right)\right)\c\dk^{\leq s}\Re\left(q \ov{q}^{3}\nab_\Z\nab_3A\right)\\
&&+\D_\mu\Big(O(a)\chi\dk^{\leq s}\psi_0\c\dk^{\leq s+2}\a\big(\Rhat^\mu, \T^\mu, \Z^\mu\big)\Big)+O(a)\chi\nab_\Rhat\dk^{\leq s}\psi\dk^{\leq s+2}\a\\
&&+ O(a)(|\chi'|+|\chi|)(\dk^{\leq s+2}\a)\c(\dk^{\leq s+2}\a),
\eeaa 
 and hence
\beaa
&&\chi ( \dk^{\leq s}\psi_0 )\c  \dk^{\leq s} \nab_{\T} N_L \\
&=& \D_\mu\Big(O(a)\chi\dk^{\leq s+2}\a\c\dk^{\leq s+2}\a\big(e_3^\mu, \Z^\mu\big)\Big) +\D_\mu\Big(O(a)\chi\dk^{\leq s}\psi_0\c\dk^{\leq s+2}\a\big(\Rhat^\mu, \T^\mu, \Z^\mu\big)\Big)\\
&&+O(a)\chi\nab_\Rhat\dk^{\leq s}\psi\dk^{\leq s+2}\a + O(a)(|\chi'|+|\chi|)(\dk^{\leq s+2}\a)\c(\dk^{\leq s+2}\a).
\eeaa  
 We deduce from the above, using integration by parts\footnote{For the boundary terms, we use again the control of $A$ provided by $E_\de[A]$ as well as the one  provided by \eqref{eq:transportA:additionalestimateformissingderivativesenergy}.},
\beaa
 ^{(En)} \NN^s[\psi_0,  N_L](\tau_1, \tau_2)  &\les&  \left|\int_{\MM(\tau_1, \tau_2)} \chi ( \dk^{\leq s}\psi_0 )\c  \dk^{\leq s} \nab_{\T} N_L\right|\\
 && + |a| \left(\sup_{\tau\in[\tau_1, \tau_2]} E_\de[\psi, A]+ B_\de^s[\psi, A](\tau_1, \tau_2)\right) +\ep_0^2\tau_1^{-2-3\dec}\\
 &\les&  |a| \left(\sup_{\tau\in[\tau_1, \tau_2]} E_\de[\psi, A]+ B_\de^s[\psi, A](\tau_1, \tau_2)\right)+\ep_0^2\tau_1^{-2-3\dec}
\eeaa
which concludes the proof of  \eqref{eq:step3-}.


\section{Transport estimates for $A$}\label{section:transport-estimates}


The goal of this section is to prove Proposition \ref{prop:MaiTransportA-steps}, i.e.  show, for $s\leq k_L$ and for all $\de\leq p \leq 2-\de$,
 \bea
\BEF^s_p[A](\tau_1, \tau_2)  &\les&  \BEF^s_p[\psi](\tau_1,\tau_2) 
 + E^s_p[A](\tau_1) +\ep_0^2\tau_1^{-2-3\dec}.
 \eea
 
 To this end, we proceed as follows:
\begin{enumerate}
\item First, we state a general lemma for transport equation in $\nab_3$ in section \ref{sec:basictransportlemmaine3:chap11}, see Lemma \ref{lemma:general-transport-estimate}.

\item Lemma \ref{lemma:general-transport-estimate} is then proved in section \ref{sec:proofoflemma:general-transport-estimate}.

\item Next, we derive estimates for $A$, $\nab_3A$ and $\nab_4A$ in section \ref{sec:controlofAnab3Anab4A:chap11}.

\item Then, we control angular derivatives of $A$ in section \ref{section:angularderivA-transport}.

\item Finally, we conclude the proof of  Proposition \ref{prop:MaiTransportA-steps} in section \ref{section:endoftheproofofprop:MaiTransportA-steps}.
\end{enumerate}


\subsection{General transport estimates}
\lab{sec:basictransportlemmaine3:chap11}


The main result of this section is the following general transport estimates.

  \begin{lemma}\lab{lemma:general-transport-estimate}
  Suppose $\Phi_1, \Phi_2 \in \sk_2(\CCC)$ satisfy the differential relation 
  \bea
  \nabc_3 \Phi_1=\Phi_2.
  \eea
  Also, let $\chi_{nt}=\chi_{nt}(r)$  a smooth cut-off function equal to 0 on $\MM_{trap}$ and equal to 1 on $r\geq 4m$. Then, for every $p\geq \de$, we have
 \bea\lab{eq:general-estimate-commuted-Rhat}
 \begin{split}
&\int_{\MM(\tau_1, \tau_2)}  r^{p-3} \big(r^2|\nab_3\Phi_1|^2+r^2|\nab_4\Phi_1|^2+|\Phi_1|^2\big)\\
&+ \int_{\pr\MM^+(\tau_1, \tau_2)} r^{p-2}\big(r^2\chi_{nt}^2|\nab_4\Phi_1|^2+|\nab_{\Rhat}\Phi_1|^2+|\Phi_1|^2\big)  \\
  \les& \int_{\MM(\tau_1, \tau_2)} r^{p-1}\Big(r^2\chi_{nt}^2|\nab_4\Phi_2|^2+|\nab_{\Rhat}\Phi_2|^2+|\Phi_2|^2\Big)\\
  & + \int_{\Si(\tau_1)} r^{p-2}\big(r^2\chi_{nt}^2|\nab_4\Phi_1|^2+|\nab_{\Rhat}\Phi_1|^2+|\Phi_1|^2\big)\\
  &   + (a^2+\ep^2)\int_{\MM(\tau_1, \tau_2)}r^{p-1} | \nab \Phi_1|^2,
\end{split}
\eea
where $\pr^{+} \MM(\tau_1, \tau_2)$ denotes the future  boundary of $\MM(\tau_1, \tau_2)$, i.e.
\beaa
\pr^{+} \MM(\tau_1, \tau_2)=\AA(\tau_1, \tau_2)\cup\Si(\tau_2)\cup\Si_*(\tau_1, \tau_2). 
\eeaa
\end{lemma}

\begin{remark} 
Observe that estimate \eqref{eq:general-estimate-commuted-Rhat} is conditional with respect to the $\nab \Phi_1$ appearing on the right hand side. 
\end{remark}

In the proof we make us of the divergence theorem in the form, see Lemma 
\ref{lemma:basicdivergenceidentitybianchipairrpweightedestimate} for a general vectorfield $X$,
\bea
\lab{eq:diverggencetheoremX}
-\int_{\pr^+\MM(\tau_1, \tau_2)} \g(X, N) + \int_{\pr^-\MM(\tau_1, \tau_2)} \g(X, N)= \int_{\MM(\tau_1, \tau_2)} \Div(X),
\eea
where  $N$ is  the normal to the boundary   such that $\g(N, e_3) =-1$. 

\begin{remark}
\lab{remark:Boundariesfortransport}
We also make use of the following properties of the boundary,  see  Definition \ref{definition:definition-oftau}, 
 \begin{itemize}
\item On $\AA$ we have
 \beaa
 \g(N_\AA, e_3)= - 1, \qquad   \g(N_\AA, e_4) \leq - \frac{1}{10}\de_\HH, \qquad \g(N_\AA, e_a) =O(\deh),
 \eeaa 
\item  On the boundary $\Si_*$ we have, with $N_*=N_{\Si_*} $,
\beaa
 \g(N_{*}, e_3)= - 1, \qquad   \g(N_{*}, e_4) \leq - 1, \qquad \g(N_{*}, e_a) =O(r^{-1}).
 \eeaa
\item  On the spacelike  boundary $\Si=\Si(\tau)$, $\g(N_{\Si}, N_{\Si}) \leq -\frac{1}{100}\frac{m^2}{r^2}$, 
\beaa
 \g(N_{\Si}, e_4)=-e_4(\tau) \leq -\frac{1}{100} \frac{m^2}{r^2}, \quad  \g(N_{\Si}, e_3)=-e_3(\tau)=-1. 
\eeaa
\end{itemize}
\end{remark}

The following basic lemma will be used in the proof of Lemma  \ref{lemma:general-transport-estimate}.
\begin{lemma}
\lab{Lemma:div-e_3}
For any function $f$, we have
\beaa
\Div \big( f e_3\big)&=& e_3(f)+\left(-\frac{2r}{|q|^2}+\Ga_b\right)f.
\eeaa
\end{lemma}

\begin{proof}
We have
\beaa
\Div ( e_3) &=& \g^{43}\g(\D_4e_3, e_3) +\g^{43}\g(\D_3e_3, e_4)+\g^{bc}\g(\D_be_3, e_c)\\
&=& -\frac{1}{2}4\omb +\trchb = -\frac{2r}{|q|^2}+\trchbc+\Ga_b\\
&=&  -\frac{2r}{|q|^2}+\Ga_b
\eeaa
and hence 
\beaa
\Div ( fe_3)= f \Div(e_3) + e_3(f) =  e_3(f)+\left(-\frac{2r}{|q|^2}+\Ga_b\right)f
\eeaa
as stated. This concludes the proof of Lemma \ref{Lemma:div-e_3}.
 \end{proof}

 
\subsection{Proof of Lemma  \ref{lemma:general-transport-estimate}}
\lab{sec:proofoflemma:general-transport-estimate}


We start with the following lemma.
 \begin{lemma}\lab{lemma:general-transport} 
 Suppose $\Phi_1, \Phi_2 \in \sk_2(\CCC)$ satisfy the differential relation 
 \bea\label{eq:relation-nab3Phi1Phi2}
 \nabc_3 \Phi_1=\Phi_2.
 \eea
  Then for every  $p \geq \de$  we have 
  \bea
\lab{eq:lemma-general-transport1}
 r |q|^{p-4} |\Phi_1|^2&\les &\frac {4}{ p^2} r^{-1}|q|^{p}|\Phi_2|^2   -\frac {2}{ p}\Div(  |q|^{p-2}|\Phi_1|^2 e_3)
\eea
and its integral form
 \bea\label{eq:general-integrated-estimate-e3}
\nn\int_{\MM(\tau_1, \tau_2)}  r^{p-3} |\Phi_1|^2+ \int_{\pr\MM^+(\tau_1, \tau_2)} r^{p-2}|\Phi_1|^2 &\les& \int_{\MM(\tau_1, \tau_2)} r^{p-1}|\Phi_2|^2+ \int_{\Si(\tau_1)} r^{p-2}|\Phi_1|^2.
\eea 
 \end{lemma}
 
\begin{proof} 
Multiplying the relation $\nabc_3\Phi_1=\Phi_2$ by $\ov{\Phi_1}$, we deduce 
\beaa
e_3(|\Phi_1|^2)  &=&2 \Re\big( (\Phi_2 +2s\omb\Phi_1)\c \ov{\Phi_1}\big)\\
&=& 2 \Re\big(\Phi_2\c \ov{\Phi_1}\big)+\Ga_b|\Phi_1|^2.
\eeaa
Multiplying by $|q|^{p-2}$, and using 
\beaa
e_3(|q|)&=& \frac{e_3(|q|^2)}{2|q|}=-\frac{r}{|q|}+O(1)\widecheck{e_3(r)}+O(r^{-1})e_3(\cos\th)=-\frac{r}{|q|}+r\Ga_b,
\eeaa
we deduce
\beaa
2|q|^{p-2} \Re\big( \Phi_2 \c \ov{\Phi_1}\big) &=& e_3( |q|^{p-2}|\Phi_1|^2)  - e_3(|q|^{p-2}) |\Phi_1|^2 +  |q|^{p-2}\Ga_b|\Phi_1|^2\\
&=& e_3( |q|^{p-2}|\Phi_1|^2)  + (p-2) r|q|^{p-4} |\Phi_1|^2 +r^{p-2}\Ga_b |\Phi_1|^2.
\eeaa
In view of Lemma \ref{Lemma:div-e_3}, we write
\beaa
\Div(  |q|^{p-2}|\Phi_1|^2 e_3)&=& e_3( |q|^{p-2}|\Phi_1|^2) +\left(-\frac{2r}{|q|^2}+\Ga_b\right) |q|^{p-2}|\Phi_1|^2\\
&=&2|q|^{p-2} \Re\big( \Phi_2 \c \ov{\Phi_1}\big) - (p-2) r|q|^{p-4} |\Phi_1|^2\\
&&  +\left(-\frac{2r}{|q|^2}+\Ga_b\right) |q|^{p-2}|\Phi_1|^2 +r^{p-2}\Ga_b|\Phi_1|^2\\
&=&2|q|^{p-2} \Re\big( \Phi_2 \c \ov{\Phi_1}\big) - p r|q|^{p-4} |\Phi_1|^2 +r^{p-2}\Ga_b|\Phi_1|^2.
\eeaa
From the above identity we deduce
\beaa
&& p r |q|^{p-4} |\Phi_1|^2 \\ 
&=&2|q|^{p-2} \Re\big( \Phi_2 \c \ov{\Phi_1}\big) -\Div(  |q|^{p-2}|\Phi_1|^2 e_3) +r^{p-2}\Ga_b|\Phi_1|^2\\
&=&2 \Re\big((\la r)^{-1/2}|q|^{p/2} \Phi_2 \c (\lambda r)^{1/2}|q|^{p/2-2}\ov{\Phi_1}\big) -\Div(  |q|^{p-2}|\Phi_1|^2 e_3) +r^{p-2}\Ga_b|\Phi_1|^2\\
&\leq & \lambda r |q|^{p-4}|\Phi_1|^2+ \lambda^{-1}r^{-1}|q|^{p}|\Phi_2|^2   -\Div(  |q|^{p-2}|\Phi_1|^2 e_3) +r^{p-2}\Ga_b|\Phi_1|^2.
\eeaa
We obtain, in view of the control of $\Ga_b$, 
\beaa
pr|q|^{p-4}|\Phi_1|^2 &\leq & \lambda r |q|^{p-4}|\Phi_1|^2+ \lambda^{-1}r^{-1}|q|^{p}|\Phi_2|^2 +O(\ep) r|q|^{p-4}|\Phi_1|^2  -\Div(  |q|^{p-2}|\Phi_1|^2 e_3).
\eeaa
Therefore, for $p \geq\de$, choosing   $\la =\frac  p 2 $, we infer, for $\ep$ sufficiently small, 
\beaa
 r |q|^{p-4} |\Phi_1|^2&\les &\frac {4}{ p^2} r^{-1}|q|^{p}|\Phi_2|^2   -\frac {2}{ p}\Div(  |q|^{p-2}|\Phi_1|^2 e_3)
\eeaa
which is precisely \eqref{eq:lemma-general-transport1}. 

The integral form \eqref{eq:general-integrated-estimate-e3} of the inequality then follows by the divergence   theorem, see \eqref{eq:diverggencetheoremX}, and  Remark \ref{remark:Boundariesfortransport}. This concludes the proof of Lemma \ref{lemma:general-transport}.
\end{proof}
 
Next, we need to control $\nab_{\Rhat}\Phi_1$ for solutions $\Phi_1$ of the transport equation $\nabc_3\Phi_1=\Phi_2$. To this end, we start with the following commutation lemma.
\begin{lemma}\lab{lemma:commutationlemmanab3nabRhat:chap11}
Let $U\in \sk_k$. Then, we have 
\beaa
\left[\nab_3, \frac{r^2+a^2}{|q|^2}\nab_\Rhat\right]U &=& O((a, \ep)r^{-1})\nab U+O(r^{-1}\ep)\nab_4U+O(r^{-2})\nab_3U+O(r^{-3})U
\eeaa
\end{lemma}

\begin{proof}
Recall that $\Rhat$ is given by 
\beaa
\Rhat &=& \frac 1 2 \left( \frac{|q|^2}{r^2+a^2} e_4-\frac{\De}{r^2+a^2}  e_3\right)
\eeaa
so that 
\beaa
\frac{r^2+a^2}{|q|^2}\Rhat &=& \frac 1 2 \left( e_4-\frac{\De}{|q|^2}  e_3\right).
\eeaa
We infer
\beaa
\left[\nab_3, \frac{r^2+a^2}{|q|^2}\nab_\Rhat\right] &=& \frac{1}{2}[\nab_3, \nab_4] -\frac{1}{2}e_3\left(\frac{\De}{|q|^2}\right) \nab_3 \\
&=& \frac{1}{2}[\nab_3, \nab_4] -\frac{1}{2}\left(\pr_r\left(\frac{\De}{|q|^2}\right)e_3(r)+O(r^{-2})e_3(\cos\th)\right) \nab_3 \\
&=& \frac{1}{2}[\nab_3, \nab_4] +\frac{1}{2}\left(\pr_r\left(\frac{\De}{|q|^2}\right)+r^{-1}\Ga_b\right) \nab_3.
\eeaa
Also,  note that the commutation formula for $ [\nab_4, \nab_3]$ of Corollary \ref{cor:comm-gen-B} implies 
\beaa
\, [\nab_4, \nab_3] U &=& \big(O(ar^{-2})+\Ga_b\big)\nab U+ 2\om\nab_3 U  +\Ga_b\nab_4 U +O(r^{-3})U\\
&=&  \big(O(ar^{-2})+\Ga_b\big)\nab U+ 2\left(-\frac{1}{2}\pr_r\left(\frac{\De}{|q|^2}\right)+\omc\right)\nab_3 U  +\Ga_b\nab_4 U +O(r^{-3})U\\
&=&  \big(O(ar^{-2})+\Ga_b\big)\nab U+ \left(-\pr_r\left(\frac{\De}{|q|^2}\right)+\Ga_g\right)\nab_3 U  +\Ga_b\nab_4 U +O(r^{-3})U,
\eeaa
where we used the definition of $\omc$ and the fact that $\omc\in\Ga_g$.  We deduce 
\beaa
\left[\nab_3, \frac{r^2+a^2}{|q|^2}\nab_\Rhat\right]U &=&  \big(O(ar^{-2})+\Ga_b\big)\nab U +(O(r^{-2})+\Ga_g)\nab_3 U  +\Ga_b\nab_4 U +O(r^{-3})U.
\eeaa
Together with the control of $\Ga_g$ and $\Ga_b$, we deduce 
\beaa
\left[\nab_3, \frac{r^2+a^2}{|q|^2}\nab_\Rhat\right]U &=& O((a, \ep)r^{-1})\nab U+O(r^{-1}\ep)\nab_4U+O(r^{-2})\nab_3U+O(r^{-3})U
\eeaa
as stated. This concludes the proof of Lemma \ref{lemma:commutationlemmanab3nabRhat:chap11}. 
\end{proof}  
 
We  next, we consider the following transport lemma.
 \begin{lemma}\lab{eq:lemmageneral-transport-Rhat} 
  Suppose $\Phi_1, \Phi_2 \in \sk_2(\CCC)$ satisfy the relation $\nabc_3 \Phi_1=\Phi_2$. Then, for every $p\geq \de$, we have 
  \bea\lab{eq:precise-transport-Rhat}
 \begin{split}
&\int_{\MM(\tau_1, \tau_2)}  r^{p-3} \big(|\nab_{\Rhat}\Phi_1|^2+|\Phi_1|^2\big)
+ \int_{\pr\MM^+(\tau_1, \tau_2)} r^{p-2}\big(|\nab_{\Rhat}\Phi_1|^2+|\Phi_1|^2\big)  \\
  \les& \int_{\MM(\tau_1, \tau_2)} r^{p-1}\Big(|\nab_{\Rhat}\Phi_2|^2+|\Phi_2|^2\Big) + \int_{\Si(\tau_1)} r^{p-2}\big(|\nab_{\Rhat}\Phi_1|^2+|\Phi_1|^2\big)\\
  & +   \ep^2\int_{\MM(\tau_1, \tau_2)}r^{p-3} | \nab_4 \Phi_1|^2+   (a^2+\ep^2)\int_{\MM(\tau_1, \tau_2)}r^{p-3} | \nab \Phi_1|^2.
\end{split}
\eea
\end{lemma}

\begin{proof}
We commute the transport equation for $\Phi_1$ with $\frac{r^2+a^2}{|q|^2}\nab_\Rhat$. Together with Lemma \ref{lemma:commutationlemmanab3nabRhat:chap11}, we infer
\beaa
\nabc_3\left(\frac{r^2+a^2}{|q|^2}\nab_{\Rhat}\Phi_1\right) &=& \left[\nab_3-2s\omb, \frac{r^2+a^2}{|q|^2}\nab_\Rhat\right]\Phi_1+\frac{r^2+a^2}{|q|^2}\nab_{\Rhat}\Phi_2\\
&=& O((a, \ep)r^{-1})\nab\Phi_1+O(r^{-1}\ep)\nab_4\Phi_1+O(r^{-2})\nab_3\Phi_1+O(r^{-3})\Phi_1\\
&& +2s\frac{r^2+a^2}{|q|^2}(\nab_\Rhat\omb)\Phi_1+O(1)\nab_{\Rhat}\Phi_2
\eeaa
and hence, since $\nab_\Rhat\omb=\dk\Ga_b=O(r^{-1})$, we obtain, using also $\nabc_3\Phi_1=\Phi_2$, 
\beaa
\nabc_3\left(\frac{r^2+a^2}{|q|^2}\nab_{\Rhat}\Phi_1\right) &=& O((a, \ep)r^{-1})\nab\Phi_1+O(r^{-1}\ep)\nab_4\Phi_1+O(r^{-1})\Phi_1\\
&&+O(r^{-2})\Phi_2+O(1)\nab_{\Rhat}\Phi_2.
\eeaa
Applying \eqref{eq:general-integrated-estimate-e3} to this transport equation, we infer
\beaa
 \bsplit
&\int_{\MM(\tau_1, \tau_2)}  r^{p-3} |\nab_{\Rhat}\Phi_1|^2
+ \int_{\pr\MM^+(\tau_1, \tau_2)} r^{p-2}|\nab_{\Rhat}\Phi_1|^2  \\
  \les& \int_{\MM(\tau_1, \tau_2)} r^{p-1}\Big((a^2+\ep^2)r^{-2}|\nab\Phi_1|^2+r^{-2}\ep^2|\nab_4\Phi_1|^2+r^{-2}|\Phi_1|^2+r^{-4}|\Phi_2|^2+|\nab_{\Rhat}\Phi_2|^2\Big)\\
  &+ \int_{\Si(\tau_1)} r^{p-2}|\nab_{\Rhat}\Phi_1|^2.
  \end{split}
\eeaa
Together with \eqref{eq:general-integrated-estimate-e3}, this yields 
\beaa
 \bsplit
&\int_{\MM(\tau_1, \tau_2)}  r^{p-3} \big(|\nab_{\Rhat}\Phi_1|^2+|\Phi_1|^2\big)
+ \int_{\pr\MM^+(\tau_1, \tau_2)} r^{p-2}\big(|\nab_{\Rhat}\Phi_1|^2+|\Phi_1|^2\big)  \\
  \les& \int_{\MM(\tau_1, \tau_2)} r^{p-1}\Big(|\nab_{\Rhat}\Phi_2|^2+|\Phi_2|^2\Big) + \int_{\Si(\tau_1)} r^{p-2}\big(|\nab_{\Rhat}\Phi_1|^2+|\Phi_1|^2\big)\\
  & +   \ep^2\int_{\MM(\tau_1, \tau_2)}r^{p-3} | \nab_4 \Phi_1|^2+   (a^2+\ep^2)\int_{\MM(\tau_1, \tau_2)}r^{p-3} | \nab \Phi_1|^2
  \end{split}
\eeaa
as stated. This concludes the proof of Lemma \ref{eq:lemmageneral-transport-Rhat}.
\end{proof}

We are now ready to prove Lemma \ref{lemma:general-transport-estimate}.

\begin{proof}[Proof of Lemma \ref{lemma:general-transport-estimate}]
We need to derive an estimate for $\chi_{nt}\nab_4\Phi_1$, where $\chi_{nt}=\chi_{nt}(r)$ denotes a smooth cut-off function equal to 0 on $\MM_{trap}$ and equal to 1 on $r\geq 4m$. Note that, in view of Lemma \ref{COMMUTATOR-NAB-C-3-DD-C-HOT},
    \beaa
\, [\nabc_3, \nabc_4] U   &=&  O\big((a+\ep)r^{-1})\nabc  U + O(r^{-3})U
 \eeaa
so that 
\beaa
\nabc_3(\chi_{nt}\nabc_4\Phi_1) &=& \chi_{nt}\nabc_4\Phi_2+\pr_r\chi_{nt} e_3(r)\nabc_4\Phi_1+O\big((a+\ep)r^{-1})\nabc\Phi_1\\
&& + O(r^{-3})\Phi_1.
\eeaa
Since $\nabc_3\Phi_1=\Phi_2$, and since $\nab_4=O(1)\nab_\Rhat+O(1)\nab_3$, we infer
\beaa
\nabc_3(\chi_{nt}\nabc_4\Phi_1) &=& \chi_{nt}\nabc_4\Phi_2+O(1)\pr_r\chi_{nt}\nab_\Rhat\Phi_1+O(1)\pr_r\chi_{nt}\Phi_2\\
&&+O\big((a+\ep)r^{-1})\nabc\Phi_1 + O(r^{-3})\Phi_1.
\eeaa
Since $p+1>p\geq\de$, applying \eqref{eq:general-integrated-estimate-e3} to this transport equation with $p$ replaced by $p+1$, and using the control of $\Phi_1$ provided by Lemma \ref{eq:lemmageneral-transport-Rhat}, we infer
\beaa
 \bsplit
&\int_{\MM(\tau_1, \tau_2)}  r^{p-3} \big(r^2\chi_{nt}^2|\nab_4\Phi_1|^2+|\nab_{\Rhat}\Phi_1|^2+|\Phi_1|^2\big)\\
&+ \int_{\pr\MM^+(\tau_1, \tau_2)} r^{p-2}\big(r^2\chi_{nt}^2|\nab_4\Phi_1|^2+|\nab_{\Rhat}\Phi_1|^2+|\Phi_1|^2\big)  \\
  \les& \int_{\MM(\tau_1, \tau_2)} r^{p-1}\Big(r^2\chi_{nt}^2|\nab_4\Phi_2|^2+|\nab_{\Rhat}\Phi_2|^2+|\Phi_2|^2\Big) \\
  &+ \int_{\Si(\tau_1)} r^{p-2}\big(r^2\chi_{nt}^2|\nab_4\Phi_1|^2+|\nab_{\Rhat}\Phi_1|^2+|\Phi_1|^2\big)\\
  & +   \ep^2\int_{\MM(\tau_1, \tau_2)}r^{p-3} | \nab_4 \Phi_1|^2+   (a^2+\ep^2)\int_{\MM(\tau_1, \tau_2)}r^{p-1} | \nab \Phi_1|^2.
  \end{split}
\eeaa
Together with the fact that $\nabc_3\Phi_1=\Phi_2$, we deduce
\beaa
 \bsplit
&\int_{\MM(\tau_1, \tau_2)}  r^{p-3} \big(r^2|\nab_3\Phi_1|^2+r^2\chi_{nt}^2|\nab_4\Phi_1|^2+|\nab_{\Rhat}\Phi_1|^2+|\Phi_1|^2\big)\\
&+ \int_{\pr\MM^+(\tau_1, \tau_2)} r^{p-2}\big(r^2\chi_{nt}^2|\nab_4\Phi_1|^2+|\nab_{\Rhat}\Phi_1|^2+|\Phi_1|^2\big)  \\
  \les& \int_{\MM(\tau_1, \tau_2)} r^{p-1}\Big(r^2\chi_{nt}^2|\nab_4\Phi_2|^2+|\nab_{\Rhat}\Phi_2|^2+|\Phi_2|^2\Big)\\
  & + \int_{\Si(\tau_1)} r^{p-2}\big(r^2\chi_{nt}^2|\nab_4\Phi_1|^2+|\nab_{\Rhat}\Phi_1|^2+|\Phi_1|^2\big)\\
  & +   \ep^2\int_{\MM(\tau_1, \tau_2)}r^{p-3} | \nab_4 \Phi_1|^2+   (a^2+\ep^2)\int_{\MM(\tau_1, \tau_2)}r^{p-1} | \nab \Phi_1|^2.
  \end{split}
\eeaa
Since $\nab_4$ is spanned by $\nab_\Rhat$ and $\nab_3$, and since $\chi_{nt}=1$ in the region $r\geq 4m$, we infer
\beaa
 \bsplit
&\int_{\MM(\tau_1, \tau_2)}  r^{p-3} \big(r^2|\nab_3\Phi_1|^2+r^2|\nab_4\Phi_1|^2+|\Phi_1|^2\big)\\
&+ \int_{\pr\MM^+(\tau_1, \tau_2)} r^{p-2}\big(r^2\chi_{nt}^2|\nab_4\Phi_1|^2+|\nab_{\Rhat}\Phi_1|^2+|\Phi_1|^2\big)  \\
  \les& \int_{\MM(\tau_1, \tau_2)} r^{p-1}\Big(r^2\chi_{nt}^2|\nab_4\Phi_2|^2+|\nab_{\Rhat}\Phi_2|^2+|\Phi_2|^2\Big)\\
  & + \int_{\Si(\tau_1)} r^{p-2}\big(r^2\chi_{nt}^2|\nab_4\Phi_1|^2+|\nab_{\Rhat}\Phi_1|^2+|\Phi_1|^2\big)\\
  & +   \ep^2\int_{\MM(\tau_1, \tau_2)}r^{p-3} | \nab_4 \Phi_1|^2+   (a^2+\ep^2)\int_{\MM(\tau_1, \tau_2)}r^{p-1} | \nab \Phi_1|^2.
  \end{split}
\eeaa
For $\ep>0$ small enough, this yields 
\beaa
 \bsplit
&\int_{\MM(\tau_1, \tau_2)}  r^{p-3} \big(r^2|\nab_3\Phi_1|^2+r^2|\nab_4\Phi_1|^2+|\Phi_1|^2\big)\\
&+ \int_{\pr\MM^+(\tau_1, \tau_2)} r^{p-2}\big(r^2\chi_{nt}^2|\nab_4\Phi_1|^2+|\nab_{\Rhat}\Phi_1|^2+|\Phi_1|^2\big)  \\
  \les& \int_{\MM(\tau_1, \tau_2)} r^{p-1}\Big(r^2\chi_{nt}^2|\nab_4\Phi_2|^2+|\nab_{\Rhat}\Phi_2|^2+|\Phi_2|^2\Big)\\
  & + \int_{\Si(\tau_1)} r^{p-2}\big(r^2\chi_{nt}^2|\nab_4\Phi_1|^2+|\nab_{\Rhat}\Phi_1|^2+|\Phi_1|^2\big)\\
  &   + (a^2+\ep^2)\int_{\MM(\tau_1, \tau_2)}r^{p-1} | \nab \Phi_1|^2
  \end{split}
\eeaa
as stated in \eqref{eq:general-estimate-commuted-Rhat}. This concludes the proof of Lemma \ref{lemma:general-transport-estimate}.   
\end{proof}


\subsection{Estimates for $A$, $\nab_3A$ and $\nab_4A$}
\lab{sec:controlofAnab3Anab4A:chap11}


Recall from Corollary \ref{cor:systemoftransportequationsforPsiandAfromqf:chap11} that $(\Psi, A)$ satisfies the following system of transport equations 
\bea\lab{eq:thesystemof2transporteuqaitonsAPsipsiactuallyused}
\bsplit
\nabc_3\Psi &=\Big(O(r^{-2})+r^{-1}\Ga_g\Big)\qf  +r^2\dk^{\leq 1}\Ga_b\nabc_3A +r\dk^{\leq 1}\Ga_b A, \\ \nabc_3\left(\frac{(\ov{\tr\Xb})^2}{(\Re(\tr\Xb))^2(\tr\Xb)^2}A\right) & = \Psi+r^2\dk^{\leq 1}(\Ga_b)\c A,
\end{split}
\eea
where $\Psi$ is given in view of Definition \ref{def:PsiforintegrationofAfromqf:chap11} by
\beaa
\Psi &=& \nabc_3\left(\frac{(\ov{\tr\Xb})^2}{(\Re(\tr\Xb))^2(\tr\Xb)^2}A\right) -\frac{r^2}{2}F_2A,
\eeaa
with $F_2$ given by \eqref{eq:defintionofF2appearingtransportequationA:chap11}.

To state the next proposition, we introduce the following partial  norms for $A$ which do not provide control for angular derivatives.
\begin{definition}
\lab{Definition:reducednormsA}
We  define
\beaa
\Bdot_p[A](\tau_1, \tau_2)&=& \int_{\MM(\tau_1, \tau_2)}   r^{p+1} \Big(r^4|\nab_3\nabc_3A|^2+r^4|\nab_4\nab_3A|^2+r^2|\nab_3A|^2+r^2|\nab_4A|^2+|A|^2\Big),\\
\Edot_p[A](\tau)&=& \int_{\Si(\tau)}   r^{p+2} \Big(r^4\chi_{nt}^2|\nab_4\nab_3A|^2+r^2|\nab_{\Rhat}\nab_3A|^2+r^2|\nab_3A|^2+r^2|\nab_4A|^2\\
&&+|A|^2 \Big),\\
\Fdot_p[A](\tau_1, \tau_2)&=&\int_{\AA\cup \Si_*(\tau_1, \tau_2)  }  r^{p+2} \Big(r^4\chi_{nt}^2|\nab_4\nab_3A|^2+r^2|\nab_{\Rhat}\nab_3A|^2+r^2|\nab_3A|^2+r^2|\nab_4A|^2\\
&&+|A|^2 \Big),
\eeaa
where $\chi_{nt}=\chi_{nt}(r)$ denotes a smooth cut-off function equal to 0 on $\MM_{trap}$ and equal to 1 on $r\geq 4m$.

We also  define the combined norms
\bea
\BEFdot_p[A](\tau_1, \tau_2)&=&\Bdot_p[A](\tau_1, \tau_2)+\sup_{\tau\in[ \tau_1, \tau_2]} \Edot_p[A] +\Fdot_p[A](\tau_1, \tau_2).
\eea
\end{definition}

\begin{remark}
Note that the norms  above  do not contain angular derivatives which will have to be recovered later. 
\end{remark}

We prove the following proposition.
\begin{proposition}
\lab{proposition:BBEEestimatesforA}
The following estimates hold true, for all $\de\leq p\le 2-\de$,
\bea\label{eq:transport-estimate-no-angular}
\begin{split}
\BEFdot_p[A](\tau_1, \tau_2)
\les & B_{p}[\qf](\tau_1, \tau_2)+ \Edot_p[A](\tau_1) \\
&+(a^2+\ep^2)\int_{\MM(\tau_1, \tau_2)}r^{p+3}\big(r^2| \nab\nab_3A|^2+| \nab A|^2\big).
\end{split}
\eea
\end{proposition}

\begin{proof}
Recall from \eqref{eq:thesystemof2transporteuqaitonsAPsipsiactuallyused} that we have
\beaa
\nabc_3\Psi &=& \Big(O(r^{-2})+r^{-1}\Ga_g\Big)\qf  +r^2\dk^{\leq 1}\Ga_b\nabc_3A +r\dk^{\leq 1}\Ga_b A.
\eeaa
By applying Lemma \ref{lemma:general-transport-estimate} with $\Phi_1=\Psi$ and 
\beaa
\Phi_2 &=& \Big(O(r^{-2})+r^{-1}\Ga_g\Big)\qf  +r^2\dk^{\leq 1}\Ga_b\nabc_3A +r\dk^{\leq 1}\Ga_b A\\
&=& O(r^{-2})\qf+O(r\ep)\nab_3A+O(\ep)A,
\eeaa
where we used the control of $\Ga_g$ and $\Ga_b$, we obtain for $p'\geq \de$,
\beaa
 \bsplit
&\int_{\MM(\tau_1, \tau_2)}  r^{p'-3} \big(r^2|\nab_3\Psi|^2+r^2|\nab_4\Psi|^2+|\Psi|^2\big)\\
&+ \int_{\pr\MM^+(\tau_1, \tau_2)} r^{p'-2}\big(r^2\chi_{nt}^2|\nab_4\Psi|^2+|\nab_{\Rhat}\Psi|^2+|\Psi|^2\big)  \\
  \les& \int_{\MM(\tau_1, \tau_2)} r^{p'-5}\Big(r^2\chi_{nt}^2|\nab_4\qf|^2+|\nab_{\Rhat}\qf|^2+|\qf|^2\Big)\\
  & + \int_{\Si(\tau_1)} r^{p'-2}\big(r^2\chi_{nt}^2|\nab_4\Psi|^2+|\nab_{\Rhat}\Psi|^2+|\Psi|^2\big) +   (a^2+\ep^2)\int_{\MM(\tau_1, \tau_2)}r^{p'-1} | \nab\Psi|^2\\
  &+\ep^2\int_{\MM(\tau_1, \tau_2)} r^{p'+1}\Big(r^2\chi_{nt}^2|\nab_4\nab_3A|^2+|\nab_{\Rhat}\nab_3A|^2+|\nab_3A|^2\Big)\\
  &+\ep^2\int_{\MM(\tau_1, \tau_2)} r^{p'-1}\Big(r^2\chi_{nt}^2|\nab_4A|^2+|\nab_{\Rhat}A|^2+|A|^2\Big)
  \end{split}
\eeaa
We choose $p'=p+2$ and obtain for $\de \leq p \leq 2-\de$, noticing that $p'\geq 2+\de>\de$ in that case, 
\beaa
 \bsplit
&\int_{\MM(\tau_1, \tau_2)}  r^{p-1} \big(r^2|\nab_3\Psi|^2+r^2|\nab_4\Psi|^2+|\Psi|^2\big)\\
&+ \int_{\pr\MM^+(\tau_1, \tau_2)} r^{p}\big(r^2\chi_{nt}^2|\nab_4\Psi|^2+|\nab_{\Rhat}\Psi|^2+|\Psi|^2\big)  \\
  \les& \int_{\MM(\tau_1, \tau_2)} r^{p-3}\Big(r^2\chi_{nt}^2|\nab_4\qf|^2+|\nab_{\Rhat}\qf|^2+|\qf|^2\Big)\\
  & + \int_{\Si(\tau_1)} r^{p}\big(r^2\chi_{nt}^2|\nab_4\Psi|^2+|\nab_{\Rhat}\Psi|^2+|\Psi|^2\big) +   (a^2+\ep^2)\int_{\MM(\tau_1, \tau_2)}r^{p+1} | \nab\Psi|^2\\
  &+\ep^2\int_{\MM(\tau_1, \tau_2)} r^{p+3}\Big(r^2\chi_{nt}^2|\nab_4\nab_3A|^2+|\nab_{\Rhat}\nab_3A|^2+|\nab_3A|^2\Big)\\
  &+\ep^2\int_{\MM(\tau_1, \tau_2)} r^{p+1}\Big(r^2\chi_{nt}^2|\nab_4A|^2+|\nab_{\Rhat}A|^2+|A|^2\Big)
  \end{split}
\eeaa
In view of the definition of the norm $\Bdot_p[A]$ and $B_p[\qf]$, we infer, for $\de \leq p \leq 2-\de$,
\bea\lab{eq:auxialiaryesitmatecontroltransportA:controlofPsi:chap11}
 \bsplit
&\int_{\MM(\tau_1, \tau_2)}  r^{p-1} \big(r^2|\nab_3\Psi|^2+r^2|\nab_4\Psi|^2+|\Psi|^2\big)\\
&+ \int_{\pr\MM^+(\tau_1, \tau_2)} r^{p}\big(r^2\chi_{nt}^2|\nab_4\Psi|^2+|\nab_{\Rhat}\Psi|^2+|\Psi|^2\big)  \\
  \les& B_p[\qf](\tau_1, \tau_2) + \int_{\Si(\tau_1)} r^{p}\big(r^2\chi_{nt}^2|\nab_4\Psi|^2+|\nab_{\Rhat}\Psi|^2+|\Psi|^2\big)\\
  & +   (a^2+\ep^2)\int_{\MM(\tau_1, \tau_2)}r^{p+1} | \nab\Psi|^2 +\ep^2\Bdot_p[A](\tau_1, \tau_2).
  \end{split}
\eea

Next, recall from \eqref{eq:thesystemof2transporteuqaitonsAPsipsiactuallyused} that we have
\beaa
\nabc_3\left(\frac{(\ov{\tr\Xb})^2}{(\Re(\tr\Xb))^2(\tr\Xb)^2}A\right) & = \Psi+r^2\dk^{\leq 1}(\Ga_b)\c A.
\eeaa
By applying Lemma \ref{lemma:general-transport-estimate} with $\Phi_1=\frac{(\ov{\tr\Xb})^2}{(\Re(\tr\Xb))^2(\tr\Xb)^2}A$ and 
\beaa
\Phi_2 &=& \Psi+r^2\dk^{\leq 1}(\Ga_b)\c A = \Psi+O(r\ep)A,
\eeaa
where we used the control of $\Ga_b$, we obtain, for $p\geq \de$,
\beaa
 \bsplit
&\int_{\MM(\tau_1, \tau_2)}  r^{p+1} \big(r^2|\nab_3A|^2+r^2|\nab_4A|^2+|A|^2\big)\\
&+ \int_{\pr\MM^+(\tau_1, \tau_2)} r^{p+2}\big(r^2\chi_{nt}^2|\nab_4A|^2+|\nab_{\Rhat}A|^2+|A|^2\big)  \\
  \les& \int_{\MM(\tau_1, \tau_2)} r^{p-1}\Big(r^2\chi_{nt}^2|\nab_4\Psi|^2+|\nab_{\Rhat}\Psi|^2+|\Psi|^2\Big)\\
  & + \int_{\Si(\tau_1)} r^{p+2}\big(r^2\chi_{nt}^2|\nab_4A|^2+|\nab_{\Rhat}A|^2+|A|^2\big)   + (a^2+\ep^2)\int_{\MM(\tau_1, \tau_2)}r^{p+3} | \nab A|^2\\
  & + \ep^2\int_{\MM(\tau_1, \tau_2)} r^{p+1}\Big(r^2\chi_{nt}^2|\nab_4A|^2+|\nab_{\Rhat}A|^2+|A|^2\Big)
 \end{split}
\eeaa
and hence, for $\ep$ small enough, 
\beaa
 \bsplit
&\int_{\MM(\tau_1, \tau_2)}  r^{p+1} \big(r^2|\nab_3A|^2+r^2|\nab_4A|^2+|A|^2\big)\\
&+ \int_{\pr\MM^+(\tau_1, \tau_2)} r^{p+2}\big(r^2\chi_{nt}^2|\nab_4A|^2+|\nab_{\Rhat}A|^2+|A|^2\big)  \\
  \les& \int_{\MM(\tau_1, \tau_2)} r^{p-1}\Big(r^2\chi_{nt}^2|\nab_4\Psi|^2+|\nab_{\Rhat}\Psi|^2+|\Psi|^2\Big)\\
  & + \int_{\Si(\tau_1)} r^{p+2}\big(r^2\chi_{nt}^2|\nab_4A|^2+|\nab_{\Rhat}A|^2+|A|^2\big)   + (a^2+\ep^2)\int_{\MM(\tau_1, \tau_2)}r^{p+3} | \nab A|^2.
 \end{split}
\eeaa
Together with \eqref{eq:auxialiaryesitmatecontroltransportA:controlofPsi:chap11}, we infer, for $\de \leq p \leq 2-\de$,
\beaa
 \bsplit
&\int_{\MM(\tau_1, \tau_2)}  r^{p-1} \big(r^2|\nab_3\Psi|^2+r^2|\nab_4\Psi|^2+|\Psi|^2\big)+\int_{\MM(\tau_1, \tau_2)}  r^{p+1} \big(r^2|\nab_3A|^2+r^2|\nab_4A|^2+|A|^2\big)\\
&+ \int_{\pr\MM^+(\tau_1, \tau_2)} r^{p}\big(r^2\chi_{nt}^2|\nab_4\Psi|^2+|\nab_{\Rhat}\Psi|^2+|\Psi|^2\big)  \\
&+ \int_{\pr\MM^+(\tau_1, \tau_2)} r^{p+2}\big(r^2\chi_{nt}^2|\nab_4A|^2+|\nab_{\Rhat}A|^2+|A|^2\big)  \\
  \les& B_p[\qf](\tau_1, \tau_2) + \int_{\Si(\tau_1)} r^{p}\big(r^2\chi_{nt}^2|\nab_4\Psi|^2+|\nab_{\Rhat}\Psi|^2+|\Psi|^2\big)\\
  &+ \int_{\Si(\tau_1)} r^{p+2}\big(r^2\chi_{nt}^2|\nab_4A|^2+|\nab_{\Rhat}A|^2+|A|^2\big)\\
  & +   (a^2+\ep^2)\int_{\MM(\tau_1, \tau_2)}r^{p+1} | \nab\Psi|^2 + (a^2+\ep^2)\int_{\MM(\tau_1, \tau_2)}r^{p+3} | \nab A|^2+\ep^2\Bdot_p[A](\tau_1, \tau_2).
  \end{split}
\eeaa
In view of the definition of $\Psi$, we have
\beaa
\Psi &=& \nabc_3\left(\frac{(\ov{\tr\Xb})^2}{(\Re(\tr\Xb))^2(\tr\Xb)^2}A\right) -\frac{r^2}{2}F_2A\\
&=& \nabc_3\left(\left(\frac{\ov{q}^4}{r^2}+r\Ga_g\right)A\right) +r^2\Ga_b\c A\\
&=& \left(\frac{\ov{q}^4}{r^2}+r\Ga_g\right)\c\nab_3A+(O(r)+r^2\dk^{\leq 1}\Ga_b)\c A
\eeaa
where we used the fact that  $F_2$ given by \eqref{eq:defintionofF2appearingtransportequationA:chap11} satisfies $F_2\in\Ga_b$. Using the control of $\Ga_g$ and $\Ga_b$, this yields 
\beaa
|\nab_3A| &\les& r^{-2}|\Psi|+r^{-1}|A|,\\
|\nab_4\nab_3A| &\les& r^{-2}|\nab_4\Psi|+r^{-3}|\Psi|+r^{-1}|\nab_4A|+r^{-2}|A|,\\
|\nab_\Rhat\nab_3A| &\les& r^{-2}|\nab_\Rhat\Psi|+r^{-2}|\Psi|+r^{-1}|\nab_\Rhat A|+r^{-1}|A|.
\eeaa
Also, using the definition of $\qf$, see \eqref{eq:definition-qf-again}, we have
\beaa
\nabc_3^2A &=& \frac{1}{q\ov{q}^3}\qf+O(r^{-1})\nabc_3A+O(r^{-2})A+\Ga_g\nabc_3A+r^{-1}\Ga_gA,
\eeaa
and hence, using also
\beaa
\nab_3\nabc_3A &=& \nabc_3^2A +\Ga_b\nabc_3A=\nabc_3^2A +\Ga_b\nab_3A+\Ga_b\Ga_bA
\eeaa
and the control of $\Ga_b$, we infer
\beaa
|\nab_3\nabc_3A| &\les& r^{-4}|\qf|+r^{-1}|\nab_3A|+r^{-2}|A|.
\eeaa
Plugging in the above estimate, we infer, for $\de \leq p \leq 2-\de$,
\beaa
 \bsplit
&\int_{\MM(\tau_1, \tau_2)}  r^{p+1} \big(r^4|\nab_3\nabc_3A|^2+r^4|\nab_4\nab_3A|^2+r^2|\nab_3A|^2+r^2|\nab_4A|^2+|A|^2\big)\\
&+ \int_{\pr\MM^+(\tau_1, \tau_2)} r^{p+2}\big(r^4\chi_{nt}^2|\nab_4\nab_3A|^2+r^2|\nab_{\Rhat}\nab_3A|^2+r^2|\nab_3A|^2+r^2|\nab_4A|^2+|A|^2\big)  \\
  \les& B_p[\qf](\tau_1, \tau_2) + \int_{\Si(\tau_1)}r^{p+2}\big(r^4\chi_{nt}^2|\nab_4\nab_3A|^2+r^2|\nab_{\Rhat}\nab_3A|^2+r^2|\nab_3A|^2+r^2|\nab_4A|^2+|A|^2\big)\\
  & +   (a^2+\ep^2)\int_{\MM(\tau_1, \tau_2)}r^{p+1} | \nab\Psi|^2 + (a^2+\ep^2)\int_{\MM(\tau_1, \tau_2)}r^{p+3} | \nab A|^2+\ep^2\Bdot_p[A](\tau_1, \tau_2).
  \end{split}
\eeaa
In view of the definition of $\BEFdot_p[A]$, we obtain, for $\de \leq p \leq 2-\de$,
\beaa
 \bsplit
 \BEFdot_p[A](\tau_1, \tau_2)   \les& B_p[\qf](\tau_1, \tau_2) + \Edot_p[A](\tau_1) +   (a^2+\ep^2)\int_{\MM(\tau_1, \tau_2)}r^{p+3}\big(r^2| \nab\nab_3A|^2+| \nab A|^2\big) \\
 &+\ep^2\Bdot_p[A](\tau_1, \tau_2).
  \end{split}
\eeaa
For $\ep$ small enough,  we deduce, for $\de \leq p \leq 2-\de$, 
\beaa
 \bsplit
 \BEFdot_p[A](\tau_1, \tau_2)   \les& B_p[\qf](\tau_1, \tau_2) + \Edot_p[A](\tau_1) +   (a^2+\ep^2)\int_{\MM(\tau_1, \tau_2)}r^{p+3}\big(r^2| \nab\nab_3A|^2+| \nab A|^2\big) 
  \end{split}
\eeaa
as stated. This concludes the proof of Proposition \ref{proposition:BBEEestimatesforA}.
\end{proof}


\subsection{Estimates for  the angular derivatives  of $A$}
\lab{section:angularderivA-transport}


In order to control angular derivatives  of $A$, we start with the following estimates for $\DDc\hot (\DDbc \c A)$ and $\DDc\hot (\DDbc \c\nab_3A)$. 
\begin{lemma}\lab{lemma:controlofDDchotDDbccAandnab3A:chap11} 
The quantity $\DDc\hot (\DDbc \c A)$ satisfies, for $\de\leq p\leq 2-\de$, 
\beaa
\int_{\MM(\tau_1, \tau_2)}r^{p+5}|\DDc\hot (\DDbc \c A)|^2 &\les& a^2\int_{\MM(\tau_1, \tau_2)}r^{p+1}|\nab A|^2+\Bdot_p[A](\tau_1, \tau_2)\\
&& +\int_{\MM(\tau_1, \tau_2)}r^{p+5}|\err_{\nab A}|^2
\eeaa
and 
\beaa
\int_{\Si(\tau)}r^{p+6}|\DDc\hot (\DDbc \c A)|^2 &\les& a^2\int_{\Si(\tau)}r^{p+2}|\nab A|^2+\Edot_p[A](\tau)\\
&& +\int_{\Si(\tau)}r^{p+6}|\err_{\nab A}|^2,
\eeaa
where
\beaa 
\err_{\nab A} := \Ga_g\c\nab_3A+r^{-1}  \dk^{\leq 1}\big( \Ga_g \c(A, B)\big)  + \Ga_b \c \Ga_g \c A.
\eeaa

Also, the quantity $\DDc\hot (\DDbc \c\nab_3A)$ satisfies, for $\de\leq p\leq 2-\de$, 
\beaa
\int_{\MM(\tau_1, \tau_2)}r^{p+7}|\DDc\hot (\DDbc \c\nab_3A)|^2 &\les& a^2\int_{\MM(\tau_1, \tau_2)}r^{p+3}\big(|\nab\nab_3A|^2+|\nab A|^2\big)+B_p[\qf](\tau_1, \tau_2)\\
&&+\Bdot_p[A](\tau_1, \tau_2) +\int_{\MM(\tau_1, \tau_2)}r^{p+7}|\err_{\nab\nab_3A}|^2
\eeaa
and 
\beaa
\int_{\Si(\tau)}r^{p+8}|\DDc\hot (\DDbc \c\nab_3A)|^2 &\les& a^2\int_{\Si(\tau)}r^{p+4}\big(|\nab\nab_3A|^2+|\nab A|^2\big)+E_p[\qf](\tau)\\
&&+\Edot_p[A](\tau) +\int_{\Si(\tau)}r^{p+8}|\err_{\nab\nab_3A}|^2,
\eeaa
where
\beaa
\err_{\nab\nab_3A} := r^{-4}\Ga_g\c\qf+r^{-1}  \dk^{\leq 1}\big( \Ga_g \c  (\nab_3A, \nab_3B)\big)+r^{-2}  \dk^{\leq 2}\big( \Ga_b \c  (A, B)\big)+\dk^{\leq 1}(\Ga_b \c \Ga_g \c A).
\eeaa
\end{lemma}

\begin{proof}
In view of Proposition \ref{TEUKOLSKY-PROPOSITION}, we have 
 \beaa
 \frac{1}{4}\DDc\hot (\DDbc \c A) &=& \nabc_4\nabc_3A -\left(- \frac 1 2 \tr X -2\ov{\tr X} \right)\nabc_3A+\frac{1}{2}\tr\Xb \nabc_4A\\
 && -\left( 4H+\Hb +\ov{\Hb} \right)\c \nabc A- \left(-\ov{\tr X} \tr \Xb +2\ov{P}\right) A-  H   \hot (\ov{\Hb} \c A)\\ 
 && +r^{-1}  \dk^{\leq 1}\big( \Ga_g \c  B\big)  +  \Ga_b \c \Ga_g \c A,
\eeaa
where we used the fact that the frame used in this chapter satisfies $\Hc\in \Ga_g$ and $\nab_3\Xi\in r^{-1}\dk^{\leq 1}\Ga_g$. We infer
 \beaa
 \frac{1}{4}\DDc\hot (\DDbc \c A) &=& \nabc_4\nabc_3A +O(r^{-1})\nabc_3A+O(r^{-1})\nabc_4A\\
 && +O(ar^{-2})\nabc A +O(r^{-2})A +\Ga_g\c\nab_3A+r^{-1}  \dk^{\leq 1}\big( \Ga_g \c(A, B)\big)\\
 &&  + \Ga_b \c \Ga_g \c A,
\eeaa
where we used again that fact $\Hc\in \Ga_g$. Thus, we obtain 
 \beaa
 \frac{1}{4}\DDc\hot (\DDbc \c A) &=& \nabc_4\nabc_3A +O(r^{-1})\nabc_3A+O(r^{-1})\nabc_4A\\
 && +O(ar^{-2})\nabc A +O(r^{-2})A +\err_{\nab A},\\
\err_{\nab A} &=& \Ga_g\c\nab_3A+r^{-1}  \dk^{\leq 1}\big( \Ga_g \c(A, B)\big)  + \Ga_b \c \Ga_g \c A.
\eeaa
This yields, for $\de\leq p\leq 2-\de$, in view of the definition of the norms $\BEFdot_p[A]$, 
\beaa
\int_{\MM(\tau_1, \tau_2)}r^{p+5}|\DDc\hot (\DDbc \c A)|^2 &\les& a^2\int_{\MM(\tau_1, \tau_2)}r^{p+1}|\nab A|^2+\Bdot_p[A](\tau_1, \tau_2)\\
&& +\int_{\MM(\tau_1, \tau_2)}r^{p+5}|\err_{\nab A}|^2
\eeaa
and 
\beaa
\int_{\Si(\tau)}r^{p+6}|\DDc\hot (\DDbc \c A)|^2 &\les& a^2\int_{\Si(\tau)}r^{p+2}|\nab A|^2+\Edot_p[A](\tau)\\
&& +\int_{\Si(\tau)}r^{p+6}|\err_{\nab A}|^2
\eeaa
which are the stated estimates for $\DDc\hot (\DDbc \c A)$.

Next, we consider the control of $\DDc\hot (\DDbc \c \nab_3A)$. Recall from above the following identity
\beaa
 \frac{1}{4}\DDc\hot (\DDbc \c A) &=& \nabc_4\nabc_3A -\left(- \frac 1 2 \tr X -2\ov{\tr X} \right)\nabc_3A+\frac{1}{2}\tr\Xb \nabc_4A\\
 && -\left( 4H+\Hb +\ov{\Hb} \right)\c \nabc A- \left(-\ov{\tr X} \tr \Xb +2\ov{P}\right) A-  H   \hot (\ov{\Hb} \c A)\\ 
 && +r^{-1}  \dk^{\leq 1}\big( \Ga_g \c  B\big)  +  \Ga_b \c \Ga_g \c A,
\eeaa
which we differentiate w.r.t. $\nabc_3$. Using the following consequences of the null structure equations, taking into account that $\Hc\in\Ga_g$ in the frame used in this chapter, 
\beaa
\nabc_3\tr X &=& O(r^{-2})+r^{-1}\dk^{\leq 1}\Ga_g,\\
\nabc_3\tr\Xb &=& O(r^{-2})+r^{-1}\dk^{\leq 1}\Ga_b, 
\eeaa
as well as 
\beaa
&&\tr X=\frac{2}{r}+O(ar^{-2})+\Ga_g, \qquad \tr\Xb=O(r^{-1})+\Ga_g, \qquad P=O(r^{-3})+r^{-1}\Ga_g, \\
&& H=O(ar^{-2})+\Ga_g, \qquad \Hb=O(ar^{-2})+\Ga_g,
\eeaa
and the fact that $\nab_3\Ga_g=r^{-1}\dk^{\leq 1}\Ga_b$, we have 
\beaa
&&\nabc_3\Bigg\{-\left(- \frac 1 2 \tr X -2\ov{\tr X} \right)\nabc_3A+\frac{1}{2}\tr\Xb \nabc_4A\\
 && -\left( 4H+\Hb +\ov{\Hb} \right)\c \nabc A- \left(-\ov{\tr X} \tr \Xb +2\ov{P}\right) A-  H   \hot (\ov{\Hb} \c A)\\ 
 && +r^{-1}  \dk^{\leq 1}\big( \Ga_g \c  B\big)  +  \Ga_b \c \Ga_g \c A\Bigg\}\\
 &=& \left(\frac{5}{r}+O(ar^{-2})\right)\nabc_3\nabc_3A +O(r^{-2})\nabc_3A +O(r^{-1})\nabc_4\nabc_3A\\
 &&+O(r^{-1})[\nabc_3, \nabc_4]A +O(r^{-2})\nabc_4A+O(ar^{-2})\nabc\nabc_3A\\
 &&+O(ar^{-2})[\nabc_3, \nabc]A+O(ar^{-2})\nab A+O(r^{-2})\nabc_3A+O(r^{-3})A+\Ga_g\c\nabc_3^2A\\
 &&+r^{-1}  \dk^{\leq 1}\big( \Ga_g \c  (\nab_3A, \nab_3B)\big)+r^{-2}  \dk^{\leq 2}\big( \Ga_b \c  (A, B)\big)+\dk^{\leq 1}(\Ga_b \c \Ga_g \c A).
\eeaa
Using the following consequences of the commutation formulas of Lemma \ref{COMMUTATOR-NAB-C-3-DD-C-HOT}, taking into account that $\Hc\in\Ga_g$ in the frame used in this chapter, 
\beaa
\, [\nabc_3, \nabc_4] U  &=& O(ar^{-2})\nab U+O(r^{-3})U+r^{-1}\Ga_g\dk^{\leq 1}U,\\ 
 \, [\nabc_3, \nabc] U  &=& O(ar^{-2})\nab_3U+O(r^{-1})\nab U+O(ar^{-3})U+\Ga_g\c  \dk^{\leq 1} U,
\eeaa
we infer
\beaa
&&\nabc_3\Bigg\{-\left(- \frac 1 2 \tr X -2\ov{\tr X} \right)\nabc_3A+\frac{1}{2}\tr\Xb \nabc_4A\\
 && -\left( 4H+\Hb +\ov{\Hb} \right)\c \nabc A- \left(-\ov{\tr X} \tr \Xb +2\ov{P}\right) A-  H   \hot (\ov{\Hb} \c A)\\ 
 && +r^{-1}  \dk^{\leq 1}\big( \Ga_g \c  B\big)  +  \Ga_b \c \Ga_g \c A\Bigg\}\\
 &=& \left(\frac{5}{r}+O(ar^{-2})\right)\nabc_3\nabc_3A +O(r^{-1})\nabc_4\nabc_3A +O(ar^{-2})\nabc\nabc_3A\\
 && +O(r^{-2})\nabc_3A +O(r^{-2})\nabc_4A+O(ar^{-2})\nab A +O(r^{-3})A\\
 &&+\Ga_g\c\nabc_3^2A+r^{-1}  \dk^{\leq 1}\big( \Ga_g \c  (\nab_3A, \nab_3B)\big)+r^{-2}  \dk^{\leq 2}\big( \Ga_b \c  (A, B)\big)+\dk^{\leq 1}(\Ga_b \c \Ga_g \c A).
\eeaa
Next, using the definition of $\qf$, see \eqref{eq:definition-qf-again}, we have
\beaa
\nabc_3^2A &=& \left(\frac{1}{r^4}+O(ar^{-5})\right)\qf+O(r^{-1})\nabc_3A+O(r^{-2})A+\Ga_g\nabc_3A+r^{-1}\Ga_gA,
\eeaa
and hence
\beaa
&&\nabc_3\Bigg\{-\left(- \frac 1 2 \tr X -2\ov{\tr X} \right)\nabc_3A+\frac{1}{2}\tr\Xb \nabc_4A\\
 && -\left( 4H+\Hb +\ov{\Hb} \right)\c \nabc A- \left(-\ov{\tr X} \tr \Xb +2\ov{P}\right) A-  H   \hot (\ov{\Hb} \c A)\\ 
 && +r^{-1}  \dk^{\leq 1}\big( \Ga_g \c  B\big)  +  \Ga_b \c \Ga_g \c A\Bigg\}\\
 &=& \left(\frac{5}{r}+O(ar^{-6})\right)\qf +O(r^{-1})\nabc_4\nabc_3A +O(ar^{-2})\nabc\nabc_3A\\
 && +O(r^{-2})\nabc_3A +O(r^{-2})\nabc_4A+O(ar^{-2})\nab A +O(r^{-3})A\\
 &&+r^{-4}\Ga_g\c\qf+r^{-1}  \dk^{\leq 1}\big( \Ga_g \c  (\nab_3A, \nab_3B)\big)+r^{-2}  \dk^{\leq 2}\big( \Ga_b \c  (A, B)\big)+\dk^{\leq 1}(\Ga_b \c \Ga_g \c A).
\eeaa
Coming back to the identity 
\beaa
 \frac{1}{4}\DDc\hot (\DDbc \c A) &=& \nabc_4\nabc_3A -\left(- \frac 1 2 \tr X -2\ov{\tr X} \right)\nabc_3A+\frac{1}{2}\tr\Xb \nabc_4A\\
 && -\left( 4H+\Hb +\ov{\Hb} \right)\c \nabc A- \left(-\ov{\tr X} \tr \Xb +2\ov{P}\right) A-  H   \hot (\ov{\Hb} \c A)\\ 
 && +r^{-1}  \dk^{\leq 1}\big( \Ga_g \c  B\big)  +  \Ga_b \c \Ga_g \c A,
\eeaa
and differentiating it w.r.t. $\nabc_3$, we infer
\beaa
 && \frac{1}{4}\nabc_3\DDc\hot (\DDbc \c A)\\ 
 &=& \nabc_3\nabc_4\nabc_3A  +\left(\frac{5}{r}+O(ar^{-2})\right)\qf +O(r^{-1})\nabc_4\nabc_3A +O(ar^{-2})\nabc\nabc_3A\\
 && +O(r^{-2})\nabc_3A +O(r^{-2})\nabc_4A+O(ar^{-2})\nab A +O(r^{-3})A\\
 &&+r^{-4}\Ga_g\c\qf+r^{-1}  \dk^{\leq 1}\big( \Ga_g \c  (\nab_3A, \nab_3B)\big)+r^{-2}  \dk^{\leq 2}\big( \Ga_b \c  (A, B)\big)+\dk^{\leq 1}(\Ga_b \c \Ga_g \c A).
\eeaa

Next, using again the commutation formula
\beaa
\, [\nabc_3, \nabc_4] U  &=& O(ar^{-2})\nab U+O(r^{-3})U+r^{-1}\Ga_g\dk^{\leq 1}U,
\eeaa
we have
\beaa
[\nabc_3, \nabc_4]\nabc_3A &=& O(ar^{-2})\nab\nabc_3A +O(r^{-3})\nabc_3A+r^{-1}\Ga_g\dk^{\leq 1}\nabc_3A
\eeaa
and hence
\beaa
 && \frac{1}{4}\nabc_3\DDc\hot (\DDbc \c A)\\ 
 &=& \nabc_4\nabc_3^2A  +\left(\frac{5}{r}+O(ar^{-2})\right)\qf +O(r^{-1})\nabc_4\nabc_3A +O(ar^{-2})\nabc\nabc_3A\\
 && +O(r^{-2})\nabc_3A +O(r^{-2})\nabc_4A+O(ar^{-2})\nab A +O(r^{-3})A\\
 &&+r^{-4}\Ga_g\c\qf+r^{-1}  \dk^{\leq 1}\big( \Ga_g \c  (\nab_3A, \nab_3B)\big)+r^{-2}  \dk^{\leq 2}\big( \Ga_b \c  (A, B)\big)+\dk^{\leq 1}(\Ga_b \c \Ga_g \c A).
\eeaa
Using again the fact that 
\beaa
\nabc_3^2A &=& \left(\frac{1}{r^4}+O(ar^{-5})\right)\qf+O(r^{-1})\nabc_3A+O(r^{-2})A+\Ga_g\nabc_3A+r^{-1}\Ga_gA,
\eeaa
we infer
\beaa
 && \frac{1}{4}\nabc_3\DDc\hot (\DDbc \c A)\\ 
 &=& \left(\frac{1}{r^4}+O(ar^{-5})\right)\nabc_4\qf+\left(\frac{1}{r^5}+O(ar^{-6})\right)\qf +O(r^{-1})\nabc_4\nabc_3A\\
 && +O(ar^{-2})\nabc\nabc_3A +O(r^{-2})\nabc_3A +O(r^{-2})\nabc_4A+O(ar^{-2})\nab A +O(r^{-3})A\\
 &&+r^{-4}\Ga_g\c\qf+r^{-1}  \dk^{\leq 1}\big( \Ga_g \c  (\nab_3A, \nab_3B)\big)+r^{-2}  \dk^{\leq 2}\big( \Ga_b \c  (A, B)\big)+\dk^{\leq 1}(\Ga_b \c \Ga_g \c A)
\eeaa
and hence
\beaa
 && \frac{1}{4}\nabc_3\DDc\hot (\DDbc \c A)\\ 
 &=& O(r^{-5})e_4(r\qf)+O(r^{-6})\qf +O(r^{-1})\nabc_4\nabc_3A +O(ar^{-2})\nabc\nabc_3A\\
 && +O(r^{-2})\nabc_3A +O(r^{-2})\nabc_4A+O(ar^{-2})\nab A +O(r^{-3})A\\
 &&+r^{-4}\Ga_g\c\qf+r^{-1}  \dk^{\leq 1}\big( \Ga_g \c  (\nab_3A, \nab_3B)\big)+r^{-2}  \dk^{\leq 2}\big( \Ga_b \c  (A, B)\big)+\dk^{\leq 1}(\Ga_b \c \Ga_g \c A).
\eeaa

Next, using the following consequences of the commutation formulas of Lemma \ref{COMMUTATOR-NAB-C-3-DD-C-HOT}, taking into account that $\Hc\in\Ga_g$ in the frame used in this chapter, 
\beaa
 \, [\nabc_3, \DD\hot] U  &=& -\frac{1}{2}\tr X\DD\hot U + O(ar^{-2})\nab_3U+O(ar^{-3})U+\Ga_g\nab_3U+r^{-1}\Ga_g\c  \dk^{\leq 1} U,\\
  \, [\nabc_3, \ov{\DD}\c] U  &=& -\frac{1}{2}\ov{\tr X}\,\ov{\DD}\c U + O(ar^{-2})\nab_3U+O(ar^{-3})U+\Ga_g\nab_3U+r^{-1}\Ga_g\c  \dk^{\leq 1} U,
\eeaa
we have
\beaa
&& \nabc_3\DDc\hot (\DDbc \c A) - \DDc\hot (\DDbc \c\nabc_3A)\\ 
&=& [\nabc_3,\DDc\hot](\DDbc \c A)+\DDc\hot[\nabc_3, \DDbc \c]A\\
&=& \left(-\frac{1}{2}\tr X\DD\hot  + O(ar^{-2})\nab_3+O(ar^{-3})+\Ga_g\nab_3+r^{-1}\Ga_g\c  \dk^{\leq 1}\right)(\DDbc \c A)\\
&&+ \DDc\hot\left(-\frac{1}{2}\ov{\tr X}\,\ov{\DD}\c A + O(ar^{-2})\nab_3A+O(ar^{-3})A+\Ga_g\nab_3A+r^{-1}\Ga_g\c  \dk^{\leq 1}A\right)\\
&=& -\trch\DDc\hot (\DDbc \c A) +O(ar^{-2})\nab\nab_3A+O(ar^{-2})[\nab_3, \nab]A +O(ar^{-3})\nab_3A\\
&&+O(ar^{-3})\nab A+O(ar^{-4})A +r^{-1}\dk^{\leq 1}(\Ga_g\nab_3A)+r^{-2}\dk^{\leq 2}(\Ga_g\c A)\\
&=& O(r^{-1})\DDc\hot (\DDbc \c A) +O(ar^{-2})\nab\nab_3A +O(ar^{-3})\nab_3A +O(ar^{-3})\nab A\\
&& +O(ar^{-4})A+r^{-1}\dk^{\leq 1}(\Ga_g\nab_3A)+r^{-2}\dk^{\leq 2}(\Ga_g\c A).
\eeaa
Together with the above, we infer
\beaa
 && \DDc\hot (\DDbc \c\nabc_3A)\\ 
 &=& O(r^{-5})e_4(r\qf)+O(r^{-6})\qf +O(r^{-1})\nab_4\nab_3A +O(ar^{-2})\nab\nab_3A\\
 && +O(r^{-1})\DDc\hot (\DDbc \c A)  +O(r^{-2})\nab_3A +O(r^{-2})\nab_4A+O(ar^{-2})\nab A +O(r^{-3})A\\
 &&+r^{-4}\Ga_g\c\qf+r^{-1}  \dk^{\leq 1}\big( \Ga_g \c  (\nab_3A, \nab_3B)\big)+r^{-2}  \dk^{\leq 2}\big( \Ga_b \c  (A, B)\big)+\dk^{\leq 1}(\Ga_b \c \Ga_g \c A).
\eeaa
Recalling 
 \beaa
 \frac{1}{4}\DDc\hot (\DDbc \c A) &=& \nabc_4\nabc_3A +O(r^{-1})\nabc_3A+O(r^{-1})\nabc_4A\\
 && +O(ar^{-2})\nabc A +O(r^{-2})A +\err_{\nab A},
\eeaa
we deduce 
\beaa
 && \DDc\hot (\DDbc \c\nabc_3A)\\ 
 &=& O(r^{-5})e_4(r\qf)+O(r^{-6})\qf +O(r^{-1})\nab_4\nab_3A +O(ar^{-2})\nab\nab_3A\\
 &&  +O(r^{-2})\nab_3A +O(r^{-2})\nab_4A+O(ar^{-2})\nab A +O(r^{-3})A\\
 &&+r^{-4}\Ga_g\c\qf+r^{-1}  \dk^{\leq 1}\big( \Ga_g \c  (\nab_3A, \nab_3B)\big)+r^{-2}  \dk^{\leq 2}\big( \Ga_b \c  (A, B)\big)+\dk^{\leq 1}(\Ga_b \c \Ga_g \c A)\\
 &&+O(r^{-1})\err_{\nab A}.
\eeaa
In view of the form of $\err_{\nab A}$, i.e. 
\beaa
\err_{\nab A} &=& \Ga_g\c\nab_3A+r^{-1}  \dk^{\leq 1}\big( \Ga_g \c(A, B)\big)  + \Ga_b \c \Ga_g \c A,
\eeaa
we obtain 
\beaa
 && \DDc\hot (\DDbc \c\nabc_3A)\\ 
 &=& O(r^{-5})e_4(r\qf)+O(r^{-6})\qf +O(r^{-1})\nab_4\nab_3A +O(ar^{-2})\nab\nab_3A\\
 &&  +O(r^{-2})\nab_3A +O(r^{-2})\nab_4A+O(ar^{-2})\nab A +O(r^{-3})A +\err_{\nab\nab_3A}
\eeaa
where 
\beaa
\err_{\nab\nab_3A} &=& r^{-4}\Ga_g\c\qf+r^{-1}  \dk^{\leq 1}\big( \Ga_g \c  (\nab_3A, \nab_3B)\big)+r^{-2}  \dk^{\leq 2}\big( \Ga_b \c  (A, B)\big)+\dk^{\leq 1}(\Ga_b \c \Ga_g \c A).
\eeaa
This yields, for $\de\leq p\leq 2-\de$, in view of the definition of the norms $\BEF_p[\qf]$ and $\BEFdot_p[A]$, 
\beaa
\int_{\MM(\tau_1, \tau_2)}r^{p+7}|\DDc\hot (\DDbc \c\nab_3A)|^2 &\les& a^2\int_{\MM(\tau_1, \tau_2)}r^{p+3}\big(|\nab\nab_3A|^2+|\nab A|^2\big)+B_p[\qf](\tau_1, \tau_2)\\
&&+\Bdot_p[A](\tau_1, \tau_2) +\int_{\MM(\tau_1, \tau_2)}r^{p+7}|\err_{\nab\nab_3A}|^2
\eeaa
and 
\beaa
\int_{\Si(\tau)}r^{p+8}|\DDc\hot (\DDbc \c\nab_3A)|^2 &\les& a^2\int_{\Si(\tau)}r^{p+4}\big(|\nab\nab_3A|^2+|\nab A|^2\big)+E_p[\qf](\tau)\\
&&+\Edot_p[A](\tau) +\int_{\Si(\tau)}r^{p+8}|\err_{\nab\nab_3A}|^2
\eeaa
which are the stated estimates for $\DDc\hot (\DDbc \c\nab_3A)$. This concludes the proof of Lemma \ref{lemma:controlofDDchotDDbccAandnab3A:chap11}.
\end{proof}

To control angular derivatives of $A$, we will also need the following lemma.
\begin{lemma}
\lab{Lemma:EstimateangderivativesA-LL(a)}
The following estimate holds true,  for any  $U\in\sk_2(\CCC)$ and $S \subset\MM$,
\beaa
\int_S  \Big(|\nab U|^2  +   r^{-2} |U|^2 \Big)  &\les&\left| \int_S   \ov{U}\c \DD\hot (\DDb \c U)  \right| +O(a^2)  \int_S r^{-2} |( \nab_3, \nab_4)  U|^2. 
\eeaa
\end{lemma}

\begin{proof}
In view of Lemma \ref{Lemma:DDhot(DDb)=lap}  we  write
\beaa
\DD \hot (\DDb \c U)&=& 2 \lap_2U - 4 \Kh U -i  \big(\atrch\nab_3+\atrchb \nab_4\big)U
\eeaa
where $\Kh=- \frac 14  \trch \trchb-\frac 1 4 \atrch \atrchb+\frac 1 2 \chih \c \chibh-  \frac 1 4 \rho$. Therefore
\beaa
 \ov{U} \c \DD \hot (\DDb \c U)  &=&  2  \ov{U} \c \lap_2U - 4 \Kh |U|^2  -i  \ov{U}\c \big(\atrch\nab_3+\atrchb \nab_4\big)U\\
 &=&  -|\nab U|^2 +  \nab^a(\ov{U}  \c \nab_aU) - 4 \Kh |U|^2  -i  \ov{U}\c \big(\atrch\nab_3+\atrchb \nab_4\big)U.
\eeaa
Proceeding as in the proof of Proposition \ref{Prop:HodgeThmM8}, with the help of Lemma  \ref{Lemma:projectionS},
\beaa
 \nab^a(\ov{U} \c \nab_a U)= \div^\S\big( U\c \nab U)+O(ar^{-1} )
   \ov{U}\c (\nab_3, \nab_4, \nab) U.
\eeaa
Thus, by integration on $S$, we deduce
\beaa
\int_S\big(|\nab U|^2  + 4 \Kh |U|^2\big)&\les &\left|\int_S  \ov{U} \c \DD \hot (\DDb \c U)\right| +  a  \int_S r^{-1 } \big|U\big| \big|  (\nab_3, \nab_4, \nab) U\big|.
\eeaa
Since $\Kh=  r^{-2} +O(a^2 r^{-4} )$, we deduce
 \beaa
\int_S\big(|\nab U|^2  +  r^{-2} |U|^2\big)&\les &\left|\int_S  \ov{U} \c \DD \hot (\DDb \c U)\right| +  a^2  \int_S r^{-2} \big|  (\nab_3, \nab_4)U\big|
\eeaa 
as stated.
\end{proof}

We now obtain the desired control of $\nab A$ and $\nab\nab_3A$ as a corollary of Lemma \ref{lemma:controlofDDchotDDbccAandnab3A:chap11} and Lemma \ref{Lemma:EstimateangderivativesA-LL(a)}. 
\begin{corollary}\lab{cor:controlofnabAandnabnab3A:chap11}
We have, for all $\de\leq p\leq 2-\de$, 
   \beaa
         \bsplit
    \int_{\MM(\tau_1, \tau_2)} r^{p+3}|\nab A|^2  \les & \Bdot_p[A](\tau_1, \tau_2) +\int_{\MM(\tau_1, \tau_2)}r^{p+5}|\err_{\nab A}|^2, 
    \end{split}
     \eeaa 
     \beaa
         \bsplit
    \int_{\Si(\tau)} r^{p+4}|\nab A|^2     \les &  \Edot_p[A](\tau) +\int_{\Si(\tau)}r^{p+6}|\err_{\nab A}|^2,
    \end{split}
     \eeaa
    \beaa
         \bsplit
    \int_{\MM(\tau_1, \tau_2)} r^{p+5}|\nab\nab_3A|^2     \les & B_p[\qf](\tau_1, \tau_2) +\Bdot_p[A](\tau_1, \tau_2) +\int_{\MM(\tau_1, \tau_2)}r^{p+7}|\err_{\nab\nab_3A}|^2,
    \end{split}
     \eeaa
     and 
     \beaa
         \bsplit
    \int_{\Si(\tau)} r^{p+6}|\nab\nab_3A|^2     \les & E_p[\qf](\tau)+\Edot_p[A](\tau) +\int_{\Si(\tau)}r^{p+8}|\err_{\nab\nab_3A}|^2,
    \end{split}
     \eeaa
     where the error terms $\err_{\nab A}$ and $\err_{\nab\nab_3A}$ are defined in Lemma \ref{lemma:controlofDDchotDDbccAandnab3A:chap11}.
\end{corollary}

\begin{proof}
According to the  elliptic type   estimates of  Lemma \ref{Lemma:EstimateangderivativesA-LL(a)}  we have
   for any  $ S \subset\MM$.
\beaa
\int_S  \Big(|\nab A|^2  +   r^{-2} |A|^2 \Big)  &\les&\left| \int_SA\c \DD\hot (\DDb \c A)  \right| +(a^2+\ep^2)  \int_S r^{-2} |( \nabc_3, \nabc_4)A|^2. 
\eeaa
We  deduce,  for all $\de\leq p\leq 2-\de$, 
     \beaa
         \bsplit
    \int_{\MM(\tau_1, \tau_2)} r^{p+3}|\nab A|^2 \les &  \int_{\MM(\tau_1, \tau_2)} r^{p+1}  |A|^2 +    \int_{\MM(\tau_1, \tau_2)} r^{p+5}  \big|\DD \hot (\DDb \c A)\big|^2     \\
    &+(a^2+\ep^2)  \int_{\MM(\tau_1, \tau_2)}   r^{p+1} |( \nab_3, \nab_4)A|^2\\
    \les & \Bdot_p[A](\tau_1, \tau_2)+\int_{\MM(\tau_1, \tau_2)} r^{p+5}  \big|\DD \hot (\DDb \c A)\big|^2, 
    \end{split}
     \eeaa
     and 
     \beaa
         \bsplit
    \int_{\Si(\tau)} r^{p+4}|\nab A|^2 \les&   \int_{\Si(\tau)} r^{p+2}  |A|^2 +    \int_{\Si(\tau)} r^{p+6}  \big|\DD \hot (\DDb \c A)\big|^2     \\
    &+(a^2+\ep^2)  \int_{\Si(\tau)}   r^{p+2} |( \nab_3, \nab_4)A|^2\\
    \les &  \Edot_p[A](\tau) +    \int_{\Si(\tau)} r^{p+6}  \big|\DD \hot (\DDb \c A)\big|^2 
    \end{split}
     \eeaa
         Next, recall from Lemma \ref{lemma:controlofDDchotDDbccAandnab3A:chap11} that we have, for $\de\leq p\leq 2-\de$, 
\beaa
\int_{\MM(\tau_1, \tau_2)}r^{p+5}|\DDc\hot (\DDbc \c A)|^2 &\les& a^2\int_{\MM(\tau_1, \tau_2)}r^{p+1}|\nab A|^2+\Bdot_p[A](\tau_1, \tau_2)\\
&& +\int_{\MM(\tau_1, \tau_2)}r^{p+5}|\err_{\nab A}|^2
\eeaa
and 
\beaa
\int_{\Si(\tau)}r^{p+6}|\DDc\hot (\DDbc \c A)|^2 &\les& a^2\int_{\Si(\tau)}r^{p+2}|\nab A|^2+\Edot_p[A](\tau)\\
&& +\int_{\Si(\tau)}r^{p+6}|\err_{\nab A}|^2,
\eeaa
see Lemma \ref{lemma:controlofDDchotDDbccAandnab3A:chap11} for the definition of $\err_{\nab A}$. 
We infer,  for all $\de\leq p\leq 2-\de$, 
     \beaa
         \bsplit
    \int_{\MM(\tau_1, \tau_2)} r^{p+3}|\nab A|^2  \les & a^2\int_{\MM(\tau_1, \tau_2)}r^{p+1}|\nab A|^2+\Bdot_p[A](\tau_1, \tau_2) +\int_{\MM(\tau_1, \tau_2)}r^{p+5}|\err_{\nab A}|^2, 
    \end{split}
     \eeaa
     and 
     \beaa
         \bsplit
    \int_{\Si(\tau)} r^{p+4}|\nab A|^2     \les &   a^2\int_{\Si(\tau)}r^{p+2}|\nab A|^2+\Edot_p[A](\tau) +\int_{\Si(\tau)}r^{p+6}|\err_{\nab A}|^2.
    \end{split}
     \eeaa
For $a$ small enough, we infer, for all $\de\leq p\leq 2-\de$, 
   \beaa
         \bsplit
    \int_{\MM(\tau_1, \tau_2)} r^{p+3}|\nab A|^2  \les & \Bdot_p[A](\tau_1, \tau_2) +\int_{\MM(\tau_1, \tau_2)}r^{p+5}|\err_{\nab A}|^2, 
    \end{split}
     \eeaa
     and 
     \beaa
         \bsplit
    \int_{\Si(\tau)} r^{p+4}|\nab A|^2     \les &  \Edot_p[A](\tau) +\int_{\Si(\tau)}r^{p+6}|\err_{\nab A}|^2,
    \end{split}
     \eeaa
as stated. 

Next, proceeding as above, with $\DD \hot (\DDb \c A)$ replaced by $\DD \hot (\DDb \c\nab_3A)$, we have, for all $\de\leq p\leq 2-\de$, 
     \beaa
         \bsplit
    \int_{\MM(\tau_1, \tau_2)} r^{p+5}|\nab\nab_3A|^2 \les &  \int_{\MM(\tau_1, \tau_2)} r^{p+3}|\nab_3A|^2 +    \int_{\MM(\tau_1, \tau_2)} r^{p+7}  \big|\DD \hot (\DDb \c\nab_3A)\big|^2     \\
    &+(a^2+\ep^2)  \int_{\MM(\tau_1, \tau_2)}   r^{p+3} |( \nab_3^2, \nab_4\nab_3)A|^2\\
    \les & \Bdot_p[A](\tau_1, \tau_2)+\int_{\MM(\tau_1, \tau_2)} r^{p+7}  \big|\DD \hot (\DDb \c \nab_3A)\big|^2, 
    \end{split}
     \eeaa
     and 
     \beaa
         \bsplit
    \int_{\Si(\tau)} r^{p+6}|\nab\nab_3A|^2 \les&   \int_{\Si(\tau)} r^{p+4}  |\nab_3A|^2 +    \int_{\Si(\tau)} r^{p+8}  \big|\DD \hot (\DDb \c\nab_3A)\big|^2     \\
    &+(a^2+\ep^2)  \int_{\Si(\tau)}   r^{p+4} |( \nab_3^2, \nab_4\nab_3)A|^2\\
    \les &  \Edot_p[A](\tau) +    \int_{\Si(\tau)} r^{p+8}  \big|\DD \hot (\DDb \c\nab_3A)\big|^2 
    \end{split}
     \eeaa
    Next, recall from Lemma \ref{lemma:controlofDDchotDDbccAandnab3A:chap11} that we have, for $\de\leq p\leq 2-\de$, 
\beaa
\int_{\MM(\tau_1, \tau_2)}r^{p+7}|\DDc\hot (\DDbc \c\nab_3A)|^2 &\les& a^2\int_{\MM(\tau_1, \tau_2)}r^{p+3}\big(|\nab\nab_3A|^2+|\nab A|^2\big)+B_p[\qf](\tau_1, \tau_2)\\
&&+\Bdot_p[A](\tau_1, \tau_2) +\int_{\MM(\tau_1, \tau_2)}r^{p+7}|\err_{\nab\nab_3A}|^2
\eeaa
and 
\beaa
\int_{\Si(\tau)}r^{p+8}|\DDc\hot (\DDbc \c\nab_3A)|^2 &\les& a^2\int_{\Si(\tau)}r^{p+4}\big(|\nab\nab_3A|^2+|\nab A|^2\big)+E_p[\qf](\tau)\\
&&+\Edot_p[A](\tau) +\int_{\Si(\tau)}r^{p+8}|\err_{\nab\nab_3A}|^2,
\eeaa
see Lemma \ref{lemma:controlofDDchotDDbccAandnab3A:chap11} for the definition of $\err_{\nab\nab_3A}$. 
We infer,  for all $\de\leq p\leq 2-\de$, 
    \beaa
         \bsplit
    \int_{\MM(\tau_1, \tau_2)} r^{p+5}|\nab\nab_3A|^2     \les &\, a^2\int_{\MM(\tau_1, \tau_2)}r^{p+3}\big(|\nab\nab_3A|^2+|\nab A|^2\big)+B_p[\qf](\tau_1, \tau_2)\\
&+\Bdot_p[A](\tau_1, \tau_2) +\int_{\MM(\tau_1, \tau_2)}r^{p+7}|\err_{\nab\nab_3A}|^2
    \end{split}
     \eeaa
     and 
     \beaa
         \bsplit
    \int_{\Si(\tau)} r^{p+6}|\nab\nab_3A|^2     \les & \, a^2\int_{\Si(\tau)}r^{p+4}\big(|\nab\nab_3A|^2+|\nab A|^2\big)+E_p[\qf](\tau)\\
&+\Edot_p[A](\tau) +\int_{\Si(\tau)}r^{p+8}|\err_{\nab\nab_3A}|^2. 
    \end{split}
     \eeaa
For $a$ small enough, we infer, for all $\de\leq p\leq 2-\de$, 
    \beaa
         \bsplit
    \int_{\MM(\tau_1, \tau_2)} r^{p+5}|\nab\nab_3A|^2     \les & B_p[\qf](\tau_1, \tau_2) +\Bdot_p[A](\tau_1, \tau_2) +\int_{\MM(\tau_1, \tau_2)}r^{p+7}|\err_{\nab\nab_3A}|^2
    \end{split}
     \eeaa
     and 
     \beaa
         \bsplit
    \int_{\Si(\tau)} r^{p+6}|\nab\nab_3A|^2     \les & E_p[\qf](\tau)+\Edot_p[A](\tau) +\int_{\Si(\tau)}r^{p+8}|\err_{\nab\nab_3A}|^2
    \end{split}
     \eeaa
as stated. This concludes the proof of Corollary \ref{cor:controlofnabAandnabnab3A:chap11}.
\end{proof}

To conclude this section, we provide the control of the error terms $\err_{\nab A}$ and $\err_{\nab\nab_3A}$ defined in  Lemma \ref{lemma:controlofDDchotDDbccAandnab3A:chap11}.
\begin{lemma}\lab{lemma:controloftheerrortermserrnabAanderrnabnab3A:chap11}
The error terms $\err_{\nab A}$ and $\err_{\nab\nab_3A}$ defined in  Lemma \ref{lemma:controlofDDchotDDbccAandnab3A:chap11} satisfy the following estimates, for all $s\leq k_L$, and for all $\de\leq p\leq 2-\de$, 
   \beaa
   \int_{\MM(\tau_1, \tau_2)}r^{p+5}|\dk^s\err_{\nab A}|^2 +\int_{\MM(\tau_1, \tau_2)}r^{p+7}|\dk^s\err_{\nab\nab_3A}|^2 \les \ep_0^2\tau_1^{-2-3\dec},
     \eeaa 
     and
     \beaa
       \int_{\Si(\tau)}r^{p+6}|\dk^s\err_{\nab A}|^2 + \int_{\Si(\tau)}r^{p+8}|\dk^s\err_{\nab\nab_3A}|^2 &\les& \ep_0^2\tau^{-2-3\dec}.
     \eeaa
\end{lemma}

\begin{proof}
Recall that the error terms $\err_{\nab A}$ and $\err_{\nab\nab_3A}$ are given respectively by
\beaa 
\err_{\nab A} = \Ga_g\c\nab_3A+r^{-1}  \dk^{\leq 1}\big( \Ga_g \c(A, B)\big)  + \Ga_b \c \Ga_g \c A
\eeaa
and 
\beaa
\err_{\nab\nab_3A} = r^{-4}\Ga_g\c\qf+r^{-1}  \dk^{\leq 1}\big( \Ga_g \c  (\nab_3A, \nab_3B)\big)+r^{-2}  \dk^{\leq 2}\big( \Ga_b \c  (A, B)\big)+\dk^{\leq 1}(\Ga_b \c \Ga_g \c A).
\eeaa
In view of the Bianchi identities for $\nab_3A$ and $\nab_3B$, we have 
\beaa
\nab_3A=O(r^{-1})\dk^{\leq 1}B+O(r^{-1})A+r^{-3}\Ga_g, \qquad \nab_3B=r^{-2}\dk^{\leq 1}\Ga_g. 
\eeaa
Together with the fact that  $\qf\in r\dk^{\leq 2}\Ga_g$, we infer
\beaa 
\err_{\nab A} = r^{-1}  \dk^{\leq 1}\big( \Ga_g \c(A, B)\big) + r^{-3}\Ga_g\c\Ga_g  + \Ga_b \c \Ga_g \c A
\eeaa
and 
\beaa
\err_{\nab\nab_3A} = r^{-3}  \dk^{\leq 2}\big( \Ga_g \c  \Ga_g\big)+r^{-2}  \dk^{\leq 2}\big( \Ga_b \c  (A, B)\big)+\dk^{\leq 1}(\Ga_b \c \Ga_g \c A).
\eeaa
Using the bootstrap assumptions for $\Ga_g$, $\Ga_b$, $A$ and $B$, we infer, for all $s\leq k_L$,  
 \beaa
 &&  \int_{\MM(\tau_1, \tau_2)}r^{7-\de}|\dk^s\err_{\nab A}|^2 +\int_{\MM(\tau_1, \tau_2)}r^{9-\de}|\dk^s\err_{\nab\nab_3A}|^2\\
   &\les& \ep^4\left(\int_{\tau_1}^{+\infty}\frac{d\tau}{\tau^{3+3\dec}}\right)\left(\int_{r\geq r_+(1-\deh)}\frac{dr}{r^{1+\de}}\right)\\
   &\les& \ep_0^2\tau_1^{-2-3\dec}
     \eeaa
and 
 \beaa
      \int_{\Si(\tau)}r^{8-\de}|\dk^s\err_{\nab A}|^2 + \int_{\Si(\tau)}r^{10-\de}|\dk^s\err_{\nab\nab_3A}|^2   &\les& \ep^4\tau^{-2-3\dec}\left(\int_{r\geq r_+(1-\deh)}\frac{dr}{r^{1+\de}}\right)\\
       &\les& \ep_0^2\tau^{-2-3\dec}
     \eeaa
     as stated.
\end{proof}


\subsection{End of the proof of Proposition    \ref{prop:MaiTransportA-steps}}
\lab{section:endoftheproofofprop:MaiTransportA-steps}


In view of Corollary \ref{cor:controlofnabAandnabnab3A:chap11} and Lemma \ref{lemma:controloftheerrortermserrnabAanderrnabnab3A:chap11}, we have, for all $\de\leq p\leq 2-\de$, 
  \beaa
         \bsplit
    \int_{\MM(\tau_1, \tau_2)} r^{p+3}\big( r^2|\nab\nab_3A|^2+|\nab A|^2\big)  \les &  B_p[\qf](\tau_1, \tau_2) +\Bdot_p[A](\tau_1, \tau_2) +\ep_0^2\tau_1^{-2-3\dec}.
    \end{split}
     \eeaa
  Also,  recall from Proposition \ref{proposition:BBEEestimatesforA} that we have, for $\de\leq p\le 2-\de$, 
\beaa
\begin{split}
\BEFdot_p[A](\tau_1, \tau_2)
\les & B_{p}[\qf](\tau_1, \tau_2)+ \Edot_p[A](\tau_1) \\
&+(a^2+\ep^2)\int_{\MM(\tau_1, \tau_2)}r^{p+3}\big(r^2| \nab\nab_3A|^2+| \nab A|^2\big).
\end{split}
\eeaa
Putting the two estimates together, and using the smallness of $\ep$ and $a$, we infer
\beaa
\begin{split}
& \BEFdot_p[A](\tau_1, \tau_2) +\int_{\MM(\tau_1, \tau_2)}r^{p+3}\big(r^2| \nab\nab_3A|^2+| \nab A|^2\big)\\
\les & B_{p}[\qf](\tau_1, \tau_2)+ \Edot_p[A](\tau_1) +\ep_0^2\tau_1^{-2-3\dec}.
\end{split}
\eeaa 
Noticing that $F_p[A](\tau_1, \tau_2) = \Fdot_p[A](\tau_1, \tau_2)$, $E_p[A](\tau) = \Edot_p[A](\tau)$ and 
\beaa
B_p[A](\tau_1, \tau_2) &=& \Bdot_p[A](\tau_1, \tau_2) +\int_{\MM(\tau_1, \tau_2)}r^{p+3}\big(r^2| \nab\nab_3A|^2+| \nab A|^2\big),
\eeaa  
 we infer, since $\qf=\psi+i\dual\psi$, 
\beaa
\BEF_p[A](\tau_1, \tau_2) &\les & B_{p}[\psi](\tau_1, \tau_2)+ \Edot_p[A](\tau_1) +\ep_0^2\tau_1^{-2-3\dec}
\eeaa   
 which is \eqref{eq:transportA} in the case $s=0$.   
  
Next, in view of Corollary \ref{cor:controlofnabAandnabnab3A:chap11} and Lemma \ref{lemma:controloftheerrortermserrnabAanderrnabnab3A:chap11}, we have, for all $\de\leq p\leq 2-\de$,
     \beaa
         \bsplit
    \int_{\Si(\tau)} r^{p+2}\big(r^4|\nab\nab_3A|^2+r^2|\nab A|^2\big)     \les &  E_p[\qf](\tau)+\Edot_p[A](\tau) +\ep_0^2\tau^{-2-3\dec}.
    \end{split}
     \eeaa
Moreover, we have in view of the above, for all $\de\leq p\leq 2-\de$,
\beaa
\sup_{\tau\in[\tau_1, \tau_2]}E_p[A](\tau)  &\les& \BEF_p[A](\tau_1, \tau_2) \les  B_{p}[\psi](\tau_1, \tau_2)+ \Edot_p[A](\tau_1) +\ep_0^2\tau_1^{-2-3\dec}
\eeaa   
and hence, we have in particular, for any $\tau\in[\tau_1, \tau_2]$, , for all $\de\leq p\leq 2-\de$,  
\beaa
&&\int_{\Si(\tau)}   r^{p+2} \Big(r^4\chi_{nt}^2|\nab_4\nab_3A|^2+r^2|\nab_{\Rhat}\nab_3A|^2+r^2|\nab_3A|^2+|A|^2\Big)\\ 
&\les&  B_{p}[\psi](\tau_1, \tau_2)+ \Edot_p[A](\tau_1) +\ep_0^2\tau_1^{-2-3\dec}.
\eeaa
Also, recall that we have
\beaa
\nabc_3^2A &=& O(r^{-4})\qf+O(r^{-1})\nabc_3A+O(r^{-2})A+\Ga_g\nabc_3A+r^{-1}\Ga_gA,
\eeaa
and hence
\beaa
\nab_3\nabc_3A &=& \nabc_3^2A+\Ga_b\nabc_3A\\
&=&\nabc_3^2A+\Ga_b\nab_3A+\Ga_b\Ga_b A\\
&=& O(r^{-4})\qf+O(r^{-1})\nabc_3A+O(r^{-2})A+\Ga_g\nabc_3A+r^{-1}\Ga_gA\\
&& +\Ga_b\nab_3A+\Ga_b\Ga_b A\\
&=& O(r^{-4})\qf+O(r^{-1})\nab_3A+O(r^{-2})A
\eeaa
where we used the control of $\Ga_b$ and $\Ga_g$. This yields, for any $\tau\in[\tau_1, \tau_2]$, for all $\de\leq p\leq 1-\de$,   
\beaa
\int_{\Si(\tau)}   r^{p+6}|\nab_3\nabc_3A|^2 &\les& E_{p}[\qf](\tau_1, \tau_2)+\int_{\Si(\tau)}   r^{p+2} \Big(r^2|\nab_3A|^2+|A|^2\Big).
\eeaa
Grouping the above estimates, we infer, for any $\tau\in[\tau_1, \tau_2]$, for all $\de\leq p\leq 2-\de$,   
    \beaa
         \bsplit
    &\int_{\Si(\tau)} r^{p+2}\big(r^{\min(4, 5-\de-p)}|\nab_3\nabc_3A|^2+ r^4\chi_{nt}^2|\nab_4\nab_3A|^2+r^2|\nab_{\Rhat}\nab_3A|^2\\
    &+r^2|\nab_3A|^2+|A|^2+r^4|\nab\nab_3A|^2+r^2|\nab A|^2\big)\\    
     \les &  EB_p[\psi](\tau_1, \tau_2)+E_p[A](\tau_1) +\ep_0^2\tau_1^{-2-3\dec}.
    \end{split}
     \eeaa
In particular, since $\chi_{nt}=1$ on $r\geq 4m$, and since $\nab_4$ is spanned by $\nab_3$ and $\nab_\Rhat$ on $r\leq 4m$, we infer, for any $\tau\in[\tau_1, \tau_2]$, for all $\de\leq p\leq 2-\de$,   
    \beaa
         \bsplit
    &\int_{\Si(\tau)} r^{p+2}\big(r^{\min(4, 5-\de-p)}|\nab_3\nabc_3A|^2+ r^4|\nab_4\nab_3A|^2+r^4|\nab\nab_3A|^2+r^2|\nab A|^2\big)\\    
     \les &  EB_p[\psi](\tau_1, \tau_2)+E_p[A](\tau_1) +\ep_0^2\tau_1^{-2-3\dec}.
    \end{split}
     \eeaa
 which is  \eqref{eq:transportA:additionalestimateformissingderivativesenergy} in the case $s=0$.

It remains to recover  \eqref{eq:transportA} and  \eqref{eq:transportA:additionalestimateformissingderivativesenergy} for $1\leq s\leq k_L$. To this end, we proceed as follows:
\begin{enumerate}
\item We argue by iteration assuming that \eqref{eq:transportA} and  \eqref{eq:transportA:additionalestimateformissingderivativesenergy} hold for some $0\leq s\leq k_L-1$. It it true for $s=0$ by the above, and our goal is to prove that \eqref{eq:transportA} and  \eqref{eq:transportA:additionalestimateformissingderivativesenergy} hold with $s$ replaced by $s+1$. 

\item We commute the system of transport equations \eqref{eq:thesystemof2transporteuqaitonsAPsipsiactuallyused}, i.e.
\beaa
\bsplit
\nabc_3\Psi &=\Big(O(r^{-2})+r^{-1}\Ga_g\Big)\qf  +r^2\dk^{\leq 1}\Ga_b\nabc_3A +r\dk^{\leq 1}\Ga_b A, \\ \nabc_3\left(\frac{(\ov{\tr\Xb})^2}{(\Re(\tr\Xb))^2(\tr\Xb)^2}A\right) & = \Psi+r^2\dk^{\leq 1}(\Ga_b)\c A,
\end{split}
\eeaa
with $\Lieb_\T$, $\ov{q}\,\ov{\DDc}\c$ and $\nabc_4$. In view of the commutation formulas of Lemma \ref{lemma:basicpropertiesLiebTfasdiuhakdisug:chap9} and Lemma \ref{COMMUTATOR-NAB-C-3-DD-C-HOT}, we have, for $U\in\sk_2$, 
\beaa
[\nabc_3, \Lieb_\T]U &=& [\nab_3, \Lieb_\T]U +4[\omb, \Lieb_\T]U = [\Lieb_\T, \dk]U +\dk^{\leq 1}(\Ga_b) U\\
&=& \dk^{\leq 1}(\Ga_b U),
\eeaa 
\beaa
 \, [\nabc_3, q\ov{\DDc}\c] U&=& q [\nabc_3, \ov{q}\,\ov{\DDc}\c] U +e_3(q)\ov{\DDc}\c U\\
 &=& - \frac 1 2q\left(\ov{\tr\Xb} -\frac{2}{q}e_3(q)\right) \ov{\DDc} \c U +2q\ov{\tr\Xb}\,\ov{H}\c U+q\ov{H} \c \nabc_3 U\\
 &&+ \Ga_b \c  \dk^{\leq 1} U\\
 &=& O(ar^{-1})\nabc_3U + O(ar^{-2})U +r\Ga_g\nabc_3U+ \Ga_b \c  \dk^{\leq 1} U,
\eeaa
and
\beaa
 \, [\nabc_3, \nabc_4] U&=& O(ar^{-2})\nab U+O(r^{-3})U+r^{-1}\Ga_g\dk^{\leq 1}U.
 \eeaa
This yields the commuted systems 
\beaa
\bsplit
\nabc_3\Lieb_\T\Psi  =&\Big(O(r^{-2})+r^{-1}\Ga_g\Big)\Lieb_\T\qf+\Big(O(r^{-3})+r^{-2}\dk^{\leq 1}\Ga_b\Big)\qf  \\
& +r^2\dk^{\leq 1}\Ga_b\nabc_3\Lieb_\T A +r\dk^{\leq 1}\Ga_b\Lieb_\T A\\
& +r^2\dk^{\leq 2}\Ga_b\nabc_3A +r\dk^{\leq 2}\Ga_b A, \\ 
\nabc_3\Lieb_\T\left(\frac{(\ov{\tr\Xb})^2}{(\Re(\tr\Xb))^2(\tr\Xb)^2}A\right)  =& \Lieb_\T\Psi+r^2\dk^{\leq 1}(\Ga_b)\c \Lieb_\T A+r^2\dk^{\leq 2}(\Ga_b)\c A,
\end{split}
\eeaa
\beaa
\bsplit
\nabc_3 q\ov{\DDc}\c\Psi &=\Big(O(r^{-1})+\Ga_g\Big)\nab\qf+\Big(O(r^{-2})+r^{-1}\dk^{\leq 1}\Ga_g\Big)\qf \\
& +r\dk^{\leq 1}\Ga_b\c\nab\nabc_3A+r^2\dk^{\leq 2}\Ga_b\nabc_3A \\
&+\dk^{\leq 1}\Ga_b\nab A+r\dk^{\leq 2}\Ga_b A,\\ 
\nabc_3 q\ov{\DDc}\c\left(\frac{(\ov{\tr\Xb})^2}{(\Re(\tr\Xb))^2(\tr\Xb)^2}A\right) & = q\ov{\DDc}\c\Psi +O(ar)\nabc_3A + O(a)A \\
&+r^2\Ga_g\nabc_3A+ r^2\Ga_b \c  \dk^{\leq 1}A +r^2\dk^{\leq 2}(\Ga_b)\c A, 
\end{split}
\eeaa
and
\beaa
\bsplit
\nabc_3\nabc_4\Psi =&\Big(O(r^{-2})+r^{-1}\Ga_g\Big)\nabc_4\qf+\Big(O(r^{-3})+r^{-2}\dk^{\leq 1}\Ga_g\Big)\qf\\
&  +r^2\dk^{\leq 1}\Ga_b\nabc_3A +r\dk^{\leq 1}\Ga_b A\\
& +O(ar^{-2})\nab \Psi+O(r^{-3})\Psi+r^{-1}\Ga_g\dk^{\leq 1}\Psi,\\ 
\nabc_3\nabc_4\left(\frac{(\ov{\tr\Xb})^2}{(\Re(\tr\Xb))^2(\tr\Xb)^2}A\right)  =& \nabc_4\Psi +O(a)\nab A+O(r^{-1})A+r^2\dk^{\leq 1}(\Ga_b)\c A.
\end{split}
\eeaa

\item Using the iteration assumption for these commuted systems, and using the original system to recover the $\nab_3$ derivative in the redshift region, we infer that \eqref{eq:transportA} and  \eqref{eq:transportA:additionalestimateformissingderivativesenergy} hold for $s$ derivatives with $A$ replaced with $(\Lieb_\T, q\ov{\DDc}\c, \nabc_4, \chi_{red}\nabc_3)\Ab$. Together with:
\begin{enumerate}
\item the link between $\Lieb_\T$ and $\nab_\T$ of Lemma \ref{lemma:basicpropertiesLiebTfasdiuhakdisug:chap9},

\item the Hodge elliptic estimates of Proposition \ref{Prop:HodgeThmM8},

\item the fact that $(\nab_\T, r\nab_4, \dkb, \chi_{red}\nab_3)$ span $\dk$,
\end{enumerate} 
and using the iteration assumption to absorb lower order terms in differentiability, we infer that \eqref{eq:transportA} and  \eqref{eq:transportA:additionalestimateformissingderivativesenergy} hold for $s$ derivatives with $A$ replaced with $\dk^{\leq 1}A$. In particular, \eqref{eq:transportA} and  \eqref{eq:transportA:additionalestimateformissingderivativesenergy} hold with $s$ replaced by $s+1$. Thus, by iteration,  \eqref{eq:transportA} and  \eqref{eq:transportA:additionalestimateformissingderivativesenergy} hold for all $s$ such that $0\leq s\leq k_L$. This end the proof of Proposition \ref{prop:MaiTransportA-steps}.  
\end{enumerate}


\section{Proof of Theorem \ref{THEOREM:GRW1-P-WEAK-PSIWC}}
\lab{section:ProofThm-gRW1-p-weak-psiwc}


  In this section we prove Theorem \ref{THEOREM:GRW1-P-WEAK-PSIWC}, i.e.  we 
   establish the estimate,  for $2\le s\le \kl-1$, for all $-1+\de\leq q \leq 1-\de$,
   \beaa
       \begin{split}
       \BEF_q^s[\psiwc](\tau_1, \tau_2)  &\les
       E_q^s[\psiwc](\tau_1)+E_{\max(q, \de)}^{s+1}[\psi, A](\tau_1)+\widecheck{\NN}_{q}^{s}[\psiwc, N_\err](\tau_1, \tau_2)\\
       &+\NN_{\max(q, \de)}^{s+1}[\psi, N_\err](\tau_1, \tau_2)  +\ep_0^2\tau_1^{-2-3\dec},
       \end{split}
       \eeaa
       where
       \beaa
       \widecheck{\NN}_{q}^{s}[\psiwc, N_\err](\tau_1, \tau_2)&=& \int_{\Mext} r^{q+2} \widecheck{\nab}_4 \dk^{\leq s}\psiwc \c \left( \nab_4 \dk^{\leq s}N_\err + \frac 3 r \dk^{\leq s} N_\err \right).
       \eeaa
       
\begin{proof}[Proof of Theorem \ref{THEOREM:GRW1-P-WEAK-PSIWC}]
 According to Theorem \ref{THEOREM:GENRW2-Q}   we have,   for solutions $\psi\in\sk_2$  of  \eqref{eq:Gen.RW} on $\MM$,   for all $-1+\de<q\le 1-\de$, $s\le\kl -1$, 
  \bea
  \lab{eq:theorem:GenRW2-q-Chap10.}
            \BEF^s_q[\psiwc](\tau_1, \tau_2)
       \les  \Et_q^s[\psiwc](\tau_1)+\NNt_q^s[\psiwc, N](\tau_1, \tau_2) + \NN_{\max\{q,\de\}}^{s+1}[\psi, N](\tau_1, \tau_2),
       \eea
 where  $\psiwc=r^2(e_4\psi+\frac{r}{|q|^2}\psi)$, $  \Et_q^s[\psiwc](\tau)=E_q^s[\psiwc](\tau)+ E^{s+1}_{\max\{q,\de\} }[\psi](\tau)$ and 
                 \beaa
         \bsplit
          \NNt_q^s[\psiwc, N](\tau_1, \tau_2)         & =   \int_{\MM_{\ge R}(\tau_1,\tau_2) } r^{q+2} \dk^{\le s} \psiwc \c \left(\nab_4 \dk^{\le s }N+ \frac 3 r \dk^{\le s} N\right),
         \end{split}
         \eeaa       
          where, see section \ref{section:FullgRWchap11}, $ N=N_0+N_L+N_\err$. In view of \eqref{eq:theorem:GenRW2-q-Chap10.}, it remains to estimate the terms $\NN_{\max(q, \de)}^{s+1}[\psi, N_0+N_L](\tau_1, \tau_2)$ and $\widecheck{\NN}_{q}^{s}[\psiwc, N_0+N_L](\tau_1, \tau_2)$. This is done in the following steps.

\bigskip

{\bf Step 1.} We first estimate the term $\NN_{\max(q, \de)}^{s+1}[\psi, N_0+N_L](\tau_1, \tau_2)$. Note that 
\beaa
&&\NN_{\max(q, \de)}^{s+1}[\psi, N_0+N_L](\tau_1, \tau_2)\\ 
&\les&  ^{(Mor)}\NN^{s+1}[\psi, N_0+N_L](\tau_1, \tau_2)+\Next_{\max(q, \de)}^{s+1}[\psi, N_0+N_L](\tau_1, \tau_2)\\
&&+^{(En)}\NN^{s+1}[\psi, N_0](\tau_1, \tau_2)+^{(En)}\NN^{s+1}[\psi, N_L](\tau_1, \tau_2),
\eeaa
which together with \eqref{eq:step1-}-\eqref{eq:step3-} implies
\bea
\NN_{\max(q, \de)}^{s+1}[\psi, N_0+N_L](\tau_1, \tau_2) &\les& \BEF^{s+1}_{\max(q, \de)}[\psi, A](\tau_1, \tau_2)    +\ep_0^2\tau_1^{-2-3\dec}.
\eea

{\bf Step 2.} Next, we estimate the term $\widecheck{\NN}_{q}^{s}[\psiwc, N_0+N_L](\tau_1, \tau_2)$. 
Notice  that 
 \beaa
\widecheck{\NN}_{q}^{s}[\psiwc, N](\tau_1, \tau_2) &=&\left|  \int_{\Mext(\tau_1, \tau_2)} r^{q+2} \widecheck{\nab}_4 \dk^{\leq s}\psiwc \c \left( \nab_4 \dk^{\leq s}N + \frac 3 r \dk^{\leq s}N\right)\right|\\
&\les&\left|  \int_{\Mext(\tau_1, \tau_2)} r^{q+1} \widecheck{\nab}_4 \dk^{\leq s}\psiwc \c \big( \dk^{\leq s+1}N\big)\right|.
\eeaa
We thus have
\beaa
\widecheck{\NN}_{q}^{s}[\psiwc, N_0+N_L](\tau_1, \tau_2) &\les& \widecheck{\NN}_{q}^{s}[\psiwc, N_0](\tau_1, \tau_2)+\widecheck{\NN}_{q}^{s}[\psiwc, N_L](\tau_1, \tau_2),\\
\widecheck{\NN}_{q}^{s}[\psiwc, N_0](\tau_1, \tau_2) &\les&\left|  \int_{\Mext(\tau_1, \tau_2)} r^{q+1} \widecheck{\nab}_4 \dk^{\leq s}\psiwc \c \big( \dk^{\leq s+1}N_0\big)\right|,\\
\widecheck{\NN}_{q}^{s}[\psiwc, N_L](\tau_1, \tau_2) &\les&\left|  \int_{\Mext(\tau_1, \tau_2)} r^{q+1} \widecheck{\nab}_4 \dk^{\leq s}\psiwc \c \big( \dk^{\leq s+1}N_L\big)\right|,
\eeaa
and we estimate below the terms $\widecheck{\NN}_{q}^{s}[\psiwc, N_0](\tau_1, \tau_2)$ and $\widecheck{\NN}_{q}^{s}[\psiwc, N_L](\tau_1, \tau_2)$ separately.

{\bf Step 2a.} We have
 \beaa
\widecheck{\NN}_{q}^{s}[\psiwc, N_0](\tau_1, \tau_2) &\les& \bigg(  \int_{\Mext(\tau_1, \tau_2)}   r^{q-3} |\dk^{\le s+1} \psiwc|^2\bigg)^{1/2} \Bigg( \int_{\Mext(\tau_1, \tau_2)} r^{q+5}   |\dk^{\le s+1 }(N_0)| ^2 \bigg)^{1/2}\\
&\les&\Big( \Bext_q^s[\psiwc](\tau_1, \tau_2) \Big)^{1/2}  \Bigg( \int_{\Mext(\tau_1, \tau_2)} r^{q+5}   |\dk^{\le s+1 }(N_0)| ^2 \bigg)^{1/2},
 \eeaa
 and, since $N_0= O(a^2 r^{-4} ) \psi$, 
 \beaa
  \int_{\Mext(\tau_1, \tau_2)} r^{q+5}   |\dk^{\le s+1 }N_0| ^2 &\les & O(a^4)  \int_{\Mext(\tau_1, \tau_2)} r^{q-3 }   |\dk^{\le s+1 }\psi| ^2\\
  &\les & O( a^4) B_{\max(q, \de)}^{s+1}[\psi](\tau_1, \tau_2)
 \eeaa
 which yields
 \bea
\widecheck{\NN}_{q}^{s}[\psiwc, N_0](\tau_1, \tau_2) &\les& \Big( \Bext_q^s[\psiwc] (\tau_1, \tau_2)\Big)^{1/2}  \Big(B_{\max(q, \de)}^{s+1}[\psi] (\tau_1, \tau_2) \Big)^{1/2}.
 \eea

 {\bf Step 2b.} Since, according to  \eqref{eq:definition-N-L-psi-again:schematic}, we have
\beaa
N_L &=& O(a)\dk^{\leq1}\nab_3 A+ O(ar^{-1} ) \dk^{\leq1} A,
\eeaa
 we infer
 \beaa
\widecheck{\NN}_{q}^{s}[\psiwc, N_L]&\les&\left|  \int_{\Mext} r^{q+1} \widecheck{\nab}_4 \dk^{\leq s}\psiwc \c \big( \dk^{\leq s+1}(O(a) \dk^{\leq1} \nab_3 A)  \big)\right|\\
&&+\left|  \int_{\Mext} r^{q} \widecheck{\nab}_4 \dk^{\leq s}\psiwc \c \big( \dk^{\leq s+1}(O(a) \dk^{\leq1}A)  \big)\right|\\
&\les& |a|\left(\int_{\Mext} r^{q-1}|\widecheck{\nab}_4 \dk^{\leq s}\psiwc|^2\right)^{\frac{1}{2}}\left(\int_{\Mext} r^{q+1}\Big(r^2|\dk^{\leq s+2}\nab_3A|^2+|\dk^{\leq s+2}A|^2\Big)\right)^{\frac{1}{2}}
\eeaa
which yields, recalling the  definition of the norms $B_q[A]$ in  Definition \ref{Definition:NormsBEF-A},
 \bea
\widecheck{\NN}_{q}^{s}[\psiwc, N_L](\tau_1, \tau_2) &\les& \Big( \Bext_q^s[\psiwc] (\tau_1, \tau_2)\Big)^{1/2}  \Big(B_{\max(q, \de)}^{s+1}[A] (\tau_1, \tau_2) \Big)^{1/2}.
 \eea

\bigskip

{\bf Step 3.} Combining the estimates of Step 1 and Step 2, we have
 \beaa
&&\NN_{\max(q, \de)}^{s+1}[\psi, N_0+N_L](\tau_1, \tau_2)+\widecheck{\NN}_{q}^{s}[\psiwc, N_0+N_L]\\
&\les& \Big( \Bext_q^s[\psiwc] (\tau_1, \tau_2)\Big)^{1/2}  \Big(B_{\max(q, \de)}^{s+1}[\psi, A] (\tau_1, \tau_2) \Big)^{1/2}\\
&&+\BEF^{s+1}_{\max(q, \de)}[\psi, A](\tau_1, \tau_2)+\ep_0^2\tau_1^{-2-3\dec}.
 \eeaa
Since $N=N_0+N_L+N_\err$, this yields, together with \eqref{eq:theorem:GenRW2-q-Chap10.}, 
\beaa
            \BEF^s_q[\psiwc](\tau_1, \tau_2)
      &\les&  \Et_q^s[\psiwc](\tau_1)+\NNt_q^s[\psiwc, N_\err](\tau_1, \tau_2) + \NN_{\max\{q,\de\}}^{s+1}[\psi, N_\err](\tau_1, \tau_2)\\
       &&+\Big( \Bext_q^s[\psiwc] (\tau_1, \tau_2)\Big)^{1/2}  \Big(B_{\max(q, \de)}^{s+1}[\psi, A] (\tau_1, \tau_2) \Big)^{1/2}\\
&&+\BEF^{s+1}_{\max(q, \de)}[\psi, A](\tau_1, \tau_2)+\ep_0^2\tau_1^{-2-3\dec}.
       \eeaa
In view of the control of $\BEF^{s+1}_{\max(q, \de)}[\psi, A](\tau_1, \tau_2)$ provided by Theorem \ref{eqtheorem:gRW1-p-chap11}, we infer
 \beaa
       \begin{split}
       \BEF_q^s[\psiwc](\tau_1, \tau_2)  &\les
       E_q^s[\psiwc](\tau_1)+E_{\max(q, \de)}^{s+1}[\psi, A](\tau_1)+\widecheck{\NN}_{q}^{s}[\psiwc, N_\err](\tau_1, \tau_2)\\
       &+\NN_{\max(q, \de)}^{s+1}[\psi, N_\err](\tau_1, \tau_2)  +\ep_0^2\tau_1^{-2-3\dec},
       \end{split}
       \eeaa
       as stated. This concludes the proof of  Theorem \ref{THEOREM:GRW1-P-WEAK-PSIWC} .
       \end{proof}


\section{Eliminating $N_\err$}


The goal of this section is to eliminate  the  error term $N_\err$ appearing in the RHS of the estimates of Theorems \ref{Thm:Nondegenerate-Morawetz} and  \ref{THEOREM:GRW1-P-WEAK-PSIWC}  and derive the following results.
\begin{theorem}
\lab{Thm:Nondegenerate-Morawetz-strong}
Under the   assumptions made in section \ref{section:Preliminaries:Transport}, 
the following  estimates hold true, for   all   $\de\le p\le  2-\de$ and  $2\leq s\le \kl$,
 \bea
        \BEF_p^s[\psi, A](\tau_1, \tau_2) \les  E_p^s[\psi, A](\tau_1) +\ep_0^2\tau_1^{p-2-3\dec}.
       \eea          
\end{theorem}

\begin{theorem}
\lab{theorem:gRW1-p-strong-psiwc}
Under the   assumptions made in section \ref{section:Preliminaries:Transport}, the following  estimates hold true, 
 for $2\le s\le \kl-1$ and $-1+\de\leq q \leq 3\dec$,
 \bea
       \begin{split}
       \BEF_q^s[\psiwc](\tau_1, \tau_2)  &\les
       E_q^s[\psiwc](\tau_1)+E_{\max(q, \de)}^{s+1}[\psi, A](\tau_1)+\ep_0^2\tau_1^{q -3\dec}.
       \end{split}
       \eea
       \end{theorem}

\begin{proof}[Proof of Theorem \ref{Thm:Nondegenerate-Morawetz-strong}]
Recall, see \eqref{eq:N_err=N_g+},
\beaa
 N_{\err}
 &=&N_g+ \nab_3(rN_g)+ N_m[\qf],\quad
  N_g= r^2 \dk^{\le 2}\big(\Ga_g\c(A,B) \big), \quad N_m[\qf]=\dk^{\le 1}\big(\Ga_g \c \qf\big).
  \eeaa
The proof of the theorem  follows step by step the   proof of Theorem 5.14 in \cite{KS}.  As in that paper, the  terms   $ N_m[\qf]$  and  $N_g$ are estimated directly  using 
  the bounds for $\Ga_g, A, B$,  making sure to absorb to the left  the  corresponding  bulk  contribution in $\qf$. The more difficult term  $\nab_3(rN_g)$  requires an integration by parts  which is  explained  in  detail in  the last  
   part  of    section 5.3.2 of   \cite{KS}.
   \end{proof}

\begin{proof}[Proof of Theorem  \ref{theorem:gRW1-p-strong-psiwc}]
The proof is similar to the proof of Theorem 5.15 in \cite{KS}, and we recall below the main steps of the proof. 

\bigskip

{\bf Step 1.} First observe that the control of $N_\err$ in the proof of Theorem \ref{Thm:Nondegenerate-Morawetz-strong} yields in particular
\beaa
\NN_{\max(q, \de)}^{s+1}[\psi, N_\err](\tau_1, \tau_2) &\les& \ep_0^2\tau_1^{\max(q, \de)-2-3\dec}.
\eeaa
We therefore only need to estimate
       \beaa
       \widecheck{\NN}_{q}^{s}[\psiwc, N_\err](\tau_1, \tau_2)&=& \int_{\Mext} r^{q+2} \widecheck{\nab}_4 \dk^{\leq s}\psiwc \c \left( \nab_4 \dk^{\leq s}N_\err + \frac 3 r \dk^{\leq s} N_\err \right)\\
       &\les&\left| \int_{\Mext} r^{q}( r\widecheck{\nab}_4 \dk^{\leq s}\psiwc)  \big(  \dk^{\leq s+1}N_\err  \big) \right|.
       \eeaa
       
       \bigskip 
       
       {\bf Step 2.} We use fact that 
        \beaa
\dk^{\leq s+1}N_\err=\dk^{\leq s+1}N_g+ e_3(\dk^{\leq s+1}r N_g) +\dk^{\leq s+1}N_m[\qf],
\eeaa
and estimate each contribution separately. 

{\bf Step 2a.} We have 
       \beaa
       \widecheck{\NN}_{q}^{s}[\psiwc, N_g](\tau_1, \tau_2)       &\les& \int_{\Mext} r^{q}| r\widecheck{\nab}_4 \dk^{\leq s}\psiwc|  |  \dk^{\leq s+1}N_g |\\
        &\les& \left( \int_{\Mext} r^{q-3}| r\widecheck{\nab}_4 \dk^{\leq s}\psiwc|^2 \right)^{1/2} \left(\int_{\Mext} r^{q+3}  |  \dk^{\leq s+1}N_g |^2 \right)^{1/2}\\
        &\les& \de_1 B_q^s[\psiwc](\tau_1, \tau_2) + \de_1^{-1} \NN_{q+2}^{s+1}[\psi, N_g](\tau_1, \tau_2),
       \eeaa
       where $\de_1>0$ is chosen sufficiently small so that we can later absorb the term $\de_1 B_q^s[\psiwc](\tau_1, \tau_2)$ on the left hand side of the main estimate.

{\bf Step 2b.} As in the proof of Theorem  \ref{Thm:Nondegenerate-Morawetz-strong}, the integral due to $\dk^{\leq s+1} N_m[\qf]$ can be bounded using the control of $\Ga_g$. Since $\de \leq q+1 \leq 2-\de$, we finally infer 
       \beaa
       \widecheck{\NN}_{q}^{s}[\psiwc, N_m[\qf]](\tau_1, \tau_2)              &\les& \de_1 B_q^s[\psiwc](\tau_1, \tau_2) + \de_1^{-1} \NN_{q+2}^{s+1}[\psi, N_g](\tau_1, \tau_2).
       \eeaa

{\bf Step 2c.} The integral due to $e_3(\dk^{\leq s+1}r N_g)$ is estimated through the integration by parts in $e_3$ and the use of the wave equation for $\psiwc$, as in the proof of Theorem  \ref{Thm:Nondegenerate-Morawetz-strong}. As in the proof of Theorem 5.15 in \cite{KS}, we obtain
       \beaa
       \widecheck{\NN}_{q}^{s}[\psiwc, e_3(r N_g)]](\tau_1, \tau_2)           &\les&\sum_{k\le s}\Big| \int_{\Mext(\tau_1,\tau_2)}  r^{q} e_3e_4 (r \dk^k \psiwc ) \c  \dk^{k+1} (r N_g)\Big|\\
 &&+ \de_1 B_q^s[\psiwc](\tau_1, \tau_2) + \de_1^{-1}  \NN^{s+1}_{q+2}[\psi,  N_g](\tau_1, \tau_2),\\
     &\les& \de_1 B_q^s[\psiwc](\tau_1, \tau_2)+ \de_1^{-1} \NN_{q+2}^{s+2}[\psi, N_g](\tau_1, \tau_2).
       \eeaa
       
{\bf Step 2d.} The estimates in Steps 2a, 2b and 2c  imply, for any $\de_1>0$, 
\beaa
       \widecheck{\NN}_{q}^{s}[\psiwc, N_\err](\tau_1, \tau_2)&\les& \de_1 B_q^s[\psiwc](\tau_1, \tau_2)+ \de_1^{-1} \NN_{q+2}^{s+2}[\psi, N_g](\tau_1, \tau_2).
 \eeaa           
 
{\bf Step 3.}  Together with Theorem \ref{THEOREM:GRW1-P-WEAK-PSIWC}, steps 1 and 2 imply, for all $-1+\de\leq q\leq 1-\de$ and any $\de_1>0$,
\beaa
       \begin{split}
       \BEF_q^s[\psiwc](\tau_1, \tau_2)  &\les
       E_q^s[\psiwc](\tau_1)+E_{\max(q, \de)}^{s+1}[\psi, A](\tau_1) +\ep_0^2\tau_1^{\max(q, \de)-2-3\dec}\\
       &+\de_1 B_q^s[\psiwc](\tau_1, \tau_2)+ \de_1^{-1} \NN_{q+2}^{s+2}[\psi, N_g](\tau_1, \tau_2).
       \end{split}
       \eeaa
Choosing $\de_1>0$ small enough to absorb the fourth term on the RHS from the LHS, we infer, for all $-1+\de\leq q\leq 1-\de$,
       \beaa
       \begin{split}
       \BEF_q^s[\psiwc](\tau_1, \tau_2)  &\les
       E_q^s[\psiwc](\tau_1)+E_{\max(q, \de)}^{s+1}[\psi, A](\tau_1) +\ep_0^2\tau_1^{\max(q, \de)-2-3\dec}\\
       &+\NN_{q+2}^{s+2}[\psi, N_g](\tau_1, \tau_2).
       \end{split}
       \eeaa       
Now, arguing as in Proposition 5.10 of \cite{KS}, we have, for all $-1+\de\leq q\leq 3\dec$,
       \beaa
       \NN_{q+2}^{s+2}[\psi, N_g](\tau_1, \tau_2) &\les& \ep_0^2\tau_1^{q -3\dec}.
       \eeaa
       We infer, for all $-1+\de\leq q\leq 3\dec$,
       \beaa
       \begin{split}
       \BEF_q^s[\psiwc](\tau_1, \tau_2)  &\les
       E_q^s[\psiwc](\tau_1)+E_{\max(q, \de)}^{s+1}[\psi, A](\tau_1) +\ep_0^2\tau_1^{q -3\dec}
       \end{split}
       \eeaa  
as stated. This concludes the proof of Theorem  \ref{theorem:gRW1-p-strong-psiwc}.
\end{proof}


\section{Proof of Theorem M1}
\lab{sec:finallyproofofThmM1}


In this section we  make use of the  results of Theorems 
 \ref{Thm:Nondegenerate-Morawetz-strong} and  \ref{theorem:gRW1-p-strong-psiwc}
 to complete  the proof of Theorem M1 of   \cite{KS:Kerr} which we   restate below.
 \begin{theorem}[Theorem M1 in \cite{KS:Kerr}]
\lab{theoremM1:Chap11}
Assume that the spacetime $\MM$  verifies the assumptions   \eqref{eq:assumptionsonMextforpartII-1}--\eqref{eq:GlobalFrame-HcinGa_g2}  as well as  the assumption \eqref{eq:controlofinitialdataforThM8-intro}  on the  initial data. 
 Then, if $\ep_0>0$ is sufficiently small, there exists $\dee>\dec$ such that we have the following   estimates in  $\MM$, for all  $ s\le \kl -10$,
\bea
\lab{eq:theoremM1-Chap11-1}
\sup_{\MM}\Big(\frac{r^2(2r+\tau)^{1+\dee}}{\log(1+\tau)}+r^3(2r+\tau)^{\frac{1}{2}+\dee}\Big)\Big(|\dk^{\le s} A|+r|\dk^{\le s-1}\nab_3 A|\Big) \les \ep_0.
\eea
Also,  for all  $ s\le \kl -10$,
\bea
\lab{eq:theoremM1-Chap11-2}
  \int_{\Si_*(\ge \tau)}|\nab_3 \dk^{s-1} \psi|^2 &\les& \ep^2_0 \tau^{-2-2\de_{extra}}.
\eea
\end{theorem}

\begin{proof} 
We proceed according to the following steps.

{\bf Step 1.}  Starting with the result of Theorem  \ref{Thm:Nondegenerate-Morawetz-strong}, 
we  run the basic mean value argument  for $\de\leq p\leq 2-\de$ as in Theorem 5.21 of \cite{KS}. We obtain\footnote{With $\qf$ replaced here by the pair $(\psi, A)$. Note that there is a loss one derivative for each application of the mean value theorem. See section  5.4.1 in \cite{KS}, and in particular  Theorem 5.21 in that paper.}  for $\tau_1\le \tau\le \tau_*$ and  $s\le \kl-2$,
\bea
\lab{eq:theoremM1-Chap11-3}
\BEF_p[\psi, A](\tau, \tau_*) &\les& \ep_0^2\tau^{-(2-p-\de)}.
\eea

{\bf Step 2.} According to Theorem \ref{theorem:gRW1-p-strong-psiwc} we have
for $2\le s\le \kl-1$, for all $-1+\de\leq q \leq 3\dec$,
 \beaa
       \begin{split}
       \BEF_q^s[\psiwc](\tau_1, \tau_2)  &\les
       E_q^s[\psiwc](\tau_1)+E_{\max(q, \de)}^{s+1}[\psi, A](\tau_1)+\ep_0^2\tau_1^{q-3\dec}.
       \end{split}
       \eeaa
  In view of Step 1, we have
  \beaa
  E_{\max(q, \de)}^{s+1}[\psi, A](\tau_1) \les \ep_0^2\tau_1^{-(2-\max(q, \de)-\de)}
  \eeaa 
  and hence, for $ s\le \kl-3$ and $-1+\de\leq q \leq 3\dec$,
  \beaa
  \BEF_q^s[\psiwc](\tau_1, \tau_2)  &\les    E_q^s[\psiwc](\tau_1) +\ep_0^2\tau_1^{q-3\dec}.
  \eeaa
  Applying the basic mean value theorem argument in the range $-1+3\dec\leq q \leq 3\dec$ as  in the proof of Theorem 5.22 in \cite{KS} we  deduce, for all $-1+3\dec\leq q \leq 3\dec$,
 \beaa
      \BEF_q^s[\psiwc](\tau_1, \tau_2) &\les& \ep_0^2\tau_1^{-(3\dec -q)}.
    \eeaa

{\bf Step 3.} In view of Step 2, we have in particular, for $ s\le \kl-3$ and $q=-\de$, 
 \beaa
     \sup_{\tau\in[\tau_1, \tau_2]}E^s_{-\de}[\psiwc](\tau_1, \tau_2) &\les& \ep_0^2\tau_1^{-3\dec-\de}.
    \eeaa
From  the relation between $\psi$ and $\widecheck{\psi}$, we deduce, as in (5.4.20) of \cite{KS},  for  $s\le \kl-3$,
 \beaa
       \sup_{\tau\in[\tau_1, \tau_2]} E_{2-\de}^s[\psi](\tau)  &\les& \ep_0^2\tau_1^{-3\dec -\de}.
       \eeaa

{\bf Step 4.} Recall from Proposition  \ref{prop:Step1-psi} that the following  estimate for solutions  $\psi$ of the  full gRW  equation hold true\footnote{Here, we crucially exploit the fact that,  in the estimate  of the Proposition  \ref{prop:Step1-psi},  the term in $A$  on the right hand  side  appears in the norm  $\BEF^s_{\de}[\psi, A]$ (rather than $\BEF^s_p[\psi, A]$).}, for all $s\le k_L$ and all $\de\le p\le 2-\de$,
  \beaa
  \BEF_p^s[\psi](\tau_1, \tau_2) \les
       E_p^s[\psi](\tau_1) +   \BEF^s_{\de}[\psi, A](\tau_1,\tau_2) +\NN_p^s[\psi, N_\err](\tau_1, \tau_2).
  \eeaa
Eliminating $N_\err$  as before  in  Theorems \ref{Thm:Nondegenerate-Morawetz-strong} and 
 Theorem \ref{theorem:gRW1-p-strong-psiwc},  see also  Proposition 5.10 and the proof of Theorem 5.14 in \cite{KS},
 we deduce, for $s\le \kl-1$ and all $\de\le p\le 2-\de$,
 \beaa
  \BEF_p^s[\psi](\tau_1, \tau_2) \les
       E_p^s[\psi](\tau_1) +   \BEF^s_{\de}[\psi, A](\tau_1,\tau_2) +\ep_0^2\tau_1^{p-2-3\dec}.
  \eeaa
  In  view of Step 1  we have
   $\BEF_\de[\psi, A](\tau, \tau_*) \les  \ep_0^2\tau^{-(2-2\de)}$.
Hence, for $s\le \kl-1$ and all $\de\le p\le 2-\de$,
 \beaa
  \BEF_p^s[\psi](\tau_1, \tau_2) \les
       E_p^s[\psi](\tau_1) +  \ep_0^2\tau_1^{-(2-2\de)} +\ep_0^2\tau_1^{-(2+3\dec -p)}.
  \eeaa
In particular, we infer for $s\le \kl-3$,
 \beaa
  \BEF_p^s[\psi](\tau_1, \tau_2) \les
       E_p^s[\psi](\tau_1)  +\ep_0^2\tau_1^{-(2+3\dec -p)}, \qquad 3\dec +2\de \leq p\leq 2-\de.
  \eeaa
Together with Step 3 and the usual mean value argument, we deduce  for $s\le \kl-4$,
 \bea\lab{eq:Flux-psiStep4:beforetolast}
  \BEF_p^s[\psi](\tau_1, \tau_2) \les
        \ep_0^2\tau_1^{-(2+3\dec -p)}, \qquad 1-\de \leq p\leq 2-\de.
  \eea
In particular,  for $s\le \kl-4$, we have 
\bea
\lab{eq:Flux-psiStep4}
  \BEF_{1-\de}^s[\psi](\tau_1, \tau_2) \les
        \ep_0^2\tau_1^{-(3\dec+1+\de)}.
  \eea

{\bf Step 5.} As in Proposition 5.12 in \cite{KS}, we interpolate the control on $E_p[\psi]$ provided \eqref{eq:Flux-psiStep4:beforetolast} between $p=1+\de$ and $p=1-\de$ and obtain 
\beaa
\tau^{1+3\dec}\int_{S_r}|\dk^{\le s} \psi|^2 &\les& \ep_0^2.  
\eeaa
Using Sobolev, we infer the following pointwise decay estimate for $\psi$, for all $s\le \kl-6$,
\bea\lab{eq:pointwisedecayforpsiatrminus1plustaudecauStep5}
|\dk^{\le s }\psi| &\les& \ep_0r^{-1}\tau^{-\frac{1+3\dec}{2}}.
\eea

{\bf Step 6.}   We can now make use of the system of transport equations, see Corollary \ref{cor:systemoftransportequationsforPsiandAfromqf:chap11},
\bea\lab{eq:weuseagainthetransportsystemforAPsiintermsofqfforproofThoeremM1}
\bsplit
\nabc_3\Psi &=\Big(O(r^{-2})+r^{-1}\Ga_g\Big)\qf  +r^2\dk^{\leq 1}\Ga_b\nabc_3A +r\dk^{\leq 1}\Ga_b A, \\ \nabc_3\left(\frac{(\ov{\tr\Xb})^2}{(\Re(\tr\Xb))^2(\tr\Xb)^2}A\right) & = \Psi+r^2\dk^{\leq 1}(\Ga_b)\c A.
\end{split}
\eea
with  $\Psi$ introduced  in  Definition  \ref{def:PsiforintegrationofAfromqf:chap11}.
Integrating the    system of transport equations\footnote{Proceeding as in sections  6.1.3 and  6.1.4  of \cite{KS}.}
 and using   the pointwise decay for $\psi=\Re(\qf)$ in \eqref{eq:pointwisedecayforpsiatrminus1plustaudecauStep5}, 
 we  infer the following pointwise decay estimate for $A$, for all  $s\le \kl-6$,
\bea
|\dk^{\le s}A| &\les& \ep_0r^{-3}\tau^{-\frac{1+3\dec}{2}}.
\eea
Integrating this estimate on $\Si(\tau)$ and recalling the definition  of the $E_p^s[A]$ 
norms in Definition \ref{Definition:NormsBEF-A}
 we  deduce, for all $s\le \kl-6$,
\bea
\lab{eq:Flux-psiStep5}
E^s_{1-\de}[A](\tau) &\les& \ep_0^2\tau^{-(1+3\dec)}.
\eea

{\bf Step 7.} The estimates  \eqref{eq:Flux-psiStep4}  and \eqref{eq:Flux-psiStep5}  imply, for all $s\le \kl-6$,
\beaa
E^s_{1-\de}[A, \psi](\tau) &\les& \ep_0^2\tau^{-(1+3\dec)}.
\eeaa
We can  thus run again  the standard mean value argument and deduce from that estimate and Theorem  \ref{Thm:Nondegenerate-Morawetz-strong},  for $s\le \kl-7$,
\beaa
\BEF^s_{\de}[A, \psi](\tau, \tau_*) &\les& \ep_0^2\tau^{-(2+3\dec-2\de)}.
\eeaa
In particular,  for $s\le \kl-7$,
\bea
\lab{eq:Flux-psiStep7}
\BEF^s_{\de}[\psi](\tau, \tau_*)  &\les& \ep_0^2\tau^{-(2+3\dec- 2\de)}
\eea
and, in view of the definition  of the flux norms $F_\de^s[\psi]$  in section \ref{subsection:basicnormsforpsi}, we deduce 
\beaa
 \int_{\Si_*(\ge \tau)}|\nab_3 \dk^{s-1} \psi|^2 &\les& \ep_0^2\tau^{-(2+3\dec-3 \de)}.
\eeaa
Hence, choosing $\de_{extra}= \frac{3\dec-2\de}{2}>\dec$, for $s\le \kl-7$,
\beaa
 \int_{\Si_*(\ge \tau)}|\nab_3 \dk^{s-1} \psi|^2 &\les& \ep_0^2 \tau^{-2-2\de_{extra} }
 \eeaa
which establishes  the   desired  estimate \eqref{eq:theoremM1-Chap11-2}.

{\bf Step 8.}  Making  use of  the estimate \eqref{eq:Flux-psiStep7} 
 and proceeding as in the  derivation of  the estimate  (5.2.7) in   Proposition 5.12 of \cite{KS},  we derive the estimate, for any $S_r\subset\Si(\tau)$, $s\le \kl-8$,
\bea
\lab{eq:Flux-psiStep8-1}
 r^{-1} \int_{S_r} | \dk^{\le s} \psi|^2 &\les &\ep_0^2  \tau^{-(2+3\dec-2\de)}.
\eea
Similarly, proceeding as in the estimate (5.2.9)  in Proposition 5.13 in \cite{KS}  we  derive
 for any $S_r\subset\Si(\tau)$, $s\le \kl-8$,
\bea
\lab{eq:Flux-psiStep8-2}
\int_{S_r} | \dk^{\le s-1} \nab_3\psi|^2 &\les &\ep_0^2  \tau^{-(2+3\dec-2\de)}.
\eea

{\bf Step 9.} In view of \eqref{eq:pointwisedecayforpsiatrminus1plustaudecauStep5}, 
 \eqref{eq:Flux-psiStep8-1} and \eqref{eq:Flux-psiStep8-2}, and since $\psi=\Re(\qf)$, we  deduce the following estimate\footnote{Note that this corresponds    to  the estimate  for $\qf$ stated in section  3.6.1 of \cite{KS}.}  for $\qf$, for all $s\le \kl-10$,
\bea
\lab{eq:theoremM1-Chap11-psi}
\sup_{\MM}\Big(r\tau^{\frac{1}{2}+\dee}+\tau^{1+\dee}\Big)|\dk^{\leq s}\qf|+\sup_{\MM}r\tau^{1+\dee}|\dk^{\leq s-1}\nab_3\qf| \les \ep_0
\eea
for $\de_{extra}= \frac{3\dec-2\de}{2}>\dec$.

{\bf Step 10.}    Using  \eqref{eq:theoremM1-Chap11-psi}, together with the system of transport equations  in \eqref{eq:weuseagainthetransportsystemforAPsiintermsofqfforproofThoeremM1}, and proceeding as in sections  6.1.3 and  6.1.4  of \cite{KS}, we derive the following pointwise estimate for $A$, for all $s\le \kl-10$,
 \beaa
\sup_{\MM}\Big(\frac{r^2(2r+\tau)^{1+\dee}}{\log(1+\tau)}+r^3(2r+\tau)^{\frac{1}{2}+\dee}\Big)\Big(|\dk^{\le s} A|+r|\dk^{\le s-1}\nab_3 A|\Big) \les \ep_0
\eeaa
as stated in \eqref{eq:theoremM1-Chap11-1}. Together with the proof of \eqref{eq:theoremM1-Chap11-2} in Step 7, this concludes the proof of Theorem M1.
\end{proof}


   \chapter{Decay Estimates  for $\protect \Ab$}
   \lab{chapter-full-RWforqfb}
 
 
 The goal of this chapter is to  provide a complete  proof for   Theorem $M2$. 
 To this end we  proceed as follows:
 \begin{enumerate}
 \item We  derive  combined  $r^p$ weighted estimates  for  the pair $(\qfb, \Ab)$  stated in 
 Theorem \ref{theorem:unconditional-result-final-psib}. This, by far the  most demanding result of the chapter, is proved  in sections \ref{section:controlofgRWforqfb} and \ref{section:transportAb}. 

 \item We then prove Theorem M2 in section
 \ref{sec:finallyproofofThmM2}, relying in particular on these combined  $r^p$ weighted estimates. 
 \end{enumerate}

  
  \section{Preliminaries}
  \lab{sec:assumptionsontheframe:Chapter12}
 

 The spacetime $\MM$ we are dealing with here is  precisely  that 
 described in  section \ref{sec:pfdoisdvhdifuhgiwhdniwbvoubwuyf}. As in Chapter \ref{chapter-full-RWforqf}, we  make stronger assumptions  on $(\Ga_g, \Ga_b)$. We assume in fact for all 
  for all $k\le \kl$, with\footnote{This is consistent with the  value of $k_L$  used in       the bootstrap assumption needed in the proof of Theorem M2 (see  section \ref{sec:bootstrapassumptions:intro}).}  $k_L=k_{small}+120$,
\bea\lab{eq:assumptionsonMextforpartIIchap12}
\bsplit
\Big(r^2\tau^{\frac{1}{2}+\dec}+r\tau^{1+\dec}\Big)|\dk^{\leq k}\Ga_g| & \leq\ep, \\ 
r\tau^{1+\dec}|\dk^{\leq k}\Ga_b| & \leq\ep.
\end{split}
\eea
 We also assume that the curvature components   $A, B$ verify,  for  $ k\leq k_L$,
 \bea\lab{eq:assumptionsonMextforpartIIchap12:moredecayinrAB}
  r^{7/2+\dec} |\dk^{\leq k} (A, B)| &\leq \ep.
 \eea
 We make the additional gauge condition 
\bea
\lab{eq:Xi-Hb-chapter12}
\Xi\in r^{-2}\Ga_g, \qquad \Hbc\in r^{-1}\Ga_g,
\eea
 condition which  played an essential role in deriving  the   gRW equation for $\qfb$  in section   \ref{section:gRW-equation-qfb}.     Recall that  this choice of frame was necessary  to  derive the correct structure of the   nonlinear terms $\err[\square_2 \qfb] $  in Theorem \ref{THEOREM:EQ-QFB}.
 
 \begin{remark}
 The additional  conditions \eqref{eq:Xi-Hb-chapter12} are verified  by the  global frame   constructed in section 3.6  of \cite{KS:Kerr}. These are crucial in
  deriving the  correct structure of the nonlinear terms  $N_{\err}$ of the gRW  equation for $\qfb$. 
 \end{remark}


\subsection{Full Regge Wheeler equation for $\qfb$}


 Recall the definition of  $\qfb$, see Definition \ref{Definition:Define-qfb},
   \beaa
\qfb&=&  \ov{q}q^3\Big(  \nabc_4\nabc_4 \Ab + \underline{C}_1  \nabc_4\Ab + \underline{C}_2  \Ab\Big),
\eeaa
with complex scalars
\bea\lab{eq:defintionofCb1andCb2fordefintionqfb:chap12}
\begin{split}
\und{C}_1&=2\trch - 2\frac {\atrch^2}{ \trch}  -4 i \atrch, \\
\und{C}_2  &= \frac 1 2 \trch^2- 4\atrch^2+\frac 3 2 \frac{\atrch^4}{\trch^2} +  i \left(-2\trch\atrch +4\frac{\atrch^3}{\trch}\right).
\end{split}
\eea

 The real part of $\qfb$, denoted $\underline{\psi}=\Re(\qfb)$,  verifies  the  following real equation, see Theorem \ref{THEOREM:EQ-QFB},
 \bea\label{eq:Gen.RW-pert-qfb}
\squared_2 \underline{\psi} -V_0\underline{\psi}= \frac{4 a\cos\th}{|q|^2}\dual \nab_T  \underline{\psi}+N, \qquad  V_0= \frac{4\De}{ (r^2+a^2) |q|^2}, 
\eea
with the right hand side $N$ given by 
 \bea\lab{eq:NequalthesumeofN0NLNerr}
  N = N_0+N_L+N_{\err}
  \eea
   where
   \begin{itemize}
   \item[-] $N_0$ denotes the zero-th order term in $\psi$, i.e.
 \bea\lab{eq:definition-N-0-psi-again}
 N_0:= \big(V- V_0\big)  \undpsi=O\left(\frac{a}{r^4}\right) \undpsi.
 \eea
 
  \item[-] $N_L$ denotes\footnote{Here $\underline{W}_4$, $\underline{W}_3$, $\underline{W}_a$, $\underline{W}_0$ are complex functions  of $(r, \th)$, all of which vanish for zero angular momentum, having the following fall-off in $r$.} 
  \beaa
 \begin{split}
 N_L&= \Re\left(q\ov{q}^3    \Big( \frac{8a^2 \De}{r^2|q|^4}\nab_\T \Ab_4 +\frac{8a  \De }{r^2|q|^4} \nab_\Z  \Ab_4
+ \underline{W}_4\Ab_4+ \underline{W}_3 \nab_3 \Ab+\underline{W}_a \nab_a \Ab +\underline{W}_0  \Ab\right)
\end{split}
 \eeaa
 where $\underline{W}_4$, $\underline{W}_3$, $\underline{W}_0$ are complex functions  of $(r, \th)$ and $\underline{W}$ is the product of a complex function of $(r, \th)$ with $\dual\Re(\Jk)$, with the following fall-off
\beaa
 q\ov{q}^3 \underline{W}_4, q\ov{q}^3 \underline{W}=O\left(\frac{a^2}{r}\right), \qquad  \qquad  q\ov{q}^3 \underline{W}_3, q\ov{q}^3 \underline{W}_0=O\left(\frac{a^2}{r^2}\right).
\eeaa
Away from the trapping, the following schematic structure will suffice
 \bea\lab{eq:definition-N-L-undpsi-again}
 N_L = O(a r^{-1})\dk^{\le 1}  \nab_4(r \Ab)+O(ar^{-2})\dk^{\leq 1}\Ab.
 \eea

 \item[-]  $N_{\err} =\err[\square_2 \qfb]$ is  the nonlinear  quadratic  error term,  given schematically by  the expression \eqref{eq:MaiThmParq-err-bar} which we recall below
 \bea\lab{eq:MaiThmParq-err-bar-chap12}\lab{eq:N_err-Ab}
N_{\err} =\widetilde{N}_{\err}+ \dk^{\le  3 } (\Ga_g \c \Ga_b),\qquad \widetilde{N}_{\err} = r^2\dk^{\le 2 }\big(\Ga_b\c(A, B)\big),
 \eea
 with $\widetilde{N}_{\err}$ the principal term  in $N_{\err}$ with respect to decay in $r$.
  \end{itemize}


\subsection{Factorizations of $\qfb$}


\begin{lemma}\lab{lemma:factorizationofqfbusingpowersofqandr}
Assume that $\Xi\in r^{-1}\Ga_g$. Then, we have
\beaa
\qfb &=& \ov{q}q^3\left(\nabc_4+2\tr X -\frac{|\tr X|^2}{2\trch}\right)\left(\nabc_4+2\tr X -\frac{3|\tr X|^2}{2\trch}\right)\Ab +r^2\dk^{\leq 1}(\Ga_g\c\Ga_b).
\eeaa 
\end{lemma}

\begin{proof}
We have
\beaa
&& \left(\nabc_4+2\tr X -\frac{|\tr X|^2}{2\trch}\right)\left(\nabc_4+2\tr X -\frac{3|\tr X|^2}{2\trch}\right)\Ab\\
&=& \nab^2_4\Ab+\left(4\tr X-\frac{2|\tr X|^2}{\trch}\right)\nabc_4\Ab +\nabc_4\left(2\tr X -\frac{3|\tr X|^2}{2\trch}\right)\Ab\\
&&+\left(2\tr X -\frac{|\tr X|^2}{2\trch}\right)\left(2\tr X -\frac{3|\tr X|^2}{2\trch}\right)\Ab
\eeaa
Since 
\beaa
4\tr X-\frac{2|\tr X|^2}{\trch} &=& 4\trch -4i\atrch -\frac{2(\trch^2+\atrch^2)}{\trch}\\
&=& 2\trch - 2\frac {\atrch^2}{ \trch}  -4 i \atrch\\
&=& \und{C}_1,
\eeaa
where $\und{C}_1$ is defined by \eqref{eq:defintionofCb1andCb2fordefintionqfb:chap12}, we deduce
\beaa
&& \left(\nabc_4+2\tr X -\frac{|\tr X|^2}{2\trch}\right)\left(\nabc_4+2\tr X -\frac{3|\tr X|^2}{2\trch}\right)\Ab\\
&=& \nab^2_4\Ab+\und{C}_1\nabc_4\Ab +\nabc_4\left(2\tr X -\frac{3|\tr X|^2}{2\trch}\right)\Ab\\
&&+\left(2\tr X -\frac{|\tr X|^2}{2\trch}\right)\left(2\tr X -\frac{3|\tr X|^2}{2\trch}\right)\Ab.
\eeaa

Next, in view of the following consequence of the null structure equations and of the fact that $\Xi\in r^{-1}\Ga_g$, 
\beaa
\nabc_4\tr X +\frac{1}{2}(\tr X)^2 &=& r^{-2}\dk^{\leq 1}\Ga_g,
\eeaa
we have
\beaa
\nabc_4\left(2\tr X -\frac{3}{2}\frac{|\tr X|^2}{\trch}\right) &=& -(\tr X)^2 -\frac{3}{2}\frac{-\frac{1}{2}(\tr X)^2\ov{\tr X}+\tr X\ov{-\frac{1}{2}(\tr X)^2}}{\trch}\\
&&+\frac{3}{2}\frac{|\tr X|^2\Re(-\frac{1}{2}(\tr X)^2)}{\trch^2}+r^{-2}\dk^{\leq 1}\Ga_g\\
&=& -(\tr X)^2 +\frac{3}{2}|\tr X|^2 -\frac{3}{4}\frac{|\tr X|^2(\trch^2-\atrch^2)}{\trch^2}+r^{-2}\dk^{\leq 1}\Ga_g\\
&=& -\frac{1}{4}\trch^2 +\frac{5}{2}\atrch^2 +\frac{3}{4}\frac{\atrch^4}{\trch^2} +2i\trch\atrch    +r^{-2}\dk^{\leq 1}\Ga_g
\eeaa
and hence
\beaa
&& \nabc_4\left(2\tr X -\frac{3}{2}\frac{|\tr X|^2}{\trch}\right)+\left(2\tr X -\frac{1}{2}\frac{|\tr X|^2}{\trch}\right)\left(2\tr X -\frac{3}{2}\frac{|\tr X|^2}{\trch}\right)\\
&=&  -\frac{1}{4}\trch^2 +\frac{5}{2}\atrch^2 +\frac{3}{4}\frac{\atrch^4}{\trch^2} +2i\trch\atrch    +r^{-2}\dk^{\leq 1}\Ga_g\\
&&+\left(\frac{3}{2}\trch -\frac{1}{2}\frac{\atrch^2}{\trch} -2i\atrch \right)\left(\frac{1}{2}\trch -\frac{3}{2}\frac{\atrch^2}{\trch} -2i\atrch\right)\\
&=&  \frac{1}{2}\trch^2 - 4\atrch^2 +\frac{3}{2}\frac{\atrch^4}{\trch^2} + i\left( -2\trch\atrch +4\frac{\atrch^3}{\trch}\right)   +r^{-2}\dk^{\leq 1}\Ga_g\\
&=& \und{C}_2 +r^{-2}\dk^{\leq 1}\Ga_g
\eeaa
where $\und{C}_2$ is defined by \eqref{eq:defintionofCb1andCb2fordefintionqfb:chap12}. We infer
\beaa
&& \left(\nabc_4+2\tr X -\frac{|\tr X|^2}{2\trch}\right)\left(\nabc_4+2\tr X -\frac{3|\tr X|^2}{2\trch}\right)\Ab\\
&=& \nab^2_4\Ab+\und{C}_1\nabc_4\Ab +\und{C}_2\Ab+r^{-2}\dk^{\leq 1}(\Ga_g\c\Ga_b).
\eeaa
In view of the definition of $\qfb$, we deduce 
\beaa
\qfb &=& \ov{q}q^3\left(\nabc_4+2\tr X -\frac{|\tr X|^2}{2\trch}\right)\left(\nabc_4+2\tr X -\frac{3|\tr X|^2}{2\trch}\right)\Ab +r^2\dk^{\leq 1}(\Ga_g\c\Ga_b)
\eeaa
as stated. This concludes the proof of Lemma \ref{lemma:factorizationofqfbusingpowersofqandr}.
\end{proof}

Next, we introduce the tensor $\Psib$. 
\begin{definition}\lab{def:PsibforintegrationofAbfromqfb:chap12}
Let $\Psib\in \sk_2(\CCC)$ given by 
\beaa
\und{\Psi} &:=& \frac{q^4}{r^2}\left(\nabc_4+2\tr X -\frac{3|\tr X|^2}{2\trch}\right)\Ab.
\eeaa
\end{definition}

We have the following corollary of Lemma \ref{lemma:factorizationofqfbusingpowersofqandr}.
\begin{corollary}\lab{cor:systemoftransportequationsforPsibandAbfromqfb:chap12}
Let $\Psib$ as in Definition \ref{def:PsibforintegrationofAbfromqfb:chap12}. Then, $\und{\Psi}\in \dk^{\leq 1}\Ga_b$,  and $(\und{\Psi}, \Ab)$ satisfies the following system of transport equations 
\beaa
\nabc_4(r\Psib)=\frac{q}{r\ov{q}}\qfb  +r\dk^{\leq 1}(\Ga_g\c\Ga_b), \qquad \nabc_4\left( \frac{q^4}{r^3}\Ab\right)= \frac{1}{r}\Psib+r\Ga_g\c\Ga_b.
\eeaa
\end{corollary}

\begin{proof}
We have 
\beaa
\und{\Psi} &=& \frac{q^4}{r^2}\left(\nabc_4+2\tr X -\frac{3|\tr X|^2}{2\trch}\right)\Ab\\
&=&  \frac{q^4}{r^2}\left(\nabc_4+\frac{1}{2}\tr X+ \frac{3(\trch\tr X-|\tr X|^2)}{2\trch}\right)\Ab\\
&=&  O(r^2)\left(\nabc_4+\frac{1}{2}\tr X\right)\Ab +O(r^2)\atrch\Ab\\
&=&  O(r^2)\left(\nabc_4+\frac{1}{2}\tr X\right)\Ab +\Ga_b.
\eeaa
Together with the Bianchi identity for $(\nabc_4\Ab$, we infer 
\beaa
\und{\Psi} &=&  O(r^2)\left(-\frac 1 2 \DDc\hot\Bb  - 2 \Hb\hot \Bb -3P\Xbh\right)\Ab +\Ga_b\\
&=& \dk^{\leq 1}\Ga_b
\eeaa
as stated. 

Next, we compute 
\beaa
\nabc_4(r\Psib) &=& \nabc_4\left(\frac{q^4}{r}\left(\nabc_4+2\tr X -\frac{3|\tr X|^2}{2\trch}\right)\Ab\right)\\
&=& \frac{q^4}{r}\left(\nabc_4+\frac{4e_4(q)}{q}-\frac{e_4(r)}{r}\right)\left(\nabc_4+2\tr X -\frac{3|\tr X|^2}{2\trch}\right)\Ab.
\eeaa
Since 
\beaa
\frac{e_4(q)}{q} = \frac{1}{2}\tr X+\Ga_g,\qquad \frac{e_4(r)}{r} = \frac{1}{2}\frac{|\tr X|^2}{\trch}+\Ga_g,
\eeaa
we have 
\beaa
\nabc_4(r\Psib) &=& \frac{q^4}{r}\left(\nabc_4+2\tr X -\frac{|\tr X|^2}{2\trch}+\Ga_g\right)\left(\nabc_4+2\tr X -\frac{3|\tr X|^2}{2\trch}\right)\Ab\\
&=& \frac{q^4}{r}\left(\nabc_4+2\tr X -\frac{|\tr X|^2}{2\trch}\right)\left(\nabc_4+2\tr X -\frac{3|\tr X|^2}{2\trch}\right)\Ab +r\Ga_g\und{\Psi}\\
&=& \frac{q^4}{r}\left(\nabc_4+2\tr X -\frac{|\tr X|^2}{2\trch}\right)\left(\nabc_4+2\tr X -\frac{3|\tr X|^2}{2\trch}\right)\Ab\\
&&  +r\dk^{\leq 1}(\Ga_g\c\Ga_b)
\eeaa
where we used the fact that $\Psib\in \dk^{\leq 1}\Ga_b$. Now, recall from Lemma \ref{lemma:factorizationofqfbusingpowersofqandr} that we have
\beaa
\qfb &=& \ov{q}q^3\left(\nabc_4+2\tr X -\frac{|\tr X|^2}{2\trch}\right)\left(\nabc_4+2\tr X -\frac{3|\tr X|^2}{2\trch}\right)\Ab +r^2\dk^{\leq 1}(\Ga_g\c\Ga_b).
\eeaa 
We deduce 
\beaa
\nabc_4(r\Psib)=\frac{q}{r\ov{q}}\qfb  +r\dk^{\leq 1}(\Ga_g\c\Ga_b)
\eeaa
as stated. 

Finally, in view of the definition of $\Psib$, we have
\beaa
\nabc_4\left( \frac{q^4}{r^3}\Ab\right) &=& \frac{q^4}{r^3}\left(\nabc_4\Ab+\left(\frac{4e_4(q)}{q} -\frac{3e_3(r)}{r}\right)\Ab\right)\\
&=& \frac{q^4}{r^3}\left(\nabc_4+2\tr X -\frac{3|\tr X|^2}{2\trch}+\Ga_g\right)\Ab\\
&=& \frac{1}{r}\Psib+r\Ga_g\c\Ga_b
\eeaa
as stated. This concludes the proof of Corollary \ref{cor:systemoftransportequationsforPsibandAbfromqfb:chap12}.
\end{proof}

  
\section{Control of the full gRW equation for $\qfb$}
\lab{section:controlofgRWforqfb}



\subsection{Norms for $\Ab$}


 \begin{definition}
  \lab{Definition:NormsBEF-Ab}
We introduce the following norms for $\Ab$ in $\MM=\MM(\tau_1, \tau_2)$
\beaa
B_p[\Ab](\tau_1, \tau_2) &=& 
\int_{\MM(\tau_1, \tau_2)}  r^{p-3} \Big( r^4 |\nab_4 \nab_4(r\Ab) |^2 +  r^4 |\nab \nab_4(r\Ab) |^2
+ r^2 |\nab_3\nab_4(r\Ab)|^2\\
&& + r^2|\nab_4(r\Ab)|^2 + r^2|\nab \Ab|^2+|\nab_3\Ab|^2 +|\Ab|^2\Big)\\
&=&\int_{\MM(\tau_1, \tau_2)}  r^{p-3}\Big( r^2 |\dk^{\le 1}   \nab_4(r\Ab) |^2 +|\dk^{\le 1 }\Ab|^2 \Big),
\\
E_p[\Ab](\tau)&=&\int_{\Si(\tau)}  r^{p-4}  \Big( r^2 |\nab_\Rhat\nab_4(r\Ab)|^2 + r^2|\nab_4(r\Ab)|^2 +\chi_{red}^2|\nab_3\nab_4(r\Ab)|^2\Big)\\
&&+\int_{\Si(\tau)}  r^{p-4}  \Big(  |\nab_3\Ab|^2 +|\Ab|^2\Big),\\
F_p[\Ab](\tau_1,\tau_2)&=&\int_{\AA\cup \Si_*(\tau_1, \tau_2)}  r^{p-2}  \Big( r^2 |\nab_{\Rhat}\nab_4(r\Ab)|^2 + r^2|\nab_4(r\Ab)|^2 +|\nab_3\Ab|^2 +|\Ab|^2\Big)\\
&& +\int_{\AA(\tau_1, \tau_2)}|\nab_3\nab_4(r\Ab)|^2,
\eeaa
where $\chi_{red}=\chi_{red}(r)$ is a smooth function such that $\chi_{red}=1$ for $r\leq r_+(1+2\de_{red})$ and $\chi_{red}=0$ for $r\geq r_+(1+2\de_{red})$,  with the  constant $\de_{red}>0$ small enough such that there holds $\MM_{trap}\subset\{r\geq r_+(1+2\de_{red}\}$.

 \begin{remark}
    \lab{Remark:BEF[A]-norms}
    Note  that the derivatives $\nab_3^2 \Ab$, $\nab\nab_3 \Ab $ and $\nab^2 \Ab$  are missing in $\BEF_p[\Ab](\tau_1, \tau_2)$ but   are fortunately not needed to close the estimates for $\qfb$. Additional derivatives in $E_p[\Ab](\tau)$ and $F_p[\Ab](\tau_1, \tau_2)$ are missing as well and are also  not needed to close the estimates for $\qfb$, with the exception of the ones recovered in \eqref{eq:transportAb:additionalestimateformissingderivativesenergy}. 
     \end{remark}   

 The higher derivative  norms are defined by the usual procedure
    \beaa
    B^s_p[\Ab]=  B_p[\dk^{\le s} \Ab],  \qquad   E^s_p[\Ab]=  E^s_p[\dk^{\le s} \Ab],  \qquad   F^s_p[\Ab]=  F^s_p[\dk^{\le s} \Ab]. 
    \eeaa
        
     For $\undpsi$, we use the norms for solutions of RW type equations introduced in section \ref{subsection:basicnormsforpsi}. We  also  define the  combined $(\Ab, \undpsi$) norms as follows
  \beaa
  E^s_p[\undpsi, \Ab](\tau)&=&  E^s_p[\undpsi](\tau) +E^s_p[\Ab](\tau),\\
  B^s_p[\undpsi, \Ab](\tau_1, \tau_2) &=&  B^s_p[\undpsi](\tau_1, \tau_2) +B^s_p[\Ab](\tau_1, \tau_2),\\
  F^s_p[\undpsi, \Ab](\tau_1, \tau_2) &=& F^s_p[\undpsi](\tau_1, \tau_2)+ F^s_p[\Ab](\tau_1, \tau_2).
  \eeaa
  We  use the short hand notation
  \beaa
  \BEF_p^s[\Ab](\tau_1, \tau_2)  &=& B_p^s[\Ab](\tau_1, \tau_2)  +\sup_{\tau\in[\tau_1,\tau_2]} E_p^s[\Ab](\tau)+
   F_p^s[\Ab](\tau_1, \tau_2), \\
    \BEF_p^s[\undpsi](\tau_1, \tau_2)  &=& B_p^s[\undpsi](\tau_1, \tau_2)  +\sup_{\tau\in[\tau_1,\tau_2]} E_p^s[\undpsi](\tau)+
   F_p^s[\undpsi](\tau_1, \tau_2), \\
  \BEF_p^s[\undpsi, \Ab] &=&\BEF_p^s[\undpsi] +\BEF_p^s[\Ab]. 
  \eeaa
     \end{definition}

  
 \subsection{Statement of the main result of Chapter \ref{chapter-full-RWforqfb}}


First, note that Theorem \ref{THEOREM:GENRW1-P} applies also  to the gRW  model
 problem   \eqref{eq:Gen.RW-pert-qfb}  for $\undpsi$. 
        \begin{theorem}[Basic $r^p$-weighted estimates for $\undpsi$]
       \lab{theorem:Gen.RW-pert-qfb}
       The following estimates hold true  for solutions $\undpsi\in\sk_2$  of   the model gRW  equation \eqref{eq:Gen.RW-pert-qfb}, for all $s\le k_L$, $ \de\le p\le 2-\de$.
        \bea
       \lab{eqtheorem:GenRW1-p-psib'}
    \BEF_p^s[\undpsi](\tau_1, \tau_2) \les
       E_p^s[\undpsi](\tau_1)+\NN_p^s[\undpsi, N](\tau_1, \tau_2).
       \eea
            \end{theorem}
                   
  \begin{proof}     
  The proof of  the estimate \eqref{eqtheorem:GenRW1-p-psib'}  is identical to the  analogous  for $\psi$  in  Theorem \ref{THEOREM:GENRW1-P}. Indeed, the only difference between the two RW model equations is the change of  sign in front of the term $\frac{4 a\cos\th}{|q|^2}\dual \nab_T$ which is never used in the proof of Theorem \ref{THEOREM:GENRW1-P}.
\end{proof}

   The main  result  of this chapter is to  extend  \eqref{eqtheorem:GenRW1-p-psib'}
    to the full gRW  system    as follows.

\begin{theorem}
\lab{Thm:Nondegenerate-Morawetz-psib}\lab{theorem:unconditional-result-final-psib}
The following holds true for   $s\le k_L$,  for all $\de\leq p \leq 2 -\de$,
 \bea
       \lab{eqtheorem:gRW1-p}
        \BEF_p^s[\undpsi, \Ab](\tau_1, \tau_2) \les  E_p^s[\undpsi, \Ab](\tau_1) +\NN_p^s[\undpsi, \Nt_\err](\tau_1, \tau_2)+\ep_0^2\tau_1^{-2-2\dec},
       \eea
       where $\Nt_\err$ is defined in \eqref{eq:N_err-Ab}.      
\end{theorem}


\subsection{Proof of Theorem \ref{Thm:Nondegenerate-Morawetz-psib}}
\lab{section:SketchProofThmNondegenerate-Morawetz-psib}


 As in the case of the analogous result for $\psi$ in Chapter \ref{chapter-full-RWforqf}, the proof is done in steps as follows.
 
 {\bf Step 1.}  Recall that $N=N_0+N_L+N_\err$, see \eqref{eq:NequalthesumeofN0NLNerr}.  Also, recall that $N_\err=\Nt_\err +\dk^{\leq 3}(\Ga_b\c\Ga_g)$, see \eqref{eq:N_err-Ab}.
 We  first eliminate $N-\Nt_\err$ from  the right hand side  of \eqref{eqtheorem:GenRW1-p-psib'}. 
  \begin{proposition}
  \lab{prop:Step1-undpsi}
  The following  estimate for solutions  $\undpsi$ of the  full gRW  equation hold true for all $s\le k_L$ and all $\de\le p\le 2-\de$.
  \bea
  \lab{eq:Step1-undpsi}
\nn  \BEF_p^s[\undpsi](\tau_1, \tau_2) &\les&
       E_p^s[\undpsi](\tau_1) +   O(a)    \BEF^s_p[\undpsi, \Ab](\tau_1,\tau_2) \\
       &&+\NN_p^s[\undpsi, \widetilde{N}_\err](\tau_1, \tau_2)+\ep_0^2 \tau_1^{-2-2\dec}.
  \eea
  \end{proposition}
  
The proof of Proposition \ref{prop:Step1-undpsi} is an immediate consequence of \eqref{eqtheorem:GenRW1-p-psib'} and the following lemma.
\begin{lemma}
\lab{lemma:DecayfortheNterm-undpsi}
 For $\de\leq p\le 2-\de$, $N$ given by \eqref{eq:NequalthesumeofN0NLNerr} satisfies  
 \bea
 \nn\NN^s _p[\undpsi,  N](\tau_1, \tau_2) 
  &\les& |a|\BEF_p^s[\undpsi, \Ab](\tau_1, \tau_2) + \NN_p^s[\undpsi, \widetilde{N}_\err](\tau_1, \tau_2)\\
&& +  \ep_0\tau_1^{-1-\dec}\Big(BEF_p[\undpsi](\tau_1, \tau_2)\Big)^{\frac{1}{2}} +\ep_0^2 \tau_1^{-2-2\dec},
  \eea   
  where $\widetilde{N}_\err $ is defined in \eqref{eq:N_err-Ab}.
\end{lemma}

The proof of Lemma \ref{lemma:DecayfortheNterm-undpsi} is given in section \ref{section:Prooflemma:DecayfortheNterm-undpsi}.

\bigskip

  {\bf Step 2.} We control the term  $ \BEF^s_p[\Ab]$   with the help of the proposition  below. 
   \begin{proposition}
 \lab{prop:MaiTransportAb-steps} 
 The following  estimates hold true, for $ s\le k_L$, for all $\de\leq p \leq 2-\de$, 
 \bea
 \lab{eq:transportAb}
\BEF_p^s[\Ab] (\tau_1, \tau_2)&\les&   B^s_{\de}[\undpsi](\tau_1,\tau_2) 
 + E^s_p[\Ab](\tau_1) +\ep_0^2 \tau_1^{-2-2\dec}.
 \eea
Also, we have the following additional control  on $\Si(\tau)$ with $\tau\in[\tau_1, \tau_2]$, for $ s\le k_L$, for all $\de\leq p \leq 1-\de$, 
\bea\lab{eq:transportAb:additionalestimateformissingderivativesenergy}
\nn&&\int_{\Si(\tau)}r^{p-2}\Big(r^2|\nab_4\nab_4(r\dk^{\leq s}\Ab)|^2+r^2|\nab\nab_4(r\dk^{\leq s}\Ab)|^2+|\nab_3\nab_4(r\dk^{\leq s}\Ab)|^2+|\nab\dk^{\leq s}\Ab|^2\Big)\\
 &\les& EB^s_p[\undpsi](\tau_1,\tau_2)   + E^s_p[\Ab](\tau_1) +\ep_0^2 \tau_1^{-2-2\dec}.
\eea
\end{proposition}

The proof of Proposition \ref{prop:MaiTransportAb-steps} is given in section  \ref{section:transportAb}.

{\bf Step 3.} As a consequence    of Proposition \ref{prop:MaiTransportAb-steps}  and  Proposition \ref{prop:Step1-undpsi}, as well as the smallness of $|a|/m$,  we deduce, for all $s\le k_L$ and and all $\de\le p\le 2-\de$,
\beaa
\BEF_p^s[\undpsi, \Ab] (\tau_1, \tau_2)&\les&  E^s_p[\undpsi, \Ab](\tau_1) +\NN_p^s[\undpsi, \widetilde{N}_\err ](\tau_1, \tau_2)+\ep_0^2 \tau_1^{-2-2\dec}
\eeaa
as stated. This  ends the proof of Theorem  \ref{Thm:Nondegenerate-Morawetz-psib}.


\subsection{Proof of Lemma \ref{lemma:DecayfortheNterm-undpsi}}
\lab{section:Prooflemma:DecayfortheNterm-undpsi}


Since $N=N_0+N_L+N_\err$, it suffices  to  prove, for $\de\leq p\le 2-\de$, the following  estimates
 \bea
 \lab{eq:LemmaDecayfortheNterm-undpsi}
 \bsplit
 \NN_p^s[\undpsi, N_0+N_L](\tau_1, \tau_2) &\les |a|\BEF_p^s[\undpsi, \Ab](\tau_1, \tau_2)+\ep_0\tau_1^{-2-2\dec},\\
  \NN_p^s[\undpsi,  N_\err -\widetilde{N}_\err](\tau_1, \tau_2)  &\les \ep_0\tau_1^{-1-\dec}\Big(BEF^s_p[\undpsi](\tau_1, \tau_2)\Big)^{\frac{1}{2}}.
 \end{split}
 \eea

We start with the control of $N_\err-\widetilde{N}_\err$. Recalling  the definition of the  norms $\NN_p[\undpsi,  N]$, see section \ref{subsection:basicnormsforpsi} we  have, for $\de\leq p\leq 2-\de$, 
 \beaa
    &&   \NN_p[\undpsi,  N](\tau_1, \tau_2) \\
       &=&\int_{\MM(\tau_1,\tau_2)} \big(|\nab_{\Rhat} \undpsi|+r^{-1}|\undpsi|\big) |N|+\left| \int_{\Mext(\tau_1,\tau_2)}  r^{p-1}  \, \nab_4 (r\undpsi ) \c  N\right|
      +\left|\int_{\MM(\tau_1,\tau_2)} \nab_{\That_\de} \undpsi \c N\right| \\
       &\les&
 \int_{\MM_{trap}}(|\dk\undpsi|+r^{-1}|\undpsi|)|N|
 +\int_{\Mntrap}\Big((|\nab_3\undpsi|+r^{-1}|\dk^{\leq 1}\undpsi|)|N|+r^{p-1}|\nab_4(r\undpsi)||N|\Big)\\
 &\les& \left(\sup_{[\tau_1, \tau_2]}E[\undpsi](\tau)\right)^{\frac{1}{2}}\int_{\tau_1}^{\tau_2}\|N\|_{L^2(\Si_{trap}(\tau))}\\
 &+&\left(\int_{\Mntrap(\tau_1, \tau_2)}\Big(r^{-1-\de}|\nab_3\undpsi|^2+r^{p-3}|\dk^{\leq 1}\undpsi|^2\Big)\right)^{\frac{1}{2}}\left(\int_{\Mntrap(\tau_1, \tau_2)}r^{p+1}|N|^2\right)^{\frac{1}{2}}\\
 &\les& \Big(BEF_p[\undpsi](\tau_1, \tau_2)\Big)^{\frac{1}{2}}\left(\int_{\tau_1}^{\tau_2}\|N\|_{L^2(\Si_{trap}(\tau))}+\left(\int_{\MM(\tau_1, \tau_2)}r^{p+1}|N|^2\right)^{\frac{1}{2}}\right).
 \eeaa      
Since  $N_\err -\widetilde{N}_\err=\dk^{\leq 3}(\Ga_b\c\Ga_g)$ we have, according to  our bootstrap assumptions,  both 
\beaa
N_\err - \widetilde{N}_\err  = \ep^2 r^{-3} \tau^{-3/2-2\dec},  \qquad N_\err -\widetilde{N}_\err=\ep^2r^{-2}\tau^{-2-2\dec}. 
\eeaa
Thus,   for $\de\leq p\leq 2-\de$,  
\beaa
&& \int_{\tau_1}^{\tau_2}\left\|N_\err -\widetilde{N}_\err \right\|_{L^2(\Si_{trap}(\tau))}+\left(\int_{\MM(\tau_1, \tau_2)}r^{p+1}\left|N_\err -\widetilde{N}_\err \right|^2\right)^{\frac{1}{2}}\\
&\les& \ep^2\int_{\tau_1}^{+\infty}\frac{d\tau}{\tau^{2+2\dec}} +\ep^2\left(\int_{\MM(\tau_1, \tau_2)}r^{p+1 -6}\tau^{-3-4\dec}\right)^{\frac{1}{2}}\\
&\les& \ep_0\tau_1^{-1-2\dec}.
\eeaa
We deduce, for $\de\leq p\leq 2-\de$, 
 \beaa
 \NN_p[\undpsi,  N_\err -\widetilde{N}_\err](\tau_1, \tau_2)  &\les& \ep_0\tau_1^{-1-2\dec}\Big(BEF_p[\undpsi](\tau_1, \tau_2)\Big)^{\frac{1}{2}}
 \eeaa  
which proves the estimate for $N_\err -\widetilde{N}_\err $ in \eqref{eq:LemmaDecayfortheNterm-undpsi} in the case $s=0$.   
  
  For higher  derivatives, $s\le k_{L} $, we write  schematically
\beaa
\dk^{\le s}(N_\err -\widetilde{N}_\err ) &=& \dk^{\le 3+s  }(\Ga_g\c \Ga_b) = \dk^{\le 3+s  }\Ga_g  \c  \dk^{\le  \frac{k_L}{2}  }\Ga_b+\dk^{\le  \frac{k_L}{2} }  \Ga_g  \c   \dk^{\le 3+s  }\Ga_b
\eeaa
 and   we use an additional $\tau_1^{\de_0} $ by   a standard interpolation argument, see Lemma 5.1 in \cite{KS}.
  Hence, for $\de\leq p\leq 2-\de$,
  \beaa
 \NN^s _p[\undpsi,  N_\err -\widetilde{N}_\err ](\tau_1, \tau_2) 
  &\les&  \ep_0\tau_1^{-1-\dec}\Big(BEF^s_p[\undpsi](\tau_1, \tau_2)\Big)^{\frac{1}{2}}
  \eeaa  
  which concludes the proof of the estimate for $N_\err -\widetilde{N}_\err $ in \eqref{eq:LemmaDecayfortheNterm-undpsi}.

It remains to prove the estimate in \eqref{eq:LemmaDecayfortheNterm-undpsi} for the linear terms $N_0+N_L$. We decompose in two parts
 \beaa
 \NN^s _p[\psi,  N_0+N_L](\tau_1, \tau_2)  &=& \NN^s_{p, r\leq 4m}[\psi,  N_0+N_L](\tau_1, \tau_2)+\NN^s_{p, r\geq 4m}[\psi,  N_0+N_L](\tau_1, \tau_2).
\eeaa
The control of the term $\NN^s_{p, r\leq 4m}[\psi,  N_0+N_L](\tau_1, \tau_2)$ can be done exactly as in section \ref{section:Them:gRW1-p-weak}. In particular, for the treatment of ${}^{(En)}\NN^s_{r\leq 4m}[\psi,  N_L](\tau_1, \tau_2)$, we proceed exactly as in section \ref{subsection:proof-step3}, with $\nabc_3$ replaced by $\nabc_4$, where the boundary terms in the integrations by parts require to control all first order derivatives of $\Ab$ and $\nab_4(r\Ab)$ on $\Si(\tau)\cap\{r\leq 4m\}$ which follows immediately from the control of $E_p[\Ab](\tau)$ and the additional control on $\Si(\tau)$ provided by \eqref{eq:transportAb:additionalestimateformissingderivativesenergy}. We infer
 \beaa
 \NN^s_{p, r\leq 4m}[\psi,  N_0+N_L](\tau_1, \tau_2) &\les& |a|\BEF^s_\de[\undpsi](\tau_1, \tau_2)+\ep_0\tau_1^{-2-2\dec}
 \eeaa
and hence
 \beaa
 \NN^s _p[\psi,  N_0+N_L](\tau_1, \tau_2)  \les \NN^s_{p, r\geq 4m}[\psi,  N_0+N_L](\tau_1, \tau_2)+|a|\BEF^s_\de[\undpsi](\tau_1, \tau_2)+\ep_0\tau_1^{-2-2\dec}.
\eeaa
It thus remains to control $\NN^s_{p, r\geq 4m}[\psi,  N_0+N_L](\tau_1, \tau_2)$. We have, for $\de\leq p\leq 2-\de$,
\beaa
&&\NN^s_{p, r\geq 4m}[\psi,  N_0+N_L](\tau_1, \tau_2)\\ 
&\les& \int_{\MM_{r\geq 4m}(\tau_1, \tau_2)}\Big(|\nab_{\Rhat}\dk^{\leq s}\undpsi|+r^{-1}|\dk^{\leq s}\undpsi|+|\nab_{\That_\de}\dk^{\leq s}\undpsi|+r^{p-1}|\nab_4(r\dk^{\leq s}\undpsi)|\Big)|\dk^{\leq s}(N_0+N_L)|\\
&\les& \int_{\MM_{r\geq 4m}(\tau_1, \tau_2)}\Big(|\nab_3\dk^{\leq s}\undpsi|+r^{p-1}|\dk^{\leq s+1}\undpsi|\Big)|\dk^{\leq s}(N_0+N_L)|\\
&\les& \left(\int_{\MM_{r\geq 4m}(\tau_1, \tau_2)}\Big(r^{-1-\de}|\nab_3\dk^{\leq s}\undpsi|^2+r^{p-3}|\dk^{\leq s+1}\undpsi|^2\right)^{\frac{1}{2}}\\
&&\times\left(\int_{\MM_{r\geq 4m}(\tau_1, \tau_2)}\Big(r^{1+\de}+r^{1+p}\Big)\Big(|\dk^{\leq s}N_0|^2+|\dk^{\leq s}N_L|^2\Big)\right)^{\frac{1}{2}}\\
&\les& \Big(BEF^s_p[\undpsi](\tau_1, \tau_2)\Big)^{\frac{1}{2}}\left(\int_{\MM_{r\geq 4m}(\tau_1, \tau_2)}r^{1+p}\Big(|\dk^{\leq s}N_0|^2+|\dk^{\leq s}N_L|^2\Big)\right)^{\frac{1}{2}}
\eeaa
and hence
\beaa
&& \NN^s _p[\psi,  N_0+N_L](\tau_1, \tau_2)\\ 
&\les& \Big(BEF^s_p[\undpsi](\tau_1, \tau_2)\Big)^{\frac{1}{2}}\left(\int_{\MM_{r\geq 4m}(\tau_1, \tau_2)}r^{1+p}\Big(|\dk^{\leq s}N_0|^2+|\dk^{\leq s}N_L|^2\Big)\right)^{\frac{1}{2}}\\
&&+|a|\BEF^s_\de[\undpsi](\tau_1, \tau_2)+\ep_0\tau_1^{-2-2\dec}.
\eeaa
Since $N_0=O(ar^{-4})\undpsi$, we have, for $\de\leq p\leq 2-\de$, 
\beaa
\int_{\MM_{r\geq 4m}(\tau_1, \tau_2)}r^{1+p}|\dk^{\leq s}N_0|^2 &\les& a^2B_p[\undpsi](\tau_1, \tau_2).
\eeaa
For the control of $N_L$, we  write, see \eqref{eq:definition-N-L-undpsi-again},    
  \beaa
 N_L = O(a r^{-1})\dk^{\le 1}  \nab_4(r \Ab)+O(ar^{-2})\dk^{\leq 1}\Ab,
 \eeaa
 and the definition of  the  $B_p[\Ab]$ norms  (see  Definition 
\ref{Definition:NormsBEF-Ab}),
\beaa
\int_{\MM_{r\geq 4m}(\tau_1, \tau_2)} r^{p+1}   |\dk^{\le s }N_L| ^2&\les&  a^2 \int_{\MM_{r\geq 4m}(\tau_1, \tau_2)} r^{p-1} \big| \dk^{\le 1}  \nab_4( r \Ab)|^2 + a^4 \int_{\MM_{r\geq 4m}(\tau_1, \tau_2)}  r^{p-3}  |\Ab|^2\\
 &\les& a^2 B^s_p[\Ab](\tau_1, \tau_2).
\eeaa
We  deduce, for all $\de\leq p\le 2-\de$,
\beaa
 \NN_p^s[\undpsi, N_0+N_L](\tau_1, \tau_2) &\les |a|\BEF_p^s[\undpsi, \Ab](\tau_1, \tau_2)+\ep_0\tau_1^{-2-2\dec}
 \eeaa
  as stated in \eqref{eq:LemmaDecayfortheNterm-undpsi}.  This concludes the proof of Lemma \ref{lemma:DecayfortheNterm-undpsi}.


\section{Transport Estimates for $\Ab$}
\lab{section:transportAb}


In this section we prove Proposition \ref{prop:MaiTransportAb-steps}.
 i.e. we prove, for $\de\leq p\le 2-\de$, $s\le k_L$,   the following estimate
 \beaa
\BEF_p^s[\Ab] (\tau_1, \tau_2)&\les&   \BEF^s_\de[\undpsi](\tau_1,\tau_2) 
 + E^s_p[\Ab](\tau_1) +\ep_0^2(1+\tau)^{-2-2\dec}.
 \eeaa

To this end, we proceed as follows:
\begin{enumerate}
\item First, we derive a basic lemma for transport equation in $\nab_4$ in section \ref{sec:basictransportlemmaine4:chap12}.

\item Next, we derive estimates for $\Ab$, $\nab_4\Ab$ and $\nab_3\Ab$ in section \ref{sec:controlofAbnab4Abnab3Ab:chap12}.

\item It remains to control angular derivatives. To this end, we fist derive some algebraic identities involving angular derivatives of $\Ab$  in section \ref{sec:usefulalgebraicidentityDDchotDDcAb:chap12}.

\item Next, we derive estimates for $\nab\Ab$ in section \ref{sec:estimatefornabAb:chap12} and for  $\nab\nab_4(r\Ab)$ in section \ref{sec:estimatefornabnab4Ab:chap12}.

\item Finally, we conclude the proof of Proposition \ref{prop:MaiTransportAb-steps} in section \ref{sec:endoftheproofofproposition:MaiTransportAb-steps}.
\end{enumerate}


\subsection{Basic transport lemma in $\nab_4$}
\lab{sec:basictransportlemmaine4:chap12}


We prove below    the $\nabc_4$-transport lemma counterpart of  Lemma
 \ref{lemma:general-transport}  and \ref{eq:lemmageneral-transport-Rhat}. 
 
  \begin{lemma}\lab{lemma:general-transport-estimate-qfb}
 Suppose $\Phi_1, \Phi_2 \in \sk_2(\CCC)$ with signature $s\leq -1$ satisfy the differential relation 
 \bea\label{eq:relation-nab3Phi1Phi2}
 \nabc_4\Phi_1=\Phi_2.
 \eea
Also, let $\chi_{red}=\chi_{red}(r)$ is a smooth function such that $\chi_{red}=1$ for $r\leq r_+(1+2\de_{red})$ and $\chi_{red}=0$ for $r\geq r_+(1+2\de_{red})$,  with the  constant $\de_{red}>0$ small enough such that there holds $\MM_{trap}\subset\{r\geq r_+(1+2\de_{red}\}$.  Then, for every {$p \leq -\de$}:
  \begin{enumerate}
  \item   The pointwise inequality holds true
  \bea
  \lab{eq:lemma-general-transport1-qfb}
 r |q|^{p-4} |\Phi_1|^2&\les &\frac {4}{ p^2} r^{-1}|q|^{p}|\Phi_2|^2   -\frac {2}{ p}\Div(  |q|^{p-2}|\Phi_1|^2 e_4)
\eea
and its integral form
 \bea\label{eq:general-integrated-estimate-e4}
 \bsplit
&\int_{\MM(\tau_1, \tau_2)}  r^{p-3} |\Phi_1|^2
+ \int_{\Si(\tau_2)} r^{p-4}|\Phi_1|^2+  \int_{\AA\cup \Si_*(\tau_1,\tau_2)} r^{p-2}|\Phi_1|^2  \\
  &\les \int_{\MM(\tau_1, \tau_2)} r^{p-1}|\Phi_2|^2+ \int_{\Si(\tau_1)} r^{p-4}|\Phi_1|^2.
  \end{split}
\eea 
 \item We also  have
 \bea
 \lab{eq:lemma-general-transport2-qfb}
\bsplit
&\int_{\MM(\tau_1, \tau_2)}  r^{p-3}\big( r^2| \nab_4 \Phi_1|^2 
+|\nab_3\Phi_1|^2 +|\Phi_1|^2 \big)\\
& + \int_{\Si(\tau_2)}  r^{p-4}\big( |\nab_\Rhat \Phi_1|^2 +|\Phi_1|^2 + \chi_{red}^2|\nab_3\Phi_1|^2 \big)  \\
&+\int_{\AA\cup \Si_*(\tau_1, \tau_2)}  r^{p-2}\big( |\nab_\Rhat \Phi_1|^2 +|\Phi_1|^2 \big) +\int_{\AA(\tau_1, \tau_2)}|\nab_3\Phi_1|^2\\
&\les     \int_{\MM(\tau_1, \tau_2)} r^{p-1}\big( |\nab_\Rhat\Phi_2|^2+  |\Phi_2|^2\big)+\int_{\MM_{r\leq r_+(1+2\de_{red})}(\tau_1, \tau_2)}|\nab_3\Phi_2|^2\\
& + \int_{\Si(\tau_1)}  r^{p-4}\big(|\nab_\Rhat \Phi_1|^2 +|\Phi_1|^2+ \chi_{red}^2|\nab_3\Phi_1|^2 \big)   +   (a^2+\ep^2)\int_{\MM(\tau_1, \tau_2)}r^{p-3} | \nab \Phi_1|^2.
 \end{split}
 \eea
 \end{enumerate}
  \end{lemma}
  
\begin{proof}
We prove  first  the following
  analogue    of Lemma \ref{Lemma:div-e_3}.
\begin{lemma}
\lab{Lemma:div-e_4}
For any  smooth scalar function $f$ on $\MM$, we have
\beaa
\Div \big( f e_4\big)&=& e_4(f)  +\left(\frac{2r\De}{|q|^4} -2\om+\Ga_g\right)f.
\eeaa
\end{lemma}

\begin{proof}
We have
\beaa
\Div ( e_4) &=& \g^{43}\g(\D_4e_4, e_3) +\g^{43}\g(\D_3e_4, e_4)+\g^{bc}\g(\D_be_4, e_c)\\
&=& -\frac{1}{2}4\om +\trch = -2\om +\frac{2\De r}{|q|^4}+\trchc\\
&=& -2\om +\frac{2\De r}{|q|^4}+\Ga_g
\eeaa
and hence 
\beaa
\Div ( fe_4)= f \Div(e_4) + e_4(f) =  e_4(f)+\left(\frac{2r\De}{|q|^4} -2\om +\Ga_g\right)f
\eeaa
as stated. This concludes the proof of Lemma \ref{Lemma:div-e_4}.
\end{proof}

 We continue the proof of \eqref{eq:lemma-general-transport1-qfb} as follows.
 Multiplying the relation $\nabc_4 \Phi_1=\Phi_2$ by $\ov{\Phi_1}$, we deduce, since $\Phi_1$ has  signature $s$, 
\beaa
e_4(|\Phi_1|^2)  &=&2 \Re\big( (\Phi_2 -2s\om\Phi_1)\c \ov{\Phi_1}\big).
\eeaa
Multiplying by $|q|^{p-2}$, and using 
\beaa
e_4(|q|)&=& \frac{e_4(|q|^2)}{2|q|}=\frac{r\De}{|q|^3}+O(1)\widecheck{e_4(r)}+O(r^{-1})e_4(\cos\th)=\frac{r\De}{|q|^3}+\Ga_g,
\eeaa
we deduce
\beaa
2|q|^{p-2} \Re\big( \Phi_2 \c \ov{\Phi_1}\big) &=& e_4( |q|^{p-2}|\Phi_1|^2)  - e_4(|q|^{p-2}) |\Phi_1|^2 +4s\om   |q|^{p-2}     |\Phi_1|^2\\
&=& e_4( |q|^{p-2}|\Phi_1|^2)  - (p-2) r\De |q|^{p-6} |\Phi_1|^2 +4s\om|q|^{p-2} |\Phi_1|^2\\
&&+|q|^{p-3}\Ga_g|\Phi_1|^2.
\eeaa
In view of Lemma \ref{Lemma:div-e_4}, we write
\beaa
\Div(  |q|^{p-2}|\Phi_1|^2 e_4)&=& e_4( |q|^{p-2}|\Phi_1|^2) +\left(\frac{2r\De}{|q|^4} -2\om+\Ga_g\right) |q|^{p-2}|\Phi_1|^2\\
&=&2|q|^{p-2} \Re\big( \Phi_2 \c \ov{\Phi_1}\big) + (p-2) r\De |q|^{p-6} |\Phi_1|^2 -4s\om |q|^{p-2} |\Phi_1|^2\\
&& +\left(\frac{2r\De}{|q|^4} -2\om\right) |q|^{p-2}|\Phi_1|^2 +|q|^{p-2}\Ga_g|\Phi_1|^2\\
&=&2|q|^{p-2} \Re\big( \Phi_2 \c \ov{\Phi_1}\big) + p r\De |q|^{p-6} |\Phi_1|^2 -2(2s+1)\om |q|^{p-2} |\Phi_1|^2\\
&& +|q|^{p-2}\Ga_g|\Phi_1|^2.
\eeaa
From the above identity we deduce
\beaa
&& - p r\De |q|^{p-6} |\Phi_1|^2 +2(2s+1)\om |q|^{p-2} |\Phi_1|^2\\ 
&=&2|q|^{p-2} \Re\big( \Phi_2 \c \ov{\Phi_1}\big) -\Div(  |q|^{p-2}|\Phi_1|^2 e_4) +|q|^{p-2}\Ga_g|\Phi_1|^2\\
&=&2 \Re\big((\la r)^{-1/2}|q|^{p/2} \Phi_2 \c (\lambda r)^{1/2}|q|^{p/2-2}\ov{\Phi_1}\big) -\Div(  |q|^{p-2}|\Phi_1|^2 e_4) +|q|^{p-2}\Ga_g|\Phi_1|^2\\
&\leq & \lambda r |q|^{p-4}|\Phi_1|^2+ \lambda^{-1}r^{-1}|q|^{p}|\Phi_2|^2   -\Div(  |q|^{p-2}|\Phi_1|^2 e_4) +|q|^{p-2}\Ga_g|\Phi_1|^2.
\eeaa
Since  $-\om\gtrsim \frac{m}{r^2}+\Ga_g$   and  $s\leq -1$, we deduce, using also the control of $\Ga_g$, 
\beaa
pr|q|^{p-4}|\Phi_1|^2 &\leq & \lambda r |q|^{p-4}|\Phi_1|^2+ \lambda^{-1}r^{-1}|q|^{p}|\Phi_2|^2 +O(\ep) |q|^{p-4}|\Phi_1|^2  -\Div(  |q|^{p-2}|\Phi_1|^2 e_4).
\eeaa
Therefore, for $p \leq -\de$, choosing   $\la =\frac  p 2 $, we  deduce, for $\ep$ sufficiently small, 
\beaa
 r |q|^{p-4} |\Phi_1|^2&\les &\frac {4}{ p^2} r^{-1}|q|^{p}|\Phi_2|^2   -\frac {2}{ p}\Div(  |q|^{p-2}|\Phi_1|^2 e_4)
\eeaa
which is precisely \eqref{eq:lemma-general-transport1-qfb}. 

The integral form \eqref{eq:general-integrated-estimate-e4} of the inequality then follows by the divergence   theorem, see \eqref{eq:diverggencetheoremX}, and  Remark\footnote{Note in particular that  $\g(N_{\Si}, e_4)= -\frac{m^2}{r^2}$,  which is responsible for   the energy integral  on $\Si(\tau)$.} \ref{remark:Boundariesfortransport}. 
 
Next, we focus on deriving \eqref{eq:lemma-general-transport2-qfb}. To this end, we start with the following commutation lemma.
\begin{lemma}\lab{lemma:commutationlemmanab4nabRhat:chap12}
Let $U\in \sk_k$. Then, we have 
\beaa
\left[\nab_4, \frac{r^2+a^2}{|q|^2}\nab_\Rhat\right]U &=& O((a, \ep)r^{-1})\nab U+O(r^{-1}\ep)\nab_4U+O(r^{-2}\ep)\nab_3U+O(r^{-3})U.
\eeaa
\end{lemma}

\begin{proof}
Recall that $\Rhat$ is given by 
\beaa
\Rhat &=& \frac 1 2 \left( \frac{|q|^2}{r^2+a^2} e_4-\frac{\De}{r^2+a^2}  e_3\right)
\eeaa
so that 
\beaa
\frac{r^2+a^2}{|q|^2}\Rhat &=& \frac 1 2 \left( e_4-\frac{\De}{|q|^2}  e_3\right).
\eeaa
We infer
\beaa
\left[\nab_4, \frac{r^2+a^2}{|q|^2}\nab_\Rhat\right] &=&  -\frac{1}{2}e_4\left(\frac{\De}{|q|^2}\right) \nab_3 -\frac{1}{2}\frac{\De}{|q|^2}[\nab_4, \nab_3]\\
&=& -\frac{1}{2}\left(\pr_r\left(\frac{\De}{|q|^2}\right)e_4(r)+O(r^{-2})e_4(\cos\th)\right) \nab_3 -\frac{1}{2}\frac{\De}{|q|^2}[\nab_4, \nab_3]\\
&=& -\frac{1}{2}\left(\pr_r\left(\frac{\De}{|q|^2}\right)\frac{\De}{|q|^2}+r^{-2}\Ga_g\right) \nab_3 -\frac{1}{2}\frac{\De}{|q|^2}[\nab_4, \nab_3].
\eeaa
Also,  note that the commutation formula for $ [\nab_4, \nab_3]$ of Corollary \ref{cor:comm-gen-B} implies 
\beaa
\, [\nab_4, \nab_3] U &=& \big(O(ar^{-2})+\Ga_b\big)\nab U+ 2\om\nab_3 U  +\Ga_b\nab_4 U +O(r^{-3})U\\
&=&  \big(O(ar^{-2})+\Ga_b\big)\nab U+ 2\left(-\frac{1}{2}\pr_r\left(\frac{\De}{|q|^2}\right)+\omc\right)\nab_3 U  +\Ga_b\nab_4 U +O(r^{-3})U\\
&=&  \big(O(ar^{-2})+\Ga_b\big)\nab U+ \left(-\pr_r\left(\frac{\De}{|q|^2}\right)+\Ga_g\right)\nab_3 U  +\Ga_b\nab_4 U +O(r^{-3})U,
\eeaa
where we used the definition of $\omc$ and the fact that $\omc\in\Ga_g$.  We deduce 
\beaa
\left[\nab_4, \frac{r^2+a^2}{|q|^2}\nab_\Rhat\right]U &=&  \big(O(ar^{-2})+\Ga_b\big)\nab U +\Ga_g\nab_3 U  +\Ga_b\nab_4 U +O(r^{-3})U.
\eeaa
Together with the control of $\Ga_g$ and $\Ga_b$, we deduce 
\beaa
\left[\nab_4, \frac{r^2+a^2}{|q|^2}\nab_\Rhat\right]U &=& O((a, \ep)r^{-1})\nab U+O(r^{-1}\ep)\nab_4U+O(r^{-2}\ep)\nab_3U+O(r^{-3})U
\eeaa
as stated. This concludes the proof of Lemma \ref{lemma:commutationlemmanab4nabRhat:chap12}. 
\end{proof} 

We commute the transport equation for $\Phi_1$ with $\frac{r^2+a^2}{|q|^2}\nab_\Rhat$. Together with Lemma \ref{lemma:commutationlemmanab4nabRhat:chap12}, we infer
\beaa
\nabc_4\left(\frac{r^2+a^2}{|q|^2}\nab_{\Rhat}\Phi_1\right) &=& \left[\nab_4+2s\om, \frac{r^2+a^2}{|q|^2}\nab_\Rhat\right]\Phi_1+\frac{r^2+a^2}{|q|^2}\nab_{\Rhat}\Phi_2\\
&=& O((a, \ep)r^{-1})\nab\Phi_1+O(r^{-1}\ep)\nab_4\Phi_1+O(r^{-2}\ep)\nab_3\Phi_1+O(r^{-3})\Phi_1\\
&& -2s\frac{r^2+a^2}{|q|^2}(\nab_\Rhat\om)\Phi_1+O(1)\nab_{\Rhat}\Phi_2
\eeaa
and hence, since $\nab_\Rhat\om=O(r^{-3})+\dk\Ga_g=O(r^{-2})$, we obtain, using also $\nabc_4\Phi_1=\Phi_2$, 
\beaa
\nabc_4\left(\frac{r^2+a^2}{|q|^2}\nab_{\Rhat}\Phi_1\right) &=& O((a, \ep)r^{-1})\nab\Phi_1+O(r^{-2}\ep)\nab_3\Phi_1+O(r^{-2})\Phi_1\\
&&+O(\ep r^{-1})\Phi_2+O(1)\nab_{\Rhat}\Phi_2.
\eeaa
Applying \eqref{eq:general-integrated-estimate-e4} to this transport equation, we infer
\beaa
 \bsplit
&\int_{\MM(\tau_1, \tau_2)}  r^{p-3} |\nab_{\Rhat}\Phi_1|^2
+ \int_{\Si(\tau_2)} r^{p-4}|\nab_{\Rhat}\Phi_1|^2+  \int_{\AA\cup \Si_*(\tau_1,\tau_2)} r^{p-2}|\nab_{\Rhat}\Phi_1|^2  \\
  &\les \int_{\MM(\tau_1, \tau_2)} r^{p-1}\Big((a^2+\ep^2)r^{-2}|\nab\Phi_1|^2+r^{-4}\ep^2|\nab_3\Phi_1|^2+r^{-4}|\Phi_1|^2+r^{-2}|\Phi_2|^2+|\nab_{\Rhat}\Phi_2|^2\Big)\\
  &+ \int_{\Si(\tau_1)} r^{p-4}|\nab_{\Rhat}\Phi_1|^2.
  \end{split}
\eeaa
Together with \eqref{eq:general-integrated-estimate-e4} and the fact that $\nabc_4\Phi_1=\Phi_2$, this yields 
\beaa
\bsplit
&\int_{\MM(\tau_1, \tau_2)}  r^{p-3}\big( r^2| \nab_4 \Phi_1|^2 
+|\nab_\Rhat \Phi_1|^2 +|\Phi_1|^2 \big) + \int_{\Si(\tau_2)}  r^{p-4}\big( |\nab_\Rhat \Phi_1|^2 +|\Phi_1|^2 \big)  \\
&+\int_{\AA\cup \Si_*(\tau_1, \tau_2)}  r^{p-2}\big( |\nab_\Rhat \Phi_1|^2 +|\Phi_1|^2 \big)\\
&\les     \int_{\MM(\tau_1, \tau_2)} r^{p-1}\big( |\nab_\Rhat\Phi_2|^2+  |\Phi_2|^2\big)  + \int_{\Si(\tau_1)}  r^{p-4}\big( |\nab_\Rhat \Phi_1|^2 +|\Phi_1|^2 \big)  \\
& +   \ep^2\int_{\MM(\tau_1, \tau_2)}r^{p-5} | \nab_3 \Phi_1|^2+   (a^2+\ep^2)\int_{\MM(\tau_1, \tau_2)}r^{p-3} | \nab \Phi_1|^2.
 \end{split}
 \eeaa

It remains to recover $\nab_3\Ab$ in the redshift region. To this end, we commute the transport equation $\nabc_4\Phi_1=\Phi_2$ with $\chi_{red}\nabc_3$ where  $\chi_{red} $ is a smooth  cut-off     function  equal to 1 in  $r\leq r_+(1+\de_{red})$ and    $0$    for $ r\ge  r_+( 1+2\de_{red} )$.  Note that, in view of Lemma \ref{COMMUTATOR-NAB-C-3-DD-C-HOT},
    \beaa
\, [\nabc_3, \nabc_4] U   &=&  O\big((a+\ep)r^{-1})\nabc  U + O(r^{-3})U
 \eeaa
so that 
\beaa
\nabc_4(\chi_{red}\nabc_3\Phi_1) &=& \chi_{red}\nab_3\Phi_2+\pr_r\chi_{red} e_4(r)\nabc_3\Phi_1+O\big((a+\ep)r^{-1})\nabc\Phi_1\\
&& + O(r^{-3})\Phi_1.
\eeaa
Since $s-1<s\leq -1$, applying \eqref{eq:general-integrated-estimate-e4} to this transport equation, and using the above control of $\Phi_1$, we infer
\beaa
\bsplit
&\int_{\MM(\tau_1, \tau_2)}  r^{p-3}\big( r^2| \nab_4 \Phi_1|^2 
+|\nab_\Rhat \Phi_1|^2 +|\Phi_1|^2 \big) +\int_{\MM(r\leq r_+(1+\de_{red}))}|\nab_3\Phi_1|^2\\
& + \int_{\Si(\tau_2)}  r^{p-4}\big( |\nab_\Rhat \Phi_1|^2 +|\Phi_1|^2 + \chi_{red}^2|\nab_3\Phi_1|^2 \big)  \\
&+\int_{\AA\cup \Si_*(\tau_1, \tau_2)}  r^{p-2}\big( |\nab_\Rhat \Phi_1|^2 +|\Phi_1|^2 \big) +\int_{\AA(\tau_1, \tau_2)}|\nab_3\Phi_1|^2\\
&\les    \int_{\MM(\tau_1, \tau_2)} r^{p-1}\big( |\nab_\Rhat\Phi_2|^2+  |\Phi_2|^2\big)+\int_{\MM_{r\leq r_+(1+2\de_{red})}(\tau_1, \tau_2)}|\nab_3\Phi_2|^2\\
& + \int_{\Si(\tau_1)}  r^{p-4}\big(|\nab_\Rhat \Phi_1|^2 +|\Phi_1|^2+ \chi_{red}^2|\nab_3\Phi_1|^2 \big)  \\
& +   \ep^2\int_{\MM(\tau_1, \tau_2)}r^{p-5} | \nab_3 \Phi_1|^2+   (a^2+\ep^2)\int_{\MM(\tau_1, \tau_2)}r^{p-3} | \nab \Phi_1|^2
 \end{split}
 \eeaa
 and hence
\beaa
\bsplit
&\int_{\MM(\tau_1, \tau_2)}  r^{p-3}\big( r^2| \nab_4 \Phi_1|^2 
+|\nab_3\Phi_1|^2 +|\Phi_1|^2 \big)\\
& + \int_{\Si(\tau_2)}  r^{p-4}\big( |\nab_\Rhat \Phi_1|^2 +|\Phi_1|^2 + \chi_{red}^2|\nab_3\Phi_1|^2 \big)  \\
&+\int_{\AA\cup \Si_*(\tau_1, \tau_2)}  r^{p-2}\big( |\nab_\Rhat \Phi_1|^2 +|\Phi_1|^2 \big) +\int_{\AA(\tau_1, \tau_2)}|\nab_3\Phi_1|^2\\
&\les     \int_{\MM(\tau_1, \tau_2)} r^{p-1}\big( |\nab_\Rhat\Phi_2|^2+  |\Phi_2|^2\big)+\int_{\MM_{r\leq r_+(1+2\de_{red})}(\tau_1, \tau_2)}|\nab_3\Phi_2|^2\\
& + \int_{\Si(\tau_1)}  r^{p-4}\big(|\nab_\Rhat \Phi_1|^2 +|\Phi_1|^2+ \chi_{red}^2|\nab_3\Phi_1|^2 \big)  \\
& +   \ep^2\int_{\MM(\tau_1, \tau_2)}r^{p-5} | \nab_3 \Phi_1|^2+   (a^2+\ep^2)\int_{\MM(\tau_1, \tau_2)}r^{p-3} | \nab \Phi_1|^2.
 \end{split}
 \eeaa 
For $\ep>0$ small enough, we deduce 
\beaa
\bsplit
&\int_{\MM(\tau_1, \tau_2)}  r^{p-3}\big( r^2| \nab_4 \Phi_1|^2 
+|\nab_3\Phi_1|^2 +|\Phi_1|^2 \big)\\
& + \int_{\Si(\tau_2)}  r^{p-4}\big( |\nab_\Rhat \Phi_1|^2 +|\Phi_1|^2 + \chi_{red}^2|\nab_3\Phi_1|^2 \big)  \\
&+\int_{\AA\cup \Si_*(\tau_1, \tau_2)}  r^{p-2}\big( |\nab_\Rhat \Phi_1|^2 +|\Phi_1|^2 \big) +\int_{\AA(\tau_1, \tau_2)}|\nab_3\Phi_1|^2\\
&\les     \int_{\MM(\tau_1, \tau_2)} r^{p-1}\big( |\nab_\Rhat\Phi_2|^2+  |\Phi_2|^2\big)+\int_{\MM_{r\leq r_+(1+2\de_{red})}(\tau_1, \tau_2)}|\nab_3\Phi_2|^2\\
& + \int_{\Si(\tau_1)}  r^{p-4}\big(|\nab_\Rhat \Phi_1|^2 +|\Phi_1|^2+ \chi_{red}^2|\nab_3\Phi_1|^2 \big)  +   (a^2+\ep^2)\int_{\MM(\tau_1, \tau_2)}r^{p-3} | \nab \Phi_1|^2
 \end{split}
 \eeaa  
as stated in \eqref{eq:lemma-general-transport2-qfb}. This concludes the proof of Lemma \ref{lemma:general-transport-estimate-qfb}.
\end{proof}


\subsection{Estimates for $\Ab$, $\nab_4\Ab$ and $\nab_3\Ab$}
\lab{sec:controlofAbnab4Abnab3Ab:chap12}


Recall from Corollary \ref{cor:systemoftransportequationsforPsibandAbfromqfb:chap12} that $(\und{\Psi}, \Ab)$ satisfies the following system of transport equations 
\bea\lab{eq:thesystemof2transporteuqaitonsAbPsibundpsiactuallyused}
\nabc_4(r\Psib)=\frac{q}{r\ov{q}}\qfb  +r\dk^{\leq 1}(\Ga_g\c\Ga_b), \qquad \nabc_4\left( \frac{q^4}{r^3}\Ab\right)= \frac{1}{r}\Psib+r\Ga_g\c\Ga_b,
\eea
where $\Psib$ is given in view of Definition \ref{def:PsibforintegrationofAbfromqfb:chap12} by
\beaa
\und{\Psi} &=& \frac{q^4}{r^2}\left(\nabc_4+2\tr X -\frac{3|\tr X|^2}{2\trch}\right)\Ab.
\eeaa

To state the next proposition, we introduce the following partial  norms for $\Ab$ which do not provide control for angular derivatives.
 \begin{definition}
 \lab{definition:dotnormsforAb}
 We define, for all $p$, 
 \begin{enumerate}
  \item  In $\MM(\tau_1, \tau_2)$
 \beaa
 \Bdot_p[\Ab](\tau_1, \tau_2) &=&\int_{\MM(\tau_1, \tau_2)} r^{p-1} \Big(   r^2| \nab_4 \nab_4(r\Ab) |^2  +  |\nab_3 \nab_4(r\Ab) |^2 + | \nab_4(r\Ab) |^2 \Big)\\
&&
 \int_{\MM(\tau_1, \tau_2)} r^{p-3 }\Big(  r^2| \nab_4 \Ab|^2+  |\nab_3 \Ab|^2+ |\Ab|^2    \Big).
 \eeaa
 
 \item On $\Si(\tau)$, $\Edot_p[\Ab](\tau)=E_p[\Ab](\tau)$, 
 
 \item On $\AA\cup \Si_*(\tau_1,\tau_2)$, $\Fdot_p[\Ab](\tau_1, \tau_2)=F_p[\Ab](\tau_1, \tau_2)$.
 \end{enumerate}
 \end{definition}

\begin{proposition}
\lab{Prop:EstimatesforBEFdotnormsAb}
The following estimates hold true for $p\le 2-\de$
 \beaa
\BEFdot_p[\Ab](\tau_1, \tau_2) &\les&  B_{\de}[\undpsi](\tau_1, \tau_2)+\Edot_{p}[\Ab] (\tau_1) \\
&&  +   (a^2+\ep^2)\int_{\MM(\tau_1, \tau_2)}r^{p-3}\Big(r^2| \nab\nab_4(r\Ab)|^2+| \nab \Ab|^2\Big)+\ep_0^2\tau^{-2-2\dec}.
\eeaa
\end{proposition}

\begin{proof}
  For $p\le -\de$  we apply  \eqref{eq:lemma-general-transport2-qfb} of Lemma \ref{lemma:general-transport-estimate-qfb} with $\Phi_1=r\Psib$ and $\Phi_2= r^{-1}\frac{q}{\ov{q} }\qfb +r\dk^{\leq 1}(\Ga_g\c\Ga_b)$ and  we infer 
\beaa
\bsplit
&\int_{\MM(\tau_1, \tau_2)}  r^{p-3}\big( r^4| \nab_4 \Psib|^2 
+ r^2 |\nab_3 \Psib|^2 + r^2 |\Psib|^2 \big) \\
&+ \int_{\Si(\tau_2)}  r^{p-4}\big( r^2 |\nab_\Rhat\Psib|^2 + r^2 |\Psib|^2 + \chi_{red}^2|\nab_3\Psib|^2\big)  \\
&+\int_{\AA\cup \Si_*(\tau_1, \tau_2)}  r^{p-2}\big(r^2 |\nab_\Rhat \Psib|^2 + r^2 |\Psib|^2\big) +\int_{\AA(\tau_1, \tau_2)}|\nab_3\Psib|^2 \\
&\les         \int_{\MM(\tau_1, \tau_2)} r^{p-1}\big( r^{-2} |\nab_\Rhat\qfb|^2+  r^{-2} |\qfb|^2 +r^2|\dk^{\leq 2}(\Ga_g\c\Ga_b)|^2\big) \\
&+\int_{\MM_{r\leq r_+(1+2\de_{red})}(\tau_1, \tau_2)}   \Big(|\nab_3 \qfb    |^2+|\dk^{\leq 2}(\Ga_g\c\Ga_b)|^2\Big) \\
&+ \int_{\Si(\tau_1)}  r^{p-4}\big( r^2 |\nab_\Rhat\Psib|^2 + r^2 |\Psib|^2 + \chi_{red}^2|\nab_3\Psib|^2 \big)  +   (a^2+\ep^2)\int_{\MM(\tau_1, \tau_2)}r^{p-3}  r^2 | \nab \Psib|^2.
 \end{split}
 \eeaa
   i.e., since $\undpsi=\Re(\qfb)$, and in view of the control of $\Ga_g$ and $\Ga_b$, using also  $p\le -\de$, 
 \bea
 \lab{eq: TransportEstmates-Psib}
\bsplit
&\int_{\MM(\tau_1, \tau_2)}  r^{p-1}\big( r^2| \nab_4 \Psib|^2 
+  |\nab_3 \Psib|^2 +  |\Psib|^2 \big) \\
& + \int_{\Si(\tau_2)}  r^{p-2}\big( |\nab_\Rhat\Psib|^2 +  |\Psib|^2 + \chi_{red}^2|\nab_3\Psib|^2 \big)  \\
&+\int_{\AA\cup \Si_*(\tau_1, \tau_2)}  r^p\big(  |\nab_\Rhat\Psib|^2 +  |\Psib|^2 \big) +\int_{\AA(\tau_1, \tau_2)}|\nab_3\Psib|^2\\
&\les B_{\de}[\undpsi](\tau_1, \tau_2) + \int_{\Si(\tau_1)}  r^{p-2}\big( |\nab_\Rhat\Psib|^2 +  |\Psib|^2 + \chi_{red}^2|\nab_3\Psib|^2 \big) \\
& + (a^2+\ep^2)\int_{\MM(\tau_1, \tau_2)}r^{p-1}| \nab \Psib|^2+\ep_0^2\tau^{-2-2\dec}.
 \end{split}
 \eea

 For $p\leq -\de$, we apply  estimate \eqref{eq:lemma-general-transport2-qfb} with $\Phi_1=\frac{q^4}{r^3}\Ab$ and $\Phi_2=r^{-1}\und{\Psi} +r\Ga_g\c\Ga_b$ and  deduce
\beaa
&&\int_{\MM(\tau_1, \tau_2)}  r^{p-3}\big( r^4| \nab_4 \Ab|^2+ r^2 |\nab_3 \Ab|^2+ r^2|\Ab|^2 \big) \\
&& + \int_{\Si(\tau_2)}  r^{p-4}\big( r^2 |\nab_\Rhat\Ab|^2 + r^2 |\Ab|^2 + \chi_{red}^2|\nab_3\Ab|^2 \big)\\
&&+\int_{\AA\cup \Si_*(\tau_1, \tau_2)}  r^{p-2}\big( r^2 |\nab_\Rhat \Ab|^2 + r^2 |\Ab|^2 \big) +\int_{\AA(\tau_1, \tau_2)}|\nab_3\Ab|^2\\
&\les&      \int_{\MM(\tau_1, \tau_2)} r^{p-1}\big( r^{-2} |\nab_\Rhat\Psib|^2+  r^{-2} |\Psib|^2 +r^2|\dk^{\leq 2}(\Ga_g\c\Ga_b)|^2\big) \\
&&+\int_{\MM_{r\leq r_+(1+2\de_{red})}(\tau_1, \tau_2)}   \Big(|\nab_3 \Psib    |^2+|\dk^{\leq 2}(\Ga_g\c\Ga_b)|^2\Big) \\
&&+ \int_{\Si(\tau_1)}  r^{p-4}\big( r^2 |\nab_\Rhat\Ab|^2 + r^2 |\Ab|^2 + \chi_{red}^2|\nab_3\Ab|^2 \big)  +   (a^2+\ep^2)\int_{\MM(\tau_1, \tau_2)}r^{p-3}  r^2 | \nab \Ab|^2
\eeaa
i.e.,  in view of the control of $\Ga_g$ and $\Ga_b$, using also  $p\le -\de$, 
\beaa
&&\int_{\MM(\tau_1, \tau_2)}  r^{p-1}\big( r^2| \nab_4 \Ab|^2+  |\nab_3 \Ab|^2+ |\Ab|^2 \big) \\
&& + \int_{\Si(\tau_2)}  r^{p-2}\big(  |\nab_\Rhat\Ab|^2 +  |\Ab|^2 + \chi_{red}^2|\nab_3\Ab|^2 \big)\\
&&+\int_{\AA\cup \Si_*(\tau_1, \tau_2)}  r^{p}\big(  |\nab_\Rhat \Ab|^2 +  |\Ab|^2 \big) +\int_{\AA(\tau_1, \tau_2)}|\nab_3\Ab|^2\\
&\les&      \int_{\MM(\tau_1, \tau_2)} r^{p-3}\big(  |\nab_\Rhat\Psib|^2+   |\Psib|^2\big) +\int_{\MM_{r\leq r_+(1+2\de_{red})}(\tau_1, \tau_2)}|\nab_3 \Psib    |^2 \\
&&+ \int_{\Si(\tau_1)}  r^{p-2}\big(  |\nab_\Rhat\Ab|^2 +  |\Ab|^2 + \chi_{red}^2|\nab_3\Ab|^2 \big)  +   (a^2+\ep^2)\int_{\MM(\tau_1, \tau_2)}r^{p-1}| \nab \Ab|^2+\ep_0^2\tau^{-2-2\dec}.
\eeaa
Combining  with \eqref{eq: TransportEstmates-Psib}   we  deduce, for all $p\le -\de$,
   \beaa
  && \int_{\MM(\tau_1, \tau_2)} r^{p-1} \Big(  r^2| \nab_4 \Psib|^2 
+  |\nab_3 \Psib|^2 +  |\Psib|^2  +  r^2| \nab_4 \Ab|^2+  |\nab_3 \Ab|^2+ |\Ab|^2    \Big)\\
&&+\int_{\Si(\tau_2) } r^{p-2} \Big(  |\nab_\Rhat\Psib|^2 +  |\Psib|^2 + \chi_{red}^2|\nab_3\Psib|^2 + |\nab_\Rhat\Ab|^2 +  |\Ab|^2 + \chi_{red}^2|\nab_3\Ab|^2    \Big)\\
&&+\int_{\AA\cup\Si_*(\tau_1, \tau_2)} r^p  \Big(   |\nab_\Rhat \Psib|^2 +  |\Psib|^2 +  |\nab_\Rhat \Ab|^2 +  |\Ab|^2    \Big) +\int_{\AA(\tau_1, \tau_2)}\Big(|\nab_3\Psib|^2+|\nab_3\Ab|^2\Big)\\
&\les& B_{\de}[\undpsi]+\int_{\Si(\tau_1) } r^{p -2} \Big(   |\nab_\Rhat\Psib|^2 +  |\Psib|^2 + \chi_{red}^2|\nab_3\Psib|^2 +  |\nab_\Rhat\Ab|^2 +  |\Ab|^2 + \chi_{red}^2|\nab_3\Ab|^2  \Big)\\
&&+   (a^2+\ep^2)\int_{\MM(\tau_1, \tau_2)}r^{p-1}\Big(| \nab \Psib|^2+| \nab \Ab|^2\Big)+\ep_0^2\tau^{-2-2\dec}.
\eeaa
Replacing $p$ with $p-2$ we  rewrite, for all $p\le 2-\de$,
   \beaa
  && \int_{\MM(\tau_1, \tau_2)} r^{p-3} \Big(  r^2| \nab_4 \Psib|^2 
+  |\nab_3 \Psib|^2 +  |\Psib|^2  +  r^2| \nab_4 \Ab|^2+  |\nab_3 \Ab|^2+ |\Ab|^2    \Big)\\
&&+\int_{\Si(\tau_2) } r^{p-4} \Big(  |\nab_\Rhat\Psib|^2 +  |\Psib|^2 + \chi_{red}^2|\nab_3\Psib|^2 + |\nab_\Rhat\Ab|^2 +  |\Ab|^2 + \chi_{red}^2|\nab_3\Ab|^2   \Big)\\
&&+\int_{\AA\cup\Si_*(\tau_1, \tau_2)} r^{p-2}  \Big(   |\nab_\Rhat \Psib|^2 +  |\Psib|^2 +  |\nab_\Rhat \Ab|^2 +  |\Ab|^2     \Big) +\int_{\AA(\tau_1, \tau_2)}\Big(|\nab_3\Psib|^2+|\nab_3\Ab|^2\Big)\\
&\les& B_{\de}[\undpsi]+\int_{\Si(\tau_1) } r^{p -4} \Big(   |\nab_\Rhat\Psib|^2 +  |\Psib|^2 + \chi_{red}^2|\nab_3\Psib|^2 + |\nab_\Rhat\Ab|^2 +  |\Ab|^2 + \chi_{red}^2|\nab_3\Ab|^2  \Big)\\
&&+   (a^2+\ep^2)\int_{\MM(\tau_1, \tau_2)}r^{p-3}\Big(| \nab \Psib|^2+| \nab \Ab|^2\Big)+\ep_0^2\tau^{-2-2\dec}.
\eeaa

We now note that
\beaa
\Psib&=& \frac{q^4}{r^2}\left(\nabc_4+2\tr X -\frac{3|\tr X|^2}{2\trch}\right)\Ab\\
&=& r^2(1+O(r^{-1}))\left(\nab_4-4\om +\frac{4}{q}-\frac{3}{r}+\Ga_g\right)\Ab\\
&=& r^2(1+O(r^{-1}))\left(\nab_4+\frac{1}{r}+O(r^{-2})+\Ga_g\right)\Ab\\
&=& r(1+O(r^{-1}))\left(\nab_4(r\Ab) +\left(-\widecheck{e_4(r)}+O(r^{-1})+r\Ga_g\right)\Ab\right)\\
&=& r(1+O(r^{-1}))\left(\nab_4(r\Ab) +\left(O(r^{-1})+r\Ga_g\right)\Ab\right)\\
&=& (1+O(r^{-1}))\big(r\nab_4(r\Ab) +O(1)\Ab\big),
\eeaa
where we have used the control of $\Ga_g$. We can thus replace $\Psib$ with $r\nab_4(r\Ab)$ in  formula above. 
Together with the fact that 
\beaa
r^2|\nab_4\Ab|^2+|\nab_3\Ab|^2+|\Ab|^2\les |\nab_4(r\Ab)|^2+|\nab_\Rhat\Ab|^2+|\Ab|^2+\chi_{red}^2|\nab_3\Ab|^2,
\eeaa
we deduce,  for $p\le 2-\de$,
  \beaa
\BEFdot_p[\Ab](\tau_1, \tau_2) &\les&  B_{\de}[\undpsi](\tau_1, \tau_2)+\Edot_{p}[\Ab] (\tau_1) \\
&&  +   (a^2+\ep^2)\int_{\MM(\tau_1, \tau_2)}r^{p-3}\Big(r^2| \nab\nab_4(r\Ab)|^2+| \nab \Ab|^2\Big)+\ep_0^2\tau^{-2-2\dec}
\eeaa
which ends the proof of Proposition \ref{Prop:EstimatesforBEFdotnormsAb}.
\end{proof}

 
\subsection{Identities for angular derivatives of $\Ab$}
\lab{sec:usefulalgebraicidentityDDchotDDcAb:chap12}


In addition to the estimates of Proposition \ref{Prop:EstimatesforBEFdotnormsAb}, we need to control angular derivatives of $\Ab$. To this end, we derive in this section several identities.  
We start with the following identity for $\DDc\hot(\DDbc \c\Ab)$. 
\begin{lemma}
\lab{Lemma:TeukoskyAb-sharp}
$\DDc\hot(\DDbc \c\Ab)$ satisfies the following identity 
\bea
\lab{eq:TeukoskyAb-sharp}
\bsplit
 \frac{1}{4}\DDc\hot(\DDbc \c\Ab) =& \nabc_3\left(\nabc_4\Ab +\frac{1}{2}\tr X\Ab\right)  +\left( \frac 1 2 \tr \Xb+ 2\ov{\tr \Xb} \right)\left(\nabc_4\Ab +\frac{1}{2}\tr X\Ab\right)\\\
  &   +(H+\ov{H})\c\nab\Ab -\frac{1}{2}(\DDc\c\ov{H})\Ab+ O(ar^{-2})  \nab\Ab +O(r^{-3} ) \Ab\\
  & +r^{-1}\dk^{\leq 1}(\Ga_b\c\Ga_g).
\end{split}
\eea
\end{lemma}

\begin{proof}
Recall that $\Ab$ verifies the Teukolsky equation
\beaa
\bsplit
&-\nabc_3\nabc_4\Ab  + \frac{1}{4}\DDc\hot(\DDbc \c\Ab) -\frac{1}{2}\tr X\nabc_3\Ab
  -\frac{1}{2}\Big(\tr\Xb+4\ov{\tr\Xb}\Big)\nabc_4\Ab \\
   &+\Big(4\Hb+H+\ov{H}\Big)\c\nab\Ab
+\Big(-\tr X\ov{\tr\Xb} +2P\Big)\Ab  + 2\Hb\hot\Big(\ov{H}\c\Ab\Big) = r^{-2}\dk^{\leq 1}(\Ga_b\c\Ga_b)
\end{split}
\eeaa
and hence
\beaa
\bsplit
 \frac{1}{4}\DDc\hot(\DDbc \c\Ab) =& \nabc_3\nabc_4\Ab +\frac{1}{2}\tr X\nabc_3\Ab
  +\frac{1}{2}\Big(\tr\Xb+4\ov{\tr\Xb}\Big)\nabc_4\Ab \\
   &-\Big(4\Hb+H+\ov{H}\Big)\c\nab\Ab
-\Big(-\tr X\ov{\tr\Xb} +2P\Big)\Ab  - 2\Hb\hot\Big(\ov{H}\c\Ab\Big) \\
&+ r^{-2}\dk^{\leq 1}(\Ga_b\c\Ga_b)\\
=& \nabc_3\nabc_4\Ab +\frac{1}{2}\tr X\nabc_3\Ab
  +\frac{1}{2}\Big(\tr\Xb+4\ov{\tr\Xb}\Big)\nabc_4\Ab +\ov{\tr\Xb}\tr X\Ab\\
  &   +(H+\ov{H})\c\nab\Ab+ O(ar^{-2})  \nab\Ab +O(r^{-3} ) \Ab+r^{-1}\dk^{\leq 1}(\Ga_b\c\Ga_g).
\end{split}
\eeaa

Next, using the null structure equation for $\nabc_3\tr X$, i.e. 
\beaa
\nabc_3\tr X +\frac{1}{2}\tr\Xb\tr X &=& \DDc\c\ov{H}+H\c\ov{H}+2P+\Xib\c\ov{\Xi}-\frac{1}{2}\Xbh\c\ov{\Xh}\\
&=&  \DDc\c\ov{H}+O(r^{-3})+r^{-1}\Ga_g+\Ga_g\c\Ga_b,
\eeaa
we compute 
\beaa
&&\nabc_3\left(\nabc_4\Ab +\frac{1}{2}\tr X\Ab\right)  +\left( \frac 1 2 \tr \Xb+ 2\ov{\tr \Xb} \right)\left(\nabc_4\Ab +\frac{1}{2}\tr X\Ab\right)\\ 
&=& \nabc_3\nabc_4\Ab +\frac{1}{2}\tr X\nabc_3\Ab+\frac{1}{2}\nabc_3(\tr X)\Ab +\left( \frac 1 2 \tr \Xb+ 2\ov{\tr \Xb} \right)\nabc_4\Ab\\
&&+\left( \frac 1 2 \tr \Xb+ 2\ov{\tr \Xb} \right)\frac{1}{2}\tr X\Ab\\
&=& \nabc_3\nabc_4\Ab +\frac{1}{2}\tr X\nabc_3\Ab-\frac{1}{4}\tr\Xb\tr X\Ab +\left( \frac 1 2 \tr \Xb+ 2\ov{\tr \Xb} \right)\nabc_4\Ab\\
&&+\left( \frac 1 2 \tr \Xb+ 2\ov{\tr \Xb} \right)\frac{1}{2}\tr X\Ab +\frac{1}{2}(\DDc\c\ov{H})\Ab+O(r^{-3} ) \Ab+r^{-1}\Ga_b\c\Ga_g\\
&=& \nabc_3\nabc_4\Ab +\frac{1}{2}\tr X\nabc_3\Ab +\left( \frac 1 2 \tr \Xb+ 2\ov{\tr \Xb} \right)\nabc_4\Ab +\ov{\tr \Xb}\tr X\Ab\\
&& +\frac{1}{2}(\DDc\c\ov{H})\Ab+O(r^{-3} ) \Ab+r^{-1}\Ga_b\c\Ga_g.
\eeaa
We deduce 
\beaa
\bsplit
 \frac{1}{4}\DDc\hot(\DDbc \c\Ab) =& \nabc_3\left(\nabc_4\Ab +\frac{1}{2}\tr X\Ab\right)  +\left( \frac 1 2 \tr \Xb+ 2\ov{\tr \Xb} \right)\left(\nabc_4\Ab +\frac{1}{2}\tr X\Ab\right)\\\
  &   +(H+\ov{H})\c\nab\Ab -\frac{1}{2}(\DDc\c\ov{H})\Ab+ O(ar^{-2})  \nab\Ab +O(r^{-3} ) \Ab\\
  & +r^{-1}\dk^{\leq 1}(\Ga_b\c\Ga_g)
\end{split}
\eeaa
and hence 
\beaa
\bsplit
 \frac{1}{4}\DDc\hot(\DDbc \c\Ab) =& \nabc_3\left(\nabc_4\Ab +\frac{1}{2}\tr X\Ab\right)  +\left( \frac 1 2 \tr \Xb+ 2\ov{\tr \Xb} \right)\left(\nabc_4\Ab +\frac{1}{2}\tr X\Ab\right)\\\
  &   +(H+\ov{H})\c\nab\Ab -\frac{1}{2}(\DDc\c\ov{H})\Ab+ O(ar^{-2})  \nab\Ab +O(r^{-3} ) \Ab\\
  & +r^{-1}\dk^{\leq 1}(\Ga_b\c\Ga_g)
\end{split}
\eeaa
as stated. This  concludes the proof of Lemma \ref{Lemma:TeukoskyAb-sharp}.
\end{proof}

We will need to differentiate the identity of Lemma \ref{Lemma:TeukoskyAb-sharp}. To this end, we first derive the following identity. 
\begin{lemma}\lab{lemma:almostfactorizationfoqfbneededforcontralangularderivativesPsib:chap12}
We have
\beaa
 &&\ov{q}q^3\left(\nabc_4 +\frac{3}{2}\tr X - 2\frac {\atrch^2}{ \trch}  -2i \atrch\right)\left(\nabc_4\Ab +\frac{1}{2}\tr X\Ab\right)\\ 
 &=& \qfb + O(a^2)\Ab+r^2\dk^{\leq 1}(\Ga_g\c\Ga_b). 
\eeaa
\end{lemma}

\begin{proof}
Let $h$ a scalar function to be chosen below. We have
\beaa
&&\left(\nabc_4 +\frac{3}{2}\tr X+h\right)\left(\nabc_4\Ab +\frac{1}{2}\tr X\Ab\right)\\
&=& \nabc_4^2\Ab+(2\tr X+h)\nabc_4\Ab+\left(\frac{3}{4}(\tr X)^2+\frac{1}{2}\tr Xh+\frac{1}{2}\nabc_4(\tr X)\right)\Ab\\
&=& \nabc_4^2\Ab+(2\tr X+h)\nabc_4\Ab+\left(\frac{3}{4}(\tr X)^2+\frac{1}{2}\tr Xh+\frac{1}{2}\left(-\frac{1}{2}(\tr X)^2+r^{-2}\dk^{\leq 1}\Ga_g\right)\right)\Ab \\
&=& \nabc_4^2\Ab+(2\tr X+h)\nabc_4\Ab+\left(\frac{1}{2}(\tr X)^2+\frac{1}{2}\tr Xh\right)\Ab+r^{-2}\dk^{\leq 1}(\Ga_g\c\Ga_b)  
\eeaa
where we used the null structure equation for $\nabc_4\tr X$ and the fact that $\Xi\in r^{-1}\Ga_g$ in this chapter. Recall from \eqref{eq:defintionofCb1andCb2fordefintionqfb:chap12} that $\und{C}_1$ is given by 
\beaa
\und{C}_1&=& 2\trch - 2\frac {\atrch^2}{ \trch}  -4 i \atrch.
\eeaa
We choose 
\beaa
h := \und{C}_1-2\tr X
\eeaa
so that 
\beaa
h &=&  2\trch - 2\frac {\atrch^2}{ \trch}  -4 i \atrch -2\trch +2i\atrch\\
&=& - 2\frac {\atrch^2}{ \trch}  -2i \atrch.   
\eeaa
With this choice of $h$, we deduce 
\beaa
&&\left(\nabc_4 +\frac{3}{2}\tr X+h\right)\left(\nabc_4\Ab +\frac{1}{2}\tr X\Ab\right)\\
&=& \nabc_4^2\Ab+\und{C}_1\nabc_4\Ab+\left(\frac{1}{2}(\tr X)^2+\frac{1}{2}\tr X\left(- 2\frac {\atrch^2}{ \trch}  -2i \atrch\right)\right)\Ab+r^{-2}\dk^{\leq 1}(\Ga_g\c\Ga_b).  
\eeaa

Next, recall from \eqref{eq:defintionofCb1andCb2fordefintionqfb:chap12} that $\und{C}_2$ is given by
\beaa
\und{C}_2  &=& \frac 1 2 \trch^2- 4\atrch^2+\frac 3 2 \frac{\atrch^4}{\trch^2} +  i \left(-2\trch\atrch +4\frac{\atrch^3}{\trch}\right).
\eeaa
We compute 
\beaa
&&\und{C}_2  -\left(\frac{1}{2}(\tr X)^2+\frac{1}{2}\tr X\left(- 2\frac {\atrch^2}{ \trch}  -2i \atrch\right)\right)\\
&=& \frac 1 2 \trch^2- 4\atrch^2+\frac 3 2 \frac{\atrch^4}{\trch^2} +  i \left(-2\trch\atrch +4\frac{\atrch^3}{\trch}\right)\\
&& -\left(\frac{1}{2}(\tr X)^2+\frac{1}{2}\tr X\left(- 2\frac {\atrch^2}{ \trch}  -2i \atrch\right)\right)\\
&=& - \frac{3}{2}\atrch^2+\frac 3 2 \frac{\atrch^4}{\trch^2} + \frac{3i\atrch^3}{\trch}
\eeaa
where both the $\trch^2$ in the real part and the $\trch\atrch$ terms in the imaginary part cancel. We infer
\beaa
&&\left(\nabc_4 +\frac{3}{2}\tr X+h\right)\left(\nabc_4\Ab +\frac{1}{2}\tr X\Ab\right)\\
&=& \nabc_4^2\Ab+\und{C}_1\nabc_4\Ab+\left(\und{C}_2  + \frac{3}{2}\atrch^2-\frac 3 2 \frac{\atrch^4}{\trch^2} - \frac{3i\atrch^3}{\trch}\right)\Ab+r^{-2}\dk^{\leq 1}(\Ga_g\c\Ga_b).  
\eeaa
Since
\beaa
 \frac{3}{2}\atrch^2-\frac 3 2 \frac{\atrch^4}{\trch^2} - \frac{3i\atrch^3}{\trch} &=& O(a^2r^{-4})+r^{-2}\Ga_g,
\eeaa
we deduce 
\beaa
&&\left(\nabc_4 +\frac{3}{2}\tr X+h\right)\left(\nabc_4\Ab +\frac{1}{2}\tr X\Ab\right)\\
&=& \nabc_4^2\Ab+\und{C}_1\nabc_4\Ab+\und{C}_2\Ab +O(a^2r^{-4})\Ab+r^{-2}\dk^{\leq 1}(\Ga_g\c\Ga_b)  
\eeaa
and hence, using the definition of $\qfb$,   
\beaa
 \ov{q}q^3\left(\nabc_4 +\frac{3}{2}\tr X+h\right)\left(\nabc_4\Ab +\frac{1}{2}\tr X\Ab\right) &=& \qfb + O(a^2)\Ab+r^2\dk^{\leq 1}(\Ga_g\c\Ga_b). 
\eeaa
In view of the definition of $h$, this concludes the proof of Lemma \ref{lemma:almostfactorizationfoqfbneededforcontralangularderivativesPsib:chap12}.
\end{proof}

\begin{lemma}\lab{Lemma:TeukoskyAb-sharp:spacialcombinationnabc4derivative} 
We have 
\beaa
\bsplit
& \frac{1}{4}\DDc\hot\left(\DDbc \c\left(\nabc_4\Ab+\left(\frac{1}{2}\tr X - 2\frac {\atrch^2}{ \trch}  -2i \atrch\right)\Ab\right)\right)\\ 
=&  O(ar^{-3})  \nab\nab_4(r\Ab)  + O(ar^{-3})  \nab\Ab +O(r^{-4})\nab_3\qfb+O(r^{-5})\qfb +O(r^{-3})\nab_4(r\Ab)\\
  &  +O(r^{-4})\nab_3\Ab +O(r^{-4} ) \Ab  +r^{-2}\dk^{\leq 2}(\Ga_b\c\Ga_g).
\end{split}
\eeaa
\end{lemma}

\begin{proof}
Recall from Lemma \ref{Lemma:TeukoskyAb-sharp} that we have
\beaa
\bsplit
 \frac{1}{4}\DDc\hot(\DDbc \c\Ab) =& \nabc_3\left(\nabc_4\Ab +\frac{1}{2}\tr X\Ab\right)  +\left( \frac 1 2 \tr \Xb+ 2\ov{\tr \Xb} \right)\left(\nabc_4\Ab +\frac{1}{2}\tr X\Ab\right)\\\
  &   +(H+\ov{H})\c\nab\Ab -\frac{1}{2}(\DDc\c\ov{H})\Ab+ O(ar^{-2})  \nab\Ab +O(r^{-3} ) \Ab\\
  & +r^{-1}\dk^{\leq 1}(\Ga_b\c\Ga_g).
\end{split}
\eeaa
We will differentiate this identity by $\nabc_4 +\frac{3}{2}\tr X - 2\frac {\atrch^2}{ \trch}  -2i \atrch$. To this end, first, note that 
\beaa
&&\left(\nabc_4 +\frac{3}{2}\tr X - 2\frac {\atrch^2}{ \trch}  -2i \atrch\right)\Bigg((H+\ov{H})\c\nab\Ab -\frac{1}{2}(\DDc\c\ov{H})\Ab\\
&&+ O(ar^{-2})  \nab\Ab +O(r^{-3} ) \Ab +r^{-1}\dk^{\leq 1}(\Ga_b\c\Ga_g)\Bigg)\\
&=& \left(\nabc_4 +\frac{3}{2}\tr X\right)\left((H+\ov{H})\c\nab\Ab -\frac{1}{2}(\DDc\c\ov{H})\Ab\right)\\
&&+ O(ar^{-2})  \nab\nab_4\Ab  + O(ar^{-3})  \nab\Ab +O(r^{-3} )\nab_4\Ab +O(r^{-4} ) \Ab +r^{-2}\dk^{\leq 2}(\Ga_b\c\Ga_g)\\
&=& \left(\left(\nabc_4 +\frac{1}{2}\tr X\right)(H+\ov{H})\right)\c\nab\Ab+(H+\ov{H})\c\nab\left(\left(\nabc_4 +\frac{1}{2}\tr X\right)\Ab\right)\\
&& -\frac{1}{2}\left(\DDc\c\left(\left(\nabc_4 +\frac{1}{2}\tr X\right)\ov{H}\right)\right)\Ab -\frac{1}{2}(\DDc\c\ov{H})\left(\nabc_4 +\frac{1}{2}\tr X\right)\Ab\\
&&+ O(ar^{-2})  \nab\nab_4\Ab  + O(ar^{-3})  \nab\Ab +O(r^{-3} )\nab_4\Ab +O(r^{-4} ) \Ab +r^{-2}\dk^{\leq 2}(\Ga_b\c\Ga_g)
\eeaa
where we used the fact that, in view of Lemma \ref{LEMMA:COMMUTATION-FORMULAS-1} and  the fact that $\Xi\in r^{-1}\Ga_g$ in this chapter, we have
\beaa
\, [\nab_4, \nab]F&=& - \frac{1}{2}\tr X\nab F +O(ar^{-2})\nab_4F+O(ar^{-3})F+ r^{-1} \Ga_g \c  \dk^{\leq 1} F.
\eeaa
In view of the following consequence of the null structure equation
\beaa
\nabc_4H +\frac{1}{2}\ov{\tr X}H &=& \nabc_3\Xi +\frac{1}{2}\ov{\tr X}\Hb -\frac{1}{2}\Xh\c(\ov{H}-\ov{\Hb}) -B\\
&=& O(ar^{-3})+r^{-1}\dk^{\leq 1}\Ga_g,\\
\nabc_4H +\frac{1}{2}\tr X H &=& \nabc_4H +\frac{1}{2}\ov{\tr X}H+O(a^2r^{-4})+r^{-2}\Ga_b\\
 &=& O(ar^{-3})+r^{-1}\dk^{\leq 1}\Ga_g,
\eeaa
where we used again the fact that $\Xi\in r^{-1}\Ga_g$ in this chapter, we infer
\beaa
&&\left(\nabc_4 +\frac{3}{2}\tr X - 2\frac {\atrch^2}{ \trch}  -2i \atrch\right)\Bigg((H+\ov{H})\c\nab\Ab -\frac{1}{2}(\DDc\c\ov{H})\Ab\\
&&+ O(ar^{-2})  \nab\Ab +O(r^{-3} ) \Ab +r^{-1}\dk^{\leq 1}(\Ga_b\c\Ga_g)\Bigg)\\
&=& (H+\ov{H})\c\nab\left(\left(\nabc_4 +\frac{1}{2}\tr X\right)\Ab\right) -\frac{1}{2}(\DDc\c\ov{H})\left(\nabc_4 +\frac{1}{2}\tr X\right)\Ab\\
&&+ O(ar^{-2})  \nab\nab_4\Ab  + O(ar^{-3})  \nab\Ab +O(r^{-3} )\nab_4\Ab +O(r^{-4} ) \Ab +r^{-2}\dk^{\leq 2}(\Ga_b\c\Ga_g).
\eeaa
Also, since $(\nabc_4 +\frac{1}{2}\tr X)\Ab\in r^{-2}\dk^{\leq 1}\Ga_b$ in view of Bianchi, we obtain
\beaa
&&\left(\nabc_4 +\frac{3}{2}\tr X - 2\frac {\atrch^2}{ \trch}  -2i \atrch\right)\Bigg((H+\ov{H})\c\nab\Ab -\frac{1}{2}(\DDc\c\ov{H})\Ab\\
&&+ O(ar^{-2})  \nab\Ab +O(r^{-3} ) \Ab +r^{-1}\dk^{\leq 1}(\Ga_b\c\Ga_g)\Bigg)\\
&=&  O(ar^{-2})  \nab\nab_4\Ab  + O(ar^{-3})  \nab\Ab +O(r^{-3} )\nab_4\Ab +O(r^{-4} ) \Ab +r^{-2}\dk^{\leq 2}(\Ga_b\c\Ga_g).
\eeaa
Recalling
\beaa
\bsplit
 \frac{1}{4}\DDc\hot(\DDbc \c\Ab) =& \nabc_3\left(\nabc_4\Ab +\frac{1}{2}\tr X\Ab\right)  +\left( \frac 1 2 \tr \Xb+ 2\ov{\tr \Xb} \right)\left(\nabc_4\Ab +\frac{1}{2}\tr X\Ab\right)\\\
  &   +(H+\ov{H})\c\nab\Ab -\frac{1}{2}(\DDc\c\ov{H})\Ab+ O(ar^{-2})  \nab\Ab +O(r^{-3} ) \Ab\\
  & +r^{-1}\dk^{\leq 1}(\Ga_b\c\Ga_g),
\end{split}
\eeaa
and  differentiating this identity by $\nabc_4 +\frac{3}{2}\tr X - 2\frac {\atrch^2}{ \trch}  -2i \atrch$, we deduce
\beaa
\bsplit
&  \frac{1}{4}\left(\nabc_4 +\frac{3}{2}\tr X - 2\frac {\atrch^2}{ \trch}  -2i \atrch\right)\DDc\hot(\DDbc \c\Ab)\\ 
=& \nabc_3\left(\left(\nabc_4 +\frac{3}{2}\tr X - 2\frac {\atrch^2}{ \trch}  -2i \atrch\right)\left(\nabc_4\Ab +\frac{1}{2}\tr X\Ab\right)\right)\\ 
& +\left[\nabc_4 +\frac{3}{2}\tr X - 2\frac {\atrch^2}{ \trch}  -2i \atrch, \nabc_3\right]\left(\nabc_4\Ab +\frac{1}{2}\tr X\Ab\right)\\ 
& +\left( \frac 1 2 \tr \Xb+ 2\ov{\tr \Xb} \right)\left(\nabc_4 +\frac{3}{2}\tr X - 2\frac {\atrch^2}{ \trch}  -2i \atrch\right)\left(\nabc_4\Ab +\frac{1}{2}\tr X\Ab\right)\\
&+\nabc_4\left( \frac 1 2 \tr \Xb+ 2\ov{\tr \Xb} \right)\left(\nabc_4\Ab +\frac{1}{2}\tr X\Ab\right)\\
  & +O(ar^{-2})  \nab\nab_4\Ab  + O(ar^{-3})  \nab\Ab +O(r^{-3} )\nab_4\Ab +O(r^{-4} ) \Ab +r^{-2}\dk^{\leq 2}(\Ga_b\c\Ga_g).
\end{split}
\eeaa
Also, recalling from Lemma \ref{lemma:almostfactorizationfoqfbneededforcontralangularderivativesPsib:chap12} that we have
\beaa
 &&\ov{q}q^3\left(\nabc_4 +\frac{3}{2}\tr X - 2\frac {\atrch^2}{ \trch}  -2i \atrch\right)\left(\nabc_4\Ab +\frac{1}{2}\tr X\Ab\right)\\ 
 &=& \qfb + O(a^2)\Ab+r^2\dk^{\leq 1}(\Ga_g\c\Ga_b),
\eeaa
we infer
\beaa
\bsplit
& \frac{1}{4}\left(\nabc_4 +\frac{3}{2}\tr X - 2\frac {\atrch^2}{ \trch}  -2i \atrch\right)\DDc\hot(\DDbc \c\Ab)\\ 
=& \nabc_3\left(\frac{1}{\ov{q}q^3}\qfb\right) +\left[\nabc_4 +\frac{3}{2}\tr X - 2\frac {\atrch^2}{ \trch}  -2i \atrch, \nabc_3\right]\left(\nabc_4\Ab +\frac{1}{2}\tr X\Ab\right)\\ 
& +\frac{1}{\ov{q}q^3}\left( \frac 1 2 \tr \Xb+ 2\ov{\tr \Xb} \right)\qfb +\nabc_4\left( \frac 1 2 \tr \Xb+ 2\ov{\tr \Xb} \right)\left(\nabc_4\Ab +\frac{1}{2}\tr X\Ab\right)\\
  & +O(ar^{-2})  \nab\nab_4\Ab  + O(ar^{-3})  \nab\Ab +O(r^{-3} )\nab_4\Ab +O(r^{-4})\nab_3\Ab +O(r^{-4} ) \Ab \\
  & +r^{-2}\dk^{\leq 2}(\Ga_b\c\Ga_g)
\end{split}
\eeaa
and hence
\beaa
\bsplit
& \frac{1}{4}\left(\nabc_4 +\frac{3}{2}\tr X - 2\frac {\atrch^2}{ \trch}  -2i \atrch\right)\DDc\hot(\DDbc \c\Ab)\\ 
=&  \left[\nabc_4 +\frac{3}{2}\tr X - 2\frac {\atrch^2}{ \trch}  -2i \atrch, \nabc_3\right]\left(\nabc_4\Ab +\frac{1}{2}\tr X\Ab\right)\\ 
&  +\nabc_4\left( \frac 1 2 \tr \Xb+ 2\ov{\tr \Xb} \right)\left(\nabc_4\Ab +\frac{1}{2}\tr X\Ab\right)\\
  & +O(ar^{-2})  \nab\nab_4\Ab  + O(ar^{-3})  \nab\Ab +O(r^{-4})\nab_3\qfb+O(r^{-5})\qfb \\
  & +O(r^{-3} )\nab_4\Ab +O(r^{-4})\nab_3\Ab +O(r^{-4} ) \Ab  +r^{-2}\dk^{\leq 2}(\Ga_b\c\Ga_g).
\end{split}
\eeaa
Since we have, in view of the null structure equations, 
\beaa
\nabc_4\tr\Xb  &=& O(r^{-2})+r^{-1}\dk^{\leq 1}\Ga_g, 
\eeaa
and since 
\beaa
\nabc_4\Ab +\frac{1}{2}\tr X\Ab &=& r^{-1}\nab_4(r\Ab)+O(r^{-2})\Ab+\Ga_g\c\Ga_b,
\eeaa
we obtain
\beaa
\bsplit
& \frac{1}{4}\left(\nabc_4 +\frac{3}{2}\tr X - 2\frac {\atrch^2}{ \trch}  -2i \atrch\right)\DDc\hot(\DDbc \c\Ab)\\ 
=&  \left[\nabc_4 +\frac{3}{2}\tr X - 2\frac {\atrch^2}{ \trch}  -2i \atrch, \nabc_3\right]\left(\nabc_4\Ab +\frac{1}{2}\tr X\Ab\right)\\ 
  & +O(ar^{-2})  \nab\nab_4\Ab  + O(ar^{-3})  \nab\Ab +O(r^{-4})\nab_3\qfb+O(r^{-5})\qfb +O(r^{-3})\nab_4(r\Ab)\\
  &  +O(r^{-4})\nab_3\Ab +O(r^{-4} ) \Ab  +r^{-2}\dk^{\leq 2}(\Ga_b\c\Ga_g).
\end{split}
\eeaa

Next, using the null structure equation for $\nabc_3\tr X$, i.e. 
\beaa
\nabc_3\tr X +\frac{1}{2}\tr\Xb\tr X &=& \DDc\c\ov{H}+H\c\ov{H}+2P+\Xib\c\ov{\Xi}-\frac{1}{2}\Xbh\c\ov{\Xh}\\
&=&  O(r^{-3})+r^{-1}\dk^{\leq 1}\Ga_b+\Ga_g\c\Ga_b,
\eeaa
as well as the commutator formula 
\beaa
\, [\nabc_3, \nabc_4] U   &=&  O(ar^{-2})\nab U  +O(r^{-3})U+\Ga_b\c\nab U+r^{-1}\Ga_g U,
 \eeaa
we infer, using again $(\nabc_4 +\frac{1}{2}\tr X)\Ab\in r^{-2}\dk^{\leq 1}\Ga_b$ in view of Bianchi,
\beaa
&&\left[\nabc_4 +\frac{3}{2}\tr X - 2\frac {\atrch^2}{ \trch}  -2i \atrch, \nabc_3\right]\left(\nabc_4\Ab +\frac{1}{2}\tr X\Ab\right)\\
&=& O(ar^{-2})\nab\left(\nabc_4\Ab +\frac{1}{2}\tr X\Ab\right)+O(r^{-2})\left(\nabc_4\Ab +\frac{1}{2}\tr X\Ab\right)+r^{-3}\dk^{\leq 2}(\Ga_b\c\Ga_b).
\eeaa
Using again 
\beaa
\nabc_4\Ab +\frac{1}{2}\tr X\Ab &=& r^{-1}\nab_4(r\Ab)+O(r^{-2})\Ab+\Ga_g\c\Ga_b,
\eeaa
this yields 
\beaa
&&\left[\nabc_4 +\frac{3}{2}\tr X - 2\frac {\atrch^2}{ \trch}  -2i \atrch, \nabc_3\right]\left(\nabc_4\Ab +\frac{1}{2}\tr X\Ab\right)\\
&=& O(ar^{-3})\nab(\nab_4(r\Ab))+O(ar^{-4})\nab\Ab+O(r^{-3})\nab_4(r\Ab)+O(r^{-4})\Ab+r^{-2}\dk^{\leq 2}(\Ga_g\c\Ga_b)
\eeaa
and hence
\beaa
\bsplit
& \frac{1}{4}\left(\nabc_4 +\frac{3}{2}\tr X - 2\frac {\atrch^2}{ \trch}  -2i \atrch\right)\DDc\hot(\DDbc \c\Ab)\\ 
=&  O(ar^{-3})  \nab\nab_4(r\Ab)  + O(ar^{-3})  \nab\Ab +O(r^{-4})\nab_3\qfb+O(r^{-5})\qfb +O(r^{-3})\nab_4(r\Ab)\\
  &  +O(r^{-4})\nab_3\Ab +O(r^{-4} ) \Ab  +r^{-2}\dk^{\leq 2}(\Ga_b\c\Ga_g).
\end{split}
\eeaa

Next, we use again the following commutation formula 
\beaa
\, [\nab_4, \nab]F&=& - \frac{1}{2}\tr X\nab F +O(ar^{-2})\nab_4F+O(ar^{-3})F+ r^{-1} \Ga_g \c  \dk^{\leq 1} F
\eeaa
twice to obtain 
\beaa
\nabc_4\DDc\hot(\DDbc \c\Ab) &=&\DDc\hot(\DDbc \c\nabc_4\Ab) -\frac{1}{2}\tr X\DDc\hot(\DDbc \c\Ab)\\
&& -\frac{1}{2}\DDc\hot(\tr X\DDbc \c\Ab) +O(ar^{-3})  \nab\nab_4(r\Ab) + O(ar^{-3})  \nab\Ab\\
&&  +O(ar^{-3})\nab_4(r\Ab)  +O(ar^{-4} ) \Ab  +r^{-2}\dk^{\leq 2}(\Ga_b\c\Ga_g)\\
&=& \DDc\hot(\DDbc \c\nabc_4\Ab) -\tr X\DDc\hot(\DDbc \c\Ab) +O(ar^{-3})  \nab\nab_4(r\Ab)\\
&& + O(ar^{-3})  \nab\Ab  +O(ar^{-3})\nab_4(r\Ab)  +O(ar^{-4} ) \Ab  +r^{-2}\dk^{\leq 2}(\Ga_b\c\Ga_g).
\eeaa
We infer
\beaa
&&\left(\nabc_4 +\frac{3}{2}\tr X - 2\frac {\atrch^2}{ \trch}  -2i \atrch\right)\DDc\hot(\DDbc \c\Ab)\\
&=& \DDc\hot(\DDbc \c\nabc_4\Ab)+\left(\frac{1}{2}\tr X - 2\frac {\atrch^2}{ \trch}  -2i \atrch\right)\DDc\hot(\DDbc \c\Ab)\\
&& +O(ar^{-3})  \nab\nab_4(r\Ab) + O(ar^{-3})  \nab\Ab  +O(ar^{-3})\nab_4(r\Ab)  +O(ar^{-4} ) \Ab  +r^{-2}\dk^{\leq 2}(\Ga_b\c\Ga_g)\\
&=& \DDc\hot\left(\DDbc \c\left(\nabc_4\Ab+\left(\frac{1}{2}\tr X - 2\frac {\atrch^2}{ \trch}  -2i \atrch\right)\Ab\right)\right)\\
&& +O(ar^{-3})  \nab\nab_4(r\Ab) + O(ar^{-3})  \nab\Ab  +O(ar^{-3})\nab_4(r\Ab)  +O(ar^{-4} ) \Ab  +r^{-2}\dk^{\leq 2}(\Ga_b\c\Ga_g).
\eeaa
This implies
\beaa
\bsplit
& \frac{1}{4}\DDc\hot\left(\DDbc \c\left(\nabc_4\Ab+\left(\frac{1}{2}\tr X - 2\frac {\atrch^2}{ \trch}  -2i \atrch\right)\Ab\right)\right)\\ 
=&  O(ar^{-3})  \nab\nab_4(r\Ab)  + O(ar^{-3})  \nab\Ab +O(r^{-4})\nab_3\qfb+O(r^{-5})\qfb +O(r^{-3})\nab_4(r\Ab)\\
  &  +O(r^{-4})\nab_3\Ab +O(r^{-4} ) \Ab  +r^{-2}\dk^{\leq 2}(\Ga_b\c\Ga_g).
\end{split}
\eeaa
as stated. This concludes the proof of Lemma \ref{Lemma:TeukoskyAb-sharp:spacialcombinationnabc4derivative}.
\end{proof}


 \subsection{Estimates for $\nab\Ab$}
 \lab{sec:estimatefornabAb:chap12}


The following lemma provides the control of $\nab\Ab$. 
\begin{lemma}
\lab{lemma:EstimatesfornabAb}
 The following estimates hold true, for all $p\le 2-\de$,
       \bea
       \lab{Estimates:nabAb-p}
    \int_{\MM(\tau_1, \tau_2)} r^{p-1}|\nab  \Ab |^2 &\les& \Bdot_p[\Ab](\tau_1, \tau_2)+\ep_0^2 \tau_1^{-2- 2\dec}.
     \eea
Also, we have, for all $p\le 2-\de$,
\bea
\bsplit
\lab{Eq:EstimatenabAb-prM}
\int_{\Si(\tau) } r^{p-2}|\nab  \Ab |^2  &\les \Edot_p[\Ab](\tau)+\int_{\Si(\tau)}  r^{p-2}|\nab_3\nab_4(r\Ab)|^2 +\ep_0^2\tau^{-2- 2\dec}.
\end{split}
\eea
\end{lemma}

\begin{proof}
According to the  elliptic type   estimates of  Lemma \ref{Lemma:EstimateangderivativesA-LL(a)}  we have
   for any  $ S \subset\MM$.
\beaa
\int_S  \Big(|\nab \Ab|^2  +   r^{-2} |\Ab|^2 \Big)  &\les&\left| \int_S   \Ab\c \DD\hot (\DDb \c \Ab)  \right| +(a^2+\ep^2)  \int_S r^{-2} |( \nabc_3, \nabc_4)  \Ab|^2. 
\eeaa
We  deduce,    on $\MM(\tau_1, \tau_2)$,  
     \bea
     \lab{eq:applyTeukoskyAb-sharp}
     \bsplit
    \int_{\MM(\tau_1, \tau_2)} r^{p-1}\big(  |\nab  \Ab |^2 + r^{-2} | \Ab|^2\big)&\les   \int_{\MM(\tau_1, \tau_2)} r^{p-3}  |\Ab|^2 +    \int_{\MM(\tau_1, \tau_2)} r^{p+1}  \big|\DD \hot (\DDb \c \Ab)\big|^2     \\
    &+(a^2+\ep^2)  \int_{\MM(\tau_1, \tau_2)}   r^{p-3} |( \nabc_3, \nabc_4)  \Ab|^2. 
    \end{split}
     \eea
    Next, recall \eqref{eq:TeukoskyAb-sharp} 
     \beaa
\bsplit
 \frac{1}{4}\DDc\hot(\DDbc \c\Ab) =& \nabc_3\left(\nabc_4\Ab +\frac{1}{2}\tr X\Ab\right)  +\left( \frac 1 2 \tr \Xb+ 2\ov{\tr \Xb} \right)\left(\nabc_4\Ab +\frac{1}{2}\tr X\Ab\right)\\\
  &   +(H+\ov{H})\c\nab\Ab -\frac{1}{2}(\DDc\c\ov{H})\Ab+ O(ar^{-2})  \nab\Ab +O(r^{-3} ) \Ab\\
  & +r^{-1}\dk^{\leq 1}(\Ga_b\c\Ga_g),
\end{split}
\eeaa
which together with 
\beaa
\nabc_4\Ab +\frac{1}{2}\tr X\Ab &=& r^{-1}\nab_4(r\Ab)+O(r^{-2})\Ab+\Ga_g\c\Ga_b,
\eeaa
yields   
\beaa
|\DDc\hot(\DDbc \c\Ab)| &\les& ar^{-2}|\nab\Ab|+r^{-1}|\nab_3(\nab_4(r\Ab))|+r^{-2}|\nab_4(r\Ab)|+r^{-2}|\nab_3\Ab|+r^{-3}|\Ab|\\
&&+r^{-1}|\dk^{\leq 1}(\Ga_b\c\Ga_b)|.
\eeaa
In view of  the definition of  the norms  $\Bdot_p[\Ab] $, see Definition \ref{definition:dotnormsforAb},   we   infer
     \beaa
      && \int_{\MM(\tau_1, \tau_2)} r^{p+1}  \big|\DD \hot (\DDb \c A)\big|^2\\
      &\les&
    \int_{\MM(\tau_1, \tau_2)}  r^{p-1} \Big( |\nab_3 \big( \nab_4(r\Ab)\big)|^2 +| \nab_4(r\Ab) |^2 \Big)+
    a^2\int_{\MM(\tau_1, \tau_2)}   r^{p-3} |\nab \Ab|^2 \\
    &+&\int_{\MM(\tau_1, \tau_2)} r^{p-3}\Big(|\nab_3\Ab|^2+ |\Ab|^2\Big) +   \int_{\MM(\tau_1, \tau_2)}    r^{p-1} |\dk^{\leq 1}(\Ga_b\c\Ga_b) |^2 \\
      &\les&\Bdot_p[\Ab]+a^2  \int_\MM r^{p-1} |\nab \Ab|^2+ \int_\MM    r^{p-1} |\dk^{\leq 1}(\Ga_b\c\Ga_g) |^2
     \eeaa     
     where we recall
 \beaa
 \Bdot_p[\Ab](\tau_1, \tau_2) &=&\int_{\MM(\tau_1, \tau_2)} r^{p-1} \Big(   r^2| \nab_4 \nab_4(r\Ab) |^2  +  |\nab_3 \nab_4(r\Ab) |^2 + | \nab_4(r\Ab) |^2 \Big)\\
&&+
 \int_{\MM(\tau_1, \tau_2)} r^{p-3 }\Big(   r^2| \nab_4 \Ab|^2+  |\nab_3 \Ab|^2+ |\Ab|^2    \Big).
 \eeaa
 Also 
\bea
\lab{Estimate:nabAb-error}
 \int_{\MM(\tau_1, \tau_2)}    r^{p-1} |\dk^{\leq 1}(\Ga_b\c\Ga_b) |^2 \les  \ep^4\tau_1^{-2- 2 \dec}
\eea
which holds true for $p\leq 2-\de$.
Hence, for all $p\leq 2-\de$,
      \beaa
     \int_{\MM(\tau_1, \tau_2)} r^{p+1}  \big|\DD \hot (\DDb \c A)\big|^2  &\les&\Bdot_p[\Ab]+a^2  \int_\MM r^{p-1} |\nab \Ab|^2+\ep_0^2 \tau_1^{-2- 2\dec}.
     \eeaa
     Back to \eqref{eq:applyTeukoskyAb-sharp}, we deduce,
     after absorbing the  term  $a^2  \int_\MM r^{p-1} |\nab \Ab|^2$ on the left for $a$  small enough, for all $p\le 2-\de$,
       \beaa
    \int_{\MM(\tau_1, \tau_2)} r^{p-1}\big(  |\nab  \Ab |^2 + r^{-2} | \Ab|^2\big)&\les& \Bdot_p[\Ab]+\ep_0^2 \tau_1^{-2- 2\dec}
     \eeaa
which yields \eqref{Estimates:nabAb-p}. In the same vein, we derive
\beaa
\bsplit
\int_{\Si(\tau) } r^{p-2} \big( |\nab  \Ab |^2 + r^{-2} | \Ab |^2\big) &\les \Edot_p[\Ab](\tau)+\int_{\Si(\tau)}  r^{p-2}|\nab_3\nab_4(r\Ab)|^2 +\ep_0^2\tau^{-2- 2\dec},
\end{split}
\eeaa
which yields \eqref{Eq:EstimatenabAb-prM}. This concludes the proof of Lemma \ref{lemma:EstimatesfornabAb}.
\end{proof}


\subsection{Estimates for $\nab\nab_4(r\Ab)$}
\lab{sec:estimatefornabnab4Ab:chap12}


 It  remains to  estimate  the terms involving $\nab \nab_4(r \Ab)$.  This is achieved in the lemma below.
 \begin{lemma}
 \lab{lemma:estimatesnabnab_4(rAb)}
 The following estimates hold true in $\MM=\MM(\tau_1,\tau_2)$ for all $p\le 2-\de$
\beaa
     \bsplit
    \int_{\MM(\tau_1, \tau_2)} r^{p+1}|\nab\nab_4(r\Ab)|^2 &\les     B_\de[\qfb](\tau_1, \tau_2)+\Bdot_p[\Ab](\tau_1, \tau_2)    +  \ep_0^2\tau_1^{-2-2\dec},\\
\int_{\Si(\tau) } r^p |\nab\nab_4(r\Ab)|^2  &\les E_\de[\qfb](\tau)+\Edot_p[\Ab](\tau)+\int_{\Si(\tau)}  r^{p-2}|\nab_3\nab_4(r\Ab)|^2  +\ep_0^2\tau^{-2- 2\dec}.     
    \end{split}
     \eeaa  
 \end{lemma}
 
\begin{proof}
Let us introduce the notation 
\beaa
Y[\Ab] &:=& \nabc_4\Ab+\left(\frac{1}{2}\tr X - 2\frac {\atrch^2}{ \trch}  -2i \atrch\right)\Ab.
\eeaa
According to the  elliptic type   estimates of  Lemma \ref{Lemma:EstimateangderivativesA-LL(a)}  we have
   for any  $ S \subset\MM$.
\beaa
\int_S  \Big(|\nab Y[\Ab]|^2  +   r^{-2} |Y[\Ab]|^2 \Big)  &\les&\left| \int_S   \Ab\c \DD\hot (\DDb \c Y[\Ab])  \right| \\
&& +(a^2+\ep^2)  \int_S r^{-2} |( \nabc_3, \nabc_4)Y[\Ab]|^2. 
\eeaa
We  deduce,    on $\MM(\tau_1, \tau_2)$,  
     \bea
     \lab{eq:applyTeukoskyAb-sharp:secondtimenabc4combination}
     \bsplit
    &\int_{\MM(\tau_1, \tau_2)} r^{p+3}\big(  |\nab Y[\Ab]|^2 + r^{-2} |Y[\Ab]|^2\big)\\
    &\les   \int_{\MM(\tau_1, \tau_2)} r^{p+1}  |Y[\Ab]|^2 +    \int_{\MM(\tau_1, \tau_2)} r^{p+5}  \big|\DD \hot (\DDb \c Y[\Ab])\big|^2     \\
    &+(a^2+\ep^2)  \int_{\MM(\tau_1, \tau_2)}   r^{p+1} |( \nabc_3, \nabc_4)Y[\Ab]|^2. 
    \end{split}
     \eea
    Next, recall from Lemma \ref{Lemma:TeukoskyAb-sharp:spacialcombinationnabc4derivative} that we have, using also the above definition of $Y[\Ab]$, 
\beaa
\bsplit
\frac{1}{4}\DDc\hot\left(\DDbc \c\left(Y[\Ab]\right)\right) =&  O(ar^{-3})  \nab\nab_4(r\Ab)  + O(ar^{-3})  \nab\Ab +O(r^{-4})\nab_3\qfb+O(r^{-5})\qfb\\
& +O(r^{-3})\nab_4(r\Ab)  +O(r^{-4})\nab_3\Ab +O(r^{-4} ) \Ab  +r^{-2}\dk^{\leq 2}(\Ga_b\c\Ga_g).
\end{split}
\eeaa
which yields 
\beaa
&&\int_{\MM(\tau_1, \tau_2)} r^{p+5}  \big|\DD \hot (\DDb \c Y[\Ab])\big|^2 \\
&\les& \int_{\MM(\tau_1, \tau_2)} r^{p-3}\Big(|\nab_3\qfb|^2+r^{-2}|\qfb|^2+|\nab_3\Ab|^2+|\Ab|^2\Big) +\int_{\MM(\tau_1, \tau_2)} r^{p-1}|\nab_4(r\Ab)|^2\\
&&+ a^2\int_{\MM(\tau_1, \tau_2)} r^{p-1}\Big(|\nab\nab_4(r\Ab)|^2+|\nab\Ab|^2\Big) +  \int_{\MM(\tau_1, \tau_2)} r^{p+1}|\dk^{\leq 2}(\Ga_b\c\Ga_g)|^2.
\eeaa
In view of the definition of $\Bdot_p[\Ab](\tau_1, \tau_2)$ and $B_p[\qfb](\tau_1, \tau_2)$, we infer, for $p\leq 2-\de$, 
\beaa
&&\int_{\MM(\tau_1, \tau_2)} r^{p+5}  \big|\DD \hot (\DDb \c Y[\Ab])\big|^2 \\
&\les& B_\de[\qfb](\tau_1, \tau_2)+\Bdot_p[\Ab](\tau_1, \tau_2)\\
&&+ a^2\int_{\MM(\tau_1, \tau_2)} r^{p-1}\Big(|\nab\nab_4(r\Ab)|^2+|\nab\Ab|^2\Big) +  \int_{\MM(\tau_1, \tau_2)} r^{p+1}|\dk^{\leq 2}(\Ga_b\c\Ga_g)|^2.
\eeaa
Together with the control of $\Ga_g$ and $\Ga_b$, we infer, for $p\leq 2-\de$, 
\beaa
&&\int_{\MM(\tau_1, \tau_2)} r^{p+5}  \big|\DD \hot (\DDb \c Y[\Ab])\big|^2 \\
&\les& B_\de[\qfb](\tau_1, \tau_2)+\Bdot_p[\Ab](\tau_1, \tau_2) + a^2\int_{\MM(\tau_1, \tau_2)} r^{p-1}\Big(|\nab\nab_4(r\Ab)|^2+|\nab\Ab|^2\Big) +  \ep^4\tau_1^{-2-2\dec}.
\eeaa
Pluggin in \eqref{eq:applyTeukoskyAb-sharp:secondtimenabc4combination}, we deduce,  for $p\leq 2-\de$, 
 \beaa
     \bsplit
    &\int_{\MM(\tau_1, \tau_2)} r^{p+3}\big(  |\nab Y[\Ab]|^2 + r^{-2} |Y[\Ab]|^2\big)\\
    &\les   \int_{\MM(\tau_1, \tau_2)} r^{p+1}  |Y[\Ab]|^2 +   B_\de[\qfb](\tau_1, \tau_2)+\Bdot_p[\Ab](\tau_1, \tau_2) \\
    &+ a^2\int_{\MM(\tau_1, \tau_2)} r^{p-1}\Big(|\nab\nab_4(r\Ab)|^2+|\nab\Ab|^2\Big)    \\
    &+(a^2+\ep^2)  \int_{\MM(\tau_1, \tau_2)}   r^{p+1} |( \nabc_3, \nabc_4)Y[\Ab]|^2 +  \ep_0^2\tau_1^{-2-2\dec}  . 
    \end{split}
     \eeaa
Noticing that 
\beaa
Y[\Ab] &=& \nabc_4\Ab+\left(\frac{1}{2}\tr X - 2\frac {\atrch^2}{ \trch}  -2i \atrch\right)\Ab\\
&=& r^{-1}\nab_4(r\Ab)+O(ar^{-2})\Ab+\Ga_g\c\Ga_b,
\eeaa
we infer, together with the control of $\Ga_g$ and $\Ga_b$ and the definition of $\Bdot_p[\Ab](\tau_1, \tau_2)$, 
for $p\leq 2-\de$, 
 \beaa
     \bsplit
    \int_{\MM(\tau_1, \tau_2)} r^{p+1}|\nab\nab_4(r\Ab)|^2 &\les     B_\de[\qfb](\tau_1, \tau_2)+\Bdot_p[\Ab](\tau_1, \tau_2) \\
    &+ a^2\int_{\MM(\tau_1, \tau_2)} r^{p-1}\Big(|\nab\nab_4(r\Ab)|^2+|\nab\Ab|^2\Big)    +  \ep_0^2\tau_1^{-2-2\dec}  . 
    \end{split}
     \eeaa
Absorbing the term $\nab\nab_4(r\Ab)$ on the RHS from $a$ small enough, and using the control provided by \eqref{Estimates:nabAb-p} for the term $\nab\Ab$, we infer
 \beaa
     \bsplit
    \int_{\MM(\tau_1, \tau_2)} r^{p+1}|\nab\nab_4(r\Ab)|^2 &\les     B_\de[\qfb](\tau_1, \tau_2)+\Bdot_p[\Ab](\tau_1, \tau_2)    +  \ep_0^2\tau_1^{-2-2\dec}  
    \end{split}
     \eeaa
as stated.

In the same vein, we derive
\beaa
\bsplit
\int_{\Si(\tau) } r^p \big( |\nab\nab_4(r\Ab)|^2 + r^{-2} | \nab_4(r\Ab) |^2\big) &\les E_\de[\qfb](\tau)+\Edot_p[\Ab](\tau)+\int_{\Si(\tau)}  r^{p-2}|\nab_3\nab_4(r\Ab)|^2\\
& +\ep_0^2\tau^{-2- 2\dec},
\end{split}
\eeaa
as stated. This ends the proof of Lemma \ref{lemma:estimatesnabnab_4(rAb)}. 
\end{proof}

   
   \subsection{End of the proof of Proposition \ref{prop:MaiTransportAb-steps}}
   \lab{sec:endoftheproofofproposition:MaiTransportAb-steps}
   
   
We combine \eqref{Estimates:nabAb-p} and the first estimate of Lemma \ref{lemma:estimatesnabnab_4(rAb)}, i.e., for any $p\le 2-\de$,
  \beaa
  \bsplit
    \int_{\MM(\tau_1, \tau_2)} r^{p-1}|\nab  \Ab |^2 &\les \Bdot_p[\Ab](\tau_1, \tau_2)+\ep_0^2 \tau_1^{-2- 2\dec},\\
         \int_{\MM(\tau_1, \tau_2)} r^{p+1}|\nab\nab_4(r\Ab)|^2 &\les     B_\de[\qfb](\tau_1, \tau_2)+\Bdot_p[\Ab](\tau_1, \tau_2)    +  \ep_0^2\tau_1^{-2-2\dec},     
    \end{split}
     \eeaa 
with those of    Proposition  \ref{Prop:EstimatesforBEFdotnormsAb}, i.e. 
\beaa
\BEFdot_p[\Ab](\tau_1, \tau_2) &\les&  B_{\de}[\undpsi](\tau_1, \tau_2)+\Edot_{p}[\Ab] (\tau_1) \\
&&  +   (a^2+\ep^2)\int_{\MM(\tau_1, \tau_2)}r^{p-3}\Big(r^2| \nab\nab_4(r\Ab)|^2+| \nab \Ab|^2\Big)+\ep_0^2\tau^{-2-2\dec}.
\eeaa
This yields, recalling the definition of the norms $\BEF_p[\Ab]$  in Definition \ref{Definition:NormsBEF-Ab} and after absorbing the terms proportional to $a$ on the left, for all $p\le 2-\de$,
  \beaa
\BEF_p[\Ab] (\tau_1, \tau_2)&\les&   B_{\de}[\undpsi](\tau_1,\tau_2)
 + E_p[\Ab](\tau_1) +\ep_0^2 \tau_1^{-2-2\dec}
 \eeaa
 which is \eqref{eq:transportAb} in the case $s=0$. 
 
Next, notice that 
\beaa
&&\int_{\Si(\tau)}r^{p-2}\Big(r^2|\nab_4\nab_4(r\Ab)|^2+|\nab_3\nab_4(r\Ab)|^2\Big)\\
&\les& \int_{\Si(\tau)}r^{p-2}\Big(r^2|\nab_4\nab_4(r\Ab)|^2+|\nab_{\Rhat}\nab_4(r\Ab)|^2+\chi_{red}^2|\nab_3\nab_4(r\Ab)|^2\Big)\\
&\les& \int_{\Si(\tau)}r^{p-2}r^2|\nab_4\nab_4(r\Ab)|^2 +E_p[\Ab](\tau).
\eeaa
Since we have, in view of the definition of $\Psib$ and $\qf$,  
\beaa
\nab_4\nab_4(r\Ab) &=& \nab_4(r^{-1}\Psib+O(ar^{-1})\Ab+r\Ga_b\c\Ga_g)\\
&=& r^{-2}\nab_4(r\Psib) -\frac{1}{r^2}\Psib +O(ar^{-2})\nab_4(r\Ab) +O(ar^{-2})\Ab+\dk^{\leq 1}(\Ga_b\c\Ga_g)\\
&=& O(r^{-1})\qfb+O(r^{-1})\nab_4(r\Ab)+O(ar^{-2})\Ab+\dk^{\leq 1}(\Ga_b\c\Ga_g),
\eeaa
we infer, for $p\leq 1-\de$, 
\beaa
\int_{\Si(\tau)}r^{p}|\nab_4\nab_4(r\Ab)|^2 &\les& E_p[\qfb](\tau)+\ep_0^2\tau^{-2-2\dec}.
\eeaa
Together with \eqref{Eq:EstimatenabAb-prM} and the second estimate of Lemma \ref{lemma:estimatesnabnab_4(rAb)}, i.e., for any $p\le 2-\de$, 
 \beaa
     \bsplit
\int_{\Si(\tau) } r^{p-2}|\nab  \Ab |^2  &\les \Edot_p[\Ab](\tau)+\int_{\Si(\tau)}  r^{p-2}|\nab_3\nab_4(r\Ab)|^2 +\ep_0^2\tau^{-2- 2\dec},\\
\int_{\Si(\tau) } r^p |\nab\nab_4(r\Ab)|^2  &\les E_\de[\qfb](\tau)+\Edot_p[\Ab](\tau)+\int_{\Si(\tau)}  r^{p-2}|\nab_3\nab_4(r\Ab)|^2  +\ep_0^2\tau^{-2- 2\dec},     
    \end{split}
     \eeaa 
and together  with \eqref{eq:transportAb} derived above, this yields, for any $p\leq 1-\de$, 
 \beaa
\nn&&\sup_{\tau\in[\tau_1,\tau_2]}\int_{\Si(\tau)}r^{p-2}\Big(r^2|\nab_4\nab_4(r\Ab)|^2+r^2|\nab\nab_4(r\Ab)|^2+|\nab_3\nab_4(r\Ab)|^2+|\nab\Ab|^2\Big)\\
 &\les& EB_{p}[\undpsi](\tau_1,\tau_2)   + E_p[\Ab](\tau_1) +\ep_0^2 \tau_1^{-2-2\dec}.
\eeaa
which is  \eqref{eq:transportAb:additionalestimateformissingderivativesenergy} in the case $s=0$. 

It remains to recover  \eqref{eq:transportAb} and  \eqref{eq:transportAb:additionalestimateformissingderivativesenergy} for $1\leq s\leq k_L$. To this end, we proceed as follows:
\begin{enumerate}
\item We argue by iteration assuming that \eqref{eq:transportAb} and  \eqref{eq:transportAb:additionalestimateformissingderivativesenergy} hold for some $0\leq s\leq k_L-1$. It it true for $s=0$ by the above, and our goal is to prove that \eqref{eq:transportAb} and  \eqref{eq:transportAb:additionalestimateformissingderivativesenergy} hold with $s$ replaced by $s+1$. 

\item We commute the system of transport equations \eqref{eq:thesystemof2transporteuqaitonsAbPsibundpsiactuallyused}, i.e.
\beaa
\nabc_4(r\Psib)=\frac{q}{r\ov{q}}\qfb  +r\dk^{\leq 1}(\Ga_g\c\Ga_b), \qquad \nabc_4\left( \frac{q^4}{r^3}\Ab\right)= \frac{1}{r}\Psib+r\Ga_g\c\Ga_b,
\eeaa
with $\Lieb_\T$, $\ov{q}\,\ov{\DDc}\c$ and $\chi_{red}\nabc_3$. In view of the commutation formulas of Lemma \ref{lemma:basicpropertiesLiebTfasdiuhakdisug:chap9} and Lemma \ref{COMMUTATOR-NAB-C-3-DD-C-HOT}, we have, for $U\in\sk_2$, 
\beaa
[\nabc_4, \Lieb_\T]U &=& [\nab_4, \Lieb_\T]U -4[\om, \Lieb_\T]U = r^{-1}[\Lieb_\T, \dk]U +r^{-1}\Ga_b U\\
&=& r^{-1}\dk^{\leq 1}(\Ga_b U),
\eeaa 
and
\beaa
 \, [\nabc_4, \ov{q}\,\ov{\DDc}\c] U&=& \ov{q}\, [\nabc_4, \ov{q}\,\ov{\DDc}\c] U +e_4(\ov{q})\ov{\DDc}\c U\\
 &=& - \frac 1 2\ov{q}\left(\ov{\tr X} -\frac{2}{\ov{q}}e_4(\ov{q})\right) \ov{\DDc} \c U +\ov{q}\,\ov{\Hb} \c \nabc_4 U+ \Ga_g \c  \dk^{\leq 1} U\\
 &=& O(ar^{-1})\nabc_4U + \Ga_g \c  \dk^{\leq 1} U.
\eeaa
This yields the commuted systems 
\beaa
\nabc_4\left(\Lieb_\T(r\Psib)\right) &=& \frac{q}{r\ov{q}}\Lieb_\T\qfb  +r\dk^{\leq 2}(\Ga_g\c\Ga_b), \\ \nabc_4\left(\Lieb_\T\left(\frac{q^4}{r^3}\Ab\right)\right) &=& \frac{1}{r}\Lieb_\T\Psib+r\dk^{\leq 1}(\Ga_g\c\Ga_b),
\eeaa
\beaa
\nabc_4\left(\ov{q}\,\ov{\DDc}\c(r\Psib)\right) &=& \frac{q}{r\ov{q}}\dkb^{\leq 1}\qfb  +r\dk^{\leq 2}(\Ga_g\c\Ga_b), \\ 
\nabc_4\left(\ov{q}\,\ov{\DDc}\c\left(\frac{q^4}{r^3}\Ab\right)\right) &=& \frac{1}{r}\dkb^{\leq 1}\Psib+r\dk^{\leq 1}(\Ga_g\c\Ga_b),
\eeaa
and
\beaa
\nabc_4\left(\chi_{red}\nabc_3(r\Psib)\right) &=& \chi_{red}\nabc_3\left(\frac{q}{r\ov{q}}\qfb\right) +\pr_r\chi_{red} e_4(r)\nabc_3(r\Psib)\\
&&+ O\big((a+\ep)r^{-1})\chi_{red}\nabc(r\Psib) + O(r^{-2})\chi_{red}\Psib, \\ 
\nabc_4\left(\chi_{red}\nabc_3\left(\frac{q^4}{r^3}\Ab\right)\right) &=& \chi_{red}\nabc_3\left(\frac{1}{r}\Psib\right)+\pr_r\chi_{red} e_4(r)\nabc_3\left(\frac{q^4}{r^3}\Ab\right)\\
&&+ O\big((a+\ep)r^{-1})\chi_{red}\nabc\left(\frac{q^4}{r^3}\Ab\right) + O(r^{-2})\chi_{red}\Ab.
\eeaa

\item Using the iteration assumption for these commuted systems, and using the original system to recover the $\nab_4$ derivative, we infer that \eqref{eq:transportAb} and  \eqref{eq:transportAb:additionalestimateformissingderivativesenergy} hold for $s$ derivatives with $\Ab$ replaced with $(\Lieb_\T, \ov{q}\,\ov{\DDc}\c, \nab_4, \chi_{red}\nabc_3)\Ab$. Together with:
\begin{enumerate}
\item the link between $\Lieb_\T$ and $\nab_\T$ of Lemma \ref{lemma:basicpropertiesLiebTfasdiuhakdisug:chap9},

\item the Hodge elliptic estimates of Proposition \ref{Prop:HodgeThmM8},

\item the fact that $(\nab_\T, r\nab_4, \dkb, \chi_{red}\nabc_3)$ span $\dk$,
\end{enumerate} 
and using the iteration assumption to absorb lower order terms in differentiability, we infer that \eqref{eq:transportAb} and  \eqref{eq:transportAb:additionalestimateformissingderivativesenergy} hold for $s$ derivatives with $\Ab$ replaced with $\dk^{\leq 1}\Ab$. In particular, \eqref{eq:transportAb} and  \eqref{eq:transportAb:additionalestimateformissingderivativesenergy} hold with $s$ replaced by $s+1$. Thus, by iteration,  \eqref{eq:transportAb} and  \eqref{eq:transportAb:additionalestimateformissingderivativesenergy} hold for all $s$ such that $0\leq s\leq k_L$. This end the proof of Proposition \ref{prop:MaiTransportAb-steps}.  
\end{enumerate}


\section{Proof of Theorem M2}
\lab{sec:finallyproofofThmM2}


In this section, we prove Theorem M2 of \cite{KS:Kerr} by relying on Theorem \ref{Thm:Nondegenerate-Morawetz-psib}.


\subsection{Statement of Theorem M2}
\lab{sec:statementofTheoremM2:cahp12}

 
In this section \ref{sec:statementofTheoremM2:cahp12}, we restate Theorem M2 on the decay of the flux of $\aa$ on $\Si_*$. Note that the global frame used for the proof of Theorem M2 is constructed in section 3.6 of  \cite{KS:Kerr} and satisfies in particular the assumptions of section \ref{sec:assumptionsontheframe:Chapter12}.


\subsubsection{Definition of the $r$-foliation of $\Si_*$}


We consider the foliation on $\Si_*$ induced by the scalar function $r$.
\begin{definition}[$r$-foliation of $\Si_*$]
\lab{def:rfoliationonSigmastar:chap12}
The foliation on $\Si_*$ induced by the scalar function $r$ is such that:
\begin{enumerate}
\item The function $r$ foliates $\Si_*$ by spheres $S_{\Si_*}(r)$.

\item We have $\tau=\tau(r)$ on $\Si_*$, i.e. the restriction to $\Si_*$ of $\tau$ is a function of $r$.

\item We consider a null pair $(e_3^{\Si_*}, e_4^{\Si_*})$ and and orthonormal basis $e^{\Si_*}_a$, $a=1,2$, of the tangent space to $S_{\Si_*}(r)$ such that $(e_3^{\Si_*}, e_4^{\Si_*}, e^{\Si_*}_1, e^{\Si_*}_2)$ forms a null frame on    $\Si_*$.

\item We denote by $\nu$ the unique vectorfield tangent to $\Si_*$ and normal to $S_{\Si_*}(r)$ such that $\nu$ is given on $\Si_*$ by 
\bea\lab{eq:decompoisitionofthetangentvectornuonthenullframeofSigmastar}
\nu &=& e_3^{\Si_*}+b^{\Si_*}e_4^{\Si_*}
\eea
for some scalar function $b^{\Si_*}$. 

\item The Ricci coefficients associated with the null frame $(e_3^{\Si_*}, e_4^{\Si_*}, e^{\Si_*}_1, e^{\Si_*}_2)$ satisfy the following transversality conditions\footnote{Note that, in view of these transversality  conditions, all Ricci coefficients associated to the null frame $(e_3^{\Si_*}, e_4^{\Si_*}, e^{\Si_*}_1, e^{\Si_*}_2)$ are defined on $\Si_*$.} 
\bea
\xi^{\Si_*}=0, \qquad \om^{\Si_*}=0, \qquad \etab^{\Si_*}=-\ze^{\Si_*}, \qquad\textrm{on}\quad\Si_*.
\eea

\item We introduce the following linearized Ricci and curvature components  associated to  the null frame $(e_3^{\Si_*}, e_4^{\Si_*}, e^{\Si_*}_1, e^{\Si_*}_2)$:
\beaa
&&\widecheck{\trch^{\Si_*}} := \trch^{\Si_*}-\frac{2}{r}, \qquad \widecheck{\trchb^{\Si_*}} := \trchb^{\Si_*}+\frac{2\left(1-\frac{2m}{r}\right)}{r},\\ 
&& \widecheck{\omb^{\Si_*}} :=\omb^{\Si_*} -\frac{m}{r^2}, \qquad \widecheck{\rho^{\Si_*}} :=\rho^{\Si_*} +\frac{2m}{r^3}.
\eeaa
We also linearize the scalar function $b^{\Si_*}$ appearing in \eqref{eq:decompoisitionofthetangentvectornuonthenullframeofSigmastar} and $\nu(r)$ as follows
\beaa
\widecheck{b^{\Si_*}} := b^{\Si_*}+1+\frac{2m}{r}, \qquad \widecheck{\nu(r)}:=\nu(r)+2.
\eeaa

\item We group the above linearized quantities as follows
\beaa
\bsplit
\Ga_g^{\Si_*} &= \Big\{\widecheck{\trch^{\Si_*}}, \, \chih^{\Si_*}, \, \eta^{\Si_*},\, \etab^{\Si_*}, \, \ze^{\Si_*}, \,  \widecheck{\trchb^{\Si_*}}, \, \chibh^{\Si_*}, \, \widecheck{\omb^{\Si_*}}, \, \xib^{\Si_*}, \, r\a^{\Si_*}, \, r\b^{\Si_*}, \, r\rhoc^{\Si_*}, \, r\!\rhod^{\Si_*}\Big\},\\
\Ga_b^{\Si_*} &= \Big\{\chibh^{\Si_*}, \, \widecheck{\omb^{\Si_*}}, \, \xib^{\Si_*}, \, \aa^{\Si_*}, \, r\bb^{\Si_*}, \, r^{-1}\widecheck{b^{\Si_*}}, \, r^{-1}\widecheck{\nu(r)}\Big\}.
\end{split}
\eeaa
\end{enumerate}
\end{definition}

\begin{remark}
Definition \ref{def:rfoliationonSigmastar:chap12} is compatible with the definition of the $r$-foliation on the last slice $\Si_*$ in \cite{KS:Kerr}, see for example section 5.1 in  \cite{KS:Kerr}. 
\end{remark}


\subsubsection{Assumptions on the $r$-foliation of $\Si_*$}


We will need the following assumptions on the $r$-foliation of $\Si_*$:
\begin{itemize}
\item The function $r$ satisfies along $\Si_*$ the following lower bound
\bea\lab{eq:cominantconditionofronSigmastar}
\min_{\Si_*}r &\gtrsim& \ep_0^{-1}\tau_*^{1+\dec}. 
\eea

\item $\Ga_g^{\Si_*}$ and $\Ga_b^{\Si_*}$ verify the following estimates on $\Si_*$, for $0\leq k\leq k_L$,  
\bea\lab{eq:controlofGagSistarandGabSistarframeSistar:chap12}
\Big(r^2u^{\frac{1}{2}+\dec}+ru^{1+\dec}\Big)|\dk^k\Ga_g^{\Si_*}|+ru^{1+\dec}|\dk^k\Ga_b^{\Si_*}| &\les& \ep.
\eea

\item The functions $\tau$ and $r$ satisfy the following estimate on $\Si_*$
\bea\lab{eq:propertiesoffirstorderderivativesofrandtauonSigmastar}
|\nu(\tau)-2|\les \frac{1}{r}, \qquad \nab^{\Si_*}(\tau)=0, \qquad \nab^{\Si_*}(r)=0,
\eea
where that last two identities follow from the fact that $\nab^{\Si_*}$ is tangent to $S_{\Si_*}(r)$ and the fact that $\tau$ is a function of $r$ on $\Si_*$. 
\end{itemize}

\begin{remark}
The assumptions above on the $r$-foliation are compatible with the ones on the last slice $\Si_*$ in \cite{KS:Kerr}, see (3.4.5) in \cite{KS:Kerr} for \eqref{eq:cominantconditionofronSigmastar} and section 5.1.4 in \cite{KS:Kerr} for \eqref{eq:controlofGagSistarandGabSistarframeSistar:chap12} and  \eqref{eq:propertiesoffirstorderderivativesofrandtauonSigmastar}.  
\end{remark}


\subsubsection{Comparison of the global frame of $\MM$ with the $r$-foliation of $\Si_*$}


We assume that the global frame of $\MM$ and the null frame $(e_3^{\Si_*}, e_4^{\Si_*}, e^{\Si_*}_1, e^{\Si_*}_2)$ adapted to the $r$-foliation of $\Si_*$ are related on $\Si_*$ by the following  formulas
\bea\lab{eq:decompositiontangentspacerfoliationSigamastaronglobalframe:chap12}
 \bsplit
   e_4 &=\la\left(e_4^{\Si_*} + f^b  e_b^{\Si_*} +\frac 1 4 |f|^2  e_3^{\Si_*}\right),\\
  e_a&= \left(\de_a^b +\frac{1}{2}\fb_af^b\right) e_b^{\Si_*} +\frac 1 2  \fb_a  e_4^{\Si_*} +\left(\frac 1 2 f_a +\frac{1}{8}|f|^2\fb_a\right)e_3^{\Si_*},\quad a=1,2, \\
 e_3&= \la^{-1}\left(\left(1+\frac{1}{2}f\c\fb  +\frac{1}{16} |f|^2  |\fb|^2\right) e_3^{\Si_*} + \left(\fb^b+\frac 1 4 |\fb|^2f^b\right) e_b^{\Si_*}  + \frac 1 4 |\fb|^2 e_4^{\Si_*}\right),
 \end{split}
 \eea
for some scalar $\la$ and some 1-forms $(f, \fb)$, see (3.2.4) in \cite{KS:Kerr}.   

To compare the global from of $\MM$ with the  null frame $(e_3^{\Si_*}, e_4^{\Si_*}, e^{\Si_*}_1, e^{\Si_*}_2)$ adapted to the $r$-foliation of $\Si_*$, we assume that $(f, \fb, \la)$ appearing in \eqref{eq:decompositiontangentspacerfoliationSigamastaronglobalframe:chap12} satisfies on $\Si_*$
\bea\lab{eq:controloftansformationcoefficientsbetweenglobalframeandframeSigamstar}
|\dk^kf| \les \frac{1}{r},\qquad  |\dk^k\fb| \les \frac{1}{r},\qquad |\dk^k(\la-1)| \les \frac{1}{r}, \quad k\leq k_L,
\eea
see (3.2.5) and (3.2.6) in \cite{KS:Kerr}.


\subsubsection{Statement of Theorem M2}

  
We are now ready to restate Theorem M2 of \cite{KS:Kerr}, see also Theorem \ref{theoremM1:intro} in the Introduction. 
\begin{theorem}[Theorem M2 in \cite{KS:Kerr}]
\lab{thm:restatementofTheoremM2}
Assume that the global frame of $\MM$ satisfies the assumptions of section \ref{sec:assumptionsontheframe:Chapter12}. Assume in addition that the assumptions \eqref{eq:cominantconditionofronSigmastar}--\eqref{eq:propertiesoffirstorderderivativesofrandtauonSigmastar} and \eqref{eq:controloftansformationcoefficientsbetweenglobalframeandframeSigamstar} 
hold on $\Si_*$. Finally, assume that the control of the flux of $\qf$ provided by Theorem M1 holds, i.e.
\bea\lab{eq:assumptionsonqfforTheoremM2compatiblewithTheoremM1:chap12}
\int_{\Si_*(\geq\tau)}|\nab_3\dk^k\qf|^2 &\les& \ep_0^2\tau^{-2-2\dec}.
\eea
Then, $\aa$ satisfies the following estimate on $\Si_*$, for all $1\leq \tau\leq \tau_*$ and $k\leq k_L-7$, 
\beaa
\int_{\Si_*(\geq \tau)}|\dk^k\aa|^2 &\les& \ep_0^2\tau^{-2-2\dec}.
\eeaa
\end{theorem}

To prove Theorem M2, restated here as Theorem \ref{thm:restatementofTheoremM2}, one starts 
with Theorem \ref{Thm:Nondegenerate-Morawetz-psib} from which,  using the structure of the  error term $\widetilde{N}_\err$ in \eqref{eq:N_err-Ab}, one   can only derive estimates  of the form  (see \eqref{eq:step1fluxSigamstarLiebTAbchap12:1} below)    for $\de\le p\le 1-\de$,
\beaa
\BEF_p^s[\undpsi, \Ab](\tau_1, \tau_2) \les  E_p^s[\undpsi, \Ab](\tau_1) +\ep_0^2\tau_1^{-2-3\dec+p+\de}.
\eeaa
These estimates can be improved  by replacing $(\undpsi, \Ab)$ with their $\Lieb_\T$  derivatives. This fact is crucial in the proof of Theorem M2. We proceed as follows:
\begin{enumerate}
\item In section \ref{sec:decayfluxLiebT2aaandLieb2TundpsionMM:chap12}, we obtain decay estimates for $\Lieb_\T^2\aa$ on $\MM$ by relying on Theorem \ref{Thm:Nondegenerate-Morawetz-psib}. 

\item In section \ref{sec:decayfluxLiebT2aaonSigmastar:chap12:bis}, we deduce from the decay estimates of $\Lieb_\T^2\aa$ on $\MM$ of section \ref{sec:decayfluxLiebT2aaandLieb2TundpsionMM:chap12} a decay estimates for the flux of $\Lieb_\T^2\aa$ on $\Si_*$. 

\item In section \ref{sec:tukolskystarobinskytypeidentityqfSigmastart:chap12}, we derive an identity on $\Si_*$ involving $\qf$ and $\aa$.

\item Finally, in section \ref{sec:endoftheproofofTheoremM2:chap12}, we rely on the control of $\Lieb_\T^2\aa$ on $\Si_*$ and the control of $\qf$ provided by Theorem M1, which together with the above mentioned identity involving $\aa$ and $\qf$ yields an elliptic equation for $\aa$ along $\Si_*$. We then rely on this  elliptic  equation for $\aa$ along $\Si_*$ to prove Theorem M2.  
\end{enumerate}


\subsection{Decay estimate for $\Lieb_\T^2\Ab$}
\lab{sec:decayfluxLiebT2aaandLieb2TundpsionMM:chap12}


The goal of this section is to prove the following improved  decay estimate for $\Lieb_\T^2\Ab$.
\begin{proposition}\lab{prop:decayofpr2tqfb:onMMfirst}
The following decay estimate holds for $\Lieb_\T^2\Ab$, for $s\leq k_L -9$,
\beaa
B_{2-\de}^s[\Lieb_\T^2\Ab] (\tau_1, \tau_*)&\les& \ep_0^2 \tau_1^{-2-2\dec}.
 \eeaa
 Furthermore, there exists  a sequence of times $\tau^{(j)}$ such that, for $s\leq k_L-9$, 
\beaa
E^s_{2}[\Lieb_\T^2\Ab](\tau^{(j)}) &\les& \ep_0^2(\tau^{(j)})^{-2-2\dec}, \qquad \tau^{(j)}\sim 2^j.
\eeaa
\end{proposition}

 We will rely on the following lemma.
\begin{lemma}
\lab{lemma:estimateforN_{err}-chap12}
We have  for all  $\de\leq p\leq 1-\de$ and $s\leq k_L$,
\bea
   \lab{eq:estimateforN_{err}-chap12}
       \NN^s_p[\undpsi,  \widetilde{N}_\err](\tau_1, \tau_2)
      &\les& \ep_0\tau_1^{-1-\frac{3}{2}\dec+\frac{p+\de}{2}}\Big(\BEF^s_p[\undpsi](\tau_1, \tau_2)\Big)^{\frac{1}{2}}.
\eea      
\end{lemma}

\begin{proof}
Recall that $\Nt_\err= r^2\dk^{\le 2 }\big(\Ga_b\c(A, B)\big)$, see \eqref{eq:N_err-Ab}. In the particular  case $s=0$, we use as in the  proof of Lemma \ref{lemma:DecayfortheNterm-undpsi}, see section \ref{section:Prooflemma:DecayfortheNterm-undpsi},  for $\de\leq p\leq 2-\de$,
 \beaa
      \NN_p[\undpsi,  N](\tau_1, \tau_2) 
      &\les& \Big(\BEF_p[\undpsi](\tau_1, \tau_2)\Big)^{\frac{1}{2}}\left(\int_{\tau_1}^{\tau_2}\|N\|_{L^2(\Si_{trap}(\tau))}+\left(\int_{\MM(\tau_1, \tau_2)}r^{p+1}|N|^2\right)^{\frac{1}{2}}\right).
 \eeaa    
 In view of $\Nt_\err= r^2\dk^{\le 2 }\big(\Ga_b\c(A, B)\big)$, we immediately have, using \eqref{eq:assumptionsonMextforpartIIchap12}, 
 \beaa
 \int_{\tau_1}^{\tau_2}\|N\|_{L^2(\Si_{trap}(\tau))} &\les& \ep^4\tau_1^{-2-3\dec}\les \ep_0^2\tau_1^{-2-3\dec}.
 \eeaa
 Also, we have, for $\de\leq p\leq 1-\de$, 
 \beaa
 \int_{\MM(\tau_1, \tau_2)}r^{p+1}|N|^2 &\les& \ep^2\int_{\MM(\tau_1, \tau_2)}\frac{r^{p+3}}{\tau^{2+2\dec}}|\dk^{\leq 2}(A, B)|^2\\
 &\les& \frac{\ep^2}{\tau_1^{1+2\dec}}\int_{\MM(\tau_1, \tau_2)}\frac{r^{p+3}}{\tau}|\dk^{\leq 2}(A, B)|^2.
 \eeaa
 Next, interpolating between \eqref{eq:assumptionsonMextforpartIIchap12} and \eqref{eq:assumptionsonMextforpartIIchap12:moredecayinrAB} to control $(A, B)$, we infer, for $\de\leq p\leq 1-\de$, 
 \beaa
 \int_{\MM(\tau_1, \tau_2)}r^{p+1}|N|^2 &\les& \frac{\ep^4}{\tau_1^{1+2\dec}}\int_{\MM(\tau_1, \tau_2)}\frac{r^{p+3}}{\tau}\left(\frac{1}{r^6\tau^{1+2\dec}}\right)^{1-\frac{p+\de}{1+2\dec}}\left(\frac{1}{r^{7+2\dec}}\right)^{\frac{p+\de}{1+2\dec}}\\
 &\les& \frac{\ep^4}{\tau_1^{1+2\dec}}\left(\int \frac{dr}{r^{1+\de}}\right)\left(\int_{\tau\geq \tau_1}\frac{d\tau}{\tau^{2+2\dec-(p+\de)}}\right)
 \eeaa
 and hence, for $\de\leq p\leq 1-\de$, 
 \beaa
 \int_{\MM(\tau_1, \tau_2)}r^{p+1}|N|^2 &\les& \ep_0^2\tau_1^{-2-3\dec -(p+\de)}.
 \eeaa 
In view of the above, we deduce, for $\de\leq p\leq 1-\de$, 
\beaa
\NN_p[\undpsi,  \widetilde{N}_\err](\tau_1, \tau_2)
      &\les& \ep_0\tau_1^{-1-\frac{3}{2}\dec+\frac{p+\de}{2}}\Big(\BEF_p[\undpsi](\tau_1, \tau_2)\Big)^{\frac{1}{2}},
\eeaa 
which proves  \eqref{eq:estimateforN_{err}-chap12} in the particular case $s=0$. The general case can be shown  in the same manner which concludes the proof of Lemma \ref{lemma:estimateforN_{err}-chap12}.
\end{proof}

\begin{proof}[Proof of Proposition \ref{prop:decayofpr2tqfb:onMMfirst}]
We proceed in steps as follows.

{\bf Step 1.}  We derive the estimate
for  $s\leq k_L$,
\bea\lab{eq:step1fluxSigamstarLiebTAbchap12:1} 
\BEF_p^s[\undpsi, \Ab](\tau_1, \tau_2) \les  E_p^s[\undpsi, \Ab](\tau_1) +\ep_0^2\tau_1^{-2-3\dec+p+\de}, \qquad  \de\leq p\leq 1-\de.
\eea
This is an  immediate consequence  of    Theorem \ref{Thm:Nondegenerate-Morawetz-psib} and Lemma \ref{lemma:estimateforN_{err}-chap12}.

{\bf Step 2.} Next, applying the standard mean value procedure, see for instance the statement and proof of Theorem 5.21  in \cite{KS}, we infer from  \eqref{eq:step1fluxSigamstarLiebTAbchap12:1}, for $s\leq k_L-1$, 
 \bea\lab{eq:step1fluxSigamstarLiebTAbchap12:2}
\BEF_\de^s[\undpsi, \Ab](\tau_1, \tau_*) \les  \ep_0^2\tau_1^{-1+2\de}.
 \eea

{\bf Step 3.} Next, note that we have\footnote{Note that  the $B_\de$ norms  are stronger than the $\Morr$ norms, see definition  of both   in  section \ref{subsection:basicnormsforpsi}. }
\beaa
\int_{\MM(\tau_1, \tau_*)}r^{-1-\de}|\Lieb_\T\dk^{\leq s}\qfb|^2 &\les& \int_{\MM(\tau_1, \tau_*)}\Big(r^{-1-\de}|\nab_3\dk^{\leq s}\qfb|^2+r^{-3-\de}|\dk^{\leq s+1}\qfb|^2\Big)\\
&\les& B^{s+1}_\de[\undpsi](\tau_1, \tau_*)
\eeaa 
which together with \eqref{eq:step1fluxSigamstarLiebTAbchap12:2} implies, for $s\leq k_L-2$,
\beaa
\int_{\MM(\tau_1, \tau_*)}r^{-1-\de}|\Lieb_\T\dk^{\leq s}\qfb|^2 &\les&  \ep_0^2\tau_1^{-1+2\de}.
\eeaa
In view of the definition of $B^s_p[\Lieb_\T\undpsi]$,    see section  \ref{subsection:basicnormsforpsi},      we infer for $s\leq k_L-3$ and for all $\de\leq p\leq 2-\de$,   
\bea\lab{eq:step1fluxSigamstarLiebTAbchap12:3}
B^s_p[\Lieb_\T\undpsi](\tau_1, \tau_*) &\les& \ep_0^2\tau_1^{-1+2\de}.
\eea

{\bf Step 4.} Next, we commute the wave equation \eqref{eq:Gen.RW-pert-qfb} for $\undpsi$ and the system of transport equations \eqref{eq:thesystemof2transporteuqaitonsAbPsibundpsiactuallyused} with $\Lieb_\T$ and obtain 
 \beaa
\squared_2 \Lieb_\T\underline{\psi} -V_0\Lieb_\T\underline{\psi}= \frac{4 a\cos\th}{|q|^2}\dual \nab_T\Lieb_\T  \underline{\psi}+N_{\Lieb_\T}
\eeaa
and 
\beaa
\nab_4(r\Lieb_\T\Psib)=\frac{q}{r\ov{q}}\Lieb_\T\qfb+F_{1,\Lieb_\T},\qquad \nab_4\left( \frac{q^4}{r^3}\Lieb_\T\Ab\right) = r^{-1}\Lieb_\T\Psib+F_{2,\Lieb_\T},
\eeaa  
where 
\beaa
N_{\Lieb_\T} &=& \Lieb_\T N+[\squared_2, \Lieb_\T]\underline{\psi}, \\ 
F_{1, \Lieb_\T} &=& [\nab_4(r\cdot), \Lieb_\T]\Psib+r\dk^{\leq 2}(\Ga_g\c\Ga_b), \\ 
F_{2, \Lieb_\T}  &=& \left[\nab_4\left(\frac{q^4}{r^3}\cdot\right), \Lieb_\T\right]\Ab+O(r^{-2})\T(r)\Psib+r\dk^{\leq 2}(\Ga_g\c\Ga_b).
\eeaa
Now, in view of Corollary \ref{cor:commutator-Lied-squared}, Lemma \ref{lemma:basicpropertiesLiebTfasdiuhakdisug:chap9} and the above definition of $N_{\Lieb_\T}$ and $F_{\Lieb_\T}$, we have 
\beaa
N_{\Lieb_\T} &=& \Lieb_\T N+\dk(\Ga_g\c\dk\underline{\psi})+\Ga_b\c\squared_2\undpsi,\\
&=&  \Lieb_\T N+\dk^{\leq 1}(\Ga_g\c\underline{\psi})+\Ga_b\c N,\\ 
F_{1, \Lieb_\T} &=& r\dk^{\leq 2}(\Ga_g\c\Ga_b), \\ 
F_{2, \Lieb_\T}  &=& r\dk^{\leq 2}(\Ga_g\c\Ga_b).
\eeaa
The above system for $(\Lieb_\T\undpsi, \Lieb_\T\Ab)$ has error terms of the same type as before, so that the estimates established for $(\undpsi, \Ab)$ hold for $(\Lieb_\T\undpsi, \Lieb_\T\Ab)$ as well. In particular, we have the  following analog of \eqref{eqtheorem:gRW1-p} for   $s\le k_L-1$ and  for all $\de\leq p \leq 2 -\de$,
 \bea
       \lab{eqtheorem:gRW1-p:Liebversion}
 \nn       \BEF_p^s[\Lieb_\T\undpsi, \Lieb_\T\Ab](\tau_1, \tau_2) &\les&  E_p^s[\Lieb_\T\undpsi, \Lieb_\T\Ab](\tau_1)+\NN_p^s[\Lieb_\T\undpsi,  \Lieb_\T\widetilde{N}_\err](\tau_1, \tau_2)\\
        && +\ep_0^2\tau_1^{-2-2\dec},
       \eea
and the following  analog of \eqref{eq:transportAb} for   $s\le k_L-1$ and  for all $\de\leq p \leq 2 -\de$,
 \bea
 \lab{eq:transportAb:Liebversion}
 B_p^s[\Lieb_\T\Ab] (\tau_1, \tau_2)&\les&   B^s_{\de}[\Lieb_\T\undpsi](\tau_1,\tau_2) 
 + E^s_p[\Lieb_\T\Ab](\tau_1) +\ep_0^2 \tau_1^{-2-2\dec}.
 \eea

{\bf Step 5.} Next, we deduce from \eqref{eq:step1fluxSigamstarLiebTAbchap12:3} and \eqref{eq:transportAb:Liebversion},  for $s\leq k_L-3$ and for all $\de\leq p\leq 2-\de$, 
 \beaa
  B_p^s[\Lieb_\T\Ab] (\tau_1, \tau_2)&\les&     E^s_p[\Lieb_\T\Ab](\tau_1) +\ep_0^2 \tau_1^{-1+2\de}.
 \eeaa
Now, in view of  the Definition  \ref{Definition:NormsBEF-Ab}  of  the norms $B_p^s[\Ab]$ and $E^s_p[\Lieb_\T\Ab]$,  we have for all $\de\leq p\leq 2-\de$
\beaa
\int_{\tau_1}^{\tau_2}E^s_p[\Lieb_\T\Ab](\tau) &\les& B_{p-1}^s[\Lieb_\T\Ab] (\tau_1, \tau_2)
\eeaa 
and hence,  for $s\leq k_L-3$,  
\beaa
\int_{\tau_1}^{\tau_2}E^s_{1+\de}[\Lieb_\T\Ab](\tau) &\les&   B_{\de}^s[\Lieb_\T\Ab] (\tau_1, \tau_2)\les   E^s_\de[\Lieb_\T\Ab](\tau_1) +\ep_0^2 \tau_1^{-1+2\de}.
\eeaa 
Together with \eqref{eq:step1fluxSigamstarLiebTAbchap12:2}, we infer,  for $s\leq k_L-3$ and for all $\de\leq p\leq 1+\de$, 
\beaa
\int_{\tau_1}^{\tau_2}E^s_p[\Lieb_\T\Ab](\tau) &\les&   \ep_0^2 \tau_1^{-1+2\de}.
\eeaa
Thus, since 
\beaa
\int_{\tau_1}^{\tau_2}E^s_p[\Lieb_\T\undpsi](\tau_1) &\les& B^{s+1}_{p+1}[\Lieb_\T\undpsi](\tau_1, \tau_2), 
\eeaa
we deduce, in view of \eqref{eq:step1fluxSigamstarLiebTAbchap12:3}, for $s\leq k_L-4$ and for all $\de\leq p\leq 1-\de$, 
\beaa
\int_{\tau_1}^{\tau_2}E^s_p[\Lieb_\T\undpsi, \Lieb_\T\Ab](\tau) &\les&   \ep_0^2 \tau_1^{-1+2\de}.
\eeaa
We infer the existence of a sequence of times $\tau^{(j)}$ such that, for $s\leq k_L-4$, 
\bea\lab{eq:step1fluxSigamstarLiebTAbchap12:4}
E^s_{1-\de}[\Lieb_\T\undpsi, \Lieb_\T\Ab](\tau^{(j)}) &\les& \ep_0^2(\tau^{(j)})^{-1+2\de}, \qquad \tau^{(j)}\sim 2^j.
\eea

{\bf Step 6.} Next, arguing as in Step 1,  using the structure of $\Lieb_\T \Nt_{\err}$ and Lemma      \ref{lemma:estimateforN_{err}-chap12}   to deduce, for $\de\leq p\leq 1-\de$ and $s\leq k_L-1$,
\beaa
      \NN^s_p[\Lieb_\T\undpsi,  \Lieb_\T \widetilde{N}_\err](\tau_1, \tau_2)   &\les& \ep_0\tau_1^{-1-\frac{3}{2}\dec+\frac{p+\de}{2}}\Big(\BEF^s_p[\undpsi](\tau_1, \tau_2)\Big)^{\frac{1}{2}}.
\eeaa  
Together with \eqref{eqtheorem:gRW1-p:Liebversion} and  the estimate \eqref{eq:step1fluxSigamstarLiebTAbchap12:1}  for $\BEF^s_p[\undpsi](\tau_1, \tau_2)$, we deduce, for   $s\le k_L-1$ and  for all $\de\leq p \leq 1 -\de$, 
 \bea\lab{eqtheorem:gRW1-p:Liebversion:afterestimatewidetildeNerrterm}
      \BEF_p^s[\Lieb_\T\undpsi, \Lieb_\T\Ab](\tau_1, \tau_2) \les  E_p^s[\Lieb_\T\undpsi, \Lieb_\T\Ab](\tau_1)+\ep_0^2\tau_1^{-2-3\dec+p+\de}.
 \eea
We then apply again  the standard mean value procedure starting with 
\eqref{eqtheorem:gRW1-p:Liebversion:afterestimatewidetildeNerrterm}  and making use 
of   \eqref{eq:step1fluxSigamstarLiebTAbchap12:4}.  We thus infer, for $s\leq k_L-5$, 
 \bea\lab{eq:step6fluxSigamstarLiebTAbchap12:666}
\BEF_\de^s[\Lieb_\T\undpsi, \Lieb_\T\Ab](\tau_1, \tau_*) \les  \ep_0^2\tau_1^{-2+4\de}.
 \eea

{\bf Step 7.} We now run the procedure of Step 3 to Step 6 again, with $(\Lieb_\T\undpsi, \Lieb_\T\Ab)$ replaced by $(\Lieb_\T^2\undpsi, \Lieb_\T^2\Ab)$. More precisely:
\begin{enumerate}
\item Starting from \eqref{eq:step6fluxSigamstarLiebTAbchap12:666},      as analogous to 
\eqref{eq:step1fluxSigamstarLiebTAbchap12:2}, 
and proceeding as in Step 3, we obtain, for $s\leq k_L-7$ and for all $\de\leq p\leq 2-\de$,   
\beaa
B^s_p[\Lieb^2_\T\undpsi](\tau_1, \tau_*) &\les& \ep_0^2\tau_1^{-2+4\de}.
\eeaa

\item Next, commuting the system for $(\Lieb_\T\undpsi, \Lieb_\T\aa)$ with  another $\Lieb_\T$, and proceeding as in Step 4, we obtain the analog of \eqref{eqtheorem:gRW1-p:Liebversion} and \eqref{eq:transportAb:Liebversion} for $(\Lieb_\T^2\undpsi, \Lieb_\T^2\aa)$, i.e. for   $s\le k_L-2$ and  for all $\de\leq p \leq 2 -\de$, we have
 \beaa
 \nn       \BEF_p^s[\Lieb_\T^2\undpsi, \Lieb_\T^2\Ab](\tau_1, \tau_2) &\les&  E_p^s[\Lieb_\T^2\undpsi, \Lieb_\T^2\Ab](\tau_1)+\NN_p^s[\Lieb_\T^2\undpsi,  \Lieb_\T^2\widetilde{N}_\err ](\tau_1, \tau_2)\\
        && +\ep_0^2\tau_1^{-2-2\dec},
\eeaa
and 
 \beaa
\nn B_p^s[\Lieb_\T^2\Ab] (\tau_1, \tau_2)&\les&   B^s_{\de}[\Lieb_\T^2\undpsi](\tau_1,\tau_2) 
 + E^s_p[\Lieb_\T^2\Ab](\tau_1)  +\ep_0^2 \tau_1^{-2-2\dec}.
 \eeaa

\item Next, proceeding as in Step 5, we obtain the analog of \eqref{eq:step1fluxSigamstarLiebTAbchap12:4}, i.e. 
we infer the existence of a sequence of times $\tau_1^{(j)}$ such that, for $s\leq k_L-8$, 
\beaa
E^s_{1-\de}[\Lieb_\T^2\undpsi, \Lieb_\T^2\Ab](\tau_1^{(j)}) &\les& \ep_0^2(\tau_1^{(j)})^{-2+4\de}, \qquad \tau_1^{(j)}\sim 2^j.
\eeaa

\item Next, proceeding as in Step 6, we obtain the analog of \eqref{eq:step6fluxSigamstarLiebTAbchap12:666}, for $s\leq k_L-9$, 
 \beaa
\BEF_\de^s[\Lieb_\T^2\undpsi, \Lieb_\T^2\Ab](\tau_1, \tau_*) \les  \ep_0^2\tau_1^{-2-2\dec}+\ep_0^2\tau_1^{-2-3\dec+2\de}+\ep_0^2\tau_1^{-3+6\de},
 \eeaa
where:
\begin{itemize}
\item the first term in the RHS comes from the contribution of all nonlinear terms except $\Lieb_\T^2\Nt_\err$, i.e. the terms of type $\dk^{\leq 5}(\Ga_b\c\Ga_g)$, 

\item the second term in the RHS comes from the contribution of $\Lieb^2_\T \Nt_\err$, 

\item the last term comes from the mean value argument. 
\end{itemize}
Choosing $\de>0$ such that $2\de\leq \dec$, we infer, for $s\leq k_L-9$ and $\dec$ small enough, 
 \beaa
\BEF_\de^s[\Lieb_\T^2\undpsi, \Lieb_\T^2\Ab](\tau_1, \tau_*) \les  \ep_0^2\tau_1^{-2-2\dec}.
 \eeaa
\end{enumerate} 

{\bf Step 8.} Recall that we have obtained in Step 7, for all $\de\leq p \leq 2 -\de$, 
 \beaa
\nn B_p^s[\Lieb_\T^2\Ab] (\tau_1, \tau_2)&\les&   B^s_{\de}[\Lieb_\T^2\undpsi](\tau_1,\tau_2) 
 + E^s_p[\Lieb_\T^2\Ab](\tau_1)+\ep_0^2 \tau_1^{-2-2\dec}.
 \eeaa
Together with the final estimate of Step 7 for $\undpsi$, this yields, for all $\de\leq p \leq 2 -\de$ and $s\leq k_L-9$, 
 \bea\lab{eq:transportAb:Lieb2Tversion}
 B_p^s[\Lieb_\T^2\Ab] (\tau_1, \tau_*)&\les& E^s_p[\Lieb_\T^2\Ab](\tau_1)+\ep_0^2 \tau_1^{-2-2\dec}.
 \eea
As in Step 5, we use the fact  that we have for all $\de\leq p\leq 2-\de$
\beaa
\int_{\tau_1}^{\tau_2}E^s_p[\Lieb_\T^2\Ab](\tau) &\les& B_{p-1}^s[\Lieb_\T^2\Ab] (\tau_1, \tau_2).
\eeaa  
Together with the  final estimate of Step 7 for $\BEF_\de^s[ \Lieb_\T^2\Ab](\tau_1, \tau_*)$, we infer the existence of a sequence of times $\tau_1^{(j)}$ such that, for $s\leq k_L-9$, 
\beaa
E^s_{1+\de}[\Lieb_\T^2\Ab](\tau_1^{(j)}) &\les& \ep_0^2(\tau_1^{(j)})^{-3-2\dec}, \qquad \tau_1^{(j)}\sim 2^j.
\eeaa
Plugging in \eqref{eq:transportAb:Lieb2Tversion}, we infer, for all $s\leq k_L-9$,
\bea\lab{eq:transportAb:Lieb2Tversion:concsequncce}
 B_{1+\de}^s[\Lieb_\T^2\Ab] (\tau_1, \tau_*)&\les& \ep_0^2 \tau_1^{-2-2\dec}.
 \eea
Running again the same argument, we deduce the existence of a sequence of times $\tau_2^{(j)}$ such that, for $s\leq k_L-9$, 
\beaa
E^s_{2}[\Lieb_\T^2\Ab](\tau_2^{(j)}) &\les& \ep_0^2(\tau_2^{(j)})^{-3-2\dec}, \qquad \tau_2^{(j)}\sim 2^j.
\eeaa
Plugging in \eqref{eq:transportAb:Lieb2Tversion}, we finally obtain, for all $s\leq k_L-9$, 
\beaa
 B_{2-\de}^s[\Lieb_\T^2\Ab] (\tau_1, \tau_*)&\les& \ep_0^2 \tau_1^{-2-2\dec}.
 \eeaa
as stated. This concludes the proof of Proposition \ref{prop:decayofpr2tqfb:onMMfirst}. 
\end{proof}


\subsection{Decay of the flux on $\Si_*$ of $\Lieb_\T^2\Ab$}
\lab{sec:decayfluxLiebT2aaonSigmastar:chap12:bis}


The goal of this section is to prove the following decay estimate for the flux on $\Si_*$ of $\Lieb_\T^2\Ab$.
\begin{proposition}\lab{prop:decayofpr2tqfb}
The following decay estimate holds for $\Lieb_\T^2\Ab$, for $s\leq k_L -10$,
\bea
F^s_{\Si_*}[\Lieb_\T^2\Ab](\tau, \tau_*) &\les& \ep_0^2\tau^{-2-2\dec},
\eea
where $F^s_{\Si_*}$ denotes the flux on $\Si_*$. 
\end{proposition}

To prove Proposition \ref{prop:decayofpr2tqfb}, we first derive the following extension of Lemma \ref{lemma:general-transport-estimate-qfb}.

  \begin{lemma}\lab{lemma:general-transport-estimate-qfb:endpoint}
 Suppose $\Phi_1, \Phi_2 \in \sk_2(\CCC)$ with signature $s\leq -1$ satisfy the differential relation 
 \beaa
 \nabc_4\Phi_1=\Phi_2,
 \eeaa
Then, we have
 \bea\label{eq:general-integrated-estimate-e4:endpoint}
 \bsplit
& \int_{\Si(\tau_2)} r^{-4}|\Phi_1|^2+  \int_{\AA\cup \Si_*(\tau_1,\tau_2)} r^{-2}|\Phi_1|^2  \\
  &\les \int_{\MM(\tau_1, \tau_2)}|\Phi_2|^2+ \int_{\Si(\tau_1)} r^{-4}|\Phi_1|^2 +\int_{\MM(\tau_1, \tau_2)}  r^{-4} |\Phi_1|^2\\
  & +\int_{\MM(\tau_1, \tau_2)}  r^{-2}|\Ga_g||\Phi_1|^2.
  \end{split}
\eea 
\end{lemma}

\begin{proof}
Recall from the proof of Lemma \ref{lemma:general-transport-estimate-qfb} that we have
\beaa
\Div(  |q|^{p-2}|\Phi_1|^2 e_4) &=&2|q|^{p-2} \Re\big( \Phi_2 \c \ov{\Phi_1}\big) + p r\De |q|^{p-6} |\Phi_1|^2 -2(2s+1)\om |q|^{p-2} |\Phi_1|^2\\
&& +|q|^{p-2}\Ga_g|\Phi_1|^2.
\eeaa
We choose $p=0$ which yields
\beaa
\Div(  |q|^{-2}|\Phi_1|^2 e_4) &=&2|q|^{-2} \Re\big( \Phi_2 \c \ov{\Phi_1}\big)  -2(2s+1)\om |q|^{-2} |\Phi_1|^2 +|q|^{-2}\Ga_g|\Phi_1|^2.
\eeaa
Since $-\om\gtrsim \frac{m}{r^2}+\Ga_g$   and  $s\leq -1$, we deduce
\beaa
\Div(  |q|^{-2}|\Phi_1|^2 e_4) &\leq & 2|q|^{-2} \Re\big( \Phi_2 \c \ov{\Phi_1}\big)  +|q|^{-2}\Ga_g|\Phi_1|^2\\
&\leq & |\Phi_2|^2+ |q|^{-4}|\Phi_1|^2 +|q|^{-2}\Ga_g|\Phi_1|^2.
\eeaa 
We now apply the divergence theorem as in the proof of Lemma \ref{lemma:general-transport-estimate-qfb} to obtain \eqref{eq:general-integrated-estimate-e4:endpoint}. This concludes the proof of Lemma \ref{lemma:general-transport-estimate-qfb:endpoint}.
\end{proof}

We have the following higher derivatives version of Lemma \ref{lemma:general-transport-estimate-qfb:endpoint}.
\begin{corollary}\lab{cor:general-transport-estimate-qfb:endpoint:higherorder}
 Suppose $\Phi_1, \Phi_2 \in \sk_2(\CCC)$ with signature $s\leq -1$ satisfy the differential relation 
 \bea\label{eq:relation-nab3Phi1Phi2}
 \nabc_4\Phi_1=\Phi_2,
 \eea
Then, for any $0\leq s\leq k_L$, we have
 \bea\label{eq:general-integrated-estimate-e4:endpoint:higherorder}
 \bsplit
& \int_{\Si(\tau_2)} r^{-4}|\dk^s\Phi_1|^2+  \int_{\AA\cup \Si_*(\tau_1,\tau_2)} r^{-2}|\dk^s\Phi_1|^2  \\
  &\les \int_{\MM(\tau_1, \tau_2)}|\dk^s\Phi_2|^2+ \int_{\Si(\tau_1)} r^{-4}|\dk^s\Phi_1|^2\\
  &+\int_{\MM(\tau_1, \tau_2)}  r^{-4} |\dk^s\Phi_1|^2+\int_{\MM(\tau_1, \tau_2)}  r^{-2}|\dk^{\leq s}(\Ga_g\c\Phi_1)||\dk^s\Phi_1|.
  \end{split}
\eea 
\end{corollary}

\begin{proof}
The case $s=0$ holds true by Lemma \ref{lemma:general-transport-estimate-qfb:endpoint}. We then argue by iteration on $s$, and pass from $s$ to $s+1$ by following the procedure outlined at the end of section \ref{sec:endoftheproofofproposition:MaiTransportAb-steps}. 
\end{proof}

We are now ready to prove Proposition \ref{prop:decayofpr2tqfb}. 
\begin{proof}[Proof of Proposition \ref{prop:decayofpr2tqfb}.] 
Recall the second equation of \eqref{eq:thesystemof2transporteuqaitonsAbPsibundpsiactuallyused}, i.e.  
\beaa
\nabc_4\left( \frac{q^4}{r^3}\Ab\right)= \frac{1}{r}\Psib+r\Ga_g\c\Ga_b,
\eeaa
where $\Psib$ is given by
\beaa
\und{\Psi} &=& \frac{q^4}{r^2}\left(\nabc_4+2\tr X -\frac{3|\tr X|^2}{2\trch}\right)\Ab.
\eeaa
We commute this equation with $\Lieb_\T^2\Ab$ and obtain, as in the end of section \ref{sec:endoftheproofofproposition:MaiTransportAb-steps},
\beaa
\nabc_4\left( \frac{q^4}{r^3}\Lieb_\T^2\Ab\right)= \frac{1}{r}\Lieb_\T^2\Psib+r\dk^{\leq 2}(\Ga_g\c\Ga_b).
\eeaa
We apply Corollary \ref{cor:general-transport-estimate-qfb:endpoint:higherorder} to this transport equation, i.e. we choose 
\beaa
\Phi_1=\frac{q^4}{r^3}\Lieb_\T^2\Ab, \qquad \Phi_2=\frac{1}{r}\Lieb_\T^2\Psib+r\dk^{\leq 2}(\Ga_g\c\Ga_b).
\eeaa
We infer, for $0\leq s\leq k_L-2$
\beaa
 \bsplit
& \int_{\Si(\tau_2)} r^{-2}|\dk^s\Lieb_\T^2\Ab|^2+  \int_{\AA\cup \Si_*(\tau_1,\tau_2)}|\dk^s\Lieb_\T^2\Ab|^2  \\
  &\les \int_{\MM(\tau_1, \tau_2)} r^{-2}|\dk^s\Lieb_\T^2\Psib|^2+ \int_{\Si(\tau_1)} r^{-2}|\dk^s\Lieb_\T^2\Ab|^2\\
  &+\int_{\MM(\tau_1, \tau_2)}  r^{-2} |\dk^s\Lieb_\T^2\Ab|^2+\int_{\MM(\tau_1, \tau_2)}|\dk^{\leq s}(\Ga_g\c\Ga_b)||\dk^s\Ga_b|.
  \end{split}
\eeaa 
Using the definition of $\Psib$, we infer
\beaa
 \bsplit
& \int_{\Si(\tau_2)} r^{-2}|\dk^s\Lieb_\T^2\Ab|^2+  \int_{\AA\cup \Si_*(\tau_1,\tau_2)}|\dk^s\Lieb_\T^2\Ab|^2  \\
  &\les \int_{\MM(\tau_1, \tau_2)}|\dk^s\Lieb_\T^2\nab_4(r\Ab)|^2+ \int_{\Si(\tau_1)} r^{-2}|\dk^s\Lieb_\T^2\Ab|^2\\
  &+\int_{\MM(\tau_1, \tau_2)}  r^{-2} |\dk^s\Lieb_\T^2\Ab|^2+\int_{\MM(\tau_1, \tau_2)}|\dk^{\leq s}(\Ga_g\c\Ga_b)||\dk^s\Ga_b|.
  \end{split}
\eeaa 
In view of the control of $\Ga_g$ and $\Ga_b$, as well as the definition of $B^s_1[\Lieb_\T^2\Ab](\tau_1, \tau_2)$ and  $E^s_2[\Lieb_\T^2\Ab](\tau)$, we deduce
\beaa
 \bsplit
& \int_{\Si(\tau_2)} r^{-2}|\dk^s\Lieb_\T^2\Ab|^2+  \int_{\AA\cup \Si_*(\tau_1,\tau_2)}|\dk^s\Lieb_\T^2\Ab|^2  \\
  &\les B^s_1[\Lieb_\T^2\Ab](\tau_1, \tau_2)+E^s_2[\Lieb_\T^2\Ab](\tau_1)+\ep_0^2\tau_1^{-2-2\dec}.
  \end{split}
\eeaa 

Finally, recall from Proposition \ref{prop:decayofpr2tqfb:onMMfirst} that we have, for $s\leq k_L -9$,
\beaa
B_{2-\de}^s[\Lieb_\T^2\Ab] (\tau_1, \tau_*)&\les& \ep_0^2 \tau_1^{-2-2\dec}.
 \eeaa
 and that there exists  a sequence of times $\tau^{(j)}$ such that, for $s\leq k_L-9$, 
\beaa
E^s_{2}[\Lieb_\T^2\Ab](\tau^{(j)}) &\les& \ep_0^2(\tau^{(j)})^{-2-2\dec}, \qquad \tau^{(j)}\sim 2^j.
\eeaa
Plugging in the above, we deduce,  for $s\leq k_L-9$, 
\beaa
 \bsplit
  \int_{\AA\cup \Si_*(\tau_1,\tau_*)}|\dk^s\Lieb_\T^2\Ab|^2   &\les \ep_0^2\tau_1^{-2-2\dec}.
  \end{split}
\eeaa 
In particular, we infer 
The following decay estimate holds for $\Lieb_\T^2\Ab$, for $s\leq k_L -10$,
\beaa
F^s_{\Si_*}[\Lieb_\T^2\Ab](\tau, \tau_*) &\les& \ep_0^2\tau^{-2-2\dec},
\eeaa
as stated. This concludes the proof of Proposition \ref{prop:decayofpr2tqfb}. 
\end{proof}


\subsection{An identity on $\Si_*$ involving $\qf$ and $\aa$}
\lab{sec:tukolskystarobinskytypeidentityqfSigmastart:chap12}


Recall from  Definition \ref{Definition:Define-qf} that the quantity $\qf$ is given by 
 \beaa
 \bsplit
 \qf&=q \ov{q}^{3} \left( \nabc_3\nabc_3 A + C_1  \nabc_3A + C_2   A\right),\\
 C_1&=2\trchb - 2\frac {\atrchb^2}{ \trchb}  -4 i \atrchb, \\
C_2  &= \frac 1 2 \trchb^2- 4\atrchb^2+\frac 3 2 \frac{\atrchb^4}{\trchb^2} +  i \left(-2\trchb\atrchb +4\frac{\atrchb^3}{\trchb}\right).
 \end{split}
\eeaa
We first derive the following identities  involving $\qf$ which are the analog of the ones of section 2.3.4 in \cite{KS}. 

\begin{proposition}\lab{prop:roughidentitybetweenqfe3qfe3e3qfandPcBbandAb:chap12}
We have
\bea\lab{eq:roughidentitybetweenqfandDhotDhotPcandXhbonSigmastar:chap12}
\qf &=& \frac{1}{2}r^4\DDc\hot\DDc\ov{\Pc}   -6m\widehat{\Xb}    +\dk^{\leq 2}\Ga_g +r^2\dk^{\leq 1}(\Ga_b\c\Ga_g),
\eea
\bea\lab{eq:roughidentitybetweene3rqfandDhotDDcBbandAbonSigmastar:chap12}
\nabc_3(r\qf) &=&  -\frac{1}{4}r^5\DDc\hot\DDc\DDc\c\ov{\Bb} +6mr\Ab    +\dk^{\leq 3}\Ga_b +r^3\dk^{\leq 2}(\Ga_b\c\Ga_g),
\eea
and 
\bea\lab{eq:roughidentitybetweene3r2e3rqfandDhotDDcDDcAbandnab3cAbonSigmastar:chap12}
\nn\nabc_3(r^2\nabc_3(r\qf)) &=&  \frac{1}{8}r^7\DDc\hot\DDc\DDc\c\DDc\c\ov{\Ab}  +6mr^3\nab_3\Ab  +r^2\dk^{\leq 4}\Ga_b \\
&& +r^5\dk^{\leq 3}(\Ga_b\c\Ga_g).
\eea
\end{proposition}

\begin{proof}
We start with the proof of \eqref{eq:roughidentitybetweenqfandDhotDhotPcandXhbonSigmastar:chap12}. Note that $C_1$ and $C_2$ satisfy 
 \beaa
 \bsplit
  C_1&=-\frac{4}{r} +O(r^{-2})+\Ga_g, \\
C_2  &= \frac{2}{r^2} +O(r^{-3})+r^{-1}\Ga_g,
 \end{split}
\eeaa
so that 
\beaa
\qf&=& q \ov{q}^{3} \left( \nabc_3\nabc_3 A + C_1  \nabc_3A + C_2   A\right)\\
&=& r^4 \nabc_3\nabc_3 A -4r^3\nabc_3A +2r^2A +O(r^3)\nabc_3\nabc_3A\\
&&+\Big(O(r^2)+r^4\Ga_g\Big)\nabc_3A+\Big(O(r)+r^3\Ga_g\Big)A.
\eeaa
Since $A\in r^{-1}\Ga_g$, and since $\nabc_3A\in r^{-2}\dk^{\leq 1}\Ga_g$ and $\nabc_3^2A\in r^{-3}\dk^{\leq 2}\Ga_g$ by Bianchi, we infer
\beaa
\qf &=& r^4 \nabc_3\nabc_3 A -4r^3\nabc_3A +2r^2A +\dk^{\leq 2}\Ga_g.
\eeaa

Next, we have in view of the Bianchi identity for $\nabc_3A$ of Proposition \ref{prop:bianchi:complex-conf} 
\beaa
\nabc_3A  &=& \frac 1 2 \DDc\hot B -\frac{1}{2}\tr\Xb A + 2 H   \hot B -3\ov{P}\Xh\\
&=& \frac 1 2 \DDc\hot B +\frac{1}{r} A+O(r^{-2})(A, B)+\frac{6m}{r^3}\Xh+r^{-4}\Ga_g
\eeaa
and hence
\beaa
\qf &=& r^4 \nabc_3\left(\frac 1 2 \DDc\hot B +\frac{1}{r} A+O(r^{-2})(A, B)+\frac{6m}{r^3}\Xh+r^{-4}\Ga_g\right)\\ && -4r^3\left(\frac 1 2 \DDc\hot B +\frac{1}{r} A+O(r^{-2})(A, B)+\frac{6m}{r^3}\Xh+r^{-4}\Ga_g\right) +2r^2A +\dk^{\leq 2}\Ga_g\\
&=& \frac{1}{2}r^4\nabc_3\DDc\hot B +r^3\nabc_3A +r^2A +6mr\nabc_3\Xh  -2r^3\DDc\hot B  -4r^2A \\
&&+2r^2A+\dk^{\leq 2}\Ga_g +r^2\dk^{\leq 1}(\Ga_b\c\Ga_g)\\
&=& \frac{1}{2}r^4\nabc_3\DDc\hot B +r^3\nabc_3A  +6mr\nabc_3\Xh  -2r^3\DDc\hot B  -r^2A +\dk^{\leq 2}\Ga_g\\
&& +r^2\dk^{\leq 1}(\Ga_b\c\Ga_g).
\eeaa
Using again the Bianchi identity for $\nabc_3A$, as well as the following consequence of the null structure equation for $\nabc_3\Xh$ of Proposition \ref{prop-nullstr:complex-conf},
\beaa
\nabc_3\Xh  &=& -\frac{1}{r}\widehat{\Xb}+r^{-1}\dk^{\leq 1}\Ga_g,
\eeaa
we infer
\beaa
\qf &=& \frac{1}{2}r^4\nabc_3\DDc\hot B +r^3\left(\frac 1 2 \DDc\hot B +\frac{1}{r} A+O(r^{-2})(A, B)+\frac{6m}{r^3}\Xh+r^{-4}\Ga_g\right)\\
&&  +6mr\left( -\frac{1}{r}\widehat{\Xb}+r^{-1}\dk^{\leq 1}\Ga_g\right)  -2r^3\DDc\hot B  -r^2A +\dk^{\leq 2}\Ga_g +r^2\dk^{\leq 1}(\Ga_b\c\Ga_g)\\
&=& \frac{1}{2}r^4\nabc_3\DDc\hot B -\frac{3}{2}r^3\DDc\hot B -6m\widehat{\Xb}    +\dk^{\leq 2}\Ga_g +r^2\dk^{\leq 1}(\Ga_b\c\Ga_g).
\eeaa
Together with the following commutator identity, see Lemma \ref{COMMUTATOR-NAB-C-3-DD-C-HOT},  
\beaa
 \, [\nabc_3, \DDc \hot ]B &=&- \frac 1 2 \tr \Xb \left( \DDc \hot B + 2H \hot B \right)  + H \hot \nabc_3 B+ r^{-1}\Ga_b  \c \dk^{\leq 1}B\\
 &=& \frac{1}{r}\DDc \hot B + r^{-4}\dk^{\leq 1}\Ga_g+r^{-2}\dk^{\leq 1}(\Ga_b\c\Ga_g)
\eeaa
where we used in particular the fact that $\Hc\in\Ga_g$ and $\nabc_3B\in r^{-2}\dk^{\leq 1}\Ga_g$. We deduce
\beaa
\qf &=& \frac{1}{2}r^4\DDc\hot\nabc_3 B+\frac{1}{2}r^4[\nabc_3,\DDc\hot] B -\frac{3}{2}r^3\DDc\hot B -6m\widehat{\Xb}    +\dk^{\leq 2}\Ga_g \\
&&+r^2\dk^{\leq 1}(\Ga_b\c\Ga_g)\\
&=& \frac{1}{2}r^4\DDc\hot\nabc_3 B - r^3\DDc\hot B -6m\widehat{\Xb}    +\dk^{\leq 2}\Ga_g +r^2\dk^{\leq 1}(\Ga_b\c\Ga_g).
\eeaa
In view of the following consequence of the Bianchi identity for $\nab_3B$ 
\beaa
\nabc_3B &=& \DDc\ov{P} +3\ov{P}H +\frac{2}{r} B +r^{-3}\dk^{\leq 1}\Ga_g +r^{-1}\Ga_b\c\Ga_g\\
&=& \DDc\ov{\Pc}  +\frac{2}{r} B +r^{-3}\dk^{\leq 1}\Ga_g +r^{-1}\Ga_b\c\Ga_g,
\eeaa
we infer
\beaa
\qf &=& \frac{1}{2}r^4\DDc\hot\left(\DDc\ov{\Pc}  +\frac{2}{r} B \right) - r^3\DDc\hot B -6m\widehat{\Xb}    +\dk^{\leq 2}\Ga_g +r^2\dk^{\leq 1}(\Ga_b\c\Ga_g)\\
&=& \frac{1}{2}r^4\DDc\hot\DDc\ov{\Pc}   -6m\widehat{\Xb}    +\dk^{\leq 2}\Ga_g +r^2\dk^{\leq 1}(\Ga_b\c\Ga_g)
\eeaa
which concludes the proof of \eqref{eq:roughidentitybetweenqfandDhotDhotPcandXhbonSigmastar:chap12}.

Next, we prove \eqref{eq:roughidentitybetweene3rqfandDhotDDcBbandAbonSigmastar:chap12}. First, we multiply by $r$ and differentiate \eqref{eq:roughidentitybetweenqfandDhotDhotPcandXhbonSigmastar:chap12} w.r.t. $\nabc_3$ and obtain, using $\nab_3\Ga_g=r^{-1}\dk^{\leq 1}\Ga_b$,  
\beaa
\nabc_3(r\qf) &=& \frac{1}{2}r^5\nabc_3\DDc\hot\DDc\ov{\Pc}  +\frac{5}{2}r^4e_3(r)\DDc\hot\DDc\ov{\Pc} -6mr\nabc_3\widehat{\Xb}    +\dk^{\leq 3}\Ga_b\\
&& +r^3\dk^{\leq 2}(\Ga_b\c\Ga_g)\\
&=& \frac{1}{2}r^5\DDc\hot\DDc\ov{\nabc_3\Pc}+ \frac{1}{2}r^5[\nabc_3, \DDc\hot\DDc]\ov{\Pc}  -\frac{5}{2}r^4\DDc\hot\DDc\ov{\Pc}\\
&&  -6mr\nabc_3\widehat{\Xb}    +\dk^{\leq 3}\Ga_b +r^3\dk^{\leq 2}(\Ga_b\c\Ga_g).
\eeaa
Now, we have, in view of the commutator identities for $[\nabc_3, \DDc\hot]$ and $[\nabc_3, \DDc]$, see Lemma \ref{COMMUTATOR-NAB-C-3-DD-C-HOT}, 
\beaa
&& [\nabc_3, \DDc\hot\DDc]\ov{\Pc}\\ 
&=& [\nabc_3, \DDc\hot]\DDc\ov{\Pc} +\DDc\hot[\nabc_3, \DDc]\ov{\Pc}\\
&=& - \frac 1 2 \tr \Xb \left( \DDc \hot\DDc\ov{\Pc} + H \hot\DDc\ov{\Pc}\right)  + H \hot \nabc_3\DDc\ov{\Pc} + r^{-1}\Ga_b  \c \dk^{\leq 1}\DDc\ov{\Pc}\\
&&+\DDc\hot\left(- \frac 1 2 \tr \Xb\DDc\ov{\Pc}  + H \hot \nabc_3\ov{\Pc} + r^{-1}\Ga_b  \c \dk^{\leq 1}\ov{\Pc}\right)\\
&=& \frac{2}{r}\DDc\hot\DDc\ov{\Pc}  +r^{-5}\dk^{\leq 2}\Ga_b +r^{-3}\dk^{\leq 2}(\Ga_b\c\Ga_g)
\eeaa 
where we used in particular the fact that $\nabc_3\Pc\in r^{-2}\dk^{\leq 1}\Ga_b$ in view of Bianchi. We deduce 
\beaa
\nabc_3(r\qf) &=& \frac{1}{2}r^5\DDc\hot\DDc\ov{\nabc_3\Pc} -\frac{3}{2}r^4\DDc\hot\DDc\ov{\Pc}  -6mr\nabc_3\widehat{\Xb}    +\dk^{\leq 3}\Ga_b\\
&& +r^3\dk^{\leq 2}(\Ga_b\c\Ga_g).
\eeaa
Next, using the following consequence of Bianchi and the null structure equations
\beaa
\nabc_3\Pc &=&  -\frac{1}{2}\DDbc \c\Bb  +\frac{3}{r}\Pc  +r^{-3}\Ga_b +\Ga_b\c\Ga_g\\
\nabc_3\Xbh &=&    -\Ab +r^{-1}\dk^{\leq 1}\Ga_b+\Ga_b\c\Ga_g,
\eeaa
we infer
\beaa
\nabc_3(r\qf) &=& \frac{1}{2}r^5\DDc\hot\DDc\ov{\left(-\frac{1}{2}\DDbc \c\Bb  +\frac{3}{r}\Pc  +r^{-3}\Ga_b +\Ga_b\c\Ga_g\right)} -\frac{3}{2}r^4\DDc\hot\DDc\ov{\Pc}\\
&&  -6mr\Big(  -\Ab +r^{-1}\dk^{\leq 1}\Ga_b+\Ga_b\c\Ga_g\Big)    +\dk^{\leq 3}\Ga_b +r^3\dk^{\leq 2}(\Ga_b\c\Ga_g)\\
&=&  -\frac{1}{4}r^5\DDc\hot\DDc\DDc\c\ov{\Bb} +6mr\Ab    +\dk^{\leq 3}\Ga_b +r^3\dk^{\leq 2}(\Ga_b\c\Ga_g)
\eeaa
which concludes the proof of \eqref{eq:roughidentitybetweene3rqfandDhotDDcBbandAbonSigmastar:chap12}.

Finally, we prove \eqref{eq:roughidentitybetweene3r2e3rqfandDhotDDcDDcAbandnab3cAbonSigmastar:chap12}. 
First, we multiply by $r^2$ and differentiate \eqref{eq:roughidentitybetweene3rqfandDhotDDcBbandAbonSigmastar:chap12} w.r.t. $\nabc_3$ and obtain
\beaa
\nabc_3(r^2\nabc_3(r\qf)) &=&  -\frac{1}{4}r^7\nabc_3\DDc\hot\DDc\DDc\c\ov{\Bb} +6mr^3\nab_3\Ab     +r^2\dk^{\leq 4}\Ga_b\\
&& +r^5\dk^{\leq 3}(\Ga_b\c\Ga_g)\\
&=&  -\frac{1}{4}r^7\DDc\hot\DDc\DDc\c\ov{\nabc_3\Bb}  -\frac{1}{4}r^7[\nabc_3, \DDc\hot\DDc\DDc\c]\ov{\Bb}\\ && +6mr^3\nab_3\Ab     +r^2\dk^{\leq 4}\Ga_b +r^5\dk^{\leq 3}(\Ga_b\c\Ga_g).
\eeaa
Now, we have, in view of the commutator identities for $[\nabc_3, \DDc\hot]$, $[\nabc_3, \DDc]$ and $[\nabc_3, \DDc\c]$, see Lemma \ref{COMMUTATOR-NAB-C-3-DD-C-HOT}, 
\beaa
&& [\nabc_3, \DDc\hot\DDc\DDc\c]\ov{\Bb}\\ 
&=& [\nabc_3, \DDc\hot]\DDc\DDc\c\ov{\Bb} +\DDc\hot[\nabc_3, \DDc]\DDc\c\ov{\Bb}+\DDc\hot\DDc[\nabc_3, \DDc\c]\ov{\Bb}\\
&=& r^{-5}\dk^{\leq 3}\Ga_b +r^{-3}\dk^{\leq 2}(\Ga_b\c\Ga_g).
\eeaa 
We deduce
\beaa
\nabc_3(r^2\nabc_3(r\qf)) &=&  -\frac{1}{4}r^7\DDc\hot\DDc\DDc\c\ov{\nabc_3\Bb}  +6mr^3\nab_3\Ab \\
&&  +r^2\dk^{\leq 4}\Ga_b  +r^5\dk^{\leq 3}(\Ga_b\c\Ga_g).
\eeaa
Next, using the following consequence of Bianchi 
\beaa
\nabc_3\Bb &=&  -\frac{1}{2}\DDbc \c\Ab  +r^{-2}\Ga_b+\Ga_b\c\Ga_g,
\eeaa
we infer
\beaa
\nabc_3(r^2\nabc_3(r\qf)) &=&  \frac{1}{8}r^7\DDc\hot\DDc\DDc\c\DDc\c\ov{\Ab}  +6mr^3\nab_3\Ab \\
&&  +r^2\dk^{\leq 4}\Ga_b  +r^5\dk^{\leq 3}(\Ga_b\c\Ga_g)
\eeaa
which concludes the proof of \eqref{eq:roughidentitybetweene3r2e3rqfandDhotDDcDDcAbandnab3cAbonSigmastar:chap12}, and hence of Proposition \ref{prop:roughidentitybetweenqfe3qfe3e3qfandPcBbandAb:chap12}.
\end{proof}

Proposition \ref{prop:roughidentitybetweenqfe3qfe3e3qfandPcBbandAb:chap12} implies the following corollary. 
\begin{corollary}\lab{cor::roughidentitybetweenqfe3qfe3e3qfandPcBbandAb:chap12:firstcor}
We have 
\bea
\nn 3m\nab_3\aa+ r^4\DDs_2\DDs_1\DDd_1\DDd_2\aa &=& \frac{1}{2}r^{-3}\Re\big(\nabc_3(r^2\nabc_3(r\qf))\big)     +r^{-1}\dk^{\leq 4}\Ga_b\\
&&  +r^2\dk^{\leq 3}(\Ga_b\c\Ga_g).
\eea
\end{corollary}

\begin{proof}
Recall \eqref{eq:roughidentitybetweene3r2e3rqfandDhotDDcDDcAbandnab3cAbonSigmastar:chap12}
\beaa
\nabc_3(r^2\nabc_3(r\qf)) &=&  \frac{1}{8}r^7\DDc\hot\DDc\DDc\c\DDc\c\ov{\Ab}  +6mr^3\nab_3\Ab \\
&&  +r^2\dk^{\leq 4}\Ga_b  +r^5\dk^{\leq 3}(\Ga_b\c\Ga_g).
\eeaa
Taking the real part, we infer
\beaa
\Re\big(\nabc_3(r^2\nabc_3(r\qf))\big) &=&  \frac{1}{8}r^7\Re\Big(\DDc\hot\DDc\DDc\c\DDc\c\ov{\Ab}\Big) +6mr^3\nab_3\aa \\
&&  +r^2\dk^{\leq 4}\Ga_b  +r^5\dk^{\leq 3}(\Ga_b\c\Ga_g).
\eeaa
Now, we have
\beaa
\DDc\c\DDc\c\ov{\Ab} &=& \DDc\c\big(2\div\aa - i2\dual(\div\aa)\big)\\
&=& 4\div\div\aa-4i\curl\div\aa,
\eeaa
and, for two scalar functions $(f, \dual f)$, 
\beaa
\Re\Big(\DDc\hot\DDc(f-i\dual f)\Big) &=& 2\nab\hot\Re(\DDc(f-i\dual f))\\
&=& 2\nab\hot\nab f +2\nab\hot\dual\nab(\dual f)
\eeaa
and hence, using the definition of the Hodge operators $\DDd_1$, $\DDd_2$, $\DDs_1$ and $\DDs_2$, see section    \ref{section:HorizontalHodgeoperators}, we infer 
\beaa
\Re\Big(\DDc\hot\DDc\DDc\c\DDc\c\ov{\Ab}\Big) &=& 8\nab\hot\nab \div\div\aa-8\nab\hot\dual\nab\curl\div\aa \\
&=& 16\DDs_2\Big(-\nab \div\div\aa+\dual\nab\curl\div\aa\Big)\\
&=& 16\DDs_2\DDs_1\Big(\div\div\aa, \curl\div\aa\Big)\\
&=& 16\DDs_2\DDs_1\DDd_1\div\aa\\
&=& 16\DDs_2\DDs_1\DDd_1\DDd_2\aa.
\eeaa
This yields 
\beaa
\Re\big(\nabc_3(r^2\nabc_3(r\qf))\big) &=&  2r^7\DDs_2\DDs_1\DDd_1\DDd_2\aa +6mr^3\nab_3\aa  +r^2\dk^{\leq 4}\Ga_b  +r^5\dk^{\leq 3}(\Ga_b\c\Ga_g)
\eeaa
and finally 
\beaa
3m\nab_3\aa+ r^4\DDs_2\DDs_1\DDd_1\DDd_2\aa &=& \frac{1}{2}r^{-3}\Re\big(\nabc_3(r^2\nabc_3(r\qf))\big)     +r^{-1}\dk^{\leq 4}\Ga_b  +r^2\dk^{\leq 3}(\Ga_b\c\Ga_g)
\eeaa
as stated. This concludes the proof of Corollary \ref{cor::roughidentitybetweenqfe3qfe3e3qfandPcBbandAb:chap12:firstcor}.
\end{proof}

The following corollary of Corollary \ref{cor::roughidentitybetweenqfe3qfe3e3qfandPcBbandAb:chap12:firstcor} is the goal of this section.
\begin{corollary}\lab{cor::roughidentitybetweenqfe3qfe3e3qfandPcBbandAb:chap12:secondcor}
We have on $\Si_*$
\bea
\nn r^4{\DDs}_2\,^{\Si_*}\DDs_1\,^{\Si_*}\DDd_1\!\,^{\Si_*}\DDd_2\!\,^{\Si_*}\Lieb_\T^j\aa &=&  O(1)\dk^{\leq 3}\nab_3\Lieb_\T^j\aa+O(1)\dk^{\leq 1}\nab_3\Lieb_\T^j\qf    +r^{-1}\dk^{\leq 4+j}\Ga_b\\
&& +r^2\dk^{\leq 3+j}(\Ga_b\c\Ga_g).
\eea
\end{corollary}

\begin{proof}
In view of Corollary \ref{cor::roughidentitybetweenqfe3qfe3e3qfandPcBbandAb:chap12:firstcor}, and since $\qf\in r\dk^{\leq 2}\Ga_g$ and $\nab_3\qf\in \dk^{\leq 3}\Ga_b$, we have
\beaa
r^4\DDs_2\DDs_1\DDd_1\DDd_2\aa &=& -3m\nab_3\aa+O(1)\dk^{\leq 1}\nab_3\qf    +r^{-1}\dk^{\leq 4}\Ga_b  +r^2\dk^{\leq 3}(\Ga_b\c\Ga_g).
\eeaa
Differentiating this identity  with respect to $\Lieb_\T$, using the fact that $\T(r)\in r\Ga_b$, and using the commutation formula $[\Lieb_\T, \dk]U = \dk(\Ga_b U)$, see Lemma \ref
{lemma:basicpropertiesLiebTfasdiuhakdisug:chap9}, we infer
\beaa
r^4\DDs_2\DDs_1\DDd_1\DDd_2\Lieb_\T^j\aa = -3m\nab_3\Lieb_\T^j\aa +O(1)\dk^{\leq 1}\nab_3\Lieb_\T^j\qf    +r^{-1}\dk^{\leq 4+j}\Ga_b +r^2\dk^{\leq 3+j}(\Ga_b\c\Ga_g).
\eeaa

Next,  in view of the comparison provided by \eqref{eq:decompositiontangentspacerfoliationSigamastaronglobalframe:chap12} and \eqref{eq:controloftansformationcoefficientsbetweenglobalframeandframeSigamstar} between the global frame of $\MM$ and the null frame $(e_3^{\Si_*}, e_4^{\Si_*}, e^{\Si_*}_1, e^{\Si_*}_2)$ adapted to the $r$-foliation of $\Si_*$, we have 
\beaa
r^4\DDs_2\DDs_1\DDd_1\DDd_2\Lieb_\T^j\aa &=& r^4{\DDs}_2\,^{\Si_*}\DDs_1\,^{\Si_*}\DDd_1\!\,^{\Si_*}\DDd_2\!\,^{\Si_*}\Lieb_\T^j\aa+O(1)\dk^{\leq 3}\nab_3\Lieb_\T^j\aa+O(r^{-1})\dk^{\leq 4}\Lieb_\T^j\aa\\
&=& r^4{\DDs}_2\,^{\Si_*}\DDs_1\,^{\Si_*}\DDd_1\!\,^{\Si_*}\DDd_2\!\,^{\Si_*}\Lieb_\T^j\aa+O(1)\dk^{\leq 3}\nab_3\Lieb_\T^j\aa+r^{-1}\dk^{\leq 4+j}\Ga_b, 
\eeaa
where ${\DDs}_2\,^{\Si_*}$, $\DDs_1\,^{\Si_*}$, $\DDd_1\!\,^{\Si_*}$ and $\DDd_2\!\,^{\Si_*}$ denote the Hodge operators tangent to the spheres $S^{\Si_*}(r)$ of the $r$-foliation of $\Si_*$. In view of the above, we infer
\beaa
r^4{\DDs}_2\,^{\Si_*}\DDs_1\,^{\Si_*}\DDd_1\!\,^{\Si_*}\DDd_2\!\,^{\Si_*}\Lieb_\T^j\aa &=&  O(1)\dk^{\leq 3}\nab_3\Lieb_\T^j\aa+O(1)\dk^{\leq 1}\nab_3\Lieb_\T^j\qf    +r^{-1}\dk^{\leq 4+j}\Ga_b \\
&& +r^2\dk^{\leq 3+j}(\Ga_b\c\Ga_g)
\eeaa
as stated. This concludes the proof of Corollary \ref{cor::roughidentitybetweenqfe3qfe3e3qfandPcBbandAb:chap12:secondcor}.
\end{proof}


\subsection{End of the proof of Theorem M2}
\lab{sec:endoftheproofofTheoremM2:chap12}


We are now ready to prove  Theorem M2. We proceed in several steps.

{\bf Step 1.} Recall from Corollary \ref{cor::roughidentitybetweenqfe3qfe3e3qfandPcBbandAb:chap12:secondcor} that we have 
\beaa
r^4{\DDs}_2\,^{\Si_*}\DDs_1\,^{\Si_*}\DDd_1\!\,^{\Si_*}\DDd_2\!\,^{\Si_*}\Lieb_\T^j\aa &=&  O(1)\dk^{\leq 3}\nab_3\Lieb_\T^j\aa+O(1)\dk^{\leq 1}\nab_3\Lieb_\T^j\qf    +r^{-1}\dk^{\leq 4+j}\Ga_b \\
&& +r^2\dk^{\leq 3+j}(\Ga_b\c\Ga_g).
\eeaa
We take the $L^2(S^{\Si_*}(r))$ norm of this identity which yields
\beaa
&&\|r^4{\DDs}_2\,^{\Si_*}\DDs_1\,^{\Si_*}\DDd_1\!\,^{\Si_*}\DDd_2\!\,^{\Si_*}\Lieb_\T^j\aa\|_{L^2(S^{\Si_*}(r))}\\
&\les& \|\dk^{\leq 3}\nab_3\Lieb_\T^j\aa\|_{L^2(S^{\Si_*}(r))}+\|\dk^{\leq 1}\nab_3\Lieb_\T^j\qf\|_{L^2(S^{\Si_*}(r))}    +r^{-1}\|\dk^{\leq 4+j}\Ga_b\|_{L^2(S^{\Si_*}(r))} \\
&& +r^2\|\dk^{\leq 3+j}(\Ga_b\c\Ga_g)\|_{L^2(S^{\Si_*}(r))}.
\eeaa
Next, note the identity 
\beaa
\DDs_1\,^{\Si_*}\DDd_1\!\,^{\Si_*} &=& \DDd_2\!\,^{\Si_*}{\DDs}_2\,^{\Si_*}+2K(S^{\Si_*})\\
&=& \DDd_2\!\,^{\Si_*}{\DDs}_2\,^{\Si_*}+\frac{2}{r^2}+r^{-1}\Ga_g
\eeaa
which yields
\beaa
r^4{\DDs}_2\,^{\Si_*}\DDs_1\,^{\Si_*}\DDd_1\!\,^{\Si_*}\DDd_2\!\,^{\Si_*}U &=& r^2{\DDs}_2\,^{\Si_*}\DDd_2\!\,^{\Si_*}\Big(r^2{\DDs}_2\,^{\Si_*}\DDd_2\!\,^{\Si_*}+2\Big)U+r\dk^{\leq 2}(\Ga_g U)
\eeaa
and hence
\beaa
&&\left\|r^2{\DDs}_2\,^{\Si_*}\DDd_2\!\,^{\Si_*}\Big(r^2{\DDs}_2\,^{\Si_*}\DDd_2\!\,^{\Si_*}+2\Big)\Lieb_\T^j\aa\right\|_{L^2(S^{\Si_*}(r))}\\
&\les& \|\dk^{\leq 3}\nab_3\Lieb_\T^j\aa\|_{L^2(S^{\Si_*}(r))}+\|\dk^{\leq 1}\nab_3\Lieb_\T^j\qf\|_{L^2(S^{\Si_*}(r))}    +r^{-1}\|\dk^{\leq 4+j}\Ga_b\|_{L^2(S^{\Si_*}(r))} \\
&& +r^2\|\dk^{\leq 3+j}(\Ga_b\c\Ga_g)\|_{L^2(S^{\Si_*}(r))}.
\eeaa
Thanks to Hodge estimates relying on Proposition \ref{prop:2D-hodge}, $r^2{\DDs}_2\,^{\Si_*}\DDd_2\!\,^{\Si_*}$ is an elliptic coercive operator and we infer, for $j\leq 2$,  
\beaa
&&\left\|(r\nab^{\Si_*})^{\leq 4}\Lieb_\T^j\aa\right\|_{L^2(S^{\Si_*}(r))}\\
&\les& \|\dk^{\leq 3}\nab_3\Lieb_\T^j\aa\|_{L^2(S^{\Si_*}(r))}+\|\dk^{\leq 1}\nab_3\Lieb_\T^j\qf\|_{L^2(S^{\Si_*}(r))}    +r^{-1}\|\dk^{\leq 4+j}\Ga_b\|_{L^2(S^{\Si_*}(r))} \\
&& +r^2\|\dk^{\leq 3+j}(\Ga_b\c\Ga_g)\|_{L^2(S^{\Si_*}(r))}.
\eeaa
Now,  in view of the comparison provided by \eqref{eq:decompositiontangentspacerfoliationSigamastaronglobalframe:chap12} and \eqref{eq:controloftansformationcoefficientsbetweenglobalframeandframeSigamstar} between the global frame of $\MM$ and the null frame $(e_3^{\Si_*}, e_4^{\Si_*}, e^{\Si_*}_1, e^{\Si_*}_2)$ adapted to the $r$-foliation of $\Si_*$, we have, for $j\leq 2$, 
\beaa
\left\|\dkb^{\leq 4}\Lieb_\T^j\aa\right\|_{L^2(S^{\Si_*}(r))} &\les& \left\|(r\nab^{\Si_*})^{\leq 4}\Lieb_\T^j\aa\right\|_{L^2(S^{\Si_*}(r))}+ \|\dk^{\leq 3}\nab_3\Lieb_\T^j\aa\|_{L^2(S^{\Si_*}(r))} \\
&& +r^{-1}\|\dk^{\leq 4+j}\Ga_b\|_{L^2(S^{\Si_*}(r))}.
\eeaa
We deduce, for $j\leq 2$, 
\beaa
&&\left\|\dkb^{\leq 4}\Lieb_\T^j\aa\right\|_{L^2(S^{\Si_*}(r))}\\
&\les& \|\dk^{\leq 3}\nab_3\Lieb_\T^j\aa\|_{L^2(S^{\Si_*}(r))}+\|\dk^{\leq 1}\nab_3\Lieb_\T^j\qf\|_{L^2(S^{\Si_*}(r))}    +r^{-1}\|\dk^{\leq 4+j}\Ga_b\|_{L^2(S^{\Si_*}(r))} \\
&& +r^2\|\dk^{\leq 3+j}(\Ga_b\c\Ga_g)\|_{L^2(S^{\Si_*}(r))}.
\eeaa

{\bf Step 2.} Next, we square the final identity of Step 1 and integrate it in $\tau$ which implies, for $j\leq 2$,  
\beaa
\int_{\Si_*(\geq \tau)}|\dkb^{\leq 4}\Lieb_\T^j\aa|^2 &\les& \int_{\Si_*(\geq \tau)}|\dk^{\leq 3}\nab_3\Lieb_\T^j\aa|^2+\int_{\Si_*(\geq \tau)}|\dk^{\leq 1}\nab_3\Lieb_\T^j\qf|^2\\
&&+\int_{\Si_*(\geq \tau)}r^{-2}|\dk^{\leq 4+j}\Ga_b|^2+\int_{\Si_*(\geq \tau)}r^4|\dk^{\leq 3+j}(\Ga_b\c\Ga_g)|^2.
\eeaa 
Together with the bootstrap assumptions on $\Ga_g$ and $\Ga_b$, the dominant condition \eqref{eq:cominantconditionofronSigmastar} for $r$ on $\Si_*$, and the control   \eqref{eq:assumptionsonqfforTheoremM2compatiblewithTheoremM1:chap12} for $\qf$ provided by Theorem M1, we infer, for $j\leq 2$,  
\beaa
\int_{\Si_*(\geq \tau)}|\dkb^{\leq 4}\Lieb_\T^j\aa|^2 &\les& \int_{\Si_*(\geq \tau)}|\dk^{\leq 3}\nab_3\Lieb_\T^j\aa|^2+\ep_0^2\tau^{-2-2\dec}\\
&&+\ep^2\left(\max_{\Si_*}r^{-2}\right)\left(\int_{\geq\tau}\frac{d\tau'}{{\tau'}^{2+2\dec}}\right)+\ep^4\left(\int_{\geq\tau}\frac{d\tau'}{{\tau'}^{3+2\dec}}\right)\\
&\les& \int_{\Si_*(\geq \tau)}|\dk^{\leq 3}\nab_3\Lieb_\T^j\aa|^2+\ep_0^2\tau^{-2-2\dec}+\ep^2\tau^{-1-2\dec}\ep_0^2\tau_*^{-2-2\dec}
\eeaa
and hence, for $j\leq 2$,
\beaa
\int_{\Si_*(\geq \tau)}|\dkb^{\leq 4}\Lieb_\T^j\aa|^2 &\les& \int_{\Si_*(\geq \tau)}|\dk^{\leq 3}\nab_3\Lieb_\T^j\aa|^2+\ep_0^2\tau^{-2-2\dec}.
\eeaa 

{\bf Step 3.} We consider the final identity of Step 2 with $j=2$, i.e. 
\beaa
\int_{\Si_*(\geq \tau)}|\dkb^{\leq 4}\Lieb_\T^2\aa|^2 &\les& \int_{\Si_*(\geq \tau)}|\dk^{\leq 3}\nab_3\Lieb_\T^2\aa|^2+\ep_0^2\tau^{-2-2\dec}\\
&\les& F^3_{\Si_*}[\Lieb_\T^2\Ab](\tau, \tau_*)+\ep_0^2\tau^{-2-2\dec}.
\eeaa 
Together with the control of $F^3_{\Si_*}[\Lieb_\T^2\Ab](\tau, \tau_*)$ provided by Proposition \ref{prop:decayofpr2tqfb}, we infer
\beaa
\int_{\Si_*(\geq \tau)}|\dkb^{\leq 4}\Lieb_\T^2\aa|^2+\int_{\Si_*(\geq \tau)}|\dk^{\leq 3}\nab_3\Lieb_\T^2\aa|^2 &\les& \ep_0^2\tau^{-2-2\dec}.
\eeaa 
Also, as a consequence of Bianchi, we have
\beaa
\nab_4\Ab &=& -\frac{1}{r}\Ab+r^{-2}\dk^{\leq 1}\Ga_b 
\eeaa
and hence, we infer, for $k\leq k_L-1$, 
\beaa
\nn\int_{\Si_*(\geq \tau)}r^2|\dk^k\nab_4\Ab|^2 &\les& \int_{\Si_*(\geq \tau)}|\dk^k\Ab|^2+\int_{\Si_*(\geq \tau)}r^{-2}|\dk^{k+1}\Ga_b|^2\\
&\les&  \int_{\Si_*(\geq \tau)}|\dk^k\Ab|^2+\ep^2\left(\max_{\Si_*}r^{-2}\right)\left(\int_{\geq\tau}\frac{d\tau'}{{\tau'}^{2+2\dec}}\right),
\eeaa
which together with the dominant condition \eqref{eq:cominantconditionofronSigmastar} for $r$ on $\Si_*$ implies 
\bea\lab{eq:controlnab4AbfromAbandlotforcontrolAbinThM2:chap12}
\int_{\Si_*(\geq \tau)}r^2|\dk^k\nab_4\Ab|^2 &\les&  \int_{\Si_*(\geq \tau)}|\dk^k\Ab|^2+\ep_0^2\tau^{-2-2\dec}.
\eea
In view of the above estimates, we easily infer
\bea\lab{eq:controldkleq4ALiebT2AbforcontrolAbinThM2:chap12}
\int_{\Si_*(\geq \tau)}|\dk^{\leq 4}\Lieb_\T^2\aa|^2 &\les& \ep_0^2\tau^{-2-2\dec}.
\eea 

Next, we consider the following iteration assumption, for $4\leq l\leq k_L-10$, that we have
\bea\lab{eq:controldkleqlALiebT2AbforcontrolAbinThM2iterationassumption:chap12}
\int_{\Si_*(\geq \tau)}|\dk^{\leq l}\Lieb_\T^2\aa|^2 &\les& \ep_0^2\tau^{-2-2\dec}.
\eea
In view of \eqref{eq:controldkleq4ALiebT2AbforcontrolAbinThM2:chap12}, \eqref{eq:controldkleqlALiebT2AbforcontrolAbinThM2iterationassumption:chap12} holds for $l=4$. We assume from now on that \eqref{eq:controldkleqlALiebT2AbforcontrolAbinThM2iterationassumption:chap12} holds for $4\leq l\leq k_L-10$. 

Next, note that \eqref{eq:controlnab4AbfromAbandlotforcontrolAbinThM2:chap12} and the iteration assumption \eqref{eq:controldkleqlALiebT2AbforcontrolAbinThM2iterationassumption:chap12} imply 
\beaa
\int_{\Si_*(\geq \tau)}r^2|\nab_4\dk^{\leq l}\Lieb_\T^2\aa|^2 &\les& \ep_0^2\tau^{-2-2\dec}.
\eeaa
Also, Proposition \ref{prop:decayofpr2tqfb} implies 
\beaa
\int_{\Si_*(\geq \tau)}|\nab_3\dk^{\leq l}\Lieb_\T^2\aa|^2 &\les& \ep_0^2\tau^{-2-2\dec}.
\eeaa
We deduce 
\beaa
\int_{\Si_*(\geq \tau)}|\dk^{\leq l+1}\Lieb_\T^2\aa|^2 &\les& \int_{\Si_*(\geq \tau)}|\dkb^{l+1}\Lieb_\T^2\aa|^2+\ep_0^2\tau^{-2-2\dec}.
\eeaa
In view of the elliptic-Hodge estimates of Proposition \ref{Prop:HodgeThmM8}, and using again the above control of $\nab_3\dk^{\leq l}\Lieb_\T^2$ and $\nab_4\dk^{\leq l}\Lieb_\T^2\aa$, as well as the and the iteration assumption \eqref{eq:controldkleqlALiebT2AbforcontrolAbinThM2iterationassumption:chap12}, we have
\beaa
\int_{\Si_*(\geq \tau)}|\dkb^{l+1}\Lieb_\T^2\aa|^2  &\les& \int_{\Si_*(\geq \tau)}|r^4\DDs_2\DDs_1\DDd_1\DDd_2\dkb^{l-3}\Lieb_\T^2\aa|^2+\int_{\Si_*(\geq \tau)}|\nab_4\dk^{\leq l}\Lieb_\T^2\aa|\\
&&+\int_{\Si_*(\geq \tau)}|\nab_3\dk^{\leq l}\Lieb_\T^2\aa|^2+\int_{\Si_*(\geq \tau)}|\dk^{\leq l}\Lieb_\T^2\aa|^2\\
&\les& \int_{\Si_*(\geq \tau)}|r^4\DDs_2\DDs_1\DDd_1\DDd_2\dkb^{l-3}\Lieb_\T^2\aa|^2 +\ep_0^2\tau^{-2-2\dec}\\
&\les& \int_{\Si_*(\geq \tau)}|\dkb^{l-3}r^4\DDs_2\DDs_1\DDd_1\DDd_2\Lieb_\T^2\aa|^2 +\int_{\Si_*(\geq \tau)}|\dk^{\leq l}\Lieb_\T^2\aa|^2\\
&&+\ep_0^2\tau^{-2-2\dec}
\eeaa
and hence
\beaa
\int_{\Si_*(\geq \tau)}|\dk^{\leq l+1}\Lieb_\T^2\aa|^2 &\les& \int_{\Si_*(\geq \tau)}|\dkb^{l-3}r^4\DDs_2\DDs_1\DDd_1\DDd_2\Lieb_\T^2\aa|^2+\ep_0^2\tau^{-2-2\dec}.
\eeaa
On the other hand, we have, in view Corollary \ref{cor::roughidentitybetweenqfe3qfe3e3qfandPcBbandAb:chap12:firstcor} and the commutator formula for $[\Lieb_\T, \dk]$ of Lemma \ref{lemma:basicpropertiesLiebTfasdiuhakdisug:chap9}, 
\beaa
r^4\DDs_2\DDs_1\DDd_1\DDd_2\Lieb_\T^2\aa &=& -3m\nab_3\Lieb_\T^2\aa+\frac{1}{2}r^{-3}\Re\big(\Lieb_\T^2\nabc_3(r^2\nabc_3(r\qf))\big)     +r^{-1}\dk^{\leq 6}\Ga_b\\
&&  +r^2\dk^{\leq 5}(\Ga_b\c\Ga_g)
\eeaa
which yields, in view of the above, 
\beaa
\int_{\Si_*(\geq \tau)}|\dk^{\leq l+1}\Lieb_\T^2\aa|^2 &\les& \int_{\Si_*(\geq \tau)}|\dk^{\leq l}\Lieb_\T^2\aa|^2+\int_{\Si_*(\geq \tau)}|\dk^{\leq l}\nab_3\qf|^2\\
&&+\int_{\Si_*(\geq \tau)}r^{-2}|\dk^{\leq l+3}\Ga_b|^2+\int_{\Si_*(\geq \tau)}r^4|\dk^{\leq l+2}(\Ga_b\c\Ga_g)|^2\\
&&+\ep_0^2\tau^{-2-2\dec}.
\eeaa
Together with the iteration assumption \eqref{eq:controldkleqlALiebT2AbforcontrolAbinThM2iterationassumption:chap12}, the bootstrap assumptions on $\Ga_g$ and $\Ga_b$, the dominant condition \eqref{eq:cominantconditionofronSigmastar} for $r$ on $\Si_*$, and the control   \eqref{eq:assumptionsonqfforTheoremM2compatiblewithTheoremM1:chap12} for $\qf$ provided by Theorem M1, we infer
\beaa
\int_{\Si_*(\geq \tau)}|\dk^{\leq l+1}\Lieb_\T^2\aa|^2 &\les& \ep_0^2\tau^{-2-2\dec}
\eeaa
which is \eqref{eq:controldkleqlALiebT2AbforcontrolAbinThM2iterationassumption:chap12} with $l$ replaced by $l+1$. This implies by iteration that \eqref{eq:controldkleqlALiebT2AbforcontrolAbinThM2iterationassumption:chap12} holds for any $4\leq l\leq k_L-9$, i.e. 
\beaa
\int_{\Si_*(\geq \tau)}|\dk^{\leq k_L-9}\Lieb_\T^2\aa|^2 &\les& \ep_0^2\tau^{-2-2\dec}.
\eeaa

{\bf Step 4.} According to the final identity of Step 2 with $j=1$, we have 
\beaa
\int_{\Si_*(\geq \tau)}|\dkb^{\leq 4}\Lieb_\T\aa|^2 &\les& \int_{\Si_*(\geq \tau)}|\dk^{\leq 3}\nab_3\Lieb_\T\aa|^2+\ep_0^2\tau^{-2-2\dec}.
\eeaa
In view of the definition of $\T$, as well as the comparison of  $\Lieb_\T$ and $\nab_\T$ provided by Lemma \ref{lemma:basicpropertiesLiebTfasdiuhakdisug:chap9}, we have on $\Si_*$
\beaa
\nab_3\Lieb_\T\aa &=& \nab_\T\Lieb_\T\aa+r^{-1}\dk\Ga_b = \Lieb_\T^2\aa+r^{-1}\dk^{\leq 1}\Ga_b
\eeaa
which yields 
\beaa
\int_{\Si_*(\geq \tau)}|\dk^{\leq k_L-9}\nab_3\Lieb_\T\aa|^2 &\les& \int_{\Si_*(\geq \tau)}|\dk^{\leq k_L-9}\Lieb_\T^2\aa|^2+\int_{\Si_*(\geq \tau)}r^{-2}|\dk^{\leq k_L-10}\Ga_b|^2
\eeaa
and hence, in view of the final estimate of Step 3, the bootstrap assumptions on $\Ga_b$ and the dominant condition \eqref{eq:cominantconditionofronSigmastar} for $r$ on $\Si_*$, we obtain 
\beaa
\int_{\Si_*(\geq \tau)}|\dk^{\leq k_L-9}\nab_3\Lieb_\T\aa|^2 &\les& \ep_0^2\tau^{-2-2\dec}.
\eeaa
In particular, in view of the above, we infer
\beaa
\int_{\Si_*(\geq \tau)}|\dkb^{\leq 4}\Lieb_\T\aa|^2+\int_{\Si_*(\geq \tau)}|\dk^{\leq k_L-9}\nab_3\Lieb_\T\aa|^2 &\les& \ep_0^2\tau^{-2-2\dec}.
\eeaa
Together with \eqref{eq:controlnab4AbfromAbandlotforcontrolAbinThM2:chap12}, we easily infer
\bea\lab{eq:controldkleq4ALiebTAbforcontrolAbinThM2:chap12}
\int_{\Si_*(\geq \tau)}|\dk^{\leq 4}\Lieb_\T\aa|^2 &\les& \ep_0^2\tau^{-2-2\dec}.
\eea
which is the analog of \eqref{eq:controldkleq4ALiebT2AbforcontrolAbinThM2:chap12} for $\Lieb_\T\aa$. 

We then proceed by iteration, exactly as in Step 3, with $\Lieb^2_\T\aa$ replaced by $\Lieb_\T\aa$, and we deduce  the following analog of the final estimate of  Step 3
\beaa
\int_{\Si_*(\geq \tau)}|\dk^{\leq k_L-8}\Lieb_\T\aa|^2 &\les& \ep_0^2\tau^{-2-2\dec}.
\eeaa

{\bf Step 5.} According to the final identity of Step 2 with $j=0$, we have 
\beaa
\int_{\Si_*(\geq \tau)}|\dkb^{\leq 4}\aa|^2 &\les& \int_{\Si_*(\geq \tau)}|\dk^{\leq 3}\nab_3\aa|^2+\ep_0^2\tau^{-2-2\dec}.
\eeaa
Proceeding as in Step 4, with $\Lieb_\T\aa$ replaced by $\aa$, and relying on the final estimate of Step 4 for $\Lieb_\T\aa$, we obtain 
\beaa
\int_{\Si_*(\geq \tau)}|\dk^{\leq 4}\aa|^2 &\les& \ep_0^2\tau^{-2-2\dec}.
\eeaa
which is the analog of \eqref{eq:controldkleq4ALiebTAbforcontrolAbinThM2:chap12} for $\aa$. 

We then proceed by iteration, exactly as in Step 3, with $\Lieb^2_\T\aa$ replaced by $\aa$, and we deduce  the following analog of the final estimate of  Step 3
\beaa
\int_{\Si_*(\geq \tau)}|\dk^{\leq k_L-7}\aa|^2 &\les& \ep_0^2\tau^{-2-2\dec}
\eeaa
as stated. This concludes the proof of Theorem M2.


\part{Top curvature estimates for Theorem M8}



\chapter{Main results of Part III}
\label{section:preliminaries-thmM8}



\section{Properties of $\MM$}


Recall, see section \ref{sec:pfdoisdvhdifuhgiwhdniwbvoubwuyf}, that the  spacetime $\MM$ consists of two regions 
\beaa
\MM=\Mint\cup\Mext, \qquad \Mint=\MM\cap\{r\leq r_0\},\qquad \Mext=\MM\cap\{r\geq r_0\},
\eeaa
and comes together with: 
\begin{enumerate}
\item a pair of constants $(a, m)$,

\item scalar functions $(r, \tau, \th)$, where  $\tau$ is a time function\footnote{The hypersurfaces $\Sigma(\tau)$ are spacelike but not  strictly spacelike.  It is tied to the horizontal structure  by properties which will be discussed below in section \ref{section:propertiesoftau-M8}.}, and with the complex scalar $q$ defined by $q= r+i a\cos\th$,

\item a global regular null pair $(e_4, e_3)$ and its associated horizontal structure,

\item a complex horizontal 1-form $\Jk$.
\end{enumerate}


\subsection{Boundaries of $\MM$}
\lab{section:BoundariesofMM-ThmM8}


The boundaries of $\MM$ are given by
\beaa
\pr\MM &=& \AA\cup\Sigma(\tau_*)\cup\Sigma_*\cup\Sigma(1)
\eeaa
where:
\begin{enumerate}
\item the  hypersurface $\AA$ is spacelike and given by  $\AA = \{r=r_+(1-\deh)\}$
for some small constant $\deh>0$,

\item the hypersurfaces  $\Sigma(1)$ and $\Sigma(\tau_*)$ are spacelike level hypersurfaces of $\tau$, with $1\leq\tau\leq\tau_*$ on $\MM$,

\item  $\Sigma_*$ is a uniformly spacelike hypersurface connecting $\Sigma(\tau_*)$ to $\Sigma(1)$.  
\end{enumerate}

With respect to the horizontal structure of $\MM$, we have for the normal of $\AA$ 
 \beaa
 \g(N_\AA, e_3)= - 1, \qquad   \g(N_\AA, e_4) \leq - \frac{1}{10}\de_\HH, \qquad \g(N_\AA, e_a) =O(\deh),
 \eeaa
and for the normal of $\Si_*$
\beaa
 \g(N_{\Si_*}, e_3)= - 1, \qquad   \g(N_{\Si_*}, e_4) \leq - 1, \qquad \g(N_{\Si_*}, e_a) =O(r^{-1}).
 \eeaa


\subsection{Linearized quantities}


We make use of the definitions  and conventions of  section 
\ref{sec:definitionoflinearizedquantities:chap4}  in Chapter 4.  
In particular, recall that the following  quantities in vanish in Kerr and are thus  linear  quantities on $\MM$.
\begin{enumerate} 
\item The quantities 
\beaa
\Xh, \quad \Xbh,\quad \Xi, \quad \Xib,\quad A, \quad B, \quad \Bb, \quad \Ab,\quad  \nab(r), \quad e_4(\th), \quad e_3(\th), \qquad \DD\hot\Jk.
\eeaa

\item The renormalization of 0-conformally invariant quantities 
\beaa
&& \Hc = H-\frac{aq}{|q|^2}\Jk, \qquad\quad \Hbc=\Hb+\frac{a\ov{q}}{|q|^2}\Jk,\qquad \qquad  \Pc = P+\frac{2m}{q^3},\\
&&\widecheck{\DD q} =\DD q+a\Jk, \qquad\quad\,\, \widecheck{\DD \ov{q}} =\DD \ov{q}-a\Jk,\qquad \widecheck{\DD(\cos\th)} = \DD(\cos\th) -i\Jk,\\
&& \widecheck{\ov{\DD}\c\Jk} = \ov{\DD}\c\Jk-\frac{4i(r^2+a^2)\cos\th}{|q|^4}.
\eeaa

\item  The remaining renormalized quantities\footnote{With respect to the  ingoing renormalization, see Definition \ref{def:renormalizationofallnonsmallquantitiesinPGstructurebyKerrvalue:3}.} 
\beaa
\bsplit
\trXc &= \tr X-\frac{2\ov{q}\De}{|q|^4}, \quad\,     \trXbc = \tr\Xb+\frac{2}{\ov{q}},\quad 
\Zc = Z-\frac{aq}{|q|^2}\Jk,\quad\,\,\,\omc  = \om  + \frac{1}{2}\pr_r\left(\frac{\De}{|q|^2} \right),\\
\widecheck{e_3(r)} &= e_3(r)+1, \qquad \widecheck{e_4(r)} = e_4(r)-\frac{\Delta}{|q|^2}, \\
\widecheck{\nab_3\Jk}&=\nab_3\Jk -\frac{1}{\ov{q}}\Jk, \qquad \widecheck{\nab_4\Jk}=\nab_4\Jk +\frac{\De \ov{q}}{|q|^4}\Jk.
\end{split}
\eeaa
 \end{enumerate}

 Finally we  provide a definition\footnote{The definition here differs from that in  section \ref{sec:definitionofGabandGagfirsttime} in that we separate  the linearized curvature quantities,  denoted $\Rc_g, \Rc_b$, from the Ricci and metric coefficients, denoted by $\Ga_g, \Ga_b$.} for $\Ga_b$, $\Ga_g$, $\Rc_b$ and $\Rc_g$. 
\begin{definition}\lab{def:GabGagRcbRcg:P3}
We define the quantities $\Ga_g, \Ga_b$
\beaa
\Ga_g&=&  \Big\{\Xi, \,  \omc, \, \trXc,\,   \Xh,\,  \Zc,\, \Hbc, \,  \trXbc,  \,  r\widecheck{e_4(r)}, \,  r^{-1}\nab(r), \,  re_4(\cos\th), \,  r^2\widecheck{\nab_4\Jk}\Big\},\\
\Ga_b&=&  \Big\{\Hc, \,  \Xbh, \,  \ombc, \,  \Xib,\, r^{-1}\widecheck{e_3(r)}, \,      \widecheck{\DD(\cos\th)}, \,  e_3(\cos\th),  \,   r\,\widecheck{\ov{\DD}\c\Jk}, \,  r\,\DD\hot\Jk, \,  r\,\widecheck{\nab_3\Jk}\Big\}.
\eeaa
We also define the   following sets of curvature quantities.
\beaa
\Rc_g =\Big\{ \Pc, \,   B, \,   A\Big\}, \qquad 
\Rc_b =\Big\{ r\Bb,  \Ab\Big\}. 
\eeaa
For any of the sets  $\Gac=\Ga_g, \Ga_b, \Rc_g, \Rc_b$ we denote by  $ \dk^{\le k}\Gac $  all derivatives  up to order $k$ with respect to  $\dk=\{\nab_3, r\nab_4, \dkb=r\nab \}$.
\end{definition}

\begin{remark}
 \lab{Remark:Ga_g=r^{-1}Ga_b}
In the norms $\Sk_k$ for Ricci and metric coefficients introduced in section \ref{subsection:MainNormsM8}, $\Ga_g$  behaves precisely as $r^{-1} \Ga_b$. Also, in the norms $\Rk_k$ for curvature components introduced in section \ref{subsection:MainNormsM8}, $\Rc_g$ behaves precisely as  $r^{-2}\Rc_b$. Since in Part   3 we are   only interested in deriving  estimates for the  top derivatives of  curvature using the norms $\Sk_k$ and $\Rk_k$, we will often identify\footnote{The sole properties for which one might have to distinguish are the low derivatives decay rates in $\tau$ that are only needed when  estimating   quadratic error terms. But, even in that case, we in fact only need $\tau^{-1-\dec}$ decay in $\MM_{trap}$ which is again consistent with these identifications.} $\Ga_g$ with $r^{-1} \Ga_b$ and $\Rc_g$ with  $r^{-2}\Rc_b$.
 \end{remark}


\subsection{The vectorfields  $\T$, $\Z$ and $\Rhat$}


We recall below the  definition of $\T, \Z$  in  the global frame of $\MM$.
\beaa
\T &=& \frac{1}{2}\left(e_4+\frac{\Delta}{|q|^2}e_3 -2a\Re(\Jk)^be_b\right),\\
\Z &=& \frac 1 2 \left(2(r^2+a^2)\Re(\Jk)^be_b -a(\sin\th)^2 e_4 -\frac{a(\sin\th)^2\De}{ |q|^2} e_3\right)\\
&=&  (r^2+a^2)\Re(\Jk)^be_b -\frac{a\sin\th}{2}\left( e_4 +\frac{\De}{|q|^2} e_3\right).
\eeaa

\begin{remark}
Recall, see \eqref{eq:atrch-e3-atrchb-e4-pert-kerr}, the formula
\bea
\lab{eq:T-e_3e_4nab}
\atrch e_3+\atrchb e_4+ 2(\eta+\etab) \c \dual \nab&=&\frac{4a\cos\th}{|q|^2} \T+ \Ga_g \c \dk.
\eea
\end{remark}

In Part III, we will use the following definition of $\Rhat$.
\begin{definition}
\lab{def:Rhat}
We define $\Rhat$ in the global frame of $\MM$ to be
\beaa
\Rhat &:=& \frac 1 2  \left( e_4-\frac{\De}{|q|^2}  e_3\right).
\eeaa
\end{definition}

\begin{remark}
Note that the above definition of $\Rhat$ differs by a factor of $\frac{|q|^2}{r^2+a^2}$ from the one used in Part II.
\end{remark}


\subsection{Conformally invariant operators}
\lab{sec:conformallyinvaariantop:part3}


Recall, see  Lemma \ref{lemma:definition-conformal-derivatives},  that  the  conformal operators  $\nabc_3$, $\nabc_4 $, $\nabc$ for an s-conformally invariant  tensor $f$ are defined by
 \begin{itemize}
 \item $\nabc_3 f:= \nab_3f-2 s \omb f$  has  signature  $(s-1)$.
 \item $\nabc_4 f:= \nab_4f+2 s \om f$  has signature  $(s+1)$.
 \item $\nabc f:= \nab_Af+ s \ze f$  has   signature  $s$. 
  \end{itemize}
  
We also introduce the following operator
\bea\lab{eq:definitionofnabcRhatconf}
\nabc_{\Rhat} := -\frac{|q|^2}{4r}\big(\trch \nabc_3+\trchb \nabc_4\big) 
\eea
which preserves the signature. 

\begin{remark}
In view of the definition of $\Rhat$, $\trchc$ and $\trchbc$, and since $\trchc, \trchbc\in\Ga_g$, we have,  for an s-conformally invariant  tensor $f$, 
\bea
\nabc_{\Rhat}f=\nab_\Rhat f+s\om f+r\Ga_g\c (\nab_3, \nab_4)f+\Ga_b f.
\eea
\end{remark}  
  
We have the following commutator lemma for conformal derivatives. 
\begin{lemma}
 \lab{commutator-nab-c-3-DD-c-hot-ThmM8}
For any tensor   horizontal   tensorfield $U$ we have
\bea
\bsplit
\, [\nabc_3, \nabc_4] U  &= O(ar^{-2} )\nabc^{\le 1 }  U+ O(r^{-3}) U + r^{-1}\Ga_b \cdot\dk^{\leq 1} U,\\
\, [\nabc_4 , \DDc]U &=  -\frac{1}{2}\tr X\DDc U +O(ar^{-3} ) \dk^{\le 1} U +\Xi\nabc_3U + r^{-2} \Ga_b  \c \dk^{\leq 1} U,\\
 \, [\nabc_3 , \DDc]U  &= -\frac{1}{2}\tr \Xb\DDc U +O(ar^{-2} ) \dk^{\le 1} U+  \Ga_b  \c \dk^{\leq 1} U.
\end{split}
\eea
\end{lemma}

\begin{proof}
This is an immediate consequence of Lemma \ref{COMMUTATOR-NAB-C-3-DD-C-HOT}. 
\end{proof}

As a consequence   of Lemma \ref{commutator-nab-c-3-DD-c-hot-ThmM8},  the definition of $\Rhat$   and  $\om =-\frac 1 2 \pr_r\Big(\frac{\De}{|q|^2}\Big) + \Ga_g$,      we derive the following lemma.
\begin{lemma}
    \lab{lemma:comm.nabc_Rhat}
    The following commutation relations hold 
    \bea
    \bsplit
     \, [  \nabc_3, \nabc_ \Rhat] U&= -\om \nabc_3U+ O(ar^{-3} )\dkb^{\le 1 }  U+ O(r^{-3}) U + r^{-1} \Ga_b \cdot\dk^{\leq 1}U,\\
       \, [  \nabc_4, \nabc_ \Rhat] U&=\frac{\De}{|q|^2} \om \nabc_3U+ O(ar^{-3} )\dkb^{\le 1 }  U+ O(r^{-3}) U + r^{-1} \Ga_b \cdot\dk^{\leq 1}U,\\
      \, [  \DDc, \nabc_ \Rhat] U&= -\frac{\De}{2|q|^2}\tr\Xb\DDc U +  O(ar^{-2}) \dk^{\le 1} U +  \Ga_b  \c \dk^{\leq 1} U.
    \end{split}
     \eea      
   \end{lemma}


\subsection{Properties  of the function $\tau$}
\lab{section:propertiesoftau-M8}


We recall below Definition \ref{definition:definition-oftau}. 
\begin{definition}[Choice of $\tau$]
\lab{definition:definition-oftau:part3}
Let  $\de_\HH>0$ small enough. We choose the smooth scalar function $\tau$ on $r\geq r_+(1-\deh)$  such that  we have on $r\geq r_+(1-\deh)$
\begin{enumerate}
\item We have  for all  $r\geq r_+(1-\de_\HH)$
\beaa
e_4(\tau)>0, \qquad e_3(\tau)>0, \qquad |\nab\tau|^2  \leq \frac{8}{9}e_4(\tau)e_3(\tau).
\eeaa
In addition, we have the following asymptotic behavior  for $r$ large
\beaa
\frac{m^2}{r^2}\les e_4(\tau)\les \frac{m^2}{r^2}, \qquad 1\les e_3(\tau)\les 1.
\eeaa

\item The unit normal $N_\Si$  to  $\Si=\Si(\tau)$, normalized by the condition $ \g(N_\Si, e_3)=-1$,  verifies 
\beaa
\g(N_{\Si}, N_{\Si}) \leq -\frac{m^2}{8r^2}.
\eeaa

\item Finally, we assume on $\MM$ 
\beaa
\T(\tau)=1+ r\Ga_b, \qquad \nab(\tau)=a \Re( \Jk)+\Ga_b.
\eeaa
\end{enumerate}
\end{definition}

 
\section{Some structure equations}
\lab{section:BianchiPairs-M8}


We recall the following null structure equations with complex notations and conformally invariant operators that will be useful later, see Proposition \ref{prop-nullstr:complex-conf}.
\beaa
\nabc_3\Xbh+\Re(\tr\Xb) \Xbh&=&\frac 1 2  \DDc\hot \Xib+  \frac 1 2  \Xib\hot(H+\Hb)-\Ab,\\
\nabc_4\Xh+\Re(\tr X)\Xh&=&\frac 1 2  \DDc\hot \Xi+\frac 1 2   \Xi\hot(\Hb+H)-A,\\
\eeaa
\beaa
\nabc_3\Hb -\nabc_4\Xib &=&  -\frac{1}{2}\ov{\tr\Xb}(\Hb-H) -\frac{1}{2}\Xbh\c(\ov{\Hb}-\ov{H}) +\Bb,\\
\nabc_4H -\nabc_3\Xi &=&  -\frac{1}{2}\ov{\tr X}(H-\Hb) -\frac{1}{2}\Xh\c(\ov{H}-\ov{\Hb}) -B.
\eeaa
Linearizing the equations for $H, \Hb$ and  writing schematically, we obtain
\bea
\lab{eqts:usefulNullstructure}
\bsplit
\nabc_3\Xbh+\Re(\tr\Xb) \Xbh&=\frac 1 2  \DDc\hot \Xib+  O(ar^{-2})\Xib -\Ab +\Ga_b\c\Ga_b,\\
\nabc_4\Xh+\Re(\tr X)\Xh&=\frac 1 2  \DDc\hot \Xi+O(ar^{-2})\Xi -A+\Ga_b\c\Ga_g,\\
\nabc_3\Hbc -\nabc_4\Xib&=\Bb +  O(r^{-1} ) \Ga_b+\Ga_b\c\Ga_b,\\
\nabc_4\Hc -\nabc_3\Xi&=-B +  O(r^{-1} ) \Ga_b+\Ga_b\c\Ga_g.
\end{split}
\eea

For  the convenience of the the reader we recall  below  the Bianchi identities  in complex notation,    see       Proposition   \ref{prop:bianchi:complex}. 
\begin{remark}  
\lab{remark:simpligyBianch-interior} For simplicity    we make  the  simplifying  convention below  $\Ga_g= r^{-1} \Ga_b$ and $\Rc_g= r^{-2}\Rc_b$, see  Remark \ref{Remark:Ga_g=r^{-1}Ga_b}.
This  is convenient for the  interior estimates  in chapter  4.
We note however that in Chapter \ref{CHAPTER:ESTIMATES-MEXTM8}   we will need to rely on the more precise structure of the equations
 as stated  in Proposition \ref{prop:bianchi:complex}.  
\end{remark}
 \beaa
 \nabc_3A +\frac{1}{2}\tr\Xb A&=& \frac 1 2 \DDc\hot B +  O(ar^{-2}) B  +O(r^{-3} ) \Xh + r^{-2}\Ga_b\c \Rc_b,\\
\nabc_4B  +2\ov{\tr X} B &=&   \frac{1}{2} \DDbc \c A+  O(ar^{-2}) A    +O(r^{-3} ) \Xi + r^{-2} \Ga_b\c \Rc_b,\\
\nabc_3B  +\tr\Xb B &=& \DDc\ov{P}+3\ov{P}H  +r^{-2} \Ga_b\c \Rc_b, \\
\nabc_4P  +\frac{3}{2}\tr X P &=&  \frac{1}{2}\DDc\c \ov{B} + O( a r^{-2}) B +r^{-2} \Ga_b\c \Rc_b, \\
\\
\nabc_3P+ \frac{3}{2}\ov{\tr\Xb} P &=& - \frac{1}{2}\DDbc \c\Bb +  O(ar^{-2} )\Bb    +r^{-1} \Ga_b\c \Rc_b, \\
\nabc_4\Bb+ \tr X\Bb  &=&-\DDc P -3P\Hb  +r^{-1} \Ga_b\c \Rc_b ,\\
\nabc_3\Bb  +2\ov{\tr\Xb}\,\Bb &=&  -\frac{1}{2}\ov{\DDc}\c\Ab  +O(a r^{-2})\Ab -O(r^{-3})  \,\Xib + \Ga_b\c \Rc_b,\\
\nabc_4\Ab  +\frac{1}{2}\tr X \Ab   &=&-\frac 1 2 \DDc \hot\Bb  +O(ar^{-2} ) \Bb +O(r^{-3} )\Xbh + \Ga_b\c \Rc_b.
\eeaa


\section{Useful commutation formulas}


We state some of the commutation formulas which will be used later.

  \begin{lemma}
  \lab{lemma:commutwithrnab=M8}
  Let $U_{A}= U_{a_1\ldots a_k} $ be a general $k$-horizontal  tensorfield. 
  \begin{enumerate}
\item  We have
  \beaa
  \, [\nab_3, r\nab_b] U_A&=&\big( O(a^2r^{-2} )+ r\Ga_b \big)\nab_b U_A -\frac 1 2  r \atrchb \dual  \nab_b U_A   + r( \eta_b-\ze_b) \nab_3 U_A  \\
&& + r\xib_b \nab_4 U_A+r\chibh_{bc}\nab_c U_A  + r\sum_{i=1}^k\Big(-\in_{a_i c} \dual\bb_b +  \frac 1 2 \B_{a_i c 3b} \Big) U_{a_1\ldots }\,^ c \,_{\ldots a_k}. 
  \eeaa
  
  \item Also, we have
   \beaa
  \, [\nab_4, r\nab_b] U_A&=& \big( O(a^2 r^{-2} ) + r\Ga_g \big)\nab_b U_A -\frac 1 2  r \atrch \dual  \nab_b U_A   + r( \etab_b+\ze_b) \nab_4 U_A  \\
&& + r\xi_b \nab_3 U_A+r\chih_{bc}\nab_c U_A  + r\sum_{i=1}^k\Big(\in_{a_i c} \dual\b_b +  \frac 1 2 \B_{a_i c 4b} \Big) U_{a_1\ldots }\,^ c \,_{\ldots a_k}. 
  \eeaa
  \end{enumerate}
    \end{lemma}
  
  \begin{proof}
  According to Lemma \ref{LEMMA:COMM-GEN-B}
   we have
   \beaa
\, [\nab_3, r\nab_b] U_A&=& r  [\nab_3, \nab_b]U_A+ e_3(r) \nab_b U_a\\
 &=&  -\frac 1 2  r  \trchb  \nab_b U_A   + e_3 (r) \nab_b U_A  -\frac 1 2  r \atrchb \dual  \nab_b U_A   + r( \eta_b-\ze_b) \nab_3 U_A +  \\
&&r\xib_b \nab_4 U_A+r\chibh_{bc}\nab_c U_A  + r\sum_{i=1}^k\Big(-\in_{a_i c} \dual\bb_b +  \frac 1 2 \B_{a_i c 3b} \Big) U_{a_1\ldots }\,^ c \,_{\ldots a_k}. 
\eeaa
Since  $e_3(r)=-1+ r\Ga_b$ and  $\trchb=-\frac{2}{r}+ O(a^2 r^{-3}) +\Ga_g$, we deduce
\beaa
 -\frac 1 2  r  \trchb  \nab_b U_A   + e_3 (r) \nab_b U_A =\big( O(ar^{-2} )+ r\Ga_b \big)\nab_b U_A
\eeaa
which implies the first commutator identity of the lemma.

Similarly, since $e_4(r)=\frac{\De}{|q|^2}+ \Ga_g $ and
$\trch =\frac{2\De}{r|q|^2}+ O(a^2 r^{-3})+\Ga_g $, we have
\beaa
 -\frac 1 2  r  \trch  \nab_b U_A   + e_4 (r) \nab_b U_A =\big( O(a^2 r^{-2} ) + r\Ga_g \big)\nab_b U_A
\eeaa
 which together  with the formula for  $ \, [\nab_4, \nab_b] U_A$
  in Lemma \ref{LEMMA:COMM-GEN-B} implies the second commutator identity of the lemma. 
  \end{proof}
  
  We infer the following corollary.
  \begin{corollary}
  \lab{corr:commwithrnab-M8}
 Given $U$ a  general $k$ horizontal tensor, we  have
 \beaa
   \, [\nab_3, r\nab_b] U_A&=& \big( O(ar^{-1})  + r\Ga_b\big)     \nab_3 U +   O(ar^{-2} )\dkb^{\le 1} U +\Ga_b\c \dk^{\le 1}  U  + \Rc_b\c U,\\
  \, [\nab_4, r\nab_b] U &=&  O(ar^{-2} )\dk^{\le 1} U +r\Ga_g\nab_3U+\Ga_g \c \dk^{\le 1}U + r^{-1} \Rc_b \c U.
 \eeaa
  \end{corollary}
  
  \begin{proof}
  Recall, see  Proposition \ref{proposition:componentsofB}, that we have
  \beaa
\B_{ a   b  c 3}&=&     -  \trchb  \big( \de_{ca}\eta_b-  \de_{cb} \eta_a\big)  -  \atrchb \big( \in_{ca}  \eta_b -  \in_{cb}  \eta_a\big) + 2 \big(- \chibh_{ca}  \eta_b + \chibh_{cb} \eta_a-  \chi_{ca} \xib_b+  \chi_{cb} \xib_a\big),\\
&=& O(r^{-1} ) ( O(a r^{-2}) +\Ga_b) +\Ga_b\c \Ga_b, 
\\
\B_{ a   b  c 4}&=&     -  \trch  \big( \de_{ca}\etab_b-  \de_{cb} \etab_a\big)  -  \atrch \big( \in_{ca}  \etab_b -  \in_{cb}  \etab_a\big) + 2 \big(- \chih_{ca}  \etab_b + \chih_{cb} \etab_a-  \chib_{ca} \xi_b+  \chib_{cb} \xi_a\big)\\
&=&  O(r^{-1} ) ( O(a r^{-2}) +\Ga_g)+\Ga_b\c\Ga_g.
\eeaa
  Hence 
  \beaa
  \, [\nab_3, r\nab_b] U_A&=&\big( O(ar^{-2} )+ r\Ga_b \big)\nab_b U_A -\frac 1 2  r \atrchb \dual  \nab_b U_A   + r( \eta_b-\ze_b) \nab_3 U_A  \\
&& + r\xib_b \nab_4 U_A+r\chibh_{bc}\nab_c U_A  + r\sum_{i=1}^k\Big(-\in_{a_i c} \dual\bb_b +  \frac 1 2 \B_{a_i c 3b} \Big) U_{a_1\ldots }\,^ c \,_{\ldots a_k} \\
&=& \big( O(ar^{-3} )+ \Ga_b \big)\dkb U+ r\eta\nab_3 U+\Ga_b \dk U+\Rc_b\c U+  ( O(a r^{-2}) +\Ga_b)\c U\\
&=&  r\eta\nab_3 U +  O(ar^{-2} )\dkb^{\le 1} U +\Ga_b\c \dk^{\le 1}  U+ \Rc_b\c U
  \eeaa
  and hence
  \beaa
    \, [\nab_3, r\nab_b] U_A&=& \big( O(ar^{-1})  + r\Ga_b\big)     \nab_3 U +   O(ar^{-2} )\dkb^{\le 1} U +\Ga_b\c \dk^{\le 1}  U + \Rc_b\c U 
  \eeaa
  as stated.
  
Similarly
 \beaa
  \, [\nab_4, r\nab_b] U_A&=& \big( O(a^2 r^{-2} ) + r\Ga_g \big)\nab_b U_A -\frac 1 2  r \atrch \dual  \nab_b U_A   + r( \etab_b+\ze_b) \nab_4 U_A  \\
&& + r\xi_b \nab_3 U_A+r\chih_{bc}\nab_c U_A  + r\sum_{i=1}^k\Big(\in_{a_i c} \dual\b_b +  \frac 1 2 \B_{a_i c 4b} \Big) U_{a_1\ldots }\,^ c \,_{\ldots a_k} \\
&=& \big( O(a^2 r^{-3} ) + \Ga_g \big)\dkb U+ O(a r^{-2} )\dk U + r^{-1}\Rc_b\c U  +  ( O(a r^{-2}) +\Ga_g) \c U \\
&=& O(ar^{-2} )\dk^{\le 1} U+r\Ga_g\nab_3U+\Ga_g \c \dk^{\le 1}U + r^{-1} \Rc_b \c U
  \eeaa
  as stated.
  \end{proof}

  
  \subsection{Commutation formulas with $\Lieb_\T$}
  

We recall below   some   the following commutation Lemma with $\Lieb_\T$. The definition of $\Lieb_X$ for an arbitrary  vectorfield $X$ was given in  section \ref{section:horizLieDerivatives}.

\begin{lemma}
\lab{lemma:commLieb_TM8}
For a horizontal covariant k-tensor $U$ we have
 \beaa
  \bsplit
  \nab_b(\Lieb_\T U_{A})-\Lieb_\T(\nab_b U_{A})&= r^{-1}\dk\Ga_b \c U,\\
   \nab_4(\Lieb_\T U_A)-\Lieb_\T(\nab_4 U_A) &=r^{-1}\dk(\Ga_b \c U),\\
   \nab_3(\Lieb_\T U_A)-\Lieb_\T(\nab_3 U_A)  &=\dk(\Ga_b\c  U).
  \end{split}
  \eeaa
\end{lemma}

\begin{proof}
The proof, based on Lemma  \ref{lemma:commutation-Lieb-nab}
and the fact that   $\piT\in \Ga_b$, see Lemma \ref{LEMMA:DEFORMATION-TENSORS-T}, appears  in the proof  of Lemma \ref{lemma:basicpropertiesLiebTfasdiuhakdisug:chap9}.
\end{proof}


\subsection{Commutation formulas with $\ov{q} \nab_4$  and  $q\nab_4$}


\begin{lemma}
\lab{Lemma:commutationwith-qnabc_4}
The following commutation formulas hold true for an arbitrary  horizontal tensor $U$
\bea
\lab{eq:Lemma-commutationwith-qnabc_4-1}
\bsplit
\,[\nabc_{\ov{q}e_4}, \nabc_3] U &= - \frac {1}{ 2} \tr \Xb \nabc_{  \ov{q} e_4} U  +  O(ar^{-1} )\nab^{\le 1 }  U+ O(r^{-2}) U + \Ga_b \cdot \dk^{\le 1} U,\\
\,[\nabc_{\ov{q}e_4}, \nabc_4] U &=-
\frac 1 2 \ov{\tr X} \nabc_{\ov{q} e_4}U +r^{-1} \Ga_b  \c \dk^{\leq 1} U,\\
\,[\nabc_{ q e_4}, \nabc_3] U &= - \frac {1}{ 2} \ov{\tr \Xb }\nabc_{  q e_4} U  +  O(ar^{-1} )\nab^{\le 1 }  U+ O(r^{-2}) U + \Ga_b \cdot \dk^{\le 1} U,\\
\,[\nabc_{ q e_4}, \nabc_4] U &=-
\frac 1 2\tr X \nabc_{ q e_4} \nabc_4 U +r^{-1} \Ga_b  \c \dk^{\leq 1} U.
\end{split}
\eea
Also,
\bea
\lab{eq:Lemma-commutationwith-qnabc_4-2}
\bsplit
[ q \nabc_4 , \DDc]U&= -\frac{1}{2}\tr X q \DDc U +O(ar^{-2} ) \dk^{\le 1} U +r\Xi\nabc_3U+ r^{-1} \Ga_b  \c \dk^{\leq 1} U,\\
[ \ov{q} \nabc_4 , \DDc]U&= -\frac{1}{2}\tr X \ov{q} \DDc U +O(ar^{-2} ) \dk^{\le 1} U +r\Xi\nabc_3U+ r^{-1} \Ga_b  \c \dk^{\leq 1} U,\\
[\ov{q} \nabc_4 , \ov{\DDc}]U&= -\frac{1}{2}\ov{\tr X}\, \ov{q}\, \ov{\DDc} U +O(ar^{-2} ) \dk^{\le 1} U +r\ov{\Xi}\nabc_3U+ r^{-1} \Ga_b  \c \dk^{\leq 1} U,\\
[q\nabc_4 , \ov{\DDc}]U&= -\frac{1}{2}\ov{\tr X}\, q\, \ov{\DDc} U +O(ar^{-2} ) \dk^{\le 1} U +r\ov{\Xi}\nabc_3U+ r^{-1} \Ga_b  \c \dk^{\leq 1} U.
\end{split}
\eea
\end{lemma}

\begin{proof}
We have
\beaa
\,[\nabc_{\ov{q}e_4}, \nabc_3] U&=&\ov{q} [\nabc_4, \nabc_3] U - e_3(\ov{q}) \nabc_4  U.
\eeaa
Making use of $\nab_3 \ov{q} =  \frac {1}{ 2} \tr \Xb  \ov{q}+ r\Ga_b $ and 
Lemma \ref{commutator-nab-c-3-DD-c-hot-ThmM8} we deduce
\beaa
\bsplit
\,[\nabc_{\ov{q}e_4}, \nabc_3] U&=\left(- \frac {1}{ 2} \tr \Xb  \ov{q}+ r\Ga_b\right)\c \nabc_4 U+  O(ar^{-1} )\nab^{\le 1 }  U+ O(r^{-2}) U + \Ga_b \cdot \dk^{\leq 1}U\\
&= - \frac {1}{ 2} \tr \Xb \nabc_{  \ov{q} e_4} U  +  O(ar^{-1} )\nab^{\le 1 }  U+ O(r^{-2}) U + \Ga_b \cdot \dk^{\le 1} U
\end{split}
\eeaa
as stated. Similarly, using  
$\nab_4 \ov{q} =  \frac {1}{ 2}\ov{ \tr X } \ov{q} + \Ga_b $, we deduce
\beaa
\,[\nabc_{\ov{q}e_4}, \nabc_4] U&=&- \nab_4(\ov{q} ) \nabc_4 U=-
\frac 1 2 \ov{\tr X} \nabc_{\ov{q} e_4}U+r^{-1} \Ga_b  \c \dk^{\leq 1} U.
\eeaa
The other results  in \eqref{eq:Lemma-commutationwith-qnabc_4-1} follow  by conjugation. 

To check \eqref{eq:Lemma-commutationwith-qnabc_4-2}  we make use of   the formula for $\, [\nabc_4 , \nabc]U $   in Lemma \ref{commutator-nab-c-3-DD-c-hot-ThmM8} and  the fact that 
 $\DD q= -a \Jk+\Ga_b+r\Ga_g= O(ar^{-1}) +\Ga_b$
we derive
\beaa
[q\nabc_4 , \DDc]U&=& q[\nabc_4, \DDc \,] U- \DD(q)\c \nabc_4  U\\
&=&   -\frac{1}{2}\tr X q \DDc U +O(ar^{-2} ) \dk^{\le 1} U +r\Xi\nabc_3U+ r^{-1} \Ga_b  \c \dk^{\leq 1} U.
\eeaa
The second identity in  \eqref{eq:Lemma-commutationwith-qnabc_4-2} is derived in the same way, and the last two identities in  \eqref{eq:Lemma-commutationwith-qnabc_4-2} follow  by conjugation. This concludes the proof of Lemma \ref{Lemma:commutationwith-qnabc_4}.
\end{proof}

 
  \section{Non-integrable Hodge estimates}
\lab{section:Hodge-ThmM8}


The following is an immediate corollary of Proposition \ref{Prop:HodgeThmM8}.
    \begin{corollary}
    \lab{Cor:HodgeThmM8}
    Given $f\in \sk_p$,  $p=0, 1, 2$, and $\DD$ any of the Hodge  operators  acting on $f$. Then, for any $q$, 
    \bea
    \int_S r^q |\nab f|^2 \les    \int_S r^q |\DD f|^2 +   \int_S r^{q-2}  | f|^2     +O(a^2, \ep^2  )    
    \int_S r^{q-2}  |\dk f|^2. 
    \eea
    Also, for higher derivatives,
     \bea
    \int_S r^q |\nab \dk^{\le k}  f|^2 \les    \int_S r^q \Big( |\dk^{\le k} \DD f|^2 +   r^{-2}   |\dk^{\le k} f|^2 \Big)    +O(a^2, \ep^2  )    
    \int_S r^{q-2}  | \dk^{\le k+1}  f   |^2. 
    \eea
    \end{corollary}


\section{Main norms}
\lab{subsection:MainNormsM8}


We  recall  here the main norms  appearing  in section 9.4.1 of \cite{KS:Kerr} in connection with the higher order curvatures estimates of Theorem 9.4.10 in \cite{KS:Kerr}.


\subsection{Norms of $\Mext$}


\begin{definition}
\lab{def:exteriorGac.norms}
We define the following norms  for the Ricci coefficients in $\Mext$ with respect to the global frame of $\MM$
\bea
\nn\Skext_k^2 &:=& \sup_{\la \ge r_0}  \int_{r=\la}\Big(r^2\big|\dk^{\le k}(\Ga_g\setminus\{\trXbc\}) |^2+ \big| \dk^{\le k}(\Ga_b\setminus\{\Xib\})\big|^2\\
&& \qquad\qquad +r^{2-\dt}|\dk^{\le k}\trXbc |^2+ r^{-\dt}|\dk^{\le k}\Xib |^2\Big).
\eea
For convenience we introduce the notation 
 \bea\lab{eq:defintionofGabprimeandGagprime:partIII}
 \Ga_b':=\Ga_b\setminus \{\Xib\}, \qquad \Ga_g':=\Ga_g\setminus \{\trXbc\}. 
 \eea
\end{definition}

\begin{definition}
\lab{Definiition:Rkext}
We define  the following norms for the  curvature coefficients in $\Mext$, with respect to the global frame of $\MM$.
\bea
\bsplit
\Rkext^2_k :=&\int_{\Mext} r^{3+\dt}|\dk^{\le k}(A, B)|^2 \\
& +r^{3-\dt}\big(|\dk^{\le k}\Pc |^2 +r^{-2} |\dk^{\le k}\Bb|^2 +r^{-4} |\dk^{\le k}\Ab|^2\big).
\end{split}
\eea
\end{definition}


\subsection{Norms of $\Mint$} 


\begin{definition}
\lab{definit:norms-SkMint'}
We define the following norms  for the Ricci coefficients  in   $\Mint$ with respect to the global frame of $\MM$
\beaa
\Skint_k^2 := \int_{\Mint}\big|\dk^{\le k}\Gac\big|^2, 
\eeaa
 where $\Gac$ denotes the set of all  linearized Ricci and metric coefficients  with respect  to the global frame of $\MM$, i.e.   
  \beaa
\Gac &:=& \Big\{\Xib,\, \omb,\,  \trXbc,\,  \Xbh,\,  \widecheck{e_3(r)},\, e_3(\cos\th),\, \widecheck{\nab_3\Jk},\, \Zc,\,  \Hbc,\, \Hc,\,  \DD r,\, \widecheck{ \DD\cos \th},\, \widecheck{ \ov{\DD}\c \Jk},\,  \DD\hot\Jk, \\
&& \trXc,\, \Xh,\, \omc,\,  \widecheck{e_4(r)},\, e_4(\cos\th),\, \widecheck{\nab_4 \Jk},\, \Xi\Big\}.
  \eeaa
 \end{definition}

For the curvature norms in $\Mint$, we rely in particular on the scalar function $\tau$  of Definition \ref{definition:definition-oftau:part3}. 
 We also   recall the definition, see Definition \ref{def:causalregions},  of  the non trapped region of $\Mint$ given by
\bea
\Mint_{\ntrap} &:=& \Mint\cap\left\{\frac{|\TT|}{r^3}\geq \de_{trap}\right\},
\eea
where  $\TT$ is  the   polynomial in $r$ defined in \eqref{definition-TT}, i.e.
    \beaa
    \TT=r^3-3mr^2+ a^2r+ma^2,
    \eeaa
and where we choose $\delta_{trap} = \frac{1}{10}$ as in Lemma \ref{rem:de-trap}.

\begin{definition}
\lab{definit:norms-RkMint'}
We define the following norms for the null  curvature  components  in $\Mint$
\beaa
\Rkint_k^2 &=&  \int_{\Mint}\Big( \big| \nab_{\Rhat} \dk^{\le k-1}\Rc\big|^2
+|\dk^{\le k-1}\Rc|^2\Big) +\int_{\Mint_{\ntrap}}  \big| \dk^{\le k}\Rc\big|^2    \\
&&+\sup_\tau\int_{\Mint\cap\Si(\tau)}  |\dk^{\le k}  \Rc|^2,
\eeaa
 where, see  Definition \ref{def:Rhat}, $\Rhat=  \frac 1 2  \left( e_4-\frac{\De}{|q|^2}  e_3\right)$, and  where $\Rc$ denotes the set of all  linearized curvature components  with respect  to the global frame of $\MM$, i.e.   
  \beaa
\Rc &:=& \Big\{\Ab,\,\, \Bb,\,\,  \Pc,\,\,  B,\,\,  A\Big\}.
  \eeaa
\end{definition}


\subsubsection{Global norms}


We   define the global  norms  of $\MM$ as follows
\bea\lab{eq:defintionofglobalnorms}
\bsplit
\Sk_k &=\Skext_k+\Skint_k, \qquad 
\Rk_k =\Rkext_k+\Rkint_k.
\end{split}
\eea

The following lemma  will be used  frequently in this chapter.
\begin{lemma}
\lab{lemma:auxilliarynormsforGa_b}
The following estimates hold true:
\begin{enumerate}
\item We have, in  $\Mext$,  for $p> 1+\dt$.
\bea
\lab{eq:auxilliarynormsforGa_b1}
\int_{\Mext}\Big(r^{-p}  | \dk^{\le k}\Ga_b|^2 +r^{-p+2}  | \dk^{\le k}\Ga_g|^2\Big) &\les  &  r_0^{1-p+\dt}    \Skext_k^2.
\eea
\item We also have the stronger  estimate for the notations $\Ga_b'$ and $\Ga_g'$ introduced in \eqref{eq:defintionofGabprimeandGagprime:partIII}, for $p>1$,
\bea
\lab{eq:auxilliarynormsforGa_b1-strong}
\int_{\Mext} \Big(r^{-p }| \dk^{\le k}\Ga' _b|^2 +r^{-p+2 }| \dk^{\le k}\Ga' _g|^2\Big) &\les  &  r_0^{1-p}    \Skext_k^2.
\eea

\item  We have
\bea
\lab{eq:auxilliarynormsforRb:boundary}
\nn\sup_{\tau\in[1, \tau_*]}\int_{\Si(\tau)} r^{3-\dt}\Big(| \dk^{\le k}\Rc_g|^2 +r^{-4}  | \dk^{\le k}\Rc_b|^2\Big)\\
+\int_{\AA\cup\Si_*} r^{3-\dt}\Big(| \dk^{\le k}\Rc_g|^2 +r^{-4}  | \dk^{\le k}\Rc_b|^2\Big) &\les& \Rk_{k+1}\Rk_k, 
\eea
and, for  $ p>1$, 
\bea
\lab{eq:auxilliarynormsforGa_b2}
\nn\sup_{\tau\in[1, \tau_*]}\int_{\Si(\tau)} \Big(r^{-p}  | \dk^{\le k}\Ga_b'|^2 +r^{-p+2}  | \dk^{\le k}\Ga_g'|^2\Big)\\
\nn +\int_{\AA\cup\Si_*} \Big(r^{-p}  | \dk^{\le k}\Ga_b'|^2 +r^{-p+2}  | \dk^{\le k}\Ga_g'|^2\Big)\\
 \nn +\sup_{\tau\in[1, \tau_*]}\int_{\Si(\tau)}  \Big(r^{-p-\dt}  | \dk^{\le k}\Xib|^2+r^{-p-\dt+2}  | \dk^{\le k}\trXbc|^2\Big)\\
  +\int_{\AA\cup\Si_*}  \Big(r^{-p-\dt}  | \dk^{\le k}\Xib|^2+r^{-p-\dt+2}  | \dk^{\le k}\trXbc|^2\Big) &\les& \Sk_{k+1}\Sk_k. 
\eea
\end{enumerate}
\end{lemma}

\begin{proof}
We write for $p>  1+\dt$, in view of the definition of $\Skext$ norms.
\beaa
\int_{\Mext} r^{-p}  | \dk^{\le k}\Ga_b|^2&\les  & \int_{r_0} ^\infty \la^{-p+\de_B} \left(\int_{r=\la}  r^{-\de_B}|\dk^{\le k}\Ga_b|^2\right) d\la  \les   r_0^{1-p+\dt}  \Skext_k^2,
\eeaa
and similarly for $\Ga_g$. Also, we have, for $p>1$,
\beaa
\int_{\Mext} r^{-p}  | \dk^{\le k}\Ga_b'|^2&\les  & \int_{r_0} ^\infty \la^{-p} \Big(\int_{r=\la}   |\dk^{\le k}\Ga_b'|^2\Big) d\la  \les   r_0^{1-p}  \Skext_k^2,
\eeaa
and similarly for $\Ga_g'$. Finally, the last  estimates of the lemma  follow from the standard trace theorem, using the definition of $\Rk_k$ for curvature and the two first estimates for $\Ga_b$.
\end{proof}


\subsubsection{Initial data norms} 


The initial data  is set on the spacelike  hypersurface $\Si_1$ as follows.

\begin{definition}\lab{def:initialdatanorm}
We define the following initial data  norms on $\Si_1$
\bea
\bsplit
\Ik_{k} := &\sup_{S\subset\Si_1 }r^{\frac{5}{2}  +\de_B} \Big(  \big\| \dk^k\, (A, B) \big\|_{L^2(S)}  + \big\| \dk^k\, B\big\|_{L^2(S)}\Big)\\
&+ \sup_{S\subset\Si_1 }\Big(  r^2 \big\| \dk^k\,  \Pc  \big\|_{L^2(S)}+  r  \big\| \dk^k\, \Bb\big\|_{L^2(S)} + \big\| \dk^k\, \Ab\big\|_{L^2(S)}\Big).
\end{split}
\eea
\end{definition}


\section{Top curvature estimates for Theorem M8}



\subsection{Main assumptions} 
\lab{section:MainAss.ThmM8}


The goal of this part, i.e. Part III, is to provide the control of high order curvature estimates needed to complete the proof of Theorem M8 of \cite{KS:Kerr}. In this section, we recall the assumptions on which this proof rests.


\subsubsection{Control of the initial data} 


We have the following control of the initial data norm of Definition \ref{def:initialdatanorm}
\bea\lab{eq:controlofinitialdataforThM8}
\Ik_{k_L}\les \ep_0.
\eea
This results has been proved in Theorem 9.4.12 of \cite{KS:Kerr}.


\subsubsection{Bootstrap assumptions} 


Relative to the global norms \eqref{eq:defintionofglobalnorms},  we make the following bootstrap assumption 
\bea\lab{eq:mainbootassforchapte13}
\Sk_k+\Rk_k &\leq& \ep, \quad k\le k_L,
\eea
 see the bootstrap assumption (9.4.32) in \cite{KS:Kerr}. 
 
In addition, we make a bootstrap assumption on decay for low derivatives, weaker that the corresponding one  in \cite{KS:Kerr}, to deal with trapping. Recall  the scalar function 
$\tau_{trap}$ defined by
\beaa
\tau_{trap} := \left\{\ba{lll}
1+\tau & \textrm{on} & \MM_{trap},\\
1& \textrm{on} & \Mntrap.
\ea\right.
\eeaa
Then, we assume that $(\Ga_g, \Ga_b)$ and $(A, B, \Pc, \Rc_b)$ satisfy the following estimates on $\MM$
\bea\lab{eq:auxlowderivativebootassdecayforchapte13}
\bsplit
r^{\frac{7}{2}+\dec}|\dk^{\leq k}(A, B)|+r^3|\dk^{\leq k}\Pc|+r^2|\dk^{\leq k}\Ga_g|+r|\dk^{\leq k}(\Rc_b, \Ga_b)| &\leq \frac{\ep}{\tau_{trap}^{1+\dec}}, \quad k\leq \frac{\kl}{2}.
\end{split}
\eea


\subsubsection{Iteration assumption} 


We make  the  iteration assumption for $J$  in the range  $\frac{k_L}{2}\leq J\leq k_L-1$,
\bea\lab{eq:iterationassumptiondiscussionThM8:bis}
\Sk_J+\Rk_J  & \les&\ep_J,
\eea
see the iteration assumption (9.4.33) in \cite{KS:Kerr}. 

\begin{remark}
We refer to (9.4.34) in \cite{KS:Kerr} for the specific choice of $\ep_J$. We do not recall it here since it is irrelevant for the statements and proofs of Part III.
\end{remark}


\subsubsection{Identities satisfied by the global frame of $\MM$} 


We assume that the following identities hold for  the global frame of $\MM$
\bea\lab{eq:specialidentityforthegloablframeofMMinpartIII}
\Xi=0, \qquad \Hbc=0, \quad\textrm{for}\quad r\geq r_0+1.
\eea

\begin{remark}\lab{rmk:specialidentityforthegloablframeofMMinpartIII}
The global frame used in Part III is constructed in section 9.4 of \cite{KS:Kerr} and indeed satisfies \eqref{eq:specialidentityforthegloablframeofMMinpartIII}. The identity $\Hbc=0$ is only used for the control of $\qfb$ in Theorem \ref{Thm:Estimates-forqf}, while the identity $\Xi=0$ is used for the control of $\qfb$ in Theorem \ref{Thm:Estimates-forqf} and implies in particular  $\Xi\in r^{-1}\Ga_g$ which is used to control the error terms of some commutators appearing in Part III.
\end{remark}


\subsection{Control of high order curvature estimates}


The following is our main result of this Part on the control of high order curvature estimate for the proof of Theorem M8 of \cite{KS:Kerr}. 

\begin{theorem}[Control of Curvature]
\lab{prop:rpweightedestimatesiterationassupmtionThM8}
Let $J$ such that $\frac{k_L}{2}\leq J\leq k_L-1$. Assume
\begin{itemize}
\item the control of initial data in \eqref{eq:controlofinitialdataforThM8},

\item the bootstrap assumptions \eqref{eq:mainbootassforchapte13} and \eqref{eq:auxlowderivativebootassdecayforchapte13},

\item the iteration assumption \eqref{eq:iterationassumptiondiscussionThM8:bis}.
\end{itemize}
Then  the following estimates hold  in $\MM$
\bea
\lab{eq:InteriorcurvEstimatesThmM8}
   \nn\Rkint_{J+1}^2&\les& r_0^{18}\Big( \ep_J(\Sk_{J+1}+\Rk_{J+1}) +\ep_J^2+\ep_0^2\Big)+|a|r_0^3\Sk^2_{J+1}\\
&&+r_0^{\frac{27}{4}}\Sk_{J+1}^{\frac{3}{2}}\Big(\ep_0+\sqrt{\ep_J}\sqrt{\Sk_{J+1}+\Rk_{J+1}}\Big)^{\frac{1}{2}},\\
\lab{eq:ExteriorcurvEstimatesThmM8}
\Rkext_{J+1}^2 &\les& r_0^{3+\de_B}\Rkint^2_{J+1}+ r_0^{-\de_B} \Skext^2_{J+1} +\ep_J^2 +\ep_0^2,
\eea
where the constant in $\les$ is   independent of $r_0$.
\end{theorem}

\begin{remark}
\eqref{eq:InteriorcurvEstimatesThmM8} and \eqref{eq:ExteriorcurvEstimatesThmM8} imply
\bea\lab{eq:TotalIntplusExtcurvEstimatesThmM8}
\nn \Rk_{J+1} &\les&  r_0^{21+\dt}\ep_J+ r_0^{\frac{21}{2}+\frac{\dt}{2}}\Big( \sqrt{\Sk_{J+1}}\sqrt{\ep_J} +\ep_0\Big)+\sqrt{|a|}r_0^{3+\frac{\dt}{2}}\Sk_{J+1}\\
\nn&& + r_0^{-\frac{\de_B}{2}} \Skext_{J+1} +r_0^{\frac{39}{8}+\frac{\dt}{2}}\Sk_{J+1}^{\frac{3}{4}}\Big(\ep_0+\sqrt{\ep_J}\sqrt{\Sk_{J+1}}\Big)^{\frac{1}{4}}\\
&& +r_0^{\frac{39}{7}+\frac{4\dt}{7}}\Sk_{J+1}^{\frac{6}{7}}\ep_J^{\frac{1}{7}}
\eea
where the constant in $\les$ is   independent of $r_0$. \eqref{eq:TotalIntplusExtcurvEstimatesThmM8} proves 
Theorem 9.4.15 of \cite{KS:Kerr}. The control of $\Rk_{J+1}$ provided by \eqref{eq:TotalIntplusExtcurvEstimatesThmM8}, together with the control for $\Sk_{J+1}$ provided by Proposition 9.4.17--9.4.20 in \cite{KS:Kerr}, allows to obtain the iteration assumption \eqref{eq:iterationassumptiondiscussionThM8:bis} with $J$ replaced by $J+1$ for $\ep_J$ given by (9.4.34) in \cite{KS:Kerr}. This iteration procedure then concludes the proof of Theorem M8 of \cite{KS:Kerr}, see section 9.4.7 in \cite{KS:Kerr} for the iteration procedure, and section 9.4.3 in \cite{KS:Kerr} for the proof of Theorem M8.
\end{remark}

\begin{remark}\lab{rmk:theinteriorestimatesareinfactglobal}
The estimate \eqref{eq:InteriorcurvEstimatesThmM8} results in fact from  global energy-Morawetz estimates on $\MM$ which are then restricted to $\Mint$. 
\end{remark}


\subsection{Structure of  the  proof of Theorem \ref{prop:rpweightedestimatesiterationassupmtionThM8}}


 Theorem \ref{prop:rpweightedestimatesiterationassupmtionThM8}  is proved  as follows: 

\begin{enumerate}
\item {\bf Energy-Morawetz estimates.}  Recall from Remark \ref{rmk:theinteriorestimatesareinfactglobal} that the interior  estimates \eqref{eq:InteriorcurvEstimatesThmM8} follow  in fact from global energy-Morawetz estimates on $\MM$ which are then restricted to $\Mint$. These energy-Morawetz estimates are given in  Chapters  \ref{Chapter:EN-MorforPc} and  \ref{CHAPTER:ENERGYMOR-INTERIORM8} as follows:
\begin{enumerate}
\item We start by deriving the following energy-Morawetz estimates for $\Pc$ in $\MM$ stated in Theorem \ref{theorem:Morawetz-EnergyPc}
\bea
\lab{eq:mainestimatePc-summaryM8}
\BEF_\de^{J}[r^2 \Pc ] &\les  r_0^{15}\Big( \Sk_{J+1} \Sk_J +\Rk_{J+1}\Rk_J +\ep_J^2+\ep_0^2\Big),
 \eea
 where the $\BEF_\de^J$ norms are defined in section \ref{subsection:recallMorawetz-Energy}.  
The proof requires   a linearized  version  of  the  wave equation verified   by  $P$, see  section \ref{section:waveEqforP},  the general Morawetz- Energy estimates  for scalar  wave equations derived in Part II and  recalled here in Proposition  \ref{Prop:scalarwavePsi-M8}, as well as  the main  energy-Morawetz estimates for $\qf$, derived  also  in Part II. 

\item Next, we derive the following energy-Morawetz estimates for $B$ and $\Bb$ in $\MM$ stated in 
Proposition \ref{proposition:EstimatesBBb-interior}
\bea
 \lab{eq:mainestimateBBb-summaryM8}
  \BEF^J_\de[r^2B]+\BEF^J_\de[\Bb] &\les&    \BEF_\de^{J}[r^2 \Pc ] + \ep_0^2+\ep_J^2+    |a|  \Sk^2_{J+1}.
 \eea
These estimates   are based  on   integral   estimates for   Bianchi pairs   stated in Proposition \ref{Prop:Bainchi-pairsEstimates-integrated}.
    
\item Finally, we derive the following energy-Morawetz estimates for $A$ and $\Ab$ in $\MM$ stated in Proposition  \ref{proposition:EstimatesAAb-interior}
     \bea
 \lab{eq:mainestimateAAb-summaryM8}
  \bsplit
   \BEF^J_\de[ A]&\les  \BEF_\de^J[ r^2 B], +\ep_0^2+\ep_J^2 + |a|\Sk^2 _{J+1},\\
    \BEF^J_\de[\Ab]&\les \BEF_\de^J[\Bb]+\ep_0^2+\ep_J^2 + |a|\Sk^2 _{J+1}.
   \end{split}
   \eea
  These estimates   are again based  on   integral   estimates for   Bianchi pairs   stated in Proposition \ref{Prop:Bainchi-pairsEstimates-integrated}.
\end{enumerate}
The combination of the estimates \eqref{eq:mainestimatePc-summaryM8}, \eqref{eq:mainestimateBBb-summaryM8} and \eqref{eq:mainestimateAAb-summaryM8} then yields  \eqref{eq:InteriorcurvEstimatesThmM8}, see section \ref{sec:proofofeq:InteriorcurvEstimatesThmM8}.

\item {\bf Exterior estimates.} The  exterior estimates  \eqref{eq:ExteriorcurvEstimatesThmM8} are proved in  Chapter \ref{CHAPTER:ESTIMATES-MEXTM8} based on $r^p$ weighted estimates  for 
  Bianchi pairs.
\end{enumerate}


\chapter{Energy-Morawetz estimates for $\Pc$}
\lab{Chapter:EN-MorforPc}



\section{Control of $\Pc$}



\subsection{Morawetz-Energy norms}
\lab{subsection:recallMorawetz-Energy}


 For the convenience of the reader  we recall  the  following  norms, see  section   \ref{subsection:basicnormsforpsi},
 \beaa
 \bsplit
B^k_\de[\psi](\tau_1, \tau_2)&=  \int_{\MM_{trap}(\tau_1,\tau_2)}\big(| \nab_\Rhat \dk^{\le k}\psi |^2 + |\dk^{\le k}\psi|^2\big)  +\int_{\Mntrap(\tau_1,\tau_2)} r^{\de-3} |\dk^{\le  k+1} \psi |^2,\\
E^k_\de[\psi](\tau)&=\int_{\Si(\tau)}  \Big(r^\de\big(|  \nab_4\dk^{\le k} \psi |^2 +r^{-2}|\dk^{\le k}\psi|^2\big) +|\nab \dk^{\le k}\psi|^2 + r^{-2}|\nab_3 \dk^{\le k}\psi|^2\Big),\\
 F^k_\de[\psi](\tau_1, \tau_2) &= \int_{\Si_*(\tau_1,\tau_2)}   \Big( r^{\de}\big(|\nab_4 \dk^{\le k}\psi|^2+|\nab \dk^{\le k}\psi|^2+r^{-2}| \dk^{\le k}\psi|^2\big)+|\nab_3 \dk^{\le k}\psi|^2\Big)\\
 &+\int_{\AA(\tau_1, \tau_2) } |\dk^{\le k+1 } \psi|^2.
 \end{split}
  \eeaa
 Throughout the chapter we  set   $\tau_1=1$ and $\tau_2=\tau_*$  and we simplify our expressions  by   writing  $\MM=\MM(1, \tau_*)$, $\AA=\AA(1, \tau_*)$, $\Si_*=\Si_*(1, \tau_*)$, $B^k_\de=B^k_\de(1, \tau_*)$,   $F^k_\de=F^k_\de(1, \tau_*)$.  We  introduce the  short hand notation
 \beaa
 \bsplit
\BEF_\de^k[\psi]:&=B_\de^k[\psi]+\sup_{\tau\in[1,\tau_*]}E_\de^k[\psi](\tau) +F_\de^k[\psi],\\
\EF_\de^k[\psi]:&=\sup_{\tau\in[1,\tau_*]}E_\de^k[\psi](\tau) +F_\de^k[\psi].
\end{split}
\eeaa
  We also recall the $\NN$  norm
 \beaa
 \NN^k[\psi,  N] &:=& \NNmor^k[\psi,  N]+ \NNred^k [\psi,  N]+ \NNen^k [\psi,  N]
 \eeaa
where
\beaa
 \NNmor^k[\psi,  N] &:=&\int_{\MM(\tau_1, \tau_2)} \big(|\dk^{\le k}  \nab_{\Rhat} \psi|+r^{-1}|\dk^{\le k}\psi|\big) | \dk^{\le k}N|,\\
 \NNred^k[\psi,  N] &:=& \int_{\MMred(\tau_1, \tau_2)} \big|\dk^{\le k+1} \psi \big| \big| \dk^{\le k} N \big|,\\
 \NNen^k [\psi,  N] &:=& \left|\int_{\MM(\tau_1, \tau_2)} \nab_{\That_\de} \dk^{\le k} \psi \c \dk^{\le k} N\right|,
 \eeaa
Here,  $\MMred$ is  the region of $\MM$ where  $r\le  r_+(1+2 \de_{red} )$ with  $\de_{red} $   a universal constant  verifying     $\de_\HH \ll \de_{red} \ll  m-|a|$, and $\That_\de$ is given in Definition \ref{definitionThat_de}. Also, the $\NN_\de$ norm is given by
  \beaa
 \NN^k_\de[\psi,  N] &:=& \NN^k [\psi,  N] +   \int_{\Mext} r^\de \Big(\big|\nab_4\dk^{\le k} \psi \big|+r^{-1}\big|\dk^{\le k} \psi \big|\Big) \big| \dk^{\le k} N \big|.
 \eeaa

These norms will be used in this chapter to estimate $\Pc$.   With these goal in mind we need to compare the  norms $\BEF$  with   the standard curvature  norms $\Rk_k$,
 see  section \ref{subsection:MainNormsM8}.  
 \begin{lemma}
 \lab{Lemma:comparisonofnormsBEF-Rk}
 We have
 \bea
 \lab{induction:hypoth-Pc}
 B^{k-1}_\de[ r^2\Rc_g] +B^{k-1}_\de[\Rc_b] \les \Rk_k^2, \qquad   \EF^{k-1}_\de[ r^2 \Rc_g] +\EF^{k-1}_\de[\Rc_b] \les \Rk_k\Rk_{k+1}. 
 \eea
 \end{lemma} 
 
 \begin{proof}
   In view of the definition of $\BEF$ norms, we have
  \beaa
  \bsplit
B^k_\de[r^2\Rc_g] =&  \int_{\MM_{trap}}\big(| \nab_\Rhat \dk^{\le k}\Rc_g|^2 + | \dk^{\le k}\Rc_g|^2\big)  +\int_{\Mntrap}\Big(r^{1+\de} |\dk^{\le k+1}\Rc_g|^2 +r^{3-\de} |\dk^{\le k}\nab_3\Rc_g|^2\Big),\\
E^k_{\de}[r^2\Rc_g] =&\sup_{\tau}  \left(\int_{\Si(\tau)}r^4\Big(r^{\de}\big(|  \nab_4  \dk^{\le k}\Rc_g |^2 +r^{-2}| \dk^{\le k}\Rc_g|^2\big)+|\nab \dk^{\le k}\Rc_g|^2 + r^{-2}|\nab_3 \dk^{\le k}\Rc_g|^2\Big)\right),\\
 F^k_{\de}[r^2\Rc_g] =& \int_{\Si_*}   r^4\Big(r^{\de}\big( |\nab_4 \dk^{\le k}\Rc_g|^2+|\nab \dk^{\le k}\Rc_g|^2 +r^{-2}| \dk^{\le k}\Rc_g|^2\big)+|\nab_3  \dk^{\le k}\Rc_g|^2\Big) \\
 &+\int_\AA   |\dk^{\le k+1}\Rc_g|^2.
 \end{split}
  \eeaa
Hence, in view of the definition of the $\Rk_k$ norms in  section
 \ref{subsection:MainNormsM8} and \eqref{eq:auxilliarynormsforRb:boundary}, and the fact that $\Rc_g=(A, B, \Pc)$, we have $B^{k-1}_\de[ r^2\Rc_g] \les \Rk_k^2$ and $\EF^{k-1}[ r^2\Rc_g]\les  \Rk_k\Rk_{k+1}$  as stated. The estimate for $\Rc_b$ can be proved similarly. 
 \end{proof} 
 
We will also need to  compare  the norms $\BEF$  for  Ricci coefficients  with those 
given by $\Sk_k$.
\begin{lemma}
\lab{Lemma:comparisonofnormsBEF[Gac)-Sk}
The following estimates hold true
\beaa
B^{k-1}_\de[\Ga_b'] + B^{k-1}_\de[r^{-\dt} \Ga_b] \les \Sk_k^2, \qquad \EF^{k-1}_\de[\Ga_b'] +\EF^{k-1}_\de[r^{-\dt} \Ga_b] \les \Sk_k  \Sk_{k+1}. 
\eeaa
\end{lemma}

\begin{proof}
The proof  follows easily from Lemma  \ref{lemma:auxilliarynormsforGa_b}.
\end{proof}


\subsection{Statement of main results on the control of $\Pc$}


 We assume that all the bootstrap, initial data assumptions and  the induction hypothesis made in section \ref{section:MainAss.ThmM8} hold true.  The goal of this chapter is to derive the following  Morawetz energy estimates   for $\Pc $ in  $\MM=\MM(1, \tau_*)$.
\begin{theorem}[Morawetz-Energy for $\Pc$] 
\lab{theorem:Morawetz-EnergyPc}
The following estimates hold true in $\MM=\MM(1, \tau_*)$, for sufficiently small $a$,
\bea
\BEF_\de^J[r^2\Pc]   &\les& r_0^{15}\Big( \Sk_{J+1} \Sk_J +\Rk_{J+1}\Rk_J +\ep_J^2+\ep_0^2\Big).  
\eea
\end{theorem}

The proof is based on the conditional weighted  estimates   for scalar  wave equations of Proposition \ref{Prop:scalarwavePsi-M8-chap10}  which  we recall below.

\begin{proposition}
\lab{Prop:scalarwavePsi-M8}
Let $\psi$ be a solution to the following scalar wave equation 
\bea
\square_\g\psi+V\psi &=& N,
\eea 
where $V$ is real and satisfies $V=O(r^{-3})$ for $r$ large, in a spacetime  $\MM=\MM(1, \tau_*)$   verifying  the assumptions of section \ref{section:MainAss.ThmM8}. Then:
  \begin{enumerate}
\item The following conditional Morawetz estimates hold true in $\MM$
\bea
\lab{eq:Estimatesforpsi-M8-1}
\bsplit
   B^k_{\de}[\psi]  &\les   \EF_\de^{k}[\psi]  +B ^{k-1}_{\de}[\psi]  +\int_{\MM_{trap}} |\dk^{\le k}\psi|^2 + \NNmor^k[\psi,  N]+ \NNred^k [\psi,  N]\\
   & +   \int_{\Mext} r^\de \Big(\big|\nab_4\dk^{\le k} \psi \big|+r^{-1}\big|\dk^{\le k} \psi \big|\Big) \big| \dk^{\le k} N \big|.
\end{split}
\eea

 \item The following conditional  Energy-Morawetz estimates hold true
 \bea
\lab{eq:Estimatesforpsi-M8-2}
\bsplit
  B ^k_{\de}[\psi] + \sup_{\tau\in[1,\tau_*]}E^k _{\de}[\psi](\tau)+ F ^k_{\de}[\psi] 
& \les   E^k _{\de}[\psi](0)       + \BEF^{k-1}_\de[\psi] \\
& +\int_{\MM_{trap}} |\dk^{k}\psi|^2  +\NN^k_\de[\psi,  N]. 
    \end{split}
\eea
 \end{enumerate}
\end{proposition}

\begin{remark}
 Note  that    both  estimates    are conditional   on the control of  the  lower order term of $\Mor^{k}[\psi]$ in $\MM_{trap}$, i.e. the term $\int_{\MM_{trap}} |\dk^{k}\psi|^2 $.  
\end{remark}

We  will also  make use of the  following control of $\qfb$ which follows from a variant of Theorem \ref{theorem:unconditional-result-final-psib}.
\begin{theorem}
\lab{Thm:Estimates-forqf}
We have,  for $k\le k_L-3$,
\bea
\lab{eq:Estimatesforqf-M8}
 \sup_{\tau\in[1,\tau_*]}E\,^{k}_{\de}[\qfb](\tau) +   B ^{k}_{\de}[\qfb]  +   F ^{k} _{\de}[\qfb] &\les \ep_0^2. 
\eea
\end{theorem} 

\begin{proof}
Recall from \eqref{eq:specialidentityforthegloablframeofMMinpartIII} that the global frame of $\MM$ satisfies in Part III the identities $\Xi=\Hbc=0$ for $r\geq r_0+1$. In particular, $\undpsi=\Re(\qfb)$ satisfies the generalized RW equation \eqref{eq:defintionofCb1andCb2fordefintionqfb:chap12}. One can then easily adapt the proof of Theorem \ref{theorem:unconditional-result-final-psib} to obtain that the following  holds true, for   $s\le k_L-3$,  
 \beaa
        \BEF_\de^s[\undpsi, \Ab](1, \tau_*) \les  E_\de^s[\undpsi, \Ab](1) +\NN_\de^s[\undpsi, N_{\err} ](1, \tau_*),
 \eeaa 
  where $N_{\err}$ is defined in \eqref{eq:N_err-Ab}, i.e.  
  \beaa
N_{\err} &=& r^2\dkb^{\leq 1}\aa\c\dk^{\leq 2}(\a, \b)+\dk^{\leq 3}(\Ga_b\c\Ga_g).
\eeaa
 Using the  control of the initial data provided by \eqref{eq:controlofinitialdataforThM8} and the definition of $\BEF_\de^s[\undpsi, \Ab](1, \tau_*)$ and of $\undpsi$, we infer, for   $s\le k_L-3$,  
 \beaa
        \BEF_\de^s[\qfb] \les  \ep_0^2 +\NN_\de^s[\undpsi, N_{\err} ](1, \tau_*).
 \eeaa  
Also, as established in the proof of Lemma \ref{lemma:DecayfortheNterm-undpsi}, see section \ref{section:Prooflemma:DecayfortheNterm-undpsi}, we have
 \beaa
    &&   \NN^s_\de[\undpsi,  N](\tau_1, \tau_2) \\
       &\les& \Big(BEF^s_\de[\undpsi](\tau_1, \tau_2)\Big)^{\frac{1}{2}}\left(\int_{\tau_1}^{\tau_2}\|\dk^{\leq s}N\|_{L^2(\Si_{trap}(\tau))}+\left(\int_{\MM(\tau_1, \tau_2)}r^{\de+1}|\dk^{\leq s}N|^2\right)^{\frac{1}{2}}\right).
 \eeaa        
Now, in view of the above form of $N_{\err}$ and the control provided by the bootstrap assumptions \eqref{eq:mainbootassforchapte13} \eqref{eq:auxlowderivativebootassdecayforchapte13}, we have, for   $s\le k_L-3$,
\beaa
&& \int_{\tau_1}^{\tau_2}\|\dk^{\leq s}N_{\err}\|_{L^2(\Si_{trap}(\tau))}+\left(\int_{\MM(\tau_1, \tau_2)}r^{\de+1}|\dk^{\leq s}N_{\err}|^2\right)^{\frac{1}{2}}\\
&\les& (\Sk_s+\Rk_s)\int_{\tau_1}^{\tau_2}\|\dk^{\leq \frac{s}{2}}(\Gac, \Rc)\|_{L^2(\Si_{trap}(\tau))} +\ep\int_{\MM(\tau_1, \tau_2)}\Big(r^{3+\de}|\dk^{\leq s+2}(\a, \b)|^2\\
&&+r^{-2-2\dec+\de}|\dk^{\leq s+1}\aa|^2+r^{-1+\de}|\dk^{\leq s+3}\Ga_g|^2+r^{-3+\de}|\dk^{\leq s+3}\Ga_b|^2\Big)\\
&\les& \ep^2 \int_{\tau_1}^{\tau_2}\frac{d\tau}{\tau^{1+\dec}}+\ep(\Sk_{s+3}+\Rk_{s+3})\les \ep^2\les\ep_0,
\eeaa
where we used the fact that $s+3\leq k_L$ and $\de<\dec<\dt$. We infer
 \beaa
  \NN^s_\de[\undpsi,  N](\tau_1, \tau_2)   &\les& \ep_0\Big(BEF^s_\de[\undpsi](\tau_1, \tau_2)\Big)^{\frac{1}{2}}
 \eeaa  
 and hence, for   $s\le k_L-3$, 
 \beaa
        \BEF_\de^s[\qfb] \les  \ep_0^2
 \eeaa
as stated. This concludes the proof of Theorem \ref{Thm:Estimates-forqf}.
\end{proof}

These results will be used together with the following lemma  relating  $\DDc\hot \DDc\Pc$ to  $\qfb$, to prove Theorem \ref{theorem:Morawetz-EnergyPc}.
 
  \begin{lemma}
  \lab{Lemma:Formula-qf-DDhotDDPc}
  The following relation between $\qfb$ and $\Pc$ holds true.
\bea
 \lab{eq-Formula-qf-DDhotDDPc:P3}
 \bsplit
 \qfb =& \frac{1}{2}\ov{q}q^3 \DDc\hot \DDc\Pc +  \dk^{\leq 1}\Ga_b'  +O(a)\dk^{\leq 1}\Rc_b + O(ar) \dk^{\le 1}\Pc +r\dk^{\leq 1}(\Ga_b\c \Rc_b),
 \end{split}
 \eea
 where we recall that $\Ga_b'=\Ga_b\setminus\{\Xib\}$.
  \end{lemma}
  
\begin{proof}
This follows immediately from Proposition \ref{prop:Formula-qf-DDhotDDPc-ch5} and the fact that $r\Bb\in \Rc_b$, $\Ga_g'=\Ga_g\setminus\{\trXbc\}$ and the fact that we identity $\Rc_g$ with $r^{-2}\Rc_b$, $\Ga_g$ with $r^{-1}\Ga_b$, and $\Ga_g'$ with $r^{-1}\Ga_b'$ in Part III.
\end{proof}

The main part in  the  proof of Theorem \ref{theorem:Morawetz-EnergyPc} is to derive the following
  result  for  the pair of scalars $\psi=q^2(\T(P), \Z(P) )$ stated below.
 \begin{proposition}
 \lab{proposition:Morawetz-Energy-psi}
 The pair of scalars  
 \beaa
 \psi:=q^2(\T(P), \Z(P) )
 \eeaa
 verify the following estimates
 \bea
 \BEF_\de^{J-1}[\psi] \les r_0^{10}\Big( \Sk_{J+1} \Sk_J +\Rk_{J+1}\Rk_J +\ep_J^2+\ep_0^2\Big).
 \eea
 \end{proposition}


\section{Proof of the Morawetz-Energy estimates for $\Pc$}


The goal of this section is to prove Theorem \ref{theorem:Morawetz-EnergyPc} on Morawetz-Energy estimates for $\Pc$. We start by proving Proposition \ref{proposition:Morawetz-Energy-psi} on Morawetz-Energy estimates for the pair of scalars $\psi=q^2(\T(P), \Z(P) )$.


\subsection{Proof of Proposition \ref{proposition:Morawetz-Energy-psi}}
\lab{sec:proofofproposition:Morawetz-Energy-psi}


The proof of Proposition \ref{proposition:Morawetz-Energy-psi} proceeds in several steps.

{\bf Step 0.}   We start with the equation  of Lemma \ref{LEMMA:WAVEEQP} for $P$ 
\bea
\lab{eq:P-WaveEq-M8}
\nn\square_\g P &=&  \tr X\nab_3P +\ov{\tr\Xb}\nab_4P  -\ov{H}\c\DD P  - \Hb\c\ov{\DD}P + V P  \\
&& + r^{-3}\dk^{\leq 1}(\Ga_b \c\Rc_b)-\Ab\c \ov{A},
\eea
see also Remark \ref{rem:error-terms-square-P}, and recall that we identify $\Rc_g$ with $r^{-2}\Rc_b$ and $\Ga_g$ with $r^{-1}\Ga_b$.

\begin{remark}\lab{rmk:problemwhenusingtriviallinearizationbycoordinates:partIII}
We will need to linearize \eqref{eq:P-WaveEq-M8}. A possibility consists in deriving a wave equation for $\Pc$, where we recall that  $\Pc=P+\frac{2m}{q^3}$. In view of \eqref{eq:P-WaveEq-M8}, this leads to 
\bea
\lab{eq:P-WaveEq-linearized}
\bsplit
\square_\g \Pc &=  \tr X\nab_3\Pc  +\ov{\tr\Xb}\nab_4\Pc  -\ov{H}\c\DD \Pc   - \Hb\c\ov{\DD}\Pc + V \Pc    + r^{-4}  \dk^{\le 1} \Ga_b \\
& + r^{-3}\dk^{\leq 1}(\Ga_b \c\Rc_b)-\Ab\c \ov{A}.
\end{split}
\eea
This linearization   leads however  to terms   of the  form $\dk^{J+1} \Gac$   in  estimates for  $\dk^{J+1} \Pc$, so that the iteration assumption \eqref{eq:iterationassumptiondiscussionThM8:bis} on $\Sk_J$ cannot be used.
\end{remark}

To avoid the linearization issue outlined in Remark \ref{rmk:problemwhenusingtriviallinearizationbycoordinates:partIII}, we  proceed by the  linearization  procedure  discussed in Step 1 below which consists in deriving  Morawetz-Energy estimates  for  the linearized quantity  
\bea
\psi:= q^2(\T P, \Z P).
\eea

 To apply the induction hypothesis we need to compare the Morawetz-Energy norms   $\BEF^{k} _\de[\psi]$   with $\Rk_{k+1}^2[\Pc]$. 
 \begin{lemma}
 \lab{Lemma:comparisonofnormsBEF-Rk-forpsi}
 The induction hypothesis   \eqref{eq:iterationassumptiondiscussionThM8:bis} implies
 \bea
 \lab{induction:hypoth-psi}
 B^{J-2} _\de[\psi] \les \ep_J^2, \qquad   \EF_\de^{J-2}[\psi] \les   \ep_J^2 +\Rk_J\Rk_{J+1}. 
 \eea
 \end{lemma}
 
 \begin{proof}
 Note that
\beaa
q^2 \T P&=& q^2 \T\left(\Pc-\frac{2m}{q^3} \right)= q^2 \T \Pc + r^{-1}\Ga_b, \\
q^2 \Z P&=& q^2 \Z\left(\Pc-\frac{2m}{q^3}\right)= q^2 \Z \Pc + r^{-1}\Ga_b.
\eeaa
We deduce, in view of Lemmas  \ref{Lemma:comparisonofnormsBEF-Rk} and  \ref{Lemma:comparisonofnormsBEF[Gac)-Sk}  and the induction hypothesis.
\beaa
B^{J-2} _\de[\psi] &\les&   B_\de^{J-2}[r^2(\T \Pc, \Z\Pc)]  +\BEF_\de^{J-2} [r^{-1}\Ga_b]\\
 &\les&   B_\de^{J-1}[r^2\Pc]  +\BEF_\de^{J-2} [r^{-1}\Ga_b]\\
&\les& \Rk_J^2+   \BEF_\de^{J-2} [r^{-1}\Ga_b]\les \ep_J^2+ \Sk_J\Sk_{J-1}\\
&\les&\ep_J^2
\eeaa 
and similarly 
\beaa
\EF_\de^{J-2}[\psi] &\les&   \EF_\de^{J-2}[r^2(\T \Pc, \Z\Pc)  ] +\BEF_\de^{J-2} [r^{-1}\Ga_b]\\
 &\les&    \EF_\de^{J-1}[r^2\Pc   ]  +\BEF_\de^{J-2} [r^{-1}\Ga_b]\\
&\les& \Rk_J\Rk_{J+1}+   \BEF_\de^{J-2} [r^{-1}\Ga_b]\les \Rk_J\Rk_{J+1}+ \Sk_J\Sk_{J-1}\\
&\les&\ep_J^2 +\Rk_J\Rk_{J+1}
\eeaa  
as stated.
\end{proof}

{\bf Step 1.} Next, we derive a wave  equation for the  pair of scalars $\psi= q^2(\T P, \Z P)$.
\begin{lemma}\lab{LEMMA:WAVEEQUATIONFORMODQ2TPTZ}
The pair of scalars $\psi= q^2(\T P, \Z P)$ satisfies a wave equation 
of the schematic  form\footnote{Recall that $\Rc_b$ is a curvature term which behaves like $\Ga_b$ in terms  of  powers of $r$.}
\bea
\lab{eq:Schematicswave|q|^2TP-TZ}
\square_\g \psi&=& W  \psi +   r^{-2} \dk^{\leq 1}\Ga_b + r^{-1}\dk^{\leq 2}(\Ga_b \c\Rc_b)-r^2\dk^{\leq 1}(\Ab\c \ov{A})
\eea
where $W$ is a complex potential of the form
\bea
\lab{eq:DefineN^{le 3}}
\Re(W)&=&O( m r^{-3}),   \qquad  \Im(W)=O(ma r^{-4}).
\eea
\end{lemma}

\begin{proof}
This is done by making use of the commutation formulas \eqref{eq:commTZsquare}
\beaa
\bsplit
 \,[ \T, \square_\g] \psi&=\dk \big(\Ga_g \c \dk \psi\big)+\Ga_b \c \square_\g\psi ,\\
   \,[ \Z, \square_\g] \psi &= \dk \big(\Ga_g \c \dk \psi\big)+r\Ga_b \c \square_\g\psi,
 \end{split}
 \eeaa
as well as the renormalization  Lemma \ref{Le:squareq^2Psi}. See section \ref{section:DetailesStep1}.
\end{proof}

{\bf Step 2.} We  apply  the first  conditional estimate of Proposition \ref{Prop:scalarwavePsi-M8}
 to equation \eqref{eq:Schematicswave|q|^2TP-TZ}
 to  derive  a conditional  Morawetz estimate  for $\psi$ of the form.
\bea
\lab{estimate:conditional-B^{J-1}[psi]}
B_\de^{J-1}[\psi]&\les&  \EF_\de^{J-1}[\psi] +\int_{\MM_{trap}} | \dk^{J-1}\psi|^2  + \ep^2_J+\ep_0^2.
\eea

 \begin{remark}
 Note that we cannot estimate the energy flux $\EF$   at this step because of the  presence of the complex potential  $W$ in  \eqref{eq:Schematicswave|q|^2TP-TZ}.
 \end{remark}

The proof of \eqref{estimate:conditional-B^{J-1}[psi]}    is  a      straightforward  application of  the estimate \eqref{eq:Estimatesforpsi-M8-1}  of   Proposition \ref{Prop:scalarwavePsi-M8} which yields, for $k=J-1$, 
\beaa
 B ^{J-1}_{\de}[\psi] 
&\les &  \EF_\de^{J-1}[\psi]  +B ^{J-2}_{\de}[\psi] +\int_{\MM_{trap}} |\dk^{J-1}\psi|^2 \\
&& + \NNmor^{J-1}[\psi,  N]+ \NNred^{J-1}[\psi,  N]\\
   && +   \int_{\Mext} r^\de \Big(\big|\nab_4\dk^{\le J-1} \psi \big|+r^{-1}\big|\dk^{\le J-1} \psi \big|\Big) \big| \dk^{\le J-1} N \big|,
\eeaa
 where we have, in view of \eqref{eq:Schematicswave|q|^2TP-TZ}, $V=-\Re(W)$ and 
\beaa 
 N :=  i\Im(W)  \psi +   r^{-2} \dk^{\leq 1}\Ga_b + r^{-1}\dk^{\leq 2}(\Ga_b \c\Rc_b)-r^2\dk^{\leq 1}(\Ab\c \ov{A}). 
 \eeaa
 We have 
\beaa
 && \NNmor^{J-1}[\psi,  N]+ \NNred^{J-1}[\psi,  N]\\
   &&+   \int_{\Mext} r^\de \Big(\big|\nab_4\dk^{\le J-1} \psi \big|+r^{-1}\big|\dk^{\le J-1} \psi \big|\Big) \big| \dk^{\le J-1} N \big|\\
       &\les& \Big(B ^{J-1}_{\de}[\psi]\Big)^{\frac{1}{2}}\left(\int_{\MM}r^{\de+1}|\dk^{\le J-1}N|^2\right)^{\frac{1}{2}}.
\eeaa
and, using the definition of $N$, the fact that $\Im(W)=O(ma r^{-4})$,  the induction hypothesis for the linear term involving $\Ga_b$   and the  bootstrap  assumptions for the nonlinear one,   and Lemma \ref{lemma:auxilliarynormsforGa_b}, we estimate 
 \beaa
 \int_{\MM}r^{\de+1}|\dk^{\le J-1}N|^2 &\les& |a|B ^{J-1}_{\de}[\psi]+\int_{\MM} r^{-3+\de} |\dk^{\le J} \Ga_b|^2 +\ep^2\int_{\MM} r^{-3+\de} |\dk^{\le J+1} \Rc_b|^2\\
 &&+\ep^2\int_{\MM} r^{-2+\de} |\dk^{\le J+1} \Ga_b|^2+\ep^2\int_{\MM} r^{3+\de} |\dk^{\le J+1}A|^2\\
 &\les& |a|B ^{J-1}_{\de}[\psi]+\Sk_{J}^2+ \ep^2(\Sk_{J+1}+\Rk_{J+1})^2\\
 &\les& |a|B ^{J-1}_{\de}[\psi]+\ep_J^2+ \ep_0^2.
 \eeaa
 We deduce
 \beaa
 B ^{J-1}_{\de}[\psi] 
&\les &  \EF^{J-1}_\de[\psi]  +B ^{J-2}_{\de}[\psi] +\int_{\MM_{trap}} |\dk^{J-1}\psi|^2 \\
&& +\Big(B ^{J-1}_{\de}[\psi]\Big)^{\frac{1}{2}}\left(|a|B ^{J-1}_{\de}[\psi]+\ep_J^2+ \ep_0^2\right)^{\frac{1}{2}},
\eeaa
and hence, using also the fact that $B^{J-2} _\de[\psi]\les   \ep_J^2$ in view of Lemma \ref{Lemma:comparisonofnormsBEF-Rk-forpsi}, we obtain
\beaa
B_\de^{J-1}[\psi] &\les& |a|B_\de^{J-1}[\psi]+ \EF^{J-1}_\de[\psi] +\int_{\MM_{trap}} | \dk^{J-1}\psi|^2  + \ep^2_J+\ep_0^2.
\eeaa
For $a$ small enough, this yields \eqref{estimate:conditional-B^{J-1}[psi]} as stated.

 {\bf Step 3.}  To control  the  energy-flux term $\EF^{J-1}_\de[\psi]$  in \eqref{estimate:conditional-B^{J-1}[psi]},   we   take a circuitous route  by   commuting the  wave equation for $\psi$  with $\nab^2_\Rhat $. We rely on on the following commutation lemma.
  
\begin{lemma}
\lab{LEMMA:SQUARENAB^2_RHAT}
Assume  $\square_\g \psi= N$. Then, we have
\bea
\bsplit
\square _\g \nab^2_\Rhat \psi&=   O(r^{-1} )\nab_\Rhat \lap\psi +O(ar^{-3})\nab_\Rhat\dk^{\leq 2} \psi +O(r^{-2}) \lap \psi +O(r^{-2}) \nab_\Rhat \psi\\
&+ O(a r^{-3}) \dk^{\le 2} \psi  + r^{-1}\dk^{\leq 3}\big( \Ga_b \c \psi\big) + \nab^2 _\Rhat N+r^{-1}\nab_\Rhat  N + r^{-2}  N.
\end{split}
\eea
\end{lemma}
 
 \begin{proof}
 See section \ref{section:DetailesStep3}.
 \end{proof} 

\begin{remark}
The reason for commuting  the  wave equation for $\psi$  with $\nab^2_\Rhat $ is to ensure that all linear terms involving top order curvature components on the RHS of the wave equation of Lemma \ref{LEMMA:SQUARENAB^2_RHAT} contain at least on $\nab_{\Rhat}$ derivative. 
\end{remark}

We rewrite the wave equation for $\psi$ in \eqref{eq:Schematicswave|q|^2TP-TZ} as 
\beaa
\square_\g \psi &=& N, \\
N&:=& W  \psi +   r^{-2} \dk^{\leq 1}\Ga_b + r^{-1}\dk^{\leq 2}(\Ga_b \c\Rc_b)-r^2\dk^{\leq 1}(\Ab\c \ov{A}),
\eeaa
and compute, using $W=O(r^{-3})$, 
\beaa
\nab^2 _\Rhat N+r^{-1}\nab_\Rhat  N + r^{-2}N  &=& W\nab_\Rhat^2\psi  +O(r^{-3})\dk^{\leq 1}\psi + r^{-2} \dk^{\leq 3}\Ga_b \\
&& + r^{-1}\dk^{\leq 4}(\Ga_b \c\Rc_b) -r^2\dk^{\leq 3}(\Ab\c \ov{A}).
\eeaa
Commuting $\square_\g \psi = N$ with $\nab^2_\Rhat $, and relying on Lemma \ref{LEMMA:SQUARENAB^2_RHAT}, we obtain the following wave equation for $\nab_\Rhat^2\psi$
\bea
\lab{eq:wave-nab^2_Rhatpsi}
\bsplit
\square_\g(\nab_\Rhat^2\psi)&=  W \nab_\Rhat^2 \psi      +O(r^{-2}) \lap \psi      +O(r^{-1})\nab_\Rhat\lap  \psi  +O(ar^{-3})\Big(\nab_\Rhat\dk^{\le  2} \psi    + \dk^{\le  2} \psi\Big)   \\
&+O(r^{-2}) \dk^{\le 1} \psi    +N^{\le 3},\\
 N^{\le 3}&= r^{-2}   \dk^{\leq 3 }\Ga_b  + r^{-1}\dk^{\leq 4}(\Ga_b \c\Rc_b)-r^2\dk^{\leq 3}(\Ab\c \ov{A}).
 \end{split}
\eea

{\bf Step 4.}  We   apply   the  second conditional estimate of Proposition  \ref{Prop:scalarwavePsi-M8}
 to derive the following lemma.
 \begin{lemma}\lab{LEMMA:ESTIMATE:BEF[NAB^2_RHAT].}
 For any $0<\de_1\leq 1$, we have
  \bea
  \lab{estimate:BEF[nab^2_Rhat].}
   \begin{split}
    \BEF ^{J-3}_{\de}[\nab^2_\Rhat \psi]& \les ( \de_1+ |a| )  B_\de^{J-1}[\psi] +\de_1^{-1}  B_\de^{J-3}[ r\lap \psi] +\de_1^{-1}\Big(\ep_J^2 +\ep_0^2\Big) +\Rk_J\Rk_{J+1}.
 \end{split}
   \eea
\end{lemma}

\begin{proof}
See section \ref{section:DetailesStep4}.
\end{proof}

\begin{remark}
In the proof of Lemma \ref{LEMMA:ESTIMATE:BEF[NAB^2_RHAT].} in section \ref{section:DetailesStep4}, when applying the   second conditional estimate of Proposition  \ref{Prop:scalarwavePsi-M8},   we will have  to take into account 
 the additional term $\int_{\MM_{trap}} |\dk^{J-3} \nab_\Rhat^2 \psi|^2 $  on the right hand side.   Fortunately, thanks to the $\nab_\Rhat$ derivative, this term  is bounded by $B_\de^{J-1}[\psi] $   and  thus  can be  controlled by  the induction hypothesis.
\end{remark}

{\bf Step 5.} Recall the following commutator formula, see Lemma \ref{lemma:commutator-triangle}, 
\beaa
\, [|q|^2\lap,|q|^2 \square_\g]\psi&=& |q|^2\Big[O(a^2r^{-3}) \dk^{\leq 2}\psi+ \dk^2 \big( \Ga_g \c \dk \psi)+ \Ddot_3 \dk \big(|q|^2 \xi \c \Ddot_a \psi \big)\Big].
\eeaa
In view of \eqref{eq:specialidentityforthegloablframeofMMinpartIII}, we have in particular $\Xi\in r^{-1}\Ga_g$, and hence
\beaa
\, [|q|^2\lap,|q|^2 \square_\g]\psi&=& |q|^2\Big[O(a^2r^{-3}) \dk^{\leq 2}\psi+ \dk^2 \big( \Ga_g \c \dk \psi)\Big].
\eeaa
We infer
\beaa
\square_\g (|q|^2\lap \psi) &=& \frac{1}{|q|^2}\Big(|q|^2\lap(|q|^2\square_\g\psi)- [|q|^2\lap,|q|^2 \square_\g]\psi\Big)\\
&=& \lap(|q|^2\square_\g\psi) + O(a^2r^{-3}) \dk^{\leq 2}\psi+ \dk^2 \big( \Ga_g \c \dk \psi).
\eeaa
Plugging the wave equation for $\psi$, see \eqref{eq:Schematicswave|q|^2TP-TZ}, in the RHS, and using $W=O(r^{-3})$, we infer
 \bea
\lab{eq:wave-lappsi}
\square_\g (|q|^2\lap \psi)&=& W( |q|^2\lap  \psi)   + O(ar^{-4 }) \dk^{\leq 2}\psi + O(r^{-3 }) \dk^{\leq 1}\psi  +N^{\le 3},
\eea
where $N^{\le 3}$ is as in \eqref{eq:wave-nab^2_Rhatpsi}. We can then proceed  as in Step 2, using the first  conditional estimate  of Proposition  \ref{Prop:scalarwavePsi-M8}, to derive  the estimate
\bea
\lab{estimate:conditional-B^{J-3}[lap-psi]}
B_\de^{J-3}[ r^2 \lap \psi] \les  |a|B_\de^{J-1}[\psi]+\EF_\de^{J-3}[ r^2 \lap\psi]+\int_{\MM_{trap}} |\dk^{J-3} \lap \psi|^2  + \ep^2_J+\ep_0^2.
\eea

 Next, using the  Hodge  type estimate of Corollary \ref{Cor:HodgeThmM8}, integrated  on  $\MM$,  we   derive
 \beaa
 B_\de^{J-3}[ r^2 \nab^2  \psi]&\les  B_\de^{J-3}[ r^2 \lap \psi]+O(a, \ep)  B_\de^{J-1}[\psi].
 \eeaa
 Combining with \eqref{estimate:conditional-B^{J-3}[lap-psi]} we  thus infer
 \bea
 \begin{split}
 \lab{estimate:conditional-B^{J-3}[nab^2-psi]}
  B_\de^{J-3}[ r^2 \nab^2  \psi]&\les  \EF_\de^{J-3}[ r^2 \lap\psi] +O(a, \ep)  B_\de^{J-1}[\psi]  +\int_{\MM_{trap}} |\dk^{J-3} \lap \psi|^2\\
  &+ \ep^2_J+\ep_0^2.
  \end{split}
 \eea

 {\bf Step 6.} Combining  \eqref{estimate:BEF[nab^2_Rhat].}  with  \eqref{estimate:conditional-B^{J-3}[nab^2-psi]}  we derive, for any $0<\de_1\leq 1$, 
   \beaa
    \BEF ^{J-3}_{\de}[\nab^2_\Rhat \psi]& \les& ( \de_1+ |a| )  B_\de^{J-1}[\psi] +\de_1^{-1}  B_\de^{J-3}[ r \lap \psi] +\de_1^{-1}\Big(\ep_J^2 +\ep_0^2\Big) +\Rk_J\Rk_{J+1}\\
    & \les & ( \de_1+ |a| )  B_\de^{J-1}[\psi]  +\de_1^{-1}  \EF_\de^{J-3}[ r^2 \lap \psi] +\de_1^{-1} O(a, \ep)  B_\de^{J-1}[\psi] \\
   &&+\de_1^{-1} \int_{\MM_{trap}} |\dk^{J-3} \lap \psi|^2+  \de_1^{-1} \big( \ep_J^2 +\ep_0^2\big) +\Rk_J\Rk_{J+1}.
    \eeaa
    Hence, for any $0<\de_1\leq 1$,
    \bea
   \lab{estimate:conditional-B^{J-3}[nab^2R-nab^2-psi]}
    \bsplit
       \BEF ^{J-3}_{\de}[(\nab^2_\Rhat,  r^2 \nab^2 )  \psi]& \les   \Big( \de_1+ O(a,  \ep) \de_1^{-1} \Big)  B_\de^{J-1}[\psi] +\de_1^{-1}  \EF_\de^{J-3}[ r^2 \lap \psi]\\
       &+\de_1^{-1} \int_{\MM_{trap}}  \ |\dk^{J-3} \lap \psi|^2 + \de_1^{-1} \big( \ep_J^2 +\ep_0^2\big) +\Rk_J\Rk_{J+1}.
       \end{split}
    \eea

{\bf Step 7.}  We  make use  of the following consequence of  Lemma \ref{lemma:squaredkpsi-perturbations}
\beaa
-\nab_\That^2\psi  +\nab^2_\Rhat\psi &=& O(1) \square \psi+O(1)\lap \psi   +O(r^{-1} )\nab_\Rhat \psi    
   +O(a r^{-2} )\nab\psi + \Ga_g \c \dk \psi
\eeaa
to derive the   estimates
\beaa
\EF_\de^{J-3}[\nab^2_\That \psi]&\les& \EF_\de^{J-3}[\nab^2_\Rhat \psi]+ \EF_\de^{J-3}[ \lap \psi]+\ep_J^2+\ep_0^2,\\
B_\de^{J-3}[\nab^2_\That \psi]&\les&   B_\de^{J-3}[\nab^2_\Rhat \psi]+ B_\de^{J-3}[ \lap \psi]+\ep_J^2+\ep_0^2.
\eeaa
Combining this with \eqref{estimate:conditional-B^{J-3}[nab^2R-nab^2-psi]}, we infer that,  for any $0<\de_1\leq 1$,
\beaa
\bsplit
\BEF_\de^{J-3}[\big( \nab^2_\That, \nab^2_\Rhat, r^2 \nab^2\big)  \psi]& \les   \Big( \de_1+ O(a,\ep) \de_1^{-1} \Big)   B_\de^{J-1}[\psi]   +\de_1^{-1}  \EF_\de^{J-3}[r^2 \lap \psi] \\
& +\de_1^{-1} \int_{\MM_{trap}}  \ |\dk^{J-3} \lap \psi|^2 +\de_1^{-1}\Big(\ep_J^2 +\ep_0^2\Big) +\Rk_J\Rk_{J+1}.
\end{split}
\eeaa
Moreover, since $J\geq\frac{k_L}{2}$ and $\kl$ is large, we may assume that $J\geq 5$. In particular, we have
\beaa
\BEF_\de^{J-5}[\big( \nab_\That, \nab_\Rhat, r\nab\big)^{\leq 4}\psi] &\les& \BEF_\de^{J-3}[\big( \nab^2_\That, \nab^2_\Rhat, r^2 \nab^2\big)  \psi]+\BEF_\de^{J-2}[\psi]
\eeaa 
and hence,  for any $0<\de_1\leq 1$,
\beaa
\bsplit
\BEF_\de^{J-5}[\big( \nab_\That, \nab_\Rhat, r\nab\big)^{\leq 4}\psi] \les&   \Big( \de_1+ O(a,\ep) \de_1^{-1} \Big)   B_\de^{J-1}[\psi]   +\de_1^{-1}  \EF_\de^{J-3}[r^2 \lap \psi] \\
& +\de_1^{-1} \int_{\MM_{trap}}  \ |\dk^{J-3} \lap \psi|^2 +\de_1^{-1}\Big(\ep_J^2 +\ep_0^2\Big) +\Rk_J\Rk_{J+1}\\
&+\BEF_\de^{J-2}[\psi].
\end{split}
\eeaa
Together with \eqref{induction:hypoth-psi}, we infer,  for any $0<\de_1\leq 1$,
\bea\lab{estimate:conditional-B^{J-3}[nab^2R-nab_T-nab^2-psi]} 
\bsplit
\BEF_\de^{J-5}[\big( \nab_\That, \nab_\Rhat, r\nab\big)^{\leq 4}\psi] \les&   \Big( \de_1+ O(a,\ep) \de_1^{-1} \Big)   B_\de^{J-1}[\psi]   +\de_1^{-1}  \EF_\de^{J-3}[r^2 \lap \psi] \\
& +\de_1^{-1} \int_{\MM_{trap}}  \ |\dk^{J-3} \lap \psi|^2 +\de_1^{-1}\Big(\ep_J^2 +\ep_0^2\Big) +\Rk_J\Rk_{J+1}.
\end{split}
\eea

{\bf Step 8.} Next, we derive an estimate for $\BEF_{\de; r\geq r_0}^{J-2}[re_4\psi]$. To this end, we commute the wave equation for $\psi$ in \eqref{eq:Schematicswave|q|^2TP-TZ} with $r\nab_4$. We obtain  
\beaa
\square_\g\big(re_4\psi\big) &=& re_4\big(\square_\g\psi\big) - [re_4, \square_\g]\psi\\
&=& O(r^{-3})\dk^{\leq 1}\psi+ r^{-2} \dk^{\leq 2}\Ga_b + r^{-1}\dk^{\leq 3}(\Ga_b \c\Rc_b)-r^2\dk^{\leq 2}(\Ab\c \ov{A}) - [re_4, \square_\g]\psi,
\eeaa
where we used the fact that $W=O(r^{-3})$. Together with the commutator formula of Lemma \ref{LEMMA:COMMUTATOR-NAB3-NAB4-SQUARE}, we infer
\beaa
\bsplit
\square_\g\big(re_4\psi\big) =& \frac{1}{r}\nab_4(r\nab_4\psi)+N_{re_4},\\
N_{re_4} :=& O(r^{-2})\dkb^2\psi+O(r^{-2})\dk^{\leq 1}\psi +O(r^{-3})\dk^{\leq 2}\psi\\
&+ r^{-2} \dk^{\leq 2}\Ga_b + r^{-1}\dk^{\leq 3}(\Ga_b \c\Rc_b)-r^2\dk^{\leq 2}(\Ab\c \ov{A}).
\end{split}
\eeaa
Applying the $r^p$ weighted estimates of Proposition \ref{Proposition:Step3-Chap10}, and noticing that the first term on the RHS of the above wave equation for $\square_\g\big(re_4\psi\big)$ has the right sign in the estimate, 
we infer 
\beaa
\BEF_{\de; r\geq r_0}^{J-2}[re_4\psi] &\les& \ep_0^2+r_0B_{\de; r_0/2\leq r\leq r_0}^{J-2}[re_4\psi] \\
&&+ \int_{\MM(r\geq r_0/2)}r^\de\big(|re_4\dk^{\leq J-2}re_4\psi|+|\dk^{\leq J-2}re_4\psi|\big)|\dk^{\leq J-2}N_{re_4}|.
\eeaa
In view of the form of $N_{re_4}$, we infer
\beaa
\BEF_{\de; r\geq r_0}^{J-2}[re_4\psi] &\les& \ep_0^2+\ep_J^2+r_0B_{\de; r_0/2\leq r\leq r_0}^{J-2}[re_4\psi] +(\ep_0+\ep_J)\sqrt{\BEF_{\de; r\geq r_0}^{J-2}[re_4\psi]}\\
&&+\sqrt{\BEF_{\de; r\geq r_0}^{J-3}[(r\nab)^2\psi]}\sqrt{\BEF_{\de; r\geq r_0}^{J-2}[re_4\psi]}+r_0^{-1}\BEF_{\de; r\geq r_0}^{J-1}[\psi].
\eeaa
For $r_0$ large enough, we may absorb the part of the last term on the RHS that has $e_4$ derivatives and obtain 
\beaa
\BEF_{\de; r\geq r_0}^{J-2}[re_4\psi] &\les& \ep_0^2+\ep_J^2+r_0B_{\de; r_0/2\leq r\leq r_0}^{J-2}[re_4\psi]\\
&&+\sqrt{\BEF_{\de; r\geq r_0}^{J-3}[(r\nab)^2\psi]}\sqrt{\BEF_{\de; r\geq r_0}^{J-2}[re_4\psi]}+\BEF_\de^{J-5}[\big( \nab_\That, r\nab\big)^{\leq 4}\psi].
\eeaa
Also, integrating by parts, one easily obtains 
\beaa
\BEF_{\de; r\geq r_0}^{J-3}[(r\nab)^2\psi] &\les& \BEF_\de^{J-5}[\big( \nab_\That, r\nab\big)^{\leq 4}\psi] +\ep_J\sqrt{\EF_{\de; r\geq r_0}^{J-2}[re_4\psi]}\\
&&+\sqrt{\BEF_\de^{J-5}[\big( \nab_\That, r\nab\big)^{\leq 4}\psi]}\sqrt{\BEF_{\de; r\geq r_0}^{J-2}[re_4\psi]}
\eeaa
and hence
\beaa
\BEF_{\de; r\geq r_0}^{J-2}[re_4\psi] &\les& \ep_0^2+\ep_J^2+r_0B_{\de; r_0/2\leq r\leq r_0}^{J-2}[re_4\psi] +\BEF_\de^{J-5}[\big( \nab_\That, r\nab\big)^{\leq 4}\psi].
\eeaa
Also, since $e_4$ is spanned by $\That$ and $\Rhat$, 
\beaa
\BEF_{\de; r\geq r_0}^{J-2}[re_4\psi] &\les& \ep_0^2+\ep_J^2+r_0B_{\de; r_0/2\leq r\leq r_0}^{J-2}[re_4\psi] +\BEF_\de^{J-5}[\big( \nab_\That, r\nab\big)^{\leq 4}\psi]\\
&\les&  \ep_0^2+\ep_J^2+r_0^5\BEF_\de^{J-5}[\big( \nab_\That, \nab_\Rhat, r\nab\big)^{\leq 4}\psi].
\eeaa
Since 
\beaa
\BEF_\de^{J-5}[\big( \nab_\That, r\nab_4, r\nab\big)^{\leq 4}\psi] &\les& r_0^5\BEF_\de^{J-5}[\big( \nab_\That, \nab_\Rhat, r\nab\big)^{\leq 4}\psi]+\BEF_{\de; r\geq r_0}^{J-2}[re_4\psi],
\eeaa
we infer
\beaa
\BEF_\de^{J-5}[\big( \nab_\That, r\nab_4, r\nab\big)^{\leq 4}\psi] &\les& \ep_0^2+\ep_J^2+r_0^5\BEF_\de^{J-5}[\big( \nab_\That, \nab_\Rhat, r\nab\big)^{\leq 4}\psi].
\eeaa
Together with \eqref{estimate:conditional-B^{J-3}[nab^2R-nab_T-nab^2-psi]}, we deduce, for any $0<\de_1\leq 1$,
\bea\lab{estimate:conditional-B^{J-3}[nab^2R-nab_T-nab^2-psi]:plusre4}
\bsplit
\BEF_\de^{J-5}[\big( \nab_\That, r\nab_4, r\nab\big)^{\leq 4}\psi] \les& r_0^5\Big( \de_1+ O(a,\ep) \de_1^{-1} \Big)   B_\de^{J-1}[\psi]   +r_0^5\de_1^{-1}  \EF_\de^{J-3}[r^2 \lap \psi] \\
& +r_0^5\de_1^{-1} \int_{\MM_{trap}}  \ |\dk^{J-3} \lap \psi|^2 +r_0^5\de_1^{-1}\Big(\ep_J^2 +\ep_0^2\Big)\\ 
& +r_0^5\Rk_J\Rk_{J+1}.
\end{split}
\eea

{\bf Step 9.} Next, we derive an estimate for $\BEF_{\de; r\leq r_+(1+\de_{red})}^{J-2}[e_3\psi]$. To this end, we commute the wave equation for $\psi$ in \eqref{eq:Schematicswave|q|^2TP-TZ} with $\nab_3$. In the region $r\leq 4m$, we obtain 
\beaa
\square_\g\big(e_3\psi\big) &=& e_3\big(\square_\g\psi\big) - [e_3, \square_\g]\psi\\
&=& O(1)\dk^{\leq 1}\psi+ r^{-2} \dk^{\leq 2}\Ga_b + r^{-1}\dk^{\leq 3}(\Ga_b \c\Rc_b)-r^2\dk^{\leq 2}(\Ab\c \ov{A}) - [e_3, \square_\g]\psi.
\eeaa
Together with the commutator formula of Lemma \ref{lemma:commutationwithe3forredshift} which states that 
we have,  for $r\leq 4m$, 
\beaa
[\nab_3, \square_\g] &=& -\pr_r\left(\frac{\De}{|q|^2}\right)\nab_3^2\psi+O(1)\nab\nab_3\psi +O(1)\nab_4 \nab_3\psi\\
&& +O(1)\square_\g \psi +O(1)\dk^{\leq 1}\psi+\Ga_b\dk^{\leq 2}\psi,
\eeaa
we infer
\beaa
\bsplit
\square_\g\big(e_3\psi\big) =& N_{e_3},\\
N_{e_3} :=& \pr_r\left(\frac{\De}{|q|^2}\right)\nab_3^2\psi+O(1)\nab\nab_3\psi +O(1)\nab_4 \nab_3\psi  +O(1)\dk^{\leq 1}\psi\\
& + r^{-2} \dk^{\leq 2}\Ga_b + r^{-1}\dk^{\leq 3}(\Ga_b \c\Rc_b)-r^2\dk^{\leq 2}(\Ab\c \ov{A})+\Ga_b\dk^{\leq 2}\psi.
\end{split}
\eeaa
Making use of the favorable sign of $\pr_r\left(\frac{|\De|}{|q|^2}\right) $  in the red shift  region $r\leq r_+(1+\de_{red})$, and proceeding as in the redshift estimates of section \ref{sec:proofofThm:HigherDerivs-Morawetz-chp3:chap9}, we easily infer
\beaa
\BEF_{\de; r\leq r_+(1+\de_{red})}^{J-2}[\nab_3\psi] &\les& \ep_0^2+\ep_J^2+\de_{red}^{-1}\BEF_{\de; r_+(1+\de_{red})\leq r\leq r_+(1+2\de_{red})}^{J-2}[\nab_3\psi].
\eeaa
This yields
\beaa
\BEF_{\de; r\leq r_+(1+\de_{red})}^{J-2}[\nab_3\psi] &\les& \ep_0^2+\ep_J^2+\de_{red}^{-5}\BEF_\de^{J-5}[\big( \nab_\That, r\nab_4, r\nab\big)^{\leq 4}\psi].
\eeaa
Since
\beaa
\BEF_\de^{J-1}[\psi] &\les& \de_{red}^{-5}\BEF_\de^{J-5}[\big( \nab_\That, r\nab_4, r\nab\big)^{\leq 4}\psi]+\BEF_{\de; r\leq r_+(1+\de_{red})}^{J-2}[\nab_3\psi]
\eeaa
we infer, fixing the value of $\de_{red}>0$ small enough for the redshift estimate used above to hold, 
\beaa
\BEF_\de^{J-1}[\psi] &\les& \ep_0^2+\ep_J^2+\BEF_\de^{J-5}[\big( \nab_\That, r\nab_4, r\nab\big)^{\leq 4}\psi].
\eeaa
Together with \eqref{estimate:conditional-B^{J-3}[nab^2R-nab_T-nab^2-psi]:plusre4}, this yields, for any $0<\de_1\leq 1$, 
\bea\lab{estimate:conditional-B^{J-3}[nab^2R-nab_T-nab^2-psi]:plusre4:pluse3}
\bsplit
\BEF_\de^{J-1}[\psi] \les& r_0^5\Big( \de_1+ O(a,\ep) \de_1^{-1} \Big)   B_\de^{J-1}[\psi]   +r_0^5\de_1^{-1}  \EF_\de^{J-3}[r^2 \lap \psi] \\
& +r_0^5\de_1^{-1} \int_{\MM_{trap}}  \ |\dk^{J-3} \lap \psi|^2 +r_0^5\de_1^{-1}\Big(\ep_J^2 +\ep_0^2\Big)  +r_0^5\Rk_J\Rk_{J+1}.
\end{split}
\eea

{\bf Step 10.}  In view of \eqref{estimate:conditional-B^{J-3}[nab^2R-nab_T-nab^2-psi]:plusre4:pluse3}, we need to  estimate $ \EF_\de^{J-3}[ r^2\Delta\psi]  +\int_{\MM_{trap}}  \ |\dk^{J-3} \lap \psi|^2 $. To  this end, we rely 
on the following identity, see  Lemma \ref{Lemma:Formula-qf-DDhotDDPc}, 
\beaa
 \qfb &=& \frac{1}{2}\ov{q}q^3 \DDc\hot \DDc\Pc +  \dk^{\leq 1}\Ga_b'  +O(a)\dk^{\leq 1}\Rc_b + O(ar) \dk^{\le 1}\Pc +r\dk^{\leq 1}(\Ga_b\c \Rc_b).
 \eeaa
Commuting with $(\Lieb_\T, \Lieb_\Z)$,  and using, as in the proof of Lemma \ref{Lemma:comparisonofnormsBEF-Rk-forpsi},
\beaa
\psi = q^2\big( \T P, \T Z) = q^2( \T \Pc, \Z\Pc)  + r^{-1}\Ga_b,
\eeaa
we easily  derive the following  identity for 
$\psi =q^2( \Lie_\T P , \Lie_\Z P)   $
\beaa
 \DDc\hot \DDc\psi &=& O(r^{-2} ) \Lieb_{\T, \Z}\qfb +  O(r^{-2} ) \dk^{\leq 2}\Ga_b' 
 +O(ar^{-2})\Lieb_{\T, \Z}\dk^{\leq 1}\Rc_b + O(ar^{-1} ) \dk^{\le 1} \Lieb_{\T, \Z}\Pc\\
  &&+ r^{-1} \dk^{\leq 2}(\Ga_b\c \Rc_b).
\eeaa
 Making use of our bootstrap  assumptions, we deduce
 \beaa
\EF_\de^{J-3}[   r^2  \DDc\hot \DDc\psi] &\les &   \EF_\de^{J-3}[ \Lieb_\T \qfb, \Lieb_\Z  \qfb]+
\EF^{J-3} _\de[\dk^{\leq 2}\Ga_b' ]+ \EF_\de^{J-3} [r\dk^{\leq2}(\Ga_b\c \Rc_b)]\\
&&+ O(a) \EF_\de^{J-3}[\Lieb_{\T, \Z}\dk^{\leq 1}\Rc_b]+ O(a) \EF_\de^{J-3}[ r \dk^{\le 1} \Lieb_{\T, \Z}\Pc]\\
&\les&\EF_\de^{J-2}[\qfb] +  \EF^{J-1} _\de[\Ga_b' ] +  \EF_\de^{J-1}[ \Rc_b] +  \EF_\de^{J-1}[ r\Pc]+\ep_0^2.
 \eeaa
 Together with Lemmas  \ref{Lemma:comparisonofnormsBEF-Rk} and  \ref{Lemma:comparisonofnormsBEF[Gac)-Sk}, and the control of $\qfb$ in Theorem \ref{Thm:Estimates-forqf}, we obtain 
 \beaa
 \EF_\de^{J-3}[   r^2  \DDc\hot \DDc\psi]&\les &\Sk_{J+1} \Sk_J +\Rk_{J+1}\Rk_J+\ep_J^2+\ep_0^2.
 \eeaa
 Using the  Hodge  type estimates  of Corollary \ref{Cor:HodgeThmM8}    over spheres  of fixed $r$ in  $\Si(\tau)\cup  \AA\cup \Si_*$, integrating  then over  these regions,  and making use of  \eqref{induction:hypoth-psi}, we deduce
 \beaa
  \EF_\de^{J-3}[ r^2\nab^2\psi]&\les&   \EF_\de^{J-3}[   r^2  \DDc\hot \DDc\psi] 
  +\EF^{J-2}_\de[\psi]+O(a+\ep) \EF^{J-1}_\de[\psi]\\
  &\les& \Sk_{J+1} \Sk_J +\Rk_{J+1}\Rk_J +\ep_J^2+\ep_0^2 +O(a+\ep) \EF^{J-1}_\de[\psi].
 \eeaa

 The term $  \int_{\MM_{trap}}   | \dk^{J-3} \lap \psi|^2$ is  estimated by  the same procedure, and is even easier to control.  We  obtain
 \beaa
 \int_{\MM_{trap}}   | \dk^{J-3}\Delta\psi|^2  &\les&  B^{J-3}_\de[\DDc\hot \DDc\psi] +O(a+\ep) B^{J-1}_\de[\psi]+B^{J-2}_\de[\psi]\\
&\les&  \int_{\MM_{trap}}|\nab^2\dk^{J-3}\psi|^2 +O(a+\ep) B^{J-1}_\de[\psi]+B^{J-2}_\de[\psi]\\
&\les& B_\de^{J-2}[\qfb] +  B^{J-1} _\de[\Ga_b' ] +  B_\de^{J-1}[ \Rc_b] +  B_\de^{J-1}[ r\Pc] +O(a+\ep) B^{J-1}_\de[\psi]\\
&& +\ep_J^2+\ep_0^2\\
 &\les&  \Sk_{J+1} \Sk_J +\ep_J^2+\ep_0^2 +O(a+\ep)B^{J-1}_\de[\psi].
 \eeaa 
 Therefore,
 \bea
 \lab{estimate:EF^{J-3}[r^2nab^2psi]}
 \bsplit
  \EF_\de^{J-3}[r^2\lap \psi] + \int_{\MM_{trap}}   | \dk^{J-3} \lap \psi|^2 \les& \Sk_{J+1} \Sk_J +\Rk_{J+1}\Rk_J +\ep_J^2+\ep_0^2\\
  &  +O(a+\ep) \BEF^{J-1}_\de[\psi].
 \end{split}
 \eea

 {\bf Step 11.}  Combining   \eqref{estimate:conditional-B^{J-3}[nab^2R-nab_T-nab^2-psi]:plusre4:pluse3} with 
 \eqref{estimate:EF^{J-3}[r^2nab^2psi]}, we infer, for any $0<\de_1\leq 1$, 
\beaa
\bsplit
\BEF_\de^{J-1}[\psi] \les& r_0^5\Big( \de_1+ O(a+\ep) \de_1^{-1} \Big)\BEF_\de^{J-1}[\psi]   +r_0^5\de_1^{-1}\Big( \Sk_{J+1} \Sk_J +\Rk_{J+1}\Rk_J +\ep_J^2+\ep_0^2\Big).
\end{split}
\eeaa
  We may now fix $\de_1$ such that $r_0^5\de_1$ is small, and then $a$ and $\ep$, such that both $r_0^5\de_1$ and $r_0^5\de_1^{-1}(a+\ep)$ are small enough to absorb the first term on the RHS from the LHS. We infer 
\beaa
\bsplit
\BEF_\de^{J-1}[\psi] \les& r_0^{10}\Big( \Sk_{J+1} \Sk_J +\Rk_{J+1}\Rk_J +\ep_J^2+\ep_0^2\Big)
\end{split}
\eeaa
  which ends the proof of Proposition  \ref{proposition:Morawetz-Energy-psi}.

 
\subsection{Proof of Theorem \ref{theorem:Morawetz-EnergyPc}}


The proof of Theorem \ref{theorem:Morawetz-EnergyPc} proceeds along the following steps. 

{\bf Step 1.} Recalling the  definition of $\psi$,  we write  as in the proof of Lemma \ref{Lemma:comparisonofnormsBEF-Rk-forpsi},
\beaa
\psi=q^2\big( \T P, \T Z) = q^2( \T \Pc, \Z\Pc)  + r^{-1}\Ga_b.
\eeaa
 According to   Proposition  \ref{proposition:Morawetz-Energy-psi} and Lemma \ref{Lemma:comparisonofnormsBEF[Gac)-Sk},  we deduce
\beaa
 \BEF^{J-1}_\de[ r^2 ( \T \Pc, \Z\Pc)]   &\les & r_0^{10}\Big( \Sk_{J+1} \Sk_J +\Rk_{J+1}\Rk_J +\ep_J^2+\ep_0^2\Big).
\eeaa
Recalling that $\That=\T+\frac{a}{r^2+a^2}\Z$, we infer
\bea\lab{eq:slidfhzdkvjzbdkjzbjhbsldiofh:ThatPc}
 \BEF^{J-1}_\de[ r^2  \That (\Pc) ] &\les &  r_0^{10}\Big( \Sk_{J+1} \Sk_J +\Rk_{J+1}\Rk_J +\ep_J^2+\ep_0^2\Big).
 \eea
 
{\bf Step 2.}  Next, renormalizing the wave equation \eqref{eq:P-WaveEq-linearized} for $\Pc$ using Lemma \ref{Le:squareq^2Psi}, we  obtain  
\beaa
\square_\g(q^2\Pc)&=& Wq^2\Pc+  r^{-2} \dk^{\leq 1}\Ga_b+ \Ga_b \c  \dk^{\leq 1}\Rc_b
\eeaa
which we rewrite as
\beaa
\square_\g(\Psi)&=& W\Psi+  r^{-2} \dk^{\leq 1}\Ga_b+ \Ga_b \c  \dk^{\leq 1}\Rc_b, \qquad \Psi:= q^2\Pc.
\eeaa
Therefore, using  the formula of Lemma  \ref{lemma:squaredkpsi-perturbations}, we infer
\beaa
\nab_\Rhat^2\Psi &=&\nab^2_\That \Psi+ O(1) \square_\g\Psi+O(1)\lap \Psi   +O(r^{-1} )\nab_\Rhat \Psi    
   +O(a r^{-2} )\nab\Psi + \Ga_g \c \dk \Psi\\
   &=& \nab^2_\That \Psi+ O(1)\lap \Psi    +O(r^{-1} )\nab_\Rhat \Psi     +O(a r^{-2} )\nab\Psi  +O(r^{-2} ) \Psi \\
   && +r^{-2} \dk^{\leq 1}\Ga_b + \Ga_b \c  \dk^{\leq 1}\Rc_b.
   \eeaa
   Arguing as in \eqref{induction:hypoth-psi}, we have
   \beaa
  \BEF^{J-1} _\de[\Psi] &\les  \ep_J^2+\Rk_J\Rk_{J+1}. 
 \eeaa
 Also,  \eqref{eq:slidfhzdkvjzbdkjzbjhbsldiofh:ThatPc} yields for $\Psi= q^2\Pc$
 \beaa
 \BEF^{J-1}_\de[ \That(\Psi) ] &\les &  r_0^{10}\Big( \Sk_{J+1} \Sk_J +\Rk_{J+1}\Rk_J +\ep_J^2+\ep_0^2\Big).
 \eeaa
 Consequently
\beaa
\BEF^{J-2}_\de[\nab^2_\Rhat  \Psi ] &\les& \BEF^{J-2}_\de[\nab^2_\That  \Psi ]  + \BEF^{J-2}_\de[\lap  \Psi ] + \BEF^{J-1} _\de[\Psi] \\
&&+ \BEF^{J-2}_\de[r^{-2} \dk^{\leq 1}\Ga_b]+\BEF^{J-2}_\de[\Ga_b \c  \dk^{\leq 1}\Rc_b ]\\
&\les&  \BEF^{J-2}_\de[\lap  \Psi  ] +r_0^{10}\Big( \Sk_{J+1} \Sk_J +\Rk_{J+1}\Rk_J +\ep_J^2+\ep_0^2\Big).
\eeaa
We thus  deduce
\beaa
\bsplit
\BEF^{J-2}_\de[\nab^2_\That \Psi,  \nab^2_\Rhat\Psi ]  \les   \BEF^{J-2}_\de[ r^2\lap  \Psi  ]   +r_0^{10}\Big( \Sk_{J+1} \Sk_J +\Rk_{J+1}\Rk_J +\ep_J^2+\ep_0^2\Big).
\end{split}
\eeaa

{\bf Step 3.}  We can then proceed exactly as in Step 10 of section \ref{sec:proofofproposition:Morawetz-Energy-psi}, with the help of  the identity \eqref{eq-Formula-qf-DDhotDDPc:P3} and  the Hodge estimate of  Corollary \ref{Cor:HodgeThmM8}    over spheres  of fixed $r$ either in $\MM$ or $\Si(\tau)\cup  \AA\cup \Si_*$, 
  to derive 
 \beaa
\BEF^{J-2}_\de[ r^2 \nab^2  \Psi  ] &\les&    \Sk_{J+1} \Sk_J +\Rk_{J+1}\Rk_J +\ep_J^2+\ep_0^2  +O(a+\ep) \BEF^{J}_\de[\Psi].
\eeaa
Together with Step 2, we infer 
\beaa
\bsplit
\BEF^{J-2}_\de[\nab^2_\That \Psi,  \nab^2_\Rhat\Psi,   r^2\nab^2   \Psi  ]  \les&   r_0^{10}\Big( \Sk_{J+1} \Sk_J +\Rk_{J+1}\Rk_J +\ep_J^2+\ep_0^2\Big)  +O(a+\ep) \BEF^{J}_\de[\Psi] .
\end{split}
\eeaa
Since $\Psi= q^2\Pc$, we have, using Lemma \ref{Lemma:comparisonofnormsBEF-Rk},
\beaa
\BEF^{J-2}_\de[r^2(\nab^2_\That,  \nab^2_\Rhat,   r^2\nab^2)\Pc] &\les& \BEF^{J-2}_\de[\nab^2_\That \Psi,  \nab^2_\Rhat\Psi,   r^2\nab^2   \Psi  ]+\BEF^{J-1}[r^2\Pc]\\
&\les&  \BEF^{J-2}_\de[\nab^2_\That \Psi,  \nab^2_\Rhat\Psi,   r^2\nab^2   \Psi  ] +\Rk_J\Rk_{J+1},
\eeaa
and hence
\beaa
\bsplit
\BEF^{J-2}_\de[r^2(\nab^2_\That,  \nab^2_\Rhat,   r^2\nab^2)\Pc]  \les&   r_0^{10}\Big( \Sk_{J+1} \Sk_J +\Rk_{J+1}\Rk_J +\ep_J^2+\ep_0^2\Big)  +O(a+\ep) \BEF^{J}_\de[r^2\Pc] .
\end{split}
\eeaa
Moreover, since $J\geq\frac{k_L}{2}$ and $\kl$ is large, we may assume that $J\geq 4$. In particular, we have
\beaa
\BEF_\de^{J-4}[r^2\big( \nab_\That, \nab_\Rhat, r\nab\big)^{\leq 4}\Pc] &\les& \BEF_\de^{J-2}[r^2\big( \nab^2_\That, \nab^2_\Rhat, r^2 \nab^2\big)\Pc]+\BEF_\de^{J-1}[\Pc]
\eeaa 
which yields
\bea\lab{eq:interemdiaryodjaedigha}
\bsplit
\BEF_\de^{J-4}[r^2\big( \nab_\That, \nab_\Rhat, r\nab\big)^{\leq 4}\Pc]   \les&   r_0^{10}\Big( \Sk_{J+1} \Sk_J +\Rk_{J+1}\Rk_J +\ep_J^2+\ep_0^2\Big)\\
&  +O(a+\ep) \BEF^{J}_\de[r^2\Pc] .
\end{split}
\eea

{\bf Step 4.} \eqref{eq:interemdiaryodjaedigha} degenerates in the redshift region. To improve on it, we next derive estimates  for $\nab_3\Pc$ in that region.   These 
can be obtained, as in  Step 9 in the proof of Proposition  \ref{proposition:Morawetz-Energy-psi},
according to the following steps\footnote{Note that we  avoid using  the direct   linearization \eqref{eq:P-WaveEq-linearized} of the wave equation for $P$ for the same reason as before, i.e. it leads to   $\dk^{J+1} \Ga_b$ on the right hand side.  We thus  first commute the  wave equation for $P$ with $\nab_3$ and then linearize.}.
\begin{enumerate}
\item  First, commute the wave equation  \eqref{eq:P-WaveEq-M8} for $P$ with $  \nab_3$ using 
Lemma \ref{lemma:commutationwithe3forredshift} and linearize it to derive a wave equation for  the linearized quantity  $\widecheck{e_3(P)}=e_3(P)+\frac{6m}{q^4}$. This yields,  for $r\leq 4m$, 
\beaa
\square_\g(\widecheck{e_3(P)}) &=& 
  \left(\pr_r\left(\frac{|\De|}{|q|^2}\right)+O\left(\frac{\De}{r^2}\right)\right)\nab_3(\widecheck{e_3(P)})
+O(1)\nab(\widecheck{e_3(P)})\\
&&+O(1)e_4(\widecheck{e_3(P)})   +O(1)\dk^{\leq 1}\Pc+ \dk^{\leq 1}\Ga_b+O(\ep)\dk^{\leq 2}\Pc.
\eeaa

\item Then, we commute further with $J-1$ non degenerate derivatives and apply the redshift estimate to the resulting commuted equation. Making use of the favorable sign of $\pr_r\left(\frac{|\De|}{|q|^2}\right) $  in the red shift  region $r\leq r_+(1+\de_{red})$, and proceeding as in the redshift estimates of section \ref{sec:proofofThm:HigherDerivs-Morawetz-chp3:chap9}, we easily infer
\beaa
\BEF_{\de; r\leq r_+(1+\de_{red})}^{J-1}[\widecheck{e_3(P)}] &\les& \ep_0^2+\ep_J^2+\de_{red}^{-1}\BEF_{\de; r_+(1+\de_{red})\leq r\leq r_+(1+2\de_{red})}^{J-1}[\widecheck{e_3(P)}].
\eeaa
This yields
\beaa
\BEF_{\de; r\leq r_+(1+\de_{red})}^{J-1}[\widecheck{e_3(P)}] &\les& \ep_0^2+\ep_J^2+\de_{red}^{-5}\BEF_\de^{J-4}[\big( \nab_\That, r\nab_4, r\nab\big)^{\leq 4}\Pc].
\eeaa

\item Using the fact that 
\beaa
\widecheck{e_3(P)} &=& e_3(P)+\frac{6m}{q^4}=e_3\left(\Pc -\frac{2m}{q^3}\right)+\frac{6m}{q^4}\\
&=& e_3(\Pc)+r^{-3}\Ga_b,
\eeaa
we infer
\beaa
\BEF_{\de; r\leq r_+(1+\de_{red})}^{J-1}[e_3(\Pc)] &\les& \ep_0^2+\ep_J^2+\de_{red}^{-5}\BEF_\de^{J-4}[\big( \nab_\That, r\nab_4, r\nab\big)^{\leq 4}\Pc].
\eeaa
\end{enumerate}
Together with \eqref{eq:interemdiaryodjaedigha}, we obtain, fixing the value of $\de_{red}>0$ small enough for the redshift estimate used above to hold, 
\bea\lab{eq:eq:intermediaryredshiftaodifa}
\bsplit
\BEF_\de^{J-4}[r^2\big( \nab_3, \nab_4, r\nab\big)^{\leq 4}\Pc]   \les&  r_0^{10}\Big( \Sk_{J+1} \Sk_J +\Rk_{J+1}\Rk_J +\ep_J^2+\ep_0^2\Big)\\
&  +O(a+\ep) \BEF^{J}_\de[r^2\Pc] .
\end{split}
\eea

 {\bf Step 5.} In view of \eqref{eq:eq:intermediaryredshiftaodifa}, it remains to recover $r\nab_4$ derivatives of $\Pc$ in $r\geq r_0$. We only sketch this step:
 \begin{enumerate}
 \item As in Step 4, we do not use the direct   linearization \eqref{eq:P-WaveEq-linearized} of the wave equation for $P$. Instead, we commute the wave equation for $P$ with $re_4$, and linearize it using $\widecheck{re_4(P)}=re_4(P)-\frac{6m}{q^4}\frac{r\De}{|q|^2}$. Finally, we  then renormalize the wave equation for $\widecheck{re_4(P)}$ using Lemma \ref{Le:squareq^2Psi}, thus yielding a wave equation for $q^2\widecheck{re_4(P)}$. 
 
\item  We then use an $r^p$ weighted estimate in the region $r\geq r_0$ with $p=\de$ for the wave equation satisfied by $q^2\widecheck{re_4(P)}$ as in Step 8 of section \ref{sec:proofofproposition:Morawetz-Energy-psi}.

\item As in Step 4, we relate $\widecheck{re_4(P)}$ and $re_4(\Pc)$ and deduce an $r^p$ weighted estimate in the region $r\geq r_0$ with $p=\de$ for  $q^2re_4(\Pc)$.
\end{enumerate}
This $r^p$ weighted procedure yields
\beaa
\BEF_{\de; r\geq r_0}^{J-1}[r^3e_4(\Pc)] &\les&  \ep_0^2+\ep_J^2+r_0^5\BEF_\de^{J-4}[r^2\big( \nab_3, \nab_4, r\nab\big)^{\leq 4}re_4(\Pc)]
\eeaa
which together with \eqref{eq:eq:intermediaryredshiftaodifa} yields
\beaa
\bsplit
\BEF_\de^J[r^2\Pc]   \les&   r_0^{15}\Big( \Sk_{J+1} \Sk_J +\Rk_{J+1}\Rk_J +\ep_J^2+\ep_0^2\Big)  +r_0^5O(a+\ep) \BEF^{J}_\de[r^2\Pc].
\end{split}
\eeaa
For $a$ and $\ep$ small enough compared to $r_0^{-5}$, we infer
\beaa
\bsplit
\BEF_\de^J[r^2\Pc]   \les&  r_0^{15}\Big( \Sk_{J+1} \Sk_J +\Rk_{J+1}\Rk_J +\ep_J^2+\ep_0^2\Big)  
\end{split}
\eeaa
as stated. This ends the proof of Theorem \ref{theorem:Morawetz-EnergyPc}


 \section{Proof of Lemmas \ref{LEMMA:WAVEEQUATIONFORMODQ2TPTZ}, \ref{LEMMA:SQUARENAB^2_RHAT} and 
 \ref{LEMMA:ESTIMATE:BEF[NAB^2_RHAT].}}


In this section, we provide the proof of Lemmas \ref{LEMMA:WAVEEQUATIONFORMODQ2TPTZ}, \ref{LEMMA:SQUARENAB^2_RHAT} and \ref{LEMMA:ESTIMATE:BEF[NAB^2_RHAT].} which are used in the proof of Proposition  \ref{proposition:Morawetz-Energy-psi} in section \ref{sec:proofofproposition:Morawetz-Energy-psi}.


\subsection{Proof of Lemma \ref{LEMMA:WAVEEQUATIONFORMODQ2TPTZ}}
\lab{section:DetailesStep1}


The proof of Lemma \ref{LEMMA:WAVEEQUATIONFORMODQ2TPTZ} follows from checking  that 
the linearized  quantities $q^2 \T P,\, q^2  \Z P$ satisfy  the following wave equations
\bea
\lab{waveeq:TZ(P)}
\bsplit
\square_\g(q^2\T P) &=  W q^2\T P+ r^{-2} \dk^{\leq 1}\Ga_b + r^{-1}\dk^{\leq 2}(\Ga_b \c\Rc_b)-r^2\dk^{\leq 1}(\Ab\c \ov{A}),\\
\square_\g(q^2\Z P) &=  W q^2 \Z P+ r^{-2} \dk^{\leq 1}\Ga_b  + r^{-1}\dk^{\leq 2}(\Ga_b \c\Rc_b)-r^2\dk^{\leq 1}(\Ab\c \ov{A}),
\end{split}
\eea
where the  complex potential  $W$  verifies 
\beaa
\Re( W )=O(mr^{-3}), \qquad  \Im(W)=O(am r^{-4}).
\eeaa

Recall the wave equation for $P$,  see  \eqref{eq:P-WaveEq-M8},
\beaa
\nn\square_\g P &=&  \tr X\nab_3P +\ov{\tr\Xb}\nab_4P  -\ov{H}\c\DD P  - \Hb\c\ov{\DD}P + V P   + r^{-3}\dk^{\leq 1}(\Ga_b \c\Rc_b)-\Ab\c \ov{A}.
\eeaa
 Using    the commutator   formula  \eqref{eq:commTZsquare} for $[ \T, \square_\g]$, i.e. 
\beaa
 \,[ \T, \square_\g]P =\dk \big(\Ga_g \c \dk P\big)+\Ga_b \c \square_\g P,
 \eeaa
 and, $[\T, e_3]=[\T, e_4]=[\T, e_a]=\Ga_b \c \dk$, we deduce
\beaa
 \square_\g(\T P)&=&  \T \square_\g P+\dk \big(\Ga_g \c \dk P\big)+\Ga_b \c \square_\g P\\
 &=&  \T\Big(\tr X\nab_3P +\ov{\tr\Xb}\nab_4P  -\ov{H}\c\DD P  - \Hb\c\ov{\DD}P + V P  + r^{-3}\dk^{\leq 1}(\Ga_b \c\Rc_b)-\Ab\c \ov{A}\Big) \\
 &&+\dk \big(\Ga_g \c \dk P\big)+\Ga_b \c \Big(r^{-1}\dk^{\leq 1}P + r^{-3}\dk^{\leq 1}(\Ga_b \c\Rc_b)-\Ab\c \ov{A}\Big)\\
 &=& \tr X\nab_3(\T P) +\ov{\tr\Xb}\nab_4(\T P)  -\ov{H}\c\DD( \T P)  - \Hb\c\ov{\DD}(\T P) +V \T P  \\
 &&      + r^{-3} \dk^{\le 1}   \Ga_g + r^{-3}\dk^{\leq 2}(\Ga_b \c\Rc_b)-\dk^{\leq 1}(\Ab\c \ov{A})\\
&+&  \T (\tr X ) \nab_3P + \T( \ov{\tr\Xb}) \nab_4P  - \T (\ov{H})\c\DD P  -\T( \Hb) \c\ov{\DD}P +\T(V)   P  \\
&+& \tr X [\T, \nab_3] P +\ov{\tr\Xb} [\T, \nab_4] P -\ov{H}\c [\T,    \DD]P   - \Hb\c  [\T,      \ov{\DD}] P\\
&=& \tr X\nab_3(\T P) +\ov{\tr\Xb}\nab_4(\T P)  -\ov{H}\c\DD( \T P)  - \Hb\c\ov{\DD}(\T P) +V (\T P)\\
&&+  r^{-3} \dk^{\le 1}   \Ga_g + r^{-3}\dk^{\leq 2}(\Ga_b \c\Rc_b)-\dk^{\leq 1}(\Ab\c \ov{A}).
\eeaa
We are  now in a position to apply  Lemma  \ref{Le:squareq^2Psi}    and deduce
\beaa
\square_\g(q^2\T P) &=&  \Big[ V+q^{-2}\square_\g(q^2)\Big]q^2\T P +r\Ga_b\c \dk(\T P)  +  r^{-1} \dk^{\le 1}   \Ga_g\\
&& + r^{-1}\dk^{\leq 2}(\Ga_b \c\Rc_b)-r^2\dk^{\leq 1}(\Ab\c \ov{A}).
\eeaa
Thus, since $\T P= \T(\Pc-\frac{2m}{q^3} )     =   \dk\Pc +   r^{-1}\dk \Ga_g $, and $\Ga_g=r^{-1}\Ga_b$ in view of Remark \ref{Remark:Ga_g=r^{-1}Ga_b}, 
\beaa
\square_\g(q^2\T P) &=W |q|^2  \T P  +r^{-2} \dk^{\le 1}   \Ga_b  + r^{-1}\dk^{\leq 2}(\Ga_b \c\Rc_b)-r^2\dk^{\leq 1}(\Ab\c \ov{A})
\eeaa
which is  the first identity in \eqref{waveeq:TZ(P)}. Also, we have $ W=V+q^{-2}\square_\g(q^2)$ which satisfies indeed $\Re( W )=O(mr^{-3})$ and  $\Im(W)=O(am r^{-4})$. Finally, the second identity in \eqref{waveeq:TZ(P)} can be derived in the same manner.  This concludes the proof of \eqref{waveeq:TZ(P)}, and hence of Lemma \ref{LEMMA:WAVEEQUATIONFORMODQ2TPTZ}.


\subsection{Proof of Lemma \ref{LEMMA:SQUARENAB^2_RHAT}}
\lab{section:DetailesStep3}


Starting with  the  equation $\square_\g \psi= N$, and using the following commutation formula, see Lemma \ref{LEMMA:COMMUTATOR-NAB-RHAT-SQUARE}, 
  \beaa
   \, [\nab_\Rhat, |q|^2\square_\g]\psi &=&O(r)\square_\g \psi+O(r)\lap\psi +O(ar^{-1})\dk^{\le 2}  \psi+ O(1) \nab_\Rhat \psi +O(ar^{-1}) \nab \psi\\
 &&+ r \dk \big( \Ga_b \c \dk \psi\big),
 \eeaa
we obtain 
\beaa
\square_\g\psi(\nab_\Rhat\psi) &=& \frac{1}{|q|^2}\nab_\Rhat(|q|^2N)+\frac{1}{|q|^2}[\nab_\Rhat, |q|^2\square_\g]\psi
\eeaa
and hence
\beaa
\square_\g \nab_\Rhat \psi= [N],
\eeaa
where
\beaa
 [N]&:=& \frac{1}{|q|^2}\nab_\Rhat(|q|^2N)+\frac{1}{|q|^2}[\nab_\Rhat, |q|^2\square_\g]\psi\\
&=&  O(r^{-1})\square_\g \psi+O(r^{-1} )\lap\psi +O(ar^{-3})\dk^{\leq 2}\psi + O(r^{-2}) \nab_\Rhat \psi+O(ar^{-3}) \nab \psi\\
&&+ r^{-1}\dk \big( \Ga_b \c \dk \psi\big) + \nab_\Rhat N+r^{-1} N\\
&=& O(r^{-1} )\lap\psi +O(ar^{-3})\dk^{\leq 2}\psi+ O(r^{-2}) \nab_\Rhat \psi+O(ar^{-4}) \nab \psi\\
&& +r^{-1}\dk \big( \Ga_b \c \dk \psi\big) + \nab_\Rhat N+r^{-1} N.
\eeaa
Repeating the process, we find 
\beaa
\square_\g \nab^2_\Rhat \psi= [[N]]
\eeaa 
where
\beaa
[[N]]&:=& \frac{1}{|q|^2}\nab_\Rhat(|q|^2 [N])+\frac{1}{|q|^2}[\nab_\Rhat, |q|^2\square_\g]\nab_\Rhat\psi\\
&=&  O(r^{-1} )\lap \Rhat \psi +O(ar^{-3})\dk^{\leq 2}\nab_\Rhat\psi +O(r^{-2}) \nab_\Rhat^2 \psi+O(ar^{-4}) \nab\nab_\Rhat  \psi\\
&& + r^{-1}\dk^{\leq 3}\big( \Ga_b \c \psi\big) + \nab_\Rhat [N]+r^{-1} [N]\\
&=& O(r^{-1} ) \nab_\Rhat  \lap \psi  +O(r^{-2})  \lap \psi          +O(ar^{-3})   \nab_\Rhat\dk^{\leq 2}\psi  +O(r^{-2}) \nab_\Rhat \psi+ O(ar^{-3} )\dk^{\le 2} \psi\\
&&+  r^{-1}\dk^{\leq 3}\big( \Ga_b \c \psi\big)   + \nab_\Rhat [N]+r^{-1} [N]
\eeaa
where we used in particular the commutation formula for $[\nab_\Rhat, \Delta]$ in Lemma \ref{LEMMA:COMMUTATOR-NAB-RHAT-SQUARE}, i.e.
\beaa
[\nab_\Rhat, \lap]\psi&=&O(r^{-1}) \lap \psi +O(ar^{-5}) \dkb \psi+ r^{-1} \dk \big( \Ga_g \c \dk \psi \big).
\eeaa
Now
\beaa
 \nab_\Rhat [N]&=&\nab_\Rhat\Big(O(r^{-1} )\lap\psi +O(ar^{-3})\dk^{\leq 2}\psi+ O(r^{-2}) \nab_\Rhat \psi+O(ar^{-4}) \nab \psi\Big)\\
 && +r^{-1}\dk^{\leq 3}\big( \Ga_b \c \psi\big) + \nab^2 _\Rhat N+r^{-1}\nab_\Rhat  N + r^{-2}  N\\
 &=& O(r^{-1} )\nab_\Rhat \lap\psi +O(ar^{-3})\nab_\Rhat \dk^{\leq 2}\psi+O(ar^{-4}) \nab_\Rhat \nab \psi\\
 &&+O(r^{-2}) \lap \psi +O(r^{-2}) \nab_\Rhat \psi       + O(a r^{-4}) \dk^{\le 2} \psi +r^{-1}\dk^{\leq 3}\big( \Ga_b \c \psi\big)\\
 && + \nab^2 _\Rhat N+r^{-1}\nab_\Rhat  N + r^{-2}  N.
 \eeaa
Therefore,
\beaa
\square_\g \nab^2_\Rhat \psi&=& O(r^{-1} ) \nab_\Rhat  \lap \psi  +O(r^{-2})  \lap \psi          +O(ar^{-3})   \nab_\Rhat\dk^{\leq 2}\psi  +O(r^{-2}) \nab_\Rhat \psi+ O(ar^{-3} )\dk^{\le 2} \psi\\
&& +O(r^{-1} )\nab_\Rhat \lap\psi +O(ar^{-3})\nab_\Rhat\dk^{\leq 2}\psi+O(ar^{-4}) \nab_\Rhat \nab \psi\\
 &&+O(r^{-2}) \lap \psi +O(r^{-2}) \nab_\Rhat \psi       + O(a r^{-4}) \dk^{\le 2} \psi +r^{-1}\dk^{\leq 3}\big( \Ga_b \c \psi\big)\\
&&+ O(r^{-2} )\lap\psi +O(ar^{-4})\dk^{\leq 2}\psi+ O(r^{-3}) \nab_\Rhat \psi+O(ar^{-5}) \nab \psi\\
&& +  r^{-1}\dk^{\leq 3}\big( \Ga_b \c \psi\big)  + \nab^2 _\Rhat N+r^{-1}\nab_\Rhat  N + r^{-2}  N
\eeaa
and hence
\beaa
\square_\g \nab^2_\Rhat \psi&=& O(r^{-1} )\nab_\Rhat \lap\psi +O(ar^{-3})\nab_\Rhat\dk^{\leq 2} \psi +O(r^{-2}) \lap \psi \\
&&+O(r^{-2}) \nab_\Rhat \psi+ O(a r^{-3}) \dk^{\le 2} \psi  + r^{-1}\dk^{\leq 3}\big( \Ga_b \c \psi\big) + \nab^2 _\Rhat N+r^{-1}\nab_\Rhat  N + r^{-2}  N
\eeaa
as  desired. This concludes the proof of Lemma \ref{LEMMA:SQUARENAB^2_RHAT}.


\subsection{Proof of Lemma \ref{LEMMA:ESTIMATE:BEF[NAB^2_RHAT].}}
\lab{section:DetailesStep4}


  According to  Lemma \ref{Lemma:comparisonofnormsBEF-Rk-forpsi}
  we have
 \bea
 \label{eq:InductionHyp-forPsi}
 \BEF_\de^{J-2}[\psi] \les \ep_J^2 +\Rk_J\Rk_{J+1}.
 \eea
 We apply   the second estimate of    Proposition \ref{Prop:scalarwavePsi-M8} to equation \eqref{eq:wave-nab^2_Rhatpsi} with $V=-\Re(W)$, $k=J-3$, and an inhomogeneous  term $N$ given by
  \beaa
 N&=& i\Im(W)\nab_\Rhat^2 \psi      +O(r^{-2}) \lap \psi      +O(r^{-1})\nab_\Rhat\lap  \psi  +O(ar^{-3})\Big(\nab_\Rhat\dk^{\le  2} \psi    + \dk^{\le  2} \psi\Big)   \\ 
&&+O(r^{-2})  \dk^{\le 1} \psi    +N^{\le 3},
 \eeaa
 where $N^{\le 3}$ is defined in \eqref{eq:wave-nab^2_Rhatpsi}. We obtain
\beaa
 \bsplit
  \BEF^{J-3}_\de[\nab^2_\Rhat \psi] &\les\ep_0^2 +\BEF^{J-4}_{\de}[\nab^2_\Rhat\psi]  +\int_{\MM_{trap}} |\dk^{J-3}\nab^2_\Rhat\psi|^2  +\NN^{J-3}_\de[\nab^2_\Rhat \psi,  N]\\
  &\les\ep_0^2 +\BEF^{J-2}_{\de}[\psi] +\int_{\MM_{trap}}\Big(|\nab_\Rhat\dk^{\leq J-2}\psi|^2+|\dk^{\leq J-2}\psi|^2\Big) \\
  & +\NN^{J-3}_\de[\nab^2_\Rhat \psi,  N]\\
  &\les\ep_0^2 +\BEF^{J-2}_{\de}[\psi]   +\NN^{J-3}_\de[\nab^2_\Rhat \psi,  N]
  \end{split}
  \eeaa
  which together with \eqref{eq:InductionHyp-forPsi} implies
\bea
 \lab{eq:En-Mor-nab_Rhat^2psi1}
  \BEF^{J-3}_\de[\nab^2_\Rhat \psi]  &\les& \ep_0^2 +\ep_J^2 +\Rk_J\Rk_{J+1} +\NN^{J-3}_\de[\nab^2_\Rhat \psi,  N].
  \eea

It remains to estimate $\NN^{J-3}_\de[\nab^2_\Rhat \psi,  N]$. We have
\beaa
 && \NN^{J-3}_\de[\nab^2_\Rhat \psi,  N]\\
&=&  \NNmor^{J-3}[\nab^2_\Rhat\psi,  N] +\NNen^{J-3}[\nab^2_\Rhat\psi,  N]+ \NNred^{J-3}[\nab^2_\Rhat\psi,  N]\\
   &&+   \int_{\Mext} r^\de \Big(\big|\nab_4\dk^{\le J-3}\nab^2_\Rhat\psi \big|+r^{-1}\big|\dk^{\le J-3} \nab^2_\Rhat\psi \big|\Big) \big| \dk^{\le J-3} N \big|\\
       &\les& \Big(B ^{J-3}_{\de}[\nab^2_\Rhat\psi]\Big)^{\frac{1}{2}}\left(\int_{\MM}r^{\de+1}|\dk^{\le J-3}N|^2\right)^{\frac{1}{2}}+\NNen^{J-3}[\nab^2_\Rhat\psi,  N]\\
        &\les& \Big(B ^{J-1}_{\de}[\psi]\Big)^{\frac{1}{2}}\left(\int_{\MM}r^{\de+1}|\dk^{\le J-3}N|^2\right)^{\frac{1}{2}}+\int_{\Mtrap}\big|\nab_{\That_\de}\nab_{\Rhat}^2\dk^{\leq J-3}\psi|\big| \dk^{\le J-3} N \big|.
\eeaa
Since $\nab_{\That_\de}\nab_{\Rhat}^2\dk^{\leq J-3}\psi=\nab_{\Rhat}\dk^{\leq J-1}\psi+O(r^{-1})\dk^{\leq J-1}\psi$, we infer
\beaa
 \NN^{J-3}_\de[\nab^2_\Rhat \psi,  N]   &\les& \Big(B ^{J-1}_{\de}[\psi]\Big)^{\frac{1}{2}}\left(\int_{\MM}r^{\de+1}|\dk^{\le J-3}N|^2\right)^{\frac{1}{2}}.
\eeaa
Now, using the definition of $N$, the fact that $\Im(W)=O(ma r^{-4})$, and \eqref{eq:InductionHyp-forPsi}, we have
\beaa
 \int_{\MM}r^{\de+1}|\dk^{\le J-3}N|^2 &\les& |a|B ^{J-1}_{\de}[\psi]+B ^{J-3}_{\de}[r\Delta\psi]+B ^{J-2}_{\de}[\psi]+\int_{\MM}r^{\de+1}|\dk^{\le J-3}N^{\le 3}|^2\\
 &\les& |a|B ^{J-1}_{\de}[\psi]+B ^{J-3}_{\de}[r\Delta\psi]+\ep_J^2+\int_{\MM}r^{\de+1}|\dk^{\le J-3}N^{\le 3}|^2.
 \eeaa
Together with the definition of $N^{\le 3}$  in \eqref{eq:wave-nab^2_Rhatpsi}, i.e.
 \beaa
 N^{\le 3}&= r^{-2}   \dk^{\leq 3 }\Ga_b  + r^{-1}\dk^{\leq 4}(\Ga_b \c\Rc_b)-r^2\dk^{\leq 3}(\Ab\c \ov{A}),
\eeaa
the induction hypothesis for the linear term involving $\Ga_b$   and the  bootstrap  assumptions for the nonlinear one,   and Lemma \ref{lemma:auxilliarynormsforGa_b}, we infer
\beaa
 \int_{\MM}r^{\de+1}|\dk^{\le J-3}N|^2 &\les& |a|B ^{J-1}_{\de}[\psi]+B ^{J-3}_{\de}[r\Delta\psi]+\ep_J^2+\int_{\MM} r^{-3+\de} |\dk^{\le J} \Ga_b|^2\\
 && +\ep^2\int_{\MM} r^{-3+\de} |\dk^{\le J+1} \Rc_b|^2 +\ep^2\int_{\MM} r^{-2+\de} |\dk^{\le J+1} \Ga_b|^2\\
 &&+\ep^2\int_{\MM} r^{3+\de} |\dk^{\le J+1}A|^2\\
 &\les& |a|B ^{J-1}_{\de}[\psi]+B ^{J-3}_{\de}[r\Delta\psi]+\ep_J^2+\Sk_{J}^2+ \ep^2(\Sk_{J+1}+\Rk_{J+1})^2\\
 &\les& |a|B ^{J-1}_{\de}[\psi]+B ^{J-3}_{\de}[r\Delta\psi]+\ep_J^2+ \ep_0^2.
 \eeaa
 This yields
 \beaa
 \NN^{J-3}_\de[\nab^2_\Rhat \psi,  N]   &\les& \Big(B ^{J-1}_{\de}[\psi]\Big)^{\frac{1}{2}}\left(\int_{\MM}r^{\de+1}|\dk^{\le J-3}N|^2\right)^{\frac{1}{2}}\\
 &\les& \Big(B ^{J-1}_{\de}[\psi]\Big)^{\frac{1}{2}}\Big(|a|B ^{J-1}_{\de}[\psi]+B ^{J-3}_{\de}[r\Delta\psi]+\ep_J^2+ \ep_0^2\Big)^{\frac{1}{2}}.
\eeaa
Together with \eqref{eq:En-Mor-nab_Rhat^2psi1}, we deduce
\beaa
  \BEF^{J-3}_\de[\nab^2_\Rhat \psi]  &\les& \ep_0^2 +\ep_J^2 +\Rk_J\Rk_{J+1} +\NN^{J-3}_\de[\nab^2_\Rhat \psi,  N]\\
  &\les&  \ep_0^2 +\ep_J^2 +\Rk_J\Rk_{J+1} +\Big(B ^{J-1}_{\de}[\psi]\Big)^{\frac{1}{2}}\Big(|a|B ^{J-1}_{\de}[\psi]+B ^{J-3}_{\de}[r\Delta\psi]+\ep_J^2+ \ep_0^2\Big)^{\frac{1}{2}}.
  \eeaa
  Hence, we have for any $0<\de_1\leq 1$
   \beaa
   \begin{split}
    \BEF ^{J-3}_{\de}[\nab^2_\Rhat \psi]& \les ( \de_1+ |a| )  B_\de^{J-1}[\psi] +\de_1^{-1}  B_\de^{J-3}[r \lap \psi] +\de_1^{-1}\Big(\ep_J^2  +\ep_0^2\Big) +\Rk_J\Rk_{J+1}
 \end{split}
   \eeaa
   as stated in \eqref{estimate:BEF[nab^2_Rhat].}. This concludes the proof of Lemma \ref{LEMMA:ESTIMATE:BEF[NAB^2_RHAT].}.

  
\chapter{Energy-Morawetz for $\protect A, B, \Bb, \Ab$}
\lab{CHAPTER:ENERGYMOR-INTERIORM8}


The goal of this chapter  is to derive  estimates Energy-Morawetz for  the curvature components  $A, B, \Bb, \Ab$   using the Bianchi identities as well as the estimates for $\Pc$ derived in Chapter \ref{Chapter:EN-MorforPc}. This will complete the proof of the estimate \eqref{eq:InteriorcurvEstimatesThmM8} in Theorem \ref{prop:rpweightedestimatesiterationassupmtionThM8}, see section \ref{sec:proofofeq:InteriorcurvEstimatesThmM8}.

 
\section{Statement of the main results of Chapter \ref{CHAPTER:ENERGYMOR-INTERIORM8}}


In order to derive Energy-Morawetz for  the curvature components  $A, B, \Bb, \Ab$, we first control $(B, \Bb)$, and then $(A, \Ab)$. 

  
\subsection{Energy-Morawetz for $\protect B, \Bb$}


The following proposition provides energy-Morawetz estimates for $(B, \Bb)$. 

  \begin{proposition}
     \lab{proposition:EstimatesBBb-interior}
    The following estimates hold true in 
     $\MM=\MM(1, \tau_*) $
         \bea
 \lab{eq:mainestimateBBb-M8}
  \nn\BEF^J_\de[r^2B]+\BEF^J_\de[\Bb] &\les&  \de_{J+1}[\Pc]  +\ep_0^2+\ep_J^2 +|a|^2\Sk_{J+1}^2\\
  && +\ep_J\Rk_{J+1} +\left(\sqrt{ \de_{J+1}[\Pc]} +\ep_0+\ep_J\right)\Sk_{J+1},
 \eea
   where
   \beaa
   \de_{J+1}[\Pc] :=\BEF_\de^{J}[r^2 \Pc ].
   \eeaa
   \end{proposition}

Proposition \ref{proposition:EstimatesBBb-interior} will be proved in section \ref{section:InteriorEstimatesBBb}.

  
\subsection{Energy-Morawetz for $\protect A, \Ab$}


The following proposition provides energy-Morawetz estimates for $(A, \Ab)$. 
    \begin{proposition}
     \lab{proposition:EstimatesAAb-interior}
    The following estimates hold true in  $\MM=\MM(1, \tau_*)$
     \bea
 \lab{eq:proposition-EstimatesAAb-interior}
  \bsplit
   \BEF^J_\de[ r^2A]&\les  \de_{J+1}[B]  +\ep_0^2+\ep_J^2  +\left(\sqrt{ \de_{J+1}[B]} +\ep_0+\ep_J\right)\Sk_{J+1}+|a|^2\Sk_{J+1}^2,\\
    \BEF^J_\de[\Ab]&\les \de_{J+1}[\Bb]  +\ep_0^2+\ep_J^2  +\left(\sqrt{ \de_{J+1}[\Bb]} +\ep_0+\ep_J\right)\Sk_{J+1}+|a|^2\Sk_{J+1}^2,
   \end{split}
   \eea
   where 
   \beaa
\de_{J+1}[ B] :=\BEF_\de^J[ r^2 B],\qquad \de_{J+1}[ \Bb] := \BEF_\de^J[\Bb].
\eeaa
    \end{proposition}

Proposition \ref{eq:proposition-EstimatesAAb-interior} will be proved in section \ref{section:InteriorEstimatesAAb}.

  
\subsection{Proof of \eqref{eq:InteriorcurvEstimatesThmM8} in Theorem \ref{prop:rpweightedestimatesiterationassupmtionThM8}}
\lab{sec:proofofeq:InteriorcurvEstimatesThmM8}


We are now ready to prove the estimate \eqref{eq:InteriorcurvEstimatesThmM8} in Theorem \ref{prop:rpweightedestimatesiterationassupmtionThM8} on the control of the curvature norm $\Rkint$. 

\begin{proof}[Proof of \eqref{eq:InteriorcurvEstimatesThmM8} of Theorem \ref{prop:rpweightedestimatesiterationassupmtionThM8}]
First,  in view of Theorem \ref{theorem:Morawetz-EnergyPc} and the definition of $\de_{J+1}[\Pc]$, we have
\beaa
  \de_{J+1}[\Pc] =   \BEF^J_\de[r^2 \Pc] &\les&  r_0^{15}\Big( \Sk_{J+1} \Sk_J +\Rk_{J+1}\Rk_J +\ep_J^2+\ep_0^2\Big).
\eeaa
Also, in view of Proposition \ref{proposition:EstimatesBBb-interior}, we have
         \beaa
  \BEF^J_\de[r^2B]+\BEF^J_\de[\Bb] &\les&    \de_{J+1}[\Pc]  +\ep_0^2+\ep_J^2 +|a|^2\Sk_{J+1}^2\\
  && +\ep_J\Rk_{J+1} +\left(\sqrt{ \de_{J+1}[\Pc]} +\ep_0+\ep_J\right)\Sk_{J+1}.
 \eeaa
Together with the definition of $\de_{J+1}[ B]$ and $\de_{J+1}[ \Bb]$, and the above control of $\de_{J+1}[\Pc]$, we infer
\beaa
\de_{J+1}[\Pc]+\de_{J+1}[B]+\de_{J+1}[\Bb] &\les& r_0^{15}\Big( \ep_J(\Sk_{J+1}+\Rk_{J+1}) +\ep_J^2+\ep_0^2\Big)+|a|^2\Sk^2_{J+1}\\
&& +r_0^{\frac{15}{2}}\Sk_{J+1}\Big(\ep_0+\sqrt{\ep_J}\sqrt{\Sk_{J+1}+\Rk_{J+1}}\Big).
\eeaa

Next, in view of Proposition \ref{eq:proposition-EstimatesAAb-interior}, we have
    \beaa
  \bsplit
    \BEF^J_\de[ r^2A]&\les  \de_{J+1}[B]  +\ep_0^2+\ep_J^2  +\left(\sqrt{ \de_{J+1}[B]} +\ep_0+\ep_J\right)\Sk_{J+1}+|a|^2\Sk_{J+1}^2,\\
    \BEF^J_\de[\Ab]&\les \de_{J+1}[\Bb]  +\ep_0^2+\ep_J^2  +\left(\sqrt{ \de_{J+1}[\Bb]} +\ep_0+\ep_J\right)\Sk_{J+1}+|a|^2\Sk_{J+1}^2.
   \end{split}
   \eeaa
We infer
\beaa
&&\de_{J+1}[\Pc]+\de_{J+1}[B]+\de_{J+1}[\Bb] +\BEF^J_\de[ r^2A] + \BEF^J_\de[\Ab] \\
&\les& r_0^{15}\Big( \ep_J(\Sk_{J+1}+\Rk_{J+1}) +\ep_J^2+\ep_0^2\Big)+|a|\Sk^2_{J+1}\\
&&+r_0^{\frac{15}{4}}\Sk_{J+1}^{\frac{3}{2}}\Big(\ep_0+\sqrt{\ep_J}\sqrt{\Sk_{J+1}+\Rk_{J+1}}\Big)^{\frac{1}{2}}.
\eeaa
Together with the definition of $\de_{J+1}[\Pc]$, $\de_{J+1}[B]$ and $\de_{J+1}[\Bb]$, this yields
\beaa
\BEF^J_\de[A, B, \Pc, \Bb, \Ab] &\les& r_0^{15}\Big( \ep_J(\Sk_{J+1}+\Rk_{J+1}) +\ep_J^2+\ep_0^2\Big)+|a|\Sk^2_{J+1}\\
&&+r_0^{\frac{15}{4}}\Sk_{J+1}^{\frac{3}{2}}\Big(\ep_0+\sqrt{\ep_J}\sqrt{\Sk_{J+1}+\Rk_{J+1}}\Big)^{\frac{1}{2}}.
\eeaa
Also, in view of the definition of $\Rkint_k$ in section \ref{subsection:MainNormsM8}, and the one of $B^k_\de$ in section \ref{subsection:recallMorawetz-Energy}, we have
\beaa
\Rkint_{J+1}^2 &\les& r_0^3 \BEF^J_\de[A, B, \Pc, \Bb, \Ab].
\eeaa
We deduce 
\beaa
\Rkint_{J+1}^2 &\les& r_0^{18}\Big( \ep_J(\Sk_{J+1}+\Rk_{J+1}) +\ep_J^2+\ep_0^2\Big)+|a|r_0^3\Sk^2_{J+1}\\
&&+r_0^{\frac{27}{4}}\Sk_{J+1}^{\frac{3}{2}}\Big(\ep_0+\sqrt{\ep_J}\sqrt{\Sk_{J+1}+\Rk_{J+1}}\Big)^{\frac{1}{2}}
\eeaa
which concludes the proof of  \eqref{eq:InteriorcurvEstimatesThmM8} of Theorem \ref{prop:rpweightedestimatesiterationassupmtionThM8}.
\end{proof}

The rest of the chapter is devoted to the proof of Proposition \ref{proposition:EstimatesBBb-interior} and 
Proposition \ref{eq:proposition-EstimatesAAb-interior}. We first exhibit useful properties of  Bianchi pairs in sections \ref{section:linearization-Bianchi-pairs}, \ref{sec:hyperbolicestimatesforBianchipairs} and \ref{sec:bianchipairsforhigherorderderivatives}. Proposition \ref{proposition:EstimatesBBb-interior} and 
Proposition \ref{eq:proposition-EstimatesAAb-interior} are then proved respectively in section \ref{section:InteriorEstimatesBBb} and \ref{section:InteriorEstimatesAAb}.


\section{Linearization of second  and  third Bianchi pairs}
\label{section:linearization-Bianchi-pairs}


Observe that the first and fourth Bianchi  pair,  see section \ref{section:BianchiPairs-M8},    involve $A$, $B$, $\Bb$,  $\Ab$, $\Xh$ and $\Xi$, as well as quadratic terms, and are therefore already in linearized form. On the other hand, in addition to $B$ and $\Bb$, the second and third Bianchi pairs   contain  also $P$, $\tr X$, $\tr\Xb$, $H$ and $\Hb$ which are not  linearized quantities.  In this section, we linearize  the second and third Bianchi pairs based on commutation\footnote{Recall the definition of $\Lieb_X$ in  section \ref{section:horizLieDerivatives}.} with $\Lieb_\T$.     
    The result is stated in the following\footnote{See  also Remark \ref{remark:simpligyBianch-interior}.}.
   \begin{lemma}
  \lab{Lemma:linearizationbyLie_T}
   The quantities  
   \beaa
   \Bdot:=\Lieb_\T B, \qquad \Pdot:=\Lieb_\T P=\T(P), \qquad \Bbdot:=\Lieb_\T \Bb,
   \eeaa 
   verify the following equations:
   \begin{enumerate}
   \item The linearization by $\Lieb_\T$ of the second Bianchi pair  can be written in the form
     \bea
     \lab{eq:linearizationbyLie_T-B2}
     \begin{split}
  \nabc_3\Bdot  +\tr\Xb \Bdot =&  \DD\ov{\Pdot} +O(ar^{-2}) \ov{\Pdot}  +O(r^{-4}) \Ga_b+O(r^{-3} )\Lieb_\T\Hc \\
  &  + r^{-2} \dk^{\le 1}(\Ga_b\c \Rc_b),\\
   \nabc_4\Pdot  +\frac{3}{2}\tr X \Pdot  =&  \frac{1}{2}\DD\c \ov{\Bdot}  +O(ar^{-2}) \ov{\Bdot}+O(r^{-3} )  \dk^{\le 1} \Ga_b   + r^{-2} \dk^{\le 1}(\Ga_b\c \Rc_b).
   \end{split}
 \eea
 
  \item The linearization by $\Lieb_\T$ of the  third  Bianchi pair  can be written in the form 
  \bea
    \lab{eq:linearizationbyLie_T-B3}
  \begin{split}
\nabc_3\Pdot +      \frac{3}{2}\ov{\tr\Xb}   \Pdot    =& - \frac{1}{2}\ov{\DD}\c\Bbdot + O(ar^{-2} ) \Bbdot+O(r^{-3} )\dk^{\le 1} \Ga_b +r^{-1} \dk^{\le 1}(\Ga_b\c \Rc_b),\\
\nabc_4\Bbdot +\tr X\Bbdot  =& -\DD \Pdot + O(ar^{-2} )\Pdot +O(r^{-4} ) \Ga_b +O(r^{-3}) \Lieb_\T \Hbc\\
& + r^{-1} \dk^{\le 1}(\Ga_b\c \Rc_b).
 \end{split}
\eea
 \end{enumerate}
   \end{lemma}
   
 \begin{remark}\lab{remark-conf2}
Note that we can replace the operators $\DD,  \DDb $ by the corresponding conformal ones $\DDc,  \DDbc $ without changing the structure of the equations. Indeed, since $\Pdot$ has signature 0, we have $\DD\ov{\Pdot}=\DDc\ov{\Pdot}$ and $\DD \Pdot=\DDc \Pdot$, while for $\Adot$, $\Bdot$, $\Bbdot$ and $\Abdot$, we use the fact that $\DDc=\DD+sZ=\DD+O(ar^{-2})+\Ga_g$ so that the extra terms  do indeed not change the structure of the equations.
\end{remark}

   \begin{proof}
   We apply  first $\Lieb_\T $ to    the Bianchi identity involving $\nabc_3B$, i.e., see section \ref{section:BianchiPairs-M8},
   \beaa
   \nabc_3B  +\tr\Xb B &=& \DD\ov{P}+3\ov{P}H  +r^{-2}\Ga_b\c \Rc_b,
   \eeaa
   where we have used the fact that $\DDc\ov{P}=\DD\ov{P}$ by definition since $P$ has signature 0. Using the commutator estimates, see Lemma \ref{lemma:commutation-Lieb-nab}, we have $ [  \Lieb_\T, \DD] P=0$ since $P$ is a scalar. Also, $[  \Lieb_\T, \nab_3] B=  r^{-2} \dk^{\le 1 }(\Ga_b\c \Rc_b)$ in view of Lemma \ref{lemma:commLieb_TM8}. Therefore
      \beaa
 \nabc_3   \Bdot  + \tr\Xb \Bdot &=& \DD\ov{\Pdot}+3\ov{\Pdot}H+  3 \ov{P} \Lieb_\T H  +r^{-2} \dk^{\leq 1}(\Ga_b\c \Rc_b)\\
 &=&\DD\ov{\Pdot}+ O(ar^{-2} )\ov{\Pdot}+ O(r^{-3} ) \Lieb_\T H  + r^{-2} \dk^{\leq 1}(\Ga_b\c \Rc_b).
   \eeaa
 We write $\Lieb_\T H=  \Lieb_\T \Hc +\Lieb_\T\big(\frac{aq}{|q|^2}\Jk )=  \Lieb_\T \Hc +O(r^{-1}) \Ga_b  $ and hence
 \beaa
  \nabc_3\Bdot  +\tr\Xb \Bdot&=&  \DD\ov{\Pdot} +O(ar^{-2}) \ov{\Pdot} + O(r^{-3} )\Lieb_\T \Hc +O(r^{-4} ) \Ga_b   +r^{-2} \dk^{\leq 1}(\Ga_b\c \Rc_b)
 \eeaa
  as stated.
  
Next, we apply $\Lieb_\T $ to    the Bianchi identity involving $\nabc_4\Bb$, i.e., see section \ref{section:BianchiPairs-M8}, and we deduce as above 
      \beaa
 \nabc_4   \Bbdot  + \tr X \Bbdot &=&-\DD \Pdot -3\Pdot \Hb -3P\Lieb_\T \Hb  +r^{-1} \dk^{\leq 1}(\Ga_b\c \Rc_b)\\
 &=& -\DD \Pdot+ O(ar^{-2} )\Pdot +O(r^{-3})\Lieb_\T \Hb +O(r^{-4} ) \Ga_b + r^{-1} \dk^{\leq 1}(\Ga_b\c \Rc_b).
   \eeaa
   The two remaining equations  can be obtained in the same manner and are in fact easier. This concludes the proof of Lemma \ref{Lemma:linearizationbyLie_T}.
    \end{proof}


\section{Hyperbolic estimates for  Bianchi pairs}
\lab{sec:hyperbolicestimatesforBianchipairs}



   \subsection{Complex Hodge operators}
   \lab{section:Complexadjointoperators}
   

Consider the following complex Hodge  operators
\bea
\lab{def:complexadjointperators.M8}
\DDs_2 &=&-\frac 1 2 \DD\hot, \qquad \DDd_2 =\frac 1  2 \DDov \c, \qquad \DDs_1= -\DD, \qquad \DDd_1=\frac 1 2 \DDov\c.
\eea

\begin{remark}
\lab{remark:Comple-realHodgeoperators}
Observe that these operators are the complexified version of the real Hodge operators $\DDd_k$ and $\DDs_k$,  see  section \ref{section:HorizontalHodgeoperators}. More precisely, for $\Psi=\psi+i \dual \psi$ a complex $\sk_p$-tensor with   $p=1,2$,  or for a complex scalar $\Psi=(\psi+i \psi_*) $,  we have 
\beaa
\Re(\DDs_2\,  \Psi ) &=&2 \DDs_2 \, \Re(\Psi)= 2 \DDs_2\psi,\\
\Re(\DDd_2\,  \Psi ) &=& \DDd_2 \, \Re(\Psi)=  \DDd_2\psi,\\
\Re(\DDs_1\,  \Psi )&=&\DDs_1\,   (\psi, \psi_*),\\
\DDd_1\, \Psi&=& \div\psi+i \curl \psi.
\eeaa
In particular, above and in the rest of Part III, the same notation $\DDd_k$ and $\DDs_k$ is used for two different operators, i.e. if $\Psi$ is complex, then $\DDd_k\Psi$ and $\DDs_k\Psi$ should be interpreted as in \eqref{def:complexadjointperators.M8}, while if $\psi$ is real, $\DDd_k\psi$ and $\DDs_k\psi$ should be interpreted as in section \ref{section:HorizontalHodgeoperators}.
\end{remark}

The following lemma gives a sense to the fact that $\DDs_p$ is a formal  adjoint to $\DDd_p$.   
   \begin{lemma}
\lab{le:nab-DDd-DDs}
Given  $\Psi_{(1)}  \in  \sk_p$,   $\Psi_{(2)}  \in  \sk_{p-1}$, or $\Psi_{(2)}  \in  \sk_p$,   $\Psi_{(1)}  \in  \sk_{p-1}$, we have, in both cases, for  $p=1,2$,
 \bea\label{relation-nab-DDd-DDs}
 \Re\Big(  \DDd_p \Psi_{(1)}  \c \, \ov{\Psi_{(2)}}\Big)-\frac 1 2 \Re\Big(\Psi_{(1)}  \c \ov{ \DDs_p\, \Psi_{(2)} }\Big)&=& \nab \c \Re \Big(\Psi_{(1)} \c \ov{\Psi_{(2)}}\Big).
\eea
\end{lemma}

\begin{proof}
We  give below  the proof for both cases $p=1,2$.

{\bf Case   $p=2$:} We have for $\Psi_{(1)}=\psi_{(1)}+i \dual \psi_{(1)} \in \sk_2$ and $\Psi_{(2)}=\psi_{(2)}+i \dual \psi_{(2)} \in \sk_1$,
\beaa
&& \Re\Big(  \DDd_2 \Psi_{(1)}  \c \, \ov{\Psi_{(2)}}\Big)-\frac 1 2 \Re\Big(\Psi_{(1)}  \c \ov{ \DDs_2\, \Psi_{(2)} }\Big)\\
&=&\frac 1 2 \Re\Big( ( \ov{\DD} \c \Psi_{(1)} ) \c \, \ov{\Psi_{(2)}}\Big)+\frac 1 4 \Re\Big(\Psi_{(1)}  \c \ov{ \DD \hot\, \Psi_{(2)} }\Big)\\
&=&\Re\Big( ( \div \psi_{(1)}+ i \dual (\div \psi_{(1)}) ) \c \, \ov{\Psi_{(2)}}\Big)+\frac1 2 \Re\Big(\Psi_{(1)}  \c ( \nab \hot \psi_{(2)}- i \dual (\nab\hot \psi_{(2)}))\Big)\\
&=& 2( \div \psi_{(1)}) \c \, \psi_{(2)}+\psi_{(1)}  \c ( \nab \hot \psi_{(2)}).
\eeaa
Using that we have, for $f\in\sk_1$ and $u\in\sk_2$, 
 \beaa
 ( \nab\hot   f) \c   u  &=& \big(\nab_a f_b+\nab_b f_a- \de_{ab} \div f \big) u_{ab} =2(\nab_a f_b) u_{ab}=2 \nab_a (u_{ab} f_b)-2(\div u) \c f, 
 \eeaa
 we obtain
 \beaa
 \Re\Big(  \DDd_2 \Psi_{(1)}  \c \, \ov{\Psi_{(2)}}\Big)-\frac 1 2 \Re\Big(\Psi_{(1)}  \c \ov{ \DDs_2\, \Psi_{(2)} }\Big) &=& 2\nab\c(\psi_{(1)} \c \psi_{(2)})= \nab\c \Re(\Psi_{(1)} \c \ov{\Psi_{(2)}})
\eeaa
as stated.

{\bf Case   $p=1$:} we have for $\Psi_{(1)}=\psi_{(1)}+i \dual \psi_{(1)} \in \sk_1$ and $\Psi_{(2)}=a+i b \in \sk_0$,
 \beaa
&& \Re\Big(  \DDd_1 \Psi_{(1)}  \c \, \ov{\Psi_{(2)}}\Big)-\frac 1 2 \Re\Big(\Psi_{(1)}  \c \ov{ \DDs_1\, \Psi_{(2)} }\Big)\\
&=&\frac 1 2 \Re\Big( ( \ov{\DD} \c \Psi_{(1)})  \c \, \ov{\Psi_{(2)}}\Big)+ \frac 1 2 \Re\Big(\Psi_{(1)}  \c \ov{ \DD\, \Psi_{(2)} }\Big)\\
&=& \Re\Big( ( \div \psi_{(1)}+ i \curl \psi_{(1)} )  \c \, \ov{\Psi_{(2)}}\Big)+ \frac 1 2 \Re\Big(\Psi_{(1)}  \c (\nab a -\dual \nab b -i(\dual \nab a + \nab b) )\Big)\\
&=&a \ \div \psi_{(1)}+ b \ \curl \psi_{(1)}+  \psi_{(1)} \c (\nab a - \dual \nab b).
\eeaa
Using that we have, for $f\in\sk_1$, 
\beaa
f \c \nab a&=& f_c \nab_c a=\nab_c (a f_c)-a (\div f),\\
f \c \dual \nab b&=& - \dual f_c  \nab_c b=-\nab_c( b \dual f_c)+b (\curl f),
\eeaa
we obtain
 \beaa
 \Re\Big(  \DDd_1 \Psi_{(1)}  \c \, \ov{\Psi_{(2)}}\Big)-\frac 1 2 \Re\Big(\Psi_{(1)}  \c \ov{ \DDs_1\, \Psi_{(2)} }\Big)&=&\nab\c (a {\psi_{(1)}}+b \dual \psi_{(1)})= \nab\c \Re(\Psi_{(1)} \c \ov{\Psi_{(2)}}),
\eeaa
as stated. This concludes the proof of Lemma \ref{le:nab-DDd-DDs}.
\end{proof}

We also derive the following complex, non-integrable  version 
 of the  identities   \eqref{eq:dcalident}.
\begin{lemma}
\lab{eq:dcalident-complex}
The following formulas hold true
\bea
\bsplit
\DDd_1\DDs_1\, \Psi&= -\lap_0 \Psi    - \frac 1 2 i    (\atrch\nab_3+\atrchb \nab_4) \Psi,\\
\DDs_1\, \DDd_1 \Psi&=  -\lap_1 \Psi + \frac 1 2 i   (\atrch\nab_3+\atrchb \nab_4) \Psi  + \Kh  \Psi,\\
\DDd_2\DDs_2\,  \Psi&= - \lap_1\Psi - \frac 1 2 i    (\atrch\nab_3+\atrchb \nab_4) \Psi - \Kh \Psi,\\
\DDs_2 \DDd_2 \Psi &=-\lap_2 \psi   + \frac 1 2 i    (\atrch\nab_3+\atrchb \nab_4) \Psi + 2 \Kh \Psi,
\end{split}
\eea
with the scalar $\Kh$ given by the formula\footnote{Recall that $\Kh$ coincides    with the   Gauss curvature $K$ in the integrable case.} \eqref{eq:definition-K-in}.
\end{lemma}

\begin{proof}
We start with the identities 
\beaa
\bsplit
\DDd_1\DDs_1\, \Psi&=-\lap_0 \Psi - \frac 1 2 i \in^{ab}[\nab_a, \nab_b] \Psi,\\
\DDs_1\, \DDd_1 \Psi&=-\lap_1 \Psi + \frac 1 2 i \in_{ab} [\nab_a, \nab_b] \Psi,\\
\DDd_2\DDs_2\,  \Psi&=- \lap_1\Psi -  \frac 1  2 i  \in_{ab} [\nab_a, \nab_b]\Psi, \\
\DDs_2 \DDd_2 \Psi &=-\lap_2\Psi  +\frac 1 2 i \in_{ab} [ \nab_a, \nab_b] \Psi.
\end{split}
\eeaa
The last three identities are immediate consequences of 
Lemma
\ref{le:Bochner-nointegrable} and Remark \ref{remark:Comple-realHodgeoperators}. It thus remains to check the first one. For $\Psi\in \sk_0$, we have
\beaa
  \DDd_1\DDs_1  \Psi&=&-\frac 1 2 \ov{\DD}^a \DD_a \Psi = -\frac 1 2 
  (\nab^a-i\dual \nab^a)  (\nab_a+i\dual \nab_a) \Psi\\
  &=&  -\frac 1 2  \lap \Psi -\frac 1 2  \dual \nab^a (\dual \nab_ a \Psi)
  + \frac 1 2  i \in^{ab}(\nab_b\nab_a \Psi -\nab_a\nab_b \Psi)\\
  &=& -\lap\Psi  - \frac 1 2 i \in_{ab} [\nab_a, \nab_b] \Psi,
  \eeaa
  as stated.

We  rewrite  the formulas  above  by   making  use of the  Gauss  type formula  of  Proposition \ref{Gauss-equation-2-tensors}   according to which we  have, for $\Psi\in \sk_p(\CCC)$ with $p=0, 1, 2$, 
\beaa
 \in^{ab}[\nab_a,  \nab_b] \Psi =  (\atrch\nab_3+\atrchb \nab_4) \Psi -    2 i  \, p   \, \Kh \Psi. 
 \eeaa
We deduce
\beaa
\bsplit
\DDd_1\DDs_1\, \Psi&=-\lap_0 \Psi - i\frac 1 2  \in^{ab}[\nab_a, \nab_b] \Psi= -\lap_0 \Psi    - \frac 1 2 i    (\atrch\nab_3+\atrchb \nab_4) \Psi,\\
\DDs_1\, \DDd_1 \Psi&=-\lap_1 \Psi +\frac 1 2 i \in^{ab} [\nab_a, \nab_b] \Psi=  -\lap_1 \Psi + \frac 1 2 i   (\atrch\nab_3+\atrchb \nab_4) \Psi  + \Kh  \Psi,\\
\DDd_2\DDs_2\,  \Psi&=- \lap_1\Psi -  \frac 1  2 i  \in_{ab} [\nab_a, \nab_b]\Psi= - \lap_1\Psi - \frac 1 2 i    (\atrch\nab_3+\atrchb \nab_4) \Psi - \Kh \Psi,\\
\DDs_2 \DDd_2 \Psi &=-\lap_2\Psi  +\frac 1 2 i \in_{ab} [ \nab_a, \nab_b] \Psi=-\lap_2 \psi   + \frac 1 2 i    (\atrch\nab_3+\atrchb \nab_4) \Psi + 2 \Kh \Psi,
\end{split}
\eeaa
as stated. This concludes the proof of Lemma \ref{eq:dcalident-complex}.
\end{proof}


\subsection{Bianchi pairs using the Hodge operators $\DDd_p, \DDs_p$}


We rewrite below  the  Bianchi   equations, see section \ref{section:BianchiPairs-M8}, using  the complex  Hodge operators  introduced  in section \ref{section:Complexadjointoperators}. We split them in pairs as follows.
 \begin{definition} 
\lab{definition:BianchiPairs}
We define the following pairs of Bianchi identities:
\begin{enumerate}
\item The first pair, involving $A$ and $B$:
\bea\label{eq:first-pair-A-B-lin}
\begin{split}
 \nabc_3A +\frac{1}{2}\tr\Xb A&=  -\DDs_2\, B+  O(ar^{-2}) B  +O(r^{-3} ) \Xh +r^{-2} \Ga_b\c \Rc_b ,\\
\nabc_4B  +2\ov{\tr X} B &=    \DDd_2 A+  O(ar^{-2}) A    +O(r^{-3} ) \Xi +r^{-2} \Ga_b\c \Rc_b.
\end{split}
\eea

\item The second pair, involving $B$ and $\ov{P}$:
\bea\label{eq:second-pair-B-P}
\begin{split}
\nabc_3B  +\tr\Xb B &=-\DDs_1\,  \ov{P}+3\ov{P}H  +r^{-2} \Ga_b\c \Rc_b, \\
\nabc_4\ov{P}   +\frac{3}{2}\ov{\tr X}\,  \ov{P}   &=  \DDd_1 B + O( a r^{-2}) \ov{B}  +r^{-2} \Ga_b\c \Rc_b.
\end{split}
\eea

\item The third pair, involving $P$ and $\Bb$:
\bea\label{eq:third-pair-P-Bb}
\begin{split}
\nabc_3P+ \frac{3}{2}\ov{\tr\Xb} P &= - \DDd_1\Bb +  O(ar^{-2} )\Bb    +r^{-1} \Ga_b\c \Rc_b, \\
\nabc_4\Bb+ \tr X\Bb  &=\DDs_1\,  P -3P\Hb  +r^{-1} \Ga_b\c \Rc_b.
\end{split}
\eea

\item The fourth pair, involving $\Bb$ and $\Ab$:
\bea\label{eq:fourth-pair-Bb-Ab-lin}
\begin{split}
\nabc_3\Bb  +2\ov{\tr\Xb}\,\Bb &=  -\DDd_2\Ab  +O(a r^{-2})\Ab -O(r^{-3})  \,\Xib + \Ga_b\c \Rc_b,\\
\nabc_4\Ab  +\frac{1}{2}\tr X \Ab   &=\DDs_2\, \Bb  +O(ar^{-2} ) \Bb +O(r^{-3} )\Xbh + \Ga_b\c \Rc_b.
\end{split}
\eea
\end{enumerate}
\end{definition}


\subsection{Main lemma for Bianchi pairs}


We start with the following definition exhibiting the general form of Bianchi pairs. 
\begin{definition}
\lab{definition:GenBianchPairs}
   We consider the following  general  Bianchi pairs in $\MM$, which generalize the Bianchi pairs written as in  Definition \ref{definition:BianchiPairs}:
   \begin{itemize}
   \item   For  $\Psi_{(1)}\in\sk_p$, $\Psi_{(2)}\in\sk_{p-1}$, and    $F_{(1)}\in\sk_{p}$, $F_{(2)}\in \sk_{p-1}$,   
   \bea\lab{eq:modelbainchipairequations11-simple}
\begin{split}
\nabc_3(\Psi_{(1)})+c_{(1) }\tr \Xb\Psi_{(1)} &= -\DDs_p\, \Psi_{(2)}  +F_{(1)},\\[2mm]
\nabc_4(\Psi_{(2)})+c_{(2)} \ov{\tr X}\Psi_{(2)} &= \DDd_p\, \Psi_{(1)}   +F_{(2)}.
\end{split}
\eea

\item   For  $\Psi_{(1)}\in\sk_{p-1}$, $\Psi_{(2)}\in\sk_{p}$, and    $F_{(1)}\in\sk_{p-1}$, $F_{(2)}\in \sk_{p}$,         
\bea\lab{eq:modelbainchipairequations12-simple}
\begin{split}
\nabc_3(\Psi_{(1)})+c_{(1)}\ov{\tr \Xb}\Psi_{(1)}&=- \DDd_p\, \Psi_{(2)} +F_{(1)},\\[2mm]
\nabc_4(\Psi_{(2)})+c_{(2)}\tr X\Psi_{(2)} &=\DDs_p\, \Psi_{(1)}  +F_{(2)}.
\end{split}
\eea
 \end{itemize}
 \end{definition}

 \begin{remark}
 \lab{remark:funnysecondpair}
 Note that   the  third and fourth Bianchi pairs,  according to Definition \ref{definition:BianchiPairs},  are of the type  \eqref{eq:modelbainchipairequations12-simple}. Also, the first and 
  second Bianch pairs are of the type  \eqref{eq:modelbainchipairequations11-simple}  provided we write
  $\Psi_{(1)}= B$, $\Psi_{(2)}=\ov{P} $  for the second Bianchi pair.
  \end{remark}
 
\begin{remark} 
We can also define the general Bianchi pairs in a conformally invariant way. 
 Recall  the  conformal operators defined for tensors $f$ of signature $s$ by $\nabc_a f=\nab_a f+  s \ze_a  f$
    and $\DDc= \nabc_a+i \dual \nabc_a$.  By introducing the conformal  operators
   \beaa
\DDsc_2 &=&-\frac 1 2 \DDc\hot, \quad \DDdc_2 =\frac 1  2 \DDcov \c, \quad \DDsc_1= -\DDc, \quad \DDdc_1=\frac 1 2 \DDcov\c,
\eeaa
we can define the above pairs with $\DDsc_p$ and $\DDdc_p$, i.e.
     \bea\lab{eq:modelbainchipairequations11-simple-c}
\begin{split}
\nabc_3(\Psi_{(1)})+c_{(1) }\tr \Xb\Psi_{(1)} &= -\DDsc_p\, \Psi_{(2)}  +F_{(1)},\\[2mm]
\nabc_4(\Psi_{(2)})+c_{(2)} \ov{\tr X}\Psi_{(2)} &= \DDdc_p\, \Psi_{(1)}   +F_{(2)},
\end{split}
\eea
and
\bea\lab{eq:modelbainchipairequations12-simple-c}
\begin{split}
\nabc_3(\Psi_{(1)})+c_{(1)}\ov{\tr \Xb}\Psi_{(1)}&=- \DDdc_p\, \Psi_{(2)} +F_{(1)},\\[2mm]
\nabc_4(\Psi_{(2)})+c_{(2)}\tr X\Psi_{(2)} &=\DDsc_p\, \Psi_{(1)}  +F_{(2)}.
\end{split}
\eea
Note that the use   of  these conformally invariant   horizontal Hodge operators 
does not modify the structure of the terms $F_{(1)}$ and $F_{(2)}$, see Remark  \ref{remark-conf2}. 
\end{remark}

\begin{lemma}
\lab{Le:BasicBianchiPairs.simple}
Let $\Psi_{(1)}, \Psi_{(2)}$ verifying either one of the equations 
   \eqref{eq:modelbainchipairequations11-simple},   \eqref{eq:modelbainchipairequations12-simple} for positive real numbers  $c_{(1)}$ and $c_{(2)}$, with $\Psi_{(1)}$ of  signature  $k$ and $\Psi_{(2)}$ of  signature  $k-1$. Then denoting
   \beaa
   \La_{(1)}:= -2c_{(1)}+ 1+\frac b 2, \qquad  \La_{(2)} :=  -2c_{(2)}+ 1+\frac b 2,
   \eeaa
 the following  pointwise  identity holds true for any real $b$:
   \begin{enumerate}
\item If $\Psi_{(1)}, \Psi_{(2)}$ verify equation \eqref{eq:modelbainchipairequations11-simple}, then
\bea\lab{eq:basicdivergenceidentitybianchipairrpweightedestimate-complex-1.simple-0}
\nn&& \Div\left(\frac 1 2|q|^b|\Psi_{(1)}|^2e_3+|q|^b|\Psi_{(2)}|^2e_4-2 |q|^b \Re ( \Psi_{(1)} \c \ov{\Psi_{(2)}})\right)\\
\nn&=&\frac 1 2 |q|^{b}  \La_{(1)} \trchb |\Psi_{(1)}|^2+|q|^{b}  \La_{(2)}\trch |\Psi_{(2)}|^2+  |q|^b(2k-1)\big( \omb |\Psi_{(1)}|^2- 2\om  |\Psi_{(2)}|^2\big)\\
\nn&&+O(ar^{b-2}) \c  \Re ( \Psi_{(1)} \c \ov{\Psi_{(2)}})+  |q|^b \Re\Big( F_{(1)} \c \ov{\Psi_{(1)}} \Big)+ 2|q|^b\Re\Big(  F_{(2)} \c \ov{\Psi_{(2)}}\Big)\\
&&+r^{b} \Ga_b\Big(|\Psi_{(1)}|^2 +|\Psi_{(1)}||\Psi_{(2)}|+r^{-1}|\Psi_{(2)}|^2\Big).
\eea

\item If $\Psi_{(1)}, \Psi_{(2)}$ verify equation \eqref{eq:modelbainchipairequations12-simple}, then 
\bea\lab{eq:basicdivergenceidentitybianchipairrpweightedestimate-complex-2.simple-0}
\nn&& \Div\left(|q|^b|\Psi_{(1)}|^2e_3+\frac 1 2|q|^b|\Psi_{(2)}|^2e_4+2 |q|^b \Re ( \Psi_{(1)} \c \ov{\Psi_{(2)}})\right)\\
\nn&=& |q|^{b} \La_{(1)} \trchb |\Psi_{(1)}|^2+\frac 1 2 |q|^{b}  \La_{(2)}\trch |\Psi_{(2)}|^2+ |q|^b(2k-1)\big(  2\omb  |\Psi_{(1)}|^2- \om  |\Psi_{(2)}|^2\big)\\
\nn&&+O(ar^{b-2}) \c  \Re ( \Psi_{(1)} \c \ov{\Psi_{(2)}})+2|q|^b  \Re\Big( F_{(1)} \c \ov{\Psi_{(1)}} \Big)+|q|^b \Re\Big(  F_{(2)} \c \ov{\Psi_{(2)}}\Big)\\
&&+r^{b} \Ga_b\Big(|\Psi_{(1)}|^2 +|\Psi_{(1)}||\Psi_{(2)}|+r^{-1}|\Psi_{(2)}|^2\Big).
\eea
\end{enumerate}
\end{lemma}

\begin{remark}  
Note that the two  identities differ by a factor of $\frac 1 2$ when interchanging $\Psi_{(1)}$ and $\Psi_{(2)}$, and by  the sign of the third term in the first line.
\end{remark}

\begin{proof}
We note  first
\beaa
\D_\ga e_4^\ga &=& -\frac{1}{2}\g(\D_4e_4, e_3)-\frac{1}{2}\g(\D_3e_4, e_4)+\ga^{ab}\g(\D_a  e_4, e_b)= \trch  -2\om,\\
\D_\ga e_3^\ga &=& -\frac{1}{2}\g(\D_4e_3, e_3)-\frac{1}{2}\g(\D_3e_3, e_4)+\ga^{ab}\g(\D_a e_3, e_b)= \trchb  -2\omb.
\eeaa
We now calculate, 
\beaa
&&\Div\Big(|q|^b|\Psi_{(1)}|^2e_3\Big)\\
&=& 2|q|^b \Re\Big(  \nab_3\Psi_{(1)}\c \ov{\Psi_{(1)}}\Big)+\frac b 2 |q|^{b-2}\big( \ov{q} e_3(q)+ q e_3(\ov{q})\big) |\Psi_{(1)}|^2+|q|^b|\Psi_{(1)}|^2\D_\ga e_3^\ga\\
&=& 2|q|^b \Re\Big(  \nab_3\Psi_{(1)}\c \ov{\Psi_{(1)}}\Big)+\frac b 2 |q|^{b-2}\big( \ov{q} e_3(q)+ q e_3(\ov{q})\big) |\Psi_{(1)}|^2+|q|^b( \trchb  -2\omb)|\Psi_{(1)}|^2\\
&=& 2|q|^b \Re\Big(  \nab_3\Psi_{(1)}\c \ov{\Psi_{(1)}} \Big)+|q|^{b} \left( 1+\frac b 2\right) \trchb |\Psi_{(1)}|^2- 2\omb  |q|^b |\Psi_{(1)}|^2\\
&&+\frac b 2 |q|^{b-2}\big( \ov{q} e_3(q)+ q e_3(\ov{q})- \trchb |q|^2 \big) |\Psi_{(1)}|^2.
\eeaa
Similarly,
\beaa
\Div\Big(|q|^b|\Psi_{(2)}|^2e_4\Big)&=& 2|q|^b \Re\Big(  \nab_4\Psi_{(2)}\c \ov{\Psi_{(2)}} \Big)+|q|^{b} \left( 1+\frac b 2\right) \trch |\Psi_{(2)}|^2- 2\om  |q|^b |\Psi_{(2)}|^2\\
&&+\frac b 2 |q|^{b-2}\big( \ov{q} e_4(q)+ q e_4(\ov{q})- \trch |q|^2 \big) |\Psi_{(2)}|^2.
\eeaa
Note that\footnote{In particular, we use the fact that $-(q+\ov{q} ) -\Re(-\frac{2}{\ov{q}})|q|^2=-(q+\ov{q} )+2 \Re(q)=0$.}
\beaa
\ov{q} e_3(q)+ q e_3(\ov{q})- \trchb |q|^2 = r^2\Ga_b, \qquad 
\ov{q} e_4(q)+ q e_4(\ov{q})- \trch |q|^2 = r^2\Ga_g.
\eeaa
Hence
\beaa
\Div\Big(|q|^b|\Psi_{(1)}|^2e_3\Big) &=& 2|q|^b \Re\Big(  \nab_3\Psi_{(1)}\c \ov{\Psi_{(1)}} \Big)+|q|^{b} \left( 1+\frac b 2\right) \trchb |\Psi_{(1)}|^2- 2\omb  |q|^b |\Psi_{(1)}|^2\\
&&+r^b\Ga_b |\Psi_{(1)}|^2,\\
\Div\Big(|q|^b|\Psi_{(2)}|^2e_4\Big)&=& 2|q|^b \Re\Big(  \nab_4\Psi_{(2)}\c \ov{\Psi_{(2)}} \Big)+|q|^{b} \left( 1+\frac b 2\right) \trch |\Psi_{(2)}|^2- 2\om  |q|^b |\Psi_{(2)}|^2 \\
&&+r^{b-1}\Ga_b |\Psi_{(2)}|^2.
\eeaa

We now consider the two cases:

{\bf Case 1.}  Consider the case when   $\Psi_{(1)}, \Psi_{(2)}$ verify equations \eqref{eq:modelbainchipairequations11-simple}, with $\Psi_{(1)}$ of  signature  $k$ and $\Psi_{(2)}$ of  signature  $k-1$.  Then
\beaa
\nabc_3\Psi_{(1)}&=& \nab_3 \Psi_{(1)}-2k \omb \Psi_{(1)},\\
\nabc_4 \Psi_{(2)}&=& \nab_4 \Psi_{(2)}+2(k-1)\om \Psi_{(2)}.
\eeaa
Hence, using the equation for $\Psi_{(1)}$,
\beaa
&& \Re\Big(  \nab_3\Psi_{(1)}\c \ov{\Psi_{(1)}} \Big)\\
&=&  \Re\Big(  \big(-c_{(1) }\tr \Xb\Psi_{(1)}+2k \omb \Psi_{(1)} -\DDs_p\, \Psi_{(2)}+F_{(1)} \big)\c \ov{\Psi_{(1)}} \Big)\\
 &=&-c_{(1)}\trchb |\Psi_{(1)}|^2+2k \omb |\Psi_{(1)}|^2-  \Re \big( \DDs_p\, \Psi_{(2)}\c \ov{\Psi_{(1)}} \big)+  \Re\Big(F _{(1)} \c \ov{\Psi_{(1)}} \Big),
\eeaa
and using the equation for $\Psi_{(2)}$,
\beaa
&&\Re\Big(  \nab_4\Psi_{(2)}\c \ov{\Psi_{(2)}} \Big) \\
&=&  \Re\Big(  \big(- c_{(2)} \ov{\tr X} \Psi_{(2)}-2(k-1)\om \Psi_{(2)} +\DDd_p  \Psi_{(1)} \c \ov{\Psi_{(2)}} +  F_{(2)}   \big)    \c \ov{\Psi_{(2)}}\Big) \\
&=& - c_{(2)} \trch  |\Psi_{(2)}|^2 - 2(k-1)\om   |\Psi_{(2)}|^2+\Re\big(\DDd_p  \Psi_{(1)} \c \ov{\Psi_{(2)} }\big)+ \Re\big(  F_{(2)} \c \ov{\Psi_{(2)}}\big),
\eeaa
we deduce,
\beaa
&&\Div\Big(|q|^b|\Psi_{(1)}|^2e_3\Big)\\
&=& 2|q|^b\Big( -c_{(1)}\trchb |\Psi_{(1)}|^2+2k \omb |\Psi_{(1)}|^2-  \Re \big( \DDs_p\, \Psi_{(2)}\c \ov{\Psi_{(1)}} \big)\Big)\\
&&+|q|^{b} \left( 1+\frac b 2\right) \trchb |\Psi_{(1)}|^2- 2\omb  |q|^b |\Psi_{(1)}|^2 + 2 |q|^b  \Re\big(F _{(1)} \c \ov{\Psi_{(1)}} \big)+r^{b}\Ga_b |\Psi_{(1)}|^2\\
&=&|q|^b \trchb |\Psi_{(1)}|^2\left(-2 c_{(1)}+1 +\frac b 2 \right)+ 2 (2k-1) |q|^b \omb  |\Psi_{(1)}|^2 - 2|q|^b \Re \big( \DDs_p\, \Psi_{(2)}\c \ov{\Psi_{(1)}} \big)\\
&&+ 2 |q|^b  \Re\big(F _{(1)} \c \ov{\Psi_{(1)}} \big) +r^b\Ga_b |\Psi_{(2)}|^2 
\eeaa
and
\beaa
&&\Div\Big(|q|^b|\Psi_{(2)}|^2e_4\Big)\\
&=&  2|q|^b\Big(- c_{(2)} \trch  |\Psi_{(2)}|^2 - 2(k-1)\om   |\Psi_{(2)}|^2+\Re\big(\DDd_p  \Psi_{(1)} \c \ov{\Psi_{(2)} }\big)\Big)\\
&&  2 |q|^b  \Re\big(F _{(2)} \c \ov{\Psi_{(2)}} \big) +|q|^{b} \left( 1+\frac b 2\right) \trch |\Psi_{(2)}|^2- 2\om  |q|^b |\Psi_{(2)}|^2
 +r^{b}|\Psi_{(2)}|^2\\
 &=&|q|^b \trch |\Psi_{(2)}|^2\left(-2  c_{(2)}+ 1 +\frac b 2 \right)- 2 (2k-1)\om |\Psi_{(2)}|^2 + 2|q|^b \Re\big(\DDd_p  \Psi_{(1)} \c \ov{\Psi_{(2)} }\big)\\
 &&+ 2 |q|^b  \Re\big(F _{(2)} \c \ov{\Psi_{(2)}} \big) +r^{b-1}\Ga_b|\Psi_{(2)}|^2.
\eeaa
We deduce, making use of   the  relation \eqref{relation-nab-DDd-DDs}, 
\beaa
&&\frac 1 2 \Div\Big(|q|^b|\Psi_{(1)}|^2e_3\Big)+\Div\Big(|q|^b|\Psi_{(2)}|^2e_4\Big)-2|q|^b \nab \c \Re \Big(\Psi_{(1)} \c \ov{\Psi_{(2)}}\Big)\\
&=& \frac 1 2  |q|^b  \La_{(1)}  \trchb |\Psi_{(1)}|^2  +   |q|^b \La_{(2)}  \trch |\Psi_{(2)}|^2+ (2k-1) |q|^b\Big( \omb  |\Psi_{(1)}|^2 - 2\om |\Psi_{(2)}|^2\Big) \\
&&+  |q|^b \Re\Big( F_{(1)} \c \ov{\Psi_{(1)}} \Big)+ 2|q|^b\Re\Big(  F_{(2)} \c \ov{\Psi_{(2)}}\Big)+r^{b} \Ga_b\Big(|\Psi_{(1)}|^2+r^{-1}|\Psi_{(2)}|^2\Big).
\eeaa
Finally, using Lemma \ref{lemma:divergence-spacetime-horizontal-ch1}, i.e. $\D^\a f_\a= \nab^a f_a + (\eta +\etab) \c f$,
\beaa
|q|^b\nab \c \Re \Big(\Psi_{(1)} \c \ov{\Psi_{(2)}}\Big)&=&\nab \c \Big( |q|^b \Re ( \Psi_{(1)} \c \ov{\Psi_{(2)}})\Big)-\nab(|q|^b) \c \Re ( \Psi_{(1)} \c \ov{\Psi_{(2)}})\\
&=&\Div \Big( |q|^b \Re ( \Psi_{(1)} \c \ov{\Psi_{(2)}})\Big)-\left( 1+\frac b2\right) |q|^b(\eta+\etab) \c  \Re ( \Psi_{(1)} \c \ov{\Psi_{(2)}})\\
&& +r^{b} \Ga_b \c \Re ( \Psi_{(1)} \c \ov{\Psi_{(2)}})
\eeaa
where we used that $\nab(|q|^b)=\frac b2 (\eta+\etab)|q|^b+ r^{b}\Ga_b$. Combining with the above, we obtain \eqref{eq:basicdivergenceidentitybianchipairrpweightedestimate-complex-1.simple-0}.

{\bf Case 2.} Consider the case when  $\Psi_{(1)}, \Psi_{(2)}$ verify equation \eqref{eq:modelbainchipairequations12-simple}, with $\Psi_{(1)}$ of  signature  $k$ and $\Psi_{(2)}$ of  signature  $k-1$. Then 
\beaa
&& \Re\Big(  \nab_3\Psi_{(1)}\c \ov{\Psi_{(1)}} \Big)\\
&=&  \Re\Big(  \big(-c_{(1) } \ov{\tr \Xb}\Psi_{(1)}+2k \omb \Psi_{(1)} -\DDd_p\, \Psi_{(2)}+ F_{(1)} \big)\c \ov{\Psi_{(1)}}  \Big)\\
 &=&-c_{(1)}\trchb |\Psi_{(1)}|^2+2k \omb |\Psi_{(1)}|^2-  \Re \Big( \DDd_p\, \Psi_{(2)}\c \ov{\Psi_{(1)}} \Big)+
   \Re\Big( F_{(1)} \c \ov{\Psi_{(1)}} \Big),
\\
\\
&&\Re\Big(  \nab_4\Psi_{(2)}\c \ov{\Psi_{(2)}} \Big) \\
&=&  \Re\Big(  \big(- c_{(2)}\tr X  \Psi_{(2)}-2(k-1)\om \Psi_{(2)} +\DDs_p  \Psi_{(1)}  \big)      \c \ov{\Psi_{(2)}} + F_{(2)}      \c \ov{\Psi_{(2)}}\Big)\\
&=& - c_{(2)} \trch  |\Psi_{(2)}|^2 - 2(k-1)\om   |\Psi_{(2)}|^2+\Re\Big(\DDs_p  \Psi_{(1)} \c \ov{\Psi_{(2)} }\Big)
+ \Re\Big(  F_{(2)} \c \ov{\Psi_{(2)}}\Big).
\eeaa
Using again relation \eqref{relation-nab-DDd-DDs} and Lemma \ref{lemma:divergence-spacetime-horizontal-ch1},  we obtain \eqref{eq:basicdivergenceidentitybianchipairrpweightedestimate-complex-2.simple-0}, which concludes  the proof of Lemma \ref{Le:BasicBianchiPairs.simple}.
\end{proof} 
  
\begin{remark}
\lab{Remark:Bianchi-conditions}
Recall that  in  Kerr, we have
\beaa
\omb=0, \quad\om=-\frac{a^2\cos^2\th (r-m)+mr^2-a^2r}{|q|^4}<0, \quad \trchb=-\frac{2r}{|q|^2}, \quad \trch= \frac{2 r\De}{|q|^4}.
\eeaa
The dominant terms on the right hand side of \eqref{eq:basicdivergenceidentitybianchipairrpweightedestimate-complex-1.simple-0} and \eqref{eq:basicdivergenceidentitybianchipairrpweightedestimate-complex-2.simple-0} are given by, modulo a factor of $2$,
\beaa
|q|^{-b}J&:=&  \La_{(1)} \trchb |\Psi_{(1)}|^2+\frac 1 2   \La_{(2)}\trch |\Psi_{(2)}|^2+ (2k-1)\big(  2\omb  |\Psi_{(1)}|^2- \om  |\Psi_{(2)}|^2\big).
\eeaa
Then:
\begin{enumerate}
\item  The first two Bianchi pairs  are applied in situations 
 where $ \Psi_{(2)}$ is already  under control and  $2k-1$ is strictly positive, see part 1 in Proposition \ref{Prop:Bainchi-pairsEstimates-integrated} below. To derive estimates 
 for  $\Psi_{(1)}$,  we  need   the coefficient  
 $ \La_{(1)} \trchb +2(2k-1)\omb$  to be strictly negative. 
   Thus, since $\omb=\Ga_b$ and $ \trchb<0$,    we need to choose $b$ such that 
 \beaa
  \La_{(1)}=-2c_{(1)}+ 1+\frac b 2>0.
 \eeaa

 \item   By contrast,  the last  two Bianchi pairs  are applied in situations where $\Psi_{(1)}$ is under control and  $2k-1$ is strictly negative, see part 2 in Proposition \ref{Prop:Bainchi-pairsEstimates-integrated} below. To derive estimates  for  $\Psi_{(2)}$,  we  need   the coefficient  
 $ \La_{(2)} \trch -2(2k-1)\om$  to be strictly negative. Since $\om<0$ and is non degenerate near $r=r_+$, and since $\trch>0$ for $r>r_+$, it thus suffices to choose $b$ such that 
  \beaa
  \La_{(2)}=-2c_{(2)}+ 1+\frac b 2<0.
 \eeaa 
\end{enumerate}
\end{remark}


\subsection{Main integrated estimates for the Bianchi pairs}


\begin{lemma}[Divergence lemma]
\lab{lemma:basicdivergenceidentitybianchipairrpweightedestimate}
Consider  a vectorfield $ X$  in $\MM(\tau_1, \tau_2) $. We have
\beaa
- \int_{\AA(\tau_1, \tau_2)}  \g(X, N)    -    \int_{\Si(\tau_2) } \g(X, N) -\int_{\Si_*(\tau_1, \tau_2)}\g(X, N) + \int_{\Si(\tau_1)} \g(X, N)  = \int_{\MM(\tau_1, \tau_2) } \Div(X),
\eeaa
where $N$ is  the normal to the boundary   such that $\g(N, e_3) =-1$. 

We rewrite it in the form
\beaa
-\int_{\pr^+\MM(\tau_1, \tau_2)} \g(X, N) + \int_{\pr^-\MM(\tau_1, \tau_2)} \g(X, N)= \int_{\MM(\tau_1, \tau_2)} \Div(X),
\eeaa
where
\beaa
\pr^+\MM(\tau_1, \tau_2)=\AA(\tau_1, \tau_2)\cup \Si(\tau_2) \cup \Si_*(\tau_1, \tau_2), \qquad  \pr^-\MM(\tau_1, \tau_2)=  \Si(\tau_1).
\eeaa
\end{lemma}

\begin{proof}
Immediate consequence of   the standard divergence lemma. 
\end{proof}

We use the divergence lemma to obtain integrated estimates for the Bianchi pairs in the following proposition.

\begin{proposition}
\lab{Prop:Bainchi-pairsEstimates-integrated}
The following estimates\footnote{Note that the roles of  $\Psi_{(1)}$ and $\Psi_{(2)}$ are inverted in the two estimates. This is due to the fact that,  in applications,  $\Psi_{(2)}$  is already under control for the first Bianchi pair,    while  $\Psi_{(1)}$ is already under control for the  second Bianchi pair.}  hold true in $\MM=\MM(\tau_1, \tau_2)$.
\begin{enumerate}
 \item Let $\Psi_{(1)}, \Psi_{(2)}$ verifying equations \eqref{eq:modelbainchipairequations11-simple}
   for positive real numbers  $c_{(1)}$ and $c_{(2)}$ with $2k-1>0$.
    Let $b$ such that   $ \La_{(1)}=-2 c_{(1)} +1+\frac b 2 > 0$. Then, the following integrated estimate holds:
   \bea\label{eq:general-proposition-first-bianchi-pairs}
\begin{split}
& \int_{\MM } r^{b-1} |\Psi_{(1)}|^2+\int_{\pr^+\MM} r^{b}|\Psi_{(1)}|^2\\
\les &\int_{\MM} r^{b-1} |\Psi_{(2)}|^2  +\left|\int_{\MM}|q|^{b}\Re(F_{(1)}\c\ov{\Psi_{(1)}})\right|+\int_{\MM}r^b|F_{(2)}||\Psi_{(2)}| \\&  + \int_{\pr^-\MM} \big( r^{b}|\Psi_{(1)}|^2+  r^{b-2}|\Psi_{(2)}|^2\big).
\end{split}
\eea 
      
\item  Let $\Psi_{(1)}, \Psi_{(2)}$ verifying equations 
 \eqref{eq:modelbainchipairequations12-simple} for positive real numbers  $c_{(1)}$ and $c_{(2)}$, with $2k-1<0$.
  Let $b$ such that $ \La_{(2)}=-2 c_{(2)} +1+\frac b 2 < 0$. Then the following integrated estimate holds:
  \bea\label{eq:general-proposition-second-bianchi-pairs}
\begin{split}
& \int_{\MM } r^{b-1}\left(1+|2k-1|\frac{m}{r}\right) | \Psi_{(2)}|^2 +  \int_{\pr^+\MM } r^{b-2} |\Psi_{(2)}|^2\\ 
\les &\int_{\MM} r^{b-1} |\Psi_{(1)}|^2    +\left|\int_{\MM}|q|^{b}\Re(F_{(2)}\c\ov{\Psi_{(2)}})\right|    +\int_{\MM}r^{b}|F_{(1)}||\Psi_{(1)}|  \\
&+  \int_{\pr^-\MM} \big( r^{b}|\Psi_{(1)}|^2+  r^{b-2}|\Psi_{(2)}|^2\big).
 \end{split}
\eea
\end{enumerate}
\end{proposition}

\begin{proof} 
We proceed as follows. 

{\bf Case 1.} We first consider $\Psi_{(1)}, \Psi_{(2)}$ verifying equation \eqref{eq:modelbainchipairequations11-simple} with $2k-1>0$, and we apply Lemma \ref{lemma:basicdivergenceidentitybianchipairrpweightedestimate}  to the vectorfield
\beaa
X:=\frac 1 2 |q|^b|\Psi_{(1)}|^2e_3+|q|^b|\Psi_{(2)}|^2e_4-2|q|^b  \Re \Big(\Psi_{(1)} \c \ov{\Psi_{(2)}}\Big).
\eeaa
Using \eqref{eq:basicdivergenceidentitybianchipairrpweightedestimate-complex-1.simple-0}, we obtain
\bea\lab{eq:basicdivergenceidentitybianchipairrpweightedestimate-complex-1.simple}
\begin{split}
 \Div X
&=\frac 1 2 |q|^{b}  \La_{(1)} \trchb |\Psi_{(1)}|^2+|q|^{b}  \La_{(2)}\trch |\Psi_{(2)}|^2+  |q|^b(2k-1)\big( \omb |\Psi_{(1)}|^2- 2\om  |\Psi_{(2)}|^2\big)\\
&+O(ar^{b-2}) \c  \Re ( \Psi_{(1)} \c \ov{\Psi_{(2)}})+  |q|^b \Re\Big( F_{(1)} \c \ov{\Psi_{(1)}} \Big)+ 2|q|^b\Re\Big(  F_{(2)} \c \ov{\Psi_{(2)}}\Big)\\
&+r^{b}\Ga_b\Big(|\Psi_{(1)}|^2+|\Psi_{(2)}|^2\Big).
\end{split}
\eea
Let $b$ such that $\La_{(1)}=-2 c_{(1)} +1+\frac b 2 > 0$. Since $\trchb=-\frac{2r}{|q|^2}+\Ga_g $ and $\omb\in\Ga_b$, we can bound the above by 
\beaa
\begin{split}
\Div (X)&\leq - r|q|^{b-2}\Big(\La_{(1)}  -r\Ga_b \Big)|\Psi_{(1)}|^2+O(r^{b-1}) |\Psi_{(2)}|^2\\
& +  |q|^b \Re\Big( F_{(1)} \c \ov{\Psi_{(1)}} \Big)+ 2|q|^b\Re\Big(  F_{(2)} \c \ov{\Psi_{(2)}}\Big). 
\end{split}
\eeaa
We deduce, since $r\leq |q|\leq \sqrt{2}r$,
\beaa
\begin{split}
\Div (X)&\leq -\frac{\La_{(1)}}{2} r^{b-1} |\Psi_{(1)}|^2+O(r^{b-1}) |\Psi_{(2)}|^2 +  |q|^b \Re\Big( F_{(1)} \c \ov{\Psi_{(1)}} \Big)+ 2|q|^b\Re\Big(  F_{(2)} \c \ov{\Psi_{(2)}}\Big).
\end{split}
\eeaa
In view of Lemma \ref{lemma:basicdivergenceidentitybianchipairrpweightedestimate},       integrated 
 in the region $\MM(\tau_1, \tau_2)$,     we obtain
\beaa
-\int_{\pr^+\MM} \g(X, N) + \int_{\pr^-\MM} \g(X, N)&\les&- \int_{\MM } r^{b-1} |\Psi_{(1)}|^2+\int_{\MM} r^{b-1} |\Psi_{(2)}|^2\\
&& +\left|\int_{\MM}|q|^{b}\Re(F_{(1)}\c\ov{\Psi_{(1)}})\right|+\int_{\MM}r^b|F_{(2)}||\Psi_{(2)}|. 
\eeaa
Hence
\beaa
 \int_{\MM } r^{b-1} |\Psi_{(1)}|^2-\int_{\pr^+\MM} \g(X, N)&\les& - \int_{\pr^-\MM} \g(X, N)+\int_{\MM} r^{b-1} |\Psi_{(2)}|^2\\
&&+\left|\int_{\MM}|q|^{b}\Re(F_{(1)}\c\ov{\Psi_{(1)}})\right|+\int_{\MM}r^b|F_{(2)}||\Psi_{(2)}| .
\eeaa
For the boundary terms, we compute
\beaa
\g(X, N)&=&  |q|^b\left(\frac 1 2 |\Psi_{(1)}|^2 \g(e_3, N)+ |\Psi_{(2)}|^2 \g(e_4, N) - 2 \g(N, e_a)   \Re \big(\Psi_{(1)} \c \ov{\Psi_{(2)}}\big)_a\right),
\eeaa
and recall that $\pr^+\MM(\tau_1, \tau_2)=\AA(\tau_1, \tau_2)\cup \Si(\tau_2) \cup \Si_*(\tau_1, \tau_2)$ and $\pr^-\MM(\tau_1, \tau_2)=  \Si(\tau_1)$.  Based on    the assumptions made in sections \ref{section:BoundariesofMM-ThmM8} and \ref{section:propertiesoftau-M8}, we  have the following lemma.

\begin{lemma}
\lab{Lemma:g(X, N)onPrMM}
The following      inequalities hold on $\AA$, $\Si_*$ and $\Si(\tau)$:
\begin{itemize}
\item On $\AA$ we have
 \beaa
 \g(N_\AA, e_3)= - 1, \qquad   \g(N_\AA, e_4) \leq - \frac{1}{10}\de_\HH, \qquad \g(N_\AA, e_a) =O(\deh).
 \eeaa
 Therefore
 \beaa
 -\g(X, N) &=& -|q|^b\left(\frac 1 2 |\Psi_{(1)}|^2 \g(e_3, N)+ |\Psi_{(2)}|^2 \g(e_4, N) - 2 \g(N, e_a)   \Re \big(\Psi_{(1)} \c \ov{\Psi_{(2)}}\big)_a\right)\\
 &\ges &    |\Psi_{(1)}|^2-   O(\de_\HH)   |\Psi_{(2)}|^2.
 \eeaa
 
\item  On the boundary $\Si_*$ we have, with $N_*=N_{\Si_*} $,
\beaa
 \g(N_{*}, e_3)= - 1, \qquad   \g(N_{*}, e_4) \leq - 1, \qquad \g(N_{*}, e_a) =O(r^{-1}).
 \eeaa
 Therefore
 \beaa
-\g(X, N_*)&=&-  |q|^b\left(\frac 1 2 |\Psi_{(1)}|^2 \g(e_3, N_*) + |\Psi_{(2)}|^2 \g(e_4, N) - 2 \g(N_*, e_a)   \Re \big(\Psi_{(1)} \c \ov{\Psi_{(2)}}\big)_a\right)\\
 &\ges& r^{b}\big(|\Psi_{(1)}|^2+ |\Psi_{(2)}|^2 \big).
\eeaa

\item  On the boundary $\Si=\Si(\tau)$
\beaa
\g(N_{\Si}, N_{\Si}) \leq -\frac{1}{100}\frac{m^2}{r^2}, \quad e_4(\tau)\geq \frac{1}{100}\frac{m^2}{r^2}, \quad e_3(\tau)=1, \quad |\nab\tau|^2  \leq \frac{99}{100}e_4(\tau)e_3(\tau).
\eeaa
Therefore\footnote{Note that 
\beaa
2|\g(N, e_a)   \Re \big(\Psi_{(1)} \c \ov{\Psi_{(2)}}\big)_a| \leq \sqrt{2}|\nab\tau||\Psi_{(1)}||\Psi_{(2)}|\leq \frac{99}{100}\left(\frac 1 2 |\Psi_{(1)}|^2e_3(\tau) + |\Psi_{(2)}|^2e_4(\tau)\right).
\eeaa}
\beaa
-\g(X, N)&=&-  |q|^b\left(\frac 1 2 |\Psi_{(1)}|^2 \g(e_3, N)+ |\Psi_{(2)}|^2 \g(e_4, N) - 2 \g(N, e_a)   \Re \big(\Psi_{(1)} \c \ov{\Psi_{(2)}}\big)_a\right)\\
&\gtrsim& r^b|\Psi_{(1)}|^2 + r^{b-2}  |\Psi_{(2)}|^2.
\eeaa
\end{itemize}
\end{lemma}

Using Lemma \ref{Lemma:g(X, N)onPrMM} to control the boundary terms, we finally obtain
\beaa
\begin{split}
 \int_{\MM } r^{b-1} |\Psi_{(1)}|^2+\int_{\pr^+\MM} r^{b}|\Psi_{(1)}|^2&\les \int_{\MM} r^{b-1} |\Psi_{(2)}|^2   +\left|\int_{\MM}|q|^{b}\Re(F_{(1)}\c\ov{\Psi_{(1)}})\right|+\int_{\MM}r^b|F_{(2)}||\Psi_{(2)}| \\
 & + \int_{\pr^-\MM} \big( r^{b}|\Psi_{(1)}|^2+  r^{b-2}|\Psi_{(2)}|^2\big).
\end{split}
\eeaa

{\bf Case 2.}  To obtain the second part of the proposition, we  
consider $\Psi_{(1)}, \Psi_{(2)}$ verifying equation \eqref{eq:modelbainchipairequations12-simple} 
 with  $2k-1<0$ and apply  Lemma \ref{lemma:basicdivergenceidentitybianchipairrpweightedestimate}  to the vectorfield
\beaa
X:= |q|^b|\Psi_{(1)}|^2e_3+\frac 12|q|^b|\Psi_{(2)}|^2e_4+2|q|^b  \Re \Big(\Psi_{(2)} \c \ov{\Psi_{(1)}}\Big).
\eeaa
According to  identity  \eqref{eq:basicdivergenceidentitybianchipairrpweightedestimate-complex-2.simple-0}, we have
\beaa
\begin{split}
 \Div X
&=|q|^{b} \La_{(1)} \trchb |\Psi_{(1)}|^2+\frac 1 2 |q|^{b}  \La_{(2)}\trch |\Psi_{(2)}|^2+ |q|^b(2k-1)\big(  2\omb  |\Psi_{(1)}|^2- \om  |\Psi_{(2)}|^2\big)\\
&+O(ar^{b-2}) \c  \Re ( \Psi_{(1)} \c \ov{\Psi_{(2)}})+2|q|^b  \Re\Big( F_{(1)} \c \ov{\Psi_{(1)}} \Big)+|q|^b \Re\Big(  F_{(2)} \c \ov{\Psi_{(2)}}\Big)\\
&+r^{b} \Ga_b\Big(|\Psi_{(1)}|^2+|\Psi_{(2)}|^2\Big).
\end{split}
\eeaa
Since $\La_{(2)}< 0$,  $\trch=\frac{2r\De}{|q|^4} + \Ga_g$, $2k-1<0$ and $\om\gtrsim -\frac{m}{r^2}+\Ga_g$, the coefficient $C$ of    $|q|^b |\Psi_{(2)}|^2$ is given by\footnote{In particular, we use the fact  that $\De\geq 0$ on $r\geq r_+$, that $\De=O(\deh)$ on $\MM(r\leq r_+)$, and that $\om\gtrsim -\frac{m}{r^2}$ on $\MM$.}  
\beaa
C &=& \frac 1 2 \La_{(2)} \trch- (2k-1)\om= \frac{r\De}{|q|^4} \La_{(2)} +|2k-1|\om +\Ga_g\\
&\gtrsim& -\frac{|\La_{(2)}|}{r} -  \frac{m}{r^2}|2k-1|.
\eeaa
The desired estimate  follows, as in the first case, by integration  on $\MM(\tau_1, \tau_2)$, with  the boundary terms being  treated  in the same manner as before. This concludes the proof of Proposition \ref{Prop:Bainchi-pairsEstimates-integrated}.
\end{proof}


\section{Bianchi pairs for higher derivatives}
\lab{sec:bianchipairsforhigherorderderivatives}


The goal of this section is to  commute the equations
\eqref{eq:modelbainchipairequations11-simple}, \eqref{eq:modelbainchipairequations12-simple}
 with  higher  derivatives in $\nabc_3, \nabc_4, \nabc_\Rhat$, see section \ref{sec:conformallyinvaariantop:part3} for the definition of these conformally invariant operators,  and prove the following two propositions. 
\begin{proposition}
\lab{Prop:BianchPairs-tilde4-k}
The following   higher order Bianchi identities  hold true:    
 \begin{itemize}
\item If $\Psi_{(1)}, \Psi_{(2)} $ verify  the equations \eqref{eq:modelbainchipairequations11-simple}, 
   then the quantities
   \beaa
   \Psit_{(1, k)}=( \ov{q}\nabc_4)^k \nabc_\Rhat^2  \Psi_{(1)}, \qquad 
   \Psit_{(2, k)}= (\ov{q}\nabc_4)^k\nabc_\Rhat^2  \Psi_{(2)},   
   \eeaa
    verify the  equations
      \bea
     \lab{eq:modelbainchipairequations11-simple-tildek}    
\bsplit
\nabc_3 \Psit_{(1, k)}+ \left(c_{(1)} -\frac k 2 \right) \tr \Xb  \Psit_{(1, k)} &=
-\DDs_p\, \Psit_{(2, k)}+\Ft_{(1,k)},\\
\nabc_4\Psit_{(2, k)} +\left(c_{(2)} -\frac k 2 \right)\ov{\tr  X}  \Psit_{(2, k)}&= \DDd_p \Psit_{(1, k)} +\Ft_{(2, k)},
\end{split}
\eea
where
\bea
\lab{eq:Ft_(1,k)Ft_(2,k)}
 \bsplit
 \Ft_{(1,k)}&=\big( \ov{q} \nabc_4 \big)^k \nabc^2_\Rhat   F_{(1)} - 2 \om \big( \ov{q} \nabc_4\big)^k \nabc_3 \nabc_\Rhat \Psi_{(1)} \\
 &+O(r^{-1} )\dk^{\le k +1}\nabc_\Rhat   \Psi_{(2)} +O(ar^{-2} ) \dk^{\le k}\dkb\nabc_\Rhat\Psi_{(1)} \\
 &+ O(r^{-1})  \dk^{\le k+1}  \big(\Psi_{(1)}, \Psi_{(2)}\big) + \dk^{\le k+2 }  \Big( \Ga_b  \c ( \Psi_{(1)}, \Psi_{(2)} )\Big),\\
  \Ft_{(2,k)}&=\big( \ov{q} \nabc_4 \big)^k \nabc^2_\Rhat   F_{(2)} + 2\frac{\De}{|q|^2}  \om \big( \ov{q} \nabc_4\big)^k \nabc_3 \nabc_\Rhat \Psi_{(2)} \\
 &+O(r^{-1} )\dk^{\le k +1}\nabc_\Rhat \Psi_{(1)} +O(ar^{-2} ) \dk^{\le k}\dkb\nabc_\Rhat\Psi_{(2)} \\
 &+ O(r^{-1})  \dk^{\le k+1}  \big(\Psi_{(1)}, \Psi_{(2)}\big) + \dk^{\le k+2 }  \Big( \Ga_b  \c ( \Psi_{(1)}, \Psi_{(2)} )\Big).
 \end{split}
 \eea

 \item  If $\Psi_{(1)}, \Psi_{(2)} $ verify  the equations \eqref{eq:modelbainchipairequations12-simple}, 
   then the quantities
   \beaa
   \Psit_{(1, k)}=( q\nabc_4)^k \nabc_\Rhat^2  \Psi_{(1)}, \qquad 
   \Psit_{(2, k)}= (q\nabc_4)^k\nabc_\Rhat^2  \Psi_{(2)},   
   \eeaa
    verify the  equations
      \bea
     \lab{eq:modelbainchipairequations12-simple-tildek}    
\bsplit
\nabc_3 \Psit_{(1, k)}+ \left(c_{(1)} -\frac k 2 \right)\ov{ \tr \Xb}  \Psit_{(1, k)} &=-
\DDd_p\, \Psit_{(2, k)}+\Ft_{(1,k)},\\
\nabc_4\Psit_{(2, k)} +\left(c_{(2)} -\frac k 2 \right) \tr  X  \Psit_{(2, k)}&= \DDs_p \, \Psit_{(1, k)} +\Ft_{(2, k)},
\end{split}
\eea
with  $ \Ft_{(1,k)},  \Ft_{(2,k)}$ as in  \eqref{eq:Ft_(1,k)Ft_(2,k)}.
\end{itemize}
\end{proposition}

\vspace{0.05cm}

\begin{proposition}
\lab{Prop:BianchPairs-tilde3-k}
The following   higher order Bianchi identities  hold true:
    \begin{itemize}
\item
   If $\Psi_{(1)}, \Psi_{(2)} $ verify  the equations \eqref{eq:modelbainchipairequations11-simple}, 
   then the quantities
   \beaa
   \Psit_{(1, k)}=( \nabc_3)^k \nabc_\Rhat^2  \Psi_{(1)}, \qquad 
   \Psit_{(2, k)}= (\nabc_3)^k\nabc_\Rhat^2  \Psi_{(2)},   
   \eeaa
    verify the  equations
      \bea
     \lab{eq:modelbainchipairequations11-simple-tilde3k}    
\bsplit
\nabc_3 \Psit_{(1, k)}+ c_{(1)}  \tr \Xb  \Psit_{(1, k)} &=
\DDs_p\, \Psit_{(2, k)}+\Ft_{(1,k)},\\
\nabc_4\Psit_{(2, k)} +c_{(2)} \ov{\tr  X}  \Psit_{(2, k)}&= \DDd_p \Psit_{(1, k)} +\Ft_{(2, k)},
\end{split}
\eea
where
\bea
\lab{eq:Ft_(1,k)Ft_(2,k)3}
 \bsplit
 \Ft_{(1,k)}&=\big( \nabc_3 \big)^k \nabc^2_\Rhat   F_{(1)} - 2 \om \big( \nabc_3\big)^k \nabc_3 \nabc_\Rhat \Psi_{(1)} \\
 &+O(r^{-1} )\dk^{\le k+1}\nabc_\Rhat   \Psi_{(2)} +O(ar^{-2} ) \dk^{\le k}\dkb\nabc_\Rhat\Psi_{(1)} \\
 &+ O(r^{-1})  \dk^{\le k+1}  \big(\Psi_{(1)}, \Psi_{(2)}\big) + \dk^{\le k+2 }  \Big( \Ga_b  \c ( \Psi_{(1)}, \Psi_{(2)} )\Big),\\
  \Ft_{(2,k)}&=\big(  \nabc_3 \big)^k \nabc^2_\Rhat   F_{(2)} + 2\frac{\De}{|q|^2}  \om \big(  \nabc_3\big)^k \nabc_3 \nabc_\Rhat \Psi_{(2)} \\
 &+O(r^{-1} )\dk^{\le k+1}\nabc_\Rhat \Psi_{(1)}  +O(ar^{-2} )\dk^{\le k}\dkb\nabc_\Rhat\Psi_{(2)} \\
 &+ O(r^{-1})  \dk^{\le k+1}  \big(\Psi_{(1)}, \Psi_{(2)}\big) + \dk^{\le k+2 }\Big( \Ga_b  \c ( \Psi_{(1)}, \Psi_{(2)} )\Big).
 \end{split}
 \eea
 
 \item If $\Psi_{(1)}, \Psi_{(2)} $ verify  the equations \eqref{eq:modelbainchipairequations12-simple}, 
   then the quantities
   \beaa
   \Psit_{(1, k)}=( \nabc_{3})^k \nabc_\Rhat^2  \Psi_{(1)}, \qquad 
   \Psit_{(2, k)}= (\nabc_3)^k\nabc_\Rhat^2  \Psi_{(2)},   
   \eeaa
    verify the  equations
      \bea
     \lab{eq:modelbainchipairequations12-simple-tilde3k}    
\bsplit
\nabc_3 \Psit_{(1, k)}+ c_{(1)} \ov{ \tr \Xb}  \Psit_{(1, k)} &=-
\DDd_p\, \Psit_{(2, k)}+\Ft_{(1,k)},\\
\nabc_4\Psit_{(2, k)} +c_{(2)}  \tr  X  \Psit_{(2, k)}&= \DDs_p \, \Psit_{(1, k)} +\Ft_{(2, k)},
\end{split}
\eea
with  $ \Ft_{(1,k)},  \Ft_{(2,k)}$ as in  \eqref{eq:Ft_(1,k)Ft_(2,k)3}.

\end{itemize}
\end{proposition}

Propositions \ref{Prop:BianchPairs-tilde4-k} and \ref{Prop:BianchPairs-tilde3-k} will be proved in section \ref{sec:proofofProp:BianchPairs-tilde4and3-k} by relying on the commutation lemmas of sections  \ref{sec:commutationofBianchipairsiwithnabcRhat} and \ref{sec:commutationofBianchipairsiwithnabc3andnabc4}.


\subsection{Commutation of Bianchi pairs with $\nabc_{\Rhat}$}
\lab{sec:commutationofBianchipairsiwithnabcRhat}


In the following lemma, we commute the Bianchi identities with $\nabc_{\Rhat}$, see \eqref{eq:definitionofnabcRhatconf} for the definition of $\nabc_{\Rhat}$. 
\begin{lemma}
\lab{Lemma:commuteBianchi-dot}
  The  following identities hold true:
     \begin{itemize}
        \item   If $\Psi_{(1)}, \Psi_{(2)} $ verify  the equations \eqref{eq:modelbainchipairequations11-simple}, 
then  $\Psid_{(1)}=\nabc_\Rhat   \Psi_{(1)},   \Psid_{(2)}=\nabc_\Rhat \Psi_{(2)} $  verify the following 
 \bea
     \lab{eq:modelbainchipairequations11-simple-dot}    
\bsplit
\nabc_3 \Psid_{(1)}+ c_{(1)}  \tr \Xb  \Psid_{(1)} &=
-\DDs_p\, \Psid_{(2)}+\Fd_{(1)},\\
\nabc_4\Psid_{(2)} +c_{(2)}\ov{\tr  X}  \Psid_{(2)}&= \DDd_p \Psid_{(1)} +\Fd_{(2)},
\end{split}
\eea
 where
 \bea
 \lab{eq:Fd1Fd2}
 \bsplit
\Fd_{(1)}&= \nabc_{\Rhat} F_{(1)}  -\om \nabc_3 \Psi_{(1)}        + O(r^{-2} )\dk^{\le 1}  \Psi_{(2)} +O(ar^{-3}) \dkb\Psi_{(1)} +O(r^{-2} )\Psi_{(1)}\\
& + \Ga_b  \c \dk^{\leq 1} \Psi_{(2)}+ r^{-1} \Ga_b \cdot \dk^{\leq 1}  \Psi_{(1)},\\
  \Fd_{(2)}  &=\nabc_{\Rhat} F_{(2)}+\frac{\De}{|q|^2} \om \nabc_3 \Psi_{(2)} +O(r^{-2})\dk^{\le 1 }\Psi_{(1)} + O(ar^{-3} )\dkb^{\le 1 }   \Psi_{(2)}\\
  &+ O(r^{-2})  \Psi_{(2)} + \Ga_b  \c \dk^{\leq 1}  \Psi_{(1)} + r^{-1} \Ga_b \cdot  \Psi_{(2)}.
   \end{split}
\eea

\item  If $\Psi_{(1)}, \Psi_{(2)} $ verify  the equations \eqref{eq:modelbainchipairequations12-simple}, 
then  $\Psid_{(1)}=\nabc_\Rhat   \Psi_{(1)},   \Psid_{(2)}=\nabc_\Rhat \Psi_{(2)} $  verify the following 
\bea\lab{eq:modelbainchipairequations12-simple-dot}
\begin{split}
\nabc_3(\Psid_{(1)})+c_{(1)}\ov{\tr \Xb}\Psid_{(1)}&=- \DDd_p\, \Psid_{(2)} +\Fd_{(1)},\\[2mm]
\nabc_4(\Psid_{(2)})+c_{(2)}\tr X\Psid_{(2)} &=\DDs_p\, \Psid_{(1)}  +\Fd_{(2)},
\end{split}
\eea
with $\Fd_{(1)}, \Fd_{(2)} $ as in \eqref{eq:Fd1Fd2}.
\end{itemize}
\end{lemma}

\begin{proof}
We start with  the equation $\nabc_3 \Psi_{(1)}+ c_{(1)}  \tr \Xb  \Psi_{(1)} = -\DDs_p\, \Psi_{(2)}+F_{(1)}$ and write 
\beaa
\nabc_3 \Psid_{(1)}&=& \nabc_3 \nabc_\Rhat   \Psi_{(1)}=\nabc_\Rhat  \nabc_3   \Psi_{(1)}+[ \nabc_3, \nabc_\Rhat]  \Psi_{(1)}\\
&=& \nabc_\Rhat\Big(- c_{(1)}  \tr \Xb  \Psi_{(1)} - \DDs_p\, \Psi_{(2)}+F_{(1)} \Big) +[ \nabc_3, \nabc_\Rhat]  \Psi_{(1)}\\
&=& - c_{(1)}  \tr \Xb   \Psid_{(1)} +O(r^{-2} )\Psi_{(1)} - \DDs_p\, \Psid_{(2)}  - [\nab_\Rhat, \DDs_p\, ] \Psi_{(2)}  +[ \nabc_3, \nabc_\Rhat]  \Psi_{(1)}\\
&&+\nabc_{\Rhat} F_{(1)}.
\eeaa
Henceforth
\beaa
\nabc_3 \Psid_{(1)} +    c_{(1)}  \tr \Xb   \Psid_{(1)}   &=& -\DDs_p\, \Psid_{(2)}+\Fd_{(1)}
\eeaa
where
\beaa
\Fd_{(1)}&=&\nab_{\Rhat} F_{(1)} - [\nab_\Rhat, \DDs_p\, ] \Psi_{(2)}  +[ \nabc_3, \nabc_\Rhat]  \Psi_{(1)}+ O(r^{-2} )\Psi_{(1)}. 
\eeaa
In view of the commutation Lemma \ref{lemma:comm.nabc_Rhat}
\beaa
\,[ \nabc_3, \nabc_\Rhat]  \Psi_{(1)}&=& -\om \nabc_3 \Psi_{(1)}+ O(ar^{-3} )\dkb^{\le 1 }   \Psi_{(1)}+ O(r^{-3})  \Psi_{(1)} + r^{-1} \Ga_b \cdot  \dk^{\leq 1} \Psi_{(1)},\\
\,[\nab_\Rhat, \DDs_p\, ] \Psi_{(2)}&=&  \frac{\De}{2|q|^2}\tr\Xb
       \DDs_p \,  \Psi_{(2)} +O(ar^{-2})\dk^{\le 1} \Psi_{(2)} + \Ga_b  \c \dk^{\leq 1} \Psi_{(2)}.
\eeaa
Hence
\beaa
\Fd_{(1)}&=& \nabc_{\Rhat} F_{(1)}  -\om \nabc_3 \Psi_{(1)}        + O(r^{-2} )\dk^{\le 1}  \Psi_{(2)} +O(ar^{-3}) \dkb\Psi_{(1)} +O(r^{-2} )\Psi_{(1)}\\
&& + \Ga_b  \c \dk^{\leq 1} \Psi_{(2)}+ r^{-1} \Ga_b \cdot \dk^{\leq 1}  \Psi_{(1)}.
\eeaa
Similarly,  starting with the  second  equation 
 $\nabc_4\Psi_{(2)} +c_{(2)} \ov{\tr  X}  \Psi_{(2)}= \DDd_p \Psi_{(1)} +F_{(2)}$, we write
\beaa
\nabc_4 \Psid_{(2)}&=& \nabc_4 \nabc_\Rhat   \Psi_{(2)}=\nabc_\Rhat  \nabc_4   \Psi_{(2)}+[ \nabc_4, \nabc_\Rhat]  \Psi_{(2)}\\
&=& \nabc_\Rhat\Big(- c_{(2)}  \ov{\tr X}  \Psi_{(2)} + \DDd_p\, \Psi_{(1)}+F_{(2)} \Big) +[ \nabc_4, \nabc_\Rhat]  \Psi_{(2)}\\
&=& - c_{(2)}  \ov{\tr X}   \Psid_{(2)} +O(r^{-2} )\Psi_{(2)} +\DDd_p\, \Psid_{(1)}  +[\nab_\Rhat, \DDd_p\, ] \Psi_{(1)}  +[ \nabc_4, \nabc_\Rhat]  \Psi_{(2)}\\
&&+\nabc_{\Rhat} F_{(2)}.
\eeaa
Hence $\nabc_4\Psid_{(2)} +c_{(2)} \ov{\tr  X}  \Psid_{(2)}= \DDd_p \Psid_{(1)} +\Fd_{(2)}$
with
\beaa
\Fd_{(2)}&=&\nabc_{\Rhat} F_{(2)} +[\nab_\Rhat, \DDd_p\, ] \Psi_{(1)}  +[ \nabc_4, \nabc_\Rhat]  \Psi_{(2)} +O(r^{-2} )\Psi_{(2)} \\
&=&\nabc_{\Rhat} F_{(2)} +\frac{\De}{2|q|^2}\tr\Xb
      \DDd_p\Psi_{(1)}+  O(ar^{-2}) \dk^{\le 1}  \Psi_{(1)} + \Ga_b  \c \dk^{\leq 1}  \Psi_{(1)}\\
      &&+\frac{\De}{|q|^2} \om \nabc_3 \Psi_{(2)}+ O(ar^{-3} )\dkb^{\le 1 }   \Psi_{(2)}+ O(r^{-2})  \Psi_{(2)} + r^{-1} \Ga_b \cdot \dk^{\leq 1} \Psi_{(2)}\\
      &=& \nabc_{\Rhat} F_{(2)}+\frac{\De}{|q|^2} \om \nabc_3 \Psi_{(2)} +O(r^{-2})\dk^{\le 1 }\Psi_{(1)} + O(ar^{-3} )\dkb^{\le 1 }   \Psi_{(2)}+ O(r^{-2})  \Psi_{(2)}\\
      && + \Ga_b  \c \dk^{\leq 1}  \Psi_{(1)} + r^{-1} \Ga_b \cdot  \Psi_{(2)}.
\eeaa
This concludes the proof of Lemma \ref{Lemma:commuteBianchi-dot}.
\end{proof}

In the next lemma we commute  the   Bianchi pairs \eqref{eq:modelbainchipairequations11-simple}, \eqref{eq:modelbainchipairequations12-simple}   once more with  $\nabc_\Rhat  $.  

\begin{lemma}\lab{Lemma:commuteBianchi-ddot}
The  following identities hold true:
\begin{itemize}  
     \item    If $\Psi_{(1)}, \Psi_{(2)} $ verify  the equations \eqref{eq:modelbainchipairequations11-simple},
then  $\Psidd_{(1)}=\nabc^2_\Rhat   \Psi_{(1)},   \Psidd_{(2)}=\nabc^2_\Rhat \Psi_{(2)} $  verify the following 
 \bea
     \lab{eq:modelbainchipairequations11-simple-ddot}    
\bsplit
\nabc_3 \Psidd_{(1)}+ c_{(1)}  \tr \Xb  \Psidd_{(1)} &=
-\DDs_p\, \Psidd_{(2)}+\Fdd_{(1)},\\
\nabc_4\Psidd_{(2)} +c_{(2)} \ov{\tr  X}  \Psidd_{(2)}&= \DDd_p \Psidd_{(1)} +\Fdd_{(2)},
\end{split}
\eea
with
\bea
\lab{eq:Fdd1Fdd2}
\bsplit
\Fdd_{(1)}&=\nabc^2_\Rhat   F_{(1)}  - 2 \om \nabc_3 \Psid_{(1)}        + O(r^{-2} )\dk^{\le 1}  \Psid_{(2)} +O(ar^{-3}) \dkb\Psid_{(1)}\\
& + O(r^{-2})  \dk^{\le 1}\big(\Psi_{(1)}, \Psi_{(2)}\big)+  \dk^{\le 2}  \big( \Ga_b  \c \Psi_{(2)}\big)+ r^{-1}\dk^{\le 2} \big( \Ga_b \cdot   \Psi_{(1)} \big),\\
\Fdd_{(2)}&= \nabc^2_\Rhat   F_{(2)}+2 \frac{\De}{|q|^2} \om \nabc_3 \Psid_{(2)} +O(r^{-2}) \dk^{\le 1} \Psid_{(1)} +O(ar^{-3} )\dkb\Psid_{(2)}\\
  & + O(r^{-2})  \dk^{\le 1}\big(\Psi_{(1)}, \Psi_{(2)}\big) + \dk^{\le 2}  \big( \Ga_b  \c \Psi_{(1)}\big)+ r^{-1}\dk^{\le 2} \big( \Ga_b \cdot   \Psi_{(2)} \big).
\end{split}
\eea

\item  If $\Psi_{(1)}, \Psi_{(2)} $ verify  the equations \eqref{eq:modelbainchipairequations12-simple},
then  $\Psidd_{(1)}=\nabc^2_\Rhat   \Psi_{(1)},   \Psidd_{(2)}=\nabc^2_\Rhat \Psi_{(2)} $  verify the following 
\bea\lab{eq:modelbainchipairequations12-simple-ddot}
\begin{split}
\nabc_3(\Psidd_{(1)})+c_{(1)}\ov{\tr \Xb}\Psidd_{(1)}&=- \DDd_p\, \Psidd_{(2)} +\Fdd_{(1)},\\[2mm]
\nabc_4(\Psidd_{(2)})+c_{(2)}\tr X\Psidd_{(2)} &=\DDs_p\, \Psidd_{(1)}  +\Fdd_{(2)},
\end{split}
\eea
with $\Fdd_{(1)}, \Fdd_{(2)} $ as in \eqref{eq:Fdd1Fdd2}.
\end{itemize}
\end{lemma}

\begin{proof}
Starting with the first  Bianchi  pair  of  Lemma \ref{Lemma:commuteBianchi-dot}, we  commute once more 
with $\nabc_\Rhat$ and deduce
\beaa
\nabc_3(\Psidd_{(1)})+ c_{(1)}  \tr \Xb  \Psidd_{(1)} &=-
\DDs_p\,\, \Psidd_{(2)}+\Fdd_{(1)},\\
\nabc_4(\Psidd_{(2)})+c_{(2)}\ov{\tr X}\Psidd_{(2)} &=\DDd_p\, \Psidd_{(1)}  +\Fdd_{(2)},
\eeaa
 where
 \beaa
 \Fdd_{(1)}&=& \nabc_\Rhat\Fd_{(1)}  -\om \nabc_3 \Psid_{(1)}        + O(r^{-2} )\dk^{\le 1}  \Psid_{(2)} +O(ar^{-3})\dkb\Psid_{(1)} +O(r^{-2} )\Psid_{(1)}  \\
&&+ \Ga_b  \c \dk^{\leq 1} \Psid_{(2)}+ r^{-1} \Ga_b \cdot   \dk^{\leq 1}\Psid_{(1)}\\
&=&\nabc_\Rhat\Big(  \nabc_{\Rhat} F_{(1)}  -\om \nabc_3 \Psi_{(1)}        + O(r^{-2} )\dk^{\le 1}  \Psi_{(2)} +O(ar^{-3})\dkb\Psi_{(1)} +O(r^{-2})\Psi_{(1)}  \Big) \\
&& + \dk^{\le 1 }  \big( \Ga_b  \c \dk^{\leq 1} \Psi_{(2)}\big)+ r^{-1}\dk^{\le 1} \big( \Ga_b \cdot\dk^{\leq 1}\Psi_{(1)} \big)\\
&&-\om \nabc_3 \Psid_{(1)}        + O(r^{-2} )\dk^{\le 1}  \Psid_{(2)} +O(ar^{-3}) \dkb\Psid_{(1)} +O(r^{-2} )\Psid_{(1)} \\
&&+ \Ga_b  \c \dk^{\leq 1} \Psid_{(2)}+ r^{-1} \Ga_b \cdot  \dk^{\leq 1}\Psid_{(1)}\\
&=& \nabc^2_\Rhat   F_{(1)}  - 2 \om \nabc_3 \Psid_{(1)}        + O(r^{-2} )\dk^{\le 1}  \Psid_{(2)} +O(ar^{-3}) \dkb\Psid_{(1)}\\
&& + O(r^{-2})  \dk^{\le 1}\big(\Psi_{(1)}, \Psi_{(2)}\big)+  \dk^{\le 2}  \big( \Ga_b  \c \Psi_{(2)}\big) + r^{-1}\dk^{\le 2} \big( \Ga_b \cdot   \Psi_{(1)} \big).
 \eeaa
 Similarly
 \beaa 
  \Fdd_{(2)}  &=& \nabc^2_\Rhat   F_{(2)}+2 \frac{\De}{|q|^2} \om \nabc_3 \Psid_{(2)} +O(r^{-2})\dk^{\le 1} \Psid_{(1)} +O(ar^{-3} )\dkb\Psid_{(2)}\\
  && + O(r^{-2})  \dk^{\le 1}\big(\Psi_{(1)}, \Psi_{(2)}\big) + \dk^{\le 2}  \big( \Ga_b  \c \Psi_{(1)}\big)+ r^{-1}\dk^{\le 2} \big( \Ga_b \cdot   \Psi_{(2)} \big).
 \eeaa
 
 The proof of \eqref{eq:modelbainchipairequations12-simple-ddot} is similar and left to the reader.
\end{proof}


\subsection{Commutation of  Bianchi pairs with $\nabc_3$ and $\nabc_4$}
\lab{sec:commutationofBianchipairsiwithnabc3andnabc4}
    

     The following  two  lemmas concern commutation of the Bianchi pairs respectively with $\ov{q}\nabc_4$,  $q\nabc_4$, and $\nabc_3$.
     \begin{lemma}
     \lab{lemma:CommBianchi-int1}
       The  following identities hold true:
     \begin{itemize}
     \item  
      If $\Psi_{(1)}, \Psi_{(2)} $ verify  the equations \eqref{eq:modelbainchipairequations11-simple}, then   
     $\Psit_{(1)}=\ov{q} \nabc_4 \Psi_{(1)}, \Psit_{(2)}=\ov{q} \nabc_4 \Psi_{(2)}  $ verify the equations
     \bea
     \lab{eq:modelbainchipairequations11-simple-tilde}    
\bsplit
\nabc_3 \Psit_{(1)}+ \left(c_{(1)} -\frac 1 2 \right) \tr \Xb \, \Psit_{(1)} &=-
\DDs_p\, \Psit_{(2)}+\Ft_{(1)},\\
\nabc_4\Psit_{(2)} +\left(c_{(2)} -\frac 1 2 \right)\ov{ \tr  X } \, \Psit_{(2)}&= \DDd_p \Psit_{(1)} +\Ft_{(2)},
\end{split}
\eea
 where
 \bea
 \lab{eq:modelbainchipairequations11-simple-Ft}
 \bsplit
 \Ft_{(1)}&= \ov{q} \nabc_4 F_{(1)} +O(ar^{-2})\dkb\Psi_{(1)}+O(r^{-1}) \Psi_{(1)}   +O(r^{-1} )\dk^{\leq 1} \Psi_{(2)}   \\
 &   +  \Ga_b \c  \dk^{\leq 1}(\Psi_{(1)},  \Psi_{(2)}),\\
 \Ft_{(2)}&=  \ov{q} \nabc_4 F_{(2)} + O(r^{-1})\Psi_{(2)} +O(r^{-1})\dk^{\leq 1}\Psi_{(1)} +r^{-1} \Ga_b  \c \dk^{\leq 1}\big(\Psi_{(1)}, \Psi_{(2)}\big).
 \end{split}
 \eea
 
 \item     If $\Psi_{(1)}, \Psi_{(2)} $ verify  \eqref{eq:modelbainchipairequations12-simple}, then   $\Psit_{(1)}=q  \nabc_4\Psi_{(1)}, \Psit_{(2)}= q  \nabc_4\Psi_{(2)}  $ 
  verify the equations
     \bea
      \lab{eq:modelbainchipairequations12-simple-tilde}  
\bsplit
\nabc_3 \Psit_{(1)}+ \left(c_{(1)} -\frac 1 2 \right)\ov{ \tr \Xb}  \Psit_{(1)} &=
-\DDd_p\, \Psit_{(2)}+\Ft_{(1)},\\
\nabc_4\Psit_{(2)} +\left(c_{(2)} -\frac 1 2 \right)\tr   X \Psit_{(2)}&= \DDs_p \, \Psit_{(1)} +\Ft_{(2)},
\end{split}
\eea
with  $\Ft_{(1)}, \Ft_{(2)}$ defined from $F_{(1)}, F_{(2)}$ as in \eqref{eq:modelbainchipairequations11-simple-Ft}.  
     \end{itemize}    
     \end{lemma}
    
     \begin{proof}
  We start with the first equation  in \eqref{eq:modelbainchipairequations11-simple-tilde}. Using  Lemma \ref{Lemma:commutationwith-qnabc_4}, we have 
      \beaa
 \bsplit
 \nabc_3 \Psit_{(1)}=& \nabc_{\ov{q}e_4}\nabc_3 \Psi_{(1)} +[\nabc_3, \nabc_{ \ov{q} e_4}]\Psi_{(1)}\\
 =& \ov{q}\nabc_4\Big( - c_{(1) }\tr \Xb\Psi_{(1)}  -\DDs_p\, \Psi_{(2)}  +F_{(1)}\Big)\\
 & + \frac {1}{ 2} \tr \Xb \nabc_{  \ov{q} e_4}\Psi_{(1)}  +  O(ar^{-1} )\nab^{\le 1 }  \Psi_{(1)}+ O(r^{-2}) \Psi_{(1)} + \Ga_b \cdot \dk^{\le 1} \Psi_{(1)}\\
=& -\left(c_{(1)} -\frac 1 2 \right) \tr \Xb  \Psit_{(1)}  -\DDs_p\, \Psit_{(2)} +O(ar^{-2})\dkb\Psi_{(1)}+O(r^{-1}) \Psi_{(1)}   +O(r^{-1} )\dk^{\leq 1} \Psi_{(2)}       \\
 &+\ov{q} \nabc_4 F_{(1)}+  \Ga_b \c  \dk^{\leq 1}(\Psi_{(1)},  \Psi_{(2)}),
 \end{split}
 \eeaa
 and hence 
 \beaa
  \Ft_{(1)} = \ov{q} \nabc_4 F_{(1)} +O(ar^{-2})\dkb\Psi_{(1)}+O(r^{-1}) \Psi_{(1)}   +O(r^{-1} )\dk^{\leq 1} \Psi_{(2)}      +  \Ga_b \c  \dk^{\leq 1}(\Psi_{(1)},  \Psi_{(2)})
 \eeaa
 as stated.
 
 Also, using again Lemma \ref{Lemma:commutationwith-qnabc_4}, we have
  \beaa
 \bsplit
 \nabc_4\Psit_{(2)}=& \nabc_{\ov{q}e_4}\nabc_4 \Psi_{(2)} +[\nabc_4, \nabc_{ \ov{q} e_4}]\Psi_{(2)}\\
 =& \nabc_{\ov{q}e_4}\Big( - c_{(2) }\ov{\tr X}\Psi_{(2)}  +\DDd_p\, \Psi_{(1)}  +F_{(2)}\Big) +\frac 1 2 \ov{\tr X} \nabc_{\ov{q} e_4}\Psi_{(2)} +r^{-1} \Ga_b  \c \dk^{\leq 1} \Psi_{(2)}\\
 =& - c_{(2) }\ov{\tr X}\Psit_{(2)}+O(r^{-1})\Psi_{(2)} +\DDd_p\, \Psit_{(1)} +[\nab_{\ov{q}e_4}, \DDd_p]\, \Psi_{(1)}+
 \nabc_{\ov{q}e_4}(F_{(2)})\\
 & +\frac 1 2 \ov{\tr X} \nabc_{\ov{q} e_4}\Psi_{(2)} +r^{-1} \Ga_b  \c \dk^{\leq 1} \Psi_{(2)}\\
 =& - \left(c_{(2) }-\frac{1}{2}\right)\ov{\tr X}\Psit_{(2)}+O(r^{-1})\Psi_{(2)} +\DDd_p\, \Psit_{(1)} \\
 &+O(r^{-1})\dk^{\leq 1}\Psi_{(1)}+ \nabc_{\ov{q}e_4}(F_{(2)}) +r^{-1} \Ga_b  \c \dk^{\leq 1}\big(\Psi_{(1)}, \Psi_{(2)}\big)
 \end{split}
 \eeaa
 and hence
\beaa
  \Ft_{(2)} &=&  \ov{q} \nabc_4 F_{(2)}+ O(r^{-1})\Psi_{(2)} +O(r^{-1})\dk^{\leq 1}\Psi_{(1)} +r^{-1} \Ga_b  \c \dk^{\leq 1}\big(\Psi_{(1)}, \Psi_{(2)}\big)
\eeaa 
as stated.
 
 \eqref{eq:modelbainchipairequations12-simple-tilde}     is derived in the same manner.  This concludes the proof of Lemma \ref{lemma:CommBianchi-int1}.
     \end{proof}

     \begin{lemma}
       \lab{lemma:CommBianchi-int2}
      The  following identities hold true:
     \begin{itemize}
     \item  
      If $\Psi_{(1)}, \Psi_{(2)} $ verify  the equations \eqref{eq:modelbainchipairequations11-simple} then   
     $\Psit_{(1)}= \nabc_3 \Psi_{(1)}, \Psit_{(2)}= \nabc_3 \Psi_{(2)}  $ verify the equations
     \bea
     \lab{eq:modelbainchipairequations11-simple-tilde-3}    
\bsplit
\nabc_3 \Psit_{(1)}+ c_{(1)}  \tr \Xb  \Psit_{(1)} &=
-\DDs_p\, \Psit_{(2)}+\Ft_{(1)},\\
\nabc_4\Psit_{(2)} +c_{(1)} \ov{\tr X}  \Psit_{(2)}&= \DDd_p \Psit_{(1)} +\Ft_{(2)},
\end{split}
\eea
 where
 \bea
 \lab{eq:modelbainchipairequations11-simple-Ft-3}
 \bsplit
 \Ft_{(1)}&= \nabc_3 F_{(1)} +O(r^{-1} )\dk^{\le 1 } \Psi_{(2)}    +O(r^{-2} )\Psi_{(1)} + \Ga_b\c  \dk^{\leq 1}(\Psi_{(1)},  \Psi_{(2)}),\\
 \Ft_{(2)}&=  \nabc_3 F_{(2)}+O(r^{-1} )\dk^{\le 1}  \Psi_{(1)}    +O(ar^{-2} )\dkb^{\le 1}\Psi_{(2)} +O(r^{-2} )\Psi_{(2)}\\
 &+  \Ga_b \c  \dk^{\leq 1}(\Psi_{(1)},  \Psi_{(2)}).
 \end{split}
 \eea
 
 \item     If $\Psi_{(1)}, \Psi_{(2)} $ verify  \eqref{eq:modelbainchipairequations12-simple} then   $\Psit_{(1)}=  \nabc_3\Psi_{(1)}, \Psit_{(2)}=   \nabc_3\Psi_{(2)}  $ 
  verify the equations
     \bea
      \lab{eq:modelbainchipairequations12-simple-tilde-3}  
\bsplit
\nabc_3 \Psit_{(1)}+ c_{(1)} \ov{ \tr \Xb}  \Psit_{(1)} &=-\DDd_p\, \Psit_{(2)}+\Ft_{(1)},\\
\nabc_4\Psit_{(2)} +c_{(2)}  \tr  X  \Psit_{(2)}&= \DDs_p\, \Psit_{(1)} +\Ft_{(2)},
\end{split}
\eea
with  $\Ft_{(1)}, \Ft_{(2)}$ defined from $F_{(1)}, F_{(2)}$ as in \eqref{eq:modelbainchipairequations11-simple-Ft-3}.  
     \end{itemize}    
     \end{lemma}  
     
     \begin{proof}
     The proof relies on the commutators in Lemma \ref{commutator-nab-c-3-DD-c-hot-ThmM8}. It is similar to that of  Lemma \ref{lemma:CommBianchi-int1} and is in fact simpler.
     \end{proof}

 
 \subsection{Proof of Propositions \ref{Prop:BianchPairs-tilde4-k} and \ref{Prop:BianchPairs-tilde3-k}}
 \lab{sec:proofofProp:BianchPairs-tilde4and3-k}


The proof of Propositions \ref{Prop:BianchPairs-tilde4-k} and \ref{Prop:BianchPairs-tilde3-k} are similar, so we focus on the one of Proposition \ref{Prop:BianchPairs-tilde4-k} and consider first  the  case $k=1$.  

To  check  $k=1$ in the case  of the Bianchi pair \eqref{eq:modelbainchipairequations11-simple}, we 
  apply the result of Lemma \ref{lemma:CommBianchi-int1}
   to the   Bianchi pair  derived  in  Lemma
   \ref{Lemma:commuteBianchi-ddot}.  More precisely we start
    with  the equation \eqref{eq:modelbainchipairequations11-simple-ddot} 
 \beaa 
\bsplit
\nabc_3 \Psidd_{(1)}+ c_{(1)}  \tr \Xb  \Psidd_{(1)} &=-
\DDs_p\, \Psidd_{(2)}+\Fdd_{(1)},\\
\nabc_4\Psidd_{(2)} +c_{(2)} \tr  X  \Psidd_{(2)}&= \DDd_p \Psidd_{(1)} +\Fdd_{(2)},
\end{split}
\eeaa
with $\Fdd_{(1)}, \Fdd_{(2)}$ given by  \eqref{eq:Fdd1Fdd2}
and apply  to it the result of  the first part of Lemma \ref{lemma:CommBianchi-int1}  to deduce
 \bea
     \lab{eq:modelbainchipairequations11-simple-tilde-1}    
\bsplit
\nabc_3 \Psit_{(1, 1)}+ \left(c_{(1)} -\frac 1 2 \right) \tr \Xb  \Psit_{(1, 1)} &=
-\DDs_p\, \Psit_{(2, 1)}+\Ft_{(1,1)},\\
\nabc_4\Psit_{(2, 1)} +\left(c_{(2)} -\frac 1 2 \right) \tr  X  \Psit_{(2, 1)}&= \DDd_p \Psit_{(1, 1)} +\Ft_{(2, 1)},
\end{split}
\eea
where
\beaa
\Ft_{(1,1)}&=&\ov{q} \nabc_4 \Fdd_{(1)} +O(ar^{-2} ) \dkb\Psidd_{(1)}+O(r^{-1} )\Psidd_{(1)}+O(r^{-1} ) \dk^{\le 1 } \Psidd_{(2)} + \Ga_b \c  \dk^{\leq 1}(\Psidd_{(1)},  \Psidd_{(2)}),\\
\Ft_{(1,2)}&=&\ov{q} \nabc_4 \Fdd_{(2)} +O(r^{-1} )\Psidd_{(2)}+O(r^{-1} ) \dk^{\le 1 } \Psidd_{(1)} + r^{-1} \Ga_b \c  \dk^{\leq 1}(\Psidd_{(1)},  \Psidd_{(2)}).
\eeaa
Thus, in view of formulas \eqref{eq:Fdd1Fdd2}  for  $ \Fdd_{(1)}, \Fdd_{(2)}$   and the definition of $\Psid_{(1)}, \Psid_{(2)}$         we deduce
\beaa
\Ft_{(1,1)}&=&\ov{q} \nabc_4 \Big\{\nabc^2_\Rhat   F_{(1)}  - 2 \om\nabc_3 \Psid_{(1)}        + O(r^{-2} )\dk^{\le 1}  \Psid_{(2)} +O(ar^{-3})\dkb\Psid_{(1)}\Big\}\\
&&+\ov{q} \nabc_4 \Big\{ O(r^{-2})  \dk^{\le 1}\big(\Psi_{(1)}, \Psi_{(2)}\big)+  \dk^{\le 2}  \big( \Ga_b  \c \Psi_{(2)}\big)+ r^{-1}\dk^{\le 2} \big( \Ga_b \cdot   \Psi_{(1)} \big)\Big\}\\
&&+O(ar^{-2} ) \dkb\Psidd_{(1)}+O(r^{-1} )\Psidd_{(1)}+O(r^{-1} ) \dk^{\le 1 } \Psidd_{(2)} + \Ga_b \c  \dk^{\leq 1}(\Psidd_{(1)},  \Psidd_{(2)})\\
&=& \ov{q} \nabc_4 \nabc^2_\Rhat   F_{(1)} - 2 \om  \ov{q} \nabc_4 \nabc_3 \nabc_\Rhat \Psi_{(1)} +O(r^{-2}) \dk^{\le 2}\nabc_\Rhat \Psi_{(2)}\\
&&+O(ar^{-3})\dk\dkb\nabc_\Rhat \Psi_{(1)}+O(r^{-2})\dk^{\le 2}\big(\Psi_{(1)}, \Psi_{(2)}\big) +\dk^{\le 3}  \big( \Ga_b  \c \Psi_{(2)}\big)\\
&&+ r^{-1}\dk^{\le 3} \big( \Ga_b \cdot   \Psi_{(1)} \big) +O(ar^{-2} ) \dkb\nabc_{\Rhat}^2\Psi_{(1)}+O(r^{-1} )\nabc_{\Rhat}^2\Psi_{(1)}\\
&&+O(r^{-1} ) \dk^{\le 1 }\nabc_{\Rhat}^2\Psi_{(2)} + \Ga_b \c  \dk^{\leq 1}(\nabc_{\Rhat}^2\Psi_{(1)},  \nabc_{\Rhat}^2\Psi_{(2)}).
\eeaa
Thus, in  simplified  form, 
 \beaa
 \Ft_{(1,1)}&=& \ov{q} \nabc_4 \nabc^2_\Rhat   F_{(1)} - 2 \om  \ov{q} \nabc_4 \nabc_3 \nabc_\Rhat \Psi_{(1)} +O(r^{-1} )\dk^{\le 2}\nabc_\Rhat\Psi_{(2)}\\
 &&+O(ar^{-2})\dk\dkb\nabc_\Rhat \Psi_{(1)}+ O(r^{-1})  \dk^{\le 2}  \big(\Psi_{(1)}, \Psi_{(2)}\big) + \dk^{\le 3 }  \Big( \Ga_b  \c ( \Psi_{(1)}, \Psi_{(2)} )\Big)
 \eeaa
 which corresponds to the first equation of \eqref{eq:Ft_(1,k)Ft_(2,k)} in the case $k=1$.
 
 In the same fashion we find
 \beaa
  \Ft_{(2,1)}&=&  \ov{q} \nabc_4 \nabc^2_\Rhat   F_{(2)} + 2\frac{\De}{|q|^2}  \om  \ov{q} \nabc_4 \nabc_3 \nabc_\Rhat \Psi_{(2)} +O(r^{-1} )\dk^{\le 2}\nabc_\Rhat\Psi_{(1)}\\
 &&+O(ar^{-2})\dk\dkb\nabc_\Rhat \Psi_{(2)}+ O(r^{-1})  \dk^{\le 2}  \big(\Psi_{(1)}, \Psi_{(2)}\big) + \dk^{\le 3 }  \Big( \Ga_b  \c ( \Psi_{(1)}, \Psi_{(2)} )\Big)
 \eeaa
 which corresponds to the second equation of \eqref{eq:Ft_(1,k)Ft_(2,k)} in the case $k=1$. The general case, for all $k$, can be easily derived in the same manner  by induction on $k$, hence concluding the proof of Proposition \ref{Prop:BianchPairs-tilde4-k}. The proof of Proposition \ref{Prop:BianchPairs-tilde3-k} is similar.


    \section{Estimates for  $B$ and $\Bb$}
    \lab{section:InteriorEstimatesBBb}
  

The goal of this section is to prove Proposition \ref{proposition:EstimatesBBb-interior} providing energy-Morawetz estimates for $(B, \Bb)$ assuming corresponding energy-Morawetz estimates for $\Pc$.


\subsection{Estimates for $\nab_3 B$, $\nab B$, $\nab_4\Bb$, $\nab\Bb$}
\label{section:est-nab3-nab-B}


In this section, we obtain the following lemma. 
\begin{lemma}
\lab{Lemma:BianchB-steps1-2} 
Recall the notation $\de_{J+1}[\Pc]=\BEF^J_\de[r^2 \Pc]$. The following hold true in $\MM$:
\begin{enumerate}
\item We have
\bea
\lab{eq:BianchB-steps1-2-B}
\bsplit
\BEF^{J-1}_\de[ r^2( \nab_3B, r\nab B) ]\les&   \de_{J+1}[\Pc] +\Sk_J \Sk_{J+1} +\Rk_{J} \Rk_{J+1}  +\ep_0^2\\
& +O(a^2, \ep^2)\BEF_\de[  r^{J+1}  \nab^J_4 B]. 
 \end{split}
\eea

\item We have 
\bea
\lab{eq:BianchB-steps1-2-Bb}
\bsplit
\BEF^{J-1}_\de[ r (\nab_4\Bb, \nab \Bb) ] \les& 
 \de_{J+1}[\Pc]  +\Sk_J \Sk_{J+1} +\Rk_{J} \Rk_{J+1} +\ep_0^2\\
 &  +O(a^2, \ep^2) \BEF_\de[\nab^J_3 \Bb].
 \end{split}
\eea
\end{enumerate}
\end{lemma} 

\begin{proof}
We rely on the standard  linearization of the second Bianchi pair which will be stated in \eqref{eq:second-pair-B-P-lin}, i.e.
\bea
\lab{eq:BianchiforB-linear}
\bsplit
\nabc_3B +  \tr\Xb B &= \DD\ov{\Pc}+O(r^{-2}) \Pc +O(r^{-4} ) \Ga_b+ r^{-2} \Ga_b\c \Rc_b,\\
\nabc_4\Pc+\frac 3 2 \tr X \Pc&= \frac{1}{2}\DD\c \ov{B} +O(r^{-2} ) B +O(r^{-4} ) \Ga_b+ r^{-2}\Ga_b\c \Rc_b.
\end{split}
\eea
From the first equation we deduce
\beaa
\BEF^{J-1}_\de[ r^2 \nabc_3B] &\les & \BEF^{J-1}_\de[ r^2 \DD \Pc]  +\BEF^{J-1}_\de[  r B]  + \BEF^{J-1}_\de[  r^{- 2 }  \Ga_b]   +\ep_0^2\\
&\les&  \BEF^{J}_\de[ r\Pc] +\BEF^{J-1}_\de[r^{-1}\Rc_b]+\Sk_{J} \Sk_{J+1}  +\ep_0^2 \\
&\les& \de_{J+1}[\Pc]  +\Sk_{J} \Sk_{J+1} +\Rk_{J} \Rk_{J+1} +\ep_0^2,
\eeaa
where we used Lemmas \ref{Lemma:comparisonofnormsBEF-Rk} and \ref{Lemma:comparisonofnormsBEF[Gac)-Sk} to control $\BEF^{J-1}_\de[r^{-1}\Rc_b]$ and $\BEF^{J-1}_\de[  r^{-2}  \Ga_b]$. Similarly, from  the second equation,
\beaa
\BEF^{J}_\de[ r^3\DDov\c  B] &\les& \de_{J+1}[\Pc] +\Sk_J \Sk_{J+1} +\Rk_{J} \Rk_{J+1} +\ep_0^2.
\eeaa
Thus,
\bea
\lab{eq:EstmatesforB-1}
\bsplit
\BEF^{J-1}_\de[ r^2 \nabc_3B]+ \BEF^{J-1}_\de[ r^3 \DDov\c B  ] &\les \de_{J+1}[\Pc] +\Sk_J \Sk_{J+1} +\Rk_{J} \Rk_{J+1} +\ep_0^2.
\end{split}
\eea

\begin{remark}
\lab{remark:EstimatesforB-1-strong}
Note that  the terms  $\Sk_{J} \Sk_{J+1}$ and $\Rk_{J} \Rk_{J+1}$ on the RHS of \eqref{eq:EstmatesforB-1} are only needed to estimate  the energy flux  terms $\EF^{J-1}_\de[r^{-1}\Ga_b]$ and $\EF^{J-1}_\de[\Rc_b]$ thanks to Lemmas \ref{Lemma:comparisonofnormsBEF-Rk} and \ref{Lemma:comparisonofnormsBEF[Gac)-Sk}. In particular, in view of the control of the norms $B^{J-1}_\de[r^{-1}\Ga_b]$ and $B^{J-1}_\de[\Rc_b]$ provided by Lemmas \ref{Lemma:comparisonofnormsBEF-Rk} and \ref{Lemma:comparisonofnormsBEF[Gac)-Sk}, we infer the following stronger analog of \eqref{eq:EstmatesforB-1} for the flux
\bea
B^{J-1}_\de[ r^2 \nabc_3B]+ B^{J-1}_\de[r^3\DDov\c B] &\les &  \de_{J+1}[\Pc]  +\ep_J^2 +\ep_0^2.
\eea
\end{remark}

Next, we apply  Corollary \ref{Cor:HodgeThmM8} to derive
\beaa
 \BEF^{J-1}_\de[ r^3  \nab B  ] &\les& \BEF^{J-1}_\de[ r^3  \DDov\c B  ] +  \BEF^{J-1}_\de[r^2B] +O(a^2, \ep^2) \BEF^{J-1}_\de[  r^2(\nab_3B, \nab_4 B)]\\
 &\les& \BEF^{J-1}_\de[ r^3  \DDov\c B  ] +\Rk_J\Rk_{J+1}+O(a^2, \ep^2) \BEF^{J-1}_\de[   r^2 (\nab_3B, \nab_4 B)].
 \eeaa
 We deduce
 \beaa
 \bsplit
 \BEF^{J-1}_\de[ r^2 \nab_3B]+ \BEF^{J-1}_\de[ r^3  \nab B  ] \les&  \de_{J+1}[\Pc] +\Sk_J \Sk_{J+1} +\Rk_{J} \Rk_{J+1}  +\ep_0^2\\
 & +O(a^2, \ep^2) \BEF^{J-1}_\de[  r^2 (\nab_3B, \nab_4 B)].
 \end{split}
 \eeaa
 Thus, for small $a$ and $\ep$,
 \beaa
 \BEF^{J-1}_\de[ r^2 \nab_3B]+ \BEF^{J-1}_\de[ r^3  \nab B  ]&\les&  \de_{J+1}[\Pc] +\Sk_J \Sk_{J+1} +\Rk_{J} \Rk_{J+1}  +\ep_0^2\\
&& +O(a^2, \ep^2) \BEF^{J-1}_\de[  r^2 \nab_4B].
 \eeaa
Since
\beaa
\BEF^{J-1}_\de[  r^2 \nab_4 B] &\les& \BEF_\de[  r^2  (r\nab_4)^{J-1} \nab_4 B] + \BEF^{J-1} _\de[r (\nab_3 B, \nab B)]+\BEF^{J-2} _\de[rB]\\
&\les&  \BEF_\de[ r^{J+1}\nab_4^JB] + \BEF^{J-1} _\de[r (\nab_3 B, \nab B)]+\ep_J^2,
\eeaa
we infer, for small $a$ and $\ep$,
 \beaa
 \BEF^{J-1}_\de[ r^2(\nab_3B, r \nab  B)  ]\les \de_{J+1}[\Pc] +\Sk_J \Sk_{J+1} +\Rk_{J} \Rk_{J+1}  +\ep_0^2 +O(a^2, \ep^2)\BEF_\de[  r^{J+1}  \nab^J_4 B]
 \eeaa
 as stated in \eqref{eq:BianchB-steps1-2-B}. Also, in view of Remark \ref{remark:EstimatesforB-1-strong}, we have the following stronger analog of \eqref{eq:BianchB-steps1-2-B} for the flux
\beaa
B^{J-1}_\de[ r^2(\nab_3B, r\nab B)] &\les &  \de_{J+1}[\Pc]  +\ep_J^2 +\ep_0^2 +O(a^2, \ep^2)B_\de[  r^{J+1}  \nab^J_4 B],
\eeaa
and hence, since $\ep^2B_\de[  r^{J+1}  \nab^J_4 B]\les \ep^2B^J_\de[rB]\les \ep^2\Rk_{J+1}^2\les\ep_0^2$ in view of Lemma \ref{Lemma:comparisonofnormsBEF-Rk}, we obtain
\bea
\lab{eq:EstmatesforB-1-strong}
B^{J-1}_\de[ r^2(\nab_3B, r\nab B)] &\les &  \de_{J+1}[\Pc]  +\ep_J^2 +\ep_0^2 +O(a^2)B_\de[  r^{J+1}  \nab^J_4 B].
\eea

The estimate \eqref{eq:BianchB-steps1-2-Bb}  follows in the same fashion from the standard  linearization of the third Bianchi pair which will be stated in \eqref{eq:third-pair-P-Bb-lin}, i.e.
 \beaa
\begin{split}
\nabc_3\Pc+ \frac{3}{2}\ov{\tr\Xb} \Pc &= - \frac{1}{2}\DDb \c\Bb +  O(ar^{-2} )\Bb  +O(r^{-3})\Ga_b  +r^{-1} \Ga_b\c \Rc_b, \\
\nabc_4\Bb+ \tr X\Bb  &=-\DD \Pc  +O(ar^{-2}) \Pc   +O(r^{-3})\Ga_b +r^{-1} \Ga_b\c \Rc_b.
\end{split}
\eeaa
This concludes the proof of Lemma \ref{Lemma:BianchB-steps1-2}. 
\end{proof} 

\begin{remark}
Note that we have the following analog of \eqref{eq:EstmatesforB-1-strong}
\bea
\lab{eq:EstmatesforBb-1-strong}
B^{J-1}_\de[ r(\nab_4\Bb, \nab\Bb)] &\les &  \de_{J+1}[\Pc]  +\ep_J^2 +\ep_0^2 +O(a^2)B_\de[\nab^J_3\Bb].
\eea
\end{remark}

\begin{remark}
\lab{remark:BtPtpBbtPtm}
 It only remains to provide estimates for the top  $\nab_4 $ derivatives of $B$ and top $\nab_3$ derivatives of $\Bb$.  To this end, we need to include  derivatives with respect to $\Rhat$ which 
    leads us to introduce  the following quantities\footnote{Note  that  the signatures of   $\Bt, \Ptp, \Ptm, \Bbt$ are, respectively, $ J-1$, $J-2$, $-J+2$, and $- J+1$. Also, recall that $\Bdot=\Lieb_\T B$, $\Pdot=\T(P)$ and $\Bbdot=\Lieb_\T\Bb$.}
   \bea
   \bsplit
  &   \Bt:=(\ov{q}\nabc_4 )^{\le J-2}  \nabc^2_\Rhat\Bdot,  \qquad       \Ptp:=(\ov{q}\nabc_4)^{\le J-2} \nabc^2_{\Rhat} \ov{\Pdot}, \\
  & \Ptm:=(\nabc_3)^{\le J-2} \nabc^2_{\Rhat}\Pdot,\qquad\,\,\, \Bbt:=(\nabc_3 )^{\le J-2} 
   \nabc^2_\Rhat\Bbdot.
   \end{split}
   \eea  
  \end{remark}


\subsection{Bianchi equations for the quantities  $\protect \Bt, \Ptp,\Ptm,\Bbt$}


The following lemma provides Bianchi equations for the quantities  $\Bt, \Ptp,\Ptm,\Bbt$.
\begin{lemma}
   \lab{lemma:EstimatesBtPtp}
   The   following equations hold true for the quantities $\Bt, \Ptp,\Ptm,\Bbt$ introduced in Remark \ref{remark:BtPtpBbtPtm}:
   \begin{enumerate}
   \item   The quantities  $ \Bt, \Ptp$ verify the Bianchi  pair
    \bea\lab{eq:linearization-tildes-B-P}
  \bsplit
   \nabc_3  \Bt +  \left(1-\frac{J-2}{2}\right)\tr\Xb \Bt &= -\DDs_1 \Ptp +\Ft_{(1, J-2)},\\
    \nabc_4 \Ptp       +\left(\frac 3 2 -\frac{J-2}{2}\right) \ov{\tr X} \Ptp &=  \DDd_1     \Bt +\Ft_{(2, J-2)},
   \end{split}
  \eea
with
\bea
\lab{eq:Ft-Bt-Ptp}
\bsplit
 \Ft_{(1, J-2)}&=   O(r^{-1} )  \dk^{J-2} \nab_\Rhat \big(\nabc_3  \Bdot, \dkb \Bdot\big)  + O(r^{-1})  \dk^{J-1} \nab_\Rhat \Pdot \\
 &+  O(r^{-3} )(\nab_3, \nab)\dk^{\le J}\Ga_b  + O(r^{-1})  \dk^{\le J-1}\Pdot+ O(r^{-1})  \dk^{\le J}B\\
 & + O(r^{-3} )\dk^{\le J}  \Ga_b  + r^{-2}   \dk^{\le J+1}(\Ga_b\c \Rc_b) + r^{-2} \dk^{\le J} \big( \Ga_b\c\Ga_b\big),\\
   \Ft_{(2, J-2)}&=   O(r^{-1})
   \dk^{J-1} \nab_\Rhat (\Pdot, \Bdot)+ O(r^{-3} )   \dk^{\le J+1}   \Ga_b + O(r^{-1})  \dk^{\le J-1}(\Pdot, \Bdot) \\
   &+ r^{-2}   \dk^{\le J+1}(\Ga_b\c \Rc_b).
\end{split}
\eea

\item The quantities  $ \Bbt, \Ptm$ verify the Bianchi  pair
   \bea\lab{eq:linearization-tildes-Bb-P}
  \begin{split}
\nabc_3  \Ptm+      \frac{3}{2}\ov{\tr\Xb}  \Ptm  &=- \DDd_1 \Bbt+\Ft_{1, J-2},  \\\
  \nabc_4 \Bbt +\tr X   \Bbt  &=\DDs_1 \Ptm +\Ft_{2, J-2},
  \end{split}
\eea
with
\bea
\lab{eq:Ft-Bbt-Ptp}
\bsplit
\Ft_{(1,J-2)} &= O(r^{-1})
   \dk^{J-1} \nab_\Rhat (\Pdot, \Bbdot)+ O(r^{-3} )   \dk^{\le J+1}   \Ga_b + O(r^{- 1 })  \dk^{\le J-1}(\Pdot, \Bbdot) \\
   &+ r^{-1}   \dk^{\le J+1}(\Ga_b\c \Rc_b),\\
   \Ft_{(2,J-2)} &=  -4\om \Bbt+  O(r^{-3})\dk^{\le J}(\nab_4, \nab)\Ga_b  +O(r^{-1} )\dk^{\le J-1} \nab_\Rhat\Pdot \\
   &+O(ar^{-2} ) \dk^{\le J-2}\dkb  \nab_\Rhat\Bbdot  + O(r^{-1})  \dk^{\le J-1}\Pdot+ O(r^{-1})  \dk^{\le J}\Bb \\
   &+ O(r^{-3}) \dk^{\le J}\Ga_b+  \dk^{\le J+1 }  \Big( \Ga_b  \c \Rc_b\Big).
   \end{split}
\eea
   \end{enumerate}
   \end{lemma}
   
   \begin{remark}
   Note that  the error terms  for the  equations which  contain $\Bt$ or $\Bbt$ 
    on the left hand side
   are more structured than those corresponding to  the equations for $\Ptp $ and $\Ptm$. 
   The reason, as it will become apparent in the next section, is that we already control
    the  quantities  $\Ptp$ and $\Ptm$.
   \end{remark} 
   
    \begin{proof}
    To prove  the first statement, we apply Proposition
    \ref{Prop:BianchPairs-tilde4-k}   with $k=J-2$   to the Bianchi pair \eqref{eq:linearizationbyLie_T-B2}
\beaa
     \begin{split}
  \nabc_3\Bdot  +\tr\Xb \Bdot&=- \DDs_1\,\, \ov{\Pdot} +O(ar^{-2}) \ov{\Pdot}  +O(r^{-3} )\Lieb_\T\Hc  +O(r^{-4}) \Ga_b  + r^{-2} \dk^{\le 1}(\Ga_b\c \Rc_b),\\
   \nabc_4 \ov{\Pdot } +\frac{3}{2}\ov{\tr X}\, \ov{ \Pdot}  &= \DDd_1\Bdot  +O(ar^{-2}) \Bdot+O(r^{-3} )  \dk^{\le 1} \Ga_b   + r^{-2} \dk^{\le 1}(\Ga_b\c \Rc_b).
   \end{split}
 \eeaa
    We deduce, with $\Psi_{(1)}=\Bdot, \, \Psi_{(2)}=\ov{\Pdot}$,
    \beaa
    \nabc_3  \Bt +  \left(1-\frac{J-2}{2}\right)\tr\Xb \Bt &=& -\DDs_1\, \ov{ \Ptp} +\Ft_{(1, J-2)},\\
    \nabc_4 \Ptp       +\left(\frac 3 2 -\frac{J-2}{2}\right) \ov{\tr X} \Ptp &=&\DDd_1     \ov\Bt +\Ft_{(2, J-2)},
    \eeaa
    with, see \eqref{eq:Ft_(1,k)Ft_(2,k)},
    \beaa
 \bsplit
 \Ft_{(1,k)}&=\big( \ov{q} \nabc_4 \big)^k \nabc^2_\Rhat   F_{(1)} - 2 \om \big( \ov{q} \nabc_4\big)^k \nabc_3 \nabc_\Rhat \Psi_{(1)} \\
 &+O(r^{-1}) \dk^{\le k+1} \nabc_\Rhat  \Psi_{(2)} 
 +O(ar^{-2} ) \dk^{\le k}\dkb  \nabc_\Rhat \Psi_{(1)} \\
 &+ O(r^{- 1 })  \dk^{\le k+1}  \big(\Psi_{(1)}, \Psi_{(2)}\big) + \dk^{\le k+2 }  \Big( \Ga_b  \c ( \Psi_{(1)}, \Psi_{(2)} )\Big),\\
  \Ft_{(2,k)}&=\big( \ov{q} \nabc_4 \big)^k \nabc^2_\Rhat   F_{(2)} + 2\frac{\De}{|q|^2}  \om \big( \ov{q} \nabc_4\big)^k \nabc_3 \nabc_\Rhat \Psi_{(2)} \\
 &+O(r^{-1} )\dk^{\le k +1}\nabc_\Rhat \Psi_{(1)}
 +O(ar^{-2} ) \dk^{\le k}\dkb  \nabc_\Rhat  \Psi_{(2)} \\
 &+ O(r^{- 1 })  \dk^{\le k+1}  \big(\Psi_{(1)}, \Psi_{(2)}\big) + \dk^{\le k+2 }  \Big( \Ga_b  \c ( \Psi_{(1)}, \Psi_{(2)} )\Big),
 \end{split}
 \eeaa
 where
    \beaa
    F_{(1)}&=&O(ar^{-2}) \ov{\Pdot} +O(r^{-3} )\Lieb_\T\Hc  +O(r^{-4}) \Ga_b
    + r^{-2} \dk^{\le 1}(\Ga_b\c \Rc_b),\\
   F_{(2)}&=&  O(ar^{-2}) \Bdot+O(r^{-3} )  \dk^{\le 1} \Ga_b   + r^{-2} \dk^{\le 1}(\Ga_b\c \Rc_b).
    \eeaa  
   Setting $k=J-2$ and since we have chosen $\Psi_{(1)}=\Bdot, \, \Psi_{(2)}=\ov{\Pdot}$, we deduce
   \beaa
   \Ft_{(1, J-2)}&=& - 2 \om \big( \ov{q} \nabc_4\big)^{J-2} \nabc_3 \nabc_\Rhat \Bdot+ O(r^{- 1 })
   \dk^{\leq J-1} \nabc_\Rhat \Pdot+ O(r^{-3} )   \dk^{\le J} \Lieb_\T \Hc\\
   && +O(r^{-4})\dk^{\leq J}\Ga_b +O(ar^{-2} )\dk^{\le J-2 }  \dkb \nabc_\Rhat \Bdot  + O(r^{- 1 })  \dk^{\le J-1}(\Pdot, \Bdot)\\
   && + r^{-2}   \dk^{\le J+1}(\Ga_b\c \Rc_b).
   \eeaa
   We rewrite in the  simplified form, using in particular $2\Lieb_\T\Hc=\nab_4\Hc+O(1)(\nab_3, \nab)\Ga_b+O(1)\Ga_b$, 
   \beaa
    \Ft_{(1, J-2)}&=&  O(r^{-1} )  \dk^{J-2} \nab_\Rhat \big(\nabc_3  \Bdot, \dkb \Bdot) + O(r^{-1})  \dk^{J-1} \nab_\Rhat \Pdot + O(r^{-3} )   \dk^{\le J}   \nab_4 \Hc\\
    && +O(r^{-3})\dk^{\leq J}(\nab_3, \nab)\Ga_b+O(r^{-3})\dk^{\leq J}\Ga_b+ O(r^{-1})  \dk^{\le J-1}(\Pdot, \Bdot) \\
    &&+ r^{-2}   \dk^{\le J+1}(\Ga_b\c \Rc_b).
   \eeaa      
   We  next rewrite the  dangerous     term   $O(r^{-3} )   \dk^{\le J}   \nab_4 \Hc$  using 
    the equation, see \eqref{eqts:usefulNullstructure},
    \beaa
   \nabc_4\Hc -\nabc_3\Xi&=-B +  O(r^{-1} ) \Ga_b+\Ga_b\c\Ga_g,
   \eeaa
   Thus
   \beaa
   O(r^{-3} )   \dk^{\le J}   \nab_4 \Hc&=&O(r^{-3} ) \nabc_3  \dk^{\le J}     \Xi+ 
   O(r^{-3} )   \dk^{\le J} B +O(r^{-4} )\dk^{\le J}  \Ga_b+ r^{-3} \dk^{\le J} \big( \Ga_b\c\Ga_b\big).
   \eeaa
   Therefore,
   \beaa
   \Ft_{(1, J-2)}  &=&  O(r^{-1} )  \dk^{J-2} \nab_\Rhat \big(\nabc_3  \Bdot, \dkb \Bdot)  + O(r^{-1})  \dk^{J-1} \nab_\Rhat \Pdot +O(r^{-3} )(\nab_3, \nab)\dk^{\le J}\Ga_b \\
   && + O(r^{-1})  \dk^{\le J-1}\Pdot+ O(r^{-1})  \dk^{\le J}B +O(r^{-3} )\dk^{\le J}  \Ga_b + r^{-2}   \dk^{\le J+1}(\Ga_b\c \Rc_b) \\
   &&  + r^{-2} \dk^{\le J} \big( \Ga_b\c\Ga_b\big)
   \eeaa
   as stated in the first equation of \eqref{eq:Ft-Bt-Ptp}. The second equation of \eqref{eq:Ft-Bt-Ptp} is derived in the same manner and is in fact simpler.

   \medskip
   
      To prove the second  statement  we apply  the second part of  Proposition \ref{Prop:BianchPairs-tilde3-k} to   the Bianchi pair \eqref{eq:linearizationbyLie_T-B3}
\beaa
     \begin{split}
\nabc_3\Pdot +      \frac{3}{2}\ov{\tr\Xb}   \Pdot      &=- \DDs_1\Bbdot + O(ar^{-2} ) \Bbdot+O(r^{-3} )\dk^{\le 1} \Ga_b +r^{-1} \dk^{\le 1}(\Ga_b\c \Rc_b),\\
\nabc_4\Bbdot +\tr X\Bbdot  &=\DDs_1\,  \Pdot + O(ar^{-2} )\Pdot+O(r^{-3}) \Lieb_\T \Hbc  +O(r^{-4})\Ga_b + r^{-1} \dk^{\le 1}(\Ga_b\c \Rc_b).
 \end{split}
\eeaa
We deduce, for $\Ptm=(\nabc_3)^{\le J-2} \nabc^2_{\Rhat}\Pdot,\, \Bbt=(\nabc_3 )^{\le J-2} 
   \nabc^2_\Rhat\Bdot$,
 \beaa
  \begin{split}
\nabc_3  \Ptm+      \frac{3}{2}\ov{\tr\Xb}  \Ptm  &=- \DDd_1 \Bbt+\Ft_{1, J-2},  \\
  \nabc_4 \Bbt +\tr X   \Bbt  &=\DDs_1 \Ptm +\Ft_{2, J-2},
  \end{split}
\eeaa
with, see  \eqref{eq:Ft_(1,k)Ft_(2,k)3},
\beaa
\Ft_{(1,k)}&=&\big( \nabc_3 \big)^k \nabc^2_\Rhat   F_{(1)} - 2 \om \big(  \nabc_3\big)^k \nabc_3 \nabc_\Rhat \Psi_{(1)} \\
 &&+O(r^{-1} )\dk^{\le k +1} \nabc_\Rhat  \Psi_{(2)} 
 +O(ar^{-2} )\dk^{\le k}\dkb  \nabc_\Rhat \Psi_{(1)} \\
 &&+ O(r^{- 1 })  \dk^{\le k+1}  \big(\Psi_{(1)}, \Psi_{(2)}\big) + \dk^{\le k+2}\Big( \Ga_b  \c ( \Psi_{(1)}, \Psi_{(2)} )\Big),\\
   \Ft_{(2,k)}&=&\big(  \nabc_3 \big)^k \nabc^2_\Rhat   F_{(2)} + 2\frac{\De}{|q|^2}  \om \big( \nabc_3\big)^k \nabc_3 \nab_\Rhat \Psi_{(2)} \\
 &&+O(r^{-1} )\dk^{\le k+1} \nabc_\Rhat \Psi_{(1)}
 +O(ar^{-2} ) \dk^{\le k}\dkb  \nabc_\Rhat \Psi_{(2)} \\
 &&+ O(r^{- 1 })  \dk^{\le k+1}  \big(\Psi_{(1)}, \Psi_{(2)}\big) + \dk^{\le k+2 }  \Big( \Ga_b  \c ( \Psi_{(1)}, \Psi_{(2)} )\Big),
 \eeaa
where $\Psi_{(1)}=\Pdot, \Psi_{(2)}= \Bbdot$ and
\beaa
F_{(1)}&=&O(ar^{-2} ) \Bbdot+O(r^{-3} )\dk^{\le 1} \Ga_b +r^{-1} \dk^{\le 1}(\Ga_b\c \Rc_b),\\
F_{(2)}&=& O(ar^{-2} )\Pdot+O(r^{-3}) \Lieb_\T \Hbc  +O(r^{-4})\Ga_b  + r^{-1} \dk^{\le 1}(\Ga_b\c \Rc_b).
\eeaa
We deduce
\beaa
\Ft_{(1,J-2)} &= &O(r^{-1})
   \dk^{J-1} \nab_\Rhat (\Pdot, \Bbdot)+ O(r^{-3} )   \dk^{\le J+1}   \Ga_b + O(r^{- 1 })  \dk^{\le J-1}(\Pdot, \Bdot) \\
   &&+ r^{-1}   \dk^{\le J+1}(\Ga_b\c \Rc_b)
\eeaa
 as stated in  \eqref{eq:Ft-Bbt-Ptp}. Also
 \beaa
 \Ft_{(2,J-2)} &= &\big(  \nabc_3 \big)^{J-2} \nabc^2_\Rhat \Big( O(ar^{-2} )\Pdot+O(r^{-3}) \Lieb_\T \Hbc  +O(r^{-4})\Ga_b  + r^{-1} \dk^{\le 1}(\Ga_b\c \Rc_b)\Big)\\
 &&  +  2    \om\big( \nabc_3\big)^{J-2} \left( \frac{\De}{|q|^2}\nabc_3 \nabc_\Rhat \Bbdot\right) +O(r^{-1} )\dk^{\le J-1} \nab_\Rhat\Pdot \\
 &&+O(ar^{-2} ) \dk^{\le J-2}\dkb  \nab_\Rhat\Bbdot + O(r^{- 1 })  \dk^{\le J-1}  \big(\Pdot, \Bbdot\big) \\
 &&+ \dk^{\le J }  \Big( \Ga_b  \c ( \Pdot, \Bbdot)\Big).
 \eeaa
 Recalling that  $\Rhat= \frac 1 2( e_4-\frac{\De}{|q|^2}  e_3)$, we write
 \beaa
  \frac{\De}{|q|^2}  \nabc_3 \nabc_\Rhat \Bbdot=\nabc_4 \nabc_\Rhat \Bbdot- 2 \nabc_\Rhat^2 \Bbdot.
 \eeaa
 We deduce
 \beaa
  \big( \nabc_3\big)^{J-2} \left( \frac{\De}{|q|^2}\nabc_3 \nabc_\Rhat \Bbdot\right) &=&-  2 \Bbt
  +\dk^{\le J-2}  \nabc_\Rhat\nabc_4 \Bbdot+ O(r^{-1} )\dk^{\leq J-1}  \Bbdot
  \eeaa
  and hence, using also $2\Lieb_\T\Hbc=\frac{\De}{|q|^2}\nab_3\Hbc+O(1)(\nab_4, \nab)\Ga_g+O(1)\Ga_g$, we obtain
 \beaa
 \Ft_{(2,J-2)} &= & O(r^{-3})    \dk^{\le J}      \nabc_3 \Hbc +O(r^{-3})\dk^{\le J}(\nab_4, \nab)\Ga_b +O(r^{-4})\dk^{\le J}\Ga_b  -  4\om \Bbt \\
 && +O(r^{-1} )\dk^{\le J-1} \nab_\Rhat\Pdot +O(ar^{-2} )\dk^{\le J-2}\dkb  \nab_\Rhat\Bbdot + O(r^{- 1 })  \dk^{\le J-1}  \big(\Pdot, \Bbdot\big) \\
 &&+ \dk^{\le J+1 }\Big( \Ga_b  \c \Rc_b\Big).
 \eeaa
 In view of  the equation  \eqref{eqts:usefulNullstructure}  for  $\nabc_3\Hbc $, we  have
\beaa
\nabc_3\Hbc -\nabc_4\Xib&=\Bb +  O(r^{-1} ) \Ga_b+\Ga_b\c\Ga_b.
\eeaa 
 Therefore,
 \beaa
 \Ft_{(2,J-2)} &= &  - 4\om \Bbt+  O(r^{-3})\dk^{\le J}(\nab_4, \nab)\Ga_b  +O(r^{-1} )\dk^{\le J-1} \nab_\Rhat\Pdot +O(ar^{-2} ) \dk^{\le J-2}\dkb  \nab_\Rhat\Bbdot\\
&&   + O(r^{-1})  \dk^{\le J-1}\Pdot+ O(r^{-1})  \dk^{\le J}\Bb + O(r^{-3}) \dk^{\le J}\Ga_b+  \dk^{\le J+1 }  \Big( \Ga_b  \c \Rc_b\Big)
 \eeaa
as stated in  \eqref{eq:Ft-Bbt-Ptp}. This concludes the proof of Lemma \ref{lemma:EstimatesBtPtp}.
         \end{proof}


\subsection{Estimates for $\Bt$} 
    \lab{subsection:Estimatesfor-Bt}
    

   We provide estimates for $\Bt$  using  \eqref{eq:linearization-tildes-B-P}  and Proposition \ref{Prop:Bainchi-pairsEstimates-integrated}.
   \begin{proposition}
    \lab{Prop:-EstimatesBtPtp3}
   The following   estimate holds true for $\Bt$, with $b=2+\de$ and for $a$ sufficiently small,  
   \bea\lab{eq:lemma-EstimatesBtPtp3}
   \bsplit
  \int_{\MM } r^{b-1} |\Bt |^2+\int_{\pr^+\MM} r^{b}|\Bt |^2 &\les \de_{J+1}[\Pc]  +\ep_J^2   +\ep_0^2 +\sqrt{ \de_{J+1}[\Pc]} \left(B^J_\de[r^2B]+\Sk_{J+1}^2\right)^{\frac{1}{2}}\\
   &+O(a^2)B_\de[  r^{J+1}  \nab^J_4 B] +\ep_J\Big(B^J_\de[ r^2B]\Big)^{\frac{1}{2}}\\
   & +\Sk_{J+1}\Big(\de_{J+1}[\Pc]  +\ep_J^2 +\ep_0^2 +O(a^2)B_\de[  r^{J+1}  \nab^J_4 B]\Big)^{\frac{1}{2}}.
   \end{split}
  \eea
    \end{proposition}
    
    \begin{proof}
   Note that the system \eqref{eq:linearization-tildes-B-P} is of  the form \eqref{eq:modelbainchipairequations11-simple} with 
   \beaa
   \Psi_{(1)}=\Bt, \qquad \Psi_{(2)}=\Ptp, \qquad c_{(1)}=1-\frac{J-2}{2}, \qquad c_{(2)}= \frac  3 2 -\frac{J-2}{2}.
   \eeaa
 Recall also  that the  signatures  of  $\Ptp$ and $\Bt$  are  respectively $J-2$ and $J-1$. We are thus in the case
  corresponding to $2k-1>0$. 
    Observe that the condition $ -2 c_{(1)} +1+\frac b 2 > 0$ is verified for $ b= 2+ \de $
   and therefore  we can apply \eqref{eq:general-proposition-first-bianchi-pairs} in $\MM=\MM(1, \tau)$, $\tau\le \tau_*$,
   to derive  
  \beaa 
\begin{split}
& \int_{\MM } r^{b-1} |\Psi_{(1)}|^2+\int_{\pr^+\MM} r^{b}|\Psi_{(1)}|^2\\
\les& \int_{\MM} r^{b-1} |\Psi_{(2)}|^2 +\left|\int_{\MM}|q|^{b}\Re(\Ft_{(1, J-2)}\c\ov{\Psi_{(1)}})\right|  +\int_{\MM}r^{b}\big| \Ft_{(2, J-2)}\big|\, \big|\Psi_{(2)}\big| \\
&+ \int_{\pr^-\MM} \big( r^{b}|\Psi_{(1)}|^2+  r^{b-2}|\Psi_{(2)}|^2\big),
\end{split}
\eeaa 
  with  $ \Ft_{(1, J-2)}, \Ft_{(2, J-2)}$ given by formula \eqref{eq:Ft-Bt-Ptp}.
  We decompose
  \beaa
  \Ft_{(1, J-2)}= O(r^{-3} ) (\nab_3, \nab)\dk^{\le J}\Ga_b + \Ft'_{(1, J-2)}.
  \eeaa
  Thus
  \beaa
   \int_{\MM } r^{b-1} |\Bt |^2+\int_{\pr^+\MM} r^{b}|\Bt |^2 &\les& \int_{\MM} r^{b-1} | \Ptp |^2 +\int_{\MM}r^{b}   \Big(\big| \Ft'_{(1, J-2)}\big|\, \big|\Bt \big| +  \big| \Ft_{(2, J-2)}\big|\, \big|\Ptp\big| \Big)\\
   &&+\big|I\big| +\ep_0^2 
  \eeaa
  where
  \bea
   I := \int_{\MM}    O(r^{-3} )    |q|^{b}\Re(    (\nab_3, \nab)\dk^{\le J}\Ga_b \c\ov{\Bt}) +\int_{\MM}O(r^{-1})|q|^{b}\Re(\dk^{\le J}B\c\ov{\Bt}).
  \eea
  Recall that
  \beaa
  \de_{J+1}[\Pc]= \BEF^{ J }_{\de}[r^2\Pc]\ge \int_{\Mtrap} |\nab_\Rhat \dk^{\le J}\Pc|^2 + | \dk^{\le J}\Pc|^2 +\int_{\Mntrap} r^{\de+1}| \dk^{\le J+1}\Pc|^2.  
  \eeaa
  Therefore, with $b=2+\de$,
 \beaa
 \int_{\MM} r^{b-1} | \Ptp |^2  = \int_{\MM} r^{1+\de} \big|(\ov{q}\nabc_4)^{\le J-2} \nabc^2_{\Rhat} \ov{\Pdot}\big|^2\les  \int_{\MM} r^{1+\de} \big| \dk^{J-1} \nabc_{\Rhat} \ov{\Pdot}\big|\les \de_{J+1}[\Pc].
 \eeaa
 By Cauchy-Schwartz and absorbing the  term in $\Bt$ to the left  we easily deduce
 \bea
 \lab{eq:EstimateBt:1}
 \bsplit
   \int_{\MM } r^{b-1} |\Bt |^2+\int_{\pr^+\MM} r^{b}|\Bt |^2 &\les  \de_{J+1}[\Pc] +  \int_{\MM} r^{b+1}   \big|\Ft_{(1, J-2 )}'\big|^2\\
   &+ \sqrt{ \de_{J+1}[\Pc]} \left(\int_{\MM} r^{b+1} \big|\Ft_{(2, J-2 )}\big|^2\right)^{\frac{1}{2}}
   +\big|I\big| +\ep_0^2.
   \end{split}
  \eea
 In view of \eqref{eq:Ft-Bt-Ptp} and the definition of $I$, we have
 \beaa
\bsplit
 \Ft_{(1, J-2)}' &=   O(r^{-1} )  \dk^{J-2} \nab_\Rhat \big(\nabc_3  \Bdot, \dkb \Bdot\big)  + O(r^{-1})  \dk^{J-1} \nab_\Rhat \Pdot +   O(r^{-1})  \dk^{\le J-1}\Pdot\\
 & + O(r^{-3} )\dk^{\le J}  \Ga_b  + r^{-2}   \dk^{\le J+1}(\Ga_b\c \Rc_b) + r^{-2} \dk^{\le J} \big( \Ga_b\c\Ga_b\big),\\
   \Ft_{(2, J-2)} &= O(r^{-1})  \dk^{J-1} \nab_\Rhat (\Pdot, \Bdot)+ O(r^{-3} )   \dk^{\le J+1}   \Ga_b + O(r^{-1})  \dk^{\le J-1}(\Pdot, \Bdot) \\
   &+ r^{-2}   \dk^{\le J+1}(\Ga_b\c \Rc_b).
\end{split}
\eeaa
 We deduce, using the  bootstrap assumptions, as well as the choice $b=2+\de$, 
 \beaa
 \int_{\MM} r^{b+1}  \big|\Ft_{(1, J-2 )}'\big|^2&\les& B^{J-2}_\de[ r^2( \nab_3 \Bdot, \nab \Bdot) ]+   B^{J-1}_\de[ r^2\Pdot] +\int_\MM r^{b-5} |\dk^{\le J} \Ga_b|^2   +\ep_0^2, 
 \eeaa
 and
 \beaa
  \int_{\MM} r^{b+1}  \big|\Ft_{(2, J-2 )}\big|^2&\les&  B^{J-1}_\de[ r^2(\Bdot, \Pdot)]+\int_\MM r^{b-5} |\dk^{\le J+1} \Ga_b|^2  +\ep_0^2. 
  \eeaa
  According to Lemma \ref{lemma:auxilliarynormsforGa_b} we have,  since $b=2+\de$,
  \beaa
   \int_{\MM} r^{b-5}| \dk^{\leq J+1}\Ga_b|^2\les    \Sk^2_{J+1}, \qquad \int_{\MM} r^{b-5}| \dk^{J}\Ga_b|^2\les    \Sk^2_{J}\les \ep_J^2,  
  \eeaa
where we have used the induction hypothesis.  Also, in view of the proof of Lemma \ref{Lemma:comparisonofnormsBEF-Rk-forpsi}, we have $\Pdot=\T\Pc+r^{-3}\Ga_g$ so that 
\beaa
B^{J-1}_\de[ r^2\Pdot] &\les& B^{J-1}_\de[ r^2\T\Pc]+B^{J-1}_\de[r^{-1}\Ga_b]\les  \de_{J+1}[\Pc]+\Sk_J^2\les  \de_{J+1}[\Pc]+\ep_J^2.
\eeaa
Therefore,
\beaa
&&\int_{\MM} r^{b+1}   \big|\Ft_{(1, J-2 )}'\big|^2 + \sqrt{ \de_{J+1}[\Pc]} \left(\int_{\MM} r^{b+1} \big|\Ft_{(2, J-2 )}\big|^2\right)^{\frac{1}{2}}\\
&\les&  B^{J-2}_\de[ r^2( \nab_3 \Bdot, \nab \Bdot) ] +\de_{J+1}[\Pc]  +\ep_J^2   +\ep_0^2\\
&& +\sqrt{ \de_{J+1}[\Pc]} \left(B^{J-1}_\de[ r^2\Bdot]+\de_{J+1}[\Pc]+\Sk_{J+1}^2+  \ep_J^2+\ep_0^2\right)^{\frac{1}{2}}\\
&\les&  B^{J-2}_\de[ r^2( \nab_3 \Bdot, \nab \Bdot) ] +\de_{J+1}[\Pc]  +\ep_J^2   +\ep_0^2 +\sqrt{ \de_{J+1}[\Pc]} \left(B^J_\de[r^2B]+\Sk_{J+1}^2\right)^{\frac{1}{2}}\\
&\les&B^{J-1}_\de[ r^2(\nab_3B, r\nab B] +\de_{J+1}[\Pc]  +\ep_J^2   +\ep_0^2 +\sqrt{ \de_{J+1}[\Pc]} \left(B^J_\de[r^2B]+\Sk_{J+1}^2\right)^{\frac{1}{2}}.
\eeaa
 Using \eqref{eq:EstmatesforB-1-strong}, i.e. 
\beaa
B^{J-1}_\de[ r^2(\nab_3B, r\nab B)] &\les &  \de_{J+1}[\Pc]  +\ep_J^2 +\ep_0^2 +O(a^2)B_\de[  r^{J+1}  \nab^J_4 B],
\eeaa
we infer 
\beaa
&&\int_{\MM} r^{b+1}   \big|\Ft_{(1, J-2 )}'\big|^2 + \sqrt{ \de_{J+1}[\Pc]} \left(\int_{\MM} r^{b+1} \big|\Ft_{(2, J-2 )}\big|^2\right)^{\frac{1}{2}}\\
&\les&  \de_{J+1}[\Pc]  +\ep_J^2   +\ep_0^2 +\sqrt{ \de_{J+1}[\Pc]} \left(B^J_\de[r^2B]+\Sk_{J+1}^2\right)^{\frac{1}{2}}+O(a^2)B_\de[  r^{J+1}  \nab^J_4 B].
\eeaa
 Back to \eqref{eq:EstimateBt:1} we infer that
\bea \lab{eq:EstimateBt:2}
\nn   \int_{\MM } r^{b-1} |\Bt |^2+\int_{\pr^+\MM} r^{b}|\Bt |^2 &\les& \big|I\big|+\de_{J+1}[\Pc]  +\ep_J^2   +\ep_0^2 +\sqrt{ \de_{J+1}[\Pc]} \left(B^J_\de[r^2B]+\Sk_{J+1}^2\right)^{\frac{1}{2}}\\
   &&+O(a^2)B_\de[  r^{J+1}  \nab^J_4 B].
  \eea
    
  It remains to estimate the term  $I$. We decompose it as 
  \beaa
  I &=& I_1+I_2,\\
  I_1 &:=& \int_{\MM}    O(r^{-3} )    |q|^{b}\Re(  (\nab_3, \nab)\dk^{\le J}\Ga_b \c\ov{\Bt}),\\
  I_2 &:=& \int_{\MM}O(r^{-1})|q|^{b}\Re(\dk^{\le J}B\c\ov{\Bt}),
  \eeaa
  and estimate $I_1$ and $I_2$ separately starting with $I_1$.  Integrating by parts,  we have
    \beaa
 I_1 &=&  \int_{\MM}    O(r^{-3} )    |q|^{b}\Re(    (\nab_3, \nab)\dk^{\le J}\Ga_b \c\ov{\Bt}) \\
    &=& -\int_{\MM}O(r^{-3})|q|^{b}\Re\left(\dk^{\leq J}\Ga_b\c\ov{(\nab_3, \nab)\Bt }\right) + \int_{\pr^+\MM}O(r^{b-3})|\dk^{\leq J}\Ga_b||\Bt|\\
    &&+\int_{\MM}O(r^{b-4})|\dk^{\leq J}\Ga_b||\Bt|+\ep_0^2.
     \eeaa
     Since $\Bt=(\ov{q}\nabc_4 )^{\le J-2}  \nabc^2_\Rhat\Bdot$, we may integrate the first term on the RHS by parts again and obtain 
    \beaa
I_1  &=& -\int_{\MM}O(r^{b-3})|\dk^{\leq J+1}\Ga_b||\nab_{\Rhat}(\nab_3, \nab)\dk^{\leq J-1}B| + \int_{\pr^+\MM}O(r^{b-3})|\dk^{\leq J}\Ga_b||\dk^{\leq J+1}B|\\
    &&+\int_{\MM}O(r^{b-4})|\dk^{\leq J}\Ga_b||\nab_{\Rhat}\dk^{\leq J}B|+\ep_0^2.
     \eeaa  
     Hence, using Lemma \ref{lemma:auxilliarynormsforGa_b}, Lemma \ref{Lemma:comparisonofnormsBEF-Rk}, and the fact that $b=2+\de$, we infer
     \beaa
     |I_1| &\les& \left(\int_{\MM}r^{b-5}|\dk^{\leq J+1}\Ga_b|^2\right)^{\frac{1}{2}}\Big(B^{J-1}_\de[ r^2(\nab_3B, r\nab B)]\Big)^{\frac{1}{2}}\\
&&+\left(\int_{\MM}r^{b-5}|\dk^{\leq J}\Ga_b|^2\right)^{\frac{1}{2}}\Big(B^{J}_\de[ r^2B]\Big)^{\frac{1}{2}}+\ep_0^2\\
     &\les& \Sk_{J+1}\Big(B^{J-1}_\de[ r^2(\nab_3B, r\nab B)]\Big)^{\frac{1}{2}}+\ep_J\Big(B^J_\de[ r^2B]\Big)^{\frac{1}{2}}+\ep_0^2.
     \eeaa
     Together with \eqref{eq:EstmatesforB-1-strong}, we deduce
      \beaa
     |I_1|      &\les& \Sk_{J+1}\Big(\de_{J+1}[\Pc]  +\ep_J^2 +\ep_0^2 +O(a^2)B_\de[  r^{J+1}  \nab^J_4 B]\Big)^{\frac{1}{2}}+\ep_J\Big(B^J_\de[ r^2B]\Big)^{\frac{1}{2}}+\ep_0^2.
     \eeaa

Next, we estimate $I_2$. Using the fact that $\Bt=(\ov{q}\nabc_4 )^{\le J-2}  \nabc^2_\Rhat\Bdot$, we have
\beaa
\Bt &=& \nab_{\Rhat}^2\dk^{\leq J-1}B+O(1)\nab_{\Rhat}\dk^{\leq J-1}B+O(1)\dk^{\leq J-1}B
\eeaa
and hence, introducing the notation 
\beaa
I_{2,1}:=\int_{\MM}O(r^{-1})|q|^{b}\Re(\dk^{\le J}B\c\ov{\nab_{\Rhat}^2\dk^{\leq J-1}B}),
\eeaa
we have, using Lemma \ref{Lemma:comparisonofnormsBEF-Rk} and the fact that $b=2+\de$,
\beaa
|I_2| &\les& |I_{2,1}|+\int_{\MM}r^{b-1}|\dk^{\le J}B||\nab_{\Rhat}\dk^{\leq J-1}B|+\int_{\MM}r^{b-1}|\dk^{\le J}B||\dk^{\leq J-1}B|\\
&\les& |I_{2,1}|+\Big(B_\de^J[r^2B]\Big)^{\frac{1}{2}}\Big(B_\de^{J-1}[r^2B]\Big)^{\frac{1}{2}}\\
&\les& |I_{2,1}|+\ep_J\Big(B_\de^{J}[r^2B]\Big)^{\frac{1}{2}}.
\eeaa
Also, integrating by parts, we have
\beaa
I_{2,1} &=& \int_{\MM}O(r^{-1})|q|^{b}\Re(\dk^{\le J}B\c\ov{\nab_{\Rhat}^2\dk^{\leq J-1}B})\\
&=& -\int_{\MM}O(r^{-1})|q|^{b}\Re(\nab_{\Rhat}\dk^{\le J}B\c\ov{\nab_{\Rhat}\dk^{\leq J-1}B})\\
&&+\int_{\pr^+\MM}O(r^{-1})|q|^{b}\Re(\dk^{\le J}B\c\ov{\nab_{\Rhat}^2\dk^{\leq J-1}B})\\
&&+\int_{\MM}O(r^{-2})|q|^{b}\Re(\dk^{\le J}B\c\ov{\nab_{\Rhat}\dk^{\leq J-1}B})+\ep_0^2
\eeaa
and hence
\beaa
|I_{2,1}| &\les& \Big(B_\de^J[r^2B]\Big)^{\frac{1}{2}}\Big(B_\de^{J-1}[r^2B]\Big)^{\frac{1}{2}}+\ep_0^2\\
&\les& \ep_J\Big(B_\de^{J}[r^2B]\Big)^{\frac{1}{2}}+\ep_0^2.
\eeaa
We deduce
\beaa
|I_2| &\les& |I_{2,1}|+\ep_J\Big(B_\de^{J}[r^2B]\Big)^{\frac{1}{2}}\\
&\les& \ep_J\Big(B_\de^{J}[r^2B]\Big)^{\frac{1}{2}}+\ep_0^2.
\eeaa
Together with the above estimate for $I_1$, this yields
 \beaa
 |I| &\les& |I_1|+|I_2|\\
 &\les&  \Sk_{J+1}\Big(\de_{J+1}[\Pc]  +\ep_J^2 +\ep_0^2 +O(a^2)B_\de[  r^{J+1}  \nab^J_4 B]\Big)^{\frac{1}{2}}+\ep_J\Big(B^J_\de[ r^2B]\Big)^{\frac{1}{2}}+\ep_0^2.
     \eeaa

 Finally, plugging the above estimate for $I$ in \eqref{eq:EstimateBt:2}, we infer
 \beaa
\nn   \int_{\MM } r^{b-1} |\Bt |^2+\int_{\pr^+\MM} r^{b}|\Bt |^2 &\les& \de_{J+1}[\Pc]  +\ep_J^2   +\ep_0^2 +\sqrt{ \de_{J+1}[\Pc]} \left(B^J_\de[r^2B]+\Sk_{J+1}^2\right)^{\frac{1}{2}}\\
   &&+O(a^2)B_\de[  r^{J+1}  \nab^J_4 B] +\ep_J\Big(B^J_\de[ r^2B]\Big)^{\frac{1}{2}}\\
   && +\Sk_{J+1}\Big(\de_{J+1}[\Pc]  +\ep_J^2 +\ep_0^2 +O(a^2)B_\de[  r^{J+1}  \nab^J_4 B]\Big)^{\frac{1}{2}}
  \eeaa
as stated. This concludes the proof of Proposition \ref{Prop:-EstimatesBtPtp3}.
   \end{proof}

   
   \subsection{Estimates for $\Bbt$}
   \lab{subsection:Estimatesfor-Bbt}
   
   
    We provide estimates for $\Bbt$  using \eqref{eq:linearization-tildes-Bb-P}  and Proposition \ref{Prop:Bainchi-pairsEstimates-integrated}.
    \begin{proposition}
    \lab{Prop:lemma-EstimatesBbtPtm3}
   The following   estimates hold true for $\Bt$, with $b=-\de$ and $a$ sufficiently small,
   \bea
   \lab{eq:lemma-EstimatesBbtPtm3}
   \begin{split}
 \int_{\MM } r^{b-1}   | \Bbt |^2 +  \int_{\pr^+\MM } r^{b-2} |\Bbt |^2 &\les  \de_{J+1}[\Pc]  +\ep_J^2   +\ep_0^2 +\sqrt{ \de_{J+1}[\Pc]} \left(B^J_\de[\Bb]+\Sk_{J+1}^2\right)^{\frac{1}{2}}\\
 & +O(a^2)B_\de[\nab^J_3\Bb] +\ep_J\Big(B^J_\de[\Bb]\Big)^{\frac{1}{2}}\\
 & +\Sk_{J+1}\Big(\de_{J+1}[\Pc]  +\ep_J^2 +\ep_0^2 +O(a^2)B_\de[\nab^J_3\Bb]\Big)^{\frac{1}{2}}.
 \end{split}
\eea
   \end{proposition}
  
   \begin{proof}
   The proof is similar to the one of Proposition \ref{Prop:-EstimatesBtPtp3}.  We  start with 
   the Bianchi  pair \eqref{eq:linearization-tildes-Bb-P}
     \beaa
  \begin{split}
\nabc_3  \Ptm+      \frac{3}{2}\ov{\tr\Xb}  \Ptm  &=- \DDd_1 \Bbt+\Ft_{1, J-2},\\
  \nabc_4 \Bbt +\tr X   \Bbt  &=\DDs_1 \Ptm +\Ft_{2, J-2},
  \end{split}
\eeaa
with
\beaa
\bsplit
\Ft_{(1,J-2)} &= O(r^{-1})
   \dk^{J-1} \nab_\Rhat (\Pdot, \Bbdot)+ O(r^{-3} )   \dk^{\le J+1}   \Ga_b + O(r^{- 1 })  \dk^{\le J-1}(\Pdot, \Bbdot) \\
   &+ r^{-1}   \dk^{\le J+1}(\Ga_b\c \Rc_b),\\
   \Ft_{(2,J-2)} &=   - 4\om \Bbt+  O(r^{-3})\dk^{\le J}(\nab_4, \nab)\Ga_b  +\Ft'_{(2,J-2)},\\
\Ft'_{(2,J-2)}   &=O(r^{-1} )\dk^{\le J-1} \nab_\Rhat\Pdot +O(ar^{-2} ) \dk^{\le J-2}\dkb  \nab_\Rhat\Bbdot  + O(r^{-1})  \dk^{\le J-1}\Pdot+ O(r^{-1})  \dk^{\le J}\Bb \\
   &+ O(r^{-3}) \dk^{\le J}\Ga_b+  \dk^{\le J+1 } \Big( \Ga_b  \c \Rc_b\Big).
   \end{split}
\eeaa
    which  can be written in the form
     \eqref{eq:modelbainchipairequations12-simple} with
     \beaa
     \Psi_{(1)}=\Ptm, \qquad  \Psi_{(2)}=\Bbt, \qquad c_{(1)}=\frac 3 2 , \qquad c_{(2)}= 1,
     \eeaa
     where the signature of  $\Ptm$ is equal to $-J+2$ and that of  $\Bbt$ is equal  to $-J+1$. 
      This corresponds to  the  case  when  $2k-1=2(-J+2) -1 = -2J+3    <0$.  To apply   the integral estimate  \eqref{eq:general-proposition-second-bianchi-pairs}, we  need $\La_{(2)}= -2 c_{(2)} +1+\frac b 2 <0$  to be  satisfied  which holds true for  the choice $b=-\de$.   Therefore    
  \bea
  \lab{Equation:EstimateBt1}
\begin{split}
& \int_{\MM } r^{b-1} \left(1+|2J-3|\frac{m}{r}\right)      | \Bbt |^2 +  \int_{\pr^+\MM } r^{ b-2 } |\Bbt |^2\\
 \les & \int_{\MM} r^{b-1} |\Ptm|^2   +\big|I_1 \big|    +\big|I_2 \big|     +\int_{\MM}r^{b}\Big(  \big| F_{(1,J-2)}\big|\, \big|\Ptm\big|+ \big| F'_{(2, J-2)}\big|\, \big|\Bbt\big|\Big)+ \ep_0^2,
 \end{split}
\eea
where
\beaa
I_1&=& \int_{\MM}|q|^{b}\Re\Big(  4 \om \Bbt    \big) \c\ov{\Bbt}\Big) =\int_{\MM}   4 |q|^b\om\big|\Bbt\big|^2,\\
I_2&=& \int_{\MM}|q|^{b}\Re\Big(   O(r^{-3}) \dk^{\le J}(\nab_4, \nab)\Ga_b\c  \ov{\Bbt}\Big).
\eeaa
Since $\om =O (m r^{-2})$, we may absorb the term $I_1$ from the LHS for a sufficiently large\footnote{Recall that the iteration assumption \eqref{eq:iterationassumptiondiscussionThM8:bis} holds for $J\geq \frac{\kl}{2}$ and that $\kl$ is chosen large enough.} choice of $J$. We deduce
\beaa
 \int_{\MM } r^{b-1}   | \Bbt |^2 +  \int_{\pr^+\MM } r^{ b-2 } |\Bbt |^2 &\les &\int_{\MM} r^{b-1} |\Ptm|^2   + \left(\int_{\MM} r^{b-1} |\Ptm|^2\right)^{\frac{1}{2}}\left(\int_{\MM } r^{b+1} | F_{(1,J-2)}  |^2\right)^{\frac{1}{2}}\\
 && + \int_{\MM } r^{b+1} |  F_{(2,J-2)}' |^2+\big|I_2 \big| + \ep_0^2
\eeaa
and hence
\beaa
 \int_{\MM } r^{b-1}   | \Bbt |^2 +  \int_{\pr^+\MM } r^{ b-2 } |\Bbt |^2 &\les &  \big|I_2 \big| + \de_{J+1}[\Pc]   + \sqrt{\de_{J+1}[\Pc]}\left(\int_{\MM } r^{b+1} | F_{(1,J-2)}  |^2\right)^{\frac{1}{2}}\\
 && + \int_{\MM } r^{b+1} |  F_{(2,J-2)}' |^2 + \ep_0^2.
\eeaa

We then estimate  the integrals $\int_{\MM } r^{b+1} | F_{(1,J-2)}  |^2 $ and  $\int_{\MM } r^{b+1} |  F'_{(2,J-2)} |^2$. Taking advantage of the structure of $\Ft_{(1,J-2)}$ and $F_{(2,J-2)}'$, proceeding as for the corresponding estimate in the proof of Proposition \ref{Prop:-EstimatesBtPtp3}, and using in particular  \eqref{eq:EstmatesforBb-1-strong} to control $B^{J-1}_\de[ r(\nab_4\Bb, \nab\Bb)] $, we obtain 
\beaa
&&\sqrt{\de_{J+1}[\Pc]}\left(\int_{\MM } r^{b+1} | F_{(1,J-2)}  |^2\right)^{\frac{1}{2}} + \int_{\MM } r^{b+1} |  F_{(2,J-2)}' |^2 \\
&\les&  \de_{J+1}[\Pc]  +\ep_J^2   +\ep_0^2 +\sqrt{ \de_{J+1}[\Pc]} \left(B^J_\de[\Bb]+\Sk_{J+1}^2\right)^{\frac{1}{2}}+O(a^2)B_\de[\nab^J_3\Bb].
\eeaa
We deduce
\beaa
 \int_{\MM } r^{b-1}   | \Bbt |^2 +  \int_{\pr^+\MM } r^{b-2} |\Bbt |^2 &\les &  \big|I_2 \big|   +\de_{J+1}[\Pc]  +\ep_J^2   +\ep_0^2 +\sqrt{ \de_{J+1}[\Pc]} \left(B^J_\de[\Bb]+\Sk_{J+1}^2\right)^{\frac{1}{2}}\\
 && +O(a^2)B_\de[\nab^J_3\Bb].
\eeaa
Then, the term $I_2$ can be integrated by parts twice,  as the corresponding estimate  in the proof of Proposition \ref{Prop:-EstimatesBtPtp3},  to obtain 
 \beaa
     |I_2|      &\les& \Sk_{J+1}\Big(\de_{J+1}[\Pc]  +\ep_J^2 +\ep_0^2 +O(a^2)B_\de[\nab^J_3\Bb]\Big)^{\frac{1}{2}}+\ep_J\Big(B^J_\de[\Bb]\Big)^{\frac{1}{2}}+\ep_0^2.
     \eeaa
We infer
 \beaa
 \int_{\MM } r^{b-1}   | \Bbt |^2 +  \int_{\pr^+\MM } r^{b-2} |\Bbt |^2 &\les & \de_{J+1}[\Pc]  +\ep_J^2   +\ep_0^2 +\sqrt{ \de_{J+1}[\Pc]} \left(B^J_\de[\Bb]+\Sk_{J+1}^2\right)^{\frac{1}{2}}\\
 && +O(a^2)B_\de[\nab^J_3\Bb] +\ep_J\Big(B^J_\de[\Bb]\Big)^{\frac{1}{2}}\\
 && +\Sk_{J+1}\Big(\de_{J+1}[\Pc]  +\ep_J^2 +\ep_0^2 +O(a^2)B_\de[\nab^J_3\Bb]\Big)^{\frac{1}{2}}
\eeaa
as stated in \eqref{eq:lemma-EstimatesBbtPtm3}. This concludes the proof of Proposition \ref{Prop:lemma-EstimatesBbtPtm3}.
   \end{proof}

      
   \subsection{Proof of the estimates for $B$ in Proposition \ref{proposition:EstimatesBBb-interior}}
      \lab{section:EstimatesforBfromBt}
   

 First, we have in view of  Lemma \ref{Lemma:comparisonofnormsBEF-Rk}
\beaa
 \BEF^{J-1}_\de[r^2B]\les \Rk_J\Rk_{J+1}
 \eeaa
 and from \eqref{eq:BianchB-steps1-2-B} 
  \beaa
\bsplit
\BEF^{J-1}_\de[ r^2( \nab_3B, r \nab B) ]\les&   \de_{J+1}[\Pc] +\Sk_J \Sk_{J+1} +\Rk_{J} \Rk_{J+1}  +\ep_0^2\\
& +O(a^2, \ep^2)\BEF_\de[  r^{J +1 }  \nab^J_4 B]. 
 \end{split}
\eeaa  
Also, since 
  \beaa
  2\Rhat=e_4+O(1)e_3+O(1)\nab, \qquad 2\T=e_4+O(1)e_3+O(1)\nab, 
  \eeaa
  and since $\Bt =(\ov{q}\nabc_4 )^{\le J-2}  \nabc^2_\Rhat\Bdot$, we have
  \beaa
  \BEF_{\de}[r^{J-1}\nab_4^JB] &\les& \int_{\MM}r^{\de+1}|\Bt|^2+\sup_{\tau\leq \tau_*}\int_{\pr^+\MM(1,\tau)}r^{\de+2}|\Bt|^2+ \BEF^{J-1}_\de[r^2(\nab_3, r\nab)B]\\
  &&+\BEF^{J-1}_\de[r^2B].
  \eeaa
  Together with the above bounds for $\BEF^{J-1}_\de[r^2B]$ and $\BEF^{J-1}_\de[ r^2( \nab_3B, r \nab B)]$, and using the control of $\Bt$ in \eqref{eq:lemma-EstimatesBtPtp3}, i.e.
  \beaa
   \bsplit
  \int_{\MM } r^{1+\de} |\Bt |^2+\sup_{\tau\leq \tau_*}\int_{\pr^+\MM(1,\tau)} r^{2+\de}|\Bt |^2 &\les \de_{J+1}[\Pc]  +\ep_J^2   +\ep_0^2 +\sqrt{ \de_{J+1}[\Pc]} \left(B^J_\de[r^2B]+\Sk_{J+1}^2\right)^{\frac{1}{2}}\\
   &+O(a^2)B_\de[  r^{J+1}  \nab^J_4 B] +\ep_J\Big(B^J_\de[ r^2B]\Big)^{\frac{1}{2}}\\
   & +\Sk_{J+1}\Big(\de_{J+1}[\Pc]  +\ep_J^2 +\ep_0^2 +O(a^2)B_\de[  r^{J+1}  \nab^J_4 B]\Big)^{\frac{1}{2}},
   \end{split}
  \eeaa
  we infer
  \bea\lab{eq:intermediaryalmostdoneesitmateforBaisdhfalif:chap15}
\nn&&  \BEF^{J-1}_\de[ r^2( \nab_3B, r \nab B)]+\BEF_{\de}[r^{J-1}\nab_4^JB]+ \BEF^{J-1}_\de[r^2B]\\
\nn&\les& \de_{J+1}[\Pc] +\Sk_J \Sk_{J+1} +\Rk_{J} \Rk_{J+1}  +\ep_0^2 +O(a^2, \ep^2)\BEF_\de[  r^{J +1 }  \nab^J_4 B]\\
\nn&& +\sqrt{ \de_{J+1}[\Pc]} \left(B^J_\de[r^2B]+\Sk_{J+1}^2\right)^{\frac{1}{2}} +\ep_J\Big(B^J_\de[ r^2B]\Big)^{\frac{1}{2}}\\
   && +\Sk_{J+1}\Big(\de_{J+1}[\Pc]  +\ep_J^2 +\ep_0^2 +O(a^2)B_\de[  r^{J+1}  \nab^J_4 B]\Big)^{\frac{1}{2}}.
  \eea  
  
  Note that $\BEF_{\de}[r^{J-1}\nab_4^JB]$ appearing on the LHS of \eqref{eq:intermediaryalmostdoneesitmateforBaisdhfalif:chap15} is not consistent in terms of powers of $r$ compared to the other terms of the LHS. We thus need to upgrade this estimate. To this end, we introduce a smooth  function $\chi_{far}(r)$ such that $\chi_{far}(r)=0$ on $\Mtrap$ and  $\chi_{far}(r)=1$ for $r\geq 5m$, and we derive the following bound 
  \bea\lab{eq:intermediaryalmostdoneesitmateforBaisdhfalif:chap15bis}
  \bsplit
    &\int_{\MM} r^{1+\de}|\chi_{far}(\ov{q}\nab_4)^{J+1}B|^2+\sup_{\tau\leq \tau_*}\int_{\pr^+\MM(1,\tau)} r^{2+\de}|\chi_{far}(\ov{q}\nab_4)^{J+1}B|^2 \\
    \les& \de_{J+1}[\Pc] +\Sk_J \Sk_{J+1} +\Rk_{J} \Rk_{J+1}  +\ep_0^2 +O(a^2, \ep^2)\BEF_\de[  r^{J +1 }  \nab^J_4 B]\\
& +\sqrt{ \de_{J+1}[\Pc]} \left(B^J_\de[r^2B]+\Sk_{J+1}^2\right)^{\frac{1}{2}} +\ep_J\Big(B^J_\de[ r^2B]\Big)^{\frac{1}{2}}\\
   & +\Sk_{J+1}\Big(\de_{J+1}[\Pc]  +\ep_J^2 +\ep_0^2 +O(a^2)B_\de[  r^{J+1}  \nab^J_4 B]\Big)^{\frac{1}{2}}.
    \end{split}
  \eea
   Assuming \eqref{eq:intermediaryalmostdoneesitmateforBaisdhfalif:chap15bis},  and since $\chi_{far}(r)=1$ for $r\geq 5m$, we have
  \beaa
  \BEF^J_\de[r^2B] &\les& \BEF^{J-1}_\de[r^2(\nab_3, r\nab)B]+ \BEF_\de[r^{J+2}\nab_4^JB]+\BEF^{J-1}_\de[r^2B]\\
  &\les& \BEF^{J-1}_\de[ r^2( \nab_3B, r \nab B)]+\BEF_{\de}[r^{J-1}\nab_4^JB]+ \BEF^{J-1}_\de[r^2B]\\
  &&+\int_{\MM} r^{1+\de}|\chi_{far}(\ov{q}\nab_4)^{J+1}B|^2+\sup_{\tau\leq \tau_*}\int_{\pr^+\MM(1,\tau)} r^{2+\de}|\chi_{far}(\ov{q}\nab_4)^{J+1}B|^2
  \eeaa
  which together with \eqref{eq:intermediaryalmostdoneesitmateforBaisdhfalif:chap15} and \eqref{eq:intermediaryalmostdoneesitmateforBaisdhfalif:chap15bis} yields
  \beaa
  \BEF^J_\de[r^2B]   &\les& \de_{J+1}[\Pc] +\Sk_J \Sk_{J+1} +\Rk_{J} \Rk_{J+1}  +\ep_0^2 +O(a^2, \ep^2)\BEF_\de[  r^{J +1 }  \nab^J_4 B]\\
\nn&& +\sqrt{ \de_{J+1}[\Pc]} \left(B^J_\de[r^2B]+\Sk_{J+1}^2\right)^{\frac{1}{2}} +\ep_J\Big(B^J_\de[ r^2B]\Big)^{\frac{1}{2}}\\
   && +\Sk_{J+1}\Big(\de_{J+1}[\Pc]  +\ep_J^2 +\ep_0^2 +O(a^2)B_\de[  r^{J+1}  \nab^J_4 B]\Big)^{\frac{1}{2}}.
  \eeaa 
 For $a$ and $\ep$ small enough, we infer
  \beaa
  \BEF^J_\de[r^2B]   &\les& \de_{J+1}[\Pc]  +\ep_0^2+\ep_J^2 +\ep_J\Rk_{J+1} +\left(\sqrt{ \de_{J+1}[\Pc]} +\ep_0+\ep_J\right)\Sk_{J+1}+|a|\Sk_{J+1}^2
  \eeaa 
  which is the stated estimate for $\BEF^J_\de[r^2B]$ in \eqref{eq:mainestimateBBb-M8}.

It thus only remains  to derive the  estimate \eqref{eq:intermediaryalmostdoneesitmateforBaisdhfalif:chap15bis} relying in particular on  \eqref{eq:intermediaryalmostdoneesitmateforBaisdhfalif:chap15}.   We sketch below the main steps.
    
{\bf Step 1. }  We commute  the second Bianchi  pair, see  \eqref{eq:second-pair-B-P}, 
\beaa
\begin{split}
\nabc_3B  +\tr\Xb B &=-\DDs_1\,  \ov{P}+3\ov{P}H  +r^{-2} \Ga_b\c \Rc_b, \\
\nabc_4\ov{P}   +\frac{3}{2}\ov{\tr X}\,  \ov{P}   &=  \DDd_1 B + O( a r^{-2}) \ov{B}  +r^{-2} \Ga_b\c \Rc_b.
\end{split}
\eeaa
 with $\ov{q}\nabc_4$ and then    linearize the quantity $\nab_4\ov{P}$ by subtracting its Kerr value\footnote{It is crucial to first commute with $\ov{q}\nabc_4$ and then linearize the quantity $\nab_4\ov{P}$. Indeed, linearizing first $\ov{P}$ and then commuting with $\ov{q}\nabc_4$ would lead to a dangerous term $\dk^{\leq J+1}\Ga_b$ in $F_{(1)}$.}, i.e.
\beaa
\widecheck{\nab_4\ov{P}}:=\nab_4\ov{P}-\frac{\De}{|q|^2}\frac{6m}{\ov{q}^4}, \qquad \widecheck{\nab_4\ov{P}}=\nab_4\Pc+O(r^{-4})\Ga_g.
\eeaa  
 We then commute with $\chi_{far}(\ov{q}\nabc_4)^J$ where we recall that $\chi_{far}(r)=0$ on $\Mtrap$ and  $\chi_{far}(r)=1$ for $r\geq 5m$. Setting  $\Psi_{(1)}=\chi_{far}(\ov{q}\nabc_4)^{J+1}B$, $\Psi_{(2)}=\chi_{far}(\ov{q}\nabc_4)^J(\ov{q}\widecheck{\nab_4\ov{P}})$, we deduce, in view of Lemma \ref{lemma:CommBianchi-int1}, 
  \bea
\begin{split}
\lab{eq:third-pair-e_4P-e_4B-new}
\nabc_3 \Psi_{(1) }+ \left(1-\frac{J+1}{2}\right)\tr\Xb \Psi_{(1)} &= - \DDs_1\,\Psi_{(2)}  +F_{(1)}, \\
\nabc_4\Psi_{(2)}+ \left(\frac{3}{2} -\frac{J+1}{2}\right)\ov{\tr X}\,\Psi_{(2)}  &=\DDd_1\Psi_{(1)}+F_{(2)},
\end{split}
\eea
with $F_{(1)}, F_{(2)}$  of the form
\beaa
\bsplit
 F_{(1)} =& \chi_{far}(r)\Big[O(r^{-1} )  \dk^J\big(\nabc_3B, \dkb B\big)  + O(r^{-1})  \dk^{\leq J+1}\Pc +  O(r^{-3} )(\nab_3, \nab)\dk^{\le J}\Ga_b  \\
 &+  O(r^{-1})  \dk^{\le J}B+ O(r^{-3} )\dk^{\le J}  \Ga_b +r^{-2}\dk^{\le J+1}(\Ga_b\c \Rc_b)\Big]  + \chi_{far}'(r)\dk^{\leq J+1}(\Pc, B),\\
  F_{(2)} =&  \chi_{far}(r)\Big[O(r^{-1})+ \dk^{\leq J+1}(\Pc, B)+ O(r^{-3} )   \dk^{\le J+1}   \Ga_b +   r^{-2}\dk^{\le J+1}(\Ga_b\c \Rc_b)\Big] \\
  &+ \chi_{far}'(r)\dk^{\leq J+1}(\Pc, B).
\end{split}
\eeaa

{\bf Step  2.}   We apply  the  integral estimate  \eqref{eq:general-proposition-first-bianchi-pairs} of  Proposition \ref{Prop:Bainchi-pairsEstimates-integrated} to the system
      \eqref{eq:third-pair-e_4P-e_4B-new}.   Note that  the signature  of $\Psi_{(1)}$ is given by $J+1$, and that we have in this case $b=2+\de$ so that  $ \La_{(1)} =  -2+(J+1)+ 1+\frac b 2>0$. We deduce
\beaa
\begin{split}
& \int_{\MM } r^{b-1} |\Psi_{(1)}|^2+\int_{\pr^+\MM} r^{b}|\Psi_{(1)}|^2\\
\les &\int_{\MM} r^{b-1} |\Psi_{(2)}|^2  +\left|\int_{\MM}|q|^{b}\Re(F_{(1)}\c\ov{\Psi_{(1)}})\right|+\int_{\MM}r^b|F_{(2)}||\Psi_{(2)}| \\&  + \int_{\pr^-\MM} \big( r^{b}|\Psi_{(1)}|^2+  r^{b-2}|\Psi_{(2)}|^2\big).
\end{split}
\eeaa

  {\bf Step 3.}    We then proceed as in the proof of Proposition \ref{Prop:-EstimatesBtPtp3}, see section \ref{subsection:Estimatesfor-Bt}, with $(\Bt, \Ptp)$ being replaced by $(\Psi_{(1)}, \Psi_{(2)})$ where 
  \beaa
  \Psi_{(1)}=\chi_{far}(\ov{q}\nabc_4)^{J+1}B, \qquad \Psi_{(2)}=\chi_{far}(\ov{q}\nabc_4)^J(\ov{q}\widecheck{\nab_4\ov{P}}). 
  \eeaa
  The main differences are the fact that the argument is now simpler since $F_{(1)}$ and $F_{(2)}$ provided by Step 1 are supported away from $\Mtrap$, and the fact that the new term, generated by $\chi_{far}'(r)\dk^{\leq J+1}B$ in $F_{(1)}$, can be controlled using \eqref{eq:intermediaryalmostdoneesitmateforBaisdhfalif:chap15} since it is compactly supported in $r$. This finally leads to \eqref{eq:intermediaryalmostdoneesitmateforBaisdhfalif:chap15bis}, hence concluding the proof of the control of 
 $\BEF^J_\de[r^2B]$ in \eqref{eq:mainestimateBBb-M8}.

   
  \subsection{Proof of the estimates for $\Bb$ in Proposition \ref{proposition:EstimatesBBb-interior}}
\lab{section:InteriorEstimates-Bb}

 
  Using the estimates for $\Bbt $ derived in Proposition \ref{Prop:lemma-EstimatesBbtPtm3} and the estimates for $\nab_4 \Bb, \nab \Bb$ derived in Lemma \ref{Lemma:BianchB-steps1-2}, 
    we  derive estimates for $\nab_\Rhat  (\nab_3)^k \Bb$.  To  this end, we   make use of the following
    lemma.
    \begin{lemma}
    \lab{lemma:identityBbt-nab_3Bb}
    We have the identity    
    \beaa
     \left( \frac{\De}{|q|^2 } \right)^2 \nab_\Rhat \nab_3^{ J }\Bb = -2\Bbt  + O(1 )  \nab_\Rhat  \dk^{\le J- 1 }(\nab_4, \nab)\Bb + O(r^{-1} )\nab_{\Rhat}\dk^{\le J-1}\Bb+ O(r^{-1} )  \dk^{\le J-1}\Bb.
    \eeaa
      \end{lemma}
      
    \begin{proof}
    We write 
  \beaa
\Bbt&=&\nabc_3^{J-2}  \nabc^2_\Rhat\Bbdot=  \nabc^2_\Rhat \nabc_3 ^{J-2} \Lieb_\T \Bb + O(r^{-1} )\nab_{\Rhat}\dk^{\le J-1}\Bb+ O(r^{-1} )  \dk^{\le J-1}\Bb\\
&=&  \frac 1 2  \nab^2_\Rhat \nab_3^{J-2} \left( \nab_4+\frac{\De}{|q|^2} \nab _3 -2a\Re(\Jk)^b\nab_b\right)\Bb + O(r^{-1} )\nab_{\Rhat}\dk^{\le J-1}\Bb+ O(r^{-1} )  \dk^{\le J-1}\Bb\\
&=& \frac 1 2 \frac{\De}{|q|^2 } \nab^2_\Rhat \nab_3^{J-1} \Bb+ \nab^2_\Rhat \nab_3^{J-2}\nab _4\Bb+ O( ar^{-2} )  \nab^2_\Rhat \nab_3^{J-2} \nab \Bb \\
&& + O(r^{-1} )\nab_{\Rhat}\dk^{\le J-1}\Bb+ O(r^{-1} )  \dk^{\le J-1}\Bb\\
&=&  \frac 1 2  \frac{\De}{|q|^2 } \  \nab^2_\Rhat \nab_3^{J-1} \Bb+O(1)\nab_\Rhat  \dk^{\le J- 1 }(\nab_4, \nab) \Bb + O(r^{-1} )\nab_{\Rhat}\dk^{\le J-1}\Bb+ O(r^{-1} )  \dk^{\le J-1}\Bb.
\eeaa
Also,
\beaa
  \nab^2_\Rhat \nab_3^{J-1} \Bb&=& \frac 1 2 \nab_\Rhat \left(\nab_4 - \frac{\De}{|q|^2 }  \nab_3\right) \nab_3^{J-1} \Bb + O(r^{-1} )\nab_{\Rhat}\dk^{\le J-1}\Bb+ O(r^{-1} )  \dk^{\le J-1}\Bb\\
  &=&- \frac 1 2 \frac{\De}{|q|^2}  \nab_\Rhat  \nab_3 ^{J} \Bb+  \nab_\Rhat  \nab_3^{J-1}\nab_4\Bb + O(r^{-1} )\nab_{\Rhat}\dk^{\le J-1}\Bb+ O(r^{-1} )  \dk^{\le J-1}\Bb.
\eeaa
We deduce
\beaa
\Bbt =   -  \frac 1 2 \left( \frac{\De}{|q|^2 } \right)^2 \nab_\Rhat \nab_3^{J} \Bb+ O(1 )  \nab_\Rhat  \dk^{\le J-1}(\nab_4, \nab) \Bb + O(r^{-1} )\nab_{\Rhat}\dk^{\le J-1}\Bb+ O(r^{-1} )  \dk^{\le J-1}\Bb
\eeaa
as stated.
\end{proof}

  We then make use  of  the  integral estimate of Proposition  \ref{Prop:lemma-EstimatesBbtPtm3} for $\Bbt$,  the estimates  for $\nab_4\Bb, \nab\Bb$  already derived, see \eqref{eq:BianchB-steps1-2-Bb},  the induction hypothesis, and Lemma \ref{lemma:identityBbt-nab_3Bb},   to deduce, with $b= -\de $, 
    \bea
    \lab{estimate:degenerate-Bb1}
    \bsplit
   &  \int_{\MM } r^{b-1}  \left( \frac{\De}{|q|^2 } \right)^4       \big|\nab_\Rhat \nab_3 ^{ J }  \Bb\big|^2        +\sup_{\tau\leq\tau_*}\int_{\pr^+ \MM(1,\tau)}  r^{ b-2 }   \left( \frac{\De}{|q|^2 } \right)^4  \big|\nab_\Rhat \nab_3^{ J }  \Bb\big|^2  \\
   \les& \de_{J+1}[\Pc]  +\Sk_J \Sk_{J+1} +\Rk_{J} \Rk_{J+1} +\ep_0^2  +O(a^2, \ep^2)B_\de[\nab^J_3\Bb]\\
 &   +\sqrt{ \de_{J+1}[\Pc]} \left(B^J_\de[\Bb]+\Sk_{J+1}^2\right)^{\frac{1}{2}} +\ep_J\Big(B^J_\de[\Bb]\Big)^{\frac{1}{2}}\\
 & +\Sk_{J+1}\Big(\de_{J+1}[\Pc]  +\ep_J^2 +\ep_0^2 +O(a^2)B_\de[\nab^J_3\Bb]\Big)^{\frac{1}{2}}.
 \end{split}
    \eea
   Note that \eqref{estimate:degenerate-Bb1} is degenerate in the redshift region $r\leq r_+(1+\deh)$. To get rid of this degeneracy, we introduce a smooth cut-off function $\chi_{red}(r)$ such that $\chi_{red}(r)=1$ for $r\le r_+(1+\deh)$ and $\chi_{red}(r)=0$ for $r\geq r_+(1+2\deh)$, and we  derive in Steps 1--4 below the following bound
    \bea
    \lab{estimate:nondegenerate-Bb2}
    \bsplit
   &  \int_{\MM }\big|\chi_{red}\nab_3^{J+1}\Bb\big|^2        +\sup_{\tau\leq\tau_*}\int_{\pr^+ \MM(1,\tau)}\big|\chi_{red}\nab_3^{J+1}  \Bb\big|^2  \\
   \les& \de_{J+1}[\Pc]  +\Sk_J \Sk_{J+1} +\Rk_{J} \Rk_{J+1} +\ep_0^2  +O(a^2, \ep^2)B_\de[\nab^J_3\Bb]\\
 &   +\sqrt{ \de_{J+1}[\Pc]} \left(B^J_\de[\Bb]+\Sk_{J+1}^2\right)^{\frac{1}{2}} +\ep_J\Big(B^J_\de[\Bb]\Big)^{\frac{1}{2}}\\
 & +\Sk_{J+1}\Big(\de_{J+1}[\Pc]  +\ep_J^2 +\ep_0^2 +O(a^2)B_\de[\nab^J_3\Bb]\Big)^{\frac{1}{2}}.
 \end{split}
    \eea
   Since $\chi_{red}(r)=1$ for $r\geq r_+(1+\deh)$, we have
  \beaa
  \BEF^J_\de[\Bb] &\les& \BEF^{J-1}_\de[r(\nab_4, \nab)\Bb]+ \BEF_\de[\nab_3^J\Bb]+\BEF^{J-1}_\de[\Bb]\\
  &\les& \BEF^{J-1}_\de[r(\nab_4, \nab)\Bb]+\BEF^{J-1}_\de[\Bb]\\
  && +\int_{\MM } r^{-\de-1}  \left( \frac{\De}{|q|^2 } \right)^4       \big|\nab_\Rhat \nab_3 ^{ J }  \Bb\big|^2        +\sup_{\tau\leq\tau_*}\int_{\pr^+ \MM(1,\tau)}  r^{-\de-2}   \left( \frac{\De}{|q|^2 } \right)^4  \big|\nab_\Rhat \nab_3^{ J }  \Bb\big|^2\\
  &&+\int_{\MM }\big|\chi_{red}\nab_3^{J+1}\Bb\big|^2        +\sup_{\tau\leq\tau_*}\int_{\pr^+ \MM(1,\tau)}\big|\chi_{red}\nab_3^{J+1}  \Bb\big|^2
  \eeaa
  which together with \eqref{estimate:degenerate-Bb1}, \eqref{estimate:nondegenerate-Bb2}, and the estimates for $\nab_4 \Bb, \nab \Bb$ derived in Lemma \ref{Lemma:BianchB-steps1-2}, yields
  \beaa
 \bsplit
  \BEF^J_\de[\Bb] \les& \de_{J+1}[\Pc]  +\Sk_J \Sk_{J+1} +\Rk_{J} \Rk_{J+1} +\ep_0^2  +O(a^2, \ep^2)B_\de[\nab^J_3\Bb]\\
 &   +\sqrt{ \de_{J+1}[\Pc]} \left(B^J_\de[\Bb]+\Sk_{J+1}^2\right)^{\frac{1}{2}} +\ep_J\Big(B^J_\de[\Bb]\Big)^{\frac{1}{2}}\\
 & +\Sk_{J+1}\Big(\de_{J+1}[\Pc]  +\ep_J^2 +\ep_0^2 +O(a^2)B_\de[\nab^J_3\Bb]\Big)^{\frac{1}{2}}.
\end{split}
\eeaa  
   For $a$ and $\ep$ small enough, we infer
  \beaa
  \BEF^J_\de[\Bb]   &\les& \de_{J+1}[\Pc]  +\ep_0^2+\ep_J^2 +\ep_J\Rk_{J+1} +\left(\sqrt{ \de_{J+1}[\Pc]} +\ep_0+\ep_J\right)\Sk_{J+1}+|a|^2\Sk_{J+1}^2
  \eeaa 
  which is the stated estimate for $\BEF^J_\de[\Bb]$ in \eqref{eq:mainestimateBBb-M8}.
 
  \medskip
    
It thus only remains  to derive the non-degenerate estimate \eqref{estimate:nondegenerate-Bb2} by relying in particular on \eqref{estimate:degenerate-Bb1}.   We sketch below the main steps in the proof of this  red shift type argument.
    
    {\bf Step 1.}   Using the fact that  $\Rhat=e_4+\frac{\De}{|q|^2} e_3 $, we obtain, as a consequence of \eqref{estimate:degenerate-Bb1} and the estimates for  $(\nab_4, \nab)\Bb$ of Lemma  \ref{Lemma:BianchB-steps1-2},
\bea
\lab{estimate:degenerateBb}
\bsplit
   &  \int_{\MM(r\leq r_+(1+2\deh))}   \left( \frac{\De}{|q|^2 } \right)^6    \big|\nab_3 ^{J+1}  \Bb\big|^2   +   \int_{\pr\MM_+(r\leq r_+(1+2\deh))}   \left( \frac{\De}{|q|^2 } \right)^6     \big|\nab_3 ^{J+1}  \Bb\big|^2
   \\
    \les&  \de_{J+1}[\Pc]  +\Sk_J \Sk_{J+1} +\Rk_{J} \Rk_{J+1} +\ep_0^2  +O(a^2, \ep^2)B_\de[\nab^J_3\Bb]\\
 &   +\sqrt{ \de_{J+1}[\Pc]} \left(B^J_\de[\Bb]+\Sk_{J+1}^2\right)^{\frac{1}{2}} +\ep_J\Big(B^J_\de[\Bb]\Big)^{\frac{1}{2}}\\
 & +\Sk_{J+1}\Big(\de_{J+1}[\Pc]  +\ep_J^2 +\ep_0^2 +O(a^2)B_\de[\nab^J_3\Bb]\Big)^{\frac{1}{2}}.
   \end{split}
\eea  

{\bf Step 2. }  We commute  the second Bianchi  pair, see  \eqref{eq:third-pair-P-Bb}, 
     \beaa
    \begin{split}
\nabc_3P+ \frac{3}{2}\ov{\tr\Xb} P &= - \DDd_1\Bb +  O(ar^{-2} )\Bb    +r^{-1} \Ga_b\c \Rc_b, \\
\nabc_4\Bb+ \tr X\Bb  &=\DDs_1  P -3P\Hb  +r^{-1} \Ga_b\c \Rc_b,
\end{split}
\eeaa
 with  $\nabc_3 $ and then    linearize the quantity $\nab_3 P $ by subtracting its Kerr value\footnote{It is crucial to first commute with $\nabc_3$ and then linearize the quantity $\nab_3P$. Indeed, linearizing first $P$ and then commuting with $\nabc_3$ would lead to a dangerous term $\dk^{\leq J+1}\Ga_b$ in $F_{(2)}$.}, i.e.
\beaa
\widecheck{\nab_3P}:=\nab_3P +\frac{6m}{\ov{q}^4}, \qquad \widecheck{\nab_3P}=\nab_3\Pc+O(r^{-3})\Ga_b.
\eeaa
We then commute with $\chi_{red}\nabc_3 ^{ J }$ where   $\chi_{red} $ is a smooth  cut-off     function  equal to $1$ in the red shift region\footnote{Recall that the red shift region  $\MM_{red} $ is  defined by 
  $\MM(r\le  r_+( 1+\de_{red} ))$ with  $\de_{red}\geq \deh$.}  $\MM_{red}$ and  $0$    for $ r\ge  r_+( 1+2\de_{red} )$.
 Setting  $\Psi_{(1)}=\chi_{red}\nabc_3^{ J }\widecheck{ \nabc_3  P}$, $ \Psi_{(2)}=\chi_{red} \nabc_3^{ J+1 } \Bb$, we deduce, in view of Lemma \ref{lemma:CommBianchi-int2},
  \bea
\begin{split}
\lab{eq:third-pair-e_3P-e_3Bb-new}
\nabc_3 \Psi_{(1) }+ \frac{3}{2}\ov{\tr\Xb} \Psi_{(1)} &= - \DDd_1\Psi_{(2)}  +F_{(1)}, \\
\nabc_4\Psi_{(2)}+ \tr X\Psi_{(2)}  &=\DDs_1\, \Psi_{(1)}+F_{(2)},
\end{split}
\eea
with $F_{(1)}, F_{(2)}$ supported in the region  $ r\leq  r_+( 1+2\de_{red} )$ and of the form
\beaa
 F_{(1)} &=&  \chi_{red}(r)\Big[O(1) \dk^{\le J+1} (\Pc,  \Bb)  + O(1)\dk^{\le J+1}\Ga_b  +\dk^{\le J+1} \big( \Ga_b\c \Rc_b\big)\Big]\\
 &&+\chi_{red}'(r)\dk^{\leq J+1}(\Pc, \Bb),\\
  F_{(2)}&=&  \chi_{red}(r)\Big[O(1) \dk^{\le J}\nab B +O(1)\dk^{\le J}(\nab_4, \nab)\Ga_b+O(1) \dk^{\le J+1}\Pc\\
  && +O(1) \dk^{\le J}\Bb+O(1) \dk^{\le J}\Ga_b  +\dk^{\le J+1} \big( \Ga_b\c \Rc_b\big)\Big]+\chi_{red}'(r)\dk^{\leq J+1}(\Pc, \Bb).
\eeaa

     {\bf Step  3.}   We apply  the  integral estimate  \eqref{eq:general-proposition-second-bianchi-pairs} of  Proposition \ref{Prop:Bainchi-pairsEstimates-integrated} to the system
      \eqref{eq:third-pair-e_3P-e_3Bb-new}.   Note that  the signature  of $\Psi_{(1)}$ is given by $-1-J$, and that we have in this case $b= - \de$ so that  $ \La_{(2)} =  -2+ 1+\frac b 2<0$. We deduce
\beaa
\begin{split}
& \int_{\MM} r^{b-1}\left(1+|2k+1|\frac{m}{r}\right) | \Psi_{(2)}|^2 +  \int_{\pr^+\MM} r^{b-2} |\Psi_{(2)}|^2\\ 
\les &\int_{\MM} r^{b-1} |\Psi_{(1)}|^2    +\left|\int_{\MM}|q|^{b}\Re(F_{(2)}\c\ov{\Psi_{(2)}})\right|    +\int_{\MM}r^{b}|F_{(1)}||\Psi_{(1)}|  \\
&+  \int_{\pr^-\MM} \big( r^{b}|\Psi_{(1)}|^2+  r^{b-2}|\Psi_{(2)}|^2\big).
 \end{split}
\eeaa

 {\bf Step 4.}    We then proceed as in the proof of Proposition \ref{Prop:lemma-EstimatesBbtPtm3}, see section \ref{subsection:Estimatesfor-Bbt}, with $(\Ptm, \Bbt)$ being replaced by $(\Psi_{(1)}, \Psi_{(2)})$ where 
 \beaa
 \Psi_{(1)}=\chi_{red}\nabc_3^{ J }\widecheck{ \nabc_3  P}, \qquad \Psi_{(2)}=\chi_{red} \nabc_3^{ J+1 } \Bb.
 \eeaa
 The main differences are the fact that the argument is now simpler since $F_{(1)}$ and $F_{(2)}$ provided by Step 2 are supported away from $\Mtrap$, and the fact that the new term, generated by $\chi_{red}'(r)\dk^{\leq J+1}\Bb$ in $F_{(2)}$, can be controlled using \eqref{estimate:degenerateBb} and  the estimates for  $(\nab_4, \nab)\Bb$ of Lemma  \ref{Lemma:BianchB-steps1-2} since it is supported in $r_+(1+\deh)\leq r\leq r_+(1+2\deh)$. This finally leads to \eqref{estimate:nondegenerate-Bb2}, hence concluding the proof of the control of 
 $\BEF^J_\de[\Bb]$ in \eqref{eq:mainestimateBBb-M8}.

   
  \section{Estimates for $A$ and $\Ab$}
  \lab{section:InteriorEstimatesAAb}
  

The goal of this section is to prove Proposition \ref{proposition:EstimatesAAb-interior} providing energy-Morawetz estimates for $(A, \Ab)$ assuming corresponding energy-Morawetz estimates for $(B, \Bb)$.


\subsection{Estimates for $\nab_3 A$, $\nab A$, $\nab_4\Ab$, $\nab\Ab$}
\label{section:est-nab3-nab-AAb}


We  derive the following lemma. 
\begin{lemma}
\lab{Lemma:BianchA-steps1-2} 
The following estimates hold true in $\MM$:
\begin{enumerate}
\item We have
\bea
\lab{eq:BianchA-steps1-2-A}
\bsplit
\BEF^{J-1}_\de[ r^2( \nab_3A, r \nab A) ]&\les   \de_{J+1}[B] +\ep_0^2+\ep_J^2+ O(a^2, \ep^2) \BEF_\de[  r^{J + 1} \nab^J_4 A]. 
 \end{split}
\eea

\item We have
\bea
\lab{eq:BianchA-steps1-2-Ab}
\bsplit
\BEF^{J-1}_\de[  r (\nab_4\Ab, \nab \Ab) ] &\les 
 \de_{J+1}[\Bb] +\ep_0^2+\ep_J^2 + O(a^2, \ep^2) \BEF_\de[\nab^J_3\Ab].
 \end{split}
\eea
\end{enumerate}
\end{lemma} 

\begin{proof}
To prove \eqref{eq:BianchA-steps1-2-Ab}, we rely on the fourth Bianchi pair   \eqref{eq:fourth-pair-Bb-Ab-lin}, i.e. 
\beaa
\nabc_3\Bb  +2\ov{\tr\Xb}\,\Bb &=&  -\frac{1}{2}\ov{\DD}\c\Ab  +O(a r^{-2})\Ab -O(r^{-3})  \,\Xib +  \Ga_b\c \Rc_b,\\
\nabc_4\Ab  +\frac{1}{2}\tr X \Ab   &=&-\frac 1 2 \DD \hot\Bb  +O(ar^{-2} ) \Bb +O(r^{-3} )\Xbh +\Ga_b\c \Rc_b.
\eeaa
From  second  equation we deduce
\beaa
\BEF^{J-1}_\de[ r \nab_4 \Ab]  &\les&\BEF_\de ^J[\Bb]+\ep_J^2+\ep_0^2,
\eeaa
while from the first equation  we deduce
\beaa
\BEF^{J-1}_\de[ r \ov{\DD}\c\Ab]  &\les & \BEF^{J}_\de[\Bb]  +\ep_J^2+\ep_0^2.
\eeaa
 Also,  making use of  the Hodge type estimates  of  Corollary \ref{Cor:HodgeThmM8},  we have
  \beaa
  \BEF^{J-1}_\de[ r \nab\Ab]  &\les& \BEF^{J-1}_\de[ r \ov{\DD}\c\Ab]   + O(a^2, \ep^2)\BEF^{J-1}_\de[( 
  \nab_3\Ab, \nab_4\Ab)].
  \eeaa
  Combining the above estimates, we infer, for $a$ and $\ep$ small enough,
  \beaa
  \BEF^{J-1}_\de[r \nab_4 \Ab] +   \BEF^{J-1}_\de[r \nab\Ab] &\les & \BEF^{J}_\de[   \Bb ]  +\ep_J^2 +\ep_0^2 + O(a^2,\ep^2) \BEF_\de[\nab^J_3 \Ab]
  \eeaa
   as desired. 
   
   The    estimates \eqref{eq:BianchA-steps1-2-A} for $A$ are derived in the same manner, by relying on the first Bianchi pair \eqref{eq:first-pair-A-B-lin}. This concludes the proof of Lemma \ref{Lemma:BianchA-steps1-2}.
   \end{proof}
   
   In view of Lemma \ref{Lemma:BianchA-steps1-2}, it remains to estimate  the top $\nab_4$ derivatives of $A$ and the top $\nab_3$ derivatives 
   of $\Ab$.  We thus introduce the quantities
   \bea
   \bsplit
   \Bbt&:= \nabc^{J-1} _3  \nabc_{\Rhat}^2 \Bb,\qquad\quad  \Abt:=\nabc_3^{J-1}  \nabc_{\Rhat}^2 \Ab,\\
    \quad  \Bt&:= (\ov{q}\nabc_4)^{J-1}\nabc_{\Rhat}^2 B, \quad\,\, \At:= (\ov{q}\nabc_4)^{J-1}   \nabc_{\Rhat}^2A.
    \end{split}
   \eea

   
\subsection{Estimates for $\Abt$}


   We first derive   equations for the Bianchi pair  $\Abt, \Bbt$. 
    Commuting  the  Bianchi pair equation \eqref{eq:fourth-pair-Bb-Ab-lin}   with  $\nabc^{J-1} _3  \nabc_{\Rhat}^2$, 
     we  derive the following   lemma,  analogous to the second part of Lemma \ref{lemma:EstimatesBtPtp}.
     \begin{lemma}\lab{Lemma:linearization-tildes-Ab-Bb}
      We have
     \bea
    \lab{eq:linearization-tildes-Ab-Bb}   
    \bsplit
\nabc_3  \Bbt+    2\ov{\tr\Xb}\, \Bbt  &=- \DDd_2 \Abt +\Ft_{(1)},\\
    \nabc_4 \Abt +\frac 1 2 \tr X  \Abt &=\DDs_2\, \Bbt  +\Ft_{(2)}, 
    \end{split}
   \eea
   with
   \bea
   \lab{eq:linearization-tildes-Ab-Bb-Ft}
   \bsplit
   \Ft_{(1)}&=O(r^{-1})
   \dk^{J} \nab_\Rhat ( \Bb, \Ab )+ O(r^{-3} )   \dk^{\le J+1}   \Ga_b + O(r^{-2})  \dk^{\le J}( \Bb, \Ab) +   \dk^{\le J+1}(\Ga_b\c \Rc_b),\\
   \Ft_{(2)}&=  - 4\om \Abt +O(r^{- 1 }) \nab_\Rhat \dk^{J-1}\big(\nab_4 \Ab, \nab \Ab)+ O(r^{- 1 }) \nab_\Rhat\dk^{\leq J} \Bb + O(r^{-3} ) \dk^{J}\nab\Ga_b
  \\
  & +O(r^{- 3  } ) \dk^{\le J}\Ga_b+   O(r^{- 1 }) \dk^{\le J}  (\Bb, \Ab) +   \dk^{\le J+1}(\Ga_b\c \Rc_b).
  \end{split}
   \eea
      \end{lemma}

   \begin{proof}
  We apply  the second part   of  Proposition \ref{Prop:BianchPairs-tilde3-k} to the fourth Bianchi pair \eqref{eq:fourth-pair-Bb-Ab-lin}, i.e. to
   \beaa
\begin{split}
\nabc_3\Bb  +2\ov{\tr\Xb}\,\Bb &=  -\DDd_2\Ab  +O(a r^{-2})\Ab -O(r^{-3})  \,\Xib + \Ga_b\c \Rc_b,\\
\nabc_4\Ab  +\frac{1}{2}\tr X \Ab   &=\DDs_2\, \Bb  +O(ar^{-2} ) \Bb +O(r^{-3} )\Xbh + \Ga_b\c \Rc_b.
\end{split}
\eeaa
This yields  the system \eqref{eq:linearization-tildes-Ab-Bb} for $(\Bbt, \Abt)$ 
with $\Ft_{(1)} $ as in \eqref{eq:linearization-tildes-Ab-Bb-Ft}
and $\Ft_{(2)}$  given by 
\beaa
\Ft_{(2)} &=&  - 4\om \Abt +O(r^{- 1 }) \nab_\Rhat \dk^{J-1}\big(\nab_4 \Ab, \nab \Ab)+ O(r^{- 1 }) \nab_\Rhat\dk^{\leq J} \Bb + O(r^{-3} ) \dk^{J} \nab_3 \Xbh
  \\
  && +O(r^{- 3  } ) \dk^{\le J}\Ga_b+   O(r^{- 1 }) \dk^{\le J}  (\Bb, \Ab) +   \dk^{\le J+1}(\Ga_b\c \Rc_b).
\eeaa
Finally,  we make use of  the following consequence of the null structure equations, see \eqref{eqts:usefulNullstructure},
\beaa
\nabc_3\Xbh+\Re(\tr\Xb) \Xbh&=& \frac 1 2  \DDc\hot \Xib+  O(ar^{-2})\Ga_b -\Ab +\Ga_b\c\Ga_b
\eeaa
and obtain  the desired expression for $F_{(2)}$.
\end{proof}

Based  on the  Bianchi pair equation in Lemma \ref{Lemma:linearization-tildes-Ab-Bb} for $(\Abt, \Bbt)$,  we  derive the following estimate for $\Abt$.
\begin{lemma}
  \lab{lemma:lemma-EstimatesabtBbtm3}
      The following   estimates hold true for $\Abt$, with $b= -\de $, 
   \bea
   \lab{eq:-EstimatesabtAbBb3}
\begin{split}
 &\int_{\MM } r^{b-1} |\Abt|^2+\int_{\pr^+\MM} r^{ b-2 }|\Abt|^2\\
 \les& \de_{J+1}[\Bb] +\ep_0^2+\ep_J^2 + O(a^2, \ep^2)\BEF_\de[\nab^J_3\Ab] +\sqrt{\de_{J+1}[\Bb]}\Big(B^{J}_\de[\Ab]+\Sk_{J+1}^2\Big)^{\frac{1}{2}}\\
&+\ep_J\Big(B^J_\de[\Ab]\Big)^{\frac{1}{2}}+\Sk_{J+1}\Big(\de_{J+1}[\Bb] +\ep_0^2+\ep_J^2 + O(a^2, \ep^2)\BEF_\de[\nab^J_3\Ab]\Big)^{\frac{1}{2}}.
 \end{split}
\eea
   \end{lemma}
   
   \begin{proof}
   We proceed as in the proof of   \eqref{eq:lemma-EstimatesBbtPtm3} starting with 
   the Bianchi  pair \eqref{eq:linearization-tildes-Ab-Bb},   
      which  can be written in the form   \eqref{eq:modelbainchipairequations12-simple} with 
     \beaa
     \Psi_{(1)}=\Bbt, \qquad  \Psi_{(2)}=\Abt, \qquad c_{(1)}= 2, \qquad  c_{(2)}= \frac 1 2, 
     \eeaa
     and $ F_{(1)}=\Ft_{(1)},  F_{(2)}= \Ft_{(2)}$.
    Note that   the signature of  $\Psi_{(1)} =\Bbt $  is  $k=-J$ and that  of  $\Psi_{(2)}=\Abt$ equal  to $k-1=-J-1$.
      This corresponds to  the  case  when  $2k-1<0$  in Proposition \ref{Prop:Bainchi-pairsEstimates-integrated}.  To apply   the integral estimate  \eqref{eq:general-proposition-second-bianchi-pairs}, we  need $\La_{(2)}= -2 c_{(2)} +1+\frac b 2 <0$  to be  satisfied  which holds true for  the choice $b= -\de $.   Therefore    
  \beaa
\begin{split}
& \int_{\MM } r^{b-1}\left(1+|2J-1|\frac{m}{r} \right)  | \Abt |^2 +  \int_{\pr^+\MM } r^{ b-2 } |\Abt |^2\\
\les& \int_{\MM} r^{b-1} |\Bbt|^2    +\int_{\MM}|\om||\Abt|^2 + \left|\int_\MM O(r^{-3}) \Re\Big(\Abt\c \ov{ \nab\dk^{\le J}\Ga_b}\Big)\right| \\
&+\int_{\MM}r^{b}\Big(  \big| F_{(1)}\big|\, \big|\Bbt\big|+ \big| F'_{(2)}\big|\, \big|\Abt\big|\Big)+ \ep_0^2,
 \end{split}
\eeaa
where 
\beaa
F_{(2)} &=& F'_{(2)} -4\om \Abt +  O(r^{-3} )\nab\dk^{J}\Ga_b,\\
F'_{(2)} &:=& O(r^{- 1 }) \nab_\Rhat \dk^{J-1}\big(\nab_4 \Ab, \nab \Ab)+ O(r^{- 1 }) \nab_\Rhat\dk^{\leq J} \Bb +O(r^{- 3  } ) \dk^{\le J}\Ga_b\\
&&+   O(r^{- 1 }) \dk^{\le J}  (\Bb, \Ab) +   \dk^{\le J+1}(\Ga_b\c \Rc_b).
\eeaa
Since $\om=O(mr^{-2})$, the term $\om|\Abt|^2$ can be absorbed from the LHS for $J$ large enough (recalling that $J\geq\frac{\kl}{2}\ll 1$), and we infer, using also Cauchy Schwartz, 
  \beaa
\begin{split}
& \int_{\MM } r^{b-1}| \Abt |^2 +  \int_{\pr^+\MM } r^{ b-2 } |\Abt |^2\\
\les& \de_{J+1}[\Bb]    + \left|\int_\MM O(r^{-3}) \Re\Big(\Abt\c \ov{ \nab\dk^{\le J}\Ga_b}\Big)\right| +\sqrt{\de_{J+1}[\Bb]}\left(\int_{\MM}r^{b+1}|F_{(1)}|^2\right)^{\frac{1}{2}}\\
&+ \int_{\MM}r^{b+1}|F_{(2)}'|^2+ \ep_0^2.
 \end{split}
\eeaa
Since
\beaa
&&\sqrt{\de_{J+1}[\Bb]}\left(\int_{\MM}r^{b+1}|F_{(1)}|^2\right)^{\frac{1}{2}} + \int_{\MM}r^{b+1}|F_{(2)}'|^2\\
&\les& \de_{J+1}[\Bb]+B^{J-1}_\de[r(\nab_4, \nab)\Ab]+\ep_J^2+\ep_0^2+\sqrt{\de_{J+1}[\Bb]}\Big(B^{J}_\de[\Ab]+\Sk_{J+1}^2\Big)^{\frac{1}{2}}, 
\eeaa
we deduce
  \beaa
\begin{split}
& \int_{\MM } r^{b-1}| \Abt |^2 +  \int_{\pr^+\MM } r^{ b-2 } |\Abt |^2\\
\les&  \left|\int_\MM O(r^{-3}) \Re\Big(\Abt\c \ov{ \nab\dk^{\le J}\Ga_b}\Big)\right| + \de_{J+1}[\Bb]+B^{J-1}_\de[r(\nab_4, \nab)\Ab]+\ep_J^2+\ep_0^2\\
&+\sqrt{\de_{J+1}[\Bb]}\Big(B^{J}_\de[\Ab]+\Sk_{J+1}^2\Big)^{\frac{1}{2}}.
 \end{split}
\eeaa
Also, by integration by parts, proceeding as for the control of the term $I_2$ in the proof of Proposition \ref{Prop:lemma-EstimatesBbtPtm3}, see section \ref{subsection:Estimatesfor-Bbt}, we have
\beaa
 \left|\int_\MM O(r^{-3}) \Re\Big(\Abt\c \ov{ \nab\dk^{\le J}\Ga_b}\Big)\right|
 &\les& \Sk_{J+1}\Big(\BEF^{J-1}_\de[r(\nab_4, \nab)\Ab]\Big)^{\frac{1}{2}}+\ep_J\Big(B^J_\de[\Ab]\Big)^{\frac{1}{2}}+\ep_0^2.
 \eeaa
We deduce
\beaa
\begin{split}
& \int_{\MM } r^{b-1}| \Abt |^2 +  \int_{\pr^+\MM } r^{b-2} |\Abt |^2\\
\les&   \de_{J+1}[\Bb]+B^{J-1}_\de[r(\nab_4, \nab)\Ab]+\ep_J^2+\ep_0^2 +\sqrt{\de_{J+1}[\Bb]}\Big(B^{J}_\de[\Ab]+\Sk_{J+1}^2\Big)^{\frac{1}{2}}+\ep_J\Big(B^J_\de[\Ab]\Big)^{\frac{1}{2}}\\
&+\Sk_{J+1}\Big(\BEF^{J-1}_\de[r(\nab_4, \nab)\Ab]\Big)^{\frac{1}{2}}.
 \end{split}
\eeaa
Together with the control of $\BEF^{J-1}_\de[r(\nab_4, \nab)\Ab]$ provided by \eqref{eq:BianchA-steps1-2-Ab}, this implies 
\beaa
\begin{split}
& \int_{\MM } r^{b-1}| \Abt |^2 +  \int_{\pr^+\MM } r^{b-2} |\Abt |^2\\
\les& \de_{J+1}[\Bb] +\ep_0^2+\ep_J^2 +  O(a^2, \ep^2)  \BEF_\de[\nab^J_3\Ab] +\sqrt{\de_{J+1}[\Bb]}\Big(B^{J}_\de[\Ab]+\Sk_{J+1}^2\Big)^{\frac{1}{2}}\\
&+\ep_J\Big(B^J_\de[\Ab]\Big)^{\frac{1}{2}}+\Sk_{J+1}\Big(\de_{J+1}[\Bb] +\ep_0^2+\ep_J^2 +  O(a^2, \ep^2)  \BEF_\de[\nab^J_3\Ab]\Big)^{\frac{1}{2}}.
 \end{split}
\eeaa
as stated in \eqref{eq:-EstimatesabtAbBb3}. This concludes the proof of Lemma \ref{lemma:lemma-EstimatesabtBbtm3}.
\end{proof}

   
   \subsection{Estimates for $\Ab$ in Proposition \ref{proposition:EstimatesAAb-interior}}
    \lab{section:InteriorEstimates-Ab}
    
   
  We proceed as in section \ref{section:InteriorEstimates-Bb}.  
  Using the estimates for $\Abt $ derived in Lemma \ref{lemma:lemma-EstimatesabtBbtm3} and the estimates for $\nab_4 \Ab, \nab \Ab$ derived in Lemma \ref{Lemma:BianchA-steps1-2},   we  derive estimates for $\nab_\Rhat  (\nab_3)^k \Ab$. To this end, we make use of the following identity
     \beaa
     \left( \frac{\De}{|q|^2 } \right)\nab_\Rhat \nab_3^{ J }\Ab =  - \Abt  + O(1 )  \nab_\Rhat  \dk^{\le J- 1 }(\nab_4, \nab)\Ab + O(r^{-1} )\nab_{\Rhat}\dk^{\le J-1}\Ab+ O(r^{-1} )  \dk^{\le J-1}\Ab
    \eeaa
    which is derived as the corresponding one in Lemma \ref{lemma:identityBbt-nab_3Bb}. We then make use of the estimates for $\Abt $ derived in Lemma \ref{lemma:lemma-EstimatesabtBbtm3} and the estimates for $\nab_4 \Ab, \nab \Ab$ derived in Lemma \ref{Lemma:BianchA-steps1-2}, the induction hypothesis, and the above identity,   to deduce the following analog of \eqref{estimate:degenerate-Bb1}, with $b= -\de $,      
      \bea\lab{estimate:degenerate-Ab1}
    \bsplit
   &  \int_{\MM } r^{b-1}  \left( \frac{\De}{|q|^2 } \right)^{ 2 }       \big|\nab_\Rhat (\nab_3 )^{J}  \Ab\big|^2        +\int_{\pr^+ \MM}  r^{ b-2 }   \left( \frac{\De}{|q|^2 } \right)^{ 2 }  \big|\nab_\Rhat (\nab_3 )^{J}  \Ab\big|^2  \\
   &\les \de_{J+1}[\Bb] +\ep_0^2+\ep_J^2 +  O(a^2, \ep^2)  \BEF_\de[\nab^J_3\Ab ] +\sqrt{\de_{J+1}[\Bb]}\Big(B^{J}_\de[\Ab]+\Sk_{J+1}^2\Big)^{\frac{1}{2}}\\
&+\ep_J\Big(B^J_\de[\Ab]\Big)^{\frac{1}{2}}+\Sk_{J+1}\Big(\de_{J+1}[\Bb] +\ep_0^2+\ep_J^2 +  O(a^2, \ep^2)  \BEF_\de[ \nab^J_3\Ab ]\Big)^{\frac{1}{2}}.
 \end{split}
    \eea
 Finally we  can repeat the procedure\footnote{Note that no additional linearization is needed in this case.}  of Steps 1--4 of section  \ref{section:InteriorEstimates-Bb} to derive the following analog of \eqref{estimate:nondegenerate-Bb2}
\bea\lab{estimate:nondegenerate-Ab2}
\bsplit
   &  \int_{\MM }   r^{b-1}   \big|\chi_{red}\nab_3 ^{J+1}  \Ab\big|^2   +   \int_{\pr^+\MM }    r^{ b-2 }    \big|\chi_{red}\nab_3 ^{J+1}  \Ab\big|^2\\
   \les& \de_{J+1}[\Bb] +\ep_0^2+\ep_J^2 +  O(a^2, \ep^2)  \BEF_\de[ \nab^J_3\Ab ] +\sqrt{\de_{J+1}[\Bb]}\Big(B^{J}_\de[\Ab]+\Sk_{J+1}^2\Big)^{\frac{1}{2}}\\
&+\ep_J\Big(B^J_\de[\Ab]\Big)^{\frac{1}{2}}+\Sk_{J+1}\Big(\de_{J+1}[\Bb] +\ep_0^2+\ep_J^2 +  O(a^2, \ep^2)  \BEF_\de[ \nab^J_3\Ab ]\Big)^{\frac{1}{2}},
   \end{split}
\eea
where the smooth cut-off function $\chi_{red}(r)$ is such that $\chi_{red}(r)=1$ for $r\le r_+(1+\deh)$ and $\chi_{red}(r)=0$ for $r\geq r_+(1+2\deh)$.
    
Since $\chi_{red}(r)=1$ for $r\geq r_+(1+\deh)$, we have
  \beaa
  \BEF^J_\de[\Ab] &\les& \BEF^{J-1}_\de[r(\nab_4, \nab)\Ab]+ \BEF_\de[\nab_3^J\Ab]+\BEF^{J-1}_\de[\Ab]\\
  &\les& \BEF^{J-1}_\de[r(\nab_4, \nab)\Ab]+\BEF^{J-1}_\de[\Ab]\\
  && +\int_{\MM } r^{-\de-1}  \left( \frac{\De}{|q|^2 } \right)^2       \big|\nab_\Rhat \nab_3 ^{ J }  \Ab\big|^2        +\sup_{\tau\leq\tau_*}\int_{\pr^+ \MM(1,\tau)}  r^{-\de-2}   \left( \frac{\De}{|q|^2 } \right)^2  \big|\nab_\Rhat \nab_3^{ J }  \Ab\big|^2\\
  &&+\int_{\MM }\big|\chi_{red}\nab_3^{J+1}\Ab\big|^2        +\sup_{\tau\leq\tau_*}\int_{\pr^+ \MM(1,\tau)}\big|\chi_{red}\nab_3^{J+1}  \Ab\big|^2
  \eeaa
  which together with \eqref{estimate:degenerate-Ab1}, \eqref{estimate:nondegenerate-Ab2}, and the estimates for $\nab_4 \Ab, \nab \Ab$ derived in Lemma \ref{Lemma:BianchA-steps1-2}, yields
  \beaa
 \bsplit
  \BEF^J_\de[\Ab] \les& \de_{J+1}[\Bb] +\ep_0^2+\ep_J^2 +  O(a^2, \ep^2)  \BEF_\de[ \nab^J_3\Ab ] +\sqrt{\de_{J+1}[\Bb]}\Big(B^{J}_\de[\Ab]+\Sk_{J+1}^2\Big)^{\frac{1}{2}}\\
&+\ep_J\Big(B^J_\de[\Ab]\Big)^{\frac{1}{2}}+\Sk_{J+1}\Big(\de_{J+1}[\Bb] +\ep_0^2+\ep_J^2 +  O(a^2, \ep^2)  \BEF_\de[ \nab^J_3\Ab ]\Big)^{\frac{1}{2}}.
\end{split}
\eeaa  
   For $a$ and $\ep$ small enough, we infer
  \beaa
  \BEF^J_\de[\Ab]   &\les& \de_{J+1}[\Bb]  +\ep_0^2+\ep_J^2  +\left(\sqrt{ \de_{J+1}[\Bb]} +\ep_0+\ep_J\right)\Sk_{J+1}+|a|^2\Sk_{J+1}^2
  \eeaa 
  which is the estimate for $\BEF^J_\de[\Ab]$ in \eqref{eq:proposition-EstimatesAAb-interior} stated in Proposition \ref{proposition:EstimatesAAb-interior}.


\subsection{Estimates for $A$ in Proposition \ref{proposition:EstimatesAAb-interior}} 
    

  We start by deriving  equations  for $(\At, \Bt)$ where we recall that 
  \beaa
  \At=( \ov{q} \nabc_4)^{J-1}   \nabc_{\Rhat}^2  A, \qquad \Bt= ( \ov{q} \nabc_4)^{J-1} \nabc_{\Rhat}^2 B.
  \eeaa
   This is done by commuting  the  first Bianchi pair  \eqref{eq:first-pair-A-B-lin}  
  \beaa
\begin{split}
 \nabc_3A +\frac{1}{2}\tr\Xb A&= -\DDs_2\, B +  O(ar^{-2}) B  +O(r^{-3} ) \Xh +r^{-2} \Ga_b\c \Rc_b,\\
\nabc_4B  +2\ov{\tr X} B &=   \DDd_2 A+  O(ar^{-2}) A    +O(r^{-3} ) \Xi + r^{-2} \Ga_b\c \Rc_b,
\end{split}
\eeaa
  with $( \ov{q} \nabc_4)^{J-1}\nabc_\Rhat^2$  to derive the following analog of the first part of  Lemma \ref{lemma:EstimatesBtPtp}.
  
  \begin{lemma}
 \lab{Lemma:linearization-tildes-A-B}
 The quantities $\At, \Bt$ verify the following  system
 \bea
 \bsplit
  \nabc_3 \At +\left(\frac 1 2- \frac {J-1}{ 2} \right) \tr \Xb  \At &=-\DDs_2\, \Bt  +\Ft_{(1)},\\
  \nabc_4   \Bt+    \left(2-\frac{J-1}{2}\right) \ov{\tr X}\Bt &= \DDd_2 \At +\Ft_{(2)}, 
  \end{split}
  \eea
  with
  \beaa
  \Ft_{(1)}&=&  O(r^{-1} )  \dk^{J-1} \nab_\Rhat \big(\nabc_3  A, \dkb A)   + O(r^{- 1 })  \dk^{J} \nab_\Rhat B +   O(r^{-3} ) \nab\dk^{\le J}\Ga_b\\
  && + O(r^{- 1 })  \dk^{\le J}(A, B) +O(r^{-3} )\dk^{\le J}  \Ga_b   + r^{-2}   \dk^{\le J+1}(\Ga_b\c \Rc_b),\\
   \Ft_{(2)}&=&  O(r^{-1})
   \dk^{J} \nab_\Rhat (A,  B)+ O(r^{-3} )   \dk^{\le J+1}   \Ga_b + O(r^{- 1 })  \dk^{\le J}(A, B) + r^{-2}   \dk^{\le J+1}(\Ga_b\c \Rc_b).
  \eeaa
  \end{lemma}
  
 \begin{remark}
 A priori, $\Ft_{(1)}$ contains the dangerous term $O(r^{-3})(\ov{q}\nabc_4)^{J-1}\nabc_\Rhat^2\Xh$ and we make use of the following null structure equation, see   \eqref{eqts:usefulNullstructure},
 \beaa
 \nabc_4\Xh+\Re(\tr X)\Xh&=& \frac 1 2  \DDc\hot \Xi+O(ar^{-2})\Xi -A+ r^{-2} \Ga_b\c\Ga_b,
\eeaa
to trade it for  suitable contributions to $\Ft_{(1)}$ such as $O(r^{-3} )\nab\dk^{\le J}\Ga_b$.
 \end{remark}

 We are now ready to derive spacetime estimates   by appealing to  the first part of  Proposition  \ref{Prop:Bainchi-pairsEstimates-integrated} with $\Psi_{(1)}= \At$, $\Psi_{(2)}=\Bt$. In this case $2k-1>0$,   $c_{(1)}=\frac 1 2(2-J)  $ and   we need $\La_{(1)} =-2 c_{(1)} +1+\frac b 2= J-1+ \frac{b}{2}   > 0$.   We can thus choose $b= 2+\de $ and proceed as in the derivation of the estimates for 
 $\Bt$ in Proposition \ref{Prop:-EstimatesBtPtp3} to derive, for small $a$ and  $ b=2+\de $,
 \beaa
   \begin{split}
 &\int_{\MM } r^{b-1} |\At |^2+\int_{\pr^+\MM} r^{b}|\At|^2\\
 \les& \de_{J+1}[B] +\ep_0^2+\ep_J^2 + O(a^2, \ep^2)\BEF_\de[r^{J+1}\nab^J_4A] +\sqrt{\de_{J+1}[B]}\Big(B^{J}_\de[A]+\Sk_{J+1}^2\Big)^{\frac{1}{2}}\\
&+\ep_J\Big(B^J_\de[A]\Big)^{\frac{1}{2}}+\Sk_{J+1}\Big(\de_{J+1}[B] +\ep_0^2+\ep_J^2 + O(a^2, \ep^2)\BEF_\de[r^{J+1}\nab^J_4A]\Big)^{\frac{1}{2}}.
\end{split}
  \eeaa 
 We can then proceed as in section \ref{section:EstimatesforBfromBt} to derive  the following estimate for $A$
 \beaa
   \BEF^J_\de[  r^2 A]&\les \de_{J+1}[B]  +\ep_0^2+\ep_J^2  +\left(\sqrt{ \de_{J+1}[B]} +\ep_0+\ep_J\right)\Sk_{J+1}+|a|^2\Sk_{J+1}^2
 \eeaa
which is the estimate for $\BEF^J_\de[r^2A]$ in \eqref{eq:proposition-EstimatesAAb-interior} stated in Proposition \ref{proposition:EstimatesAAb-interior}.


\chapter{Curvature  estimates in $\Mext$}
\lab{CHAPTER:ESTIMATES-MEXTM8}



\section{Statement of the main result of Chapter \ref{CHAPTER:ESTIMATES-MEXTM8}} 


The goal of this chapter is to prove the following theorem.
\begin{theorem}
\lab{THM:MAINRESULTMEXT}
The following estimate holds  
\bea
\Rkext_{J+1}^2 &\les& r_0^{3+\de_B}\Rkint^2_{J+1} + r_0^{-\de_B}\Skext^2_{J+1} +\ep_J^2 +\ep_0^2,
\eea
with $\Rkext_k$, $\Rkint_k$ and $\Skext_k$  defined as in  section \ref{subsection:MainNormsM8}.
\end{theorem}

\begin{remark}
Theorem \ref{THM:MAINRESULTMEXT} implies \eqref{eq:ExteriorcurvEstimatesThmM8}. Together with the proof of \eqref{eq:InteriorcurvEstimatesThmM8} in section \ref{sec:proofofeq:InteriorcurvEstimatesThmM8}, this concludes the proof of Theorem \ref{prop:rpweightedestimatesiterationassupmtionThM8}. 
\end{remark}

\begin{remark}
The proof of Theorem \ref{THM:MAINRESULTMEXT} follows the  main steps    in the proof of the corresponding result   in \cite{KS}, see the second estimate in Proposition 8.10 of  \cite{KS} for the corresponding statement, and section 8.7 in \cite{KS} for the corresponding proof. 
\end{remark}

The rest of Chapter \ref{CHAPTER:ESTIMATES-MEXTM8} focuses on the proof of Theorem \ref{THM:MAINRESULTMEXT}. Before moving to $J+1$ derivatives, we first control $\Rkext_0$ in section \ref{sec:proofofThm:MainResultMext:k=0case}, by relying on the linearization of the Bianchi pairs of section \ref{sec:proofofThm:MainResultMext:k=0case:prelim2} and on the estimates of section \ref{sec:proofofThm:MainResultMext:k=0case:prelim1}. Then, we prove Theorem \ref{THM:MAINRESULTMEXT} in section \ref{sec:proofofThm:MainResultMext:generalcase} by relying on the Bianchi pairs for higher order derivatives  derived first for angular derivatives in section \ref{sec:proofofThm:MainResultMext:generalcase:prelim1} and then for general derivatives in section \ref{sec:proofofThm:MainResultMext:generalcase:prelim2}.

 
\section{Standard linearization of the Bianchi pairs}
\lab{sec:proofofThm:MainResultMext:k=0case:prelim2} 


In order to prove Theorem \ref{THM:MAINRESULTMEXT}, it suffices to make use of the standard linearization  of the second and third Bianchi pairs\footnote{See  also Definition
\ref{definition:BianchiPairs},  but note that we need a more precise form of the the error terms here.} in Proposition \ref{prop:bianchi:complex}.

\begin{lemma}\label{lemma:standard-linearization-B-Pc}
The standard linearization of the second Bianchi pair \eqref{eq:second-pair-B-P} has the following  form\footnote{Recall  that we  split $\Ga_b$ into $\Ga_b'=\Ga_b\setminus\{\Xib\}$ and $\Xib$  as the latter  behaves somewhat worse in powers of $r$, see Definition \ref{def:exteriorGac.norms}, and that we similarly also split $\Ga_g$ into $\Ga_g'=\Ga_g\setminus\{\trXbc\}$ and $\trXbc$. Finally, recall that, in part III, we often identity $\Rc_g$ with $r^{-2}\Rc_b$, $\Ga_g$ with $r^{-1}\Ga_b$ and $\Ga_g'$ with $r^{-1}\Ga_b'$, see Remark \ref{Remark:Ga_g=r^{-1}Ga_b}.}
\bea
\lab{eq:BianchiPair2-linearized.M}\label{eq:second-pair-B-P-lin}
\bsplit
\nabc_3B  +\tr\Xb B &= -\DDs_1\,   \ov{ \Pc } + O(ar^{-2} ) \Pc +O(r^{-3} )  \Ga'_b + r^{-2} \Ga' _b\c \Rc_b + \frac{1}{2}\ov{\Xib} \c  A,\\
 \nabc_4 \ov{\Pc}   +\frac{3}{2} \ov{\tr X} \,  \ov{\Pc}   &=  \DDd_1B + O( a r^{-2}) B  +O(r^{-4})\Ga_b  + r^{-2}\Ga'_b\c \Rc_b.
 \end{split}
\eea
The standard linearization of the  third  Bianchi pair \eqref{eq:third-pair-P-Bb} has the following form
\bea\lab{eq:BianchiPair3-linearized.M}\label{eq:third-pair-P-Bb-lin}
\begin{split}
\nabc_3\Pc+ \frac{3}{2}\ov{\tr\Xb} \Pc &=  - \DDd_1\Bb +  O(ar^{-2} )\Bb  +O(r^{-3})\Ga_b  + r^{-1}\Ga_b\c \Rc_b, \\
\nabc_4\Bb+ \tr X\Bb  &=\DDs_1 \, \Pc  +O(ar^{-2}) \Pc   +O(r^{-3})\Ga_b +r^{-1} \Ga_b\c \Rc_b.
\end{split}
\eea
\end{lemma}

\begin{proof}
The proof follows  easily from the form of the second and third Bianchi pairs in \eqref{eq:second-pair-B-P}  \eqref{eq:third-pair-P-Bb} by writing $P=\Pc-\frac{2m}{q^3} $,  $H=\frac{aq}{|q|^2}\Jk +\Ga_b'$,  $\Hb=-\frac{a\ov{q}}{|q|^2}\Jk +\Ga_g$,  $\widecheck{e_4(q)}=\Ga_g$, $\widecheck{e_3(q)}=r\Ga_b$, and $\DD q =-a\Jk+ r \Ga_g$. 
\end{proof}

\begin{remark}
We note that the smallness of $a$ plays no role  in this  chapter  
 and we may  thus ignore it in what follows. We  will keep it however 
  whenever it is convenient to  take the scaling into account.
\end{remark}

We recall below Definition \ref{definition:GenBianchPairs} exhibiting the general form of Bianchi pairs.  
\begin{definition}
\lab{definition:GenBianchPairs:gen}
   We consider the following  general  Bianchi pairs in $\MM$:
   \begin{itemize}
   \item   For  $\Psi_{(1)}\in\sk_p$, $\Psi_{(2)}\in\sk_{p-1}$, and    $F_{(1)}\in\sk_{p}$, $F_{(2)}\in \sk_{p-1}$,  
   \bea\lab{eq:modelbainchipairequations11-simple:gen}
\begin{split}
\nabc_3(\Psi_{(1)})+c_{(1) }\tr \Xb\Psi_{(1)} &= -\DDs_p\, \Psi_{(2)}  +F_{(1)},\\[2mm]
\nabc_4(\Psi_{(2)})+c_{(2)} \ov{\tr X}\Psi_{(2)} &= \DDd_p\, \Psi_{(1)}   +F_{(2)}.
\end{split}
\eea

\item   For  $\Psi_{(1)}\in\sk_{p-1}$, $\Psi_{(2)}\in\sk_{p}$, and    $F_{(1)}\in\sk_{p-1}$, $F_{(2)}\in \sk_{p}$,         
\bea\lab{eq:modelbainchipairequations12-simple:gen}
\begin{split}
\nabc_3(\Psi_{(1)})+c_{(1)}\ov{\tr \Xb}\Psi_{(1)}&=- \DDd_p\, \Psi_{(2)} +F_{(1)},\\[2mm]
\nabc_4(\Psi_{(2)})+c_{(2)}\tr X\Psi_{(2)} &=\DDs_p\, \Psi_{(1)}  +F_{(2)}.
\end{split}
\eea
 \end{itemize}
 \end{definition}

We rewrite below all the linearized Bianchi identities, see Proposition \ref{prop:bianchi:complex}  and Lemma \ref{lemma:standard-linearization-B-Pc},
 in the form of Definition \ref{definition:GenBianchPairs:gen}.
  
 \begin{proposition}
 \lab{prop:structureBainchPairs}
 The linearized Bianchi identities have the following structure:
 \begin{itemize}
 \item The first Bianchi  pair for $(\Psi_{(1)}=A, \Psi_{(2)}=B)$
 can be written in the form \eqref{eq:modelbainchipairequations11-simple:gen}
  \beaa
\begin{split}
\nabc_3(\Psi_{(1)})+c_{(1) }\tr \Xb\Psi_{(1)} &= -\DDs_2\, \Psi_{(2)}   +F_{(1)},\\[2mm]
\nabc_4(\Psi_{(2)})+c_{(2)} \ov{\tr X}\Psi_{(2)} &= \DDd_2\, \Psi_{(1)}  +F_{(2)},
\end{split}
\eeaa
  with $c_{(1)}=\frac 1 2, \, c_{(2)}= 2$ and
  \beaa
  F_{(1)}&=& O(ar^{-2} ) B+ O(r^{-3} ) \Ga_g + \Ga_b \c B, \\
  F_{(2)}&=& O(a r^{-2} ) A + O(r^{-3} ) \Ga_g+ \Ga_b\c A.
  \eeaa

  \item The second Bianchi  pair for\footnote{See Remark \ref{remark:funnysecondpair}.}  $(\Psi_{(1)}=B, \Psi_{(2)}=\ov{\Pc})$  can be written in the form
  \eqref{eq:modelbainchipairequations11-simple:gen}
  \beaa
\begin{split}
\nabc_3(\Psi_{(1)})+c_{(1) }\tr \Xb\Psi_{(1)} &= -\DDs_1\, \Psi_{(2)}   +F_{(1)},\\[2mm]
\nabc_4(\Psi_{(2)})+c_{(2)} \ov{\tr X}\Psi_{(2)} &= \DDd_1\, \Psi_{(1)}   +F_{(2)},
\end{split}
\eeaa
  with $c_{(1)}=1, \, c_{(2)}= \frac{3}{2}$  and
  \beaa
   F_{(1)}&=&  O(ar^{-2} ) \Pc +O(r^{-3} )  \Ga'_b + r^{-2} \Ga'_b\c \Rc_b+ \Ga_b \c  A,\\
   F_{(2)}&=&   O( a r^{-2}) B  +O(r^{-4})\Ga_b  + \Ga'_b\c (A, B) +r^{-1}\Xi\c\Rc_b,
  \eeaa
where,   see \eqref{eq:defintionofGabprimeandGagprime:partIII},  $\Ga_b'=\Ga_b\setminus \Xib $.

 \item The third Bianchi pair for $( \Psi_{(1)}=\Pc,\Psi_{(2)}= \Bb)$ 
can be written in the form \eqref{eq:modelbainchipairequations12-simple:gen}
\beaa
\begin{split}
\nabc_3(\Psi_{(1)})+c_{(1)}\ov{\tr \Xb}\Psi_{(1)}&=- \DDd_1\,  \Psi_{(2)}+F_{(1)},\\[2mm]
\nabc_4(\Psi_{(2)})+c_{(2)}\tr X\Psi_{(2)} &=\DDs_1\, \Psi_{(1)}+ F_{(2)},
\end{split}
\eeaa
with  $c_{(1)}=\frac{3}{2}, \, c_{(2)}= 1$  and
\beaa
 F_{(1)}&=& O(ar^{-2} )\Bb  +O(r^{-3})\Ga_b  + r^{-1} \Ga_b' \c \Rc_b, \\
    F_{(2)}&=&O(ar^{-2}) \Pc   +O(r^{-3})\Ga_b +r^{-2} \Ga_b\c \Rc_b +\Xi\c \Rc_b.
\eeaa

\item The fourth Bianchi pair for $( \Psi_{(1)}=\Bb,\Psi_{(2)}= \Ab)$  can be written in the form  \eqref{eq:modelbainchipairequations12-simple:gen}
\beaa
\begin{split}
\nabc_3(\Psi_{(1)})+c_{(1)}\ov{\tr \Xb}\Psi_{(1)}&=- \DDd_2\, \Psi_{(2)}+ O(ar^{-2}) \Psi_{(2)} +E_{(1)},\\[2mm]
\nabc_4(\Psi_{(2)})+c_{(2)}\tr X\Psi_{(2)} &=\DDs_2\, \Psi_{(1)} + O(ar^{-2}) \Psi_{(1)}  +E_{(2)},
\end{split}
\eeaa
with 
 $c_{(1)}=2, \, c_{(2)}= \frac 1 2 $  and
\beaa
 F_{(1)}&=&   O(ar^{-2}) \Ab+  O(r^{-3} ) \Ga_b+ \Ga_b' \c\Rc_b, 
 \\
    F_{(2)}&=&  O(ar^{-2}) \Bb + O(r^{-3} ) \Ga_b+ r^{-1} \Ga_b\c\Rc_b.
\eeaa

 \end{itemize}
 \end{proposition}


\section{Exterior estimates for generalized Bianchi pairs} 
\lab{sec:proofofThm:MainResultMext:k=0case:prelim1}


We now state the main  integral exterior estimates for  generalized Bianchi pairs. 

 \begin{proposition}
 \lab{Prop:BasicWeighteBianch-Ext}
 For a given $b$ real we define, as in  Lemma \ref{Le:BasicBianchiPairs.simple},
  \beaa
   \La_{(1)}:= -2c_{(1)}+ 1+\frac b 2, \qquad  \La_{(2)} :=  -2c_{(2)}+ 1+\frac b 2.
   \eeaa
   Then in  both cases \eqref{eq:modelbainchipairequations11-simple:gen} and \eqref{eq:modelbainchipairequations12-simple:gen}, the following holds\footnote{Recall that $\Mext$     is defined by $r\ge r_0$.}:

{\bf Case 1.} If  $\La_{(1)}>0$ and $\La_{(2)} <0$, we have 
 \bea\lab{eq:breakingrpweightedestimatesforBinachipairsin3cases:case1}
\nn && \int_\Mext  r^{b-1} \Big(  |\Psi_{(1)}|^2 +  |\Psi_{(2)}|^2\Big)+\sup_\tau     \int_{\Sext(\tau)}  
 r^b \big( |\Psi_{(1)}|^2  + r^{-2} |\Psi_{(2)}|^2 \big) \\
 \nn&&+\int_{\Si_*}r^b \Big(  |\Psi_{(1)}|^2 +  |\Psi_{(2)}|^2\Big)\\
\nn &\les &   \int_{\Sext(1)}  
 r^b \big( |\Psi_{(1)}|^2  + r^{-2} |\Psi_{(2)}|^2 \big) + r_0^b   \Int(\Psi_{(1)}, \Psi_{(2)}) \\
 && +\int_{\MM( r\ge r_0/2)} r^{b+1}\Big(|F_{(1)}|^2+|F_{(2)}|^2 \Big),
 \eea
where the  quantity 
\beaa
  \Int(\Psi_{(1)}, \Psi_{(2)})&:=&  r_0^{-1} \int_{\MM( r_0/ 2, r_0)  }\big(|\Psi_{(1)}|^2 + |\Psi_{(2)}|^2\big)
  \eeaa
   depends only on the control of  $\Psi_{(1)}, \Psi_{(2)}$ in $\Mint$, i.e. in $r\leq r_0$.

{\bf Case 2.}  If  $\La_{(1)}\leq 0$ and $\La_{(2)} <0$, we have 
 \bea\lab{eq:breakingrpweightedestimatesforBinachipairsin3cases:case2}
\nn && \int_\Mext  r^{b-1}|\Psi_{(2)}|^2 +\sup_\tau     \int_{\Sext(\tau)}  
 r^b \big( |\Psi_{(1)}|^2  + r^{-2} |\Psi_{(2)}|^2 \big) \\
  \nn&&+\int_{\Si_*}r^b \Big(  |\Psi_{(1)}|^2 +  |\Psi_{(2)}|^2\Big)\\
\nn &\les &   \int_{\Sext(1)}  
 r^b \big( |\Psi_{(1)}|^2  + r^{-2} |\Psi_{(2)}|^2 \big) + r_0^b   \Int(\Psi_{(1)}, \Psi_{(2)}) +\int_\Mext  r^{b-1}|\Psi_{(1)}|^2\\
 && +\int_{\MM( r\ge r_0/2)} r^{b+1}\Big(|F_{(1)}|^2+|F_{(2)}|^2 \Big).
 \eea

{\bf Case 3.}  If  $\La_{(2)}=0$, we have 
 \bea\lab{eq:breakingrpweightedestimatesforBinachipairsin3cases:case3}
\nn &&  \sup_\tau     \int_{\Sext(\tau)}  
 r^b \big( |\Psi_{(1)}|^2  + r^{-2} |\Psi_{(2)}|^2 \big) +\int_{\Si_*}r^b \Big(  |\Psi_{(1)}|^2 +  |\Psi_{(2)}|^2\Big)\\
\nn &\les &   \int_{\Sext(1)}  
 r^b \big( |\Psi_{(1)}|^2  + r^{-2} |\Psi_{(2)}|^2 \big) + r_0^b   \Int(\Psi_{(1)}, \Psi_{(2)}) +\int_\Mext  r^{b-1+\dt}|\Psi_{(1)}|^2\\
 && +\int_\Mext  r^{b-1-\dt}|\Psi_{(2)}|^2 +\int_{\MM( r\ge r_0/2)} r^{b+1}\Big(r^{-\dt}|F_{(1)}|^2+r^{\dt}|F_{(2)}|^2 \Big).
 \eea

{\bf Case 4.}  If  $\La_{(1)}>0$, we have 
 \bea\lab{eq:breakingrpweightedestimatesforBinachipairsin3cases:case4}
\nn && \int_\Mext  r^{b-1}|\Psi_{(1)}|^2 +\sup_\tau     \int_{\Sext(\tau)}  
 r^b \big( |\Psi_{(1)}|^2  + r^{-2} |\Psi_{(2)}|^2 \big) \\
  \nn&&+\int_{\Si_*}r^b \Big(  |\Psi_{(1)}|^2 +  |\Psi_{(2)}|^2\Big)\\
\nn &\les &   \int_{\Sext(1)}  
 r^b \big( |\Psi_{(1)}|^2  + r^{-2} |\Psi_{(2)}|^2 \big) + r_0^b   \Int(\Psi_{(1)}, \Psi_{(2)}) +\int_\Mext  r^{b-1}|\Psi_{(2)}|^2\\
 && +\int_{\MM( r\ge r_0/2)} r^{b+1}\Big(|F_{(1)}|^2+|F_{(2)}|^2 \Big).
 \eea
 \end{proposition}

 \begin{proof}
 The proof, similar to the proof of Proposition \ref{Prop:Bainchi-pairsEstimates-integrated},  is based on the following  steps.
 
 {\bf Step 1.} The following is an immediate corollary of Lemma \ref{Le:BasicBianchiPairs.simple} and the fact that, in $\Mext$, we have $\om=O( r^{-2})$, $\omb=\Ga_b$, and $\Ga_b=O(\ep r^{-1})$. 
 \begin{corollary}
\lab{Cor:BasicBianchiPairs.simple}
Let $\Psi_{(1)}, \Psi_{(2)}$, verifying either one of the equations \eqref{eq:modelbainchipairequations11-simple:gen} and \eqref{eq:modelbainchipairequations12-simple:gen}, for positive real numbers  $c_{(1)}$ and $c_{(2)}$, with $\Psi_{(1)}$ of  signature  $k$ and $\Psi_{(2)}$ of  signature  $k-1$. Then denoting
   \beaa
   \La_{(1)}= -2c_{(1)}+ 1+\frac b 2, \qquad  \La_{(2)} =  -2c_{(2)}+ 1+\frac b 2,
   \eeaa
 the following  pointwise  identity holds true for any real $b$ in $\Mext$:
   \begin{enumerate}
\item If $\Psi_{(1)}, \Psi_{(2)}$, verify equation \eqref{eq:modelbainchipairequations11-simple:gen}, then
\bea\lab{eq:basicdivergenceidentitybianchipairrpweightedestimate-complex-1.simple-ext}
\begin{split}
& \Div\left(\frac 1 2|q|^b|\Psi_{(1)}|^2e_3+|q|^b|\Psi_{(2)}|^2e_4-2 |q|^b \Re ( \Psi_{(1)} \c \ov{\Psi_{(2)}})\right)\\
=& \frac 1 2 |q|^{b}  \La_{(1)} \trchb |\Psi_{(1)}|^2+|q|^{b}  \La_{(2)}\trch |\Psi_{(2)}|^2 +  O(r^b)\big( |F_{(1)}||\Psi_{(1)}| +  |F_{(2)}||\Psi_{(2)}|\big) \\
&+  O(r^{b-2})\big( |\Psi_{(1)}|^2 + |\Psi_{(2)}|^2\big) +  O(\ep r^{b-1})\big( |\Psi_{(1)}|^2 +|\Psi_{(1)}||\Psi_{(2)}|\big).
\end{split}
\eea

\item If $\Psi_{(1)}, \Psi_{(2)}$, verify equation \eqref{eq:modelbainchipairequations12-simple:gen}, then 
\bea\lab{eq:basicdivergenceidentitybianchipairrpweightedestimate-complex-2.simple-ext}
\begin{split}
& \Div\left(|q|^b|\Psi_{(1)}|^2e_3+\frac 1 2|q|^b|\Psi_{(2)}|^2e_4+2 |q|^b \Re ( \Psi_{(1)} \c \ov{\Psi_{(2)}})\right)\\
=& |q|^{b} \La_{(1)} \trchb |\Psi_{(1)}|^2+\frac 1 2 |q|^{b}  \La_{(2)}\trch |\Psi_{(2)}|^2 +  O(r^b)\big( |F_{(1)}||\Psi_{(1)}| +  |F_{(2)}||\Psi_{(2)}|\big)\\
&+  O(r^{b-2})\big( |\Psi_{(1)}|^2 + |\Psi_{(2)}|^2\big) +  O(\ep r^{b-1})\big( |\Psi_{(1)}|^2 +|\Psi_{(1)}||\Psi_{(2)}|\big).
\end{split}
\eea
\end{enumerate}
\end{corollary}

{\bf Step 2.} We multiply   the  divergence identities of Corollary \ref{Cor:BasicBianchiPairs.simple}   
by a smooth  cut-off  function $\chi_0(r)$  supported for $r\ge r_0/2$ and identically   $1$ for $r\ge r_0$. Integrating, applying the divergence Lemma  \ref{lemma:basicdivergenceidentitybianchipairrpweightedestimate},   and  using Lemma
\ref{Lemma:g(X, N)onPrMM} to treat the boundary terms, we deduce
 \beaa
&& \int_\Mext  |q|^{b}\left(  - \La_{(1)} \trchb |\Psi_{(1)}|^2  -  \frac 1 2 |q|^{b}  \La_{(2)}\trch |\Psi_{(2)}|^2\right) +\int_{\Sext}r^b \big( |\Psi_{(1)}|^2  + r^{-2} |\Psi_{(2)}|^2 \big)\\
&&+\int_{\Si_*}r^b\big( |\Psi_{(1)}|^2  + |\Psi_{(2)}|^2 \big)\\
 & \les&   \int_{\Sext(1)}  
 r^b \big( |\Psi_{(1)}|^2  + r^{-2} |\Psi_{(2)}|^2 \big)+  r_0^b   \Int(\Psi_{(1)}, \Psi_{(2)})\\
 && +\int_{\MM( \ge r_0/2)}r^b\Big( |F_{(1)}||\Psi_{(1)}| +  |F_{(2)}||\Psi_{(2)}|\Big) +\int_{\Mext}r^{b-2}\big( |\Psi_{(1)}|^2 + |\Psi_{(2)}|^2\big)\\
 && +\ep\int_{\Mext}r^{b-1}\big( |\Psi_{(1)}|^2 + |\Psi_{(1)}||\Psi_{(2)}|\big),
 \eeaa
where, from now on, we carry out the proof in the case \eqref{eq:modelbainchipairequations12-simple:gen}, the case \eqref{eq:modelbainchipairequations11-simple:gen} being completely analogous.

 {\bf Step 3.} Since $\trch=\frac{2}{r} +O(r^{-2} )$ and $\trchb =-\frac{2}{r}+O(r^{-2} )$,  
we  deduce
 \beaa
&& \int_\Mext  |q|^{b}r^{-1}\left( \La_{(1)}  |\Psi_{(1)}|^2 - \frac 1 2 |q|^{b}  \La_{(2)}|\Psi_{(2)}|^2\right) +\int_{\Sext}r^b \big( |\Psi_{(1)}|^2  + r^{-2} |\Psi_{(2)}|^2 \big)\\
&&+\int_{\Si_*}r^b\big( |\Psi_{(1)}|^2  + |\Psi_{(2)}|^2 \big)\\
 & \les&   \int_{\Sext(1)}  
 r^b \big( |\Psi_{(1)}|^2  + r^{-2} |\Psi_{(2)}|^2 \big)+  r_0^b   \Int(\Psi_{(1)}, \Psi_{(2)})\\
 && +\int_{\MM( \ge r_0/2)}r^b\Big( |F_{(1)}||\Psi_{(1)}| +  |F_{(2)}||\Psi_{(2)}|\Big) +\int_{\Mext}r^{b-2}\big( |\Psi_{(1)}|^2 + |\Psi_{(2)}|^2\big)\\
 && +\ep\int_{\Mext}r^{b-1}\big( |\Psi_{(1)}|^2 + |\Psi_{(1)}||\Psi_{(2)}|\big).
 \eeaa

  {\bf Step 4.}  We now conclude the proof of Proposition \ref{Prop:BasicWeighteBianch-Ext} in the cases 1, 2 and 4, i.e. we prove \eqref{eq:breakingrpweightedestimatesforBinachipairsin3cases:case1}, \eqref{eq:breakingrpweightedestimatesforBinachipairsin3cases:case2} and \eqref{eq:breakingrpweightedestimatesforBinachipairsin3cases:case4}. By Cauchy- Schwartz, for any  $\la>0$, we have
  \beaa
  \int_{\MM( \ge r_0/2)} \Big( |F_{(1)}||\Psi_{(1)}| +  |F_{(2)}||\Psi_{(2)}|\Big) &\les&
  \la^{-1}  \int_{\MM( \ge r_0/2)} r^{b-1}\big(|\Psi_{(1)}|^2 + |\Psi_{(2)}|^2 \big)
  \\
  &&+\la  \int_{\MM( \ge r_0/2)} r^{b+1} \big(|F_{(1)}|^2 + |F_{(2)}|^2 \big).
  \eeaa
 Together with Step 3, we infer, using in particular the fact that $r\geq r_0$ on $\Mext$, 
 \beaa
&& \int_\Mext |q|^{b}r^{-1}\left( \La_{(1)}  |\Psi_{(1)}|^2 - \frac 1 2 |q|^{b}  \La_{(2)}|\Psi_{(2)}|^2\right) +\int_{\Sext}r^b \big( |\Psi_{(1)}|^2  + r^{-2} |\Psi_{(2)}|^2 \big)\\
&&+\int_{\Si_*}r^b\big( |\Psi_{(1)}|^2  + |\Psi_{(2)}|^2 \big)\\
 & \les&   \int_{\Sext(1)}  
 r^b \big( |\Psi_{(1)}|^2  + r^{-2} |\Psi_{(2)}|^2 \big)+  r_0^b   \Int(\Psi_{(1)}, \Psi_{(2)})+ \la  \int_{\MM( \ge r_0/2)} r^{b+1} \big(|F_{(1)}|^2 + |F_{(2)}|^2 \big)\\
 &&  +\Big(r_0^{-1}+\la^{-1}+\ep\Big)\int_{\Mext}r^{b-1}\big( |\Psi_{(1)}|^2 + |\Psi_{(2)}|^2\big).
 \eeaa
 Choosing $\la>0$ large enough, and for $r_0$ sufficiently large and $\ep$ sufficiently small, since $r\leq |q|\leq 2r$, we immediately conclude the proof of the stated estimates \eqref{eq:breakingrpweightedestimatesforBinachipairsin3cases:case1}, \eqref{eq:breakingrpweightedestimatesforBinachipairsin3cases:case2} and \eqref{eq:breakingrpweightedestimatesforBinachipairsin3cases:case4}.  
  
{\bf Step 5.} It remains to treat case 3, i.e. to prove  \eqref{eq:breakingrpweightedestimatesforBinachipairsin3cases:case3}. In this case, by Cauchy- Schwartz, we have
  \beaa
  \int_{\MM( \ge r_0/2)} \Big( |F_{(1)}||\Psi_{(1)}| +  |F_{(2)}||\Psi_{(2)}|\Big) &\les& \int_{\MM( \ge r_0/2)} r^{b-1}\big(r^{\dt}|\Psi_{(1)}|^2 + r^{-\dt}|\Psi_{(2)}|^2 \big)
  \\
  &&+  \int_{\MM( \ge r_0/2)} r^{b+1} \big(r^{-\dt}|F_{(1)}|^2 + r^{\dt}|F_{(2)}|^2 \big).
  \eeaa
 Together with Step 3, since $\La_{(2)}=0$ in this case, we immediately infer the stated estimate \eqref{eq:breakingrpweightedestimatesforBinachipairsin3cases:case3}. 
This concludes the proof of Proposition \ref{Prop:BasicWeighteBianch-Ext}.
 \end{proof}


\section{Basic curvature estimate in $\Mext$}
\lab{sec:proofofThm:MainResultMext:k=0case}


In this section, we provide basic curvature estimate in $\Mext$ on the control of $\Rkext_0$, i.e. we prove the following estimate 
\bea\lab{eq:controlofRkext0asawarmup}
\Rkext_0^2  &\les    r_0^{-\dt}\Skext_0^2+ r_0^{ 3 +\dt}  \Rkint^2_0+\ep_0^2.
\eea
The procedure leading to \eqref{eq:controlofRkext0asawarmup} will be repeated in section \ref{sec:proofofThm:MainResultMext:generalcase} in order to prove Theorem \ref{THM:MAINRESULTMEXT} on the control of $\Rkext_{J+1}$.

In order to prove \eqref{eq:controlofRkext0asawarmup}, we consider successively each of the four Bianchi pairs starting with the first one.


\subsection{First Bianchi pair}
\lab{sec:proofofThm:MainResultMext:k=0case:firstBianchipair}

 
 In the case  of the first Bianchi pair we have $c_{(1)}=\frac 1 2$, $c_{(2)}= 2$. We choose
 $b= 4+\dt$  so that we have in this case
 \beaa
 \La_{(1)}=2+\frac{\dt}{2}>0, \qquad  \La_{(2)}= -1+\frac{\dt}{2}<0.
 \eeaa
We may thus apply case 1 of Proposition \ref{Prop:BasicWeighteBianch-Ext}, i.e. estimate \eqref{eq:breakingrpweightedestimatesforBinachipairsin3cases:case1}.  We infer
 \beaa
&& \int_{\Mext} r^{3+\dt}|(A, B)|^2+  \sup_{\tau} \int_{\Si(\tau)\cap\Mext}   r^{4+\dt}  \Big( |A|^2  + r^{-2} |B|^2\Big)\\
&&+\int_{\Si_*}  r^{4+\de_B} \Big( |A|^2  + |B|^2\Big)\\
&\les&  \int_{\Si(1)} r^{4+\de_B} \Big( |A|^2  + r^{-2} |B|^2\Big)  +r_0^{4+\dt}   \Int(A, B)  +\int_{\MM( r\ge r_0/2)} r^{5+\dt}\Big(|F_{(1)}|^2+|F_{(2)}|^2 \Big).
 \eeaa
  In view of the control of our initial conditions  $\int_{\Si(1)} r^{4+\de_B}( |A|^2  + r^{-2} |B|^2)\les \ep_0^2$,  the interior estimates $  \Int(A, B) \les   r_0^{-1} \Rkint^2_0[A, B]$, and the definition of the $\Rkext$ norms in  Definition \ref{Definiition:Rkext},  we deduce
  \beaa
  \Rkext^2_0[A]+\Rkext_0'[B]^2  &\les& r_0^{ 3 +\dt}  \Rkint^2_0[A, B]+\int_{\MM( r\ge r_0/2)} r^{5+\dt}\Big(|F_{(1)}|^2+|F_{(2)}|^2 \Big)+\ep_0^2
  \eeaa
  where
  \bea
  \bsplit
  \Rkext_0'[B]^2:&=\int_{\Mext} r^{3+\dt}| B |^2 +\sup_{\tau} \int_{\Si(\tau)\cap\Mext}  r^{2+\de_B} |B|^2 
   +\int_{\Si_*} r^{4+\de_B} |B|^2
  \end{split}
  \eea
  agrees with the norm $\Rk_0[B]$ in  Definition \ref{Definiition:Rkext} on $\Mext$ and $\Si_*$, but has a weaker weight in $r$ on $\Si(\tau)$.
   
   Now,
  \beaa
  \int_{\MM( r\ge r_0/2)} r^{5+\dt} |F_{(1)}|^2 &\les&
  \int_{\MM( r\ge r_0/2)}  r^{5+\dt} \Big(    r^{-4 } | B  |^2 + r^{-6} |\Ga_g|^2 + |\Ga_b|^2 |B|^2 \Big)\\
  &\les&  \int_{\MM( r\ge r_0/2)}  r^{1+\dt} |B|^2 +\int_{ \MM( r\ge r_0/2)}  r^{-1+\dt}|\Ga_g|^2 +\ep^4 \\
  &\les&  r_0^{-2} \Rkext_0[B] ^2  +\int_{ \MM( r\ge r_0/2)}  r^{-1+\dt}|\Ga_g|^2 +\ep_0^2. 
  \eeaa
  Hence, in view of Lemma \ref{lemma:auxilliarynormsforGa_b},
  \beaa
    \int_{\MM( r\ge r_0/2)} r^{5+\dt} |F_{(1)}|^2 &\les&   r_0^{-2} \Rkext_0[B] ^2 + r_0^{-2+ 2\dt}\Skext_0^2  +\ep_0^2. 
  \eeaa
  Similarly
  \beaa
    \int_{\MM( r\ge r_0/2)} r^{5+\dt} |F_{(2)}|^2 &\les&   r_0^{-2} \Rkext_0[A] ^2 + r_0^{-2+ 2\dt}\Skext_0^2  +\ep_0^2 
  \eeaa
  and we deduce
 \medskip 
  \bea
  \lab{Estimate:firstPair0}
  \Rkext_0^2[A] +  \Rkext' _0[B]^2 \les  r_0^{-2} \Rkext_0 ^2 + r_0^{-2+ 2\dt}\Skext_0^2+ r_0^{ 3 +\dt}  \Rkint^2+\ep_0^2.
  \eea


\subsection{Second Bianchi pair}



\subsubsection{First estimate for the second Bianchi pair}


 In the case of the second Bianchi pair, we have   $ c_{(1)}=1$, $c_{(2)}=\frac 3 2 $. 
  We choose $b= 4-\de_B$ so that we have in this case
  \beaa
  \La_{(1)} =  1  +\frac{\dt}{2}>0, \qquad   \La_{(2)}=-\frac{\dt}{2}<0.
  \eeaa
 We may thus apply case 1 of Proposition \ref{Prop:BasicWeighteBianch-Ext}, i.e. estimate \eqref{eq:breakingrpweightedestimatesforBinachipairsin3cases:case1}.  We infer
 \beaa
&& \int_{\Mext} r^{3-\dt}|(B, \Pc )|^2+  \sup_{\tau} \int_{\Si(\tau)\cap\Mext}   r^{4-\dt}  \Big( |B|^2  + r^{-2} |\Pc |^2\Big)\\
&&+\int_{\Si_*}  r^{4-\de_B} \Big( |B|^2  + |\Pc |^2\Big)\\
&\les&  \int_{\Si(1)} r^{4-\de_B} \Big( |B|^2  + r^{-2} |\Pc|^2\Big)  +r_0^{4-\dt}   \Int(B, \Pc)  +\int_{\MM( r\ge r_0/2)} r^{5- \dt}\Big(|F_{(1)}|^2+|F_{(2)}|^2 \Big).
 \eeaa
    In view of the control of our initial conditions $\int_{\Si(1)} r^{4-\de_B}( |B|^2  + r^{-2} |\Pc|^2)\les \ep_0^2$,  the interior estimates $\Int(B, \Pc) \les   r_0^{-1} \Rkint^2_0[B, \Pc]$, and the definition of the $\Rkext$ norms in  Definition \ref{Definiition:Rkext},  we deduce
  \beaa
  \bsplit
   \Rkext_0''[B]^2+\Rkext_0'[\Pc]^2\les& r_0^{ 3 -\dt}  \Rkint^2_0[B, \Pc]+   \int_{\MM( r\ge r_0/2)} r^{5 - \dt} \Big(|F_{(1)}|^2+|F_{(2)}|^2 \Big)  +\ep_0^2,\\
    \Rkext_0''[B]^2:=&  \int_{\Mext} r^{3-\dt}|B|^2+  \sup_{\tau} \int_{\Si(\tau)\cap\Mext}   r^{4-\dt}   |B|^2 +\int_{\Si_*}  r^{4-\de_B}  |B|^2,\\
  \Rkext_0'[\Pc]^2:=&    \int_{\Mext} r^{3-\dt}|\Pc |^2+  \sup_{\tau} \int_{\Si(\tau)\cap\Mext}   r^{2-\dt}  |\Pc |^2+\int_{\Si_*}  r^{4-\de_B}   |\Pc |^2.
  \end{split}
  \eeaa
 In this case
  \beaa
   F_{(1)}&=&  O(ar^{-2} ) \Pc +O(r^{-3} )  \Ga'_b + r^{-2} \Ga'_b\c \Rc_b+ \Ga_b \c  A,\\
   F_{(2)}&=&   O( a r^{-2}) B  +O(r^{-4})\Ga_b  + \Ga'_b\c (A, B)+ r^{-1}\Xi\c\Rc_b.
  \eeaa
  Now,  writing  $\Pc, B= r^{-2} \Rc _b$,
  \beaa
    \int_{\MM( r\ge r_0/2)} r^{5- \dt}\Big(|F_{(1)}|^2 +|F_{(2)}|^2\Big) &\les&  \int_{\MM( r\ge r_0/2)} r^{-3-\dt} |\Rc_b|^2     +\int_{\MM( r\ge r_0/2)}  r^{-1-\dt} |\Ga_b'|^2\\
    &&+\int_{\MM( r\ge r_0/2)}  r^{1-\dt}\Big( |\Ga_b'|^2 |\Rc_b|^2 +r^4| \Ga_b |^2 |(A, B)|^2\Big).
  \eeaa
  In view of \eqref{eq:auxilliarynormsforGa_b1-strong}, we have
  \beaa
   \int_{\MM( r\ge r_0/2)}  r^{-1-\dt} |\Ga_b'|^2&\les & r_0^{-\dt}\Skext_0^2.
  \eeaa
  Also, we have
  \beaa
  \int_{\MM( r\ge r_0/2)}  r^{1-\dt} |\Ga_b'|^2 |\Rc_b|^2  +\int_{\MM( r\ge r_0/2)}  r^{5-\dt} |\Ga_b|^2 |(A, B)|^2 &\les& \\ep^4\les\ep_0^2.
  \eeaa
  Hence
  \bea
  \lab{Estimate:secondPair0}
    \Rkext_0''[B]^2+\Rkext_0'[\Pc]^2\les  r_0^{3-\dt}  \Rkint^2+  r_0^{-\dt}\Skext_0^2 + r_0^{-2} \Rkext_0^2 +\ep_0^2.
  \eea


\subsubsection{Second estimate for the second Bianchi pair}

  
In order to control  the norm  $ \Rk_0[B]^2$, we still need to recover the correct weight in $r$ for the part of the norm on $\Si(\tau)$. To do this, we choose $b=4$ so that we have $\La_{(2)}=0$ in this case. We may thus apply case 3 of Proposition \ref{Prop:BasicWeighteBianch-Ext}, i.e. estimate \eqref{eq:breakingrpweightedestimatesforBinachipairsin3cases:case3}.  We infer
  \beaa
  &&  \sup_{\tau} \int_{\Si(\tau)\cap\Mext}   r^{4}  \Big( |B|^2  + r^{-2} |\Pc |^2\Big)+\int_{\Si_*}  r^{4} \Big( |B|^2  + |\Pc |^2\Big)\\
 &\les&    \ep_0^2+ r_0^{ 3 }  \Rkint^2_0[B, \Pc]+
  \int_{\MM( \ge r_0/2)} r^5\Big(r^{-\dt}|F_{(1)}|^2 + r^{\dt}|F_{(2)}|^2\Big)\\
  && +\int_{\Mext} r^{3+\dt} |B|^2+\int_{\Mext} r^{3-\dt} |\Pc|^2.
   \eeaa
  In view of the definition of $\Rkext_0'[B]$ and $\Rkext_0'[\Pc]$, we infer
  \beaa
  && \sup_{\tau} \int_{\Si(\tau)\cap\Mext}\Big( r^4|B|^2  + r^2|\Pc |^2\Big)+\int_{\Si_*}  r^4 |\Pc |^2\\
  &\les&    \ep_0^2+ r_0^3  \Rkint^2_0[B, \Pc] +\Rkext_0'[B]^2+\Rkext_0'[\Pc]^2\\
 &&+\int_{\MM( \ge r_0/2)} r^5\Big(r^{-\dt}|F_{(1)}|^2 + r^{\dt}|F_{(2)}|^2\Big). 
   \eeaa
   Together with \eqref{Estimate:firstPair0} and \eqref{Estimate:secondPair0}, this yields
    \beaa
    &&\sup_{\tau} \int_{\Si(\tau)\cap\Mext}\Big( r^4|B|^2  + r^2|\Pc |^2\Big)+\int_{\Si_*}  r^4 |\Pc |^2\\
     &\les&   r_0^{3+\dt}  \Rkint^2_0 +  r_0^{-\dt}\Skext_0^2 +r_0^{-2}\Rkext_0^2 +\ep_0^2 +\int_{\MM( \ge r_0/2)} r^5\Big(r^{-\dt}|F_{(1)}|^2 + r^{\dt}|F_{(2)}|^2\Big). 
   \eeaa
   Also, in the proof of \eqref{Estimate:secondPair0}, we have obtained the following estimate  
    \beaa
    \int_{\MM( r\ge r_0/2)} r^{5- \dt}|F_{(1)}|^2 &\les&  r_0^{-\dt}\Skext_0^2 +r_0^{-2}\Rkext_0^2 +\ep_0^2
  \eeaa
and hence 
\beaa
 &&\sup_{\tau} \int_{\Si(\tau)\cap\Mext}\Big( r^4|B|^2  + r^2|\Pc |^2\Big)+\int_{\Si_*}  r^4 |\Pc |^2\\
  &\les&   r_0^{3+\dt}  \Rkint^2_0 +  r_0^{-\dt}\Skext_0^2 +r_0^{-2}\Rkext_0^2 +\ep_0^2+\int_{\MM( \ge r_0/2)} r^{5+\dt}|F_{(2)}|^2. 
   \eeaa
   It thus remain to control the last term on the RHS. Recalling that 
   \beaa
   F_{(2)}&=&   O( a r^{-2}) B  +O(r^{-4})\Ga_b  +  \Ga'_b\c (A, B)+ r^{-1}\Xi\c\Rc_b,
  \eeaa
 we have 
  \beaa
    \int_{\MM( r\ge r_0/2)} r^{5+\dt}|F_{(2)}|^2 &\les&  \int_{\MM( r\ge r_0/2)} r^{1+\dt} |B|^2     +\int_{\MM( r\ge r_0/2)}  r^{-3+\dt} |\Ga_b|^2\\
    &&+\int_{\MM( r\ge r_0/2)} \Big( r^{5+\dt}|\Ga_b'|^2 |(A, B)|^2+r^{3+\dt}|\Xi|^2|\Rc_b|^2\Big)\\
    &\les& r_0^{-2}\Rkext_0^2 +r_0^{-2+2\dt}\Skext_0^2 +\ep^4\\
    &\les& r_0^{-2}\Rkext_0^2 +r_0^{-2+2\dt}\Skext_0^2 +\ep^2_0,
  \eeaa
where we have used in particular the fact that $\Xi\in r^{-1}\Ga_g$ in Part III in view of  \eqref{eq:specialidentityforthegloablframeofMMinpartIII}. Hence, we infer
  \bea\lab{Estimate:secondPair1}
 \nn&&\sup_{\tau} \int_{\Si(\tau)\cap\Mext}\Big( r^4|B|^2  + r^2|\Pc |^2\Big)+\int_{\Si_*}  r^4 |\Pc |^2\\
  &\les&   r_0^{3+\dt}  \Rkint^2_0 +  r_0^{-\dt}\Skext_0^2 +r_0^{-2}\Rkext_0^2 +\ep_0^2. 
   \eea
   Together with \eqref{Estimate:firstPair0} and \eqref{Estimate:secondPair0}, we deduce
    \bea
      \lab{Estimate:secondPair2}
    \Rkext_0[B, \Pc]^2 \les  r_0^{3+\dt}  \Rkint^2_0+r_0^{-\dt}\Skext_0^2 + r_0^{-2} \Rkext_0^2 +\ep_0^2.
    \eea


\subsection{Third Bianchi pair}



\subsubsection{First estimate for the third Bianchi pair}


In the case of the third Bianchi pair, we have $ c_{(1)}=3/2$,  $c_{(2)}= 1$.  We choose $b= 2-\de_B$  so that we have in this case 
 \beaa
  \La_{(1)}&=& -1 -\frac{\dt}{2}<0, \qquad   \La_{(2)}= -\frac{\dt}{2}<0.
  \eeaa
 We may thus apply case 2 of Proposition \ref{Prop:BasicWeighteBianch-Ext}, i.e. estimate \eqref{eq:breakingrpweightedestimatesforBinachipairsin3cases:case2}.  We infer
 \beaa
&& \int_{\Mext} r^{1-\dt}|\Bb|^2+  \sup_{\tau} \int_{\Si(\tau)\cap\Mext}   r^{2-\dt}  \Big( |\Pc|^2  + r^{-2} |\Bb|^2\Big)+\int_{\Si_*}  r^{2-\dt} \Big( |\Pc|^2  + |\Bb|^2\Big)\\
&\les&  \int_{\Si(1)} r^{2-\dt}  \Big( |\Pc|^2  + r^{-2} |\Bb|^2\Big)+r_0^{2-\dt}   \Int(\Pc, \Bb)  +\int_{\MM( r\ge r_0/2)} r^{3-\dt}\Big(|F_{(1)}|^2+|F_{(2)}|^2 \Big)\\
&&+\int_{\Mext} r^{1-\dt}|\Pc|^2\\
&\les& r_0^{3 +\dt}  \Rkint^2_0+r_0^{-\dt}\Skext_0^2 +r_0^{-2}\Rkext_0^2 +\ep_0^2 +\int_{\MM( r\ge r_0/2)} r^{3-\dt}\Big(|F_{(1)}|^2+|F_{(2)}|^2 \Big),
 \eeaa 
 where we used in particular the control of $\Rkext_0[\Pc]$ in \eqref{Estimate:secondPair2}. Given that 
 \beaa
 F_{(1)}&=& O(ar^{-2} )\Bb  +O(r^{-3})\Ga_b  + r^{-1} \Ga_b' \c \Rc_b, \\
    F_{(2)}&=&O(ar^{-2}) \Pc   +O(r^{-3})\Ga_b +r^{- 2 } \Ga_b\c \Rc_b +\Xi\c \Rc_b,
    \eeaa
  we easily deduce, using in particular the fact that $\Xi\in r^{-1}\Ga_g$ in Part III in view of  \eqref{eq:specialidentityforthegloablframeofMMinpartIII}, 
  \beaa
  \int_{\MM( r\ge r_0/2)} r^{3-\dt}\Big(|F_{(1)}|^2+|F_{(2)}|^2 \Big)
  &\les&    r_0^{-2} \Rkext_0 ^2+ r_0^{-2-\dt} \Skext_0+\ep_0^2,
  \eeaa
 and hence
 \bea
  \lab{Estimate:ThirdPair0}
 \bsplit
 \Rkext'_0[\Bb]^{ 2 } &\les   r_0^{3 +\dt}\Rkint^2_0+r_0^{-\dt}\Skext_0^2 +r_0^{-2}\Rkext_0^2  +\ep_0^2,\\
  \Rkext'_0[\Bb]^{ 2 }  &:= \int_{\Mext} r^{1-\dt}|\Bb|^2+  \sup_{\tau} \int_{\Si(\tau)\cap\Mext}   r^{-\dt}    |\Bb|^2+\int_{\Si_*}  r^{2-\dt}  |\Bb|^2.
  \end{split}
 \eea


\subsubsection{Second estimate for the third Bianchi pair}

 
In order to control  the norm  $ \Rk_0[\Bb]^2$, we still need to recover the correct weight in $r$ for the part of the norm on $\Si(\tau)$ and on $\Si_*$. To do this, we choose $b=2$ so that we have $\La_{(2)}=0$ in this case. We may thus apply case 3 of Proposition \ref{Prop:BasicWeighteBianch-Ext}, i.e. estimate \eqref{eq:breakingrpweightedestimatesforBinachipairsin3cases:case3}.  We infer 
   \beaa
&&  \sup_{\tau} \int_{\Si(\tau)\cap\Mext}   r^{2}  \Big( |\Pc|^2  + r^{-2} |\Bb|^2\Big)+\int_{\Si_*}  r^{2} \Big( |\Pc|^2  + |\Bb|^2\Big)\\
&\les&  \int_{\Si(1)} r^{2}  \Big( |\Pc|^2  + r^{-2} |\Bb|^2\Big) +r_0^2\Int(\Pc, \Bb)  +\int_{\MM( r\ge r_0/2)} r^{3}\Big(|F_{(1)}|^2+|F_{(2)}|^2 \Big)\\
&& +\int_{\Mext}r^{1+\dt}|\Pc|^2+\int_{\Mext}r^{1-\dt}|\Bb|^2\\
&\les& r_0^{3 +\dt}\Rkint^2_0+r_0^{-\dt}\Skext_0^2 +r_0^{-2}\Rkext_0^2 +\ep_0^2 \\
&& +\int_{\MM( r\ge r_0/2)} r^{3}\Big( r^{-\dt} |F_{(1)}|^2+ r^{\dt} |F_{(2)}|^2 \Big),
 \eeaa
 where we used in particular the control of $\Rkext_0[\Pc]$ and $\Rkext'_0[\Bb]$ provided respectively by   \eqref{Estimate:secondPair2} and \eqref{Estimate:ThirdPair0}.
 
Next, recalling that we have in this case
 \beaa
 F_{(1)}&=& O(ar^{-2} )\Bb  +O(r^{-3})\Ga_b  + r^{-1}\Ga_b'\c \Rc_b, \\
    F_{(2)}&=&O(ar^{-2}) \Pc   +O(r^{-3})\Ga_b +r^{-2}\Ga_b\c \Rc_b +\Xi\c \Rc_b,
    \eeaa 
 we derive
 \beaa
 \int_{\MM( r\ge r_0/2)} r^{3}\Big(r^{-\dt}|F_{(1)}|^2+r^{\dt}|F_{(2)}|^2\Big) &\les& r_0^{-2}\Rkext_0^2+ r_0^{-2+2\dt} \Skext_0^2+\ep_0^2,
 \eeaa
 where we used again the fact that $\Xi\in r^{-1}\Ga_g$ in Part III in view of  \eqref{eq:specialidentityforthegloablframeofMMinpartIII}. We infer
  \bea\lab{Estimate:ThirdPair01}
\nn\sup_{\tau} \int_{\Si(\tau)\cap\Mext}|\Bb|^2 +\int_{\Si_*}  r^2|\Bb|^2 &\les& r_0^{3 +\dt}\Rkint^2_0+r_0^{-\dt}\Skext_0^2 \\
&&+r_0^{-2}\Rkext_0^2 +\ep_0^2.
 \eea
 Combining with \eqref{Estimate:ThirdPair0}, we thus deduce
  \bea
  \lab{Estimate:ThirdPair1-0}
  \Rkext_0[\Bb]^{ 2 }  &\les r_0^{3 +\dt}\Rkint^2_0+r_0^{-\dt}\Skext_0^2 +r_0^{-2}\Rkext_0^2+\ep_0^2.
  \eea


\subsection{Fourth Bianchi pair}



\subsubsection{First estimate for the fourth Bianchi pair}

 
In the case of the third Bianchi pair, we have $ c_{(1)}=2$,  $c_{(2)}= 1/ 2$.  We choose $b= -\de_B$  so that we have in this case 
 \beaa
  \La_{(1)}&=& -3 -\frac{\dt}{2}<0, \qquad   \La_{(2)}= -\frac{\dt}{2}<0.
  \eeaa
 We may thus apply case 2 of Proposition \ref{Prop:BasicWeighteBianch-Ext}, i.e. estimate \eqref{eq:breakingrpweightedestimatesforBinachipairsin3cases:case2}.  We infer 
 \beaa
&& \int_{\Mext} r^{-1-\dt}|\Ab|^2+  \sup_{\tau} \int_{\Si(\tau)\cap\Mext}   r^{-\dt}  \Big( |\Bb|^2  + r^{-2} |\Ab|^2\Big)+\int_{\Si_*}  r^{-\dt} \Big( |\Bb|^2  + |\Ab|^2\Big)\\
&\les&  \int_{\Si(1)} r^{-\dt}  \Big( |\Bb|^2  + r^{-2} |\Ab|^2\Big)+r_0^{-\dt}   \Int(\Bb, \Ab)  +\int_{\MM( r\ge r_0/2)} r^{1-\dt}\Big(|F_{(1)}|^2+|F_{(2)}|^2 \Big)\\
&& +\int_{\Mext} r^{-1-\dt}|\Bb|^2\\
&\les& r_0^{3 +\dt}\Rkint^2_0+r_0^{-\dt}\Skext_0^2 +r_0^{-2}\Rkext_0^2 +\ep_0^2  +\int_{\MM( r\ge r_0/2)} r^{1-\dt}\Big(|F_{(1)}|^2+|F_{(2)}|^2 \Big) 
 \eeaa 
 where we used in particular the control of $\Rkext_0[\Bb]$ in \eqref{Estimate:ThirdPair1-0}. Given that 
 \beaa
 F_{(1)}&=&   O(ar^{-2}) \Ab+  O(r^{-3} ) \Ga_b+ \Ga_b' \c\Rc_b,\\
    F_{(2)}&=&  O(ar^{-2}) \Bb + O(r^{-3} ) \Ga_b+ r^{-1} \Ga_b\c\Rc_b,
\eeaa
and proceeding as before, we have
 \beaa
 \int_{\MM( r\ge r_0/2)} r^{1-\dt}\Big(|F_{(1)}|^2+|F_{(2)}|^2 \Big)&\les&   r_0^{-2} \Rkext_0 ^2+  r_0^{-4}\Skext_0+\ep_0^2.
 \eeaa
 We deduce
 \bea
  \lab{Estimate:FourthPair0}
 \bsplit
 \Rk'_0[\Ab]  &\les   r_0^{3 +\dt}\Rkint^2_0+r_0^{-\dt}\Skext_0^2 +r_0^{-2}\Rkext_0^2 +\ep_0^2,\\
  \Rk'_0[\Ab]  &:= \int_{\Mext} r^{-1-\dt}|\Ab|^2+  \sup_{\tau} \int_{\Si(\tau)\cap\Mext}   r^{-2 -\dt}    |\Ab|^2+\int_{\Si_*}  r^{-\dt}  |\Ab|^2.
  \end{split}
 \eea


\subsubsection{Second estimate for the fourth Bianchi pair}


In order to control  the norm  $ \Rk_0[\Ab]^2$, we still need to recover the correct weight in $r$ for the part of the norm on $\Si(\tau)$ and on $\Si_*$. To do this, we choose $b=0$ so that we have $\La_{(2)}=0$ in this case. We may thus apply case 3 of Proposition \ref{Prop:BasicWeighteBianch-Ext}, i.e. estimate \eqref{eq:breakingrpweightedestimatesforBinachipairsin3cases:case3}.  We infer  
 \beaa
&& \sup_{\tau} \int_{\Si(\tau)\cap\Mext}    \Big( |\Bb|^2  + r^{-2} |\Ab|^2\Big)+\int_{\Si_*}   \Big( |\Bb|^2  + |\Ab|^2\Big)\\
&\les&  \int_{\Si(1)} r^{2}  \Big( |\Bb|^2  + r^{-2} |\Ab|^2\Big)  +   \Int(\Bb, \Ab)  +\int_{\MM( r\ge r_0/2)} r\Big( r^{-\dt} |F_{(1)}|^2+ r^{\dt} |F_{(2)}|^2 \Big)\\
&& +\int_{\Mext}r^{-1+\dt}|\Bb|^2+\int_{\Mext}r^{-1-\dt}|\Ab|^2\\
&\les& r_0^{3 +\dt}\Rkint^2_0+r_0^{-\dt}\Skext_0^2 +r_0^{-2}\Rkext_0^2 +\ep_0^2  +\int_{\MM( r\ge r_0/2)} r\Big( r^{-\dt} |F_{(1)}|^2+ r^{\dt} |F_{(2)}|^2 \Big),
 \eeaa
 where we used in particular the control of $\Rkext_0[\Bb]$ and $\Rkext'_0[\Ab]$ provided respectively by  \eqref{Estimate:ThirdPair1-0} and \eqref{Estimate:FourthPair0}.

Next, recalling that we have in this case
 \beaa
 F_{(1)}&=&   O(ar^{-2}) \Ab+  O(r^{-3} ) \Ga_b+ \Ga_b' \c\Rc_b,\\
    F_{(2)}&=&  O(ar^{-2}) \Bb + O(r^{-3} ) \Ga_b+ r^{-1} \Ga_b\c\Rc_b,
\eeaa
we derive 
 \beaa
  \int_{\MM}r\Big( r^{-\dt} |F_{(1)}|^2+ r^{\dt} |F_{(2)}|^2 \Big)
  &\les&    r_0^{-2} \Rkext_0 ^2+ r_0^{-4+2\dt} \Skext_0+\ep_0^2.
  \eeaa
  We thus derive
   \bea
  \lab{Estimate:FourthPair01}
  \nn &&\sup_{\tau} \int_{\Si(\tau)\cap\Mext}    r^{-2} |\Ab|^2 +\int_{\Si_*}|\Ab|^2\\
   &\les &  r_0^{3 +\dt}\Rkint^2_0+r_0^{-\dt}\Skext_0^2 +r_0^{-2}\Rkext_0^2+\ep_0^2.
  \eea 
  Combining with \eqref{Estimate:FourthPair0}, we thus deduce
  \bea
  \lab{Estimate:FourthPair1-0}
  \Rkext_0[\Ab]  &\les  r_0^{3 +\dt}\Rkint^2_0+r_0^{-\dt}\Skext_0^2 +r_0^{-2}\Rkext_0^2+\ep_0^2.  
  \eea


\subsubsection{Conclusion}


Finally, combining \eqref{Estimate:firstPair0}, \eqref{Estimate:secondPair2}, \eqref{Estimate:ThirdPair1-0} and  \eqref{Estimate:FourthPair1-0},  we deduce
 \beaa
 \Rkext_0^2  &\les     r_0^{-2} \Rkext_0 ^2 + r_0^{-\dt}\Skext_0^2+ r_0^{ 3 +\dt}  \Rkint^2_0+\ep_0^2.
 \eeaa
 Thus,  for $r_0$ sufficiently  large, 
 \beaa
 \Rkext_0^2  &\les    r_0^{-\dt}\Skext_0^2+ r_0^{ 3 +\dt}  \Rkint^2_0+\ep_0^2,
 \eeaa
 which concludes the proof of \eqref{eq:controlofRkext0asawarmup}.

\begin{remark}\lab{rmk:summaryofbasicruvatureestimatesMextusedlater}
For convenience, we summarize the main steps in the proof of the basic curvature estimates \eqref{eq:controlofRkext0asawarmup}:
\begin{enumerate}
\item First, we control the first Bianchi pair by applying case 1 of Proposition \ref{Prop:BasicWeighteBianch-Ext} with the choice $b= 4+\dt$.

\item Next, we control the second Bianchi pair in two steps:
\begin{enumerate}
\item First, we apply case 1 of Proposition \ref{Prop:BasicWeighteBianch-Ext} with the choice $b= 4-\dt$.

\item Then, we apply case 3 of Proposition \ref{Prop:BasicWeighteBianch-Ext} with the choice $b= 4$. 
\end{enumerate}

\item Next, we control the third Bianchi pair in two steps:
\begin{enumerate}
\item First, we apply case 2 of Proposition \ref{Prop:BasicWeighteBianch-Ext} with the choice $b= 2-\dt$.

\item Then, we apply case 3 of Proposition \ref{Prop:BasicWeighteBianch-Ext} with the choice $b= 2$. 
\end{enumerate}

\item Finally, we control the fourth Bianchi pair in two steps:
\begin{enumerate}
\item First, we apply case 2 of Proposition \ref{Prop:BasicWeighteBianch-Ext} with the choice $b= -\dt$.

\item Then, we apply case 3 of Proposition \ref{Prop:BasicWeighteBianch-Ext} with the choice $b= 0$. 
\end{enumerate}
\end{enumerate}
\end{remark}


\section{Hodge systems for angular derivatives}
\lab{sec:proofofThm:MainResultMext:generalcase:prelim1}


Having outlined the procedure for deriving basic curvature estimates in $\Mext$ in section \ref{sec:proofofThm:MainResultMext:k=0case}, we focus in the rest of this chapter on higher order derivatives curvature estimates in order to prove Theorem \ref{THM:MAINRESULTMEXT} on the control of $\Rkext_{J+1}$. To this end, we first commute in this section the Bianchi pairs with weighted angular derivatives.

\begin{definition}
\lab{Def:angularDerivatives-M8}
Let $\La_p$ the operator acting on $\sk_p$ and defined by 
\bea\lab{eq:definitionoftheoperatorLambdapasr2Deltap}
\La_p := r^2\De_p.
\eea
We define weighted angular derivatives $\dkb^j$ as follows:
\begin{enumerate}
\item 
If  $f \in \sk_p$ is  a complex  $p$-tensor, $p=1, 2 $,   and  $j$ is  a positive integer,   we define
\bea\lab{def:angularderivativesonreducedksclars}
\dkb^j f=
\begin{cases}
 \La_p^{\frac{j}{2}} f , \qquad\qquad  \mbox{if $j$ is even,}\\
 r\DDd_p  \La_p^{\frac{j-1}{2}}  f,\qquad \mbox{if $j$ is odd.} 
\end{cases}
\eea

\item  If  $f \in \sk_0$   is a  complex scalar  and    $j$ is a positive integer,   we define
\bea\lab{def:angularderivativesonreducedksclars}
\dkb^j f=
\begin{cases}
 \La_0^{\frac{j}{2}} f , \qquad\qquad\,\,  \mbox{if $j$ is even,}\\
 r \DDs_1\, \,  \La_0^{\frac{j-1}{2}}  f,\qquad \mbox{if $j$ is odd.} 
\end{cases}
\eea
\end{enumerate}
\end{definition}

\begin{remark}
Note that  if  $f \in \sk_p$, for $p=1,2$ and $j$ is even, then  $\dkb^j f\in \sk_p$, and if  $j$ is odd then $\dkb^j f\in \sk_{p-1}$. Also, if $f \in \sk_0$   is a  complex scalar and $j$ is even, then  $\dkb^j f\in \sk_0$, and if  $j$ is odd then $\dkb^j f\in \sk_{1}$. 
\end{remark}


\subsection{Commutation formulas   with $\nabc_3, \nabc_4$}


\begin{lemma}\lab{lemma:commutation-complexM6}    
The following holds true for any  $\psi\in \sk_{p}(\CCC)$ 
\beaa
\, [\nabc_4, r  \DDd_p] \psi&=& O(r^{-2}) \dk^{\le 1} \psi + \Ga_g   \c \dk^{\leq 1} \psi,\\
 \, [\nabc_3, r  \DDd_p] \psi&=&O(ar^{-1}) \nab_3 \psi 
+O(r^{-2}) \dk^{\le 1} \psi  + r \Ga_b' \c \nab_3\psi +     \Ga_b   \c \dk^{\leq 1} \psi. 
\eeaa
Also, The following holds true for any  $\psi\in \sk_{p-1}(\CCC)$
\beaa
\, [\nabc_4, r  \DDs_p\, ] \psi&=& O(r^{-2}) \dk^{\le 1} \psi + \Ga_g   \c \dk^{\leq 1} \psi,\\
 \, [\nabc_3, r  \DDs_p\, ] \psi&=&O(ar^{-1}) \nab_3 \psi 
+O(r^{-2}) \dk^{\le 1} \psi  + r \Ga_b' \c \nab_3 \psi +     \Ga_b   \c \dk^{\leq 1} \psi. 
\eeaa
Finally, the following holds true for any  $\psi\in \sk_{p}(\CCC)$ 
\beaa
\, [\nabc_4, \La_p] \psi&=& O(r^{-2} )  \dk^{\le 2 } \psi+ \dk^{\le 1 }\big( \Ga_g   \c \dk^{\leq 1} \psi \big),\\
\, [\nabc_3, r^2  \La_p] \psi&=& O(ar^{-1}) \dk\nab_3 \psi     + O(r^{-2} )  \dk^{\le 2 } \psi +  \dk^{\le 1 }\big(r \Ga_b' \c\nab_3\psi \big)  + \dk^{\le 1 }\big( \Ga_b   \c \dk^{\leq 1} \psi \big),
\eeaa
where $\La_p = r^2\De_p$.
\end{lemma}

\begin{proof} 
Straightforward verification based on the commutation formulas  obtained in Lemma \ref{COMMUTATOR-NAB-C-3-DD-C-HOT}, and using the fact that $\Xi\in r^{-1}\Ga_g$ in Part III in view of  \eqref{eq:specialidentityforthegloablframeofMMinpartIII} so that the error terms $r\Xi\c\nab_3\psi$ in commutations with $\nabc_4$ satisfy $r\Xi\c\nab_3\psi=\Ga_g\c\dk\psi$.
\end{proof}

The following corollary follows from Lemma \ref{lemma:commutation-complexM6}  and Definition \ref{Def:angularDerivatives-M8}.
\begin{corollary}
\lab{cor:commutation-complexM6}   
We have, for $j\geq 1$, 
\beaa
\,[\nabc_4, \dkb^j] \psi&=&  O(r^{-2} )  \dk^{\le j } \psi+ \dk^{\le j-1 }\big( \Ga_g   \c \dk^{\leq 1} \psi \big),\\
\,[\nabc_3, \dkb^j] \psi &=& O(ar^{-1} )\dk^{\le j-1  }\nab_3 \psi      +O(r^{-2} )  \dk^{\le j } \psi +         \dk^{\le j-1 }\big(r \Ga_b' \c \nab_3  \psi \big)\\
&&+
    \dk^{\le j-1 }\big( \Ga_b   \c \dk^{\leq 1} \psi \big).
\eeaa
\end{corollary}


\subsection{Commutation formulas  with Hodge operators}
\lab{sec:commutationformulasHodgeoperators:chap16}

 
 We start with the following lemma.
 \begin{lemma}
 \lab{Lemma:HodgeFormulas-M8}
 The following identities  hold true for complex anti selfadjoint   tensors $ \psi\in \sk_2(\CCC)$
\beaa
\DDd_2 \DDs_2 \,\psi = \DDs_1\,\DDd_1\psi  - \frac{4a i \cos \th }{|q|^2}\Lieb_\T \psi -\frac{2}{r^2} + O(ar^{-3} ) \dk^{\leq 1}\psi  +  r^{-1}\Ga_b' \c \dk^{ \leq 1 }\psi.
\eeaa
 \end{lemma}
 
 \begin{proof}
 In view of Lemma \ref{eq:dcalident-complex}, we   have
 \beaa
\DDd_2 \DDs_2 \,\psi = \DDs_1\,\DDd_1\psi - i   (\atrch\nab_3+\atrchb \nab_4) \psi - 2 \Kh\psi.
\eeaa
On the other hand, in view of \eqref{eq:T-e_3e_4nab}, 
\beaa
\atrch e_3+\atrchb e_4+ 2(\eta+\etab) \c \dual \nab&=& \frac{4a\cos\th}{|q|^2} \T+  r^{-1}  \Ga_b'   \c \dk,
\eeaa
where we used in particular the fact that $\Hc\in \Ga_b'$, $\trXc\in \Ga_g'$, and the identification $\Ga_g'=r^{-1}\Ga_b'$ valid in Part III. Thus,
\beaa
\atrch e_3+\atrchb e_4&=&\frac{4a\cos\th}{|q|^2} \T+ O(ar^{-3} ) \dkb  +  r^{-1}  \Ga_b'   \c \dk.
\eeaa
Also, in view of \eqref{eq:definition-K-in}, we have 
\beaa
\Kh&=&- \frac 14  \trch \trchb-\frac 1 4 \atrch \atrchb+\frac 1 2 \chih \c \chibh-  \rho\\
&=& \frac{1}{r^2}+O(a^2r^{-4})+r^{-1}\Ga_g+\Ga_b\c\Ga_g.
\eeaa
Hence
\beaa
\DDd_2 \DDs_2 \,\psi = \DDs_1\,\DDd_1\psi  - \frac{4a i \cos \th }{|q|^2}\Lieb_\T \psi -\frac{2}{r^2} + O(ar^{-3} ) \dk^{ \leq 1 }\psi  +  r^{-1} \Ga_b' \c \dk^{ \leq 1 }\psi
\eeaa
as stated.
 \end{proof}
 
 \begin{lemma}
 \lab{Lemma:commHodge-LapM8}
  The following commutation formulas hold true
 \beaa
 -\DDd_p \lap_p +\lap_{p-1} \DDd_p&=&  O(ar^{-3} ) \dkb^{\le 1 } \Lieb_\T +  O(r^{-4} ) \dkb^2 +       O(r^{-3} ) \dkb^{\le 1}  +  r^{- 2 } \dk^{\le 1} \big( \Ga_b' \c \dk \big),\\
 -\DDs_p \, \lap_{p-1}  +\lap_{p-1} \DDs_p\, &=&  O(ar^{-3} ) \dkb^{\le 1}  \Lieb_\T +     O(r^{-4} ) \dkb^2 +       O(r^{-3} ) \dkb^{\le 1} + r^{- 2 } \dk^{\le 1} \big( \Ga_b' \c \dk \big).
 \eeaa
 \end{lemma}
 
 \begin{proof}
 The proof can be derived in the same manner as  the proof of Lemma \ref{Lemma:HodgeFormulas-M8}.
 \end{proof}
 
We deduce the following corollary. 
 \begin{corollary}
 \lab{Corr:commHodge-LapM8}
 We have, for $p=0,1,2$, and for any positive integer $k$,
 \bea
  \lab{eq:commLap^kHodge-even} 
  \bsplit
   \La_{p-1}^k   \DDd_{p}   =&\DDd_p\La_p^k+
 O(a r^{-1} )    \dkb^{\le 2k -1}  \Lieb_\T +O(r^{-2} ) \dkb^{2k}  +O(r^{-1} ) \dkb^{\le 2k-1}\\
 &   +  \dk^{\le 2k-1}\big( \Ga_b' \c \dk \big),\\
 \La_p^k\DDs_{p} \,  =& \DDs_p\,  \La_{p-1}^k+
 O(a r^{-1} )   \dkb^{\le 2k -1}  \Lieb_\T  +O(r^{-2} ) \dkb^{2k}  +O(r^{-1} ) \dkb^{\le 2k-1}\\
 & +  \dk^{\le 2k-1} \big( \Ga_b' \c \dk \big).
 \end{split}
 \eea 
 \end{corollary}
 
 \begin{proof} 
 According to Lemma \ref{Lemma:commHodge-LapM8} and the fact that $\nab(r)\in r\Ga_g$, we have
 \beaa
  (r^2  \lap_{p-1} )   \DDd_{p}   &=&\DDd_p  (r^2 \lap_p)+
 O(a r^{-1} ) \dkb^{\le 1}  \Lieb_\T +  O(r^{-2} ) \dkb^2 +  O(r^{-1} ) \dkb^{\le 1}  +  \dk^{\le 1}\big(  \Ga_b'  \c \dk \big),\\
   (r^2  \lap_{p-1} )   \DDs_{p} \,  &=&\DDs_p\,  (r^2 \lap_p)+
 O(a r^{-1} ) \dkb^{\le 1}  \Lieb_\T +  O(r^{-2} ) \dkb^2 +  O(r^{-1} ) \dkb^{\le 1} +  \dk^{\le 1} \big(  \Ga_b'  \c \dk \big),
 \eeaa
 which corresponds to the case $k=1$ since $\La_p = r^2\De_p$. The general case can be proved recursively by taking also into account the commutation formula $[ \Lieb_\T, \dkb]=  \Ga_b\c  \dk$ and the fact that $\T(r)\in r\Ga_b$.
  \end{proof}
 
 \begin{lemma}
 \lab{Lemma:commdkb^j-Hodge}
For an integer $j\geq 1$ and $\psi\in\sk_p$, $p=0,1,2$, let  $\err[\dkb^j, \psi]$ denote a schematic expression of the following type
 \beaa
 \err[\dkb^j, \psi] &=& O(a r^{-1} ) 
     \dkb^{\le j-1 }\Lieb_\T \psi +O(r^{-2} ) \dkb^j\psi +O(r^{-1} ) \dkb^{\le j-1} \psi  +  \dk^{\le j-1}\big(  \Ga_b'  \c \dk \psi \big).
 \eeaa
 Then, the following commutation identities hold true:
 \begin{enumerate}
 \item  
 If $\psi \in \sk_2$  and $j=2k$, we have  $\dkb^j \psi\in \sk_2$
  and 
 \beaa
 \dkb^j (\DDd_2 \psi) =\DDd_2 (\dkb^j \psi )+ \err[\dkb^j, \psi].
 \eeaa
 
 \item If $\psi \in \sk_2$  and $j=2k+1$, we have $\dkb^j\psi\in \sk_1$ and  
 \beaa
 \dkb^j (\DDd_2 \psi) &=&  \DDd_1 (\dkb^j \psi ) +  \err[\dkb^j, \psi].
 \eeaa
 
 \item  If  $\psi \in \sk_1$  and $j=2k$, we have  $\dkb^j \psi\in \sk_1$
  and
  \beaa
  \bsplit
   \dkb^j (\DDd_1 \psi) & = \DDd_1(\dkb^j\psi) +  \err[\dkb^j, \psi],\\
     \dkb^j (\DDs_2 \psi) & = \DDs_2\, (\dkb^j\psi) +  \err[\dkb^j, \psi].
     \end{split}
  \eeaa
  
   \item  If  $\psi \in \sk_1$  and $j=2k+1$, we have $\dkb^j \psi \in \sk_0$ and 
    \beaa
    \bsplit
    \dkb^j (\DDd_1 \psi) &=  \DDs_1\, (\dkb^j \psi)  +  \err[\dkb^j, \psi],\\
     \dkb^j (\DDs_2 \psi) &= \DDs_1(\dkb^j  \psi ) +  \err[\dkb^j, \psi].
     \end{split}
    \eeaa
    
    \item  If  $\psi \in \sk_0$  and $j=2k$,  we have 
     $\dkb^j \psi \in \sk_0$ and
     \beaa
 \dkb^j (\DDs_1 \psi)    = \DDs_1 \, ( \dkb^j \psi ) +  \err[\dkb^j, \psi].
     \eeaa
     
   \item If  $\psi \in \sk_0$  and $j=2k+1$, we have 
    $\dkb^j \psi \in \sk_1$ and 
    \beaa
    \dkb^j( \DDs_1\, \psi)&=& \DDd_1( \dkb^j \psi) +  \err[\dkb^j, \psi].
    \eeaa
    \end{enumerate}
 \end{lemma} 
 
 \begin{proof}
 The case of even $j$ is an immediate consequence of Corollary 
 \ref{Corr:commHodge-LapM8}.    The case $j=2k+1$   can be checked as follows.
 
  {\bf Case 1.}  Let  $\psi\in \sk_{p}(\CCC)$ with $p=2$      and consider  the expression $ \dkb^j(\DDd_2 \psi)=r \DDd_1  \La^k_{p-1} ( \DDd_2\psi )$ with $j= 2k+1$. Commuting  with the help of Corollary 
 \ref{Corr:commHodge-LapM8}, we find
  \beaa
 \dkb^j( \DDd_2\psi)&=& r \DDd_1  \La^k_{p-1}  \DDd_2\psi \\
  &=& r \DDd_1 \Big(  \DDd_2 \La_p^k \psi +
 O(a r^{-1} )    \dkb^{\le 2k -1}  \Lieb_\T \psi +O(r^{-2} ) \dk^{2k} \psi +O(r^{-1} ) \dk^{\le 2k-1} \psi \\
 && +\dk^{\le 2k-1}\big(  \Ga_b'  \c \dk \psi \big)\Big)\\
 &=&  r \DDd_1( r^{-1} \dkb^j \psi) +  \err[\dkb^j, \psi] \\
 &=&  \DDd_1\dkb^j \psi +   \err[\dkb^j, \psi] 
 \eeaa
 as stated, where we used the fact that $\nab(r)\in r\Ga_g$.
  
{\bf  Case 2.}   If   $\psi\in \sk_{p}(\CCC)$ with $p=1$,    we derive, using Corollary 
 \ref{Corr:commHodge-LapM8} and Lemma \ref{Lemma:HodgeFormulas-M8},    
      \beaa
   \dkb^j( \DDs_2\, \psi)&=& (  r \DDd_2)  \La^k_{p+1}   \DDs_2 \psi = r   \DDd_2 \DDs_2 \,  \La^k_{p} \psi +  \err[\dkb^j, \psi] \\
     &=&  r  \DDs_1\DDd_1  \La^k_{p} \psi +  \err[\dkb^j, \psi]  =  \DDs_1\big( r\DDd_1  \La^k_{p} \psi \big) +  \err[\dkb^j, \psi] \\
     &=&  \DDs_1(\dkb^j  \psi ) +  \err[\dkb^j, \psi] 
   \eeaa
 as stated.   Also, using again Corollary \ref{Corr:commHodge-LapM8},
 \beaa
  \dkb^j(  \DDd_1\, \psi)&=&r\DDs_1\,  \La_0^k  ( \DDd_1\, \psi) = r\DDs_1\,   ( \DDd_1\, \La_1^k  \psi) +  \err[\dkb^j, \psi] \\
     &=& \DDs_1\, ( r\DDd_1\, \La_1^k  \psi) +  \err[\dkb^j, \psi]  = \DDs_1\, (\dkb^j \psi) +  \err[\dkb^j, \psi] 
 \eeaa
 as stated.
 
 {\bf Case 3.}  If $\psi\in \sk_0(\CCC)$, we derive, using again Corollary \ref{Corr:commHodge-LapM8}, 
 \beaa
 \dkb^j( \DDs_1\, \psi)&=& r\DDd_1\La_1^k  \DDs_1\, \psi =  r \DDd_1   (\DDs_1 \La_0^k\psi) +  \err[\dkb^j, \psi] \\
      &= & \DDd_1 (r\DDs_1 \La_0^k\psi) +  \err[\dkb^j, \psi]  = \DDd_1( \dkb^j \psi) +  \err[\dkb^j, \psi] 
 \eeaa
 as stated. This concludes the proof of Lemma \ref{Lemma:commdkb^j-Hodge}.
 \end{proof}

 
 \subsection{Bianchi pairs for  higher angular derivatives}
 
 
 We apply the results  obtained in section \ref{sec:commutationformulasHodgeoperators:chap16} to the main Bianchi pairs verified by the curvature components $A, B, \Pc, \Bb, \Ab$.

 \begin{proposition}
 \lab{Proposition:BianchPairsAngulaDerivs}
 The following assertions hold true for\footnote{\lab{footnote:dkb^{-1}}Note that  for $j=0$ the terms    $O(ar^{-1} )\dkb^{j-1} \Lieb_\T$ are  not present. By convention we set  $\dkb^{-1} =0$.}       $j \ge 1$:

{\bf First Bianchi Pair.}  \,
  Consider the first Bianchi pair  \eqref{eq:first-pair-A-B-lin}  verified by  $A, B$ and  define 
   the     new quantities $ \At=  \dkb ^{j}  A,   \,    \Bt =  \dkb ^{j} B$.  
 \begin{itemize}
 \item If $j $ is  even,  then  $\At\in \sk_2$, $\Bt\in \sk_1$,  and  we have
 \beaa
 \nabc_3 \At+ \frac 1 2  \tr \Xb \At &=&-\DDs_2\, \Bt+ 
 O(ar^{-1} )\dkb^{j-1}  \Lieb_\T  B +\Ft_{(1)},\\
 \nabc_4 \Bt+  2    \ov{\tr X}  \Bt  &=&\DDd_2\, \At+ 
 O(ar^{-1} ) \dkb^{j-1} \Lieb_\T A +\Ft_{(2)}.
 \eeaa
 Note that  this   can be written in the form of  the generalized Bianchi pair \eqref{eq:modelbainchipairequations11-simple:gen} with 
 $\Psi_{(1)}=\At, \Psi_{(2)}=\Bt$.
  
 \item If $j $ is  odd,   $\At\in \sk_1$, $\Bt\in \sk_0$,  and  we have
 \beaa
 \nabc_3 \At + \frac 1 2  \tr \Xb \At&=&-\DDs_1\,\Bt+  O(ar^{-1} ) \dkb^{j-1}  \Lieb_\T  B+\ \Ft_{(1)},\\
 \nabc_4 \Bt+  2   \ov{\tr X}  \Bt  &=&\DDd_1\, \At +  O(ar^{-1} ) \dkb^{j-1}  \Lieb_\T  A+\Ft_{(2)}.
 \eeaa
  Note that  this   can also  be written in the form of  the generalized Bianchi pair \eqref{eq:modelbainchipairequations11-simple:gen} with 
 $\Psi_{(1)}=\At, \Psi_{(2)}=\Bt$.
 \end{itemize} 
 In both cases we have
\beaa
\Ft_{(1)} &=& O(r^{-3} )\dk^{\leq j}\Ga_g  +O(r^{-2} )  \dk^{\le j }\big(A, B\big) +O(r^{-1} ) \dk^{\le j-1}B  +\dk^{\leq j}\big(\Ga_b\c(A, B)\big), \\
\Ft_{(2)} &=& O(r^{-3} )\dk^{\leq j}\Ga_g +O(r^{-2} )  \dk^{\le j }\big(A, B\big) +O(r^{-1} ) \dk^{\le j-1}A  +  \dk^{\le j} \big( \Ga_b \c (A, B)\big).
\eeaa
 
 {\bf Second  Bianchi Pair.}  \, Consider the  second  Bianchi pair   \eqref{eq:BianchiPair2-linearized.M} verified by  $B, \Pc$, and  define 
   the     new quantities\footnote{See remark \ref{remark:funnysecondpair} for the  definition of  $\Pt$ here.}  $ \Bt=  \dkb ^{j}  B,   \,    \Pt =  \dkb ^{j} \ov{\Pc} $.  
 \begin{itemize}
 \item If $j $ is  even,   $\Bt\in \sk_1$, $\Pt\in \sk_0$,  and  we have
 \beaa
 \nabc_3 \Bt+   \tr \Xb \Bt &=&-\DDs_1\, \Pt+ 
 O(ar^{-1} )\dkb^{j-1}  \Lieb_\T  \Pc  +\Ft_{(1)},\\
 \nabc_4 \Pt+  \frac 3 2    \ov{\tr  X} \Pt  &=&\DDd_1\, \Bt+ 
 O(ar^{-1} ) \dkb^{j-1} \Lieb_\T B +\Ft_{(2)}.
 \eeaa
  Note that  this   can be written in the form of  the generalized Bianchi pair \eqref{eq:modelbainchipairequations11-simple:gen} with 
 $\Psi_{(1)}=\Bt, \Psi_{(2)}=\Pt$.

 \item If $j $ is  odd,   $\Bt\in \sk_0$, $\Pt\in \sk_1$,  and  we have
 \beaa
 \nabc_3 \Bt +   \ov{\tr \Xb} \Bt&=&-\DDd_1\,\Pt+  O(ar^{-1} ) \dkb^{j-1}  \Lieb_\T  \Pc +\ \Ft_{(1)},\\
 \nabc_4 \Pt+  \frac 3 2  \tr  X \Pt  &=&\DDs_1\, \Bt +  O(ar^{-1} ) \dkb^{j-1}  \Lieb_\T  B   +\Ft_{(2)}.
 \eeaa
Note that  this   can be written in the form of  the generalized Bianchi pair \eqref{eq:modelbainchipairequations12-simple:gen} with 
  $\Psi_{(1)}=\Bt, \Psi_{(2)}=\Pt$.
 \end{itemize} 
  In both cases we have
\beaa
\Ft_{(1)} &=& O(r^{-3} )\dk^{\leq j}\Ga'_b+O(r^{-2} )  \dk^{\le j }\big(B, \Pc\big) +O(r^{-1} ) \dk^{\le j-1}\Pc   +r^{-2}\dk^{\leq j}(\Ga_b'\c\Rc_b)\\
&& +  \dk^{\le j} \big( \Ga_b \c(A, B)\big), \\
\Ft_{(2)} &=& O(r^{-4})\dk^{\leq j}\Ga_b+O(r^{-2} )  \dk^{\le j }\big(B, \Pc\big) +O(r^{-1} ) \dk^{\le j-1}B  +  \dk^{\le j} \big( \Ga_b' \c (A, B)\big)\\
&& +r^{-1}\dk^{\leq j}(\Ga_b\c\Pc) +r^{-1}\dk^{\leq j}(\Xi\c\Rc_b).
\eeaa    
    
   {\bf Third   Bianchi Pair.}  Consider the   third   Bianchi pair  
   \eqref{eq:BianchiPair3-linearized.M} verified by $\Pc, \Bb$,  
    and  define   the     new quantities $ \Pt=  \dkb ^{j}  \Pc ,   \,    \Bbt =  \dkb ^{j} \Bb$. 
    \begin{itemize}
 \item If $j $ is  even,   $\Pt\in \sk_0$,  $\Bbt\in \sk_1$,  and  we have
 \beaa
 \nabc_3 \Pt+ \frac 3 2 \ov{ \tr \Xb} \, \Pt &=&-\DDd_1\, \Bbt+ 
 O(ar^{-1} )\dkb^{j-1}  \Lieb_\T   \Bb   +\Ft_{(1)},\\
 \nabc_4 \Bbt+     \tr  X\,  \Bbt  &=&\DDs_1\, \Pt+ 
 O(ar^{-1} ) \dkb^{j-1} \Lieb_\T \Pc  +\Ft_{(2)}.
 \eeaa
 Note that  this   can be written in the form of  the generalized Bianchi pair \eqref{eq:modelbainchipairequations12-simple:gen} with 
 $\Psi_{(1)}=\Pt, \Psi_{(2)}=\Bbt$. 
 
 \item If $j $ is  odd,  $\Pt\in \sk_1$, $\Bbt\in \sk_0$,  and  we have
 \beaa
 \nabc_3  \Pt  + \frac 3 2 \tr \Xb  \Pt  &=&- \DDs_1\, \Bbt +  O(ar^{-1} ) \dkb^{j-1}  \Lieb_\T \Bb  +\ \Ft_{(1)},\\
 \nabc_4  \Bbt +  \ov{\tr X \Bbt}  &=& \DDd_1\, \Pt  +  O(ar^{-1} ) \dkb^{j-1}  \Lieb_\T \Pc    +\Ft_{(2)}.
 \eeaa
Note that  this   can be written in the form of  the generalized Bianchi pair \eqref{eq:modelbainchipairequations11-simple:gen} with 
  $\Psi_{(1)}=\Pt, \Psi_{(2)}=\Bbt$.
  \end{itemize} 
 In both cases we have
\beaa
\Ft_{(1)} &=& O(r^{-3})\dk^{\leq j}\Ga_b +O(r^{-2} )   \dk^{\le j }\big(\Pc, \Bb\big) +O(r^{-1} ) \dk^{\le j-1}\Bb    +r^{-1}\dk^{\leq j}\big( \Ga_b' \c \Rc_b\big), \\
\Ft_{(2)} &=& O(r^{-3})\dk^{\leq j}\Ga_b+O(r^{-2} )    \dk^{\le j }\big(\Pc, \Bb\big) +O(r^{-1} ) \dk^{\le j-1}\Pc +r^{-2}\dk^{\leq j}\big(\Ga_b\c \Rc_b\big)\\
&& +\dk^{\leq j}\big(\Xi\c \Rc_b\big).
\eeaa
 
{\bf Fourth    Bianchi Pair.} 
   Consider the    fourth  Bianchi pair  
   \eqref{eq:fourth-pair-Bb-Ab-lin} verified by $\Bb, \Ab$,  
    and  define   the     new quantities $ \Bbt=  \dkb ^{j}  \Bb ,   \,    \Abt =  \dkb ^{j} \Ab$. 
    \begin{itemize}
    \item   If $j $ is  even,   $ \Bbt \in \sk_1$,  $\Abt\in \sk_2$, 
     and  we have
     \beaa
     \nabc_3 \Bbt + 2\ov{\tr \Xb}\,\Bbt&=&-\DDd_2 \Abt +O(ar^{-1})  \dkb^{j-1} \Lieb_\T \Ab +F_{(1)}, \\
     \nabc_4\Abt+\frac 1 2 \tr  X \Abt&=&\DDs_2\,  \Bbt +O(ar^{-1}) \dkb^{j-1}      \Lieb_\T \Bb+F_{(2)}.
     \eeaa
      Note that  this   can   be written in the form of  the generalized Bianchi pair \eqref{eq:modelbainchipairequations12-simple:gen} with 
 $\Psi_{(1)}=\Bbt, \Psi_{(2)}=\Abt$. 
 
     \item   If $j $ is  even,   $ \Bbt \in \sk_0$,  $\Abt\in \sk_1$,  and  we have
       \beaa
     \nabc_3 \Bbt + 2\ov{\tr \Xb}\,\Bbt&=&-\DDd_1 \Abt +O(ar^{-1})  \dkb^{j-1} \Lieb_\T \Ab +F_{(1)}, \\
     \nabc_4\Abt+\frac 1 2 \tr  X \Abt&=&\DDs_1\,  \Bbt +O(ar^{-1}) \dkb^{j-1}      \Lieb_\T \Bb+F_{(2)}. 
     \eeaa
     Note that  this   can  also   be written in the form of  the generalized Bianchi pair \eqref{eq:modelbainchipairequations12-simple:gen} with 
 $\Psi_{(1)}=\Bbt, \Psi_{(2)}=\Abt$. 
    \end{itemize}
    In both cases
  \beaa
\Ft_{(1)} &=& O(r^{-3})\dk^{\leq j}\Ga_b +O(r^{-2} )\dk^{\le j }\big(\Bb, \Ab\big) +O(r^{-1} ) \dk^{\le j-1}\Bb  +\dk^{\leq j}( \Ga_b' \c\Rc_b), \\
\Ft_{(2)} &=& O(r^{-3})\dk^{\leq j}\Ga_b +O(r^{-2} )  \dk^{\le j }\big(\Bb, \Ab\big) +O(r^{-1} ) \dk^{\le j-1}\Bb   +r^{-1}\dk^{\leq j}(\Ga_b\c\Rc_b).
\eeaa
 \end{proposition}

\begin{proof}
We start with the first and the second Bianchi pairs, which, in view of Proposition \ref{prop:structureBainchPairs}, are of the general form \eqref{eq:modelbainchipairequations11-simple:gen}, i.e.
\beaa
\begin{split}
\nabc_3(\Psi_{(1)})+c_{(1) }\tr \Xb\Psi_{(1)} &= -\DDs_p\, \Psi_{(2)}  +F_{(1)},\\[2mm]
\nabc_4(\Psi_{(2)})+c_{(2)} \ov{\tr X}\Psi_{(2)} &= \DDd_p\, \Psi_{(1)}   +F_{(2)}.
\end{split}
\eeaa
We consider the case $j$ even and commute with $\dkb^j$. We obtain
\beaa
\begin{split}
\nabc_3(\dkb^j\Psi_{(1)})+c_{(1) }\tr \Xb\dkb^j\Psi_{(1)} &= -\DDs_p\, \dkb^j\Psi_{(2)}  +\Ft_{(1)}',\\[2mm]
\nabc_4(\dkb^j\Psi_{(2)})+c_{(2)} \ov{\tr X}\dkb^j\Psi_{(2)} &= \DDd_p\, \dkb^j\Psi_{(1)}   +\Ft_{(2)}',
\end{split}
\eeaa
where 
\beaa
\Ft'_{(1)} &=& \dkb^j F_{(1)}+[\nabc_3, \dkb^j]\Psi_{(1)} -  [\dkb^j, \DDs_p\,]\Psi_{(2)}+  O(r^{-2}) \dkb^{\leq j-1}\Psi_{(1)} +\dkb^{\leq j}(\Ga_g\c\Psi_{(1)}), \\
\Ft'_{(2)} &=& \dkb^j F_{(2)}+[\nabc_4, \dkb^j]\Psi_{(2)} -  [\dkb^j, \DDd_p\,]\Psi_{(1)} +  O(r^{-2}) \dkb^{\leq j-1}\Psi_{(2)}+\dkb^{\leq j}(\Ga_g\c\Psi_{(2)}).
\eeaa
Using  Corollary \ref{cor:commutation-complexM6}   and   Lemma \ref{Lemma:commdkb^j-Hodge}, we infer
\beaa
\Ft'_{(1)} &=& \dkb^j F_{(1)} +O(ar^{-1} )\dk^{\le j-1  }\nab_3\Psi_{(1)}     +O(r^{-2} )  \dk^{\le j }\Psi_{(1)}+         \dk^{\le j-1 }\big(r\Ga_b'\c \nab_3\Psi_{(1)}\big)\\
&&+ \dk^{\le j-1 }\big( \Ga_b   \c \dk^{\leq 1}\Psi_{(1)}\big) +O(a r^{-1} )\dkb^{\le j-1 }\Lieb_\T\Psi_{(2)} +O(r^{-2} ) \dkb^j\Psi_{(2)}\\
&&+O(r^{-1} ) \dkb^{\le j-1}\Psi_{(2)} +  \dk^{\le j-1} \big( \Ga_b' \c\dk\Psi_{(2)}\big) +  O(r^{-2}) \dkb^{\leq j-1}\Psi_{(1)} +\dkb^{\leq j}(\Ga_g\c\Psi_{(1)}), \\
\Ft'_{(2)} &=& \dkb^j F_{(2)} +O(r^{-2} )  \dk^{\le j }\Psi_{(2)}+ \dk^{\le j-1 }\big( \Ga_g   \c \dk^{\leq 1}\Psi_{(2)} \big)  +O(a r^{-1} )\dkb^{\le j-1 }\Lieb_\T\Psi_{(1)}\\
&& +O(r^{-2} ) \dkb^j\Psi_{(1)}+O(r^{-1} ) \dkb^{\le j-1}\Psi_{(1)}+  \dk^{\le j-1} \big( \Ga_b' \c \dk\Psi_{(1)}\big) +  O(r^{-2}) \dkb^{\leq j-1}\Psi_{(2)}\\
&&+\dkb^{\leq j}(\Ga_g\c\Psi_{(2)}).
\eeaa
We obtain
\beaa
\begin{split}
\nabc_3(\dkb^j\Psi_{(1)})+c_{(1) }\tr \Xb\dkb^j\Psi_{(1)} &= -\DDs_p\, \dkb^j\Psi_{(2)} +O(a r^{-1} )\dkb^{\le j-1 }\Lieb_\T\Psi_{(2)} +\Ft_{(1)},\\[2mm]
\nabc_4(\dkb^j\Psi_{(2)})+c_{(2)} \ov{\tr X}\dkb^j\Psi_{(2)} &= \DDd_p\, \dkb^j\Psi_{(1)} +O(a r^{-1} )\dkb^{\le j-1 }\Lieb_\T\Psi_{(1)}  +\Ft_{(2)},
\end{split}
\eeaa
where 
\beaa
\Ft_{(1)} &=& \dkb^j F_{(1)} +O(ar^{-1} )\dk^{\le j-1  }\nab_3\Psi_{(1)}     +O(r^{-2} )  \dk^{\le j }\big(\Psi_{(1)}, \Psi_{(2)}\big) +O(r^{-1} ) \dk^{\le j-1}\Psi_{(2)}\\
&&+         \dk^{\le j-1 }\big(r\Ga_b'\c \nab_3\Psi_{(1)}\big)   +\dk^{\leq j}\big(\Ga_b\c\Psi_{(1)}\big)+\dk^{\leq j}\big(\Ga_b'\c\Psi_{(2)}\big), \\
\Ft_{(2)} &=& \dkb^j F_{(2)} +O(r^{-2} )  \dk^{\le j }\big(\Psi_{(1)}, \Psi_{(2)}\big) +O(r^{-1} ) \dk^{\le j-1}\Psi_{(1)}  +  \dk^{\le j} \big( \Ga_b' \c\Psi_{(1)}\big)\\
&& +r^{-1}\dk^{\leq j}(\Ga_b\c\Psi_{(2)}).
\eeaa
In view of $\nabc_3\Psi_{(1)} = -c_{(1) }\tr \Xb\Psi_{(1)}  -\DDs_p\, \Psi_{(2)}  +F_{(1)}$, we have
\beaa
\nab_3\Psi_{(1)} &=& O(r^{-1})\Psi_{(1)} +O(r^{-1})\dk^{\leq 1}\Psi_{(2)} +\Ga_b\c\Psi_{(1)} +F_{(1)}
\eeaa 
which we use to remove the $\nab_3\Psi_{(1)}$ terms from $\Ft_{(1)}$. This yields
\beaa
\Ft_{(1)} &=& \dkb^j F_{(1)} +O(ar^{-1} )\dk^{\le j-1}F_{(1)} +O(r^{-2} )  \dk^{\le j }\big(\Psi_{(1)}, \Psi_{(2)}\big) +O(r^{-1} ) \dk^{\le j-1}\Psi_{(2)}\\
&&+         \dk^{\le j-1 }\big(r\Ga_b\c F_{(1)}\big)   +  \dk^{\le j} \big( \Ga_b \c\Psi_{(1)}\big) +\dk^{\leq j}(\Ga_b'\c\Psi_{(2)}), \\
\Ft_{(2)} &=& \dkb^j F_{(2)} +O(r^{-2} )  \dk^{\le j }\big(\Psi_{(1)}, \Psi_{(2)}\big) +O(r^{-1} ) \dk^{\le j-1}\Psi_{(1)}  +  \dk^{\le j} \big( \Ga_b' \c\Psi_{(1)}\big)\\
&& +r^{-1}\dk^{\leq j}(\Ga_b\c\Psi_{(2)}).
\eeaa

We proceed in the same manner for $j$ odd, and for the other Bianchi pairs, and we obtain in all cases that the Bianchi equations have the stated structure, with $\Ft_{(1)}$, $\Ft_{(2)}$ given as above by
\beaa
\Ft_{(1)} &=& \dkb^j F_{(1)} +O(ar^{-1} )\dk^{\le j-1}F_{(1)} +O(r^{-2} )  \dk^{\le j }\big(\Psi_{(1)}, \Psi_{(2)}\big) +O(r^{-1} ) \dk^{\le j-1}\Psi_{(2)}\\
&&+         \dk^{\le j-1 }\big(r\Ga_b\c F_{(1)}\big)   +  \dk^{\le j} \big( \Ga_b \c\Psi_{(1)}\big) +\dk^{\leq j}(\Ga_b'\c\Psi_{(2)}), \\
\Ft_{(2)} &=& \dkb^j F_{(2)} +O(r^{-2} )  \dk^{\le j }\big(\Psi_{(1)}, \Psi_{(2)}\big) +O(r^{-1} ) \dk^{\le j-1}\Psi_{(1)}  +  \dk^{\le j} \big( \Ga_b' \c\Psi_{(1)}\big)\\
&& +r^{-1}\dk^{\leq j}(\Ga_b\c\Psi_{(2)}).
\eeaa

\begin{remark}
In the case where $j$ is odd, for the second and third Bianchi pairs, when implementing the above procedure, we need in addition to use the fact that 
\beaa
\ov{\tr X} &=& \tr X+O(r^{-2})+\Ga_g= \tr X+O(r^{-2})+r^{-1}\Ga_b,\\ 
\ov{\tr\Xb} &=& \tr\Xb+O(r^{-2})+\Ga_g= \tr\Xb+O(r^{-2})+r^{-1}\Ga_b, 
\eeaa
in order to put the LHS in the general form of Definition \ref{definition:GenBianchPairs:gen}. Note that this generates extra terms of the type $O(r^{-2})\dk^{\leq j}(\Psi_{(1)}, \Psi_{(2)})$ and $r^{-1}\dk^{\leq j}\big(\Ga_b\c(\Psi_{(1)}, \Psi_{(2)})\big)$ which are incorporated in $\Ft_{(1)}$, $\Ft_{(2)}$. 
\end{remark}

To conclude, it remains to plug the structure of $F_{(1)}$, $F_{(2)}$ in the above formula for  $\Ft_{(1)}$ and $\Ft_{(2)}$. We do this for all Bianchi pairs starting with the first one for which we have, according to Proposition \ref{prop:structureBainchPairs},
 \beaa
  F_{(1)}&=& O(ar^{-2} ) B+ O(r^{-3} ) \Ga_g + \Ga_b \c B, \\
  F_{(2)}&=& O(a r^{-2} ) A + O(r^{-3} ) \Ga_g+ \Ga_b\c A.
  \eeaa
  In this case, we infer, using also $\Psi_{(1)}=A$, $\Psi_{(2)}=B$, we obtain
  \beaa
\Ft_{(1)} &=& O(r^{-3} )\dk^{\leq j}\Ga_g  +O(r^{-2} )  \dk^{\le j }\big(A, B\big) +O(r^{-1} ) \dk^{\le j-1}B  +\dk^{\leq j}\big(\Ga_b\c(A, B)\big), \\
\Ft_{(2)} &=& O(r^{-3} )\dk^{\leq j}\Ga_g +O(r^{-2} )  \dk^{\le j }\big(A, B\big) +O(r^{-1} ) \dk^{\le j-1}A  +  \dk^{\le j} \big( \Ga_b \c (A, B)\big),
\eeaa
as stated.

Next, we consider the second Bianchi pair for which we have, according to Proposition \ref{prop:structureBainchPairs},
\beaa
   F_{(1)}&=&  O(ar^{-2} ) \Pc +O(r^{-3} )  \Ga'_b + r^{-2} \Ga'_b\c \Rc_b+\Ga_b\c  A,\\
   F_{(2)}&=&   O( a r^{-2}) B  +O(r^{-4})\Ga_b  + \Ga'_b\c (A, B) +r^{-1}\Xi\c\Rc_b.
  \eeaa
    In this case, we infer, using also $\Psi_{(1)}=B$, $\Psi_{(2)}=\ov{\Pc}$, we obtain
\beaa
\Ft_{(1)} &=& O(r^{-3} )\dk^{\leq j}\Ga'_b+O(r^{-2} )  \dk^{\le j }\big(B, \Pc\big) +O(r^{-1} ) \dk^{\le j-1}\Pc   +r^{-2}\dk^{\leq j}(\Ga_b'\c\Rc_b)\\
&& +  \dk^{\le j} \big( \Ga_b \c(A, B)\big), \\
\Ft_{(2)} &=& O(r^{-4})\dk^{\leq j}\Ga_b+O(r^{-2} )  \dk^{\le j }\big(B, \Pc\big) +O(r^{-1} ) \dk^{\le j-1}B  +  \dk^{\le j} \big( \Ga_b' \c (A, B)\big)\\
&& +r^{-1}\dk^{\leq j}(\Ga_b\c\Pc) +r^{-1}\dk^{\leq j}(\Xi\c\Rc_b)
\eeaa
as stated.

Next, we consider the third Bianchi pair for which we have, according to Proposition \ref{prop:structureBainchPairs},
\beaa
 F_{(1)}&=& O(ar^{-2} )\Bb  +O(r^{-3})\Ga_b  + r^{-1}\Ga_b'\c \Rc_b, \\
    F_{(2)}&=&O(ar^{-2}) \Pc   +O(r^{-3})\Ga_b +r^{-2} \Ga_b\c \Rc_b +\Xi\c \Rc_b.
\eeaa
    In this case, we infer, using also $\Psi_{(1)}=\Pc$, $\Psi_{(2)}=\Bb$, we obtain
    \beaa
\Ft_{(1)} &=& O(r^{-3})\dk^{\leq j}\Ga_b +O(r^{-2} )  \dk^{\le j }\big(\Pc, \Bb\big) +O(r^{-1} ) \dk^{\le j-1}\Bb    +r^{-1}\dk^{\leq j}\big(\Ga_b'\c \Rc_b\big), \\
\Ft_{(2)} &=& O(r^{-3})\dk^{\leq j}\Ga_b+O(r^{-2} )  \dk^{\le j }\big(\Pc, \Bb\big) +O(r^{-1} ) \dk^{\le j-1}\Pc    +r^{-2}\dk^{\leq j}\big(\Ga_b\c \Rc_b\big)\\
&& +\dk^{\leq j}\big(\Xi\c \Rc_b\big)
\eeaa
as stated.

Finally, we consider the fourth Bianchi pair for which we have, according to Proposition \ref{prop:structureBainchPairs},
\beaa
 F_{(1)}&=&   O(ar^{-2}) \Ab+  O(r^{-3} ) \Ga_b+\Ga_b'\c\Rc_b, 
 \\
    F_{(2)}&=&  O(ar^{-2}) \Bb + O(r^{-3} ) \Ga_b+r^{-1}\Ga_b\c\Rc_b.
\eeaa
In this case, we infer, using also $\Psi_{(1)}=\Bb$, $\Psi_{(2)}=\Ab$, we obtain
\beaa
\Ft_{(1)} &=& O(r^{-3})\dk^{\leq j}\Ga_b +O(r^{-2} )  \dk^{\le j }\big(\Bb, \Ab\big) +O(r^{-1} ) \dk^{\le j-1}\Bb  +\dk^{\leq j}(\Ga_b'\c\Rc_b), \\
\Ft_{(2)} &=& O(r^{-3})\dk^{\leq j}\Ga_b +O(r^{-2} )  \dk^{\le j }\big(\Bb, \Ab\big) +O(r^{-1} ) \dk^{\le j-1}\Bb   +r^{-1}\dk^{\leq j}(\Ga_b\c\Rc_b)
\eeaa
as stated. This concludes the proof of Proposition \ref{Proposition:BianchPairsAngulaDerivs}. 
\end{proof}

 
\section{Bianchi pairs for general higher derivatives}
\lab{sec:proofofThm:MainResultMext:generalcase:prelim2}

 
  To derive   estimates for all derivatives  $\le J+1$,  we  commute the Bianchi  equations of Proposition \ref{prop:structureBainchPairs}  with  $( \ov{q}\nabc_4)^{j} \dkb ^{j_2} \Lieb_{ \T}^{j_1} $ or\footnote{Depending on the type of Bianchi pair in Definition \ref{definition:GenBianchPairs:gen}.}  $(  q \nabc_4)^{j} \dkb ^{j_2} \Lieb_{ \T}^{j_1} $   for all multi-indices $(j, j_1, j_2)$    with $j+j_{1}+j_{2} = J+1$. The result is stated in the following proposition.  We also make the convention $\dkb^{-1}=0$, see footnote \ref{footnote:dkb^{-1}}.
  
  \medskip
  
  \begin{proposition}
  \lab{Proposition:BianchiPairs-Higher Derivatives}
   The following  higher derivatives Bianchi pairs   hold:
  \begin{enumerate}
  \item Consider the first Bianchi pair in $A, B$,  and   set, for $j+j_{1}+j_{2} = J+1$,
 \beaa
 \At =( \ov{q}\nabc_4)^{j}  \dkb ^{j_2}  \Lieb_{ \T}^{j_1} A, \qquad   \Bt =( \ov{q} \nabc_4)^{j}   \dkb ^{j_2}  \Lieb_{\T}^{j_1}B.
 \eeaa
 Then, with  $c_{(1)}=\frac 1 2, \,  c_{(2)}= 2 $:
  \begin{itemize}
 \item If $j_2 $ is  even,     $\At\in \sk_2$, $\Bt\in \sk_1$,  and  we have\footnote{\lab{footnote:dkb^{-1}qe4j-1}Note that  for $j_2=0$,   the terms in  $O(ar^{-1} ) (\ov{q} \nabc_4)\dkb^{j_2-1}  \Lieb^{j_1+1} _\T $  do not appear at all. Also, for $j=0$, the terms $O(r^{-1} ) (\ov{q} \nabc_4)^{j-1}\dkb^{j_2+1}  \Lieb^{j_1} _\T$ do not appear at all. This holds true for all Bianchi pairs  below.}
 \beaa
 \nabc_3 \At+ \left(c_{(1) } - \frac{j}{2}\right) \tr \Xb \At &=&-\DDs_2\, \Bt 
+  O(ar^{-1} ) (\ov{q} \nabc_4)^{ j }\dkb^{j_2-1}  \Lieb^{j_1+1} _\T  B\\
&& +  O(r^{-1} ) (\ov{q} \nabc_4)^{j-1}\dkb^{j_2+1}  \Lieb^{j_1} _\T(A,B) +\Ft_{(1)},\\
 \nabc_4 \Bt+ \left(c_{(2) } - \frac{j}{2}\right)  \ov{\tr X}  \Bt  &=&\DDd_2\, \At+ 
  O(ar^{-1} ) (\ov{q} \nabc_4)^{ j }\dkb^{j_2-1}  \Lieb^{j_1+1} _\T A\\
  && +  O(r^{-1} ) (\ov{q} \nabc_4)^{j-1}\dkb^{j_2+1}  \Lieb^{j_1} _\T A +\Ft_{(2)}.
 \eeaa
 Note that  this   can be written in the form of  the generalized Bianchi pair \eqref{eq:modelbainchipairequations11-simple:gen} with 
 $\Psi_{(1)}=\At, \Psi_{(2)}=\Bt$.
 
 \item If $j_2 $ is  odd,   $\At\in \sk_1$, $\Bt\in \sk_0$,  and  we have
 \beaa
 \nabc_3 \At +\left(c_{(1) } - \frac{j}{2}\right)  \tr \Xb \At&=&-\DDs_1\,\Bt+  O(ar^{-1} ) (\ov{q} \nabc_4) ^j\dkb^{j_2-1}  \Lieb^{j_1+1} _\T  B\\
 &&+  O(r^{-1} ) (\ov{q} \nabc_4)^{j-1}\dkb^{j_2+1}  \Lieb^{j_1} _\T(A,B) +\ \Ft_{(1)},\\
 \nabc_4 \Bt+  \left(c_{(2) } - \frac{j}{2}\right)   \ov{\tr X}  \Bt  &=&\DDd_1\, \At +  O(ar^{-1} ) (\ov{q} \nabc_4)^j\dkb^{j_2-1}  \Lieb^{j_1+1} _\T   A\\
 && +  O(r^{-1} ) (\ov{q} \nabc_4)^{j-1}\dkb^{j_2+1}  \Lieb^{j_1} _\T A +\Ft_{(2)}.
 \eeaa
  Note that  this   can also  be written in the form of  the generalized Bianchi pair \eqref{eq:modelbainchipairequations11-simple:gen} with 
 $\Psi_{(1)}=\At, \Psi_{(2)}=\Bt$.
 \end{itemize} 
  In both cases we have
 \beaa
\Ft_{(1)} &=& O(r^{-3} )\dk^{\leq J+1}\Ga_g  +O(r^{-2} )  \dk^{\le J+1}\big(A, B\big) +O(r^{-1} ) \dk^{\le J}B  +\dk^{\leq J+1}\big(\Ga_b\c(A, B)\big), \\
\Ft_{(2)} &=& O(r^{-3} )\dk^{\leq J+1}\Ga_g +O(r^{-2} )  \dk^{\le J+1}\big(A, B\big) +O(r^{-1} ) \dk^{\le J}A  +  \dk^{\le J+1} \big( \Ga_b \c (A, B)\big).
\eeaa
 
 {\bf Second  Bianchi Pair.}\, Consider the  second  Bianchi pair in $B, \Pc$,  and set\footnote{See remark \ref{remark:funnysecondpair} for the  definition of  $\Pt$ here.}, for $j+j_1+j_2=J+1$, 
 \beaa
 \Bt =( \ov{q}\nabc_4)^{j}  \dkb ^{j_2}  \Lieb_{ \T}^{j_1} B, \qquad   \Pt =( \ov{q} \nabc_4)^{j}   \dkb ^{j_2}  \Lieb_{\T}^{j_1} \ov{\Pc} .
 \eeaa
Then, with  $c_{(1)}=1, \,  c_{(2)}= 3/2 $:
\begin{itemize}
 \item If $j_2$ is  even,  $\Bt\in \sk_1$, $\Pt\in \sk_0$,  and  we have
 \beaa
 \nabc_3 \Bt+ \left(c_{(1) } - \frac{j}{2}\right)  \tr \Xb \Bt &=&-\DDs_1\, \Pt+  O(ar^{-1} ) (\ov{q} \nabc_4)^j\dkb^{j_2-1}  \Lieb^{j_1+1} _\T  \Pc\\
 && +  O(r^{-1} ) (\ov{q} \nabc_4)^{j-1}\dkb^{j_2+1}  \Lieb^{j_1} _\T(B, \Pc)  +\Ft_{(1)},\\
 \nabc_4 \Pt+  \left(c_{(2) } - \frac{j}{2}\right)    \ov{\tr  X} \Pt  &=&\DDd_1\, \Bt+ 
  O(ar^{-1} ) (\ov{q} \nabc^j_4)\dkb^{j_2-1}  \Lieb^{j_1+1} _\T B\\
  && +  O(r^{-1} ) (\ov{q} \nabc_4)^{j-1}\dkb^{j_2+1}  \Lieb^{j_1} _\T B +\Ft_{(2)}.
 \eeaa
  Note that  this   can be written in the form of  the generalized Bianchi pair \eqref{eq:modelbainchipairequations11-simple:gen} with 
 $\Psi_{(1)}=\Bt, \Psi_{(2)}=\Pt$.

 \item If $j_2$ is  odd,  $\Bt\in \sk_0$, $\Pt\in \sk_1$,  and  we have
 \beaa
 \nabc_3 \Bt + \left(c_{(1) } - \frac{j}{2}\right) \ov{\tr \Xb} \Bt&=&-\DDd_1\,\Pt+  O(ar^{-1} )( q \nabc_4)^j\dkb^{j_2-1}  \Lieb^{j_1+1} _\T \Pc\\
 && +  O(r^{-1} ) (\ov{q} \nabc_4)^{j-1}\dkb^{j_2+1}  \Lieb^{j_1} _\T(B, \Pc) + \Ft_{(1)},\\
 \nabc_4 \Pt+\left(c_{(2) } - \frac{j}{2}\right)   \tr X \Pt  &=&\DDs_1\, \Bt +  O(ar^{-1} ) (q \nabc_4)^j\dkb^{j_2-1}  \Lieb^{j_1+1} _\T  B\\
 && +  O(r^{-1} ) (\ov{q} \nabc_4)^{j-1}\dkb^{j_2+1}  \Lieb^{j_1} _\T B    +\Ft_{(2)}.
 \eeaa
 Note that  this   can be written in the form of  the generalized Bianchi pair \eqref{eq:modelbainchipairequations12-simple:gen} with 
 $\Psi_{(1)}=\Bt, \Psi_{(2)}=\Pt$.
 \end{itemize} 
  In both cases we have
 \beaa
\Ft_{(1)} &=& O(r^{-3} )\dk^{\leq J+1}\Ga'_b+O(r^{-2} )  \dk^{\le J+1}\big(B, \Pc\big) +O(r^{-1} ) \dk^{\le J}\Pc   +r^{-2}\dk^{\leq J+1}(\Ga_b'\c\Rc_b)\\
&& +  \dk^{\le J+1} \big( \Ga_b \c(A, B)\big), \\
\Ft_{(2)} &=& O(r^{-4})\dk^{\leq J+1}\Ga_b+O(r^{-2} )  \dk^{\le J+1}\big(B, \Pc\big) +O(r^{-1} ) \dk^{\le J}B  +  \dk^{\le J+1} \big( \Ga_b' \c (A, B)\big)\\
&& +r^{-1}\dk^{\leq J+1}(\Ga_b\c\Pc) +r^{-1}\dk^{\leq J+1}(\Xi\c\Rc_b).
\eeaa
    
      {\bf Third   Bianchi Pair.}   Consider the third  Bianchi pair
       in $\Pc, \Bb$,  and   set, for $j+j_{1}+j_{2} = J+1$,
 \beaa
     \Pt =( q \nabc_4)^{j}   \dkb ^{j_2}  \Lieb_{\T}^{j_1}\Pc ,   \qquad     \Bbt = ( q \nabc_4)^{j}   \dkb ^{j_2}  \Lieb_{\T}^{j_1} \Bb .
     \eeaa
     Then, with  $c_{(1)}=3/2, \,  c_{(2)}= 1 $:
    \begin{itemize}
 \item If $ j_2 $ is even,   
  $ \Pt \in \sk_0$,   $\Bbt\in \sk_1$,  and we have
 \beaa
 \nabc_3 \Pt+\left(c_{(1) } - \frac{j}{2}\right)\ov{ \tr \Xb} \, \Bbt &=&-\DDd_1\, \Bbt+ 
  O(ar^{-1} ) (q \nabc_4)^j\dkb^{j_2-1}  \Lieb^{j_1+1} _\T  \Bbt  \\
  && +  O(r^{-1} ) (\ov{q} \nabc_4)^{j-1}\dkb^{j_2+1}  \Lieb^{j_1} _\T(\Pc, \Bb) +\Ft_{(1)},\\
 \nabc_4 \Bbt+   \left(c_{(2) } - \frac{j}{2}\right)  \tr  X\,  \Bbt  &=&\DDs_1\, \Pt+ 
 O(ar^{-1} ) (q \nabc_4)^j\dkb^{j_2-1}  \Lieb^{j_1+1} _\T \Pc \\
 && +  O(r^{-1} ) (\ov{q} \nabc_4)^{j-1}\dkb^{j_2+1}  \Lieb^{j_1} _\T\Pc  +\Ft_{(2)}.
 \eeaa
 Note that  this   can be written in the form of  the generalized Bianchi pair \eqref{eq:modelbainchipairequations12-simple:gen} with 
 $\Psi_{(1)}=\Pt, \Psi_{(2)}=\Bbt$. 
 
 \item If $ j_2 $ is  odd,  $\Pt\in \sk_1$, $\Bbt\in \sk_0$,  and  we have
 \beaa
 \nabc_3  \Pt  +\left(c_{(1) } - \frac{j}{2}\right) \ov{ \tr \Xb}  \Pt &=&- \DDs_1\, \Bbt +   O(ar^{-1} ) (q \nabc_4)^j\dkb^{j_2-1}  \Lieb^{j_1+1} _\T \Bb \\
 && +  O(r^{-1} ) (\ov{q} \nabc_4)^{j-1}\dkb^{j_2+1}  \Lieb^{j_1} _\T(\Pc, \Bb) + \Ft_{(1)},\\
 \nabc_4  \Bb + \left(c_{(2) } - \frac{j}{2}\right)   \ov{\tr X}   \Bb  &=& \DDd_1\,  \Pt  +  O(ar^{-1} ) (q \nabc_4)^j\dkb^{j_2-1}  \Lieb^{j_1+1} _\T \Pc \\
 && +  O(r^{-1} ) (\ov{q} \nabc_4)^{j-1}\dkb^{j_2+1}  \Lieb^{j_1} _\T\Pc  +\Ft_{(2)}.
 \eeaa
  Note that  this   can  also be written in the form of  the generalized Bianchi pair  \eqref{eq:modelbainchipairequations11-simple:gen} with 
 $\Psi_{(1)}=\Pt, \Psi_{(2)}=\Bbt$. 
  \end{itemize} 
 In both cases we have
\beaa
\Ft_{(1)} &=& O(r^{-3})\dk^{\leq J+1}\Ga_b +O(r^{-2} )    \dk^{\le J+1}\big(\Pc, \Bb\big) +O(r^{-1} ) \dk^{\le J}\Bb     +r^{-1}\dk^{\leq J+1}\big( \Ga_b' \c \Rc_b\big), \\
\Ft_{(2)} &=& O(r^{-3})\dk^{\leq J+1}\Ga_b+O(r^{-2} )    \dk^{\le J+1}\big(\Pc, \Bb\big) +O(r^{-1} ) \dk^{\le J}\Pc     +r^{-2}\dk^{\leq J+1}\big(\Ga_b\c \Rc_b\big)\\
&& +\dk^{\leq J+1}\big(\Xi\c \Rc_b\big).
\eeaa

{\bf Fourth    Bianchi Pair.} 
 Consider the fourth  Bianchi pair
       in $\Bb, \Ab$,  and   set, for $j+j_{1}+j_{2} = J+1$,
 \beaa
     \Bbt =( q \nabc_4)^{j}   \dkb ^{j_2}  \Lieb_{\T}^{j_1}\Bb ,   \qquad     \Abt = ( q \nabc_4)^{j}   \dkb ^{j_2}  \Lieb_{\T}^{j_1} \Ab.
     \eeaa
  Then,  with  $c_{(1)}=2, \,  c_{(2)}= 1/ 2 $:
    \begin{itemize}
    \item   If $ j_2 $ is  even,   $ \Bbt \in \sk_1$,  $\Abt\in \sk_2$, 
     and  we have
     \beaa
     \nabc_3 \Bbt + \left(c_{(1) } - \frac{j}{2}\right)\ov{\tr \Xb}\,\Bbt&=&-\DDd_2 \Abt +O(ar^{-1}) ( q \nabc_4)^{j}  \dkb^{j_2-1} \Lieb^{j_1+1}_\T \Ab\\
     &&+  O(r^{-1} ) (\ov{q} \nabc_4)^{j-1}\dkb^{j_2+1}  \Lieb^{j_1} _\T(\Bb, \Ab) +\Ft_{(1)}, \\
     \nabc_4\Abt+\left(c_{(2) } - \frac{j}{2}\right) \tr  X \Abt&=&\DDs_2\,  \Bbt +O(ar^{-1}) ( q \nabc_4)^{j}  \dkb^{j_2-1}      \Lieb^{j_1+1}_\T \Bb\\
     && +  O(r^{-1} ) (\ov{q} \nabc_4)^{j-1}\dkb^{j_2+1}  \Lieb^{j_1} _\T\Bb +\Ft_{(2)}. 
     \eeaa
      Note that  this   can   be written in the form of  the generalized Bianchi pair \eqref{eq:modelbainchipairequations12-simple:gen} with 
 $\Psi_{(1)}=\Bbt, \Psi_{(2)}=\Abt$. 
 
     \item   If $ j_2 $ is  even,   $ \Bbt \in \sk_0$,  $\Abt\in \sk_1$,  and  we have
       \beaa
     \nabc_3 \Bbt + \left(c_{(1) } - \frac{j}{2}\right)\ov{\tr \Xb}\,\Bbt&=&-\DDd_1 \Abt +O(ar^{-1})  ( q \nabc_4)^{j} \dkb^{j_2-1} \Lieb^{j_1+1}_\T \Ab\\
     &&+  O(r^{-1} ) (\ov{q} \nabc_4)^{j-1}\dkb^{j_2+1}  \Lieb^{j_1} _\T(\Bb, \Ab) +\Ft_{(1)}, \\
     \nabc_4\Abt+\left(c_{(2) } - \frac{j}{2}\right) \tr  X \Abt&=&\DDs_1\,  \Bbt +O(ar^{-1})( q \nabc_4)^{j}  \dkb^{j_2-1}      \Lieb^{j_1+1}_\T \Bb\\
     && +  O(r^{-1} ) (\ov{q} \nabc_4)^{j-1}\dkb^{j_2+1}  \Lieb^{j_1} _\T\Bb +\Ft_{(2)}. 
     \eeaa
     Note that  this   can  also   be written in the form of  the generalized Bianchi pair \eqref{eq:modelbainchipairequations12-simple:gen} with 
 $\Psi_{(1)}=\Bbt, \Psi_{(2)}=\Abt$. 
    \end{itemize}
    In both cases
 \beaa
\Ft_{(1)} &=& O(r^{-3})\dk^{\leq J+1}\Ga_b +O(r^{-2} )  \dk^{\le J+1}\big(\Bb, \Ab\big) +O(r^{-1} ) \dk^{\le J}\Bb  +\dk^{\leq J+1}( \Ga_b' \c\Rc_b), \\
\Ft_{(2)} &=& O(r^{-3})\dk^{\leq J+1}\Ga_b +O(r^{-2} )  \dk^{\le J+1}\big(\Bb, \Ab\big) +O(r^{-1} ) \dk^{\le J}\Bb   +r^{-1}\dk^{\leq J+1}(\Ga_b\c\Rc_b).
\eeaa
\end{enumerate}
  \end{proposition}
  
 \begin{proof}
 The case $j=j_2=0$  can be easily  be derived  by commuting
   the Bianchi pairs  in Proposition \ref{prop:structureBainchPairs}    with  $\Lieb_\T $,  in view of Lemma  \ref{lemma:commLieb_TM8}. The case  $j=0$  can then be derived 
       by   proceeding exactly as Proposition \ref{Proposition:BianchPairsAngulaDerivs}, starting  with  the results already derived  by commutation with $\Lieb_\T$.  Finally
        we derive the results  for all $(j, j_1, j_2)$ by  making use of the commutators in Lemma \ref{Lemma:commutationwith-qnabc_4}. 
 \end{proof}

 
 \section{Proof of Theorem  \ref{THM:MAINRESULTMEXT}}
\lab{sec:proofofThm:MainResultMext:generalcase}
 
  
We are now in position to prove Theorem  \ref{THM:MAINRESULTMEXT} on the control of $\Rkext_{J+1}$. To this end, we rely on the procedure for deriving basic curvature estimates in $\Mext$ outlined in section \ref{sec:proofofThm:MainResultMext:k=0case} and on the structure of Bianchi pairs for higher derivatives in Proposition \ref{Proposition:BianchiPairs-Higher Derivatives}. The proof proceeds in the following steps.
 
{\bf Step 1.} We first derive estimates for $\Lieb_\T^{J+1}$ derivatives of curvature in $\Mext$. To this end, we consider the Bianchi pairs in Proposition \ref{Proposition:BianchiPairs-Higher Derivatives} in the particular case 
\beaa
j=0, \qquad j_2=0, \qquad j_1=J+1. 
\eeaa
In that case,  the terms  
      $O(ar^{-1} ) (q \nabc_4)^{ j }\dkb^{j_2-1}  \Lieb^{j_1+1} _\T$ and $O(r^{-1} ) (q \nabc_4)^{j-1}\dkb^{j_2+1}  \Lieb^{j_1} _\T$ are not present, see footnote \ref{footnote:dkb^{-1}qe4j-1}, and hence the structure of the Bianchi pairs in Proposition \ref{Proposition:BianchiPairs-Higher Derivatives} are the exact analog for $J+1$ derivatives of the ones in Proposition \ref{prop:structureBainchPairs} used for the proof of the basic curvature estimates of section \ref{sec:proofofThm:MainResultMext:k=0case}, with the exception of 
\begin{itemize}
\item the new terms of type $O(r^{-1})\dk^{\leq J}\Psi_{(1)}$ in $\Ft_{(1)}$,
\item the new terms of type $O(r^{-1})\dk^{\leq J}\Psi_{(2)}$ in $\Ft_{(2)}$.
\end{itemize}      
Following the same steps as in section section \ref{sec:proofofThm:MainResultMext:k=0case}, that are summarized in Remark \ref{rmk:summaryofbasicruvatureestimatesMextusedlater}, we thus obtain the following analog of  \eqref{eq:controlofRkext0asawarmup} for $\Lieb_\T^{J+1}$ derivatives of curvature 
\bea\lab{eq:controlofLiebTJ+1curvatureinMextasStep1}
\nn\Rkext_0^2[\Lieb_\T^{J+1}(A, B, \Pc, \Bb, \Ab)]  &\les&    r_0^{-\dt}\Skext_0^2+ r_0^{ 3 +\dt}  \Rkint^2_0+r_0^{-2}\Rkext_{J+1}^2\\
&&+\ep_J^2+\ep_0^2,
\eea
where the extra terms $r_0^{-2}\Rkext_{J+1}^2$ and $\ep_J^2$ are respectively due to the fact that $r_0^{-2}\Rkext_{J+1}^2$ cannot yet be absorbed by the LHS and from the new terms of type $O(r^{-1})\dk^{\leq J}\Psi_{(1)}$ and $O(r^{-1})\dk^{\leq J}\Psi_{(2)}$ in $\Ft_{(1)}$ and $\Ft_{(2)}$ mentioned above.

{\bf Step 2.} Next, we derive estimates for $(\dkb, \Lieb_\T)^{J+1}$ derivatives of curvature in $\Mext$. To this end, we consider the Bianchi pairs in Proposition \ref{Proposition:BianchiPairs-Higher Derivatives} in the particular case 
\beaa
j=0, \qquad j_2\geq 1, \qquad j_1+j_2=J+1. 
\eeaa
In that case,  the terms $O(r^{-1} ) (q \nabc_4)^{j-1}\dkb^{j_2+1}  \Lieb^{j_1} _\T$ are not present, see footnote \ref{footnote:dkb^{-1}qe4j-1}, and hence the structure of the Bianchi pairs in Proposition \ref{Proposition:BianchiPairs-Higher Derivatives} are the exact analog for $J+1$ derivatives of the ones in Proposition \ref{prop:structureBainchPairs} used for the proof of the basic curvature estimates of section \ref{sec:proofofThm:MainResultMext:k=0case}, with the exception of
\begin{itemize}
\item the new terms of type $O(r^{-1})\dk^{\leq J}\Psi_{(1)}$ and $O(ar^{-1} )\dkb^{j_2-1}  \Lieb^{j_1+1} _\T\Psi_{(2)}$ in $\Ft_{(1)}$,
\item the new terms of type $O(r^{-1})\dk^{\leq J}\Psi_{(2)}$ and $O(ar^{-1} )\dkb^{j_2-1}  \Lieb^{j_1+1} _\T\Psi_{(1)}$ in $\Ft_{(2)}$.
\end{itemize}  
Following the same steps as in section section \ref{sec:proofofThm:MainResultMext:k=0case}, that are summarized in Remark \ref{rmk:summaryofbasicruvatureestimatesMextusedlater}, we thus obtain the following analog of  \eqref{eq:controlofRkext0asawarmup} for $\dkb^{j_2}\Lieb^{j_1}_\T$ derivatives of curvature 
\beaa
\nn\Rkext_0^2[\dkb^{j_2}\Lieb^{j_1}_\T(A, B, \Pc, \Bb, \Ab)]  &\les&    r_0^{-\dt}\Skext_0^2+ r_0^{ 3 +\dt}  \Rkint^2_0+r_0^{-2}\Rkext_{J+1}^2\\
&&+\Rkext_0^2[\dkb^{j_2-1}\Lieb^{j_1+1}_\T(A, B, \Pc, \Bb, \Ab)]+\ep_J^2+\ep_0^2,
\eeaa
where the new extra term $\Rkext_0^2[\dkb^{j_2-1}\Lieb^{j_1+1}_\T(A, B, \Pc, \Bb, \Ab)]$ is due to the new terms of type $O(ar^{-1} )\dkb^{j_2-1}  \Lieb^{j_1+1} _\T\Psi_{(2)}$ in $\Ft_{(1)}$, and of 
type $O(ar^{-1} )\dkb^{j_2-1}  \Lieb^{j_1+1} _\T\Psi_{(1)}$ in $\Ft_{(2)}$. Together with the estimate in the case $j_2=0$, i.e. \eqref{eq:controlofLiebTJ+1curvatureinMextasStep1}, we immediately infer, by iteration on $j_2$,
\bea\lab{eq:controlofLiebTj1dkbj2curvatureinMextasStep2}
\nn\Rkext_0^2[(\dkb, \Lieb_\T)^{J+1}(A, B, \Pc, \Bb, \Ab)]  &\les&    r_0^{-\dt}\Skext_0^2+ r_0^{ 3 +\dt}  \Rkint^2_0+r_0^{-2}\Rkext_{J+1}^2\\
&&+\ep_J^2+\ep_0^2.
\eea
 
{\bf Step 3.}  We now recover weighted $e_4$ derivatives and start with the ones for $(B, \Pc, \Bb, \Ab)$. To this end, we use the following simple consequence of the second Bianchi identity of a given general Bianchi pair as in Definition \ref{definition:GenBianchPairs:gen} 
\beaa
\nab_4\Psi_{(2)} &=& O(r^{-1})\dkb\Psi_{(1)}+O(r^{-1})\Psi_{(2)} +F_{(2)},
\eeaa
which implies
\beaa
\Rkext_0[r\nab_4\Psi_{(2)}] &\les& \Rkext_0[\dkb\Psi_{(1)}]+\Rkext_0[\Psi_{(2)}]+\Rkext_0[F_{(2)}]. 
\eeaa
Applying this estimate to the second Bianchi identity in each of the four Bianchi pairs in Proposition \ref{Proposition:BianchiPairs-Higher Derivatives} in the particular case
\beaa
j\geq 1, \qquad j+j_1+j_2=J+1,
\eeaa 
we immediately obtain
\bea\lab{eq:controlofre4jLiebTj1dkbj2curvatureinMextasStep3}
\nn\Rkext_0^2[(r\nab_4)^j\dkb^{j_2}\Lieb^{j_1}_\T(B, \Pc, \Bb, \Ab)]  &\les& \Rkext_0^2[(r\nab_4)^{j-1}\dkb^{j_2+1}\Lieb^{j_1}_\T(A, B, \Pc, \Bb)]\\
\nn&&+\Rkext_0^2[(r\nab_4)^j\dkb^{j_2-1}\Lieb^{j_1+1}_\T(A, B, \Pc, \Bb)]\\
\nn&&+   r_0^{-\dt}\Skext_0^2+ r_0^{ 3 +\dt}  \Rkint^2_0+r_0^{-2}\Rkext_{J+1}^2\\
&&+\ep_J^2+\ep_0^2,
\eea
where we recall the convention $\dkb^{-1}=0$. 
 
{\bf Step 4.}  We now recover weighted $e_4$ derivatives for $A$. To this end, we consider again the particular case
\beaa
j\geq 1, \qquad j+j_1+j_2=J+1
\eeaa 
and focus on the first Bianchi pair in Proposition \ref{Proposition:BianchiPairs-Higher Derivatives}. In that case, we choose $b=4+\dt$, and we have then  
  \beaa
   \La_{(1)}=j+2+\frac{\dt}{2}>0.
   \eeaa
We may thus apply case 4 in Proposition \ref{Prop:BasicWeighteBianch-Ext} to the first Bianchi pair in Proposition \ref{Proposition:BianchiPairs-Higher Derivatives} which yields
\beaa
\nn\Rkext_0^2[(r\nab_4)^j\dkb^{j_2}\Lieb^{j_1}_\T A]  &\les& \Rkext_0^2[(r\nab_4)^j\dkb^{j_2}\Lieb^{j_1}_\T B]+\Rkext_0^2[(r\nab_4)^{j-1}\dkb^{j_2+1}\Lieb^{j_1}_\T(A, B)]\\
\nn&&+\Rkext_0^2[(r\nab_4)^j\dkb^{j_2-1}\Lieb^{j_1+1}_\T(A, B)]\\
\nn&&+   r_0^{-\dt}\Skext_0^2+ r_0^{ 3 +\dt}  \Rkint^2_0+r_0^{-2}\Rkext_{J+1}^2+\ep_J^2+\ep_0^2.
\eeaa
 Together with \eqref{eq:controlofre4jLiebTj1dkbj2curvatureinMextasStep3}, we infer
 \bea\lab{eq:controlofre4jLiebTj1dkbj2curvatureinMextasStep4:intermediary}
\nn\Rkext_0^2[(r\nab_4)^j\dkb^{j_2}\Lieb^{j_1}_\T(A, B, \Pc, \Bb, \Ab)]  &\les& \Rkext_0^2[(r\nab_4)^{j-1}\dkb^{j_2+1}\Lieb^{j_1}_\T(A, B, \Pc, \Bb)]\\
\nn&&+\Rkext_0^2[(r\nab_4)^j\dkb^{j_2-1}\Lieb^{j_1+1}_\T(A, B, \Pc, \Bb)]\\
\nn&&+   r_0^{-\dt}\Skext_0^2+ r_0^{ 3 +\dt}  \Rkint^2_0+r_0^{-2}\Rkext_{J+1}^2\\
&&+\ep_J^2+\ep_0^2.
\eea
In view of \eqref{eq:controlofLiebTj1dkbj2curvatureinMextasStep2} and \eqref{eq:controlofre4jLiebTj1dkbj2curvatureinMextasStep4:intermediary}, and arguing by iteration on $j$, we infer\footnote{Notice that $j=0$ holds in view of \eqref{eq:controlofLiebTj1dkbj2curvatureinMextasStep2}, and for each $j$ with $1\leq j \leq J+1$, start with the fact that $\dkb^{-1}=0$ by convention and argue by iteration on $j_2$ from $j_2=0$ to $j_2=J+1-j$.}
 \beaa
\nn\Rkext_0^2[(r\nab_4, \dkb, \Lieb_\T)^{J+1}(A, B, \Pc, \Bb, \Ab)]  &\les&  r_0^{-\dt}\Skext_0^2+ r_0^{ 3 +\dt}  \Rkint^2_0+r_0^{-2}\Rkext_{J+1}^2\\
&&+\ep_J^2+\ep_0^2.
\eeaa
 Comparing $\Lieb_\T$ and $\nab_\T$, and using $\ep_J$ to absorb the corresponding lower term, we infer
 \beaa
\nn\Rkext_0^2[(r\nab_4, \dkb, \nab_\T)^{J+1}(A, B, \Pc, \Bb, \Ab)]  &\les&  r_0^{-\dt}\Skext_0^2+ r_0^{ 3 +\dt}  \Rkint^2_0+r_0^{-2}\Rkext_{J+1}^2\\
&&+\ep_J^2+\ep_0^2
\eeaa
and hence 
 \beaa
\nn\Rkext_{J+1}^2 &\les&  r_0^{-\dt}\Skext_0^2+ r_0^{ 3 +\dt}  \Rkint^2_0+r_0^{-2}\Rkext_{J+1}^2+\ep_J^2+\ep_0^2.
\eeaa
For $r_0$ large enough, we deduce 
\beaa
\nn\Rkext_{J+1}^2 &\les&  r_0^{-\dt}\Skext_0^2+ r_0^{ 3 +\dt}  \Rkint^2_0+\ep_J^2+\ep_0^2
\eeaa
as stated. This ends the proof of  Theorem  \ref{THM:MAINRESULTMEXT}.


\appendix



\chapter{Complement for Chapter \ref{CHAPTER-NON-INTEGRABLE-STRUCTURES}}



\section{Corollary to Lemma \ref{LEMMA:COMM-GEN-B}}


\begin{corollary} 
\lab{cor:comm-gen-B}
Let $U_{A}= U_{a_1\ldots a_k} $ be a $k$-horizontal  tensorfield symmetric traceless in all indices, i.e. $U\in \sk_k$.
\begin{enumerate}
\item We have
\bea
\bsplit
\,[\nab_3, \nab_b] U_A &=-\frac 1 2 \big( \trchb \nab_b U_A +\atrchb \dual \nab_b U_A\big)+( \eta_b-\ze_b) \nab_3 U_A \\
  &+\frac 1 2 \sum_{i=1}^k\big(\de_{a_i b} \trchb+\in_{ b a_i} \atrchb\big) \eta_c  U_{a_1\ldots }\,^ c \,_{\ldots a_k}\\
  &-\frac 1 2 \sum_{i=1}^k \eta_{a_i}  \big(\trchb U_{a_1\ldots b\ldots a_k} +\atrchb \dual U_{a_1\ldots b\ldots a_k}  \big)\\
&+ \err_{3bA}[U],\\
\err_{3bA}[U]&=\sum_{i=1}^k\Big(-\in_{a_i c} \dual\bb_b +   \chibh_{ba_i}  \eta_c - \chibh_{bc} \eta_{a_i}+  \chi_{ba_i} \xib_c-  \chi_{bc} \xib_{a_i} \Big) U_{a_1\ldots }\,^ c \,_{\ldots a_k} \\
&-\chibh_{bc}   \nab_c U_A+\xib_b \nab_4 U_A .
\end{split}
\eea
\item We have
\bea
\bsplit
\,[\nab_4, \nab_b] U_A &=-\frac 1 2 \big( \trch \nab_b U_A +\atrch \dual \nab_b U_A\big)+( \etab_b+\ze_b) \nab_4 U_a \\
  &+\frac 1 2 \sum_{i=1}^k\big(\de_{a_i b} \trch+\in_{ b a_i} \atrch\big) \etab_c  U_{a_1\ldots }\,^ c \,_{\ldots a_k}\\
  &-\frac 1 2 \sum_{i=1}^k \etab_{a_i}  \big(\trch U_{a_1\ldots b\ldots a_k} +\atrch \dual U_{a_1\ldots b\ldots a_k}  \big)\\
&+ \err_{4bA}[U],\\
\err_{4bA}[U]&=\sum_{i=1}^k\Big(\in_{a_i c} \dual\b_b +   \chih_{ba_i}  \etab_c - \chih_{bc} \etab_{a_i}+  \chib_{ba_i} \xi_c-  \chib_{bc} \xi_{a_i} \Big) U_{a_1\ldots }\,^ c \,_{\ldots a_k} \\
&-\chih_{bc}   \nab_c U_A+\xi_b \nab_3 U_A.
\end{split}
\eea
\item We have,
\bea
\bsplit
\, [\nab_4, \nab_3] U_A&= 2(\etab_b-\eta_b ) \nab_b U_A+ 2 \om \nab_3 U_A -2\omb \nab_4 U_A \\
&+ 2\sum_{i=1}^k\big( \eta_{a_i} \etab_b-\etab_{a_i} \eta_b- \in_{a_i b}\dual \rho) U_{a_1\ldots}\,^b\,_{\ldots a_k}  +\err_{43A},\\
\err_{43A}&= 2\sum_{i=1}^k \big( \xib_{a_i}  \xi_b- \xi_{a_i}  \xib_b )U_{a_1\ldots} \,^b\,_{\ldots a_k}.
\end{split}
\eea
\end{enumerate}
\end{corollary}


\section{Proof of Lemma \ref{LEMMA:COMMUTATIONNAB_XSQUARED}}
\label{sec:proof-lemma-comm-NabXsquared}


We write
\beaa
&& \squared   (  X^\b \Db_\b  U_a) =
 \g^{\mu\nu} \Db_\mu\Db_\nu(  X^\b \Db_\b  U_a)\\
 &=& ( \g^{\mu\nu} \Db_\mu \Db_\nu   X^\b) \Db_\b U_a + \g^{\mu\nu}\big(\Db_\mu X^\b  \Db_\nu \Db_\b   U_a+ \Db_\nu X^\b  \Db_\mu  \Db_\b  U_a \big)+
   X^\b \g^{\mu\nu}\Db_\mu\Db_\nu \Db_\b U_a \\
   &=& ( \g^{\mu\nu} \Db_\mu\Db_\nu  X^\b)\Db_\b U_a+ X^\b  \Db_\b \big(  \g^{\mu\nu}   \Db_\mu\Db_\nu ) U_a\\
   &&  + \g^{\mu\nu}\big(\Db_\mu X^\b  \Db_\nu   \Db_\b U_a+ \Db_\nu X^\b  \Db_\mu \Db_\b  U_a \big) + X^\b \g^{\mu\nu}\Big(\Db_\mu\Db_\nu \Db_\b - \Db_\b 
\Db_\mu\Db_\nu \Big) U_a .
\eeaa
Hence
\beaa
\squared  (  X^\b \Db_\b  U_a)&=& X^\b \Db_\b \squared  U_a+(\squared  X^\b)\Db_\b  U_a 
 + \g^{\mu\nu}\big(\Db_\mu X^\b  \Db_\nu \Db_\b   U_a+ \Db_\nu X^\b  \Db_\mu   \Db_\b U_a \big)\\
 &&+ X^\b  \g^{\mu\nu} \Big(\Db_\mu \big(\Db_\nu\Db_\b-\Db_\b \Db_\nu\big)U_a +\big(\Db_\mu\Db_\b -\Db_\b \Db_\mu \big)\Db_\nu  U_a \Big).
\eeaa
Using Lemma \ref{Lemma:R-PiXformula},  \eqref{definition:GaX} and the fact that $\R_{\mu\nu}=0$, we obtain
\beaa
\g^{\mu\nu} \D_\mu \D_\nu X_{\b}&=&\g^{\mu\nu}\R_{\b  \mu  \nu  \ga  }X^{\ga }+ \g^{\mu\nu}\GaX_{\mu \nu \b}\\
&=&\frac 1 2 \g^{\mu\nu} (\D_{\mu} \piX_{\nu \b}+\D_{\nu} \piX_{\mu \b}-\D_\b\piX_{\mu \nu })= \D^\mu  \pi_{\mu \b}-\frac 1 2  \D_\b  \tr \pi .
\eeaa
Consider  the term 
\beaa
\g^{\mu\nu}\big(\Db_\mu X^\b  \Db_\nu \Db_\b  U_a+ \Db_\nu X^\b  \Db_\mu\Db_\b   U_a \big)&=&\Db^\nu X^\b \Db_\nu \Db_\b U_a + \Db^\mu X^\b \Db_\mu\Db_\b U_a \\
&=& \Db^\nu X^\b \Db_\nu \Db_\b U_a+  \Db^\b X^\nu\Db_\b \Db_\nu U_a\\
&=& (\Db^\nu X^\b+\Db^\b X^\nu) \Db_\nu \Db_\b U_a\\
&&+  
 \Db^\b X^\nu( \Db_\b \Db_\nu -\Db_\nu \Db_\b )      U_a.
 \eeaa
 Using  the commutator formula of Proposition \ref{Proposition:commutehorizderivatives} and definition of $\pi$ we deduce
 \beaa
  \g^{\mu\nu}\big(\Db_\mu X^\b  \Db_\nu \Db_\b  U_a+ \Db_\nu X^\b  \Db_\mu\Db_\b   U_a \big)&=& \pi^{\mu\nu}  \Db_\mu \Db_\nu  U_a +  \Rdot_{ac \b\nu}  U^c  \D^\b X^\nu.
\eeaa
Therefore,
\beaa
&&\squared  (  X^\b \Db_\b  U_a)- X^\b \Db_\b \squared U_a\\
&=& \pi^{\mu\nu}  \Db_\mu \Db_\nu  U_a+
\big( \D^\mu  \pi_\mu \,^\b-\frac 1 2  \D^\b  \tr \pi) \Db_\b U_a +  \Rdot_{ac \b\nu}  U^c  \D^\b X^\nu\\
&&+ X^\b  \g^{\mu\nu} \Big(\Db_\mu \big(\Db_\nu\Db_\b-\Db_\b \Db_\nu\big)U_a +\big(\Db_\mu\Db_\b -\Db_\b \Db_\mu \big)\Db_\nu  U_a \Big).
\eeaa
Now,
\beaa
\g^{\mu\nu}\Db_\mu \big(\Db_\nu\Db_\b-\Db_\b \Db_\nu\big)U_a &=&\g^{\mu\nu}
\Db_\mu\big( \Rdot_{a  c\nu \b} U^c \big)= (\D^\mu \Rdot_{a  c\mu \b})  U^c         +\g^{\mu\nu}   \Rdot_{a  c\nu \b} \Db_\mu U_c, \\
\g^{\mu\nu}\big(\Db_\mu\Db_\b -\Db_\b \Db_\mu \big)\Db_\nu  U_a &=&\g^{\mu\nu}\big(
\R_{\nu \la \mu\b} \Db_\la U_a + \Rdot_{a c  \mu\b} \Db_\nu U_c \big)=\g^{\mu\nu} \Rdot_{a c  \mu\b} \Db_\nu U_c.
\eeaa
Using the fact that $\D^\mu \R_{a  c\mu \b}=0$, we finally have
\beaa
&&\squared  (  X^\b \Db_\b  U_a)- X^\b \Db_\b \squared U_a\\
&=& \pi^{\mu\nu}  \Db_\mu \Db_\nu  U_a+
\big( \D^\mu  \pi_\mu \,^\b-\frac 1 2  \D^\b  \tr \pi) \Db_\b U_a +  \Rdot_{ac \b\nu}  U^c  \D^\b X^\nu\\
&&+\frac 12   X^\b  (\D^\mu \B_{a  c\mu \b})  U^c  + 2 X^\b \Rdot_{a c  \mu\b} \Db^\mu U^c .
\eeaa
Writing that $\Rdot=\R+\frac 12 \B$, we obtain the stated identity.


\section{Proof of Proposition \ref{COMMUTATION-KK-SQUARE}}\label{sec:proof-comm-KK}


We first prove the following.

\begin{lemma}\label{lemma:carter-tensor-1} Let $K$ be a symmetric tensor which satisfies \eqref{definition-Pi}. Then
\bea
\D_\mu K^{\mu\nu}&=&-\frac 1 2 \D^\nu (\tr K)+\frac 3 2 {\Pi_{\mu}}^{\mu\nu} \label{divergence-K-trace}, \\
 \D^\mu\D_\a K^{\a\nu}-\D^\nu\D_\a K^{\a\mu}&=&\frac 3 2 \D^\mu {\Pi_{\a}}^{\a\nu} -\frac 3 2 \D^\nu {\Pi_{\a}}^{\a\mu} \label{relation-K-2}, \\
 \D^\mu\D^\a \D_{\a} K_{\mu\nu}&=& - \left(\frac 1 3 {{\R^{\a}}_{\nu\mu}}^\ep- {{\R^{\a}}_{\mu \nu}}^\ep \right) \D^\mu K_{\a \ep} - \D^\mu\D^\a \Pi_{\a\nu\mu}\nn \\
 &&+\frac 1 2 \D^\mu\D_\mu {\Pi^{\a}}_{\a\nu} -\frac 1 2 \D^\mu\D_\nu {\Pi^{\a}}_{\a\mu}. \label{relation-K-5}
\eea
Let $\phi$ be a scalar function. Then
\bea
(\D_\a K^{\mu\nu}) \D^\a \D_\nu \phi&=& - \frac 1 2 \D^{\mu} {K^{\nu}}_{\a} \D^\a \D_\nu \phi+\frac 3 2 {\Pi_{\a}}^{\mu\nu}\D^\a \D_\nu \phi, \label{relation-K-3} \\
(\D^\mu {K^{\nu}}_{\a}) \D_\mu \D^{\a}\D_\nu \phi&=&-\frac 2 3 \D^{\a} {K^{\nu}}_{\mu}{ {\R^\mu}_{\a \nu}}^\delta \D_\delta \phi+ \Pi_{\mu\a\nu} \D^\mu \D_{\a}\D_\nu \phi. \label{relation-K-4}
\eea
Let $\Psi \in \sk_2$. Then 
\bea\label{relation-K-6}
\begin{split}
(\D_\a K^{\mu\nu}) \Ddot^\a \Ddot_\nu \Psi_{ab}&= - \frac 1 2 \D^{\mu} {K^{\nu}}_{\a} \Ddot^\a \Ddot_\nu \Psi_{ab}\\
&+\frac 1 2 \D^{\nu} {K^{\mu} }_{\a}({{{\R_{\nu}}^\a}_{ a}}^c \Psi_{c b}+{ {{\R_{\nu}}^{\a}}_b }^c\Psi_{ac})+\frac 3 2 {\Pi_{\a}}^{\mu\nu}\Ddot^\a \Ddot_\nu \Psi_{ab},
\end{split}
\eea
and
\bea\label{relation-K-7}
\begin{split}
(\D^\mu {K^{\nu}}_{\a}) \Ddot_\mu \Ddot^{\a}\Ddot_\nu \Psi_{ab}&=-\frac 2 3 \D^{\a} {K^{\nu}}_{\mu}{ {\R^\mu}_{\a \nu}}^\delta \Ddot_\delta \Psi_{ab}+\Pi_{\a\mu\nu} \Ddot^\a \Ddot_{\mu}\Ddot_\nu \Psi_{ab}\\
&-\frac 23\D^{\mu} {K^{\nu}}_{\a}({ {\R^\a}_{\mu a} }^c \Ddot_\nu \Psi_{c b}+{ {\R^\a}_{\mu b} }^c \Ddot_\nu \Psi_{ac}) \\
&-\frac 1 3 \D^{\nu} {K^{\mu} }_{\a} \Ddot^\a({ \R_{\nu\mu a} }^c \Psi_{c b}+{ \R_{\nu\mu b }}^c \Psi_{a c}).
\end{split}
\eea
\end{lemma}
\begin{proof} 
 From contracting $3 \Pi_{\mu\nu\rho}=\D_{\mu} K_{\nu\rho}+\D_{\nu} K_{\rho\mu}+\D_{\rho} K_{\mu\nu} $ with $\g^{\mu\nu}$ we obtain \eqref{divergence-K-trace}. 
Using \eqref{divergence-K-trace}, we can write
\beaa
\D^\mu\D_\a K^{\a\nu}-\D^\nu\D_\a K^{\a\mu}&=&-\frac 1 2 \D^\mu \D^\nu (\tr K)+\frac 1 2 \D^\nu \D^\mu (\tr K)+\frac 3 2 \D^\mu {\Pi_{\a}}^{\a\nu} -\frac 3 2 \D^\nu {\Pi_{\a}}^{\a\mu} \\
&=&\frac 3 2 \D^\mu {\Pi_{\a}}^{\a\nu} -\frac 3 2 \D^\nu {\Pi_{\a}}^{\a\mu},
\eeaa
which proves \eqref{relation-K-2}. Applying $\D^\mu$ to $3 \Pi_{\mu\nu\rho}=\D_{\mu} K_{\nu\rho}+\D_{\nu} K_{\rho\mu}+\D_{\rho} K_{\mu\nu} $, we have
\beaa
3 \D^\mu \Pi_{\mu\nu\rho}&=&\D^\mu \D_{\mu} K_{\nu\rho}+\D^\mu \D_{\nu} K_{\rho\mu}+\D^\mu \D_{\rho} K_{\mu\nu} \\
&=&\D^\mu \D_{\mu} K_{\nu\rho}+ \D_{\nu} \D^\mu K_{\rho\mu}+ \D_{\rho} \D^\mu K_{\mu\nu} + {{\R^{\mu}}_{\nu\rho}}^\ep K_{\ep \mu}+ {{\R^{\mu}}_{\rho \nu}}^\ep K_{\mu \ep},
\eeaa
which gives
\beaa
\D^\a \D_{\a} K_{\mu\nu}&=&- \D_{\nu} \D^\a K_{\mu\a}- \D_{\mu} \D^\a K_{\a\nu} - {{\R^{\a}}_{\nu\mu}}^\ep K_{\ep \a}- {{\R^{\a}}_{\mu \nu}}^\ep K_{\a \ep}-3 \D^\a \Pi_{\a\nu\mu}.
\eeaa
Using \eqref{relation-K-2}, this gives
\beaa
\D^\a \D_{\a} K_{\mu\nu}&=&- 2\D_\mu\D^\a K_{\a\nu} - \left({{\R^{\a}}_{\nu\mu}}^\ep+ {{\R^{\a}}_{\mu \nu}}^\ep \right) K_{\a \ep}-3 \D^\a \Pi_{\a\nu\mu}\\
&&+\frac 3 2 \D_\mu {\Pi^{\a}}_{\a\nu} -\frac 3 2 \D_\nu {\Pi^{\a}}_{\a\mu}.
\eeaa
Applying $\D^\mu$ to the above we have 
\beaa
\D^\mu\D^\a \D_{\a} K_{\mu\nu}+2\D^\mu\D_\mu\D^\a K_{\a\nu}&=& - \left({{\R^{\a}}_{\nu\mu}}^\ep+ {{\R^{\a}}_{\mu \nu}}^\ep \right) \D^\mu K_{\a \ep}\\
&&-3 \D^\mu\D^\a \Pi_{\a\nu\mu}+\frac 3 2 \D^\mu\D_\mu {\Pi^{\a}}_{\a\nu} -\frac 3 2 \D^\mu\D_\nu {\Pi^{\a}}_{\a\mu}.
\eeaa
The left hand side is given by
\beaa
\D^\mu\D^\a \D_{\a} K_{\mu\nu}+2\D^\mu\D_\mu\D^\a K_{\a\nu}&=&3\D^\a \D_{\a}  \D^\mu K_{\mu\nu}+[\D^\mu, \D^\a \D_{\a}] K_{\mu\nu}\\
&=&3\D^\a \D_{\a}  \D^\mu K_{\mu\nu}+2{{\R_{\a}}^{\mu}}_{\nu\ep} \D_\mu K^{\a\ep },
\eeaa
which gives
\beaa
\D^\a \D_{\a}  \D^\mu K_{\mu\nu}&=& - \left(\frac 1 3 {{\R^{\a}}_{\nu\mu}}^\ep+ {{\R^{\a}}_{\mu \nu}}^\ep \right) \D^\mu K_{\a \ep}\\
&&- \D^\mu\D^\a \Pi_{\a\nu\mu}+\frac 1 2 \D^\mu\D_\mu {\Pi^{\a}}_{\a\nu} -\frac 1 2 \D^\mu\D_\nu {\Pi^{\a}}_{\a\mu} .
\eeaa
Writing again $\D^\mu\D^\a \D_{\a} K_{\mu\nu}=\D^\a \D_{\a}  \D^\mu K_{\mu\nu}+2{{\R_{\a}}^{\mu}}_{\nu\ep} \D_\mu K^{\a\ep }$,
we obtain \eqref{relation-K-5}.

From $3 {\Pi_{\a}}^{\mu\nu}=\D_{\a} K^{\mu\nu}+\D^{\mu} {K^{\nu}}_{\a}+\D^{\nu} {K^{\mu} }_{\a}$,
we have
\beaa
(\D_\a K^{\mu\nu}) \D^\a \D_\nu \phi&=& (-\D^{\mu} {K^{\nu}}_{\a}-\D^{\nu} {K^{\mu} }_{\a}+3 {\Pi_{\a}}^{\mu\nu}) \D^\a \D_\nu \phi\\
&=& -\D^{\mu} {K^{\nu}}_{\a} \D^\a \D_\nu \phi-\D^{\nu} {K^{\mu} }_{\a} \D^\a \D_{\nu} \phi+3 {\Pi_{\a}}^{\mu\nu}\D^\a \D_\nu \phi\\
&=& -\D^{\mu} {K^{\nu}}_{\a} \D^\a \D_\nu \phi-\D_{\a} K^{\mu\nu}  \D^\a \D_\nu \phi+3 {\Pi_{\a}}^{\mu\nu}\D^\a \D_\nu \phi,
\eeaa
which implies \eqref{relation-K-3}.

We similarly compute
\beaa
(\D_\a K^{\mu\nu}) \D^\a \D_{\mu}\D_\nu \phi&=& (-\D^{\mu} {K^{\nu}}_{\a}-\D^{\nu} {K^{\mu} }_{\a}+3 \Pi_{\a\mu\nu}) \D^\a \D_{\mu}\D_\nu \phi\\
&=& -\D^{\mu} {K^{\nu}}_{\a} \D^\a \D_{\mu}\D_\nu \phi-\D^{\nu} {K^{\mu} }_{\a} \D^\a \D_{\nu}\D_\mu \phi+3 \Pi_{\a\mu\nu} \D^\a \D_{\mu}\D_\nu \phi\\
&=&-2\D^{\mu} {K^{\nu}}_{\a} \D^\a \D_{\mu}\D_\nu \phi+3 \Pi_{\a\mu\nu} \D^\a \D_{\mu}\D_\nu \phi\\
&=&-2\D^{\mu} {K^{\nu}}_{\a}( \D_\mu \D^{\a}\D_\nu \phi+[\D^\a, \D_{\mu}]\D_\nu \phi ) +3 \Pi_{\a\mu\nu} \D^\a \D_{\mu}\D_\nu \phi\\
&=&-2\D^{\mu} {K^{\nu}}_{\a}( \D_\mu \D^{\a}\D_\nu \phi+{ {\R^\a}_{\mu \nu}}^\delta \D_\delta \phi) +3 \Pi_{\a\mu\nu} \D^\a \D_{\mu}\D_\nu \phi.
\eeaa
Observe that the first term on the right hand side is the same as the term on the left hand side. We therefore obtain \eqref{relation-K-4}.

From $3 {\Pi_{\a}}^{\mu\nu}=\D_{\a} K^{\mu\nu}+\D^{\mu} {K^{\nu}}_{\a}+\D^{\nu} {K^{\mu} }_{\a}$,
we have
\beaa
(\D_\a K^{\mu\nu}) \D^\a \D_\nu \Psi_{ab}&=& (-\D^{\mu} {K^{\nu}}_{\a}-\D^{\nu} {K^{\mu} }_{\a}+3 {\Pi_{\a}}^{\mu\nu}) \D^\a \D_\nu \Psi_{ab}\\
&=& -\D^{\mu} {K^{\nu}}_{\a} \D^\a \D_\nu \Psi_{ab}-\D^{\nu} {K^{\mu} }_{\a} \D^\a \D_{\nu} \Psi_{ab}+3 {\Pi_{\a}}^{\mu\nu}\D^\a \D_\nu \Psi_{ab}\\
&=& -\D^{\mu} {K^{\nu}}_{\a} \D^\a \D_\nu \Psi_{ab}-\D^{\nu} {K^{\mu} }_{\a} \D_\nu \D^\a \Psi_{ab}\\
&&+\D^{\nu} {K^{\mu} }_{\a}({{{\R_{\nu}}^\a}_{ a}}^\ep \Psi_{\ep b}+{ {{\R_{\nu}}^{\a}}_b }^\ep\Psi_{a\ep})+3 {\Pi_{\a}}^{\mu\nu}\D^\a \D_\nu \Psi_{ab}\\
&=& -\D^{\mu} {K^{\nu}}_{\a} \D^\a \D_\nu \Psi_{ab}-\D_{\a} K^{\mu\nu}  \D^\a \D_\nu \Psi_{ab}\\
&&+\D^{\nu} {K^{\mu} }_{\a}({{{\R_{\nu}}^\a}_{ a}}^\ep \Psi_{\ep b}+{ {{\R_{\nu}}^{\a}}_b }^\ep\Psi_{a\ep})+3 {\Pi_{\a}}^{\mu\nu}\D^\a \D_\nu \Psi_{ab},
\eeaa
which implies \eqref{relation-K-6}. 

We similarly compute
\beaa
&&(\D_\a K^{\mu\nu}) \D^\a \D_{\mu}\D_\nu \Psi_{ab}\\
&=& (-\D^{\mu} {K^{\nu}}_{\a}-\D^{\nu} {K^{\mu} }_{\a}+3 \Pi_{\a\mu\nu}) \D^\a \D_{\mu}\D_\nu \Psi_{ab}\\
&=& -\D^{\mu} {K^{\nu}}_{\a} \D^\a \D_{\mu}\D_\nu \Psi_{ab}-\D^{\nu} {K^{\mu} }_{\a} \D^\a \D_{\nu}\D_\mu \Psi_{ab}\\
&&-\D^{\nu} {K^{\mu} }_{\a} \D^\a({ \R_{\nu\mu a} }^c \Psi_{c b}+{ \R_{\nu\mu b }}^c \Psi_{a c})+3 \Pi_{\a\mu\nu} \D^\a \D_{\mu}\D_\nu \Psi_{ab}\\
&=&-2\D^{\mu} {K^{\nu}}_{\a} \D^\a \D_{\mu}\D_\nu \Psi_{ab}-\D^{\nu} {K^{\mu} }_{\a} \D^\a({ \R_{\nu\mu a} }^c \Psi_{c b}+{ \R_{\nu\mu b }}^c \Psi_{a c})+3 \Pi_{\a\mu\nu} \D^\a \D_{\mu}\D_\nu \Psi_{ab}\\
&=&-2\D^{\mu} {K^{\nu}}_{\a}( \D_\mu \D^{\a}\D_\nu \Psi_{ab}+[\D^\a, \D_{\mu}]\D_\nu \Psi_{ab} ) \\
&&-\D^{\nu} {K^{\mu} }_{\a} \D^\a({ \R_{\nu\mu a} }^c \Psi_{c b}+{ \R_{\nu\mu b }}^c \Psi_{a c})+3 \Pi_{\a\mu\nu} \D^\a \D_{\mu}\D_\nu \Psi_{ab}\\
&=&-2\D^{\mu} {K^{\nu}}_{\a}( \D_\mu \D^{\a}\D_\nu \Psi_{ab}+{ {\R^\a}_{\mu \nu}}^\delta \D_\delta \Psi_{ab}+{ {\R^\a}_{\mu a} }^c \D_\nu \Psi_{c b}+{ {\R^\a}_{\mu b} }^c \D_\nu \Psi_{ac}) \\
&&-\D^{\nu} {K^{\mu} }_{\a} \D^\a({ \R_{\nu\mu a} }^c \Psi_{c b}+{ \R_{\nu\mu b }}^c \Psi_{a c})+3 \Pi_{\a\mu\nu} \D^\a \D_{\mu}\D_\nu \Psi_{ab}.
\eeaa
Observe that the first term on the right hand side is the same as the term on the left hand side. We therefore obtain \eqref{relation-K-7}. 
\end{proof}

We now prove the following.

\begin{lemma}\label{preliminaries-commutators-lemma} 
In a vacuum spacetime, the following commutation formulas hold for a scalar $\phi$, a horizontal 1-tensor $\psi$ and a horizontal 2-tensor $\Psi$:
\beaa
\, [\D_\nu, \square_\g]\phi&=&0, \qquad \, [\Ddot_\mu, \squared_1 ]\psi^\mu =0\qquad \, [\Ddot_\nu, \squared_2]{\Psi^\nu}_{\de}=2\R_{\nu\a\de\ep} \Ddot^\a \Psi^{\nu\ep }.
\eeaa
Also, for $\Phi$ a 3-tensor which is symmetric in the last two indices, we have,
\beaa
\, [\Ddot_\nu, \square]{\Phi^\nu}_{\de\la}&=&2 \R_{\nu\a \de\ep} \Ddot^\a {\Phi^{\nu\ep}}_{\la} +2\R_{\nu\a \la\ep }\Ddot^\a {\Phi^{\nu\ep}}_{\de}.
\eeaa
\end{lemma}

\begin{proof} 

We have for a scalar function $\phi$:
\beaa
\, [\D_\a, \D_\b] \phi &=&0, \qquad  [\D_\a, \D_\b] \D_\gamma \phi = {\R_{\a\b\gamma}}^\delta \D_\delta \phi.
\eeaa
We therefore compute 
\beaa
\, [\D_\nu, \square_\g]\phi&=&  [\D_\nu, \D^\a \D_\a]\phi=[\D_\nu, \D^\a] \D_\a\phi + \D^\a[\D_\nu,  \D_\a]\phi\\
&=&\g^{\a\mu}{\R_{\nu\mu \a}}^\de \D_\de\phi = - {\R_{\nu}}^\de \D_\de \phi=0.
\eeaa
For a 1-tensor $X$ we have
\beaa
\, [\D_\a, \D_\b] X_\gamma  &=& {\R_{\a\b\gamma\epsilon }} X^\epsilon, \qquad  [\D_\a, \D_\b] \D_\delta X_\gamma  = {\R_{\a\b\delta \epsilon }} \D^\epsilon X_\gamma+ {\R_{\a\b\gamma\epsilon }} \D_\delta X^\epsilon.
\eeaa
This gives
\beaa
\, [\D_\mu, \square_\g ]X_\b &=&  [\D_\mu, \D^\a \D_\a]X_\b =[\D_\mu, \D^\a] \D_\a X_\b+ \D^\a ([\D_\mu, \D_\a]X_\b )\\
&=&\g^{\a\ze} ({\R_{\mu\ze\a}}^\ep \D_\ep X_\b+ {\R_{\mu\ze\b}}^\ep \D_\a  X_\ep )+ \D^\a ({\R_{\mu\a\b}}^\ep X_\ep )\\
&=& -{\R_{\mu}}^\ep \D_\ep X_\b +{{{\R_{\mu}}^\a}_{\b}}^\ep \D_\a  X_\ep + \D^\a {\R_{\mu\a\b}}^\ep X_\ep + {\R_{\mu\a\b}}^\ep \D^\a X_\ep.
\eeaa
Since by the Bianchi identities $\D^\a\R_{\a\mu\nu\gamma}=\D_\nu \R_{\mu\gamma}-\D_\gamma \R_{\mu\nu}$, the divergence of the Riemann tensor vanishes in the case of a vacuum spacetime. We therefore have
\beaa
\, [\D_\mu, \square_\g ]X_\be &=&2{\R_{\mu\a\b}}^\ep \D^\a X_\ep, \qquad  [\D_\mu, \square_\g ]X^\mu =2{\R_{\mu\a}}^{\mu\ep} \D^\a X_\ep=0 .
\eeaa
For a 2-tensor $\Psi$ we have
\beaa
\,[ \D_\nu , \D_\a ] \Psi_{\ga\de} &=&{ \R_{\nu\a \ga} }^\ep \Psi_{\ep\de}+{ \R_{\nu\a\de }}^\ep\Psi_{\gamma\ep}, \\
\, [\D_\nu, \D_\ze ] \D_\a\Psi_{\ga\de}&=& { \R_{\nu\ze \a} }^\ep \D_\ep \Psi_{\ga\de}+{ \R_{\nu\ze \gamma} }^\ep \D_\a \Psi_{\ep\de}+{ \R_{\nu\ze \de} }^\ep \D_\a \Psi_{\gamma\ep}.
\eeaa
In the case of a horizontal 2-tensor 
\beaa
\,[ \Ddot_\nu , \Ddot_\mu ] \Psi_{ab} &=&{ \R_{\nu\mu a} }^c \Psi_{c b}+{ \R_{\nu\mu b }}^c \Psi_{a c}, \\
\, [\Ddot_\nu, \Ddot_\ze ] \D_\a\Psi_{ab}&=& { \R_{\nu\ze \a} }^\ep \Ddot_\ep \Psi_{ab}+{ \R_{\nu\ze a} }^c \Ddot_\a \Psi_{c b}+{ \R_{\nu\ze b} }^c \Ddot_\a \Psi_{ac}.
\eeaa
We therefore compute
\beaa
\, [\D_\nu, \square_\g]\Psi_{\ga\de}&=& \, [\D_\nu, \D^\a \D_\a]\Psi_{\ga\de}= [\D_\nu, \D^\a ] \D_\a\Psi_{\ga\de}+ \D^\a [\D_\nu,  \D_\a]\Psi_{\ga\de}\\
&=&\g^{\a\ze} \left( { \R_{\nu\ze \a} }^\ep \D_\ep \Psi_{\ga\de}+{ \R_{\nu\ze \gamma} }^\ep \D_\a \Psi_{\ep\de}+{ \R_{\nu\ze \de} }^\ep \D_\a \Psi_{\gamma\ep} \right)\\
&&+ \D^\a({ \R_{\nu\a \ga} }^\ep \Psi_{\ep\de}+{ \R_{\nu\a\de }}^\ep\Psi_{\gamma\ep})\\
&=& -{ \R_{\nu} }^\ep \D_\ep \Psi_{\ga\de}+{{{ \R_{\nu}}^\a}_{ \gamma} }^\ep \D_\a \Psi_{\ep\de}+{{{ \R_{\nu}}^\a}_{ \de} }^\ep \D_\a \Psi_{\gamma\ep} \\
&&+ \D^\a{ \R_{\nu\a \ga} }^\ep \Psi_{\ep\de}+{ \R_{\nu\a \ga} }^\ep \D^\a\Psi_{\ep\de}+\D^\a{ \R_{\nu\a\de }}^\ep\Psi_{\gamma\ep}+{ \R_{\nu\a\de }}^\ep \D^\a \Psi_{\gamma\ep}\\
&=&2{{{ \R_{\nu}}^\a}_{ \gamma} }^\ep \D_\a \Psi_{\ep\de}+2{{{ \R_{\nu}}^\a}_{ \de} }^\ep \D_\a \Psi_{\gamma\ep}.
\eeaa
This gives
\beaa
\, [\D_\nu, \square_\g]{\Psi^\nu}_{\de}&=&2 {\R_{\nu}}^{\a\nu\ep} \D_\a \Psi_{\ep\de}+2{{{ \R_{\nu}}^\a}_{ \de} }^\ep \D_\a {\Psi^{\nu}}_{\ep} =2{{{ \R_{\nu}}^\a}_{ \de} }^\ep \D_\a {\Psi^{\nu}}_{\ep}.
\eeaa

For a $3$-tensor $\Phi$ we have 
\beaa
\,[ \D_\nu , \D_\a ] \Phi_{\ga\de\la} &=&{ \R_{\nu\a \ga} }^\ep \Phi_{\ep\de\la}+{ \R_{\nu\a\de }}^\ep\Phi_{\gamma\ep\la}+{ \R_{\nu\a \la} }^\ep \Phi_{\gamma \de\ep}, \\
\, [\D_\nu, \D_\ze ] \D_\a\Phi_{\ga\de\la}&=& { \R_{\nu\ze \a} }^\ep \D_\ep \Phi_{\ga\de\la}+{ \R_{\nu\ze \gamma} }^\ep \D_\a \Phi_{\ep\de\la}+{ \R_{\nu\ze \de} }^\ep \D_\a \Phi_{\gamma\ep\la }+{ \R_{\nu\ze \la} }^\ep \D_\a \Phi_{\gamma\de\ep }.
\eeaa
We therefore compute
\beaa
\, [\D_\nu, \square_\g]\Phi_{\ga\de\la}&=& \, [\D_\nu, \D^\a \D_\a]\Phi_{\ga\de\la}= [\D_\nu, \D^\a ] \D_\a\Phi_{\ga\de\la}+ \D^\a [\D_\nu,  \D_\a]\Phi_{\ga\de\la}\\
&=&\g^{\a\ze} \left( { \R_{\nu\ze \a} }^\ep \D_\ep \Phi_{\ga\de\la}+{ \R_{\nu\ze \gamma} }^\ep \D_\a \Phi_{\ep\de\la}+{ \R_{\nu\ze \de} }^\ep \D_\a \Phi_{\gamma\ep\la }+{ \R_{\nu\ze \la} }^\ep \D_\a \Phi_{\gamma\de\ep }\right)\\
&&+ \D^\a({ \R_{\nu\a \ga} }^\ep \Phi_{\ep\de\la}+{ \R_{\nu\a\de }}^\ep\Phi_{\gamma\ep\la}+{ \R_{\nu\a \la} }^\ep \Phi_{\gamma \de\ep})\\
&=& {{{ \R_{\nu}}^\a}_{ \gamma} }^\ep \D_\a \Phi_{\ep\de\la}+{{{ \R_{\nu}}^\a}_{ \de} }^\ep \D_\a \Phi_{\gamma\ep\la} +{{{ \R_{\nu}}^\a}_{ \la} }^\ep \D_\a \Phi_{\gamma\de\ep} \\
&&+{ \R_{\nu\a \ga} }^\ep \D^\a\Phi_{\ep\de\la}+{ \R_{\nu\a\de }}^\ep \D^\a \Phi_{\gamma\ep\la}+{ \R_{\nu\a \la} }^\ep \D^\a\Phi_{\gamma \de\ep}\\
&=&2{{{ \R_{\nu}}^\a}_{ \gamma} }^\ep \D_\a \Phi_{\ep\de\la}+2{{{ \R_{\nu}}^\a}_{ \de} }^\ep \D_\a \Phi_{\gamma\ep\la} +2{{{ \R_{\nu}}^\a}_{ \la} }^\ep \D_\a \Phi_{\gamma\de\ep}.
\eeaa
This gives
\beaa
\, [\D_\nu, \square]{\Phi^\nu}_{\de\la}&=&2 \R_{\nu\a \de\ep} \D^\a {\Phi^{\nu\ep}}_{\la} +2\R_{\nu\a \la\ep }\D^\a {\Phi^{\nu\ep}}_{\de}.
\eeaa
We can similarly adapt the above computations to the case of horizontal tensors.
\end{proof}

We can finally prove Proposition \ref{COMMUTATION-KK-SQUARE}.
We first consider the case of a $k$-tensor $\Psi$:
\beaa
[\KK, \square_\g]\Psi &=& [ \D_{\mu} K^{\mu\nu} \D_\nu, \square_\g]\Psi=\D_{\mu} [K^{\mu\nu} \D_\nu, \square_\g]\Psi+  [ \D_{\mu}, \square_\g] K^{\mu\nu} \D_\nu\Psi\\
&=&\D_{\mu}( K^{\mu\nu} \D_\nu \square_\g\Psi- \square_\g (K^{\mu\nu} \D_\nu \Psi))+  [ \D_{\mu}, \square_\g] K^{\mu\nu} \D_\nu\Psi.
\eeaa
Writing $\square_\g(K^{\mu\nu} \D_\nu \Psi)= \square_\g K^{\mu\nu} \D_\nu \Psi+2 \D^\a K^{\mu\nu} \D_\a \D_\nu \Psi+ K^{\mu\nu} \square_\g \D_\nu \Psi $, we have
\beaa
[\KK, \square_\g]\Psi &=&\D_{\mu}( K^{\mu\nu}[ \D_\nu,  \square_\g] \Psi- \square_\g K^{\mu\nu} \D_\nu \Psi-2 \D^\a K^{\mu\nu} \D_\a \D_\nu \Psi)+  [ \D_{\mu}, \square_\g] K^{\mu\nu} \D_\nu\Psi.
\eeaa
We now specialize to the case of a scalar function $\phi$. 
By Lemma \ref{preliminaries-commutators-lemma}, we have $ [\D_\nu, \square_\g]\phi=0$,  and $[\D_\mu, \square ]X^\mu =0$ applied to $X^\mu=K^{\mu\nu} \D_\nu \phi$ gives $ [ \D_{\mu}, \square_k] K^{\mu\nu} \D_\nu\phi=0$. We are left with
\beaa
[\KK, \square_\g]\phi &=&\D_{\mu}( K^{\mu\nu}[ \D_\nu,  \square_\g] \phi- \square_\g K^{\mu\nu} \D_\nu \phi-2 \D^\a K^{\mu\nu} \D_\a \D_\nu \phi)+  [ \D_{\mu}, \square_\g] K^{\mu\nu} \D_\nu\phi\\
&=&\D_{\mu}( - \square_\g K^{\mu\nu} \D_\nu \phi-2 \D_\a K^{\mu\nu} \D^\a \D_\nu \phi).
\eeaa
Using \eqref{relation-K-3}, we obtain
\beaa
[\KK, \square_\g]\phi &=&\D_{\mu}( - \square_2 K^{\mu\nu} \D_\nu \phi+ \D^{\mu} {K^{\nu}}_{\a} \D^\a \D_\nu \phi)-3 \D_\mu({\Pi_{\a}}^{\mu\nu}\D^\a \D_\nu \phi).
\eeaa
Observe that in expanding the derivative there is a cancellation of the term $(\square_\g K^{\mu\nu})\D_\mu\D_\nu \phi$. Indeed,
\beaa
[\KK, \square_\g]\phi &=& - \D_{\mu}(\square_\g K^{\mu\nu}) \D_\nu \phi- \square_\g K^{\mu\nu} \D_{\mu}\D_\nu \phi+\D_{\mu} \D^{\mu} {K^{\nu}}_{\a}  \D_\nu \D^\a\phi\\
&&+ \D^{\mu} {K^{\nu}}_{\a} \D_{\mu}\D^\a \D_\nu \phi-3 \D_\mu({\Pi_{\a}}^{\mu\nu}\D^\a \D_\nu \phi)\\
&=& - \D_{\mu}(\square_\g K^{\mu\nu}) \D_\nu \phi+ \D^{\mu} {K^{\nu}}_{\a} \D_{\mu}\D^\a \D_\nu \phi-3 \D_\mu({\Pi_{\a}}^{\mu\nu}\D^\a \D_\nu \phi).
\eeaa
Using \eqref{relation-K-5} and \eqref{relation-K-4}, we obtain
\beaa
[\KK, \square_\g]\phi  &=& \D^\mu K_{\a \ep}  \left(\frac 1 3 {{\R^{\a\de}}_{\mu}}^\ep- {{\R^{\a}}_{\mu }}^{\de\ep} \right) \D_\de \phi -\frac 2 3 \D^{\mu} {K^{\ep}}_{\a}{ {\R^\a}_{\mu \ep}}^\delta \D_\delta \phi\\
&&+ (\D^\mu\D^\a \Pi_{\a\nu\mu}-\frac 1 2 \D^\mu\D_\mu {\Pi^{\a}}_{\a\nu} +\frac 1 2 \D^\mu\D_\nu {\Pi^{\a}}_{\a\mu} ) \D_\nu \phi\\
&&+ \Pi_{\mu\a\nu} \D^\mu \D_{\a}\D_\nu \phi-3 \D_\mu({\Pi_{\a}}^{\mu\nu}\D^\a \D_\nu \phi).
\eeaa
Using that ${ {\R^\a}_{\mu \ep}}^\delta=-{{ {\R^\a}_{\mu}}^\delta}_{\ep}$ we have
\beaa
[\KK, \square_\g]\phi  &=&\frac 1 3 \D^\mu K_{\a \ep}  \left( {{\R^{\a\de}}_{\mu}}^\ep- {{\R^{\a}}_{\mu }}^{\de\ep} \right) \D_\de \phi \\
&&+ \left(\D^\mu\D^\a \Pi_{\a\nu\mu}-\frac 1 2 \D^\mu\D_\mu {\Pi^{\a}}_{\a\nu} +\frac 1 2 \D^\mu\D_\nu {\Pi^{\a}}_{\a\mu} \right) \D_\nu \phi\\
&&+ \Pi_{\mu\a\nu} \D^\mu \D_{\a}\D_\nu \phi-3 \D_\mu({\Pi_{\a}}^{\mu\nu}\D^\a \D_\nu \phi).
\eeaa
Writing 
\beaa
\D_{\mu} K_{\ep\a}&=& -\D_{\ep} K_{\a\mu}-\D_{\a} K_{\mu\ep} +3 \Pi_{\mu\ep\a}
\eeaa
we can observe that the first term is symmetric in $\a\mu$ while the second Riemann tensor term in antisymmetric in $\a\mu$, and the second term in symmetric in $\mu\ep$ while the first Riemann tensor is antisymmetric in $\mu\nu$. We therefore are left with
\beaa
[\KK, \square_\g]\phi  &=&\frac 1 3 (-\D_{\ep} {K_{\a}}^{\mu}{{\R^{\a\de}}_{\mu}}^\ep+\D_{\a} {K^{\mu}}_{\ep}{{\R^{\a}}_{\mu }}^{\de\ep}) \D_\de \phi \\
&&+ \left(\D^\mu\D^\a \Pi_{\a\nu\mu}-\frac 1 2 \D^\mu\D_\mu {\Pi^{\a}}_{\a\nu} +\frac 1 2 \D^\mu\D_\nu {\Pi^{\a}}_{\a\mu} \right) \D_\nu \phi\\
&&+ \Pi_{\mu\a\nu} \D^\mu \D_{\a}\D_\nu \phi-3 \D_\mu({\Pi_{\a}}^{\mu\nu}\D^\a \D_\nu \phi).
\eeaa
Writing $\D^{\nu} {K_{\a}}^{\mu} {{\R}^{\a \de}}_{\mu\nu}=\D^{\nu} {K_{\a}}^{\mu} {{\R}_{\mu\nu}}^{\a \de}=\D^{\a} {K_{\nu}}^{\mu} {{\R}_{\mu\a}}^{\nu \de}=-\D^{\a} {K_{\nu}}^{\mu} {{\R}_{\a\mu}}^{\nu \de}$ we obtain the cancellation of the first line, which gives
\beaa
[\KK, \square_\g]\phi  &=&\err[\Pi](\phi),
\eeaa
with
\beaa
\err[\Pi](\phi)&=& \left(\D^\mu\D^\a \Pi_{\a\nu\mu}-\frac 1 2 \D^\mu\D_\mu {\Pi^{\a}}_{\a\nu} +\frac 1 2 \D^\mu\D_\nu {\Pi^{\a}}_{\a\mu} \right) \D_\nu \phi\\
&&+ \Pi_{\mu\a\nu} \D^\mu \D_{\a}\D_\nu \phi-3 \D_\mu({\Pi_{\a}}^{\mu\nu}\D^\a \D_\nu \phi).
\eeaa
 We now simplify $\err[\Pi](\phi)$.
By writing 
\beaa
 \Pi_{\mu\a\nu}  \D^\mu\D^{\a}\D^\nu \phi= \D^\mu( \Pi_{\mu\a\nu} \D^{\a}\D^\nu \phi)-  \D^\mu \Pi_{\mu\a\nu} \D^{\a}\D^\nu \phi,
\eeaa
we obtain
\beaa
\err[\Pi](\phi)&=&\D^\mu  \left(\D^\a \Pi_{\a\nu\mu}-\frac 1 2 \D_\mu {\Pi^{\a}}_{\a\nu} +\frac 1 2 \D_\nu {\Pi^{\a}}_{\a\mu} \right) \D^\nu \phi\\
&&-  \D^\mu \Pi_{\mu\a\nu} \D^{\a}\D^\nu \phi-2 \D^\mu( \Pi_{\mu\a\nu} \D^{\a}\D^\nu \phi).
\eeaa
By writing
\beaa
&&\D^\mu \Big[ \left(\D^\a \Pi_{\a\nu\mu}-\frac 1 2 \D_\mu {\Pi^{\a}}_{\a\nu} +\frac 1 2 \D_\nu {\Pi^{\a}}_{\a\mu} \right) \D^\nu \phi \Big]\\
&=&\D^\mu \left(\D^\a \Pi_{\a\nu\mu}-\frac 1 2 \D_\mu {\Pi^{\a}}_{\a\nu} +\frac 1 2 \D_\nu {\Pi^{\a}}_{\a\mu} \right) \D^\nu \phi \\
&&+  \left(\D^\a \Pi_{\a\nu\mu}-\frac 1 2 \D_\mu {\Pi^{\a}}_{\a\nu} +\frac 1 2 \D_\nu {\Pi^{\a}}_{\a\mu} \right)\D^\mu \D^\nu \phi.
\eeaa
Observe that 
\beaa
\left(-\frac 1 2 \D_\mu {\Pi^{\a}}_{\a\nu} +\frac 1 2 \D_\nu {\Pi^{\a}}_{\a\mu} \right)\D^\mu \D^\nu \phi =0
\eeaa
since the first term is antisymmetric in $\mu\nu$ and the second term is symmetric in $\mu\nu$. We can therefore write
\beaa
&&\D^\mu \left(\D^\a \Pi_{\a\nu\mu}-\frac 1 2 \D_\mu {\Pi^{\a}}_{\a\nu} +\frac 1 2 \D_\nu {\Pi^{\a}}_{\a\mu} \right) \D^\nu \phi \\
&=&\D^\mu \Big[ \left(\D^\a \Pi_{\a\nu\mu}-\frac 1 2 \D_\mu {\Pi^{\a}}_{\a\nu} +\frac 1 2 \D_\nu {\Pi^{\a}}_{\a\mu} \right) \D^\nu \phi \Big]- (\D^\a \Pi_{\a\nu\mu} )\D^\mu \D^\nu \phi 
\eeaa
and we obtain
\beaa
\err[\Pi](\phi)&=&\D^\mu \Big[ (\D^\a \Pi_{\a\nu\mu}-\frac 1 2 \D_\mu {\Pi^{\a}}_{\a\nu} +\frac 1 2 \D_\nu {\Pi^{\a}}_{\a\mu} ) \D^\nu \phi -2 \Pi_{\mu\a\nu} \D^{\a}\D^\nu \phi \Big]\\
&&- 2(\D^\a \Pi_{\a\nu\mu} )\D^\mu \D^\nu \phi 
\eeaa
which proves the Proposition.


\section{Proof of Proposition  \ref{PROP:BOCHNERFOR-HORIZONTAL-LAP}}
\label{subsection:appendix-prop-lap2}


We write, for any $\psi\in \O_k$
\beaa
 \nab_a\nab_b\psi\c \nab_c\nab_d \psi &=&\nab_a\Big(\nab_b\psi\c \nab_c\nab_d \psi\Big) - \nab_b\psi\c   \nab_a \nab_c\nab_d \psi\\
 &=& - \nab_b\psi\c  \nab_c   \nab_a \nab_d \psi -  \nab_b\psi\c [  \nab_a, \nab_c ] \nab_d\psi +\nab_a\Big(\nab_b\psi\c \nab_c\nab_d \psi\Big)\\
 &=&  \nab_c  \nab_b\psi\c   \nab_a \nab_d \psi -  \nab_b\psi\c [  \nab_a, \nab_c ] \nab_d\psi +\nab_a\Big(\nab_b\psi\c \nab_c\nab_d \psi \Big)\\
&& -  \nab_c  \Big(\nab_b\psi\c    \nab_a \nab_d \psi \Big).
\eeaa
We deduce, for $\lap=\lap_k$,
\beaa
\big|\lap \psi|^2 &=& g^{ab} g^{cd}  \nab_a\nab_b\psi\c \nab_c\nab_d \psi \\
&=&\nab^a\nab^c\psi\c \nab_c\nab_a \psi- g^{ab} g^{cd} \nab_b\psi\c [  \nab_a, \nab_c ] \nab_d\psi  \\
&&+\nab_a\Big(\nab^a\psi \c\lap\psi \Big)-\nab_c\Big(\nab_a \psi\c \nab^a\nab^c\psi\Big)\\
&=&|\nab^2 \psi|^2 - \nab^a\nab^c\psi\c [\nab_a, \nab_c] \psi - g^{ab} g^{cd} \nab_b\psi\c [  \nab_a, \nab_c ] \nab_d\psi + \div_k[\lap\psi] \\
&=&|\nab^2 \psi|^2- \frac 1 2   [\nab^a, \nab^c] \psi  \c  [\nab_a, \nab_c] \psi   - g^{ab} g^{cd} \nab_b\psi\c [  \nab_a, \nab_c ] \nab_d\psi + \div_k[\lap\psi],
\eeaa
 with divergence term of the form
 \beaa
\div_k[\lap\psi]:&=&\nab_a\Big(\nab^a \lap\psi \Big)-\nab_c\Big(\nab_a\c \nab^a\nab^c\psi\Big)= \nab_a\Big(\nab^a\psi\c\lap\psi \Big)-\nab_a\Big(\nab_c\psi \c \nab^c\nab^a\psi\Big).
\eeaa
Hence
\bea
\lab{eq:Bochner-lap1}
\bsplit
\big|\lap \psi|^2 &= \big|\nab^2 \psi\big|^2  -A-  \frac 1 2 B+\div_{k}[\De\psi]\\
A&:=  g^{ab} g^{cd} \nab_b\psi\c [  \nab_a, \nab_c ] \nab_d\psi\\
B&:=   [\nab^a, \nab^c] \psi  \c  [\nab_a, \nab_c] \psi  \\
\div_{k}[\lap\psi]&:=  \nab_a\Big(\nab^a\psi\c\lap\psi \Big)-\nab_a\Big(\nab_c\psi \c \nab^c\nab^a\psi\Big).
\end{split}
\eea
Now,  for horizontal indices $I= i_1\ldots i_k$, we deduce using Proposition \ref{Gauss-equation-k-tensors} 
\beaa
 [  \nab_a, \nab_c ] \nab_d\psi_I &=& \frac 1 2 \big(\atrch\nab_3+\atrchb \nab_4\big) \nab_d\psi_I\in_{ac} + \Kh \big( g_{da} g_{tc}- g_{dc} g_{ta}\big) \nab^t \psi_I\\
 &&+\Kh \Big[ \big( g_{i_1 a }  g_{tc} - g_{i_1 c}  g_{ta} \big) \nab_d U^t\,_{i_2\ldots i_k}+\cdots  \big( g_{i_k a }  g_{tc} - g_{i_k c}  g_{ta} \big)     \nab_d  U_{ i_1\ldots i_{k-1}}\,^t\Big].
\eeaa
We deduce,  assuming that $\psi\in\O_k $ is symmetric
\bea
\bsplit
A&= \frac 1 2 \nab \psi\c   \big(\atrch\nab_3+\atrchb \nab_4\big)\dual \nab\psi -\Kh |\nab\psi|^2\\
&+ k\Kh \Big(  \nab_a\psi^{a\ldots}  \nab^a \psi_{a\ldots} - \nab_a \psi_{c\cdots}\nab^c\psi^{a \cdots} \Big).
\end{split}
\eea
Similarly,
\beaa
  [\nab_a, \nab_c] \psi _I &=&
   \frac 1 2 \big(\atrch\nab_3+\atrchb \nab_4\big) \psi_I\in_{ac} \\
   &&+\Kh \Big[ \big( g_{i_1 a }  g_{tc} - g_{i_1 c}  g_{ta} \big)  U^t\,_{i_2\ldots i_k}+\cdots  \big( g_{i_k a }  g_{tc} - g_{i_k c}  g_{ta} \big)    U_{ i_1\ldots i_{k-1}}\,^t\Big],
\eeaa
and
\beaa
B&=&\frac 1 2 \Big| \big(\atrch\nab_3+\atrchb \nab_4 \big)\psi\Big|^2  + 2  k\Kh  
\big(\atrch\nab_3+\atrchb \nab_4\big) \psi\c \dual \psi  + 2 k \Kh^2 |\psi|^2.
\eeaa
We thus  deduce,
\bea
\lab{eq:Bochner-lap2}
\bsplit
\big|\lap \psi|^2 &=|\nab^2 \psi|^2   +\Kh\Big( |\nab\psi|^2- k \Kh |\psi|^2\Big)\\
& + k\Kh \Big( \nab_a \psi_{c\cdots}\nab^c\psi^{a \cdots} -  \nab_a\psi^{a\ldots}  \nab^a \psi_{a\ldots}\Big) \\
& +\err_{k,1}[\lap \psi] +\div_{k}[\psi],
\end{split}
\eea
with error and divergence terms
\beaa
\err_{k,1}[\lap \psi]&=&-\frac 1 2 \nab \psi\c   \big(\atrch\nab_3+\atrchb \nab_4\big)\dual \nab\psi +\frac 1 2 \Big| \big(\atrch\nab_3+\atrchb \nab_4 \big)\psi\Big|^2\\
&&- k \Kh  
\big(\atrch\nab_3+\atrchb \nab_4\big) \psi\c \dual \psi 
\\
\div_{k}[\lap\psi]&=&\nab_a\Big(\nab^a\psi\c\lap\psi -\nab_c\psi \c \nab^c\nab^a\psi\Big).
\eeaa
It remains to calculate the term
\beaa
J_k &:=&k \Big( \nab_a \psi_{c\cdots}\nab^c\psi^{a \cdots} -  \nab_a\psi^{a\ldots}  \nab^a \psi_{a\ldots}\Big).
\eeaa

\textit{Scalar case.} In the particular case when $\psi$ is a scalar, $k=0$,
we deduce
\beaa
\big|\lap \psi|^2 &=& \big|\nab^2 \psi\big|^2 +\Kh |\nab\psi|^2+\err_0[\lap \psi] +\div_0[\lap \psi]\\
\err_0[\psi]&=& -\frac 1 2 \nab \psi\c   \big(\atrch\nab_3+\atrchb \nab_4\big)\dual \nab\psi +\frac 1 2 \Big| \big(\atrch\nab_3+\atrchb \nab_4 \big)\psi\Big|^2\\
\div[\lap\psi]&=& \nab_a\Big(\nab^a\psi\c  \lap\psi  - \frac 1 2 \nab^a|\nab\psi|^2 \Big)  ,
\eeaa
 as stated.

 \bigskip
 
 \textit{Case $k=1$.} We consider  now the case $k=1$  in which case  the term  $J_1[\psi]$ takes the form
 \beaa
 J_1&=&\Big( \nab_a \psi_{c}\nab^c\psi^{a } -  \nab_a\psi^{a}  \nab^a \psi_{a}\Big)=  |\nab\psi|^2     +  \nab_a\psi_c\big( \nab^c\psi^{a} -\nab^a\psi^c \big) -|\div\psi|^2\\
 &=&|\nab\psi|^2 +\frac   1 2 \big( \nab_a\psi_c-\nab_c\psi_a)\big( \nab^c\psi^{a} -\nab^a\psi^c \big)-|\div\psi|^2 \\
 &=& |\nab\psi|^2- |\curl\psi|^2 -|\div \psi|^2 =  |\nab\psi|^2-|\DDd_1\psi|^2.
 \eeaa
 We now recall \eqref{eq:hodgeident1-nonint-spacet-fin-div}  in  Proposition \ref{prop:2D-hodge-non-integrable-pert-Kerr-div}.
We deduce
\beaa
J_1&=&  |\nab\psi|^2-|\DDd_1\psi|^2 \\
&=&- \Kh |\psi|^2 +\frac 1 2 \bigg(\Big(\atrch\nab_3+\atrchb \nab_4\Big) \dual \psi  \bigg) \c \psi +\div[\DDd_1\psi]\\
div[\DDd_1\psi]&=&  \nab_a \Big( \nab^a \psi \c \psi -             (\div \psi ) \psi ^{a}-(\curl \psi ) (\dual \psi)^a \Big).
\eeaa
Hence, \eqref{eq:Bochner-lap2} becomes
\beaa
\big|\lap \psi|^2 &=&|\nab^2 \psi|^2   +\Kh\Big( |\nab\psi|^2-  \Kh |\psi|^2\Big)+ \Kh J_1
 +\err_{1,1}[\lap \psi] +\div_{1}[\psi]\\
&=& |\nab^2 \psi|^2   +\Kh \big(|\nab\psi|^2-2\Kh |\psi|^2 \big) +\err_1[\lap \psi]+\div_1[\lap\psi],
\eeaa
where
\beaa
\err_1[\lap \psi]&=&\err_{1,1}[\lap \psi]  +\frac 1 2 \Kh\bigg(\Big(\atrch\nab_3+\atrchb \nab_4\Big) \dual \psi  \bigg) \c \psi+\Kh\div[\DDd_1\psi]  \\
&=&-\frac 1 2 \nab \psi\c   \big(\atrch\nab_3+\atrchb \nab_4\big)\dual \nab\psi +\frac 1 2 \Big| \big(\atrch\nab_3+\atrchb \nab_4 \big)\psi\Big|^2\\
&&-  \Kh  
\big(\atrch\nab_3+\atrchb \nab_4\big) \psi\c \dual \psi  \\
&&+\frac 1 2 \Kh\bigg(\Big(\atrch\nab_3+\atrchb \nab_4\Big) \dual \psi  \bigg) \c \psi+\Kh\div[\DDd_1\psi] \\
&=&-\frac 1 2 \nab \psi\c   \big(\atrch\nab_3+\atrchb \nab_4\big)\dual \nab\psi +\frac 1 2 \Big| \big(\atrch\nab_3+\atrchb \nab_4 \big)\psi\Big|^2\\
&&-\frac 3 2  \Kh  
\big(\atrch\nab_3+\atrchb \nab_4\big) \psi\c \dual \psi +\Kh\div[\DDd_1\psi],
\eeaa
and
\beaa
\div_1[\lap\psi]&=& \div_1[\lap\psi]  + \div[\DDd_1\psi]\\
&=& \nab_a\Big(\nab^a\psi\c  \lap\psi  - \frac 1 2 \nab^a|\nab\psi|^2 \Big) +\nab_a \Big( \nab^a \psi \c \psi -             (\div \psi ) \psi ^{a}-(\curl \psi ) (\dual \psi)^a \Big),
\eeaa
as stated in part 2 of the proposition.

\bigskip

\textit{Case $k=2$.}
 It remains to consider  the case $k=2$  in which case  the term 
 \beaa
  J_2=2\Big( \nab_a \psi_{ci}\nab^c\psi^{a i} -  \nab_a\psi^{i}  \nab^a \psi_{a i }\Big)
  \eeaa
    in \eqref{eq:Bochner-lap2} becomes
 \beaa
 \frac 1 2 J_2&=& \nab_a \psi_{ci}\nab^c\psi^{ai} -  |\div\psi|^2  =  |\nab\psi|^2     +  \nab_a\psi_{ci}\big( \nab^c\psi^{ai } -\nab^a\psi^{ci} \big) -|\div\psi|^2\\
 &=&|\nab\psi|^2 +\frac   1 2 \big( \nab_a\psi_{ci}-\nab_c\psi_{ai})\big( \nab^c\psi^{ai} -\nab^a\psi^{ci} \big)-|\div\psi|^2.
 \eeaa
Note that $ \nab_a\psi_{ci}-\nab_c\psi_{ai}=\in_{ac} (\dual \div \psi)_i$. Hence
\beaa
 J_2&=& 2 |\nab\psi|^2- 4 \div \psi=2 |\nab\psi|^2- 4 \DDd_2 \psi.
\eeaa
Recalling the identity  \ref{eq:hodgeident2-nonint-spacet-fin-div}  of Proposition \ref{prop:2D-hodge-non-integrable-pert-Kerr-div} we deduce
\beaa
J_2&=& 2\big( |\nab\psi|^2- 2 \DDd_2 \psi\big)=- 4 \Kh|\psi|^2 + \bigg(\Big(\atrch\nab_3+\atrchb \nab_4\Big) \dual f\psi\bigg) \c \psi+ 2\div[\DDd_2\psi]
\eeaa
and, in view of \eqref{eq:Bochner-lap2},
\beaa
\big|\lap \psi|^2 &=&|\nab^2 \psi|^2   +\Kh\Big( |\nab\psi|^2- 2 \Kh |\psi|^2\Big) + \Kh J_2 +\err_{2,1}[\lap \psi] +\div_{2}[\psi]\\
&=& |\nab^2 \psi|^2   +\Kh\Big( |\nab\psi|^2- 6 \Kh |\psi|^2\Big) +\Kh \bigg(\Big(\atrch\nab_3+\atrchb \nab_4\Big) \dual \psi\bigg) \c \psi\\
&&+ 2\Kh\div[\DDd_2\psi]  +\err_{2,1}[\lap \psi] +\div_{2}[\psi]\
\eeaa
or,
\beaa
\big|\lap \psi|^2 &=& |\nab^2 \psi|^2   +\Kh\Big( |\nab\psi|^2- 6 \Kh |\psi|^2\Big) +\err_2[\psi]+\div_2\psi,
\eeaa
where
\beaa
\err_2[\lap \psi]&=&-\frac 1 2 \nab \psi\c   \big(\atrch\nab_3+\atrchb \nab_4\big)\dual \nab\psi +\frac 1 2 \Big| \big(\atrch\nab_3+\atrchb \nab_4 \big)\psi\Big|^2\\
&&- 3\Kh  
\big(\atrch\nab_3+\atrchb \nab_4\big) \psi\c \dual \psi + 2\Kh\div[\DDd_2\psi] \\
\div_{k}[\lap\psi]&=&\nab_a\Big(\nab^a\psi\c\lap\psi -\nab_c\psi \c \nab^c\nab^a\psi\Big).
\eeaa
This ends the proof of the Proposition.


\section{Proof of Lemma \ref{SIMPLIFICATION-ANGULAR}}
\lab{Appendix:ProofLemma-simplificationangular}


The identities in \eqref{DD-hot-hF} are proved by straightforward computation. We show here some of the computations. 

To prove the fourth formula in \eqref{DD-hot-hF}, we write
\beaa
 \DD\hot(\ov{F}\c U)_{ab}&=& \DD_a(\ov{F}\c U)_b+\DD_b(\ov{F}\c U)_a- \de_{ab}\DD^c (\ov{F}\c U)_c\\
&=& \DD_a ( \ov{F} ^c  U_{cb}  ) +\DD_b( \ov{F} ^c  U_{ca}  ) -\de_{ab} \DD^d ( \ov{F}^c U_{cd}) \\
&=& \DD_a  \ov{F} ^c  U_{cb}   +\DD_b  \ov{F} ^c  U_{ca}   -\de_{ab} \DD^d  \ov{F}^c U_{cd}\\
&+& \ov{F}^c\big(   \DD_a  U_{cb} +  \DD_b  U_{ca}- \de_{ab} \DD^d U_{cd}\big).
\eeaa
Now, in view of Lemma \ref{le:nonsym-product},
\beaa
 \DD_a  \ov{F} ^c  U_{cb}   +\DD_b  \ov{F} ^c  U_{ca} &=&\de_{ab}   (\DD^d  \ov{F} ^c) U_{cd}+(\DD\c \ov{F} ) U_{ab} 
 \\
 &&+\frac 1 2\Big(\big( \DD_a \ov{F} _c-\DD_c \ov{F}_a \big) U_{cb}+ \big( \DD_b \ov{F}_c -\DD_c \ov{F}_b \big) U_{ca}\Big).
\eeaa
Hence
\beaa
 \DD_a  \ov{F} ^c  U_{cb}   +\DD_b  \ov{F} ^c  U_{ca}   -\de_{ab} \DD^d  \ov{F}^c U_{cd}&=&(\DD\c \ov{F} ) U_{ab} \\
 &&+\frac 1 2\Big(\big( \DD_a \ov{F} _c-\DD_c \ov{F}_a \big) U_{cb}+ \big( \DD_b \ov{F}_c -\DD_c \ov{F}_b \big) U_{ca}\Big).
 \eeaa
 Recall that $\dual F=-i F, \,  \dual U=-i U, \,  \dual \DD=-i\DD $. We deduce,
 \beaa
 \DD_a \ov{F}_b -\DD_b \ov{F} _a&=& i\in_{ab}( \DD\c \ov{F} ).
 \eeaa
To check note that 
 \beaa
   \DD_1 \ov{F}_2 -\DD_2 \ov{F} _1&=& 2\Big[(\nab_1 f_2 -\nab_2 f_1) + i (\nab_1 f_1+\nab_2 f_2)\Big],\\
   ( \DD\c \ov{F} )&=& 2 \Big[  (\nab_1 f_1+\nab_2 f_2)- i (\nab_1 f_2 -\nab_2 f_1)\Big].
 \eeaa
 We deduce,
 \beaa
 \DD_a  \ov{F} ^c  U_{cb}   +\DD_b  \ov{F} ^c  U_{ca}   -\de_{ab} \DD^d  \ov{F}^c U_{cd}&=&(\DD\c \ov{F} ) U_{ab}+\frac 1 2 i  (\DD\c \ov{F} )\big( \in_{ac} U_{cb}+\in_{bc} U_{ca}\big)\\
 &=& (\DD\c \ov{F} ) U_{ab} +\frac 1 2 i  (\DD\c \ov{F} )\big( - 2 i U_{ab}\big)\\
&=& 2  (\DD\c \ov{F} ) U_{ab}.
 \eeaa
 Therefore,
 \beaa
  \DD\hot(\ov{F}\c U)_{ab}&=&  2  (\DD\c \ov{F} ) U_{ab}+\ov{F}^c\big(   \DD_a  U_{cb} +  \DD_b  U_{ca}- \de_{ab} \DD^d U_{cd}\big).
 \eeaa
 It remains to re-express the tensor 
 \beaa
   \DD_a  U_{cb} +  \DD_b  U_{ca}- \de_{ab} \DD^d U_{cd}.
   \eeaa
 Note also that  $ \DD^d U_{cd}=0$. We claim
 \beaa
  \DD_a  U_{cb} +  \DD_b  U_{ca}&=& 2 \DD_c U_{ab}.
 \eeaa
Indeed, for $a=b=1$, $c=2$,
 \beaa
 2\DD_1  U_{21} &=&2(\nab_1+i \dual \nab_1)  U_{21}=2(\nab_1+i  \nab_2)  U_{12}=-2i(\nab_1+i  \nab_2)  U_{11}=2( \nab_2-i\nab_1)  U_{11},\\
 2\DD_2 U_{11}&=& 2(\nab_2+i \dual \nab_2)  U_{11}=2(\nab_2-i \nab_1)  U_{11}.
 \eeaa
 For $a=c=1$, $b=2$,
 \beaa
 \DD_1  U_{12} +  \DD_2  U_{11}&=&(\nab_1+i\nab_2) U_{12} + (\nab_2-i\nab_1 )  U_{11}=(\nab_1+i\nab_2) U_{12} + i(\nab_2-i\nab_1 )  U_{12}\\
&=&2(\nab_1+i\nab_2) U_{12} =2\DD_1 U_{12}.
 \eeaa
 We deduce,
 \beaa
 \DD\hot(\ov{F}\c U)_{ab}&=&  2  (\DD\c \ov{F} ) U_{ab}+ 2 \ov{F}^c \DD_c U_{ab},
 \eeaa
as stated.

To prove the last formula in \eqref{DD-hot-hF}, we write
\beaa
( \DDb \c F )&=& \DDb^a F_{a}=(\nab - i \dual \nab)^a F_{a}=(\nab_1 - i \dual \nab_1) F_{1}+(\nab_2 - i \dual \nab_2) F_{2}\\
&=&(\nab_1 - i  \nab_2) F_{1}+(\nab_2 + i  \nab_1) (-i F_{1})=2(\nab_1 - i  \nab_2) F_{1},
\eeaa
which gives
\beaa
(U (\ov{\DD} \c F))_{ab}&=& 2(\nab_1 - i  \nab_2) F_{1} U_{ab}.
\eeaa
On the other hand,
\beaa
(U \c \ov{\DD} F)_{ab}&=&U_{ad}  \ov{\DD}_d F_{b}=U_{a1}  \ov{\DD}_1 F_{b}+U_{a2}  \ov{\DD}_2 F_{b}\\
&=&(u_{a1}+ i \dual u_{a1}) (\nab_1 - i \dual\nab_1) F_{b}+(u_{a2}+ i \dual u_{a2}) (\nab_2 - i \dual\nab_2)  F_{b}\\
&=&(u_{a1}+ i \dual u_{a1}) (\nab_1 - i\nab_2) F_{b}+(u_{a2}+ i \dual u_{a2}) (\nab_2 + i \nab_1)  F_{b},
\eeaa
and therefore 
\beaa
(\ov{U} \c \DD F)_{11}&=&(u_{11}+ iu_{12}) (\nab_1 - i\nab_2) F_{1}+(u_{12}- i  u_{11}) (\nab_2 + i \nab_1)  F_{1}\\
&=&2(u_{11}\nab_1-i u_{11} \nab_2+u_{12}\nab_2 + i u_{12}\nab_1)F_{1}\\
&=& 2(\nab_1 - i  \nab_2) F_{1} U_{ab}.
\eeaa
Similarly for the other components. We obtain the stated identity.

To prove the first formula in \eqref{Leib-eq-DDb}, we write
\beaa
( \DDb \c U )_1&=& \DDb^a U_{1a}=(\nab - i \dual \nab)^a U_{1a}\\
&=&(\nab_1 - i \dual \nab_1) U_{11}+(\nab_2 - i \dual \nab_2) U_{12}\\
&=&(\nab_1 - i  \nab_2) U_{11}+(\nab_2 + i  \nab_1) (-i U_{11})\\
&=&2(\nab_1 - i  \nab_2) U_{11},
\eeaa
and
\beaa
( \DDb \c U )_2&=& \DDb^a U_{2a}=(\nab - i \dual \nab)^a U_{2a}\\
&=&(\nab_1 - i \dual \nab_1) U_{21}+(\nab_2 - i \dual \nab_2) U_{22}\\
&=&(\nab_1 - i  \nab_2)(-i U_{11})+(\nab_2 + i  \nab_1) (-U_{11})\\
&=&-2i(\nab_1 -i  \nab_2) U_{11}.
\eeaa
Therefore 
\beaa
(F \hot ( \DDb \c U ))_{11}&=& 2F_1( \DDb \c U )_1-\de_{11} F \c ( \DDb \c U )\\
&=& F_1( \DDb \c U )_1- F_2 ( \DDb \c U )_2\\
&=& (f_1+i \dual f_1)2(\nab_1 - i  \nab_2) U_{11}- (f_2+i\dual f_2)(-2i(\nab_1 -i  \nab_2) U_{11})\\
&=& 4f_1\nab_1U_{11} - 4i  f_1\nab_2U_{11} +4i f_2\nab_1U_{11} +4 f_2   \nab_2U_{11}. 
\eeaa
On the other hand,
\beaa
(F \hot ( \DDb \c U ))_{12}&=& F_1( \DDb \c U )_2+F_2( \DDb\c U )_1\\
&=& (f_1+i \dual f_1)(2(\nab_1 -i  \nab_2) U_{12})+ (f_2+i\dual f_2)2(\nab_1 - i  \nab_2) iU_{12}\\
&=& (f_1+i  f_2)(2(\nab_1 -i  \nab_2) U_{12})+ (f_2-i f_1)2(\nab_1 - i  \nab_2) iU_{12}\\
&=& 4f_1\nab_1U_{12} - 4i  f_1\nab_2U_{12} +4i f_2\nab_1U_{12} +4 f_2   \nab_2U_{12} ,
\eeaa
which therefore gives
\beaa
\frac 1 2 (F \hot ( \DDb \c U ))_{ab}&=& 2f_1\nab_1U_{ab} - 2i  f_1\nab_2U_{ab} +2i f_2\nab_1U_{ab} +2 f_2   \nab_2U_{ab}. 
\eeaa
On the other hand,
\beaa
(F\c \DDb U)_{ab}&=& F^c\DDb_c U_{ab}= F_1\DDb_1 U_{ab}+F_2\DDb_2 U_{ab}\\
&=&(f_1+ i \dual f_1)(\nab_1- i \dual \nab_1) U_{ab}+(f_2+ i \dual f_2)(\nab_2- i \dual \nab_2)  U_{ab}\\
&=&(f_1+ i  f_2)(\nab_1- i  \nab_2) U_{ab}+(f_2-i  f_1)(\nab_2+ i  \nab_1)  U_{ab}\\
&=&2f_1\nab_1U_{ab} - 2i  f_1\nab_2U_{ab} +2i f_2\nab_1U_{ab} +2 f_2   \nab_2U_{ab}.
\eeaa
Therefore $\frac 1 2 (F \hot ( \DDb \c U ))_{ab}= (F\c \DDb U)_{ab}$
as stated. Finally,
\beaa
 F^c\ov{ \DD_c}U &=&  f^c\ov{ \DD_c}U +i( \dual f^c) \ov{ \DD_c}U = f^c\ov{ \DD_c}U - i f^c( \dual \ov{ \DD_c}U) = 2f^c\ov{ \DD_c}U =2 F^c \nab_c U.
 \eeaa

To prove the second formula in \eqref{Leib-eq-DDb}, we write from the above
\beaa
(\ov{F}\c \DD U)_{ab}&=&2f_1\nab_1U_{ab} + 2i  f_1\nab_2U_{ab} -2i f_2\nab_1U_{ab} +2 f_2   \nab_2U_{ab}
\eeaa
which implies
\beaa
(F\c \DDb U)_{ab}+(\ov{F}\c \DD U)_{ab}&=&4f_1\nab_1U_{ab} +4 f_2   \nab_2U_{ab}=4f \c \nab U
\eeaa
as stated.


\chapter{Complement for Chapter \ref{SECTION:KERR}}



\section{Proof of Proposition \ref{PROP:COMPUTATIONS-PI}}\label{sec:proof-prop-Pi}


We decompose the symmetric 3-tensor $\D_{(\mu} K_{\nu \rho)}$ in the null frame. We have
\beaa
3 \Pi_{\mu\nu\rho}&=&\D_{\mu} K_{\nu\rho}+\D_{\nu} K_{\rho\mu}+\D_{\rho} K_{\mu\nu} .
\eeaa
 We have
\beaa
3 \Pi_{abc}&=&\D_{a} K_{bc}+\D_{b} K_{c a}+\D_{c} K_{ab} \\
&=&\nab_{a} K_{bc}- K_{\D_a b c}-K_{b \D_ac}+\nab_{b} K_{c a}-K_{\D_b c a}- K_{c \D_b a}+\nab_{c} K_{ab} - K_{\D_c a b} - K_{a \D_c b}.
\eeaa
Since $\D_a e_b=\nab_a e_b+\frac 1 2 \chi_{ab} e_3+\frac 1 2  \chib_{ab}e_4$ and $K_{a3}=K_{a4}=0$, we have
\beaa
3 \Pi_{abc}&=&\nab_{a} K_{bc}+\nab_{b} K_{c a}+\nab_{c} K_{ab} \\
&=&\nab_{a}(r^2 \de_{bc})+\nab_{b} (r^2 \de_{c a})+\nab_{c}(r^2 \de_{ab} )=3\nab_{(a}r^2 \de_{bc)}.
\eeaa
Since in Kerr $\nab_a(r)=0$, $\Pi_{abc}=0$. 
 We have
\beaa
3 \Pi_{ab3}&=&\D_{a} K_{b3}+\D_{b} K_{3 a}+\D_{3} K_{ab} \\
&=&\nab_{a} K_{b3}-K_{\D_a b 3} - K_{b \D_a 3}+\nab_{b} K_{3 a}- K_{\D_b 3 a}-K_{3 \D_b a}+\nab_{3} K_{ab} -K_{\D_3 a b } - K_{a \D_3 b}\\
&=&- \frac 1 2 \chib_{ab} K_{4 3} - \chib_{ac} K_{b c}- \chib_{bc} K_{c a}-\frac 1 2 \chib_{ba} K_{3 4}+\nab_{3} K_{ab} \\
&=&- ( \chib_{ab}+\chib_{ba}) (a^2\cos^2\th)  - r^2 \chib_{ac} \de_{b c}- r^2\chib_{bc} \de_{c a}+\nab_{3} (r^2 \de_{ab} )\\
&=&- ( \chib_{ab}+\chib_{ba})(r^2+ a^2\cos^2\th)  +\nab_{3} (r^2  )\de_{ab}\\
&=& (2r \nab_{3} (r )-|q|^2 \trchb  ) \de_{ab} -2 \chibh_{ab}|q|^2,
\eeaa
which can also be written as
\beaa
3 \Pi_{ab3}&=&( \nab_{3} (r^2 )-|q|^2 \trchb  ) \de_{ab} -2 \chibh_{ab}|q|^2,
\eeaa
and similarly for $3 \Pi_{ab4}$. Using Lemma \ref{lemma:equations-q}, we see that $\Pi_{ab3}=\Pi_{ab4}=0$ in Kerr.
We have
\beaa
3 \Pi_{a 34}&=&\D_{a} K_{34}+\D_{3} K_{4 a}+\D_{4} K_{a3} \\
&=&\nab_{3} K_{4a}-K_{\D_3 4 a}-K_{4 \D_3 a}+\nab_{4} K_{a 3}-K_{\D_4 a 3}- K_{a \D_4 3}+\nab_{a} K_{34} -K_{\D_a 3 4}- K_{3 \D_a 4}\\
&=&-2\eta_b K_{ ab }-\eta_a  K_{34}-\etab_a K_{34}-2\etab_b K_{a b}+\nab_{a} K_{34} \\
&=&\nab_{a} K_{34} -2( \eta_b+\etab_b) K_{ ab }-(\eta_a +\etab_a) K_{34}\\
&=&\nab_{a} 2(a^2\cos^2\th) -2( \eta_b+\etab_b) r^2 \de_{ab}-(\eta_a +\etab_a)  2(a^2\cos^2\th)\\
&=&2\nab_{a} (a^2\cos^2\th) -2( \eta_a+\etab_a) |q|^2.
\eeaa
Using Lemma \ref{lemma:equations-q}, we see that $\Pi_{a34}=0$ in Kerr.
We have
\beaa
3 \Pi_{a 33}&=&\D_{a} K_{33}+2\D_{3} K_{3 a}\\
&=&\nab_{a} K_{33}- 2 K_{\D_a 3 3}+2\nab_{3} K_{3 a}-2 K_{\D_3 3 a}- 2 K_{3 \D_3 a}\\
&=&-4\xib_b K_{ b a}- 2\xib_a K_{3  4}\\
&=&-4\xib_b r^2 \de_{ b a}- 4\xib_a (a^2 \cos^2\th)= - 4 |q|^2 \xib_a,
\eeaa
and similarly for $3\Pi_{a44}$.  
We have
\beaa
3 \Pi_{34 3}&=&2\D_{3} K_{43}+\D_{4} K_{33}\\
&=&2\nab_{3} K_{43}-2 K_{\D_3 43} -2 K_{4 \D_3 3}+\nab_{4} K_{33}-2K_{\D_4 3 3}\\
&=&2\nab_{3} K_{43}-4\omb  K_{ 43} +4\omb K_{43}\\
&=&4\nab_{3}(a^2\cos^2\th),
\eeaa
and similarly for $\Pi_{43 4}$. 
Consider the component $\Pi_{333}$. We have
\beaa
\Pi_{333}&=& \D_3 K_{33}=\nab_3 K_{33} - 2 K_{\D_3 3 3}=0,
\eeaa
and similarly for $\Pi_{444}$.


\section{Proof of Proposition \ref{PROP:CARTER-OPERATOR}}\label{sec:proof-carter-operator}


The operator $\KK$ applied to $\psi \in \sk_k$ is given by
\beaa
\KK(\psi)&=& \Ddot_{\mu}( K^{\mu\nu} \Ddot_\nu(\psi))=\Ddot_{\mu} K^{\mu\nu} \Ddot_\nu(\psi)+ K^{\mu\nu} \Ddot_{\mu}\Ddot_\nu(\psi).
\eeaa
Using the general computations in \eqref{divergence-K-trace} for $\Pi=0$, we have
\beaa
\KK(\psi)&=&-\frac 1 2 \Ddot^\nu (\tr K) \Ddot_\nu(\psi)+ K^{\mu\nu} \Ddot_{\mu}\Ddot_\nu(\psi) \\
&=&-\frac 1 2\de^{ab} \Ddot_a (\tr K) \Ddot_b(\psi)-\frac 1 2\g^{34} \Ddot_4 (\tr K) \Ddot_3(\psi)-\frac 1 2\g^{43} \Ddot_3 (\tr K) \Ddot_4(\psi)\\
&&+ (-(a^2\cos^2\th)  \g^{\mu\nu} +|q|^2\ga^{ab} e_a^{\mu} e_b^{\nu}) \Ddot_{\mu}\Ddot_\nu(\psi),
\eeaa
which gives
\beaa
\KK(\psi)&=&-(a^2\cos^2\th) \squared_k \psi -\frac 1 2 \nab^c (\tr K) \nab_c(\psi)+\frac 1 4 \nab_4 (\tr K) \nab_3(\psi)+\frac 1 4 \nab_3 (\tr K) \nab_4(\psi)\\
&&+ |q|^2 \ga^{ab} \Ddot_a \Ddot_b \psi.
\eeaa
Using that $\tr K= 2(r^2 -a^2\cos^2\th)$, we compute
\beaa
 \nab^c (\tr K)&=& 2\nab^c(r^2)-2\nab^c ( a^2\cos^2\th)=-2( \eta^c+\etab^c) |q|^2.
\eeaa
 Similarly we compute
 \beaa
  \nab_3(\tr K)&=& 4r\nab_3r -2\nab_3(a^2\cos^2\th)=2 |q|^2 \trchb, \\
    \nab_4(\tr K)&=&2 |q|^2 \trch.
\eeaa
The above gives 
\beaa
\KK(\psi)&=&-(a^2\cos^2\th) \squared_k \psi -\frac 1 2 \left(  -2 ( \eta^c+\etab^c) |q|^2 -3 {\Pi^{c}}_{ 34} \right) \nab_c\psi\\
&&+\frac 1 4\left(2|q|^2 \trch \right) \nab_3\psi+\frac 1 4 \left(2 |q|^2 \trchb \right) \nab_4\psi+ |q|^2 \de^{ab}( \nab_b \nab_a\psi- \nab_{\D_b e_a} \psi)\\
&=&-(a^2\cos^2\th) \squared_k \psi +   |q|^2 ( \eta+\etab) \c  \nab\psi+\frac 1 2 |q|^2 \trch  \nab_3\psi+\frac 1 2 |q|^2 \trchb \nab_4\psi\\
&&+ |q|^2\de^{ab}( \nab_b \nab_a\psi -\frac 1 2 \chi_{ba} \nab_{3} \psi -\frac 1 2 \chib_{ba} \nab_{4}\psi )\\
&=&-(a^2\cos^2\th) \squared_k \psi + |q|^2\lap_k \psi +   |q|^2    ( \eta+\etab) \c \nab \psi,
\eeaa
as stated.


\section{Proof of Proposition \ref{LEMMA:MOD-LAPLACIAN-KERR}}\label{sec:proof-mod-lap}


Recall from Lemma \ref{lemma:expression-wave-operator} that we have 
\beaa
|q|^2\squared_2 \psi&=&-\frac 1 2|q|^2 \big(\nab_3\nab_4\psi+\nab_4 \nab_3 \psi\big)+|q|^2\left(\omb -\frac 1 2 \trchb\right) \nab_4\psi+|q|^2\left(\om -\frac 1 2 \trch\right) \nab_3\psi\\
&& +\OO(\psi).
\eeaa
Using the ingoing frame to write $e^{(in)}_4=\frac{r^2+a^2}{|q|^2} \big(\That+\Rhat\big)$, see \eqref{eq-incoming-3-4-T-R}, we have
\beaa
\nab_3\nab_4+\nab_4 \nab_3 &=&\nab_3\Big(\frac{r^2+a^2}{|q|^2} \big(\nab_\That+\nab_\Rhat\big)\Big)+\frac{r^2+a^2}{|q|^2} \big(\nab_\That+\nab_\Rhat\big) \nab_3 \\
&=&\frac{r^2+a^2}{|q|^2} \Big( \nab_\That \nab_3+ \nab_3 \nab_\That+\nab_\Rhat \nab_3+ \nab_3\nab_\Rhat  \Big) \\
&&+e^{(in)}_3\big(\frac{r^2+a^2}{|q|^2}\big) \big(\nab_\That+\nab_\Rhat\big)\\
&=&\frac{r^2+a^2}{|q|^2} \Big( \nab_\That \nab_3+ \nab_3 \nab_\That+\nab_\Rhat \nab_3+ \nab_3\nab_\Rhat  \Big) +O(a^2 r^{-3}) \big(\nab_\That+\nab_\Rhat\big).
\eeaa
We can also write in the ingoing frame, using \eqref{eq:ThatRhat-e_3e_4-Kerr}, 
\beaa
&&|q|^2\left(\omb -\frac 1 2 \trchb\right) \nab_4+|q|^2\left(\om -\frac 1 2 \trch\right) \nab_3\\
&=&|q|^2\left(\frac{r}{|q|^2}\right) \nab_4+|q|^2\left(-\frac 1 2 \pr_r \big(\frac{\De}{|q|^2} \big) -\frac{r\De}{|q|^4}\right) \nab_3 \\
&=& r \big(  \nab_4 -\frac{\De}{|q|^2} \nab_3 \big) -\frac 1 2 |q|^2\pr_r \big(\frac{\De}{|q|^2} \big) \nab_3\\
&=& 2\frac{r(r^2+a^2)}{|q|^2} \nab_\Rhat  -\frac 1 2 |q|^2\pr_r \big(\frac{\De}{|q|^2} \big) \nab_3.
\eeaa
We deduce in the ingoing frame
\beaa
|q|^2\squared_2 \psi&=&-\frac 1 2 (r^2+a^2) \Big( \nab_\That \nab_3+ \nab_3 \nab_\That+\nab_\Rhat \nab_3+ \nab_3\nab_\Rhat  \Big)\psi \\
&&+2\frac{r(r^2+a^2)}{|q|^2} \nab_\Rhat \psi -\frac 1 2 |q|^2\pr_r \big(\frac{\De}{|q|^2} \big) \nab_3\psi+O(a^2 r^{-1}) \big(\nab_\That+\nab_\Rhat\big)\psi+\OO(\psi),
\eeaa
and therefore
\beaa
[\OO, |q|^2\squared_2] \psi&=&-\frac 1 2 (r^2+a^2)[\OO, \big( \nab_\That \nab_3+ \nab_3 \nab_\That+\nab_\Rhat \nab_3+ \nab_3\nab_\Rhat  \big)]\psi \\
&&+ \frac{2r(r^2+a^2)}{|q|^2} [\OO,\nab_\Rhat] \psi -\frac 1 2|q|^2\pr_r \big(\frac{\De}{|q|^2} \big) [\OO,  \nab_3]\psi\\
&&+O(a^2 r^{-1})[\OO, \big(\nab_\That+\nab_\Rhat\big)]\psi+ O(a^2r^{-2}) \dk^{\leq 1} \psi .
\eeaa
Using Lemma \ref{LEMMA:COMMUTATOR-NAB3-NAB4-LAP}, Lemma \ref{lemma:commutationpropertiesofthesymmetryoperatorsafsoiudf:chap9} and  Corollary \ref{cor:commutatornabRhatwith|q|nab}, i.e.
\beaa
\, [ \OO, \nab_3]\psi&=& O(ar^{-2}) \dk^{\leq 1} \psi, \\
\,[\OO, \nab_\Rhat]\psi &=& O(ar^{-2})\dk^{\leq 1}\psi, \\
\,[\OO, \nab_\That] &=& O(ar^{-2})\dkb^{\leq 1}\psi,
\eeaa
we obtain
\beaa
[\OO, |q|^2\squared_2] \psi&=&-\frac 1 2 (r^2+a^2) \Big([\OO,  \nab_\That ]\nab_3+ \nab_\That[\OO, \nab_3]+ [\OO, \nab_3] \nab_\That + \nab_3[\OO, \nab_\That]  \big)\psi \Big)\\
&&-\frac 1 2 (r^2+a^2) \Big([\OO,  \nab_\Rhat ]\nab_3+ \nab_\Rhat[\OO, \nab_3]+ [\OO, \nab_3] \nab_\Rhat + \nab_3[\OO, \nab_\Rhat]  \big)\psi \Big)\\
&&+ O(ar^{-1})\dk^{\leq 1}\psi\\
&=&-\frac 1 2 (r^2+a^2) \Big([\OO,  \nab_\That ]\nab_3+ \nab_\That[\OO, \nab_3]+ [\OO, \nab_3] \nab_\That + \nab_3[\OO, \nab_\That]  \Big)\psi \\
&&-\frac 1 2 (r^2+a^2) \Big([\OO,  \nab_\Rhat ]\nab_3+ \nab_3[\OO, \nab_\Rhat]  \Big)\psi \\
&&+O(a)\nab_\Rhat \dk^{\leq1}\psi+ O(ar^{-1})\dk^{\leq 1}\psi.
\eeaa
Observe that, since the commutators $[ \OO, \nab_3]$, $[ \OO, \nab_\Rhat]$, $[ \OO, \nab_\That]$ involve only one derivative, then
\beaa
\nab_3[\OO, \nab_\That]=[\OO, \nab_\That]\nab_3+O(ar^{-3})\dk^{\leq1}
\eeaa
and similarly for the others. We then have
\beaa
[\OO, |q|^2\squared_2] \psi&=&- (r^2+a^2) \Big( [\OO,  \nab_\That ]\nab_3+ [\OO, \nab_3] \nab_\That +[\OO,  \nab_\Rhat ]\nab_3  \Big)\psi \\
&&+O(a)\nab_\Rhat \dk^{\leq1}\psi+ O(ar^{-1})\dk^{\leq 1}\psi.
\eeaa
Now using \eqref{eq-incoming-3-4-T-R} to write $\Rhat=  \That-\frac{\De}{r^2+a^2} e_3^{(in)}$, we obtain
\beaa
[\OO, |q|^2\squared_2] \psi&=&- (r^2+a^2) \Big( [\OO,  \nab_\That ]\nab_3+ [\OO, \nab_3] \nab_\That +[\OO,   \nab_\That-\frac{\De}{r^2+a^2} \nab_3]\nab_3  \Big)\psi \\
&&+O(a)\nab_\Rhat \dk^{\leq1}\psi+ O(ar^{-1})\dk^{\leq 1}\psi\\
&=&- (r^2+a^2) \Big(2 [\OO,  \nab_\That ]\nab_3+ [\OO, \nab_3]\big( \nab_\That -\frac{\De}{r^2+a^2}\nab_3 \big) \Big)\psi \\
&&+O(a)\nab_\Rhat \dk^{\leq1}\psi+ O(ar^{-1})\dk^{\leq 1}\psi\\
&=&- 2(r^2+a^2)  [\OO,  \nab_\That ]\nab_3\psi +O(a)\nab_\Rhat \dk^{\leq1}\psi+ O(ar^{-1})\dk^{\leq 1}\psi.
\eeaa

Next we compute $[\OO,  \nab_\That ]$.
Recall that we have from Lemma \ref{lemma:basicpropertiesLiebTfasdiuhakdisug:chap9},
\beaa
\nab_\T\psi &=& \Lieb_\T\psi + \frac{4amr\cos\th}{|q|^4}\dual\psi.
\eeaa
We infer
\beaa
[\nab_\T, \OO]\psi &=& \left[ \Lieb_\T\psi + \frac{4amr\cos\th}{|q|^4}\dual, \OO\right]\psi\\
&=& 4amr\left[ \frac{\cos\th}{|q|^4}, \OO\right]\dual\psi\\
&=& 4amr|q|^2\left[ \frac{\cos\th}{|q|^4}, \Delta+(\eta+\etab)\c\nab\right]\dual\psi\\
&=& -8amr|q|^2\nab\left(\frac{\cos\th}{|q|^4}\right)\c\nab\dual\psi+O(a r^{-3})\psi.
\eeaa
Also, recall that we have 
\beaa
\nab_\Z\psi &=& \Lieb_\Z\psi - \frac{2\cos\th((r^2+a^2)^2-a^2(\sin\th)^2\De)}{|q|^4}\dual\psi.
\eeaa
We infer
\beaa
[\nab_\Z, \OO]\psi &=& \left[ \Lieb_\Z\psi - \frac{2\cos\th((r^2+a^2)^2-a^2(\sin\th)^2\De)}{|q|^4}\dual, \OO\right]\psi\\
&=&- 2\left[ \frac{\cos\th((r^2+a^2)^2-a^2(\sin\th)^2\De)}{|q|^4}, \OO\right]\dual\psi\\
&=& - 2|q|^2\left[ \frac{\cos\th((r^2+a^2)^2-a^2(\sin\th)^2\De)}{|q|^4}, \Delta+(\eta+\etab)\c\nab\right]\dual\psi\\
&=&  4|q|^2\nab\left(\frac{\cos\th((r^2+a^2)^2-a^2(\sin\th)^2\De)}{|q|^4}\right)\c\nab\dual\psi+O(1)\psi.
\eeaa

Next, recall that $\That=\T+\frac{a}{r^2+a^2}\Z$ which together with the above yields 
\beaa
[\nab_\That, \OO]\psi &=& [\nab_\T, \OO]\psi+\frac{a}{r^2+a^2}[\nab_\Z, \OO]\psi\\
&=&- 8amr|q|^2\nab\left(\frac{\cos\th}{|q|^4}\right)\c\nab\dual\psi\\
&&+ \frac{4a|q|^2}{r^2+a^2}\nab\left(\frac{\cos\th((r^2+a^2)^2-a^2(\sin\th)^2\De)}{|q|^4}\right)\c\nab\dual\psi +O(a r^{-2})\psi.
\eeaa
Observe that we can simplify
\beaa
f(r, \cos\th)&:=& - 8amr\frac{\cos\th}{|q|^4}+ \frac{4a}{r^2+a^2}\frac{\cos\th((r^2+a^2)^2-a^2(\sin\th)^2\De)}{|q|^4}\\
&=&\frac{4a\cos\th}{|q|^4(r^2+a^2)} \Big(  - 2mr(r^2+a^2)+ (r^2+a^2)^2-a^2(\sin\th)^2\De \Big) \\
&=&\frac{4a\cos\th}{|q|^4(r^2+a^2)} \Big(  (r^2+a^2)\De-a^2(\sin\th)^2\De \Big) \\
&=&\frac{4a\cos\th \De}{|q|^2(r^2+a^2)}  .
\eeaa
We deduce
\beaa
[\nab_\That, \OO]\psi &=&\frac{|q|^2\De}{r^2+a^2}\nab\left(\frac{4a\cos\th }{|q|^2}  \right)\c\nab\dual\psi +O(a r^{-2})\psi.
\eeaa
This gives
\beaa
[\OO, |q|^2\squared_2] \psi&=& |q|^2\nab\left(\frac{8a\cos\th }{|q|^2}  \right)\c\nab(\De \nab_3)  \dual \psi +O(a)\nab^{\leq 1}_\Rhat \dk^{\leq1}\psi.
\eeaa
By writing $\De e_3^{(in)}=(r^2+a^2) (\That - \Rhat)$, we finally obtain
\beaa
[\OO, |q|^2\squared_2] \psi&=& |q|^2\nab\left(\frac{8a(r^2+a^2) \cos\th }{|q|^2}  \right)\c\nab \nab_\That  \dual \psi +O(a)\nab^{\leq 1}_\Rhat \dk^{\leq1}\psi,
\eeaa
as stated.
This concludes the proof of the lemma.


\chapter{Complement for Chapter \ref{PERTURBATIONS-SECTION}}



\section{Proof of Lemma \ref{LEMMA:COMMUTATION-FORMULAS-1}}
\label{sec:proof-lemma-comm-form-1}


Formula \eqref{eq:comm-nab4-nab3-DD-h-precise} is straightforward from \eqref{eq:comm-nab3-nab4-naba-f-general}.

According to Corollary \ref{corr:comm}, from \eqref{last-statement-item1} we have
\beaa
\, [\nab_4, \nab \hot] f&=&- \frac 1 2 \trch \left( \nab\hot f +\etab \hot f\right)- \frac 1 2 \atrch\,  \dual  \left(\nab \hot  f+\etab \hot  f\right)+ (\etab+\ze) \hot \nab_4 f, \\
&&+\dual \b \hot \dual f+\xi \hot \nab_3 f-\xi \hot ( \chib \c f)+\chibh \,(\xi\c f)-\chih \c \nab f\\
&&-\etab \hot (\chih \c f)+\chih (\etab\c f).
\eeaa
Hence for $F=f+i \dual f$,
\beaa
\, [\nab_4, \nab \hot] F&=&- \frac 1 2 \trch \left( \nab\hot F +\etab \hot F\right)- \frac 1 2 \atrch\,  \dual  \left(\nab \hot  F+\etab \hot  F\right)+ (\etab+\ze) \hot \nab_4 F, \\
&&+\dual \b \hot \dual F+\xi \hot \nab_3 F-\xi \hot ( \chib \c F)+\chibh \,(\xi\c F)-\chih \c \nab F\\
&&-\etab \hot (\chih \c F)+\chih (\etab\c F).
\eeaa
Recalling that  $ \dual F=-iF$  and $\dual (\nab \hot f)=\dual \nab\hot f=\nab\hot \dual f$ we deduce,
 \beaa
 \, [\nab_4, \nab \hot] F&=&- \frac 1 2 \trch \left( \nab\hot F +\etab \hot F\right)+ i \frac 1 2 \atrch  \left(\nab \hot  F+\etab \hot F\right)
 + (\etab+\ze) \hot \nab_4 F\\
 && -i \dual \b \hot F+\xi \hot \nab_3 F-\xi \hot ( \chib \c F)+\chibh \,(\xi\c F)-\chih \c \nab F\\
 &&-\etab \hot (\chih \c F)+\chih (\etab\c F)\\
 &=&-\frac 1 2 \tr X\left(\nab \hot  F+\etab \hot F\right)  + (\etab+\ze) \hot \nab_4 F\\
  && -i \dual \b \hot F+\xi \hot \nab_3 F-\xi \hot ( \chib \c F)+\chibh \,(\xi\c F)-\chih \c \nab F\\
  &&-\etab \hot (\chih \c F)+\chih (\etab\c F).
 \eeaa
Taking the dual
\beaa
\, [\nab_4,  \dual \nab \hot] F &=&-\frac 1 2 \tr X\left( \dual \nab \hot  F+\dual \etab \hot F\right)  +\dual (\etab+\ze) \hot \nab_4 F +i  \b \hot F+\dual \xi \hot \nab_3 F\\
&&-\dual \xi \hot ( \chib \c F)+\dual\chibh \,(\xi\c F)+\chih \c \dual \nab F-\dual \etab \hot (\chih \c F)+\dual \chih (\etab\c F).
\eeaa
Finally adding the above we derive
\beaa
 [\nab_4,  \DD \hot] F &=&-\frac 1 2 \tr X\left( \DD \hot  F+\Hb\hot F\right)+(\Hb+Z)\hot\nab_4 F+ \Xi \hot \nab_3 F \\
 &&-B \hot F -\Xi \hot ( \chib \c F)-\chih \c \ov{\DD} F-\Hb \hot (\chih \c F)+\Xh (\etab\c F)+\Xbh \,(\xi\c F)\\
 &=&-\frac 1 2 \tr X\left( \DD \hot  F+\Hb\hot F\right)+(\Hb+Z)\hot\nab_4 F+ \Xi \hot \nab_3 F \\
 &&-B \hot F - \frac 1 2 \tr \Xb \Xi \hot  F-\frac 1 2\Xh \c \ov{\DD} F+\frac12\Xh (\ov{\Hb}\c F)+ (\Ga_b \c \Ga_g) F,
\eeaa
where observe that $\Hb \hot (\chih \c F)=0$. 
By symmetry we also have
\beaa
 [\nab_3,  \DD \hot] F  &=&-\frac 1 2 \tr \Xb\left( \DD \hot  F+H \hot F\right)+(H-Z)\hot\nab_3 F+ \Xib \hot \nab_4 F \\
 &&+\Bb \hot F - \frac 1 2 \tr X \Xib \hot  F-\frac 1 2\Xbh \c \ov{\DD} F+\frac12\Xbh (\ov{H}\c F)+ (\Ga_b \c \Ga_g) F.
\eeaa
This proves \eqref{eq:comm:nab4-nab3-DDhot-precise}.

According to Corollary \ref{corr:comm}, from \eqref{eq:comm-nab3-nab4-div-all-errors} we have
\beaa
\,[\nab_4,  \div] u &=&-\frac 1 2 \trch \big(  \div u - 2\etab \c u\big) +\frac 1 2 \atrch\big(  \div\dual u -2 \etab\c \dual u\big) +(\etab+\ze)\c\nab_4 u \\
&&+2\dual \b \c \dual u +\xi\c \nab_3 u -\xi\c  \chib\c  u -(\chib\c u)\xi +\xi\c u\c\chib -\chih\c \nab u\\
&& -\etab\c \chih \c u-(\chih\c u)\etab+\etab\c u\c\chih.
\eeaa
Hence for $U= u + i \dual u$,
\beaa
\,[\nab_4,  \div] U &=&-\frac 1 2 \trch \big(  \div U - 2\etab \c U\big) +\frac 1 2 \atrch\big(  \div\dual U -2 \etab\c \dual U\big) +(\etab+\ze)\c\nab_4 U \\
&&+2\dual \b \c \dual U +\xi\c \nab_3 U -\xi\c  \chib\c  U -(\chib\c U)\xi +\xi\c U\c\chib -\chih\c \nab U\\
&& -\etab\c \chih \c U-(\chih\c U)\etab+\etab\c U\c\chih.
\eeaa
Recalling that $\dual U= - i U$ and $\dual (\div U)=\div \dual U$ we deduce
\beaa
\,[\nab_4,  \div] U &=&-\frac 1 2 \trch \big(  \div U - 2\etab \c U\big) -\frac 1 2 i\atrch\big(  \div U -2 \etab\c U\big) +(\etab+\ze)\c\nab_4 U \\
&&-2i \dual \b \c U +\xi\c \nab_3 U -\xi\c  \chib\c  U -(\chib\c U)\xi +\xi\c U\c\chib -\chih\c \nab U\\
&& -\etab\c \chih \c U-(\chih\c U)\etab+\etab\c U\c\chih\\
&=&-\frac 1 2\ov{\tr X} \big(  \div U - 2\etab \c U\big) +(\etab+\ze)\c\nab_4 U \\
&&-2i \dual \b \c U +\xi\c \nab_3 U -\xi\c  \chib\c  U -(\chib\c U)\xi +\xi\c U\c\chib -\chih\c \nab U\\
&& -\etab\c \chih \c U-(\chih\c U)\etab+\etab\c U\c\chih.
\eeaa
Taking the dual
\beaa
\,[\nab_4, \dual \div] U &=&-\frac 1 2\ov{\tr X} \big( \dual \div U - 2\dual \etab \c U\big) +\dual (\etab+\ze)\c\nab_4 U +2i  \b \c U +\dual \xi\c \nab_3 U\\
&& -\dual \xi\c  \chib\c  U -\dual (\chib\c U)\xi +\dual \xi\c U\c\chib -\chih\c \dual  \nab U\\
&& -\dual \etab\c \chih \c U-\dual(\chih\c U)\etab+\dual\etab\c U\c\chih.
\eeaa
Finally we derive for $\ov{\DD}\c=\div - i \dual \div$, 
\beaa
\,[\nab_4,  \ov{\DD} \c ] U &=&-\frac 1 2\ov{\tr X} \big(  \ov{\DD}\c U - 2\ov{\Hb} \c U\big) +(\Hb + Z)\c\nab_4 U+ \ov{\Xi} \c \nab_3 U \\
&&+2\ov{B} \c U  -\ov{\Xi}\c  \chib\c  U -(\chib\c U)\ov{\Xi} +\ov{\Xi}\c U\c\chib -\chih\c \ov{\DD} U\\
&& -\ov{\Hb} \c \chih \c U-(\chih\c U)\ov{\Hb}+\ov{\Hb} \c U\c\chih\\
&=&-\frac 1 2\ov{\tr X} \big(  \ov{\DD}\c U - 2\ov{\Hb} \c U\big) +(\Hb + Z)\c\nab_4 U+ \ov{\Xi} \c \nab_3 U \\
&&+2\ov{B} \c U  -\frac 1 2 \ov{\tr \Xb} \ov{\Xi}\c  U  -\frac 1 2 \Xh \c \ov{\DD} U-\frac 1 2 (\ov{\Xh}\c U)\ov{\Hb}+ (\Ga_b \c \Ga_g) U.
\eeaa
By symmetry we also have
\beaa
\,[\nab_3,  \ov{\DD} \c ] U &=&-\frac 1 2\ov{\tr \Xb} \big(  \ov{\DD}\c U - 2\ov{H} \c U\big) +(H- Z)\c\nab_3 U+ \ov{\Xib} \c \nab_4 U \\
&&-2\ov{\Bb} \c U  -\frac 1 2 \ov{\tr X} \ov{\Xib}\c  U  -\frac 1 2 \Xbh \c \ov{\DD} U-\frac 1 2 (\ov{\Xbh}\c U)\ov{H}+ (\Ga_b \c \Ga_g) U.
\eeaa
This proves  \eqref{eq:comm:nab4-nab3-DDc-precise}.

   Recall from \eqref{commutator-3-a-u-bc}, 
\beaa
\,[\nab_3, \nab_4]U_{ab} &=& - 2\om \nab_3 U_{ab} + 2\omb \nab_4 U_{ab}  + 2 (\eta_c-\etab_c) \nab_c U_{ab}\\
&&+2\etab_a\eta_c U_{bc}+2\etab_b \eta_c U_{ac}- 2\eta_a\etab_c U_{bc}-2\eta_b \etab_c U_{ac}+4\rhod \dual U_{ab}+\err_{34ab}[U].
\eeaa
Applying it to $U \in \sk_k(\mathbb{C})$ for which $\dual U=-i U$, we have
\beaa
\,[\nab_3, \nab_4]U_{ab} &=& - 2\om \nab_3 U_{ab} + 2\omb \nab_4 U_{ab}  + 2 (\eta_c-\etab_c) \nab_c U_{ab}-4i \rhod U_{ab} +C^{\eta, \etab}_{3, 4}(U)_{ab}\\
&&+\err_{34ab}[U],
\eeaa
where
\beaa
C^{\eta, \etab}_{3, 4}(U)_{ab}&=& 2\etab_a\eta_c U_{bc}+2\etab_b \eta_c U_{ac}- 2\eta_a\etab_c U_{bc}-2\eta_b \etab_c U_{ac}.
\eeaa
We have that 
\beaa
C^{\eta, \etab}_{3, 4}(U)&=&4i (\eta \wedge \etab  ) U.
\eeaa
Indeed, we compute
\beaa
C^{\eta, \etab}_{3, 4}(U)_{11}&=& 2\etab_1\eta_c U_{1c}+2\etab_1 \eta_c U_{1c}- 2\eta_1\etab_c U_{1c}-2\eta_1 \etab_c U_{1c}\\
&=& 4\etab_1\eta_c U_{1c}- 4\eta_1\etab_c U_{1c}\\
&=& 4\etab_1(\eta_1 U_{11}+\eta_2 U_{12})- 4\eta_1(\etab_1 U_{11}+\etab_2 U_{12})\\
&=&4( \etab_1\eta_2 - \eta_1\etab_2 )U_{12}=4i(  \eta_1\etab_2 -\etab_1\eta_2)U_{11}.
\eeaa
Also,
\beaa
C^{\eta, \etab}_{3, 4}(U)_{12}&=& 2\etab_1\eta_c U_{2c}+2\etab_2 \eta_c U_{1c}- 2\eta_1\etab_c U_{2c}-2\eta_2 \etab_c U_{1c}\\
&=& 2\etab_1(\eta_1 U_{21}+\eta_2 U_{22})+2\etab_2( \eta_1 U_{11}+ \eta_2 U_{12})\\
&&- 2\eta_1(\etab_1 U_{21}+\etab_2 U_{22})-2\eta_2( \etab_1 U_{11}+ \etab_2 U_{12})\\
&=& 2\etab_1\eta_2 U_{22}+2\etab_2 \eta_1 U_{11}- 2\eta_1\etab_2 U_{22}-2\eta_2 \etab_1 U_{11}\\
&=& -4\etab_1\eta_2 U_{11}+ 4\eta_1\etab_2 U_{11}=4i ( \eta_1\etab_2 -\etab_1\eta_2  ) U_{12}.
\eeaa
We conclude by writing $\eta_1\etab_2 -\etab_1\eta_2 =\eta \wedge \etab$.

Starting with \eqref{commutator-3-a-u-bc}
\beaa
\,  [\nab_3,\nab_a] u_{bc}    &=&-\frac 1 2 \trchb\, (\nab_a u_{bc}+\eta_bu_{ac}+\eta_c u_{ab}-\de_{a b}(\eta \c u)_c-\de_{a c}(\eta \c u)_b )\\
&&-\frac 1 2 \atrchb\, (\dual \nab_a u_{bc} +\eta_b\dual u_{ac}+\eta_c\dual u_{ab}- \in_{a b}(\eta \c u)_c- \in_{a c}(\eta \c u)_b )\\
&&+(\eta_a-\ze_a)\nab_3 u_{bc}+\err_{3abc}[u],
\eeaa
and adding the  same expression  for $u$ replaced with $i \dual u$ we derive
\beaa
\,  [\nab_3,\nab_a] U_{bc}    &=&-\frac  1 2   \trchb\, \Big(\nab_a U_{bc}+\eta_bU_{ac}+\eta_c U_{ab}-\de_{a b}(\eta \c U)_c-\de_{a c}(\eta \c U)_b \Big)\\
&-&\frac 1 2 \atrchb\, \Big(\dual \nab_a  U_{bc} +\eta_b \dual U_{ac}+\eta_c \dual U_{ab}- \in_{a b}(\eta \c  U)_c- \in_{a c}(\eta \c  U)_b \Big)\\
 &+&(\eta_a-\ze_a)\nab_3 U_{bc}+ \err_{3abc}[U]
\eeaa
as stated.


\section{Proof of Lemma \ref{COMMUTATOR-NAB-C-3-DD-C-HOT}}\label{sec:proof-lemma-comm-2}


Using \eqref{eq:comm-nab4-nab3-DD-h-precise}, we have 
\beaa
 \, [\nabc_4 , \DDc]h &=& (\nab_4+2s\om)(\DD+sZ)h- (\DD+(s+1)Z)(\nab_4+2s\om)h\\
 &=& [\nab_4, \DD]h+2s\om\DDc h+s\nab_4(Zh) -(s+1)Z\nabc_4h -2s\DD(\om h)\\
 &=& -\frac{1}{2}\tr X\DD h+(\Hb+Z)\nab_4 h -\chih\c\ov{\DD} h+\Xi \nab_3h +2s\om\DDc h+s(\nab_4Z)h\\
 &&+sZ \nab_4h -(s+1)Z\nabc_4h -2s\DD(\om)h -2s\om\DD h.
\eeaa
Using the null structure equation for $\nab_4Z$, i.e.
\beaa
\nab_4Z +\frac{1}{2}\tr X(Z-\Hb)-2\om(Z+\Hb) &=& 2\DD\om +\frac{1}{2}\widehat{X}\c(-\ov{Z}+\ov{\Hb})\\
&&-\frac{1}{2}\tr\Xb\Xi-2\omb\Xi -B-\frac{1}{2}\ov{\Xi}\c\Xbh,
\eeaa
we infer
\beaa
 \, [\nabc_4 , \DDc]h  &=& -\frac{1}{2}\tr X\DDc h+(\Hb+Z)\nab_4 h -\frac{1}{2}\Xh\c\ov{\DDc} h+\Xi \nab_3h\\
&& +2s\om\DDc h +s\left(\frac{1}{2}\tr X\Hb+2\om\Hb  +\frac{1}{2}\widehat{X}\c\ov{\Hb}-\frac{1}{2}\tr\Xb\Xi-2\omb\Xi -B\right)h\\
&& +sZ\nabc_4h -(s+1)Z\nabc_4h  -2s\om\DD h+ (\Ga_b \c \Ga_g) h \\
&=& -\frac{1}{2}\tr X\DDc h+(\Hb+Z)\nab_4 h -\frac{1}{2}\Xh\c\ov{\DDc} h+\Xi \nab_3h\\
&& +2s\om\DDc h +s\left(\frac{1}{2}\tr X\Hb+2\om\Hb  +\frac{1}{2}\widehat{X}\c\ov{\Hb}-\frac{1}{2}\tr\Xb\Xi-2\omb\Xi -B\right) h\\
&&  - Z(\nab_4+2s\om)h  -2s\om\DD h+ (\Ga_b \c \Ga_g) h\\
&=& -\frac{1}{2}\tr X\DDc h+\Hb\nabc_4 h -\frac{1}{2}\Xh\c\ov{\DDc} h+\Xi\nabc_3h\\
&&  +s\left(\frac{1}{2}\tr X\Hb  +\frac{1}{2}\widehat{X}\c\ov{\Hb}-\frac{1}{2}\tr\Xb\Xi -B\right)h+ (\Ga_b \c \Ga_g) h,
\eeaa
as stated. Similarly for $[\nabc_3, \DDc]h$. This proves \eqref{eq:comm-nabc4-DDc-h-precise}.

For a  $s$-conformally invariant scalar $h$,  we  have
\beaa
\,[ \nabc_3, \nabc_4] h&=& \nabc_3 \nabc_ 4  h - \nabc_4 \nabc_3  h\\
 &=& \big( e_3   -2 (s+1) \omb \big) \big ( e_4 h +2 s \om h\big)-
 \big( e_4 + 2(s-1) \om \big)( e_3 h - 2 s \omb h \big)\\
&=&  e_3 e_4 h- 2 (s+1) \omb e_4 h + 2 s e_3(\om h )   -  4 s  (s+1)\omb \om  h \\
&-& e_4 e_3 h - 2 (s-1) \om e_3 h + 2 s e_4(\omb h) + 4 s (s-1) \om \omb h \\
&=&[e_3, e_4] h - 2\omb e_4 h + 2\om   e_3h +2 s \Big( e_3 \om + e_4 \omb  - 4 \om \omb \Big) h. 
\eeaa
Using \eqref{eq:comm-nab3-nab4-naba-f-general} and the null structure equation $\nab_3\om+\nab_4\omb=   \rho  +4\om\omb +\xi\c \xib +(\eta-\etab)\c\ze -\eta\c\etab$, we obtain
\beaa
\,[ \nabc_3, \nabc_4] h&=& 2(\eta-\etab)\c \nab h+2s\Big(\rho+ (\eta-\etab)\c\ze -\eta\c\etab+\xi\c \xib \Big) h\\
&=&2(\eta-\etab)\c\nabc h  +2s\Big(\rho -\eta\c\etab+\xi\c \xib \Big) h,
\eeaa
as stated. This proves \eqref{eq:comm-nabc3-nabc4-h}.

For $F \in \sk_1 (\mathbb{C})$ $s$-conformally invariant we have
 \beaa
 \nabc_4 (\DDc \hot F )&=&\nabc_4 (\DD\hot F+s Z \hot F )\\
 &=&\nab_4 (\DD\hot F+s Z \hot F )+2s \om (\DD\hot F+s Z \hot F ) \\
  &=&\nab_4 (\DD\hot F)+s \nab_4(Z) \hot F +s Z \hot \nab_4F +2s \om (\DD\hot F+s Z \hot F ).
 \eeaa
 Using \eqref{eq:comm:nab4-nab3-DDhot-precise} and the null structure equation for $\nab_4 Z$,  we obtain
 \beaa
 \nabc_4 (\DDc \hot F )  &=&\mathcal{D}\hot(\nab_4F)-\frac 1 2 \tr X\left( \DD \hot  F+\Hb\hot F\right)+(\Hb+Z)\hot\nab_4 F+ \Xi \hot \nab_3 F \\
 &&-B \hot F - \frac 1 2 \tr \Xb \Xi \hot  F-\frac 1 2\Xh \c \ov{\DD} F +\frac12\Xh (\ov{\Hb}\c F)\\
 &&+s \Big[-\frac{1}{2}\tr X(Z-\Hb)+2\om(Z+\Hb)+ 2\DD\om\\
 &&+\frac{1}{2}\widehat{X}\c(-\ov{Z}+\ov{\Hb})-\frac{1}{2}\tr\Xb\Xi-2\omb\Xi -B \Big] \hot F \\
 &&+s Z \hot \nab_4F +2s \om (\DD\hot F+s Z \hot F )+ (\Ga_b \c \Ga_g) F.
 \eeaa
 We can rewrite the above as
  \beaa
 \nabc_4 (\DDc \hot F )  &=&\mathcal{D}\hot(\nab_4F)+s Z \hot \nab_4F+(\Hb+Z)\hot\nab_4 F+2s\om(Z+\Hb)\hot F\\
 &&+ \Xi \hot \nab_3 F-2s\omb\Xi \hot F -\frac 1 2 \tr X\left( \DD \hot  F+\Hb\hot F+ s (Z- \Hb) \hot F\right)\\
 &&+2s \om \DD\hot F+2s^2\om Z \hot F +2s \DD\om  \hot F\\
 &&-(s+1)B \hot F - (s+1)\frac 1 2 \tr \Xb \Xi \hot  F-\frac 1 2\Xh \c \ov{\DD} F-\frac 1 2 s (\widehat{X}\c\ov{Z} ) \hot F\\
 &&+\frac12\Xh (\ov{\Hb}\c F) +s  \frac{1}{2}(\widehat{X}\c\ov{\Hb})  \hot F + (\Ga_b \c \Ga_g) F,
 \eeaa
 which gives
   \beaa
 \nabc_4 (\DDc \hot F )  &=&\mathcal{D}\hot(\nab_4F)+s Z \hot \nab_4F+(\Hb+Z)\hot\nab_4 F+2s\om(Z+\Hb)\hot F\\
 &&+ \Xi \hot \nab_3 F-2s\omb\Xi \hot F -\frac 1 2 \tr X\left( \DD \hot  F+\Hb\hot F+ s (Z- \Hb) \hot F\right)\\
 &&+2s \om \DD\hot F+2s^2\om Z \hot F +2s \DD\om  \hot F\\
 &&-(s+1)B \hot F - (s+1)\frac 1 2 \tr \Xb \Xi \hot  F-\frac 1 2\Xh \c \ov{\DDc} F\\
 &&+\frac12\Xh (\ov{\Hb}\c F) +s  \frac{1}{2}(\widehat{X}\c\ov{\Hb})  \hot F  + (\Ga_b \c \Ga_g) F.
 \eeaa
 Recall that $\nabc_4F$ is conformal of type $s+1$, therefore $$\DDc \hot(\nabc_4F)= \mathcal{D}\hot(\nabc_4F)+(s+1) Z \hot \nabc_4F.$$ 
 We finally obtain 
  \beaa
[ \nabc_4, \DDc \hot] F   &=& -\frac 1 2 \tr X\left( \DDc \hot  F+(1-s)\Hb\hot F\right)+\Hb\hot\nabc_4 F+ \Xi \hot \nabc_3 F\\
 &&-(s+1)B \hot F - (s+1)\frac 1 2 \tr \Xb \Xi \hot  F-\frac 1 2\Xh \c \ov{\DDc} F\\
 && +\frac12\Xh (\ov{\Hb}\c F) +s  \frac{1}{2}(\widehat{X}\c\ov{\Hb})  \hot F+ (\Ga_b \c \Ga_g) F,
 \eeaa
 as stated. Similarly for $[\nabc_3, \DDc\hot] F$. This proves \eqref{eq:comm-nabc4nabc3DDchot-precise}.

For $U\in \sk_2(\mathbb{C})$ $s$-conformally invariant, we have
\beaa
\nabc_4 \ov{\DDc}\c U&=&\nabc_4 ( \ov{\DD}\c U + s \ov{Z} \c U)\\
&=& \nab_4 ( \ov{\DD}\c U + s \ov{Z} \c U)+2s\om  ( \ov{\DD}\c U + s \ov{Z} \c U)\\
&=& \nab_4 ( \ov{\DD}\c U) + s \nab_4(\ov{Z}) \c U+ s \ov{Z} \c \nab_4U+2s\om  ( \ov{\DD}\c U + s \ov{Z} \c U).
\eeaa
Using the null structure equation for $ \nab_4\ov{Z}$ and \eqref{eq:comm:nab4-nab3-DDc-precise}, 
we have
\beaa
\nabc_4 \ov{\DDc}\c U&=&\ov{\DD} \c \nab_4 U -\frac 1 2\ov{\tr X} \big(  \ov{\DD}\c U - 2\ov{\Hb} \c U\big) +\ov{(\Hb + Z)}\c\nab_4 U+ \ov{\Xi} \c \nab_3 U \\
&&+2\ov{B} \c U  -\frac 1 2 \ov{\tr \Xb} \ov{\Xi}\c  U  -\frac 1 2 \Xh \c \ov{\DD} U-\frac 1 2 (\ov{\Xh}\c U)\ov{\Hb}\\
&&+ s \Big[-\frac{1}{2}\ov{\tr X}(\ov{Z-\Hb})+2\om(\ov{Z+\Hb})+ 2\ov{\DD}\om\\
 &&+\frac{1}{2}\ov{\widehat{X}}\c(-Z+\Hb)-\frac{1}{2}\ov{\tr\Xb}\ov{\Xi}-2\omb\ov{\Xi} -\ov{B} \Big] \c U\\
 &&+ s \ov{Z} \c \nab_4U+2s\om  ( \ov{\DD}\c U + s \ov{Z} \c U)+ (\Ga_b \c \Ga_g) U.
\eeaa
We can rewrite the above as
\beaa
\nabc_4 \ov{\DDc}\c U&=&\ov{\DDc} \c \nabc_4 U -\frac 1 2\ov{\tr X} \big(  \ov{\DDc}\c U - (s+2)\ov{\Hb} \c U\big) +\ov{\Hb}\c\nabc_4 U\\
&&+ \ov{\Xi} \c \nabc_3 U -(s-2)\ov{B} \c U  -(s+1)\frac 1 2 \ov{\tr \Xb} \ov{\Xi}\c  U  -\frac 1 2 \Xh \c \ov{\DDc} U\\
&&-\frac 1 2 (\ov{\Xh}\c U)\ov{\Hb}+ s \frac{1}{2}(\ov{\widehat{X}}\c(\Hb)  )\c U+ (\Ga_b \c \Ga_g) U.
\eeaa
Similarly for $[\nabc_3, \ov{\DDc} \c]$. This proves \eqref{eq:comm-nabc4nabc3-ovDDc-U-precise}.

Using \eqref{commutator-nab4-ov-DDcF}, we write
 \beaa
 \nabc_4 (\ov{\DDc} \c F )&=&\nabc_4 (\ov{\DD}\c F+s \ov{Z}\c  F )\\
 &=&\nab_4 (\ov{\DD}\c F+s \ov{Z}\c  F )+2s \om (\ov{\DD}\c F+s \ov{Z}\c  F )\\
  &=&\nab_4 (\ov{\DD}\c F)+s \nab_4(\ov{Z}) \c F +s \ov{Z} \c \nab_4F +2s \om (\ov{\DD}\c F+s \ov{Z}\c  F )\\
    &=&\ov{\DD} \c (\nab_4F)- \frac 1 2\ov{\tr X}\, ( \ov{\DD} \c F -  \ov{\Hb} \c F)+\ov{(\Hb+Z)} \c \nab_4 F\\
 &&+s (-\frac{1}{2}\ov{\tr X}(\ov{Z-\Hb})+2\om\ov{(Z+\Hb)}+ 2\ov{\DD}\om ) \c F \\
 &&+s \ov{Z} \c \nab_4F +2s \om (\ov{\DD}\c F+s \ov{Z}\c  F ) + r^{-1} \Ga_g  \c \dk^{\leq 1} F.
 \eeaa
 We rewrite the above as
  \beaa
 \nabc_4 (\ov{\DDc} \c F )    &=&\ov{\DD} \c (\nabc_4F)+s \ov{Z} \c \nabc_4F+\ov{(\Hb+Z)} \c \nabc_4 F\\
 &&- \frac 1 2\ov{\tr X}\, ( \ov{\DD} \c F -  \ov{\Hb} \c F+ s(\ov{Z-\Hb}) \c F ) + r^{-1} \Ga_g  \c \dk^{\leq 1} F.
 \eeaa
Since $\ov{\DDc} \c (\nabc_4F)= \ov{\mathcal{D}}\c(\nabc_4F)+(s+1) \ov{Z} \c \nabc_4F$, 
we obtain the stated identity.

 We have 
 \beaa
 [\nabc_3, \nabc_4] U&=&\nabc_3 \nabc_4 U -\nabc_4 \nabc_3 U\\
 &=&(\nab_3-2(s+1) \omb ) (\nab_4 U+2s \om U) -(\nab_4+2(s-1)\om) (\nab_3 U -2s \omb U)\\
 &=&\nab_3\nab_4 U+2s\nab_3(\om U)-2(s+1) \omb (\nab_4 U)-4s(s+1) \om  \omb U \\
 &&-\nab_4 \nab_3 U +2s \nab_4(\omb U)-2(s-1)\om\nab_3 U+4s(s+1) \om  \omb U\\
 &=&[\nab_3, \nab_4] U-2 \omb \nab_4 U+2\om \nab_3 U  +2s(\nab_3 \om+ \nab_4\omb -4\om\omb)U.
 \eeaa
 Using \eqref{commutator-nab-4-D-c}   and
 \beaa
 \nab_3\om+\nab_4\omb -4\om\omb &=&   \rho  +(\eta-\etab)\c\ze -\eta\c\etab+\xi\c \xib, 
 \eeaa
 we obtain
  \beaa
 [\nabc_3, \nabc_4] U &=& - 2\om \nab_3 U+ 2\omb \nab_4 U  + 2 (\eta_c-\etab_c) \nab_c U +4i \left(- \rhod+ \eta \wedge \etab  \right) U\\
&&-2 \omb \nab_4 U+2\om \nab_3 U  +2s(\rho  +(\eta-\etab)\c\ze -\eta\c\etab)U + \err_{34}[U]\\
 &=&  2(\eta_c-\etab_c )  \nabc_c  U +2s(\rho   -\eta\c\etab)U +4i \left(- \rhod+ \eta \wedge \etab  \right) U+ \err_{34}[U]
 \eeaa
 as stated.

We have
\beaa
\nabc_3 \nabc_a U_{bc}&=&\nabc_3 (\nab_a U_{bc}+s \ze_a U_{bc})\\
&=&\nab_3 (\nab_a U_{bc}+s \ze_a U_{bc})-2s\omb (\nab_a U_{bc}+s \ze_a U_{bc})\\
&=&\nab_3 \nab_a U_{bc}+s \nab_3\ze_a U_{bc}+s \ze_a \nab_3U_{bc}-2s\omb (\nab_a U_{bc}+s \ze_a U_{bc}).
\eeaa
Using \eqref{commutator-nab-3-nab-a-U} and the null structure equation 
\beaa
\nab_3 \ze+2\nab\omb&=& -\frac{1}{2}\trchb(\ze+\eta)-\frac{1}{2}\atrchb(\dual\ze+\dual\eta)+ 2 \omb(\ze-\eta)\\
&&+\hch\c\xib+\frac{1}{2}\trch\,\xib+\frac{1}{2}\atrch\dual\xib +2 \om \xib-\chibh\c(\ze+\eta) -\bb
\eeaa
we obtain
\beaa
\nabc_3 \nabc_a U_{bc}&=&\nab_a \nab_3 U_{bc}-\frac  1 2   \trchb\, \Big(\nab_a U_{bc}+\eta_bU_{ac}+\eta_c U_{ab}-\de_{a b}(\eta \c U)_c-\de_{a c}(\eta \c U)_b \Big)\\
&&-\frac 1 2 \atrchb\, \Big(\dual \nab_a  U_{bc} +\eta_b \dual U_{ac}+\eta_c \dual U_{ab}- \in_{a b}(\eta \c  U)_c- \in_{a c}(\eta \c  U)_b \Big)\\
 &&+(\eta_a-\ze_a)\nab_3 U_{bc}\\
 &&+s (-2\nab_a\omb -\frac{1}{2}\trchb(\ze_a+\eta_a)-\frac{1}{2}\atrchb(\dual\ze_a+\dual\eta_a)+ 2 \omb(\ze_a-\eta_a)) U_{bc}\\
 &&+s \ze_a \nab_3U_{bc}-2s\omb (\nab_a U_{bc}+s \ze_a U_{bc})+\Ga_g \c \dk^{\leq 1} U.
\eeaa
We can rewrite the above as
\beaa
&&\nabc_3 \nabc_a U_{bc}\\
&=&\nab_a \nab_3 U_{bc}+s \ze_a \nab_3U_{bc}+(\eta_a-\ze_a)\nab_3 U_{bc}+ 2s \omb(\ze_a-\eta_a) U_{bc}\\
&&-\frac  1 2   \trchb\, \Big(\nab_a U_{bc}+s(\ze_a+\eta_a) U_{bc}+\eta_bU_{ac}+\eta_c U_{ab}-\de_{a b}(\eta \c U)_c-\de_{a c}(\eta \c U)_b \Big)\\
&&-\frac 1 2 \atrchb\, \Big(\dual \nab_a  U_{bc}+s (\dual\ze_a+\dual\eta_a) U_{bc} +\eta_b \dual U_{ac}+\eta_c \dual U_{ab}\\
&&- \in_{a b}(\eta \c  U)_c- \in_{a c}(\eta \c  U)_b \Big)\\
 &&-2s\omb \nab_a U_{bc}-2s^2\omb  \ze_a U_{bc}-2s\nab_a\omb  U_{bc}+\Ga_g \c \dk^{\leq 1} U,
\eeaa
which gives
\beaa
&&\nabc_3 \nabc_a U_{bc}\\
&=&\nab_a \nabc_3 U_{bc}+s \ze_a \nabc_3U_{bc}+(\eta_a-\ze_a)\nabc_3 U_{bc}\\
&&-\frac  1 2   \trchb\, \Big(\nabc_a U_{bc}+s(\eta_a) U_{bc}+\eta_bU_{ac}+\eta_c U_{ab}-\de_{a b}(\eta \c U)_c-\de_{a c}(\eta \c U)_b \Big)\\
&&-\frac 1 2 \atrchb\, \Big(\dual \nabc_a  U_{bc}+s (\dual\eta_a) U_{bc} +\eta_b \dual U_{ac}+\eta_c \dual U_{ab}\\
&&- \in_{a b}(\eta \c  U)_c- \in_{a c}(\eta \c  U)_b \Big)+\Ga_g \c \dk^{\leq 1} U.
\eeaa
Recall that $\nabc_3 U$ is of conformal type $s-1$, therefore
\beaa
\nabc_a \nabc_3 U_{bc}&=& \nab_a \nabc_3 U_{bc}+(s-1) \ze_a \nabc_3 U_{bc}.
\eeaa
We finally obtain
\beaa
&&\nabc_3 \nabc_a U_{bc}\\
&=&\nabc_a \nabc_3 U_{bc}+\eta_a\nabc_3 U_{bc}\\
&&-\frac  1 2   \trchb\, \Big(\nabc_a U_{bc}+s(\eta_a) U_{bc}+\eta_bU_{ac}+\eta_c U_{ab}-\de_{a b}(\eta \c U)_c-\de_{a c}(\eta \c U)_b \Big)\\
&&-\frac 1 2 \atrchb\, \Big(\dual \nabc_a  U_{bc}+s (\dual\eta_a) U_{bc} +\eta_b \dual U_{ac}+\eta_c \dual U_{ab}\\
&&- \in_{a b}(\eta \c  U)_c- \in_{a c}(\eta \c  U)_b \Big)+\Ga_g \c \dk^{\leq 1} U,
\eeaa
as stated.


\section{Proof of Lemma \ref{LEMMA:DEFORMATION-TENSORS-T}}\label{appendix-proof-deformation-tensor}


For $\piT$, one can easily adapt the proof of Proposition 2.6.10 in \cite{KS:Kerr}, which is done in the particular case of outgoing geodesic frame. Also,
 \beaa
 ( \Div{}^{(\T)}\pi)_a &=& -\frac{1}{2}\D_3{}^{(\T)}\pi_{4a} -\frac{1}{2}\D_4{}^{(\T)}\pi_{3a}+\g^{bc}\D_b{}^{(\T)}\pi_{ca}=\dk^{\leq 1}\Ga_g\\
   (\Div \piT)_4&=& -\frac{1}{2}\D_3{}^{(\T)}\pi_{44} -\frac{1}{2}\D_4{}^{(\T)}\pi_{34}+\g^{ab}\D_a{}^{(\T)}\pi_{b4}\\
 &=&  -\frac{1}{2}\D_3{}^{(\T)}\pi_{44}+\g^{ab}\D_a{}^{(\T)}\pi_{b4}= \dk^{\leq 1}\Ga_g,\\
 (\Div \piT)_3&=& -\frac{1}{2}\D_3{}^{(\T)}\pi_{43} -\frac{1}{2}\D_4{}^{(\T)}\pi_{33}+\g^{ab}\D_a{}^{(\T)}\pi_{b3}\\
 &=& -\frac{1}{2}\D_4{}^{(\T)}\pi_{33}+\g^{ab}\D_a{}^{(\T)}\pi_{b3}=r^{-1} \dk^{\leq 1}\Ga_b.
 \eeaa

For $\piZ$, since the only difficulty is to track the weights in $r$, we provide the proof in the case where the  normalization of $(e_3, e_4)$ is outgoing, in which case we have
\beaa
\Z &=& \frac 1 2 \left(2(r^2+a^2)\Re(\Jk)^be_b -a(\sin\th)^2 e_3 -\frac{a(\sin\th)^2\De}{ |q|^2} e_4\right).
\eeaa
We write
\beaa
\Z=\Z^{(1)}+\Z^{(2)}, \qquad \Z^{(1)}:=(r^2+a^2)\Re(\Jk)^be_b, \qquad \Z^{(2)}:=\frac 1 2 \left( -a(\sin\th)^2 e_3 -\frac{a(\sin\th)^2\De}{ |q|^2} e_4\right).
\eeaa
Using $\widecheck{\nab\Jk}\in r^{-1}\Ga_b$, we compute
\beaa
\widecheck{{}^{(\Z^{(1)})}\pi_{44}} &=& 0,\\
\widecheck{{}^{(\Z^{(1)})}\pi_{43}} &=& r\widecheck{e_3(r)}\Re(\Jk)^b\g(e_b, e_4)+r^2\widecheck{\nab_3\Re(\Jk)^b}=r^2\widecheck{\nab_3\Re(\Jk)^b}=r\Ga_b,\\
\widecheck{{}^{(\Z^{(1)})}\pi_{4a}} &=& r\widecheck{\nab(r)}\Re(\Jk)^b\g(e_b, e_4)=r\Ga_g,\\
\widecheck{{}^{(\Z^{(1)})}\pi_{ab}} &=& r\widecheck{\nab(r)}\Re(\Jk)^b\g(e_b, e_a)+r^2\widecheck{\nab\Re(\Jk)^b}=r\Ga_b,\\
\widecheck{{}^{(\Z^{(1)})}\pi_{3a}} &=& r\widecheck{e_3(r)}\Re(\Jk)^b\g(e_b, e_a) +r\widecheck{\nab(r)}\Re(\Jk)^b\g(e_b, e_3)+r^2\widecheck{\nab_3\Re(\Jk)^b}+r^2\widecheck{\nab\Re(\Jk)^b}=r\Ga_b,\\
\widecheck{{}^{(\Z^{(1)})}\pi_{33}} &=& r\widecheck{e_3(r)}\Re(\Jk)^b\g(e_b, e_3) +r\widecheck{\nab_3\Re(\Jk)^b}=r^2\widecheck{\nab_3\Re(\Jk)^b}=r\Ga_b.
\eeaa 
Since $\Z^{(2)}$ is at the same level as $\T$, then $\widecheck{{}^{(\Z^{(1)})}\pi}$ has the same decay properties as $\piT$. We deduce
\beaa
{}^{(\Z)}\pi_{44}\in \Ga_g, \qquad {}^{(\Z)}\pi_{4a}\in r\Ga_g, \qquad {}^{(\Z)}\pi_{ab},\,\, {}^{(\Z)}\pi_{43}, \,\,{}^{(\Z)}\pi_{33}, \,\,   {}^{(\Z)}\pi_{3a}\in r\Ga_b,
\eeaa
as stated.  
 In particular, we have
   \beaa
\tr({}^{(\Z)}\pi) &=& -{}^{(\Z)}\pi_{43}+\g^{ab}{}^{(\Z)}\pi_{ab}\in r\Ga_b.
 \eeaa
 Also we compute
 \beaa
(  \Div{}^{(\Z)}\pi)_a &=& -\frac{1}{2}\D_3({}^{(\Z)}\pi)_{4a} -\frac{1}{2}\D_4({}^{(\Z)}\pi)_{3a}+\g^{bc}\D_b({}^{(\Z)}\pi)_{ca}= r\dk^{\leq 1}\Ga_g, \\
 (\Div \piZ)_4 &=&  -\frac{1}{2}\D_3({}^{(\Z)}\pi)_{44}+\g^{ab}\D_a({}^{(\Z)}\pi)_{b4}= \dk^{\leq 1}(\Ga_g),\\
 (\Div \piZ)_3 &=& -\frac{1}{2}\D_4({}^{(\Z)}\pi)_{33}+\g^{ab}\D_a({}^{(\Z)}\pi)_{b3}= \dk^{\leq 1}(\Ga_b),
 \eeaa
 as stated.


\section{Proof of Proposition \ref{LE:COMMTZSQUARE}}\label{proof:comm-T-Z-square}


Recall the following  general commutation  formula 
  for a scalar $\psi$ and vectorfield $X$ in a Ricci flat manifold, as implied by Lemma \ref{LEMMA:COMMUTATIONNAB_XSQUARED}:
\bea
\label{eq:generalcommutation-scalarwave}
 \,[ X, \square_\g] \psi&=&-2\piX^{\mu\a} \D_{\mu}\D_\a \psi+\big( \D^{\a}(\tr \piX )-2  (\Div \piX)^\a\big) \D_\a \psi.
\eea
From \eqref{eq:generalcommutation-scalarwave}, we infer
\beaa
 \,[ X, \square_\g] \psi &=& r^{-2}\Big(\piX_{33}, \, \piX_{3a},\, \piX_{ab}\Big)\dk^{\leq 2}\psi + r^{-1}\piX_{4a}\dk^{\leq 2}\psi+\piX_{44}\dk^{\leq 2}\psi\\
 &&+\piX_{34}\D_3\D_4\psi +\Big(\D_4(\tr \piX )-2  (\Div \piX)_4\Big)\dk^{\leq 1}\psi\\
 &&+ r^{-1}\Big(\D_3(\tr \piX )-2  (\Div \piX)_3,\,\, \D_a(\tr \piX )-2(\Div \piX)_a\Big)\dk^{\leq 1}\psi.
\eeaa
Writing
\beaa
\D_3\D_4\psi &=& -\square_\g\psi+\g^{ab}\D_a\D_b\psi,
\eeaa
we obtain
\bea\lab{cor:commutationbetweenvectofieldXandscalarwave}
\begin{split}
 \,[ X, \square_\g] \psi &= r^{-2}\Big(\piX_{33}, \, \piX_{3a},\, \piX_{ab}, \, \piX_{34}\Big)\dk^{\leq 2}\psi + r^{-1}\piX_{4a}\dk^{\leq 2}\psi+\piX_{44}\dk^{\leq 2}\psi\\
 &-\piX_{34}\square_\g\psi +\Big(\D_4(\tr \piX )-2  (\Div \piX)_4\Big)\dk^{\leq 1}\psi\\
 &+ r^{-1}\Big(\D_3(\tr \piX )-2  (\Div \piX)_3,\,\, \D_a(\tr \piX )-2(\Div \piX)_a\Big)\dk^{\leq 1}\psi.
 \end{split}
\eea

Using Lemma \ref{LEMMA:DEFORMATION-TENSORS-T} and \eqref{cor:commutationbetweenvectofieldXandscalarwave}, we obtain for $X=\T$:
\beaa
 \,[ \T, \square_\g] \psi &=& r^{-2}\Ga_b \c \dk^{\leq 2}\psi + r^{-1}\Ga_g \c \dk^{\leq 2}\psi+\Ga_g \c \dk^{\leq 2}\psi\\
 &&+\Ga_b \c \square_\g\psi +\Big( \dk^{\leq 1} \Ga_g \Big)\dk^{\leq 1}\psi+ r^{-1}\Big(\dk^{\leq 1} \Ga_g+r^{-1} \dk^{\leq 1} \Ga_b\Big)\dk^{\leq 1}\psi\\
 &=& \Ga_g \c \dk^{\leq 2}\psi+ \dk^{\leq 1} \Ga_g \c \dk^{\leq 1}\psi+\Ga_b \c \square_\g\psi \\
 &=&\dk \big(\Ga_g \c \dk \psi\big)+\Ga_b \c \square_\g\psi,
\eeaa
as stated.
 
 Similarly, using Lemma \ref{LEMMA:DEFORMATION-TENSORS-T} and applying \eqref{cor:commutationbetweenvectofieldXandscalarwave} to $X=\Z$, we obtain
\beaa
 \,[ \Z, \square_\g] \psi &=& r^{-1}  \Ga_b  \c \dk^{\leq 2}\psi + \Ga_g \c \dk^{\leq 2}\psi+\Ga_g \c \dk^{\leq 2}\psi\\
 &&+ r \Ga_b \c \square_\g\psi +\Big(\D_4(\tr \piZ )+ \dk^{\leq1} \Ga_g \Big)\dk^{\leq 1}\psi\\
 &&+ r^{-1}\Big(\D_3(\tr \piZ )+ \dk^{\leq 1} \Ga_b ,\,\, \Ga_b + r \dk^{\leq 1} \Ga_g\Big)\dk^{\leq 1}\psi\\
 &=&  \Ga_g \c \dk^{\leq 2}\psi+ r \Ga_b \c \square_\g\psi \\
 &&+\Big(\D_4(\tr \piZ )+ r^{-1} \D_3(\tr \piZ )+ \dk^{\leq1} \Ga_g \Big)\dk^{\leq 1}\psi.
\eeaa

We now show that 
\bea\label{eq:D4-D3-trpiZ}
\D_4(\tr \piZ )+ r^{-1} \D_3(\tr \piZ ) &=&  \dk^{\leq1} \Ga_g.
\eea
We have 
 \beaa
 \D_4(\tr \piZ )&=& \D_4\Big( -{}^{(\Z)}\pi_{43}+\g^{ab}{}^{(\Z)}\pi_{ab}\Big)= \D_4\Big(\g^{ab}{}^{(\Z)}\pi_{ab}\Big) \\
  &=& \nab_4(r^2\widecheck{\div(\Re(\Jk))})= \dk^{\leq 1}\Ga_g,
 \eeaa
 where we used the good transport equation in $e_4$ for $\widecheck{\div(\Re(\Jk))}$. Also,
 \beaa
 \D_3(\tr \piZ ) &=& \D_3\Big( -{}^{(\Z)}\pi_{43}+\g^{ab}{}^{(\Z)}\pi_{ab}\Big) = \D_3\Big(\g^{ab}{}^{(\Z)}\pi_{ab}\Big)\\
 &=& \nab_3(r^2\widecheck{\nab\Jk})= \dk^{\leq 1}\Ga_b,
 \eeaa
 where we use Lemma \ref{lemma:improved-nab3-r^2nabJ} below. Since $r^{-1} \Ga_b$ decays faster than $\Ga_g$, this proves \eqref{eq:D4-D3-trpiZ}, and therefore the Proposition.

\begin{lemma}\label{lemma:improved-nab3-r^2nabJ}
We have  the following improved estimate:
\beaa
\nab_3(r^2\widecheck{\nab\Jk})\in \dk^{\leq 1}\Ga_b. 
\eeaa 
\end{lemma}
\begin{proof}
Note first from Proposition 5.6.16 in  \cite{KS:Kerr} that 
\beaa
\nab_\nu(r^2\widecheck{\nab\Jk})\in \dk^{\leq 1}\Ga_b\textrm{ on }\Si_*.
\eeaa
Also, we have in view of Lemma 6.1.15 in  \cite{KS:Kerr}: 
\beaa
\nab_4 (\DD\hot\Jk)+\frac{2}{q}\DD\hot\Jk &=& O(r^{-1})B +O(r^{-2})\trXc+O(r^{-2})\Xh\\
&&+O(r^{-2})\Zc +O(r^{-3})\widecheck{\DD(\cos\th)}+r^{-1}\Ga_b\c\Ga_g,\\
\nab_4\big(\widecheck{\ov{\DD}\c\Jk} \big)+\Re\left(\frac{2}{q}\right)\widecheck{\ov{\DD}\c\Jk} &=& O(r^{-1})B+O(r^{-2})\trXc+O(r^{-2})\Xh +O(r^{-2})\Zc\\
&&+O(r^{-3})\widecheck{\DD(\cos\th)}+r^{-1}\Ga_b\c\Ga_g.
\eeaa
In particular, we infer
\beaa
\nab_4\Big(q^2\DD\hot\Jk\Big) &=& \Big(rB,\, \trXc,\, \Xh,\, \Zc\Big) +r^{-1}\Ga_b,\\
\nab_4\Big(|q|^2\widecheck{\ov{\DD}\c\Jk}\Big) &=& \Big(rB,\, \trXc,\, \Xh,\, \Zc\Big) +r^{-1}\Ga_b.
\eeaa
Commuting with $\nab_3$, we infer
\beaa
\nab_4\nab_3\Big(q^2\DD\hot\Jk\Big) &=& r^{-1}\dk^{\leq 1}\Ga_b,\\
\nab_4\nab_3\Big(|q|^2\widecheck{\ov{\DD}\c\Jk}\Big) &=& r^{-1}\dk^{\leq 1}\Ga_b,
\eeaa
which immediately implies by integration from $\Si_*$ the estimate $\nab_3(r^2\widecheck{\nab\Jk})\in \dk^{\leq 1}\Ga_b$ as desired. 
\end{proof}


\section{Proof of Proposition \ref{LEMMA:MOD-LAPLACIAN-PERT-KERR}}
\label{proof-lemma:commutator-OO-|q|2square-scalar}


Using Lemma \ref{lemma:expression-wave-operator}, we have 
\beaa
|q|^2\squared_2 \psi&=&-\frac 1 2|q|^2 \big(\nab_3\nab_4\psi+\nab_4 \nab_3 \psi\big)+|q|^2\left(\omb -\frac 1 2 \trchb\right) \nab_4\psi+|q|^2\left(\om -\frac 1 2 \trch\right) \nab_3\psi\\
&& +|q|^2( \triangle \psi +(\eta+\etab)\c \nab \psi) \\
&=&-\frac 1 2|q|^2 \big(\nab_3\nab_4\psi+\nab_4 \nab_3 \psi\big)+|q|^2\left(\omb -\frac 1 2 \trchb\right) \nab_4\psi+|q|^2\left(\om -\frac 1 2 \trch\right) \nab_3\psi\\
&& +\OO(\psi)+ r \Ga_b \c \dk \psi ,
\eeaa
and therefore
\beaa
[\OO, |q|^2 \squared_2]\psi &=& -\frac 1 2[\OO, |q|^2 \big(\nab_3\nab_4+\nab_4 \nab_3 \big)]\psi\\
&&+\big[\OO, |q|^2\left(\omb -\frac 1 2 \trchb\right) \nab_4\big]\psi+\big[\OO, |q|^2\left(\om -\frac 1 2 \trch\right) \nab_3\big]\psi+ r \dk^2\big(  \Ga_b \c \dk \psi \big).
\eeaa

 We compute $[\OO, |q|^2 \big(\nab_3\nab_4+\nab_4 \nab_3 \big)]\psi$. We write,
\beaa
[\OO, |q|^2 \nab_3\nab_4]\psi&=& |q|^2[\OO, \nab_3\nab_4]\psi+2|q|^2\nab (|q|^2) \c \nab \nab_3 \nab_4 \psi+\OO(|q|^2) \nab_3 \nab_4\psi\\
&=& |q|^2[\OO, \nab_3]\nab_4\psi+|q|^2\nab_3([\OO, \nab_4]\psi)\\
&&+2|q|^2\nab (|q|^2) \c \nab \nab_3 \nab_4 \psi+\OO(|q|^2) \nab_3 \nab_4\psi.
\eeaa
Using the explicit expressions for the commutators $[\OO, \nab_3]$ and $[\OO, \nab_4]$ in Lemma \ref{LEMMA:COMMUTATOR-NAB3-NAB4-LAP}, we obtain
\beaa
[\OO, \nab_3]\nab_4\psi&=& -(\eta-\ze) \c  |q|^2\nab_3 \nab \nab_4\psi-(\eta-\ze) \c  |q|^2\nab\nab_3 \nab_4\psi \\
&&-\big( \div(\eta-\ze)+(\eta+\etab) \c (\eta - \ze) \big) |q|^2 \nab_3 \nab_4\psi \\
&&+O(ar^{-3}) \dk^{\leq 2} \psi+r^{-1}\dk \big( \Ga_b \c \dk^2\psi),
\eeaa
and 
\beaa
\nab_3([ \OO, \nab_4 ]\psi)&=&- (\etab+\ze) \c  |q|^2 \nab_3\nab_4 \nab \psi-(\etab+\ze) \c  |q|^2\nab_3 \nab\nab_4 \psi \\
&&- \nab_3 (|q|^2(\etab+\ze) ) \c \nab_4 \nab \psi-\nab_3 (|q|^2(\etab+\ze)) \c \nab\nab_4 \psi \\
&&-\big( \div(\etab+\ze)+(\eta+\etab) \c (\etab + \ze) \big) |q|^2 \nab_3\nab_4 \psi \\
&&- \nab_3 \Big(\big( \div(\etab+\ze)+(\eta+\etab) \c (\etab + \ze) \big) |q|^2\Big) \nab_4 \psi \\
&&+O(ar^{-2}) \dk^{\leq 2} \psi +\Ddot_3 \Ddot_3 \big(|q|^2 \xi \c \Ddot_a \psi \big)+ \dk^2 \big( \Ga_g \c \dk \psi).
\eeaa
By considering for each term in the above its value in Kerr plus the error terms we obtain
\beaa
|q|^{-2} [\OO, |q|^2 \big(\nab_3\nab_4+\nab_4 \nab_3 \big)]\psi&=&|q|^{-2} [\OO, |q|^2 \big(\nab_3\nab_4+\nab_4 \nab_3 \big)]_{Kerr}\psi\\
&&+\Ddot_3 \dk \big(|q|^2 \xi \c \Ddot_a \psi \big)+ \dk^2 \big( \Ga_g \c \dk \psi),
\eeaa

We compute $\Big[\OO, |q|^2\left(\omb -\frac 1 2 \trchb\right) \nab_4\Big]\psi+\Big[\OO, |q|^2\left(\om -\frac 1 2 \trch\right) \nab_3\Big]\psi$.

We write 
\beaa
&&\Big[\OO, |q|^2\left(\omb -\frac 1 2 \trchb\right) \nab_4\Big]\psi\\
&=&|q|^2\left(\omb -\frac 1 2 \trchb\right)  [\OO, \nab_4]\psi\\
&&+ 2 |q|^2 \nab \big(|q|^2\big(\omb -\frac 1 2 \trchb\big)  \big) \c \nab \nab_4 \psi + \OO\Big(|q|^2\big(\omb -\frac 1 2 \trchb\big) \Big) \nab_4 \psi,
\eeaa
which gives
\beaa
&&|q|^{-2} \Big[\OO, |q|^2\left(\omb -\frac 1 2 \trchb\right) \nab_4\Big]\psi\\
&=&\big(\omb -\frac 1 2 \trchb\big) \Big(- (\etab+\ze) \c  |q|^2\nab_4 \nab \psi-(\etab+\ze) \c  |q|^2\nab\nab_4 \psi \Big)\\
&&+ 2 \nab \left(|q|^2\big(\omb -\frac 1 2 \trchb\big)  \right) \c \nab \nab_4 \psi+O(ar^{-3}) \dk \psi+ \dk \big( \Ga_g \c \dk \psi).
\eeaa
By symmetry we obtain
\beaa
&&|q|^{-2} \Big[\OO, |q|^2\left(\om -\frac 1 2 \trch\right) \nab_3\Big]\psi\\
&=&\big(\om -\frac 1 2 \trch\big) \Big(- (\eta-\ze) \c  |q|^2\nab_3 \nab \psi-(\eta-\ze) \c  |q|^2\nab\nab_3 \psi \Big)\\
&&+ 2 \nab \left(|q|^2\big(\om -\frac 1 2 \trch\big)  \right) \c \nab \nab_3 \psi+O(ar^{-3}) \dk \psi+r^{-1} \dk \big( \Ga_b \c \dk \psi).
\eeaa
We therefore deduce 
\beaa
&&|q|^{-2} \Big([\OO, |q|^2\big(\omb -\frac 1 2 \trchb\big) \nab_4]\psi+[\OO, |q|^2\big(\om -\frac 1 2 \trch\big) \nab_3]\psi\Big)\\
&=& |q|^{-2} \Big([\OO, |q|^2\big(\omb -\frac 1 2 \trchb\big) \nab_4]\psi+[\OO, |q|^2\big(\om -\frac 1 2 \trch\big) \nab_3]\psi\Big)_{Kerr}+ \dk \big( \Ga_g \c \dk \psi).
\eeaa

 We sum the above two contributions, and obtain
 \beaa
|q|^{-2} [\OO, |q|^2 \square_\g]\psi &=&|q|^{-2} [\OO, |q|^2 \square_\g]_{Kerr}\psi +\Ddot_3 \dk \big(|q|^2 \xi \c \Ddot_a \psi \big)+ \dk^2 \big( \Ga_g \c \dk \psi).
\eeaa
Using Proposition \ref{LEMMA:MOD-LAPLACIAN-KERR} to write 
\beaa
[\OO, |q|^2 \square_\g]_{Kerr}\psi&=&|q|^2 \Big[\nab\left(\frac{8a(r^2+a^2)\cos\th}{|q|^2}\right)\c\nab\nab_\That\dual\psi +O(ar^{-2})\nab^{\leq 1}_{\Rhat}\dk^{\leq 1}\psi \Big],
\eeaa
 we prove the Proposition.


\section{Proof of Lemma \ref{LEMMA:SYMM-OPERATORS}}\label{appendix-proof-symmetry-operators}


We write, for two vectorfields $X$ and $Y$:
\beaa
&&|q|^2 \Ddot_\a(|q|^{-2}X^{(\a} Y^{\b)} \Ddot_\b \psi)\\
&=&\nab_{(X} \nab_{Y)} \psi +\frac 12 \Big[ |q|^2 \nab_X(|q|^{-2}) \nab_Y \psi+ |q|^2 \nab_Y(|q|^{-2}) \nab_X \psi\\
&&+ \tr \piX \nab_Y \psi+ \tr {}^{(Y)}\pi \nab_X \psi+ {}^{(Y)}\pi_{\a\b} X^{\a} \Ddot^\b \psi+\piX_{\a\b} Y^{\a} \Ddot^\b \psi\Big].
\eeaa
For $X=Y=\T$, we obtain
\beaa
\SS_1(\psi)&=&\nab_\T \nab_\T \psi + |q|^2 \nab_\T(|q|^{-2}) \nab_\T \psi+ \tr \piT \nab_\T \psi+\piT_{\a\b} \T^{\a} \Ddot^\b \psi.
\eeaa
Writing $\nab_\T=\dk$ and using that $\nab_\T (|q|) \in r \Ga_b$, $\tr \piT \in \Ga_g$, and 
\beaa
\piT_{\a\b} \T^{\a} \Ddot^\b \psi&=& \piT_{3\b} \T^{3} \Ddot^\b \psi+\piT_{4\b} \T^{4} \Ddot^\b \psi+\piT_{a\b} \T^{a} \Ddot^\b \psi\\
&=& \piT_{34} \T^{3} \Ddot_3 \psi+\Ga_g \c \dk \psi=\Ga_b \c \dk \psi,
\eeaa
we have
\beaa
\SS_1(\psi)&=&\nab_\T \nab_\T \psi +\Ga_b \c \dk \psi.
\eeaa
For $X=Y=\Z$, we obtain
\beaa
\SS_3(\psi)&=&a^2\nab_\Z \nab_\Z \psi + |q|^2 \nab_\Z(|q|^{-2}) \nab_\Z \psi+ \tr \piZ \nab_\Z \psi+\piZ_{\a\b} \Z^{\a} \Ddot^\b \psi.
\eeaa
Writing $\nab_\Z=\dk$ and using that $\nab_\Z (|q|) \in r \Ga_b$, $\tr \piZ \in r\Ga_b$, and 
\beaa
\piZ_{\a\b} \Z^{\a} \Ddot^\b \psi&=& \piZ_{3\b} \Z^{3} \Ddot^\b \psi+\piZ_{4\b} \Z^{4} \Ddot^\b \psi+\piZ_{a\b} \Z^{a} \Ddot^\b \psi\\
&=& \piZ_{34} \Z^{3} \Ddot_3 \psi+\Ga_b \c \dk \psi=r\Ga_b \c \dk \psi,
\eeaa
we have
\beaa
\SS_3(\psi)=a^2\nab_\Z \nab_\Z \psi +r \Ga_b \c \dk \psi.
\eeaa
Combining the above, for $X=\T$ and $Y=\Z$ we have
\beaa
\SS_2(\psi)&=&a\nab_\T \nab_\Z \psi +r \Ga_b \c \dk \psi+\piT_{\a\b} \Z^{\a} \Ddot^\b \psi+\piZ_{\a\b} \T^{\a} \Ddot^\b \psi= a\nab_\T \nab_\Z \psi +r \Ga_b \c \dk \psi.
\eeaa
Finally, we write
\beaa
\SS_4(\psi)&=& |q|^2\Ddot_\b(|q|^{-2}O^{\a\b}\Ddot_\a \psi)=|q|^2\Ddot_\b( \ga^{ab} e_a^\a e_b^\b\Ddot_\a \psi)\\
&=& |q|^2 \ga^{ab} e_a^\a e_b^\b \Ddot_\a \Ddot_\b \psi+  \ga^{ab}|q|^2\Ddot_\b( e_a^\a) e_b^\b\Ddot_\a \psi+ \ga^{ab}|q|^2e_a^\a \Ddot_\b( e_b^\b)\Ddot_\a \psi \\
&=& |q|^2 \ga^{ab} \Ddot_a \Ddot_b \psi+  \ga^{ab}|q|^2\Ddot_b( e_a^\a)\Ddot_\a \psi+ \ga^{ab}|q|^2 \Ddot_\b( e_b^\b)\Ddot_a \psi \\
&=& |q|^2 \ga^{ab} \nab_a \nab_b \psi - \frac 1 2 |q|^2\trch \nab_3 \psi -\frac 1 2|q|^2 \trchb\nab_4 \psi
\\
&&+  \ga^{ab}|q|^2\big( \frac 12 \chi_{ba} e_3^\a +\frac 1 2 \chib_{ba} e_4^\a) \Ddot_\a \psi+ \ga^{ab}|q|^2\big((\Ddot_3 e_b)^3+ (\Ddot_4 e_b)^4\big)\Ddot_a \psi\\
&=&  |q|^2 \ga^{ab} \nab_a \nab_b \psi + \ga^{ab}|q|^2(\eta_b+ \etab_b)\nab_a \psi= \OO(\psi)+r\Ga_b\c\dk, 
\eeaa
as stated.


\section{Proof of Lemma \ref{LEMMA:SQUARED-K-THAT-RHAT-KERR}}\label{proof:alternative-repr-square}


Using the expression for the metric given by \eqref{inverse-metric-vfs-perturbations}, we have 
\beaa
|q|^2 \squared_k \psi&=& |q|^2 \g^{\a\b} \Ddot_\a \Ddot_\b \psi =\left(\frac{(r^2+a^2)^2}{\De} \big( -\That^\a\That^\b +      \Rhat^\a \Rhat^\b\big) + O^{\a\b} \right) \Ddot_\a \Ddot_\b \psi \\
&=&\frac{(r^2+a^2)^2}{\De} \big( - \Ddot_{\That} \Ddot_{\That} \psi +    \Ddot_{\Rhat} \Ddot_{\Rhat} \psi  \big) + O^{\a\b} \Ddot_\a \Ddot_\b \psi .
\eeaa
We write
\beaa
 \Ddot_{\That} \Ddot_{\That} \psi&=& \nab_\That \nab_\That \psi - \Ddot_{\D_\That \That} \psi, \qquad  \Ddot_{\Rhat} \Ddot_{\Rhat} \psi= \nab_\Rhat \nab_\Rhat \psi - \Ddot_{\D_\Rhat \Rhat} \psi.
\eeaa
Using the definitions of $\That$ and $\Rhat$ in the outgoing frame, we obtain
\beaa
2\D_\That \That&=& \frac{\De}{r^2+a^2} \D_4 \That +\frac{|q|^2}{r^2+a^2} \D_3 \That\\
&=&\frac 1 2 \frac{\De}{r^2+a^2} \D_4 \left( \frac{\De}{r^2+a^2} e_4+\frac{|q|^2}{r^2+a^2}  e_3\right) + \frac 1 2 \frac{|q|^2}{r^2+a^2} \D_3  \left( \frac{\De}{r^2+a^2} e_4+\frac{|q|^2}{r^2+a^2}  e_3\right)\\
&=&\frac 1 2 \frac{\De}{r^2+a^2} \left( \frac{\De}{r^2+a^2}  \D_4 e_4+\frac{|q|^2}{r^2+a^2}   \D_4 e_3\right)\\
&&+ \frac 1 2 \frac{|q|^2}{r^2+a^2}  \left( \frac{\De}{r^2+a^2} \D_3 e_4+\frac{|q|^2}{r^2+a^2}  \D_3 e_3\right)\\
&& +\frac 1 2 \frac{\De}{r^2+a^2}  \left( e_4\big(\frac{\De}{r^2+a^2}\big) e_4+e_4\big(\frac{|q|^2}{r^2+a^2}\big)  e_3\right) \\
&&+ \frac 1 2 \frac{|q|^2}{r^2+a^2}  \left( e_3\big(\frac{\De}{r^2+a^2}\big) e_4+e_3\big(\frac{|q|^2}{r^2+a^2} \big) e_3\right). 
\eeaa
Similarly, 
\beaa
2\D_\Rhat \Rhat&=&\frac 1 2 \frac{\De}{r^2+a^2} \left( \frac{\De}{r^2+a^2}  \D_4 e_4-\frac{|q|^2}{r^2+a^2}   \D_4 e_3\right) \\
&&- \frac 1 2 \frac{|q|^2}{r^2+a^2}  \left( \frac{\De}{r^2+a^2} \D_3 e_4-\frac{|q|^2}{r^2+a^2}  \D_3 e_3\right)\\
&&+\frac 1 2 \frac{\De}{r^2+a^2}  \left( e_4\big(\frac{\De}{r^2+a^2}\big) e_4-e_4\big(\frac{|q|^2}{r^2+a^2}\big)  e_3\right) \\
&&- \frac 1 2 \frac{|q|^2}{r^2+a^2}  \left( e_3\big(\frac{\De}{r^2+a^2}\big) e_4-e_3\big(\frac{|q|^2}{r^2+a^2} \big) e_3\right). 
\eeaa
This gives
\beaa
2 \big( \D_\That \That- \D_\Rhat \Rhat\big)&=& \frac{\De |q|^2}{(r^2+a^2)^2} \left( \D_4 e_3+ \D_3e_4\right)\\
&& + \frac{\De}{r^2+a^2}  e_4\big(\frac{|q|^2}{r^2+a^2}\big)  e_3+ \frac{|q|^2}{r^2+a^2}  e_3\big(\frac{\De}{r^2+a^2}\big) e_4.
\eeaa
Using \eqref{eq:expressions-Riccif-formula}, we obtain
\beaa
2 \big( \D_\That \That- \D_\Rhat \Rhat\big)&=&  \Big( \frac{ |q|^2}{r^2+a^2}  2\om +  e_4\big(\frac{|q|^2}{r^2+a^2}\big) \Big)  \frac{\De }{r^2+a^2} e_3\\
&& +\Big( \frac{\De }{r^2+a^2} 2\omb +   e_3\big(\frac{\De}{r^2+a^2}\big)\Big)  \frac{ |q|^2}{r^2+a^2} e_4\\
&&+ \frac{\De |q|^2}{(r^2+a^2)^2} 2 (\eta_b+\etab_b) e_b.
\eeaa
Computing in the outgoing frame
\bea\label{eq:computations-e3e4-out}
\begin{split}
e_4\left(\frac{|q|^2}{r^2+a^2} \right)&= \frac{2a^2r\sin^2\th}{(r^2+a^2)^2} + r^{-3} \Ga_g \\
e_3\left(\frac{\De}{r^2+a^2} \right)&= -\frac{\De}{|q|^2} \frac{\pr_r (\De)}{r^2+a^2} + \frac{2r\De^2}{|q|^2(r^2+a^2)^2} + \Ga_b,
\end{split}
\eea
we obtain
\beaa
\frac{ |q|^2}{r^2+a^2}  2\om +  e_4\big(\frac{|q|^2}{r^2+a^2}\big)&=&  \frac{2ra^2\sin^2\th}{(r^2+a^2)^2} + \Ga_g\\
\frac{\De }{r^2+a^2} 2\omb +   e_3\big(\frac{\De}{r^2+a^2}\big)&=& \frac{\De }{r^2+a^2} \pr_r\left(\frac{\De}{|q|^2} \right)  -\frac{\De}{|q|^2} \frac{\pr_r (\De)}{r^2+a^2} + \frac{2r\De^2}{|q|^2(r^2+a^2)^2}+\Ga_b\\
&=&- \frac{2r\De^2 }{|q|^4(r^2+a^2)} + \frac{2r\De^2}{|q|^2(r^2+a^2)^2}+\Ga_b\\
&=&- \frac{2ra^2\sin^2\th \De^2 }{|q|^4(r^2+a^2)^2}+\Ga_b.
\eeaa
We write
\beaa
2 \big( \D_\That \That- \D_\Rhat \Rhat\big)&=&    \frac{2ra^2\sin^2\th \De}{|q|^2(r^2+a^2)^2} \big( \frac{|q|^2 }{r^2+a^2} e_3 -  \frac{\De }{r^2+a^2} e_4\big) \\
&&+ \frac{\De |q|^2}{(r^2+a^2)^2} 2 (\eta_b+\etab_b) e_b + \Ga_g \c \dk \\
&=& -   \frac{4ra^2\sin^2\th \De}{|q|^2(r^2+a^2)^2} \Rhat+ \frac{\De |q|^2}{(r^2+a^2)^2} 2 (\eta_b+\etab_b) e_b + \Ga_g \c \dk .
\eeaa
We therefore obtain
\beaa
|q|^2 \squared_k \psi&=&\frac{(r^2+a^2)^2}{\De} \big( -  \nab_\That \nab_\That \psi+   \nab_\Rhat \nab_\Rhat \psi + \Ddot_{\D_\That \That-\D_\Rhat \Rhat} \psi \big) + O^{\a\b} \Ddot_\a \Ddot_\b \psi \\
&=&\frac{(r^2+a^2)^2}{\De} \big( -  \nab_\That \nab_\That \psi+   \nab_\Rhat \nab_\Rhat \psi \big) -   \frac{2ra^2\sin^2\th }{|q|^2} \nab_\Rhat \psi\\
&&+ O^{\a\b} \Ddot_\a \Ddot_\b \psi   + |q|^2  (\eta_b+\etab_b) \nab_b\psi  + r^2 \Ga_g \c \dk \psi.
\eeaa
Finally, using the computations in Lemma \ref{lemma:expression-wave-operator}, i.e.
\beaa
\Ddot_c \Ddot_d \psi&=&\nab_c\nab_d \psi - \frac 1 2 \chi_{cd}\nab_3 \psi -\frac 1 2 \chib_{cd}\nab_4 \psi,
\eeaa
we write
\beaa
O^{\a\b} \Ddot_\a \Ddot_\b \psi &=& |q|^2 \ga^{ab} e_a^\a e_b^\b \Ddot_\a \Ddot_\b \psi \\
&=& |q|^2 \ga^{ab} \nab_a \nab_b \psi - \frac 1 2 |q|^2\trch \nab_3 \psi -\frac 1 2|q|^2 \trchb\nab_4 \psi\\
&=& |q|^2 \lap_k \psi - r \nab_3 \psi + \frac{r\De}{|q|^2}\nab_4 \psi + r^2\Ga_g \c \dk \psi \\
&=& |q|^2 \lap_k \psi +\frac{2r(r^2+a^2)}{|q|^2} \nab_\Rhat \psi + r^2\Ga_g \c \dk \psi.
\eeaa
This finally gives
\beaa
|q|^2 \squared_k \psi&=&\frac{(r^2+a^2)^2}{\De} \big( -  \nab_\That \nab_\That \psi+   \nab_\Rhat \nab_\Rhat \psi \big) -   \frac{2ra^2\sin^2\th }{|q|^2} \nab_\Rhat \psi\\
&&+  |q|^2 \lap_k \psi +\frac{2r(r^2+a^2)}{|q|^2} \nab_\Rhat \psi   + |q|^2  (\eta_b+\etab_b) \nab_b\psi  + r^2 \Ga_g \c \dk \psi\\
&=&\frac{(r^2+a^2)^2}{\De} \big( -  \nab_\That \nab_\That \psi+   \nab_\Rhat \nab_\Rhat \psi \big) +2r \nab_\Rhat \psi\\
&&+  |q|^2 \lap_k \psi   + |q|^2  (\eta_b+\etab_b) \nab_b\psi  + r^2 \Ga_g \c \dk \psi,
\eeaa
as stated.


\section{Proof of Lemma \ref{COMPLEXWAVE-DECOMP}}\label{sec:proof-complex-wave}


We define $Y_{ab}:=\big(\DD\hot( \DDb \c \psi)\big)_{ab}$, and express $Y$ in frames. 
We have
\beaa
Y_{ab}&=& \DD_a   \DDb^c \psi_{cb} +\DD_b    \DDb^c \psi_{c a } -\de_{ab} \DD^d    \DDb^c \psi_{dc}.
\eeaa
By construction, $Y$ is symmetric and traceless. 
We calculate first $Y_{11}=-Y_{22}$. For $a=b=1$ we derive
\beaa
Y_{11}&=& \DD_1   \DDb^c \psi_{c1} +\DD_1    \DDb^c \psi_{c 1 } -\de_{11} \DD^d    \DDb^c \psi_{dc}\\
&=& 2\DD_1   \DDb_c \psi_{c1}  - \DD^d    \DDb^c \psi_{dc}\\
&=& 2\DD_1   \DDb_1 \psi_{11}+2\DD_1   \DDb_2 \psi_{21}  - \DD^d    \DDb^1 \psi_{d1}- \DD^d    \DDb^2 \psi_{d2}\\
&=& 2\DD_1   \DDb_1 \psi_{11}+2\DD_1   \DDb_2 \psi_{21}  - \DD_1    \DDb_1 \psi_{11}- \DD_2    \DDb_1 \psi_{21}- \DD_1    \DDb_2 \psi_{12}- \DD_2    \DDb_2 \psi_{22}\\
&=& \DD_1   \DDb_1 \psi_{11}+\DD_1   \DDb_2 \psi_{21} - \DD_2    \DDb_1 \psi_{21}- \DD_2    \DDb_2 \psi_{22}.
\eeaa
Writing $\psi_{22}=-\psi_{11}$, we have
\beaa
Y_{11}&=& (\DD_1   \DDb_1 + \DD_2    \DDb_2) \psi_{11}+(\DD_1   \DDb_2  - \DD_2    \DDb_1) \psi_{12}.
\eeaa
We compute
\beaa
Y_{11}&=& (\big(\nab_1+i\dual\nab_1\big)   \big(\nab_1-i\dual\nab_1\big) + \big(\nab_2+i\dual\nab_2\big)   \big(\nab_2-i\dual\nab_2\big)) \psi_{11}\\
&&+(\big(\nab_1+i\dual\nab_1\big)   \big(\nab_2-i\dual\nab_2\big)  - \big(\nab_2+i\dual\nab_2\big)   \big(\nab_1-i\dual\nab_1\big)) \psi_{12}\\
&=& (\big(\nab_1+i\nab_2\big)   \big(\nab_1-i\nab_2\big) + \big(\nab_2-i\nab_1\big)   \big(\nab_2+i\nab_1\big)) \psi_{11}\\
&&+(\big(\nab_1+i\nab_2\big)   \big(\nab_2+i\nab_1\big)  - \big(\nab_2-i\nab_1\big)   \big(\nab_1-i\nab_2\big)) \psi_{12}\\
&=& 2\lap  \Psi_{11} -2 i (\nab_1\nab_2-\nab_2\nab_1) \Psi_{11}+2 \big(\nab_1\nab_2-\nab_2\nab_1\big) \psi_{12} + 2 i \lap  \psi_{12}.
\eeaa
Using that $\psi_{12}=-i \psi_{11}$, we obtain
\beaa
Y_{11}&=& 2\lap  \psi_{11} -2 i (\nab_1\nab_2-\nab_2\nab_1) \psi_{11}-2i \big(\nab_1\nab_2-\nab_2\nab_1\big)  \psi_{11} + 2  \lap \psi_{11}\\
&=& 4\lap  \psi_{11} -4 i (\nab_1\nab_2-\nab_2\nab_1) \psi_{11}.
\eeaa
Similarly we compute $Y_{12}=Y_{21}$.
We therefore obtain
\beaa
Y_{ab}&=& 4\lap  \psi_{ab} -4 i (\nab_1\nab_2-\nab_2\nab_1) \psi_{ab}.
\eeaa
Using the Gauss formulas \eqref{Gauss-eq-first-com}, we prove the desired formula.


\section{Proof of Lemma \ref{LEMMA:CONFORMAL-DD-LAP}}\label{sec:proof-lemma-conformal-DD}


We have for a $s$-conformally invariant horizontal tensor $\psi$,
\beaa
\DDc\hot( \DDbc \c \psi)&=&\DD\hot( \DDbc \c \psi)+ s Z \hot ( \DDbc \c \psi)\\
&=&\DD\hot( \DDb \c \psi + s \ov{Z} \c \psi )+ s Z \hot ( \DDb \c \psi+ s \ov{Z} \c \psi )\\
&=&\DD\hot( \DDb \c \psi) + s  \DD\hot (\ov{Z} \c \psi )+ s Z \hot (\DDb \c \psi)+  s^2 Z \hot ( \ov{Z} \c \psi ).
\eeaa
Using \eqref{simil-Leibniz}, \eqref{Leibniz-hot}, and \eqref{Leib-eq-DDb-DD-nab}, we have
\beaa
\DDc\hot( \DDbc \c \psi)&=&\DD\hot( \DDb \c \psi) +2 s  (\DD\c\ov{Z} ) \psi+2s(\ov{Z} \c \DD) \psi +2 s (Z \c \DDb )\psi+  2s^2 (Z \c \ov{Z} ) \psi \\
&=&\DD\hot( \DDb \c \psi)+8s \ze \c \nab \psi  +   4s(\div \ze - i \curl \ze) \psi+  4s^2|\ze|^2  \psi,
\eeaa
where we wrote $Z \c \ov{Z}=2 |\ze|^2$ and $\DD\c\ov{Z} =2(\div \ze - i \curl \ze)$. On the other hand,
\beaa
^{(c)}\lap_2\psi&=& \nabc^a \nabc_a \psi= \nab^a \nabc_a \psi +s\ze^a \nabc_a \psi= \nab^a( \nab_a \psi +s\ze_a \psi )+s\ze^a (\nab_a  \psi +s \ze_a \psi)\\
&=& \nab^a \nab_a \psi +s \nab^a \ze_a \psi +s\ze_a \nab^a\psi +s\ze^a \nab_a  \psi +s^2 |\ze|^2 \psi\\
&=&\lap_2 \psi +2s\ze \c \nab\psi  +s \nab^a \ze_a \psi +s^2 |\ze|^2 \psi.
\eeaa
Using \eqref{eq:expression-DD-laplacian-real}, 
we obtain
\beaa
&&\DDc\hot( \DDbc \c \psi)\\
&=&\DD\hot( \DDb \c \psi)+8s \ze \c \nab \psi  +   4s(\div \ze - i \curl \ze) \psi+  4s^2|\ze|^2  \psi \\
&=&4\lap_2  \psi -2 i (\atrch\nab_3+\atrchb \nab_4) \psi  +  2 \left( \trch\trchb+ \atrch\atrchb+4\rho\right) \psi \\
&&+8s \ze \c \nab \psi  +   4s(\div \ze - i \curl \ze) \psi+  4s^2|\ze|^2  \Psi +(\Ga_g \c \Ga_b) \c \psi \\
&=&4\lapc_2  \psi - 2i (\atrch\nab_3+\atrchb \nab_4) \psi  +2   \left( \trch\trchb+ \atrch\atrchb+4\rho\right) \psi \\
&&  -   4s  i \curl \ze \psi+(\Ga_g \c \Ga_b) \c \psi.
\eeaa
Writing $\nab_3=\nabc_3 + 2 s \omb$, $\nab_4=\nabc_4-2s\om$, we have
\beaa
&&\DDc\hot( \DDbc \c \psi)\\
&=&4\lapc_2  \psi -2 i \left(\atrch(\nabc_3\psi + 2 s \omb \psi)+\atrchb (\nabc_4\psi -2s\om\psi)\right)  \\
&&+ 2  \left( \trch\trchb+ \atrch\atrchb+4\rho\right) \psi   -   4s  i \curl \ze \psi+(\Ga_g \c \Ga_b) \c \psi\\
&=&4\lapc_2  \psi -2 i (\atrch \nabc_3+\atrchb \nabc_4)\psi  + 2  \left( \trch\trchb+ \atrch\atrchb+4\rho\right) \psi  \\
&&- 4si \big(  \curl \ze+ \omb  \atrch -\om \atrchb \big)\psi+(\Ga_g \c \Ga_b) \c \psi.
\eeaa
Using the null structure equation
\beaa
\curl\ze&=&-\frac 1 2 \chih\wedge\chibh   +\frac 1 4 \big(  \trch\atrchb-\trchb\atrch   \big)+\om \atrchb -\omb\atrch+\dual \rho,
\eeaa
we finally obtain the stated relation.


\section{Proof of Lemma \ref{LEMMA:COMMUTATOR-NAB3-NAB4-LAP}}
\label{proof-lemma-comm-nab3-nab4-lap}


Writing that $\lap=\ga^{ab} \nab_a \nab_b$, and using Lemma \ref{lemma:comm}, i.e. for a scalar $\psi$, 
       \beaa
        \,[\nab_3, \nab_a] \psi &=&-\frac 1 2 \left(\trchb \nab_a \psi+\atrchb \dual \nab_a \psi\right)+(\eta_a-\ze_a) \nab_3 \psi-\chibh_{ab}\nab_b \psi  +\xib_a \nab_4 \psi,
        \eeaa
we obtain
\beaa
[\nab_3, \lap]\psi&=& \ga^{ab} [\nab_3, \nab_a]\nab_b \psi+\ga^{ab} \nab_a [\nab_3, \nab_b]\psi\\
&=& -\frac 1 2 \left(\trchb \ga^{ab}\nab_a \nab_b \psi+\atrchb \ga^{ab}\dual \nab_a \nab_b \psi\right)+\ga^{ab}(\eta_a-\ze_a) \nab_3 \nab_b \psi\\
&&+\ga^{ab} \nab_a \big(-\frac 1 2 \left(\trchb \nab_b \psi+\atrchb \dual \nab_b \psi\right)+(\eta_b-\ze_b) \nab_3 \psi \big)\\
&&-\ga^{ab}\chibh_{ac}\nab_c \nab_b \psi +\ga^{ab}\xib_a \nab_4\nab_b \psi +\ga^{ab} \nab_a \big(-\chibh_{bc}\nab_c \psi  +\xib_b \nab_4 \psi \big)\\
&=& -\frac 1 2 \left(\trchb \ga^{ab}\nab_a \nab_b \psi+\atrchb \ga^{ab}\dual \nab_a \nab_b \psi\right)+\ga^{ab}(\eta_a-\ze_a) \nab_3 \nab_b \psi\\
&&-\frac 1 2 \left(\trchb \ga^{ab} \nab_a\nab_b \psi+\atrchb \ga^{ab} \nab_a\dual \nab_b \psi\right)+ \ga^{ab}(\eta_b-\ze_b) \nab_a\nab_3 \psi \\
&&-\frac 1 2 \left(\ga^{ab} \nab_a\trchb \nab_b \psi+\ga^{ab} \nab_a\atrchb \dual \nab_b \psi\right)+\ga^{ab} \nab_a(\eta_b-\ze_b) \nab_3 \psi +\err_3,
\eeaa
where $\err_3:=-\ga^{ab}\chibh_{ac}\nab_c \nab_b \psi +\ga^{ab}\xib_a \nab_4\nab_b \psi +\ga^{ab} \nab_a \big(-\chibh_{bc}\nab_c \psi  +\xib_b \nab_4 \psi \big)$. Writing $\dual \nab_a= \in_{ac} \nab_c$, we obtain 
\beaa
[\nab_3, \lap]\psi&=& -\frac 1 2 \left(\trchb \lap \psi+\atrchb \in_{ba} \nab_a \nab_b \psi\right)+\ga^{ab}(\eta_a-\ze_a) \nab_3 \nab_b \psi\\
&&-\frac 1 2 \left(\trchb \lap \psi+\atrchb \in_{ab} \nab_a \nab_b \psi\right)+ \ga^{ab}(\eta_b-\ze_b) \nab_a\nab_3 \psi \\
&&-\frac 1 2 \left( \nab \trchb  \c \nab \psi+\nab \atrchb \c  \dual \nab \psi\right)+\div(\eta-\ze) \nab_3 \psi +\err_3\\
&=& -\trchb \lap \psi +(\eta-\ze) \c  \nab_3 \nab \psi+(\eta-\ze) \c  \nab\nab_3 \psi +\div(\eta-\ze) \nab_3 \psi \\
&&-\frac 1 2 \left( \nab \trchb  \c \nab \psi+\nab \atrchb \c  \dual \nab \psi\right)+\err_3.
\eeaa
Similarly, we have
\beaa
[\nab_4, \lap]\psi&=& -\trch \lap \psi +(\etab+\ze) \c \nab_4 \nab\psi+(\etab+\ze) \c \nab \nab_4 \psi +\div(\etab+\ze) \nab_4 \psi \\
&&-\frac 1 2 \left( \nab \trch  \c \nab \psi+\nab \atrch \c  \dual \nab \psi\right)+\err_4,
\eeaa
where $\err_4=-\ga^{ab}\chih_{ac}\nab_c \nab_b \psi +\ga^{ab}\xi_a \nab_3\nab_b \psi +\ga^{ab} \nab_a \big(-\chih_{bc}\nab_c \psi  +\xi_b \nab_3 \psi \big)$. In particular, $\err_3=r^{-2} \dk \big( \Ga_b \c \dk \psi)$ and $\err_4=\Ddot_3 \big( \xi \c \Ddot_a \psi \big)+ r^{-2} \dk \big( \Ga_g \c \dk \psi)$, as stated.

Writing $\nab_3(|q|^2)=\trchb |q|^2 + \big( r \widecheck{e_3(r)} +\Ga_b\big)$, we have
\beaa
[\nab_3, |q|^2 \lap]\psi&=& |q|^2 [\nab_3, \lap]\psi+\nab_3(|q|^2) \lap\psi\\
&=&(\eta-\ze) \c  |q|^2\nab_3 \nab \psi+(\eta-\ze) \c  |q|^2\nab\nab_3 \psi +\div(\eta-\ze)|q|^2 \nab_3 \psi \\
&&-\frac 1 2|q|^2 \left( \nab \trchb  \c \nab \psi+\nab \atrchb \c  \dual \nab \psi\right)\\
&&+r^2\err_3 +\big( r \widecheck{e_3(r)} +\Ga_b\big)  \lap\psi.
\eeaa
The error terms are given by
\beaa
r^2\err_3 +\big( r \widecheck{e_3(r)} +\Ga_b\big)  \lap\psi&=& \dk \big( \Ga_b \c \dk \psi)+r^2\Ga_b \c r^{-2} \dk^2\psi= \dk \big( \Ga_b \c \dk \psi).
\eeaa
Similarly we obtain
\beaa
[\nab_4, |q|^2 \lap]\psi&=&(\etab+\ze) \c |q|^2\nab_4 \nab\psi+(\etab+\ze) \c |q|^2\nab \nab_4 \psi +\div(\etab+\ze) |q|^2\nab_4 \psi \\
&&-\frac 1 2|q|^2 \left( \nab \trch  \c \nab \psi+\nab \atrch \c  \dual \nab \psi\right)\\
&&+r^2\Ddot_3 \big( \xi \c \Ddot_a \psi \big)+ \dk \big( \Ga_g \c \dk \psi),
\eeaa
as stated.

Using the above, we deduce for $\OO= |q|^2 \lap \psi - 2a^2\cos\th \Im(\Jk)^b \nab_b \psi   $, 
\beaa
[\nab_3, \OO]\psi&=& [\nab_3, |q|^2 \lap]\psi - [\nab_3, 2a^2\cos\th \Im(\Jk)^b \nab_b]\psi\\
&=& (\eta-\ze) \c  |q|^2\nab_3 \nab \psi+(\eta-\ze) \c  |q|^2\nab\nab_3 \psi +\div(\eta-\ze)|q|^2 \nab_3 \psi \\
&&-\frac 1 2|q|^2 \left( \nab \trchb  \c \nab \psi+\nab \atrchb \c  \dual \nab \psi\right)\\
&&- \nab_3\big(2a^2\cos\th \Im(\Jk) \big) \c \nab\psi-2a^2\cos\th \Im(\Jk)\c  [\nab_3, \nab]\psi+\dk \big( \Ga_b \c \dk \psi)\\
&=& (\eta-\ze) \c  |q|^2\nab_3 \nab \psi+(\eta-\ze) \c  |q|^2\nab\nab_3 \psi +\div(\eta-\ze)|q|^2 \nab_3 \psi \\
&&-\frac 1 2|q|^2 \left( \nab \trchb  \c \nab \psi+\nab \atrchb \c  \dual \nab \psi\right)\\
&&- \nab_3\big(2a^2\cos\th \Im(\Jk) \big) \c \nab\psi\\
&&- 2a^2\cos\th \Im(\Jk)\c \left(-\frac 1 2 \left(\trchb \nab \psi+\atrchb \dual \nab \psi\right)+(\eta-\ze) \nab_3 \psi \right)+\dk \big( \Ga_b \c \dk \psi),
\eeaa
which gives
\beaa
[\nab_3, \OO]\psi&=& (\eta-\ze) \c  |q|^2\nab_3 \nab \psi+(\eta-\ze) \c  |q|^2\nab\nab_3 \psi \\
&&+\big( \div(\eta-\ze)+(\eta+\etab) \c (\eta - \ze) \big) |q|^2 \nab_3 \psi +O(ar^{-2}) \dk^{\leq 1} \psi +\dk \big( \Ga_b \c \dk \psi).
\eeaa
Similarly, 
\beaa
[\nab_4, \OO]\psi&=& (\etab+\ze) \c  |q|^2\nab_4 \nab \psi+(\etab+\ze) \c  |q|^2\nab\nab_4 \psi \\
&&+\big( \div(\etab+\ze)+(\eta+\etab) \c (\etab + \ze) \big) |q|^2 \nab_4 \psi+O(ar^{-2}) \dk^{\leq 1} \psi \\
&&+r^2\Ddot_3 \big( \xi \c \Ddot_a \psi \big)+ \dk \big( \Ga_g \c \dk \psi),
\eeaa
as stated.

In the case of a 2-tensor, we have to modify the above using the formulas \eqref{commutator-3-a-u-bc}       \beaa
        \,[\nab_3, \nab_a] \psi &=&-\frac 1 2 \left(\trchb \nab_a \psi+\atrchb \dual \nab_a \psi\right)+(\eta_a-\ze_a) \nab_3 \psi-\chibh_{ab}\nab_b \psi  +\xib_a \nab_4 \psi\\
        &&+ O(ar^{-3})\psi+ r^{-1} \Ga_b \c \psi.
        \eeaa
        Following the same steps as above, we obtain for $\psi \in \sk_2$
\beaa
[\nab_3, \lap_2]\psi&=& -\trchb \lap \psi +(\eta-\ze) \c  \nab_3 \nab \psi+(\eta-\ze) \c  \nab\nab_3 \psi +\div(\eta-\ze) \nab_3 \psi \\
&&-\frac 1 2 \left( \nab \trchb  \c \nab \psi+\nab \atrchb \c  \dual \nab \psi\right)+O(ar^{-4}) \dk^{\leq 1} \psi +r^{-2} \dk \big( \Ga_b \c \dk \psi),
\eeaa
and therefore
\beaa
[\nab_3, |q|^2 \lap_2]\psi&=&(\eta-\ze) \c  |q|^2\nab_3 \nab \psi+(\eta-\ze) \c  |q|^2\nab\nab_3 \psi +\div(\eta-\ze)|q|^2 \nab_3 \psi \\
&&-\frac 1 2|q|^2 \left( \nab \trchb  \c \nab \psi+\nab \atrchb \c  \dual \nab \psi\right)+O(ar^{-4}) \dk^{\leq 1} \psi+\dk \big( \Ga_b \c \dk \psi),
\eeaa
and
\beaa
[\nab_3, \OO]\psi&=& (\eta-\ze) \c  |q|^2\nab_3 \nab \psi+(\eta-\ze) \c  |q|^2\nab\nab_3 \psi \\
&&+\big( \div(\eta-\ze)+(\eta+\etab) \c (\eta - \ze) \big) |q|^2 \nab_3 \psi +O(ar^{-2}) \dk^{\leq 1} \psi +\dk \big( \Ga_b \c \dk \psi).
\eeaa
Similarly for $\nab_4$.


\section{Proof of Lemma \ref{LEMMA:COMMUTATOR-NAB3-NAB4-SQUARE}}
\label{proof:lemma-comm-nab3-nab4-square}


 Using Lemma \ref{lemma:expression-wave-operator} and the commutator in Lemma \ref{lemma:comm}, we deduce for a scalar $\psi$,
\beaa
\square_\g \psi&=&-\nab_3\nab_4\psi-\frac 1 2 \trch \nab_3\psi+\left(2\omb -\frac 1 2 \trchb\right) \nab_4\psi+\lap \psi + 2\eta \c\nab \psi.
\eeaa
We therefore compute
\beaa
[\nab_3, \square_\g]\psi&=& -\nab_3 [\nab_3, \nab_4]\psi-\frac 1 2 \nab_3 \trch \nab_3 \psi \\
&&+\nab_3\left(2\omb -\frac 1 2 \trchb\right) \nab_4\psi+\left(2\omb -\frac 1 2 \trchb\right)[\nab_3, \nab_4]\psi\\
&&+[\nab_3,\lap] \psi + 2\nab_3\eta \c\nab \psi+ 2\eta \c[ \nab_3,\nab] \psi.
\eeaa
Using Lemma \ref{lemma:comm} and Lemma \ref{LEMMA:COMMUTATOR-NAB3-NAB4-LAP}, we obtain
\beaa
&&[\nab_3, \square_\g]\psi\\
&=& -\nab_3 \big(2(\eta-\etab ) \c \nab \psi - 2 \om \nab_3 \psi +2\omb \nab_4 \psi \big)-\frac 1 2 \nab_3 \trch \nab_3 \psi \\
&&+\left(2\nab_3\omb -\frac 1 2\nab_3 \trchb\right) \nab_4\psi+\left(2\omb -\frac 1 2 \trchb\right)\big(2(\eta-\etab ) \c \nab \psi - 2 \om \nab_3 \psi +2\omb \nab_4 \psi \big)\\
&& -\trchb \lap \psi +2(\eta-\ze) \c \nab \nab_3\psi +\big(\div(\eta-\ze)+ |\eta-\ze|^2 \big) \nab_3 \psi \\
&&-\frac 1 2 \big(\nab \trchb  + \trchb (\eta-\ze) \big)  \c \nab \psi-\frac 1 2  \big( \nab \atrchb+\atrchb (\eta-\ze) \big) \c  \dual \nab \psi \\
&&+ 2\nab_3\eta \c\nab \psi+ 2\eta \c \left(-\frac 1 2 \left(\trchb \nab \psi+\atrchb \dual \nab \psi\right)+(\eta-\ze) \nab_3 \psi \right)+\err_3,
\eeaa
which gives
\beaa
&&[\nab_3, \square_\g]\psi\\
&=&2 \om \nab_3\nab_3 \psi -2\omb \nab_3\nab_4 \psi -\trchb \lap \psi \\
&&+ \left( 2 \nab_3\om -\frac 1 2 \nab_3 \trch -2\om \left(2\omb -\frac 1 2 \trchb\right)+\div(\eta-\ze)+ |\eta-\ze|^2+ 2 \eta \c (\eta-\ze)  \right) \nab_3 \psi \\
&&+\left(-\frac 1 2\nab_3 \trchb+2\omb (2\omb -\frac 1 2 \trchb) \right)\nab_4\psi \\
&& - 2(\eta-\etab ) \c \nab_3\nab \psi+2(\eta-\ze) \c \nab \nab_3\psi  \\
&&+\left( 2\nab_3\etab  +2\left(2\omb -\frac 1 2 \trchb\right)(\eta-\etab )-\frac 1 2 \big(\nab \trchb  + \trchb (3\eta-\ze) \big)\right) \c \nab \psi \\
&&-\frac 1 2  \big( \nab \atrchb+\atrchb (3\eta-\ze) \big) \c  \dual \nab \psi +\err_3.
\eeaa
Using again the expression for $\square_\g$, writing
\beaa
\lap \psi &=&\square_\g \psi+\nab_3\nab_4\psi+\frac 1 2 \trch \nab_3\psi-\left(2\omb -\frac 1 2 \trchb\right) \nab_4\psi-2\eta \c\nab \psi,
\eeaa
we finally obtain
\beaa
&&[\nab_3, \square_\g]\psi\\
&=&2 \om \nab_3\nab_3 \psi -(\trchb +2\omb) \nab_3\nab_4 \psi -\trchb\square_\g \psi \\
&&+ \Big[ 2 \nab_3\om -\frac 1 2 \nab_3 \trch-\frac 1 2 \trch\trchb -2\om \left(2\omb -\frac 1 2 \trchb\right)+\div(\eta-\ze)+ |\eta-\ze|^2\\
&&+ 2 \eta \c (\eta-\ze)  \Big] \nab_3 \psi +\left(-\frac 1 2\nab_3 \trchb+2\omb \left(2\omb -\frac 1 2 \trchb\right) +\trchb (2\omb -\frac 1 2 \trchb)\right)\nab_4\psi \\
&& - 2(\eta-\etab ) \c \nab_3\nab \psi+2(\eta-\ze) \c \nab \nab_3\psi  \\
&&+\left( 2\nab_3\etab  +2\left(2\omb -\frac 1 2 \trchb\right)(\eta-\etab )-\frac 1 2 \big(\nab \trchb  + \trchb (3\eta-\ze) \big)+2\trchb \eta \right) \c \nab \psi \\
&&-\frac 1 2  \big( \nab \atrchb+\atrchb (3\eta-\ze) \big) \c  \dual \nab \psi +\err_3,
\eeaa
and similarly for $\nab_4$. Schematically, the commutator can be written as
\beaa
[\nab_3, \square_\g]\psi&=&2 \om \nab_3\nab_3 \psi -(\trchb +2\omb) \nab_3\nab_4 \psi -\trchb\square_\g \psi \\
&&+r^{-2}  \dk \psi +O(ar^{-2})  \nab_3\nab \psi   +r^{-2} \dk \big( \Ga_b \c \dk \psi),
\eeaa
and 
\beaa
[\nab_4, \square_\g]\psi&=&2 \omb \nab_4\nab_4 \psi -(\trch +2\om) \nab_4\nab_3 \psi -\trch\square_\g \psi \\
&&+r^{-2}  \dk \psi +O(ar^{-2})  \nab_4\nab \psi   +r^{-1} \dk \big( \Ga_g \c \dk \psi),
\eeaa
as stated.

For $\psi \in \sk_2$, using \eqref{first-equation-square}, we compute
\beaa
[\nab_4, \squared_2]\psi&=& -\nab_4 [\nab_4, \nab_3]\psi-\frac 1 2 \nab_4 \trchb \nab_4 \psi \\
&&+\nab_4\left(2\om -\frac 1 2 \trch\right) \nab_3\psi+\left(2\om -\frac 1 2 \trch\right)[\nab_4, \nab_3]\psi\\
&&+[\nab_4,\lap_2] \psi + 2\nab_4\etab \c\nab \psi+ 2\etab \c[ \nab_4,\nab] \psi+ 2i \nab_4\left( \rhod- \eta \wedge \etab \right) \psi.
\eeaa
Writing that
\beaa
 \, [\nab_4, \nab_3] \psi &=  2\om \nab_3 \psi  - 2\omb \nab_4 \psi  + 2 (\etab_c-\eta_c) \nab_c \psi +O(ar^{-4}) \psi ,
\eeaa
and using Lemma \ref{LEMMA:COMMUTATOR-NAB3-NAB4-LAP}, to write
\beaa
[\nab_4, \lap_2]\psi&=& -\trch \lap_2 \psi +2(\etab+\ze) \c  \nab_4 \nab \psi+O(ar^{-4}) \dk^{\leq 1} \psi \\
&&+\Ddot_3 \big( \xi \c \Ddot_a \psi \big)+ r^{-2}\dk \big( \Ga_g \c \dk \psi),
\eeaa
we obtain
\beaa
[\nab_4, \squared_2]\psi&=& -\nab_4 \big(2\om \nab_3 \psi  - 2\omb \nab_4 \psi  + 2 (\etab_c-\eta_c) \nab_c \psi  \big)+O(r^{-3}) \dk^{\leq 1} \psi  \\
&& -\frac 1 2 \nab_4\trch \nab_3\psi-\trch \lap_2 \psi +2(\etab+\ze) \c  \nab_4 \nab \psi \\
&&+\Ddot_3 \big( \xi \c \Ddot_a \psi \big)+ r^{-2}\dk \big( \Ga_g \c \dk \psi),
\eeaa
which gives
\beaa
[\nab_4, \squared_2]\psi&=&2\omb \nab_4 \nab_4 \psi  -2\om \nab_4 \nab_3 \psi  -\trch \lap_2 \psi + 2 (\eta+\zeta) \c \nab_4 \nab \psi   \\
&& -\frac 1 2 \nab_4\trch \nab_3\psi +O(r^{-3}) \dk^{\leq 1} \psi  +\Ddot_3 \big( \xi \c \Ddot_a \psi \big)+ r^{-2}\dk \big( \Ga_g \c \dk \psi) .
\eeaa
Using again the expression for $\squared_2$, writing
\beaa
\lap_2 \psi&=&\squared_2 \psi+\nab_4 \nab_3 \psi   +\frac 1 2 \trch \nab_3\psi+O(r^{-2}) \dk^{\leq 1}\psi,
\eeaa
we finally obtain
\beaa
[\nab_4, \squared_2]\psi&=&2\omb \nab_4 \nab_4 \psi -(\trch +2\om) \nab_4 \nab_3 \psi + 2 (\eta+\zeta) \c \nab_4 \nab \psi  -\trch \squared_2 \psi  \\
&&-\frac 1 2\big( \nab_4\trch+( \trch)^2  \big) \nab_3\psi +O(r^{-3}) \dk^{\leq 1} \psi  +\Ddot_3 \big( \xi \c \Ddot_a \psi \big)+ r^{-2}\dk \big( \Ga_g \c \dk \psi).
\eeaa
Using the null structure equation for $\nab_4 \trch$ we obtain the stated.

Using the above commutator,
we compute
\beaa
[r \nab_4, \squared_2]\psi &=& r[\nab_4, \squared_2]\psi +(\nab_3 r) \nab_4 \nab_4\psi + (\nab_4 r) \nab_3\nab_4 \psi -(\square_\g r) \nab_4 \psi\\
&=&-r(\trch +2\om) \nab_4 \nab_3 \psi - \nab_4 \nab_4 \psi -r\trch \squared_2 \psi -\frac 1 4r( \trch)^2  \nab_3\psi \\
&& + \frac 1 2 \trch r \nab_3\nab_4 \psi +O(r^{-2}) \dk^{\leq 1} \psi+O(r^{-3}) \dk^{\leq 2} \psi\\
&&+r\Ddot_3 \big( \xi \c \Ddot_a \psi \big)+ r^{-1}\dk \big( \Ga_b \c \dk \psi),
\eeaa
where we used $\nab_4 r=\frac 1 2 \trch r +O(r^{-2})+\Ga_g$, and $e_3(r)=-\frac{\De}{|q|^2}+r\Ga_b$ (and thus $e_3(r)=-1$ at main order). We therefore obtain
\beaa
[r \nab_4, \squared_2]\psi &=&-r\left(\frac 1 2 \trch +2\om\right) \nab_4 \nab_3 \psi-\nab_4 \nab_4\psi  -\frac 1 4r( \trch)^2  \nab_3\psi -r\trch \squared_2 \psi \\
&&+O(r^{-2}) \dk^{\leq 1} \psi+O(r^{-3}) \dk^{\leq 2} \psi+r\Ddot_3 \big( \xi \c \Ddot_a \psi \big)+ r^{-1}\dk \big( \Ga_g \c \dk \psi).
\eeaa
Using once again that
\beaa
\nab_4 \nab_3 \psi &=&-\squared_2 \psi+\lap_2 \psi-\frac 1 2 \trch \nab_3\psi+O(r^{-2}) \dk^{\leq 1}\psi  ,
\eeaa
we have
\beaa
&&[r \nab_4, \squared_2]\psi\\
 &=&-r\left(\frac 1 2 \trch +2\om\right) \big(-\squared_2 \psi+\lap_2 \psi-\frac 1 2 \trch \nab_3\psi \big)-\nab_4\nab_4\psi -\frac 1 4r( \trch)^2  \nab_3\psi \\
&&-r\trch \squared_2 \psi +O(r^{-2}) \dk^{\leq 1} \psi+O(r^{-3}) \dk^{\leq 2} \psi+r\Ddot_3 \big( \xi \c \Ddot_a \psi \big)+ r^{-1}\dk \big( \Ga_g \c \dk \psi)\\
&=&r\left(-\frac 1 2 \trch +2\om\right) \squared_2 \psi -r\left(\frac 1 2 \trch +2\om\right) \lap_2 \psi -\nab_4\nab_4\psi  \\
&&+O(r^{-2}) \dk^{\leq 1} \psi+O(r^{-3}) \dk^{\leq 2} \psi+r\Ddot_3 \big( \xi \c \Ddot_a \psi \big)+ r^{-1}\dk \big( \Ga_g \c \dk \psi),
\eeaa
as stated.


\section{Proof of Lemma \ref{LEMMA:COMMUTATOR-NAB-RHAT-SQUARE}}
\label{proof:lemma-comm-nab-Rhat-square}


 As a consequence of Lemma \ref{lemma:expression-wave-operator}, for scalar functions $f$ and $\psi$, we have
 \beaa
\square_\g( f \psi)&=& \square_\g(f) \psi+f \square_\g \psi - \nab_3 f \nab_4\psi- \nab_4f \nab_3 \psi +2\nab f \c \nab \psi,
\eeaa
which implies
\beaa
[f X, \square_\g]\psi&=& f [X, \square_\g]\psi+\nab_3 f \nab_4 X \psi +\nab_4f \nab_3 X \psi -2\nab f \c \nab X \psi-\square_\g (f) X\psi.
\eeaa
Using the definition in the outgoing frame $2\Rhat=\frac{\De}{r^2+a^2} e_4-\frac{|q|^2}{r^2+a^2}  e_3$, we compute
\beaa
2[\nab_\Rhat, \square_\g]\psi&=&\frac{\De}{r^2+a^2} [ \nab_4, \square_\g]\psi-\frac{|q|^2}{r^2+a^2} [  \nab_3, \square_\g]\psi\\
&&+e_4\left(\frac{\De}{r^2+a^2}\right) \nab_3 \nab_4 \psi-e_3\left(\frac{|q|^2}{r^2+a^2}\right)\nab_4 \nab_3 \psi\\
&&+ e_3 \left(\frac{\De}{r^2+a^2}\right)\nab_4 \nab_4 \psi -e_4\left(\frac{|q|^2}{r^2+a^2}\right) \nab_3 \nab_3 \psi -2 \nab\left(\frac{\De}{r^2+a^2}\right) \c \nab \nab_4 \psi \\
&&+ 2 \nab \left(\frac{|q|^2}{r^2+a^2}\right) \c \nab \nab_3 \psi -\square_\g\left(\frac{\De}{r^2+a^2}\right) \nab_4 \psi +\square_\g \left(\frac{|q|^2}{r^2+a^2}\right) \nab_3\psi.
\eeaa
Using Lemma \ref{LEMMA:COMMUTATOR-NAB3-NAB4-SQUARE} we deduce
\beaa
&&2[\nab_\Rhat, \square_\g]\psi\\
&=&\left( -\trch \frac{\De}{r^2+a^2}+\trchb \frac{|q|^2}{r^2+a^2}  \right) \square_\g \psi \\
&&+\Big( e_4\left(\frac{\De}{r^2+a^2}\right)+ (\trchb +2\omb) \frac{|q|^2}{r^2+a^2} \Big) \nab_3 \nab_4 \psi\\
&&-\Big( e_3\left(\frac{|q|^2}{r^2+a^2}\right)+(\trch +2\om)\frac{\De}{r^2+a^2}\Big) \nab_4 \nab_3 \psi\\
&&+\Big(  e_3 \left(\frac{\De}{r^2+a^2}\right) + 2\omb\frac{\De}{r^2+a^2} \Big) \nab_4 \nab_4 \psi - \Big( e_4\left(\frac{|q|^2}{r^2+a^2}\right)+2\om \frac{|q|^2}{r^2+a^2}  \Big)  \nab_3 \nab_3 \psi \\
&&+r^{-2} \dk^{\leq 1} \psi +O(ar^{-2}) \big( \nab_3\nab \psi +\nab_4 \nab \psi\big) +r^{-1} \dk \big( \Ga_g \c \dk \psi).
\eeaa
Using \eqref{eq:computations-e3e4-out} we write
\beaa
&&\Big(  e_3 \left(\frac{\De}{r^2+a^2}\right) + 2\omb\frac{\De}{r^2+a^2} \Big) \nab_4 \nab_4 \psi - \Big( e_4\left(\frac{|q|^2}{r^2+a^2}\right)+2\om \frac{|q|^2}{r^2+a^2}  \Big)  \nab_3 \nab_3 \psi \\
&=& \Big( - \frac{2ra^2\sin^2\th \De^2 }{|q|^4(r^2+a^2)^2}+\Ga_b \Big) \nab_4 \nab_4 \psi - \Big( \frac{2ra^2\sin^2\th}{(r^2+a^2)^2} + \Ga_g \Big)  \nab_3 \nab_3 \psi \\
&=& - \frac{2ra^2\sin^2\th  }{(r^2+a^2)^2} \Big( \frac{\De^2}{|q|^4} \nab_4 \nab_4 \psi+ \nab_3 \nab_3 \psi\Big) +\Ga_g \c \dk^{2} \psi.
\eeaa
Also, writing $\nab_3 \nab_4\psi=-\square_\g \psi -\lap \psi +r^{-1} \dk\psi $ and $\nab_4 \nab_3\psi=-\square_\g \psi -\lap \psi +r^{-1} \dk\psi $, we deduce 
 \beaa
[\nab_\Rhat, \square_\g]\psi&=&O(r^{-1})  \square_\g \psi+O(r^{-1})  \lap \psi +r^{-2} \dk^{\leq 1} \psi\\
&&+O(a^2r^{-3})  \big( \nab_4 \nab_4 \psi+ \nab_3 \nab_3 \psi\big)  +O(ar^{-2}) \big( \nab_3\nab \psi +\nab_4 \nab \psi\big) \\
&&+ \dk \big( \Ga_g \c \dk \psi),
\eeaa
as stated. The last commutations with $\triangle$ are obtained in a similar way as a consequence of Lemma \ref{LEMMA:COMMUTATOR-NAB3-NAB4-LAP}.


\section{Proof of Lemma \ref{LEMMA:COMMUTATIONOFHODGEELLIPTICORDER1WITHSQAURED2FDILUHS}}\label{proof:lemma:commofHodgewithsquare}


Combining Lemma \ref{le:Bochner-nointegrable} and Proposition \ref{Gauss-equation-2-tensors}, we have for $\xi \in \sk_1$, $u \in \sk_2$
\beaa
\DDd_2\DDs_2 \xi &=&-\f12\lap_1\xi -  \frac 1 2 \, \Kh \xi + \frac 1 4 (\atrch\nab_3+\atrchb \nab_4)\dual \xi  ,\\
\DDs_2 \DDd_2 u &=&-\f12\lap_2u+  \, \Kh  u  - \frac 1 4 (\atrch\nab_3+\atrchb \nab_4)\dual  u .
\eeaa
We then have for $\psi \in \sk_2$:
\beaa
\DDd_2 \lap_2\psi - \lap_1 \DDd_2\psi&=& \DDd_2\big( -2\DDs_2 \DDd_2 \psi+ 2 \Kh \psi -\frac 1 2(\atrch\nab_3+\atrchb \nab_4)\dual \psi \big)\\
&& - \big(-2\DDd_2\DDs_2- \Kh +\frac 1 2 (\atrch\nab_3+\atrchb \nab_4)\dual \big) \DDd_2\psi\\
&=&  2 \Kh \DDd_2\psi +2\nab  \Kh \c  \psi -\frac 1 2\DDd_2((\atrch\nab_3+\atrchb \nab_4)\dual \psi) \\
&& + \Kh \DDd_2 \psi -\frac 1 2 (\atrch\nab_3+\atrchb \nab_4)\dual \DDd_2\psi,
\eeaa
which gives
\beaa
\DDd_2 \lap_2\psi - \lap_1 \DDd_2\psi&=&  3 \Kh \DDd_2\psi   - (\atrch\nab_3+\atrchb \nab_4)\dual \DDd_2\psi\\
&&+ O(ar^{-3})\dk^{\leq 1}  \psi + r^{-1} \dk^{\leq 2} (\Ga_g \c \psi).
\eeaa

From Corollary \ref{corr:comm}, i.e.
    \beaa
  \,    [\nab_3,\DDd_2] u&=&-\frac 1 2 \trchb \,\big( \DDd_2 u-\eta\c u\big) +
  \frac 1 2 \atrchb\,\big(  \DDd_2\dual u  -\eta \c \dual u\big)      +(\eta-\ze)\c\nab_3 u \\
  &&+r^{-1} \Ga_b \c \dk^{\leq 1} u ,\\
\,   [\nab_4,\DDd_2] u&=&-\frac 1 2 \trch \, \big(\DDd_2 u-\etab\c u\big) +\frac 1 2 \atrch\, \big(\DDd_2\dual u -   \etab\c \dual u\big) 
+(\etab+\ze)\c\nab_4 u \\
&&+\xi \c \nabc_3 f +r^{-1}\Ga_g \c \dk^{\leq 1} u,
    \eeaa
    we deduce, see also Corollary \ref{corr:commwithrnab-M8},
        \beaa
  \,    [\nab_3, |q|\DDd_2] u&=&\frac 1 2 |q|\atrchb \DDd_2 \dual u  + |q|( \eta- \ze) \c \nab_3 u + O(ar^{-3} ) u  + r \Hc \nab_3 u +\Ga_b \c \dk^{\leq 1} u ,\\
\,   [\nab_4,|q| \DDd_2] u&=&\frac 1 2 |q| \atrch \DDd_2 \dual u + |q| (\etab+ \ze) \c \nab_4 u + O(ar^{-3} ) u +r\xi \c \nab_3 u +\Ga_g \c \dk^{\leq 1} u.
    \eeaa
We can therefore compute, using the null decomposition in Lemma \ref{lemma:expression-wave-operator},
\beaa
&&|q|\DDd_2\squared_2 - \squared_1|q|\DDd_2\\
 &=&-\frac 1 2 \big([|q| \DDd_2,\nab_3]\nab_4\psi+\nab_3[|q| \DDd_2,\nab_4]\psi+[|q| \DDd_2, \nab_4] \nab_3 \psi+\nab_4[|q|\DDd_2, \nab_3] \psi\big)\\
&&+\left(\omb -\frac 1 2 \trchb\right)[|q|\DDd_2, \nab_4]\psi+\left(\om -\frac 1 2 \trch\right)[|q|\DDd_2, \nab_3]\psi \\
&&+|q|\DDd_2 \lap_2\psi - \lap_1 |q|\DDd_2\psi+O(ar^{-2} )\dk^{\leq 1}\psi,
\eeaa
which gives
\beaa
&&|q|\DDd_2\squared_2 - \squared_1|q|\DDd_2\\
 &=&\frac 1 2 \Big[\left(\frac 1 2 |q|\atrchb \DDd_2 \dual   +|q|( \eta- \ze) \c \nab_3  \right)\nab_4\psi+\nab_3\left(\frac 1 2 |q| \atrch \DDd_2 \dual \psi + |q| (\etab+ \ze) \c \nab_4 \psi \right)\\
 &&+\left(\frac 1 2 |q| \atrch \DDd_2 \dual  + |q| (\etab+ \ze) \c \nab_4  \right)\nab_3 \psi+\nab_4\left(\frac 1 2 |q|\atrchb \DDd_2 \dual \psi  + |q|( \eta- \ze) \c \nab_3 \psi \right)\Big]\\
&&+  3 \Kh |q|\DDd_2\psi   - |q|(\atrch\nab_3+\atrchb \nab_4)\dual \DDd_2\psi\\
&&+ O(ar^{-2})\dk^{\leq 1}  \psi + \dk^{\leq 2} (\Ga_g \c \psi)+O(ar^{-2} )\dk^{\leq 1}\psi\\
&=& 3 \Kh |q|\DDd_2\psi -\frac 1 2 |q|\atrchb  \dual \nab_4\DDd_2\psi -\frac 1 2  |q| \atrch  \dual \nab_3\DDd_2 \psi  +|q|( \eta+\etab) \c \nab_3\nab_4\psi   \\
&&+ O(ar^{-2})\dk^{\leq 1}  \psi + \dk^{\leq 2} (\Ga_g \c \psi)+O(ar^{-2} )\dk^{\leq 1}\psi\\
&&+r \Hc \left( -\nab_3 \nab_4\psi- \frac 12 \trch   \nab_3 \psi\right)  -\frac 1 2\nab_4( r \Hc) \nab_3 \psi + \Ddot_3 (r\xi \c \nab_3 u).
\eeaa
Writing $\nab_3\nab_4=-\squared_2+\lap_2+O(r^{-1}) \dk^{\leq1}=-\squared_2-2\DDs_2\DDd_2+O(r^{-1}) \dk^{\leq1}$,
and recalling that, see \eqref{eq:atrch-e3-atrch-e4-etaetab-kerr},
\beaa
 \atrch e_3+\atrchb e_4+ 2(\eta+\etab) \c \dual \nab&=&\frac{4a\cos\th}{|q|^2} \T,
\eeaa
we have
\beaa
|q|\DDd_2\squared_2 - \squared_1|q|\DDd_2&=& 3 \Kh |q|\DDd_2\psi -\frac{2a\cos\th}{|q|}\dual \nab_\T\DDd_2\psi  -|q|( \eta+\etab) \c \squared_2\psi   \\
&&+ O(ar^{-2})\dk^{\leq 1}  \psi + \dk^{\leq 2} (\Ga_g \c \psi)+O(ar^{-2} )\dk^{\leq 1}\psi+O(ar^{-3})\dk^{\leq 2}\psi\\
&&+r \Hc \left( -\nab_3 \nab_4\psi- \frac 12 \trch   \nab_3 \psi\right)  -\frac 1 2\nab_4( r \Hc) \nab_3 \psi + \Ddot_3 (r\xi \c \nab_3 u).
\eeaa
Using that,  see Lemma \ref{lemma:basicpropertiesLiebTfasdiuhakdisug:chap9},
\beaa
\nab_\T\psi &=& \Lieb_\T\psi + \frac{4amr\cos\th}{|q|^4}\dual\psi,
\eeaa
we obtain
\beaa
|q|\DDd_2\squared_2 - \squared_1|q|\DDd_2&=& 3 \Kh |q|\DDd_2\psi -\frac{2a\cos\th}{|q|}\dual \DDd_2 \Lieb_\T\psi  -|q|( \eta+\etab) \c \squared_2\psi   \\
&&+ O(ar^{-2})\dk^{\leq 1}  \psi + \dk^{\leq 2} (\Ga_g \c \psi)+O(ar^{-2} )\dk^{\leq 1}\psi+O(ar^{-3})\dk^{\leq 2}\psi\\
&&+r \Hc \left( -\nab_3 \nab_4\psi- \frac 12 \trch   \nab_3 \psi\right)  -\frac 1 2\nab_4( r \Hc) \nab_3 \psi + \Ddot_3 (r\xi \c \nab_3 u).
\eeaa
Finally, using the linearized null structure equation \eqref{eq:linearized-nabc4Hc}, we deduce that $\nab_4(r \Hc) =r^{-1} \Ga_g$ and by writing $-\nab_3 \nab_4\psi- \frac 12 \trch   \nab_3 \psi=\squared_2 \psi + r^{-2} \dk^{\leq 2} \psi $, we obtain the stated identity.

Similarly, combining Lemma \ref{le:Bochner-nointegrable} and Proposition \ref{Gauss-equation-2-tensors}, we have for $\xi \in \sk_1$, $f \in \sk_0$
\beaa
\DDs_1\DDd_1\xi&=&-\lap_1\xi+ \, \Kh  \xi - \frac 1 2 (\atrch\nab_3+\atrchb \nab_4) \dual \xi, \\
\DDd_1 \DDs_1f&=&-\lap_0f.
\eeaa
We then have for $\psi \in \sk_1$:
\beaa
&&\DDd_1 \lap_1 \psi - \lap_0 \DDd_1 \psi \\
&=& \DDd_1 \big(- \DDs_1\DDd_1 \psi+ \, \Kh  \psi - \frac 1 2 (\atrch\nab_3+\atrchb \nab_4) \dual \psi  \big)+\DDd_1 \DDs_1\DDd_1 \psi\\
&=&\Kh  \DDd_1  \psi +O(ar^{-3}) \dk^{\leq 2} \psi + r^{-1} \dk^{\leq 2} (\Ga_g \c \psi),
\eeaa
from which we deduce the stated formula.


\chapter{Complement for Chapter \ref{CHAPTER-DERIVATION-MAIN-EQS}}\label{section:appendix-chap5}



\section{Proof of Proposition \ref{TEUKOLSKY-PROPOSITION}}
\label{section:proof-teukolsky-eq}


According to Proposition \ref{prop:bianchi:complex}, we have the following Bianchi identity for $A$:
\bea\label{Bianchi-identity-A}
 \nabc_3A -\frac 1 2 \DDc\hot B &=& -\frac{1}{2}\tr\Xb A +  2H   \hot B -3\ov{P}\Xh.
\eea
We apply $\nabc_4$ to \eqref{Bianchi-identity-A}:
\beaa
\nabc_4\nabc_3A &=&\frac 1 2 \nabc_4(\DDc\hot B )-\frac{1}{2}\nabc_4\tr\Xb A-\frac{1}{2}\tr\Xb \nabc_4A \\
&&+ 2\nabc_4( H )  \hot B + 2 H   \hot\nabc_4( B )-3\nabc_4\ov{P}\Xh-3\ov{P}\nabc_4\Xh.
\eeaa
According to Lemma \ref{COMMUTATOR-NAB-C-3-DD-C-HOT} applied to $B$, which is $1$-conformally invariant, we have
 \beaa
 \, [\nabc_4 , \DDc \hot ]B &=&- \frac 1 2 \tr X( \DDc\hot B )+ \underline{H} \hot \nabc_4 B+r^{-1} \Ga_g \c \dk^{\leq 1} B.
 \eeaa
We therefore obtain
 \beaa
&&  \nabc_4(\DDc\hot B )+ 4 H   \hot\nabc_4( B )\\
&=&\DDc\hot(\nabc_4 B )+ 4 H   \hot\nabc_4( B ) - \frac 1 2 \tr X \DDc\hot B + \underline{H} \hot \nabc_4 B+r^{-1} \Ga_g \c \dk^{\leq 1} B\\
 &=&\DDc\hot(\nabc_4 B )  + \left(4 H  + \underline{H} \right)\hot \nabc_4 B- \frac 1 2 \tr X \DDc\hot B  +r^{-1} \Ga_g \c \dk^{\leq 1} B.
 \eeaa
According to Proposition \ref{prop:bianchi:complex}, we have the following Bianchi identity for $B$: 
\beaa
\nabc_4B -\frac{1}{2} \DDbc \c A &=& -2\ov{\tr X} B +\frac{1}{2}A\c \ov{\Hb}+3 \ov{P} \  \Xi.
\eeaa
This gives
 \bea\label{term1}
 \begin{split}
& \nabc_4(\DDc\hot B )+ 4 H   \hot\nabc_4( B )\\
 &=\DDc\hot\left(\frac{1}{2} \DDbc \c A -2\ov{\tr X} B +\frac{1}{2}A\c \ov{\Hb}+3 \ov{P} \  \Xi\right)\\
 &  + \left(4 H  + \underline{H} \right)\hot \left(\frac{1}{2} \DDbc \c A -2\ov{\tr X} B +\frac{1}{2}A\c \ov{\Hb}+3 \ov{P} \  \Xi\right)- \frac 1 2 \tr X \DDc\hot B +r^{-1} \Ga_g \c \dk^{\leq 1} B\\
  &=\frac{1}{2}\DDc\hot (\DDbc \c A + \ov{\Hb} \c A) +\left(- \frac 1 2 \tr X -2\ov{\tr X} \right)\DDc \hot B \\
 &  + \left( 2H  +\frac 1 2 \underline{H} \right)\hot ( \DDbc \c A + \ov{\Hb} \c A)+2\left(-\DDc\ov{\tr X}-\ov{\tr X} \left(4 H  + \underline{H} \right) \right)\hot  B \\
 &+3 \ov{P}\left( \DDc \hot \Xi+(4H+\underline{H})\hot\Xi\right)+ 3 \DDc \ov{P} \hot \Xi+r^{-1} \Ga_g \c \dk^{\leq 1} B.
 \end{split}
 \eea
According to Proposition \ref{prop-nullstr:complex-conf}, we have the following null structure equation:
\beaa
\nabc_4\tr\Xb +\frac{1}{2}\tr X\tr\Xb &=& \DDc\c\ov{\Hb}+\Hb\c\ov{\Hb}+2\ov{P}+\Ga_b \c  \Ga_g,
\eeaa
which gives
\bea\label{term2}
\begin{split}
&-\frac{1}{2}\nabc_4\tr\Xb A-\frac{1}{2}\tr\Xb \nabc_4A\\
&=\left(\frac{1}{4}\tr X\tr\Xb -\frac 1 2 ( \DDc\c\ov{\Hb}+\Hb\c\ov{\Hb})-\ov{P}\right) A-\frac{1}{2}\tr\Xb \nabc_4A+\Ga_b \c  \Ga_g \c A.
\end{split}
\eea
According to Proposition \ref{prop-nullstr:complex-conf}, we have the following null structure equation:
\beaa
\nabc_4H-\nabc_3 \Xi &=&  -\frac{1}{2}\ov{\tr X}(H-\Hb) -\frac{1}{2}\Xh\c(\ov{H}-\ov{\Hb}) -B,
\eeaa
which gives
\bea\label{term3}
4\nabc_4( H )  \hot B&=&  -2\ov{\tr X}(H-\Hb)  \hot B+\left(\nabc_3\Xi +r^{-1} \Ga_g\right)\c B.
\eea
Finally using again Proposition \ref{prop-nullstr:complex-conf} and Proposition \ref{prop:bianchi:complex},
\beaa
\nabc_4\Xh+\Re(\tr X)\Xh&=&\frac 1 2 \DDc \hot \Xi +\frac 1 2 \Xi \hot (\Hb +H) -A, \\
\nabc_4P -\frac{1}{2}\DDc\c \ov{B} &=& -\frac{3}{2}\tr X P + \Hb \c\ov{B} -\ov{\Xi} \c \Bb +\Ga_b  \c A,
\eeaa
we obtain
\bea\label{term4}
\begin{split}
&-3\nabc_4\ov{P}\Xh-3\ov{P}\nabc_4\Xh\\
&=-3\left(-\frac{3}{2}\ov{\tr X} \ \ov{ P}\right)\Xh-3\ov{P}\left(-\frac 1 2 (\tr X+\ov{\tr X})\Xh -A+ \frac 1 2 \DDc \hot \Xi +\frac 1 2  \Xi \hot (\Hb +H)\right)\\
&+r^{-1} \Ga_g \c \dk^{\leq 1} B+ \Ga_b \c \Ga_g \c A \\
&= \left(\frac 3 2\tr X+6\ov{\tr X}\right)\ov{P}\Xh +3\ov{P}A-\frac 3 2 \ov{P}\left( \DDc \hot \Xi + \Xi \hot (\Hb +H)\right)+r^{-1} \Ga_g \c \dk^{\leq 1} B.
\end{split}
\eea
Summing $\frac 1 2 $ \eqref{term1}, \eqref{term2}, \eqref{term3} and \eqref{term4}, we obtain
\beaa
\nabc_4\nabc_3A &=& \left(- \frac 1 2 \tr X -2\ov{\tr X} \right)\frac 1 2 \DDc \hot B +\left(\frac 3 2\tr X+6\ov{\tr X}\right)\ov{P}\Xh\\
 &&  +\left(-\DDc\ov{\tr X}-5\ov{\tr X}  H  \right)\hot  B \\
&&+ \left(\frac{1}{4}\tr X\tr\Xb -\frac 1 2 ( \DDc\c\ov{\Hb}+\Hb\c\ov{\Hb})+2\ov{P}\right) A-\frac{1}{2}\tr\Xb \nabc_4A\\
&&+\frac{1}{4}\DDc\hot (\DDbc \c A + \ov{\Hb} \c A)+ \left( H  +\frac 1 4 \underline{H} \right)\hot ( \DDbc \c A + \ov{\Hb} \c A)\\
&&+r^{-1}  \dk^{\leq 1}\big( \Ga_g \c  B\big) + \nabc_3\Xi  \c B +  \Ga_b \c \Ga_g \c A.
\eeaa
Using again \eqref{Bianchi-identity-A} we have 
\beaa
&&\left(- \frac 1 2 \tr X -2\ov{\tr X} \right)\frac 1 2 \DDc \hot B +\left(\frac 3 2\tr X+6\ov{\tr X}\right)\ov{P}\Xh\\
&=&\left(- \frac 1 2 \tr X -2\ov{\tr X} \right)\left(\frac 1 2 \DDc \hot B -3\ov{P}\Xh\right)\\
&=&\left(- \frac 1 2 \tr X -2\ov{\tr X} \right)\left(\nabc_3A+\frac{1}{2}\tr\Xb A-2 H   \hot B \right).
\eeaa
This gives
\beaa
\nabc_4\nabc_3A &=& \left(- \frac 1 2 \tr X -2\ov{\tr X} \right)\nabc_3A-\frac{1}{2}\tr\Xb \nabc_4A\\
 &&  +\left(-\DDc\ov{\tr X}+(\tr X-\ov{\tr X} ) H  \right)\hot  B \\
&&+ \left(-\ov{\tr X} \tr \Xb -\frac 1 2 ( \DDc\c\ov{\Hb}+\Hb\c\ov{\Hb})+2\ov{P}\right) A\\
&&+\frac{1}{4}\DDc\hot (\DDbc \c A + \ov{\Hb} \c A)+ \left( H  +\frac 1 4 \underline{H} \right)\hot ( \DDbc \c A + \ov{\Hb} \c A)\\
&&+r^{-1}  \dk^{\leq 1}\big( \Ga_g \c  B\big) + \nabc_3\Xi  \c B +  \Ga_b \c \Ga_g \c A.
\eeaa
By the Codazzi equation
\beaa
\frac{1}{2}\DDbc\c\Xh  &=& \frac{1}{2}\DDc\ov{\tr X} -i\Im(\tr X)(H+\Xi)-B,
\eeaa
 the second line is absorbed by the quadratic terms in the last line. 
We also simplify
\beaa
&&\frac{1}{4}\DDc\hot (\DDbc \c A + \ov{\Hb} \c A)+ \left( H  +\frac 1 4 \underline{H} \right)\hot ( \DDbc \c A + \ov{\Hb} \c A)\\
&=& \frac{1}{4}\DDc\hot (\DDbc \c A)+\frac{1}{4}\DDc\hot ( \ov{\Hb} \c A)+ \left( H  +\frac 1 4 \underline{H} \right)\hot ( \DDbc \c A )\\
&&+ \left( H  +\frac 1 4 \underline{H} \right)\hot (\ov{\Hb} \c A).
\eeaa
Applying \eqref{simil-Leibniz} and \eqref{Leibniz-hot}, we write
\beaa
\DDc\hot ( \ov{\Hb} \c A)&=&  2(\DDc\c \ov{\Hb} ) A + 2(\ov{\Hb}\c \DDc) A,\\
\Hb \hot ( \ov{\Hb} \c A)&=& 2( \Hb \c \ov{\Hb}) \ A,
\eeaa
which implies
\beaa
\nabc_4\nabc_3A &=& \frac{1}{4}\DDc\hot (\DDbc \c A)+\left(- \frac 1 2 \tr X -2\ov{\tr X} \right)\nabc_3A-\frac{1}{2}\tr\Xb \nabc_4A\\
&&+ \left(-\ov{\tr X} \tr \Xb -\frac 1 2 ( \DDc\c\ov{\Hb}+\Hb\c\ov{\Hb})+2\ov{P}\right) A\\
&&+\frac{1}{2}\left( (\DDc\c \ov{\Hb} ) A + (\ov{\Hb}\c \DDc) A\right)+ \left( H  +\frac 1 4 \underline{H} \right)\hot ( \DDbc \c A )\\
&&+ H  \hot (\ov{\Hb} \c A)+  \frac 1 2( \Hb \c \ov{\Hb}) \ A+r^{-1}  \dk^{\leq 1}\big( \Ga_g \c  B\big) + \nabc_3\Xi  \c B +  \Ga_b \c \Ga_g \c A\\
&=& \frac{1}{4}\DDc\hot (\DDbc \c A)+\left(- \frac 1 2 \tr X -2\ov{\tr X} \right)\nabc_3A-\frac{1}{2}\tr\Xb \nabc_4A\\
&&+\frac{1}{2}  (\ov{\Hb}\c \DDc) A+ \left( H  +\frac 1 4 \underline{H} \right)\hot ( \DDbc \c A )+ \left(-\ov{\tr X} \tr \Xb +2\ov{P}\right) A\\
&&+  H  \hot (\ov{\Hb} \c A)+r^{-1}  \dk^{\leq 1}\big( \Ga_g \c  B\big) + \nabc_3\Xi  \c B +  \Ga_b \c \Ga_g \c A.
\eeaa
Using \eqref{Leib-eq-DDb}, we further simplify the angular part writing 
\beaa
\frac 12 (\ov{\Hb}\c \DDc )A+ \left( H  +\frac 1 4 \underline{H} \right)\hot ( \DDbc \c A )&=& \frac 12 (\ov{\Hb}\c \DDc )A+ \left( 2H  +\frac 1 2 \underline{H} \right)\c \DDbc A \\
&=& \frac 12 (\ov{\Hb}\c \DDc+\Hb \c \DDbc )A+(2H \c \DDbc ) A \\
&=& \left( 4H+\Hb +\ov{\Hb} \right)\c \nabc A.
\eeaa
This proves the proposition.


\section{Proof of Corollary \ref{COROLLARY-REAL-TEUK}}\label{proof-corollary-teukolsky}


Applying Lemma \ref{LEMMA:CONFORMAL-DD-LAP} to $A$, which is $2$-conformally invariant, we obtain
\beaa
\DDc\hot( \DDbc \c A)&=&4\lapc_2  A - 2i (\atrch \nabc_3+\atrchb \nabc_4)A\\
&&  +2   \left( \trch\trchb+ \atrch\atrchb+4\rho\right) A  \\
&&- 2i  \big(   \trch\atrchb-\trchb\atrch  +4\dual \rho \big)A.
\eeaa
From \eqref{Teukolsky-operator-ch5}, we then obtain
\beaa
\LL(A) &=&-\nabc_4\nabc_3A+\lapc_2  A -\frac 1 2  i (\atrch \nabc_3+\atrchb \nabc_4)A \\
&& +  \frac 1 2 \left( \trch\trchb+ \atrch\atrchb+4\rho\right) A  - i  \Big( \frac 1 2 \big(  \trch\atrchb-\trchb\atrch   \big)+2\dual \rho \Big)A\\
&&+\left(- \frac 1 2 \tr X -2\ov{\tr X} \right)\nabc_3A-\frac{1}{2}\tr\Xb \nabc_4A\\
&&+\left( 4H+\Hb +\ov{\Hb} \right)\c \nabc A+ \left(-\ov{\tr X} \tr \Xb +2\ov{P}\right) A+  H   \hot (\ov{\Hb} \c A).
\eeaa
Writing 
\beaa
-\frac{1}{2}\tr\Xb&=& -\frac{1}{2}\trchb +\frac 1 2  i \atrchb \\
- \frac 1 2 \tr X -2\ov{\tr X} &=&\left(-\frac 5 2 \trch\right) +i \left(- \frac 3 2 \atrch\right)\\
-\ov{\tr X} \tr \Xb +2\ov{P}&=& -\trch \trchb - \atrch \atrchb+2\rho + i\big(\trch \atrchb-\trchb\atrch-2\dual \rho\big),
\eeaa
we obtain
\beaa
\LL(A) &=&-\nabc_4\nabc_3A+\lapc_2  A +(4\eta+2\etab) \c \nabc A -\frac{1}{2}\trchb  \nabc_4A-\frac 5 2 \trch \nabc_3 A\\
&& +   \left( -\frac 1 2\trch\trchb-\frac 1 2 \atrch\atrchb+4\rho\right) A +  H   \hot (\ov{\Hb} \c A)\\
&&+ i \Big[-2\atrch \nabc_3A+  4  \dual \eta \c \nabc A+ \left( \frac 1 2  \trch \atrchb- \frac 1 2 \trchb\atrch-4\dual \rho\right) A\Big]
\eeaa
Using \eqref{simil-Leib-half}, we write
\bea\label{eq:alternative-HovHb}
H   \hot (\ov{\Hb} \c A)=\left( 4\eta \c \etab -4i \eta \wedge \etab \right)A,
\eea
and this completes the proof.


\section{Proof of Proposition \ref{PROP:FACTORIZATION-QF}}
\label{proof:prop-factorization-qf}


We first prove the relations in Kerr.
 For simplicity we choose  the normalization $e_3= e_3^{(in)}$ for which  $\omb=0$ and 
thus $\nabc_3=\nab_3$.
Let $I$ be the expression
\bea
I:=r \nab_3\Big( r^2   \big(   \nab_3  ( r f  A)\big)\Big),
\eea
with 
\beaa
f=\frac{\ov{q}^4}{r^4}= f_1 f^2 _2, \qquad f_1=\frac{|q|^4}{r^4}, \qquad f_2=\frac{\ov {q}}{q}.
\eeaa
 We calculate
\beaa
I&=& r \nab_3\Big( r^3  f\nab_3 A+    r^2  e_3( rf)  A \Big)\\
&=& r^4  f\nab_3 \nab_3 A+\big(  re_3(r^3 f)   + r^3 e_3( r f) \big)     \nab_3 A+ r\nab_3(r^2 e_3(r f)) A\\
&=& r^4  f\nab_3 \nab_3 A+ \big(2  r^4   e_3 f- 4 r^3  f\big) \nab_3 A+  
\big(r^4 \nab_3 e_3 (f)  - 4 r^3 e_3 f +2 r^2 f\big) A.
\eeaa
We deduce
\bea
\lab{eq:Prop-factorization-qf-I_1I_2}
\bsplit
I&= r^4  f\Big( \nab_3 \nab_3 A+ I_1    \nab_3 A+ I_2 A\Big),\\
 I_1&= 2 f^{-1}( e_3 f -2 r^{-1}  f),  \\
 I_2&= f^{-1} \Big(   \nab_3 e_3 (f)-  4 r^{-1}   e_3 f + 2 r^{-2} f\Big).
 \end{split}
\eea
We rewrite $I_2$  in the form
\beaa
 I_2&=& \nab_3\big( e_3 (f)- 2 r^{-1} f) \big) + 2 r^{-1} e_3 f + 2 r^{-2} f  -  4 r^{-1}   e_3 f + 2 r^{-2} f\\
  &=&  \nab_3\big( e_3 (f)- 2 r^{-1} f) \big) -2 r^{-1} \big( e_3 f - 2  r^{-1} f \big)= \frac 1 2  \nab_3( f I_1) - r^{-1}  f I_1.
\eeaa
Hence
\beaa
\bsplit
I&= r^4  f\Big( \nab_3 \nab_3 A_{11}+ I_1    \nab_3 A_{11}+ I_2 A\Big),\\
 I_1&= 2 f^{-1}( e_3 f -2 r^{-1}  f),  \\
 I_2&=  \frac 1 2  \nab_3( f I_1) - r^{-1}  f I_1.
 \end{split}
\eeaa
Now, in view of our choices of the scalar functions  $f_1, f_2$, 
\beaa
 e_3 f_1&=&\frac{ r^4 e_3(|q|^4)- |q|^4  e_3(r^4)}{ r^8}=\frac{ -4 r^5  |q|^2 + 4 |q|^4  r^3 }{ r^8}=
 -\frac{4|q|^2}{r^5}\Big( r^2 -|q|^2\Big)\\
 &=&\frac{ 4 a^2 \cos^2 |q|^2}{r^5}= \frac{ 4a^2 \cos^2 \th}{ r|q|^2 }  f_1,   \\
 e_3 f_2&=&- i \atrchb \frac{\ov{q}}{ q}=-i \atrchb  f_2 .
 \eeaa
 Hence
 \beaa
 e_3 f&=& e_3( f_1) f^2_2+  2 f_1 f_2  e_3( f_2)=\frac{ 4a^2 \cos^2 \th}{ r|q|^2 }  f_1  f^2_2  - 2 i  \atrchb f_1 f^2_2\\
 &=& \Big(\frac{ 4a^2 \cos^2 \th}{ r|q|^2 }-  2i \atrchb \Big) f
 \eeaa
 i.e.
 \beaa
 e_3 f&=&\Big(\frac{ 4a^2 \cos^2 \th}{ r|q|^2 }- i \atrchb \Big) f.
 \eeaa
We deduce
 \beaa
 I_1&=&  2 f^{-1}( e_3 f - 2 r^{-1} f) =  2 \Big( \frac{ 4a^2 \cos^2 \th}{ r|q|^2 }  - 2 i \atrchb  -\frac 2  r\Big)\\
 &=&  2 \Big( \frac{ 2a^2 \cos^2 \th}{ r|q|^2 } +\frac 2 r  -   \frac{2r}{|q|^2}      - \frac 2  r - 2 i \atrchb\Big)= 2 \trchb -2 \frac { \atrchb^2}{\trchb}-  4  i \atrchb .
 \eeaa
 Hence, recalling the definition of $C_1$ in \eqref{eq:C1-C2-comparison-Ma},
 $I_1= C_1$.
 
 It remains to calculate $I_2$.
\beaa
I_2&=& \frac 1 2  \nab_3( f I_1) - r^{-1}  f I_1=\frac 1 2  I_1   \nab_3 f+ \frac 1 2  f \nab_3 I_1 - r^{-1}  f I_1\\
 &=&  \frac 1 2  f \nab_3 I_1 + \frac 1 2  I_1\big(    \nab_3 f- 2r^{-1} f\big)=\frac 1 2  f \nab_3 I_1 +  \frac 1 4  f  I^2 _1  .
\eeaa
Therefore,
\beaa
I_2&=&\frac 1 2  f \big(  \nab_3 I_1 +   \frac 1 2   I^2 _1\big).
\eeaa
Now,
\beaa
 \nab_3 I_1&=& \nab_3\Big(  2 \trchb -2 \frac { \atrchb^2}{\trchb}-  4  i \atrchb \Big)\\
  &=& 2 \nab_3  \trchb  -2 \frac {    \trchb\nab_3  (  \atrchb^2)  -   \atrchb^2  \nab_3 \trchb  }{\trchb^2} -4 i \nab_3 \atrchb.
\eeaa
Recall, in Kerr,
   \beaa
   \nab_3\trchb&=&-\frac 1 2 \big( \trchb^2-\atrchb^2\big), \qquad   \nab_3\atrchb= -\trchb\atrchb.
   \eeaa
   Therefore,
   \beaa
    \nabc_3 E_1&=& - \trchb^2+ \atrchb^2  -2 \frac {    - 2 \trchb^2 \atrchb^2  +\frac 1 2   \atrchb^2  \big( \trchb^2-\atrchb^2\big)   }{\trchb^2}+ 4 i \trchb \atrchb\\
    &=&  - \trchb^2+ \atrchb^2 +\frac{3 \trchb^2 \atrchb^2+  \atrchb^4}{\trchb^2} + 4 i \trchb \atrchb\\
    &=&  - \trchb^2+ \atrchb^2+ 3\atrchb^2 + \frac{ \atrchb^4}{\trchb^2} + 4 i \trchb \atrchb\\
    &=&   - \trchb^2+ 4 \atrchb^2  + \frac{ \atrchb^4}{\trchb^2} + 4 i \trchb \atrchb.
    \eeaa
  Adding, simplifying  and  recalling  the definition of $C_2$ we deduce
   \beaa
   \nab_3 I_1+  \frac 1 2 I_1^2 &=& - \trchb^2+ 4 \atrchb^2  + \frac{ \atrchb^4}{\trchb^2} + 4 i \trchb \atrchb\\
  &+&2\Big(   \trchb - \frac { \atrchb^2}{\trchb}-  2  i \atrchb\Big)^2 \\
  &=&  - \trchb^2+ 4 \atrchb^2  + \frac{ \atrchb^4}{\trchb^2} + 4 i \trchb \atrchb\\
  &+&  2\Big( \trchb^2 +   \frac{ \atrchb^4}{\trchb^2} - 4 \atrchb^2- 2 \atrch^2 - 4 i \trchb \atrchb+ 4i \frac{\atrchb^3}{\trchb}\Big)\\
  &=& \trchb^2 - 8\atrchb^2+ 3  \frac{ \atrchb^4}{\trchb^2} -  4 i \trchb \atrchb + 8 i \trchb \atrchb\\
  &=& 2 C_2.
   \eeaa
   Hence $I_2 = C_2$ and, recalling the definition of  $f=\frac{\ov{q}^4}{r^4} $, 
 \beaa
 I&= &r^4  f\Big( \nab_3 \nab_3 A+ I_1    \nab_3 A+ I_2 A\Big)=\frac{\ov{q} }{q }   q\ov{q}^3  \Big( \nab_3 \nab_3 A+ C_1    \nab_3 A+ C_2 A\Big)= \frac{\ov{q} }{q } \qf,
\eeaa
as stated in \eqref{eq:Prop-factorization-qf}.

It remains to check \eqref{eq:secondfactorization-qf}. To do this 
it  helps to write  the expression 
$ E=r^2\nab_3\nab_3\big(\frac{\ov{q}^4}{r^2} A \big) $  using  $ f= \frac{\ov{q}^4}{r^4}$,
\beaa
E&=& r^2\nab_3\nab_3\big( r^2 f A  \big)= r^4 f \Big( \nab_3\nab_3  A+
\frac{2}{r^2 f} e_3(r^2f) \nab_3 A+\frac{1}{r^2 f}  e_3 e_3 (r^2f)   A\Big)\\
&=&r^4  f\Big( \nab_3 \nab_3 A_{11}+ I_1    \nab_3 A_{11}+ I_2 A\Big),
\eeaa
with $I_1, I_2$ as in  \eqref{eq:Prop-factorization-qf-I_1I_2}. Therefore
\beaa
r^2\nab_3\nab_3\big( r^2 f A  \big)=r \nab_3\Big( r^2   \big(   \nab_3  ( r f  A)\big)\Big)=\frac{\ov{q} }{q } \qf,
\eeaa
as stated.

Using that $\widecheck{e_3(r)}= r \Ga_b$ and $e_3(\cos\th)=\Ga_b$ we deduce the relations in perturbations of Kerr.


\section{Proof of Theorem \ref{MAIN-THEOREM-PART1}}\label{proof:theorem-main-Part1}


In this section we present the proof of Theorem \ref{MAIN-THEOREM-PART1}. The computations are obtained in the outgoing frame.


\subsection{Step 1. Compute the commutator $[Q,\LL ]$}\label{proof:step1-main-thm-part1}
\label{proof-appendix}


 Recall the Teukolsky equation as in Proposition \ref{TEUKOLSKY-PROPOSITION}, i.e.
\bea\label{Teuk-repeat}
 \LL(A)&=& \err[\LL(A)].
\eea
We apply the Chandrasekhar transformation, i.e. the operator $Q$ as defined in \eqref{definition-Q(A)}, to the above. We obtain
 \bea\label{Teuk-repeat=2}
 \LL(Q(A))+[Q, \LL](A)&=& Q( \err[\LL(A)]).
 \eea

We compute the commutator $[Q,\LL ]$ between $\LL$ and the second order differential operator $Q$ for any scalar functions $C_1$ and $C_2$. In order to obtain cancellations in the commutator, we impose conditions on the functions $C_1$ and $C_2$, which for $a=0$ coincide with the ones in Schwarzschild as in \cite{KS}.  We obtain the following.

\begin{proposition}\label{first-intermediate-step-main-theorem} 
Let $Q(A)= \nabc_3 \nabc_3 A+ C_1 \nabc_3 A + C_2  A$ with $C_1$ and $C_2$ given by
\bea\label{eq:first-assumptions-C1-C2}
C_1= 2\trchb +\widetilde{C_1}, \qquad C_2=\frac 1 2 \trchb^2 + \widetilde{C_2}
\eea
where $\widetilde{C_1}$ and $\widetilde{C_2}$ are complex functions satisfying $\widetilde{C_1}=O(|a|r^{-2})$, $\widetilde{C_2}=O(|a|r^{-3})$.
 Then the commutator between $Q$ and $\LL$ is given by 
 \bea\label{final-commutator}
 \begin{split}
 [Q, \LL](A) &= 4\etab  \c \nabc Q(A)- 2\trchb \nabc_4Q(A)+\hat{V} Q(A)+ L_Q(A)+\err[ [Q, \LL]A],
 \end{split}
 \eea
 where 
 \begin{itemize}
\item the potential 
$\hat{V}$ is given by
\bea
\hat{V}&=& I_{33}+J_{33}+K_{33}+M_{33}
\eea
where $ I_{33}$, $J_{33}$, $K_{33}$, $M_{33}$ are given in Section \ref{proof-appendix},
\item  $ L_{Q}[A ]$ is  a  second order linear  operator  in $A$,   given by
  \bea\label{definition-L_QA}
  \begin{split}
 L_{Q}[A]&=Z_{43} \  \nabc_4\nabc_3A+Z_4 \  \nabc_4A+Z_{a3}  \nabc_a \nabc_3A \\
&+ Z_3\nabc_3A+Z_a \nabc_a A +Z_0  A,
\end{split}
\eea
where $Z_{a3}$, $Z_a$ are complex one-forms and $Z_{43}$, $Z_4$, $Z_3$ and $Z_0$ are complex functions  of $(r, \th)$, all of which vanish for zero angular momentum, having the following fall-off in $r$ 
\beaa
Z_{43},  Z_{a3}=O\left(\frac{|a|}{r^3}\right), \qquad  Z_4, Z_3, Z_a=O\left(\frac{|a|}{r^4}\right), \qquad   Z_0=O\left(\frac{|a|}{r^5}\right).
\eeaa
More precisely,
\beaa
Z_{43}&=&\nabc_3 \widetilde{C_1} + \trchb \widetilde{C_1}, \\
Z_4&=&\nabc_3 \widetilde{C_2}+ 2\trchb  \widetilde{C_2}-\frac 1 4 (\trchb^2+ \atrchb^2) \widetilde{C_1},
\eeaa
and 
\beaa
 Z_{a3}&=&I_{a3}+J_{a3}+L_{a3}+M_{a3}-4C_1 \etab,
\eeaa
where $ I_{a3}$, $J_{a3}$, $L_{a3}$, $M_{a3}$ are given in Section \ref{proof-appendix},
\item the error terms $\err[ [Q, \LL]A]$ are given by
\beaa
\err[ [Q, \LL]A]&=&r^{-2} \dk^{\leq 3}(\Ga_g \c (A, B))+\nab_3 \big( r^{-2} \Ga_b  \c \frak{d}^{\leq 2} A \big)+r^{-1}\dk\big( (\Ga_b \c \Ga_g) A \big).
\eeaa
\end{itemize}
 \end{proposition}
 
The proof of Proposition \ref{first-intermediate-step-main-theorem} is obtained as a result of the computations below.
 We first collect the following commutators for $Q$.

\begin{proposition}\label{commutation-general-Q} 
Let $U \in \sk_2(\CCC)$ $s$-conformally invariant. We have:
\begin{itemize}
\item the following commutators with $\nabc_3$ and $\nabc_4$:
\beaa
\, [Q, \nabc_3] U&=& (- \nabc_3 C_1) \  \nabc_3U+( - \nabc_3 C_2 ) \  U,\\
\, [Q, \nabc_4] U  &=&4 (\eta-\etab ) \c \nabc \nabc_3 U\\
 && +\Big( 2\nabc_3(\eta-\etab )+(2C_1 -   \trchb ) (\eta-\etab )+\atrchb  \dual (\eta-\etab )\Big) \c \nabc U\\
 &&+\Big(\mathcal{V}^s_{[3, 4]} +\mathcal{V}^{s-1}_{[3, 4]}  +2 \eta \c (\eta-\etab ) -\nabc_4 C_1 \Big) \nabc_3 U\\
 &&+\Big(\nabc_3 \mathcal{V}^s_{[3, 4]} + C_1  (\mathcal{V}^s_{[3, 4]} ) -\nabc_4 C_2 \Big) U +2 (\eta-\etab ) \c \mathcal{V}^s_{[3,a]}(U)\\
   &&+ r^{-3} \Ga_b  \c \frak{d}^{\leq 1} U+\dk\big( (\Ga_b \c \Ga_g) U \big),
 \eeaa
 where $\mathcal{V}^s_{[3, 4]}$ and $\mathcal{V}^s_{[3,a]}$ are given by \eqref{expression-C-0-3-4} and \eqref{expression-C-0-3-a} respectively. 
 \item the following commutator with $\nabc_a$:
  \beaa
&&\, [Q, \nabc_a] U\\
&=&-   \trchb\,  \nabc_a\nabc_3 U- \atrchb\,  \dual \nabc_a\nabc_3 U+2\eta_a \nabc_3\nabc_3 U\\
&&+\left(-\frac  1 2   \trchb (C_1-\trchb)-\frac 1 2 \atrchb^2 \right) \, \nabc_a U-\frac 1 2 \atrchb \big(C_1-2\trchb\big) \, \dual \nabc_a  U\\
&&+\left(\nabc_3\eta_a+(C_1-\frac  1 2   \trchb )\eta_a-\frac 1 2 \atrchb \dual  \eta_a -\nabc_a  C_1  \right) \nabc_3 U \\
&&+\nabc_3(\mathcal{V}^s_{[3,a]}(U))+\mathcal{V}^{s-1}_{[3,a]}(\nabc_3 U)-\frac  1 2   \trchb\mathcal{V}^s_{[3,a]}(U)-\frac 1 2 \atrchb \dual \mathcal{V}^s_{[3,a]}(U)\\
&&   -\nabc_a (C_2)  U+ C_1\mathcal{V}^s_{[3,a]}(U) +\nab_3 \big( r^{-1} \Ga_b  \c \frak{d}^{\leq 1} U \big) .
\eeaa
  
 \item the following commutator with $\DDbc \c$:
   \beaa
\, [Q, \DDbc \c] U  &=&2\ov{H} \c \nabc_3\nabc_3 U- \ov{\tr\Xb}\,\ov{\DDc} \c  \nabc_3U+\frac 1 2\ov{\tr\Xb} (\ov{\tr\Xb}-C_1)\,  \ov{\DDc} \c U \\
&&+\left(\nabc_3\ov{H} +(- (s-2)\ov{\tr\Xb}+C_1)\,  \ov{H}-\DDbc ( C_1 )\right)\c \nabc_3U\\
&&+\left(\frac 1 2 (s-2) \ov{\tr\Xb}\left(-  \nabc_3\ov{H}+ (\ov{\tr\Xb}-C_1)\,  \ov{H}\right)-\DDbc (C_2) \right)\c U\\
&&+\nab_3 \big( r^{-1} \Ga_b  \c \frak{d}^{\leq 1} U \big).
 \eeaa
 \end{itemize}

Let $F\in \sk_1(\CCC)$ of conformal type $s$. Then we have the following commutator with $\DDc \hot$:
\beaa
 \, [Q, \DDc \hot] F  &=&2H \hot \nabc_3\nabc_3 F-  \tr \Xb  \DDc\hot \nabc_3F+\frac 1 2 (\tr \Xb)(\tr\Xb-C_1)  \DDc\hot F \\
 &&+\left(\nabc_3H+(- (s+1) \tr \Xb +C_1 )H -\DDc ( C_1 )\right) \hot \nabc_3 F\\
 &&+\left( \frac 1 2 (s+1) \tr \Xb \left(-\nabc_3H+(\tr\Xb-C_1)H\right)  -\DDc (C_2) \right)\hot F \\
 &&+\nab_3 \big( r^{-1} \Ga_b  \c \frak{d}^{\leq 1} U \big).
 \eeaa
\end{proposition}

\begin{proof} 
We compute 
\beaa
\, [Q, \nabc_3] U&=&(\nabc_3\nabc_3  + C_1   \nabc_3 +C_2 )\nabc_3U\\
&& - \nabc_3(\nabc_3\nabc_3U  + C_1  \nabc_3U +C_2 U)\\
&=& (- \nabc_3 C_1)   \nabc_3U+( - \nabc_3C_2 )   U.
\eeaa
 Similarly, we compute
 \beaa
 [Q, \nabc_4] U  &=&\nabc_3([\nabc_3, \nabc_4]U)+ [\nabc_3,\nabc_4]\nabc_3 U \\
 &&+ C_1   [\nabc_3, \nabc_4]U  -\nabc_4( C_1 )  \nabc_3U  -\nabc_4(C_2) \ U.
 \eeaa
Recall from Lemma \ref{commutator-nab-3-nab-4}, 
\beaa
\, [\nabc_3, \nabc_4] U   &=& 2(\eta-\etab ) \c \nabc U +\mathcal{V}^s_{[3, 4]}  U+ (\Ga_b \c \Ga_g) U ,
 \eeaa   
 where 
 \bea\label{expression-C-0-3-4}
\mathcal{V}^s_{[3, 4]} =2s\left(\rho   -\eta\c\etab\right) +4i \left(- \rhod+ \eta \wedge \etab \right).
 \eea 
In particular, since $\nabc_3U$ is conformal of type $s-1$, we have 
\beaa
[\nabc_3, \nabc_4] \nabc_3U  &=& 2(\eta-\etab ) \c \nabc \nabc_3U +\mathcal{V}^{s-1}_{[3, 4]}  \nabc_3U + (\Ga_b \c \Ga_g) \dk U .
\eeaa
 On the other hand we compute
 \beaa
  \nabc_3([\nabc_3, \nabc_4]U) &=&2 (\eta-\etab ) \c \nabc_3\nabc U+2\nabc_3(\eta-\etab ) \c \nabc U\\
  &&+\mathcal{V}^s_{[3, 4]}  \nabc_3 U+\nabc_3 \mathcal{V}^s_{[3, 4]}  U+\dk\big( (\Ga_b \c \Ga_g) U \big).
 \eeaa
Recall from Lemma \ref{commtator-3-a}, 
 \beaa
 [\nabc_3, \nabc_a] U_{bc} &=&-\frac  1 2   \trchb\, \nabc_a U_{bc}-\frac 1 2 \atrchb\, \dual \nabc_a  U_{bc}+\eta_a \nabc_3 U_{bc}+\mathcal{V}^s_{[3,a]}(U)\\
 &&+ r^{-1} \Ga_b  \c \frak{d}^{\leq 1} U,
\eeaa
where
\bea\label{expression-C-0-3-a}
\begin{split}
\mathcal{V}^s_{[3,a]}(U)&= -\frac  1 2   \trchb\, \Big(s(\eta_a) U_{bc}+\eta_bU_{ac}+\eta_c U_{ab}-\de_{a b}(\eta \c U)_c-\de_{a c}(\eta \c U)_b \Big)\\
&-\frac 1 2 \atrchb\, \Big(s (\dual\eta_a) U_{bc} +\eta_b \dual U_{ac}+\eta_c \dual U_{ab}- \in_{a b}(\eta \c  U)_c- \in_{a c}(\eta \c  U)_b \Big).
\end{split}
\eea
We therefore obtain
\beaa
&&\,  \nabc_3([\nabc_3, \nabc_4]U)\\
 &=&2 (\eta-\etab ) \c \Big(\nabc \nabc_3 U-\frac  1 2   \trchb\, \nabc U-\frac 1 2 \atrchb\, \dual \nabc  U+\eta \nabc_3 U+\mathcal{V}^s_{[3,a]}(U)\Big) \\
  &&+2\nabc_3(\eta-\etab ) \c \nabc U+\mathcal{V}^s_{[3, 4]}  \nabc_3 U+\nabc_3 \mathcal{V}^s_{[3, 4]}  U\\
  &&+ r^{-3} \Ga_b  \c \frak{d}^{\leq 1} U+\dk\big( (\Ga_b \c \Ga_g) U \big)\\
   &=&2 (\eta-\etab ) \c \nabc \nabc_3 U +\Big( 2\nabc_3(\eta-\etab ) -   \trchb  (\eta-\etab )+\atrchb  \dual (\eta-\etab )\Big) \c \nabc U\\
  &&+\Big(\mathcal{V}^s_{[3, 4]} +2 \eta \c (\eta-\etab )\Big) \nabc_3 U+\nabc_3 \mathcal{V}^s_{[3, 4]}  U +2 (\eta-\etab ) \c \mathcal{V}^s_{[3,a]}(U)\\
  &&+ r^{-3} \Ga_b  \c \frak{d}^{\leq 1} U+\dk\big( (\Ga_b \c \Ga_g) U \big).
\eeaa
 Putting the above together, we obtain the stated expression.

We compute 
\beaa
\, [Q, \nabc_a] U&=& \nabc_3([\nabc_3, \nabc_a]U)+ [\nabc_3,\nabc_a]\nabc_3 U + C_1  [\nabc_3, \nabc_a]U \\
  && -\nabc_a ( C_1 )  \nabc_3U  -\nabc_a (C_2)  U.
  \eeaa
We have, again from Lemma \ref{commtator-3-a}, 
\beaa
[\nabc_3,\nabc_a]\nabc_3 U&=&-\frac  1 2   \trchb\, \nabc_a \nabc_3U-\frac 1 2 \atrchb\, \dual \nabc_a  \nabc_3U+\eta_a \nabc_3 \nabc_3 U\\
&&+\mathcal{V}^{s-1}_{[3,a]}(\nabc_3 U)+ r^{-1} \Ga_b  \c \frak{d}^{\leq 2} U,
\eeaa
and 
  \beaa
\nabc_3( [\nabc_3, \nabc_a] U) &=&-\frac  1 2   \trchb\, \nabc_3\nabc_a U-\frac 1 2 \atrchb\, \dual \nabc_3\nabc_a  U+\eta_a \nabc_3\nabc_3 U\\
&&-\frac  1 2  (\nabc_3 \trchb)\, \nabc_a U-\frac 1 2 (\nabc_3\atrchb)\, \dual \nabc_a  U\\
&&+(\nabc_3\eta_a) \nabc_3 U+\nabc_3(\mathcal{V}^s_{[3,a]}(U))+r^{-1}  \dk (\Ga_b  \c \frak{d}^{\leq 1} U).
\eeaa
Applying once again the commutator, the above becomes
  \beaa
&&\nabc_3( [\nabc_3, \nabc_a] U_{bc}) \\
&=&-\frac  1 2   \trchb\, \Big( \nabc_a\nabc_3 U_{bc}-\frac  1 2   \trchb\, \nabc_a U_{bc}-\frac 1 2 \atrchb\, \dual \nabc_a  U_{bc}+\eta_a \nabc_3 U_{bc}+\mathcal{V}^s_{[3,a]}(U)\Big)\\
&&-\frac 1 2 \atrchb\, \dual \Big[ \nabc_a\nabc_3 U_{bc}-\frac  1 2   \trchb\, \nabc_a U_{bc}-\frac 1 2 \atrchb\, \dual \nabc_a  U_{bc}\\
&&+\eta_a \nabc_3 U_{bc}+\mathcal{V}^s_{[3,a]}(U)\Big]\\
&&+\eta_a \nabc_3\nabc_3 U_{bc}-\frac  1 2  (\nabc_3 \trchb)\, \nabc_a U_{bc}-\frac 1 2 (\nabc_3\atrchb)\, \dual \nabc_a  U_{bc}\\
&&+(\nabc_3\eta_a) \nabc_3 U_{bc}+\nabc_3(\mathcal{V}^s_{[3,a]}(U))+r^{-1}  \dk (\Ga_b  \c \frak{d}^{\leq 1} U)\\
&=&-\frac  1 2   \trchb\,  \nabc_a\nabc_3 U_{bc}-\frac 1 2 \atrchb\,  \dual \nabc_a\nabc_3 U_{bc}+\eta_a \nabc_3\nabc_3 U_{bc}\\
&&+\big(\frac  1 2   \trchb^2-\frac 1 2 \atrchb^2 \big) \, \nabc_a U_{bc}+\big(\trchb \atrchb \big)\, \dual \nabc_a  U_{bc}\\
&&+(\nabc_3\eta_a-\frac  1 2   \trchb \eta_a-\frac 1 2 \atrchb \dual  \eta_a )  \nabc_3 U_{bc}\\
&&+\nabc_3(\mathcal{V}^s_{[3,a]}(U))-\frac  1 2   \trchb\mathcal{V}^s_{[3,a]}(U)-\frac 1 2 \atrchb \dual \mathcal{V}^s_{[3,a]}(U)+r^{-1}  \dk (\Ga_b  \c \frak{d}^{\leq 1} U).
\eeaa
Putting the above together, we obtain the stated expression.

We compute  
\beaa
\, [Q, \DDbc \c] U  &=&\nabc_3([\nabc_3, \DDbc \c]U)+ [\nabc_3,\DDbc \c]\nabc_3 U + C_1 \  [\nabc_3, \DDbc \c]U \\
  && -\DDbc ( C_1 )\c  \nabc_3U  -\DDbc (C_2) \c  U.
 \eeaa
 Recall from Lemma \ref{commtator-3-a}, 
  \beaa
 [\nabc_3 , \ov{\DDc} \c] U &=&- \frac 1 2\ov{\tr\Xb}\, \left( \ov{\DDc} \c U + (s-2) \ov{H}\c U\right) +\ov{H} \c \nabc_3 U+ r^{-1}\Ga_b \c \frak{d}^{\leq 1} U.
  \eeaa
  In particular, since $\nabc_3U$ is conformal of type $s-1$, we have 
  \beaa
 [\nabc_3 , \ov{\DDc} \c] \nabc_3 U &=&- \frac 1 2\ov{\tr\Xb}\, \left( \ov{\DDc} \c \nabc_3 U + (s-3) \ov{H}\c \nabc_3 U\right) +\ov{H} \c \nabc_3 \nabc_3 U\\
 &&+r^{-1} \Ga_b  \c \frak{d}^{\leq 2} U.
  \eeaa
  On the other hand we compute 
  \beaa
&& \nabc_3( [\nabc_3 , \ov{\DDc} \c] U )\\
&=&- \frac 1 2\nabc_3\ov{\tr\Xb}\, \left( \ov{\DDc} \c U + (s-2) \ov{H}\c U\right)\\
 &&- \frac 1 2\ov{\tr\Xb}\, \left( \nabc_3\ov{\DDc} \c U + (s-2) \nabc_3\ov{H}\c U+ (s-2) \ov{H}\c \nabc_3U\right)\\
 && +\nabc_3\ov{H} \c \nabc_3 U +\ov{H} \c \nabc_3\nabc_3 U+r^{-1}\dk(  \Ga_b \c \frak{d}^{\leq 1} U)\\
 &=& \frac 1 2(\ov{\tr\Xb})^2\, \left( \ov{\DDc} \c U + (s-2) \ov{H}\c U\right)\\
 &&- \frac 1 2\ov{\tr\Xb}\, \left(\ov{\DDc} \c  \nabc_3U  + (s-2) \nabc_3\ov{H}\c U+ (s-1) \ov{H}\c \nabc_3U\right)\\
 && +\nabc_3\ov{H} \c \nabc_3 U +\ov{H} \c \nabc_3\nabc_3 U+r^{-1}\dk(  \Ga_b \c \frak{d}^{\leq 1} U).
  \eeaa
  Putting the above together we obtain the desired formula.

  We compute  
\beaa
\, [Q, \DDc \hot] F  &=&\nabc_3([\nabc_3, \DDc \hot]F)+ [\nabc_3,\DDc \hot]\nabc_3 F + C_1 \  [\nabc_3, \DDc \hot]F \\
  && -\DDc ( C_1 )\hot  \nabc_3F  -\DDc (C_2) \hot  F.
 \eeaa
Recall from Lemma \ref{commtator-3-a}, 
  \beaa
 \, [\nabc_3 , \DDc \hot ]F &=&- \frac 1 2 \tr \Xb \left( \DDc\hot F + (s+1)H \hot F \right)  + H \hot \nabc_3 F+r^{-1}  \Ga_b  \c \dk^{\leq 1} F .
 \eeaa
In particular, since $\nabc_3F$ is conformal of type $s-1$, we have 
  \beaa
 \, [\nabc_3 , \DDc \hot ]\nabc_3F &=&- \frac 1 2 \tr \Xb \left( \DDc\hot \nabc_3F + (s)H \hot \nabc_3F \right)  + H \hot \nabc_3 \nabc_3F\\
 &&+r^{-1}  \Ga_b  \c \dk^{\leq 2} F.
 \eeaa
 On the other hand we compute 
 \beaa
&& \nabc_3(\, [\nabc_3 , \DDc \hot ]F )\\
&=&- \frac 1 2 \nabc_3\tr \Xb \left( \DDc\hot F + (s+1)H \hot F \right)\\
 &&- \frac 1 2 \tr \Xb \left( \nabc_3\DDc\hot F + (s+1)\nabc_3H \hot F+ (s+1)H \hot \nabc_3F \right)  \\
 &&+ \nabc_3H \hot \nabc_3 F+ H \hot \nabc_3\nabc_3 F+r^{-1}\dk(  \Ga_b  \c \dk^{\leq 1} F) \\
 &=& \frac 1 2 (\tr \Xb)^2 \left( \DDc\hot F + (s+1)H \hot F \right)\\
 &&- \frac 1 2 \tr \Xb \left( \DDc\hot \nabc_3F   + (s+1)\nabc_3H \hot F+ (s+2)H \hot \nabc_3F \right)  \\
 &&+ \nabc_3H \hot \nabc_3 F+ H \hot \nabc_3\nabc_3 F+r^{-1}\dk(  \Ga_b  \c \dk^{\leq 1} F).
 \eeaa
 Putting the above together we obtain the desired formula. 
\end{proof}

 We now compute the commutator between $Q$ and $\LL$.  
Using \eqref{Teukolsky-operator-ch5}, we separate the commutator into 
\beaa
[Q, \LL]A&=& I + J+K+L+M+N
\eeaa
where 
\beaa
I&=& - [Q, \nabc_4\nabc_3]A, \qquad\qquad\qquad\quad\,\,\, J= \frac{1}{4}[Q, \DDc\hot\DDbc \c ]A,\\
K&=&  \left[Q, \left(- \frac 1 2 \tr X -2\ov{\tr X} \right)\nabc_3\right]A, \qquad L= [Q,-\frac{1}{2}\tr\Xb \nabc_4]A,\\
M&=&[Q, \left( 4H+\Hb +\ov{\Hb} \right)\c \nabc] A, \qquad\quad\,\, N= [Q, \left(-\ov{\tr X} \tr \Xb +2\ov{P}\right)] A+ [Q, H   \hot \ov{\Hb} \c ]A.
\eeaa


\subsubsection{Expression for $I$}


We have
\beaa
 I&=& -[Q, \nabc_4] \nabc_3A-\nabc_4( [Q,  \nabc_3]A).
\eeaa
Using Proposition \ref{commutation-general-Q}, applied to $U=\nabc_3A$ of conformal type $s=1$, we obtain
\beaa
&&\, [Q, \nabc_4] \nabc_3A  \\
&=&4 (\eta-\etab ) \c \nabc \nabc_3 \nabc_3A\\
 && +\Big( 2\nabc_3(\eta-\etab )+(2C_1 -   \trchb ) (\eta-\etab )+\atrchb  \dual (\eta-\etab )\Big) \c \nabc \nabc_3A \\
 &&+\Big(\mathcal{V}^{s=1}_{[3, 4]} +\mathcal{V}^{s-1=0}_{[3, 4]}  +2 \eta \c (\eta-\etab ) -\nabc_4( C_1 )\Big) \nabc_3 \nabc_3A\\
 &&+\Big(\nabc_3 \mathcal{V}^{s=1}_{[3, 4]} + C_1  (\mathcal{V}^{s=1}_{[3, 4]} ) -\nabc_4(C_2) \Big) \nabc_3A +2 (\eta-\etab ) \c \mathcal{V}^{s=1}_{[3,a]}(\nabc_3A)\\
 &&+ r^{-3} \Ga_b  \c \frak{d}^{\leq 2} A+\dk\big( (\Ga_b \c \Ga_g) \dk^{\leq 1} A \big).
\eeaa
We also deduce
 \beaa
 \nabc_4( [Q,  \nabc_3]A)&=& \nabc_4( (- \nabc_3 C_1) \  \nabc_3A+( - \nabc_3C_2 ) \  A)\\
 &=& (- \nabc_3 C_1) \  \nabc_4\nabc_3A+ (- \nabc_4\nabc_3 C_1) \  \nabc_3A\\
 &&+( - \nabc_3C_2 ) \  \nabc_4A+( - \nabc_4\nabc_3C_2 ) \  A.
 \eeaa
 We therefore obtain
 \beaa
  I&=& -4 (\eta-\etab ) \c \nabc \nabc_3 \nabc_3A+I_{43}  \nabc_4\nabc_3A+I_{33} \nabc_3 \nabc_3A\\
  && +I_{a3} \c \nabc \nabc_3A +I_{4}  \nabc_4A+I_{3} \nabc_3A +I_{0}  A\\
  &&+ r^{-3} \Ga_b  \c \frak{d}^{\leq 2} A+\dk\big( (\Ga_b \c \Ga_g) \dk^{\leq 1} A \big),
 \eeaa
 where
 \beaa
 I_{43}&=& \nabc_3 C_1 \\
 I_{33}&=& -\mathcal{V}^{s=1}_{[3, 4]} -\mathcal{V}^{s-1=0}_{[3, 4]}  -2 \eta \c (\eta-\etab ) +\nabc_4( C_1 )\\
 I_{a3}&=& - 2\nabc_3(\eta-\etab )-(2C_1 -   \trchb ) (\eta-\etab )-\atrchb  \dual (\eta-\etab )\\
 I_{4}&=&  \nabc_3C_2\\
 I_{3}&=& -\nabc_3 \mathcal{V}^{s=1}_{[3, 4]} - C_1  (\mathcal{V}^{s=1}_{[3, 4]} ) +\nabc_4(C_2)+ \nabc_4\nabc_3 C_1-2 (\eta-\etab ) \c \mathcal{V}^{s=1}_{[3,a]}\\
 I_{0}&=& \nabc_4\nabc_3C_2.
 \eeaa
 Using \eqref{expression-C-0-3-4}, we have
  \beaa
\mathcal{V}^{s=1}_{[3, 4]} &=&2\rho-2\eta\c\etab +i \left( - 4 \rhod  +4 \eta \wedge \etab \right), \\
\mathcal{V}^{s-1=0}_{[3, 4]} &=&i \left( - 4 \rhod  +4 \eta \wedge \etab \right)
\eeaa
and therefore we can write for $I_{33}$
 \bea\label{eq:expression-I33}
 I_{33}&=& -2\rho-2 \eta \c (\eta-2\etab )  +i \left(  8 \rhod  -8 \eta \wedge \etab \right)  +\nabc_4( C_1 ).
 \eea


\subsubsection{Expression for $J$}


We have
\beaa
J&=& \frac 1 4 [Q, \DDc\hot] (\DDbc \c A)+\frac 1 4 \DDc\hot ([ Q, \DDbc \c ]A).
\eeaa
Using Proposition \ref{commutation-general-Q}, applied to $F=\DDbc \c A$ of conformal type $s=2$, we obtain
\beaa
[Q, \DDc\hot] (\DDbc \c A)&=&2H \hot \nabc_3\nabc_3 \DDbc \c A-  \tr \Xb  \DDc\hot \nabc_3\DDbc \c A\\
&&+\frac 1 2 (\tr \Xb)(\tr\Xb-C_1)  \DDc\hot \DDbc \c A \\
 &&+\left(\nabc_3H+(- 3 \tr \Xb +C_1 )H -\DDc ( C_1 )\right) \hot \nabc_3 \DDbc \c A\\
 &&+\left( \frac 3 2  \tr \Xb \left(-\nabc_3H+(\tr\Xb-C_1)H\right)  -\DDc (C_2) \right)\hot \DDbc \c A\\
 &&+\nab_3 \big( r^{-2} \Ga_b  \c \frak{d}^{\leq 2} A \big).
\eeaa
We write
  \beaa
 \nabc_3  \ov{\DDc} \c A &=&\ov{\DDc} \c \nabc_3 A- \frac 1 2\ov{\tr\Xb}\, \left( \ov{\DDc} \c A \right) +\ov{H} \c \nabc_3 A+ r^{-1} \Ga_b \c \dk^{\leq 1} A, 
 \eeaa
 \beaa
 \DDc\hot \nabc_3\DDbc \c A&=& \nabc_3\DDc\hot \DDbc \c A-[\nabc_3, \DDc\hot]\DDbc \c A\\
 &=& \nabc_3\DDc\hot \DDbc \c A+ \frac 1 2 \tr \Xb \DDc\hot \DDbc \c A \\
 &&- H \hot \nabc_3 \DDbc \c A+  \frac 3 2 \tr \Xb H \hot \DDbc \c A + r^{-2} \Ga_b \c \dk^{\leq 2} A \\
  &=& \nabc_3\DDc\hot \DDbc \c A+ \frac 1 2 \tr \Xb \DDc\hot \DDbc \c A - H \hot (\ov{\DDc} \c \nabc_3 A)\\
  &&- H \hot (\ov{H} \c \nabc_3 A)+  \left(\frac 3 2 \tr \Xb +\frac 1 2\ov{\tr\Xb}\right)H \hot \DDbc \c A+ r^{-2} \Ga_b \c \dk^{\leq 2} A,   
  \eeaa
  and 
 \beaa
 \nabc_3\nabc_3 \DDbc \c A&=& \nabc_3 \DDbc \c\nabc_3 A+\nabc_3([\nabc_3,  \DDbc \c] A)\\
 &=& \DDbc \c\nabc_3 \nabc_3 A+[\nabc_3 ,\DDbc \c]\nabc_3 A+\nabc_3([\nabc_3,  \DDbc \c] A)\\
  &=& \DDbc \c\nabc_3 \nabc_3 A\\
  &&- \frac 1 2\ov{\tr\Xb}\, \left( \ov{\DDc} \c \nabc_3 A - \ov{H}\c \nabc_3 A\right) +\ov{H} \c \nabc_3 \nabc_3 A\\
  &&+\frac 1 2(\ov{\tr\Xb})^2\, \left( \ov{\DDc} \c A \right)- \frac 1 2\ov{\tr\Xb}\, \left(\ov{\DDc} \c  \nabc_3A  +  \ov{H}\c \nabc_3A\right)\\
 && +\nabc_3\ov{H} \c \nabc_3 A +\ov{H} \c \nabc_3\nabc_3 A+ r^{-1} \Ga_b \c \dk^{\leq 2} A\\
 &=& \DDbc \c\nabc_3 \nabc_3 A +2\ov{H} \c \nabc_3 \nabc_3 A- \ov{\tr\Xb}\,  \ov{\DDc} \c \nabc_3 A \\
  &&+\frac 1 2(\ov{\tr\Xb})^2\,  \ov{\DDc} \c A  +\nabc_3\ov{H} \c \nabc_3 A+ r^{-1} \Ga_b \c \dk^{\leq 2} A, 
 \eeaa
where we used the intermediate computations in Proposition \ref{commutation-general-Q}. Putting the above together we obtain 
 \beaa
&&[Q, \DDc\hot] (\DDbc \c A)\\
&=&-  \tr \Xb  \nabc_3\DDc\hot \DDbc \c A+2H \hot \DDbc \c\nabc_3 \nabc_3 A -\frac 1 2 (\tr \Xb) C_1 \DDc\hot \DDbc \c A\\
&&+4H \hot(\ov{H} \c \nabc_3 \nabc_3 A)+2H \hot (\nabc_3\ov{H} \c \nabc_3 A)\\
&&+\left(\nabc_3H+(-2\tr\Xb- 2\ov{\tr\Xb}+C_1)\, H-\DDc ( C_1 )\right) \hot (\ov{\DDc} \c \nabc_3 A) \\
  && +\left(\nabc_3H+(- 2 \tr \Xb +C_1 )H -\DDc ( C_1 )\right) \hot (\ov{H} \c \nabc_3 A )\\
 &&+\Bigg[ \left(\frac 3 2  \tr \Xb +\frac 1 2 \ov{\tr\Xb}\right)\left(-\nabc_3H-C_1 H\right)+\left(\tr\Xb \ov{\tr\Xb}+(\ov{\tr\Xb})^2\right)\,  H\\
 &&+\frac 1 2\ov{\tr\Xb}\DDc ( C_1 )  -\DDc (C_2) \Bigg]\hot \DDbc \c A+\nab_3 \big( r^{-2} \Ga_b  \c \frak{d}^{\leq 2} A \big).
\eeaa
Using Lemma \ref{SIMPLIFICATION-ANGULAR},
the above becomes
 \beaa
&&\,[Q, \DDc\hot] (\DDbc \c A)\\
&=&-  \tr \Xb  \nabc_3\DDc\hot \DDbc \c A+8H \c \nabc \nabc_3 \nabc_3 A \\
&&-\frac 1 2 (\tr \Xb) C_1 \DDc\hot \DDbc \c A+8 (H \c \ov{H} ) \nabc_3 \nabc_3 A\\
&&+4\left(\nabc_3H+(-2\tr\Xb- 2\ov{\tr\Xb}+C_1)\, H-\DDc ( C_1 )\right) \c \nabc  \nabc_3 A \\
  &&+2H \hot (\nabc_3\ov{H} \c \nabc_3 A) \\
 &&+\left(\nabc_3H+(- 2 \tr \Xb +C_1 )H -\DDc ( C_1 )\right) \hot (\ov{H} \c \nabc_3 A )\\
 &&+4\Bigg[ \left(\frac 3 2  \tr \Xb +\frac 1 2 \ov{\tr\Xb}\right)\left(-\nabc_3H-C_1 H\right)+\left(\tr\Xb \ov{\tr\Xb}+(\ov{\tr\Xb})^2\right)\,  H\\
 &&+\frac 1 2\ov{\tr\Xb}\DDc ( C_1 )  -\DDc (C_2) \Bigg]\c \nabc A+\nab_3 \big( r^{-2} \Ga_b  \c \frak{d}^{\leq 2} A \big).
\eeaa
Using Proposition \ref{commutation-general-Q}, we deduce
\beaa
&&\DDc\hot ([ Q, \DDbc \c ]A)\\
&=& \DDc\hot \Bigg[2\ov{H} \c \nabc_3\nabc_3 A- \ov{\tr\Xb}\,\ov{\DDc} \c  \nabc_3A\\
&&+\frac 1 2\ov{\tr\Xb} (\ov{\tr\Xb}-C_1)\,  \ov{\DDc} \c A +\left(\nabc_3\ov{H} +C_1  \ov{H}-\DDbc ( C_1 )\right)\c \nabc_3A\\
&&+\left(-\DDbc (C_2) \right)\c A+\nab_3 \big( r^{-2} \Ga_b  \c \frak{d}^{\leq 1} A \big)\Bigg],
\eeaa
which gives
\beaa
\DDc\hot ([ Q, \DDbc \c ]A)&=& 2\DDc\hot (\ov{H} \c \nabc_3\nabc_3 A)- \ov{\tr\Xb}\,\DDc\hot\ov{\DDc} \c  \nabc_3A\\
&&- \DDc\ov{\tr\Xb}\,\hot \ov{\DDc} \c  \nabc_3A+\frac 1 2\ov{\tr\Xb} (\ov{\tr\Xb}-C_1)\, \DDc\hot \ov{\DDc} \c A\\
&&+\DDc(\frac 1 2\ov{\tr\Xb} (\ov{\tr\Xb}-C_1))\hot\,  \ov{\DDc} \c A \\
&&+\DDc\hot \big[\left(\nabc_3\ov{H} +C_1  \ov{H}-\DDbc ( C_1 )\right)\c \nabc_3A \big]\\
&&-\DDc\hot\left(\DDbc (C_2) \c A\right)+\nab_3 \big( r^{-2} \Ga_b  \c \frak{d}^{\leq 2} A \big).
\eeaa
Writing
\beaa
\DDc\hot\ov{\DDc} \c  \nabc_3A&=&\DDc\hot\nabc_3\ov{\DDc} \c  A-\DDc\hot[\nabc_3, \ov{\DDc} \c ]A \\
&=&\nabc_3\DDc\hot \DDbc \c A+ \frac 1 2 (\tr \Xb +\ov{\tr\Xb})\DDc\hot \DDbc \c A \\
&&- H \hot (\ov{\DDc} \c \nabc_3 A)- H \hot (\ov{H} \c \nabc_3 A)\\
&&+  \left(\frac 3 2 \tr \Xb +\frac 1 2\ov{\tr\Xb}\right)H \hot \DDbc \c A  +\frac 1 2\DDc (\ov{\tr\Xb})\hot(  \ov{\DDc} \c A )\\
&& -\DDc \hot (\ov{H} \c \nabc_3 A)+\nab_3 \big( r^{-2} \Ga_b  \c \frak{d}^{\leq 2} A \big),
\eeaa
we have 
\beaa
\DDc\hot ([ Q, \DDbc \c ]A)&=&- \ov{\tr\Xb}\,\nabc_3\DDc\hot \DDbc \c A +\frac 1 2\ov{\tr\Xb} (-\tr\Xb-C_1)\, \DDc\hot \ov{\DDc} \c A\\
&& + 2\DDc\hot (\ov{H} \c \nabc_3\nabc_3 A)+\ov{\tr\Xb}\DDc \hot (\ov{H} \c \nabc_3 A)\\
&&+\ov{\tr\Xb} H \hot (\ov{H} \c \nabc_3 A)+\left(\ov{\tr\Xb} H - \DDc\ov{\tr\Xb}\right)\,\hot (\ov{\DDc} \c  \nabc_3A)\\
&+ & \Big[\left(-\frac 3 2 \tr \Xb \ov{\tr\Xb}-\frac 1 2\ov{\tr\Xb}^2\right)H +\frac 1 2\ov{\tr\Xb}\DDc (\ov{\tr\Xb})\\
&&-\frac 1 2 \DDc(\ov{\tr\Xb} C_1)\Big]\hot \DDbc \c A \\
&&+\DDc\hot((\nabc_3\ov{H} +C_1  \ov{H}-\DDbc ( C_1 ))\c \nabc_3A)\\
&&-\DDc\hot(\DDbc (C_2) \c A)+\nab_3 \big( r^{-2} \Ga_b  \c \frak{d}^{\leq 2} A \big).
\eeaa
Using Lemma \ref{SIMPLIFICATION-ANGULAR}, the above simplifies to
\beaa
&&\DDc\hot ([ Q, \DDbc \c ]A)\\
&=&- \ov{\tr\Xb}\,\nabc_3\DDc\hot \DDbc \c A +  8 \ov{H} \c \nabc \nabc_3\nabc_3 A\\
&&+\frac 1 2\ov{\tr\Xb} (-\tr\Xb-C_1)\, \DDc\hot \ov{\DDc} \c A+ 4(\DDc\c\ov{H} )\nabc_3\nabc_3 A\\
&&+2\left(4\ov{\tr\Xb}\eta  - 2\DDc\ov{\tr\Xb}+2(\nabc_3\ov{H} +C_1  \ov{H}-\DDbc ( C_1 ))\right)\c \nabc\nabc_3A\\
&&+2\left(\ov{\tr\Xb}(\DDc \c\ov{H}+ H \c\ov{H}+\DDc\c(\nabc_3\ov{H} +C_1  \ov{H}-\DDbc ( C_1 )) \right)\nabc_3 A\\
&&+ 4\Big[\left(-\frac 3 2 \tr \Xb \ov{\tr\Xb}-\frac 1 2\ov{\tr\Xb}^2\right)H +\frac 1 2\ov{\tr\Xb}\DDc (\ov{\tr\Xb})\\
&&-\frac 1 2 \DDc(\ov{\tr\Xb} C_1)-\DDbc (C_2) \Big]\c \nabc  A-2(\DDc \c\DDbc (C_2))  A+\nab_3 \big( r^{-2} \Ga_b  \c \frak{d}^{\leq 2} A \big).
\eeaa
Putting the above together we finally obtain 
\beaa
J&=&-  \left(\tr \Xb+\ov{\tr \Xb}\right)  \nabc_3\left(\frac 1 4\DDc\hot \DDbc \c A\right)\\
&&+\left(- \frac 1 2\ov{\tr\Xb} \tr\Xb-\frac 1 2 (\tr \Xb+\ov{\tr\Xb})C_1\right)\, \left(\frac 1 4 \DDc\hot \ov{\DDc} \c A\right)\\
&&+4 \eta \c \nabc \nabc_3 \nabc_3 A +\tilde{J}_{33} \nabc_3 \nabc_3 A+ \tilde{J}_{a3}\c  \nabc \nabc_3 A\\
&&+\tilde{J}_{3} \nabc_3 A +\tilde{J}_a \c \nabc A+\tilde{J}_0  A+\nab_3 \big( r^{-2} \Ga_b  \c \frak{d}^{\leq 2} A \big),
\eeaa
where
\beaa
\tilde{J}_{33}&=&\DDc\c\ov{H}+2H \c\ov{H}, \\
\tilde{J}_{a3}&=&\nabc_3H+(-2\tr\Xb- 2\ov{\tr\Xb}+\C)\, H-\DDc ( C_1 )\\
&&+2\ov{\tr\Xb}\eta  - \DDc\ov{\tr\Xb}+\nabc_3\ov{H} +C_1  \ov{H}-\DDbc ( C_1 ),\\
\tilde{J}_3&=&\frac 1 2 \ov{\tr\Xb}(\DDc \c\ov{H}+ H \c\ov{H})+\frac 1 2 \DDc\c(\nabc_3\ov{H} +C_1  \ov{H}-\DDbc ( C_1 ))\\
&&+H \hot (\nabc_3\ov{H} \c ) +\frac 1 2 \left(\nabc_3H+(- 2 \tr \Xb +C_1 )H -\DDc ( C_1 )\right) \hot (\ov{H} \c  ),\\
 \tilde{J}_a&=& \left(\frac 3 2  \tr \Xb +\frac 1 2 \ov{\tr\Xb}\right)\left(-\nabc_3H-C_1 H\right)+\left(\tr\Xb \ov{\tr\Xb}+(\ov{\tr\Xb})^2\right)\,  H\\
 &&+\frac 1 2\ov{\tr\Xb}\DDc ( C_1 )  -\DDc (C_2)+\left(-\frac 3 2 \tr \Xb \ov{\tr\Xb}-\frac 1 2\ov{\tr\Xb}^2\right)H \\
 &&+\frac 1 2\ov{\tr\Xb}\DDc (\ov{\tr\Xb})-\frac 1 2 \DDc(\ov{\tr\Xb} C_1)-\DDbc (C_2), \\
\tilde{J}_0&=&-\frac 1 2 (\DDc \c\DDbc (C_2)).
  \eeaa
We can simplify $\tilde{J}_{a3}$ as
  \beaa
  \tilde{J}_{a3}&=&\nabc_3H+(-2\tr\Xb- 2\ov{\tr\Xb}+C_1)\, H-\DDc ( C_1 )\\
  &&+2\ov{\tr\Xb}\eta  - \DDc\ov{\tr\Xb}+\nabc_3\ov{H} +C_1  \ov{H}-\DDbc ( C_1 )\\
  &=&2\nabc_3\eta-2 \nabc ( C_1 ) +2C_1  \eta- \DDc\ov{\tr\Xb}+(-2\tr\Xb- \ov{\tr\Xb})\, H+\ov{\tr\Xb}\ov{H}  \\
   &=&2\nabc_3\eta-2 \nabc ( C_1) +2C_1  \eta- (\tr\Xb-\ov{\tr\Xb})\Hb+(-2\tr\Xb- \ov{\tr\Xb})\, H+\ov{\tr\Xb}\ov{H} .
  \eeaa

Recalling the definition \eqref{Teukolsky-operator-ch5} of $\LL(A)$, we write 
\beaa
 \frac{1}{4}\DDc\hot (\DDbc \c A) &=&\nabc_4\nabc_3A+\left( \frac 1 2 \tr X +2\ov{\tr X} \right)\nabc_3A+\frac{1}{2}\tr\Xb \nabc_4A\\
&&-\left( 4H+\Hb +\ov{\Hb} \right)\c \nabc A+ \left(\ov{\tr X} \tr \Xb -2\ov{P}-4\eta \c \etab +4i \eta \wedge \etab\right) A
\eeaa
which gives
\beaa
&&\nabc_3\left( \frac{1}{4}\DDc\hot (\DDbc \c A) \right)\\
&=&\nabc_3\nabc_4\nabc_3A+\left( \frac 1 2 \tr X +2\ov{\tr X} \right)\nabc_3\nabc_3A\\
&&+\nabc_3\left( \frac 1 2 \tr X +2\ov{\tr X} \right)\nabc_3A+\frac{1}{2}\tr\Xb \nabc_3\nabc_4A+\frac{1}{2}\nabc_3\tr\Xb \nabc_4A\\
&&-\left( 4H+\Hb +\ov{\Hb} \right)\c \nabc_3\nabc A-\nabc_3\left( 4H+\Hb +\ov{\Hb} \right)\c \nabc A\\
&&+ \left(\ov{\tr X} \tr \Xb -2\ov{P}-4\eta \c \etab +4i \eta \wedge \etab\right) \nabc_3A\\
&&+ \nabc_3\left(\ov{\tr X} \tr \Xb -2\ov{P}-4\eta \c \etab +4i \eta \wedge \etab\right) A.
\eeaa
Writing 
\beaa
\nabc_3 \nabc_4 A&=& \nabc_4 \nabc_3 A+2(\eta-\etab ) \c \nabc A +(\mathcal{V}^{s=2}_{[3,4]} )A,\\
\nabc_3 \nabc_4 \nabc_3 A&=& \nabc_4 \nabc_3 \nabc_3 A+2(\eta-\etab ) \c \nabc \nabc_3A +(\mathcal{V}^{s=1}_{[3,4]} ) \nabc_3 A,\\
\,  \nabc_3\nabc A    &=&\nabc \nabc_3 A- \frac  1 2   \trchb\, \nabc A-\frac 1 2 \atrchb\, \dual \nabc  A+\eta \nabc_3 A+\mathcal{V}^{s=2}_{[3,a]} (A),
\eeaa
we have 
\beaa
\nabc_3\left( \frac{1}{4}\DDc\hot (\DDbc \c A) \right)&=&\nabc_4 \nabc_3 \nabc_3 A+\frac{1}{2}\tr\Xb \nabc_4 \nabc_3 A-\frac{1}{4}\tr\Xb^2 \nabc_4A\\
&&+\hat{J}_{a3}\c \nabc \nabc_3A+\hat{J}_{33}\nabc_3\nabc_3A+\hat{J}_3 \nabc_3A+\hat{J}^a_3(A)\\
&&+\hat{J}_a\c \nabc A  +\hat{J}_{\dual a }\c  \dual \nabc  A+ \hat{J}_0A+\hat{J}^a_0(A)
\eeaa
where
\beaa
\hat{J}_{a3}&=& 2(\eta-\etab ) -\left( 4H+\Hb +\ov{\Hb} \right),\\
\hat{J}_{33}&=& \frac 1 2 \tr X +2\ov{\tr X}, \\
\hat{J}_3&=&\nabc_3\left( \frac 1 2 \tr X +2\ov{\tr X} \right)+\ov{\tr X} \tr \Xb -P+\ov{P}  -2\eta\c\etab-\left( 4H+\Hb +\ov{\Hb} \right)\c \eta, \\
\hat{J}^a_3(A)&=&  - 4  \eta \hot (\etab \c \nabc_3A)+4 \etab \hot (\eta \c \nabc_3A)-  2 H   \hot (\ov{\Hb} \c \nabc_3A),\\
\hat{J}_a&=&\tr\Xb (\eta-\etab )-\nabc_3\left( 4H+\Hb +\ov{\Hb} \right) +\frac  1 2   \trchb\left( 4H+\Hb +\ov{\Hb} \right),\\
\hat{J}_{\dual a}&=&\frac 1 2 \atrchb\left( 4H+\Hb +\ov{\Hb} \right),\\
\hat{J}_0&=&\nabc_3\left(\ov{\tr X} \tr \Xb -2\ov{P}\right) +\frac{1}{2}\tr\Xb(4\ov{P}  -4\eta\c\etab)\\
&&+(\trchb \left( 4H+\Hb +\ov{\Hb} \right) -\atrchb \dual \left( 4H+\Hb +\ov{\Hb} \right) )\c \eta,\\
\hat{J}^a_0(A)&=&\frac  1 2   \trchb\, \Big(\eta \hot (\left( 4H+\Hb +\ov{\Hb} \right) \c A)-\left( 4H+\Hb +\ov{\Hb} \right) \hot (\eta \c A) \Big)\\
&&+\frac 1 2 \atrchb\, \Big(-\eta \hot (\dual \left( 4H+\Hb +\ov{\Hb} \right) \c A)+\dual  \left( 4H+\Hb +\ov{\Hb} \right) \hot(\eta \c  A) \Big)\\
&&-  2\nabc_3H   \hot (\ov{\Hb} \c A)-  2H   \hot (\nabc_3 \ov{\Hb} \c A)+\frac{1}{2}\tr\Xb ( - 4  \eta \hot (\etab \c A)+4 \etab \hot (\eta \c A)).
\eeaa
We therefore finally obtain
\beaa
J&=&-  \left(\tr \Xb+\ov{\tr \Xb}\right)  \nabc_4 \nabc_3 \nabc_3 A+4 \eta \c \nabc \nabc_3 \nabc_3 A\\
&&+J_{43} \nabc_4 \nabc_3 A +J_4 \nabc_4 A +J_{a3}\c \nabc \nabc_3A+J_{33} \nabc_3 \nabc_3 A\\
&&+J_3 \nabc_3A+J^a_3(A)+J_a\c \nabc A +J_{\dual a }\c  \dual \nabc  A+J_0 A+J^a_0(A)\\
&&+\nab_3 \big( r^{-2} \Ga_b  \c \frak{d}^{\leq 2} A \big),
\eeaa
where
\beaa
J_{43}&=& - \frac 1 2\ov{\tr\Xb} \tr\Xb-\frac {C_1 +\tr\Xb}{2} (\tr \Xb+\ov{\tr\Xb}),\\
J_4&=&\frac{1}{4}(\tr\Xb)^2 \left(\tr \Xb\right) -\frac 1 4\tr\Xb(\tr \Xb+\ov{\tr\Xb})C_1,
\eeaa
and
\beaa
J_{a3}&=&-  \left(\tr \Xb+\ov{\tr \Xb}\right)  \hat{J}_{a3} +\tilde{J}_{a3},\\
J_{33}&=& -  \left(\tr \Xb+\ov{\tr \Xb}\right)  \hat{J}_{33}+\tilde{J}_{33},\\
J_3&=& -  \left(\tr \Xb+\ov{\tr \Xb}\right)  \hat{J}_3 +\left(- \frac 1 2\ov{\tr\Xb} \tr\Xb-\frac 1 2 C_1(\tr \Xb+\ov{\tr\Xb})\right)\, \left( \frac 1 2 \tr X +2\ov{\tr X} \right)+\tilde{J}_{3},\\
J^a_3(A)&=&-  \left(\tr \Xb+\ov{\tr \Xb}\right)  (\hat{J}^a_3(A)+ \tilde{J}^a_3(A),\\
J_a&=&-  \left(\tr \Xb+\ov{\tr \Xb}\right)  \hat{J}_a -\left(- \frac 1 2\ov{\tr\Xb} \tr\Xb-\frac 1 2 C_1 (\tr \Xb+\ov{\tr\Xb})\right)\, \left( 4H+\Hb +\ov{\Hb} \right)+\tilde{J}_a, \\
J_{\dual a}&=& -  \left(\tr \Xb+\ov{\tr \Xb}\right)  \hat{J}_{\dual a },\\
J_0&=& -  \left(\tr \Xb+\ov{\tr \Xb}\right)  \hat{J}_0+\left(- \frac 1 2\ov{\tr\Xb} \tr\Xb-\frac 1 2 C_1 (\tr \Xb+\ov{\tr\Xb})\right)\,  \left(\ov{\tr X} \tr \Xb -2\ov{P}\right)+\tilde{J}_0,\\
J^a_0(A)&=&-  \left(\tr \Xb+\ov{\tr \Xb}\right)  \hat{J}^a_0(A)-\left(- \frac 1 2\ov{\tr\Xb} \tr\Xb-\frac 1 2 C_1(\tr \Xb+\ov{\tr\Xb})\right)\, (   2H   \hot (\ov{\Hb} \c A)).
\eeaa
Using the above computation for $\tilde{J}_{a3}$, we simplify
\beaa
J_{a3}&=&-  2\trchb (2(\eta-\etab ) -\left( 4(\eta+ i \dual \eta )+2\etab \right)) \\
&&+2\nabc_3\eta-4\trchb\, H-2\nabc ( C_1 )+2\ov{\tr\Xb}\eta  - \DDc\ov{\tr\Xb} +2C_1  \eta\\
&=&-  4\trchb (\eta-\etab )+  2\trchb \left( 4(\eta+ i \dual \eta )+2\etab \right) \\
&&+2\nabc_3\eta-4\trchb\, (\eta+ i \dual \eta)-2\nabc ( C_1 )\\
&&+2(\trchb + i \atrchb)\eta  +2i\atrchb (\etab+ i \dual \etab) +2C_1  \eta\\
&=&2\nabc_3\eta-2\nabc ( C_1 ) +2C_1  \eta+  \trchb (8\etab+2\eta) -2\atrchb   \dual \etab \\
&& +  i  \left(4 \trchb \dual \eta+2\atrchb \eta +2\atrchb \etab \right) 
\eeaa
Also, we simplify
\bea\label{eq:J33-final}
J_{33}&=& -  2\trchb   \left(\frac 1 2 \tr X +2\ov{\tr X} \right)+\DDc\c\ov{H}+2H \c\ov{H}.
\eea


\subsubsection{Expression for $K$}

Observe that
 \bea\label{Q-f-g}
 \begin{split}
 Q(f U) &= Q(f) U+ fQ(U)+2\nabc_3f \nabc_3U  -C_2  f U \\
 &= fQ(U)+2\nabc_3f \nabc_3U +\left( \nabc_3\nabc_3 f + C_1\nabc_3 f\right) U.
 \end{split}
 \eea 
Using Proposition \ref{commutation-general-Q}, we obtain
 \beaa
&&[Q, \FF \ \nabc_3]A\\
&=& \left( \nabc_3\nabc_3 \FF + C_1\nabc_3 \FF\right) \nabc_3A+ \FF [Q, \nabc_3]A+2\nabc_3\FF  \nabc_3\nabc_3A \\
&=& \left(2\nabc_3\FF  \right)\nabc_3\nabc_3A +\left(\nabc_3\nabc_3 \FF + C_1\nabc_3 \FF- \FF \nabc_3 C_1\right) \nabc_3A\\
&&+  ( - \FF \nabc_3C_2 )   A. 
\eeaa
In particular, 
\beaa
K&=& [Q, \FF \ \nab_3^{(c)}]A \qquad \text{for} \qquad \FF= - \frac 1 2 \tr X -2\ov{\tr X}.
\eeaa
We therefore obtain
\beaa
K&=& K_{33} \nabc_3 \nabc_3 A+K_3 \nabc_3 A +K_0 A,
\eeaa
where
\bea
K_{33}&=& 2\nabc_3\left(- \frac 1 2 \tr X -2\ov{\tr X} \right), \label{eq:expression-K33}\\
K_3&=& \nabc_3\nabc_3\left(- \frac 1 2 \tr X -2\ov{\tr X} \right) + C_1\nabc_3 \left(- \frac 1 2 \tr X -2\ov{\tr X} \right) \nn\\
&& - \left(- \frac 1 2 \tr X -2\ov{\tr X} \right) \nabc_3 C_1, \nn \\
K_0&=& -\left(- \frac 1 2 \tr X -2\ov{\tr X} \right)\nabc_3C_2. \nn
\eea


\subsubsection{Expression for $L$}


Using \eqref{Q-f-g}, we obtain for a scalar $\EE$, 
 \beaa
&&[Q, \EE \ \nabc_4]A\\
&=& \left( \nabc_3\nabc_3 \EE + C_1\nabc_3 \EE\right)\nabc_4A+ \EE [Q, \nabc_4]A+2\nabc_3\EE  \nabc_3\nabc_4A\\
&=& \left( \nabc_3\nabc_3 \EE + C_1\nabc_3 \EE\right)\nabc_4A+ \EE [Q, \nabc_4]A+2\nabc_3\EE  \nabc_4\nabc_3A\\
&&+2\nabc_3\EE[\nabc_3, \nabc_4]A.
\eeaa
Using Proposition \ref{commutation-general-Q} applied to $U=A$ of conformal type $s=2$, we obtain
 \beaa
[Q, \EE \ \nabc_4]A&=& 2\nabc_3\EE  \nabc_4\nabc_3A+\left( \nabc_3\nabc_3 \EE + C_1\nabc_3 \EE\right)\nabc_4A\\
&&+ \EE \Big[4 (\eta-\etab ) \c \nabc \nabc_3 A\\
 && +\Big( 2\nabc_3(\eta-\etab )+(2C_1 -   \trchb ) (\eta-\etab )+\atrchb  \dual (\eta-\etab )\Big) \c \nabc A\\
 &&+\Big(\mathcal{V}^{s=2}_{[3, 4]} +\mathcal{V}^{s-1=1}_{[3, 4]}  +2 \eta \c (\eta-\etab ) -\nabc_4( C_1 )\Big) \nabc_3 A\\
 &&+\Big(\nabc_3 \mathcal{V}^{s=2}_{[3, 4]} + C_1  (\mathcal{V}^{s=2}_{[3, 4]} ) -\nabc_4(C_2) \Big) A +2 (\eta-\etab ) \c \mathcal{V}^{s=2}_{[3,a]}(A)\\
 &&+ r^{-3} \Ga_b  \c \frak{d}^{\leq 1} U+\dk\big( (\Ga_b \c \Ga_g) U \big) \Big]\\
&&+2\nabc_3\EE \Big(2(\eta-\etab ) \c \nabc A +(\mathcal{V}^{s=2}_{[3, 4]} ) A \Big).
\eeaa
In particular, 
\beaa
L&=& [Q, \EE \ \nabc_4]A \qquad \text{for} \qquad \EE= -\frac{1}{2}\tr\Xb.
\eeaa
We therefore obtain
 \beaa
L&=& L_{43}  \nabc_4\nabc_3A+L_{a3} \c \nabc\nabc_3A+ L_4 \nabc_4A+L_a \c \nabc A \\
&&+ r^{-4} \Ga_b  \c \frak{d}^{\leq 1} A+r^{-1}\dk\big( (\Ga_b \c \Ga_g) A \big),
\eeaa
where
\beaa
L_{43}&=& 2\nabc_3\EE=\frac 1 2 (\tr\Xb)^2\\
L_{a3}&=& 4 \EE(\eta-\etab )=-2\tr\Xb(\eta-\etab ),\\
L_4&=& \nabc_3\nabc_3 \EE + C_1\nabc_3 \EE =-\frac 1 4 (\tr\Xb)^3+\frac 1 4 (\tr\Xb)^2C_1,\\
L_{a}&=& \EE \Big( 2\nabc_3(\eta-\etab )+(2C_1 -   \trchb ) (\eta-\etab )+\atrchb  \dual (\eta-\etab )\Big)+ 4\nabc_3\EE (\eta-\etab )\\
L_{3}&=&\EE \Big(\mathcal{V}^{s=2}_{[3, 4]} +\mathcal{V}^{s-1=1}_{[3, 4]}  +2 \eta \c (\eta-\etab ) -\nabc_4( C_1 )\Big)\\
L_0&=& \EE\Big(\nabc_3 \mathcal{V}^{s=2}_{[3, 4]} + C_1  (\mathcal{V}^{s=2}_{[3, 4]} ) -\nabc_4(C_2) \Big) +2\nabc_3\EE (\mathcal{V}^{s=2}_{[3, 4]} ) +2\EE (\eta-\etab ) \c \mathcal{V}^{s=2}_{[3,a]}.
\eeaa
Using \eqref{expression-C-0-3-4}, we compute
\beaa
L_{3}&=& -\frac{1}{2}\tr\Xb \Big(6\left(\rho   -\eta\c\etab\right) +8i \left(- \rhod+ \eta \wedge \etab \right) +2 \eta \c (\eta-\etab ) -\nabc_4( C_1 )\Big),
\eeaa
and
\beaa
L_0&=& -\frac{1}{2}\tr\Xb\Big[\nabc_3\big(4\left(\rho   -\eta\c\etab\right) +4i \left(- \rhod+ \eta \wedge \etab \right)\big) \\
&&+ C_1  (4\left(\rho   -\eta\c\etab\right) +4i \left(- \rhod+ \eta \wedge \etab \right)) -\nabc_4(C_2) \Big] +\frac 1 2 (\tr\Xb)^2 (4\left(\rho   -\eta\c\etab\right)\\
&& +4i \left(- \rhod+ \eta \wedge \etab \right) ) +2\EE (\eta-\etab ) \c \mathcal{V}^{s=2}_{[3,a]}.
\eeaa


\subsubsection{Expression for $M$}


Observe that
 \beaa
 Q(F \c U) &= &F \c Q(U)+2\nabc_3F \c \nabc_3U +\left( \nabc_3\nabc_3 F + C_1\nabc_3 F\right)\c  U.
 \eeaa
We therefore obtain
 \beaa
 M&=& \left( 4H+\Hb +\ov{\Hb} \right) \c [Q, \nabc ]A+ 2 \nabc_3 \left( 4H+\Hb +\ov{\Hb} \right) \c \nabc_3 \nabc A\\
 &&+\left( \nabc_3\nabc_3 (4H+\Hb +\ov{\Hb} ) + C_1\nabc_3 (4H+\Hb +\ov{\Hb} )\right)\c  \nabc A\\
 &=& \left( 4H+\Hb +\ov{\Hb} \right) \c [Q, \nabc ]A+ 2 \nabc_3 \left( 4H+\Hb +\ov{\Hb} \right) \c \nabc \nabc_3 A\\
 &&+ 2 \nabc_3 \left( 4H+\Hb +\ov{\Hb} \right) \c [\nabc_3, \nabc] A\\
 &&+\left( \nabc_3\nabc_3 (4H+\Hb +\ov{\Hb} ) + C_1\nabc_3 (4H+\Hb +\ov{\Hb} )\right)\c  \nabc A.
 \eeaa
  Using Proposition \ref{commutation-general-Q} applied to $U=A$ of conformal type $s=2$, we obtain
  \beaa
 M &=& \left( 4H+\Hb +\ov{\Hb} \right)_a  \Bigg[-   \trchb\,  \nabc_a\nabc_3 A_{bc}- \atrchb\,  \dual \nabc_a\nabc_3 A_{bc}\\
 &&+2\eta_a \nabc_3\nabc_3 A_{bc}+\left(-\frac  1 2   \trchb (C_1-\trchb)-\frac 1 2 \atrchb^2 \right) \, \nabc_a A_{bc}\\
 &&-\frac 1 2 \atrchb \big(C_1-2\trchb\big) \, \dual \nabc_a  A_{bc}\\
&&+\left(\nabc_3\eta_a+(C_1-\frac  1 2   \trchb )\eta_a-\frac 1 2 \atrchb \dual  \eta_a -\nabc_a  C_1  \right) \nabc_3 A_{bc} \\
&&+\nabc_3(\mathcal{V}^{s=2}_{[3,a]}(A))+\mathcal{V}^{s-1=1}_{[3,a]}(\nabc_3 A)-\frac  1 2   \trchb\mathcal{V}^{s=2}_{[3,a]}(A)-\frac 1 2 \atrchb \dual \mathcal{V}^{s=2}_{[3,a]}(A)\\
&&   -\nabc_a (C_2)  A_{bc}+ C_1\mathcal{V}^{s=2}_{[3,a]}(A) +r^{-1}  \dk (\Ga_b  \c \frak{d}^{\leq 1} A) \Bigg]\\
&&+ 2 \nabc_3 \left( 4H+\Hb +\ov{\Hb} \right) \c \nabc \nabc_3 A\\
 &&+ 2 \nabc_3 \left( 4H+\Hb +\ov{\Hb} \right)_a \Big[ -\frac  1 2   \trchb\, \nabc_a A_{bc}-\frac 1 2 \atrchb\, \dual \nabc_a  A_{bc}+\eta_a \nabc_3 A_{bc}\\
 &&+\mathcal{V}^{s=2}_{[3,a]}(A)+r^{-1} \Ga_b \c  \frak{d}^{\leq 1} A\Big]\\
 &&+\left( \nabc_3\nabc_3 (4H+\Hb +\ov{\Hb} ) + C_1\nabc_3 (4H+\Hb +\ov{\Hb} )\right)\c  \nabc A,
 \eeaa
 which gives
 \beaa
  M &=&M_{a3} \c \nabc \nabc_3 A+M_{33} \nabc_3\nabc_3 A+M_{a} \c \nabc A+M_{3} \nabc_3 A+M_0 A\\
  &&+r^{-3}  \dk (\Ga_b  \c \frak{d}^{\leq 1} A),
 \eeaa
 where
 \bea
 M_{a3}&=& 2 \nabc_3 \left( 4H+\Hb +\ov{\Hb} \right)-   \trchb\left( 4H+\Hb +\ov{\Hb} \right)+\atrchb \dual  \left( 4H+\Hb +\ov{\Hb} \right) \nn\\
 M_{33}&=& 2\eta \c \left( 4H+\Hb +\ov{\Hb} \right) \label{eq:expression-M33}\\
 M_{a}&=& -\trchb  \nabc_3 \left( 4H+\Hb +\ov{\Hb} \right)+\atrchb \dual \nabc_3 \left( 4H+\Hb +\ov{\Hb} \right) \nn\\
 &&+\left(-\frac  1 2   \trchb (C_1-\trchb)-\frac 1 2 \atrchb^2 \right)\left( 4H+\Hb +\ov{\Hb} \right) \nn\\
 &&+\frac 1 2 \atrchb \big(C_1-2\trchb\big)\dual  \left( 4H+\Hb +\ov{\Hb} \right)\nn\\
 &&+ \nabc_3\nabc_3 (4H+\Hb +\ov{\Hb} ) + C_1\nabc_3 (4H+\Hb +\ov{\Hb} )\nn \\
 M_{3}&=&  \left( 4H+\Hb +\ov{\Hb} \right)_a  \Big(\left(\nabc_3\eta_a+(C_1-\frac  1 2   \trchb )\eta_a-\frac 1 2 \atrchb \dual  \eta_a -\nabc_a  C_1  \right) \nn \\
 &&+\nabc_3(\mathcal{V}^{s=2}_{[3,a]})+\mathcal{V}^{s-1=1}_{[3,a]} \Big)+ 2 \nabc_3 \left( 4H+\Hb +\ov{\Hb} \right)_a \Big[ \eta_a \nabc_3 A_{bc}\Big] \nn\\
 M_0&=&  \left( 4H+\Hb +\ov{\Hb} \right)_a  \Bigg[-\frac  1 2   \trchb\mathcal{V}^{s=2}_{[3,a]}(A)-\frac 1 2 \atrchb \dual \mathcal{V}^{s=2}_{[3,a]}(A) \nn\\
 && -\nabc_a (C_2)  A_{bc}+ C_1\mathcal{V}^{s=2}_{[3,a]}(A)  \Bigg]+ 2 \nabc_3 \left( 4H+\Hb +\ov{\Hb} \right)_a \Big[ \mathcal{V}^{s=2}_{[3,a]}(A)\Big]. \nn
 \eea

  
 \subsubsection{Expression  for $N$}
 
 
 Using \eqref{eq:alternative-HovHb} we write
 \beaa
 N&=& [Q, \left(-\ov{\tr X} \tr \Xb +2\ov{P}\right)] A+ [Q, H   \hot \ov{\Hb} \c ]A\\
 &=& [Q, \left(-\ov{\tr X} \tr \Xb +2\ov{P}+4\eta \c \etab -4i \eta \wedge \etab\right)] A.
 \eeaa
 Using \eqref{Q-f-g}, we obtain
 \beaa
 N&=& N_3 \nabc_3A +N_0 A,
 \eeaa
 where
 \beaa
 N_3&=& 2\nabc_3\left(-\ov{\tr X} \tr \Xb +2\ov{P}+4\eta \c \etab -4i \eta \wedge \etab\right)\\
 N_0&=& \nabc_3\nabc_3 \left(-\ov{\tr X} \tr \Xb +2\ov{P}+4\eta \c \etab -4i \eta \wedge \etab\right)\\
 &&+ C_1\nabc_3 \left(-\ov{\tr X} \tr \Xb +2\ov{P}+4\eta \c \etab -4i \eta \wedge \etab\right).
 \eeaa

  
\subsubsection{The sum}


 Putting the above expressions together we obtain
 \beaa
[Q, \LL]A&=& I + J+K+L+M+N\\
&=&  -4(\eta-\etab ) \c \nabc\nabc_3\nabc_3A+4 \eta \c \nabc \nabc_3 \nabc_3 A \\
&&-  \left(\tr \Xb+\ov{\tr \Xb}\right)  \nabc_4 \nabc_3 \nabc_3 A\\
&&+\left( I_{43}+J_{43}+L_{43}\right) \  \nabc_4\nabc_3A+\left(I_4+J_4+L_4 \right) \  \nabc_4A\\
 &&+\left(I_{33}+J_{33}+K_{33}+M_{33}\right) \nabc_3 \nabc_3A \\
 &&+\left(I_{a3}+J_{a3}+L_{a3}+M_{a3}\right)\c \nabc \nabc_3A \\
 &&+ \left(I_3+J_3 +K_3+L_3+N_3\right)\nabc_3A+J^a_3(A)\\
 &&+\left(J_a+L_a+M_a\right)\c \nabc A +\left(J_{\dual a } \right)\c  \dual \nabc  A\\
 &&+\left(I_0+J_0+K_0+L_0+M_0+N_0\right) \  A+J^a_0(A) \\
 &&+\nab_3 \big( r^{-2} \Ga_b  \c \frak{d}^{\leq 2} A \big)+r^{-3}\dk^{\leq2}(  \Ga_b \c \frak{d}^{\leq 1} A)+r^{-1}\dk\big( (\Ga_b \c \Ga_g) A \big),
 \eeaa
 which gives
 \beaa
 [Q, \LL]A&=&  4\etab  \c \nabc\nabc_3\nabc_3A -  \left(\tr \Xb+\ov{\tr \Xb}\right)  \nabc_4 \nabc_3 \nabc_3 A\\
&&+\left( I_{43}+J_{43}+L_{43}\right) \  \nabc_4\nabc_3A+\left(I_4+J_4+L_4 \right) \  \nabc_4A\\
 &&+\left(I_{33}+J_{33}+K_{33}+M_{33}\right) \nabc_3 \nabc_3A \\
 &&+\left(I_{a3}+J_{a3}+L_{a3}+M_{a3}\right)\c \nabc \nabc_3A \\
 &&+ \left(I_3+J_3 +K_3+L_3+N_3\right)\nabc_3A+J^a_3(A)\\
 &&+\left(J_a+L_a+M_a\right)\c \nabc A +\left(J_{\dual a } \right)\c  \dual \nabc  A\\
 &&+\left(I_0+J_0+K_0+L_0+M_0+N_0\right) \  A+J^a_0(A)+\err,
\eeaa
with $\err=\nab_3 \big( r^{-2} \Ga_b  \c \frak{d}^{\leq 2} A \big)+r^{-3}\dk^{\leq2}(  \Ga_b \c \frak{d}^{\leq 1} A)+r^{-1}\dk\big( (\Ga_b \c \Ga_g) A \big)$.
Recalling the definition of $Q(A)$, we write
  \beaa
\nabc_3\nabc_3 A &=& Q(A)- C_1  \nabc_3A -C_2  A,
 \eeaa
and therefore
 \beaa
 \nabc \nabc_3 \nabc_3 A &=& \nabc Q(A)- C_1  \nabc \nabc_3A- (\nabc C_1)  \nabc_3A \\
 &&-C_2 \nabc A-(\nabc C_2) A,\\
 \nabc_4\nabc_3\nabc_3A&=& \nabc_4Q(A)- C_1  \nabc_4\nabc_3A -C_2  \nabc_4A\\
 &&- \nabc_4C_1  \nabc_3A -(\nabc_4C_2)  A.
 \eeaa
 Hence,
  \beaa
[Q, \LL]A &=&  4\etab  \c \nabc Q(A)  -  2\trchb  \nabc_4Q(A)+\hat{V}  Q(A) \\
&&+ Z_{43} \  \nabc_4\nabc_3A+Z_4 \  \nabc_4A+Z_{a3}  \nabc_a \nabc_3A \\
&&+ Z_3\nabc_3A+Z_a \nabc_a A +Z_0  A+\err,
  \eeaa
  where
  \beaa
  \hat{V}&=&I_{33}+J_{33}+K_{33}+M_{33},
  \eeaa
  and
  \beaa
  Z_{43}&=&  I_{43}+J_{43}+L_{43}+C_1 \left(\tr \Xb+\ov{\tr \Xb}\right)\\
  Z_4&=& I_4+J_4+L_4 +C_2 \left(\tr \Xb+\ov{\tr \Xb}\right), \\
  Z_{a3}&=&I_{a3}+J_{a3}+L_{a3}+M_{a3}-4C_1 \etab, \\
   Z_{3}&=&I_3+J_3 +K_3+L_3+N_3-4\etab \nabc C_1+\nabc_4C_1 \left(\tr \Xb+\ov{\tr \Xb}\right)\\
   &&-C_1 \left(I_{33}+J_{33}+K_{33}+M_{33}\right)+J^a_3,\\
  Z_a&=& J_a+L_a+M_a-4C_2 \etab-\dual J_{\dual a }, \\
   Z_0&=& I_0+J_0+K_0+L_0+N_0-4\etab \c \nabc C_2+\nabc_4C_2 \left(\tr \Xb+\ov{\tr \Xb}\right)\\
   &&-C_2 \left(I_{33}+J_{33}+K_{33}+M_{33}\right)+J^a_0.
               \eeaa

  
\subsubsection{The $Z$ coefficients}

               
    We now show that with the choices of $C_1$ and $C_2$ given by \eqref{eq:first-assumptions-C1-C2}, i.e. $C_1= 2\trchb +\widetilde{C_1}$, $C_2=\frac 1 2 \trchb^2 + \widetilde{C_2}$ all the $Z$ coefficients are $O(|a|)$.  We denote by $O(|a| r^{-c})$ any function which vanishes in Schwarzschild, such as multiples of $\eta$, $\etab$, $\atrchb$, $\dual \rho$, and has a $r^{-c}$ fall-off in $r$.
       
       We start with $Z_{43}$ and $Z_4$. We have         
\beaa
Z_{43}&=&I_{43} +J_{43}+L_{43}+C_1 \left(\tr \Xb+\ov{\tr \Xb}\right)\\
&=& \nabc_3 C_1+\frac {1 }{2} (\tr \Xb+\ov{\tr\Xb}) C_1- \ov{\tr\Xb} \tr\Xb \\
&=&  \nabc_3(2\trchb + \widetilde{C_1})+(2\trchb + \widetilde{C_1}) \trchb  - (\trchb^2+\atrchb^2)\\
&=&  - \big( \trchb^2-\atrchb^2\big) +\nabc_3 \widetilde{C_1}+(2 \trchb + \widetilde{C_1}) \trchb  - (\trchb^2+ \atrchb^2)+r^{-1}\Ga_b\\
&=&\nabc_3 \widetilde{C_1} + \trchb \widetilde{C_1}+r^{-1} \Ga_b,
\eeaa
and
\beaa
Z_4&=& I_4 +J_4  + L_4 +C_2 \left(\tr \Xb+\ov{\tr \Xb}\right)\\
&=& \nabc_3 C_2+ \left(\tr \Xb+\ov{\tr \Xb}\right)C_2-\frac 1 4\tr\Xb\ov{\tr\Xb}C_1\\
&=& \nabc_3 \left( \frac 1 2 \trchb^2 +\widetilde{C_2}\right) +\left( \frac 1 2 \trchb^2 +\widetilde{C_2} \right) 2\trchb-\frac 1 4\big( 2 \trchb +\widetilde{C_1}\big) (\trchb^2+ \atrchb^2)\\
&=&  \trchb  \left(-\frac 1 2 \big( \trchb^2-\atrchb^2\big)+r^{-1}\Ga_b\right) +\nabc_3 \widetilde{C_2}\\
&&+\left( \frac 1 2 \trchb^2 + \tilde{C_2} \right) 2\trchb-\frac 1 4( 2 \trchb + \widetilde{C_1}) (\trchb^2+ \atrchb^2)\\
&=&\nabc_3 \widetilde{C_2}+ 2\trchb  \widetilde{C_2}-\frac 1 4 (\trchb^2+ \atrchb^2) \widetilde{C_1}  +r^{-2} \Ga_b.
\eeaa
From the above expressions we clearly see that $Z_{43}=O(|a|r^{-3})$ and $Z_4=O(|a|r^{-4})$, and the error terms are of the same form as $\err$. Indeed,
\beaa
r^{-1} \Ga_b \c \nabc_4 \nabc_3 A+ r^{-2} \Ga_b \c \nabc_4 A=r^{-2} \Ga_b \c \dk^{\leq1} \nabc_3A=r^{-3} \Ga_b \c \dk^{\leq 2}(A, B),
\eeaa
where we used the Bianchi identity for $A$.

From the expressions above, we immediately have that $Z_{a3}$ and $Z_a$ are $O(|a|)$, and are given by 
\beaa
 Z_{a3}&=&I_{a3}+J_{a3}+L_{a3}+M_{a3}-4C_1 \etab= O(|a|r^{-3}) +\dk^{\leq 1} \Ga_g, \\
 Z_a&=& J_a+L_a+M_a-4C_2 \etab-\dual J_{\dual a }=O(|a|r^{-4}) + r^{-1}\dk^{\leq 1}\Ga_g.
\eeaa
The error terms are then given by
\beaa
\dk^{\leq 1}\Ga_g \c \nabc_a \nabc_3 A+ r^{-1} \dk^{\leq 1}\Ga_g \c \nabc_a A=r^{-1} \dk^{\leq 1}\Ga_g \c \dk^{\leq1} \nabc_3A=r^{-2}\dk^{\leq 1} \Ga_g \c \dk^{\leq 2}(A, B),
\eeaa
as above.

Finally, we consider the coefficients $Z_3$ and $Z_0$. 
   In particular observe that, according to \eqref{eq:first-assumptions-C1-C2}, we can write
  \beaa
  C_1&=& 2 \trchb +O(|a| r^{-2}), \qquad C_2=\frac 1 2 \trchb^2 + O(|a| r^{-3}).
  \eeaa
  We therefore have
  \beaa
  \nabc_4 C_1&=& - \trch\trchb+4\rho +O(|a| r^{-3})+ r^{-1} \Ga_g, \\
  \nabc_3 C_1&=& -\trchb^2 +O(|a| r^{-3}) + r^{-1} \Ga_b, \\
  \nabc_4 \nabc_3 C_1&=&  \trch\trchb^2 -4\trchb \rho+O(|a| r^{-4})+ r^{-2}\Ga_b, \\
  \nabc_4C_2&=& -\frac 1 2 \trch\trchb^2+2\trchb \rho +O(|a| r^{-4})+r^{-2} \Ga_g, \\
  \nabc_3 C_2&=& -\frac 1 2 \trchb^3 +O(|a| r^{-4})+r^{-2} \Ga_b,\\
  \nabc_4 \nabc_3 C_2&=& \frac 3 4  \trch\trchb^3- 3  \trchb^2 \rho  +O(|a| r^{-5})+r^{-3} \Ga_b.
  \eeaa
  In what follows we omit the error terms as they are of the form $\err$ given above.
We compute
\beaa
I_{3}&=&-\nabc_3(2\rho)- 2\trchb \ (2\rho) -\frac 1 2 \trch\trchb^2+2\trchb \rho+  \trch\trchb^2 -4\trchb \rho +O(|a| r^{-4})\\
&=&\frac 1 2 \trch\trchb^2-3\trchb \rho+O(|a| r^{-4}), \\
J_3&=& -  2\trchb \left(\nabc_3\left(\frac 5 2 \trch\right)+(\trch) (\trchb ) \right) +\left(- \frac 52\trchb^2\right)\, \left( \frac 5 2 \trch\right)+O(|a| r^{-4})\\
&=&-\frac{23}{4} \trch\trchb^2 -  10\trchb  \rho+O(|a| r^{-4}),\\
K_3&=& \nabc_3\nabc_3\left(-\frac 5 2 \trch \right)+2\trchb \nabc_3\left(-\frac 5 2 \trch \right)+\frac 5 2 \trch(-\trchb^2) +O(|a| r^{-4})\\
&=& -\frac 5 2\nabc_3\left(-\frac 1 2 \trchb\trch+2\rho\right)-5\trchb \left(-\frac 1 2 \trchb\trch+2\rho\right)+\frac 5 2 \trch(-\trchb^2)\\
&& +O(|a| r^{-4})\\
&=& -\frac 5 2\left(\frac 1 4 \trchb^2\trch-\frac 1 2 \trchb\left(-\frac 1 2 \trchb\trch+2\rho\right)-3\trchb\rho\right)-10\trchb \rho+O\left(\frac{|a|}{r^4}\right)\\
&=& -\frac 5 4 \trchb^2\trch +O(|a| r^{-4}),
\eeaa
\beaa
 L_3 &=& -\frac{1}{2}\trchb\left(6\rho  -(- \trch\trchb+4\rho)\right)+O(|a| r^{-4})= -\frac{1}{2}\trch\trchb^2-\trchb \rho+O(|a| r^{-4}),\\
N_3&=& 2\left(-\nabc_3\trch\trchb-\trch\nabc_3\trchb+2\nabc_3\rho\right) +O(|a| r^{-4})\\
&=& 2 \trch\trchb^2-10\rho\trchb+O(|a| r^{-4}).
\eeaa
We compute 
\beaa
I_{33}&=& 3P-5\ov{P} +4\eta\c\etab-2|\eta|^2  +\nabc_4( C_1 )=  - \trch\trchb+2\rho +O(|a| r^{-3}),\\
J_{33}&=& -  \left(\tr \Xb+\ov{\tr \Xb}\right) (\frac 1 2 \tr X +2\ov{\tr X})+\DDc\c\ov{H}+2H \c\ov{H}= - 5\trchb  \trch+O(|a| r^{-3}),\\
K_{33}&=& 2\nabc_3\left(- \frac 1 2 \tr X -2\ov{\tr X} \right)= \frac 5 2 \trch\trch-10\rho +O(|a| r^{-3}),\\
 M_{33}&=& 2\left( 4H+\Hb +\ov{\Hb} \right) \c \eta=O(|a| r^{-4}).
\eeaa
We compute
\beaa
I_0&=&  \nabc_4\nabc_3C_2= \frac 3 4  \trch\trchb^3- 3  \trchb^2 \rho  +O(|a| r^{-5}),\\
J_0&=& -  2\trchb \Big(\nabc_3\left(\trch\trchb -2\rho\right) +\frac{1}{2}\trchb(4\rho  ) \Big)+\left(-\frac 5 2\trchb^2  \right)\,  \left(\trch\trchb -2\rho\right)\\
&&+O(|a| r^{-5})\\
&=& -\frac 1 2  \trch\trchb^3-9\rho\trchb^2+O(|a| r^{-5}),\\
K_0&=& -\left(- \frac 1 2 \tr X -2\ov{\tr X} \right)\nabc_3C_2= -\frac 5 4 \trch  \trchb^3+O(|a| r^{-5}),\\
 L_0  &=& -2 \trchb^2  \rho +\frac{1}{2}\trchb(-\frac 1 2 \trch\trchb^2+2\trchb \rho)\   -2\trchb \nabc_3(\rho)+O(|a| r^{-5})\\
   &=&-\frac{1}{4} \trch\trchb^3+ 2 \trchb^2  \rho \  +O(|a| r^{-5})\\
N_0&=& \nabc_3 \nabc_3\left(-\trch\trchb + 2\rho \right)+2\trchb  \nabc_3\left(-\trch\trchb + 2\rho\right)\\
&=& ( -\frac 1 2 \trch\trchb+2\rho)\trchb^2+ 2\trch\trchb (-\frac 1 2 \trchb^2)-5(-\frac 3 2 \trchb\rho)\trchb-5\rho(-\frac 1 2 \trchb^2)\\
&&+2\trchb  \left( \trch\trchb^2-5\rho\trchb\right)+O(|a| r^{-5})= \frac 1 2 \trch \trchb^3+ 2\rho\trchb^2+O(|a| r^{-5}).
\eeaa
We finally obtain
\beaa
 Z_{3} &=&\frac 1 2 \trch\trchb^2-3\trchb \rho-\frac{23}{4} \trch\trchb^2 -  10\trchb  \rho -\frac 5 4 \trchb^2\trch-\frac{1}{2}\trch\trchb^2-\trchb \rho\\
 &&+2 \trch\trchb^2-10\rho\trchb+ 2\trchb(- \trch\trchb+4\rho) \\
 &&-2\trchb \left(- \trch\trchb+2\rho - 5\trchb  \trch+\frac 5 2 \trch\trch-10\rho\right)+O(|a| r^{-4})=O(|a| r^{-4}).
 \eeaa
and
\beaa
 Z_0 &=&  \frac 3 4  \trch\trchb^3- 3  \trchb^2 \rho -\frac 1 2  \trch\trchb^3-9\rho\trchb^2-\frac 5 4 \trch  \trchb^3+2 \trchb^2  \rho -\frac{1}{4} \trch\trchb^3\\
 &&+2\rho\trchb^2+ \frac 1 2 \trch \trchb^3+2\trchb (-\frac 1 2 \trch\trchb^2+2\trchb \rho) \\
 &&-\frac 1 2 \trchb^2 \left(- \trch\trchb+2\rho - 5\trchb  \trch+\frac 5 2 \trch\trch-10\rho\right)+O(|a| r^{-5})= O(|a| r^{-5}).
\eeaa
as stated. This completes the proof of Proposition \ref{first-intermediate-step-main-theorem}.


\subsection{Step 2. Derive the wave equation for $Q(A)$ and $\qf$}
\label{proof:step2-main-thm-Part1}


We use the previous two steps to derive the wave equation for $Q=Q(A)$ from the Teukolsky equation for $A$ and the commutator $[Q, \LL]$.

\begin{proposition}\label{prop:wave-eq-Q}
Let $Q=Q(A) \in \sk_2$ as in Proposition \ref{first-intermediate-step-main-theorem}. Then $Q$ satisfies the following wave equation:
\bea\label{wave-eq-Q}
\begin{split}
\squared_2 Q&= ( 2\ov{\tr X} )\nab_3Q+\left( \tr \Xb+\ov{\tr \Xb}\right)  \nab_4Q-\left( 4H+2\Hb +2\ov{\Hb}\right)\c \nab Q+\widetilde{V}  Q\\
&-L_Q(A) + \err[\square_2 Q],
\end{split}
\eea
where 
\bea\label{definition-tilde-V}
\begin{split}
\widetilde{V}&=   \frac 1 2\trch\trchb+ \frac 1 2 \atrch\atrchb-4\rho-4\eta \c \etab\\
& + i \left(-\trch  \atrchb+ \atrch \trchb+4 \rhod+ 2\eta \wedge \etab   \right)  -\hat{V} ,
\end{split}
\eea
and 
\bea
\err[\squared_2 Q]&=&  Q( \err[\LL(A)])+\err[ [Q, \LL]A]+(\Ga_b \c \Ga_g) \c Q.
\eea
\end{proposition}

\begin{proof}
 Recall \eqref{Teuk-repeat=2}, i.e.
 \beaa
 \LL(Q(A))+[Q, \LL](A)&=& Q( \err[\LL(A)]).
 \eeaa
Applying the definition of the operator $\LL$ given by \eqref{Teukolsky-operator-ch5} to $Q=Q(A)$, we obtain\footnote{Recall that $Q(A)$ is of conformal type $0$, therefore all conformal derivatives coincide with the non-conformal ones, see also \eqref{first-equation-square-conformal-sign0}.}
 \beaa
  \LL(Q) &=&-\nab_4\nab_3Q+ \frac{1}{4}\DD\hot (\DDb \c Q)+\left(- \frac 1 2 \tr X -2\ov{\tr X}-2\om \right)\nab_3Q-\frac{1}{2}\tr\Xb \nab_4Q\\
&&+\left( 4H+\Hb +\ov{\Hb} \right)\c \nab Q+ \left(-\ov{\tr X} \tr \Xb +2\ov{P}+4\eta \c \etab -4i \eta \wedge \etab\right) Q.
 \eeaa
Using \eqref{final-commutator} of Proposition \ref{first-intermediate-step-main-theorem}, equation \eqref{Teuk-repeat=2} gives
\beaa
&&-\nab_4\nab_3Q+ \frac{1}{4}\DD\hot (\DDb \c Q)+\left(- \frac 1 2 \tr X -2\ov{\tr X} \right)\nab_3Q-\left(\frac 3 2 \tr \Xb+\ov{\tr \Xb}\right)  \nab_4Q\\
&&+\left( 4H+\Hb +\ov{\Hb} +4\etab\right)\c \nab Q+ \left(-\ov{\tr X} \tr \Xb +2\ov{P}+\hat{V}+4\eta \c \etab -4i \eta \wedge \etab\right) Q\\
&=&-L_Q(A) + Q( \err[\LL(A)])+\err[ [Q, \LL]A].
\eeaa
Using the formula for the wave equation \eqref{wave-equation-Psi} applied to $Q\in \sk_2(\CCC)$, i.e.
\beaa
\squared_2 Q&=&-\nab_4 \nab_3 Q +\frac 1 4  \DD\hot( \DDb \c Q)+\left(2\om -\frac 1 2 \tr X\right) \nab_3Q- \frac 1 2 \tr\Xb \nab_4Q+2\etab \c\nab Q \\
&& +  \left( -\frac 1 2\trch\trchb- \frac 1 2 \atrch\atrchb-2\rho\right) Q+ 2i \left( \rhod-\eta \wedge \etab  \right) Q +(\Ga_b \c \Ga_g) \c Q,
\eeaa
we obtain \eqref{wave-eq-Q}, with
\beaa
\widetilde{V}&=&   -\frac 1 2\trch\trchb- \frac 1 2 \atrch\atrchb-2\rho + 2i \left( \rhod- \eta \wedge \etab   \right) \\
&&- \left(-\ov{\tr X} \tr \Xb +2\ov{P}+\hat{V}+4\eta \c \etab -4i \eta \wedge \etab\right), 
\eeaa
as stated.
\end{proof}

  We now want to rescale $Q$ in order to absorb the real parts of the first order terms in \eqref{wave-eq-Q} into the wave operator. The rescaling is obtained through a scalar function of $q$ and $\ov{q}$, i.e.
  \beaa
  \qf&=& f Q(A)=f \left(\nabc_3\nabc_3 A + C_1  \nabc_3A + C_2  A\right) \in \sk_2(\CCC),
  \eeaa
  with $C_1$ and $C_2$ given as in Proposition \ref{first-intermediate-step-main-theorem}.

\begin{proposition}\label{prop:rescaling-f} 
Let $f$ be given by
\beaa
f=q \ov{q}^3.
\eeaa
Then $\qf =f Q(A) \in \sk_2$ satisfies the following wave equation:
    \bea\label{final-eq}\label{eq:squared-2-qf-almost-end}
  \squared_2 \qf  - i   \frac{4a\cos\th}{|q|^2} \nab_\T  \qf- V_1 \qf   &=&  \widetilde{L_{\qf}[A]} + f \Big( \err[\squared_2 Q]+ \Ga_g \c \dk^{\leq 1} Q \Big),
 \eea
 where
 \bea
 V_1&:=& f^{-1}\square_\g(f) +\widetilde{V}, \label{definition-V-final}\\
 \widetilde{ L_{\qf}[A]}&:=& - f L_Q(A) , \nn \\
   \err[\squared_2 \qf]&:=&f \Big( \err[\squared_2 Q]+ \Ga_g \c \dk^{\leq 1} Q \Big). \nn
 \eea
\end{proposition}
\begin{proof} 
Let $f$ be given by 
 \beaa
 f&=& (q)^n (\ov{q})^{m}.
 \eeaa
Using Lemma \ref{lemma:equations-q}, we deduce
 \beaa
 \nab_3(f)&=& \left(\frac n 2 \ov{\tr \Xb} +\frac m 2   \tr \Xb \right) f+r^{n+m} \Ga_b, \\ 
  \nab_4(f) &=& \left(\frac n 2 \tr X +\frac m 2   \ov{\tr X} \right) f+r^{n+m-1} \Ga_g, \\
    \nab f    &=&  \left(\frac{m}{2} H+\frac n2\ov{H}+\frac n2   \Hb+ \frac{m }{2}  \ov{\Hb}  \right)f+r^{n+m} \Ga_g.
    \eeaa
 We obtain for $\qf=f Q$ from \eqref{wave-eq-Q}, 
 \beaa
  \squared_2 \qf &=& \square_\g(f) Q+f \squared_2(Q)- \nab_3 f \nab_4Q- \nab_4f \nab_3 Q +2\nab f \c \nab Q\\
  &=& \square_\g(f) Q\\
  &&+ \Big(2\ov{\tr X} f\nab_3Q+\left( \tr \Xb+\ov{\tr \Xb}\right) f \nab_4Q-\left( 4H+2\Hb +2\ov{\Hb} \right)\c f\nab Q+f \widetilde{V} Q \Big)\\
&&- \left(\frac n 2 \ov{\tr \Xb} +\frac m 2   \tr \Xb + \Ga_b \right) f\nab_4Q-  \left(\frac n 2 \tr X +\frac m 2   \ov{\tr X} +r^{-1} \Ga_g \right) f \nab_3 Q \\
  &&+2 \left(\frac{m}{2} H+\frac n2\ov{H}+\frac n2   \Hb+ \frac{m }{2}  \ov{\Hb} +\Ga_g \right)f \c \nab Q+ f \Big(-L_Q(A) + \err[\square_2 Q] \Big),
 \eeaa
which gives
 \beaa
  \squared_2 \qf &=& \left((1-\frac n 2) \ov{\tr \Xb}+(1 -\frac m 2 )  \tr \Xb \right) f\nab_4Q+  \left(-\frac n 2 \tr X +(2-\frac m 2)   \ov{\tr X} \right) f \nab_3 Q \\
  &&+ \left((m-4) H+n\ov{H}+(n-2)   \Hb+ (m-2) \ov{\Hb}  \right)f \c \nab Q+ \left(f^{-1}\square(f) +\widetilde{V} \right)\qf\\
   &&+ f \Big(-L_Q(A) + \err[\square_2 Q]+r^{-1} \Ga_b \c \dk Q \Big).
 \eeaa
 Observe that the real part of the coefficients of all the first derivatives are multiple of $m+n-4$. To cancel the real part of those coefficients we then take $m=4-n$, therefore for $f=(q)^n (\ov{q})^{4-n}$, we have
  \beaa
  \squared_2 \qf    &=& \left((1-\frac n 2) \ov{\tr \Xb}-(1 -\frac {n}{ 2} ) \tr \Xb \right) f\nab_4Q+  \left(-\frac n 2 \tr X +\frac {n}{ 2}   \ov{\tr X} \right) f \nab_3 Q \\
  &&+ \left((-n) H+n\ov{H}+(n-2)   \Hb+ (2-n) \ov{\Hb}  \right)f \c \nab Q+ \left(f^{-1}\square_\g(f) +\widetilde{V} \right)\qf\\
  &&+ f \Big(-L_Q(A) + \err[\square_2 Q]+r^{-1} \Ga_b \c \dk Q \Big)\\
  &=&i  f \Big( (2- n )\atrchb \nab_4+ n \atrch  \nab_3 + \left((-2n) \dual \eta+2(n-2) \dual \etab \right) \c \nab \Big)  Q\\
  &&+ \left(f^{-1}\square_\g(f) +\widetilde{V} \right)\qf+ f \Big(-L_Q(A) + \err[\square_2 Q] +r^{-1} \Ga_b \c \dk Q\Big).
 \eeaa
 Using, see \eqref{eq:atrch-e3-atrchb-e4-pert-kerr}, that
 \beaa
 \atrch e_3+\atrchb e_4+ 2(\eta+\etab) \c \dual \nab&=&\frac{4a\cos\th}{|q|^2} \T+ \Ga_g \c \dk,
 \eeaa
 we obtain for $n=1$, 
   \beaa
  \squared_2 \qf      &=&i  f \frac{4a\cos\th}{|q|^2} \nab_\T  Q+ \left(f^{-1}\square_\g(f) +\widetilde{V} \right)\qf\\
  &&+ f \Big(-L_Q(A) + \err[\square_2 Q]+ \Ga_g \c \dk^{\leq 1} Q \Big),
 \eeaa
 as stated. 

\end{proof}

We now analyze the error terms. We have
\beaa
   \err[\squared_2 \qf]&:=&f \Big( \err[\squared_2 Q]+ \Ga_g \c \dk^{\leq 1} Q \Big)\\
   &=&f \Big(  Q( \err[\LL(A)])+\err[ [Q, \LL]A]+(\Ga_b \c \Ga_g) \c Q+ \Ga_g \c \dk^{\leq 1} Q \Big)\\
      &=&f \Big(  Q\big( r^{-1}  \dk^{\leq 1}\big( \Ga_g \c  B\big) + \nabc_3\Xi  \c B + r^{-1} \Ga_b \c \Ga_g \c \Ga_g\big)\\
      &&+r^{-2} \dk^{\leq 3}(\Ga_g \c (A, B))+\nab_3 \big( r^{-2} \Ga_b  \c \frak{d}^{\leq 2} A \big)+r^{-1}\dk\big( (\Ga_b \c \Ga_g) A \big)\\
      &&+(\Ga_b \c \Ga_g) \c Q+ \Ga_g \c \dk^{\leq 1} Q \Big),
\eeaa
which gives
   \beaa
 \err[\squared_2 \qf]&=& r^2 \frak{d}^{\leq 3} (\Ga_g \c (A, B))+ \nab_3 (r^2 \frak{d}^{\leq 2}( \Ga_b \c (A, B)))\\
&&+\frak{d}^{\leq 1} (\Ga_g \c \qf) + r^{3}\dk^{\leq 2} \big( \Ga_b \c \Ga_g \c \Ga_g\big),
 \eeaa
 as stated in Theorem \ref{MAIN-THEOREM-PART1}.


\subsection{Step 3. Reality of the potential and lower order terms}\label{proof:step-3-main-thm-Part1}


Observe that in order to obtain equation \eqref{final-eq} for $\qf$ we only needed so far to impose that $C_1$ and $C_2$ are given by \eqref{eq:first-assumptions-C1-C2}, i.e.
\beaa
C_1= 2\trchb +\widetilde{C_1}, \qquad C_2=\frac 1 2 \trchb^2 + \widetilde{C_2},
\eeaa
with $\widetilde{C_1}$ and $\widetilde{C_2}$ are complex functions satisfying $\widetilde{C_1}=O(|a|r^{-2})$, $\widetilde{C_2}=O(|a|r^{-3})$.

 To complete the proof of Theorem \ref{MAIN-THEOREM-PART1}, we need to show that there exists a choice of complex functions $\widetilde{C_1}$ and $\widetilde{C_2}$ for which the potential $V$ is real, the scalar $Z_{43}$ and the one-form $Z_{a3}$ are real and the one-form $Z_{13}$ vanishes.

We have the following.

 \begin{proposition}\label{prop:potential-real}\label{prop:lot-real}\label{lemma-lot}  
Let $\widetilde{C_1}$ be the complex function with $C_1=2\trchb + \widetilde{C_1}$. Then
\begin{enumerate}
\item if $\Im(\widetilde{C_1}) = -4 \atrchb$, then 
\begin{itemize}
\item  the potential $V_1$ as given in \eqref{definition-V-final} is real, and is given by 
\beaa
V_1=-\trch\trchb+ O(\frac{|a|}{r^4});
\eeaa
\item the scalar $Z_{43}$ and the one-form $Z_{a3}$ are real;
\end{itemize}
\item if $\Re(\widetilde{C_1})=-2\frac {(\atrchb)^2}{ \trchb} $, then the one-form $Z_{13}$ vanishes. Moreover,  $V_1$, $Z_{43}$ and $Z_{23}$ are given (in the outgoing frame) by
\beaa
 V_1&=& \frac{4}{|q|^2}\frac{r^2-2mr+2a^2}{r^2}\\
&&-\frac{4a^2\cos^2\th}{r^2|q|^6}( 2r^4+4mr^3+a^2r^2+2a^2\cos^2\th r^2-2mr a^2\cos^2\th+a^4\cos^2\th), \\
Z_{43}&=& \atrchb^2\big( 1 +\frac{\atrchb^2}{\trchb^2}\big)=\atrchb^2 \frac{|q|^2}{r^2}=\frac{4a^2\De^2\cos^2\th}{r^2|q|^6}, \\
Z_{23}&=& \frac{8a\sin\th \De}{|q|^5}.
\eeaa
\end{enumerate} 
\end{proposition}

\begin{proof} 
Here we compute $\Im(V_1)=\Im(f^{-1}\square_\g(f) +\widetilde{V})$. We start with the following.

\begin{lemma}\label{lemma:squaref-1f} Let $f= q \ov{q}^{3}$. Then
\beaa
\Im(f^{-1}\square_\g(f))&=&   -3\trchb \atrch+2 \trch \atrchb +2 \dual \rho + \div(\dual \eta+\dual \etab) - 10\eta \wedge \etab +r^{-1}\dk^{\leq1 } \Ga_b,\\
\Re(f^{-1}\square_\g f)&=& - 5 \trch\trchb-2\atrch\atrchb -4\rho\\
&&+2\div(\eta-\etab)+ 3( |\eta|^2+ |\etab|^2)+14\eta \c \etab +r^{-1}\dk^{\leq1 } \Ga_g.
\eeaa
\end{lemma}
\begin{proof}
Using that
 \beaa
 \nab_3f&=& \left(\frac 1 2 \ov{\tr \Xb} +\frac 3 2   \tr \Xb \right) f+r^{4} \Ga_b= \left(2 \trchb -i\atrchb \right) f+r^{4} \Ga_b, \\ 
  \nab_4f &=& \left(\frac 1 2 \tr X +\frac 3 2   \ov{\tr X} \right) f+r^{3} \Ga_g= \left(2\trch +i \atrch \right) f+r^{3} \Ga_g, \\
    \nab f    &=&  \left(\frac{3}{2} H+\frac 12\ov{H}+\frac 12   \Hb+ \frac{3 }{2}  \ov{\Hb}  \right)f+r^{4} \Ga_g=  \left(2(\eta+ \etab) +i (\dual \eta-\dual \etab) \right)f+r^{4} \Ga_g,
    \eeaa
    we deduce
    \beaa
    \nab_4\nab_3 f&=&\left(2 \nab_4\trchb -i\nab_4\atrchb \right) f+\left(2 \trchb -i\atrchb \right) \nab_4f+r^{3}\dk \Ga_b \\
    &=&\left(2 \nab_4\trchb -i\nab_4\atrchb \right) f+\left(2 \trchb -i\atrchb \right)  \left(2\trch +i \atrch \right) f+r^{3}\dk \Ga_b \\
    \lap f&=&   \left(2\div(\eta+ \etab) +i \div(\dual \eta-\dual \etab) \right)f+\left(2(\eta+ \etab) +i (\dual \eta-\dual \etab) \right)\c \nab f+r^{3}\dk \Ga_g\\
    &=&   \left(2\div(\eta+ \etab) +i \div(\dual \eta-\dual \etab) \right)f\\
    &&+\left(2(\eta+ \etab) +i (\dual \eta-\dual \etab) \right)\c \left(2(\eta+ \etab) +i (\dual \eta-\dual \etab) \right)f+r^{3}\dk \Ga_g.
    \eeaa
    Using Lemma \ref{lemma:expression-wave-operator-pert} applied to a scalar function, we obtain 
\beaa
\Im(f^{-1}\square_\g f)&=&-\Im(f^{-1}\nab_4 \nab_3 f) -\frac 1 2 \trchb \Im(f^{-1}\nab_4f)+\left(2\om -\frac 1 2 \trch\right) \Im(f^{-1}\nab_3f)\\
&&+\Im(f^{-1}\lap_2 f)+2\etab \c \Im(f^{-1}\nab f )\\
&=&\nab_4\atrchb  -2 \trchb \atrch+2 \trch \atrchb   -\frac 1 2 \trchb \atrch\\
&&-\left(2\om -\frac 1 2 \trch\right)\atrchb+ \div(\dual \eta-\dual \etab) \\
&&+ 4(\eta+ \etab) \c (\dual \eta-\dual \etab) +2\etab \c (\dual \eta-\dual \etab)+r^{-1}\dk^{\leq1 } \Ga_b\\
&=&   -3\trchb \atrch+2 \trch \atrchb +2 \dual \rho + \div(\dual \eta+\dual \etab) - 10\eta \wedge \etab \\
&&+r^{-1}\dk^{\leq1 } \Ga_b,
\eeaa
where we used the null structure equation for $\nab_4\atrchb$. Similarly we have
\beaa
\Re(f^{-1}\square_\g f)&=&-\Re(f^{-1}\nab_4 \nab_3 f) -\frac 1 2 \trchb \Re(f^{-1}\nab_4f)+\left(2\om -\frac 1 2 \trch\right) \Re(f^{-1}\nab_3f)\\
&&+\Re(f^{-1}\lap_2 f)+2\etab \c \Re(f^{-1}\nab f )\\
&=&-2\nab_4 \trchb - 4 \trch\trchb-\atrch\atrchb -\frac 1 2 \trchb (2\trch)\\
&&+\left(2\om -\frac 1 2 \trch\right) (2\trchb)+2\div(\eta+\etab)+ 4 |\eta+\etab|^2-|\dual \eta - \dual \etab|^2\\
&&+2\etab \c (2(\eta+ \etab) )\\
&=& - 5 \trch\trchb-2\atrch\atrchb -4\rho\\
&&+2\div(\eta-\etab)+ 3( |\eta|^2+ |\etab|^2)+14\eta \c \etab +r^{-1}\dk^{\leq1 } \Ga_g, 
\eeaa
where we used the null structure equation for $\nab_4 \trchb$. 
\end{proof}

We now compute $\Im(\widetilde{V})$ using \eqref{definition-tilde-V}. We have
\beaa
\Im(\widetilde{V})&= & -\trch  \atrchb+ \atrch \trchb+4 \rhod+ 2\eta \wedge \etab   -\Im(\hat{V} ),
\eeaa
and 
\beaa
\Im(\hat{V})&=&\Im(I_{33})+\Im(J_{33})+\Im(K_{33})+\Im(M_{33}).
\eeaa
Using \eqref{eq:expression-I33}, \eqref{eq:J33-final}, \eqref{eq:expression-K33}, and \eqref{eq:expression-M33}, we compute
\beaa
\Im( I_{33})&=& 8 \rhod  -8 \eta \wedge \etab   +\nabc_4\Im( \widetilde{C_1}), \\
\Im(J_{33})&=& -  3\trchb  \atrch-2\div \dual \eta, \\
\Im(K_{33})&=&-3 \nabc_3\atrch=\frac 3 2(\atrchb \trch+\trchb\atrch) +6 \dual \rho-6 \div \dual \eta , \\
\Im(M_{33})&=&\Im( 2\left( 4(\eta+ i \dual \eta)+2\etab \right) \c \eta)=0,
\eeaa
which gives
\beaa
\Im(\hat{V})&=&\frac 3 2(\atrchb \trch-\trchb\atrch)+14 \rhod -8\div \dual \eta  -8 \eta \wedge \etab   +\nabc_4\Im( \widetilde{C_1}),
\eeaa
and 
\beaa
\Im(\widetilde{V})&= & -\frac 5 2 (\trch  \atrchb- \trchb\atrch )-10 \rhod   +8\div \dual \eta  +10 \eta \wedge \etab   -\nabc_4\Im(\widetilde{ C_1}).
\eeaa
Combining the above with Lemma \ref{lemma:squaref-1f}, we finally deduce
\beaa
\Im(V_1)&=&\Im(f^{-1}\square_\g(f) +\widetilde{V})\\
&=&-\frac 1 2 (\trch  \atrchb+ \trchb\atrch )  -8 \dual \rho + \div(\dual \eta+\dual \etab)   +8\div \dual \eta   -\nabc_4\Im( \widetilde{C_1})\\
&&+r^{-1}\dk^{\leq1 } \Ga_b.
\eeaa
Using that, see \eqref{eq:vanishing-trch-atrch} and \eqref{eq:vanishing-div-dual-eta-etab},
\beaa
\trch  \atrchb+ \trchb\atrch &=&r^{-1} \Ga_g, \\
\div(\dual \eta+\dual \etab)&=&r^{-1} \dk \Ga_g,
\eeaa
we deduce that in order to obtain $\Im(V_1)=0$ we need to have
\bea\label{eq:condition-on-Im-V}
\nabc_4\Im( \widetilde{C_1})&=&  -8 \dual \rho    -8\div \dual \etab  +r^{-1}\dk^{\leq1 } \Ga_b.
\eea
 Therefore, if $\Im(\widetilde{C_1})=-4 \atrchb$, relation \eqref{eq:condition-on-Im-V} is satisfied, and $V_1$ is a real function.

We now explicitly compute $V_1$ for $\Re(C_1)=2\trchb -2\frac {(\atrchb)^2}{ \trchb} $.
We have
\beaa
\Re(\widetilde{V})&=&   \frac 1 2\trch\trchb+ \frac 1 2 \atrch\atrchb-4\rho-4\eta \c \etab  -\Re(\hat{V} ),
\eeaa
and 
\beaa
\Re(\hat{V} )&=& \Re(I_{33})+\Re(J_{33})+\Re(K_{33})+\Re(M_{33}).
\eeaa
Using \eqref{eq:expression-I33}, \eqref{eq:J33-final}, \eqref{eq:expression-K33}, and \eqref{eq:expression-M33}, we compute
\beaa
\Re( I_{33})&=& -2\rho-2 \eta \c (\eta-2\etab ) +\nabc_4 \Re(C_1 ),\\
\Re(J_{33})&=& -  5\trch \trchb  +2\div \eta+4|\eta|^2, \\
\Re(K_{33})&=& -5\nabc_3\trch= \frac 5 2 \trchb\trch-\frac 5 2 \atrchb\atrch    -10   \div \eta -10 |\eta|^2-10\rho, \\
\Re(M_{33})&=& \Re(2\eta \c \left( 4(\eta+ i \dual \eta)+2\etab \right) )= 8 |\eta|^2 +4 \eta \c \etab,
\eeaa
which gives
\beaa
\Re(\hat{V} )&=&\nabc_4 \Re(C_1 )  - \frac 5 2 \trchb\trch-\frac 5 2 \atrchb\atrch  -12\rho  -8   \div \eta +8 \eta \c \etab,
\eeaa
and
\beaa
\Re(\widetilde{V})&=&  3\trch\trchb+ 3 \atrch\atrchb+8\rho+8   \div \eta-12\eta \c \etab  - \nabc_4 \Re(C_1 ).
\eeaa
Combining the above with Lemma \ref{lemma:squaref-1f}, we finally deduce
\beaa
V_1&=& \Re(V_1)=\Re(f^{-1}\square_\g(f) +\widetilde{V})\\
&=&  - 2 \trch\trchb+\atrch\atrchb +4\rho+2\div(\eta-\etab)+ 3( |\eta|^2+ |\etab|^2)+2\eta \c \etab\\
&&+8   \div \eta  - \nabc_4 \Re(C_1 ).
\eeaa
Using \eqref{eq:vanishing-div-eta-etab} and \eqref{eq:vanishing-eta2-etab2}, we obtain
\beaa
V_1&=&  - 2 \trch\trchb+\atrch\atrchb +4\rho+8   \div \etab+ 6 |\etab|^2+2\eta \c \etab  - \nabc_4 \Re(C_1 ).
\eeaa
For $\Re(C_1)=2\trchb -2\frac {(\atrchb)^2}{ \trchb} $, we have
\beaa
\nabc_4\Re(C_1)&=&2\nabc_4\trchb -2\nabc_4\left(\frac {\atrchb^2}{ \trchb} \right)\\
&=&2\nabc_4\trchb -4\frac{\atrchb}{\trchb} \nabc_4 \atrchb +2\frac{\atrchb^2}{\trchb^2} \nabc_4 \trchb\\
&=&- \trch\trchb+ \atrch\atrchb    +  4   \div \etab+ 4 |\etab|^2+4\rho -4\frac{\atrchb}{\trchb} \big(2 \curl\etab  +2 \dual \rho \big) \\
&&+2\frac{\atrchb^2}{\trchb^2} \left( -\frac 1 2 \trch\trchb+\frac 1 2 \atrch\atrchb    +  2   \div \etab+ 2 |\etab|^2+2\rho \right).
\eeaa
We therefore obtain
\beaa
V_1&=&  -  \trch\trchb +4   \div \etab+ 2 |\etab|^2+2\eta \c \etab \\
&& +\frac{\atrchb}{\trchb} \Big[  \atrchb\trch-\frac{\atrchb^2}{\trchb} \atrch    -\frac{\atrchb}{\trchb}\big( 4   \div \etab+ 4 |\etab|^2+4\rho\big)\\
&&+8 \curl\etab  +8 \dual \rho \Big].
\eeaa
Using the values in Kerr given in Section \ref{section:values-Kerr}, we have
\beaa
V_1&=&  \frac{4r^2\Delta}{|q|^6} + \frac{4a^4\sin^2\th \cos^2\th }{|q|^6}+4   \div \etab \\
&& -\frac{a\cos\th}{r} \left( \frac{4a\Delta\cos\th}{r|q|^4}   +\frac{a\cos\th}{r}\big( 4   \div \etab+ 4|\etab|^2 +4\rho\big)+8 \curl\etab  +8 \dual \rho \right)\\
&=&  \frac{4r^2\Delta}{|q|^6} +\frac{r^2-a^2\cos^2\th}{r^2}4   \div \etab -\frac{a\cos\th}{r}8 \curl\etab\\
&& -\frac{4a^2\cos^2\th}{r^2} \left( \frac{\Delta}{|q|^4}   + \frac{a^4\sin^2\th \cos^2\th }{|q|^6}  \right)-\frac{4a\cos\th}{r} \left( \frac{a\cos\th}{r} \rho  +2 \dual \rho \right).
\eeaa
By writing
\beaa
\div\etab&=& e_1(\etab_1) + \Lambda \etab_1= \frac{1}{|q|^6}((-3\cos^2\th+1)a^2r^2+a^4\cos^2\th(-3+\cos^2\th)), \\
\curl \etab&=&e_1(\etab_2)+\La \etab_2= \frac{a\cos\th}{|q|^6}(- 2r^3+2a^2\cos^2\th r-4a^2 r ),
\eeaa
we obtain
\beaa
\frac 1 4 V_1&=&  \frac{r^2(r^2-2mr+a^2)}{|q|^6} + \frac{1}{|q|^6}((-3\cos^2\th+1)a^2r^2+a^4\cos^2\th(-3+\cos^2\th))\\
&&-\frac{a^2\cos^2\th}{r^2|q|^6} (- 4r^4+(\cos^2\th-7)a^2r^2-2a^4\cos^2\th)\\
&& -\frac{a^2\cos^2\th}{r^2|q|^4 }  \Delta   -\frac{a^2\cos^2\th}{r |q|^6} \big( 10mr^2+2m a^2\cos^2\th   \big)\\
&=&  \frac{r^2(r^2-2mr+2a^2)}{|q|^6} - \frac{a^2\cos^2 \th}{r^2|q|^6}(8mr^3-3a^2r^2+a^2\cos^2\th r^2-a^4\cos^2\th).
\eeaa
Writing
\beaa
&&\frac{1}{|q|^2}\frac{r^2-2mr+2a^2}{r^2}\\
&=& \frac{(r^4+2a^2r^2\cos^2\th+a^4\cos^4\th)}{|q|^6}\frac{(r^2-2mr+2a^2)}{r^2}\\
&=& \frac{r^2(r^2-2mr+2a^2)}{|q|^6}\\
&&+\frac{a^2\cos^2\th}{r^2|q|^6}( 2r^4-4mr^3+4a^2r^2+a^2\cos^2\th r^2-2mr a^2\cos^2\th+2a^4\cos^2\th),
\eeaa
we obtain 
\beaa
\frac 1 4 V_1&=& \frac{1}{|q|^2}\frac{r^2-2mr+2a^2}{r^2}\\
&&-\frac{a^2\cos^2\th}{r^2|q|^6}( 2r^4+4mr^3+a^2r^2+2a^2\cos^2\th r^2-2mr a^2\cos^2\th+a^4\cos^2\th),
\eeaa
as stated.

Recall from Proposition \ref{first-intermediate-step-main-theorem}, that
\beaa
Z_{43}&=&\nabc_3 \widetilde{C_1} + \trchb \widetilde{C_1}.
\eeaa
In particular we have
\beaa
\Im(Z_{43})&=&\nabc_3 \Im(\widetilde{C_1}) + \trchb \Im(\widetilde{C_1})=r^{-1}\Ga_b
\eeaa
if $\Im(\widetilde{C_1})=n \atrchb$ for any $n$. We can also explicitly compute for $\Re(\widetilde{C_1})=-2\frac {(\atrchb)^2}{ \trchb} $,
\beaa
Z_{43}&=& -2\nabc_3 \left(\frac {(\atrchb)^2}{ \trchb}\right) + \trchb \left(-2\frac {(\atrchb)^2}{ \trchb}\right)\\
&=& -4\frac{\atrchb}{\trchb} \nabc_3 \atrchb +2\frac{\atrchb^2}{\trchb^2} \nabc_3 \trchb  -2(\atrchb)^2\\
&=& -4\frac{\atrchb}{\trchb} (-\trchb\atrchb ) -\frac{\atrchb^2}{\trchb^2}(  \trchb^2-\atrchb^2)  -2(\atrchb)^2+r^{-2} \Ga_b\\
&=& \atrchb^2 +\frac{\atrchb^4}{\trchb^2} +r^{-2} \Ga_b.
\eeaa

We now consider $Z_{a3}$. Recall from Proposition \ref{first-intermediate-step-main-theorem}, that
\beaa
 Z_{a3}&=&I_{a3}+J_{a3}+L_{a3}+M_{a3}-4C_1 \etab,
\eeaa
with
\beaa
 I_{a3}&=& - 2\nabc_3(\eta-\etab )-(2C_1 -   \trchb ) (\eta-\etab )-\atrchb  \dual (\eta-\etab )\\
 J_{a3}&=&2\nabc_3\eta-2\nabc ( C_1 ) +2C_1  \eta+  \trchb (8\etab+2\eta) -2\atrchb   \dual \etab \\
&& +  i  \left(4 \trchb \dual \eta+2\atrchb \eta +2\atrchb \etab \right) \\
L_{a3}&=& -2(\trchb - i \atrchb )(\eta-\etab )  \\
M_{a3}&=& 2 \nabc_3 \left( 4H+\Hb +\ov{\Hb} \right)-   \trchb\left( 4H+\Hb +\ov{\Hb} \right)+\atrchb \dual  \left( 4H+\Hb +\ov{\Hb} \right).
\eeaa

We compute
\beaa
\Im(L_{a3})&=& 2 \atrchb (\eta-\etab ),\\
\Im(M_{a3})&=& \Im(2 \nabc_3 \left( 4\eta+ 4i \dual \eta \right)-   \trchb\left( 4\eta+ 4i \dual \eta \right)+\atrchb \dual  \left(  4\eta+ 4i \dual \eta \right))\\
&=&8 \nabc_3 \dual \eta-  4\trchb  \dual\eta-  4  \atrchb  \eta.
\eeaa
and, recalling that $\Im(C_1)=\Im(\widetilde{C_1})$,
\beaa
\Im(I_{a3})&=& -2 \Im(\widetilde{C_1}) (\eta-\etab), \\
\Im(J_{a3})&=&-2\nabc \Im(\widetilde{C_1}) +2\Im(\widetilde{C_1})  \eta +  4 \trchb \dual \eta+2\atrchb \eta +2\atrchb \etab.
\eeaa
We therefore obtain
\beaa
\Im( Z_{a3})&=&-2\nabc \Im(\widetilde{C_1})  -2 \Im(\widetilde{C_1}) \etab +8 \nabc_3 \dual \eta.
\eeaa
Using that, see \eqref{eq:vanishing-relations-Kerr-real}, 
\beaa
\nab_3 \dual\eta+\trchb \dual \eta -\atrchb \eta= \dk^{\leq 1} \Ga_g,
\eeaa
we deduce that in order to have $\Im(Z_{a3})=0$, we need to have
\bea\label{eq:condition-on-Im-Za3}
\nabc \Im(\widetilde{C_1})  + \Im(\widetilde{C_1}) \etab =-4 (\trchb \dual \eta -\atrchb \eta)+\dk^{\leq 1} \Ga_g.
\eea
Observe that, using \eqref{eq:vanishing-relations-Kerr-real-2}, i.e.
\beaa
\nab(\atrchb)&=& -\frac 3 2  \atrchb  (\etab+ \eta) +\frac 1 2\trchb  (\dual \eta-  \dual \etab )+r^{-1} \Ga_g,
\eeaa
we obtain
\beaa
\nabc \atrchb  + \atrchb \etab &=& \nab \atrchb  + \atrchb (\etab-\ze) \\
&=&-\frac 3 2  \atrchb  (\etab+ \eta) +\frac 1 2\trchb  (\dual \eta-  \dual \etab ) + \atrchb (\etab-\ze) +r^{-1} \Ga_g\\
&=&\frac 1 2\atrchb \etab -\frac 3 2  \atrchb   \eta +\frac 1 2\trchb  (\dual \eta-  \dual \etab )  +r^{-1} \Ga_g\\
&=& \trchb  \dual \eta - \atrchb   \eta  +r^{-1} \Ga_g,
\eeaa
where we used, see \eqref{eq:vanishing-relations-Kerr-real}, that 
$\trchb (\dual \eta + \dual\etab) - \atrchb ( \etab - \eta)=r^{-1} \Ga_g$. Therefore, if $\Im(\widetilde{C_1})=-4 \atrchb$, relation \eqref{eq:condition-on-Im-Za3} is satisfied, and $Z_{a3}$ is a real one-form.

We now consider $\Re(Z_{a3})$. We compute
\beaa
\Re(L_{a3})&=& -2\trchb  (\eta-\etab )  \\
\Re(M_{a3})&=& 4 \nabc_3 \left( 2\eta+\etab \right)-   \trchb\left( 4\eta+2\etab\right)+\atrchb \dual  \left( 4\eta+2\etab \right),
\eeaa
and, recalling that $\Re(C_1)=2\trchb +\Re(\widetilde{C_1})$,
\beaa
\Re( I_{a3})&=& - 2\nabc_3(\eta-\etab )-(2\Re(C_1) -   \trchb ) (\eta-\etab )-\atrchb  \dual (\eta-\etab )\\
&=& - 2\nabc_3(\eta-\etab )-(2\Re(\widetilde{C_1}) +3  \trchb ) (\eta-\etab )-\atrchb  \dual (\eta-\etab )\\
\Re( J_{a3})&=&2\nabc_3\eta-2\nabc  \Re(C_1)  +2\Re(C_1)  \eta+  \trchb (8\etab+2\eta) -2\atrchb   \dual \etab \\
&=&2\nabc_3\eta-4\nabc \trchb  -2\nabc \Re(\widetilde{C_1})  +2(2\trchb +\Re(\widetilde{C_1}))  \eta\\
&&+  \trchb (8\etab+2\eta) -2\atrchb   \dual \etab\\
&=&2\nabc_3\eta-4\nabc \trchb  -2\nabc \Re(\widetilde{C_1})  +2\Re(\widetilde{C_1})  \eta+  \trchb (8\etab+6\eta) -2\atrchb   \dual \etab.
\eeaa
We therefore obtain
\beaa
\Re( Z_{a3})&=&- 2\nabc_3(\eta-\etab )-(2\Re(\widetilde{C_1}) +3  \trchb ) (\eta-\etab )-\atrchb  \dual (\eta-\etab )\\
&&+2\nabc_3\eta-4\nabc \trchb  -2\nabc \Re(\widetilde{C_1})  +2\Re(\widetilde{C_1})  \eta+  \trchb (8\etab+6\eta) -2\atrchb   \dual \etab\\
&& -2\trchb  (\eta-\etab )+ 4 \nabc_3 \left( 2\eta+\etab \right)-   \trchb\left( 4\eta+2\etab\right)+\atrchb \dual  \left( 4\eta+2\etab \right)\\
&&-4(2\trchb +\Re(\widetilde{C_1})) \etab\\
&=& 6\nabc_3\etab + 8 \nabc_3  \eta-4\nabc \trchb  -2\nabc \Re(\widetilde{C_1}) - 2\Re(\widetilde{C_1})   \etab \\
&& +  \trchb (3\etab-3\eta) +\atrchb \dual  \left( 3\eta+\etab \right).
\eeaa
Using the null structure equation for $\nabc_3 \etab$,  \eqref{eq:vanishing-relations-Kerr-real} and \eqref{eq:vanishing-relations-Kerr-real-2}, we obtain
\beaa
\Re( Z_{a3})&=& 6\big( -\frac{1}{2}\trchb(\etab-\eta)+\frac{1}{2}\atrchb(\dual\etab-\dual\eta)  \big) + 8 (-\trchb \eta -\atrchb \dual \eta)\\
&&-4\big( -\frac 3 2 \trchb   \left(\etab + \eta\right) -\frac 1 2 \atrchb  \left(  \dual \eta-   \dual \etab\right)+ \trchb \etab \big)  -2\nabc \Re(\widetilde{C_1}) - 2\Re(\widetilde{C_1})   \etab \\
&& +  \trchb (3\etab-3\eta) +\atrchb \dual  \left( 3\eta+\etab \right)+r^{-1} \Ga_b\\
&=&-2\nabc \Re(\widetilde{C_1}) - 2\Re(\widetilde{C_1})   \etab+ 2\trchb   \left(\etab - \eta\right) +2\atrchb(\dual\etab-3\dual\eta) +r^{-1} \Ga_b.
\eeaa

Observe that 
\beaa
\nabc \left( \frac {\atrchb^2}{ \trchb}\right) &=& \nab \left( \frac {\atrchb^2}{ \trchb}\right)-\ze \frac {\atrchb^2}{ \trchb} \\
&=& \nab \left( \frac {\atrchb^2}{ \trchb}\right)+\etab \frac {\atrchb^2}{ \trchb} \\
&=&-\frac {\atrchb^2}{ \trchb^2}\nab \trchb +\frac {2\atrchb}{ \trchb}\nab \atrchb+\etab \frac {\atrchb^2}{ \trchb}.
\eeaa
Using \eqref{eq:vanishing-relations-Kerr-real-2}, we obtain
\beaa
&&\nabc \left( \frac {\atrchb^2}{ \trchb}\right) \\
&=&-\frac {\atrchb^2}{ \trchb^2}\left(-\frac 3 2 \trchb   \left(\etab + \eta\right) -\frac 1 2 \atrchb  \left(  \dual \eta-   \dual \etab\right) \right)\\
&& +\frac {2\atrchb}{ \trchb}\left(-\frac 3 2  \atrchb  (\etab+ \eta) +\frac 1 2\trchb  (\dual \eta-  \dual \etab ) \right)+\etab \frac {\atrchb^2}{ \trchb} + r^{-2} \Ga_b\\
&=&- \frac {\atrchb^2}{ \trchb} \left(\frac 1 2\etab + \frac 3 2\eta\right)  +\frac 1 2\frac {\atrchb^3}{ \trchb^2} \left(  \dual \eta-   \dual \etab\right)  +\atrchb  (\dual \eta-  \dual \etab ) + r^{-2} \Ga_b,
\eeaa
which implies
\beaa
&&\nabc \left( \frac {\atrchb^2}{ \trchb}\right) +\etab  \frac {\atrchb^2}{ \trchb}\\
&=&  \frac {\atrchb^2}{ \trchb} \left(\frac 1 2\etab - \frac 3 2\eta\right)  +\atrchb  (\dual \eta-  \dual \etab ) +\frac 1 2\frac {\atrchb^3}{ \trchb^2} \left(  \dual \eta-   \dual \etab\right)  + r^{-2} \Ga_b.
\eeaa
Therefore, if $\Re(\widetilde{C_1})=-2\frac {(\atrchb)^2}{ \trchb} $, we have
\beaa
&&\Re( Z_{a3})\\
&=&4\Big(\nabc \left(\frac {(\atrchb)^2}{ \trchb}\right) +\frac {(\atrchb)^2}{ \trchb}   \etab\Big)+ 2\trchb   \left(\etab - \eta\right) +2\atrchb(\dual\etab-3\dual\eta) +r^{-1} \Ga_b\\
&=&4\Big( \frac {\atrchb^2}{ \trchb} \left(\frac 1 2\etab - \frac 3 2\eta\right)  +\atrchb  (\dual \eta-  \dual \etab ) +\frac 1 2\frac {\atrchb^3}{ \trchb^2} \left(  \dual \eta-   \dual \etab\right)\Big)\\
&&+ 2\trchb   \left(\etab - \eta\right) +2\atrchb(\dual\etab-3\dual\eta) +r^{-1} \Ga_b\\
&=&2 \frac {\atrchb^2}{ \trchb} \left(\etab -  3 \eta\right)  +2\frac {\atrchb^3}{ \trchb^2} \left(  \dual \eta-   \dual \etab\right)+ 2\trchb   \left(\etab - \eta\right) -2\atrchb(\dual\etab+\dual\eta) +r^{-1} \Ga_b.
\eeaa
Evaluating the above to $e_a=e_1$ and using that $\eta_1-\etab_1=\Ga_g$, $\dual \eta_1+\dual\etab_1=\Ga_g$, we obtain
\beaa
\Re( Z_{13})&=&-4 \frac {\atrchb^2}{ \trchb}  \etab_1  -4\frac {\atrchb^3}{ \trchb^2}   \dual \etab_1+r^{-1} \Ga_b\\
&=&-4 \frac {\atrchb^2}{ \trchb} \Big(  -\frac{a^2\sin\th \cos\th }{|q|^3}  +\frac {a\cos\th}{ r}  \frac{ar\sin\th }{|q|^3} \big)+r^{-1} \Ga_b=r^{-1} \Ga_b.
\eeaa
Evaluating the above to $e_a=e_2$ and using that $\eta_2+\etab_2=\Ga_g$, $\dual \eta_2-\dual \etab_2=\Ga_g$
\beaa
\Re( Z_{23})&=&8 \frac {\atrchb^2}{ \trchb} \etab_2   + 4\trchb  \etab_2 -4\atrchb \dual\etab_2 +r^{-1} \Ga_b\\
&=&8 \frac {a\cos\th}{ r}  \frac{2a\Delta\cos\th}{|q|^4} \frac{ar\sin\th }{|q|^3}   +4\frac{2r\Delta}{|q|^4}  \frac{ar\sin\th }{|q|^3} -4\frac{2a\Delta\cos\th}{|q|^4} \frac{a^2\sin\th \cos\th}{|q|^3} +r^{-1} \Ga_b\\
&=&\frac{8a\sin\th \De}{|q|^5} +r^{-1} \Ga_b,
\eeaa
as stated.
\end{proof}

We therefore deduce the following intermediate theorem for the gRW equation.

   \begin{theorem}\label{main-theorem-intermediate} 
For the following choices of complex scalar functions $C_1$, $C_2$, i.e.
\bea\label{finalchoicefordefinition-C}
C_1&=&2\trchb - 2\frac {\atrchb^2}{ \trchb}  -4 i \atrchb, \qquad C_2=\frac 1 2 \trchb^2 + \widetilde{C_2},
\eea
where $\widetilde{C_2}$ is any complex function satisfying $\widetilde{C_2}=O(|a|r^{-3})$, the invariant symmetric traceless $2$-tensor $\qf \in \sk_2(\CCC)$ satisfies the equation
 \bea\label{wave-equation-qf-intermediate}
 \squared_2 \qf   -i \frac{4 a\cos\th}{|q|^2} \nab_\T \qf   - V_1  \qf &=&   \widetilde{L_{\qf}[A]} + \err[\squared_2 \qf],
 \eea
 where
 \begin{itemize}
 \item The potential $V_1$ is a \textbf{real} scalar function given by $V_1=-\trch\trchb+ O(\frac{|a|}{r^4})$, 
 \item $ \widetilde{L_{\qf}[A ]}$ is  a  linear second order   operator  in $A$ vanishing for zero angular momentum, given by
 \beaa
   - \widetilde{L_{\qf}[A]}&=&  q \ov{q}^3 \Big[Z_{43} \  \nabc_4\nabc_3A+ Z_{23}\nabc_2 \nabc_3A \\
&&+Z_4 \  \nabc_4A+ Z_3\nabc_3A+Z_a \nabc_a A +Z_0  A \Big],
  \eeaa
  where in the outgoing frame
  \beaa
  Z_{43}=\frac{4a^2\De^2\cos^2\th}{r^2|q|^6}, \qquad Z_{23}=\frac{8a\sin\th \De}{|q|^5} 
  \eeaa
and $Z_4$, $Z_3$, $Z_a$, $Z_0$ have the following fall-off in $r$:
\beaa
Z_4= O(|a| r^{-4}), \qquad Z_3=O(|a| r^{-4}), \qquad Z_a=O(|a| r^{-4}), \qquad Z_0=O(|a| r^{-5}).
\eeaa

 \item $\err[\squared_2 \qf]$ is  the nonlinear correction term, which is given schematically by  the expression
   \beaa
 \err[\squared_2 \qf]&=& r^2 \frak{d}^{\leq 3} (\Ga_g \c (A, B))+ \nab_3 (r^2 \frak{d}^{\leq 2}( \Ga_b \c (A, B)))\\
&&+\frak{d}^{\leq 1} (\Ga_g \c \qf) + r^{3}\dk^{\leq 2} \big( \Ga_b \c \Ga_g \c \Ga_g\big).
 \eeaa
 \end{itemize}
 \end{theorem}

Observe that, by expanding the conformal derivatives, the linear second order operator $ \widetilde{L_{\qf}[A ]}$ can be written as
 \beaa
   - \widetilde{L_{\qf}[A]}&=&  q \ov{q}^3 \Big[Z_{43} \  \nab_4\nab_3A+Z_{23} \nab_2 \nab_3A +W_4 \  \nab_4A+ W_3\nab_3A+W_a \nab_a A +W_0  A \Big],
  \eeaa
where $W_4$, $W_3$, $W_a$, $W_0$ have the same fall-off in $r$ as $Z_4$, $Z_3$, $Z_a$, $Z_0$ respectively.

We now combine parts of the potential with the lower order terms $ \widetilde{L_{\qf}[A ]}$.
In particular, we can write
\beaa
&&-\frac{|q|^2}{\Delta} Z_{43}  \ \qf + \widetilde{L_{\qf}[A]}\\
&=& -\frac{|q|^2}{\Delta} Z_{43}  q\ov{q}^3 \big[\nab_3 \nab_3 A + O(r^{-1}) \nab_3 A + O(r^{-2}) A \big] \\
&&-q \ov{q}^3 \Big[Z_{43} \  \nab_4\nab_3A+Z_{23} \nab_2 \nab_3A +W_4 \  \nab_4A+ W_3\nab_3A+W_a \nab_a A +W_0  A \Big]\\
&=& -q\ov{q}^3    \big[Z_{43} \big( \frac{|q|^2}{\Delta}  \nab_3  + \nab_4\big) \nab_3A+ Z_{23} \nab_2 \nab_3A\\
&&+ W_4 \  \nab_4A+ W_3\nab_3A+W_a \nab_a A +W_0  A \Big].
\eeaa

Using that with the outgoinng normalization of $(e_3, e_4)$,
\beaa
 \frac 1 2 \left( \frac{\De}{r^2+a^2} e_4+\frac{|q|^2}{r^2+a^2}  e_3\right)=\T+\frac{a}{r^2+a^2}\Z,
\eeaa
and $e_2=\frac{a\sin\th}{|q|}\nab_\T +\frac{1}{|q|\sin\th}\nab_\Z$, 
we obtain
\beaa
-\frac{|q|^2}{\Delta} Z_{43}  \ \qf + \widetilde{L_{\qf}[A]}&=& -q\ov{q}^3    \big[ \big( 2Z_{43}\frac{r^2+a^2}{\De}+ Z_{23}\frac{a\sin\th}{|q|}\big) \nab_\T \nab_3 A \\
&&+ \big( 2Z_{43} \frac{a}{\Delta}+ Z_{23}\frac{1}{|q|\sin\th} \big)  \nab_\Z \nab_3A\\
&&+ W_4 \  \nab_4A+ W_3\nab_3A+W_a \nab_a A +W_0  A \Big]\\
&=& -q\ov{q}^3    \big[ \frac{8a^2 \De}{r^2|q|^4}\nab_\T \nab_3 A +\frac{8a  \De }{r^2|q|^4} \nab_\Z \nab_3A\\
&&+ W_4 \  \nab_4A+ W_3\nab_3A+W_a \nab_a A +W_0  A \Big],
\eeaa
where we used the values of $Z_{43}$ and $Z_{23}$ obtained in Proposition \ref{prop:potential-real}. Combining the above with equation \eqref{wave-equation-qf-intermediate}, we finally obtain
\beaa
 \squared_2 \qf   -i \frac{4 a\cos\th}{|q|^2} \nab_\T \qf   -\big( V_1 +\frac{|q|^2}{\Delta} Z_{43} \big) \ \qf &=&L_\qf[A]+ \err[\squared_2 \qf],
 \eeaa
 where $L_\qf[A]:=-\frac{|q|^2}{\Delta} Z_{43}  \ \qf  + \widetilde{L_{\qf}[A]}$, is given by 
 \beaa
 -L_\qf[A]&=& q\ov{q}^3    \big[ \frac{8a^2 \De}{r^2|q|^4}\nab_\T \nab_3 A +\frac{8a  \De }{r^2|q|^4} \nab_\Z \nab_3A\\
&&+ W_4 \  \nab_4A+ W_3\nab_3A+W_a \nab_a A +W_0  A\big].
 \eeaa
By defining
\beaa
V&:=& V_1+\frac{|q|^2}{\Delta} Z_{43}\\
&=& \frac{4}{|q|^2}\frac{r^2-2mr+2a^2}{r^2}-\frac{4a^2\cos^2\th}{|q|^6}( r^2+6mr+a^2\cos^2\th ),
\eeaa
we finally completed the proof of Theorem \ref{MAIN-THEOREM-PART1}.


\section{Proof of Proposition \ref{PROP:TEUK-AB}}
\label{proof-teukolsky-Ab}



\subsection{Preliminaries}


We make use of  our  gauge conditions $\Xi=0$, $\Hbc=0$   to derive  the following   linearized
 equations.

\begin{lemma}\lab{Lemma:DDc{trX}}
We have 
\bea\label{eq:DDc-ov-trX}
\bsplit
 \DDc \ov{\tr X}&=     2 i\Im(\tr  X) (H-\Hc)+ r^{-1}\Ga_g,\\
 \DDc \tr X&= -2\tr X \Hb + r^{-1} \Ga_g.
\end{split}
\eea
Also,
\bea\label{eq:DDc-ov-trXb}
\bsplit
 \DDc \ov{\tr \Xb }&=    2 i\Im(\tr  \Xb) \Hb  + r^{-1}\Ga_g,\\
 \DDc \tr \Xb&=-2\tr\Xb (H- \Hc)+r^{-1} \Ga_g.
\end{split}
\eea
\end{lemma}
\begin{proof}      
 Since $\tr X=\frac{2\De\ov{q}}{|q|^4}+\Ga_g $, 
       and $Z=\frac{aq}{|q|^2}\Jk+\Ga_g$, we have
       \beaa
\DDc (\ov{\tr X})&=& \DD \ov{\tr X} + \ov{\tr X} Z= \DD \left(\frac{2\De q}{|q|^4}+\Ga_g \right) + \left(\frac{2\De q}{|q|^4}+\Ga_g \right) \left(\frac{aq}{|q|^2}\Jk+\Ga_g \right)\\
&=& \DD \left(\frac{2\De q}{|q|^4} \right) + \frac{2\De q}{|q|^4}\frac{aq}{|q|^2}\Jk+r^{-1} \Ga_g.
\eeaa
       Since  $\DD q=- a \Jk + r\Ga_g$, $\DD \ov{q}=a \Jk + r\Ga_g$, $\DD(\De) = r^2 \Ga_g$, we have
 \beaa
\DD\left(  \frac{2 q \De}{|q|^4} \right)&=&2\De \DD\left(\frac{1} { \ov{q}^2q}\right)+ r^{-1} \Ga_g=  2\De\big( - \ov{q}^{-2}  \,  q^{-2} \DD(q) - 2 \ov{q}^{-3}\,  q^{-1} \DD( \ov{q}) \big)+ r^{-1} \Ga_g\\
&=&-\frac{2\De}{|q|^4} \Big(\DD(q)+\frac{ 2 q }{\ov{q}} \DD(\ov{q}) \Big)+ r^{-1} \Ga_g=\frac{2a\De}{|q|^4} \Big(  \Jk-\frac{ 2 q }{\ov{q}}  \Jk\Big)+ r^{-1} \Ga_g\\
&=&
\frac{2aq\De}{|q|^4} \Big(  \frac{\ov{q}- 2 q }{|q|^2} \Big)\Jk+ r^{-1} \Ga_g,
\eeaa
which gives, using that $H=\Hc+ \frac{a q}{|q|^2} \Jk $,
\beaa
\DDc (\ov{\tr X})&=&\frac{2aq\De}{|q|^4} \Big(  \frac{\ov{q}- 2 q }{|q|^2} \Big)\Jk+ \frac{2\De q}{|q|^4}\frac{aq}{|q|^2}\Jk+r^{-1} \Ga_g\\
&=&\frac{2aq\De}{|q|^4} \Big(  \frac{\ov{q}-  q }{|q|^2} \Big)\Jk+r^{-1} \Ga_g=2i \Im(\tr X)(H-\Hc)+r^{-1} \Ga_g.
\eeaa
   Similarly, we have
\beaa
\DDc (\tr X)&=& \DD \tr X + \tr X Z= \DD \left(\frac{2\De\ov{q}}{|q|^4}+\widecheck{\tr X} \right) + \left(\frac{2\De\ov{q}}{|q|^4}+\widecheck{\tr X} \right) \left(\frac{aq}{|q|^2}\Jk+\widecheck{Z} \right)\\
&=& \DD \left(\frac{2\De\ov{q}}{|q|^4} \right) + \frac{2\De\ov{q}}{|q|^4}\frac{aq}{|q|^2}\Jk+\DD \widecheck{\tr X}+\tr X \Zc+\widecheck{\tr X} Z,
\eeaa
and 
 \beaa
\DD\left(  \frac{2 \ov{q} \De}{|q|^4} \right)&=&2\De \DD\left(\frac{1} { q^2 \ov{q}}\right)+  \frac{2 \ov{q} }{|q|^4}\DD (\De)=  2\De\big( - q^{-2}  \,  \ov{q}^{-2} \DD(\ov{q}) - 2 q^{-3}\,  \ov{q}^{-1} \DD( q) \big)+ r^{-1} \Ga_g\\
&=&-\frac{2\De}{|q|^4} \Big(\DD(\ov{q})+\frac{ 2 \ov{q} }{q} \DD(q) \Big)+ r^{-1} \Ga_g=
-\frac{2a\De}{|q|^4}\left(\Jk- \frac{ 2 \ov{q} }{q}\Jk\right) + r^{-1} \Ga_g\\
&=& -\frac{2a \ov{q} \De}{|q|^4}\,\frac{q-2\ov{q} }{|q|^2}\Jk+ r^{-1} \Ga_g,
\eeaa
which gives
\beaa
\DDc (\tr X)&=&\frac{4a \ov{q}^2 \De}{|q|^6}\, \Jk +r^{-1} \Ga_g=-\frac{4 \De\ov{q}}{|q|^4}\,\left( -\frac{ a}{q}\Jk\right) +r^{-1} \Ga_g\\
&=&-2\tr X \Hb +r^{-1}\Ga_g,
\eeaa
as stated.

We next calculate, using that $\tr \Xb=-\frac{2}{\ov{q}}+ \Ga_g$, and  $H=\Hc+ \frac{a q}{|q|^2} \Jk $
\beaa
\DDc (\tr \Xb) &=& \DD \tr \Xb -\tr \Xb Z= \DD\left(-\frac{2}{\ov{q}}+ \Ga_g \right) -\left(-\frac{2}{\ov{q}}+ \Ga_g \right) \left(\frac{aq}{|q|^2}\Jk+\Ga_g \right)\\
&=&- \DD\left(\frac{2}{\ov{q}} \right) +\frac{2}{\ov{q}}\frac{aq}{|q|^2}\Jk+r^{-1} \Ga_g=\frac{2}{ \ov{q}^{2}} \DD\ov{q} +\frac{2}{\ov{q}}\frac{aq}{|q|^2}\Jk+r^{-1} \Ga_g\\
&=&\frac{2a}{ \ov{q}^{2}} \Jk +\frac{2a}{\ov{q}^2}\Jk+r^{-1} \Ga_g=-2\tr\Xb (H -\Hc)+r^{-1} \Ga_g.
\eeaa
Similarly, 
\beaa
\DDc (\ov{\tr \Xb}) &=& \DD \ov{\tr \Xb} -\ov{\tr \Xb} Z= \DD\left(-\frac{2}{q}+ \Ga_g \right) -\left(-\frac{2}{q}+ \Ga_g \right) \left(\frac{aq}{|q|^2}\Jk+\Ga_g \right)\\
&=&- \DD\big(\frac{2}{q} \big) +\frac{2}{q}\frac{aq}{|q|^2}\Jk+r^{-1} \Ga_g=\frac{2}{ q^{2}} \DD q +\frac{2a}{|q|^2}\Jk+r^{-1} \Ga_g\\
&=&-\frac{2a}{ q^{2}}  \Jk  +\frac{2a}{|q|^2}\Jk+r^{-1} \Ga_g= \big(\frac{2}{\ov{q}}-\frac{2}{q} \big) \frac{a}{q} \Jk+r^{-1} \Ga_g=2 i\Im(\tr  \Xb) \Hb  + r^{-1}\Ga_g,
\eeaa
as stated.
\end{proof}

\begin{lemma}
\lab{Lemma:nabc_4Hb}
In the frame  for which $ \Hb= -\frac{a}{q}\Jk$ we have
\beaa
\nab_4\Hb+\tr X\Hb=r^{-2} \Ga_g.
\eeaa
\end{lemma}
\begin{proof}
Recall that we have, see   Definitions \ref{def:renormalizationofallnonsmallquantitiesinPGstructurebyKerrvalue:3} and  \ref{definition.Ga_gGa_b},
\beaa
\nab_4\Jk =-\frac{\De \ov{q}}{|q|^4}\Jk+ r^{-1} \Ga_g, \,\, \, 
\nab_4 q =\frac{\De}{|q|^2} +\Ga_g, \,\,\,  \tr X=\frac{2\ov{q}\De }{|q|^4}+\Ga_g, \,\,\, \nabc_4 q=\frac 1 2 \tr X q + r\Ga_g.
\eeaa
Hence
\beaa
\nab_4 \Hb&=& \nab_4 \left( -\frac{a}{q}\Jk\right)= \frac{a\nab_4 q }{q^2}\Jk  -\frac{a}{q}\left(-\frac{\De \ov{q}}{|q|^4}\Jk+ r^{-1} \Ga_g\right), \\
&=&\Big(\frac{a\De}{q^2|q|^2}+ \frac{a\ov{q} \De}{q|q|^4 }\Big)\Jk+ r^{-2} \Ga_g=\frac{a\De}{|q|^2 }      \Big(\frac{1}{q^2}+ \frac{\ov{q} }{q|q|^2 }\Big)\Jk+ r^{-2} \Ga_g\\
&=&    2 \frac{a\De}{|q|^2 }    \frac{\ov{q}^2 }{|q|^4}\Jk +           r^{-2} \Ga_g \\
&=&-\frac{2\De \ov{q} }{|q|^4} \Hb+r^{-2} \Ga_g =-\tr X  \Hb+ r^{-2} \Ga_g,
\eeaa
as stated.
\end{proof}

\begin{proposition}\label{prop:preliminaries-qfb-null-structure}
The following  equations hold true.
\bea\label{eq:linearized-nabc4Hc}\label{equation:nabc_4Hc}
\bsplit
\nabc_4\Hc+\frac{1}{2}\ov{\tr X}\Hc&=  - B-\frac 1 2 \Xh\c \ov{\Hc}  + r^{-2} \Ga_g,\\
\nabc_4 \Xib&=-\frac 1 2 \ov{\tr X}\Hc-\Bb -\frac 1 2 \Xbh \c \Hc  + r^{-2}\Ga_b.
\end{split}
\eea
Also
\bea
\lab{eq:linearizedCoazzi-qfb}
\bsplit
 \DDc(\ov{\trXc})&=   2 i \Im( \tr \Xb )\Hc+2  B    + \ov{\DDc}\c\Xh + r^{-1} \Ga_g, \\
 \DDc\ov{\trXbc}&=  -   2 \Bb    + \ov{\DDc}\c\Xbh+2  i\Im(\tr X)\Xib + r^{-1} \Ga_g .
 \end{split}
 \eea
\end{proposition}
\begin{proof}
Recall that $H$ verifies the equation
\beaa
\nabc_4H  &=&  -\frac{1}{2}\ov{\tr X}(H-\Hb) -\frac{1}{2}\Xh\c(\ov{H}-\ov{\Hb}) -B.
\eeaa
Writing $H=\Hc+ \frac{a q}{|q|^2} \Jk $ and $\Hb=-\frac{a\ov{q}}{|q|^2} \Jk$ we deduce
\beaa
\nabc_4\Hc+\nabc_4\left( \frac{a q}{|q|^2} \Jk\right)=  -\frac{1}{2}\ov{\tr X}\Hc 
 -\frac{1}{2}\ov{\tr X}\left(\frac{a q}{|q|^2} \Jk  + \frac{a\ov{q}}{|q|^2} \Jk\right)- \frac 1 2 \Xh \c \ov{\Hc}  - B+ r^{-2} \Ga_g,
\eeaa
or,
\beaa
\nabc_4\Hc+\frac{1}{2}\ov{\tr X}\Hc&=&  - B-\frac 1 2 \Xh\c \ov{\Hc}  + r^{-2} \Ga_g- E\\
E&=&\nabc_4\left( \frac{a q}{|q|^2} \Jk\right)+ \frac{1}{2}\ov{\tr X}\left(\frac{a q}{|q|^2} \Jk + \frac{a\ov{q}}{|q|^2} \Jk\right).
\eeaa
Using
\beaa
\nab_4\Jk &=&-\frac{\De \ov{q}}{|q|^4}\Jk+ r^{-1} \Ga_g, \qquad 
\nab_4 q =\frac{\De}{|q|^2} +\Ga_g, \qquad \tr X=\frac{2\ov{q}\De }{|q|^2}+\Ga_g,
\eeaa
 we   calculate 
\beaa
E    &=& \nabc_4 \big( \frac{a}{\ov{q} }\big)\Jk+ \frac{a}{\ov{q} }\Big(-\frac{\De \ov{q}}{|q|^4}\Jk+ r^{-1} \Ga_g\Big)+ \Big(\frac{ q \De }{|q|^2}+\Ga_g\Big) \big(\frac{a q}{|q|^2} +\frac{a\ov{q}}{|q|^2} \big) \Jk \\
    &=&\Big(-\frac{a}{\ov{q}^2} \frac{\De}{|q|^2} -\frac{a\De}{|q|^4} + \frac{ q \De }{|q|^2}\big(\frac{a q}{|q|^2} +\frac{a\ov{q}}{|q|^2} \big)\Big)\Jk+ r^{-2} \Ga_g= r^{-2} \Ga_g.
\eeaa
Similarly
the equation
\beaa
   \nabc_3\Hb -\nabc_4\Xib &=&  -\frac{1}{2}\ov{\tr\Xb}(\Hb-H) -\frac{1}{2}\Xbh\c(\ov{\Hb}-\ov{H}) +\Bb,
 \eeaa
takes the form
 \beaa
 \bsplit
 \nabc_4 \Xib&=-\frac 1 2 \ov{\tr X}\Hc-\Bb + r^{-2}\Ga_b -E\\
 E&=  \nabc_3 \left(\frac{a}{q} \Jk\right)+\frac 1 2 
 \ov{\tr \Xb} \left(\frac{a}{q} \Jk+  \frac{a q}{|q|^2} \Jk\right).
 \end{split}
 \eeaa
 Using the relations 
\beaa
 \tr \Xb =-\frac{2}{\ov{q}} +\Ga_g,\,\,  \nabc_3 \Jk=\frac{1}{\ov{q}}\Jk+r^{-1} \Ga_b, \,\, e_3 (q) = -1 + r\Ga_b,
\eeaa
we  deduce
\beaa
 E &=&-\frac{a}{q^2} \big(-1+ r\Ga_b\big)\Jk+\frac{a}{q} \left( \frac{1}{\ov{q}}\Jk+r^{-1} \Ga_b\right)
 +\frac 1 2 \left( -\frac{2}{q} +\Ga_g\right) \big(\frac{a}{q}   +\frac{a q}{|q|^2} \big)\Jk= r^{-2} \Ga_b,
\eeaa
 and thus
 \beaa
 \nabc_4 \Xib&=&-\frac 1 2 \ov{\tr X}\Hc-\Bb -\frac 1 2 \Xbh \c \Hc  + r^{-2}\Ga_b,
 \eeaa
  as stated.
  To derive the   equations for $\trXbc, \trXc $ we start with the corresponding  Codazzi equation  written in the form\footnote{Recall that $\DDc  \ov{\tr\Xb}= \DDc  \ov{\tr\Xb}-   Z  \ov{\tr\Xb} $, $ \DDc\ov{\tr X} =  \DD\ov{\tr X}+\ov{\tr X }Z$.}
  \bea
  \bsplit
   \DD\ov{\tr X}+\ov{\tr X}Z&=   2 i\Im(\tr  X) H     +2  B    + \ov{\DDc}\c\Xh +2 i\Im(\tr X)\Xi+ r^{-2} \Ga_g,  \\
    \DD\ov{\tr\Xb}-\ov{\tr\Xb}Z&=   2 i \Im( \tr \Xb )\Hb-   2 \Bb    + \ov{\DDc}\c\Xbh+2  i\Im(\tr X)\Xib + r^{-2} \Ga_g  .
    \end{split}
     \eea       
     For the first equation, since $\tr X=\frac{2\De\ov{q}}{|q|^4}+\Ga_g $, 
      $2 i \Im(\tr X) = \tr X-\ov{\tr X}=-\frac{2\De(q-\ov{q})  }{|q|^4} +\Ga_g $, $H=\Hc+ \frac{a q}{|q|^2} \Jk$
       and $Z=\frac{aq}{|q|^2}+\Ga_g$
     we write 
     \beaa
      &&  \DD\ov{\tr X}+\ov{\tr X}Z-   2 i\Im(\tr  X) H \\
      &=&\DD(\trXc + \frac{2 q \De}{|q|^4}) + (\trXc + \frac{2 q \De}{|q|^4})( \frac{aq}{|q|^2}\Jk +\Ga_g)
        -  2i  \Im(\tr X)(\Hc+ \frac{aq}{|q|^2}\Jk\big) \\
        &=&\DD\trXc  -   2 i\Im(\tr  X) \Hc +\DD\big(  \frac{2 q \De}{|q|^4} \big)+  \frac{2 q \De}{|q|^4}  \frac{aq}{|q|^2}\Jk -   2 i\Im(\tr  X) \Hc +\frac{2\De(q-\ov{q})  }{|q|^4} \frac{aq}{|q|^2}\Jk \\
        &&+ r^{-2} \Ga_g \\
        &=&\DD\trXc-   2 i\Im(\tr  X) \Hc + \DD\big(  \frac{2 q \De}{|q|^4} \big)+\frac{2 a q\De}{|q|^4}  \Jk\big( \frac{q}{|q|^2}+
        \frac{q-\ov{q} } {|q|^2}\big)+  r^{-1} \Ga_g.
          \eeaa
      Since  $\DD q=- a \Jk + r\Ga_g$, $\DD \ov{q}=a \Jk + r\Ga_g$ , 
       \beaa
       \bsplit
        \DD\big(  \frac{2 q \De}{|q|^4} \big)&= 2 \De \DD\big(\frac{1} { q \ov{q}^2}\big)= 2 \De\big( - q^{-2}   \ov{q}^{-2} \DD(q) - 2 q^{-1}  \ov{q}^{-3} \DD(\ov{q}) \big) \\
       &=-\frac{2\De}{|q|^4} \Big(\DD(q)+\frac{ 2 q}{\ov{q}} \DD(\ov{q}) \Big)=\frac{2a\De}{|q|^4}\Big( \Jk - \frac{ 2 q}{\ov{q}}\Jk \Big)= \frac{2a q \De}{|q|^4}\Big(\frac{\ov{q}}{|q|^2}- \frac{2 q}{|q|^2 }  \big)\Jk.
       \end{split}
       \eeaa
       Hence
       \beaa
        \DD\ov{\tr X}+\ov{\tr X}Z-   2 i\Im(\tr  X) H&=&    \DD(\ov{\trXc}) -   2 i\Im(\tr  X) \Hc +r^{-1} \Ga_g\\
        &=& \DDc(\ov{\trXc}) -   2 i\Im(\tr  X) \Hc+ r^{-1}\Ga_g,
\eeaa
and therefore 
\beaa
  \DDc(\ov{\trXc})&=   2 i \Im( \tr \Xb )\Hc+2  B    + \ov{\DDc}\c\Xh + r^{-1} \Ga_g  .
\eeaa
The second equation in  \eqref{eq:linearizedCoazzi-qfb} is proved in the same manner.
 \end{proof}

    \begin{lemma}
    \lab{remark:DDcP-DDcov{P}}
    The following identity  holds true
    \bea\label{eq:error-DDP}
    \begin{split}
    \DDc P&=- 3 P \Hb+\DDc \Pc + r^{-3} \Ga_g, \\
    \DDc \ov{P}&=- 3  \ov{P}(H-\Hc) +\DDc \ov{\Pc } + r^{-3} \Ga_g.
    \end{split}
    \eea
    Also,
\beaa
 \nabc_4P+ \frac{3}{2}\tr X P  &=&\nab_4 \Pc+\frac 3 2 \tr X \Pc+ r^{-3}\Ga_g.
 \eeaa
    \end{lemma}
    \begin{proof}
    We write  $P=-\frac{2m}{q^3} +\Pc$.  Thus, since  
$\DD(q) =- a \Jk+ r\Ga_g$ and $\Hb=-\frac{a}{q}\Jk$,
\beaa
\DDc P=\DD\big(-\frac{2m}{q^3} +\Pc\big)= \frac{6m}{q^4}\DD(q)+\DD  \Pc  = \frac{6m}{q^3}\Hb +\DD \Pc +r^{-3} \Ga_g.
\eeaa
Similarly, since  $\DD \ov{q}=a \Jk + r\Ga_g$ and  $H=\Hc+ \frac{a q}{|q|^2} \Jk $,
\beaa
\DDc \ov{P}&=&  \frac{6m}{\ov{q}^4}\DD(\ov{q} )+\DD \ov{ \Pc}= \frac{6m}{\ov{q}^4}(a \Jk + r\Ga_g) +\DD \ov{ \Pc}=-3 \ov{P}(H-\Hc) +\DD \ov{ \Pc}+ r^3 \Ga_g.
\eeaa
 
Writing $P=\frac{2m}{q^3} +\Pc$, $ \nab_4 q =\frac{\De}{|q|^2} +\Ga_g $  and  $\tr X=\frac{2\ov{q}\De }{|q|^2}+\Ga_g $ we derive
    \beaa
     \nabc_4P+ \frac{3}{2}\tr X P&=&\nab_4\big(\frac{2m}{q^3} +\Pc\big)+  \frac{3}{2}\tr X \big(\frac{2m}{q^3} +\Pc\big)\\
     &=&\nab_4 \Pc+\frac 3 2 \tr X \Pc -\frac{6m}{q^4} \nabc_4 q +\frac 3 2 \big(  \frac{2\ov{q}\De }{|q|^2}+\Ga_g \big)\frac{2m}{q^3} \\
     &=& \nab_4 \Pc+\frac 3 2 \tr X \Pc  -\frac{6m}{q^4} \big( \frac{\De}{|q|^2} +\Ga_g\big)+\frac{6m}{q^3}  \big(\frac{\ov{q}\De }{|q|^2}+\Ga_g \big)\\
     &=&\nab_4 \Pc+\frac 3 2 \tr X \Pc+ r^{-3}\Ga_g,
    \eeaa
as stated.
    \end{proof} 
    
\begin{proposition}\label{prop:preliminaries-qfb-bianchi-lin}
The following  linearized Bianchi equations  hold true.
\bea
\bsplit
 \nabc_4\Bb + \tr X\Bb   &= -  \DDc \Pc +\ov{B}\c \Xbh + r^{-3} \Ga_g\\
 \nabc_4\Pc + \frac{3}{2}\tr X \Pc &= \frac{1}{2}\DDc\c \ov{B}+ \Hb \c\ov{B}  -\frac{1}{4}\Xbh\c \ov{A}
 +r^{-3} \Ga_g.
 \end{split}
\eea
\end{proposition}
\begin{proof}
In view of  Lemma \ref{remark:DDcP-DDcov{P}},   $ \DDc P=- 3 P \Hb+\DDc \Pc + r^{-3} \Ga_g $, 
       the equation
    \beaa
    \nabc_4\Bb +\tr X\Bb
     &=& -  \DDc P     +\ov{B}\c \Xbh-3P\Hb -\frac{1}{2}\Ab\c\ov{\Xi}
    \eeaa
    becomes
     \beaa
    \nabc_4\Bb + \tr X\Bb   &=& -  \DDc \Pc +\ov{B}\c \Xbh + r^{-3} \Ga_g.
    \eeaa
    Using  equation $ \nabc_4P+ \frac{3}{2}\tr X P =\nab_4 \Pc+\frac 3 2 \tr X \Pc+ r^{-3}\Ga_g$, the equation
     \beaa
    \nabc_4P -\frac{1}{2}\DDc\c \ov{B} &=& -\frac{3}{2}\tr X P + \Hb \c\ov{B} -\ov{\Xi}\c\Bb -\frac{1}{4}\Xbh\c \ov{A}
\eeaa 
becomes
\beaa
 \nabc_4\Pc + \frac{3}{2}\tr X \Pc &=& \frac{1}{2}\DDc\c \ov{B}+ \Hb \c\ov{B}  -\frac{1}{4}\Xbh\c \ov{A}
 +r^{-3} \Ga_g.
\eeaa
as stated.   
\end{proof}


\subsection{Proof of Proposition \ref{PROP:TEUK-AB}}


We start with the Bianchi identity
\beaa
\Ab_4&=& -\frac 1 2 \DDc\hot \Bb -  2\Hb    \hot \Bb -3 P\Xbh.
\eeaa
To the above, we apply the operator $\nabc_3+ \big( 2 \ov{\tr \Xb}+\frac 1 2 \tr \Xb\big)$ and we deduce
\beaa
- \nabc_3\Ab_4-\big( 2 \ov{\tr \Xb}+\frac 1 2 \tr \Xb\big) \Ab_4&=&\frac 1 2  \nabc_3 \DDc\hot \Bb+\frac 1 2 \big( 2 \ov{\tr \Xb}+\frac 1 2 \tr \Xb\big) \DDc\hot \Bb\\
 &&+2\nabc_3(\Hb    \hot \Bb)+ 2\big( 2 \ov{\tr \Xb}+\frac 1 2 \tr \Xb\big)\Hb    \hot \Bb\\
 &&+3\nabc_3(P\Xbh)+ 3\big( 2 \ov{\tr \Xb}+\frac 1 2 \tr \Xb\big)P\Xbh\\
&=& I+J+K.
\eeaa

{\bf Step 1:   Calculation of $I$.} Using the commutation formula \eqref{eq:comm-nabc4nabc3DDchot-precise} applied to $\Bb$ of signature $-1$, we have
\beaa
 [ \nabc_3 , \DDc \hot] \Bb   &=&-\frac 1 2 \tr \Xb \DDc \hot  \Bb+H\hot\nabc_3 \Bb + \Xib \hot \nabc_4 \Bb \\
 &&-2\Bb \hot \Bb- \tr X \Xib \hot  \Bb-\frac 1 2\Xbh \c \ov{\DDc} \Bb\\
 &&+\frac12\Xbh (\ov{H}\c \Bb) +\frac{1}{2}(\widehat{\Xb}\c\ov{H}) \hot \Bb +\Xh \c \Xib \c \Bb,
\eeaa
which can be written as
\beaa
 [ \nabc_3 , \DDc \hot] \Bb   &=&-\frac 1 2 \tr \Xb \DDc \hot  \Bb+H\hot\nabc_3 \Bb\\
 && + \Xib \hot \nabc_4 \Bb -2\Bb \hot \Bb- \tr X \Xib \hot  \Bb-\frac 1 2\Xbh \c \ov{\DDc} \Bb+  r^{-2} \dk^{\leq 1}(\Ga_g\c\Ga_b)\\
 &&+ (\widehat{\Xb}\c\ov{\Hc}) \Bb,
\eeaa
where we can write $(\Xh \c \Xib) \Bb= r^{-2} \dk^{\leq 1}(\Ga_g\c\Ga_b)$.
In view of the equation for $\nabc_4 \Bb$, $\Xi=0$  and, see \eqref{eq:error-DDP},  $\DDc P= -3P \Hb + r^{-2} \dk^{\le 1} \Ga_g $ we deduce
\beaa
\Xib\hot  \nabc_4 \Bb&=&\Xib\hot\big( -\DDc P  -\tr X\Bb+\ov{B}\c \Xbh-3P\Hb\big)\\
&=& -\tr X \Xib\hot \Bb + r^{-2} \dk^{\le 1}(\Ga_g\c \Ga_b)+ ( \Ga_b\c \Ga_b) \c B .
\eeaa
This gives
         \beaa
          \nabc_3 \DDc\hot \Bb&=&   \DDc\hot \nabc_3\Bb- \frac 1 2 \tr \Xb \DDc \hot \Bb  + H \hot \nabc_3 \Bb\\
         &&+\err_1+ r^{-2} \dk^{\le 1} (\Ga_g \c \Ga_b),
         \eeaa
         where
         \beaa
         \err_1&=& -2\Bb \hot \Bb- 2\tr X \Xib \hot  \Bb-\frac 1 2\Xbh \c \ov{\DDc} \Bb+ (\widehat{\Xb}\c\ov{\Hc}) \Bb+ ( \Ga_b\c \Ga_b) \c B.
         \eeaa
Using the definition $\Bb_3=\nabc_3\Bb+2 \ov{\tr \Xb} \Bb$, we have
         \beaa
       &&   \nabc_3 \DDc\hot \Bb\\
       &=&   \DDc\hot \big( \Bb_3-2 \ov{\tr \Xb} \Bb \big) - \frac 1 2 \tr \Xb \DDc \hot \Bb  + H \hot \big( \Bb_3-2 \ov{\tr \Xb} \Bb\big)\\
         &&+  \err_1 + r^{-2} \dk^{\le 1} (\Ga_g \c \Ga_b)\\
         &=&   \DDc\hot  \Bb_3 + H \hot  \Bb_3-\big( 2 \ov{\tr \Xb}+ \frac 1 2 \tr \Xb\big) \DDc\hot \Bb-\big( 2 \DDc\ov{\tr \Xb} +2 \ov{\tr \Xb} H \big)\hot  \Bb\\
         &&+  \err_1+ r^{-2} \dk^{\le 1} (\Ga_g \c \Ga_b).
         \eeaa
This implies
\beaa
I&=& \frac 1 2  \nabc_3 \DDc\hot \Bb+\frac 1 2 \big( 2 \ov{\tr \Xb}+\frac 1 2 \tr \Xb\big) \DDc\hot \Bb\\
&=& \frac 1 2 \DDc\hot  \Bb_3 +\frac 1 2  H \hot  \Bb_3-\big(  \DDc\ov{\tr \Xb} + \ov{\tr \Xb} H \big)\hot  \Bb+  \err_1+ r^{-2} \dk^{\le 1} (\Ga_g \c \Ga_b).
\eeaa
Making use of the equation 
\beaa
\Bb_3  &=&-\frac{1}{2} \big( \DDbc \c\Ab  + \ov{H} \c \Ab\big)-3P \,\Xib,
       \eeaa
        we deduce
        \beaa
I&=&  -  \frac 1 4   \DDc \hot \big(\DDbc \c\Ab +  \ov{H}  \c  \Ab\big)- \frac 1 4 H \hot \big(\DDbc \c\Ab +  \ov{H} \c \Ab \big)\\
    &&-  \frac 3 2 P  \DDc \hot  \,\Xib -  \frac 3 2   (\DDc P) \hot  \,\Xib - \frac 3 2 P H \hot  \,\Xib -\big(  \DDc\ov{\tr \Xb} + \ov{\tr \Xb} H \big)\hot  \Bb\\
&&+  \err_1+ r^{-2} \dk^{\le 1} (\Ga_g \c \Ga_b).
\eeaa
Using again \eqref{eq:error-DDP} to write $\DD P= -3P \Hb + r^{-2} \Ga_g $, we obtain
        \beaa
I&=&  -  \frac 1 4   \DDc \hot \big(\DDbc \c\Ab +  \ov{H}  \c  \Ab\big)- \frac 1 4 H \hot \big(\DDbc \c\Ab +  \ov{H} \c \Ab \big)\\
    &&-  \frac 3 2 P  \DDc \hot  \,\Xib +\big(   \frac 9 2   P \Hb - \frac 3 2 P H\big) \hot  \,\Xib -\big(  \DDc\ov{\tr \Xb} + \ov{\tr \Xb} H \big)\hot  \Bb\\
&&+  \err_1+ r^{-2} \dk^{\le 1} (\Ga_g \c \Ga_b).
\eeaa
                
     {\bf Step 2: Calculation of $J$. }  
     We have
     \beaa
     J&=& 2\nabc_3(\Hb    \hot \Bb)+ 2\big( 2 \ov{\tr \Xb}+\frac 1 2 \tr \Xb\big)\Hb    \hot \Bb\\
    &=& 2\nabc_3 \Hb    \hot \Bb+ 2\Hb    \hot \nabc_3\Bb+ 4 \ov{\tr \Xb}  (\Hb    \hot \Bb) +\tr \Xb \Hb \hot \Bb\\
        &=& 2\nabc_3 \Hb    \hot \Bb+ 2\Hb    \hot \big(\Bb_3-2 \ov{\tr \Xb} \Bb \big)+ 4 \ov{\tr \Xb}  (\Hb    \hot \Bb)+\tr\Xb \Hb \hot \Bb\\
                &=& 2\nabc_3 \Hb    \hot \Bb+ 2\Hb    \hot \Bb_3+\tr\Xb \Hb \hot \Bb.
     \eeaa
Making use of the equations
       \beaa
     \nabc_3\Hb&=&  -\frac{1}{2}\ov{\tr\Xb}(\Hb-H) +\Bb +r^{-2} \Ga_b,
     \eeaa
     and
     \beaa
\Bb_3  &=&-\frac{1}{2} \big( \DDbc \c\Ab  + \ov{H} \c \Ab\big)-3P \,\Xib,
       \eeaa
we deduce
     \beaa
    J                &=&  -\ov{\tr\Xb}(\Hb-H)    \hot \Bb+ 2\Hb    \hot \Bb_3+\tr\Xb \Hb \hot \Bb+2 \Bb\hot \Bb+ r^{-3} (\Ga_b \c \Ga_b)        \\
       &=&  2\Hb    \hot \big(-\frac{1}{2} \big( \DDbc \c\Ab  + \ov{H} \c \Ab\big)-3P \,\Xib \big) -\ov{\tr\Xb}(\Hb-H)    \hot \Bb+\tr\Xb \Hb \hot \Bb\\
       &&+ 2\Bb\hot \Bb+ r^{-3} (\Ga_b \c \Ga_b)   \\
       &=& -\Hb    \hot \big(\DDbc \c\Ab  + \Ab\c \ov{H} \big) +\big(  -\ov{\tr\Xb}(\Hb-H)    +\tr\Xb \Hb \big)\hot \Bb-6P\Hb    \hot  \,\Xib\\
       &&+2 \Bb\hot \Bb+  r^{-3} (\Ga_b \c \Ga_b)   .
     \eeaa

{ \bf Step 3. Calculation of $K$.}
We have
\beaa
K&=& 3\nabc_3(P\Xbh)+ 3\big( 2 \ov{\tr \Xb}+\frac 1 2 \tr \Xb\big)P\Xbh\\
&=& 3 \nabc_3P \Xbh+3 P \nabc_3\Xbh+ 6  \ov{\tr \Xb} P \Xbh+ \frac 3 2 \tr \Xb P\Xbh.
\eeaa
Making use of the equations,
\beaa
\nabc_3P  &=& -\frac{3}{2}\ov{\tr\Xb} P-\frac{1}{2}\DDbc \c\Bb - \ov{H} \c\Bb +\Xib \c \ov{B}-\frac 1  4 \Xbh\c \ov{A} , \\
\nabc_3\Xbh&=& -\frac 1 2 (\tr\Xb+\ov{\tr\Xb}) \Xbh+\frac 1 2  \DDc\hot \Xib+  \frac 1 2  \Xib\hot(H+\Hb)-\Ab,
\eeaa
we obtain
\beaa
K&=& 3\big( -\frac{3}{2}\ov{\tr\Xb} P-\frac{1}{2}\DDbc \c\Bb - \ov{H} \c\Bb +(A, B)  \c \Ga_b \big)\Xbh\\
&&+3 P\big( -\frac 1 2 (\tr\Xb+\ov{\tr\Xb}) \Xbh+\frac 1 2  \DDc\hot \Xib+  \frac 1 2  \Xib\hot(H+\Hb)-\Ab \big)\\
&&+ 6  \ov{\tr \Xb} P \Xbh+ \frac 3 2 \tr \Xb P\Xbh\\
&=&\frac 3 2 P \DDc\hot \Xib+  \frac 3 2P   \Xib\hot(H+\Hb)-3P\Ab\\
&&+ \big(-\frac 3 2 \DDbc \c\Bb-3 \ov{H} \c\Bb   \big) \Xbh+(\Ga_b \c \Ga_b)\c (A, B)
\eeaa

{ \bf Step 4. Final sum. } By summing the above three terms we obtain the cancellation of the terms in $\Xib$, and we deduce
\beaa
- \nabc_3\Ab_4-\big( 2 \ov{\tr \Xb}+\frac 1 2 \tr \Xb\big) \Ab_4&=& I+J+K\\
&=&  -  \frac 1 4   \DDc \hot \big(\DDbc \c\Ab +  \ov{H}  \c  \Ab\big)\\
&&- \frac 1 4 (H+4\Hb) \hot \big(\DDbc \c\Ab +  \ov{H} \c \Ab \big)-3P\Ab\\
    && -\big(  \DDc\ov{\tr \Xb} +(\ov{\tr\Xb}   -\tr\Xb) \Hb \big)\hot \Bb\\
&&+\err_1+ r^{-2} \dk^{\le 1} (\Ga_g \c \Ga_b)+2 \Bb\hot \Bb+  r^{-3} (\Ga_b \c \Ga_b)\\
&&+ \big(-\frac 3 2 \DDbc \c\Bb-3 \ov{H} \c\Bb   \big) \Xbh+(\Ga_b \c \Ga_b)\c (A, B).
\eeaa
Using \eqref{eq:DDc-ov-trXb}, i.e.
     \beaa
    \DDc \ov{\tr \Xb }&=&    2 i\Im(\tr  \Xb) \Hb  + r^{-1}\Ga_g , 
     \eeaa
       we finally obtain
\beaa
- \nabc_3\Ab_4-\big( 2 \ov{\tr \Xb}+\frac 1 2 \tr \Xb\big) \Ab_4&=&  -  \frac 1 4   \DDc \hot \big(\DDbc \c\Ab +  \ov{H}  \c  \Ab\big)\\
&&- \frac 1 4 (H+4\Hb) \hot \big(\DDbc \c\Ab +  \ov{H} \c \Ab \big)-3P\Ab\\
&&+\err_1+ r^{-2} \dk^{\le 1} (\Ga_g \c \Ga_b)+2 \Bb\hot \Bb+  r^{-3} (\Ga_b \c \Ga_b)\\
&&+ \left(-\frac 3 2 \DDbc \c\Bb-3 \ov{H} \c\Bb   \right) \Xbh+(\Ga_b \c \Ga_b)\c (A, B).
\eeaa

Therefore
\beaa
\bsplit
- \nabc_3\Ab_4-\big( 2 \ov{\tr \Xb}+\frac 1 2 \tr \Xb\big) \Ab_4&= -  \frac 1 4 \big(  \DDc+ H+ 4\Hb\big)  \hot \big(\DDbc \c\Ab +  \ov{H}  \c  \Ab\big)\\
&-3P\Ab+ \err_{TE}[\Ab],
\end{split}
\eeaa
with  error term $ \err_{TE}[\Ab]$ given by
\beaa
\err_{TE}&=& \err_1+ r^{-2} \dk^{\le 1} (\Ga_g \c \Ga_b)+2 \Bb\hot \Bb+  r^{-3} (\Ga_b \c \Ga_b)\\
&&+ \left(-\frac 3 2 \DDbc \c\Bb-3 \ov{H} \c\Bb   \right) \Xbh+(\Ga_b \c \Ga_b)\c (A, B)\\
&=& -2\Bb \hot \Bb- 2\tr X \Xib \hot  \Bb-\frac 1 2\Xbh \c \ov{\DDc} \Bb+ (\widehat{\Xb}\c\ov{\Hc}) \Bb+ ( \Ga_b\c \Ga_b) \c B\\
&&+ r^{-2} \dk^{\le 1} (\Ga_g \c \Ga_b)+2 \Bb\hot \Bb+  r^{-3} (\Ga_b \c \Ga_b)\\
&&+ \left(-\frac 3 2 \DDbc \c\Bb-3 \ov{H} \c\Bb   \right) \Xbh+(\Ga_b \c \Ga_b)\c (A, B),
\eeaa
which gives
\beaa
\err_{TE}&=& - 2\tr X \Xib \hot  \Bb-\frac 1 2\Xbh \c \ov{\DDc} \Bb-\frac 3 2( \DDbc \c\Bb) \Xbh\\
 &&+ (\widehat{\Xb}\c\ov{\Hc}) \Bb+(\Ga_b \c \Ga_b)\c (A, B)+ r^{-2} \dk^{\le 1} (\Ga_g \c \Ga_b),
\eeaa
as stated.


\section{Proof of Theorem \ref{THEOREM:EQ-QFB}}
\label{proof:thm-eq-qfb}


We take  a  modified $\nabc_4$ derivative  of the Teukolsky equation.

\begin{proposition}\label{prop:appendix-nabc43Ab4} We have
\bea\label{eq:Ab-434}
\begin{split}
&\Big( \nabc_4+  \tr X+ \frac 1 2\ov{\tr X}\Big) \Big( \nabc_3+ 2 \ov{\tr \Xb}+\frac 1 2 \tr \Xb\Big) \Ab_4\\
&=    \frac 1 4\big(  \DDc+H+5\Hb \big) \hot \Big( \DDbc \c \Ab_4+(\ov{H}+\ov{\Hb}) \c \Ab_4\Big)+3P\Ab_4\\
& +3\Big( \frac 1 2\ov{\tr X}-\tr X\Big) P\Ab+\mathcal{J}_{434} \hot \big(\DDbc \c\Ab +  \ov{H} \c \Ab \big)+  \err_{434},
\end{split}
\eea
  where $ \mathcal{J}_{434}$ is a one-form given by
  \bea
 \mathcal{J}_{434}&=& \frac 1 4 \Big(-\frac 1 2  \DDc( \tr X+\ov{\tr X})+ \frac 1 2 (\tr X -\ov{\tr X}) H   -4 \tr X\Hb  \Big)-\frac 1 4 B \label{eq:mathcal-J-434-r-2Gag}\\
 &&+r^{-2} \Ga_g+ \Xh \c \Hc \nn\\
 &=&-\frac 3 4 \tr X \Hb+ r^{-1}\dk^{\le 1 } \Ga_g, \label{eq:mathcal-J-434-r-1Gag}
  \eea
 and the error terms are schematically given by
    \bea
    \begin{split}
  \err_{434}&= \nabc_4\err_{TE} +\big( \frac 1 2 \tr X+ \ov{\tr X}\big)\err_{TE}\\
  &+\DDc\hot( \Xh \c \ov{\DDc} \Ab)+ \DDc \hot ((  \DDc\c\ov{\Xh})\c \Ab)\\
 &+ r^{-1} \dk^{\le 1 } ((A,B)\c \Ga_b) + r^{-3}\dk^{\leq1}( \Ga_g  \c \Ga_b)  .
 \end{split}
  \eea 
\end{proposition}

\begin{proof} We apply the operator $\nabc_4+\big(  \tr X+ \frac 1 2\ov{\tr X}\big)$ to the Teukolsky equation \eqref{eq:Teuk-Ab-Ab4}, and we deduce
\beaa
&&\Big( \nabc_4+  \tr X+ \frac 1 2\ov{\tr X}\Big) \Big( \nabc_3\Ab_4+\big( 2 \ov{\tr \Xb}+\frac 1 2 \tr \Xb\big) \Ab_4\Big)\\
 &=&   I+J+K \\
&&+\nabc_4\err_{TE} +\big( \tr X+ \frac 1 2 \ov{\tr X}\big)\err_{TE} + r^{-3} \dk^{\le 1} (\Ga_g \c \Ga_b)
\eeaa
where
\beaa
I&=&  \frac 1 4 \nabc_4  \DDc \hot \big(\DDbc \c\Ab +  \ov{H}  \c  \Ab\big)+  \frac 1 4 \big(  \tr X+ \frac 1 2\ov{\tr X}\big)  \DDc \hot \big(\DDbc \c\Ab +  \ov{H}  \c  \Ab\big)\\
J&=& \frac 1 4 \nabc_4\Big( (H+4\Hb) \hot \big(\DDbc \c\Ab +  \ov{H} \c \Ab \big)\Big)\\
&&+\frac 1 4\big( \tr X+\frac 1 2 \ov{\tr X}\big) (H+4\Hb) \hot \big(\DDbc \c\Ab +  \ov{H} \c \Ab \big)\\
K&=& 3\nabc_4(P\Ab)+3\big(  \tr X+ \frac 1 2\ov{\tr X}\big) P\Ab.
\eeaa

  { \bf Step 1. Calculation of $I$.}
  Here we use  the commutation formula \eqref{eq:comm-nabc4nabc3DDchot-precise} applied to $F=\big(\DDbc \c\Ab +  \ov{H}  \c  \Ab\big)$ of signature $s=-2$, and we obtain\footnote{Because of the gauge conditions $\Xi=0$, $\Hbc=0$ there is no cubic term involving $\Ab$. }
  \beaa
 &&\nabc_4 \DDc \hot \big(\DDbc \c\Ab +  \ov{H}  \c  \Ab\big) \\
 &=& \DDc \hot \nabc_4\big(\DDbc \c\Ab +  \ov{H}  \c  \Ab\big)+\Hb\hot\nabc_4 \big(\DDbc \c\Ab +  \ov{H}  \c  \Ab\big)\\
  && -\frac 1 2 \tr X\left( \DDc \hot  \big(\DDbc \c\Ab +  \ov{H}  \c  \Ab\big)+3\Hb\hot \big(\DDbc \c\Ab +  \ov{H}  \c  \Ab\big)\right)\\
 &&+B \hot \big(\DDbc \c\Ab \big) -\frac 1 2\Xh \c \ov{\DDc} \big(\DDbc \c\Ab \big)+ r^{-3}\dk^{\leq1}( \Ga_g  \c \Ga_b).
 \eeaa
 We now compute $ \nabc_4 \big(\DDbc \c\Ab +  \ov{H}  \c  \Ab\big)$. We first
apply \eqref{eq:comm-nabc4nabc3-ovDDc-U-precise} to $U=\Ab$ of signature $s=-2$ and we obtain
\beaa
\,[\nabc_4, \ov{\DDc}\c] \Ab&=& -\frac 1 2\ov{\tr X} \big(  \ov{\DDc}\c \Ab\big) +\ov{\Hb}\c\nabc_4 \Ab\\
&& +4\ov{B} \c \Ab   -\frac 1 2 \Xh \c \ov{\DDc} \Ab+ r^{-2} \dk^{\leq 1}(\Ga_g \c  \Ga_b).
\eeaa
The above gives
\beaa
 \nabc_4 \big(\DDbc \c\Ab +  \ov{H}  \c  \Ab\big)&=& \DDbc \c( \nabc_4\Ab) -\frac 1 2\ov{\tr X} \big(  \ov{\DDc}\c \Ab\big) +(\ov{H}+\ov{\Hb}) \c \nabc_4 \Ab\\
&&+ \nabc_4 \ov{H}  \c  \Ab +4\ov{B} \c \Ab   -\frac 1 2 \Xh \c \ov{\DDc} \Ab+ r^{-2} \dk^{\leq 1}(\Ga_g \c  \Ga_b)\\
&=& \DDbc \c( \nabc_4\Ab) -\frac 1 2\ov{\tr X} \big(  \ov{\DDc}\c \Ab\big) +(\ov{H}+\ov{\Hb}) \c \nabc_4 \Ab\\
&&-\frac{1}{2}\tr X(\ov{H}-\ov{\Hb})   \c  \Ab +3\ov{B} \c \Ab   -\frac 1 2 \Xh \c \ov{\DDc} \Ab+ r^{-2} \dk^{\leq 1}(\Ga_g \c  \Ga_b),
\eeaa
where we used the null structure equation 
  \beaa
  \nabc_4H &=&  -\frac{1}{2}\ov{\tr X}(H-\Hb)  -B+r^{-2} \Ga_g.
  \eeaa
Now writing  $\nabc_4 \Ab=\Ab_4-\frac 1 2 \tr X \Ab $, we obtain
\beaa
 \nabc_4 \big(\DDbc \c\Ab +  \ov{H}  \c  \Ab\big)&=& \DDbc \c( \Ab_4-\frac 1 2 \tr X \Ab) -\frac 1 2\ov{\tr X} \big(  \ov{\DDc}\c \Ab\big) \\
 &&+(\ov{H}+\ov{\Hb}) \c ( \Ab_4-\frac 1 2 \tr X \Ab)-\frac{1}{2}\tr X(\ov{H}-\ov{\Hb})   \c  \Ab \\
 &&+3\ov{B} \c \Ab   -\frac 1 2 \Xh \c \ov{\DDc} \Ab+ r^{-2} \dk^{\leq 1}(\Ga_g \c  \Ga_b)\\
    &=& \DDbc \c \Ab_4+(\ov{H}+\ov{\Hb}) \c \Ab_4-\frac 1 2( \tr X+\ov{\tr X})  \DDbc \c \Ab\\
         &&- \big( \frac 1 2  \DDbc \tr X +\tr X\ov{H}\big)  \c  \Ab \\
         &&+3\ov{B} \c \Ab   -\frac 1 2 \Xh \c \ov{\DDc} \Ab+ r^{-2} \dk^{\leq 1}(\Ga_g \c  \Ga_b).
\eeaa
 Using the Codazzi equation
 \beaa
     \ov{ \DDc}\tr X&=& - ( \tr X-\ov{\tr X})\ov{H}+    2 \ov{B}    + \DDc\c\ov{\Xh} + r^{-2} \Ga_g  ,
      \eeaa
we obtain
    \bea\label{eq:nabc4-DDbcAb}
    \begin{split}
& \nabc_4 \big(\DDbc \c\Ab +  \ov{H}  \c  \Ab\big)     \\
   &= \DDbc \c \Ab_4+(\ov{H}+\ov{\Hb}) \c \Ab_4-\frac 1 2( \tr X+\ov{\tr X}) \big( \DDbc \c \Ab+\ov{H}  \c  \Ab \big) \\
            &+2\ov{B} \c \Ab   -\frac 1 2 \Xh \c \ov{\DDc} \Ab-\frac 1 2 (  \DDc\c\ov{\Xh})\c \Ab + r^{-2} \dk^{\leq 1}(\Ga_g \c  \Ga_b).
         \end{split}
 \eea
 By applying the operator $\DDc\hot$ to \eqref{eq:nabc4-DDbcAb}, we obtain
 \beaa
&& \DDc\hot \nabc_4 \big(\DDbc \c\Ab +  \ov{H}  \c  \Ab\big)     \\
  &=&\DDc\hot \Big( \DDbc \c \Ab_4+(\ov{H}+\ov{\Hb}) \c \Ab_4\Big)-\frac 1 2( \tr X+\ov{\tr X}) \DDc\hot \big( \DDbc \c \Ab+\ov{H}  \c  \Ab \big)\\
  &&-\frac 1 2\DDc( \tr X+\ov{\tr X}) \hot \big( \DDbc \c \Ab+\ov{H}  \c  \Ab \big) \\
  &&+2\DDc\hot( \ov{B} \c \Ab)   -\frac 1 2\DDc\hot( \Xh \c \ov{\DDc} \Ab)-\frac 1 2 \DDc \hot ((  \DDc\c\ov{\Xh})\c \Ab)  + r^{-3} \dk^{\leq 1}(\Ga_g \c  \Ga_b).
 \eeaa
Also, from \eqref{eq:nabc4-DDbcAb} we obtain
\beaa
 \underline{H} \hot \nabc_4 \big(\DDbc \c\Ab +  \ov{H}  \c  \Ab\big)&=&  \underline{H} \hot \Big( \DDbc \c \Ab_4+(\ov{H}+\ov{\Hb}) \c \Ab_4 \Big)\\
 && -\frac 1 2( \tr X+\ov{\tr X})  \Hb \hot \big( \DDbc \c \Ab+\ov{H}  \c  \Ab \big)+ r^{-3} \dk^{\leq 2}(\Ga_g \c  \Ga_b).
\eeaa
We therefore deduce
\beaa
4I&=&\DDc\hot \Big( \DDbc \c \Ab_4+(\ov{H}+\ov{\Hb}) \c \Ab_4\Big)-\frac 1 2( \tr X+\ov{\tr X}) \DDc\hot \big( \DDbc \c \Ab+\ov{H}  \c  \Ab \big)\\
  &&-\frac 1 2\DDc( \tr X+\ov{\tr X}) \hot \big( \DDbc \c \Ab+\ov{H}  \c  \Ab \big) \\
         &&+   \underline{H} \hot \Big( \DDbc \c \Ab_4+(\ov{H}+\ov{\Hb}) \c \Ab_4 \Big) -\frac 1 2( \tr X+\ov{\tr X})  \Hb \hot \big( \DDbc \c \Ab+\ov{H}  \c  \Ab \big) \\
 &&- \frac 1 2 \tr X\left( \DDc\hot \big(\DDbc \c\Ab +  \ov{H}  \c  \Ab\big) + 3\Hb\hot \big(\DDbc \c\Ab +  \ov{H}  \c  \Ab\big)\right)\\
   &&+  \big(  \tr X+ \frac 1 2\ov{\tr X}\big)  \DDc \hot \big(\DDbc \c\Ab +  \ov{H}  \c  \Ab\big)\\
   &&+\frac 1 4 B \hot \big(\DDbc \c\Ab \big)+\frac 1 2\DDc\hot( \ov{B} \c \Ab)  -\frac 1 8\Xh \c \ov{\DDc} \big(\DDbc \c\Ab \big)\\
   &&  -\frac 1 8\DDc\hot( \Xh \c \ov{\DDc} \Ab)-\frac 1 8 \DDc \hot ((  \DDc\c\ov{\Xh})\c \Ab)+ r^{-3}\dk^{\leq1}( \Ga_g  \c \Ga_b),
\eeaa
which gives
\beaa
I&=&\frac 1 4 \big( \DDc+\Hb\big)\hot \Big( \DDbc \c \Ab_4+(\ov{H}+\ov{\Hb}) \c \Ab_4\Big)\\
  &&-\frac 1 8\big( \DDc( \tr X+\ov{\tr X})+( 4\tr X+\ov{\tr X})  \Hb \big)  \hot \big( \DDbc \c \Ab+\ov{H}  \c  \Ab \big) \\
  &&+\frac 1 4 B \hot \big(\DDbc \c\Ab \big)+\frac 1 2\DDc\hot( \ov{B} \c \Ab)  -\frac 1 8\Xh \c \ov{\DDc} \big(\DDbc \c\Ab \big)\\
   &&  -\frac 1 8\DDc\hot( \Xh \c \ov{\DDc} \Ab)-\frac 1 8 \DDc \hot ((  \DDc\c\ov{\Xh})\c \Ab)+ r^{-3}\dk^{\leq1}( \Ga_g  \c \Ga_b).
\eeaa

    { \bf Step 2. Calculation of $J$.} We have, using \eqref{eq:nabc4-DDbcAb}, 
    \beaa
    J&=& \frac 1 4 (H+4\Hb) \hot  \nabc_4 \big(\DDbc \c\Ab +  \ov{H} \c \Ab \big)+\frac 1 4  \nabc_4(H+4\Hb) \hot \big(\DDbc \c\Ab +  \ov{H} \c \Ab \big)\\
    &&+\frac 1 4\big(  \tr X+ \frac 1 2\ov{\tr X}\big) (H+4\Hb) \hot \big(\DDbc \c\Ab +  \ov{H} \c \Ab \big)\\
    &=& \frac 1 4 (H+4\Hb) \hot  \Big[\DDbc \c \Ab_4+(\ov{H}+\ov{\Hb}) \c \Ab_4-\frac 1 2( \tr X+\ov{\tr X}) \big( \DDbc \c \Ab+\ov{H}  \c  \Ab \big) \Big]\\
         &&+\frac 1 4  \nabc_4(H+4\Hb) \hot \big(\DDbc \c\Ab +  \ov{H} \c \Ab \big)\\
    &&+\frac 1 4\big( \tr X+ \frac 1 2\ov{\tr X}\big) (H+4\Hb) \hot \big(\DDbc \c\Ab +  \ov{H} \c \Ab \big) + r^{-3} \dk^{\leq 1}(\Ga_g \c  \Ga_b),
    \eeaa
    which gives
        \beaa
    J    &=& \frac 1 4 (H+4\Hb) \hot  \Big(\DDbc \c \Ab_4+(\ov{H}+\ov{\Hb}) \c \Ab_4 \Big)\\
         &&+\Big( \frac 1 4  \nabc_4(H+4\Hb)+\frac 1 8\tr X (H+4\Hb)  \Big) \hot \big(\DDbc \c\Ab +  \ov{H} \c \Ab \big)        + r^{-3} \dk^{\leq 1}(\Ga_g \c  \Ga_b).
    \eeaa

    { \bf Step 3. Calculation of $K$.} Writing  $\nabc_4 \Ab=\Ab_4-\frac 1 2 \tr X \Ab $, and   making use of the equation
  \beaa
  \nabc_4P -\frac{1}{2}\DDc\c \ov{B} &=& -\frac{3}{2}\tr X P + \Hb \c\ov{B} -\frac{1}{4}\Xbh\c \ov{A},
  \eeaa
  we have
    \beaa
     K&=& 3P\nabc_4\Ab+3\nabc_4(P)\Ab+3\big(  \tr X+\frac 1 2 \ov{\tr X}\big) P\Ab\\
     &=& 3P\big(\Ab_4-\frac 1 2 \tr X \Ab  \big)+3\big(-\frac{3}{2}\tr X P  \big)\Ab+3\big( \tr X+\frac 1 2 \ov{\tr X}\big) P\Ab\\
     &&+\frac 3 2 (\DDc\c \ov{B}) \Ab + r^{-3} (\Ga_g \c \Ga_b)+ \Xbh \c \Ab \c A \\
          &=& 3P\Ab_4 +3\big( \frac 1 2\ov{\tr X}-\tr X\big) P\Ab+\frac 3 2 (\DDc\c \ov{B}) \Ab + r^{-3} (\Ga_g \c \Ga_b)+ r^{-1} (\Ga_b \c A),
     \eeaa
    where we wrote $\Xbh \c \Ab \c A=r^{-1} (\Ga_b \c A)$.

  { \bf Step 4. Final sum.} We deduce
  \beaa
&&\Big( \nabc_4+  \tr X+ \frac 1 2\ov{\tr X}\Big) \Big( \nabc_3\Ab_4+\big( 2 \ov{\tr \Xb}+\frac 1 2 \tr \Xb\big) \Ab_4\Big)\\
&=&    \frac 1 4\big(  \DDc+H+5\Hb \big) \hot \Big( \DDbc \c \Ab_4+(\ov{H}+\ov{\Hb}) \c \Ab_4\Big)+3P\Ab_4\\
&& +3\big( \frac 1 2\ov{\tr X}-\tr X\big) P\Ab+\mathcal{J}_{434} \hot \big(\DDbc \c\Ab +  \ov{H} \c \Ab \big)+  \err_{434},
  \eeaa
  where $\mathcal{J}_{434}$ is the one form given by
  \beaa
  \mathcal{J}_{434}&:=& \frac 1 4  \nabc_4(H+4\Hb)+\frac 1 8 \tr X (H+4\Hb) -\frac 1 8\big( \DDc( \tr X+\ov{\tr X})+(4 \tr X+\ov{\tr X})  \Hb  \big).
  \eeaa
  and $\err_{434}$ are the error terms explicitly given by
  \beaa
  \err_{434}&=& \nabc_4\err_{TE}[\Ab] +\big( \frac 1 2 \tr X+ \ov{\tr X}\big)\err_{TE}[\Ab]  \\
 &&+\frac 1 4 B \hot \big(\DDbc \c\Ab \big)+\frac 1 2\DDc\hot( \ov{B} \c \Ab)+\frac 3 2 (\DDc\c \ov{B}) \Ab \\
 && -\frac 1 8\Xh \c \ov{\DDc} \big(\DDbc \c\Ab \big)  -\frac 1 8\DDc\hot( \Xh \c \ov{\DDc} \Ab)-\frac 1 8 \DDc \hot ((  \DDc\c\ov{\Xh})\c \Ab)\\
 &&+ r^{-3}\dk^{\leq1}( \Ga_g  \c \Ga_b)+ r^{-1} (\Ga_b \c A) .
  \eeaa
  Observe that the above error terms can schematically be written as
    \beaa
  \err_{434}&=& \nabc_4\err_{TE} +\big( \frac 1 2 \tr X+ \ov{\tr X}\big)\err_{TE}+\DDc\hot( \Xh \c \ov{\DDc} \Ab)+ \DDc \hot ((  \DDc\c\ov{\Xh})\c \Ab)\\
 &&+ r^{-1} \dk^{\le 1 } ((A,B)\c \Ga_b) + r^{-3}\dk^{\leq1}( \Ga_g  \c \Ga_b).
  \eeaa

We now compute $\mathcal{J}_{434}$. Recall, see Lemma  \ref{Lemma:DDc{trX}},
\beaa
\DDc \ov{\tr X}&=&     2 i\Im(\tr  X) H+ r^{-1}\Ga_g\\
 \DDc \tr X&=& -2\tr X \Hb + r^{-1} \Ga_g.
 \eeaa
Hence
\beaa
 \DDc( \tr X+\ov{\tr X})&=& (\tr  X-\ov{\tr X}) H-2\tr X \Hb + r^{-1} \Ga_g.
\eeaa
 Recalling also  the equations, see Lemma \ref{Lemma:nabc_4Hb},
 \beaa
 \nab_4\Hb+\tr X\Hb&=&r^{-2} \Ga_g\\
 \nabc_4H  &=&  -\frac{1}{2}\ov{\tr X}(H-\Hb)  -B+ \Xh \c \Hc
 \eeaa
 We deduce
\beaa
 \mathcal{J}_{434}&:=& \frac 1 4  \nabc_4(H+4\Hb)+\frac 1 8 \tr X (H+4\Hb) -\frac 1 8\big( \DDc( \tr X+\ov{\tr X})+(4 \tr X+\ov{\tr X})  \Hb  \big)\\
 &=&\frac 1 4 \Big(\big(\nabc_4 H+\frac 1 2 \tr X H\big)+ 4 
 \big(\nabc_4 \Hb +\frac 1 2 \tr X \Hb\big)\Big)-\frac 1  8((\tr  X-\ov{\tr X}) H-2\tr X \Hb )\\
 &&-\frac 1 8 (4 \tr X+\ov{\tr X})  \Hb+ r^{-1}\dk^{\le 1 } \Ga_g\\
 &=&\frac 1 4 \Big(\big( -\frac{1}{2}\ov{\tr X}(H-\Hb)  -B+ \Ga_g\c \Ga_b+\frac 1 2 \tr X H\big)+ 4 
 \big(-\tr X\Hb +\frac 1 2 \tr X \Hb\big)\Big)\\
 &&-\frac 1  8((\tr  X-\ov{\tr X}) H-2\tr X \Hb )-\frac 1 8 (4 \tr X+\ov{\tr X})  \Hb+ r^{-1}\dk^{\le 1 } \Ga_g\\
  &=&-\frac 3 4 \tr X \Hb+ r^{-1}\dk^{\le 1 } \Ga_g.
  \eeaa
Hence 
\beaa
 \mathcal{J}_{434}&=&-\frac 3 4 \tr X \Hb+ r^{-1}\dk^{\le 1 } \Ga_g.
\eeaa
We also have, by keeping only error terms that behave like $r^{-2} \Ga_g$ or better,
\beaa
 \mathcal{J}_{434}&:=& \frac 1 4  \nabc_4(H+4\Hb)+\frac 1 8 \tr X (H+4\Hb) -\frac 1 8\big( \DDc( \tr X+\ov{\tr X})+(4 \tr X+\ov{\tr X})  \Hb  \big)\\
 &=&\frac 1 4 \Big(  -\frac{1}{2}\ov{\tr X}(H-\Hb)  -B+  \Xh \c \Hc+\frac 1 2 \tr X H+ 4 
 \big(-\frac 1 2 \tr X\Hb+r^{-2} \Ga_g \big)\Big)\\
 && -\frac 1 8\big( \DDc( \tr X+\ov{\tr X})+(4 \tr X+\ov{\tr X})  \Hb  \big)\\
  &=&\frac 1 4 \Big(-\frac 1 2  \DDc( \tr X+\ov{\tr X})+ \frac 1 2 (\tr X -\ov{\tr X}) H   -4 \tr X\Hb  \Big)\\
  &&-\frac 1 4 B+r^{-2} \Ga_g+ \Xh \c \Hc,
  \eeaa
  as stated. This ends the proof of Proposition \ref{prop:appendix-nabc43Ab4}. 
\end{proof}

We now take another derivative of equation \eqref{eq:Ab-434}.

\begin{proposition}\label{prop:appendix-nabc4nabc4nabc3Ab} We have
  \bea\label{eq:nabc4-nabc4-nabc4-Ab4}
  \begin{split}
&\Big(\nabc_4 + 3  \tr X-\frac 1 2 \ov{\tr X}- 2\frac {\atrch^2}{ \trch}\Big)\Big( \nabc_4+  \tr X+ \frac 1 2\ov{\tr X}\Big) \Big( \nabc_3+ 2 \ov{\tr \Xb}+\frac 1 2 \tr \Xb\Big) \Ab_4\\
         &+ \big(2\tr X- \ov{\tr X}  \big)3P\Ab_4 \\
&=   \frac 1 4  \big( \DDc+H+6\Hb\big)\hot \Big( \DDbc \c \widetilde{\underline{Q}(\Ab)}+(\ov{H}+2\ov{\Hb}) \c \widetilde{\underline{Q}(\Ab)} \Big)+3P \ \widetilde{\underline{Q}(\Ab)}\\
   &+\mathscr{L}[\Ab, \Ab_4, \DD \Ab_4, \DD \Ab]+\err_{4434}
   \end{split}
\eea
where $\mathscr{L}[\Ab, \Ab_4, \DD \Ab_4, \DD \Ab]$ denote $O(|a|)$ linear order terms in $\Ab, \Ab_4, \DD \Ab, \DD \Ab_4$ explicitly given by \eqref{eq:mathscr-L-Ab}, schematically given by
\beaa
\mathscr{L}[\Ab, \Ab_4,  \DD \Ab, \DD \Ab_4]&=&O(ar^{-3}) \big( \DDbc \c \Ab_4+(\ov{H} +\ov{\Hb}) \c \Ab_4\big)\\
&&+O(a^2r^{-4})\big(\DDbc \c\Ab +  \ov{H} \c \Ab \big) +O(a^2r^{-7}) \Ab,
\eeaa
 and the error terms are given by 
\beaa
\err_{4434}&=& \nabc_4  \err_{434}+ \big(   \tr X+\frac 3 2 \ov{\tr X}+2\frac {\atrch^2}{ \trch} \big)  \err_{434}\\
&&+ r^{-2} \dk^{\leq 1}( (A,B) \c \Ga_b)+ r^{-4} \dk^{\leq 2}(\Ga_g \c  \Ga_b)\\
&&+ r^{-1}(\nabc_4 (\Xh \c \Hc)+\frac 3 2 \tr X (\Xh \c \Hc))\c  \dk^{\leq 1} \Ab. 
\eeaa
\end{proposition}
\begin{proof} We apply the operator $\nabc_4 + 2  \tr X+\frac 1 2 \ov{\tr X}$ to equation \eqref{eq:Ab-434}, and we obtain
\beaa
&&\Big(\nabc_4 + 2  \tr X+\frac 1 2 \ov{\tr X}\Big)\Big( \nabc_4+  \tr X+ \frac 1 2\ov{\tr X}\Big) \Big( \nabc_3+ 2 \ov{\tr \Xb}+\frac 1 2 \tr \Xb\Big) \Ab_4\\
&=&   I+J+K+L+M\\
&&+\nabc_4  \err_{434}+ \big( 2  \tr X+\frac 1 2 \ov{\tr X}\big)  \err_{434}+ r^{-4} \dk^{\le 1} (\Ga_g \c \Ga_b),
\eeaa
where
\beaa
I&=&   \frac 1 4 \nabc_4 \DDc  \hot \Big( \DDbc \c \Ab_4+(\ov{H}+\ov{\Hb}) \c \Ab_4\Big)\\
&&+  \frac 1 4 \big( 2  \tr X+\frac 1 2 \ov{\tr X}\big) \DDc  \hot \Big( \DDbc \c \Ab_4+(\ov{H}+\ov{\Hb}) \c \Ab_4\Big),\\
J&=&  \frac 1 4\nabc_4\Big( \big( H+5\Hb \big) \hot \big( \DDbc \c \Ab_4+(\ov{H}+\ov{\Hb}) \c \Ab_4\big) \Big)\\
&&+  \frac 1 4\big( 2  \tr X+\frac 1 2 \ov{\tr X}\big) \big( H+5\Hb \big) \hot \Big( \DDbc \c \Ab_4+(\ov{H}+\ov{\Hb}) \c \Ab_4\Big),\\
K&=& 3\nabc_4(P\Ab_4)+3\big( 2  \tr X+\frac 1 2 \ov{\tr X}\big) P\Ab_4, \\
L&=& 3\nabc_4 \Big(\big( \frac 1 2 \ov{\tr X}-\tr X\big) P\Ab\Big)+3 \big( 2  \tr X+\frac 1 2 \ov{\tr X}\big)\big( \frac 1 2 \ov{\tr X}-\tr X\big) P\Ab,\\
M&=&\nabc_4\Big( \mathcal{J}_{434} \hot \big(\DDbc \c\Ab +  \ov{H} \c \Ab \big)\Big)+ \big( 2  \tr X+\frac 1 2 \ov{\tr X}\big) \mathcal{J}_{434} \hot \big(\DDbc \c\Ab +  \ov{H} \c \Ab \big).
\eeaa

  { \bf Step 1. Calculation of $I$.}
  Here we use  the commutation formula \eqref{eq:comm-nabc4nabc3-DDchot-err}  applied to $F=\big( \DDbc \c \Ab_4+(\ov{H}+\ov{\Hb}) \c \Ab_4\big)$ of signature $s=-1$, where observe that $F = r^{-3} \Ga_b$. We therefore obtain
  \beaa
&&  \nabc_4 \DDc  \hot \Big( \DDbc \c \Ab_4+(\ov{H}+\ov{\Hb}) \c \Ab_4\Big)\\
&=& \DDc  \hot \nabc_4\big( \DDbc \c \Ab_4+(\ov{H}+\ov{\Hb}) \c \Ab_4\big)+ \underline{H} \hot \nabc_4 \big( \DDbc \c \Ab_4+(\ov{H}+\ov{\Hb}) \c \Ab_4\big)\\
&&- \frac 1 2 \tr X\left( \DDc\hot \big( \DDbc \c \Ab_4+(\ov{H}+\ov{\Hb}) \c \Ab_4\big) + 2\Hb\hot \big( \DDbc \c \Ab_4+(\ov{H}+\ov{\Hb}) \c \Ab_4\big)\right)\\
   &&+ r^{-4}\dk^{\leq1}( \Ga_g  \c \Ga_b).
  \eeaa
 We now compute $ \nabc_4 \big(\DDbc \c \Ab_4+(\ov{H}+\ov{\Hb}) \c \Ab_4\big)$ by 
applying \eqref{eq:comm-nabc4nabc3-ovDDc-U-err} to $U=\Ab_4$ of signature $s=-1$ and we obtain
\beaa
 \nabc_4 \big(\DDbc \c \Ab_4+(\ov{H}+\ov{\Hb}) \c \Ab_4\big)&=& \DDbc \c(  \nabc_4 \Ab_4)- \frac 1 2\ov{\tr X}\, ( \ov{\DDc} \c \Ab_4 - \ov{\Hb} \c \Ab_4)\\
         &&+(\ov{H}+2\ov{\Hb}) \c \nabc_4\Ab_4+\nabc_4(\ov{H}+\ov{\Hb}) \c \Ab_4\\
         &&+ r^{-3} \dk^{\leq 1}(\Ga_g \c  \Ga_b).
\eeaa
  Making use of the equations
  \beaa
  \nabc_4H &=&  -\frac{1}{2}\ov{\tr X}(H-\Hb) +r^{-1} \Ga_g, \\
  \nab_4\Hb&=&-\tr X\Hb+r^{-2} \Ga_g,
  \eeaa
  we obtain
  \beaa
 \nabc_4 \big(\DDbc \c \Ab_4+(\ov{H}+\ov{\Hb}) \c \Ab_4\big)&=& \DDbc \c(  \nabc_4 \Ab_4)- \frac 1 2\ov{\tr X}\, ( \ov{\DDc} \c \Ab_4+\ov{\Hb} \c \Ab_4 )\\
         &&+(\ov{H}+2\ov{\Hb}) \c \nabc_4\Ab_4 -\frac{1}{2}\tr X(\ov{H}-\ov{\Hb})  \c \Ab_4\\
         &&+ r^{-3} \dk^{\leq 1}(\Ga_g \c  \Ga_b).
\eeaa
Using \eqref{eq:widetilde-Qb-Ab} to write $ \nabc_4\Ab_4=\widetilde{\underline{Q}(\Ab)}-\big(\frac{5}{2}\tr X-\ov{\tr X} - 2\frac {\atrch^2}{ \trch} \big)\Ab_4 $, we obtain
    \beaa
&& \nabc_4 \big(\DDbc \c \Ab_4+(\ov{H}+\ov{\Hb}) \c \Ab_4\big)\\
&=& \DDbc \c\Big( \widetilde{\underline{Q}(\Ab)}-\big(\frac{5}{2}\tr X-\ov{\tr X} - 2\frac {\atrch^2}{ \trch} \big)\Ab_4\Big)\\
 &&- \frac 1 2\ov{\tr X}\, ( \ov{\DDc} \c \Ab_4+\ov{\Hb} \c \Ab_4 )-\frac{1}{2}\tr X(\ov{H}-\ov{\Hb})  \c \Ab_4\\
         &&+(\ov{H}+2\ov{\Hb}) \c \Big( \widetilde{\underline{Q}(\Ab)}-\left(\frac{5}{2}\tr X-\ov{\tr X} - 2\frac {\atrch^2}{ \trch} \right)\Ab_4\Big)+ r^{-3} \dk^{\leq 1}(\Ga_g \c  \Ga_b),
         \eeaa
         which gives
         \bea\label{eq:nabc4-DDbc-c-Ab4}
         \begin{split}
   & \nabc_4 \big(\DDbc \c \Ab_4+(\ov{H}+\ov{\Hb}) \c \Ab_4\big)  \\
   &= \DDbc \c \widetilde{\underline{Q}(\Ab)}+(\ov{H}+2\ov{\Hb}) \c \widetilde{\underline{Q}(\Ab)}\\
         &-\left(\frac{5}{2}\tr X-\frac 1 2 \ov{\tr X} - 2\frac {\atrch^2}{ \trch} \right) \big( \DDbc \c \Ab_4+(\ov{H} +\ov{\Hb}) \c \Ab_4\big)\\
         &+\ov{\mathcal{I}} \c \Ab_4+ r^{-3} \dk^{\leq 1}(\Ga_g \c  \Ga_b),
         \end{split}
\eea
where $\ov{\mathcal{I}}$ is the one-form given by
\beaa
\ov{\mathcal{I}}&=&  -\DDbc \left(\frac{5}{2}\tr X-\ov{\tr X} - 2\frac {\atrch^2}{ \trch} \right)+ \frac 1 2\ov{\tr X}\,  \ov{H} -\frac{1}{2}\tr X(\ov{H}-\ov{\Hb}) \\
         && - \left(\frac{5}{2}\tr X-\ov{\tr X} - 2\frac {\atrch^2}{ \trch} \right)\ov{\Hb}.
\eeaa
Using Lemma \eqref{eq:DDc-ov-trX}, we obtain for $\mathcal{I}$:
\bea\label{eq:ov-mathcal-I}
\begin{split}
\ov{\mathcal{I}}&= 4i \Im(tr X) \ov{H} - \left(2\tr X+\ov{\tr X} - 2\frac {\atrch^2}{ \trch} \right)\ov{\Hb}+ 2\DDbc\left(\frac {\atrch^2}{ \trch} \right)+r^{-1} \Ga_g\\
&= O(ar^{-3})+r^{-1} \Ga_g.
\end{split}
\eea

By applying the operator $\DDc\hot$ to \eqref{eq:nabc4-DDbc-c-Ab4}, we obtain
\beaa
   &&\DDc\hot  \nabc_4 \big(\DDbc \c \Ab_4+(\ov{H}+\ov{\Hb}) \c \Ab_4\big)     \\
    &=&\DDc\hot \Big( \DDbc \c \widetilde{\underline{Q}(\Ab)}+(\ov{H}+2\ov{\Hb}) \c \widetilde{\underline{Q}(\Ab)} \Big)\\
         &&-\big(\frac{5}{2}\tr X-\frac 1 2 \ov{\tr X} - 2\frac {\atrch^2}{ \trch} \big) \DDc\hot \big( \DDbc \c \Ab_4+(\ov{H} +\ov{\Hb}) \c \Ab_4\big)\\
                  &&-\DDc \big(\frac{5}{2}\tr X-\frac 1 2 \ov{\tr X} - 2\frac {\atrch^2}{ \trch} \big)\hot  \big( \DDbc \c \Ab_4+(\ov{H} +\ov{\Hb}) \c \Ab_4\big)\\
         &&+2\ov{\mathcal{I}} \c \DDc \Ab_4+2(\DDc \c\ov{\mathcal{I}}) \Ab_4+ r^{-4} \dk^{\leq 2}(\Ga_g \c  \Ga_b).
\eeaa
Also, from \eqref{eq:nabc4-DDbc-c-Ab4} we obtain
\beaa
&&\underline{H} \hot \nabc_4 \big( \DDbc \c \Ab_4+(\ov{H}+\ov{\Hb}) \c \Ab_4\big)\\
&=&\underline{H} \hot \Big(\DDbc \c \widetilde{\underline{Q}(\Ab)}+(\ov{H}+2\ov{\Hb}) \c \widetilde{\underline{Q}(\Ab)}\Big) \\
         &&-\big(\frac{5}{2}\tr X-\frac 1 2 \ov{\tr X} - 2\frac {\atrch^2}{ \trch} \big) \Hb \hot \big( \DDbc \c \Ab_4+(\ov{H} +\ov{\Hb}) \c \Ab_4\big)+\Hb \hot \big(\ov{\mathcal{I}} \c \Ab_4\big)\\
         &&+ r^{-4} \dk^{\leq 2}(\Ga_g \c  \Ga_b).
\eeaa
We therefore finally deduce
\beaa
I&=&\frac 1 4  \big( \DDc+\Hb\big)\hot \Big( \DDbc \c \widetilde{\underline{Q}(\Ab)}+(\ov{H}+2\ov{\Hb}) \c \widetilde{\underline{Q}(\Ab)} \Big)\\
         &&-\frac 1 4 \big(\tr X- \ov{\tr X} - 2\frac {\atrch^2}{ \trch} \big) \DDc\hot \big( \DDbc \c \Ab_4+(\ov{H} +\ov{\Hb}) \c \Ab_4\big)\\
                  &&+\frac 1 4 \mathcal{J} \hot  \big( \DDbc \c \Ab_4+(\ov{H} +\ov{\Hb}) \c \Ab_4\big)+\frac1 2\ov{\mathcal{I}} \c \DDc \Ab_4+\frac 12(\DDc \c\ov{\mathcal{I}}) \Ab_4+\frac 1 4\Hb \hot \big(\ov{\mathcal{I}} \c \Ab_4\big)\\
   &&+ r^{-4} \dk^{\leq 2}(\Ga_g \c  \Ga_b),
\eeaa
where $\mathcal{J}$ is the one-form given by
\beaa
\mathcal{J}&=& -\DDc \left(\frac{5}{2}\tr X-\frac 1 2 \ov{\tr X} - 2\frac {\atrch^2}{ \trch} \right)-\left(\frac{5}{2}\tr X-\frac 1 2 \ov{\tr X} - 2\frac {\atrch^2}{ \trch} \right) \Hb- \tr X\Hb.
\eeaa
Using Lemma \eqref{eq:DDc-ov-trX}, we obtain for $\mathcal{J}$:
\bea\label{eq:mathcal-J}
\begin{split}
\mathcal{J}&= - \frac{5}{2}( -2\tr X \Hb)+\frac 1 2 ( 2 i\Im(\tr  X) (H-\Hc)) + 2\DDc\left(\frac {\atrch^2}{ \trch}\right) \\
&-\left(\frac{5}{2}\tr X-\frac 1 2 \ov{\tr X} - 2\frac {\atrch^2}{ \trch} \right) \Hb- \tr X\Hb+r^{-1} \Ga_g\\
&=  i\Im(\tr  X) (H-\Hc) + 2\DDc\left(\frac {\atrch^2}{ \trch}\right) +\left(\frac{3}{2}\tr X+\frac 1 2 \ov{\tr X} + 2\frac {\atrch^2}{ \trch} \right) \Hb+r^{-1} \Ga_g.
\end{split}
\eea

  { \bf Step 2. Calculation of $J$.} We have, using \eqref{eq:nabc4-DDbc-c-Ab4},
  \beaa
  J&=&  \frac 1 4 \big( H+5\Hb \big) \hot \nabc_4\big( \DDbc \c \Ab_4+(\ov{H}+\ov{\Hb}) \c \Ab_4\big)\\
  &&+ \frac 1 4\nabc_4 \big( H+5\Hb \big) \hot \big( \DDbc \c \Ab_4+(\ov{H}+\ov{\Hb}) \c \Ab_4\big)\\
&&+  \frac 1 4\big( 2  \tr X+\frac 1 2 \ov{\tr X}\big) \big( H+5\Hb \big) \hot \Big( \DDbc \c \Ab_4+(\ov{H}+\ov{\Hb}) \c \Ab_4\Big)\\
&=&  \frac 1 4 \big( H+5\Hb \big) \hot \Big[\DDbc \c \widetilde{\underline{Q}(\Ab)}+(\ov{H}+2\ov{\Hb}) \c \widetilde{\underline{Q}(\Ab)}\\
         &&-\big(\frac{5}{2}\tr X-\frac 1 2 \ov{\tr X} - 2\frac {\atrch^2}{ \trch} \big) \big( \DDbc \c \Ab_4+(\ov{H} +\ov{\Hb}) \c \Ab_4\big)+\ov{\mathcal{I}} \c \Ab_4 \Big]\\
  &&+ \frac 1 4\nabc_4 \big( H+5\Hb \big) \hot \big( \DDbc \c \Ab_4+(\ov{H}+\ov{\Hb}) \c \Ab_4\big)\\
&&+  \frac 1 4\big( 2  \tr X+\frac 1 2 \ov{\tr X}\big) \big( H+5\Hb \big) \hot \Big( \DDbc \c \Ab_4+(\ov{H}+\ov{\Hb}) \c \Ab_4\Big)+ r^{-4} \dk^{\leq 1}(\Ga_g \c  \Ga_b),
  \eeaa
  which gives
    \beaa
  J&=&  \frac 1 4 \big( H+5\Hb \big) \hot \Big(\DDbc \c \widetilde{\underline{Q}(\Ab)}+(\ov{H}+2\ov{\Hb}) \c \widetilde{\underline{Q}(\Ab)}\Big) \\
         &&\mathcal{K}\hot \big( \DDbc \c \Ab_4+(\ov{H} +\ov{\Hb}) \c \Ab_4\big)+ \frac 1 4 \big( H+5\Hb \big) \hot (\ov{\mathcal{I}} \c \Ab_4)+ r^{-4} \dk^{\leq 1}(\Ga_g \c  \Ga_b)
  \eeaa
  where $\mathcal{K}$ is given by
  \beaa
  \mathcal{K}&=& - \frac 1 4\big(\frac{5}{2}\tr X-\frac 1 2 \ov{\tr X} - 2\frac {\atrch^2}{ \trch} \big)  \big( H+5\Hb \big) + \frac 1 4\nabc_4 \big( H+5\Hb \big)\\
  &&+  \frac 1 4\big( 2  \tr X+\frac 1 2 \ov{\tr X}\big) \big( H+5\Hb \big),
  \eeaa
which gives, making use of equations $ \nabc_4H =  -\frac{1}{2}\ov{\tr X}(H-\Hb)+r^{-1} \Ga_g$ and $ \nab_4\Hb=-\tr X\Hb+r^{-2} \Ga_g$, 
  \bea\label{eq:mathcal-K}
  \begin{split}
  \mathcal{K}&= - \frac 1 4\big(\frac{5}{2}\tr X-\frac 1 2 \ov{\tr X} - 2\frac {\atrch^2}{ \trch} \big)  \big( H+5\Hb \big) + \frac 1 4( -\frac{1}{2}\ov{\tr X}(H-\Hb))+ \frac 5 4(-\tr X\Hb)\\
  &+  \frac 1 4\big( 2  \tr X+\frac 1 2 \ov{\tr X}\big) \big( H+5\Hb \big)+r^{-1} \Ga_g\\
  &= - \frac 1 4 \Big[  \frac{1}{2} \big(\tr X- \ov{\tr X} - 4\frac {\atrch^2}{ \trch} \big)   H + \big(\frac{15}{2}\tr X-\frac{11}{2} \ov{\tr X} - 10\frac {\atrch^2}{ \trch} \big)  \Hb  \Big]+r^{-1} \Ga_g.
  \end{split}
  \eea

    { \bf Step 3. Calculation of $K$.} Writing $ \nabc_4\Ab_4=\widetilde{\underline{Q}(\Ab)}-\big(\frac{5}{2}\tr X-\ov{\tr X} - 2\frac {\atrch^2}{ \trch} \big)\Ab_4 $, we have
\beaa
K&=& 3P\nabc_4\Ab_4+3\nabc_4(P)\Ab_4+3\big( 2  \tr X+\frac 1 2 \ov{\tr X}\big) P\Ab_4\\
&=& 3P\Big(\widetilde{\underline{Q}(\Ab)}-\big(\frac{5}{2}\tr X-\ov{\tr X} - 2\frac {\atrch^2}{ \trch} \big)\Ab_4 \Big)+3\big(  -\frac{3}{2}\tr X P\big) \Ab_4\\
&&+3\big( 2  \tr X+\frac 1 2 \ov{\tr X}\big) P\Ab_4+ r^{-4} \dk^{\leq1}(\Ga_g \c \Ga_b)\\
&=& 3P \ \widetilde{\underline{Q}(\Ab)}-3 \big(2\tr X-\frac 3 2 \ov{\tr X} - 2\frac {\atrch^2}{ \trch} \big)P \Ab_4+ r^{-4} (\Ga_g \c \Ga_b).
\eeaa

     { \bf Step 4. Calculation of $L$.} Writing  $\nabc_4 \Ab=\Ab_4-\frac 1 2 \tr X \Ab $, and   making use of the equations
  \beaa
  \nabc_4P -\frac{1}{2}\DDc\c \ov{B} &=& -\frac{3}{2}\tr X P + \Hb \c\ov{B}  -\frac{1}{4}\Xbh\c \ov{A},\\
\nabc_4\tr X +\frac{1}{2}(\tr X)^2 &=&-\frac{1}{2}\Xh\c\ov{\Xh},
  \eeaa
we obtain
\beaa
L&=& 3 \big( \frac 1 2 \ov{\tr X}-\tr X\big) P\nabc_4\Ab+ 3\big( \frac 1 2 \ov{\tr X}-\tr X\big) \nabc_4(P)\Ab\\
&&+ 3\big( \frac 1 2 \nabc_4\ov{\tr X}-\nabc_4\tr X\big) P\Ab+3 \big( 2  \tr X+\frac 1 2 \ov{\tr X}\big)\big( \frac 1 2 \ov{\tr X}-\tr X\big) P\Ab\\
&=& 3 \big( \frac 1 2 \ov{\tr X}-\tr X\big) P\big(\Ab_4-\frac 1 2 \tr X \Ab  \big)+ 3\big( \frac 1 2 \ov{\tr X}-\tr X\big)\big(-\frac{3}{2}\tr X P \big)\Ab\\
&&+ 3\big( -\frac 1 2\frac{1}{2}(\ov{\tr X})^2+\frac{1}{2}(\tr X)^2\big) P\Ab+3 \big( 2  \tr X+\frac 1 2 \ov{\tr X}\big)\big( \frac 1 2 \ov{\tr X}-\tr X\big) P\Ab\\
&&+ r^{-2} \dk^{\leq 1}((A, B) \c \Ga_b),
\eeaa
which gives
 \beaa
L&=&  \big( \frac 1 2 \ov{\tr X}-\tr X\big) 3P\Ab_4 +\frac 1 2 \tr X  \big( \tr X-  \ov{\tr X} \big)3 P\Ab+ r^{-2} \dk^{\leq 1}( (A,B) \c \Ga_b).
\eeaa

         { \bf Step 5. Calculation of $M$.} Using \eqref{eq:nabc4-DDbcAb}, we have
         \beaa
   M&=&\mathcal{J}_{434} \hot \nabc_4\big(\DDbc \c\Ab +  \ov{H} \c \Ab \big)+\nabc_4\big( \mathcal{J}_{434}\big) \hot \big(\DDbc \c\Ab +  \ov{H} \c \Ab \big)\\
   &&+ \big( 2  \tr X+\frac 1 2 \ov{\tr X}\big) \mathcal{J}_{434} \hot \big(\DDbc \c\Ab +  \ov{H} \c \Ab \big)\\
   &=&\mathcal{J}_{434} \hot \Big( \DDbc \c \Ab_4+(\ov{H}+\ov{\Hb}) \c \Ab_4 \Big)\\
   && + \Big( \nabc_4\big( \mathcal{J}_{434}\big)+  \frac 3 2   \tr X \mathcal{J}_{434} \Big) \hot \big(\DDbc \c\Ab +  \ov{H} \c \Ab \big)\\
   &&+ r^{-4} \dk^{\leq 1}(\Ga_g \c  \Ga_b)+ r^{-2} \dk^{\leq 1}( B \c \Ga_b).
   \eeaa

  { \bf Step 6. Final sum.} We deduce
  \beaa
&&\Big(\nabc_4 + 2  \tr X+\frac 1 2 \ov{\tr X}\Big)\Big( \nabc_4+  \tr X+ \frac 1 2\ov{\tr X}\Big) \Big( \nabc_3+ 2 \ov{\tr \Xb}+\frac 1 2 \tr \Xb\Big) \Ab_4\\
&=&   \frac 1 4  \big( \DDc+H+6\Hb\big)\hot \Big( \DDbc \c \widetilde{\underline{Q}(\Ab)}+(\ov{H}+2\ov{\Hb}) \c \widetilde{\underline{Q}(\Ab)} \Big)+3P \ \widetilde{\underline{Q}(\Ab)}\\
         &&-\frac 1 4 \big(\tr X- \ov{\tr X} - 2\frac {\atrch^2}{ \trch} \big) \DDc\hot \big( \DDbc \c \Ab_4+(\ov{H} +\ov{\Hb}) \c \Ab_4\big)\\
                  &&+\big( \frac 1 4 \mathcal{J}+\mathcal{K}+\mathcal{J}_{434} \big)\hot  \big( \DDbc \c \Ab_4+(\ov{H} +\ov{\Hb}) \c \Ab_4\big)+\frac1 2\ov{\mathcal{I}} \c \DDc \Ab_4+\frac 12(\DDc \c\ov{\mathcal{I}}) \Ab_4\\
         &&+ \frac 1 4 \big( H+6\Hb \big) \hot (\ov{\mathcal{I}} \c \Ab_4) + \big(-3\tr X+2\ov{\tr X} + 2\frac {\atrch^2}{ \trch} \big)3P \Ab_4\\
         &&+\frac 1 2 \tr X  \big( \tr X-  \ov{\tr X} \big)3 P\Ab + \Big( \nabc_4\big( \mathcal{J}_{434}\big)+  \frac 3 2   \tr X \mathcal{J}_{434} \Big) \hot \big(\DDbc \c\Ab +  \ov{H} \c \Ab \big)\\
         &&+\widetilde{\err_{4434}},
\eeaa
where the error terms are given by
\beaa
\widetilde{\err_{4434}}&:=& \nabc_4  \err_{434}+ \big( 2  \tr X+\frac 1 2 \ov{\tr X}\big)  \err_{434}\\
&&+ r^{-2} \dk^{\leq 1}((A,B) \c \Ga_b)+ r^{-4} \dk^{\leq 2}(\Ga_g \c  \Ga_b).
\eeaa
Finally, we use \eqref{eq:Ab-434} to substitute the term $ \DDc\hot \big( \DDbc \c \Ab_4+(\ov{H} +\ov{\Hb}) \c \Ab_4\big)$, given by
\beaa
&&  \frac 1 4  \DDc \hot \Big( \DDbc \c \Ab_4+(\ov{H}+\ov{\Hb}) \c \Ab_4\Big)\\
&=&\Big( \nabc_4+  \tr X+ \frac 1 2\ov{\tr X}\Big) \Big( \nabc_3+ 2 \ov{\tr \Xb}+\frac 1 2 \tr \Xb\Big) \Ab_4\\
&&-  \frac 1 4\big( H+5\Hb \big) \hot \Big( \DDbc \c \Ab_4+(\ov{H}+\ov{\Hb}) \c \Ab_4\Big)-3P\Ab_4\\
&& -3\Big( \frac 1 2\ov{\tr X}-\tr X\Big) P\Ab-\mathcal{J}_{434} \hot \big(\DDbc \c\Ab +  \ov{H} \c \Ab \big)+  \err_{434},
\eeaa
and we obtain
  \beaa
&&\Big(\nabc_4 + 2  \tr X+\frac 1 2 \ov{\tr X}\Big)\Big( \nabc_4+  \tr X+ \frac 1 2\ov{\tr X}\Big) \Big( \nabc_3+ 2 \ov{\tr \Xb}+\frac 1 2 \tr \Xb\Big) \Ab_4\\
         &&+ \big(\tr X- \ov{\tr X} - 2\frac {\atrch^2}{ \trch} \big)\Big( \nabc_4+  \tr X+ \frac 1 2\ov{\tr X}\Big) \Big( \nabc_3+ 2 \ov{\tr \Xb}+\frac 1 2 \tr \Xb\Big) \Ab_4\\
         &&+ \big(2\tr X- \ov{\tr X}  \big)3P\Ab_4 \\
&=&   \frac 1 4  \big( \DDc+H+6\Hb\big)\hot \Big( \DDbc \c \widetilde{\underline{Q}(\Ab)}+(\ov{H}+2\ov{\Hb}) \c \widetilde{\underline{Q}(\Ab)} \Big)+3P \ \widetilde{\underline{Q}(\Ab)}\\
   &&+\mathscr{L}[\Ab, \Ab_4, \DDc \Ab_4, \DD \Ab]+\err_{4434},
\eeaa
where
\beaa
\err_{4434}&=& \nabc_4  \err_{434}+ \big(   \tr X+\frac 3 2 \ov{\tr X}+2\frac {\atrch^2}{ \trch} \big)  \err_{434}\\
&&+ r^{-2} \dk^{\leq 1}((A,B) \c \Ga_b)+ r^{-4} \dk^{\leq 2}(\Ga_g \c  \Ga_b),
\eeaa
and where $\mathscr{L}[\Ab, \Ab_4, \DDc \Ab_4, \DD \Ab]$ are $O(|a|)$ linear order terms in $\Ab, \Ab_4, \DDc \Ab_4, \DD \Ab$ explicitly given by
\bea\begin{split}\label{eq:mathscr-L-Ab}
&\mathscr{L}[\Ab, \Ab_4, \DDc \Ab_4, \DD \Ab]\\
&= \big( \frac 1 4 \mathcal{J}+\mathcal{K}+\mathcal{J}_{434}+\frac 1 4 \big(\tr X- \ov{\tr X} - 2\frac {\atrch^2}{ \trch} \big)\big( H+5\Hb \big)  \big)\hot  \big( \DDbc \c \Ab_4+(\ov{H} +\ov{\Hb}) \c \Ab_4\big)\\
&+ \Big( \nabc_4\big( \mathcal{J}_{434}\big)+  \frac 3 2   \tr X \mathcal{J}_{434}+ \big(\tr X- \ov{\tr X} - 2\frac {\atrch^2}{ \trch} \big)\mathcal{J}_{434}  \Big) \hot \big(\DDbc \c\Ab +  \ov{H} \c \Ab \big)\\
&+\Big(  \big(\tr X- \ov{\tr X} - 2\frac {\atrch^2}{ \trch} \big)\big( \frac 1 2\ov{\tr X}-\tr X\big)+\frac 1 2 \tr X  \big( \tr X-  \ov{\tr X} \big) \Big)3 P\Ab\\
                                 &+\frac1 2\ov{\mathcal{I}} \c \DDc \Ab_4+\frac 12(\DDc \c\ov{\mathcal{I}}) \Ab_4+ \frac 1 4 \big( H+6\Hb \big) \hot (\ov{\mathcal{I}} \c \Ab_4). 
                                 \end{split}
\eea

  { \bf Step 7. Analysis of lower order terms.} We now analyze each term in $\mathscr{L}[\Ab, \Ab_4, \DDc \Ab_4, \DD \Ab]$. 
  
  The term in $\big( \DDbc \c \Ab_4+(\ov{H} +\ov{\Hb}) \c \Ab_4\big)$ is given by, using \eqref{eq:mathcal-J}, \eqref{eq:mathcal-K} and \eqref{eq:mathcal-J-434-r-1Gag},
  \beaa
 &&  \frac 1 4 \mathcal{J}+\mathcal{K}+\mathcal{J}_{434}+\frac 1 4 \big(\tr X- \ov{\tr X} - 2\frac {\atrch^2}{ \trch} \big)\big( H+5\Hb \big) \\
 &=&     \frac 1 4\Big(  i\Im(\tr  X) (H-\Hc) + 2\DDc(\frac {\atrch^2}{ \trch}) +\big(\frac{3}{2}\tr X+\frac 1 2 \ov{\tr X} + 2\frac {\atrch^2}{ \trch} \big) \Hb \Big)\\
 && - \frac 1 4 \Big[  \frac{1}{2} \big(\tr X- \ov{\tr X} - 4\frac {\atrch^2}{ \trch} \big)   H + \big(\frac{15}{2}\tr X-\frac{11}{2} \ov{\tr X} - 10\frac {\atrch^2}{ \trch} \big)  \Hb  \Big]\\
 &&-\frac 3 4 \tr X \Hb+\frac 1 4 \big(\tr X- \ov{\tr X} - 2\frac {\atrch^2}{ \trch} \big)\big( H+5\Hb \big) +r^{-1} \Ga_g\\
  &=&     \frac 1 4\Big[  2\DDc(\frac {\atrch^2}{ \trch})-3  \tr X \Hb -\big(\tr X- \ov{\tr X} - 2\frac {\atrch^2}{ \trch} \big) \Hb  + \big(\tr X- \ov{\tr X}  \big)H   \Big] \\
  &&+\Im(\tr X) \Hc+r^{-1} \Ga_g.
  \eeaa
  In particular, this is given by
  \beaa
&& \big( \frac 1 4 \mathcal{J}+\mathcal{K}+\mathcal{J}_{434}+\frac 1 4 \big(\tr X- \ov{\tr X} - 2\frac {\atrch^2}{ \trch} \big)\big( H+5\Hb \big)  \big)\hot  \big( \DDbc \c \Ab_4+(\ov{H} +\ov{\Hb}) \c \Ab_4\big) \\
&=&O(ar^{-3}) \big( \DDbc \c \Ab_4+(\ov{H} +\ov{\Hb}) \c \Ab_4\big) + r^{-4}( \Ga_g \c \Ga_b).
  \eeaa

  To compute the term in $\big(\DDbc \c\Ab +  \ov{H} \c \Ab \big)$, we first compute using  \eqref{eq:mathcal-J-434-r-2Gag},
  \beaa
  \nabc_4 \big( \mathcal{J}_{434}\big)&=&  \frac 1 4 \Big[-\frac 1 2 \nabc_4 \DDc( \tr X+\ov{\tr X})+ \frac 1 2 (\tr X -\ov{\tr X})\nabc_4 H\\
  &&+ \frac 1 2 \nabc_4(\tr X -\ov{\tr X}) H   -4 \tr X\nabc_4\Hb  -4 \nabc_4\tr X\Hb  \Big]\\
 &&-\frac 1 4 \nabc_4B+r^{-3} \Ga_g+\nabc_4 (\Xh \c \Hc)\\
 &=&  \frac 1 4 \Big[-\frac 1 2 \nabc_4 \DDc( \tr X+\ov{\tr X}) -\frac{1}{4}\ov{\tr X}  (\tr X -\ov{\tr X}) (H-\Hb)  \\
  &&+ \frac 1 2 \nabc_4(\tr X -\ov{\tr X}) H   +4 (\tr X)^2\Hb  -4 \nabc_4\tr X\Hb  \Big]\\
 &&-\frac 1 4 \nabc_4B+r^{-3} \Ga_g+\nabc_4 (\Xh \c \Hc).
  \eeaa
  Using the commutator for a scalar function $h$ of signature $s$
  \beaa
 \, [\nabc_4 , \DDc]h  &=& -\frac{1}{2}\tr X\DDc h+\Hb\nabc_4 h -\frac{1}{2}\Xh\c\ov{\DDc}h\\
&&  +s\left(\frac{1}{2}\tr X\Hb  +\frac{1}{2}\widehat{X}\c\ov{\Hb} -B\right)h  + (\Ga_b \c \Ga_g) h,
 \eeaa
  and making use of the equation
  \beaa
  \nabc_4\tr X  &=&-\frac{1}{2}(\tr X)^2+\Ga_g \c \Ga_g, 
  \eeaa
  we have
   \beaa
\nabc_4 \big( \mathcal{J}_{434}\big)  &=&  \frac 1 4 \Big[-\frac 1 2 \DDc(  \nabc_4\tr X+ \nabc_4\ov{\tr X})\\
  &&-\frac 1 2\big[ -\frac{1}{2}\tr X\DDc ( \tr X+\ov{\tr X}) +\Hb\nabc_4 ( \tr X+\ov{\tr X})\\
  &&  -\frac{1}{2}\Xh\c\ov{\DDc}( \tr X+\ov{\tr X})   +( \tr X+\ov{\tr X}) \big(\frac{1}{2}\tr X\Hb  -B\big)   \big]\\
  &&-\frac{1}{4}\ov{\tr X}  (\tr X -\ov{\tr X}) (H-\Hb)  \\
  &&+ \frac 1 2 (-\frac{1}{2}(\tr X)^2 +\frac{1}{2}(\ov{\tr X})^2) H   +4 (\tr X)^2\Hb  -4 (-\frac{1}{2}(\tr X)^2)\Hb  \Big]\\
 &&-\frac 1 4 \nabc_4B+r^{-3} \Ga_g+\nabc_4 (\Xh \c \Hc).
  \eeaa
  which gives
    \beaa
  \nabc_4 \big( \mathcal{J}_{434}\big) &=&  \frac 1 4 \Big[\frac 1 2 \tr X \DDc (\tr X) + \frac 1 2 \ov{\tr X} \DDc(\ov{\tr X}) + \frac{1}{4}\tr X\DDc ( \tr X+\ov{\tr X}) \\
  && + \big(-\frac{1}{4}(\tr X)^2 -\frac{1}{4}\ov{\tr X}  \tr X+\frac{1}{2}(\ov{\tr X})^2\big) H   + 6 (\tr X)^2 \Hb  \Big]\\
 &&-\frac 1 4 \nabc_4B+\frac 1 8 ( \tr X+\ov{\tr X}) B+\frac{1}{16}\Xh\c\ov{\DDc}( \tr X+\ov{\tr X}) \\
 && +r^{-3} \Ga_g+\nabc_4 (\Xh \c \Hc).
 \eeaa
We therefore obtain for the term in $\big(\DDbc \c\Ab +  \ov{H} \c \Ab \big)$,
  \beaa
&& \nabc_4\big( \mathcal{J}_{434}\big)+  \frac 3 2   \tr X \mathcal{J}_{434}+ \big(\tr X- \ov{\tr X} - 2\frac {\atrch^2}{ \trch} \big)\mathcal{J}_{434} \\
&=& \frac 1 4 \Big[ \frac 1 2( \ov{\tr X}-\tr X) \DDc(\ov{\tr X})+ \big(\frac{1}{2}(\tr X)^2 -\ov{\tr X}  \tr X+\frac{1}{2}(\ov{\tr X})^2\big) H     \Big]\\
 &&+ \big(\tr X- \ov{\tr X} - 2\frac {\atrch^2}{ \trch} \big)\mathcal{J}_{434}  +r^{-1} \dk^{\leq 1} B+r^{-3} \Ga_g+\nabc_4 (\Xh \c \Hc)+\frac 3 2 \tr X (\Xh \c \Hc).
   \eeaa
Observe that  $\frac12 (\tr X)^2-\ov{\tr X}\tr X +\frac 1 2 (\ov{\tr X})^2=\frac{1}{2}(\tr X - \ov{\tr X})^2 = -2i(\atrch)^2=O(r^{-4})$, and therefore
    \beaa
&& \nabc_4\big( \mathcal{J}_{434}\big)+  \frac 3 2   \tr X \mathcal{J}_{434}+ \big(\tr X- \ov{\tr X} - 2\frac {\atrch^2}{ \trch} \big)\mathcal{J}_{434} \\
&=&O(a^2r^{-5}) +r^{-1} \dk^{\leq 1} B +r^{-3} \Ga_g+\nabc_4 (\Xh \c \Hc)+\frac 3 2 \tr X (\Xh \c \Hc).
   \eeaa

  The term in $\Ab$ is given by
\beaa
&&\Big(  \big(\tr X- \ov{\tr X} - 2\frac {\atrch^2}{ \trch} \big)\big( \frac 1 2\ov{\tr X}-\tr X\big)+\frac 1 2 \tr X  \big( \tr X-  \ov{\tr X} \big) \Big)3 P\Ab\\
&=&\Big(-\frac 1 2   \big(\tr X- \ov{\tr X} \big)^2- 2\frac {\atrch^2}{ \trch} \big( \frac 1 2\ov{\tr X}-\tr X\big) \Big)3 P\Ab\\
&=&\Big(2  \atrch^2- 2\frac {\atrch^2}{ \trch} \big(-\frac 1 2 \trch+ \frac 3 2  i \atrch\big) \Big)3 P\Ab\\
&=&3 \atrch^2\Big(1 - \frac {  i \atrch}{ \trch} \Big)3 P\Ab=O(a^2r^{-7}) \Ab+ r^{-4}( \Ga_g \c \Ga_b).
\eeaa

Finally, using \eqref{eq:ov-mathcal-I} and putting together all the above we obtain from \eqref{eq:mathscr-L-Ab},
\beaa
&&\mathscr{L}[\Ab, \Ab_4, \DDc \Ab_4, \DD \Ab]\\
&=& O(ar^{-3}) \big( \DDbc \c \Ab_4+(\ov{H} +\ov{\Hb}) \c \Ab_4\big)+O(a^2r^{-4})\big(\DDbc \c\Ab +  \ov{H} \c \Ab \big) +O(a^2r^{-7}) \Ab\\
&&+\big( r^{-2} \dk^{\leq 1} B , r^{-1}(\nabc_4 (\Xh \c \Hc)+\frac 3 2 \tr X (\Xh \c \Hc))\big)  \dk^{\leq 1} \Ab+ r^{-4}( \Ga_g \c \Ga_b). 
\eeaa
This completes the proof of Proposition \ref{prop:appendix-nabc4nabc4nabc3Ab}.
\end{proof}

We finally show that the coefficient $\underline{W}_3=O(a^2 r^{-5})$. The coefficient comes from the symmetric of $W_4$ as in Theorem \ref{MAIN-THEOREM-PART1}, whose lowest decaying term is given by $Z_4$ in Proposition \ref{first-intermediate-step-main-theorem}, i.e. 
\beaa
\underline{Z}_3&=&\nabc_4 \widetilde{\underline{C}_2}+ 2\trch  \widetilde{\underline{C}_2}-\frac 1 4 (\trch^2+ \atrch^2) \widetilde{\underline{C}_1}.
\eeaa
Using the definition of $\underline{C}_1$ and $\underline{C}_2$ given by \eqref{eq:Cb1-Cb2-comparison-Ma}, we have
\beaa
\underline{Z}_3&=&\nabc_4 (-2i \trch \atrch)+ 2\trch (-2i \trch \atrch)-\frac 1 4 (\trch^2+ \atrch^2) (-4i \atrch)\\
&&+O(a^2 r^{-5})\\
&=&-2i \left(-\frac 1 2  \trch^2 \atrch+ \trch (-\trch \atrch)\right)+ 2\trch (-2i \trch \atrch)\\
&&+i \trch^2 \atrch+O(a^2 r^{-5})\\
&=& O(a^2 r^{-5}),
\eeaa
as stated.


\subsection{Treatment of the nonlinear terms}


The nonlinear terms we have encountered so far are as follows.

\begin{enumerate}
\item Nonlinear terms in the Teukolsky equation,  see Proposition \ref{PROP:TEUK-AB}, schematically given by
\beaa
\err_{TE}&=& \tr X \Xib \hot  \Bb+( \DDbc \c\Bb) \Xbh+ (\widehat{\Xb}\c\ov{\Hc}) \Bb \\
&&+(\Ga_b \c \Ga_b)\c (A, B)+ r^{-2} \dk^{\le 1} (\Ga_g \c \Ga_b).
\eeaa
 \item Nonlinear terms after the first commutation, see   Proposition   \ref{prop:appendix-nabc43Ab4},   written schematically   in the form
\beaa
  \err_{434}&=& \nabc_4\err_{TE} +\big( \frac 1 2 \tr X+ \ov{\tr X}\big)\err_{TE}\\
  &&+\DDc\hot( \Xh \c \ov{\DDc} \Ab)+ \DDc \hot ((  \DDc\c\ov{\Xh})\c \Ab)\\
 &&+ r^{-1} \dk^{\le 1 } ( (A,B)\c \Ga_b) + r^{-3}\dk^{\leq1}( \Ga_g  \c \Ga_b).
\eeaa
\item Nonlinear terms after the  second  commutation, see   Proposition \ref{prop:appendix-nabc4nabc4nabc3Ab},
\beaa
\err_{4434}&=& \nabc_4  \err_{434}+ \big(   \tr X+\frac 3 2 \ov{\tr X}+2\frac {\atrch^2}{ \trch} \big)  \err_{434}\\
&&+ r^{-2} \dk^{\leq 1}( (A,B) \c \Ga_b)+ r^{-4} \dk^{\leq 2}(\Ga_g \c  \Ga_b)\\
&&+ r^{-1}(\nabc_4 (\Xh \c \Hc)+\frac 3 2 \tr X (\Xh \c \Hc))\c  \dk^{\leq 1} \Ab.
\eeaa
\end{enumerate}

We can therefore separate the final error terms, $\err_{4434}$, as follows.
\beaa
\err_{4434}&=&\err_1 +\err_2+ r^{-2} \dk^{\leq 2}((A, B) \c \Ga_b)+ r^{-4} \dk^{\leq 2}(\Ga_g \c  \Ga_b)\\
&&+ r^{-1}\left(\nabc_4 (\Xh \c \Hc)+\frac 3 2 \tr X (\Xh \c \Hc)\right)\c  \dk^{\leq 1} \Ab.
\eeaa
where
\beaa
\err_1&=& \nabc_4  \big( \nabc_4\err_{TE} +\big( \frac 1 2 \tr X+ \ov{\tr X}\big)\err_{TE} \big)\\
&&+ \big(   \tr X+\frac 3 2 \ov{\tr X} \big) \big( \nabc_4\err_{TE} +\big( \frac 1 2 \tr X+ \ov{\tr X}\big)\err_{TE} \big), \\
\err_2&=& \nabc_4  \big(\DDc\hot( \Xh \c \ov{\DDc} \Ab)+ \DDc \hot ((  \DDc\c\ov{\Xh})\c \Ab) \big)\\
  &&+ \big(   \tr X+\frac 3 2 \ov{\tr X} \big)\big( \DDc\hot( \Xh \c \ov{\DDc} \Ab)+ \DDc \hot ((  \DDc\c\ov{\Xh})\c \Ab) \big).
\eeaa

We have the following for $\err_{4434}$. 

\begin{proposition}\label{prop:appendix-error-terms-qfb}
The following holds true: 
\beaa
\err_{4434}&=&r^{-2} \dk^{\leq 2}((A, B) \c \Ga_b)+ r^{-4} \dk^{\leq 3}(\Ga_g \c  \Ga_b) .
\eeaa
\end{proposition}
\begin{proof} We start with $\err_1$, which contains the terms coming from $\err_{TE}$.
We consider each term in $\err_{TE}$ separately. 

{\bf Step 1.} We consider the term $ \tr X \Xib \hot  \Bb$.
We first compute, 
\beaa
I_1&:=&\nabc_4( \tr X \Xib \hot  \Bb) +\big( \frac 1 2 \tr X+ \ov{\tr X}\big)( \tr X \Xib \hot  \Bb)\\
&=&\big( \nabc_4( \tr X)+\frac 1 2 (\tr X)^2 \big) \Xib \hot  \Bb + \tr X \nabc_4\Xib \hot  \Bb+ \tr X \Xib \hot \big(  \nabc_4\Bb+ \tr X \Bb\big) \\
&&+\big( \ov{\tr X}-\tr X\big)( \tr X \Xib \hot  \Bb) \\
&=& \tr X\big(-\frac 1 2 \ov{\tr X}\Hc-\Bb -\frac 1 2 \Xbh \c \Hc   \big) \hot  \Bb +r^{-3} \dk^{\leq 1} (\Ga_g \c \Ga_b),
\eeaa
where we used Proposition \ref{prop:preliminaries-qfb-null-structure} and Proposition \ref{prop:preliminaries-qfb-bianchi-lin} to write 
    \beaa
    \nabc_4 \Xib&=&-\frac 1 2 \ov{\tr X}\Hc-\Bb -\frac 1 2 \Xbh \c \Hc  + r^{-2}\Ga_b\\
    \nabc_4\Bb + \tr X\Bb   &=&r^{-2} \Ga_g.
    \eeaa
    We write schematically in the form 
    \beaa
I_1&=&(\tr X)^2 \Hc\hot \Bb      +\tr X (\Bb\hot \Bb)+    \tr X ( \Xbh \c \Hc)\hot \Bb+r^{-3} \dk^{\leq 2} (\Ga_g \c \Ga_b).
\eeaa
We then obtain that 
\beaa
I_2&:=& \nabc_4 I_1+  \big(   \tr X+\frac 3 2 \ov{\tr X} \big) I_1+r^{-4}\dk^{\leq 3} (\Ga_g \c \Ga_b),
\eeaa
is given as the sum of the following terms
  \beaa
 S_1&=& \nabc_4   \big((\tr X)^2 \Hc\hot \Bb \big)+  \big(   \tr X+\frac 3 2 \ov{\tr X} \big)  \big((\tr X)^2 \Hc\hot \Bb \big)\\
 S_2&=&  \nabc_4   \big(\tr X (\Bb\hot \Bb) \big)+ \big(   \tr X+\frac 3 2 \ov{\tr X} \big)   \big(\tr X (\Bb\hot \Bb) \big)\\
 S_3&=&\nabc_4  \big( \tr X ( \Xbh \c \Hc)\hot \Bb\big) +  \big(   \tr X+\frac 3 2 \ov{\tr X} \big)  \big( \tr X ( \Xbh \c \Hc)\hot \Bb\big).
  \eeaa
  
  We now analyze each of the above.
  
  {\bf Step 1a.} We have 
  \beaa
   S_1&=& \big( \nabc_4  (\tr X)^2+ (\tr X)^3\big) \Hc\hot \Bb + (\tr X)^2\big(  \nabc_4  \Hc + \frac 1 2 \ov{\tr X} \Hc \big) \hot \Bb\\
   &&+(\tr X)^2 \Hc\hot \big( \nabc_4  \Bb + \tr X \Bb\big)+r^{-4}\dk^{\leq 3} (\Ga_g \c \Ga_b).
  \eeaa
Recall that $H$ verifies the equation, see \eqref{equation:nabc_4Hc},
\beaa
\nabc_4\Hc+\frac{1}{2}\ov{\tr X}\Hc &=&
- \frac 1 2 \Xh \c \ov{\Hc}  - B+ r^{-2} \Ga_g,
\eeaa
and using again Proposition \ref{prop:preliminaries-qfb-bianchi-lin} and the null structure equation for $\nabc_4 \tr X$ we obtain 
    \beaa
   S_1&=&r^{-4} \dk^{\leq 3} (\Ga_g \c \Ga_b).
  \eeaa
  
    {\bf Step 1b.} We have 
    \beaa
     S_2&=&  \nabc_4   \big(\tr X (\Bb\hot \Bb) \big)+ \big(   \tr X+\frac 3 2 \ov{\tr X} \big)   \big(\tr X (\Bb\hot \Bb) \big)\\
     &=& \big( \nabc_4  (\tr X)+ \frac 1 2 (\tr X)^2\big) (\Bb\hot \Bb) +2  \tr X ( \nabc_4 \Bb+ \tr X \Bb) \hot \Bb+r^{-4}\dk^{\leq 3}  (\Ga_g \c \Ga_b)\\
     &=& r^{-4}\dk^{\leq 3}  (\Ga_g \c \Ga_b),
    \eeaa
    as above.
    
        {\bf Step 1c.} We have 
        \beaa
        S_3&=& \nabc_4  \big( \tr X ( \Xbh \c \Hc)\hot \Bb\big) +  \big(   \tr X+\frac 3 2 \ov{\tr X} \big)  \big( \tr X ( \Xbh \c \Hc)\hot \Bb\big)\\
        &=&\big( \nabc_4( \tr X)+\frac 1 2 (\tr X)^2\big)  ( \Xbh \c \Hc)\hot \Bb+ \big( \tr X (( \nabc_4 \Xbh+\frac 1 2 \tr X \Xbh) \c \Hc)\hot \Bb\big)\\
        &&+  \big( \tr X ( \Xbh \c \big(  \nabc_4 \Hc+\frac{1}{2}\ov{\tr X}\Hc \big) )\hot \Bb\big)+ \big( \tr X ( \Xbh \c \Hc)\hot (\nabc_4 \Bb+ 2 \tr X \Bb) \big) \\
        &&+ r^{-3} (\Ga_b \c \Ga_b) \c \Bb\\
                &=& r^{-4} \dk^{\leq 3} (\Ga_g \c \Ga_b),
        \eeaa
as above.
This completes the term involving $ \tr X \Xib \hot  \Bb$, by showing that 
\beaa
I_2=r^{-4} \dk^{\leq 3} (\Ga_g \c \Ga_b).
\eeaa

{\bf Step 2.} We consider the term $( \DDbc \c\Bb) \Xbh$.
We first compute 
\beaa
J_1&:=& \nabc_4( ( \DDbc \c\Bb) \Xbh) +\big( \frac 1 2 \tr X+ \ov{\tr X}\big)( ( \DDbc \c\Bb) \Xbh)\\
&=& \nabc_4 ( \DDbc \c\Bb) \Xbh +  ( \DDbc \c\Bb)\big(  \nabc_4\Xbh+\frac 1 2 \tr X \Xbh \big)+ \ov{\tr X} ( \DDbc \c\Bb) \Xbh\\
&=&   \DDbc \c(\nabc_4\Bb) \Xbh+([\nabc_4 , \DDbc \c]\Bb) \Xbh + \ov{\tr X} ( \DDbc \c\Bb) \Xbh\\
&&+  ( \DDbc \c\Bb)\big(  \nabc_4\Xbh+\frac 1 2 \tr X \Xbh \big).
\eeaa
Using  the  equation
\beaa
 \nabc_4\widehat{\Xb} +\frac{1}{2}\tr X\, \widehat{\Xb} &=& r^{-1} \Ga_g,
 \eeaa
  and the commutator  formula \eqref{eq:comm-nabc4nabc3-ovDDc-U-err},  written schematically
    \beaa
         \, [\nabc_4, \ov{\DDc}\c] \Bb&=&- \frac 1 2\ov{\tr X}\, \big( \ov{\DDc} \c  \Bb)+ r^{-3}\Bb  +  r^{-1}( \dk^{\leq 1}\Ga_g) \c  \Bb,
         \eeaa
         we have, using also from Lemma \ref{Lemma:DDc{trX}},  $\DDc\ov{\tr X}= a^2  r^{-3}  + r^{-2}\Ga_b      +r^{-1} \dk^{\le 1}\Ga_g$, 
         \beaa
J_1&=&   \DDbc \c(\nabc_4\Bb) \Xbh +\frac 1 2  \ov{\tr X} ( \DDbc \c\Bb) \Xbh+ r^{-3} (\Ga_g \c \Ga_b)\\
&=&   \DDbc \c(\nabc_4\Bb+ \tr X \Bb) \Xbh -\frac 1 2  \ov{\tr X} ( \DDbc \c\Bb) \Xbh+ r^{-3} (\Ga_g \c \Ga_b)\\
&=&    -\frac 1 2  \ov{\tr X} ( \DDbc \c\Bb) \Xbh+ r^{-3} \dk^{\leq 2} (\Ga_g \c \Ga_b).
\eeaa
We write schematically in the form
\beaa
J_1&=& \ov{\tr X} ( \DDbc \c\Bb) \Xbh+ r^{-3} \dk^{\leq 2} (\Ga_g \c \Ga_b), 
\eeaa
and we then obtain that 
\beaa
J_2&=&  \nabc_4 J_1+  \big(   \tr X+\frac 3 2 \ov{\tr X} \big) J_1+r^{-4} \dk^{\leq 3} (\Ga_g \c \Ga_b).
\eeaa
Using the equations stated above, we obtain
\beaa
J_2&=&  \nabc_4 (\ov{\tr X} ( \DDbc \c\Bb) \Xbh)+  \big(   \tr X+\frac 3 2 \ov{\tr X} \big) \ov{\tr X} ( \DDbc \c\Bb) \Xbh+r^{-4} \dk^{\leq 3} (\Ga_g \c \Ga_b)\\
&=&\big(  \nabc_4 (\ov{\tr X})+\frac 12 (\ov{\tr X})^2\big) ( \DDbc \c\Bb) \Xbh+  \ov{\tr X} \big( \nabc_4 ( \DDbc \c\Bb)+\frac 3 2 \ov{\tr X} \DDbc \c \Bb\big) \Xbh\\
&&+  \ov{\tr X} ( \DDbc \c\Bb)\big( \nabc_4 ( \Xbh)+\frac 1 2 \tr X \Xbh \big) +r^{-4} \dk^{\leq 3} (\Ga_g \c \Ga_b)\\
&=& r^{-4} \dk^{\leq 3} (\Ga_g \c \Ga_b),
\eeaa
which completes the term involving $ ( \DDbc \c\Bb) \Xbh$.

{\bf Step 3.} We consider the term $(\widehat{\Xb}\c\ov{\Hc}) \Bb$. Following a similar pattern as above we obtain
\beaa
K_1&:=& \nabc_4 ((\widehat{\Xb}\c\ov{\Hc}) \Bb) +\big( \frac 1 2 \tr X+ \ov{\tr X}\big)( (\widehat{\Xb}\c\ov{\Hc}) \Bb)\\
&=& ( (\nabc_4 \widehat{\Xb}+\frac 1 2 \tr X \Xbh) \c\ov{\Hc}) \Bb+ (\widehat{\Xb}\c (\nabc_4\ov{\Hc}+\frac 1 2 \tr X \ov{\Hc} )) \Bb\\
&&+ (\widehat{\Xb}\c\ov{\Hc})(\nabc_4 \Bb+ \frac 1 2 \tr X \Bb) -\frac 1 2  \ov{\tr X} (\widehat{\Xb}\c\ov{\Hc}) \Bb+ r^{-1} \dk^{\leq 1} (\Ga_b \c \Ga_g) \c \Bb\\
&=&  -\frac 1 2  \ov{\tr X} (\widehat{\Xb}\c\ov{\Hc}) \Bb+ r^{-1} \dk^{\leq 1} (\Ga_b \c \Ga_g) \c \Bb.
\eeaa
By taking the second derivative, we finally get 
\beaa
J_2&=&  \nabc_4 K_1+  \big(   \tr X+\frac 3 2 \ov{\tr X} \big) K_1+r^{-4} \dk^{\leq 3} (\Ga_g \c \Ga_b)\\
&=&  \nabc_4( \ov{\tr X} (\widehat{\Xb}\c\ov{\Hc}) \Bb)+  \big(   \tr X+\frac 3 2 \ov{\tr X} \big)  \ov{\tr X} (\widehat{\Xb}\c\ov{\Hc}) \Bb+r^{-4} \dk^{\leq 3} (\Ga_g \c \Ga_b)\\
&=& \big( \nabc_4( \ov{\tr X})+ \frac 1 2 (\ov{\tr X})^2\big) (\widehat{\Xb}\c\ov{\Hc}) \Bb+ \ov{\tr X} ((\nabc_4\widehat{\Xb}+ \frac 12 \tr X \Xbh)\c\ov{\Hc}) \Bb\\
&&+\ov{\tr X} (\widehat{\Xb}\c( \nabc_4\ov{\Hc}+\frac 1 2 \tr X \ov{\Hc})) \Bb+  \ov{\tr X} (\widehat{\Xb}\c\ov{\Hc})(\nabc_4 \Bb+\tr X \Bb)+r^{-4} \dk^{\leq 3} (\Ga_g \c \Ga_b)\\
&=& r^{-4} \dk^{\leq 3} (\Ga_g \c \Ga_b).
\eeaa

{\bf Step 4.} Combining the three steps above with the remaining terms in $\err_{TE}$, we finally obtain
\beaa
\err_1&=& r^{-4} \dk^{\leq 3} (\Ga_g \c \Ga_b) +r^{-3} \dk^{\leq 2}  \Ga_b \c (A, B).
\eeaa

{\bf Step 5.} We now prove the bound for $\err_2$. We first treat the term $\DDc\hot( \Xh \c \ov{\DDc} \Ab)$, while the other term can be treated in a similar manner. We have
\beaa
L&:=& \nabc_4  \big(\DDc\hot( \Xh \c \ov{\DDc} \Ab) \big)+ \big(   \tr X+\frac 3 2 \ov{\tr X} \big)\big( \DDc\hot( \Xh \c \ov{\DDc} \Ab) \big)\\
&=& \DDc\hot\nabc_4  ( \Xh \c \ov{\DDc} \Ab)+ [\nabc_4 , \DDc\hot]( \Xh \c \ov{\DDc} \Ab)\\
&&+ \big(   \tr X+\frac 3 2 \ov{\tr X} \big)\big( \DDc\hot( \Xh \c \ov{\DDc} \Ab) \big)\\
&=& \DDc\hot  ( \nabc_4\Xh \c \ov{\DDc} \Ab)+ \DDc\hot ( \Xh \c \nabc_4\ov{\DDc} \Ab)-\frac 1 2 \tr X \DDc\hot( \Xh \c \ov{\DDc} \Ab)\\
&&+ \big(   \tr X+\frac 3 2 \ov{\tr X} \big)\big( \DDc\hot( \Xh \c \ov{\DDc} \Ab) \big)+ r^{-4} (\Ga_g \c \Ga_b).
\eeaa
Using the null structure equation 
\beaa
\nabc_4\Xh+\Re(\tr X)\Xh&=&-A,
\eeaa
we then obtain
\beaa
L&=& \DDc\hot  ( (\nabc_4\Xh+\Re(\tr X)\Xh)   \c \ov{\DDc} \Ab)+ \DDc\hot ( \Xh \c \ov{\DDc}\nabc_4 \Ab)\\
&&+ \DDc\hot ( \Xh \c [\nabc_4,\ov{\DDc}] \Ab)+ \big(   \ov{\tr X} \big)\big( \DDc\hot( \Xh \c \ov{\DDc} \Ab) \big)+ r^{-4} (\Ga_g \c \Ga_b)\\
&=&- \DDc\hot  ( A  \c \ov{\DDc} \Ab)+ \DDc\hot ( \Xh \c \ov{\DDc}\nabc_4 \Ab)\\
&&-\frac 1 2 \ov{\tr X}  \DDc\hot ( \Xh \c\DDc \Ab)+ \big(   \ov{\tr X} \big)\big( \DDc\hot( \Xh \c \ov{\DDc} \Ab) \big)+ r^{-4} (\Ga_g \c \Ga_b)\\
&=&  r^{-2} \dk^{\leq 2}( A \c \Ga_b)+ \DDc\hot ( \Xh \c \ov{\DDc}\Ab_4)+ r^{-3} (\Ga_g \c \Ga_b)\\
&=&  r^{-4} \dk^{\leq 3} (\Ga_g \c \Ga_b)+  r^{-2} \dk^{\leq 2}( A \c \Ga_b).
\eeaa
  This shows 
  \beaa
  \err_2&=&  r^{-4} \dk^{\leq 3} (\Ga_g \c \Ga_b)+  r^{-2} \dk^{\leq 2}( A \c \Ga_b).
  \eeaa

{\bf Step 6.} We are now left to treat the following term in $\err_{4434}$:
\beaa
M&:=&  r^{-1}(\nabc_4 (\Xh \c \Hc)+\frac 3 2 \tr X (\Xh \c \Hc))\c  \dk^{\leq 1} \Ab\\
&=&  r^{-1}(\big( \nabc_4 (\Xh)+ \Re(\tr X) \Xh\big) \c \Hc+\Xh \c \big( \nabc_4 (\Hc)+\frac 1 2 \ov{\tr X}  \Hc\big))\c  \dk^{\leq 1} \Ab\\
&&+  r^{-4} \dk^{\leq 3} (\Ga_g \c \Ga_b)\\
&=&  r^{-1}( A  \c \Hc+\Xh \c \big( \Xh \c \Hc + B  \big))\c  \dk^{\leq 1} \Ab+  r^{-4} \dk^{\leq 3} (\Ga_g \c \Ga_b),
\eeaa
which gives
\beaa
M&=&   r^{-4} \dk^{\leq 3} (\Ga_g \c \Ga_b)+  r^{-2} \dk^{\leq 2}( (A, B) \c \Ga_b).
\eeaa

  {\bf Step 7.} Finally, we combine
  \beaa
\err_{4434}&=&\err_1 +\err_2+ r^{-2} \dk^{\leq 2}((A, B) \c \Ga_b)+ r^{-4} \dk^{\leq 2}(\Ga_g \c  \Ga_b)\\
&&+ r^{-1}\left(\nabc_4 (\Xh \c \Hc)+\frac 3 2 \tr X (\Xh \c \Hc)\right)\c  \dk^{\leq 1} \Ab\\
&=& r^{-2} \dk^{\leq 2}( (A, B) \c \Ga_b)+  r^{-4} \dk^{\leq 3} (\Ga_g \c \Ga_b) ,
  \eeaa
which concludes the proof of the Proposition. 
\end{proof}


\section{Proof of Proposition \ref{THEOREM:TEUK-STAR}}
\label{proof:teukolsky-starobinski}



\subsection{Preliminaries}


We use \eqref{eq:DDP-Kerr} and Lemma \ref{lemma:additional-relations-Kerr} to define the following $O(\ep)$ quantities in perturbations of Kerr.

\begin{definition} We define the following $O(\ep)$ quantities in $O(\ep)$ perturbations of Kerr:
\beaa
\mathcal{A}_1&=& \DDc P+3P\Hb \in \sk_1(\CCC) \\
\mathcal{A}_2&=& \DDc\hot\Hb  +\Hb\hot\Hb \in \sk_2( \CCC)\\
\AA_3&=& 2\nabc_4(\Hb)  -\DDc\left(\tr X\right) \in \sk_1(\CCC)\\
\AA_4&=& \DDc\hot \DDc\left(\tr X\right)  +3\Hb\hot \DDc\left(\tr X\right) \in \sk_2(\CCC).
\eeaa
We also define
\beaa
\BB&=& \DDc\c \ov{B}+ 2\Hb \c\ov{B} \in \sk_0.
\eeaa
\end{definition}

We collect here some useful derivatives of the above.

\begin{lemma}\label{lemma:derivatives-AA1} We have, modulo quadratic terms,
\beaa
\nabc_4 \mathcal{A}_1&=&  -2\tr X\AA_1 +\frac{1}{2}\DDc\BB +2\Hb\BB+\frac 3 2P\left(\ov{H} \c \Xh-\ov{\tr\Xb} \Xi+\AA_3\right),
\eeaa
and
\beaa
\nabc_4 \DDc\hot \AA_1&=& -\frac 5 2 \tr X \DDc\hot \AA_1 + \frac 3 2 \tr X\Hb \hot\AA_1 +\frac{1}{2}\DDc\hot\DDc\BB +\frac 5 2 \Hb \hot \DDc\BB\\
&&+\DDc\hot\Big(\frac{3}{2} P \left(\ov{H} \c \Xh-\ov{\tr\Xb} \Xi+\AA_3\right) \Big)\\
&&+ \underline{H} \hot \Big(\frac{3}{2} P \left(\ov{H} \c \Xh-\ov{\tr\Xb} \Xi+\AA_3\right) \Big).
\eeaa
\end{lemma}
\begin{proof} Applying Lemma \ref{comm:scalar}, \eqref{eq:comm-nabc4-DDc-h-precise}, to $h=P$ and $s=0$, we obtain
 \beaa
 \, [\nabc_4 , \DDc]P  &=& -\frac{1}{2}\tr X\DDc P+\Hb\nabc_4 P -\frac{1}{2}\Xh\c\ov{\DDc}P+\Xi\nabc_3P\\
 &=& -\frac{1}{2}\tr X\Big(\DDc P+3P\Hb\Big)+\Hb\left(\frac{1}{2}\DDc\c \ov{B}+ \Hb \c\ov{B}\right)\\
 && -\frac{1}{2}\Xh\c\ov{\DDc}P -\frac{3}{2}\ov{\tr\Xb} P\Xi\\
  &=& -\frac{1}{2}\tr X\mathcal{A}_1+\Hb\left(\frac{1}{2}\DDc\c \ov{B}+ \Hb \c\ov{B}\right) +\frac{3}{2} P \ov{H} \c \Xh -\frac{3}{2}\ov{\tr\Xb} P\Xi
 \eeaa
 where we used $\ov{\DDc}P=-3P\ov{H}+O(\ep)$. This gives
\beaa
\nabc_4 \AA_1&=& \nabc_4\Big(\DDc P+3P\Hb\Big)\\
&=& \DDc\nabc_4P + [\nabc_4 , \DDc]P+3\nabc_4(P)\Hb+3P\nabc_4(\Hb)\\
&=& \DDc\left(-\frac{3}{2}\tr X P +\frac{1}{2}\DDc\c \ov{B}+ \Hb \c\ov{B}\right) \\
&&+3\left(-\frac{3}{2}\tr X P +\frac{1}{2}\DDc\c \ov{B}+ \Hb \c\ov{B}\right)\Hb+3P\nabc_4(\Hb)\\
&& -\frac{1}{2}\tr X\mathcal{A}_1+\Hb\left(\frac{1}{2}\DDc\c \ov{B}+ \Hb \c\ov{B}\right) +\frac{3}{2} P \ov{H} \c \Xh -\frac{3}{2}\ov{\tr\Xb} P\Xi\\
 &=& -\frac{3}{2}\tr X\Big(\DDc P+3P\Hb\Big) +\DDc\left(\frac{1}{2}\DDc\c \ov{B}+ \Hb \c\ov{B}\right) \\
&&+3\Hb\left(\frac{1}{2}\DDc\c \ov{B}+ \Hb \c\ov{B}\right)+3\left(\nabc_4(\Hb)  -\frac{1}{2}\DDc\left(\tr X\right)\right)P\\
&& -\frac{1}{2}\tr X\mathcal{A}_1+\Hb\left(\frac{1}{2}\DDc\c \ov{B}+ \Hb \c\ov{B}\right) +\frac{3}{2} P \ov{H} \c \Xh -\frac{3}{2}\ov{\tr\Xb} P\Xi,
\eeaa
which gives
\beaa
\nabc_4 \AA_1 &=& -2\tr X\AA_1 +\DDc\left(\frac{1}{2}\DDc\c \ov{B}+ \Hb \c\ov{B}\right) +4\Hb\left(\frac{1}{2}\DDc\c \ov{B}+ \Hb \c\ov{B}\right)\\
&&+3\left(\nabc_4(\Hb)  -\frac{1}{2}\DDc\left(\tr X\right)\right)P+\frac{3}{2} P \ov{H} \c \Xh -\frac{3}{2}\ov{\tr\Xb} P\Xi.
\eeaa
Recalling the definition of $\BB$ and $\AA_3$ we obtain the first identity.

To obtain the second identity we apply Lemma \ref{comm:1form}, \eqref{eq:comm-nabc4nabc3-DDchot-err}, to $F=\AA_1$ and $s=0$,  and we have
\beaa
[\nabc_4,\DDc\hot]\AA_1 &=&  - \frac 1 2 \tr X\left( \DDc\hot\AA_1 + \Hb\hot\AA_1\right)+ \underline{H} \hot \nabc_4\AA_1.
\eeaa
This gives
\beaa
&&\nabc_4 \DDc\hot \AA_1\\
&=& \DDc\hot  \nabc_4 \AA_1+[\nabc_4, \DDc\hot] \AA_1\\
&=& \DDc\hot  \nabc_4 \AA_1- \frac 1 2 \tr X\left( \DDc\hot\AA_1 + \Hb\hot\AA_1\right)+ \underline{H} \hot \nabc_4\AA_1\\
&=& \DDc\hot \Big( -2\tr X\AA_1 +\frac{1}{2}\DDc\BB +2\Hb\BB+\frac{3}{2} P \left(\ov{H} \c \Xh-\ov{\tr\Xb} \Xi+\AA_3\right) \Big)\\
&&- \frac 1 2 \tr X\left( \DDc\hot\AA_1 + \Hb\hot\AA_1\right)\\
&&+ \underline{H} \hot \Big( -2\tr X\AA_1 +\frac{1}{2}\DDc\BB_1 +2\Hb\BB_1+\frac{3}{2} P \left(\ov{H} \c \Xh-\ov{\tr\Xb} \Xi+\AA_3\right) \Big)\\
&=& -\frac 5 2 \tr X \DDc\hot \AA_1 - \big(2 \DDc \tr X + \frac 5 2 \tr X\Hb\big)\hot\AA_1  \\
&&+\DDc\hot\Big(\frac{1}{2}\DDc\BB +2\Hb\BB\Big)+ \underline{H} \hot \Big(\frac{1}{2}\DDc\BB +2\Hb\BB \Big)\\
&&+\DDc\hot\Big(\frac{3}{2} P \left(\ov{H} \c \Xh-\ov{\tr\Xb} \Xi+\AA_3\right) \Big)+ \underline{H} \hot \Big(\frac{3}{2} P \left(\ov{H} \c \Xh-\ov{\tr\Xb} \Xi+\AA_3\right) \Big)\\
&=& -\frac 5 2 \tr X \DDc\hot \AA_1 - \big(2 \DDc \tr X + \frac 5 2 \tr X\Hb\big)\hot\AA_1  \\
&&+\frac{1}{2}\DDc\hot\DDc\BB +\frac 5 2 \Hb \hot \DDc\BB+2\big(\DDc\hot\Hb+ \Hb \hot \Hb \big)\BB\\
&&+\DDc\hot\Big(\frac{3}{2} P\left(\ov{H} \c \Xh-\ov{\tr\Xb} \Xi+\AA_3\right) \Big)+ \underline{H} \hot \Big(\frac{3}{2} P \left(\ov{H} \c \Xh-\ov{\tr\Xb} \Xi+\AA_3\right) \Big).
\eeaa
In view of $\DDc\hot\Hb+\Hb\hot\Hb=O(\ep)$ and $\DDc(\tr X) = -2\tr X\Hb+O(\ep)$ we obtain the stated.
\end{proof}

\begin{lemma}\label{lemma:DDcAA3} We have, modulo quadratic terms, 
\beaa
\DDc\hot \AA_3  -2\nabc_4\mathcal{A}_2&=&  \tr X \AA_2-3\Hb\hot\AA_3-\AA_4+2B \hot \Hb  \\
&&+(\tr \Xb+ \overline{\tr\Xb} ) \Hb \hot \Xi- \overline{\tr\Xb}  H \hot \Xi  +\Xh \c \ov{\DDc} \Hb  -\Xh (\ov{\Hb}\c \Hb).
\eeaa
\end{lemma}
\begin{proof} Recalling the definitions of $\AA_2$ and $\AA_3$ we obtain
\beaa
\DDc\hot \AA_3  -2\nabc_4\mathcal{A}_2&=& \DDc\hot \big( 2\nabc_4(\Hb)  -\DDc\left(\tr X\right) \big)  -2\nabc_4\big(\DDc\hot\Hb+\Hb\hot\Hb \big)\\
&=&  2[\DDc\hot ,\nabc_4]\Hb  -\DDc\hot \DDc\left(\tr X\right)  -4\nabc_4\Hb\hot\Hb. 
\eeaa
Applying Lemma \ref{comm:1form}, \eqref{eq:comm-nabc4nabc3DDchot-precise},  to $F=\Hb$ and $s=0$, we obtain, up to quadratic terms,
  \beaa
[ \nabc_4, \DDc \hot] \Hb   &=& -\frac 1 2 \tr X\left( \DDc \hot  \Hb+\Hb\hot \Hb \right)+\Hb\hot\nabc_4 \Hb+ \Xi \hot \nabc_3 \Hb\\
 &&-B \hot \Hb - \frac 1 2 \tr \Xb \Xi \hot  \Hb-\frac 1 2\Xh \c \ov{\DDc} \Hb+\frac12\Xh (\ov{\Hb}\c \Hb) \\
   &=& -\frac 1 2 \tr X \AA_2+\Hb\hot\nabc_4 \Hb-\frac 1 2(\tr \Xb+ \overline{\tr\Xb} ) \Hb \hot \Xi+\frac 1 2 \overline{\tr\Xb}  H \hot \Xi \\
 &&-B \hot \Hb -\frac 1 2\Xh \c \ov{\DDc} \Hb +\frac12\Xh (\ov{\Hb}\c \Hb) 
 \eeaa
where we used that $\nab_3\Hb=-\frac 1 2 \overline{\tr\Xb}(\Hb- H)+O(\ep)$. 
This gives
\beaa
\DDc\hot \AA_3  -2\nabc_4\mathcal{A}_2&=&  \tr X \AA_2-6\Hb\hot\nabc_4 \Hb+2B \hot \Hb -\DDc\hot \DDc\left(\tr X\right) \\
&&+(\tr \Xb+ \overline{\tr\Xb} ) \Hb \hot \Xi- \overline{\tr\Xb}  H \hot \Xi  +\Xh \c \ov{\DDc} \Hb -\Xh (\ov{\Hb}\c \Hb).
\eeaa
Writing $ 2\nabc_4(\Hb)=\AA_3+\DDc\left(\tr X\right)$, we obtain
\beaa
&&\DDc\hot \AA_3  -2\nabc_4\mathcal{A}_2\\
&=&  \tr X \AA_2-3\Hb\hot\AA_3+2B \hot \Hb -\DDc\hot \DDc\left(\tr X\right)-3\Hb \hot\DDc\left(\tr X\right) \\
&&+(\tr \Xb+ \overline{\tr\Xb} ) \Hb \hot \Xi- \overline{\tr\Xb}  H \hot \Xi +\Xh \c \ov{\DDc} \Hb -\Xh (\ov{\Hb}\c \Hb) .
\eeaa
Recalling the definition of $\AA_4$, this concludes the lemma.
\end{proof}

\begin{lemma}\label{lemma:nabc4AA4} We have, modulo quadratic terms,
    \beaa
&&\nabc_4\AA_4  \\
&=& -2(\tr X)\AA_4  - 3\tr X \Hb \hot   \AA_3 -\tr X \DDc\hot B +2\tr X \Hb \hot B +\MM[\Xi]\\
 &&-\frac{1}{2}(\ov{\tr X}-\tr X)\DDc\hot(\Xh\c\ov{H})  +\frac{1}{2}\tr X \DDc\hot(\widehat{X}\c\ov{\Hb})\\
    &&-\frac{1}{2}(\tr X-\ov{\tr X})(H \c \ov{H})\Xh  + \tr X (  \Hb \c \ov{\Hb}) \Xh+ \tr X\underline{H} \hot (\Xh\c\ov{H})+\Xh \c \ov{\DDc} (\tr X\Hb)
\eeaa
where $\MM[\Xi]$ only depends on $\Xi$, and therefore vanishes if $\Xi=0$. 
\end{lemma}
\begin{proof}
Recalling the definition of $\AA_4$ we obtain
\beaa
&&\nabc_4\AA_4\\
&=& \nabc_4 \big( \DDc\hot \DDc\left(\tr X\right)  +3\Hb\hot \DDc\left(\tr X\right)\big)\\
&=&\DDc\hot  \nabc_4  \DDc\left(\tr X\right)  + [\nabc_4 , \DDc\hot] \DDc\left(\tr X\right) \\
&&+3\nabc_4\Hb\hot \DDc\left(\tr X\right) +3\Hb\hot \nabc_4\DDc\left(\tr X\right)\\
&=&\DDc\hot   \DDc\left(\nabc_4  \tr X\right) +\DDc\hot \big( [\nabc_4 , \DDc]\tr X \big) + [\nabc_4 , \DDc\hot] \DDc\left(\tr X\right) \\
&&+3\nabc_4\Hb\hot \DDc\left(\tr X\right) +3\Hb\hot \DDc\left(\nabc_4\tr X\right)+3\Hb\hot [\nabc_4,\DDc]\left(\tr X\right)
\eeaa
Applying Lemma \ref{comm:scalar} to $h= \tr X$ and $s=1$, we obtain
 \beaa
&& \, [\nabc_4 , \DDc] \tr X  \\
&=& -\frac{1}{2}\tr X\DDc  \tr X+\Hb\nabc_4  \tr X -\frac{1}{2}\Xh\c\ov{\DDc} \tr X+\Xi\nabc_3 \tr X\\
&&  +\left(\frac{1}{2}\tr X\Hb  +\frac{1}{2}\widehat{X}\c\ov{\Hb}-\frac{1}{2}\tr\Xb\Xi -B\right) \tr X\\
  &=& -\frac{1}{2}\tr X\DDc  \tr X+\Hb\big( -\frac{1}{2}(\tr X)^2 + \DDc\c\ov{\Xi}+\Xi\c\ov{H}+\ov{\Xi}\c H \big) \\
  &&-\frac{1}{2}(\ov{\tr X}-\tr X)\Xh\c\ov{H}+\Xi\big(-\frac{1}{2}\tr\Xb\tr X + \DDc\c\ov{H}+H\c\ov{H}+2P \big)\\
&&  +\left(\frac{1}{2}\tr X\Hb  +\frac{1}{2}\widehat{X}\c\ov{\Hb}-\frac{1}{2}\tr\Xb\Xi -B\right) \tr X\\
&=& -\frac{1}{2}\tr X\DDc  \tr X  -\tr X B\\
&&+\Hb\big( \DDc\c\ov{\Xi}+\Xi\c\ov{H}+\ov{\Xi}\c H \big) +\big(-\tr\Xb\tr X + \DDc\c\ov{H}+H\c\ov{H}+2P \big)\Xi\\
  &&-\frac{1}{2}(\ov{\tr X}-\tr X)\Xh\c\ov{H}  +\frac{1}{2}\tr X\widehat{X}\c\ov{\Hb}.
 \eeaa
 where we used that $\nabc_3\tr X +\frac{1}{2}\tr\Xb\tr X = \DDc\c\ov{H}+H\c\ov{H}+2P+O(\ep^2)$, $\nabc_4\tr X +\frac{1}{2}(\tr X)^2 = \DDc\c\ov{\Xi}+\Xi\c\ov{H}+\ov{\Xi}\c H$. 
 
 Using that  $\ov{\tr X} \Hb + \tr X H=O(\ep)$, we also simplify
 \beaa
 -\frac{1}{2}(\ov{\tr X}-\tr X)\Xh\c\ov{H}  +\frac{1}{2}\tr X\widehat{X}\c\ov{\Hb}&=& -\frac{1}{2}\ov{\tr X}\Xh\c\ov{H} +\frac{1}{2}\tr X\Xh\c\ov{H}  +\frac{1}{2}\tr X\widehat{X}\c\ov{\Hb}\\
 &=& \tr X\Xh\c\ov{\Hb} +\frac{1}{2}\tr X\Xh\c\ov{H} , 
 \eeaa
 which gives
  \beaa
 \, [\nabc_4 , \DDc] \tr X &=& -\frac{1}{2}\tr X\DDc  \tr X  -\tr X B\\
&&+\Hb\big( \DDc\c\ov{\Xi}+\Xi\c\ov{H}+\ov{\Xi}\c H \big) \\
&&+\big(-\tr\Xb\tr X + \DDc\c\ov{H}+H\c\ov{H}+2P \big)\Xi\\
  &&+ \tr X\Xh\c\ov{\Hb} +\frac{1}{2}\tr X\Xh\c\ov{H} .
 \eeaa

Applying Lemma \ref{comm:1form},  \eqref{eq:comm-nabc4nabc3DDchot-precise}, to $F=\DDc \tr X$ and $s=1$ we obtain
  \beaa
&&[ \nabc_4, \DDc \hot] \DDc \tr X  \\
 &=& -\frac 1 2 \tr X \DDc \hot  \DDc \tr X+\Hb\hot\nabc_4 \DDc \tr X+ \Xi \hot \nabc_3 \DDc \tr X\\
 &&-2B \hot \DDc \tr X - \tr \Xb \Xi \hot  \DDc \tr X-\frac 1 2\Xh \c \ov{\DDc} \DDc \tr X\\
 &&+\frac12\Xh (\ov{\Hb}\c \DDc \tr X) +  \frac{1}{2}(\widehat{X}\c\ov{\Hb})  \hot \DDc \tr X\\
 &=& -\frac 1 2 \tr X \DDc \hot  \DDc \tr X+ \underline{H} \hot  \DDc \nabc_4 \tr X+ \underline{H} \hot [\nabc_4, \DDc] \tr X\\
 &&+ \Xi \hot \nabc_3( -2\tr X\Hb)-2B \hot ( -2\tr X\Hb) - \tr \Xb \Xi \hot ( -2\tr X\Hb)\\
 &&-\frac 1 2\Xh \c \ov{\DDc} ( -2\tr X\Hb)\\
 &&+\frac12\Xh (\ov{\Hb}\c ( -2\tr X\Hb)) +  \frac{1}{2}(\widehat{X}\c\ov{\Hb})  \hot ( -2\tr X\Hb)\\
  &=& -\frac 1 2 \tr X \DDc \hot  \DDc \tr X+ \underline{H} \hot  \DDc \nabc_4 \tr X+ \underline{H} \hot [\nabc_4, \DDc] \tr X\\
 &&+4\tr X B \hot \Hb-2 \Xi \hot \big(  \nabc_3( \tr X\Hb) - \tr X \tr \Xb \Hb\big)\\
 &&+\Xh \c \ov{\DDc} (\tr X\Hb)-2 \tr X (  \Hb \c \ov{\Hb}) \Xh 
 \eeaa
  where we used that $\DDc(\tr X) = -2\tr X\Hb+O(\ep)$.

  We therefore obtain
  \beaa
\nabc_4\AA_4&=&\DDc\hot   \DDc\left(\nabc_4  \tr X\right) +\DDc\hot \big( [\nabc_4 , \DDc]\tr X \big)+ 4\underline{H} \hot [\nabc_4, \DDc] \tr X\\
&&-\frac 1 2 \tr X \DDc \hot  \DDc \tr X+4 \underline{H} \hot  \DDc \nabc_4 \tr X+3\nabc_4\Hb\hot \DDc\left(\tr X\right) \\
 &&+4\tr X B \hot \Hb-2 \Xi \hot \big(  \nabc_3( \tr X\Hb) - \tr X \tr \Xb \Hb\big)\\
 &&+\Xh \c \ov{\DDc} (\tr X\Hb) -2 \tr X (  \Hb \c \ov{\Hb}) \Xh .
\eeaa

Using the Ricci identity $\nabc_4  \tr X=-\frac{1}{2}(\tr X)^2 + \DDc\c\ov{\Xi}+\Xi\c\ov{H}+\ov{\Xi}\c H+O(\ep^2)$ we compute
\beaa
\DDc\left(\nabc_4  \tr X\right) &=&\DDc\left(-\frac{1}{2}(\tr X)^2 + \DDc\c\ov{\Xi}+\Xi\c\ov{H}+\ov{\Xi}\c H\right)\\
&=&-(\tr X) \DDc \tr X  +\DDc \DDc\c\ov{\Xi}+\DDc \Xi\c\ov{H}+\DDc \ov{\Xi}\c H\\
&&+\Xi \DDc\c\ov{H}+\ov{\Xi} \DDc\c H
\eeaa
and
\beaa
&&\DDc\hot   \DDc\left(\nabc_4  \tr X\right) \\
&=& \DDc\hot   \big(-(\tr X) \DDc \tr X  +\DDc \DDc\c\ov{\Xi}+\DDc \Xi\c\ov{H}+\DDc \ov{\Xi}\c H\\
&&+\Xi \DDc\c\ov{H}+\ov{\Xi} \DDc\c H \big)\\
&=& -(\tr X)  \DDc\hot \DDc \tr X - \DDc(\tr X) \hot \DDc \tr X  + \DDc\hot \DDc \DDc\c\ov{\Xi}\\
&&+ \DDc\hot \big(\DDc \Xi\c\ov{H}+\DDc \ov{\Xi}\c H+\Xi \DDc\c\ov{H}+\ov{\Xi} \DDc\c H\big).
\eeaa

Also, using the above we compute
 \beaa
\DDc\hot \big( \, [\nabc_4 , \DDc] \tr X\big) &=& -\frac{1}{2}\tr X \DDc \hot \DDc  \tr X -\frac{1}{2}\DDc \tr X\hot \DDc  \tr X\\
&&  -\tr X \DDc\hot B -\DDc\tr X \hot B\\
  &&-\frac{1}{2}(\ov{\tr X}-\tr X)\DDc\hot(\Xh\c\ov{H})  +\frac{1}{2}\tr X \DDc\hot(\widehat{X}\c\ov{\Hb})\\
    &&-\frac{1}{2}(\DDc \ov{\tr X}-\DDc\tr X)\hot (\Xh\c\ov{H})  +\frac{1}{2}\DDc \tr X\hot (\widehat{X}\c\ov{\Hb}).
 \eeaa
 Using that $\DDc(\tr X) = -2\tr X\Hb+O(\ep)$, $\DDc(\ov{\tr X})=(\tr X-\ov{\tr X})H+O(\ep)$ and $\DDc\hot \Hb=-\Hb \hot \Hb+O(\ep)$, we obtain
  \beaa
\DDc\hot \big( \, [\nabc_4 , \DDc] \tr X\big) &=& -\frac{1}{2}\tr X \DDc \hot \DDc  \tr X -\frac{1}{2}\DDc \tr X\hot \DDc  \tr X\\
&&  -\tr X \DDc\hot B +2\tr X \Hb \hot B\\
&&-\frac{1}{2}(\ov{\tr X}-\tr X)\DDc\hot(\Xh\c\ov{H})  +\frac{1}{2}\tr X \DDc\hot(\widehat{X}\c\ov{\Hb})\\
    &&-\frac{1}{2}((\tr X-\ov{\tr X})H+2\tr X\Hb)\hot (\Xh\c\ov{H}) -\tr X\Hb\hot (\widehat{X}\c\ov{\Hb}).
 \eeaa
 We finally put all together to obtain
   \beaa
&&\nabc_4\AA_4\\
&=& -2(\tr X)  \DDc\hot \DDc \tr X - \frac 3 2\big( \DDc(\tr X)-2 \nabc_4 \Hb\big)  \hot \DDc \tr X-6\tr X\underline{H} \hot \DDc  \tr X\\
&&-\tr X \DDc\hot B +2\tr X \Hb \hot B+\MM[\Xi] \\
  &&-\frac{1}{2}(\ov{\tr X}-\tr X)\DDc\hot(\Xh\c\ov{H})  +\frac{1}{2}\tr X \DDc\hot(\widehat{X}\c\ov{\Hb})\\
    &&-\frac{1}{2}((\tr X-\ov{\tr X})H+2\tr X\Hb)\hot (\Xh\c\ov{H}) -\tr X\Hb\hot (\widehat{X}\c\ov{\Hb})\\
    &&+ 4\underline{H} \hot ( \tr X\Xh\c\ov{\Hb} +\frac{1}{2}\tr X\Xh\c\ov{H})+\Xh \c \ov{\DDc} (\tr X\Hb) -2 \tr X (  \Hb \c \ov{\Hb}) \Xh  \\
 &=& -2(\tr X)  \DDc\hot \DDc \tr X - \frac 3 2\big( \DDc(\tr X)-2 \nabc_4 \Hb\big)  \hot \DDc \tr X-6\tr X\underline{H} \hot \DDc  \tr X\\
&&-\tr X \DDc\hot B +2\tr X \Hb \hot B+\MM[\Xi] \\
  &&-\frac{1}{2}(\ov{\tr X}-\tr X)\DDc\hot(\Xh\c\ov{H})  +\frac{1}{2}\tr X \DDc\hot(\widehat{X}\c\ov{\Hb})\\
    &&-\frac{1}{2}(\tr X-\ov{\tr X})(H \c \ov{H})\Xh  + \tr X(  \Hb \c \ov{\Hb}) \Xh+ \tr X\underline{H} \hot (\Xh\c\ov{H})+\Xh \c \ov{\DDc} (\tr X\Hb) .
\eeaa

 Finally writing $\AA_3= 2\nabc_4 \Hb - \DDc \tr X$ and $\AA_4= \DDc\hot \DDc\left(\tr X\right)  +3\Hb\hot \DDc\left(\tr X\right)$, we obtain
    \beaa
&&\nabc_4\AA_4  \\
&=& -2(\tr X)\AA_4  - 3\tr X \Hb \hot   \AA_3 -\tr X \DDc\hot B +2\tr X \Hb \hot B+\MM[\Xi] \\
 &&-\frac{1}{2}(\ov{\tr X}-\tr X)\DDc\hot(\Xh\c\ov{H})  +\frac{1}{2}\tr X \DDc\hot(\widehat{X}\c\ov{\Hb})\\
    &&-\frac{1}{2}(\tr X-\ov{\tr X})(H \c \ov{H})\Xh  + \tr X (  \Hb \c \ov{\Hb}) \Xh+ \tr X\underline{H} \hot (\Xh\c\ov{H})+\Xh \c \ov{\DDc} (\tr X\Hb),
\eeaa
as stated.
\end{proof}

\begin{lemma}\label{lemma:derivatives-nabc4-BB} We have, modulo quadratic terms,
\beaa
\nabc_4 \BB&=& -\frac 5 2 \tr X\BB +2\tr X\Hb \c\ov{B} + 3P \,\DDc\c \ov{\Xi}\\
&&+ \frac{1}{2}\big( \DDc + 3 \Hb\big) \c \big( \DDc\c\ov{A}+\ov{A}\c\Hb\big),\\
\nabc_4 \DDc\BB&=&  -3 \tr X \DDc \BB + 5 \tr X \Hb \BB -8 \tr X \Hb( \Hb \c\ov{B})   \\
&&+ 3P \, \DDc\DDc\c \ov{\Xi}-6\Hb P \,\DDc\c \ov{\Xi} \\
&&+ \frac{1}{2} (\DDc +\Hb)\Big( \big( \DDc + 3 \Hb\big) \c \big( \DDc\c\ov{A}+\ov{A}\c\Hb\big)\Big),
\eeaa
and 
 \beaa
&&\nabc_4\DDc\hot\DDc\BB\\
&=& -\frac 7 2  \tr X \DDc\hot\DDc \BB  +8\tr X\Hb\hot \DDc \BB -10\tr X\Hb \hot \Hb \BB\\
&& (\DDc +\Hb)\hot \Big[  -8 \tr X \Hb( \Hb \c\ov{B}) + 3P \, \DDc\DDc\c \ov{\Xi}-6\Hb P \,\DDc\c \ov{\Xi}  \Big] \\
&&+ \frac{1}{2}(\DDc+\Hb)\hot \Big(  (\DDc +\Hb)\big( \big( \DDc + 3 \Hb\big) \c \big( \DDc\c\ov{A}+\ov{A}\c\Hb\big)\big)\Big).
\eeaa
\end{lemma}
\begin{proof} Recalling the definition of $\BB$ we obtain
\beaa
\nabc_4 \BB&=& \nabc_4 ( \DDc\c \ov{B}+ 2\Hb \c\ov{B})\\
&=&\DDc\c ( \nabc_4  \ov{B})+[\nabc_4 , \DDc\c ]\ov{B}+ 2\Hb \c \nabc_4\ov{B}+ 2\nabc_4\Hb \c\ov{B}.
\eeaa
From Lemma \ref{comm:1form-div}, \eqref{eq:comm-nabc4nabc3-ovDDc-F-err},  applied to $F=B$ and $s=1$, we have
   \beaa
[ \nabc_4 ,\DDc \c] \ov{B}     &=&- \frac 1 2\tr X\, ( \DDc \c \ov{B} - 2 \Hb \c \ov{F} )+\Hb \c \nabc_4 \ov{B},
 \eeaa
and therefore 
\beaa
\nabc_4 \BB&=&\DDc\c ( \nabc_4  \ov{B})- \frac 1 2\tr X\, ( \DDc \c \ov{B} - 2 \Hb \c \ov{B} )+ 3\Hb \c \nabc_4\ov{B}+ 2\nabc_4\Hb \c\ov{B}.
\eeaa
Using the Bianchi identity 
\beaa
\nabc_4\ov{B}  &=&-2\tr X \ov{B} + 3P \,\ov{\Xi}+ \frac{1}{2}\big( \DDc\c\ov{A}+\ov{A}\c\Hb\big).
\eeaa
we obtain
\beaa
&&\nabc_4 \BB\\
&=&\DDc\c \Big( -2\tr X \ov{B} + 3P \,\ov{\Xi}+ \frac{1}{2}\big( \DDc\c\ov{A}+\ov{A}\c\Hb\big)\Big)\\
&&- \frac 1 2\tr X\, ( \DDc \c \ov{B} - 2 \Hb \c \ov{B} )+ 3\Hb \c \Big(-2\tr X \ov{B} + 3P \,\ov{\Xi}+ \frac{1}{2}\big( \DDc\c\ov{A}+\ov{A}\c\Hb\big)\Big)\\
&&+ 2\nabc_4\Hb \c\ov{B}\\
&=& -2\tr X\DDc\c  \ov{B} -2\DDc\tr X \c\ov{B} + 3P \,\DDc\c \ov{\Xi}+ 3\DDc P \c \,\ov{\Xi}+ \frac{1}{2}\DDc\c \big( \DDc\c\ov{A}+\ov{A}\c\Hb\big)\\
&&- \frac 1 2\tr X\, ( \DDc \c \ov{B} - 2 \Hb \c \ov{B} )+ 3\Hb \c\left(-2\tr X \ov{B} + 3P \,\ov{\Xi}+ \frac{1}{2}\big( \DDc\c\ov{A}+\ov{A}\c\Hb\big)\right)\\
&&+ 2\nabc_4\Hb \c\ov{B}.
\eeaa
Using that $\DDc(\tr X) = -2\tr X\Hb+O(\ep)$, $\nabc_4(\Hb)=-\tr X\Hb+O(\ep)$ and $\DDc P=-3\Hb P+O(\ep)$, we obtain
\beaa
\nabc_4 \BB&=& -\frac 5 2 \tr X\DDc\c  \ov{B} -3\tr X\Hb \c\ov{B} + 3P \,\DDc\c \ov{\Xi}\\
&&+ \frac{1}{2}\big( \DDc + 3 \Hb\big) \c \big( \DDc\c\ov{A}+\ov{A}\c\Hb\big).
\eeaa
Writing $ \DDc\c \ov{B}=\BB_1-2\Hb \c\ov{B}$ we obtain the first identity.

We now compute $\nabc_4 \DDc \BB$. 
Applying Lemma \ref{comm:scalar}, \eqref{eq:comm-nabc4-DDc-h-err}, to $h=\BB$ and $s=1$ we have, up to quadratic terms
\beaa
 \, [\nabc_4 , \DDc]\BB  &=& -\frac{1}{2}\tr X\DDc \BB+\Hb\nabc_4 \BB+  \frac{1}{2}\tr X\Hb  \BB,
 \eeaa
and therefore we obtain, using the previous obtained identity,
 \beaa
\nabc_4 \DDc\BB&=& \DDc (\nabc_4 \BB)  -\frac{1}{2}\tr X\DDc \BB+\Hb\nabc_4 \BB+  \frac{1}{2}\tr X\Hb  \BB\\
&=& \DDc \left( -\frac 5 2 \tr X\BB +2\tr X\Hb \c\ov{B} + 3P \,\DDc\c \ov{\Xi}\right)  -\frac{1}{2}\tr X\DDc \BB+  \frac{1}{2}\tr X\Hb  \BB\\
&&+\Hb\left( -\frac 5 2 \tr X\BB +2\tr X\Hb \c\ov{B} + 3P \,\DDc\c \ov{\Xi}\right)\\
&&+ \frac{1}{2} (\DDc +\Hb)\Big( \big( \DDc + 3 \Hb\big) \c \big( \DDc\c\ov{A}+\ov{A}\c\Hb\big)\Big).
\eeaa
Using Lemma \ref{SIMPLIFICATION-ANGULAR} we simplify the above to
 \beaa
\nabc_4 \DDc\BB&=&  -\frac 5 2 \tr X \DDc \BB -\frac 5 2 \DDc \tr X \BB +2\DDc \tr X\Hb \c\ov{B} +2\tr X(\DDc \Hb)\c  \ov{B} \\
&&+2\tr X\Hb ( \DDc \c \ov{B}) + 3 \DDc P \,\DDc\c \ov{\Xi}+ 3P \, \DDc\DDc\c \ov{\Xi}  \\
&&-\frac{1}{2}\tr X\DDc \BB+  \frac{1}{2}\tr X\Hb  \BB+\Hb( -\frac 5 2 \tr X\BB +2\tr X\Hb \c\ov{B} + 3P \,\DDc\c \ov{\Xi})\\
&&+ \frac{1}{2} (\DDc +\Hb)\Big( \big( \DDc + 3 \Hb\big) \c \big( \DDc\c\ov{A}+\ov{A}\c\Hb\big)\Big)\\
&=&  -3 \tr X \DDc \BB + 3 \tr X \Hb \BB -2 \tr X \Hb( \Hb \c\ov{B})\\
&& +2\tr X(\DDc \Hb)\c  \ov{B} +2\tr X\Hb ( \DDc \c \ov{B}) + 3P \, \DDc\DDc\c \ov{\Xi}-6\Hb P \,\DDc\c \ov{\Xi} \\
&&+ \frac{1}{2} (\DDc +\Hb)\Big( \big( \DDc + 3 \Hb\big) \c \big( \DDc\c\ov{A}+\ov{A}\c\Hb\big)\Big),
\eeaa
where we used that $(\DDc \Hb)\c  \ov{B}=-\Hb (\Hb\c  \ov{B}) +O(\ep)$.
By writing $ \DDc\c \ov{B}=\BB-2\Hb \c\ov{B}$, we obtain the second identity.

We now compute $\nabc_4\DDc\hot\DDc\BB_1$. Applying Lemma \ref{comm:1form}, \eqref{eq:comm-nabc4nabc3-DDchot-err}, to $F=\DDc \BB$ and $s=1$ we have, up to quadratic terms,
\beaa
 [\nabc_4 , \DDc \hot ]\DDc \BB &=&- \frac 1 2 \tr X \DDc\hot \DDc \BB + \underline{H} \hot \nabc_4 \DDc \BB,
 \eeaa
and therefore we obtain, using the previous obtained identity,
 \beaa
&&\nabc_4\DDc\hot\DDc\BB\\
&=& \DDc\hot(\nabc_4\DDc\BB)- \frac 1 2 \tr X \DDc\hot \DDc\BB + \underline{H} \hot \nabc_4 \DDc\BB\\
&=& \DDc\hot \Big[ -3 \tr X \DDc \BB + 5 \tr X \Hb \BB -8 \tr X \Hb( \Hb \c\ov{B}) \\
&&+ 3P \, \DDc\DDc\c \ov{\Xi}-6\Hb P \,\DDc\c \ov{\Xi} \\
&&+ \frac{1}{2} (\DDc +\Hb)\big( \big( \DDc + 3 \Hb\big) \c \big( \DDc\c\ov{A}+\ov{A}\c\Hb\big)\big) \Big]- \frac 1 2 \tr X \DDc\hot \DDc\BB_1 \\
&&+ \underline{H} \hot \Big[ -3 \tr X \DDc \BB + 5 \tr X \Hb \BB -8 \tr X \Hb( \Hb \c\ov{B})  \\
&&+ 3P \, \DDc\DDc\c \ov{\Xi}-6\Hb P \,\DDc\c \ov{\Xi} \\
&&+ \frac{1}{2} (\DDc +\Hb)\big( \big( \DDc + 3 \Hb\big) \c \big( \DDc\c\ov{A}+\ov{A}\c\Hb\big)\big) \Big].
\eeaa
By denoting 
\beaa
\mbox{Expr}_1(A)&:=& \frac{1}{2}(\DDc+\Hb)\hot \Big(  (\DDc +\Hb)\big( \big( \DDc + 3 \Hb\big) \c \big( \DDc\c\ov{A}+\ov{A}\c\Hb\big)\big)\Big)
\eeaa
we obtain
 \beaa
&&\nabc_4\DDc\hot\DDc\BB_1\\
&=& -3 \tr X \DDc\hot\DDc \BB  -3 \DDc \tr X\hot \DDc \BB \\
&&+ 5\DDc \tr X \hot \Hb \BB+ 5 \tr X \DDc\hot \Hb \BB_1+ 5 \tr X \Hb \hot \DDc\BB\\
&& \DDc\hot \Big[  -8 \tr X \Hb( \Hb \c\ov{B}) + 3P \, \DDc\DDc\c \ov{\Xi}-6\Hb P \,\DDc\c \ov{\Xi}  \Big]\\
&&- \frac 1 2 \tr X \DDc\hot \DDc\BB + \underline{H} \hot \Big[ -3 \tr X \DDc \BB + 5 \tr X \Hb \BB -8 \tr X \Hb( \Hb \c\ov{B})  \\
&&+ 3P \, \DDc\DDc\c \ov{\Xi}-6\Hb P \,\DDc\c \ov{\Xi}  \Big]+\mbox{Expr}_1(A)\\
&=& -\frac 7 2  \tr X \DDc\hot\DDc \BB  +8\tr X\Hb\hot \DDc \BB -10\tr X\Hb \hot \Hb \BB\\
&& (\DDc +\Hb)\hot \Big[  -8 \tr X \Hb( \Hb \c\ov{B}) + 3P \, \DDc\DDc\c \ov{\Xi}-6\Hb P \,\DDc\c \ov{\Xi}  \Big]\\
&& +\mbox{Expr}_1(A)
\eeaa
where we used $\DDc\hot \Hb= - \Hb \hot \Hb+O(\ep)$. This proves the lemma.
\end{proof}


\subsection{The derivation of the Teukolsky-Starobinski identity}


The Teukolsky-Starobinski identity is obtained by taking three $\nabc_4$ derivatives of appropriately rescalings of the Bianchi identity for $\Ab$, i.e.
\bea
\nabc_4\Ab +\frac{1}{2}\tr X \Ab&=&-\DDc\hot\Bb   - 4 \Hb\hot \Bb -3P\Xbh. \label{Bianchi1}
\eea


\subsubsection{The second $e_4$ derivative of $\Ab$}


\begin{lemma}\label{lemma:first-derivative} The quantity $\mathfrak{F} := -\nabc_4\Ab -\frac{1}{2}\tr X\Ab \in \sk_2(\CCC)$ satisfies, modulo quadratic terms,
\beaa
\nabc_4 \mathfrak{F}+\frac 3 2 \tr X  \mathfrak{F}&=& -\DDc\hot \mathcal{A}_1- 5\underline{H} \hot \mathcal{A}_1 -2 \tr X\Hb\hot  \Bb\\
&& -\frac32 P\Big(\tr X\, \widehat{\Xb} +\ov{\tr\Xb} \widehat{X}-2\mathcal{A}_2 \Big).
\eeaa
\end{lemma} 
\begin{proof}
We infer from \eqref{Bianchi1} that
\beaa
\mathfrak{F} &=&\DDc\hot\Bb   + 4 \Hb\hot \Bb +3P\Xbh.
\eeaa
We compute
\beaa
\nabc_4 \mathfrak{F}&=& \nabc_4 \DDc\hot\Bb   + 4 \Hb\hot \nabc_4 \Bb  + 4 \nabc_4\Hb\hot \Bb +3P\nabc_4\Xbh+3\nabc_4P\Xbh.
\eeaa
Applying Lemma \ref{comm:1form}, \eqref{eq:comm-nabc4nabc3-DDchot-err}, to $F=\Bb$ and $s=-1$,
we have, modulo quadratic terms,
\beaa
\nabc_4 \mathfrak{F}&=&  \DDc\hot \nabc_4\Bb- \frac 1 2 \tr X\left( \DDc\hot \Bb + 2\Hb\hot \Bb\right)+ \underline{H} \hot \nabc_4 \Bb \\
&&  + 4 \Hb\hot \nabc_4 \Bb  + 4 \nabc_4\Hb\hot \Bb +3P\nabc_4\Xbh+3\nabc_4P\Xbh\\
&=&  \DDc\hot \nabc_4\Bb+ 5\underline{H} \hot \nabc_4 \Bb- \frac 1 2 \tr X\left( \DDc\hot \Bb + 2\Hb\hot \Bb\right)  + 4 \nabc_4\Hb\hot \Bb \\
&&+3P\nabc_4\Xbh+3\nabc_4P\Xbh.
\eeaa
Using the Bianchi identity
\bea
\nabc_4\Bb &=& -\tr X\Bb- (\DDc P+3P\Hb)= -\tr X\Bb- \mathcal{A}_1
\eea
we obtain
\beaa
\nabc_4 \mathfrak{F}&=&  \DDc\hot \big( -\tr X\Bb- \mathcal{A}_1 \big)+ 5\underline{H} \hot \big( -\tr X\Bb- \mathcal{A}_1 \big)\\
&&- \frac 1 2 \tr X\left( \DDc\hot \Bb + 2\Hb\hot \Bb\right)  + 4 \nabc_4\Hb\hot \Bb +3P\nabc_4\Xbh+3\nabc_4P\Xbh\\
&=& -\tr X \DDc\hot  \Bb - \DDc\tr X \hot  \Bb -\DDc\hot \mathcal{A}_1+ 5\underline{H} \hot \big( -\tr X\Bb- \mathcal{A}_1\big)\\
&&- \frac 1 2 \tr X\left( \DDc\hot \Bb + 2\Hb\hot \Bb\right)  + 4 \nabc_4\Hb\hot \Bb +3P\nabc_4\Xbh+3\nabc_4P\Xbh\\
&=& -\frac 3 2 \tr X \DDc\hot  \Bb - 6 \tr X\Hb\hot \Bb -\DDc\hot \mathcal{A}_1- 5\underline{H} \hot \mathcal{A}_1\\
&&+\big( 4 \nabc_4\Hb - \DDc\tr X\big) \hot  \Bb +3P\nabc_4\Xbh+3\nabc_4P\Xbh.
\eeaa
By writing $\DDc\hot\Bb   + 4 \Hb\hot \Bb = \mathfrak{F}-3P\Xbh$
we obtain
\beaa
\nabc_4 \mathfrak{F}&=& -\frac 3 2 \tr X\big(  \mathfrak{F}-3P\Xbh\big) -\DDc\hot \mathcal{A}_1- 5\underline{H} \hot\mathcal{A}_1\\
&&+\big( 4 \nabc_4\Hb - \DDc\tr X\big) \hot  \Bb +3P\nabc_4\Xbh+3\big(-\frac 3 2 \tr XP \big)\Xbh,
\eeaa
where we used $\nabc_4 P=-\frac 3 2 \tr XP+O(\ep)$,
and therefore 
\beaa
\nabc_4 \mathfrak{F}+\frac 3 2 \tr X  \mathfrak{F}&=& -\DDc\hot \mathcal{A}_1- 5\underline{H} \hot \mathcal{A}_1\\
&&+\big( 4 \nabc_4\Hb - \DDc\tr X\big) \hot  \Bb +3P\nabc_4\Xbh.
\eeaa

Using the Ricci identity
\beaa
\nabc_4\widehat{\Xb} +\frac{1}{2}\tr X\, \widehat{\Xb} &=& \DDc\hot\Hb  +\Hb\hot\Hb -\frac{1}{2}\ov{\tr\Xb} \widehat{X}=\mathcal{A}_2 -\frac{1}{2}\ov{\tr\Xb} \widehat{X},
\eeaa
we finally have
\beaa
\nabc_4 \mathfrak{F}+\frac 3 2 \tr X  \mathfrak{F}&=& -\DDc\hot \mathcal{A}_1- 5\underline{H} \hot \mathcal{A}_1+3P\mathcal{A}_2\\
&&+\big( 4 \nabc_4\Hb - \DDc\tr X\big) \hot  \Bb -\frac32 P\big(\tr X\, \widehat{\Xb} +\ov{\tr\Xb} \widehat{X} \big).
\eeaa
Using that 
\beaa
 4 \nabc_4\Hb - \DDc\tr X&=&  4 (-\tr X\Hb) - (-2\tr X\Hb)+O(\ep)=  -2 \tr X\Hb+O(\ep),
\eeaa
we obtain the stated.
\end{proof}


\subsubsection{The third $e_4$ derivative of $\Ab$}


\begin{lemma}\label{lemma:third-derivative} The quantity $\mathfrak{G} := \nabc_4\mathfrak{F}+\frac{3}{2}\tr X\mathfrak{F} \in \sk_2(\CCC)$ satisfies, modulo quadratic terms,
 \beaa
\nabc_4 \mathfrak{G}+ \frac 5 2 \tr X \mathfrak{G} &=&-\frac{1}{2}\DDc\hot\DDc\BB -5 \Hb \hot \DDc\BB- 10(\underline{H} \hot \Hb)\BB+ 3 \tr X\Hb \hot\AA_1\\
&&+\frac 3 2 P \ov{\tr\Xb}A +\frac{3}{2} P \Big[\AA_4-2B \hot \Hb -\DDc\hot\big(\ov{H} \c \Xh\big)\\
&&+\big(\ov{\tr X\tr\Xb}-2\ov{ \DDc}\c\Hb-2P\big)\Xh-3   \Hb \hot\big(\ov{H} \c \Xh\big) \Big] .
\eeaa
\end{lemma}
\begin{proof} We infer from Lemma \ref{lemma:first-derivative} that 
\beaa
\mathfrak{G}= -\DDc\hot \mathcal{A}_1- 5\underline{H} \hot \mathcal{A}_1  -2 \tr X\Hb\hot  \Bb -\frac32 P\big(\tr X\, \widehat{\Xb} +\ov{\tr\Xb} \widehat{X}-2\mathcal{A}_2 \big).
\eeaa
We compute 
\beaa
\nabc_4 \mathfrak{G}&=& - \nabc_4\DDc\hot \mathcal{A}_1- 5\underline{H} \hot \nabc_4\mathcal{A}_1- 5\nabc_4\underline{H} \hot \mathcal{A}_1\\
&& -2 \tr X\Hb\hot  \nabc_4\Bb -2 \tr X\nabc_4\Hb\hot  \Bb -2 \nabc_4\tr X\Hb\hot  \Bb\\
&& -\frac32 P\nabc_4\big(\tr X\, \widehat{\Xb} +\ov{\tr\Xb} \widehat{X}-2\mathcal{A}_2 \big)  -\frac32 \nabc_4P\big(\tr X\, \widehat{\Xb} +\ov{\tr\Xb} \widehat{X}-2\mathcal{A}_2 \big).
\eeaa
Using that $\nabc_4 P=-\frac 3 2 \tr X P+O(\ep)$, $\nabc_4 \tr X=-\frac 1 2 (\tr X)^2+O(\ep)$, $\nabc_4 \Hb=-\tr X \Hb+O(\ep)$ and $\nabc_4\Bb = -\tr X\Bb- \mathcal{A}_1$, we obtain
\beaa
\nabc_4 \mathfrak{G}&=& - \nabc_4\DDc\hot \mathcal{A}_1- 5\underline{H} \hot \nabc_4\mathcal{A}_1+7 \tr X\Hb  \hot  \mathcal{A}_1+ 5(\tr X)^2\Hb\hot  \Bb \\
&& -\frac32 P\nabc_4\big(\tr X\, \widehat{\Xb} +\ov{\tr\Xb} \widehat{X}-2\mathcal{A}_2 \big)  +\frac94 \tr X P\big(\tr X\, \widehat{\Xb} +\ov{\tr\Xb} \widehat{X}-2\mathcal{A}_2 \big).
\eeaa
Using Lemma \ref{lemma:derivatives-AA1} for the derivatives of $\AA_1$, we obtain
\beaa
&&\nabc_4 \mathfrak{G}\\
&=& - \Big[ -\frac 5 2 \tr X \DDc\hot \AA_1 + \frac 3 2 \tr X\Hb \hot\AA_1 +\frac{1}{2}\DDc\hot\DDc\BB +\frac 5 2 \Hb \hot \DDc\BB\\
&&+\DDc\hot\Big(\frac{3}{2} P\left(\ov{H} \c \Xh-\ov{\tr\Xb} \Xi+\AA_3\right)  \Big)+ \underline{H} \hot \Big(\frac{3}{2} P \left(\ov{H} \c \Xh-\ov{\tr\Xb} \Xi+\AA_3\right)  \Big)\Big]\\
&&- 5\underline{H} \hot \Big[ -2\tr X\AA_1 +\frac{1}{2}\DDc\BB +2\Hb\BB+\frac{3}{2} P \left(\ov{H} \c \Xh-\ov{\tr\Xb} \Xi+\AA_3\right)  \Big]\\
&&+7 \tr X\Hb  \hot  \mathcal{A}_1 + 5(\tr X)^2\Hb\hot  \Bb -\frac32 P\nabc_4\big(\tr X\, \widehat{\Xb} +\ov{\tr\Xb} \widehat{X}-2\mathcal{A}_2 \big) \\
&& +\frac94 \tr X P\big(\tr X\, \widehat{\Xb} +\ov{\tr\Xb} \widehat{X}-2\mathcal{A}_2 \big) \\
&=& \frac 5 2 \tr X \DDc\hot \AA_1 + \frac {31}{ 2} \tr X\Hb \hot\AA_1+ 5(\tr X)^2\Hb\hot  \Bb \\
&&-\frac{1}{2}\DDc\hot\DDc\BB -5 \Hb \hot \DDc\BB- 10(\underline{H} \hot \Hb)\BB \\
&&-\DDc\hot\Big(\frac{3}{2} P \left(\ov{H} \c \Xh-\ov{\tr\Xb} \Xi+\AA_3\right)  \Big)- 6\underline{H} \hot \Big(\frac{3}{2} P \left(\ov{H} \c \Xh-\ov{\tr\Xb} \Xi+\AA_3\right)  \Big)\\
&& -\frac32 P \Big[\nabc_4\big(\tr X\, \widehat{\Xb} +\ov{\tr\Xb} \widehat{X}\big)-2\nabc_4\mathcal{A}_2  - \frac 3 2 \tr X \big(\tr X\, \widehat{\Xb} +\ov{\tr\Xb} \widehat{X}-2\mathcal{A}_2 \big)\Big] .
\eeaa
Writing 
\beaa
\DDc\hot \mathcal{A}_1+2 \tr X\Hb\hot  \Bb=-\mathfrak{G} - 5\underline{H} \hot \mathcal{A}_1   -\frac32 P\big(\tr X\, \widehat{\Xb} +\ov{\tr\Xb} \widehat{X}-2\mathcal{A}_2 \big),
\eeaa
we obtain
\beaa
&&\nabc_4 \mathfrak{G}+ \frac 5 2 \tr X \mathfrak{G} \\
&=&-\frac{1}{2}\DDc\hot\DDc\BB -5 \Hb \hot \DDc\BB- 10(\underline{H} \hot \Hb)\BB+ 3 \tr X\Hb \hot\AA_1 \\
&&-\DDc\hot\Big(\frac{3}{2} P \left(\ov{H} \c \Xh-\ov{\tr\Xb} \Xi+\AA_3\right)  \Big)- 6\underline{H} \hot \Big(\frac{3}{2} P \left(\ov{H} \c \Xh-\ov{\tr\Xb} \Xi+\AA_3\right)  \Big)\\
&& -\frac32 P \Big[\nabc_4\big(\tr X\, \widehat{\Xb} +\ov{\tr\Xb} \widehat{X}\big)-2\nabc_4\mathcal{A}_2  + \tr X \big(\tr X\, \widehat{\Xb} +\ov{\tr\Xb} \widehat{X}-2\mathcal{A}_2 \big)\Big] .
\eeaa
We now write for the second line, using that $\DDc P=-3  P \Hb+O(\ep)$,
\beaa
&&-\DDc\hot\Big(\frac{3}{2} P \left(\ov{H} \c \Xh-\ov{\tr\Xb} \Xi+\AA_3\right)  \Big)- 6\underline{H} \hot \Big(\frac{3}{2} P \left(\ov{H} \c \Xh-\ov{\tr\Xb} \Xi+\AA_3\right)  \Big)\\
&=&-\frac{3}{2} P \DDc\hot\big(\ov{H} \c \Xh-\ov{\tr\Xb} \Xi +\AA_3\big)-\frac{3}{2} (-3  P \Hb) \hot\big(\ov{H} \c \Xh-\ov{\tr\Xb} \Xi +\AA_3\big)\\
&&- 6\underline{H} \hot \Big(\frac{3}{2} P \left(\ov{H} \c \Xh-\ov{\tr\Xb} \Xi+\AA_3\right)  \Big)\\
&=&-\frac{3}{2} P \Big[\DDc\hot\big(\ov{H} \c \Xh-\ov{\tr\Xb} \Xi +\AA_3\big) +3   \Hb \hot\big(\ov{H} \c \Xh-\ov{\tr\Xb} \Xi +\AA_3\big)\Big]\\
&=&-\frac{3}{2} P \Big[\DDc \hot \AA_3 + 3 \Hb \hot \AA_3  -\ov{\tr\Xb} \DDc \hot \Xi -\DDc \ov{\tr\Xb} \hot \Xi  -3\ov{\tr\Xb}    \Hb \hot\Xi\\
&&+\DDc\hot\big(\ov{H} \c \Xh\big)+3   \Hb \hot\big(\ov{H} \c \Xh\big)\Big]\\
&=&-\frac{3}{2} P \Big[\DDc \hot \AA_3 + 3 \Hb \hot \AA_3  -\ov{\tr\Xb} \DDc \hot \Xi -(\tr\Xb+2\ov{\tr \Xb})\Hb \hot \Xi  \\
&&+\DDc\hot\big(\ov{H} \c \Xh\big)+3   \Hb \hot\big(\ov{H} \c \Xh\big)\Big]
\eeaa
where we used $\DDc\ov{\tr\Xb} =(\tr\Xb-\ov{\tr \Xb})\Hb+O(\ep)$.

Also, we have
\beaa
\nabc_4\left(\tr X\, \widehat{\Xb} +\ov{\tr\Xb} \widehat{X}\right) &=& -\frac{1}{2}(\tr X)^2 \Xbh+\ov{\nabc_4\tr\Xb}\Xh\\
&&+\tr X\left(-\frac{1}{2}\tr X\, \widehat{\Xb} +\mathcal{A}_2 -\frac{1}{2}\ov{\tr\Xb} \widehat{X}\right)\\
&&+\ov{\tr\Xb}\left(-\frac 1 2 (\tr X+\ov{\tr X})\Xh + \DDc\hot \Xi+  \Xi\hot(\Hb+H)-A\right)\\
 &=& -(\tr X)^2 \Xbh+\left(\ov{\nabc_4\tr\Xb} -\left(\tr X+\frac{1}{2}\ov{\tr X}\right)\ov{\tr\Xb}\right)\Xh\\
&& -\ov{\tr\Xb}A+\tr X \mathcal{A}_2+\ov{\tr\Xb}\left(\DDc\hot \Xi+  \Xi\hot(\Hb+H)\right).
\eeaa

By putting all together we obtain
\beaa
&&\nabc_4 \mathfrak{G}+ \frac 5 2 \tr X \mathfrak{G} \\
&=&-\frac{1}{2}\DDc\hot\DDc\BB_1 -5 \Hb \hot \DDc\BB_1- 10(\underline{H} \hot \Hb)\BB_1+ 3 \tr X\Hb \hot\AA_1+\frac 3 2 P \ov{\tr\Xb}A \\
&&-\frac{3}{2} P \Big[\DDc \hot \AA_3 + 3 \Hb \hot \AA_3 -2\nabc_4\mathcal{A}_2-\tr X\mathcal{A}_2\\
&&  -(\tr\Xb+\ov{\tr \Xb})\Hb \hot \Xi +\ov{\tr\Xb}  H \hot\Xi \\
&&+\DDc\hot\big(\ov{H} \c \Xh\big)+\big(\ov{\nabc_4\tr\Xb} -\frac{1}{2}\ov{\tr X}\ov{\tr\Xb}\big)\Xh+3   \Hb \hot\big(\ov{H} \c \Xh\big) \Big] .
\eeaa
Finally, using Lemma \ref{lemma:DDcAA3},
we obtain
\beaa
&&\nabc_4 \mathfrak{G}+ \frac 5 2 \tr X \mathfrak{G} \\
&=&-\frac{1}{2}\DDc\hot\DDc\BB_1 -5 \Hb \hot \DDc\BB_1- 10(\underline{H} \hot \Hb)\BB_1+ 3 \tr X\Hb \hot\AA_1+\frac 3 2 P \ov{\tr\Xb}A \\
&&-\frac{3}{2} P \Big[-\AA_4+2B \hot \Hb +\Xh \c \ov{\DDc} \Hb+\DDc\hot\big(\ov{H} \c \Xh\big)\\
&&+\big(\ov{\nabc_4\tr\Xb} -\frac{1}{2}\ov{\tr X}\ov{\tr\Xb}-\ov{\Hb}\c \Hb\big)\Xh+3   \Hb \hot\big(\ov{H} \c \Xh\big) \Big] .
\eeaa
Finally using \eqref{rule-1} to write $\Xh \c \ov{\DD} \Hb= \Xh (\ov{\DD} \c\Hb)$, and  that
\beaa
\nabc_4\ov{\tr\Xb}  &=&-\frac{1}{2}\ov{\tr X\tr\Xb}+\ov{ \DDc}\c\Hb+\Hb\c\ov{\Hb}+2P+O(\ep^2)
\eeaa
 we conclude the lemma.
\end{proof}


\subsubsection{The fourth $e_4$ derivative of $\Ab$}


Define 
\beaa
\mathfrak{H} &:=& \nabc_4\mathfrak{G} +\frac{5}{2}\tr X\mathfrak{G}.
\eeaa
We infer from Lemma \ref{lemma:third-derivative}, 
\beaa
\mathfrak{H}&=&\mathfrak{H}_1+\frac{3}{2} P     \mathfrak{H}_2+\frac 3 2 P \ov{\tr\Xb}A.
\eeaa
with 
\beaa
 \mathfrak{H}_1&=& -\frac{1}{2}\DDc\hot\DDc\BB -5 \Hb \hot \DDc\BB- 10(\underline{H} \hot \Hb)\BB+ 3 \tr X\Hb \hot\AA_1 \\
    \mathfrak{H}_2&=&\AA_4-2B \hot \Hb-\DDc\hot\big(\ov{H} \c \Xh\big)+\big(\ov{\tr X\tr\Xb}-2\ov{ \DDc}\c\Hb-2P\big)\Xh-3   \Hb \hot\big(\ov{H} \c \Xh\big).
\eeaa
This gives, using that $\nabc_4 P=-\frac 3 2 \tr X P$, 
\bea\label{nabc4mathfrakH}
\nabc_4\mathfrak{H} &=& \nabc_4\mathfrak{H}_1+\frac 3 2 P \big[\nabc_4 \mathfrak{H}_2-\frac 3 2 \tr X  \mathfrak{H}_2\big]+ \nabc_4 \Big( \frac 3 2 P \ov{\tr\Xb}A\Big).
\eea

Here we compute $\nabc_4 \mathfrak{H}_1$.

\begin{lemma}\label{lemma:nabc4mathfrak1} We have, modulo quadratic terms,
\beaa
 \nabc_4\mathfrak{H}_1+ \frac 7 2  \tr X\mathfrak{H}_1  &=&\frac{3}{2} P \Big[  - \DDc \hot( \DDc\DDc\c \ov{\Xi})-6\Hb  \hot \DDc  \DDc\c \ov{\Xi}\\
 &&-6\Hb\hot \Hb  \,\DDc\c \ov{\Xi}- 3 \tr X\ov{\tr\Xb}\Hb \hot  \Xi \\
 &&+ 3 \tr X\Hb \hot ( \ov{H} \c \Xh) + 3 \tr X\Hb \hot  \AA_3 \Big]+\mbox{Expr}_2(A)
\eeaa
where 
\beaa
\mbox{Expr}_2(A)&:=& - \frac{1}{4}(\DDc+11\Hb)\hot \Big(  (\DDc +\Hb)\big( \big( \DDc + 3 \Hb\big) \c \big( \DDc\c\ov{A}+\ov{A}\c\Hb\big)\big)\Big)\\
&&- 5 (\underline{H} \hot \Hb)\big( \DDc + 3 \Hb\big) \c \big( \DDc\c\ov{A}+\ov{A}\c\Hb\big).
\eeaa
\end{lemma}
\begin{proof} From the definition of $\mathfrak{H}_1$ we obtain
\beaa
 \nabc_4\mathfrak{H}_1&=& -\frac{1}{2}\nabc_4\DDc\hot\DDc\BB -5 \Hb \hot \nabc_4\DDc\BB -5 \nabc_4\Hb \hot \DDc\BB\\
 &&- 10(\underline{H} \hot \Hb)\nabc_4\BB- 20(\nabc_4\underline{H} \hot \Hb)\BB\\
 &&+ 3 \tr X\Hb \hot \nabc_4 \AA_1+ 3 \tr X\nabc_4 \Hb \hot\AA_1+3\nabc_4  \tr X\Hb \hot\AA_1\\
 &=& -\frac{1}{2}\nabc_4\DDc\hot\DDc\BB -5 \Hb \hot \nabc_4\DDc\BB +5 \tr X \Hb \hot \DDc\BB\\
 &&- 10(\underline{H} \hot \Hb)\nabc_4\BB+ 20 \tr X( \Hb \hot \Hb)\BB+ 3 \tr X\Hb \hot \nabc_4 \AA_1-\frac 9 2 (\tr X)^2\Hb \hot\AA_1
\eeaa
where we used  $\nabc_4 \Hb=-\tr X \Hb+O(\ep)$ and $\nabc_4 \tr X=-\frac 1 2 (\tr X)^2+O(\ep)$.

Using Lemma \ref{lemma:derivatives-nabc4-BB} and Lemma \ref{lemma:derivatives-AA1}, we obtain
\beaa
 \nabc_4\mathfrak{H}_1 &=& \frac 7 4  \tr X \DDc\hot\DDc \BB  -4\tr X\Hb\hot \DDc \BB +5\tr X\Hb \hot \Hb \BB\\
&&-\frac{1}{2} (\DDc +\Hb)\hot \Big[  -8 \tr X \Hb( \Hb \c\ov{B})+ 3P \, \DDc\DDc\c \ov{\Xi}-6\Hb P \,\DDc\c \ov{\Xi}  \Big] \\
&& -5 \Hb \hot \Big[  -3 \tr X \DDc \BB + 5 \tr X \Hb \BB -8 \tr X \Hb( \Hb \c\ov{B})  \\
&&+ 3P \, \DDc\DDc\c \ov{\Xi}-6\Hb P \,\DDc\c \ov{\Xi}\Big] +5 \tr X \Hb \hot \DDc\BB\\
 &&- 10(\underline{H} \hot \Hb)\Big[ -\frac 5 2 \tr X\BB +2\tr X\Hb \c\ov{B} + 3P \,\DDc\c \ov{\Xi} \Big]+ 20 \tr X( \Hb \hot \Hb)\BB\\
 &&+ 3 \tr X\Hb \hot \Big[  -2\tr X\AA_1 +\frac{1}{2}\DDc\BB +2\Hb\BB+\frac 3 2P\left(\ov{H} \c \Xh-\ov{\tr\Xb} \Xi+\AA_3\right)\Big]\\
 &&-\frac 9 2 (\tr X)^2\Hb \hot\AA_1+\mbox{Expr}_2(A) 
\eeaa
where 
\beaa
\mbox{Expr}_2(A)&:=& - \frac{1}{4}(\DDc+11\Hb)\hot \Big(  (\DDc +\Hb)\big( \big( \DDc + 3 \Hb\big) \c \big( \DDc\c\ov{A}+\ov{A}\c\Hb\big)\big)\Big)\\
&&- 5 (\underline{H} \hot \Hb)\big( \DDc + 3 \Hb\big) \c \big( \DDc\c\ov{A}+\ov{A}\c\Hb\big).
\eeaa
We simplify the above to
\beaa
 \nabc_4\mathfrak{H}_1 &=& \frac 7 4  \tr X \DDc\hot\DDc \BB  +\frac{35}{2}\tr X\Hb\hot \DDc \BB +31\tr X\Hb \hot \Hb \BB-\frac{21}{ 2} (\tr X)^2\Hb \hot\AA_1\\
&&-\frac{1}{2} (\DDc +11\Hb)\hot \Big[  -8 \tr X \Hb( \Hb \c\ov{B})+ 3P \, \DDc\DDc\c \ov{\Xi}-6\Hb P \,\DDc\c \ov{\Xi}  \Big]\\
&& - 10(\underline{H} \hot \Hb)\Big[2\tr X\Hb \c\ov{B} + 3P \,\DDc\c \ov{\Xi} \Big]\\
 &&+ \frac 3 2P 3 \tr X\Hb \hot \Big[ \ov{H} \c \Xh-\ov{\tr\Xb} \Xi+\AA_3\Big]+\mbox{Expr}_2(A).
\eeaa
By writing $\frac{1}{2}\DDc\hot\DDc\BB=-\mathfrak{H}_1  -5 \Hb \hot \DDc\BB- 10(\underline{H} \hot \Hb)\BB+ 3 \tr X\Hb \hot\AA_1 $, and re-organizing the above we obtain
\beaa
 \nabc_4\mathfrak{H}_1+ \frac 7 2  \tr X\mathfrak{H}_1 &=& -4\tr X\Hb \hot \Hb \BB\\
&&+4 (\DDc +11\Hb)\hot \Big[  \tr X \Hb( \Hb \c\ov{B}) \Big] - 20 \tr X(\underline{H} \hot \Hb)(\Hb \c\ov{B}) \\
&&-\frac{3}{2} (\DDc +11\Hb)\hot \Big[  P \, \DDc\DDc\c \ov{\Xi}-2\Hb P \,\DDc\c \ov{\Xi}  \Big]\\
&&- 30P(\underline{H} \hot \Hb) \,\DDc\c \ov{\Xi} \\
 &&+ \frac 3 2P 3 \tr X\Hb \hot \Big[ \ov{H} \c \Xh-\ov{\tr\Xb} \Xi+\AA_3\Big]+\mbox{Expr}_2(A).
\eeaa
We now simplify the second line. We obtain
\beaa
&&4 (\DDc +11\Hb)\hot \Big[  \tr X \Hb( \Hb \c\ov{B}) \Big] - 20 \tr X(\underline{H} \hot \Hb)(\Hb \c\ov{B}) \\
&=&4 \DDc\hot \Big[  \tr X \Hb( \Hb \c\ov{B}) \Big]+ 24 \tr X(\underline{H} \hot \Hb)(\Hb \c\ov{B}) \\
&=&4  \Big[ \DDc \tr X\hot  \Hb( \Hb \c\ov{B})+ \tr X \DDc \hot \Hb( \Hb \c\ov{B})+ \tr X \Hb \hot \DDc( \Hb \c\ov{B}) \Big]\\
&&+ 24 \tr X(\underline{H} \hot \Hb)(\Hb \c\ov{B}) \\
&=&4  \Big[  -2\tr X\Hb\hot  \Hb( \Hb \c\ov{B})-\tr X   \Hb \hot \Hb( \Hb \c\ov{B})+ \tr X \Hb \hot \DDc \Hb \c\ov{B}\\
&&+ \tr X \Hb \hot\Hb  (\DDc  \c\ov{B}) \Big]+ 24 \tr X(\underline{H} \hot \Hb)(\Hb \c\ov{B}) \\
&=&4  \tr X \Hb \hot\Hb  (\DDc  \c\ov{B})+ 8 \tr X(\underline{H} \hot \Hb)(\Hb \c\ov{B}) \\
&=&4 \tr X \Hb \hot \Hb \BB,
\eeaa
where we used $\DDc(\tr X) = -2\tr X\Hb+O(\ep)$, $\DDc \hot \Hb= - \Hb \hot \Hb+O(\ep)$. 

We finally obtain
\beaa
&& \nabc_4\mathfrak{H}_1+ \frac 7 2  \tr X\mathfrak{H}_1 \\
&=&-\frac{3}{2} P \Big[   \DDc \hot( \DDc\DDc\c \ov{\Xi}-2\Hb  \,\DDc\c \ov{\Xi} )\\
 &&-3\Hb   \hot \,( \DDc\DDc\c \ov{\Xi}-2\Hb  \,\DDc\c \ov{\Xi} ) \Big]\\
 &&-\frac{3}{2}P (11\Hb)\hot \Big[  \, \DDc\DDc\c \ov{\Xi}-2\Hb  \,\DDc\c \ov{\Xi}  \Big]- 3P10(\underline{H} \hot \Hb) \,\DDc\c \ov{\Xi} \\
 &&+ \frac 3 2P 3 \tr X\Hb \hot \Big[ \ov{H} \c \Xh-\ov{\tr\Xb} \Xi+\AA_3\Big]+\mbox{Expr}_2(A)\\
 &=&\frac{3}{2} P \Big[  - \DDc \hot( \DDc\DDc\c \ov{\Xi})-6\Hb  \hot \DDc  \DDc\c \ov{\Xi}-6\Hb\hot \Hb  \,\DDc\c \ov{\Xi}\\
 &&- 3 \tr X\ov{\tr\Xb}\Hb \hot  \Xi + 3 \tr X\Hb \hot ( \ov{H} \c \Xh) + 3 \tr X\Hb \hot  \AA_3 \Big]+\mbox{Expr}_2(A)
\eeaa
as stated.
\end{proof}

Here we compute $\nabc_4 \mathfrak{H}_2$.

\begin{lemma}\label{lemma:nabc4mathfrak2} We have, modulo quadratic terms,
\beaa
\nabc_4\mathfrak{H}_2  &=&- 2(\tr X)\mathfrak{H}_2 - 3\tr X \Hb \hot   \AA_3 +\widetilde{\MM}[\Xi]  \\
&&+ \Big[\frac 1 2 (3\tr X-3 \ov{\tr X})\ov{\tr X\tr\Xb}+ ( 3\ov{\tr X})\ov{ \DDc}\c\Hb-3\tr X(\Hb \c \ov{\Hb})\\
&& -3 \tr X \ov{H} \c  \Hb+3 (\ov{\tr X} P-\tr X \ov{P}) \Big]\Xh+\mbox{Expr}_4(A),
\eeaa
with
\beaa
\mbox{Expr}_4(A)&=&  -\Hb \hot \big( \DDbc \c A +A\c \ov{\Hb} \big)+\DDc\hot (\ov{H} \c A ) +4 \Hb \hot( \ov{H} \c A)\\
&&- \big(\ov{\tr X\tr\Xb}-2\ov{ \DDc}\c\Hb-2P \big) A+\tr X\big(- \nabc_3A -\frac{1}{2}\tr\Xb A \big),
\eeaa
and 
\beaa
\widetilde{\MM}[\Xi]&:=& \MM[\Xi] -6\ov{P} \,\Xi  - \ov{H} \c \DDc \DDc\hot \Xi-(\DDc \c \ov{H}) \DDc\hot \Xi\\
&&- \ov{H} \c \DDc ( \Xi\hot(\Hb+H)) -\Hb\hot\big(\ov{H} \c \DDc\hot \Xi \big)\\
&&- (\DDc \c \ov{H} )( \Xi\hot(\Hb+H))  -\Hb\hot\big(\ov{H} \c( \Xi\hot(\Hb+H)) \big) \\
&& \big(\ov{\tr X\tr\Xb}-\ov{ \DDc}\c\Hb-2P \big)\big( \DDc\hot \Xi+  \Xi\hot(\Hb+H) \big)\\
&&-3   \Hb \hot\big(\ov{H} \c \big( \DDc\hot \Xi+  \Xi\hot(\Hb+H) \big)\big).
\eeaa
\end{lemma}
\begin{proof} 
Recall that
\beaa
  \mathfrak{H}_2&=&\AA_4-2B \hot \Hb -\DDc\hot\big(\ov{H} \c \Xh\big)+\big(\ov{\tr X\tr\Xb}-2\ov{ \DDc}\c\Hb-2P\big)\Xh-3   \Hb \hot\big(\ov{H} \c \Xh\big).
\eeaa
We therefore obtain
\beaa
\nabc_4\mathfrak{H}_2&=& I_1+I_2+I_3+I_4+I_5.
\eeaa
with 
\beaa
I_1&=& \nabc_4 \AA_4, \\
I_2&=&-2 \nabc_4(B \hot \Hb), \\
I_3&=& -\nab_4\DDc\hot\big(\ov{H} \c \Xh\big)\\
I_4&=& \nabc_4 \Big(\big(\ov{\tr X\tr\Xb}-2\ov{ \DDc}\c\Hb-2P\big)\Xh\Big)\\
I_5&=& -3 \nabc_4 \Big(  \Hb \hot\big(\ov{H} \c \Xh\big)\Big).
\eeaa

We now compute each term above.

Using Lemma \ref{lemma:nabc4AA4}, we deduce
    \beaa
I_1 &=& -2(\tr X)\AA_4  - 3\tr X \Hb \hot   \AA_3 -\tr X \DDc\hot B +2\tr X \Hb \hot B+\MM[\Xi] \\
 &&-\frac{1}{2}(\ov{\tr X}-\tr X)\DDc\hot(\Xh\c\ov{H})  +\frac{1}{2}\tr X \DDc\hot(\widehat{X}\c\ov{\Hb})\\
    &&-\frac{1}{2}(\tr X-\ov{\tr X})(H \c \ov{H})\Xh  + \tr X (  \Hb \c \ov{\Hb}) \Xh+ \tr X\underline{H} \hot (\Xh\c\ov{H})+\Xh \c \ov{\DDc} (\tr X\Hb).
\eeaa
Using \eqref{rule-1} to write $\Xh \c \ov{\DDc} (\tr X\Hb)= \Xh (\ov{\DDc} \c(\tr X\Hb))$ and using that  $  \ov{\tr X}\ov{ H}=-\tr X \ov{\Hb}+O(\ep)$, 
\beaa
 \ov{ \DDc}\c(\tr X \Hb)&=&  \ov{ \DDc}(\tr X) \c  \Hb+ \tr X \ov{ \DDc}\c\Hb\\
 &=& (\ov{\tr X} - \tr X) \ov{H} \c  \Hb+ \tr X \ov{ \DDc}\c\Hb\\
  &=& \ov{\tr X}  \ov{H} \c  \Hb - \tr X \ov{H} \c  \Hb+ \tr X \ov{ \DDc}\c\Hb\\
    &=& -\tr X  \ov{\Hb} \c  \Hb - \tr X \ov{H} \c  \Hb+ \tr X \ov{ \DDc}\c\Hb,
\eeaa
we finally obtain
    \beaa
I_1 &=& -2(\tr X)\AA_4  - 3\tr X \Hb \hot   \AA_3 -\tr X \DDc\hot B +2\tr X \Hb \hot B+\MM[\Xi] \\
 &&-\frac{1}{2}(\ov{\tr X}-\tr X)\DDc\hot(\Xh\c\ov{H})  +\frac{1}{2}\tr X \DDc\hot(\widehat{X}\c\ov{\Hb})\\
    &&-\frac{1}{2}(\tr X-\ov{\tr X})(H \c \ov{H})\Xh  + \tr X\underline{H} \hot (\Xh\c\ov{H})\\
    &&+ \big(  - \tr X \ov{H} \c  \Hb+ \tr X \ov{ \DDc}\c\Hb \big)\Xh.
\eeaa

We have 
\beaa
I_2&=& -2\nabc_4 B \hot \Hb-2B \hot \nabc_4\Hb\\
&=& -2\big( -2\ov{\tr X} B+3\ov{P} \,\Xi +\frac{1}{2} \DDbc \c A +\frac{1}{2}A\c \ov{\Hb} \big) \hot \Hb+2\tr X B \hot  \Hb\\
&=&(2\tr X+4\ov{\tr X})  \Hb \hot B -6\ov{P} \,\Xi  -\Hb \hot \big( \DDbc \c A +A\c \ov{\Hb} \big),
\eeaa
where we used $\nabc_4 \Hb=-\tr X \Hb+O(\ep)$  and the Bianchi identity $\nabc_4B = -2\ov{\tr X} B+3\ov{P} \,\Xi +\frac{1}{2} \DDbc \c A +\frac{1}{2}A\c \ov{\Hb}$.

Using Lemma \ref{comm:1form}, \eqref{eq:comm-nabc4nabc3-DDchot-err}, applied to $F=\ov{H} \c \Xh$ and $s=1$, we have, up to quadratic terms,
  \beaa
[ \nabc_4, \DDc \hot] \ov{H} \c \Xh   &=& -\frac 1 2 \tr X \DDc \hot  (\ov{H} \c \Xh)+\Hb\hot\nabc_4 (\ov{H} \c \Xh),
 \eeaa
and therefore 
\beaa
I_3&=&   -\DDc\hot\nabc_4\big(\ov{H} \c \Xh\big) -[\nabc_4,\DDc\hot]\big(\ov{H} \c \Xh\big)\\
&=& -\DDc\hot\nabc_4\big(\ov{H} \c \Xh\big)  +\frac 1 2 \tr X \DDc \hot  (\ov{H} \c \Xh)-\Hb\hot\nabc_4 (\ov{H} \c \Xh)\\
&=& -\DDc\hot\big(\ov{H} \c \nabc_4\Xh+\nabc_4\ov{H} \c \Xh\big)  +\frac 1 2 \tr X \DDc \hot  (\ov{H} \c \Xh)\\
&&-\Hb\hot\big(\ov{H} \c \nabc_4\Xh+\nabc_4\ov{H} \c \Xh\big) .
\eeaa
Using $\nabc_4\ov{H} =  -\frac{1}{2}\tr X(\ov{H}-\ov{\Hb})+O(\ep)$ and the Ricci identity $\nabc_4\Xh=-\frac 1 2 (\tr X+ \ov{\tr X})\Xh+ \DDc\hot \Xi+  \Xi\hot(\Hb+H)-A$, we obtain
\beaa
\ov{H} \c \nabc_4\Xh+\nabc_4\ov{H} \c \Xh&=& \ov{H} \c \big(-\frac 1 2 (\tr X+ \ov{\tr X})\Xh+ \DDc\hot \Xi+  \Xi\hot(\Hb+H)-A \big) \\
&& -\frac{1}{2}\tr X(\ov{H}-\ov{\Hb}) \c \Xh\\
&=&-\frac 1 2 (2\tr X+ \ov{\tr X})\ov{H} \c\Xh+\frac{1}{2}\tr X\ov{\Hb}\c \Xh\\
&&+\ov{H} \c \DDc\hot \Xi+ \ov{H} \c( \Xi\hot(\Hb+H)) -\ov{H} \c A.
\eeaa
Writing that $  \ov{\tr X}\ov{ H}=-\tr X \ov{\Hb}+O(\ep)$, we obtain
\beaa
\ov{H} \c \nabc_4\Xh+\nabc_4\ov{H} \c \Xh&=&-\tr X\ov{H} \c\Xh+ \tr X \ov{\Hb}\c\Xh\\
&&+\ov{H} \c \DDc\hot \Xi+ \ov{H} \c( \Xi\hot(\Hb+H)) -\ov{H} \c A.
\eeaa
We therefore obtain
\beaa
I_3&=& -\DDc\hot\big(-\frac 1 2 (2\tr X+ \ov{\tr X})\ov{H} \c\Xh+\frac{1}{2}\tr X\ov{\Hb}\c \Xh -\ov{H} \c A \big)  +\frac 1 2 \tr X \DDc \hot  (\ov{H} \c \Xh)\\
&&-\Hb\hot\big( -\tr X\ov{H} \c\Xh+ \tr X \ov{\Hb}\c\Xh -\ov{H} \c A \big) \\
&=& \frac 1 2 (3\tr X+ \ov{\tr X})\DDc\hot(\ov{H} \c\Xh)-\frac{1}{2}\tr X\DDc\hot(\ov{\Hb}\c \Xh)\\
&&+\frac 1 2 (2\DDc\tr X+ \DDc\ov{\tr X})\hot (\ov{H} \c\Xh)-\frac{1}{2}\DDc\tr X\hot (\ov{\Hb}\c \Xh) \\
&&+\Hb\hot\big(\tr X\ov{H} \c\Xh- \tr X \ov{\Hb}\c\Xh \big) +\DDc\hot (\ov{H} \c A ) + \Hb \hot \ov{H} \c A.
\eeaa

Using that $\DDc \tr X =-2 \tr X \Hb+O(\ep)$ and  $\DDc(\ov{\tr X})=(\tr X-\ov{\tr X})H+O(\ep)$ we obtain
\beaa
I_3&=& \frac 1 2 (3\tr X+ \ov{\tr X})\DDc\hot(\ov{H} \c\Xh)-\frac{1}{2}\tr X\DDc\hot(\ov{\Hb}\c \Xh)\\
&&- \tr X \Hb \hot (\ov{H} \c\Xh)+\frac 1 2 (\tr X-\ov{\tr X})(H \c\ov{H}) \Xh \\
&&+\DDc\hot (\ov{H} \c A ) + \Hb \hot \ov{H} \c A.
\eeaa

We now compute 
\beaa
I_4&=& \big(\ov{\tr X\tr\Xb}-2\ov{ \DDc}\c\Hb-2P\big)\nabc_4\Xh+\nabc_4\big(\ov{\tr X\tr\Xb}-2\ov{ \DDc}\c\Hb-2P\big)\Xh\\
&=& \big(\ov{\tr X\tr\Xb}-2\ov{ \DDc}\c\Hb-2P\big)\big(-\frac 1 2 (\tr X+ \ov{\tr X})\Xh+ \DDc\hot \Xi+  \Xi\hot(\Hb+H)-A \big)\\
&&+\big(\nabc_4(\ov{\tr X})\ov{\tr\Xb}+\ov{\tr X}\nabc_4(\ov{\tr\Xb})-2\nabc_4\ov{ \DDc}\c\Hb-2\nabc_4P \big)\Xh\\
&=& \Big[-\frac 1 2 (\tr X+ \ov{\tr X})\ov{\tr X\tr\Xb}+ (\tr X+ \ov{\tr X})\ov{ \DDc}\c\Hb+(\tr X+ \ov{\tr X})P \Big]\Xh \\
&&+ \big(\ov{\tr X\tr\Xb}-2\ov{ \DDc}\c\Hb-2P\big)\big( \DDc\hot \Xi+  \Xi\hot(\Hb+H) \big)\\
&&- \big(\ov{\tr X\tr\Xb}-2\ov{ \DDc}\c\Hb-2P\big) A +\Big[-\frac 1 2 (\ov{\tr X})^2\ov{\tr\Xb}\\
&&+\ov{\tr X}\big(-\frac{1}{2}\ov{\tr X\tr\Xb}+\ov{ \DDc}\c\Hb+\Hb\c\ov{\Hb}+2P \big)-2\nabc_4\ov{ \DDc}\c\Hb+3 \tr X P \Big]\Xh,
\eeaa
where we used $\nabc_4 \ov{\tr X}=-\frac 1 2 (\ov{\tr X})^2+O(\ep)$, $\nabc_4\ov{\tr\Xb}  =-\frac{1}{2}\ov{\tr X\tr\Xb}+\ov{ \DDc}\c\Hb+\Hb\c\ov{\Hb}+2P+O(\ep^2)$,  and $\nabc_4 P=-\frac 3 2 \tr X P+O(\ep)$. This gives
\beaa
I_4&=& \Big[-\frac 1 2 (\tr X+3 \ov{\tr X})\ov{\tr X\tr\Xb}+ (\tr X+ 2\ov{\tr X})\ov{ \DDc}\c\Hb+ \ov{\tr X}(\Hb \c \ov{\Hb})\\
&&-2\nabc_4\ov{ \DDc}\c\Hb+(4\tr X+3 \ov{\tr X})P \Big]\Xh \\
&&+ \big(\ov{\tr X\tr\Xb}-2\ov{ \DDc}\c\Hb-2P\big)\big( \DDc\hot \Xi+  \Xi\hot(\Hb+H) \big)\\
&&- \big(\ov{\tr X\tr\Xb}-2\ov{ \DDc}\c\Hb-2P\big) A .
\eeaa
Finally, using Lemma \ref{comm:1form-div}, \eqref{eq:comm-nabc4nabc3-ovDDc-F-err}, applied to $F=\Hb$ and $s=0$,
  \beaa
[ \nabc_4 ,\ov{\DDc} \c] \Hb     &=&- \frac 1 2\ov{\tr X}\, \left( \ov{\DDc} \c \Hb -  \ov{\Hb} \c \Hb \right)+\ov{\Hb} \c \nabc_4 \Hb+O(\ep)\\
&=&- \frac 1 2\ov{\tr X}\, \left( \ov{\DDc} \c \Hb -  \ov{\Hb} \c \Hb \right) -\tr X\ov{\Hb} \c \Hb+O(\ep),
 \eeaa
we write
\beaa
&&-2\nabc_4\ov{ \DDc}\c\Hb= -2\ov{ \DDc}\c(\nabc_4\Hb)-2[\nabc_4,\ov{ \DDc}\c]\Hb\\
&=& 2\ov{ \DDc}\c(\tr X \Hb)+\ov{\tr X}\, \left( \ov{\DDc} \c \Hb -  \ov{\Hb} \c \Hb \right) +2\tr X\ov{\Hb} \c \Hb\\
&=& 2\ov{ \DDc}\c(\tr X \Hb)+\ov{\tr X}\, \left( \ov{\DDc} \c \Hb -  \ov{\Hb} \c \Hb \right) +2\tr X\ov{\Hb} \c \Hb\\
&=& 2 \ov{ \DDc}(\tr X) \c  \Hb+ 2\tr X \ov{ \DDc}\c\Hb+\ov{\tr X}\, \left( \ov{\DDc} \c \Hb -  \ov{\Hb} \c \Hb \right) +2\tr X\ov{\Hb} \c \Hb\\
&=& 2 (\ov{\tr X} - \tr X) \ov{H} \c  \Hb+ (2\tr X+\ov{\tr X}) \ov{ \DDc}\c\Hb  +(2\tr X-\ov{\tr X})\ov{\Hb} \c \Hb\\
&=& -2 \tr X \ov{\Hb} \c  \Hb -2 \tr X \ov{H} \c  \Hb+ (2\tr X+\ov{\tr X}) \ov{ \DDc}\c\Hb  +(2\tr X-\ov{\tr X})\ov{\Hb} \c \Hb\\
&=& -2 \tr X \ov{H} \c  \Hb+ (2\tr X+\ov{\tr X}) \ov{ \DDc}\c\Hb  -\ov{\tr X}\ov{\Hb} \c \Hb,
\eeaa
where we used that $  \ov{\tr X}\ov{ H}=-\tr X \ov{\Hb}+O(\ep)$. We finally obtain
\beaa
I_4&=& \Big[-\frac 1 2 (\tr X+3 \ov{\tr X})\ov{\tr X\tr\Xb}+ (3\tr X+ 3\ov{\tr X})\ov{ \DDc}\c\Hb\\
&& -2 \tr X \ov{H} \c  \Hb+(4\tr X+3 \ov{\tr X})P \Big]\Xh \\
&&+ \big(\ov{\tr X\tr\Xb}-2\ov{ \DDc}\c\Hb-2P\big)\big( \DDc\hot \Xi+  \Xi\hot(\Hb+H) \big)\\
&&- \big(\ov{\tr X\tr\Xb}-2\ov{ \DDc}\c\Hb-2P\big) A.
\eeaa

We compute
\beaa
I_5&=& -3  \nabc_4 \Hb \hot\big(\ov{H} \c \Xh\big)-3   \Hb \hot\big(\nabc_4\ov{H} \c \Xh\big)-3   \Hb \hot\big(\ov{H} \c \nabc_4\Xh\big)\\
&=& 3 \tr X \Hb \hot\big(\ov{H} \c \Xh\big)+\frac{3}{2}\tr X \Hb \hot\big((\ov{H}-\ov{\Hb}) \c \Xh\big)\\
&&-3   \Hb \hot\big(\ov{H} \c \big(-\frac 1 2 (\tr X+ \ov{\tr X})\Xh+ \DDc\hot \Xi+  \Xi\hot(\Hb+H)-A \big)\big)\\
&=& 6 \tr X \Hb \hot\big(\ov{H} \c \Xh\big)-\frac{3}{2}\tr X (\Hb \c \ov{\Hb} ) \Xh+\frac 3 2  \ov{\tr X}  \Hb \hot\big(\ov{H} \c\Xh\big) \\
&&-3   \Hb \hot\big(\ov{H} \c \big( \DDc\hot \Xi+  \Xi\hot(\Hb+H) \big)\big)+3   \Hb \hot\big(\ov{H} \c A \big)\\
&=& 6 \tr X \Hb \hot\big(\ov{H} \c \Xh\big)-3\tr X (\Hb \c \ov{\Hb} ) \Xh \\
&&-3   \Hb \hot\big(\ov{H} \c \big( \DDc\hot \Xi+  \Xi\hot(\Hb+H) \big)\big)+3   \Hb \hot\big(\ov{H} \c A \big),
\eeaa
where we used $\nabc_4 \Hb=-\tr X \Hb+O(\ep)$, $\nabc_4\ov{H} =  -\frac{1}{2}\tr X(\ov{H}-\ov{\Hb})+O(\ep)$ and the Ricci identity $\nabc_4\Xh=-\frac 1 2 (\tr X+ \ov{\tr X})\Xh+ \DDc\hot \Xi+  \Xi\hot(\Hb+H)-A$. We also used that $  \ov{\tr X}\ov{ H}=-\tr X \ov{\Hb}+O(\ep)$. 

We obtain for the sum
\beaa
&&\nabc_4\mathfrak{H}_2\\
&=& -2(\tr X)\AA_4  - 3\tr X \Hb \hot   \AA_3 -\tr X \DDc\hot B +(4\tr X+4\ov{\tr X})  \Hb \hot B \\
&&+2\tr X\DDc\hot(\Xh\c\ov{H})+ + \Big[-\frac 1 2 (\tr X+3 \ov{\tr X})\ov{\tr X\tr\Xb}+ (4\tr X+ 3\ov{\tr X})\ov{ \DDc}\c\Hb\\
&&-3\tr X(\Hb \c \ov{\Hb}) -3 \tr X \ov{H} \c  \Hb+(4\tr X+3 \ov{\tr X})P \Big]\Xh+6 \tr X \Hb \hot\big(\ov{H} \c \Xh\big) \\
&&+\mbox{Expr}_3(A)+\widetilde{\MM}[\Xi],
\eeaa
where 
\beaa
\mbox{Expr}_3(A)&=&  -\Hb \hot \big( \DDbc \c A +A\c \ov{\Hb} \big)+\DDc\hot (\ov{H} \c A ) +4 \Hb \hot( \ov{H} \c A)\\
&&- \big(\ov{\tr X\tr\Xb}-\ov{ \DDc}\c\Hb-2P \big) A.
\eeaa
Writing 
\beaa
  \mathfrak{H}_2&=&\AA_4-2B \hot \Hb -\DDc\hot\big(\ov{H} \c \Xh\big)+\big(\ov{\tr X\tr\Xb}-2\ov{ \DDc}\c\Hb-2P\big)\Xh-3   \Hb \hot\big(\ov{H} \c \Xh\big),
\eeaa
we obtain
\beaa
\nabc_4\mathfrak{H}_2&=&- 2(\tr X)\mathfrak{H}_2 - 3\tr X \Hb \hot   \AA_3 -\tr X \DDc\hot B +4\ov{\tr X}  \Hb \hot B  \\
&&+ \Big[\frac 1 2 (3\tr X-3 \ov{\tr X})\ov{\tr X\tr\Xb}+ ( 3\ov{\tr X})\ov{ \DDc}\c\Hb\\
&&-3\tr X(\Hb \c \ov{\Hb}) -3 \tr X \ov{H} \c  \Hb+3 \ov{\tr X} P \Big]\Xh+\mbox{Expr}_3(A)+\widetilde{\MM}[\Xi].
\eeaa
Observe that $\ov{\tr X} \Hb =- \tr X H+O(\ep)$ and therefore on the right hand side the terms depending on $B$ can be written in terms of $A$ and $\Xh$ using the Bianchi identity:
\beaa
 -\tr X \DDc\hot B + 4\ov{\tr X}  B \hot  \Hb&=& -\tr X \DDc\hot B - 4\tr X  B \hot  H\\
 &=&\tr X\big(- \nabc_3A -\frac{1}{2}\tr\Xb A   -3\ov{P}\Xh\big).
\eeaa
This gives
\beaa
\nabc_4\mathfrak{H}_2  &=&- 2(\tr X)\mathfrak{H}_2 - 3\tr X \Hb \hot   \AA_3   \\
&&+ \Big[\frac 1 2 (3\tr X-3 \ov{\tr X})\ov{\tr X\tr\Xb}+ ( 3\ov{\tr X})\ov{ \DDc}\c\Hb-3\tr X(\Hb \c \ov{\Hb})\\
&& -3 \tr X \ov{H} \c  \Hb+3 (\ov{\tr X} P-\tr X \ov{P}) \Big]\Xh+\mbox{Expr}_4(A)+\widetilde{\MM}[\Xi],
\eeaa
with
\beaa
\mbox{Expr}_4(A)&=& \mbox{Expr}_3(A)+\tr X\big(- \nabc_3A -\frac{1}{2}\tr\Xb A \big).
\eeaa
This concludes the proof.
\end{proof}

We can finally conclude the derivation. For simplicity, we consider the above computations in a gauge where $\Xi=0$. Then from \eqref{nabc4mathfrakH}, i.e.
\beaa
\nabc_4\mathfrak{H} &=& \nabc_4\mathfrak{H}_1+\frac 3 2 P \big[\nabc_4 \mathfrak{H}_2-\frac 3 2 \tr X  \mathfrak{H}_2\big]+ \nabc_4 \Big( \frac 3 2 P \ov{\tr\Xb}A\Big),
\eeaa
and using Lemma \ref{lemma:nabc4mathfrak1} and \ref{lemma:nabc4mathfrak2}, we obtain
\beaa
\nabc_4\mathfrak{H} &=&- \frac 7 2  \tr X\mathfrak{H}_1  +\frac{3}{2} P \Big[  3 \tr X\Hb \hot ( \ov{H} \c \Xh) + 3 \tr X\Hb \hot  \AA_3 \Big]+\mbox{Expr}_2(A)\\
&&+\frac 3 2 P \Big[-\frac 7 2 (\tr X)\mathfrak{H}_2  - 3\tr X \Hb \hot   \AA_3  \\
&&+ \Big(\frac 1 2 (3\tr X-3 \ov{\tr X})\ov{\tr X\tr\Xb}+ ( 3\ov{\tr X})\ov{ \DDc}\c\Hb-3\tr X(\Hb \c \ov{\Hb})\\
&& -3 \tr X \ov{H} \c  \Hb+3 (\ov{\tr X} P-\tr X \ov{P}) \Big)\Xh+\mbox{Expr}_4(A)\Big]+ \nabc_4 \Big( \frac 3 2 P \ov{\tr\Xb}A\Big)\\
  &=&- \frac 7 2  \tr X \big( \mathfrak{H}_1+ \frac 3 2 P \mathfrak{H}_2\big)  \\
&&+\frac 3 2 P \Big[ 3 \tr X\Hb \hot ( \ov{H} \c \Xh) + \Big(\frac 1 2 (3\tr X-3 \ov{\tr X})\ov{\tr X\tr\Xb}+ ( 3\ov{\tr X})\ov{ \DDc}\c\Hb\\
&&-3\tr X(\Hb \c \ov{\Hb}) -3 \tr X \ov{H} \c  \Hb+3 (\ov{\tr X} P-\tr X \ov{P}) \Big)\Xh \Big]\\
  &&+\mbox{Expr}_2(A)+\frac 3 2 P\mbox{Expr}_4(A) + \nabc_4 \Big( \frac 3 2 P \ov{\tr\Xb}A\Big).
\eeaa
Recalling that $\mathfrak{H}=\mathfrak{H}_1+\frac{3}{2} P     \mathfrak{H}_2+\frac 3 2 P \ov{\tr\Xb}A$, we write
\beaa
&&\nabc_4\mathfrak{H} + \frac 7 2  \tr X \mathfrak{H}   \\
 &=&\frac 3 2 P \big[  3 \tr X\Hb \hot ( \ov{H} \c \Xh) + \Big(\frac 1 2 (3\tr X-3 \ov{\tr X})\ov{\tr X\tr\Xb}+ ( 3\ov{\tr X})\ov{ \DDc}\c\Hb\\
&&-3\tr X(\Hb \c \ov{\Hb}) -3 \tr X \ov{H} \c  \Hb+3 (\ov{\tr X} P-\tr X \ov{P}) \Big)\Xh \big]\\
  &&+\mbox{Expr}_2(A)+\frac 3 2 P\mbox{Expr}_4(A) + \nabc_4 \Big( \frac 3 2 P \ov{\tr\Xb}A\Big)+ \frac{21}{4} P \tr X \ov{\tr \Xb} A.
\eeaa

Using \eqref{simil-Leibniz}, i.e. $\Hb \hot ( \ov{H} \c \Xh)+H \hot ( \ov{\Hb} \c \Xh)= ( \Hb \c \ov{H}+ \ov{\Hb} \c H) \ \Xh$ and that $\ov{\tr X} \Hb =- \tr X H +O(\ep)$, we have
\beaa
3 \tr X\Hb \hot ( \ov{H} \c \Xh)&=& 3 \tr X\big(-H \hot ( \ov{\Hb} \c \Xh)+ ( \Hb \c \ov{H}+ \ov{\Hb} \c H) \ \Xh \big)\\
&=& 3\ov{\tr X} \Hb \hot ( \ov{\Hb} \c \Xh)+ (3 \tr X \Hb \c \ov{H}- 3 \ov{\tr X}\ov{\Hb} \c  \Hb) \ \Xh\\
&=&  3 \tr X \Hb \c \ov{H} \ \Xh
\eeaa
and therefore obtain
\beaa
&&\nabc_4\mathfrak{H} + \frac 7 2  \tr X \mathfrak{H}   \\
 &=&\frac 3 2 P \big[ \Big(\frac 1 2 (3\tr X-3 \ov{\tr X})\ov{\tr X\tr\Xb}+ ( 3\ov{\tr X})\ov{ \DDc}\c\Hb\\
&&-3\tr X(\Hb \c \ov{\Hb})+3 (\ov{\tr X} P-\tr X \ov{P}) \Big)\Xh \big]\\
  &&+\mbox{Expr}_2(A)+\frac 3 2 P\mbox{Expr}_4(A) + \nabc_4 \Big( \frac 3 2 P \ov{\tr\Xb}A\Big)+ \frac{21}{4} P \tr X \ov{\tr \Xb} A.
\eeaa

We now show that the terms in the parenthesis vanish in Kerr, i.e.
\beaa
\frac 1 2 (3\tr X-3 \ov{\tr X})\ov{\tr X\tr\Xb}+3 (\ov{\tr X} P-\tr X \ov{P}) +3\ov{\tr X}\ov{ \DDc}\c\Hb-3 \tr X (\Hb \c \ov{\Hb})=O(\ep).
\eeaa

Recall that in Kerr 
\beaa
 \trch=\frac{2r}{|q|^2},\quad \atrch=\frac{2a\cos\th}{|q|^2}, \qquad \tr X=\frac{2}{q}, \qquad\qquad \tr\Xb=-\frac{2\Delta q}{|q|^4}, \qquad P=-\frac{2m}{q^3}\\
 H_1=\frac{ai\sin\th\, q}{|q|^3}, \qquad H_2=\frac{a\sin\th \, q}{|q|^3}, \qquad \Hb_1=-\frac{ai\sin\th\, \ov{q}}{|q|^3},\qquad\,\, \Hb_2=-\frac{a\sin\th\,\ov{q}}{|q|^3}.
\eeaa
We compute
\beaa
\ov{\tr X\tr\Xb}&=&- \frac{2}{\ov{q}}  \ov{\frac{2\De q}{|q|^4}}=-  \frac{4\De }{|q|^4}\\
\frac 1 2 (3\tr X-3 \ov{\tr X})\ov{\tr X\tr\Xb}&=& -  \frac{6\De }{|q|^4} \left(\frac {2}{q}- \frac{2}{\ov{q}}\right)=  i\frac{24a\cos\th \De }{|q|^6} \\
\ov{\tr X} P-\tr X \ov{P}&=& \frac{2}{\ov{q}} \left(-\frac{2m}{q^3}\right)-\frac 2 q  \left(-\frac{2m}{\ov{q}^3}\right)=-  \frac{4m}{\ov{q}q^3}+ \frac{4m}{q\ov{q}^3}= \frac{4m(q^2- \ov{q}^2)}{|q|^6}\\
&=& i \frac{16 a \cos\th mr}{|q|^6}.
\eeaa
Therefore
\beaa
&&\frac 1 2 (3\tr X-3 \ov{\tr X})\ov{\tr X\tr\Xb}+3 (\ov{\tr X} P-\tr X \ov{P})\\
 &=&  i\frac{24a\cos\th (r^2-2mr+a^2) }{|q|^6} + i \frac{48 a \cos\th mr}{|q|^6}=  i\frac{24a\cos\th (r^2+a^2) }{|q|^6}\\
&=&  i\frac{4a\cos\th}{|q|^8} (6r^4+6a^2\cos^2\th r^2+6a^2r^2+6a^2a^2\cos^2\th).
\eeaa
Also,
\beaa
\ov{ \DDc}\c\Hb&=&2 \div \etab+ 2 i \curl \etab\\
&=& 2 \frac{a^2}{|q|^6} (-2\cos^2\th r^2+\sin^2\th r^2-a^2\sin^2\th\cos^2\th-2a^2\cos^2\th)\\
&&+ i 2 \frac{a\cos\th}{|q|^6}(- 2r^3+2a^2\cos^2\th r-4a^2 r ).
\eeaa
Also,
\beaa
\Hb \c \ov{\Hb}&=&2\frac{a^2}{|q|^6} (  \sin^2\th r^2+a^2\sin^2\th \cos^2\th ).
\eeaa
We therefore obtain in Kerr
\beaa
&&\frac 1 2 (3\tr X-3 \ov{\tr X})\ov{\tr X\tr\Xb}+3 (\ov{\tr X} P-\tr X \ov{P}) +3\ov{\tr X}\ov{ \DDc}\c\Hb-3 \tr X (\Hb \c \ov{\Hb})\\
&=& i\frac{4a\cos\th}{|q|^8} (6r^4+6a^2\cos^2\th r^2+6a^2r^2+6a^2a^2\cos^2\th)\\
&&+\frac{4a\cos\th}{|q|^8}\big(   - 6r^4-12a^2 r^2 +6   a^2 \sin^2\th r^2+3a^2(-2a^2\cos^2\th) \big) =0.
\eeaa

We therefore have proved, that modulo quadratic terms,
\beaa
\nabc_4\mathfrak{H} + \frac 7 2  \tr X \mathfrak{H}    &=&\mbox{Expr}_2(A)+\frac 3 2 P\mbox{Expr}_4(A) + \nabc_4 \Big( \frac 3 2 P \ov{\tr\Xb}A\Big)+ \frac{21}{4} P \tr X \ov{\tr \Xb} A.
\eeaa
By denoting 
\beaa
\PP(A)&:=& \mbox{Expr}_2(A)+\frac 3 2 P\mbox{Expr}_4(A) + \nabc_4 \Big( \frac 3 2 P \ov{\tr\Xb}A\Big)+ \frac{21}{4} P \tr X \ov{\tr \Xb} A\\
&=&  - \frac{1}{4}(\DDc+11\Hb)\hot \Big(  (\DDc +\Hb)\big( \big( \DDc + 3 \Hb\big) \c \big( \DDc\c\ov{A}+\ov{A}\c\Hb\big)\big)\Big)\\
&&- 5 (\underline{H} \hot \Hb)\big( \DDc + 3 \Hb\big) \c \big( \DDc\c\ov{A}+\ov{A}\c\Hb\big)\\
&&+\frac 3 2 P\Big(  -\Hb \hot \big( \DDbc \c A +A\c \ov{\Hb} \big)+\DDc\hot (\ov{H} \c A ) +4 \Hb \hot( \ov{H} \c A)\\
&&- \big(\ov{\tr X\tr\Xb}-2\ov{ \DDc}\c\Hb-2P \big) A\\
&&+\tr X\big(- \nabc_3A -\frac{1}{2}\tr\Xb A \big). \Big) + \nabc_4 \Big( \frac 3 2 P \ov{\tr\Xb}A\Big)+ \frac{21}{4} P \tr X \ov{\tr \Xb} A,
\eeaa
we have, for $r\leq r_0$
\beaa
\nabc_4\mathfrak{H} + \frac 7 2  \tr X \mathfrak{H}    &=&\PP(A)+\dk^{\leq 3}(\Ga_b \c \Ga_g).
\eeaa
Observe that restricting the validity of the above expression to the region $r \leq r_0$, we can denote the error terms by $\dk^{\leq 3}(\Ga_b \c \Ga_g)$.

Finally we show that $\nabc_4\mathfrak{H} + \frac 7 2  \tr X \mathfrak{H} =\Big(\nabc_4+2\tr X\Big)^4\Ab $.
Recall that 
\beaa
\mathfrak{F} &=& -\nabc_4\Ab -\frac{1}{2}\tr X \Ab,\\
\mathfrak{G} &=& \nabc_4\mathfrak{F}+\frac{3}{2}\tr X\mathfrak{F},\\
\mathfrak{H} &=& \nabc_4\mathfrak{G} +\frac{5}{2}\tr X\mathfrak{G},
\eeaa
and 
\beaa
\nab_4\mathfrak{H}+\frac{7}{2}\tr X\mathfrak{H} &=& \PP(A).
\eeaa
Choosing a normalization such that $\om=O(\ep)$, and hence $\tr X=\frac{2}{q}+O(\ep)$, we obtain
\beaa
\mathfrak{F} =-\frac{1}{q}\nabc_4(q\Ab), \qquad \mathfrak{G}=\frac{1}{q^3}\nab_4(q^3\mathfrak{F}), \qquad  \mathfrak{H}=\frac{1}{q^5}\nab_4(q^5\mathfrak{G}), \qquad \frac{1}{q^7}\nab_4(q^7\mathfrak{H})=\PP(A).
\eeaa
We infer
\beaa
 \frac{1}{q^7}\nab_4\left(q^2\nab_4\left(q^2\nab_4\left(q^2\nabc_4(q\Ab)\right)\right)\right)=\PP(A).
\eeaa
The above can be written as
\beaa
\frac{1}{q^3}\nab_4\left(q^2\nab_4\right)^3(q\Ab) &=& \frac{1}{q^2}\nab_4\left(q^2\nab_4\right)^2\Big(q^3\nab_4\Ab+q^2\Ab\Big)\\
&=& \frac{1}{q^3}\nab_4\left(q^2\nab_4\right)\Big(q^5\nab_4^2\Ab +4q^4\nab_4\Ab+2q^3\Ab\Big)\\
&=& \frac{1}{q^3}\nab_4\Big(q^7\nab_4^3\Ab+9q^6\nab_4^2\Ab+18q^5\nab_4\Ab+6q^4\Ab\Big)\\
&=& q^4\nab_4^4\Ab+16q^3\nab_4^3\Ab+72q^2\nab_4^2\Ab+96q\nab_4\Ab+24\Ab,
\eeaa
which is equivalent to 
\beaa
\frac{1}{q^4}\nabc_4^4(q^4\Ab) &=& q^4\nab_4^4\Ab+4\nab_4(q^4)\nab_4^3\Ab+6\nab_4^2(q^4)\nab_4^2\Ab+4\nab_4^3(q^4)\nab_4\Ab+\nab_4^4(q^4)\Ab\\
&=& q^4\nab_4^4\Ab+16q^3\nab_4^3\Ab+72q^2\nab_4^2\Ab+96q\nab_4\Ab+24\Ab.
\eeaa
which can be written as
\beaa
\Big(\nab_4+2\tr X\Big)^4\Ab =\PP(A).
\eeaa
This concludes the proof of Proposition \ref{THEOREM:TEUK-STAR}.


\section{Proof of Lemma \ref{LEMMA:WAVEEQP}}\label{proof-wave-P}


We make use  of the Bianchi identity 
\beaa
\nabc_4P -\frac{1}{2}\DDc\c \ov{B} &=& -\frac{3}{2}\tr X P + \Hb \c\ov{B} -\ov{\Xi}\c\Bb -\frac{1}{4}\Xbh\c \ov{A}.
\eeaa
We differentiate w.r.t. $\nabc_3$ and obtain 
\beaa
&&\nabc_3\nabc_4P -\frac{1}{2}\DDc\c \nabc_3\ov{B} -\frac{1}{2}[\nabc_3,\DDc\c]\ov{B}  \\
&=& -\frac{3}{2}P\nabc_3\tr X -\frac{3}{2}\tr X \nabc_3P +\nabc_3\Hb \c\ov{B}+\Hb \c\nabc_3\ov{B}\\
&&  -\nabc_3(\ov{\Xi}\c\Bb) -\frac{1}{4}\nabc_3(\Xbh\c \ov{A}).
\eeaa
Next we make use of 
\beaa
\nabc_3B-\DDc\ov{P} &=& -\tr\Xb B+3\ov{P}H +\ov{\Bb}\c \Xh+\frac{1}{2}A\c\ov{\Xib},
\eeaa
and taking the complex conjugate, we infer
\beaa
\nabc_3\ov{B}-\ov{\DDc} P &=& -\ov{\tr\Xb} \ov{B}+3P \ov{H} +\Bb\c \ov{\Xh}+\frac{1}{2}\ov{A}\c \Xib.
\eeaa
We deduce 
\beaa
&&\nabc_3\nabc_4P -\frac{1}{2}\DDc\c \big(\ov{\DDc} P -\ov{\tr\Xb} \ov{B}+3P \ov{H} \big) -\frac{1}{2}[\nabc_3,\DDc\c]\ov{B}  \\
&=& -\frac{3}{2}P\nabc_3\tr X -\frac{3}{2}\tr X \nabc_3P +\nabc_3\Hb \c\ov{B}+\Hb \c\big(\ov{\DDc} P -\ov{\tr\Xb} \ov{B}+3P \ov{H} \big)\\
&&  -\nabc_3(\ov{\Xi}\c\Bb) -\frac{1}{4}\nabc_3(\Xbh\c \ov{A})+\frac 1 2 \DDc \c \big( \Bb\c \ov{\Xh}+\frac{1}{2}\ov{A}\c \Xib\big) \\
&&+\Hb \c \big( \Bb\c \ov{\Xh}+\frac{1}{2}\ov{A}\c \Xib \big),
\eeaa
and hence
\beaa
&&\nabc_3\nabc_4P+\frac{3}{2}P\nabc_3\tr X +\frac{3}{2}\tr X \nabc_3P\\
&& -\frac{1}{2}\DDc\c (\ov{\DDc} P +3P \ov{H}) -\Hb \c (\ov{\DDc} P +3P \ov{H})  \\
&=&-\frac{1}{2}\DDc\c \big( \ov{\tr\Xb} \ov{B}\big) +\frac{1}{2}[\nabc_3,\DDc\c]\ov{B} +\nabc_3\Hb \c\ov{B}-\ov{\tr\Xb} \Hb \c \ov{B} \\
&&  -\nabc_3(\ov{\Xi}\c\Bb) -\frac{1}{4}\nabc_3(\Xbh\c \ov{A})+\frac 1 2 \DDc \c \big( \Bb\c \ov{\Xh}+\frac{1}{2}\ov{A}\c \Xib\big) \\
&&+\Hb \c \big( \Bb\c \ov{\Xh}+\frac{1}{2}\ov{A}\c \Xib \big).
\eeaa
Next  we make use of    the commutator, see Lemma \ref{COMMUTATOR-NAB-C-3-DD-C-HOT}, 
 \beaa
 [\nabc_3 ,\ov{\DDc} \c] F      &=&- \frac 1 2\ov{\tr\Xb}\, ( \ov{\DDc} \c F +(s-1) \ov{H}  \c F )+\ov{H} \c \nabc_3 F\\
 &&-\Bb \c  \ov{F} +\ov{\Xib} \c \nab_4 F -\Xbh \c \ov{\DDc} F-H \c \Xbh \c \ov{F},
 \eeaa
 which for $B$ and $s=1$, we deduce
\beaa
[\nabc_3,\DDc\c]\ov{B}  &=&- \frac 1 2\tr\Xb\, \DDc \c \ov{B} +H \c \nabc_3\ov{B}\\
&&-\ov{\Bb} \c  B +\Xib \c \nab_4 \ov{B} -\ov{\Xbh} \c \DDc \ov{B}-\ov{H} \c \ov{\Xbh} \c B.
\eeaa
Hence,
\beaa
&&\nabc_3\nabc_4P+\frac{3}{2}P\nabc_3\tr X +\frac{3}{2}\tr X \nabc_3P\\
&& -\frac{1}{2}\DDc\c (\ov{\DDc} P +3P \ov{H}) -\Hb \c (\ov{\DDc} P +3P \ov{H})  \\
&=&-\frac{1}{2}\ov{\tr\Xb}\DDc\c  \ov{B}-\frac{1}{2}\DDc \ov{\tr\Xb}\c \ov{B}  +\frac{1}{2}\big(- \frac 1 2\tr\Xb\, \DDc \c \ov{B} +H \c \nabc_3\ov{B} \big) \\
&&+\nabc_3\Hb \c\ov{B}-\ov{\tr\Xb} \Hb \c \ov{B} \\
&&  -\nabc_3(\ov{\Xi}\c\Bb) -\frac{1}{4}\nabc_3(\Xbh\c \ov{A})+\frac 1 2 \DDc \c \big( \Bb\c \ov{\Xh}+\frac{1}{2}\ov{A}\c \Xib\big) \\
&&+\Hb \c \big( \Bb\c \ov{\Xh}+\frac{1}{2}\ov{A}\c \Xib \big)\\
&&+\frac 12 \big(-\ov{\Bb} \c  B +\Xib \c \nab_4 \ov{B} -\ov{\Xbh} \c \DDc \ov{B}-\ov{H} \c \ov{\Xbh} \c B \big).
\eeaa
Observe that the LHS of the above becomes
\beaa
LHS&=& \nabc_3\nabc_4P+\frac{3}{2}\tr X \nabc_3P -\frac{1}{2}\DDc\c \ov{\DDc} P-\frac{3}{2}\ov{H} \c \DDc P-\Hb \c \ov{\DDc} P\\
&&  +\frac{3}{2}\Big[ -\DDc\c  \ov{H}  +\nabc_3\tr X  -2 \Hb \c  \ov{H}\Big] P,
\eeaa
while the RHS, using
\beaa
 \frac{1}{2}\DDc\c \ov{B} &=&\nabc_4P +\frac{3}{2}\tr X P - \Hb \c\ov{B}  +\ov{\Xi}\c\Bb+\frac{1}{4}\Xbh\c \ov{A}, \\
 \nabc_3\ov{B} &=&\ov{\DDc} P +3P \ov{H} -\ov{\tr\Xb} \ov{B}+\Bb\c \ov{\Xh}+\frac{1}{2}\ov{A}\c \Xib\\
- \frac{1}{2}\DDc\ov{\tr\Xb} &=&-i\Im(\tr\Xb)\Hb-\frac{1}{2}\ov{\DDc}\c\Xbh-i\Im(\tr X)\Xib+\Bb\\
\nabc_3\Hb  &=&  -\frac{1}{2}\ov{\tr\Xb}(\Hb-H)+\nabc_4\Xib -\frac{1}{2}\Xbh\c(\ov{\Hb}-\ov{H}) +\Bb,
\eeaa
gives the following
\beaa
RHS&=&- \frac 1 4\big(2\ov{\tr\Xb}+ \tr\Xb \big) \, \DDc \c \ov{B}+\frac 1 2 H \c \nabc_3\ov{B} \\
&& +\big(-\frac{1}{2}\DDc \ov{\tr\Xb}+\nabc_3\Hb -\ov{\tr\Xb} \Hb\big) \c \ov{B}  -\nabc_3(\ov{\Xi}\c\Bb)\\
&& -\frac{1}{4}\nabc_3(\Xbh\c \ov{A})+\frac 1 2 \DDc \c \big( \Bb\c \ov{\Xh}+\frac{1}{2}\ov{A}\c \Xib\big) +\Hb \c \big( \Bb\c \ov{\Xh}+\frac{1}{2}\ov{A}\c \Xib \big)\\
&&+\frac 12 \big(-\ov{\Bb} \c  B +\Xib \c \nab_4 \ov{B} -\ov{\Xbh} \c \DDc \ov{B}-\ov{H} \c \ov{\Xbh} \c B \big)\\
&=&- \frac 1 2\big(2\ov{\tr\Xb}+ \tr\Xb \big) \, \big(\nabc_4P +\frac{3}{2}\tr X P - \Hb \c\ov{B}  \big)\\
&&+\frac 1 2 H \c \big(\ov{\DDc} P +3P \ov{H} -\ov{\tr\Xb} \ov{B} \big) \\
&&+\big(-i\Im(\tr\Xb)\Hb -\frac{1}{2}\ov{\tr\Xb}(\Hb-H) -\ov{\tr\Xb} \Hb\big) \c \ov{B} +\err\\
&=&- \frac 1 2\big(2\ov{\tr\Xb}+ \tr\Xb \big) \, \big(\nabc_4P +\frac{3}{2}\tr X P  \big)+\frac 1 2 H \c \big(\ov{\DDc} P +3P \ov{H} \big)+\err
\eeaa
where 
\beaa
\err&=&  -\nabc_3(\ov{\Xi}\c\Bb) -\frac{1}{4}\nabc_3(\Xbh\c \ov{A})+\frac 1 2 \DDc \c \big( \Bb\c \ov{\Xh}+\frac{1}{2}\ov{A}\c \Xib\big) \\
&&+\Hb \c \big( \Bb\c \ov{\Xh}+\frac{1}{2}\ov{A}\c \Xib \big)\\
&&+\frac 12 \big(-\ov{\Bb} \c  B +\Xib \c \nab_4 \ov{B} -\ov{\Xbh} \c \DDc \ov{B}-\ov{H} \c \ov{\Xbh} \c B \big)\\
&&- \frac 1 2\big(2\ov{\tr\Xb}+ \tr\Xb \big)(\ov{\Xi}\c\Bb+\frac{1}{4}\Xbh\c \ov{A})+\frac 12 H \c \big( \Bb\c \ov{\Xh}+\frac{1}{2}\ov{A}\c \Xib \big)\\
&&+\big(-\frac{1}{2}\ov{\DDc}\c\Xbh-i\Im(\tr X)\Xib+\Bb+\nabc_4\Xib -\frac{1}{2}\Xbh\c(\ov{\Hb}-\ov{H}) +\Bb\big)\c\ov{B}.
\eeaa
By putting the two together we obtain
\beaa
&& \nabc_3\nabc_4P+\frac{3}{2}\tr X \nabc_3P+ \frac 1 2\big(2\ov{\tr\Xb}+ \tr\Xb \big) \, \nabc_4P -\frac{1}{2}\DDc\c \ov{\DDc} P\\
&&-\frac{3}{2}\ov{H} \c \DDc P-\big( \Hb+\frac 1 2 H\big) \c \ov{\DDc} P   +\frac{3}{2}\Big[ \ov{\tr\Xb} \tr X+2P  -2 \Hb \c  \ov{H}  \Big] P \\
&=& \err-\frac 32 \big(\Xib\c\ov{\Xi}-\frac{1}{2}\Xbh\c\ov{\Xh}\big)P,
\eeaa
where we also used
\beaa
\nabc_3\tr X +\frac{1}{2}\tr\Xb\tr X &=& \DDc\c\ov{H}+H\c\ov{H}+2P+\Xib\c\ov{\Xi}-\frac{1}{2}\Xbh\c\ov{\Xh}.
\eeaa

From Lemma \ref{lemma:expression-wave-operator}, we deduce
\beaa
\square_\g P &=& -e_3(e_4(P))+\left(2\omb -\frac{1}{2}\trchb\right)\nab_4P -\frac{1}{2}\trch\nab_3P+2\eta\c\nab P+\Delta P.
\eeaa
Also, we have
\beaa
\frac{1}{2}\DD\c\ov{\DD}P &=& \frac{1}{2}(\nab+i\dual\nab)\c(\nab-i\dual\nab)P = \Delta P+\frac{i}{2}(\dual\nab\c \nab-\nab\c\dual\nab)P\\
&=& \Delta P+i(\nab_2\nab_1-\nab_1\nab_2)P\\
&=&\Delta P+i\Big((\D_2\D_1-\D_1\D_2)P +\D_{\D_2e1-\nab_2e_1-\D_1e_2+\nab_1e_2}P\Big)\\
&=&\Delta P+\frac{i}{2}\Big((\chi_{21}-\chi_{12})e_3P+(\chib_{21}-\chib_{12})e_4P\Big)\\
&=&\Delta P-\frac{i}{2}\Big(\atrch e_3P+\atrchb e_4P\Big),
\eeaa
and
\beaa
\ov{H}\c\DD P+H\c\ov{\DD} P  &=& (\eta-i\dual\eta)\c(\nab +i\dual\nab)P+(\eta+i\dual\eta)\c(\nab -i\dual\nab)P=4\eta\c\nab P.
\eeaa
We therefore deduce for $P$ of conformal type $0$:
\beaa
\square_\g P &=& -\nabc_3\nabc_4P+\frac{1}{2}\DDc\c\ov{\DDc}P -\frac{1}{2}\tr \Xb \nabc_4P -\frac{1}{2}\tr X\nabc_3P\\
&&+\frac 1 2 \ov{H}\c\DDc P+\frac 1 2 H\c\ov{\DDc} P.
\eeaa

We then obtain from the above computations
\beaa
\square_\g P&=&\tr X \nabc_3P+ \ov{\tr\Xb} \, \nabc_4P-\ov{H} \c \DDc P- \Hb\c \ov{\DDc} P  \\
&& +\frac{3}{2}\Big[ \ov{\tr\Xb} \tr X+2P  -2 \Hb \c  \ov{H}  \Big] P +  \err[\square_\g P],
\eeaa
where 
\beaa
\err[\square_\g P]&=& \err-\frac 32 \big(\Xib\c\ov{\Xi}-\frac{1}{2}\Xbh\c\ov{\Xh}\big)P
\eeaa
as stated.


\section{Proof of Proposition \ref{PROP:RELATION-QF-P}}
\label{proof-relation-qf-P}


From the Bianchi identity
 \beaa
 \nabc_4\Ab +\frac{1}{2}\tr X \Ab &=&-\frac 1 2 \DDc\hot \Bb-  2 \Hb   \hot \Bb -3P \Xbh,
\eeaa
we infer
\beaa
\nabc_4 \big( \nabc_4\Ab +\frac{1}{2}\tr X \Ab \big) &=&-\frac 1 2 \DDc\hot \nabc_4\Bb-\frac 1 2[\nabc_4, \DDc\hot] \Bb\\
&&-  2 \Hb   \hot \nabc_4\Bb-  2 \nabc_4\Hb   \hot \Bb -3P \nabc_4\Xbh-3\nabc_4P \Xbh.
\eeaa
Using the commutator \eqref{eq:comm-nabc4nabc3-DDchot-err}, i.e.
\beaa
[\nabc_4 , \DDc \hot ]\Bb &=&- \frac 1 2 \tr X\left( \DDc\hot \Bb + 2\Hb\hot \Bb\right)+ \underline{H} \hot \nabc_4 \Bb\\
  &&+\Xi \hot \nabc_3 \Bb+ r^{-1} \Ga_g  \c \dk^{\leq 1} \Bb,
\eeaa
we have
\beaa
\nabc_4 \big( \nabc_4\Ab +\frac{1}{2}\tr X \Ab \big) &=&-\frac 1 2 \DDc\hot \nabc_4\Bb-  \frac 5 2  \Hb   \hot \nabc_4\Bb\\
&&+ \frac 1 4 \tr X\left( \DDc\hot \Bb + 2\Hb\hot \Bb\right) \\
&&-  2 \nabc_4\Hb   \hot \Bb -3P \nabc_4\Xbh-3\nabc_4P \Xbh\\
&&+\Xi \hot \nabc_3 \Bb+ r^{-1} \Ga_g  \c \dk^{\leq 1} \Bb.
\eeaa
Next, using 
\beaa
\nabc_4\Bb+\DDc P &=& -\tr X\Bb-3P\Hb + \Ga_b \c B+\Xi \c \Ab, \\
\nabc_4P -\frac{1}{2}\DDc\c \ov{B} &=& -\frac{3}{2}\tr X P + \Hb \c\ov{B} +\Xi \c \Bb+ \Ga_b \c A, \\
\nabc_4\widehat{\Xb} +\frac{1}{2}\tr X\, \widehat{\Xb} &=&\frac 1 2  \DDc\hot\Hb  +\frac 1 2 \Hb\hot\Hb -\frac{1}{2}\ov{\tr\Xb} \widehat{X}+\frac 1 4 \Xi\hot\Xib
\eeaa
we deduce 
\beaa
&&\nabc_4 \big( \nabc_4\Ab +\frac{1}{2}\tr X \Ab \big) \\
&=&\frac 1 2 \DDc\hot \big(\DDc P  +\tr X\Bb+3P\Hb  \big)+  \frac 5 2  \Hb   \hot \big(\DDc P  +\tr X\Bb+3P\Hb  \big)\\
&&+ \frac 1 4 \tr X\left( \DDc\hot \Bb + 2\Hb\hot \Bb\right) -  2 \nabc_4\Hb   \hot \Bb\\
&& -3P \big(-\frac{1}{2}\tr X\, \widehat{\Xb} +\frac 1 2  \DDc\hot\Hb  +\frac 1 2 \Hb\hot\Hb -\frac{1}{2}\ov{\tr\Xb} \widehat{X} \big)-3\big( -\frac{3}{2}\tr X P  \big) \Xbh\\
&&+\Xi \hot \nabc_3 \Bb+r^{-1} \dk^{\leq 1} ( \Xi \c \Ab) + r^{-1} \Ga_g  \c \dk^{\leq 1} \Bb+ r^{-1} \dk^{\leq 1} (\Ga_b \c B )+ r^{-3} \Xi \c \Xib,
\eeaa
which gives
\beaa
&&\nabc_4 \big( \nabc_4\Ab +\frac{1}{2}\tr X \Ab \big) \\
&=&\frac 1 2 \DDc\hot \DDc P -\frac 32\tr X \big(- \frac 1 2 \DDc\hot  \Bb-2\Hb \hot \Bb-3 P  \Xbh\big) +3P \big(\frac{1}{2}\tr X\, \widehat{\Xb}  +\frac{1}{2}\ov{\tr\Xb} \widehat{X} \big)\\
&& +\big( \frac 1 2 \DDc \tr X  -  2 \nabc_4\Hb  \big) \hot \Bb+\big( 4 \DDc P   + 6P  \Hb   \big)\hot \Hb  \\
&&+\Xi \hot \nabc_3 \Bb+r^{-1} \dk^{\leq 1} ( \Xi \c \Ab) + r^{-1} \Ga_g  \c \dk^{\leq 1} \Bb+ r^{-1} \dk^{\leq 1} (\Ga_b \c B )+ r^{-3} \Xi \c \Xib.
\eeaa
Using that $-\frac 1 2 \DDc\hot \Bb-  2 \Hb   \hot \Bb -3P \Xbh=\nabc_4\Ab +\frac{1}{2}\tr X \Ab $, we have 
\beaa
&&\nabc_4 \big( \nabc_4\Ab +\frac{1}{2}\tr X \Ab \big) +\frac 32\tr X \big(\nabc_4\Ab +\frac{1}{2}\tr X \Ab\big)\\
&=&\frac 1 2 \DDc\hot \DDc P  +\frac 3 2 P \big(\tr X\, \widehat{\Xb}  +\ov{\tr\Xb} \widehat{X} \big)\\
&& +\big( \frac 1 2 \DDc \tr X  -  2 \nabc_4\Hb  \big) \hot \Bb+\big( 4 \DDc P   + 6P  \Hb   \big)\hot \Hb  \\
&&+\Xi \hot \nabc_3 \Bb+r^{-1} \dk^{\leq 1} ( \Xi \c \Ab) + r^{-1} \Ga_g  \c \dk^{\leq 1} \Bb+ r^{-1} \dk^{\leq 1} (\Ga_b \c B )+ r^{-3} \Xi \c \Xib.
\eeaa
Using \eqref{eq:vanishing-relations-perturbations}, we deduce that $\frac 1 2 \DDc \tr X  -  2 \nabc_4\Hb=\tr X\Hb+r^{-1} \dk^{\leq 1} \Ga_g$, which gives
\bea\label{eq:intermediate-relation-qfb-P}
\begin{split}
&\nabc_4 \big( \nabc_4\Ab +\frac{1}{2}\tr X \Ab \big) +\frac 32\tr X \big(\nabc_4\Ab +\frac{1}{2}\tr X \Ab\big)\\
&=\frac 1 2 \DDc\hot \DDc P+\big( 4 \DDc P   + 6P  \Hb + \tr X \Bb  \big)\hot \Hb   +\frac 3 2 P \big(\tr X\, \widehat{\Xb}  +\ov{\tr\Xb} \widehat{X} \big) \\
&+\Xi \hot \nabc_3 \Bb+r^{-1} \dk^{\leq 1} ( \Xi \c \Ab) + r^{-1} \dk^{\leq 1} (\Ga_g  \c \Bb)+ r^{-1} \dk^{\leq 1} (\Ga_b \c B ).
\end{split}
\eea

Consider the LHS of \eqref{eq:intermediate-relation-qfb-P}. Using that $\nabc_4 \tr X+ \frac 1 2 ( \tr X)^2=r^{-1} \dk^{\leq 1} \Xi+ \Ga_g \c \Ga_g$, we obtain
\beaa
LHS&=& \nabc_4 \big( \nabc_4\Ab +\frac{1}{2}\tr X \Ab \big) +\frac 32\tr X \big(\nabc_4\Ab +\frac{1}{2}\tr X \Ab\big)\\
&=& \nabc_4 \nabc_4\Ab +\frac{1}{2}\tr X  \nabc_4\Ab+\frac{1}{2} \nabc_4\tr X \Ab  +\frac 32\tr X \big(\nabc_4\Ab +\frac{1}{2}\tr X \Ab\big)\\
&=& \nabc_4 \nabc_4\Ab +2\tr X  \nabc_4\Ab  +\frac 12(\tr X)^2 \Ab+r^{-1} \dk^{\leq 1} \Xi\c  \Ab.
\eeaa
Using the definition \eqref{eq:definition-qfb}, we deduce
\beaa
LHS&=&\underline{Q}(\Ab)  +\left( 2\frac {\atrch^2}{ \trch} + 2 i \atrch\right) \nabc_4\Ab  \\
&&+\left(-\frac 3 2 \frac{\atrch^4}{\trch^2}+  \frac 72\atrch^2+i \trch \atrch-4i\frac{\atrch^3}{\trch}\right) \Ab+r^{-1} \dk^{\leq 1} \Xi\c  \Ab\\
&=&\underline{Q}(\Ab)  +O(ar^{-3}) \dk^{\leq 1} \Ab+r^{-2}  \Ga_g \c  \dk^{\leq 1} \Ab+r^{-1} \dk^{\leq 1} \Xi\c  \Ab.
\eeaa

Writing $\nabc_4 \Ab=\Ab_4 - \frac 1 2 \tr X \Ab$, we have
\beaa
LHS&=&\underline{Q}(\Ab)  +\left( 2\frac {\atrch^2}{ \trch} + 2 i \atrch\right) \left(\Ab_4 - \frac 1 2 (\trch - i \atrch) \Ab\right)  \\
&&+\left(-\frac 3 2 \frac{\atrch^4}{\trch^2}+  \frac 72\atrch^2+i \trch \atrch-4i\frac{\atrch^3}{\trch}\right) \Ab+r^{-1} \dk^{\leq 1} \Xi\c  \Ab\\
&=&\underline{Q}(\Ab)  +\left( 2\frac {\atrch^2}{ \trch} + 2 i \atrch\right) \Ab_4 \\
&&+\left(-\frac 3 2 \frac{\atrch^4}{\trch^2}+  \frac 32\atrch^2-3i\frac{\atrch^3}{\trch}\right) \Ab+r^{-1} \dk^{\leq 1} \Xi\c  \Ab\\
&=&\underline{Q}(\Ab)  +O(ar^{-2}) \Ab_4+O(a^2 r^{-4})\Ab+r^{-1} \dk^{\leq 1} \Xi\c  \Ab.
\eeaa

Using again the Bianchi identity for $\Ab_4$ we also have
\beaa
LHS&=&\underline{Q}(\Ab)+O(ar^{-5}) \Ga_b +O(ar^{-3}) \dk^{\leq 1} \Bb  +O(ar^{-4} ) \Ab\\
&&+r^{-2}\Ga_b \c \Pc+r^{-1} \dk^{\leq 1} \Xi\c  \Ab+ r^{-2} \Ga_b \c \Bb.
\eeaa

Now consider the RHS of \eqref{eq:intermediate-relation-qfb-P}. Using \eqref{eq:error-DDP}, i.e. $\DDc P=- 3 P \Hb+\DDc \Pc + r^{-3} \Ga_g$, we obtain
\beaa
RHS&=& \frac 1 2 \DDc\hot (- 3 P \Hb+\DDc \Pc )+\big( 4 (- 3 P \Hb+\DDc \Pc ) + 6P  \Hb + \tr X \Bb  \big)\hot \Hb \\
&&  +\frac 3 2 P \big(\tr X\, \widehat{\Xb}  +\ov{\tr\Xb} \widehat{X} \big) \\
&&+r^{-4} \dkb^{\leq 1} \Ga_g +\Xi \hot \nabc_3 \Bb+r^{-1} \dk^{\leq 1} ( \Xi \c \Ab) + r^{-1} \dk^{\leq 1} (\Ga_g  \c \Bb)+ r^{-1} \dk^{\leq 1} (\Ga_b \c B )\\
&=& \frac 1 2 \DDc\hot \DDc \Pc - \frac 3 2 \DDc P \hot  \Hb - \frac 3 2 P  \DDc\hot \Hb\\
&&+\big( 4 (\DDc \Pc ) - 6P  \Hb + \tr X \Bb  \big)\hot \Hb  +\frac 3 2 P \big(\tr X\, \widehat{\Xb}  +\ov{\tr\Xb} \widehat{X} \big) \\
&&+r^{-4} \dkb^{\leq 1} \Ga_g +\Xi \hot \nabc_3 \Bb+r^{-1} \dk^{\leq 1} ( \Xi \c \Ab) + r^{-1} \dkb^{\leq 1} (\Ga_g  \c \Bb)+ r^{-1} \dk^{\leq 1} (\Ga_b \c B ).
\eeaa
Using, see \eqref{eq:vanishing-relations-perturbations}, that $\DDc\hot \Hb= - \Hb \hot \Hb + r^{-1} \dk^{\leq 1} \Ga_g$, we obtain
\beaa
&&RHS\\
&=& \frac 1 2 \DDc\hot \DDc \Pc - \frac 3 2(- 3 P \Hb+\DDc \Pc) \hot  \Hb - \frac 3 2 P  (- \Hb \hot \Hb )\\
&&+\big( 4 (\DDc \Pc ) - 6P  \Hb + \tr X \Bb  \big)\hot \Hb +\frac 3 2 P \big(\tr X\, \widehat{\Xb}  +\ov{\tr\Xb} \widehat{X} \big)+ r^{-1} \dk^{\leq 1} \Ga_g \c \Pc \\
&&+r^{-4} \dkb^{\leq 1} \Ga_g +\Xi \hot \nabc_3 \Bb+r^{-1} \dk^{\leq 1} ( \Xi \c \Ab) + r^{-1} \dk^{\leq 1} (\Ga_g  \c \Bb)+ r^{-1} \dk^{\leq 1} (\Ga_b \c B )\\
&=& \frac 1 2 \DDc\hot \DDc \Pc + \frac 5 2(\DDc \Pc) \hot  \Hb + \tr X \Bb \hot \Hb +\frac 3 2 P \big(\tr X\, \widehat{\Xb}  +\ov{\tr\Xb} \widehat{X} \big)+ r^{-1} \dk^{\leq 1} \Ga_g \c \Pc \\
&&+r^{-4} \dkb^{\leq 1} \Ga_g +\Xi \hot \nabc_3 \Bb+r^{-1} \dk^{\leq 1} ( \Xi \c \Ab) + r^{-1} \dk^{\leq 1} (\Ga_g  \c \Bb)+ r^{-1} \dk^{\leq 1} (\Ga_b \c B ).
\eeaa
The above can be further simplified, and we obtain
\beaa
RHS&=& \frac 1 2 \DDc\hot \DDc \Pc +r^{-4} \dkb^{\leq 1}( \Ga_b, r \Ga_g)+ O(ar^{-3}) \Bb+O(ar^{-3}) \dk^{\leq 1} \Pc  \\
&&+ r^{-1} \dk^{\leq 1}( \Ga_b \c (\Pc, B)) + r^{-1} \dk^{\leq 1} (\Ga_g  \c \Bb)  +\Xi \hot \nabc_3 \Bb+r^{-1} \dk^{\leq 1} ( \Xi \c \Ab) .
\eeaa

By combining the above with the LHS, we obtain
\beaa
\underline{Q}(\Ab) &=& \frac 1 2 \DDc\hot \DDc \Pc +r^{-4} \dkb^{\leq 1} ( \Ga_b, r \Ga_g)+O(ar^{-2}) \Ab_4+ O(ar^{-4})\Ab +O(ar^{-3}) \dk^{\leq 1} \Pc  \\
&&+ r^{-1} \dk^{\leq 1}( \Ga_b \c (\Pc, B)) + r^{-2} \dk^{\leq 1} (\Ga_g  \c (r\Bb, \Ab))  +\Xi \hot \nabc_3 \Bb+r^{-1} \dk^{\leq 1} ( \Xi \c \Ab),
\eeaa
or 
\beaa
\underline{Q}(\Ab)&=& \frac 1 2 \DDc\hot \DDc \Pc +r^{-4} \dkb^{\leq 1} ( \Ga_b, r \Ga_g)+ O(ar^{-3})  \dk^{\leq 1} \Bb+ O(ar^{-4}) \Ab+O(ar^{-3}) \dk^{\leq 1} \Pc  \\
&&+ r^{-1} \dk^{\leq 1}( \Ga_b \c (\Pc, B)) + r^{-2} \dk^{\leq 1} (\Ga_g  \c (r\Bb, \Ab)) +\Xi \hot \nabc_3 \Bb+r^{-1} \dk^{\leq 1} ( \Xi \c \Ab),
\eeaa
as stated.

\end{document}